\documentclass[numbers=enddot,12pt,final,onecolumn,notitlepage]{scrartcl}%
\usepackage[headsepline,footsepline,manualmark]{scrlayer-scrpage}
\usepackage[all,cmtip]{xy}
\usepackage{amssymb}
\usepackage{amsmath}
\usepackage{amsthm}
\usepackage{framed}
\usepackage{comment}
\usepackage{color}
\usepackage[breaklinks=True]{hyperref}
\usepackage[sc]{mathpazo}
\usepackage[T1]{fontenc}
\usepackage{needspace}
\usepackage{tabls}
\usepackage{tikz}
\usepackage[type={CC}, modifier={zero}, version={1.0},]{doclicense}
\providecommand{\U}[1]{\protect\rule{.1in}{.1in}}
\newcounter{exer}

\numberwithin{exer}{subsection}
\theoremstyle{definition}
\newtheorem{theo}{Theorem}[subsection]
\newenvironment{theorem}[1][]
{\begin{theo}[#1]\begin{leftbar}}
{\end{leftbar}\end{theo}}
\newtheorem{lem}[theo]{Lemma}
\newenvironment{lemma}[1][]
{\begin{lem}[#1]\begin{leftbar}}
{\end{leftbar}\end{lem}}
\newtheorem{prop}[theo]{Proposition}
\newenvironment{proposition}[1][]
{\begin{prop}[#1]\begin{leftbar}}
{\end{leftbar}\end{prop}}
\newtheorem{defi}[theo]{Definition}
\newenvironment{definition}[1][]
{\begin{defi}[#1]\begin{leftbar}}
{\end{leftbar}\end{defi}}
\newtheorem{remk}[theo]{Remark}
\newenvironment{remark}[1][]
{\begin{remk}[#1]\begin{leftbar}}
{\end{leftbar}\end{remk}}
\newtheorem{coro}[theo]{Corollary}
\newenvironment{corollary}[1][]
{\begin{coro}[#1]\begin{leftbar}}
{\end{leftbar}\end{coro}}
\newtheorem{conv}[theo]{Convention}
\newenvironment{convention}[1][]
{\begin{conv}[#1]\begin{leftbar}}
{\end{leftbar}\end{conv}}
\newtheorem{quest}[theo]{Question}
\newenvironment{question}[1][]
{\begin{quest}[#1]\begin{leftbar}}
{\end{leftbar}\end{quest}}
\newtheorem{warn}[theo]{Warning}
\newenvironment{warning}[1][]
{\begin{warn}[#1]\begin{leftbar}}
{\end{leftbar}\end{warn}}
\newtheorem{conj}[theo]{Conjecture}

\newtheorem{exam}[theo]{Example}
\newenvironment{example}[1][]
{\begin{exam}[#1]\begin{leftbar}}
{\end{leftbar}\end{exam}}
\newtheorem{exmp}[exer]{Exercise}
\newenvironment{exercise}[1][]
{\begin{exmp}[#1]\begin{leftbar}\begin{small}}
{\end{small}\end{leftbar}\end{exmp}}
\newenvironment{statement}{\begin{quote}}{\end{quote}}
\newenvironment{fineprint}{\begin{small}}{\end{small}}
\let\sumnonlimits\sum
\let\prodnonlimits\prod
\let\cupnonlimits\bigcup
\let\capnonlimits\bigcap
\renewcommand{\sum}{\sumnonlimits\limits}
\renewcommand{\prod}{\prodnonlimits\limits}
\renewcommand{\bigcup}{\cupnonlimits\limits}
\renewcommand{\bigcap}{\capnonlimits\limits}
\newenvironment{noncompile}{}{}
\excludecomment{verlong}
\includecomment{vershort}
\excludecomment{noncompile}

\newcommand{\ZZ}{\mathbb{Z}}

\newcommand{\set}[1]{\left\{ #1 \right\}}

\newcommand{\tup}[1]{\left( #1 \right)}

\newcommand{\arinj}{\ar@{_{(}->}}
\newcommand{\arinjrev}{\ar@{^{(}->}}
\newcommand{\arsurj}{\ar@{->>}}
\newcommand{\symd}{\mathbin{\bigtriangleup}}
\newcommand{\Ker}{\operatorname{Ker}}
\usetikzlibrary{arrows.meta}
\usetikzlibrary{calc}
\usetikzlibrary{chains}
\usetikzlibrary{shapes}
\usetikzlibrary{decorations.pathmorphing}
\usetikzlibrary{lindenmayersystems}
\definecolor{darkgreen}{rgb}{0,.5,0}
\ihead{Rings and fields (Math 332 Winter 2023), version 2026-09-20}
\ohead{page \thepage}
\begin{document}

\title{An introduction to the algebra of rings and fields\\{\normalsize (Text for Math 332 Winter 2023 and Winter 2025 at Drexel
University)}}
\author{Darij Grinberg}
\date{draft, September 20, 2026}
\maketitle

\begin{abstract}
\textbf{Abstract.} This is an introduction to rings and fields, written for a
quarter-long undergraduate course. It includes the basic properties of ideals,
modules, algebras and polynomials, the constructions of ring extensions and
finite fields, some number-theoretical applications (such as a proof of
quadratic reciprocity and Jacobsthal's formulas for $p=x^{2}+y^{2}$), and
tastes of Gr\"{o}bner bases and the Smith normal form. Familiarity with groups
and vector spaces is assumed, though no deep results from either theory are used.

Over 250 exercises are included (mostly without solutions).

\end{abstract}
\tableofcontents

\doclicenseThis

\newpage

\section{Preface}

\subsection{What is this?}

These notes are written for the second part of a groups-first undergraduate
abstract algebra sequence, or for an introductory graduate course on rings and
fields. They cover the basic properties of \textbf{rings}, \textbf{modules}
and \textbf{fields}, in particular \textbf{polynomial rings} and
\textbf{finite fields}, while assuming that the reader is fluent in the basic
language of groups (and in elementary number theory, such as the properties of
prime numbers and greatest common divisors). The content is mostly
introductory, and the main results obtained are

\begin{itemize}
\item the Chinese Remainder Theorem for rings and ideals,

\item the construction of monoid algebras (including polynomial rings as a
particular case),

\item the main properties of univariate polynomial rings,

\item the existence and uniqueness of finite fields $\mathbb{F}_{p^{n}}$ of
all prime-power orders,

\item the law of quadratic reciprocity (with a proof in the odd/odd case), and

\item two proofs of Fermat's theorem about writing primes $p$ as sums of two
squares (one giving an \textquotedblleft explicit\textquotedblright%
\ expression in terms of Jacobsthal sums).
\end{itemize}

The last sections briefly introduce Gr\"{o}bner bases (without proofs) and the
Smith normal form (over $\mathbb{Z}$, with an outlined proof). Some properties
of the Fibonacci sequence are explored as applications. Thus, the text is
suited to a quarter-long course, less to a semester-long one.

This text is written rather informally and sometimes tersely. I assume that
the reader has encountered proofs before, as she will have to fill in some
details and understand some hints. Unlike my notes on combinatorics, this text
is not trying to fill any expository gaps, since the literature on abstract
algebra is already vast and includes some rather detailed and rigorous texts
(e.g., Warner \cite{Warner90}, Jacobson \cite{Jacobs1-3}, Zariski/Samuel
\cite{ZarSam86}).

On occasion, I have tried to mildly innovate, e.g., by constructing the
polynomial ring as a monoid algebra, by involving the Fibonacci numbers in a
few places as a \textquotedblleft grass-touching\textquotedblright\ example,
or by avoiding complex numbers in the construction of the cyclotomic
polynomials. I also attempt to view the subject through a more constructive
lens than usual (Section \ref{sec.polys2.factor} is a noticeable example),
although I am nowhere as consistent about it as a text dedicated to
constructive algebra (such as \cite{Edward22}) would be.

A quarter is not much time, and this text reflects the necessary tradeoffs. I
would have loved to cover some Galois theory, more about quadratic number
rings, more about multivariate polynomials, Gr\"{o}bner bases with proofs,
linear algebra with proofs, the Smith normal form over non-$\mathbb{Z}$ rings,
tensor products, determinants, exact sequences, ..., but I have not managed to
fit this into a quarter-long course. A reader whose appetite is whetted by
this text will almost surely have to satisfy it elsewhere (e.g.,
\cite{Aluffi16}, \cite{Bosch18}, \cite{ClaRos25}, \cite{CoLiOs15},
\cite{Cox12}, \cite{DumFoo04}, \cite{Edward22}, \cite{Elman22},
\cite{Fulton08}, \cite{Goodman}, \cite{Gorode17}, \cite{Grifo26}\footnote{See
\url{https://eloisagrifo.github.io/teaching.html} for further notes by Grifo,
some of which are natural continuations of \cite{Grifo26}.}, \cite{Isaacs94},
\cite{Jacobs1-3}, \cite{Klopsc08}, \cite{Knapp1}, \cite{Laurit09},
\cite{Leinst24}, \cite{LidNie00}, \cite{Lorsch20a}, \cite{Lorsch20b},
\cite{McNult16}, \cite{Milne22}, \cite{Reid13}, \cite{Rotman3e},
\cite{Sharif22}, \cite{Siksek19}, \cite{Stein10}, \cite{Steinb06},
\cite{Stewar15}, \cite{Waerde91}, \cite{Weintr22}).

The original template for the structure of this text was the book
\textit{Abstract Algebra} by Dummit and Foote (\cite{DumFoo04} in the
bibliography), specifically a subset of \cite[Chapters 7--14]{DumFoo04}.
However, I have ended up deviating from \cite{DumFoo04} in the presentation,
in the ordering, in the exercises and digressions, and even in some of the terminology.

This text was originally written as lecture notes for my Math 332 course at
Drexel University in Winter 2023, and much of it came from earlier notes for
my Math 533 course in Winter 2021. Both courses have websites containing
homework sets (with a few solutions):%
\begin{align*}
&  \text{\url{https://www.cip.ifi.lmu.de/~grinberg/t/23wa/}}\\
&  \text{\url{https://www.cip.ifi.lmu.de/~grinberg/t/21w/}}%
\end{align*}

I thank Bogdan Nica and Tommy Pett for reporting mistakes in the text that
follows, and Keith Conrad for pedagogical and thematic suggestions. Several
more corrections have been made as the text was proofread by GPT-5.5. Other
than this, GPT-5.5 played a major role in convincing me that Theorem
\ref{thm.gcdseq.fac} was true\footnote{Originally, I could only prove the easy
case of Theorem \ref{thm.gcdseq.fac} in which all of $a_{1},a_{2},a_{3}%
,\ldots$ were assumed regular. Out of curiosity, I asked GPT whether the
general case was true, and it responded readily with a proof. The proof given
in these notes, however, is my own.}.

\subsection{Plan}

This text is split into $6$ chapters:

\begin{enumerate}
\item[2.] \textbf{Rings and ideals.} Just like the notion of a group is an
abstract model for a set of symmetries or invertible operations in general,
the notion of a ring models a set of numbers or things made out of numbers
(such as polynomials or matrices). More formally, a ring is a set equipped
with two operations called \textquotedblleft addition\textquotedblright\ and
\textquotedblleft multiplication\textquotedblright\ and two elements called
\textquotedblleft zero\textquotedblright\ and \textquotedblleft
unity\textquotedblright\ that satisfy certain axioms. We will explore both
specific examples and general properties of rings, and we will study features
of rings such as subrings and ideals, and certain classes of rings such as
integral domains and fields.

\item[3.] \textbf{Modules.} Modules are the natural generalization of vector
spaces when the underlying number system is replaced by a ring. In particular,
we will learn some of the basics of abstract linear algebra (the theory of
vector spaces over fields) here.

\item[4.] \textbf{Monoid algebras and polynomials.} This is a generalization
of the classical notion of polynomials, which replaces the monomials by
something much more general (the elements of a monoid). Among other things,
this will give us a precise definition of polynomials. We will study
univariate polynomials (i.e., polynomials in one indeterminate) in more
detail, establishing in particular some of their unique features (division
with remainder and the Euclidean algorithm). This will help us
\textquotedblleft adjoin\textquotedblright\ a root of a polynomial to a
commutative ring or field.

\item[5.] \textbf{Finite fields.} Finite fields are \textquotedblleft
miniature versions\textquotedblright\ of our familiar number systems; they are
sets with extremely well-behaved \textquotedblleft addition\textquotedblright%
\ and \textquotedblleft multiplication\textquotedblright\ operations but only
finitely many elements. We will build up some of their basic theory and see a
few of their many applications.

\item[6.] \textbf{Polynomials II.} We will study polynomials in more detail,
focussing now mostly on multivariate polynomials. We will give a very
introductory treatment of Gr\"{o}bner bases (explaining their simplest uses,
but not proving any of their nontrivial properties). We will explain how
polynomials with integer coefficients can be factored (into irreducibles), and
address some parts of the ancient question of \textquotedblleft how do you
solve a system of polynomial equations?\textquotedblright. We also define and
explore the cyclotomic polynomials, a special family of polynomials beloved by
number theorists.

\item[7.] \textbf{Modules over a PID.} To be specific, we will be studying
modules over $\mathbb{Z}$ only. In particular, we will prove the structure
theorem for finite abelian groups, and construct the Smith normal form of a
matrix. As an application, we will prove the existence of primitive roots in a
finite field.
\end{enumerate}

\begin{noncompile}
TODO:

- Fields of fractions and localization (including of modules and algebras).

- Once this is done: M\"{o}bius function formula for strong gcd-sequence
factorization and for cyclotomic and dynatomic polynomials.

- Gr\"{o}bner bases with proofs and details. Maybe even extend Gr\"{o}bner
bases to monoid rings (then including noncommutative polynomials as well).

- The Smith normal form and its consequences (including the Jordan normal form
you might know from linear algebra).

- Tensor products.

- A bit of Galois theory.
\end{noncompile}

\subsection{Notations}

We shall use the following notations:

\begin{itemize}
\item We let $\mathbb{N}=\left\{  0,1,2,\ldots\right\}  $.

\item The notation $\left\vert S\right\vert $ denotes the size (i.e., the
number of elements) of a set $S$.

\item Unlike algebraic geometers, we do accept noncommutative rings as rings
(see below for the definition). Unlike \cite{DumFoo04}, we don't accept
nonunital rings (i.e., rings without a $1$) as rings. This will be discussed
in more detail below.

\item The shorthand \textquotedblleft LTTR\textquotedblright\ means
\textquotedblleft left to the reader\textquotedblright. I will try to use it
sparingly and only for arguments that are similar to ones given before (and no
harder than the latter).
\end{itemize}

\subsection{\label{sec.intro.mod}Refresher on modular arithmetic}

We will use modular arithmetic (specifically, the notion of residue classes
modulo $n$, and the algebraic operations on these classes). An introduction to
modular arithmetic can be found in almost any textbook on abstract algebra
(see, e.g., \cite[\S 3.4]{19s}), and I assume that you have seen it at least
in some form, since it underlies the standard definition of cyclic groups. Let
me nevertheless give a summary as a reminder.

For the rest of this section, we fix an integer $n$. Two integers $a$ and $b$
are said to be \textbf{congruent} (to each other) \textbf{modulo }$n$ if and
only if $n\mid a-b$. The short notation for this is \textquotedblleft$a\equiv
b\operatorname{mod}n$\textquotedblright, but we shall shorten this even
further to \textquotedblleft$a\underset{n}{\equiv}b$\textquotedblright\ in
this subsection, so that $\underset{n}{\equiv}$ becomes a binary relation on
the set $\mathbb{Z}$ of all integers.

For example, $5\underset{2}{\equiv}9$ (since $2\mid5-9$) but
$5\underset{2}{\not \equiv }8$ (since $2\nmid5-8$). (As usual,
\textquotedblleft$a\underset{n}{\not \equiv }b$\textquotedblright\ means
\textquotedblleft not $a\underset{n}{\equiv}b$\textquotedblright.)

The binary relation $\underset{n}{\equiv}$ is an equivalence relation. It is
called \textbf{congruence modulo }$n$. Its equivalence classes are called the
\textbf{residue classes of integers modulo }$n$. Explicitly, for every integer
$a$, the residue class that contains $a$ is the set%
\begin{align*}
&  \left\{  \text{all integers that are congruent to }a\text{ modulo
}n\right\} \\
&  =\left\{  \text{all integers that differ from }a\text{ by a multiple of
}n\right\} \\
&  =\left\{  \ldots,\ a-2n,\ a-n,\ a,\ a+n,\ a+2n,\ a+3n,\ \ldots\right\}  .
\end{align*}
We denote this class by $\overline{a}$. Two integers $a$ and $b$ satisfy
$\overline{a}=\overline{b}$ if and only if $a\underset{n}{\equiv}b$.

In particular, the residue class $\overline{0}$ of $0$ consists of all
integers that are multiples of $n$. That is:%
\[
\overline{0}=\left\{  \ldots,\ -3n,\ -2n,\ -n,\ 0,\ n,\ 2n,\ 3n,\ \ldots
\right\}  .
\]
The residue class $\overline{n}$ of $n$ consists of the very same integers,
since an integer is congruent to $n$ modulo $n$ if and only if it is congruent
to $0$ modulo $n$. In other words, $\overline{n}=\overline{0}$. Likewise,
$\overline{2n}=\overline{0}$ and $\overline{3n}=\overline{0}$ and so on.
Likewise, $\overline{n+1}=\overline{1}$ and $\overline{n+2}=\overline{2}$ and
so on. On the other hand, the $n$ residue classes $\overline{0},\overline
{1},\ldots,\overline{n-1}$ are all distinct (if $n$ is positive), since no two
of the integers $0,1,\ldots,n-1$ are congruent modulo $n$.

Here are some examples:

\begin{itemize}
\item For $n=2$, the only two residue classes modulo $n$ are%
\begin{align*}
\overline{0}  &  =\left\{  \text{all even integers}\right\}  =\left\{
\ldots,-6,-4,-2,0,2,4,\ldots\right\}  \ \ \ \ \ \ \ \ \ \ \text{and}\\
\overline{1}  &  =\left\{  \text{all odd integers}\right\}  =\left\{
\ldots,-5,-3,-1,1,3,5,\ldots\right\}  .
\end{align*}
For any other integer $a$, the residue class $\overline{a}$ of $a$ modulo $2$
is either $\overline{0}$ or $\overline{1}$, depending on whether $a$ is even
or odd. For instance, $\overline{2}=\overline{4}=\overline{6}=\cdots
=\overline{0}$ whereas $\overline{1}=\overline{3}=\overline{5}=\cdots
=\overline{1}$.

\item For $n=5$, there are five residue classes modulo $n$, namely%
\begin{align*}
\overline{0}  &  =\left\{  \ldots,-10,-5,0,5,10,\ldots\right\}  ,\\
\overline{1}  &  =\left\{  \ldots,-9,-4,1,6,11,\ldots\right\}  ,\\
\overline{2}  &  =\left\{  \ldots,-8,-3,2,7,12,\ldots\right\}  ,\\
\overline{3}  &  =\left\{  \ldots,-7,-2,3,8,13,\ldots\right\}  ,\\
\overline{4}  &  =\left\{  \ldots,-6,-1,4,9,14,\ldots\right\}  .
\end{align*}

\item For $n=1$, there is only one residue class modulo $n$, namely
\[
\overline{0}=\left\{  \ldots,-3,-2,-1,0,1,2,3,\ldots\right\}  =\mathbb{Z}.
\]

\item The case $n=0$ is special: Here, no two distinct integers $a$ and $b$
are congruent modulo $n$ (because only $0$ is divisible by $0$). Thus, for
each integer $a$, the residue class $\overline{a}$ of $a$ modulo $0$ is just
the singleton set $\left\{  a\right\}  $. Hence, there are infinitely many
residue classes modulo $0$, one for each integer.
\end{itemize}

As you see, the residue classes $\overline{a}$ modulo $n$ differ from their
underlying integers $a$ in that different integers $a,b$ lead to the same
residue class $\overline{a}=\overline{b}$ when their difference is a multiple
of $n$. Thus, working with residue classes of integers modulo $n$ can be
viewed as working with integers but pretending that $n$ equals $0$ (so that
two integers that differ by a multiple of $n$ are equal).

The set of all residue classes of integers modulo $n$ will be called
$\mathbb{Z}/n$ or $\mathbb{Z}/n\mathbb{Z}$ (or sometimes $\mathbb{Z}_{n}$, but
this symbol is unfortunately also used for other purposes). This set
$\mathbb{Z}/n$ has size $n$ if $n$ is positive\footnote{Indeed, when $n$ is
positive, the $n$ distinct residue classes modulo $n$ are%
\begin{align*}
\overline{0}  &  =\left\{  \ldots,\ -3n,\ -2n,\ -n,\ 0,\ n,\ 2n,\ 3n,\ \ldots
\right\}  ,\\
\overline{1}  &  =\left\{  \ldots
,\ -3n+1,\ -2n+1,\ -n+1,\ 1,\ n+1,\ 2n+1,\ 3n+1,\ \ldots\right\}  ,\\
\overline{2}  &  =\left\{  \ldots
,\ -3n+2,\ -2n+2,\ -n+2,\ 2,\ n+2,\ 2n+2,\ 3n+2,\ \ldots\right\}  ,\\
&  \ldots,\\
\overline{n-1}  &  =\left\{  \ldots
,\ -2n-1,\ -n-1,\ -1,\ n-1,\ 2n-1,\ 3n-1,\ 4n-1,\ \ldots\right\}  .
\end{align*}
}, size $-n$ if $n$ is negative, and infinite size if $n=0$ (indeed,
$\mathbb{Z}/0$ is just $\mathbb{Z}$ \textquotedblleft with its elements
relabelled\textquotedblright\footnote{You may be unused to this; some
textbooks carefully avoid the $n=0$ case when considering $\mathbb{Z}/n$. And
indeed, $\mathbb{Z}/0$ behaves unlike the \textquotedblleft
other\textquotedblright\ $\mathbb{Z}/n$'s in some regard (for example,
$\mathbb{Z}/0$ is infinite, whereas $\mathbb{Z}/n$ is finite for each nonzero
$n$). But the underlying idea is still the same: Two integers $a$ and $b$ are
congruent modulo $0$ if and only if $0$ divides $a-b$; but $0$ only divides
$0$ itself (since the only multiple of $0$ is $0$), so this means that $a$ and
$b$ are congruent modulo $0$ if and only if $a$ and $b$ are equal. Hence, each
residue class modulo $0$ just consists of a single number. Thus, the elements
of $\mathbb{Z}/0$ are the one-element sets $\ldots,\left\{  -2\right\}
,\left\{  -1\right\}  ,\left\{  0\right\}  ,\left\{  1\right\}  ,\left\{
2\right\}  ,\ldots$.}).

We note that, from the viewpoint of group theory, the residue classes modulo
$n$ are the cosets of the subgroup $n\mathbb{Z}=\left\{  \text{all multiples
of }n\right\}  $ in the group $\left(  \mathbb{Z},+,0\right)  $. Thus, the set
$\mathbb{Z}/n\mathbb{Z}$ of these residue classes is the quotient of the group
$\left(  \mathbb{Z},+,0\right)  $ by this subgroup $n\mathbb{Z}$. This is
where the notation $\mathbb{Z}/n\mathbb{Z}$ comes from. (The notation
$\mathbb{Z}/n$ is just shorthand for that.)

We can furthermore define the sum, the difference and the product of any two
residue classes modulo $n$. Namely, if $\overline{a}$ and $\overline{b}$ are
two residue classes modulo $n$, then we set%
\begin{align*}
\overline{a}+\overline{b}  &  =\overline{a+b};\\
\overline{a}-\overline{b}  &  =\overline{a-b};\\
\overline{a}\cdot\overline{b}  &  =\overline{ab}.
\end{align*}
In other words, to add (or subtract, or multiply) two residue classes, we just
add (or subtract, or multiply) the underlying numbers and take the residue
class of the result. It takes a bit of thought to prove that this is
well-defined (i.e., that the values $\overline{a+b}$, $\overline{a-b}$ and
$\overline{ab}$ really depend only on the residue classes $\overline{a}$ and
$\overline{b}$ and not on their chosen representatives $a$ and $b$), but this
is fairly easy and well-known. These operations (addition, subtraction and
multiplication) on residue classes are called \textbf{modular arithmetic}.

Here are some examples:

\begin{itemize}
\item If $n=24$, then
\[
\overline{23}+\overline{5}=\overline{23+5}=\overline{28}=\overline
{4}\ \ \ \ \ \ \ \ \ \ \left(  \text{since }28\underset{24}{\equiv}4\right)
.
\]
This is actually the known fact that \textquotedblleft$5$ hours after $23$
o'clock is $4$ o'clock\textquotedblright\ (although here in the US, you would
usually say \textquotedblleft$11$ PM\textquotedblright\ instead of
\textquotedblleft$23$ o'clock\textquotedblright, and \textquotedblleft$4$
AM\textquotedblright\ instead of \textquotedblleft$4$
o'clock\textquotedblright). The $24$ hours of a day thus naturally correspond
to the residue classes modulo $24$, and reckoning with time is a matter of
modular arithmetic.

\item If $n=12$, then%
\[
\overline{4}\cdot\overline{9}=\overline{4\cdot9}=\overline{36}=\overline
{0}\ \ \ \ \ \ \ \ \ \ \left(  \text{since }36\underset{12}{\equiv}0\right)
.
\]

\item If $n=0$, then (as we have seen above) the residue classes $\overline
{a}\in\mathbb{Z}/0$ are just the sets $\left\{  a\right\}  $, each consisting
of a single integer. They are added and multiplied just as the corresponding
integers: $\left\{  a\right\}  +\left\{  b\right\}  =\left\{  a+b\right\}  $
and $\left\{  a\right\}  \cdot\left\{  b\right\}  =\left\{  a\cdot b\right\}
$.
\end{itemize}

The addition of residue classes that we defined above turns the set
$\mathbb{Z}/n\mathbb{Z}$ of all these residue classes into a group $\left(
\mathbb{Z}/n\mathbb{Z},\ +,\ \overline{0}\right)  $. When $n$ is positive,
this group is known as the \textbf{cyclic group of order }$n$. However, the
multiplication is of interest too, and in fact is one of the main protagonists
of this course.

\subsection{Homework set \#0: (de)motivating questions}

The following exercises should be viewed as food for thought. Some of them are
easy, some hard, some close to impossible at the current state. Just think
about each of them for a little while (5 minutes? 15 minutes? an hour if you
like them?). These exercises are illustrative of some of the elementary
applications of abstract algebra. This course will teach you to solve some of
them. Solution sketches can be found in \cite[solutions to Homework set
0]{21w}.

\begin{exercise}
\label{exe.21hw0.3}\ \ 

\begin{enumerate}
\item[\textbf{(a)}] Factor the polynomial $a^{3}+b^{3}+c^{3}-3abc$.

\item[\textbf{(b)}] Factor the polynomial $bc\left(  b-c\right)  +ca\left(
c-a\right)  +ab\left(  a-b\right)  $.

\item[\textbf{(c)}] How general have your methods been? Did you use tricks
specific to the given polynomials, or do you have an algorithm for factoring
any polynomial (say, with integer coefficients)?
\end{enumerate}
\end{exercise}

\begin{exercise}
\label{exe.21hw0.4}Simplify $\sqrt[3]{2+\sqrt{5}}+\sqrt[3]{2-\sqrt{5}}$.
\end{exercise}

\begin{exercise}
\label{exe.21hw0.5}Let $n\in\mathbb{N}$. Let $a_{1},a_{2},\ldots,a_{n}$ be $n$
integers, and let $b_{1},b_{2},\ldots,b_{n}$ be $n$ further integers. The
Gaussian elimination tells you how to solve the system
\begin{align*}
a_{1}x_{1}+a_{2}x_{2}+\cdots+a_{n}x_{n}  &  =0;\\
b_{1}x_{1}+b_{2}x_{2}+\cdots+b_{n}x_{n}  &  =0
\end{align*}
for $n$ unknowns $x_{1},x_{2},\ldots,x_{n}\in\mathbb{Q}$. The answer, in
general, will have the form \textquotedblleft all $\mathbb{Q}$-linear
combinations (i.e., linear combinations with rational coefficients) of a
certain bunch of vectors\textquotedblright. (More precisely, \textquotedblleft
a certain bunch of vectors\textquotedblright\ are $n-2$ or $n-1$ or $n$
vectors with rational coordinates, depending on the rank of the $2\times
n$-matrix $%
\begin{pmatrix}
a_{1} & a_{2} & \cdots & a_{n}\\
b_{1} & b_{2} & \cdots & b_{n}%
\end{pmatrix}
$.)

Now, how can you solve the above system for $n$ unknowns $x_{1},x_{2}%
,\ldots,x_{n}\in\mathbb{Z}$ ? Will the answer still be \textquotedblleft all
$\mathbb{Z}$-linear combinations (i.e., linear combinations with integer
coefficients) of a certain bunch of vectors\textquotedblright?

What about more general systems of linear equations to be solved for integer unknowns?
\end{exercise}

\Needspace{9pc}

\begin{exercise}
\label{exe.21hw0.6}You are given a $5\times5$-grid of lamps, each of which is
either on or off. For example, writing $1$ for \textquotedblleft
on\textquotedblright\ and $0$ for \textquotedblleft off\textquotedblright, it
may look as follows:
\[%
\begin{tabular}
[c]{|c|c|c|c|c|}\hline
1 & 0 & 0 & 1 & 1\\\hline
1 & 1 & 0 & 0 & 1\\\hline
1 & 0 & 0 & 1 & 0\\\hline
0 & 1 & 1 & 1 & 1\\\hline
0 & 1 & 0 & 0 & 0\\\hline
\end{tabular}
\
\]
In a single move, you can toggle any lamp (i.e., turn it on if it was off, or
turn it off if it was on); however, this will also toggle every lamp adjacent
to it. (\textquotedblleft Adjacent to it\textquotedblright\ means
\textquotedblleft having a grid edge in common with it\textquotedblright;
thus, a lamp will have $2$ or $3$ or $4$ adjacent lamps.) For example, if we
toggle the second lamp (from the left) in the topmost row in the above example
grid, then we obtain
\[%
\begin{tabular}
[c]{|c|c|c|c|c|}\hline
\textbf{0} & \textbf{1} & \textbf{1} & 1 & 1\\\hline
1 & \textbf{0} & 0 & 0 & 1\\\hline
1 & 0 & 0 & 1 & 0\\\hline
0 & 1 & 1 & 1 & 1\\\hline
0 & 1 & 0 & 0 & 0\\\hline
\end{tabular}
\
\]
(where the boldfaced numbers correspond to the lamps that have been affected
by the move).

Assume that all lamps are initially off. Can you (by a strategically chosen
sequence of moves) achieve a state in which all lamps are on?

[\textit{Remark:} You can play this game on
\url{https://codepen.io/wintlu/pen/ZJJLGz} .]
\end{exercise}

\begin{exercise}
\label{exe.21hw0.7}\ \ 

\begin{enumerate}
\item[\textbf{(a)}] How many of the numbers $0,1,\ldots,6$ appear as
remainders of a perfect square divided by $7$ ?

\item[\textbf{(b)}] How many of the numbers $0,1,\ldots,13$ appear as
remainders of a perfect square divided by $14$ ?
\end{enumerate}

What about replacing $7$ or $14$ by $n$? Can you do better than just squaring
them all?

[For example, $3$ of the numbers $0,1,\ldots,4$ appear as remainders of a
perfect square divided by $5$ -- namely, the three numbers $0,1,4$.]
\end{exercise}

\begin{exercise}
\label{exe.21hw0.8}Solve the following system of equations:%
\begin{align*}
a^{2}+b+c  &  =1;\\
b^{2}+c+a  &  =1;\\
c^{2}+a+b  &  =1
\end{align*}
for three complex numbers $a,b,c$.
\end{exercise}

The next exercise requires some preliminary discussion.

The following triangular table shows the binomial coefficients $\dbinom{n}{k}$
for $n\in\left\{  0,1,\ldots,7\right\}  $ and $k\in\left\{  0,1,\ldots
,n\right\}  $:
\[
\hspace{-0.5095pc}\begingroup\setlength{\tabcolsep}{4pt}
\renewcommand{\arraystretch}{0.5}
\begin{tabular}
[c]{ll|ccccccccccccccccc}%
\vphantom{$\overset{k=0}{\swarrow}$} &  & $\phantom{20}$ & $\phantom{20}$ &
$\phantom{20}$ & $\phantom{20}$ & $\phantom{20}$ & $\phantom{20}$ &  &  &  &
$\overset{k=0}{\swarrow}$ & $\phantom{20}$ & $\phantom{20}$ & $\phantom{20}$ &
$\phantom{20}$ & $\phantom{20}$ & $\phantom{20}$ & $\phantom{20}$\\
\vphantom{$\overset{k=0}{\swarrow}$}$n=0$ & $\rightarrow$ &  &  &  &  &  &  &
&  & $1$ &  & $\overset{k=1}{\swarrow}$ &  &  &  &  &  & \\
\vphantom{$\overset{k=0}{\swarrow}$}$n=1$ & $\rightarrow$ &  &  &  &  &  &  &
& $1$ &  & $1$ &  & $\overset{k=2}{\swarrow}$ &  &  &  &  & \\
\vphantom{$\overset{k=0}{\swarrow}$}$n=2$ & $\rightarrow$ &  &  &  &  &  &  &
$1$ &  & $2$ &  & $1$ &  & $\overset{k=3}{\swarrow}$ &  &  &  & \\
\vphantom{$\overset{k=0}{\swarrow}$}$n=3$ & $\rightarrow$ &  &  &  &  &  & $1$
&  & $3$ &  & $3$ &  & $1$ &  & $\overset{k=4}{\swarrow}$ &  &  & \\
\vphantom{$\overset{k=0}{\swarrow}$}$n=4$ & $\rightarrow$ &  &  &  &  & $1$ &
& $4$ &  & $6$ &  & $4$ &  & $1$ &  & $\overset{k=5}{\swarrow}$ &  & \\
\vphantom{$\overset{k=0}{\swarrow}$}$n=5$ & $\rightarrow$ &  &  &  & $1$ &  &
$5$ &  & $10$ &  & $10$ &  & $5$ &  & $1$ &  & $\overset{k=6}{\swarrow}$ & \\
\vphantom{$\overset{k=0}{\swarrow}$}$n=6$ & $\rightarrow$ &  &  & $1$ &  & $6$
&  & $15$ &  & $20$ &  & $15$ &  & $6$ &  & $1$ &  & $\overset{k=7}{\swarrow}%
$\\
\vphantom{$\overset{k=0}{\swarrow}$}$n=7$ & $\rightarrow$ &  & $1$ &  & $7$ &
& $21$ &  & $35$ &  & $35$ &  & $21$ &  & $7$ &  & $1$ &
\end{tabular}
\ \ \ \endgroup
\]
(This is part of what is known as \emph{Pascal's triangle}.)

Now, in this table, let us replace each even number by a $0$ and each odd
number by a $1$. We obtain
\[
\hspace{-0.5057pc}\begingroup\setlength{\tabcolsep}{4pt}
\renewcommand{\arraystretch}{0.5}
\begin{tabular}
[c]{ll|ccccccccccccccccc}%
\vphantom{$\overset{k=0}{\swarrow}$} &  & $\phantom{20}$ & $\phantom{20}$ &
$\phantom{20}$ & $\phantom{20}$ & $\phantom{20}$ & $\phantom{20}$ &
$\phantom{20}$ & $\phantom{20}$ & $\phantom{20}$ & $\overset{k=0}{\swarrow}$ &
$\phantom{20}$ & $\phantom{20}$ & $\phantom{20}$ & $\phantom{20}$ &
$\phantom{20}$ & $\phantom{20}$ & $\phantom{20}$\\
\vphantom{$\overset{k=0}{\swarrow}$}$n=0$ & $\rightarrow$ &  &  &  &  &  &  &
&  & $1$ &  & $\overset{k=1}{\swarrow}$ &  &  &  &  &  & \\
\vphantom{$\overset{k=0}{\swarrow}$}$n=1$ & $\rightarrow$ &  &  &  &  &  &  &
& $1$ &  & $1$ &  & $\overset{k=2}{\swarrow}$ &  &  &  &  & \\
\vphantom{$\overset{k=0}{\swarrow}$}$n=2$ & $\rightarrow$ &  &  &  &  &  &  &
$1$ &  & $0$ &  & $1$ &  & $\overset{k=3}{\swarrow}$ &  &  &  & \\
\vphantom{$\overset{k=0}{\swarrow}$}$n=3$ & $\rightarrow$ &  &  &  &  &  & $1$
&  & $1$ &  & $1$ &  & $1$ &  & $\overset{k=4}{\swarrow}$ &  &  & \\
\vphantom{$\overset{k=0}{\swarrow}$}$n=4$ & $\rightarrow$ &  &  &  &  & $1$ &
& $0$ &  & $0$ &  & $0$ &  & $1$ &  & $\overset{k=5}{\swarrow}$ &  & \\
\vphantom{$\overset{k=0}{\swarrow}$}$n=5$ & $\rightarrow$ &  &  &  & $1$ &  &
$1$ &  & $0$ &  & $0$ &  & $1$ &  & $1$ &  & $\overset{k=6}{\swarrow}$ & \\
\vphantom{$\overset{k=0}{\swarrow}$}$n=6$ & $\rightarrow$ &  &  & $1$ &  & $0$
&  & $1$ &  & $0$ &  & $1$ &  & $0$ &  & $1$ &  & $\overset{k=7}{\swarrow}$\\
\vphantom{$\overset{k=0}{\swarrow}$}$n=7$ & $\rightarrow$ &  & $1$ &  & $1$ &
& $1$ &  & $1$ &  & $1$ &  & $1$ &  & $1$ &  & $1$ &
\end{tabular}
\ \ \endgroup
\]
This looks rather similar to the third evolutionary stage of
\href{https://en.wikipedia.org/wiki/Sierpinski triangle}{Sierpinski's
triangle}:
\[
\pgfdeclarelindenmayersystem{Sierpinski triangle}{
\symbol{X}{\pgflsystemdrawforward}
\symbol{Y}{\pgflsystemdrawforward}
\rule{X -> X-Y+X+Y-X}
\rule{Y -> YY}
}\tikzset{
l-system={step=2cm/(2^3), order=3, angle=-120}
} \begin{tikzpicture}
\fill [black] (0,0) -- ++(0:2) -- ++(120:2) -- cycle;
\draw [draw=none] (0,0) l-system [l-system={Sierpinski triangle, axiom=X},fill=white];
\end{tikzpicture}
\]
(Each $0$ in the above table corresponds to a white $\triangle$ triangle, and
each $1$ corresponds to a black $\blacktriangle$ triangle.)

\begin{exercise}
\label{exe.21hw0.9}Where does this similarity come from?
\end{exercise}

\begin{exercise}
\label{exe.21hw0.10}A \emph{conic} means a curve of the form
\[
\left\{  \left(  x,y\right)  \in\mathbb{R}^{2}\mid ax^{2}+bxy+cy^{2}%
+dx+ey+f=0\right\}  ,
\]
where $a,b,c,d,e,f$ are six real numbers such that $\left(  a,b,c\right)
\neq\left(  0,0,0\right)  $. Examples of conics are

\begin{itemize}
\item any circle, e.g., the unit circle $\left\{  \left(  x,y\right)
\in\mathbb{R}^{2}\mid x^{2}+y^{2}=1\right\}  $;

\item more generally, any ellipse;

\item any parabola, e.g., $\left\{  \left(  x,y\right)  \in\mathbb{R}^{2}\mid
x^{2}+y=0\right\}  $;

\item any hyperbola, e.g., $\left\{  \left(  x,y\right)  \in\mathbb{R}^{2}\mid
xy=1\right\}  $ or $\left\{  \left(  x,y\right)  \in\mathbb{R}^{2}\mid
x^{2}-y^{2}=1\right\}  $;

\item the union of any two lines, e.g., $\left\{  \left(  x,y\right)
\in\mathbb{R}^{2}\mid xy=0\right\}  $.
\end{itemize}

\noindent A conic is said to be \emph{nondegenerate} if it is not the union of
two lines.

\begin{enumerate}
\item[\textbf{(a)}] What is the maximum number of points in which a
nondegenerate conic can intersect a line?

\item[\textbf{(b)}] What is the maximum number of points in which two distinct
nondegenerate conics can intersect each other?
\end{enumerate}
\end{exercise}

\newpage

\section{\label{chp.rings}Rings and ideals}

\subsection{\label{sec.rings.def}Defining rings (\cite[\S 7.1]{DumFoo04})}

\subsubsection{The definition}

You may have seen rings before, but beware: There are at least $4$
non-equivalent notions of a \textquotedblleft ring\textquotedblright, and the
one you know might be different from the one we'll use. Let us define the one
we want:\footnote{Let us recall the notion of a \textbf{monoid}, which will be
briefly used in this definition.
\par
Essentially, a monoid is just \textquotedblleft a group without
inverses\textquotedblright. The formal definition is as follows: A
\textbf{monoid} is a set $S$ equipped with a binary operation $\ast$ (that is,
a map from $S\times S$ to $S$) and a specified element $e\in S$ such that
\par
\begin{itemize}
\item the operation $\ast$ is associative (i.e., we have $a\ast\left(  b\ast
c\right)  =\left(  a\ast b\right)  \ast c$ for any $a,b,c\in S$, where we are
using the notation $x\ast y$ for the image of a given pair $\left(
x,y\right)  $ under the operation $\ast$), and
\par
\item the element $e$ is a neutral element for this operation $\ast$ (i.e., we
have $a\ast e=e\ast a=a$ for any $a\in S$).
\end{itemize}
\par
We denote this monoid by $\left(  S,\ast,e\right)  $.}

\begin{definition}
A \textbf{ring} means a set $R$ equipped with

\begin{itemize}
\item two binary operations (i.e., maps from $R\times R$ to $R$) that are
called \textbf{addition} and \textbf{multiplication} and are denoted by $+$
and $\cdot$, and

\item two elements of $R$ that are called \textbf{zero} and \textbf{unity} and
are denoted by $0$ and $1$,
\end{itemize}

\noindent such that the following properties (the \textquotedblleft%
\textbf{ring axioms}\textquotedblright) hold:

\begin{itemize}
\item $\left(  R,+,0\right)  $ is an abelian group. In other words:

\begin{itemize}
\item The operation $+$ is associative (i.e., we have $a+\left(  b+c\right)
=\left(  a+b\right)  +c$ for any $a,b,c\in R$).

\item The element $0$ is a neutral element for the operation $+$ (i.e., we
have $a+0=0+a=a$ for any $a\in R$).

\item Each element $a\in R$ has an inverse for the operation $+$ (i.e., an
element $b\in R$ satisfying $a+b=b+a=0$).

\item The operation $+$ is commutative (i.e., we have $a+b=b+a$ for any
$a,b\in R$).
\end{itemize}

\item $\left(  R,\cdot,1\right)  $ is a monoid. In other words:

\begin{itemize}
\item The operation $\cdot$ is associative (i.e., we have $a\cdot\left(
b\cdot c\right)  =\left(  a\cdot b\right)  \cdot c$ for any $a,b,c\in R$).

\item The element $1$ is a neutral element for the operation $\cdot$ (i.e., we
have $a\cdot1=1\cdot a=a$ for any $a\in R$).
\end{itemize}

Note that we \textbf{do not} require that the operation $\cdot$ be
commutative; nor do we require elements to have inverses for it.

\item The \textbf{distributive laws} hold in $R$: That is, for all $a,b,c\in
R$, we have%
\[
\left(  a+b\right)  \cdot c=a\cdot c+b\cdot c\ \ \ \ \ \ \ \ \ \ \text{and}%
\ \ \ \ \ \ \ \ \ \ a\cdot\left(  b+c\right)  =a\cdot b+a\cdot c.
\]

\item We have $0\cdot a=a\cdot0=0$ for each $a\in R$.
\end{itemize}

The zero of $R$ and the unity of $R$ don't necessarily have to be the numbers
$0$ and $1$; we just call them $0$ and $1$ because they behave similarly to
said numbers. If things can get ambiguous (i.e., if they actually differ from
the numbers $0$ and $1$), then we will call them $0_{R}$ and $1_{R}$ instead
(see below for some examples of this).

The unity $1$ of $R$ is also known as the \textbf{identity} or the
\textbf{one} of $R$ (but beware the ambiguity of the latter words).

The product $a\cdot b$ of two elements $a,b\in R$ is often denoted $ab$ (so we
omit the $\cdot$ sign) or occasionally $a\times b$ (we will avoid the latter notation).

The inverse of an element $a\in R$ in the abelian group $\left(  R,+,0\right)
$ will be called the \textbf{additive inverse} of $a$, and is denoted $-a$.

If $a,b\in R$, then the \textbf{difference} $a-b$ is defined to be the element
$a+\left(  -b\right)  \in R$.
\end{definition}

\begin{definition}
A ring $R$ is said to be \textbf{commutative} if its multiplication is
commutative (i.e., if $ab=ba$ for all $a,b\in R$).
\end{definition}

\subsubsection{\label{subsec.rings.def.exas}Some examples of rings}

You have probably seen various rings in your mathematical life. Here are some examples:

\begin{itemize}
\item The sets $\mathbb{Z}$, $\mathbb{Q}$ and $\mathbb{R}$ (endowed with the
usual addition, the usual multiplication, the usual $0$ and the usual $1$) are
commutative rings. The same holds for the set $\mathbb{C}$ of complex
numbers\footnote{This course will not rely overly much on complex numbers, as
we will be working in more abstract settings most of the time. Thus, if you
ignore everything I say about complex numbers and $\mathbb{C}$, you will miss
out on some examples and applications, but still understand the core of this
course.
\par
However, it will still be rather helpful to understand the construction of the
complex numbers, since we will imitate this construction later on. This
construction is covered in detail in \cite[\S 4.1]{19s} or in \cite[\S A.5]%
{BeaBla19}. See also Grant Sanderson's video
\url{https://www.youtube.com/watch?v=5PcpBw5Hbwo} on the geometric meaning of
complex numbers.}.

(Notice that existence of \textbf{multiplicative} inverses -- i.e., inverses
for the operation $\cdot$ -- is not required.)

\item The set $\mathbb{N}$ of nonnegative integers\footnote{Recall once again
that $\mathbb{N}$ is defined to be $\left\{  0,1,2,\ldots\right\}  $.} (again
endowed with the usual addition, the usual multiplication, the usual $0$ and
the usual $1$) is \textbf{not} a ring, since $\left(  \mathbb{N},+,0\right)  $
is not a group (only a monoid). It's what is called a \textbf{semiring}.

(Don't be fooled by the existence of negative numbers: The number $2$ has no
additive inverse in $\mathbb{N}$, even though $-2$ is an additive inverse for
it in $\mathbb{Z}$.)

\item We can define a commutative ring $\mathbb{Z}^{\prime}$ as follows: We
define a binary operation $\mathbin{\widetilde{\times}}$ on the set
$\mathbb{Z}$ by setting%
\[
a\mathbin{\widetilde{\times}}b=-ab\ \ \ \ \ \ \ \ \ \ \text{for all }\left(
a,b\right)  \in\mathbb{Z}\times\mathbb{Z}.
\]
Now, let $\mathbb{Z}^{\prime}$ be the \textbf{set} $\mathbb{Z}$, endowed with
the usual addition $+$ and the (unusual) multiplication
$\mathbin{\widetilde{\times}}$ and with the (usual) zero $0_{\mathbb{Z}%
^{\prime}}=0$ and the (unusual) unity $1_{\mathbb{Z}^{\prime}}=-1$. It is easy
to check that $\mathbb{Z}^{\prime}$ is a commutative ring. It is an example of
a ring whose unity is \textbf{not} equal to the integer $1$; the two
\textquotedblleft$1$\textquotedblright s in the equality $1_{\mathbb{Z}%
^{\prime}}=-1$ mean different things (the first \textquotedblleft%
$1$\textquotedblright\ is the unity of $\mathbb{Z}^{\prime}$, while the second
\textquotedblleft$1$\textquotedblright\ is the number $1$). This is why it is
important to never omit the subscript $\mathbb{Z}^{\prime}$ in
\textquotedblleft$1_{\mathbb{Z}^{\prime}}$\textquotedblright.

Note that I am denoting this new ring by $\mathbb{Z}^{\prime}$ rather than by
$\mathbb{Z}$, even though \textbf{as a set} it is identical with $\mathbb{Z}$.
I do this because I want to refer to a ring by just one single letter instead
of having to specify the addition and multiplication every time; but this
cannot go well if we use the same letter for different rings. A ring is not
just a set, but rather the entire package consisting of the set, the addition,
the multiplication, the zero and the unity. The rings $\mathbb{Z}$ and
$\mathbb{Z}^{\prime}$ have the same underlying set, but they differ in the
rest of the package (specifically, in the multiplication and the unity).

This all said, $\mathbb{Z}^{\prime}$ is not a very interesting ring: It is
essentially \textquotedblleft a copy of $\mathbb{Z}$, except that every
integer $n$ has been renamed as $-n$\textquotedblright. To formalize this
intuition, we would need to introduce the notion of a \textbf{ring
isomorphism}, which I will do soon (Definition \ref{def.ringmor.ringmor}
\textbf{(b)}); the main idea is that the bijection%
\[
\varphi:\mathbb{Z}\rightarrow\mathbb{Z}^{\prime},\ \ \ \ \ \ \ \ \ \ n\mapsto
-n
\]
\footnote{This notation means \textquotedblleft the map $\varphi$ from
$\mathbb{Z}$ to $\mathbb{Z}^{\prime}$ that sends each $n$ to $-n$%
\textquotedblright.} satisfies%
\begin{align*}
\varphi\left(  a+b\right)   &  =\varphi\left(  a\right)  +\varphi\left(
b\right)  \ \ \ \ \ \ \ \ \ \ \text{for all }\left(  a,b\right)  \in
\mathbb{Z}\times\mathbb{Z};\\
\varphi\left(  a\cdot b\right)   &  =\varphi\left(  a\right)
\mathbin{\widetilde{\times}}\varphi\left(  b\right)
\ \ \ \ \ \ \ \ \ \ \text{for all }\left(  a,b\right)  \in\mathbb{Z}%
\times\mathbb{Z};\\
\varphi\left(  0\right)   &  =0_{\mathbb{Z}^{\prime}};\\
\varphi\left(  1\right)   &  =1_{\mathbb{Z}^{\prime}}%
\end{align*}
(where the \textquotedblleft$0$\textquotedblright\ and the \textquotedblleft%
$1$\textquotedblright\ without subscripts are the usual numbers $0$ and $1$),
and thus the ring $\mathbb{Z}^{\prime}$ can be viewed as the ring $\mathbb{Z}$
with its elements \textquotedblleft relabelled\textquotedblright\ using this bijection.

\item The polynomial rings
\begin{align*}
\mathbb{Z}\left[  x\right]   &  =\left\{  \text{all polynomials in one
indeterminate }x\text{ with integer coefficients}\right\}  ,\\
\mathbb{Q}\left[  x\right]   &  =\left\{  \text{all polynomials in one
indeterminate }x\text{ with rational coefficients}\right\}  ,\\
\mathbb{R}\left[  x,y\right]   &  =\left\{  \text{all polynomials in two
indeterminates }x,y\text{ with real coefficients}\right\}
\end{align*}
and%
\begin{align*}
&  \mathbb{R}\left[  z_{1},z_{2},\ldots,z_{n}\right] \\
&  =\left\{  \text{all polynomials in }n\text{ indeterminates }z_{1}%
,z_{2},\ldots,z_{n}\text{ with real coefficients}\right\}
\end{align*}
(and many others, such as $\mathbb{Z}\left[  a,b\right]  $ or $\mathbb{C}%
\left[  u,p,q\right]  $) are commutative rings. (We won't give a formal
definition of polynomials until Section \ref{sec.polys1.polyrings}, but you
probably already have a rough idea of what polynomials are, and this idea
should suffice for now.)

\item The set of all functions $\mathbb{Q}\rightarrow\mathbb{Q}$ is a
commutative ring, where addition and multiplication are defined pointwise
(i.e., addition is defined by
\[
\left(  f+g\right)  \left(  x\right)  =f\left(  x\right)  +g\left(  x\right)
\ \ \ \ \ \ \ \ \ \ \text{for all }f,g:\mathbb{Q}\rightarrow\mathbb{Q}\text{
and }x\in\mathbb{Q},
\]
and multiplication is defined by
\[
\left(  fg\right)  \left(  x\right)  =f\left(  x\right)  \cdot g\left(
x\right)  \ \ \ \ \ \ \ \ \ \ \text{for all }f,g:\mathbb{Q}\rightarrow
\mathbb{Q}\text{ and }x\in\mathbb{Q},
\]
), where the zero is the \textquotedblleft constant-$0$\textquotedblright%
\ function (sending every $x\in\mathbb{Q}$ to $0$), and where the unity is the
\textquotedblleft constant-$1$\textquotedblright\ function (sending every
$x\in\mathbb{Q}$ to $1$). Of course, the same construction works if we
consider functions $\mathbb{R}\rightarrow\mathbb{C}$, or functions
$\mathbb{C}\rightarrow\mathbb{Q}$, or many other kinds of functions.

More generally, if $R$ is a ring, and if $S$ is any set, then the set of all
functions $S\rightarrow R$ is a ring (with $+$, $\cdot$, $0$ and $1$ defined
as above). If $R$ is commutative, then so is this new ring. (For some reason,
\cite{DumFoo04} requires $S$ to be nonempty here, but this is unnecessary.)

\item For comparison, here is a non-example: Consider the set of all functions
$\mathbb{Q}\rightarrow\mathbb{Q}$, where addition is defined pointwise (that
is, $\left(  f+g\right)  \left(  x\right)  =f\left(  x\right)  +g\left(
x\right)  $) but multiplication is composition (that is, $\left(  fg\right)
\left(  x\right)  =f\left(  g\left(  x\right)  \right)  $; in other words,
$fg=f\circ g$). This structure satisfies most axioms of a ring, but
\textbf{fails} one of the distributivity laws: namely, the second
distributivity law $a\cdot\left(  b+c\right)  =a\cdot b+a\cdot c$. Indeed,
three functions $a,b,c:\mathbb{Q}\rightarrow\mathbb{Q}$ usually do not satisfy
$a\circ\left(  b+c\right)  =a\circ b+a\circ c$. Such a structure is called a
\href{https://en.wikipedia.org/wiki/Near-ring}{\textbf{near-ring}}, but
(despite its name) does not really resemble rings in most of its features;
thus it shall only serve us here as a warning sign rather than as an actual example.
\end{itemize}

When we specify a ring, we don't need to provide its zero $0$ and its unity
$1$ (although, of course, they need to exist); they are uniquely determined by
the operations $+$ and $\cdot$. This is because they are neutral elements for
the operations $+$ and $\cdot$; but the neutral element of an operation is
always unique.\footnote{To wit: If $\ast$ is a binary operation on a set $S$,
and if $u$ and $v$ are two neutral elements for $\ast$, then $u\ast v=u$ (by
the neutrality of $v$) and $u\ast v=v$ (by the neutrality of $u$), so that
$u=u\ast v=v$. You have probably seen this argument in group theory, but it
does not require a group, just an arbitrary binary operation.}

Here are some more examples of rings:

\begin{itemize}
\item The set $\mathbb{S}$ of all real numbers of the form $a+b\sqrt{5}$ with
$a,b\in\mathbb{Q}$ (endowed with the usual notions of \textquotedblleft
addition\textquotedblright\ and \textquotedblleft
multiplication\textquotedblright\ defined for real numbers) is a commutative
ring. The \textquotedblleft hard\textquotedblright\ part of proving this is
showing that the product of two numbers of this form is again a number of this
form; but this is just a matter of computation:%
\begin{align*}
\left(  a+b\sqrt{5}\right)  \left(  c+d\sqrt{5}\right)   &  =ac+bc\sqrt
{5}+ad\sqrt{5}+bd\cdot5\\
&  =\underbrace{\left(  ac+5bd\right)  }_{\in\mathbb{Q}}+\underbrace{\left(
bc+ad\right)  }_{\in\mathbb{Q}}\sqrt{5}.
\end{align*}
Associativity, distributivity, etc. come for \textquotedblleft
free\textquotedblright, or, as we say, are \textbf{inherited from }%
$\mathbb{R}$ (meaning that we already know that they hold for $\mathbb{R}$, so
they must automatically hold for $\mathbb{S}$). Only the existence of additive
inverses (i.e., of inverses for the operation $+$) does not come for free
(sure, every element of $\mathbb{S}$ has an additive inverse in $\mathbb{R}$,
but we must show that it has an additive inverse in $\mathbb{S}$), but it is
easy to check (the additive inverse of $a+b\sqrt{5}\in\mathbb{S}$ is $\left(
-a\right)  +\left(  -b\right)  \sqrt{5}\in\mathbb{S}$).

The standard notation for this ring is $\mathbb{Q}\left[  \sqrt{5}\right]  $,
not $\mathbb{S}$. We will eventually see it as a particular case of a general construction.

\item We could define a different ring structure on the set $\mathbb{S}$ (that
is, a ring which, as a set, is identical with $\mathbb{S}$, but has a
different choice of operations) as follows: We define a binary operation
$\ast$ on $\mathbb{S}$ by setting%
\begin{align*}
\left(  a+b\sqrt{5}\right)  \ast\left(  c+d\sqrt{5}\right)   &  =ac+bd\sqrt
{5}\\
\ \ \ \ \ \ \ \ \ \ \text{for all }\left(  a,b\right)   &  \in\mathbb{Q}%
\times\mathbb{Q}\text{ and }\left(  c,d\right)  \in\mathbb{Q}\times\mathbb{Q}.
\end{align*}
This is well-defined, because every element of $\mathbb{S}$ can be written in
the form $a+b\sqrt{5}$ for a \textbf{unique} pair $\left(  a,b\right)
\in\mathbb{Q}\times\mathbb{Q}$. This is a consequence of the irrationality of
$\sqrt{5}$. You could not do this with $\sqrt{4}$ instead of $\sqrt{5}$ !

Now, let $\mathbb{S}^{\prime}$ be the set $\mathbb{S}$, endowed with the usual
addition $+$ and the (unusual) multiplication $\ast$, with the (usual) zero
$0_{\mathbb{S}^{\prime}}=0$ and with the (unusual) unity $1_{\mathbb{S}%
^{\prime}}=1+\sqrt{5}$ (not the integer $1$). It is easy to check that
$\mathbb{S}^{\prime}$ is a commutative ring. The \textbf{sets} $\mathbb{S}$
and $\mathbb{S}^{\prime}$ are identical, but the \textbf{rings} $\mathbb{S}$
and $\mathbb{S}^{\prime}$ are not: For example, the ring $\mathbb{S}^{\prime}$
has two nonzero elements whose product is $0$ (namely, $1\ast\sqrt{5}=0$),
whereas the ring $\mathbb{S}$ has no such things. Thus, we don't just have
$\mathbb{S}^{\prime}\neq\mathbb{S}$ as rings, but there is also no way to
regard $\mathbb{S}^{\prime}$ as \textquotedblleft a copy of $\mathbb{S}$ with
its elements renamed\textquotedblright\ (like we did with $\mathbb{Z}^{\prime
}$ and $\mathbb{Z}$). So a ring is much more than just a set; the $+$, $\cdot
$, $0$ and $1$ matter.

\item The set $\mathbb{S}_{3}$ of all real numbers of the form $a+b\sqrt[3]%
{5}$ with $a,b\in\mathbb{Q}$ (endowed with the usual addition, etc.) is
\textbf{not} a ring. Indeed, multiplication is not a binary operation on this
set $\mathbb{S}_{3}$, as you can see by noticing that%
\[
\underbrace{\left(  1+1\sqrt[3]{5}\right)  }_{\in\mathbb{S}_{3}}%
\underbrace{\left(  1+1\sqrt[3]{5}\right)  }_{\in\mathbb{S}_{3}}%
=1+2\sqrt[3]{5}+\left(  \sqrt[3]{5}\right)  ^{2}\notin\mathbb{S}_{3}.
\]
(Strictly speaking, this requires some work to prove -- how can we be sure
there are no $a,b\in\mathbb{Q}$ that satisfy $1+2\sqrt[3]{5}+\left(
\sqrt[3]{5}\right)  ^{2}=a+b\sqrt[3]{5}$ ? -- but I'm just making a point
about how not everything that looks like a ring is a ring. For a complete
proof of this claim, see Exercise \ref{exe.polring.xxx-r} \textbf{(d)}.)

\item For any $n\in\mathbb{N}$, the set $\mathbb{R}^{n\times n}$ of all
$n\times n$-matrices with real entries (endowed with matrix addition, matrix
multiplication, the zero matrix and the identity matrix) is a ring. It is not
commutative unless $n\leq1$, since we don't usually have $AB=BA$ for matrices.

More generally: If $R$ is any ring, and if $n\in\mathbb{N}$, then the set
$R^{n\times n}$ of all $n\times n$-matrices with entries in $R$ (endowed with
matrix addition, matrix multiplication, the zero matrix and the identity
matrix) is a ring. This is called the $n\times n$\textbf{-matrix ring} over
$R$; it is denoted by $R^{n\times n}$ or $M_{n}\left(  R\right)  $ or
$\operatorname*{M}\nolimits_{n}\left(  R\right)  $. Of course, the matrix
addition is defined in terms of the addition of $R$, and the matrix
multiplication is defined in terms of both $+$ and $\cdot$ operations of $R$.
Matrix rings are one of the main reasons people study noncommutative rings.

Note that if $R$ is not commutative, then this ring $R^{n\times n}$ is not
commutative even for $n=1$.

[You may ask what a $0\times0$-matrix is. Well, it pays off to be literal: It
is a table with $0$ rows, $0$ columns and $0$ entries.]
\end{itemize}

At this point, the \textquotedblleft endowed with...\textquotedblright\ phrase
has become somewhat of a ritual incantation: Most of our rings are endowed
with the exact operations ($+$ and $\cdot$) and special elements ($0$ and $1$)
you would guess if I just told you the set. Thus, in future, I will omit this
phrase unless I actually mean to endow the ring with some unexpected
operations. In particular, if I say that a set of numbers is a ring, then I
automatically understand it to be endowed with the usual addition, the usual
multiplication, the usual zero $0$ and the usual unity $1$, unless I say otherwise.

We continue with our litany of examples:

\begin{itemize}
\item Another famous noncommutative ring is the ring of \textbf{Hamiltonian
quaternions} $\mathbb{H}$. Its elements are the \textquotedblleft formal
expressions\textquotedblright\ of the form $a+bi+cj+dk$ with $a,b,c,d\in
\mathbb{R}$. (To be more rigorous, you can define them to be $4$-tuples
$\left(  a,b,c,d\right)  $ with $a,b,c,d\in\mathbb{R}$; the \textquotedblleft
formal\textquotedblright\ sum $a+bi+cj+dk$ can be viewed as just a fancy way
to write such a $4$-tuple.) Addition is defined by the distributive law:%
\begin{align*}
&  \left(  a+bi+cj+dk\right)  +\left(  a^{\prime}+b^{\prime}i+c^{\prime
}j+d^{\prime}k\right) \\
&  =\left(  a+a^{\prime}\right)  +\left(  b+b^{\prime}\right)  i+\left(
c+c^{\prime}\right)  j+\left(  d+d^{\prime}\right)  k.
\end{align*}
Multiplication is also defined by the distributive law using the formulas%
\[
i^{2}=j^{2}=k^{2}=-1,\qquad ij=-ji=k,\qquad jk=-kj=i,\qquad ki=-ik=j
\]
(and the rule that any real number should commute with any of $i,j,k$). For
example, the distributive law yields%
\[
\left(  1+i\right)  \left(  2+k\right)  =2+k+\underbrace{i\cdot2}%
_{=2i}+\underbrace{ik}_{=-j}=2+k+2i+\left(  -j\right)  =2+2i-j+k.
\]

We will see in Section \ref{sec.modules.H} that this $\mathbb{H}$ is indeed a
ring. It is not commutative. It is used in computer graphics (quaternions
encode rotations in 3D space), physics and number theory(!). In particular,
Lagrange's four-squares theorem, which says that any positive integer can be
written as a sum of four perfect squares, can be proved using quaternions!

\item If you like the empty set, you will enjoy the \textbf{zero ring}. This
is the ring which is defined as the one-element set $\left\{  0\right\}  $,
endowed with the only possible operations $+$ and $\cdot$ and its only
possible $0$ and $1$ (there is only one possibility for each of these, since
the ring only has one element!). So its zero and its unity are both $0$
(nobody said that they have to be distinct!), and it has $0+0=0$ and
$0\cdot0=0$.

The zero ring is, of course, commutative. It plays the same role in the world
of rings as the empty set does in the world of sets: It contains no
interesting information whatsoever, but its existence is important for things
to work.

Generally, a \textbf{trivial ring} is defined to be a ring containing only one
element. Every trivial ring can be viewed as the zero ring with its element
$0$ relabelled.

\item For every integer $n$, the residue classes\footnote{See Section
\ref{sec.intro.mod} for a refresher on residue classes.} of integers modulo
$n$ form a commutative ring, which is called $\mathbb{Z}/n\mathbb{Z}$ or
$\mathbb{Z}/n$ or $\mathbb{Z}_{n}$ (depending on the author; beware that
$\mathbb{Z}_{n}$ has two different meanings). You already know its additive
group $\left(  \mathbb{Z}/n,+,0\right)  $, which is classically called the
\textbf{cyclic group of order }$n$ (at least when $n$ is positive). The
multiplication is defined just as addition is: namely, we set
\[
\overline{a}\cdot\overline{b}=\overline{a\cdot b}\ \ \ \ \ \ \ \ \ \ \text{for
any }a,b\in\mathbb{Z}.
\]
(where the overline means \textquotedblleft residue class modulo
$n$\textquotedblright). This is all known as \textbf{modular arithmetic}.

The ring $\mathbb{Z}/n$ has $n$ elements when $n>0$. In particular,
$\mathbb{Z}/1\mathbb{Z}$ is a trivial ring. In contrast, as we already
mentioned in Section \ref{sec.intro.mod}, the ring $\mathbb{Z}/0\mathbb{Z}$ is
just $\mathbb{Z}$ with its elements relabelled (since a residue class modulo
$0$ only contains a single integer).

\item Here is yet another very small ring: Let $F_{4}$ be a set consisting of
four distinct elements $0,1,a,b$. Define two binary operations $+$ and $\cdot$
on $F_{4}$ by the following tables:%
\begin{align*}
&
\begin{tabular}
[c]{|c||c|c|c|c|}\hline
$x+y$ & $y=0$ & $y=1$ & $y=a$ & $y=b$\\\hline\hline
$x=0$ & $0$ & $1$ & $a$ & $b$\\\hline
$x=1$ & $1$ & $0$ & $b$ & $a$\\\hline
$x=a$ & $a$ & $b$ & $0$ & $1$\\\hline
$x=b$ & $b$ & $a$ & $1$ & $0$\\\hline
\end{tabular}
\\
& \\
&
\begin{tabular}
[c]{|c||c|c|c|c|}\hline
$x\cdot y$ & $y=0$ & $y=1$ & $y=a$ & $y=b$\\\hline\hline
$x=0$ & $0$ & $0$ & $0$ & $0$\\\hline
$x=1$ & $0$ & $1$ & $a$ & $b$\\\hline
$x=a$ & $0$ & $a$ & $b$ & $1$\\\hline
$x=b$ & $0$ & $b$ & $1$ & $a$\\\hline
\end{tabular}
\ \ .
\end{align*}
Then, I claim that $F_{4}$ is a ring (with zero $0$ and unity $1$). This can
be proved by meticulously checking that all the ring axioms are satisfied.
Arguably, this is rather boring\footnote{For example, checking associativity
of multiplication requires proving $\left(  ab\right)  c=a\left(  bc\right)  $
for $4^{3}=64$ different triples $\left(  a,b,c\right)  \in\left(
F_{4}\right)  ^{3}$.}. Eventually, we will find a way around this busywork by
constructing this ring $F_{4}$ in a different (more conceptual) way.

This ring $F_{4}$ is easily seen to be commutative (since the table for
$x\cdot y$ is symmetric across the diagonal). Its additive group $\left(
F_{4},+,0\right)  $ is the famous
\href{https://en.wikipedia.org/wiki/Klein_four-group}{Klein four-group},
characterized by the fact that it has four elements and each element $x$
satisfies $x+x=0$.

Note that both rings $F_{4}$ and $\mathbb{Z}/4$ are commutative rings with $4$
elements each. But they are rather different; in particular, $F_{4}$ is not
just \textquotedblleft$\mathbb{Z}/4$ with its labels taken
off\textquotedblright\footnote{Here is one difference: Every element $x\in
F_{4}$ satisfies $x+x=0$, but not every element $x\in\mathbb{Z}/4$ satisfies
this.}.

\item Here is one more ring with $4$ elements: Let $D_{4}$ be a set consisting
of four distinct elements $0,1,a,b$. Define a binary operation $+$ on $D_{4}$
in the same way as for $F_{4}$ in the previous example (i.e., using the exact
same table). Define a new binary operation $\cdot$ on $D_{4}$ by the following
table:%
\[%
\begin{tabular}
[c]{|c||c|c|c|c|}\hline
$x\cdot y$ & $y=0$ & $y=1$ & $y=a$ & $y=b$\\\hline\hline
$x=0$ & $0$ & $0$ & $0$ & $0$\\\hline
$x=1$ & $0$ & $1$ & $a$ & $b$\\\hline
$x=a$ & $0$ & $a$ & $0$ & $a$\\\hline
$x=b$ & $0$ & $b$ & $a$ & $1$\\\hline
\end{tabular}
\ \ \ \ .
\]
Then, $D_{4}$ is a commutative ring (with zero $0$ and unity $1$). This ring
differs crucially from both $\mathbb{Z}/4$ and $F_{4}$.

\item Yet another ring with $4$ elements is the ring $B_{4}$, which again
consists of four distinct elements $0,1,a,b$ and again has the same binary
operation $+$ as $F_{4}$ and $D_{4}$, but now has a multiplication $\cdot$
given by the table%
\[%
\begin{tabular}
[c]{|c||c|c|c|c|}\hline
$x\cdot y$ & $y=0$ & $y=1$ & $y=a$ & $y=b$\\\hline\hline
$x=0$ & $0$ & $0$ & $0$ & $0$\\\hline
$x=1$ & $0$ & $1$ & $a$ & $b$\\\hline
$x=a$ & $0$ & $a$ & $a$ & $0$\\\hline
$x=b$ & $0$ & $b$ & $0$ & $b$\\\hline
\end{tabular}
\ \ \ \ .
\]
This is again a commutative ring with zero $0$ and unity $1$.
\end{itemize}

More examples of rings can be found below (e.g., in Exercises
\ref{exe.21hw1.1}, \ref{exe.21hw1.4} and \ref{exe.21hw1.6}) and in the next
few exercises.

\begin{exercise}
\label{exe.rings.F8}There is a certain ring $F_{8}$ consisting of eight
distinct elements $0,1,a,b,c,d,e,f$. Its addition $+$ and its multiplication
$\cdot$ are given by the following tables:%
\begin{align*}
&
\begin{tabular}
[c]{|c||c|c|c|c|c|c|c|c|}\hline
$x+y$ & $y=0$ & $y=1$ & $y=a$ & $y=b$ & $y=c$ & $y=d$ & $y=e$ & $y=f$%
\\\hline\hline
$x=0$ & $0$ &  &  &  &  &  &  & \\\hline
$x=1$ &  & $0$ &  &  &  &  &  & \\\hline
$x=a$ &  & $b$ & $0$ &  &  &  &  & \\\hline
$x=b$ &  & $a$ & $1$ &  &  &  &  & \\\hline
$x=c$ &  & $d$ & $e$ &  & $0$ &  &  & \\\hline
$x=d$ &  & $c$ & $f$ &  & $1$ &  &  & \\\hline
$x=e$ &  & $f$ & $c$ &  & $a$ &  &  & \\\hline
$x=f$ &  & $e$ & $d$ &  & $b$ &  &  & \\\hline
\end{tabular}
\\
& \\
&
\begin{tabular}
[c]{|c||c|c|c|c|c|c|c|c|}\hline
$x\cdot y$ & $y=0$ & $y=1$ & $y=a$ & $y=b$ & $y=c$ & $y=d$ & $y=e$ &
$y=f$\\\hline\hline
$x=0$ &  &  &  &  &  &  &  & \\\hline
$x=1$ &  & $1$ &  &  &  &  &  & \\\hline
$x=a$ &  & $a$ & $c$ &  & $b$ &  &  & \\\hline
$x=b$ &  & $b$ &  &  &  &  &  & \\\hline
$x=c$ &  & $c$ & $b$ &  & $e$ &  &  & \\\hline
$x=d$ &  &  &  &  &  &  &  & \\\hline
$x=e$ &  &  &  &  &  &  &  & \\\hline
$x=f$ &  &  &  &  &  &  &  & \\\hline
\end{tabular}
\ \ \ \ \ .
\end{align*}
Oops, I lost most of the entries! Reconstruct all missing entries in the
tables. (You can take it for granted that $F_{8}$ really is a ring.) \medskip

[\textbf{Hint:} Start by showing that \textquotedblleft$0$\textquotedblright%
\ and \textquotedblleft$1$\textquotedblright\ really are the zero and the unity.]
\end{exercise}

Recall that complex numbers were defined as pairs $\left(  a,b\right)  $ of
real numbers, with entrywise addition
\[
\left(  a,b\right)  +\left(  c,d\right)  =\left(  a+c,\ b+d\right)
\]
and a certain weird-looking multiplication%
\[
\left(  a,b\right)  \cdot\left(  c,d\right)  =\left(  ac-bd,\ ad+bc\right)  .
\]
By setting $i=\left(  0,1\right)  $, and identifying each real number $r$ with
the complex number $\left(  r,0\right)  $, we can then write each complex
number $\left(  a,b\right)  $ in the familiar form $a+bi$.

In the next exercise, we will define a different kind of \textquotedblleft
numbers\textquotedblright: the \textbf{dual numbers}\footnote{Do they really
deserve to be called \textquotedblleft numbers\textquotedblright? Matter of
taste. But with the adjective \textquotedblleft dual\textquotedblright, the
terminology is unambiguous.}. They, too, are defined as pairs $\left(
a,b\right)  $ of real numbers, and again they are added entrywise, but their
multiplication is different from the multiplication of complex numbers:

\begin{exercise}
\label{exe.dualnums.ring}We define a \textbf{dual number} to be a pair
$\left(  a,b\right)  $ of two real numbers $a$ and $b$.

We let $\mathbb{D}$ be the set of all dual numbers.

Define an addition $+$ and a multiplication $\cdot$ on ${\mathbb{D}}$ by
setting
\begin{align*}
\left(  a,b\right)  +\left(  c,d\right)   &  =\left(  a+c,\ b+d\right)
\ \ \ \ \ \ \ \ \ \ \text{and}\\
\left(  a,b\right)  \cdot\left(  c,d\right)   &  =\left(  ac,\ ad+bc\right)
\end{align*}
for all $\left(  a,b\right)  \in{\mathbb{D}}$ and $\left(  c,d\right)
\in{\mathbb{D}}$. (Note that the only difference from complex numbers is the
definition of $\cdot$, which is lacking a $-bd$ term.)

\begin{enumerate}
\item[\textbf{(a)}] Prove that $\mathbb{D}$ becomes a commutative ring when
equipped with these two operations and with the zero $\left(  0,0\right)  $
and the unity $\left(  1,0\right)  $.
\end{enumerate}

This ring $\mathbb{D}$ will be called the \textbf{ring of dual numbers}.

We shall identify each real number $r$ with the dual number $\left(
r,0\right)  $.

We let $\varepsilon$ denote the dual number $\left(  0,1\right)  $.

\begin{enumerate}
\item[\textbf{(b)}] Prove that, for any $a,b\in\mathbb{R}$, we have
$a+b\varepsilon=\left(  a,b\right)  $ in $\mathbb{D}$ (where $a$, of course,
means the dual number $\left(  a,0\right)  $).

\item[\textbf{(c)}] Prove that $\varepsilon^{2}=0$ in $\mathbb{D}$.
\end{enumerate}
\end{exercise}

\begin{fineprint}
The usefulness of dual numbers stems from the fact that the dual number
$\varepsilon$ is a sort of \textquotedblleft algebraic
infinitesimal\textquotedblright\ (in the sense that $\varepsilon\neq0$ but
$\varepsilon^{2}=0$). We will eventually see (in Exercise
\ref{exe.dualnums.deriv}) that by evaluating a polynomial at a dual number of
the form $\left(  a,1\right)  $, we obtain not only the value of the
polynomial at $a$ but also its derivative at $a$.
\end{fineprint}

\subsubsection{Notes on the definitions}

\begin{remark}
Our above definition of a ring has some redundancies:

First of all, the $0\cdot a=a\cdot0=0$ axiom follows from distributivity and
the groupness of $\left(  R,+,0\right)  $. This is why it appears in
\cite{DumFoo04} as a theorem (\cite[Proposition 1 on page 226]{DumFoo04}), not
as an axiom.

Second, we can drop the \textquotedblleft abelian\textquotedblright\ in the
axiom \textquotedblleft$\left(  R,+,0\right)  $ is an abelian
group\textquotedblright; in other words, we can drop the requirement that
addition be commutative. This is because this requirement can be derived from
the remaining axioms (see \cite[page 223]{DumFoo04}). But this is a bit
artificial. I am aiming not for a minimal set of axioms, but for a reasonable
set of axioms that strikes the balance between usefulness (i.e., important
things are easy to derive from the axioms) and verifiability (i.e., it is easy
to check these axioms in meaningful cases).
\end{remark}

The kind of rings we defined above aren't the same kind of rings
\cite{DumFoo04} defines. The latter differ in that they are not required to
have a unity. I will call them \textbf{nonunital rings}:

\begin{definition}
A \textbf{nonunital ring} is defined in the same way as we defined a ring,
except we no longer require a unity, and we replace the axiom
\textquotedblleft$\left(  R,\cdot,1\right)  $ is a monoid\textquotedblright%
\ by \textquotedblleft the operation $\cdot$ is associative\textquotedblright.
In particular, any ring is a nonunital ring, but not vice versa.
\end{definition}

Note that the word \textquotedblleft nonunital\textquotedblright\ means
\textquotedblleft we don't require a unity\textquotedblright, not
\textquotedblleft the ring must not have a unity\textquotedblright.

For example, the set $2\mathbb{Z}$ of all even integers (i.e., the set
$\left\{  \ldots,-4,-2,0,2,4,\ldots\right\}  $) is a nonunital ring (when
equipped with the usual operations), but not a ring in our sense.

Beware:

\begin{itemize}
\item What we call a ring is called a \textquotedblleft ring with
identity\textquotedblright\ (or \textquotedblleft ring with $1$%
\textquotedblright) in \cite{DumFoo04}.

\item What we call a nonunital ring is just called a \textquotedblleft
ring\textquotedblright\ in \cite{DumFoo04}.
\end{itemize}

For an enlightening polemic about why rings in our sense (i.e., rings with a
unity) are a more important concept than nonunital rings, see \cite{Poonen18}.

Historically, the concept of a ring originated in the late 19th century in
number-theoretical considerations of Dedekind, Kronecker and Hilbert, and
emerged gradually from particular cases. Until the late 1930s, its definition
was rather fluid: Different authors imposed fewer or more axioms depending on
their specific needs. By the 1970s perhaps, the definition had stabilized
(thanks to the work of Noether and the textbook \cite{Waerde91} by van der
Waerden), except for the questions as to whether a ring should always have a
unity (i.e., the very point on which we disagree with \cite{DumFoo04}) and
occasionally as to whether a ring should always be commutative (we believe
they should not, but a number of mathematicians whose entire career is built
on the study of commutative rings prefer to have fewer words to type).

More on the history of rings can be found in
\url{https://mathshistory.st-andrews.ac.uk/HistTopics/Ring_theory/} .

\subsection{\label{sec.rings.calc}Calculating in rings}

\subsubsection{What works in arbitrary rings}

You can think of a commutative ring as a \textquotedblleft generalized number
system\textquotedblright. In particular, all computations that can be
performed with the operations $+$, $-$ and $\cdot$ on integers can be
similarly performed in any commutative ring. To some extent, this holds also
for general (noncommutative) rings.

For instance, if $a_{1},a_{2},\ldots,a_{n}$ are $n$ elements of a ring, then
the sum $a_{1}+a_{2}+\cdots+a_{n}$ is well-defined, and can be computed by
adding the elements $a_{1},a_{2},\ldots,a_{n}$ together in any
order\footnote{This means that, for example, the three sums $\left(  \left(
a_{1}+a_{2}\right)  +a_{3}\right)  +a_{4}$ and $\left(  a_{3}+\left(
a_{2}+a_{4}\right)  \right)  +a_{1}$ and $\left(  a_{2}+a_{3}\right)  +\left(
a_{4}+a_{1}\right)  $ are equal (for fixed elements $a_{1},a_{2},a_{3},a_{4}$
of a ring).
\par
This fact is known as \textbf{general commutativity} (or \textbf{generalized
commutativity}), and is true not just for rings but also (more generally) for
arbitrary abelian monoids (where $+$ is the operation of the monoid). For a
proof, see (e.g.) \cite[Theorem 2.118 \textbf{(a)}]{detnotes} (which
superficially only discusses real numbers, but gives a proof that applies
verbatim to any ring) or
\url{https://proofwiki.org/wiki/General_Commutativity_Theorem} (where the
operation we call $+$ is called $\circ$).}. More generally, finite sums of the
form $\sum_{s\in S}a_{s}$ are defined when the $a_{s}$ belong to a
ring\footnote{It should be kept in mind that empty sums (i.e., sums of the
form $\sum_{s\in\varnothing}a_{s}$) are defined to equal the zero of the
ring.}, and these sums behave just like finite sums of numbers.

The same holds for finite products when the ring is commutative. If the ring
is not commutative, then finite products in a specified order -- like
$a_{1}a_{2}\cdots a_{n}$ -- are still well-defined\footnote{This means that
the product is the same no matter \textquotedblleft where the parentheses are
placed\textquotedblright. For example, a product $abcde$ of five elements can
be computed as $\left(  \left(  \left(  ab\right)  c\right)  d\right)  e$ or
as $a\left(  b\left(  c\left(  de\right)  \right)  \right)  $ or as $\left(
a\left(  bc\right)  \right)  \left(  de\right)  $ or in several other ways,
and all these ways lead to the same result.
\par
This fact is known as \textbf{general associativity} (or \textbf{generalized
associativity}), and is true not just for rings but also (more generally) for
arbitrary monoids. For a proof, see (e.g.) \cite[Lemma 2.1.4]{Ford22} or
\url{https://groupprops.subwiki.org/wiki/Associative_implies_generalized_associative}
or
\url{https://math.stackexchange.com/questions/2459697/prove-generalized-associative-law}
.}, but unordered finite products -- like $\prod_{s\in S}a_{s}$ -- are not,
unless you have \textquotedblleft local commutativity\textquotedblright%
\ (i.e., the $a_{s}$ commute with each other).\footnote{It should be kept in
mind that empty products (i.e., products of the form $\prod_{s\in\varnothing
}a_{s}$) are defined to equal the unity of the ring.}

In any ring, subtraction satisfies the rules you would expect: For any two
elements $a,b$ of a ring, we have%
\begin{align*}
\left(  -a\right)  b  &  =a\left(  -b\right)  =-\left(  ab\right)  ;\\
\left(  -a\right)  \left(  -b\right)   &  =ab;\\
\left(  -1\right)  a  &  =-a.
\end{align*}
See \cite[\S 7.1, Proposition 1]{DumFoo04} for the easy proofs. Furthermore,
any three elements $a,b,c$ of a ring satisfy the \textquotedblleft subtractive
distributivity laws\textquotedblright%
\[
a\left(  b-c\right)  =ab-ac\ \ \ \ \ \ \ \ \ \ \text{and}%
\ \ \ \ \ \ \ \ \ \ \left(  a-b\right)  c=ac-bc.
\]
(These follow easily from the standard distributivity laws that are part of
the ring axioms.) Moreover, distributivity holds for finite sums as well:
\[
\sum_{s\in S}a_{s}b=\left(  \sum_{s\in S}a_{s}\right)
b\ \ \ \ \ \ \ \ \ \ \text{and}\ \ \ \ \ \ \ \ \ \ \sum_{s\in S}%
ba_{s}=b\left(  \sum_{s\in S}a_{s}\right)  .
\]

If $n$ is an integer and $a$ is an element of a ring $R$, then we define an
element $na$ of $R$ by%
\[
na=%
\begin{cases}
\underbrace{a+a+\cdots+a}_{n\text{ addends}}, & \text{if }n\geq0;\\
-\left(  \underbrace{a+a+\cdots+a}_{-n\text{ addends}}\right)  , & \text{if
}n<0
\end{cases}
.
\]
This defines the product of an integer with an element of $R$. This new
\textquotedblleft multiplication\textquotedblright\ operation is usually
called \textquotedblleft\textbf{scaling}\textquotedblright\ rather than
\textquotedblleft multiplication\textquotedblright, since its two inputs are
of different types: The first is an integer, while the second is an element of
$R$. In general, it is unrelated to the product of two elements of $R$,
although these operations usually agree when $\mathbb{Z}$ is a subset of $R$
(unless the multiplication of $R$ is defined in a particularly pathological
way\footnote{If $R$ is the ring $\mathbb{Z}^{\prime}$ defined in Subsection
\ref{subsec.rings.def.exas}, then the two operations do not agree, i.e., the
expression \textquotedblleft$na$\textquotedblright\ has different values
depending on whether you are viewing it as a product of two elements of $R$ or
as a product of an integer with an element of $R$. But this is no surprise,
since our definition of $\mathbb{Z}^{\prime}$ relied on deliberately altering
the multiplication.}).

We note that $0a=0$ for any $a\in R$, where the \textquotedblleft%
$0$\textquotedblright\ on the left hand side is the integer $0$. This is
because an empty sum is defined to be $0$.

If $n$ is a nonnegative integer and $a$ is an element of a ring $R$, then we
define an element $a^{n}$ of $R$ (called the $n$\textbf{-th power} of $a$) by%
\[
a^{n}=\underbrace{a\cdot a\cdot\cdots\cdot a}_{n\text{ factors}}.
\]
In particular, applying this definition to $n=0$, we obtain%
\[
a^{0}=\left(  \text{empty product}\right)  =1_{R}\ \ \ \ \ \ \ \ \ \ \text{for
each }a\in R.
\]
Furthermore, $a^{1}=a$ for each $a\in R$.

Thus we can scale elements of a ring by integers, and take them to nonnegative
integer powers. These operations satisfy the identities you would expect them
to satisfy: For example, for any $a,b\in R$ (with $R$ being a ring) and any
$n,m\in\mathbb{Z}$, we have%
\begin{align*}
\left(  n+m\right)  a  &  =na+ma;\\
n\left(  a+b\right)   &  =na+nb;\\
\left(  nm\right)  a  &  =n\left(  ma\right)  ;\\
\left(  -1\right)  a  &  =-a.
\end{align*}
Furthermore, for any $a\in R$ and any $n,m\in\mathbb{N}$, we have%
\begin{align*}
a^{n+m}  &  =a^{n}a^{m};\\
a^{nm}  &  =\left(  a^{n}\right)  ^{m}.
\end{align*}
Also,%
\begin{align*}
1^{n}  &  =1\ \ \ \ \ \ \ \ \ \ \text{for }n\in\mathbb{N};\\
0^{n}  &  =%
\begin{cases}
0, & \text{if }n>0;\\
1, & \text{if }n=0
\end{cases}
\ \ \ \ \ \ \ \ \ \ \text{for }n\in\mathbb{N}%
\end{align*}
(where, of course, the \textquotedblleft$1$\textquotedblright\ and
\textquotedblleft$0$\textquotedblright\ stand for $1_{R}$ and $0_{R}$, except
for the two \textquotedblleft$0$\textquotedblright s in \textquotedblleft%
$n>0$\textquotedblright\ and in \textquotedblleft$n=0$\textquotedblright.)

Moreover, if $a,b\in R$ satisfy $ab=ba$, then we have%
\begin{equation}
a^{i}b^{j}=b^{j}a^{i}\ \ \ \ \ \ \ \ \ \ \text{for }i,j\in\mathbb{N}
\label{eq.sec.rings.calc.aibj}%
\end{equation}
and
\begin{equation}
\left(  ab\right)  ^{n}=a^{n}b^{n}\ \ \ \ \ \ \ \ \ \ \text{for }%
n\in\mathbb{N} \label{eq.sec.rings.calc.abn}%
\end{equation}
and (the binomial formula)%
\begin{equation}
\left(  a+b\right)  ^{n}=\sum_{k=0}^{n}\dbinom{n}{k}a^{k}b^{n-k}%
\ \ \ \ \ \ \ \ \ \ \text{for }n\in\mathbb{N}. \label{eq.sec.rings.calc.binom}%
\end{equation}
All of this is proved just as for numbers.

\begin{exercise}
Actually prove these equalities (\ref{eq.sec.rings.calc.aibj}),
(\ref{eq.sec.rings.calc.abn}) and (\ref{eq.sec.rings.calc.binom}).
\end{exercise}

\begin{exercise}
Prove that the three equalities (\ref{eq.sec.rings.calc.aibj}),
(\ref{eq.sec.rings.calc.abn}) and (\ref{eq.sec.rings.calc.binom}) can fail if
we don't require $ab=ba$. For instance, find two $2\times2$-matrices
$a,b\in\mathbb{Q}^{2\times2}$ that violate (\ref{eq.sec.rings.calc.aibj}) for
$i=j=2$, violate (\ref{eq.sec.rings.calc.abn}) for $n=2$ and violate
(\ref{eq.sec.rings.calc.binom}) for $n=2$.
\end{exercise}

\begin{fineprint}
We note that even when two elements $a$ and $b$ of a ring $R$ don't satisfy
$ab=ba$, the $n$-th power $\left(  a+b\right)  ^{n}$ can be expanded using
distributivity; the result will just not usually be as nice as
(\ref{eq.sec.rings.calc.binom}). For example,%
\[
\left(  a+b\right)  ^{3}=a^{3}+a^{2}b+aba+ab^{2}+ba^{2}+bab+b^{2}a+b^{3}.
\]

\end{fineprint}

The following exercise generalizes the well-known \textquotedblleft geometric
sum\textquotedblright\ formula $1+q+q^{2}+\cdots+q^{n-1}=\dfrac{1-q^{n}}{1-q}$
(for $q\neq1$):

\begin{exercise}
Let $R$ be any ring. Let $a,b\in R$ satisfy $ab=ba$. Let $n\in\mathbb{N}$.
Prove that%
\[
a^{n}-b^{n}=\left(  a-b\right)  \cdot\underbrace{\sum_{k=0}^{n-1}%
a^{k}b^{n-1-k}}_{=a^{n-1}+a^{n-2}b+\cdots+ab^{n-2}+b^{n-1}}.
\]
(Note that when $n=0$, the sum $\sum_{k=0}^{n-1}a^{k}b^{n-1-k}$ is an empty
sum and thus equals $0$ by definition.)
\end{exercise}

The following exercises provide some practice with calculating in rings:

\begin{exercise}
\ \ 

\begin{enumerate}
\item[\textbf{(a)}] Without computing the integer $7^{4}$, prove that
$\overline{7}^{4}=\overline{1}$ in the ring $\mathbb{Z}/10$.

\item[\textbf{(b)}] Find a simple rule for the $k$-th power $\overline{7}^{k}$
of the element $\overline{7}$ in the ring $\mathbb{Z}/10$. Specifically, this
rule should express $\overline{7}^{k}$ in terms of the remainder that $k$
leaves when divided by $4$.

\item[\textbf{(c)}] What is the units digit of the number $7^{9999}$ ?
\end{enumerate}
\end{exercise}

\begin{exercise}
\ \ 

\begin{enumerate}
\item[\textbf{(a)}] Prove that every element $x\in\mathbb{Z}/7$ satisfies
$x^{7}=x$ in $\mathbb{Z}/7$.

\item[\textbf{(b)}] In the ring $\mathbb{H}$ of Hamiltonian quaternions (see
Subsection \ref{subsec.rings.def.exas}), compute $ijk$ and $\left(
1+i+j+k\right)  ^{2}$.
\end{enumerate}

Next, recall the ring $F_{4}$ constructed in Subsection
\ref{subsec.rings.def.exas}, with its four elements $0,1,a,b$.

\begin{enumerate}
\item[\textbf{(c)}] Prove that $a^{4}=a$ in this ring.

\item[\textbf{(d)}] What is $b^{4}$ ?
\end{enumerate}

[Part \textbf{(a)} is generalized in Proposition \ref{prop.ent.flt.Z/p} below,
whereas parts \textbf{(c)} and \textbf{(d)} are generalized in Proposition
\ref{prop.finfield.flt}.]
\end{exercise}

\begin{exercise}
\label{exe.quaternions.roots-of--1}\ \ 

\begin{enumerate}
\item[\textbf{(a)}] In the ring $\mathbb{H}$ of Hamiltonian quaternions (see
Subsection \ref{subsec.rings.def.exas}), prove that $\left(  ai+bj+ck\right)
^{2}=-\left(  a^{2}+b^{2}+c^{2}\right)  $ for any $a,b,c\in\mathbb{R}$.

\item[\textbf{(b)}] Conclude that there are infinitely many quaternions
$w\in\mathbb{H}$ satisfying $w^{2}=-1$.
\end{enumerate}
\end{exercise}

\subsubsection{What doesn't work in arbitrary rings}

Here are some things that might feel less familiar. Again, we let $R$ be a
ring, and we let $a,b$ be two elements of $R$.

\begin{itemize}
\item It is not always true that $a\neq0$ and $b\neq0$ implies $ab\neq0$. This
fails in the ring $\mathbb{Z}/6$ (for example, you can pick $a=\overline{2}$
and $b=\overline{3}$ to get $ab=\overline{2}\cdot\overline{3}=\overline
{2\cdot3}=\overline{6}=\overline{0}$, even though $a$ and $b$ are
$\neq\overline{0}$) and in matrix rings like $\mathbb{Z}^{2\times2}$ (here you
can pick $a=\left(
\begin{array}
[c]{cc}%
1 & 0\\
0 & 0
\end{array}
\right)  $ and $b=\left(
\begin{array}
[c]{cc}%
0 & 0\\
0 & 1
\end{array}
\right)  $ to get $ab=\left(
\begin{array}
[c]{cc}%
0 & 0\\
0 & 0
\end{array}
\right)  $, even though neither $a$ nor $b$ is the zero matrix).

\item It is not always true that $ab=1$ implies $ba=1$. This would be true in
the classical matrix rings $\mathbb{R}^{n\times n}$ and $\mathbb{C}^{n\times
n}$, in any commutative ring (for obvious reasons), and in any finite ring
(for less obvious reasons), but may fail in arbitrary rings. (Counterexamples
are not easy to find; see \cite[\S 7.1, exercise 30 (a)]{DumFoo04} for one.)
\end{itemize}

\subsubsection{Idempotents}

The next few exercises are concerned with the notion of \textquotedblleft
idempotent elements\textquotedblright\ of a ring. This is a useful notion, but
more importantly, the exercises should provide some practice with calculations
in rings.

\begin{definition}
\label{def.ring.idempotent}Let $R$ be a ring.

\begin{enumerate}
\item[\textbf{(a)}] An element $a$ of $R$ is said to be \textbf{idempotent} if
it satisfies $a^{2}=a$.

\item[\textbf{(b)}] An element $a$ of $R$ is said to be \textbf{involutive} if
it satisfies $a^{2}=1$.
\end{enumerate}
\end{definition}

Some examples first: The idempotent elements of $\mathbb{R}$ are $0$ and $1$.
The involutive elements of $\mathbb{R}$ are $1$ and $-1$. The ring
$\mathbb{Z}/6$ has four idempotent elements ($\overline{0}$, $\overline{1}$,
$\overline{3}$ and $\overline{4}$) and two involutive elements ($\overline{1}$
and $\overline{5}$). A matrix ring like $\mathbb{R}^{n\times n}$ usually has
infinitely many idempotent elements (viz., all projection matrices on
subspaces of $\mathbb{R}^{n}$) and infinitely many involutive elements (viz.,
all matrices $A$ satisfying $A^{2}=I_{n}$; for instance, all reflections
across hyperplanes are represented by such matrices).

\begin{exercise}
Let $p$ be a prime number, and let $k$ be a positive integer.

\begin{enumerate}
\item[\textbf{(a)}] Prove that the only idempotent elements of the ring
$\mathbb{Z}/p^{k}$ are $\overline{0}$ and $\overline{1}$.

\item[\textbf{(b)}] Now assume furthermore that $p\neq2$. Prove that the only
involutive elements of the ring $\mathbb{Z}/p^{k}$ are $\overline{1}$ and
$\overline{-1}$.
\end{enumerate}
\end{exercise}

\begin{exercise}
\label{exe.idemp.basics}Let $R$ be a ring. Prove the following:

\begin{enumerate}
\item[\textbf{(a)}] If $a$ is an idempotent element of $R$, then $1-a\in R$ is
again idempotent.

\item[\textbf{(b)}] If $a$ is an involutive element of $R$, then $-a\in R$ is
again involutive.

\item[\textbf{(c)}] If $a$ is an idempotent element of $R$, then $a^{n}=a$ for
each positive $n\in\mathbb{N}$.

\item[\textbf{(d)}] If $a$ is an idempotent element of $R$, then $\left(
1+a\right)  ^{n}=1+\left(  2^{n}-1\right)  a$ for each $n\in\mathbb{N}$.

\item[\textbf{(e)}] If $a$ is an involutive element of $R$, then $\left(
1+a\right)  ^{n}=2^{n-1}\left(  1+a\right)  $ for each positive $n\in
\mathbb{N}$.
\end{enumerate}
\end{exercise}

\begin{exercise}
\label{exe.21hw1.2}Let $R$ be a ring.

\begin{enumerate}
\item[\textbf{(a)}] Let $a\in R$. Prove that if $a$ is idempotent, then $1-2a$
is involutive.

\item[\textbf{(b)}] Now, assume that $2$ is \textbf{cancellable} in $R$; this
means that if $u$ and $v$ are two elements of $R$ satisfying $2u=2v$, then
$u=v$. Prove that the converse of the claim of part \textbf{(a)} holds: If
$a\in R$ is such that $1-2a$ is involutive, then $a$ is idempotent.

\item[\textbf{(c)}] Now, let $R=\mathbb{Z}/4\mathbb{Z}$. Find an element $a\in
R$ such that $1-2a$ is involutive, but $a$ is not idempotent.
\end{enumerate}
\end{exercise}

Exercise \ref{exe.21hw1.2} \textbf{(a)} assigns an involutive element to each
idempotent element of $R$. If $2$ is invertible in $R$ (that is, if the
element $2\cdot1_{R}$ has a multiplicative inverse), then this assignment is a
bijection (as can be easily derived from Exercise \ref{exe.21hw1.2}
\textbf{(b)}). Note that this assignment, when applied to a matrix ring
$\mathbb{R}^{n\times n}$, is exactly the assignment you would expect from the
geometric point of view: To the orthogonal projection on a hyperplane $H$, it
assigns the reflection in the perpendicular axis $H^{\perp}$. Exercise
\ref{exe.21hw1.2} \textbf{(c)} shows that we cannot drop the \textquotedblleft%
$2$ is cancellable\textquotedblright\ condition in Exercise \ref{exe.21hw1.2}
\textbf{(b)}.

\subsection{\label{sec.rings.subrings}Subrings (\cite[\S 7.1]{DumFoo04})}

\subsubsection{Definition}

Groups have subgroups; vector spaces have subspaces (and so do topological
spaces, although the two notions have little in common). Not surprisingly, the
same is true for rings, and you can guess the definition:

\begin{definition}
\label{def.ring.subring}Let $R$ be a ring. A \textbf{subring} of $R$ is a
subset $S$ of $R$ such that

\begin{itemize}
\item we have $a+b\in S$ for any $a,b\in S$;

\item we have $ab\in S$ for any $a,b\in S$;

\item we have $-a\in S$ for any $a\in S$;

\item we have $0\in S$ (where the $0$ means the zero of $R$);

\item we have $1\in S$ (where the $1$ means the unity of $R$).
\end{itemize}
\end{definition}

The five conditions in Definition \ref{def.ring.subring} are called the
\textquotedblleft\textbf{subring axioms}\textquotedblright. The first of these
five axioms is often reformulated as \textquotedblleft$S$ is closed under
addition\textquotedblright; the second then becomes \textquotedblleft$S$ is
closed under multiplication\textquotedblright; the third becomes
\textquotedblleft$S$ is closed under negation\textquotedblright. Thus, a
subring of a ring is a subset that is closed under addition, closed under
multiplication, closed under negation, and contains the zero and the unity.

The following is essentially obvious:

\begin{proposition}
\label{prop.ring.subring.is-ring}Let $S$ be a subring of a ring $R$. Then, $S$
automatically is a ring in its own right (with its operations $+$ and $\cdot$
obtained by restricting the corresponding operations of $R$, and with its
elements $0$ and $1$ passed down from $R$).
\end{proposition}

\subsubsection{\label{subsec.rings.subrings.exas}Examples}

Here are some examples of subrings:

\begin{itemize}
\item From the classical construction of the number systems, you know that
$\mathbb{Z}\subseteq\mathbb{Q}\subseteq\mathbb{R}\subseteq\mathbb{C}$. Each of
these three \textquotedblleft$\subseteq$\textquotedblright\ signs can be
strengthened to \textquotedblleft is a subring of\textquotedblright\ (for
example, $\mathbb{Z}$ is a subring of $\mathbb{Q}$).

\item We can extend this chain further to the right: $\mathbb{C}$ is a subring
of $\mathbb{H}$ (the quaternions).

\item However, we \textbf{cannot} extend this chain to the left: The only
subring of $\mathbb{Z}$ is $\mathbb{Z}$ itself. Indeed, a subring of
$\mathbb{Z}$ would have to contain $0$ and $1$ (by definition), thus also any
sum of the form $1+1+\cdots+1$ (since a subring is closed under addition),
i.e., any positive integer, and therefore also any negative integer (since it
is closed under negation), and thus any integer. But this means that it is
$\mathbb{Z}$.

\item There are lots of rings between $\mathbb{Z}$ and $\mathbb{Q}$ (that is,
rings $\mathbb{B}$ such that $\mathbb{Z}$ is a subring of $\mathbb{B}$ and
$\mathbb{B}$ in turn is a subring of $\mathbb{Q}$). You will see some of these
in Exercise \ref{exe.21hw1.1}. Here is another: Let $\mathbb{Q}%
_{\operatorname*{odd}}$ be the ring of all rational numbers of the form%
\[
\dfrac{a}{b}\ \ \ \ \ \ \ \ \ \ \text{with }a\in\mathbb{Z}\text{ being
arbitrary and }b\in\mathbb{Z}\text{ being odd.}%
\]
Then, $\mathbb{Q}_{\operatorname*{odd}}$ is a subring of $\mathbb{Q}$ (this is
pretty easy to check\footnote{For example, in order to check that
$\mathbb{Q}_{\operatorname*{odd}}$ is closed under addition, we need to verify
that the sum of two numbers of the form $\dfrac{a}{b}$ (with $a\in\mathbb{Z}$
being arbitrary and $b\in\mathbb{Z}$ being odd) is again a number of this
form. But this follows easily from the formula $\dfrac{a}{b}+\dfrac{c}%
{d}=\dfrac{ad+bc}{bd}$, along with the fact that a product of two odd numbers
is odd.}), and $\mathbb{Z}$ is a subring of $\mathbb{Q}_{\operatorname*{odd}}$.

\item There are myriad rings between $\mathbb{Q}$ and $\mathbb{R}$. For
example, the ring $\mathbb{S}$ from Subsection \ref{subsec.rings.def.exas} is
one of these.

\item There are no rings between $\mathbb{R}$ and $\mathbb{C}$. That is, if a
subring of $\mathbb{C}$ contains $\mathbb{R}$ as a subring, then this subring
must be either $\mathbb{R}$ or $\mathbb{C}$ itself. This is not hard to prove
(but I won't do so here).

\item There are rings between $\mathbb{Z}$ and $\mathbb{C}$ that are neither
subrings nor \textquotedblleft superrings\textquotedblright\ of $\mathbb{R}$.
A particularly important one is the ring $\mathbb{Z}\left[  i\right]  $ of
\textbf{Gaussian integers}. A \textbf{Gaussian integer} is a complex number of
the form $a+bi$ where $a$ and $b$ are integers (and where $i$ is the imaginary
unit $\sqrt{-1}$). For example, $3+5i$ and $-7+8i$ are Gaussian integers. It
is easy to see that $\mathbb{Z}\left[  i\right]  $ is indeed a subring of
$\mathbb{C}$, and of course $\mathbb{Z}$ is a subring of $\mathbb{Z}\left[
i\right]  $. But $\mathbb{Z}\left[  i\right]  $ is not an intermediate stage
on the $\mathbb{Z}\subseteq\mathbb{Q}\subseteq\mathbb{R}\subseteq\mathbb{C}$
\textquotedblleft chain\textquotedblright; it is a \textquotedblleft
detour\textquotedblright.

Likewise, there is a ring $\mathbb{Q}\left[  i\right]  $ of \textbf{Gaussian
rationals}, which are defined just as Gaussian integers but using rational
numbers (instead of integers) for $a$ and $b$. This ring $\mathbb{Q}\left[
i\right]  $ is sandwiched between $\mathbb{Q}$ and $\mathbb{C}$.

\item Recall the ring of functions from $\mathbb{Q}$ to $\mathbb{Q}$.
Similarly, there is a ring of functions from $\mathbb{R}$ to $\mathbb{R}$. The
latter has a subring consisting of all \textbf{continuous} functions from
$\mathbb{R}$ to $\mathbb{R}$. To see that this is indeed a subring, you need
to show that the sum and the product of two continuous functions are
continuous, that the negation $-f$ of a continuous function $f$ is continuous,
and that the constant-$0$ and constant-$1$ functions are continuous.

\item Let $n\in\mathbb{N}$, and let $R$ be any ring. Recall the matrix ring
$R^{n\times n}$, consisting of all $n\times n$-matrices with entries in $R$.
Then,%
\begin{align*}
R^{n\leq n}:=  &  \ \left\{  \text{all upper-triangular }n\times
n\text{-matrices with entries in }R\right\} \\
=  &  \ \left\{  \left(
\begin{array}
[c]{cccc}%
a_{1,1} & a_{1,2} & \cdots & a_{1,n}\\
0 & a_{2,2} & \cdots & a_{2,n}\\
\vdots & \vdots & \ddots & \vdots\\
0 & 0 & \cdots & a_{n,n}%
\end{array}
\right)  \ \mid\ a_{i,j}\in R\text{ for all }i\leq j\right\}
\end{align*}
is a subring of $R^{n\times n}$ (because the sum and the product of two
upper-triangular matrices are again upper-triangular, and because the zero
matrix and the identity matrix are upper-triangular). Similarly,%
\begin{align*}
R^{n\geq n}:=  &  \ \left\{  \text{all lower-triangular }n\times
n\text{-matrices with entries in }R\right\} \\
=  &  \ \left\{  \left(
\begin{array}
[c]{cccc}%
a_{1,1} & 0 & \cdots & 0\\
a_{2,1} & a_{2,2} & \cdots & 0\\
\vdots & \vdots & \ddots & \vdots\\
a_{n,1} & a_{n,2} & \cdots & a_{n,n}%
\end{array}
\right)  \ \mid\ a_{i,j}\in R\text{ for all }i\geq j\right\}
\end{align*}
is a subring of $R^{n\times n}$. The intersection $R^{n\leq n}\cap R^{n\geq
n}$ of these two subrings is again a subring of $R^{n\times n}$ (and its
elements are the diagonal $n\times n$-matrices). However,%
\begin{align*}
R_{\operatorname*{symm}}^{n\times n}:=  &  \ \left\{  \text{all symmetric
}n\times n\text{-matrices with entries in }R\right\} \\
=  &  \ \left\{  \left(
\begin{array}
[c]{cccc}%
a_{1,1} & a_{1,2} & \cdots & a_{1,n}\\
a_{2,1} & a_{2,2} & \cdots & a_{2,n}\\
\vdots & \vdots & \ddots & \vdots\\
a_{n,1} & a_{n,2} & \cdots & a_{n,n}%
\end{array}
\right)  \ \mid\ a_{i,j}=a_{j,i}\in R\text{ for all }i,j\right\}
\end{align*}
is \textbf{not} a subring of $R^{n\times n}$ unless $n\leq1$ or $R$ is trivial
(since in all other cases, it is easy to find two symmetric matrices whose
product is not symmetric\footnote{For example, if $n=2$, then the two
symmetric matrices $\left(
\begin{array}
[c]{cc}%
1 & 0\\
0 & 0
\end{array}
\right)  $ and $\left(
\begin{array}
[c]{cc}%
0 & 1\\
1 & 0
\end{array}
\right)  $ have product $\left(
\begin{array}
[c]{cc}%
1 & 0\\
0 & 0
\end{array}
\right)  \left(
\begin{array}
[c]{cc}%
0 & 1\\
1 & 0
\end{array}
\right)  =\left(
\begin{array}
[c]{cc}%
0 & 1\\
0 & 0
\end{array}
\right)  $, which is not symmetric.}).
\end{itemize}

\begin{warning}
Beware that our Definition \ref{def.ring.subring} does not agree with the
definition of a \textquotedblleft subring\textquotedblright\ in
\cite{DumFoo04}.

Indeed, \cite{DumFoo04} does \textbf{not} require $1\in S$ for a subring,
because \cite{DumFoo04} does not require rings to have a $1$ in the first
place. Thus, for example, the nonunital ring $2\mathbb{Z}$ is a subring of
$\mathbb{Z}$ in \cite{DumFoo04}'s sense (but not in our sense, since we don't
even count $2\mathbb{Z}$ as a ring). Even more confusingly, it can happen that
$S$ and $R$ are two rings in our sense (i.e., they both have unities), and $S$
is a subring of $R$ in \cite{DumFoo04}'s sense (i.e., $S$ satisfies our
definition of a subring, minus the \textquotedblleft$1\in S$\textquotedblright%
\ axiom), but not a subring of $R$ in our sense (because its unity is not the
unity of $R$). For example, the zero ring is a subring of $\mathbb{Z}$ in
\cite{DumFoo04}'s sense, but not in ours (since the unity of the zero ring is
the number $0$). Alas, there are less pathological examples, too, so this
isn't something you can ignore. For example, you can pretend that each
$2\times2$-matrix is secretly a $3\times3$-matrix by inserting a zero row at
the bottom and a zero column at the right (i.e., identifying each $2\times
2$-matrix $\left(
\begin{array}
[c]{cc}%
a & b\\
c & d
\end{array}
\right)  $ with the $3\times3$-matrix $\left(
\begin{array}
[c]{ccc}%
a & b & 0\\
c & d & 0\\
0 & 0 & 0
\end{array}
\right)  $; note that I am not saying you \textbf{should} do that), and this
makes $\mathbb{R}^{2\times2}$ a subring of $\mathbb{R}^{3\times3}$ in
\cite{DumFoo04}'s sense, but not in ours. Of course, this is one of the
situations where you really need subscripts under the \textquotedblleft%
$1$\textquotedblright\ to avoid confusing different unities.
\end{warning}

The following exercises provide several more examples of subrings:

\begin{exercise}
\label{exe.ring.subring.quadratic}Let $c$, $d$ and $g$ be three integers with
$g\neq0$. Assume that $d$ is not a perfect square (i.e., not the square of an integer).

Let $\zeta=\dfrac{c+\sqrt{d}}{g}$. (This is a real number if $d\geq0$, and a
complex number if $d<0$.) Set%
\begin{align*}
X  &  :=\left\{  a+b\zeta\ \mid\ a,b\in\mathbb{Z}\right\}
\ \ \ \ \ \ \ \ \ \ \text{and}\\
Y  &  :=\left\{  a+b\zeta\ \mid\ a,b\in\mathbb{Q}\right\}  .
\end{align*}

\begin{enumerate}
\item[\textbf{(a)}] Prove that $Y$ is always a subring of $\mathbb{C}$.

\item[\textbf{(b)}] Prove that $X$ is a subring of $\mathbb{C}$ if and only if
we have $g\mid2c$ and $c^{2}\equiv d\operatorname{mod}g^{2}$.
\end{enumerate}

[\textbf{Hint:} Show that $g^{2}\zeta^{2}=2cg\zeta+\left(  d-c^{2}\right)  $.]
\end{exercise}

\begin{fineprint}
Note that the ring $Y$ in Exercise \ref{exe.ring.subring.quadratic}
\textbf{(a)} generalizes the ring $\mathbb{Q}\left[  i\right]  $ of Gaussian
rationals (obtained by setting $c=0$ and $d=-1$ and $g=1$), whereas the ring
$X$ in Exercise \ref{exe.ring.subring.quadratic} \textbf{(b)} generalizes the
ring $\mathbb{Z}\left[  i\right]  $ of Gaussian integers (obtained in the same way).
\end{fineprint}

\begin{exercise}
\label{exe.21hw1.1}Fix an integer $m$. An $m$\textbf{-integer} shall mean a
rational number $r$ such that there exists a $k\in\mathbb{N}$ satisfying
$m^{k}r\in\mathbb{Z}$.

For example:

\begin{itemize}
\item Each integer $r$ is an $m$-integer (since $m^{k}r\in\mathbb{Z}$ for
$k=0$).

\item The rational number $\dfrac{5}{12}$ is a $6$-integer (since $6^{k}%
\cdot\dfrac{5}{12}\in\mathbb{Z}$ for $k=2$), but neither a $2$-integer nor a
$3$-integer (since multiplying it by a power of $2$ will not \textquotedblleft
get rid of\textquotedblright\ the prime factor $3$ in the denominator, and
vice versa\footnotemark).

\item The $1$-integers are the integers (since $1^{k}r=r$ for all $r$).

\item Every rational number $r$ is a $0$-integer (since $0^{k}r\in\mathbb{Z}$
for $k=1$).
\end{itemize}

Let $R_{m}$ denote the set of all $m$-integers. Prove that $R_{m}$ is a
subring of $\mathbb{Q}$.
\end{exercise}

\footnotetext{To make this more rigorous: If we had $2^{k}\cdot\dfrac{5}%
{12}\in\mathbb{Z}$ for some $k\in\mathbb{N}$, then we would have $12\mid
2^{k}\cdot5$, which would entail that $3\mid12\mid2^{k}\cdot5$, and thus $3$
would appear as a factor in the prime factorization of $2^{k}\cdot5$. But this
is absurd. Hence, $2^{k}\cdot\dfrac{5}{12}\in\mathbb{Z}$ cannot hold.
Similarly, $3^{k}\cdot\dfrac{5}{12}\in\mathbb{Z}$ cannot hold.}The ring
$R_{m}$ in Exercise \ref{exe.21hw1.1} is an example of a ring
\textquotedblleft between $\mathbb{Z}$ and $\mathbb{Q}$\textquotedblright\ (in
the sense that $\mathbb{Z}$ is a subring of $R_{m}$, while $R_{m}$ is a
subring of $\mathbb{Q}$). Note that $R_{1}=\mathbb{Z}$ and $R_{0}=\mathbb{Q}$,
whereas $R_{2}=R_{4}=R_{8}=\cdots$ is the ring of all rational numbers that
can be written in the form $a/2^{k}$ with $a\in\mathbb{Z}$ and $k\in
\mathbb{N}$.

Here is another example of a subring of a matrix ring $R^{n\times n}$:

\begin{exercise}
Let $n\in\mathbb{N}$. Let $R$ be any ring. An $n\times n$-matrix $A=\left(
\begin{array}
[c]{cccc}%
a_{1,1} & a_{1,2} & \cdots & a_{1,n}\\
a_{2,1} & a_{2,2} & \cdots & a_{2,n}\\
\vdots & \vdots & \ddots & \vdots\\
a_{n,1} & a_{n,2} & \cdots & a_{n,n}%
\end{array}
\right)  \in R^{n\times n}$ will be called \textbf{centrosymmetric} if it
satisfies%
\[
a_{i,j}=a_{n+1-i,\ n+1-j}\ \ \ \ \ \ \ \ \ \ \text{for all }i,j\in\left\{
1,2,\ldots,n\right\}  .
\]
(Visually, this means that $A$ is preserved under \textquotedblleft%
$180^{\circ}$-rotation\textquotedblright, i.e., that any two cells of $A$ that
are mutually symmetric across the center of the matrix have the same entry.
For example, a centrosymmetric $4\times4$-matrix has the form $\left(
\begin{array}
[c]{cccc}%
a & b & c & d\\
e & f & g & h\\
h & g & f & e\\
d & c & b & a
\end{array}
\right)  $ for $a,b,\ldots,h\in R$.)

Prove that the set $\left\{  \text{all centrosymmetric }n\times
n\text{-matrices with entries in }R\right\}  $ is a subring of $R^{n\times n}%
$. \medskip

[\textbf{Hint:} This can be done in a particularly slick way as follows: Let
$W$ be the $n\times n$-matrix obtained from the identity matrix $I_{n}$ by a
horizontal reflection (or, equivalently, a vertical reflection). For example,
if $n=4$, then $W=\left(
\begin{array}
[c]{cccc}%
0 & 0 & 0 & 1\\
0 & 0 & 1 & 0\\
0 & 1 & 0 & 0\\
1 & 0 & 0 & 0
\end{array}
\right)  $. Now, show that an $n\times n$-matrix $A$ is centrosymmetric if and
only if it satisfies $AW=WA$.]
\end{exercise}

\begin{fineprint}
Subrings can be used to answer a curious question: How small can a
noncommutative ring be? A moment of thought leads us to the ring $\left(
\mathbb{Z}/2\right)  ^{2\leq2}$ of upper-triangular matrices with entries in
$\mathbb{Z}/2$; this ring is noncommutative and has size $8$ (since there are
$2$ choices for each of the three entries of such a matrix, not counting the
bottom-left entry because that entry must be $\overline{0}$). Thus, a
noncommutative ring can have size $8$. A smaller size is not possible, as
follows from the following exercise:
\end{fineprint}

\begin{exercise}
\label{exe.rings.size-8-comm}Let $R$ be a finite ring. Assume that its size
$\left\vert R\right\vert $ is either a prime number $p$ or a product $pq$ of
two (not necessarily distinct!) prime numbers $p$ and $q$. Our goal is to show
that $R$ is commutative.

Consider the abelian group $\left(  R,+,0\right)  $. If $u_{1},u_{2}%
,\ldots,u_{k}$ are any elements of $R$, then $\left\langle u_{1},u_{2}%
,\ldots,u_{k}\right\rangle $ shall denote the subgroup of this abelian group
$\left(  R,+,0\right)  $ generated by $u_{1},u_{2},\ldots,u_{k}$. (Explicitly,
this subgroup consists of all sums of the form $a_{1}u_{1}+a_{2}u_{2}%
+\cdots+a_{k}u_{k}$ with $a_{1},a_{2},\ldots,a_{k}\in\mathbb{Z}$.)

Let $x,y\in R$. Consider the following chain of subgroups of $\left(
R,+,0\right)  $:%
\[
0\leq\left\langle 1\right\rangle \leq\left\langle x,1\right\rangle \leq R.
\]
(The symbol $\leq$ means \textquotedblleft subgroup of\textquotedblright.)

\begin{enumerate}
\item[\textbf{(a)}] Prove that at least one of the three \textquotedblleft%
$\leq$\textquotedblright\ signs in this chain must be an \textquotedblleft%
$=$\textquotedblright\ sign.

\item[\textbf{(b)}] Prove that $xy=yx$ if the first \textquotedblleft$\leq
$\textquotedblright\ sign is a \textquotedblleft$=$\textquotedblright\ sign.

\item[\textbf{(c)}] Prove that $xy=yx$ if the second \textquotedblleft$\leq
$\textquotedblright\ sign is a \textquotedblleft$=$\textquotedblright\ sign.

\item[\textbf{(d)}] Prove that $xy=yx$ if the third \textquotedblleft$\leq
$\textquotedblright\ sign is a \textquotedblleft$=$\textquotedblright\ sign.

\item[\textbf{(e)}] Conclude that $R$ is commutative.
\end{enumerate}

[\textbf{Hint:} In part \textbf{(a)}, recall Lagrange's theorem about
subgroups, and observe that a number $m$ of the form $p$ or $pq$ cannot have a
nontrivial chain of three divisors $1\mid d\mid e\mid m$. Parts \textbf{(b)},
\textbf{(c)} and \textbf{(d)} are easy in their own ways.]
\end{exercise}

\begin{fineprint}
See Exercise \ref{exe.CRT-finrings} for a generalization of Exercise
\ref{exe.rings.size-8-comm} \textbf{(e)}.
\end{fineprint}

The following exercise gives a way to construct new subrings out of old ones:

\begin{exercise}
\ \ 

\begin{enumerate}
\item[\textbf{(a)}] Let $R$ be a ring. Let $S$ and $T$ be two subrings of $R$.
Prove that $S\cap T$ is again a subring of $R$.

\item[\textbf{(b)}] For each integer $m$, define the subring $R_{m}$ of
$\mathbb{Q}$ as in Exercise \ref{exe.21hw1.1}. Prove that $R_{m}\cap
R_{n}=R_{\gcd\left(  m,n\right)  }$ for all $m,n\in\mathbb{Z}$.
\end{enumerate}

[\textbf{Hint:} Part \textbf{(a)} is easy. Part \textbf{(b)} requires a bit of
elementary number theory.]
\end{exercise}

\subsubsection{A first application}

We haven't proved much so far, but we are already able to reap some first
rewards. Namely, we shall prove two properties of the famous Fibonacci
sequence. We recall its definition:

\begin{definition}
\label{def.fibonacci.fib}The \textbf{Fibonacci sequence} is the sequence of
integers defined recursively by
\[
f_{0}=0,\qquad f_{1}=1,\qquad\text{and}\qquad f_{n}=f_{n-1}+f_{n-2}\text{ for
all }n\geq2.
\]

\end{definition}

The first entries of this sequence are
\[%
\begin{tabular}
[c]{|c||c|c|c|c|c|c|c|c|c|c|c|c|c|}\hline
$n$ & $0$ & $1$ & $2$ & $3$ & $4$ & $5$ & $6$ & $7$ & $8$ & $9$ & $10$ & $11$
& $12$\\\hline
$f_{n}$ & $0$ & $1$ & $1$ & $2$ & $3$ & $5$ & $8$ & $13$ & $21$ & $34$ & $55$
& $89$ & $144$\\\hline
\end{tabular}
\ \ \ \ \ .
\]
Much more about this sequence can be found (e.g.) in \cite{Vorobi02} or
\cite{mps}. The entries $f_{0},f_{1},f_{2},\ldots$ of this sequence are known
as the \textbf{Fibonacci numbers}.

We shall prove the following two facts:

\begin{proposition}
\label{prop.fibonacci.add}The Fibonacci sequence $\left(  f_{0},f_{1}%
,f_{2},\ldots\right)  $ satisfies
\[
f_{n+m}=f_{n}f_{m+1}+f_{n-1}f_{m}%
\]
for all positive integers $n$ and all nonnegative integers $m$.
\end{proposition}

\begin{proposition}
\label{prop.fibonacci.div}The Fibonacci sequence $\left(  f_{0},f_{1}%
,f_{2},\ldots\right)  $ satisfies
\[
f_{d}\mid f_{dn}\ \ \ \ \ \ \ \ \ \ \text{for any nonnegative integers
}d\text{ and }n.
\]

\end{proposition}

There are many proofs of these propositions (see, e.g., \cite[Exercise
4.9.3]{mps} and \cite[Exercise 4.9.7]{mps} for generalizations proved in a
very elementary way). We will give a proof that uses a certain commutative
subring $\mathcal{F}$ of the matrix ring $\mathbb{Z}^{2\times2}$ as a tool:

\begin{exercise}
\label{exe.21hw1.6}Let $A$ be the $2\times2$-matrix $%
\begin{pmatrix}
0 & 1\\
1 & 1
\end{pmatrix}
\in\mathbb{Z}^{2\times2}$. Consider also the identity matrix $I_{2}%
\in\mathbb{Z}^{2\times2}$.

\begin{enumerate}
\item[\textbf{(a)}] Prove that $A^{2}=A+I_{2}$.
\end{enumerate}

Now, let $\mathcal{F}$ be the subset
\[
\left\{  aA+bI_{2}\mid a,b\in\mathbb{Z}\right\}  =\left\{
\begin{pmatrix}
b & a\\
a & a+b
\end{pmatrix}
\mid a,b\in\mathbb{Z}\right\}
\]
of the matrix ring $\mathbb{Z}^{2\times2}$.

\begin{enumerate}
\item[\textbf{(b)}] Prove that the set $\mathcal{F}$ is a \textbf{commutative}
subring of $\mathbb{Z}^{2\times2}$.
\end{enumerate}

Next, let $\left(  f_{0},f_{1},f_{2},\ldots\right)  $ be the Fibonacci sequence.

\begin{enumerate}
\item[\textbf{(c)}] Prove that $A^{n}=f_{n}A+f_{n-1}I_{2}$ for all positive
integers $n$.

\item[\textbf{(d)}] Prove that $f_{n+m}=f_{n}f_{m+1}+f_{n-1}f_{m}$ for all
positive integers $n$ and all $m\in\mathbb{N}$. (This is Proposition
\ref{prop.fibonacci.add}.)
\end{enumerate}

Now, define a further matrix $B\in\mathcal{F}$ by $B=\left(  -1\right)
A+1I_{2}=I_{2}-A$.

\begin{enumerate}
\item[\textbf{(e)}] Prove that $B^{2}=B+I_{2}$ and $B^{n}=f_{n}B+f_{n-1}I_{2}$
for all positive integers $n$.

\item[\textbf{(f)}] Prove that $A^{n}-B^{n}=f_{n}\left(  A-B\right)  $ for all
$n\in\mathbb{N}$.

\item[\textbf{(g)}] Prove that $f_{d}\mid f_{dn}$ for any nonnegative integers
$d$ and $n$. (This is Proposition \ref{prop.fibonacci.div}.)
\end{enumerate}

[\textbf{Hint:} In part \textbf{(b)}, don't forget to check commutativity! It
is not inherited from $\mathbb{Z}^{2\times2}$, since $\mathbb{Z}^{2\times2}$
is not commutative.

One way to prove part \textbf{(d)} is by comparing the $\left(  1,1\right)
$-th entries of the two (equal) matrices $A^{n}A^{m+1}$ and $A^{n+m+1}$, after
first using part \textbf{(c)} to compute these matrices.

For part \textbf{(g)}, compare the $\left(  1,1\right)  $-th entries of the
matrices $A^{d}-B^{d}$ and $A^{dn}-B^{dn}$, after first proving that
$A^{d}-B^{d}\mid A^{dn}-B^{dn}$ in the commutative ring $\mathcal{F}$. Note
that divisibility is a tricky concept in general rings, but $\mathcal{F}$ is a
commutative ring, which lets many arguments from the integer setting go
through unchanged in $\mathcal{F}$.]
\end{exercise}

\subsubsection{More computational exercises}

\begin{exercise}
\label{exe.jacobi-id}Let $R$ be any ring. For any two elements $a,b\in R$, we
define the element $\left[  a,b\right]  $ of $R$ by%
\[
\left[  a,b\right]  :=ab-ba.
\]
This element $\left[  a,b\right]  $ is called the \textbf{commutator} of $a$
and $b$ (as it \textquotedblleft measures\textquotedblright\ how much $a$ and
$b$ violate the commutative law $ab=ba$). Don't confuse it with the
group-theoretical commutator $aba^{-1}b^{-1}$, which is also denoted by
$\left[  a,b\right]  $ (but is defined for groups rather than rings).

Prove that every three elements $a,b,c\in R$ satisfy the \textbf{Leibniz
identity}%
\[
\left[  a,bc\right]  =\left[  a,b\right]  c+b\left[  a,c\right]
\]
and the \textbf{Jacobi identity}%
\[
\left[  a,\left[  b,c\right]  \right]  +\left[  b,\left[  c,a\right]  \right]
+\left[  c,\left[  a,b\right]  \right]  =0.
\]

\end{exercise}

\begin{exercise}
Let $R$ be any ring. The determinant of a $2\times2$-matrix $A\in R^{2\times
2}$ is usually defined only when $R$ is commutative, but let us (for this
specific exercise) define it in general by the formula%
\[
\det\left(
\begin{array}
[c]{cc}%
a & b\\
c & d
\end{array}
\right)  :=ad-bc\ \ \ \ \ \ \ \ \ \ \text{for any }\left(
\begin{array}
[c]{cc}%
a & b\\
c & d
\end{array}
\right)  \in R^{2\times2}.
\]

Prove that the equality $\det\left(  AB\right)  =\det A\cdot\det B$ holds for
every pair of two matrices $A,B\in R^{2\times2}$ if and only if $R$ is
commutative. \medskip

[\textbf{Hint:} The \textquotedblleft if\textquotedblright-direction can be
considered well-known from linear algebra.]
\end{exercise}

\begin{exercise}
Let $R$ be a commutative ring in which $2\cdot1_{R}=0_{R}$. (Examples of such
rings are $\mathbb{Z}/2$ or polynomial rings over $\mathbb{Z}/2$.)

Prove that the set of all idempotent elements $a\in R$ is a subring of $R$.
\end{exercise}

\begin{exercise}
Let $R$ be a commutative ring in which $8\cdot1_{R}=0_{R}$. (Examples of such
rings are $\mathbb{Z}/2$, $\mathbb{Z}/4$ and $\mathbb{Z}/8$, but there are
also many others, such as polynomial rings over $\mathbb{Z}/8$.)

Prove that the set of all elements $a\in R$ satisfying $\left(  1-2a\right)
^{2}=1$ is a subring of $R$.
\end{exercise}

\subsubsection{The center of a ring, and the centralizer of a subset}

Here is yet another way to construct subrings of a ring:

\begin{definition}
\label{def.center}Let $R$ be a ring.

\begin{enumerate}
\item[\textbf{(a)}] An element $a\in R$ is said to be \textbf{central} if all
$b\in R$ satisfy $ab=ba$. (In other words, $a$ is central if and only if $a$
commutes with every element of $R$.)

\item[\textbf{(b)}] The \textbf{center} of $R$ is the set of all central
elements of $R$. This set is denoted by $Z\left(  R\right)  $.
\end{enumerate}
\end{definition}

\begin{exercise}
\label{exe.center.subring}Let $R$ be a ring. Prove that:

\begin{enumerate}
\item[\textbf{(a)}] The center $Z\left(  R\right)  $ of $R$ is a commutative
subring of $R$.

\item[\textbf{(b)}] We have $Z\left(  R\right)  =R$ if and only if $R$ is commutative.

\item[\textbf{(c)}] All elements of the form $n\cdot1_{R}$ for $n\in
\mathbb{Z}$ belong to $Z\left(  R\right)  $.
\end{enumerate}
\end{exercise}

\begin{exercise}
\ \ 

\begin{enumerate}
\item[\textbf{(a)}] Prove that $Z\left(  \mathbb{C}\right)  =\mathbb{C}$ and
$Z\left(  \mathbb{H}\right)  =\mathbb{R}$ (where $\mathbb{H}$ is the ring of quaternions).

\item[\textbf{(b)}] Compute $Z\left(  \mathbb{R}^{2\times2}\right)  $ and
$Z\left(  \mathbb{R}^{2\leq2}\right)  $. (In other words, find the $2\times
2$-matrices that commute with all $2\times2$-matrices, and find the
upper-triangular $2\times2$-matrices that commute with all upper-triangular
$2\times2$-matrices.)
\end{enumerate}
\end{exercise}

The previous exercise illustrates a somewhat slippery point: If $R$ is a
subring of a ring $S$, then $Z\left(  R\right)  $ doesn't have to be a subring
of $Z\left(  S\right)  $. An element of $R$ that commutes with all elements of
$R$ might still fail to commute with some elements of $S$.

\begin{exercise}
Let $R$ be a ring. Let $a,b\in R$ be such that $a+b$ is central. Prove that
$ab=ba$.
\end{exercise}

\begin{exercise}
Let $R$ be a ring. Let $a,b\in R$ be such that $ab$ is central. Prove that
$\left(  ab\right)  ^{n}=a^{n}b^{n}$ for all $n\in\mathbb{N}$.
\end{exercise}

A generalization of the center is the \textquotedblleft
centralizer\textquotedblright\ of a subset of a ring:\footnote{If you have
seen centralizers in groups, you'll recognize this as an analogous notion.}

\begin{definition}
\label{def.centralizer.centralizer}Let $R$ be a ring. Let $S$ be a subset of
$R$.

\begin{enumerate}
\item[\textbf{(a)}] An element $a\in R$ is said to \textbf{centralize} $S$ if
and only if all $b\in S$ satisfy $ab=ba$. (In other words, $a$ centralizes $S$
if and only if $a$ commutes with every element of $S$.)

\item[\textbf{(b)}] The \textbf{centralizer} of $S$ in $R$ is the set of all
elements of $R$ that centralize $S$. This set is denoted by $Z_{R}\left(
S\right)  $.
\end{enumerate}
\end{definition}

Note that $Z_{R}\left(  R\right)  =Z\left(  R\right)  $ is the center of $R$.

\begin{exercise}
\label{exe.centralizer.subring}Let $R$ be a ring. Let $S$ be a subset of $R$.
Prove that:

\begin{enumerate}
\item[\textbf{(a)}] The centralizer $Z_{R}\left(  S\right)  $ is a subring of
$R$.

\item[\textbf{(b)}] We have $Z_{R}\left(  \varnothing\right)  =Z_{R}\left(
\left\{  0\right\}  \right)  =Z_{R}\left(  \left\{  1\right\}  \right)  =R$.
(This shows, in particular, that $Z_{R}\left(  S\right)  $ is not always commutative.)

\item[\textbf{(c)}] If $T$ is a subset of $S$, then $Z_{R}\left(  S\right)
\subseteq Z_{R}\left(  T\right)  $.

\item[\textbf{(d)}] If $T$ is a subset of $Z_{R}\left(  S\right)  $, then $S$
is (in turn) a subset of $Z_{R}\left(  T\right)  $.
\end{enumerate}
\end{exercise}

\begin{exercise}
Let $R$ be the matrix ring $\mathbb{R}^{2\times2}$. In this ring $R$, consider
the two matrices
\[
A:=\left(
\begin{array}
[c]{cc}%
0 & 1\\
0 & 0
\end{array}
\right)  \ \ \ \ \ \ \ \ \ \ \text{and}\ \ \ \ \ \ \ \ \ \ B:=\left(
\begin{array}
[c]{cc}%
0 & 0\\
1 & 0
\end{array}
\right)  .
\]

\begin{enumerate}
\item[\textbf{(a)}] Describe the centralizer $Z_{R}\left(  \left\{
A,B\right\}  \right)  $.

\item[\textbf{(b)}] Describe the centralizer $Z_{R}\left(  \left\{
A+B\right\}  \right)  $.

\item[\textbf{(c)}] Describe the centralizer $Z_{R}\left(  \left\{
A-B\right\}  \right)  $.
\end{enumerate}
\end{exercise}

\begin{exercise}
Let $R$ be a ring. Let $S$ be a subset of $R$. Prove the following:

\begin{enumerate}
\item[\textbf{(a)}] We have $S\subseteq Z_{R}\left(  Z_{R}\left(  S\right)
\right)  $.

\item[\textbf{(b)}] If $R=\mathbb{Q}$ and $S=\mathbb{Z}$, then $S$ is a proper
subset of $Z_{R}\left(  Z_{R}\left(  S\right)  \right)  $.

\item[\textbf{(c)}] However, we always have $Z_{R}\left(  S\right)
=Z_{R}\left(  Z_{R}\left(  Z_{R}\left(  S\right)  \right)  \right)  $.
\end{enumerate}
\end{exercise}

\subsection{\label{sec.rings.idomains}Zero divisors and integral domains
(\cite[\S 7.1]{DumFoo04})}

Here comes a rather unsurprising definition:

\begin{definition}
An element of a ring $R$ is said to be \textbf{nonzero} if it is $\neq0$.
(Here, $0$ means $0_{R}$.)
\end{definition}

As we saw above, it can happen that a product of two nonzero elements of a
ring $R$ is zero. Let us give this phenomenon a name (at least in a
commutative setting):

\begin{definition}
Let $R$ be a commutative ring. A nonzero element $a\in R$ is called a
\textbf{zero divisor} if there is a nonzero $b\in R$ such that $ab=0$.
\end{definition}

This definition is slightly controversial: Some people don't require $a$ to be
nonzero. Thus, to them, $0$ is a zero divisor unless $R$ is trivial. (The
notion of a \textquotedblleft zero divisor\textquotedblright\ is somewhat
awkward and will eventually be deprecated in favor of the notion of a
\textquotedblleft regular element\textquotedblright, which is essentially its
negative. See Definition \ref{def.rings.reg.reg}.)

Here are some examples:

\begin{itemize}
\item The elements $\overline{2}$, $\overline{3}$ and $\overline{4}$ of the
ring $\mathbb{Z}/6$ are zero divisors, since they are nonzero but satisfy
$\overline{2}\cdot\overline{3}=\overline{6}=\overline{0}=0_{\mathbb{Z}/6}$ and
$\overline{3}\cdot\overline{2}=\overline{6}=\overline{0}=0_{\mathbb{Z}/6}$ and
$\overline{4}\cdot\overline{3}=\overline{12}=\overline{0}=0_{\mathbb{Z}/6}$.
The element $\overline{0}$ is not a zero divisor (since our definition
requires a zero divisor to be nonzero). The elements $\overline{1}$ and
$\overline{-1}$ are not zero divisors either; indeed, it is easy to see that
for any commutative ring $R$, neither $1_{R}$ nor $-1_{R}$ is a zero divisor.

\item If $a$ is an idempotent element of a commutative ring $R$ (see
Definition \ref{def.ring.idempotent} \textbf{(a)}), but equals neither $0$ nor
$1$, then $a$ is a zero divisor, since $a\left(  1-a\right)
=a-\underbrace{a^{2}}_{=a}=a-a=0$.
\end{itemize}

Zero divisors themselves aren't very useful, but their non-existence (in some
rings) is:

\begin{definition}
\label{def.idomain}Let $R$ be a commutative ring. Assume that $0\neq1$ in $R$.
(By this, we mean $0_{R}\neq1_{R}$; that is, the zero and the unity of $R$ are
distinct. In other words, we assume that the ring $R$ is not trivial.) We say
that $R$ is an \textbf{integral domain} if all nonzero $a,b\in R$ satisfy
$ab\neq0$.
\end{definition}

Equivalently, a commutative ring $R$ with $0\neq1$ (in $R$, that is) is an
integral domain if and only if $R$ has no zero divisors.

Here are some examples:

\begin{itemize}
\item The rings $\mathbb{Z}$, $\mathbb{Q}$, $\mathbb{R}$ and $\mathbb{C}$ are
integral domains.

\item Any subring of an integral domain is clearly an integral domain as well.
Thus, e.g., the ring $\mathbb{S}$ from Subsection \ref{subsec.rings.def.exas}
(i.e., the subring of $\mathbb{R}$ consisting of the numbers of the form
$a+b\sqrt{5}$ with $a,b\in\mathbb{Q}$) is an integral domain.

\item The ring $\mathbb{Z}/n$ is an integral domain if and only if $n$ is $0$
or a prime or minus a prime. We will prove this later.

\item The ring $\mathbb{S}^{\prime}$ from Subsection
\ref{subsec.rings.def.exas} (i.e., the ring whose elements are numbers of the
form $a+b\sqrt{5}$ with $a,b\in\mathbb{Q}$, with multiplication $\ast$ given
by $\left(  a+b\sqrt{5}\right)  \ast\left(  c+d\sqrt{5}\right)  =ac+bd\sqrt
{5}$) is not an integral domain, since it has $1\ast\sqrt{5}=0$.

\item The ring of all functions from $\mathbb{Q}$ to $\mathbb{Q}$ is not an
integral domain, since any two functions with disjoint supports will multiply
to $0$. (For a specific example, we have $\delta_{0}\cdot\delta_{1}=0$, where
$\delta_{y}$ (for $y\in\mathbb{Q}$) is the function that sends $y$ to $1$ and
all other rational numbers to $0$.)

\item We required an integral domain to be commutative in Definition
\ref{def.idomain}. If we dropped this requirement, then the ring $\mathbb{H}$
of quaternions would be an integral domain, but the matrix ring $\mathbb{R}%
^{2\times2}$ would not be.
\end{itemize}

\begin{exercise}
\label{exe.idomain.PM}Consider the ring $\mathbb{R}^{2\times2}$ of all
$2\times2$-matrices with real entries.

Define two subsets $\mathcal{P}$ and $\mathcal{M}$ of $\mathbb{R}^{2\times2}$
by%
\begin{align*}
\mathcal{P}  &  :=\left\{  \left(
\begin{array}
[c]{cc}%
a & b\\
b & a
\end{array}
\right)  \ \mid\ a,b\in\mathbb{R}\right\}  \ \ \ \ \ \ \ \ \ \ \text{and}\\
\mathcal{M}  &  :=\left\{  \left(
\begin{array}
[c]{cc}%
a & b\\
-b & a
\end{array}
\right)  \ \mid\ a,b\in\mathbb{R}\right\}  .
\end{align*}

\begin{enumerate}
\item[\textbf{(a)}] Show that $\mathcal{P}$ and $\mathcal{M}$ are commutative
subrings of $\mathbb{R}^{2\times2}$.

\item[\textbf{(b)}] Prove that $\mathcal{P}$ is not an integral domain.

\item[\textbf{(c)}] Prove that for every nonzero $A \in\mathcal{M}$, there
exists some $B \in\mathcal{M}$ such that $AB = BA = I_{2}$. (In the parlance
we will use later, this property -- along with the nontriviality of
$\mathcal{M}$ -- shall be described as ``$\mathcal{M}$ is a field''.)
\end{enumerate}
\end{exercise}

\begin{exercise}
\ \ 

\begin{enumerate}
\item[\textbf{(a)}] Recall the commutative ring $\mathcal{F}$ from Exercise
\ref{exe.21hw1.6} \textbf{(b)}. Prove that $\mathcal{F}$ is an integral domain.

\item[\textbf{(b)}] Let $\mathcal{F}_{\mathbb{Q}}$ be the ring defined just
like $\mathcal{F}$, but using $\mathbb{Q}$ instead of $\mathbb{Z}$ (that is,
it is the subring $\left\{  aA+bI_{2}\mid a,b\in\mathbb{Q}\right\}  $ of
$\mathbb{Q}^{2\times2}$). Is $\mathcal{F}_{\mathbb{Q}}$ an integral domain?

\item[\textbf{(c)}] Let $\mathcal{F}_{\mathbb{R}}$ be the ring defined just
like $\mathcal{F}$, but using $\mathbb{R}$ instead of $\mathbb{Z}$ (that is,
it is the subring $\left\{  aA+bI_{2}\mid a,b\in\mathbb{R}\right\}  $ of
$\mathbb{R}^{2\times2}$). Is $\mathcal{F}_{\mathbb{R}}$ an integral domain?
\end{enumerate}

[\textbf{Hint:} For part \textbf{(a)}, argue that the determinant%
\[
\det\left(  aA+bI_{2}\right)  =\det%
\begin{pmatrix}
b & a\\
a & a+b
\end{pmatrix}
=-a^{2}+ab+b^{2}=\dfrac{5b^{2}-\left(  2a-b\right)  ^{2}}{4}%
\]
of any nonzero matrix $aA+bI_{2}\in\mathcal{F}\setminus\left\{  0\right\}  $
is nonzero, since $\sqrt{5}$ is irrational. Now recall that $\det\left(
AB\right)  =\det A\cdot\det B$ for any $A,B\in\mathbb{R}^{2\times2}$.]
\end{exercise}

\begin{warning}
\label{warn.nzd.subring}Let $R$ be a commutative ring, and let $S$ be a
subring of $R$. It can happen that some nonzero element $a\in S$ is a zero
divisor in $R$ (that is, there is a nonzero $b\in R$ such that $ab=0$) but not
a zero divisor in $S$ (that is, there exists no nonzero $b\in S$ such that
$ab=0$). This should not be too surprising ($R$ has more elements than $S$, so
it should be \textquotedblleft easier\textquotedblright\ to find the required
$b$ in $R$ than in $S$), although an explicit example is not easy to construct
at this point. (Using the concept of quotient rings we will learn later, we
can take $R=\mathbb{Z}\left[  x\right]  /\left(  2x\right)  $ and
$S=\mathbb{Z}$ and $a=2$, where we view $S$ as a subring of $R$ in the
\textquotedblleft obvious\textquotedblright\ way by identifying each integer
$n\in\mathbb{Z}$ with the corresponding residue class $\overline{n}\in R$.)
\end{warning}

\subsection{\label{sec.rings.units}Units and fields (\cite[\S 7.1]{DumFoo04})}

\subsubsection{Units and inverses}

By definition, any ring $R$ has an addition, a subtraction and a
multiplication. Division, on the other hand, is not guaranteed: Even the ring
$\mathbb{Z}$ doesn't really have division (unless you count division with
remainder, which is a different story). However, any ring $R$ has
\textbf{some} elements that can be divided by; the simplest such element is
its unity $1$. Let us introduce a name for these elements:

\begin{definition}
\label{def.rings.units.unit}Let $R$ be a ring.

\begin{enumerate}
\item[\textbf{(a)}] An element $a\in R$ is said to be a \textbf{unit} of $R$
(or \textbf{invertible} in $R$) if there exists a $b\in R$ such that
$ab=ba=1$. In this case, $b$ is unique and is known as the \textbf{inverse}
(or \textbf{multiplicative inverse}, or \textbf{reciprocal}) of $a$, and is
denoted by $a^{-1}$.

\item[\textbf{(b)}] We let $R^{\times}$ denote the set of all units of $R$.
\end{enumerate}
\end{definition}

A few comments:

\begin{itemize}
\item It goes without saying that the \textquotedblleft$1$\textquotedblright%
\ refers to the unity of the ring $R$.

\item We required $ab=ba=1$ rather than merely $ab=1$ because $R$ is not
necessarily commutative. When $R$ is commutative, of course, $ab=1$ suffices.

\item Why is $b$ unique in Definition \ref{def.rings.units.unit} \textbf{(a)}?
Because if $b_{1}$ and $b_{2}$ are two such $b$'s (for the same $a$), then
$ab_{1}=b_{1}a=1$ and $ab_{2}=b_{2}a=1$, so that $b_{1}\underbrace{ab_{2}%
}_{=1}=b_{1}1=b_{1}$ and thus $b_{1}=\underbrace{b_{1}a}_{=1}b_{2}%
=1b_{2}=b_{2}$. This is the exact same argument that proves the uniqueness of
inverses in a group.

\item Don't confuse \textquotedblleft unit\textquotedblright\ (= invertible
element) with \textquotedblleft unity\textquotedblright\ (= neutral element
for multiplication). The unity is always a unit, but not vice versa!

\item Some people write $R^{\ast}$ or $R^{x}$ for $R^{\times}$.
\end{itemize}

Here are some examples of units:

\begin{itemize}
\item The units of the ring $\mathbb{Q}$ are all nonzero elements of
$\mathbb{Q}$. (This is because every nonzero element of $\mathbb{Q}$ has a
reciprocal, and this reciprocal again lies in $\mathbb{Q}$.) The same holds
for $\mathbb{R}$ and for $\mathbb{C}$.

\item The units of the ring $\mathbb{Z}$ are $1$ and $-1$ (with inverses $1$
and $-1$, respectively). No other integer is a unit of $\mathbb{Z}$. For
example, $2$ has an inverse $\dfrac{1}{2}$ in $\mathbb{Q}$, but not in
$\mathbb{Z}$.

\item The units of the matrix ring $\mathbb{R}^{n\times n}$ are the invertible
$n\times n$-matrices. You have seen many ways to characterize them in your
linear algebra class. You might even remember that the set $\left(
\mathbb{R}^{n\times n}\right)  ^{\times}$ of these units is known as the
$n$\textbf{-th general linear group} of $\mathbb{R}$, and is called
$\operatorname*{GL}\nolimits_{n}\left(  \mathbb{R}\right)  $ or
$\operatorname*{GL}\left(  n,\mathbb{R}\right)  $.

\item In the ring of all functions from $\mathbb{Q}$ to $\mathbb{Q}$, the
units are the functions that never vanish (i.e., that don't take $0$ as a
value). Inverses can be computed pointwise.

\item Recall the ring $\mathbb{Z}\left[  i\right]  $ of Gaussian integers. Its
only units are $1,i,-1,-i$. This is Corollary \ref{cor.Zi.units} further below.

\item A truly trivial example: The zero of a ring $R$ is a unit if and only if
$R$ is a trivial ring.
\end{itemize}

Our next example we state as a proposition:\footnote{Two integers $a$ and $b$
are said to be \textbf{coprime} (to each other) if and only if $\gcd\left(
a,b\right)  =1$. Some authors say \textquotedblleft relatively
prime\textquotedblright\ instead of \textquotedblleft
coprime\textquotedblright.}

\begin{proposition}
\label{prop.ringunits.Z/n}Let $n\in\mathbb{Z}$.

\begin{enumerate}
\item[\textbf{(a)}] The units of the ring $\mathbb{Z}/n$ are precisely the
residue classes $\overline{a}\in\mathbb{Z}/n$ where $a\in\mathbb{Z}$ is
coprime to $n$.

\item[\textbf{(b)}] Let $a\in\mathbb{Z}$. Then, $\overline{a}\in\mathbb{Z}/n$
is a unit of $\mathbb{Z}/n$ if and only if $a$ is coprime to $n$.
\end{enumerate}
\end{proposition}

\begin{proof}
We begin by proving part \textbf{(b)}, which is the stronger claim. (Part
\textbf{(a)} will then easily follow.) \medskip

\textbf{(b)} This is an \textquotedblleft if and only if\textquotedblright%
\ statement. We shall prove its \textquotedblleft if\textquotedblright\ (i.e.,
\textquotedblleft$\Longleftarrow$\textquotedblright) and \textquotedblleft
only if\textquotedblright\ (i.e., \textquotedblleft$\Longrightarrow
$\textquotedblright) parts separately:

$\Longleftarrow:$ Assume that $a\in\mathbb{Z}$ is coprime to $n$. Bezout's
theorem\footnote{\textbf{Bezout's theorem} (from elementary number theory)
states that for any two integers $a$ and $b$, there exist two integers $x$ and
$y$ satisfying $xa+yb=\gcd\left(  a,b\right)  $. In other words, the greatest
common divisor of two integers $a$ and $b$ can always be written as a linear
combination of $a$ and $b$ with integer coefficients.
\par
Proofs of Bezout's theorem can be found in most texts on elementary number
theory (e.g., \cite[Theorem 2.3]{Burton11}) and in several introductory texts
on abstract algebra (e.g., \cite[Theorem 2.9.12]{19s}). We shall also prove a
general result about Euclidean rings (Theorem \ref{thm.rings.euclid.ea.main}
\textbf{(a)}) that easily yields Bezout's theorem (see Collary
\ref{cor.gcd-lcm.eucl-bez} and Proposition \ref{prop.gcd-lcm.ZZ}).} tells us
that there exist $x,y\in\mathbb{Z}$ with $xa+yn=\gcd\left(  a,n\right)  $.
Consider these $x,y$. We have $xa+yn=\gcd\left(  a,n\right)  =1$ (since $a$ is
coprime to $n$). Thus, $xa\equiv xa+yn=1\operatorname{mod}n$. Translating this
into the language of residue classes, we obtain $\overline{xa}=\overline{1}$
in $\mathbb{Z}/n$. Hence, $\overline{x}\cdot\overline{a}=\overline
{xa}=\overline{1}$ in $\mathbb{Z}/n$. Since the ring $\mathbb{Z}/n$ is
commutative, this shows that $\overline{a}$ is invertible (with inverse
$\overline{x}$). In other words, $\overline{a}$ is a unit of $\mathbb{Z}/n$.

$\Longrightarrow:$ Conversely, assume that $\overline{a}$ is a unit of
$\mathbb{Z}/n$. Thus, $\overline{a}$ has an inverse $\overline{b}\in
\mathbb{Z}/n$. This inverse $\overline{b}$ satisfies $\overline{ab}%
=\overline{1}$; in other words, $ab\equiv1\operatorname{mod}n$. But this
easily yields that\footnote{We are using the fact that if $u$ and $v$ are two
integers satisfying $u\equiv v\operatorname{mod}n$, then $\gcd\left(
u,n\right)  =\gcd\left(  v,n\right)  $. This is just a restatement of the
classical result that the gcd of two integers does not change if we add a
multiple of one to the other.} $\gcd\left(  ab,n\right)  =\gcd\left(
1,n\right)  =1$. In other words, $ab$ is coprime to $n$. Hence, $a$ is coprime
to $n$ as well (since any common divisor of $a$ and $n$ must be a common
divisor of $ab$ and $n$). \medskip

\textbf{(a)} This follows easily from part \textbf{(b)}.
\end{proof}

Here are some examples of Proposition \ref{prop.ringunits.Z/n}:

\begin{itemize}
\item The units of the ring $\mathbb{Z}/12$ are $\overline{1},\overline
{5},\overline{7},\overline{11}$ (because among the integers $0,1,\ldots,11$,
it is the four numbers $1,5,7,11$ that are coprime to $12$).

\item The units of the ring $\mathbb{Z}/5$ are $\overline{1},\overline
{2},\overline{3},\overline{4}$.

\item The only unit of the ring $\mathbb{Z}/2$ is $\overline{1}$.
\end{itemize}

Next, we shall show some general properties of units in rings:

\begin{theorem}
\label{thm.ringunits.group}Let $R$ be a ring. Then, the set $R^{\times}$ is a
multiplicative group. More precisely: $\left(  R^{\times},\cdot,1\right)  $ is
a group.
\end{theorem}

\begin{proof}
It suffices to show the following facts:

\begin{enumerate}
\item The unity $1$ of $R$ belongs to $R^{\times}$.

\item If $a,b\in R^{\times}$, then $ab\in R^{\times}$.

\item If $a\in R^{\times}$, then $a$ has an inverse in $R^{\times}$.
\end{enumerate}

All other group axioms for $R^{\times}$ follow from the ring axioms of $R$. So
let us prove these three facts.

\textit{Proof of Fact 1:} Fact 1 is obvious (as $1$ has inverse $1$).

\textit{Proof of Fact 2:} Let $a,b\in R^{\times}$. Thus, the elements $a,b$
are units, and thus have inverses $a^{-1},b^{-1}$, respectively. These satisfy
$aa^{-1}=a^{-1}a=1$ and $bb^{-1}=b^{-1}b=1$. Now, $a\underbrace{bb^{-1}}%
_{=1}a^{-1}=aa^{-1}=1$ and $b^{-1}\underbrace{a^{-1}a}_{=1}b=b^{-1}b=1$, so
that $ab$ is invertible as well (with inverse $b^{-1}a^{-1}$). That is, $ab\in
R^{\times}$. This proves Fact 2.

\textit{Proof of Fact 3:} Let $a\in R^{\times}$. Thus, $a$ has an inverse
$a^{-1}$ in $R$. This inverse $a^{-1}$, in turn, has an inverse (namely, $a$),
and thus also lies in $R^{\times}$. Hence, $a$ has an inverse in $R^{\times}$.
This proves Fact 3.
\end{proof}

The group $R^{\times}$ from Theorem \ref{thm.ringunits.group} is known as the
\textbf{group of units} of $R$. Thus, every ring $R$ produces \textbf{two}
groups: the additive group $\left(  R,+,0\right)  $ (which contains all
elements of $R$) and the multiplicative group of units $\left(  R^{\times
},\cdot,1\right)  $ (which only contains the units). We record two important
consequences of our proof of Theorem \ref{thm.ringunits.group}:

\begin{theorem}
[Shoe-sock theorem]\label{thm.ringunits.shoesock}Let $R$ be a ring. Let $a,b$
be two units of $R$. Then, $ab$ is a unit of $R$, and its inverse is $\left(
ab\right)  ^{-1}=b^{-1}a^{-1}$.
\end{theorem}

\begin{proof}
See the proof of Fact 2 in the proof of Theorem \ref{thm.ringunits.group}.
\end{proof}

\begin{theorem}
Let $R$ be a ring. Let $a$ be a unit of $R$. Then, $a^{-1}$ is a unit of $R$,
and its inverse is $\left(  a^{-1}\right)  ^{-1}=a$.
\end{theorem}

\begin{proof}
See the proof of Fact 3 in the proof of Theorem \ref{thm.ringunits.group}.
\end{proof}

Clearly, if $R$ is a commutative ring, then its group of units $R^{\times}$ is
abelian. Sometimes, more can be said about it. In particular, we will
eventually see (Theorem \ref{thm.finfield.gauss1}) that the group $\left(
\mathbb{Z}/n\right)  ^{\times}$ is cyclic whenever $n$ is prime. For other
$n$, it sometimes is and sometimes isn't (see Exercise \ref{exe.primroot.n}).
Here is a simple example:

\begin{example}
As we have seen, the ring $\mathbb{Z}/12$ has $4$ units: $\overline
{1},\overline{5},\overline{7},\overline{11}$. Thus, the group of units
$\left(  \mathbb{Z}/12\right)  ^{\times}$ is $\left\{  \overline{1}%
,\overline{5},\overline{7},\overline{11}\right\}  $. It is easy to see that
each of these units squares to $\overline{1}$ (since $\overline{1}%
^{2}=\overline{5}^{2}=\overline{7}^{2}=\overline{11}^{2}=\overline{1}$);
hence, this group is isomorphic to
\href{https://en.wikipedia.org/wiki/Klein_four-group}{the Klein four-group} (a
direct product of two cyclic groups of order $2$).
\end{example}

Here are a few exercises on units and inverses in some special rings. The
first exercise (\cite[homework set \#2, Exercise 5 \textbf{(a)}]{21w}) will
come in useful later:

\begin{exercise}
\label{exe.21hw2.5a}Let $p$ be a prime. Prove that the only units of the ring
$\mathbb{Z}/p$ that are their own inverses (i.e., the only $m\in\left(
\mathbb{Z}/p\right)  ^{\times}$ that satisfy $m^{-1}=m$) are $\overline{1}$
and $\overline{-1}$.
\end{exercise}

\begin{exercise}
\label{exe.units.Z/pk}Let $p$ be a prime. Let $k$ be a positive integer. Prove
that the number of units of the ring $\mathbb{Z}/p^{k}$ is $p^{k}-p^{k-1}$.
\end{exercise}

\begin{exercise}
Let $R$ be the ring $\left(  \mathbb{Z}/2\right)  ^{2\times2}$ of all
$2\times2$-matrices with entries in $\mathbb{Z}/2$. This ring has size
$2^{4}=16$, since each such matrix has $4$ entries and there are $2$ options
for each entry.

\begin{enumerate}
\item[\textbf{(a)}] Find the group of units $R^{\times}$ of this ring.

\item[\textbf{(b)}] Prove that this group $R^{\times}$ is isomorphic to the
symmetric group $S_{3}$ (that is, the group of all permutations of the set
$\left\{  1,2,3\right\}  $).
\end{enumerate}
\end{exercise}

\begin{exercise}
\label{exe.dualnums.inverse}Let $\mathbb{D}$ be the ring of dual numbers, as
defined in Exercise \ref{exe.dualnums.ring}. Prove the following:

\begin{enumerate}
\item[\textbf{(a)}] A dual number $a+b\varepsilon$ (with $a,b\in\mathbb{R}$)
is a unit of $\mathbb{D}$ if and only if $a\neq0$.

\item[\textbf{(b)}] If $a,b\in\mathbb{R}$ satisfy $a\neq0$, then the inverse
of the dual number $a+b\varepsilon$ is $\dfrac{1}{a}-\dfrac{b}{a^{2}%
}\varepsilon$.
\end{enumerate}
\end{exercise}

\begin{exercise}
Recall the Fibonacci sequence $\left(  f_{0},f_{1},f_{2},\ldots\right)  $ from
Definition \ref{def.fibonacci.fib}, and recall the matrix $A$ and the
commutative ring $\mathcal{F}$ from Exercise \ref{exe.21hw1.6} \textbf{(b)}.

We extend the Fibonacci sequence $\left(  f_{0},f_{1},f_{2},\ldots\right)  $
to an infinite-in-both-directions \textquotedblleft sequence\textquotedblright%
\ $\left(  \ldots,f_{-2},f_{-1},f_{0},f_{1},f_{2},\ldots\right)  $ by
requiring that it satisfy the original recursive equation $f_{n}%
=f_{n-1}+f_{n-2}$ for all $n\in\mathbb{Z}$ (not just for $n\geq2$). Thus, the
negatively indexed Fibonacci numbers $f_{-1},f_{-2},f_{-3},\ldots$ are
computed recursively by solving this recursive equation $f_{n}=f_{n-1}%
+f_{n-2}$ for $f_{n-2}$. For instance, $f_{-1}=f_{1}-f_{0}=1-0=1$ and
$f_{-2}=f_{0}-f_{-1}=0-1=-1$.

\begin{enumerate}
\item[\textbf{(a)}] Prove that the matrix $A$ is a unit of $\mathcal{F}$ (that
is, it has an inverse in $\mathcal{F}$).

\item[\textbf{(b)}] Prove that $A^{n}=f_{n}A+f_{n-1}I_{2}$ for all
$n\in\mathbb{Z}$.

\item[\textbf{(c)}] Prove that $f_{-n}=\left(  -1\right)  ^{n-1}f_{n}$ for
each $n\in\mathbb{Z}$.

\item[\textbf{(d)}] Prove that the units of $\mathcal{F}$ are precisely the
powers $A^{k}$ of the matrix $A$ (with $k\in\mathbb{Z}$) as well as their
negatives (i.e., the matrices $-A^{k}$ with $k \in\mathbb{Z}$). (Note that
``powers''\ includes the positive powers $A^{1},A^{2},A^{3},\ldots$, the
negative powers $A^{-1},A^{-2},A^{-3},\ldots$ and the zeroth power
$A^{0}=I_{2}$.)
\end{enumerate}

[\textbf{Hint:} Part \textbf{(d)} is surprisingly tricky! Recall again that
$\det\left(  aA+bI_{2}\right)  =-a^{2}+ab+b^{2}$ for any $a,b\in\mathbb{Z}$.
Show that this determinant $\det\left(  aA+bI_{2}\right)  $ has to be $1$ or
$-1$ if $aA+bI_{2}$ is a unit of $\mathcal{F}$. Thus, we must have
$-a^{2}+ab+b^{2}\in\left\{  1,-1\right\}  $ if $aA+bI_{2}$ is a unit. But this
means that the pair $\left(  a,b\right)  $ is a \textquotedblleft golden
pair\textquotedblright\ in the terminology of \cite[Exercise 5.4.10]{mps}, and
the set of all \textquotedblleft golden pairs\textquotedblright\ can be
described explicitly in terms of the Fibonacci sequence \cite[Exercise
5.4.10]{mps}.]
\end{exercise}

\begin{fineprint}

\begin{remark}
Let $R$ be a ring, and let $S$ be a subring of $R$. Then, any unit $u$ of $S$
is also a unit of $R$ (since its inverse belongs to $S$ and therefore to $R$,
and thus $u$ has an inverse in $R$). This is in stark contrast to the
situation for non-zero-divisors (which we discussed in Warning
\ref{warn.nzd.subring}).
\end{remark}
\end{fineprint}

\subsubsection{Some exercises on inverses}

The following exercises prove surprising results and make for good practice
with the definition of an inverse\footnote{Keep in mind that a ring is not
always commutative!}:

\begin{exercise}
\label{exe.21hw1.9}Let $R$ be a ring. Let $a$ and $b$ be two elements of $R$.

Prove that if $1-ab$ is invertible, then so is $1-ba$.

Better yet, prove the following: If $c$ is an inverse of $1-ab$, then $1+bca$
is an inverse of $1-ba$.
\end{exercise}

Note that Exercise \ref{exe.21hw1.9} yields a well-known result in functional
analysis (see \url{https://math.stackexchange.com/questions/79217} ).

\begin{exercise}
Let $R$ be a ring. Let $a$ and $b$ be two units of $R$ such that $a+b$ is a
unit as well.

\begin{enumerate}
\item[\textbf{(a)}] Prove that $a^{-1}+b^{-1}$, too, is a unit, and its
inverse is%
\[
\left(  a^{-1}+b^{-1}\right)  ^{-1}=a\cdot\left(  a+b\right)  ^{-1}\cdot
b=b\cdot\left(  a+b\right)  ^{-1}\cdot a.
\]

\item[\textbf{(b)}] Show on an example that $\left(  a^{-1}+b^{-1}\right)
^{-1}$ can be different from $ab\cdot\left(  a+b\right)  ^{-1}$.
\end{enumerate}
\end{exercise}

Let us next define some weaker variants of inverses:

\begin{definition}
\label{def.left-right-inverses}Let $R$ be a ring. Let $a$ be an element of $R$.

\begin{enumerate}
\item[\textbf{(a)}] A \textbf{left inverse} of $a$ shall mean an element $b\in
R$ satisfying $ba=1$.

\item[\textbf{(b)}] A \textbf{right inverse} of $a$ shall mean an element
$b\in R$ satisfying $ab=1$.
\end{enumerate}
\end{definition}

Thus, an inverse of $a$ is the same as a left inverse of $a$ that
simultaneously is a right inverse of $a$. It is clear that the notions of
\textquotedblleft left inverse\textquotedblright\ and \textquotedblleft right
inverse\textquotedblright\ can be defined in any monoid, not just in a ring,
since they rely only on the multiplication and the unity. As already
mentioned, a left or right inverse doesn't have to be a (proper) inverse in
general, although it is hard to find examples where it isn't. The following
exercise (a result of Jacobson) might give a hint as to why:

\begin{exercise}
Let $R$ be a ring. Let $a$ and $b$ be two elements of $R$ such that $ab=1$ but
$ba\neq1$. Let $w=1-ba$.

\begin{enumerate}
\item[\textbf{(a)}] Prove that $aw=wb=0$.

\item[\textbf{(b)}] Conclude that $a\left(  b+wa^{k}\right)  =1$ for all
$k\in\mathbb{N}$.

\item[\textbf{(c)}] Prove that the elements $b+wa^{k}$ for all $k\in
\mathbb{N}$ are distinct.

\item[\textbf{(d)}] Conclude that $a$ has infinitely many right inverses.

\item[\textbf{(e)}] Conclude that $R$ cannot be finite.
\end{enumerate}

[\textbf{Hint:} The only hard part here is \textbf{(c)}. Show first that
$wa^{i}\neq0$ for all $i\in\mathbb{N}$; then show that $wa^{i}\neq w$ for all
positive integers $i$.]
\end{exercise}

The next exercise (\cite[homework set \#2, Exercise 1]{21w}) provides another
source of units in certain rings:

\begin{exercise}
\label{exe.21hw2.1}Let $R$ be a ring. An element $a\in R$ will be called
\textbf{nilpotent} if there exists some $n\in\mathbb{N}$ such that $a^{n}=0$.
(For instance, the element $\overline{18}\in\mathbb{Z}/24$ is nilpotent, since
$\overline{18}^{3}=\overline{0}$. Note that the zero $0$ is nilpotent in any
ring, but other nilpotent elements may or may not exist. For another example,
the element $\varepsilon\in\mathbb{D}$ in Exercise \ref{exe.dualnums.ring} is nilpotent.)

Let $a\in R$ be a nilpotent element.

\begin{enumerate}
\item[\textbf{(a)}] Prove that $1-a\in R$ is a unit.

\item[\textbf{(b)}] Let $u\in R$ be a unit satisfying $ua=au$. Prove that
$u-a\in R$ is a unit.
\end{enumerate}

[\textbf{Hint:} Treat the geometric series $\dfrac{1}{1-x}=1+x+x^{2}+\cdots$
as an inspiration, noting that the infinite sum on the right hand side will
become a finite sum if the nilpotent element $a$ is substituted for $x$.]
\end{exercise}

\subsubsection{Fields}

As we saw, some rings (such as $\mathbb{Z}$) have few units, while other rings
(such as $\mathbb{Q}$) have many. The rings with the most units are the
\textquotedblleft fields\textquotedblright:

\begin{definition}
\label{def.fields.field}Let $R$ be a commutative ring. Assume that $0\neq1$ in
$R$. We say that $R$ is a \textbf{field} if every nonzero element of $R$ is a unit.
\end{definition}

Examples:

\begin{itemize}
\item The rings $\mathbb{Q}$, $\mathbb{R}$ and $\mathbb{C}$ are fields. The
ring $\mathbb{Z}$ is not (since $2$ is not a unit).

\item The ring $\mathbb{S}$ of all real numbers of the form $a+b\sqrt{5}$ with
$a,b\in\mathbb{Q}$ (as defined in Subsection \ref{subsec.rings.def.exas}) is a
field. Indeed, the inverse of a nonzero element $a+b\sqrt{5}$ is%
\[
\left(  a+b\sqrt{5}\right)  ^{-1}=\dfrac{1}{a+b\sqrt{5}}=\dfrac{a-b\sqrt{5}%
}{a^{2}-b^{2}\cdot5}=\dfrac{a}{a^{2}-5b^{2}}+\dfrac{-b}{a^{2}-5b^{2}}\sqrt{5}%
\]
(the denominators here are nonzero because $a+b\sqrt{5}\neq0$ entails
$a^{2}-5b^{2}\neq0$). So this is why they taught you rationalizing
denominators in high school!

\item The Hamiltonian quaternions $\mathbb{H}$ are not a field, but for a
stupid reason: they are noncommutative. Otherwise, they would be a field. A
noncommutative ring in which each nonzero element is invertible is called a
\textbf{division ring} or \textbf{skew-field}.

\item Let $n$ be a positive integer. The ring $\mathbb{Z}/n$ is a field if and
only if $n$ is prime. (We will prove this below.)

\item The ring $F_{4}$ constructed in Subsection \ref{subsec.rings.def.exas}
as well as the ring $F_{8}$ defined in Exercise \ref{exe.rings.F8} are fields.
\end{itemize}

\subsection{\label{sec.rings.fields2}Fields and integral domains: some
connections (\cite[\S 7.1]{DumFoo04})}

\subsubsection{Fields vs. integral domains}

The notions of fields and integral domains are closely related:

\begin{proposition}
\label{prop.fields.intdom}\ \ 

\begin{enumerate}
\item[\textbf{(a)}] Every field is an integral domain.

\item[\textbf{(b)}] Every \textbf{finite} integral domain is a field. (Here,
of course, \textquotedblleft finite\textquotedblright\ means \textquotedblleft
finite as a set\textquotedblright.)
\end{enumerate}
\end{proposition}

\begin{proof}
\textbf{(a)} Let $F$ be a field. We must show that $F$ is an integral domain.

Let $a,b\in F$ be nonzero. We must show that $ab$ is nonzero.

Indeed, $a$ and $b$ are nonzero, and thus are units (since $F$ is a field).
Thus, they have inverses $a^{-1}$ and $b^{-1}$.

Now, if we had $ab=0$, then we would have $\underbrace{ab}_{=0}b^{-1}a^{-1}%
=0$, which would yield $0=a\underbrace{bb^{-1}}_{=1}a^{-1}=aa^{-1}=1$, which
would contradict the fact that $0\neq1$ in $F$ (since $F$ is a field). Thus,
we cannot have $ab=0$. In other words, $ab$ is nonzero. This proves that $F$
is an integral domain. Thus, Proposition \ref{prop.fields.intdom} \textbf{(a)}
is proved. \medskip

\textbf{(b)} Let $R$ be a \textbf{finite} integral domain. We must show that
$R$ is a field.

Let $a\in R$ be nonzero. We must show that $a$ is a unit.

Since $R$ is an integral domain, we know that $ab\neq0$ for any $b\neq0$.
Thus, $ax\neq ay$ for any two distinct elements $x$ and $y$ of $R$ (because if
$x$ and $y$ are two distinct elements of $R$, then $x-y\neq0$, and thus the
previous sentence yields $a\left(  x-y\right)  \neq0$; but this rewrites as
$ax-ay\neq0$, so that $ax\neq ay$). In other words, the map%
\[
R\rightarrow R,\ \ \ \ \ \ \ \ \ \ x\mapsto ax
\]
is injective. Hence, this map is also bijective (since any injective map
between two \textbf{finite} sets of the \textbf{same size} is bijective --
this is one of the Pigeonhole Principles\footnote{To be specific, this is what
I call the \textquotedblleft Pigeonhole Principle for
Injections\textquotedblright. See \cite[Theorem 6.1.3]{mps}, for example.}).
Thus, in particular, this map is surjective, and hence takes $1$ as a value.
In other words, there exists an $x\in R$ such that $ax=1$. Since $R$ is
commutative, this $x$ must be an inverse of $a$, and thus we conclude that $a$
is a unit. This finishes the proof of Proposition \ref{prop.fields.intdom}
\textbf{(b)}.
\end{proof}

Without the word \textquotedblleft finite\textquotedblright, Proposition
\ref{prop.fields.intdom} \textbf{(b)} would not be true; for instance,
$\mathbb{Z}$ is an integral domain but no field. The polynomial ring
$\mathbb{R}\left[  x\right]  $ (consisting of univariate polynomials with real
coefficients) is another example of an integral domain that is not a field.
(We will prove this later.)

\subsubsection{When is $\mathbb{Z}/n$ a field?}

Our above study of units of $\mathbb{Z}/n$ lets us now easily obtain the following:

\begin{corollary}
\label{cor.field.Z/n}Let $n$ be a positive integer. Then, the following chain
of equivalences holds:%
\[
\left(  \mathbb{Z}/n\text{ is an integral domain}\right)  \Longleftrightarrow
\left(  \mathbb{Z}/n\text{ is a field}\right)  \Longleftrightarrow\left(
n\text{ is prime}\right)  .
\]

\end{corollary}

\begin{proof}
The first of the two $\Longleftrightarrow$ signs follows from Proposition
\ref{prop.fields.intdom} (since $\mathbb{Z}/n$ is finite). Let's now prove the second.

$\Longrightarrow:$ Assume that $\mathbb{Z}/n$ is a field. Then, every nonzero
element of $\mathbb{Z}/n$ is a unit. Hence, the $n-1$ residue classes
$\overline{1},\overline{2},\ldots,\overline{n-1}$ are units of $\mathbb{Z}/n$
(since they are nonzero). Therefore, the $n-1$ integers $1,2,\ldots,n-1$ are
coprime to $n$ (by Proposition \ref{prop.ringunits.Z/n} \textbf{(b)}). Hence,
$n$ is either $1$ or prime. However, if $n$ was $1$, then we would have
$\overline{0}=\overline{1}$, which would mean that $0=1$ in $\mathbb{Z}/n$;
but this is forbidden for a field. Thus, $n$ cannot be $1$, and therefore must
be prime.

$\Longleftarrow:$ Assume that $n$ is prime. Then, $n>1$, so that $\overline
{0}\neq\overline{1}$. That is, $0\neq1$ in $\mathbb{Z}/n$. Furthermore, if
$\overline{a}$ (for some integer $a$) is a nonzero element of $\mathbb{Z}/n$,
then the integer $a$ is not divisible by $n$ (since $\overline{a}$ is
nonzero), so that $a$ is coprime to $n$ (since $n$ is prime), and this entails
(by Proposition \ref{prop.ringunits.Z/n} \textbf{(b)}) that $\overline{a}$ is
a unit of $\mathbb{Z}/n$. So we have shown that every nonzero element of
$\mathbb{Z}/n$ is a unit. In other words, $\mathbb{Z}/n$ is a field.
\end{proof}

Note that the positivity of $n$ in Corollary \ref{cor.field.Z/n} is important:
The ring $\mathbb{Z}/0$ is an integral domain but not a field. (In fact, this
ring is essentially $\mathbb{Z}$, except that its elements are the singleton
sets $\left\{  a\right\}  $ instead of the integers $a$ themselves.)

\subsubsection{Application: Fermat's Little Theorem}

We can use Corollary \ref{cor.field.Z/n} to obtain an important result in
elementary number theory:

\begin{theorem}
[Fermat's little theorem, short F$\ell$T]\label{thm.ent.flt}Let $p$ be a prime
number. Let $a\in\mathbb{Z}$. Then, $a^{p}\equiv a\operatorname{mod}p$.
\end{theorem}

For example, $a^{3}\equiv a\operatorname{mod}3$ and $a^{5}\equiv
a\operatorname{mod}5$ for every $a\in\mathbb{Z}$.

Before we prove Theorem \ref{thm.ent.flt}, let us first show the following
property of the field $\mathbb{Z}/p$:

\begin{proposition}
[Fermat's little theorem in $\mathbb{Z}/p$ form]\label{prop.ent.flt.Z/p}Let
$p$ be a prime number. Let $u\in\mathbb{Z}/p$. Then, $u^{p}=u$.
\end{proposition}

\begin{proof}
[Proof of Proposition \ref{prop.ent.flt.Z/p}.]We know that $p$ is prime. Thus,
Corollary \ref{cor.field.Z/n} (applied to $n=p$) yields that $\mathbb{Z}/p$ is
a field. Hence, every nonzero element of $\mathbb{Z}/p$ is a unit.

We must prove that $u^{p}=u$. If $u=0$, then this is obvious (since
$u^{p}=0^{p}=0=u$ in this case). So let us WLOG assume that $u\neq0$. Hence,
the element $u\in\mathbb{Z}/p$ is nonzero. Therefore, $u$ is a unit of the
ring $\mathbb{Z}/p$ (since every nonzero element of $\mathbb{Z}/p$ is a unit).
In other words, $u\in\left(  \mathbb{Z}/p\right)  ^{\times}$.

However, the units of the ring $\mathbb{Z}/p$ are $\overline{1},\overline
{2},\ldots,\overline{p-1}$ (again because every nonzero element of
$\mathbb{Z}/p$ is a unit). Thus, in particular, there are $p-1$ of them. This
shows that the group $\left(  \mathbb{Z}/p\right)  ^{\times}$ has order $p-1$.
Hence, Lagrange's theorem (from group theory)\footnote{Recall that this
theorem says the following: If $G$ is a finite group of order $m$ (for some
$m\in\mathbb{N}$), then $g^{m}=1$ for each $g\in G$ (where we are writing $G$
multiplicatively, so that $1$ denotes the neutral element of $G$).} shows that
$g^{p-1}=1$ for each $g\in\left(  \mathbb{Z}/p\right)  ^{\times}$. Applying
this to $g=u$, we obtain $u^{p-1}=1$. Hence, $u^{p}=u\underbrace{u^{p-1}}%
_{=1}=u$. This proves Proposition \ref{prop.ent.flt.Z/p}.
\end{proof}

We can now easily derive Theorem \ref{thm.ent.flt} from Proposition
\ref{prop.ent.flt.Z/p}:

\begin{proof}
[Proof of Theorem \ref{thm.ent.flt}.]Consider the residue class $\overline
{a}\in\mathbb{Z}/p$. Applying Proposition \ref{prop.ent.flt.Z/p} to
$u=\overline{a}$, we obtain $\overline{a}^{p}=\overline{a}$. Thus,
$\overline{a^{p}}=\overline{a}^{p}=\overline{a}$. In other words, $a^{p}\equiv
a\operatorname{mod}p$. Theorem \ref{thm.ent.flt} is thus proven.
\end{proof}

We also observe:

\begin{corollary}
[Fermat's little theorem in the non-divisible case]\label{cor.ent.flt.p-1}Let
$p$ be a prime number. Let $a\in\mathbb{Z}$ satisfy $p\nmid a$. Then,
$a^{p-1}\equiv1\operatorname{mod}p$.
\end{corollary}

\begin{proof}
Consider the residue class $\overline{a}\in\mathbb{Z}/p$. Applying Proposition
\ref{prop.ent.flt.Z/p} to $u=\overline{a}$, we obtain $\overline{a}%
^{p}=\overline{a}$. However, from $p\nmid a$, we obtain $\overline{a}%
\neq\overline{0}$. In other words, $\overline{a}$ is nonzero. Since
$\mathbb{Z}/p$ is a field (by Corollary \ref{cor.field.Z/n}), we know that
every nonzero element of $\mathbb{Z}/p$ is a unit. Thus, $\overline{a}$ is a
unit (since $\overline{a}$ is nonzero). Hence, we can divide both sides of the
equality $\overline{a}^{p}=\overline{a}$ by $\overline{a}$. As a result, we
obtain $\overline{a}^{p-1}=\overline{1}$. In other words, $\overline{a^{p-1}%
}=\overline{1}$. In other words, $a^{p-1}\equiv1\operatorname{mod}p$. This
proves Corollary \ref{cor.ent.flt.p-1}.
\end{proof}

The following proposition generalizes Proposition \ref{prop.ent.flt.Z/p} to
arbitrary finite fields:\footnote{Recall that $\mathbb{Z}/p$ is a finite field
of size $p$ whenever $p$ is a prime. Moreover, the rings $F_{4}$ and $F_{8}$
are finite fields of sizes $4$ and $8$, respectively. We will see more finite
fields in Chapter \ref{chp.finfields}.}

\begin{proposition}
\label{prop.finfield.flt}Let $F$ be a finite field (i.e., a field with
finitely many elements). Let $u\in F$. Then, $u^{\left\vert F\right\vert }=u$.
\end{proposition}

\begin{proof}
If $u=0$, then this is obvious (since $u^{\left\vert F\right\vert
}=0^{\left\vert F\right\vert }=0=u$ in this case). So let us WLOG assume that
$u\neq0$. Hence, $u\in F\setminus\left\{  0\right\}  $.

However, $F$ is a field, so that every nonzero element of $F$ is a unit. In
other words, $F\setminus\left\{  0\right\}  \subseteq F^{\times}$. Conversely,
$F^{\times}\subseteq F\setminus\left\{  0\right\}  $, since $0$ is not a unit
of $F$ (because if $0$ were a unit, then $0\cdot0^{-1}$ would be $1$, which
contradicts the axiom $0a=0$ for all $a\in F$). Combining these two
inclusions, we find $F\setminus\left\{  0\right\}  =F^{\times}$. Hence,
$\left\vert F\setminus\left\{  0\right\}  \right\vert =\left\vert F^{\times
}\right\vert $, so that $\left\vert F^{\times}\right\vert =\left\vert
F\setminus\left\{  0\right\}  \right\vert =\left\vert F\right\vert -1$.

In other words, the group $F^{\times}$ has order $\left\vert F\right\vert -1$.
Hence, Lagrange's theorem (from group theory) shows that $g^{\left\vert
F\right\vert -1}=1$ for each $g\in F^{\times}$. Applying this to $g=u$, we
obtain $u^{\left\vert F\right\vert -1}=1$ (since $u\in F\setminus\left\{
0\right\}  =F^{\times}$). Hence, $u^{\left\vert F\right\vert }%
=u\underbrace{u^{\left\vert F\right\vert -1}}_{=1}=u$. This proves Proposition
\ref{prop.finfield.flt}.
\end{proof}

We note in passing that the converse of Fermat's little theorem does not hold
in general: There are some non-prime positive integers $p>1$ such that all
$a\in\mathbb{Z}$ satisfy $a^{p}\equiv a\operatorname{mod}p$. These integers
$p$ are called \textbf{Carmichael numbers}, and the smallest of them is $561$.

\begin{exercise}
Actually prove that $561$ is a Carmichael number, i.e., that every
$a\in\mathbb{Z}$ satisfies $a^{561}\equiv a\operatorname{mod}561$. \medskip

[\textbf{Hint:} This is not as laborious as it sounds! It is not necessary to
try all $561$ elements of $\mathbb{Z}/561$. Instead, use the prime
factorization $561=3\cdot11\cdot17$.]
\end{exercise}

\subsubsection{Division in a commutative ring}

Back to the general case. Rings have addition, subtraction and multiplication;
but we can also divide two elements of a ring, as long as the denominator
(i.e., the element we are dividing by) is a unit. If the ring is
noncommutative, this is somewhat complicated by the fact that there are two
kinds of division (\textquotedblleft left\textquotedblright\ and
\textquotedblleft right\textquotedblright\ division); however, for commutative
rings, it is as simple as for numbers:

\begin{definition}
\label{def.ringunits.div}Let $R$ be a commutative ring. Let $a\in R$ and $b\in
R^{\times}$. Then, $\dfrac{a}{b}$ means the element $ab^{-1}=b^{-1}a\in R$.
This element is also written $a/b$, and is called the \textbf{quotient} of $a$
by $b$. The operation $\left(  a,b\right)  \mapsto a/b$ is called
\textbf{division}.
\end{definition}

In particular, in a field, we can divide by any nonzero element.

Division satisfies the rules you would expect:

\begin{proposition}
\label{prop.ringunits.div.props}Let $R$ be a commutative ring. Then:

\begin{enumerate}
\item[\textbf{(a)}] For any $a,c\in R$ and $b,d\in R^{\times}$, we have%
\begin{equation}
\dfrac{a}{b}+\dfrac{c}{d}=\dfrac{ad+bc}{bd} \label{pf.ringunits.div.sum}%
\end{equation}
and%
\begin{equation}
\dfrac{a}{b}\cdot\dfrac{c}{d}=\dfrac{ac}{bd}. \label{pf.ringunits.div.prod}%
\end{equation}

\item[\textbf{(b)}] For any $a\in R$ and $b,c,d\in R^{\times}$, we have%
\[
\dfrac{a}{b}/\dfrac{c}{d}=\dfrac{ad}{bc}.
\]

\item[\textbf{(c)}] Division undoes multiplication: Three elements $a\in R$,
$b\in R^{\times}$ and $c\in R$ satisfy
\begin{equation}
\dfrac{a}{b}=c\ \ \ \ \ \ \ \ \ \ \text{if and only if }%
\ \ \ \ \ \ \ \ \ \ a=bc. \label{pf.ringunits.div.undo}%
\end{equation}

\end{enumerate}
\end{proposition}

\begin{exercise}
Prove Proposition \ref{prop.ringunits.div.props}.
\end{exercise}

\begin{exercise}
Let $p$ be a prime such that $p>3$. Prove that $2^{p-2}+3^{p-2}+6^{p-2}%
\equiv1\operatorname{mod}p$. \medskip

[\textbf{Hint:} First, show that $u^{p-2}=\dfrac{1}{u}$ for every nonzero
$u\in\mathbb{Z}/p$. Then, recall (\ref{pf.ringunits.div.sum}).]
\end{exercise}

\subsection{\label{sec.rings.mors}Ring morphisms (\cite[\S 7.3]{DumFoo04})}

\subsubsection{Definition and examples}

Groups have group homomorphisms; vector spaces have vector space homomorphisms
(= linear maps); topological spaces have topological space homomorphisms (=
continuous maps). No wonder that an analogous concept exists for
rings:\footnote{We follow the modern convention of abbreviating the word
\textquotedblleft homomorphism\textquotedblright\ as \textquotedblleft
morphism\textquotedblright. Thus, for example, a \textquotedblleft group
morphism\textquotedblright\ is the same as a group homomorphism.}

\begin{definition}
\label{def.ringmor.ringmor}Let $R$ and $S$ be two rings.

\begin{enumerate}
\item[\textbf{(a)}] A \textbf{ring homomorphism} (or, for short, \textbf{ring
morphism}, or, more informally, \textbf{ring map}) from $R$ to $S$ means a map
$f:R\rightarrow S$ that

\begin{itemize}
\item \textbf{respects addition} (i.e., satisfies $f\left(  a+b\right)
=f\left(  a\right)  +f\left(  b\right)  $ for all $a,b\in R$);

\item \textbf{respects multiplication} (i.e., satisfies $f\left(  ab\right)
=f\left(  a\right)  \cdot f\left(  b\right)  $ for all $a,b\in R$);

\item \textbf{respects the zero} (i.e., satisfies $f\left(  0_{R}\right)
=0_{S}$);

\item \textbf{respects the unity} (i.e., satisfies $f\left(  1_{R}\right)
=1_{S}$).
\end{itemize}

\item[\textbf{(b)}] A \textbf{ring isomorphism} (or, informally, \textbf{ring
iso}) from $R$ to $S$ means an invertible ring morphism $f:R\rightarrow S$
whose inverse $f^{-1}:S\rightarrow R$ is also a ring morphism.

\item[\textbf{(c)}] The rings $R$ and $S$ are said to be \textbf{isomorphic}
(this is written $R\cong S$) if there exists a ring isomorphism from $R$ to
$S$.
\end{enumerate}
\end{definition}

Here are some examples:

\begin{itemize}
\item Let $n\in\mathbb{Z}$. The map
\begin{align*}
\pi:\mathbb{Z}  &  \rightarrow\mathbb{Z}/n,\\
a  &  \mapsto\overline{a}%
\end{align*}
that sends each integer $a$ to its residue class $\overline{a}=a+n\mathbb{Z}$
is a ring morphism, because any $a,b\in\mathbb{Z}$ satisfy the equalities%
\[
\overline{a+b}=\overline{a}+\overline{b},\ \ \ \ \ \ \ \ \ \ \overline{a\cdot
b}=\overline{a}\cdot\overline{b},\ \ \ \ \ \ \ \ \ \ \overline{0}%
=0_{\mathbb{Z}/n},\ \ \ \ \ \ \ \ \ \ \overline{1}=1_{\mathbb{Z}/n}.
\]
(These equalities directly follow from the definition of the ring structure on
$\mathbb{Z}/n$.)

\item The map $\mathbb{Z}\rightarrow\mathbb{Z},\ a\mapsto2a$ is \textbf{not} a
ring morphism. It respects addition and the zero, but not multiplication and
the unity.

\item The map $\mathbb{Z}\rightarrow\mathbb{Z},\ a\mapsto0$ is \textbf{not} a
ring morphism. It respects addition, multiplication and the zero, but not the unity.

\item However, if $T$ is the zero ring (i.e., the $1$-element ring $\left\{
0\right\}  $), then the map $\mathbb{Z}\rightarrow T,\ a\mapsto0$ is a ring
morphism. Comparing this example with the preceding one, we see that the ring
structure (even a trivial-looking part like the unity) matters to whether a
given map is a ring morphism or not.

\item The map $\mathbb{Z}\rightarrow\mathbb{Z},\ a\mapsto a^{2}$ is
\textbf{not} a ring morphism. It respects multiplication, the zero and the
unity, but not addition (since $\left(  a+b\right)  ^{2}$ is usually not the
same as $a^{2}+b^{2}$).

\item Let $S$ be a subring of a ring $R$. Let $i:S\rightarrow R$ be the
\textbf{canonical inclusion}; this is simply the map that sends each $a\in S$
to itself. (You can view it as the restriction of the identity map
$\operatorname*{id}\nolimits_{R}:R\rightarrow R$ to $S$.) Then, $i$ is a ring
morphism. Indeed, it respects multiplication because the multiplication of $S$
is inherited from $R$ (so that any $a,b\in S$ satisfy $i\left(  ab\right)
=\underbrace{a}_{=i\left(  a\right)  }\underbrace{b}_{=i\left(  b\right)
}=i\left(  a\right)  i\left(  b\right)  $); for similar reasons, it satisfies
the other axioms in the definition of a ring morphism.

\item Consider the map%
\begin{align*}
f  &  :\mathbb{C}\rightarrow\mathbb{R}^{2\times2},\\
a+bi  &  \mapsto\left(
\begin{array}
[c]{cc}%
a & b\\
-b & a
\end{array}
\right)  \ \ \ \ \ \ \ \ \ \ \text{(for }a,b\in\mathbb{R}\text{).}%
\end{align*}
This map $f$ is a ring morphism. Indeed, it is easy to see that it respects
addition, the zero and the unity. To see that it respects multiplication, you
need to check that $f\left(  zw\right)  =f\left(  z\right)  \cdot f\left(
w\right)  $ for any $z,w\in\mathbb{C}$. But this is straightforward: Write
$z=a+bi$ and $w=c+di$ and multiply out\footnote{In more detail: Writing
$z=a+bi$ and $w=c+di$, we have $zw=\left(  a+bi\right)  \left(  c+di\right)
=\left(  ac-bd\right)  +\left(  ad+bc\right)  i$ and thus%
\[
f\left(  zw\right)  =\left(
\begin{array}
[c]{cc}%
ac-bd & ad+bc\\
-\left(  ad+bc\right)  & ac-bd
\end{array}
\right)  .
\]
However,%
\[
f\left(  z\right)  \cdot f\left(  w\right)  =\left(
\begin{array}
[c]{cc}%
a & b\\
-b & a
\end{array}
\right)  \left(
\begin{array}
[c]{cc}%
c & d\\
-d & c
\end{array}
\right)  =\left(
\begin{array}
[c]{cc}%
ac-bd & ad+bc\\
-\left(  ad+bc\right)  & ac-bd
\end{array}
\right)  .
\]
Comparing these two equalities yields $f\left(  zw\right)  =f\left(  z\right)
\cdot f\left(  w\right)  $.}.

This can also be proved using linear algebra: The $\mathbb{R}$-vector space
$\mathbb{C}$ has basis $\left(  1,i\right)  $. If $z\in\mathbb{C}$, then
$f\left(  z\right)  $ is the $2\times2$-matrix that represents the
\textquotedblleft multiply by $z$\textquotedblright\ operator (i.e., the map
$\mathbb{C}\rightarrow\mathbb{C},\ u\mapsto zu$) in this basis, where we
represent the elements of $\mathbb{C}$ as row vectors (a somewhat unusual
choice). Since the \textquotedblleft multiply by $zw$\textquotedblright%
\ operator is the composition of the \textquotedblleft multiply by
$z$\textquotedblright\ operator with the \textquotedblleft multiply by
$w$\textquotedblright\ operator, it thus follows that $f\left(  zw\right)
=f\left(  w\right)  \cdot f\left(  z\right)  $ (because composition of linear
maps corresponds to multiplication of their representing matrices). Of course,
since the multiplication in $\mathbb{C}$ is commutative, this can be rewritten
as $f\left(  wz\right)  =f\left(  w\right)  \cdot f\left(  z\right)  $, and
thus the map $f$ is a ring morphism (since it is easily seen to respect
addition, $0$ and $1$).

Note that the image of the map $f$ is precisely the ring $\mathcal{M}$ defined
in Exercise \ref{exe.idomain.PM}.

The ring morphism $f$ is injective, and therefore you can use the matrix
$f\left(  z\right)  \in\mathbb{R}^{2\times2}$ as a \textquotedblleft
stand-in\textquotedblright\ for any complex number $z\in\mathbb{C}$. Complex
numbers can thus be \textquotedblleft represented\textquotedblright\ by
$2\times2$-matrices with real entries. In particular, if you believe that the
complex numbers are a work of the devil\footnote{Historically,
\href{https://www.cut-the-knot.org/arithmetic/algebra/HistoricalRemarks.shtml}{mistrust
of complex numbers was widespread centuries after they had been first
introduced}. This mistrust was eventually overcome once Hamilton defined them
rigorously as pairs of real numbers (defining $a+bi$ as the pair $\left(
a,b\right)  $).}, then you can \textquotedblleft exorcise\textquotedblright%
\ them out of your mathematical work by replacing every complex number $z$
with the corresponding $2\times2$-matrix $f\left(  z\right)  $. Since $f$ is
injective, this replacement does not cause any information to be lost.
Furthermore, since $f$ is a ring morphism, addition and multiplication of
complex numbers are reflected perfectly in the addition and the multiplication
of $2\times2$-matrices, so that any calculation involving complex numbers
$z,w,u,\ldots$ can be immediately reproduced with the corresponding matrices
$f\left(  z\right)  ,f\left(  w\right)  ,f\left(  u\right)  ,\ldots$ instead.
Only the commutativity of multiplication is less clear when you work with
matrices: Two arbitrary $2\times2$-matrices don't usually commute, but of
course two $2\times2$-matrices that are values of $f$ always commute.

\item Just like the ring morphism $f:\mathbb{C}\rightarrow\mathbb{R}%
^{2\times2}$ can be used to represent complex numbers as $2\times2$-matrices
with real entries, there is another ring morphism $g:\mathbb{H}\rightarrow
\mathbb{R}^{4\times4}$ that helps represent Hamiltonian quaternions as
$4\times4$-matrices with real entries. This morphism is%
\begin{align*}
g:\mathbb{H}  &  \rightarrow\mathbb{R}^{4\times4},\\
a+bi+cj+dk  &  \mapsto\left(
\begin{array}
[c]{cccc}%
a & -b & -c & -d\\
b & a & -d & c\\
c & d & a & -b\\
d & -c & b & a
\end{array}
\right)  .
\end{align*}
Proving that this is a ring morphism is a tedious but doable exercise in calculation.

\item The map $\mathbb{R}^{2\times2}\rightarrow\mathbb{R},\ A\mapsto\det A$ is
\textbf{not} a ring morphism. It respects multiplication\footnote{This is a
particular case of the famous formula%
\[
\det\left(  AB\right)  =\det A\cdot\det B
\]
whenever $A,B\in R^{n\times n}$ are any two $n\times n$-matrices with entries
in any commutative ring $R$.} but not addition.

\item The map $\mathbb{C}\rightarrow\mathbb{C}$ that sends each complex number
$z=a+bi$ (with $a,b\in\mathbb{R}$) to its complex conjugate $\overline
{z}=a-bi$ is a ring isomorphism.

\item Let $R$ be a ring. Let $S$ be any set. Let $R^{S}$ be the ring of all
functions from $S$ to $R$ (with pointwise addition and multiplication). Fix
any $s\in S$. Then, the map
\begin{align*}
R^{S}  &  \rightarrow R,\\
f  &  \mapsto f\left(  s\right)
\end{align*}
is a ring morphism.\footnote{This is just a roundabout way of saying that any
maps $g,h\in R^{S}$ satisfy%
\begin{align*}
\left(  g+h\right)  \left(  s\right)   &  =g\left(  s\right)  +h\left(
s\right)  ;\\
\left(  gh\right)  \left(  s\right)   &  =g\left(  s\right)  \cdot h\left(
s\right)  ;\\
0\left(  s\right)   &  =0;\\
1\left(  s\right)   &  =1.
\end{align*}
But these equalities follow from our definition of the ring structure on
$R^{S}$ (namely: addition is pointwise; multiplication is pointwise; the zero
is the constant-$0$ function; the unity is the constant-$1$ function).} This
map is known as the \textbf{evaluation morphism} at $s$, since all it does is
evaluate a function at the input $s$.
\end{itemize}

Time for another warning:

\begin{warning}
Our Definition \ref{def.ringmor.ringmor} \textbf{(a)} again differs from
\cite{DumFoo04} in how it treats unities. Namely, \cite{DumFoo04} does not
require a ring morphism to respect the unity. Thus, the map $\mathbb{Z}%
\rightarrow\mathbb{Z},\ a\mapsto0$ is a ring morphism according to
\cite{DumFoo04}, but not according to us.
\end{warning}

\begin{exercise}
\label{exe.ringmor.F.omega}Define $A$ and $\mathcal{F}$ as in Exercise
\ref{exe.21hw1.6}.

Let $\omega:\mathcal{F}\rightarrow\mathcal{F}$ be the map that sends each
$aA+bI_{2}$ (with $a,b\in\mathbb{Z}$) to $-aA+\left(  a+b\right)  I_{2}$.
(This is well-defined, since each element of $\mathcal{F}$ can be
\textbf{uniquely} written as $aA+bI_{2}$ with $a,b\in\mathbb{Z}$.)

\begin{enumerate}
\item[\textbf{(a)}] Prove that $\omega$ is a ring morphism.

\item[\textbf{(b)}] Prove that $\omega\circ\omega=\operatorname*{id}$.

\item[\textbf{(c)}] Conclude that $\omega$ is a ring isomorphism.
\end{enumerate}
\end{exercise}

\begin{exercise}
\ \ 

\begin{enumerate}
\item[\textbf{(a)}] Is the map%
\begin{align*}
\mathbb{Z}^{2\times2}  &  \rightarrow\mathbb{Z}^{2\times2},\\
A  &  \mapsto A^{T}%
\end{align*}
(which sends each $2\times2$-matrix $A=\left(
\begin{array}
[c]{cc}%
a & b\\
c & d
\end{array}
\right)  $ to its transpose $A^{T}=\left(
\begin{array}
[c]{cc}%
a & c\\
b & d
\end{array}
\right)  $) a ring morphism?

\item[\textbf{(b)}] Is the map%
\begin{align*}
\mathbb{Z}^{2\times2}  &  \rightarrow\mathbb{Z}^{2\times2},\\
\left(
\begin{array}
[c]{cc}%
a & b\\
c & d
\end{array}
\right)   &  \mapsto\left(
\begin{array}
[c]{cc}%
d & c\\
b & a
\end{array}
\right)
\end{align*}
a ring morphism?
\end{enumerate}
\end{exercise}

\begin{exercise}
Let $R$ be any ring, and let $u$ be any unit of $R$. Consider the map%
\begin{align*}
f:R  &  \rightarrow R,\\
a  &  \mapsto uau^{-1}.
\end{align*}
Prove that this map $f$ is a ring isomorphism. (This map $f$ is called
\textbf{conjugation by }$u$. Despite its name, it has nothing to do with
conjugation of complex numbers.)
\end{exercise}

\begin{exercise}
\label{exe.ringmor.matrix-diagonal}Let $R$ be any ring, and let $n\in
\mathbb{N}$. Recall that $R^{n\leq n}$ is the subring of $R^{n\times n}$
consisting of the upper-triangular matrices.

Consider the map%
\begin{align*}
\delta:R^{n\times n}  &  \rightarrow R^{n\times n},\\
\left(
\begin{array}
[c]{cccc}%
a_{1,1} & a_{1,2} & \cdots & a_{1,n}\\
a_{2,1} & a_{2,2} & \cdots & a_{2,n}\\
\vdots & \vdots & \ddots & \vdots\\
a_{n,1} & a_{n,2} & \cdots & a_{n,n}%
\end{array}
\right)   &  \mapsto\left(
\begin{array}
[c]{cccc}%
a_{1,1} & 0 & \cdots & 0\\
0 & a_{2,2} & \cdots & 0\\
\vdots & \vdots & \ddots & \vdots\\
0 & 0 & \cdots & a_{n,n}%
\end{array}
\right)  ,
\end{align*}
which replaces all off-diagonal entries of a matrix by $0$ (but keeps the
diagonal entries unchanged).

\begin{enumerate}
\item[\textbf{(a)}] Is this map $\delta$ a ring morphism?

\item[\textbf{(b)}] Now, consider the restriction $\delta\mid_{R^{n\leq n}}$
of the map $\delta$ to the subring $R^{n\leq n}$. Is this restriction
$\delta\mid_{R^{n\leq n}}$ a ring morphism?
\end{enumerate}
\end{exercise}

\begin{exercise}
Let $R$ be any ring, and let $n\in\mathbb{N}$. For any $n\times n$-matrix
$A\in R^{n\times n}$, we consider the \textquotedblleft
block-diagonal\textquotedblright\ matrix $\left(
\begin{array}
[c]{cc}%
A & 0\\
0 & A
\end{array}
\right)  \in R^{2n\times2n}$, which is obtained by arranging two copies of $A$
and two zero matrices in the form suggested by the notation (for example, if
$A=\left(
\begin{array}
[c]{cc}%
a & b\\
c & d
\end{array}
\right)  $, then $\left(
\begin{array}
[c]{cc}%
A & 0\\
0 & A
\end{array}
\right)  =\left(
\begin{array}
[c]{cccc}%
a & b & 0 & 0\\
c & d & 0 & 0\\
0 & 0 & a & b\\
0 & 0 & c & d
\end{array}
\right)  $). (See, e.g., \href{https://en.wikipedia.org/wiki/Block_matrix}{the
Wikipedia} for more details about block matrices.)

\begin{enumerate}
\item[\textbf{(a)}] Prove that the map%
\begin{align*}
R^{n\times n}  &  \rightarrow R^{2n\times2n},\\
A  &  \mapsto\left(
\begin{array}
[c]{cc}%
A & 0\\
0 & A
\end{array}
\right)
\end{align*}
is an injective ring morphism.

\item[\textbf{(b)}] Generalize this to find an injective ring morphism from
$R^{n\times n}$ to $R^{kn\times kn}$ for every positive integer $k$.
\end{enumerate}
\end{exercise}

The following exercise (\cite[homework set \#2, Exercise 2]{21w}) assigns to
each ring $R$ a \textquotedblleft mirror version\textquotedblright\ (called
the \textbf{opposite ring} of $R$, and denoted by $R^{\operatorname*{op}}$).
In general, this mirror version is not isomorphic to $R$, but often enough it is.

\begin{exercise}
\label{exe.21hw2.2abc}Let $R$ be a ring. We define a new binary operation
$\mathbin{\widetilde{\cdot}}$ on $R$ by setting
\[
a\mathbin{\widetilde{\cdot}}b=ba\qquad\text{for all }a,b\in R.
\]
(Thus, $\mathbin{\widetilde{\cdot}}$ is the multiplication of $R$, but with
the two arguments switched.)

\begin{enumerate}
\item[\textbf{(a)}] Prove that the set $R$, equipped with the addition $+$,
the multiplication $\mathbin{\widetilde{\cdot}}$, the zero $0_{R}$ and the
unity $1_{R}$, is a ring.
\end{enumerate}

This new ring is called the \textbf{opposite ring} of $R$, and is denoted by
$R^{\operatorname{op}}$.

Note that the \textbf{sets} $R$ and $R^{\operatorname{op}}$ are identical (so
a map from $R$ to $R$ is the same as a map from $R$ to $R^{\operatorname{op}}%
$); but the \textbf{rings} $R$ and $R^{\operatorname{op}}$ are generally not
the same (so a ring morphism from $R$ to $R$ is not the same as a ring
morphism from $R$ to $R^{\operatorname{op}}$).

\begin{enumerate}
\item[\textbf{(b)}] Prove that the identity map $\operatorname{id}%
:R\rightarrow R$ is a ring isomorphism from $R$ to $R^{\operatorname{op}}$ if
and only if $R$ is commutative.

\item[\textbf{(c)}] Now, assume that $R$ is the matrix ring $S^{n\times n}$
for some commutative ring $S$ and some $n\in\mathbb{N}$. Prove that the map
\[
R\rightarrow R^{\operatorname{op}},\qquad A\mapsto A^{T}%
\]
(where $A^{T}$, as usual, denotes the transpose of a matrix $A$) is a ring isomorphism.
\end{enumerate}
\end{exercise}

\begin{exercise}
Let $\mathbb{H}$ be the ring of Hamiltonian quaternions.

\begin{enumerate}
\item[\textbf{(a)}] Prove that the map%
\begin{align*}
\mathbb{H}  &  \rightarrow\mathbb{H}^{\operatorname*{op}},\\
a+bi+cj+dk  &  \mapsto a+bi+dj+ck\ \ \ \ \ \ \ \ \ \ \left(  \text{for
}a,b,c,d\in\mathbb{R}\right)
\end{align*}
is a ring isomorphism.

\item[\textbf{(b)}] Prove that the map%
\begin{align*}
\mathbb{H}  &  \rightarrow\mathbb{H}^{\operatorname*{op}},\\
a+bi+cj+dk  &  \mapsto a-bi-cj-dk\ \ \ \ \ \ \ \ \ \ \left(  \text{for
}a,b,c,d\in\mathbb{R}\right)
\end{align*}
is a ring isomorphism as well.
\end{enumerate}
\end{exercise}

The last two exercises might suggest that every ring $R$ is somehow isomorphic
to its opposite ring $R^{\operatorname*{op}}$. This is not the case, but a
counterexample is tricky to find; one such counterexample is constructed in
Exercise \ref{exe.rings.opposite.not-iso}.

\subsubsection{\label{subsec.ringmor.mors.basic}Basic properties of ring
morphisms}

Let us show some basic properties of ring morphisms. We start with the fact
that a composition of two ring morphisms is again a ring morphism:

\begin{proposition}
\label{prop.ringmor.compose}Let $R$, $S$ and $T$ be three rings. Let
$f:S\rightarrow T$ and $g:R\rightarrow S$ be two ring morphisms. Then, $f\circ
g:R\rightarrow T$ is a ring morphism.
\end{proposition}

\begin{proof}
This is proved in the same way as the analogous result about groups.
\end{proof}

The next proposition shows that the \textquotedblleft respects the
zero\textquotedblright\ condition in the definition of a ring morphism is
redundant (even though the \textquotedblleft respects the
unity\textquotedblright\ condition is not):

\begin{proposition}
\label{prop.ringmor.resp0}Let $R$ and $S$ be two rings. Let $f:R\rightarrow S$
be a map that respects addition. Then, $f$ automatically respects the zero.
\end{proposition}

\begin{proof}
Since $f$ respects addition, we have $f\left(  0_{R}+0_{R}\right)  =f\left(
0_{R}\right)  +f\left(  0_{R}\right)  $. Rewrite this as $f\left(
0_{R}\right)  =f\left(  0_{R}\right)  +f\left(  0_{R}\right)  $ (since
$0_{R}+0_{R}=0_{R}$). Now, subtract $f\left(  0_{R}\right)  $ from both sides
to get $0_{S}=f\left(  0_{R}\right)  $. In other words, $f$ respects the zero.
\end{proof}

Note that we can restate our definition of a ring morphism as follows:

\begin{statement}
A \textit{ring morphism} is a map $f:R\rightarrow S$ between two rings $R$ and
$S$ that is a group homomorphism from the additive group $\left(
R,+,0\right)  $ to the additive group $\left(  S,+,0\right)  $ and
simultaneously a monoid homomorphism from the multiplicative monoid $\left(
R,\cdot,1\right)  $ to the multiplicative monoid $\left(  S,\cdot,1\right)  $.
\end{statement}

It is easy to see that ring morphisms respect all sorts of operations
constructed from $+$, $\cdot$, $0$ and $1$:

\begin{proposition}
\label{prop.ringmor.respects}Let $R$ and $S$ be two rings. Let $f:R\rightarrow
S$ be a ring morphism. Then:

\begin{enumerate}
\item[\textbf{(a)}] The map $f$ respects finite sums; i.e., we have $f\left(
a_{1}+a_{2}+\cdots+a_{n}\right)  =f\left(  a_{1}\right)  +f\left(
a_{2}\right)  +\cdots+f\left(  a_{n}\right)  $ for any $a_{1},a_{2}%
,\ldots,a_{n}\in R$.

\item[\textbf{(b)}] The map $f$ respects finite products; i.e., we have
$f\left(  a_{1}a_{2}\cdots a_{n}\right)  =f\left(  a_{1}\right)  \cdot
f\left(  a_{2}\right)  \cdot\cdots\cdot f\left(  a_{n}\right)  $ for any
$a_{1},a_{2},\ldots,a_{n}\in R$.

\item[\textbf{(c)}] The map $f$ respects differences; i.e., we have $f\left(
a-b\right)  =f\left(  a\right)  -f\left(  b\right)  $ for any $a,b\in R$.

\item[\textbf{(d)}] The map $f$ respects inverses; i.e., if $a$ is a unit of
$R$, then $f\left(  a\right)  $ is a unit of $S$, with inverse $\left(
f\left(  a\right)  \right)  ^{-1}=f\left(  a^{-1}\right)  $.

\item[\textbf{(e)}] The map $f$ respects integer multiples; i.e., if $a\in R$
and $n\in\mathbb{Z}$, then $f\left(  na\right)  =nf\left(  a\right)  $.

\item[\textbf{(f)}] The map $f$ respects powers; i.e., if $a\in R$ and
$n\in\mathbb{N}$, then $f\left(  a^{n}\right)  =\left(  f\left(  a\right)
\right)  ^{n}$.
\end{enumerate}
\end{proposition}

\begin{proof}
This is pretty straightforward, and you have probably seen the idea in group
theory already. Details LTTR\footnote{\textquotedblleft LTTR\textquotedblright%
\ means \textquotedblleft left to the reader\textquotedblright.}.
\end{proof}

\subsubsection{The image of a ring morphism}

Recall that the \textbf{image} of a map $f:R\rightarrow S$ is defined to be
the set $f\left(  R\right)  =\left\{  f\left(  r\right)  \ \mid\ r\in
R\right\}  $; it is often denoted $\operatorname{Im}f$. This makes sense for
arbitrary maps $f$ between arbitrary sets $R$ and $S$, not just for ring
morphisms between rings. However, the image of a ring morphism has a special property:

\begin{proposition}
\label{prop.ringmor.Im-subring}Let $R$ and $S$ be two rings. Let
$f:R\rightarrow S$ be a ring morphism. Then, $\operatorname{Im}f=f\left(
R\right)  $ is a subring of $S$.
\end{proposition}

\begin{proof}
Just check the axioms for a subring. For example, let us show that $f\left(
R\right)  $ is closed under multiplication:

Let $x,y\in f\left(  R\right)  $. We must show that $xy\in f\left(  R\right)
$. Since $x\in f\left(  R\right)  $, there exists some $a\in R$ such that
$x=f\left(  a\right)  $. Similarly, there exists some $b\in R$ such that
$y=f\left(  b\right)  $. Consider these $a$ and $b$. From $x=f\left(
a\right)  $ and $y=f\left(  b\right)  $, we obtain
\begin{align*}
xy  &  =f\left(  a\right)  \cdot f\left(  b\right)  =f\left(  ab\right)
\ \ \ \ \ \ \ \ \ \ \left(  \text{since }f\text{ respects multiplication}%
\right) \\
&  \in f\left(  R\right)  .
\end{align*}

This completes the proof that $f\left(  R\right)  $ is closed under
multiplication. The other ring axioms can be verified similarly. Thus, we
conclude that $f\left(  R\right)  $ is a subring of $S$. Proposition
\ref{prop.ringmor.Im-subring} is proved.
\end{proof}

\begin{exercise}
\ \ 

\begin{enumerate}
\item[\textbf{(a)}] Let $R$ and $S$ be two rings. Let $f:R\rightarrow S$ and
$g:R\rightarrow S$ be two ring morphisms. Let $\operatorname*{Eq}\left(
f,g\right)  $ be the subset%
\[
\left\{  r\in R\ \mid\ f\left(  r\right)  =g\left(  r\right)  \right\}
\]
of $R$. This subset is called the \textbf{equalizer} of $f$ and $g$. Prove
that this subset $\operatorname*{Eq}\left(  f,g\right)  $ is a subring of $R$.

\item[\textbf{(b)}] Let $\omega:\mathbb{C}\rightarrow\mathbb{C}$ be the map
sending each complex number $z=a+bi$ to its complex conjugate $\overline
{z}=a-bi$. (Recall that this is a ring morphism.) Prove that the equalizer
$\operatorname*{Eq}\left(  \operatorname*{id}\nolimits_{\mathbb{C}}%
,\omega\right)  $ is $\mathbb{R}$.

\item[\textbf{(c)}] Find a specific example where the equalizer subring
$\operatorname*{Eq}\left(  f,g\right)  $ is interesting (i.e., ideally a ring
you have not seen before).
\end{enumerate}

[\textbf{Hint:} There are many good examples for part \textbf{(c)}. For
instance, using the notation $R^{n\leq n}$ from Subsection
\ref{subsec.rings.subrings.exas}, consider the two ring morphisms%
\begin{align*}
f  &  :\mathbb{Q}^{3\leq3}\rightarrow\mathbb{Q}^{2\leq2}%
,\ \ \ \ \ \ \ \ \ \ \left(
\begin{array}
[c]{ccc}%
a & b & c\\
0 & d & e\\
0 & 0 & h
\end{array}
\right)  \mapsto\left(
\begin{array}
[c]{cc}%
a & b\\
0 & d
\end{array}
\right)  \ \ \ \ \ \ \ \ \ \ \text{and}\\
g  &  :\mathbb{Q}^{3\leq3}\rightarrow\mathbb{Q}^{2\leq2}%
,\ \ \ \ \ \ \ \ \ \ \left(
\begin{array}
[c]{ccc}%
a & b & c\\
0 & d & e\\
0 & 0 & h
\end{array}
\right)  \mapsto\left(
\begin{array}
[c]{cc}%
d & e\\
0 & h
\end{array}
\right)  .
\end{align*}
What is their equalizer?]
\end{exercise}

\subsubsection{\label{subsec.ringmor.iso.basic-props}Basic properties of ring
isomorphisms}

We shall now show some fundamental facts about ring isomorphisms.

First, let us give a somewhat simplified characterization of ring
isomorphisms. According to Definition \ref{def.ringmor.ringmor} \textbf{(b)},
if you want to prove that a map $f$ is a ring isomorphism, you have to check
(1) that $f$ is a ring morphism, (2) that $f$ has an inverse, and (3) that
this inverse $f^{-1}$ is a ring morphism. However, it turns out that step (3)
is unnecessary, since it follows from steps (1) and (2). Let us state this
fact and prove it:

\begin{proposition}
\label{prop.ringmor.invertible-iso}Let $R$ and $S$ be two rings. Let
$f:R\rightarrow S$ be an invertible ring morphism. Then, $f$ is a ring isomorphism.
\end{proposition}

\begin{proof}
This is proved using the same reasoning as for groups (but not for topological
spaces): You need to show that $f^{-1}$ is a ring morphism. Let me just show
that $f^{-1}$ respects addition (the proofs of the other axioms are similar).
So let $c,d\in S$; we must show that $f^{-1}\left(  c+d\right)  =f^{-1}\left(
c\right)  +f^{-1}\left(  d\right)  $.

It is clearly sufficient to check that $f\left(  f^{-1}\left(  c+d\right)
\right)  =f\left(  f^{-1}\left(  c\right)  +f^{-1}\left(  d\right)  \right)
$. Indeed, if we can show this equality, then we can apply $f^{-1}$ to it and
obtain $f^{-1}\left(  c+d\right)  =f^{-1}\left(  c\right)  +f^{-1}\left(
d\right)  $, which is what we want to prove.

Recall that $f$ respects addition. Thus,%
\[
f\left(  f^{-1}\left(  c\right)  +f^{-1}\left(  d\right)  \right)  =f\left(
f^{-1}\left(  c\right)  \right)  +f\left(  f^{-1}\left(  d\right)  \right)
=c+d=f\left(  f^{-1}\left(  c+d\right)  \right)  .
\]
Hence, $f\left(  f^{-1}\left(  c+d\right)  \right)  =f\left(  f^{-1}\left(
c\right)  +f^{-1}\left(  d\right)  \right)  $ is proved.
\end{proof}

Incidentally, \cite{DumFoo04} defines ring isomorphisms as invertible ring
morphisms. Proposition \ref{prop.ringmor.invertible-iso} shows that this is
equivalent to our definition.

We continue with some more straightforward results:

\begin{proposition}
\label{prop.ringiso.compose}Let $R$, $S$ and $T$ be three rings. Let
$f:S\rightarrow T$ and $g:R\rightarrow S$ be two ring isomorphisms. Then,
$f\circ g:R\rightarrow T$ is a ring isomorphism.
\end{proposition}

\begin{proof}
This is proved in the same way as for groups.
\end{proof}

\begin{proposition}
\label{prop.ringiso.inverse}Let $R$ and $S$ be two rings. Let $f:R\rightarrow
S$ be a ring isomorphism. Then, $f^{-1}:S\rightarrow R$ is a ring isomorphism.
\end{proposition}

\begin{proof}
This is proved in the same way as for groups.
\end{proof}

\begin{corollary}
\label{cor.ringiso.equiv}The relation $\cong$ for rings is an equivalence relation.
\end{corollary}

\begin{proof}
Transitivity follows from Proposition \ref{prop.ringiso.compose}. Reflexivity
follows from the obvious fact that $\operatorname*{id}:R\rightarrow R$ is a
ring isomorphism whenever $R$ is a ring. Symmetry follows from Proposition
\ref{prop.ringiso.inverse}.
\end{proof}

The most useful property of ring isomorphisms is the following
\textquotedblleft meta-theorem\textquotedblright:

\begin{statement}
\textbf{Isomorphism principle for rings:} Let $R$ and $S$ be two isomorphic
rings. Then, any \textquotedblleft ring-theoretic\textquotedblright\ property
of $R$ (that is, any property that does not refer to what the elements of $R$
are, but can be stated entirely in terms of how they \textquotedblleft
interact\textquotedblright\ via the ring operations) that holds for $R$ must
hold for $S$ as well.
\end{statement}

This is somewhat nebulous: What exactly makes a property of a ring
\textquotedblleft ring-theoretic\textquotedblright? In lieu of a formal
definition, let us give some examples of \textquotedblleft
ring-theoretic\textquotedblright\ properties of $R$ (which may or may not hold):

\begin{itemize}
\item The ring $R$ has $15$ elements.

\item The ring $R$ is commutative.

\item The ring $R$ is a field.

\item For any $a,b,c\in R$, we have $3abc\left(  a+b+c\right)  =0_{R}$ (where
the \textquotedblleft$3$\textquotedblright\ is the integer $3$).

\item The center of $R$ has $10$ elements.

\item There exist two nonzero elements $a,b\in R$ satisfying $a^{2}%
+b^{2}=0_{R}$.
\end{itemize}

Thus, all of these properties can be automatically transferred from any ring
to any isomorphic ring.

In contrast, here are some examples of properties of $R$ that are \textbf{not}
\textquotedblleft ring-theoretic\textquotedblright:

\begin{itemize}
\item The elements of $R$ are matrices.

\item The ring $R$ is a subset of $\mathbb{R}$.

\item The ring $R$ is a subring of $\mathbb{R}$. (\textbf{But} the more
professional variant \textquotedblleft there is an injective ring morphism
from $R$ to $\mathbb{R}$\textquotedblright\ would be ring-theoretic!)

\item The set $R$ is disjoint from $\mathbb{C}$.

\item There exist two nonzero elements $a,b\in R$ satisfying $a^{2}%
+b^{2}=0_{\mathbb{Z}}$ (the zero of $\mathbb{Z}$).

\item The set $R$ contains the complex number $i=\sqrt{-1}$.
\end{itemize}

Clearly, an isomorphism can destroy these properties, since it can send
elements to different elements.

To make sure you understand the meaning of ring isomorphisms, pick any of the
above \textquotedblleft ring-theoretic\textquotedblright\ properties of $R$,
and show that it is preserved by isomorphisms (i.e., if it holds for a ring
$R$, then it holds for any ring $S$ isomorphic to $R$). The proof is analogous
to the similar argument for groups.

The following exercise shows an example of a non-obvious ring isomorphism:

\begin{exercise}
Define a subring $\mathcal{M}$ of $\mathbb{R}^{2\times2}$ as in Exercise
\ref{exe.idomain.PM}. Consider the map%
\begin{align*}
f  &  :\mathbb{C}\rightarrow\mathcal{M},\\
a+bi  &  \mapsto\left(
\begin{array}
[c]{cc}%
a & b\\
-b & a
\end{array}
\right)  \ \ \ \ \ \ \ \ \ \ \text{(for }a,b\in\mathbb{R}\text{).}%
\end{align*}

\begin{enumerate}
\item[\textbf{(a)}] Prove that $f$ is a ring isomorphism.

\item[\textbf{(b)}] Use this to solve Exercise \ref{exe.idomain.PM}
\textbf{(c)} again.
\end{enumerate}
\end{exercise}

\subsubsection{Injective morphisms and their images}

If $f:R\rightarrow S$ is an injective map from some set $R$ to some set $S$,
then its image $f\left(  R\right)  $ is in one-to-one correspondence with its
domain $R$ (via the map $R\rightarrow f\left(  R\right)  $ that sends each $r$
to $f\left(  r\right)  $). The same holds for ring morphisms, except that the
one-to-one correspondence is now a ring isomorphism:

\begin{proposition}
\label{prop.ringmor.inj-image}Let $R$ and $S$ be two rings. Let
$f:R\rightarrow S$ be an injective ring morphism. Then:

\begin{enumerate}
\item[\textbf{(a)}] The subring $f\left(  R\right)  $ of $S$ (known from
Proposition \ref{prop.ringmor.Im-subring}) is isomorphic to $R$.

\item[\textbf{(b)}] More specifically: The map%
\begin{align*}
R  &  \rightarrow f\left(  R\right)  ,\\
r  &  \mapsto f\left(  r\right)
\end{align*}
is a ring isomorphism.
\end{enumerate}
\end{proposition}

\begin{proof}
The map%
\begin{align*}
R  &  \rightarrow f\left(  R\right)  ,\\
r  &  \mapsto f\left(  r\right)
\end{align*}
is clearly well-defined (since $f\left(  r\right)  \in f\left(  R\right)  $
for each $r\in R$). Let us denote it by $f^{\prime}$. This map $f^{\prime}$
differs from the map $f$ only in that it goes to $f\left(  R\right)  $ rather
than to $S$. Hence, this map $f^{\prime}$ is injective (since $f$ is
injective) and surjective (since each element of $f\left(  R\right)  $ has the
form $f\left(  r\right)  $ for some $r\in R$ by definition, and thus equals
$f^{\prime}\left(  r\right)  $ for the same $r\in R$). Hence, it is bijective,
i.e., invertible.

Moreover, $f\left(  R\right)  $ is a subring of $S$, so that its addition,
multiplication, zero and unity are inherited from $S$. Hence, from the fact
that $f$ is a ring morphism, we conclude immediately that the map $f^{\prime}$
(which differs from $f$ only in that it goes to $f\left(  R\right)  $ rather
than to $S$) is a ring morphism as well. Thus, $f^{\prime}$ is an invertible
ring morphism, hence (by Proposition \ref{prop.ringmor.invertible-iso}) a ring
isomorphism. In other words, the map%
\begin{align*}
R  &  \rightarrow f\left(  R\right)  ,\\
r  &  \mapsto f\left(  r\right)
\end{align*}
is a ring isomorphism (since this map is $f^{\prime}$). Hence, $f\left(
R\right)  $ is isomorphic to $R$. This proves Proposition
\ref{prop.ringmor.inj-image}.
\end{proof}

\subsubsection{Advanced exercises on ring isomorphisms}

Here are some further exercises on ring isomorphisms.

\begin{exercise}
\label{exe.rings.opposite.not-iso}Let $F$ be a nontrivial commutative ring.
Let $R$ be the set of all $3\times3$-matrices of the form $\left(
\begin{array}
[c]{ccc}%
a & b & c\\
0 & d & 0\\
0 & 0 & e
\end{array}
\right)  \in F^{3\times3}$ with $a,b,c,d,e\in F$. It is not hard to see that
$R$ is a subring of $F^{3\times3}$.

\begin{enumerate}
\item[\textbf{(a)}] Prove that if $A,B,C\in R$ are any three matrices in $R$,
then the matrix $C\left(  AB-BA\right)  \in R$ is a scalar multiple of the
matrix $AB-BA$. (A \textquotedblleft scalar multiple\textquotedblright\ of a
matrix $M$ means a matrix of the form $\lambda M$ with $\lambda\in F$.)

\item[\textbf{(b)}] Prove that it is not always true that if $A,B,C\in R$ are
any three matrices in $R$, then the matrix $\left(  BA-AB\right)  C\in R$ is a
scalar multiple of the matrix $BA-AB$.

\item[\textbf{(c)}] Conclude that $R$ is not isomorphic to
$R^{\operatorname*{op}}$ when $F=\mathbb{Z}$ or $F=\mathbb{Z}/2$.
\end{enumerate}
\end{exercise}

The next two exercises characterize rings of certain small sizes:

\begin{exercise}
\label{exe.rings.size-p}Let $p$ be a prime. Let $R$ be a ring of size
$\left\vert R\right\vert =p$. Prove that $R$ is isomorphic to the ring
$\mathbb{Z}/p$. \medskip

[\textbf{Hint:} Prove that the additive group $\left(  R,+,0\right)  $ must be
generated by $1$. Argue that this uniquely determines the multiplication of
$R$.]
\end{exercise}

\begin{exercise}
\label{exe.rings.size-4}Let $R$ be a ring of size $\left\vert R\right\vert
=4$. Prove that $R$ is isomorphic to one of the four rings $\mathbb{Z}/4$,
$F_{4}$, $D_{4}$ and $B_{4}$ we have seen in Subsection
\ref{subsec.rings.def.exas}. \medskip

[\textbf{Hint:} As in Exercise \ref{exe.rings.size-8-comm}, consider the
subgroup $\left\langle 1\right\rangle $ of $\left(  R,+,0\right)  $. If this
subgroup is the whole $R$, then argue that $R\cong\mathbb{Z}/4$. If not,
choose an arbitrary $x\in R\setminus\left\langle 1\right\rangle $, and
distinguish cases based on what $x^{2}$ is.]
\end{exercise}

\subsection{\label{sec.rings.ideals}Ideals and kernels (\cite[\S 7.1]%
{DumFoo04})}

\subsubsection{Kernels}

In linear algebra, the kernel (aka nullspace) of a linear map
\textquotedblleft measures how non-injective it is\textquotedblright. The same
can be done for ring morphisms:

\begin{definition}
Let $R$ and $S$ be two rings. Let $f:R\rightarrow S$ be a ring morphism. Then,
the \textbf{kernel} of $f$ (denoted $\ker f$ or $\operatorname*{Ker}f$) is
defined to be the subset%
\[
\operatorname*{Ker}f:=\left\{  a\in R\ \mid\ f\left(  a\right)  =0_{S}%
\right\}
\]
of $R$.
\end{definition}

Some examples:

\begin{itemize}
\item Let $n\in\mathbb{Z}$. The kernel of the ring morphism $\pi
:\mathbb{Z}\rightarrow\mathbb{Z}/n,\ a\mapsto\overline{a}$ is $n\mathbb{Z}%
=\left\{  \text{all multiples of }n\right\}  $.

\item Let $R$ be a ring. Let $S$ be any set. Recall the ring $R^{S}$ of all
functions from $S$ to $R$. Fix an element $s\in S$. Then, the kernel of the
ring morphism $R^{S}\rightarrow R,\ f\mapsto f\left(  s\right)  $ is the set
of all functions $f\in R^{S}$ that vanish on $s$.

\item The kernel of an injective ring morphism $f:R\rightarrow S$ is always
$\left\{  0_{R}\right\}  $. Indeed, if $f:R\rightarrow S$ is an injective ring
morphism, then $f$ sends $0_{R}$ to $0_{S}$ (since $f$ is a ring morphism),
and therefore $f$ cannot send any other element to $0_{S}$ (since $f$ is injective).
\end{itemize}

\subsubsection{Ideals}

As we saw in the above example, the kernel of a ring morphism is not usually a
subring of $R$, since it normally does not contain $1_{R}$. However, it
satisfies all the other axioms for a subring (which is why \cite{DumFoo04}
considers it a subring of $R$). We can say more, however. The type of a subset
that kernels of ring morphisms are has its own name:

\begin{definition}
\label{def.ideal}Let $R$ be a ring. An \textbf{ideal} of $R$ is a subset $I$
of $R$ such that

\begin{itemize}
\item we have $a+b\in I$ for any $a,b\in I$;

\item we have $ab\in I$ and $ba\in I$ for any $a\in R$ and $b\in I$;

\item we have $0\in I$ (where the $0$ means the zero of $R$).
\end{itemize}
\end{definition}

When $R$ is commutative, of course, the \textquotedblleft$ab\in I$%
\textquotedblright\ and \textquotedblleft$ba\in I$\textquotedblright%
\ conditions are equivalent.

The three conditions in Definition \ref{def.ideal} are called the
\textquotedblleft\textbf{ideal axioms}\textquotedblright. The first and the
third of them are familiar (they already appeared in the definition of a
subring). The second is new -- it is saying that if a factor in a product
belongs to $I$, then the whole product belongs to $I$, no matter what the
other factors are.\footnote{This second axiom is sometimes called the
\textquotedblleft\textbf{absorption axiom}\textquotedblright, referring to the
idea that the ideal $I$ \textquotedblleft absorbs\textquotedblright\ every
product as long as even one factor of the product is in the ideal. I prefer to
think of it as \textquotedblleft contagiousness\textquotedblright. Another
picture in my mind is that $I$ is some kind of ditch which you can enter but
never escape through multiplication with elements of $R$.}

Here are some easy consequences of Definition \ref{def.ideal}:

\begin{proposition}
\label{prop.ideal.subgroup}Let $R$ be a ring. Let $I$ be an ideal of $R$.
Then, $I$ is a subgroup of the additive group $\left(  R,+,0\right)  $.
\end{proposition}

\begin{proof}
The first and third \textquotedblleft ideal axioms\textquotedblright\ reveal
that $I$ is closed under addition and contains $0$. It remains to show that
$I$ is closed under negation -- i.e., that we have $-b\in I$ for each $b\in
I$. But this is easy: If $b\in I$, then the second \textquotedblleft ideal
axiom\textquotedblright\ (applied to $a=-1$) yields $\left(  -1\right)  b\in
I$ and $b\left(  -1\right)  \in I$. But this rewrites as $-b\in I$, qed.
\end{proof}

\begin{theorem}
\label{thm.ringmor.ker-ideal}Let $R$ and $S$ be two rings. Let $f:R\rightarrow
S$ be a ring morphism. Then, the kernel $\operatorname*{Ker}f$ of $f$ is an
ideal of $R$.
\end{theorem}

\begin{proof}
We need to prove the three \textquotedblleft ideal axioms\textquotedblright.
Let me only show the second, as the other two are similar. So let $a\in R$ and
$b\in\operatorname*{Ker}f$. We must prove that $ab\in\operatorname*{Ker}f$ and
$ba\in\operatorname*{Ker}f$.

We have $b\in\operatorname*{Ker}f$, so that $f\left(  b\right)  =0$ (by the
definition of $\operatorname*{Ker}f$). Now, the map $f$ is a ring morphism and
thus respects multiplication. Hence, $f\left(  ab\right)  =f\left(  a\right)
\cdot\underbrace{f\left(  b\right)  }_{=0}=f\left(  a\right)  \cdot0=0$, so
that $ab\in\operatorname*{Ker}f$ (by the definition of $\operatorname*{Ker}%
f$). Similarly, $ba\in\operatorname*{Ker}f$. Thus we have shown the second
ideal axiom.
\end{proof}

We will soon see a converse of this theorem: Every ideal of a ring $R$ is the
kernel of some ring morphism from $R$. (Namely, this follows from Theorem
\ref{thm.quotring.canproj} below.)

\subsubsection{Principal ideals}

The simplest way to construct ideals of a commutative ring is by fixing an
element and taking all its multiples:

\begin{proposition}
\label{prop.ideal.princid}Let $R$ be a commutative ring.

\begin{enumerate}
\item[\textbf{(a)}] Let $u\in R$. We define $uR$ to be the set $\left\{
ur\ \mid\ r\in R\right\}  $. The elements of this set $uR$ are called the
\textbf{multiples} of $u$ (in $R$).

Then, $uR$ is an ideal of $R$. This ideal is known as a \textbf{principal
ideal} of $R$.

\item[\textbf{(b)}] In particular, $0R=\left\{  0_{R}\right\}  $ and $1R=R$
are therefore principal ideals of $R$.
\end{enumerate}
\end{proposition}

\begin{proof}
\textbf{(a)} The only thing to prove is that $uR$ is an ideal of $R$. But this
can be easily achieved by checking that it satisfies all three ideal axioms:

\begin{itemize}
\item We have $a+b\in uR$ for any $a\in uR$ and $b\in uR$. (Indeed, if $a\in
uR$ and $b\in uR$, then there exist $x,y\in R$ satisfying $a=ux$ and $b=uy$
(since $a\in uR$ and $b\in uR$), and therefore we have $a+b=ux+uy=u\left(
x+y\right)  \in uR$.)

\item We have $ab\in uR$ and $ba\in uR$ for any $a\in R$ and $b\in uR$.
(Indeed, if $a\in R$ and $b\in uR$, then there exists an $r\in R$ satisfying
$b=ur$ (since $b\in uR$), and thus we have $ab=aur=u\left(  ar\right)  \in uR$
and therefore $ba=ab\in uR$.)

\item We have $0\in uR$ (since $0=u\cdot0$).
\end{itemize}

Thus, Proposition \ref{prop.ideal.princid} \textbf{(a)} is proved. \medskip

\textbf{(b)} The equalities $0R=\left\{  0_{R}\right\}  $ and $1R=R$ are
obvious (since each $r\in R$ satisfies $0r=0$ and $1r=r$). The rest follows
from part \textbf{(a)}.
\end{proof}

For example, $2\mathbb{Z}=\left\{  \text{all even integers}\right\}  $ is an
ideal of $\mathbb{Z}$.

\begin{exercise}
\label{exe.ideal.field-has-2}\ \ 

\begin{enumerate}
\item[\textbf{(a)}] Let $F$ be a field. Prove that the only ideals of $F$ are
$0F=\left\{  0_{F}\right\}  $ and $1F=F$.

\item[\textbf{(b)}] Conversely, let $R$ be a nontrivial commutative ring that
has only two ideals. Prove that $R$ is a field.
\end{enumerate}
\end{exercise}

The requirement that $R$ be commutative in Proposition
\ref{prop.ideal.princid} was not gratuitous; the set $uR$ is not always an
ideal when $R$ is not commutative. Nevertheless, principal ideals can also be
defined for noncommutative rings, but this is more complicated\footnote{Some
details:
\par
If $R$ is a noncommutative ring, then \textbf{in general} neither $uR=\left\{
ur\ \mid\ r\in R\right\}  $ nor its mirror analogue $Ru=\left\{
ru\ \mid\ r\in R\right\}  $ are ideals of $R$. (For example, $uR$ may fail the
\textquotedblleft$ab\in uR$ for any $a\in R$ and $b\in uR$\textquotedblright%
\ requirement, because there is no way to move the $u$ to the left of the
$a$.) This suggests considering the set $\left\{  rus\ \mid\ r,s\in R\right\}
$, but this is still not an ideal (in general), since it is not always closed
under addition.
\par
However, one can define the \textquotedblleft principal
ideal\textquotedblright\ $RuR$ to be%
\[
\left\{  \text{all finite sums of the form }r_{1}us_{1}+r_{2}us_{2}%
+\cdots+r_{n}us_{n}\text{ with }r_{i},s_{i}\in R\right\}  .
\]
This is always an ideal of $R$.}. However, we don't need all of $R$ to be
commutative in order for $uR$ to be an ideal; we can get by with a more local assumption:

\begin{exercise}
\label{exe.ideal.princid-central}Let $R$ be a ring (not necessarily
commutative). Let $u$ be a central element of $R$. (See Definition
\ref{def.center} \textbf{(a)} for the meaning of \textquotedblleft
central\textquotedblright.)

Prove that the set $uR$ (as defined in Proposition \ref{prop.ideal.princid})
is an ideal of $R$.
\end{exercise}

We can also generalize the ideal $uR$ of Proposition \ref{prop.ideal.princid}
\textbf{(a)} by replacing $R$ with an arbitrary ideal:

\begin{proposition}
\label{prop.ideal.uI}Let $R$ be a commutative ring. Let $u\in R$. Let $I$ be
an ideal of $R$. We define $uI$ to be the set $\left\{  ur\ \mid\ r\in
I\right\}  $. Then, $uI$ is an ideal of $R$ satisfying $uI\subseteq I$.
\end{proposition}

\begin{proof}
The claim that $uI$ is an ideal of $R$ is proved in the same way as
Proposition \ref{prop.ideal.princid} \textbf{(a)}, where we have to use the
fact that $I$ is an ideal (so that $r\in I$ entails $ar\in I$). The claim that
$uI\subseteq I$ also follows from this same fact.
\end{proof}

\subsubsection{Other examples of ideals}

In general, not all ideals of a ring need to be principal. However, in order
to find non-principal ideals, we need to venture beyond the classical number
rings $\mathbb{Z}$, $\mathbb{Q}$, $\mathbb{R}$ and $\mathbb{C}$, because the
latter rings have the property that all their ideals are principal (this will
follow from Proposition \ref{prop.eucldom.PID} further below). One way to
construct non-principal ideals is to work with polynomials in several
variables over a field, or even with univariate polynomials over $\mathbb{Z}$.
For example:

\begin{itemize}
\item The set of all polynomials $f\in\mathbb{Q}\left[  x,y\right]  $ that
have constant term $0$ is an ideal of $\mathbb{Q}\left[  x,y\right]  $ that is
not principal.

\item The set of all polynomials $f\in\mathbb{Z}\left[  x\right]  $ whose
constant term is even is an ideal of $\mathbb{Z}\left[  x\right]  $ that is
not principal.
\end{itemize}

We will come back to this later when we actually have defined polynomials.

For further examples of ideals, one might look into noncommutative rings.
However, matrix rings like $\mathbb{R}^{n\times n}$ are rather disappointing
in this regard:

\begin{exercise}
\label{exe.ideals.matrix-rings.F}Let $F$ be a field. Let $n$ be a positive
integer. Prove that the matrix ring $F^{n\times n}$ has only two ideals,
namely $\left\{  0\right\}  $ and the whole $F^{n\times n}$ (where $0$ stands
for the zero matrix). \medskip

[\textbf{Hint:} For each $i,j\in\left\{  1,2,\ldots,n\right\}  $, let
$E_{i,j}\in F^{n\times n}$ be the $\left(  i,j\right)  $\textbf{-th elementary
matrix} -- i.e., the $n\times n$-matrix whose $\left(  i,j\right)  $-th entry
is $1$ and whose all remaining entries are $0$. What happens when you multiply
a given matrix $A\in F^{n\times n}$ by $E_{i,j}$ from the left or from the
right? I.e., how can you describe the matrices $E_{i,j}A$ and $AE_{i,j}$ ?]
\end{exercise}

Considering matrix rings over a ring $R$ (instead of over a field $F$)
ameliorates the lack of ideals only slightly -- namely, to the extent that $R$
itself has interesting ideals. For example, for each integer $m$, the matrix
ring $\mathbb{Z}^{2\times2}$ has an ideal consisting of all the matrices whose
entries are divisible by $m$. More generally, any ideal of a ring $R$ yields
an ideal of the matrix ring $R^{n\times n}$:

\begin{exercise}
\label{exe.ideals.matrix-rings.R}Let $R$ be a ring. Let $n\in\mathbb{N}$.

For each subset $I$ of $R$, let $I^{n\times n}$ be the subset%
\[
\left\{  A\in R^{n\times n}\ \mid\ \text{all entries of }A\text{ belong to
}I\right\}
\]
of the matrix ring $R^{n\times n}$. Prove the following:

\begin{enumerate}
\item[\textbf{(a)}] If $I$ is an ideal of $R$, then $I^{n\times n}$ is an
ideal of the matrix ring $R^{n\times n}$.

\item[\textbf{(b)}] Any ideal of $R^{n\times n}$ has the form $I^{n\times n}$
for some ideal $I$ of $R$.
\end{enumerate}
\end{exercise}

But we can go beyond matrix rings in search of interesting ideals. Going
beyond matrix rings doesn't mean extending matrix rings; instead, it suffices
to consider some of their subrings. As we know, the upper-triangular matrices
form a subring of the matrix ring, as do the lower-triangular ones. Here, a
plethora of ideals appears. Some examples follow:

\begin{exercise}
\label{exe.ideals.triangular-matrices.1}Let $R$ be any ring. Recall that
$R^{n\leq n}$ denotes the ring of all upper-triangular $n\times n$-matrices
with entries in $R$. In particular,%
\[
R^{2\leq2}=\left\{  \left(
\begin{array}
[c]{cc}%
a & b\\
0 & d
\end{array}
\right)  \ \mid\ a,b,d\in R\right\}  .
\]

\begin{enumerate}
\item[\textbf{(a)}] Define four subsets $I,J,K,L$ of $R^{2\leq2}$ by%
\begin{align*}
I  &  :=\left\{  \left(
\begin{array}
[c]{cc}%
0 & b\\
0 & d
\end{array}
\right)  \ \mid\ b,d\in R\right\}  ;\\
J  &  :=\left\{  \left(
\begin{array}
[c]{cc}%
a & b\\
0 & 0
\end{array}
\right)  \ \mid\ a,b\in R\right\}  ;\\
K  &  :=\left\{  \left(
\begin{array}
[c]{cc}%
0 & b\\
0 & 0
\end{array}
\right)  \ \mid\ b\in R\right\}  ;\\
L  &  :=\left\{  \left(
\begin{array}
[c]{cc}%
a & 0\\
0 & d
\end{array}
\right)  \ \mid\ a,d\in R\right\}  .
\end{align*}
Prove that $I$, $J$ and $K$ are ideals of $R^{2\leq2}$, but $L$ is not (unless
$R$ is trivial).

\item[\textbf{(b)}] For any positive integer $n$, prove that the subset%
\begin{align*}
&  \left\{  A\in R^{n\times n}\ \mid\ \text{all nonzero entries of }A\text{
lie in the first row}\right\} \\
&  =\left\{  \left(
\begin{array}
[c]{cccc}%
a_{1,1} & a_{1,2} & \cdots & a_{1,n}\\
0 & 0 & \cdots & 0\\
\vdots & \vdots & \ddots & \vdots\\
0 & 0 & \cdots & 0
\end{array}
\right)  \ \mid\ a_{1,1},a_{1,2},\ldots,a_{1,n}\in R\right\}
\end{align*}
is an ideal of $R^{n\leq n}$, but the subset%
\[
\left\{  A\in R^{n\leq n}\ \mid\ \text{all nonzero entries of }A\text{ lie in
the second row}\right\}
\]
is not (for $R$ nontrivial and $n\geq2$).

\item[\textbf{(c)}] For any positive integer $n$, prove that the subset%
\begin{align*}
&  \left\{  A\in R^{n\times n}\ \mid\ \text{all nonzero entries of }A\text{
lie in the first row or the last column}\right\} \\
&  =\left\{  \left(
\begin{array}
[c]{ccccc}%
a_{1,1} & a_{1,2} & a_{1,3} & \cdots & a_{1,n}\\
0 & 0 & 0 & \cdots & a_{2,n}\\
0 & 0 & 0 & \cdots & a_{3,n}\\
\vdots & \vdots & \vdots & \ddots & \vdots\\
0 & 0 & 0 & \cdots & a_{n,n}%
\end{array}
\right)  \ \mid\ a_{1,1},a_{1,2},\ldots,a_{1,n},a_{2,n},a_{3,n},\ldots
,a_{n,n}\in R\right\}
\end{align*}
is an ideal of $R^{n\leq n}$.

\item[\textbf{(d)}] For any $n\in\mathbb{N}$, prove that the subset%
\begin{align*}
&  \ R^{n<n}\\
:=  &  \ \left\{  A\in R^{n\leq n}\ \mid\ \text{all diagonal entries of
}A\text{ equal }0\right\} \\
=  &  \ \left\{  A\in R^{n\times n}\ \mid\ \text{all nonzero entries of
}A\text{ lie above the main diagonal}\right\} \\
=  &  \ \left\{  \left(
\begin{array}
[c]{ccccc}%
0 & a_{1,2} & a_{1,3} & \cdots & a_{1,n}\\
0 & 0 & a_{2,3} & \cdots & a_{2,n}\\
0 & 0 & 0 & \cdots & a_{3,n}\\
\vdots & \vdots & \vdots & \ddots & \vdots\\
0 & 0 & 0 & \cdots & 0
\end{array}
\right)  \ \mid\ a_{i,j}\in R\text{ for all }i<j\right\}
\end{align*}
is an ideal of $R^{n\leq n}$. This subset $R^{n<n}$ is called the set of
\textbf{strictly upper-triangular} $n\times n$-matrices over $R$.

\item[\textbf{(e)}] (For combinatorialists familiar with partially ordered
sets:) Consider an $n\in\mathbb{N}$. Assume that some cells of an (unfilled)
$n\times n$-matrix are colored red. What combinatorial properties must our
coloring satisfy in order for the set%
\[
\left\{  A\in R^{n\times n}\ \mid\ \text{all nonzero entries of }A\text{ lie
in red cells}\right\}
\]
to be an ideal of $R^{n\leq n}$ ?
\end{enumerate}
\end{exercise}

The next exercise assigns a certain important ideal to every commutative ring
$R$:

\begin{exercise}
\label{exe.21hw1.7}Let $R$ be a ring. An element $a\in R$ is said to be
\textbf{nilpotent} if there exists an $n\in\mathbb{N}$ such that $a^{n}=0$.
(For example, the residue class $\overline{6}$ in $\mathbb{Z}/8\mathbb{Z}$ is
nilpotent, since its $3$-rd power is $\overline{0}$.)

\begin{enumerate}
\item[\textbf{(a)}] If $a$ and $b$ are two nilpotent elements of $R$
satisfying $ab=ba$, then prove that $a+b$ is nilpotent as well.

\item[\textbf{(b)}] Find a counterexample to part \textbf{(a)} if we don't
assume $ab=ba$.

\item[\textbf{(c)}] Assume that the ring $R$ is commutative. Let $N$ be the
set of all nilpotent elements of $R$. Prove that $N$ is an ideal of $R$.
\end{enumerate}
\end{exercise}

The ideal $N$ in Exercise \ref{exe.21hw1.7} \textbf{(c)} is known as the
\textbf{nilradical} of $R$.

\begin{exercise}
Let $n$ be a positive integer with prime factorization $n=p_{1}^{a_{1}}%
p_{2}^{a_{2}}\cdots p_{k}^{a_{k}}$, where $p_{1},p_{2},\ldots,p_{k}$ are
distinct primes and $a_{1},a_{2},\ldots,a_{k}$ are positive integers. Let $R$
be the ring $\mathbb{Z}/n$. Prove that the nilradical of $R$ is the principal
ideal $\overline{p_{1}p_{2}\cdots p_{k}}R$.
\end{exercise}

\begin{exercise}
Let $R$ be a ring. Let $N$ be the set of all nilpotent elements of $R$ (as in
Exercise \ref{exe.21hw1.7} \textbf{(c)}, but now without assuming that $R$ is
commutative). Describe this set $N$

\begin{enumerate}
\item[\textbf{(a)}] in the case when $R$ is the matrix ring $\mathbb{Q}%
^{2\times2}$;

\item[\textbf{(b)}] in the case when $R$ is the upper-triangular matrix ring
$\mathbb{Q}^{2\leq2}$ (see Section \ref{sec.rings.subrings} for its definition).
\end{enumerate}

Both of these rings $R$ are noncommutative, so that Exercise \ref{exe.21hw1.7}
\textbf{(c)} does not apply. Nevertheless, is $N$ an ideal of $R$ in one of
these two cases?
\end{exercise}

\begin{exercise}
Let $R$ be a ring. If $A,B,C,D$ are four subsets of $R$, then the notation
$\left(
\begin{array}
[c]{cc}%
A & B\\
C & D
\end{array}
\right)  $ shall denote the set of all $2\times2$-matrices $\left(
\begin{array}
[c]{cc}%
a & b\\
c & d
\end{array}
\right)  \in R^{2\times2}$ with $a\in A$, $b\in B$, $c\in C$ and $d\in D$.
(For instance, $\left(
\begin{array}
[c]{cc}%
\mathbb{N} & 2\mathbb{Z}\\
2\mathbb{Z} & \mathbb{N}%
\end{array}
\right)  $ is the set of all $2\times2$-matrices whose diagonal entries are
nonnegative integers and whose off-diagonal entries are even integers.)

\begin{enumerate}
\item[\textbf{(a)}] Let $I$ be a subset of $R$. Prove that $I$ is an ideal of
$R$ if and only if $\left(
\begin{array}
[c]{cc}%
R & I\\
\left\{  0\right\}  & R
\end{array}
\right)  $ is a subring of $R^{2\times2}$.

\item[\textbf{(b)}] Does the same claim hold for $\left(
\begin{array}
[c]{cc}%
R & I\\
I & R
\end{array}
\right)  $ instead of $\left(
\begin{array}
[c]{cc}%
R & I\\
\left\{  0\right\}  & R
\end{array}
\right)  $ ?

\item[\textbf{(c)}] Does the same claim hold for $\left(
\begin{array}
[c]{cc}%
R & I\\
R & R
\end{array}
\right)  $ instead of $\left(
\begin{array}
[c]{cc}%
R & I\\
\left\{  0\right\}  & R
\end{array}
\right)  $ ?

\item[\textbf{(d)}] Does the same claim hold for $\left(
\begin{array}
[c]{cc}%
R & R\\
R & I
\end{array}
\right)  $ instead of $\left(
\begin{array}
[c]{cc}%
R & I\\
\left\{  0\right\}  & R
\end{array}
\right)  $ ?
\end{enumerate}
\end{exercise}

The following exercise gives some examples of principal and non-principal ideals:

\begin{exercise}
Let $R$ be a ring, and let $S$ be a set. Let $R^{S}$ be the ring of all
functions from $S$ to $R$ (with pointwise addition and multiplication).

The \textbf{support} of a function $f:S\rightarrow R$ is defined to be the set
of all $x\in S$ such that $f\left(  x\right)  \neq0$. This support is denoted
by $\operatorname*{Supp}f$.

\begin{enumerate}
\item[\textbf{(a)}] Let $T$ be any subset of $S$. Prove that the set%
\[
R_{T}^{S}:=\left\{  f:S\rightarrow R\ \mid\ \operatorname*{Supp}f\subseteq
T\right\}
\]
is an ideal of $R^{S}$, and is in fact a principal ideal if $R$ is commutative.

\item[\textbf{(b)}] Prove that the set%
\[
R_{\operatorname*{fin}}^{S}:=\left\{  f:S\rightarrow R\ \mid
\ \operatorname*{Supp}f\text{ is a finite set}\right\}
\]
is an ideal of $R^{S}$.

\item[\textbf{(c)}] Now, assume that $R=\mathbb{Q}$ and $S=\mathbb{Q}$. Prove
that $R_{\operatorname*{fin}}^{S}$ is not a principal ideal of $R^{S}$.
\end{enumerate}
\end{exercise}

Another example of a non-principal ideal comes from real analysis:

\begin{exercise}
Let $R$ be the ring of all functions from $\mathbb{R}$ to $\mathbb{R}$ (with
pointwise addition and pointwise multiplication).

The \textbf{support} of a function $f:\mathbb{R}\rightarrow\mathbb{R}$ is
defined to be the set of all $x\in\mathbb{R}$ such that $f\left(  x\right)
\neq0$.

A subset $S$ of $\mathbb{R}$ is said to be
\href{https://en.wikipedia.org/wiki/Null_set}{\textbf{null}} (or to have
\textbf{Lebesgue measure zero}) if for every positive real $\varepsilon$,
there exists a countable union of intervals $I_{1}\cup I_{2}\cup I_{3}%
\cup\cdots$ in $\mathbb{R}$ such that $S\subseteq I_{1}\cup I_{2}\cup
I_{3}\cup\cdots$ and such that the sum of the lengths of these intervals
$I_{1},I_{2},I_{3},\ldots$ is smaller than $\varepsilon$. (In particular, any
finite or countable subset of $\mathbb{R}$ is null.)

We let $I$ be the set of all functions $f:\mathbb{R}\rightarrow\mathbb{R}$
whose support is null. (For example, the function that sends every rational
number to $1$ and every irrational number to $0$ belongs to $I$, since its
support is $\mathbb{Q}$, which is null.)

\begin{enumerate}
\item[\textbf{(a)}] Prove that $I$ is an ideal of $R$.

\item[\textbf{(b)}] Prove that this ideal $I$ is not principal.
\end{enumerate}
\end{exercise}

\subsection{\label{sec.rings.quots}Quotient rings (\cite[\S 7.3]{DumFoo04})}

We now come to one of the most abstract sections of this course: the
definition and the basic properties of quotient rings.

Before we define this notion rigorously, let me outline what it is meant to achieve.

Recall the idea behind modular arithmetic (Section \ref{sec.intro.mod}): By
passing from the integers to their residue classes modulo a given integer $n$,
we are essentially equating $n$ with $0$ (so that two integers become
\textquotedblleft equal\textquotedblright\ if they differ by a multiple of
$n$). Thus, these residue classes are \textquotedblleft what
remains\textquotedblright\ of the integers if we equate $n$ with $0$.

The same passage can be made in greater generality: We can start with any ring
$R$ and any ideal $I$ of $R$, and we can equate all elements of $I$ with $0$
(so that two elements of $R$ become \textquotedblleft equal\textquotedblright%
\ if they differ by an element of $I$). What remains is again called
\textquotedblleft residue classes\textquotedblright\ (now modulo $I$ instead
of modulo $n$), and we can again define addition and multiplication on these
residue classes. The result is a new ring, which is called the
\textbf{quotient ring} of $R$ by the ideal $I$, and is denoted by $R/I$.
Working in this quotient ring is a natural generalization of modular
arithmetic to things that aren't integers. For instance, we can start with the
Gaussian integers and equate $5$ with $0$, or we can start with the
upper-triangular $2\times2$-matrices\footnote{i.e., the matrices of the form
$\left(
\begin{array}
[c]{cc}%
a & b\\
0 & c
\end{array}
\right)  $} and equate the strictly upper-triangular $2\times2$%
-matrices\footnote{i.e., the matrices of the form $\left(
\begin{array}
[c]{cc}%
0 & b\\
0 & 0
\end{array}
\right)  $} with zero. This gives us a new way to build new rings from
old.\footnote{For comparison: When you take a subring of a ring $R$, you are
throwing away some elements of $R$. In contrast, when you take a quotient ring
of $R$, you are equating some elements of $R$ with one another. In either
case, you end up with a ring smaller than $R$.}

So much for the idea; let us now define the quotient ring $R/I$ formally.

In rigorous mathematics, you cannot just take two distinct elements and
declare them to be equal. Thus, \textquotedblleft equating\textquotedblright%
\ two elements of $R$ is easier said than done. The right way to do it is by
passing from elements to equivalence classes (just as we did in modular
arithmetic, back in Section \ref{sec.intro.mod}). Let us see how this can be done.

\subsubsection{Quotient groups}

It turns out that we don't need to reinvent the wheel: You have already seen
these equivalence classes in group theory, under the name \textquotedblleft%
\textbf{cosets}\textquotedblright. Let me recall how these were defined and used:

\begin{itemize}
\item If $H$ is a subgroup of a group $G$, then the \textbf{left cosets} of
$H$ in $G$ are the subsets $gH:=\left\{  gh\ \mid\ h\in H\right\}  $ for all
$g\in G$. There is one left coset $gH$ for each $g\in G$; but different $g\in
G$ often lead to the same left coset $gH$, so there are usually fewer left
cosets than elements of $G$. The set of all left cosets of $H$ is denoted by
$G/H$.

\item If $H$ is merely a subgroup of a group $G$, then $G/H$ is merely a
\textquotedblleft$G$-set\textquotedblright\ (i.e., a set with an action of
$G$). However, when $H$ is a \textbf{normal} subgroup of $G$ (that is, a
subgroup of $G$ satisfying $gng^{-1}\in H$ for each $g\in G$ and $n\in H$),
then $G/H$ becomes a \textbf{group} as well, with group operation defined by%
\begin{equation}
\left(  g_{1}H\right)  \left(  g_{2}H\right)  =g_{1}g_{2}%
H\ \ \ \ \ \ \ \ \ \ \text{for all }g_{1},g_{2}\in G.
\label{eq.qgroups.g1Hg2H}%
\end{equation}
This group $G/H$ is called the \textbf{quotient group} of $G$ by $H$. The left
cosets of $H$ in $G$ are just called the \textbf{cosets} of $H$ in $G$ in this case.

\item If $G$ is an abelian group, then any subgroup $H$ of $G$ is normal, so
$G/H$ always is a group.

\item Now, assume that $G$ is an \textbf{additive} group (which means that its
binary operation is written as $+$ rather than as $\cdot$). This presupposes
that $G$ is abelian, as it is considered gauche to write a non-abelian group
additively. Let $H$ be a subgroup of $G$. Then, the cosets of $H$ in $G$ are
denoted by $g+H$ instead of $gH$ (in order to match the additive notation for
the group operation). Likewise, we write $+$ instead of $\cdot$ for the binary
operation of the quotient group $G/H$. The equality (\ref{eq.qgroups.g1Hg2H})
therefore rewrites as%
\[
\left(  g_{1}+H\right)  +\left(  g_{2}+H\right)  =\left(  g_{1}+g_{2}\right)
+H\ \ \ \ \ \ \ \ \ \ \text{for all }g_{1},g_{2}\in G.
\]
Note that the quotient group $G/H$ is an abelian group.

\item The most famous example of quotient groups is when $G=\mathbb{Z}$ and
$H=n\mathbb{Z}=\left\{  \text{all multiples of }n\right\}  $ for some fixed
integer $n$. (Here, the group operation on $G$ is addition of integers.) In
this case, the quotient group $\mathbb{Z}/n\mathbb{Z}$ is the cyclic group
$\mathbb{Z}/n$, also known as $Z_{n}$. See \cite[Chapter XII]{Siksek20} for
this and other examples.
\end{itemize}

\subsubsection{Quotient rings}

Now, piggybacking on the construction of quotient groups we just recalled, we
shall define a similar quotient structure for rings instead of groups. Instead
of normal subgroups, we will use ideals this time:

\begin{definition}
\label{def.quotring}Let $I$ be an ideal of a ring $R$. Thus, $I$ is a subgroup
of the additive group $\left(  R,+,0\right)  $, hence a normal subgroup (since
$\left(  R,+,0\right)  $ is abelian). Therefore, the quotient group $R/I$ is a
well-defined abelian group. Its elements are the cosets $r+I$ of $I$ in $R$.
(Note that, since our groups are additive, we are writing $r+I$ for what would
normally be written $rI$ in group theory.)

Note that the addition on $R/I$ is given by
\begin{equation}
\left(  a+I\right)  +\left(  b+I\right)  =\left(  a+b\right)
+I\ \ \ \ \ \ \ \ \ \ \text{for all }a,b\in R. \label{eq.def.quotring.+}%
\end{equation}

We now define a multiplication operation on $R/I$ by setting%
\begin{equation}
\left(  a+I\right)  \left(  b+I\right)  =ab+I\ \ \ \ \ \ \ \ \ \ \text{for all
}a,b\in R. \label{eq.def.quotring.*}%
\end{equation}
(See below for a proof that this is well-defined.)

The set $R/I$, equipped with the addition and the multiplication we just
defined and with the elements $0+I$ and $1+I$ (as zero and unity), is a ring
(as we will show in a moment). This ring is called the \textbf{quotient ring}
of $R$ by the ideal $I$; it is also pronounced \textquotedblleft$R$
\textbf{modulo} $I$\textquotedblright. It is denoted $R/I$ (so when you hear
\textquotedblleft the ring $R/I$\textquotedblright, it always means the set
$R/I$ equipped with the structure just mentioned).

The cosets $r+I$ are called \textbf{residue classes} modulo $I$, and are often
denoted $r\operatorname{mod}I$ or $\left[  r\right]  _{I}$ or $\left[
r\right]  $ or $\overline{r}$. (The last two notations are used when $I$ is
clear from the context. We will mostly be using the notations $r+I$ and
$\overline{r}$.) Thus, the equalities (\ref{eq.def.quotring.+}) and
(\ref{eq.def.quotring.*}) can be restated as
\begin{equation}
\overline{a}+\overline{b}=\overline{a+b}\ \ \ \ \ \ \ \ \ \ \text{for all
}a,b\in R \label{eq.def.quotring.+bar}%
\end{equation}
and%
\begin{equation}
\overline{a}\cdot\overline{b}=\overline{ab}\ \ \ \ \ \ \ \ \ \ \text{for all
}a,b\in R, \label{eq.def.quotring.*bar}%
\end{equation}
respectively.
\end{definition}

\begin{theorem}
\label{thm.quotring.welldef}Let $R$ and $I$ be as in Definition
\ref{def.quotring}. Then, the multiplication on $R/I$ is well-defined, and
$R/I$ does indeed become a ring when endowed with the operations and elements
just described.
\end{theorem}

Before we prove this theorem, let us see some examples:

\begin{itemize}
\item Let $n\in\mathbb{Z}$. The set $n\mathbb{Z}=\left\{  \text{all multiples
of }n\right\}  $ is a principal ideal of $\mathbb{Z}$. The quotient ring
$\mathbb{Z}/n\mathbb{Z}$ is precisely the ring $\mathbb{Z}/n$ we discussed
above. Its elements $r+n\mathbb{Z}$ are precisely the residue classes
$\overline{r}$ defined in Section \ref{sec.intro.mod}. The equalities
(\ref{eq.def.quotring.+bar}) and (\ref{eq.def.quotring.*bar}) are precisely
the standard definitions of addition and multiplication in $\mathbb{Z}/n$.

Thus, the notion of a quotient ring generalizes the familiar concept of
modular arithmetic. (More precisely, modular arithmetic is arithmetic in $R/I$
where $R=\mathbb{Z}$ and $I=n\mathbb{Z}$.)

\item Each ring $R$ has two obvious ideals $\left\{  0_{R}\right\}  $ and $R$.
The corresponding quotient rings $R/\left\{  0_{R}\right\}  $ and $R/R$ are
fairly boring:

\begin{itemize}
\item The quotient ring $R/\left\{  0_{R}\right\}  $ is isomorphic to $R$
(since each coset $r+\left\{  0_{R}\right\}  $ is just a $1$-element set
$\left\{  r\right\}  $, and thus the elements of $R/\left\{  0_{R}\right\}  $
are just the elements of $R$ \textquotedblleft clothed in set
braces\textquotedblright).

\item The quotient ring $R/R$ is trivial (since there is only one coset,
$r+R=0+R=1+R=R$, and it contains all elements of $R$).
\end{itemize}

This generalizes the facts that $\mathbb{Z}/0\cong\mathbb{Z}$ and that
$\mathbb{Z}/1$ is trivial.

\item Let $R$ be the ring $\mathbb{Z}\left[  i\right]  =\left\{
a+bi\ \mid\ a,b\in\mathbb{Z}\right\}  $ of Gaussian integers. Consider its
principal ideal
\begin{align*}
3R  &  =\left\{  3r\ \mid\ r\in R\right\}  =\left\{  3r\ \mid\ r\in
\mathbb{Z}\left[  i\right]  \right\} \\
&  =\left\{  3a+3bi\ \mid\ a,b\in\mathbb{Z}\right\} \\
&  =\left\{  c+di\ \mid\ c,d\in\mathbb{Z}\text{ are multiples of }3\right\}  .
\end{align*}
What is the quotient ring $R/\left(  3R\right)  $ ? Each element of this
quotient ring can be written in the form\footnote{We are using $\overline{z}$
to denote the residue class of a Gaussian integer $z\in R$. This should not be
confused with the complex conjugate of $z$ (which is commonly denoted
$\overline{z}$ as well). Fortunately, this confusion will be avoided in this
example, since we will not use complex conjugates.}%
\[
\overline{a+bi}\ \ \ \ \ \ \ \ \ \ \text{with }a,b\in\left\{  0,1,2\right\}
\]
(since any Gaussian integer can be reduced to an $a+bi$ with $a,b\in\left\{
0,1,2\right\}  $ by subtracting an appropriate Gaussian-integer multiple of
$3$\ \ \ \ \footnote{\textit{Proof.} Let $c+di$ be any Gaussian integer (with
$c,d\in\mathbb{Z}$). Let $q_{c}$ and $r_{c}$ be the quotient and the remainder
obtained when we divide $c$ by $3$. Let $q_{d}$ and $r_{d}$ be the quotient
and the remainder obtained when we divide $d$ by $3$. Then, $c=3q_{c}+r_{c}$
and $r_{c}\in\left\{  0,1,2\right\}  $ and $d=3q_{d}+r_{d}$ and $r_{d}%
\in\left\{  0,1,2\right\}  $. Hence,%
\[
c+di=\left(  3q_{c}+r_{c}\right)  +\left(  3q_{d}+r_{d}\right)  i=3\left(
q_{c}+q_{d}i\right)  +\left(  r_{c}+r_{d}i\right)  .
\]
Thus, by subtracting $3\left(  q_{c}+q_{d}i\right)  $ (which is a
Gaussian-integer multiple of $3$), we can reduce $c+di$ to the Gaussian
integer $r_{c}+r_{d}i$, which has the form $a+bi$ with $a,b\in\left\{
0,1,2\right\}  $.
\par
For example, $5+7i$ can be reduced to $2+i$ by subtracting $3\left(
1+2i\right)  $. Thus, $\overline{5+7i}=\overline{2+i}$.}). In other words,%
\[
R/\left(  3R\right)  =\left\{  \overline{0},\ \ \overline{1},\ \ \overline
{2},\ \ \overline{i},\ \ \overline{1+i},\ \ \overline{2+i},\ \ \overline
{2i},\ \ \overline{1+2i},\ \ \overline{2+2i}\right\}  .
\]
It is easy to see that this is a $9$-element ring (i.e., the residue
classes\newline$\overline{0},\ \ \overline{1},\ \ \overline{2},\ \ \overline
{i},\ \ \overline{1+i},\ \ \overline{2+i},\ \ \overline{2i},\ \ \overline
{1+2i},\ \ \overline{2+2i}$ are distinct\footnote{In order to verify this, you
must show that no two of the Gaussian integers
$0,\ \ 1,\ \ 2,\ \ i,\ \ 1+i,\ \ 2+i,\ \ 2i,\ \ 1+2i,\ \ 2+2i$ differ by a
Gaussian-integer multiple of $3$ (so that their residue classes are
distinct).}), and a field (i.e., all the nonzero residue classes are
invertible\footnote{Checking this is Exercise \ref{exe.quotring.Zi/n}
\textbf{(a)}.}). So we have found a finite field with $9$ elements.

Let us do some computations in this field: We have%
\[
\overline{2+i}+\overline{2+2i}=\overline{\left(  2+i\right)  +\left(
2+2i\right)  }=\overline{4+3i}=\overline{1}%
\]
(since $4+3i$ and $1$ belong to the same coset of the ideal $3R$, as the
difference $\left(  4+3i\right)  -1=3\left(  1+i\right)  $ lies in this
ideal). We also have%
\[
\overline{2+i}\cdot\overline{2+2i}=\overline{\left(  2+i\right)  \cdot\left(
2+2i\right)  }=\overline{2+6i}=\overline{2}%
\]
(since $2+6i$ and $2$ belong to the same coset of the ideal $3R$). A similar
computation proves that%
\[
\overline{2+i}\cdot\overline{1+i}=\overline{1},
\]
which reveals that the elements $\overline{2+i}$ and $\overline{1+i}$ of the
ring $R/\left(  3R\right)  $ are inverses of each other (and thus are units of
this ring).

For the curious: If we replace $3$ by any other positive integer $n$, then
$R/\left(  nR\right)  $ will be a finite ring with $n^{2}$ elements, but not
always a field. Understanding when it will be a field is a fruitful question
in elementary number theory. (It is a field for some, but not for all, primes
$n$.)

We can also consider quotient rings of the form $R/\left(  zR\right)  $ for
non-real $z\in R$. For example, one such quotient ring is $R/\left(  \left(
1+i\right)  R\right)  $. It is much less obvious how many elements this
quotient ring has! (See Exercise \ref{exe.quotring.1+i} below for the answer.)
\end{itemize}

\begin{exercise}
\label{exe.quotring.Zi/n}Let $R=\mathbb{Z}\left[  i\right]  $, as in the
example we just did.

\begin{enumerate}
\item[\textbf{(a)}] Confirm that the quotient ring $R/\left(  3R\right)  $ is
a field by finding inverses for all its eight nonzero elements.

\item[\textbf{(b)}] Confirm that the quotient ring $R/\left(  5R\right)  $ is
not a field by checking that $\overline{1+2i}\cdot\overline{1-2i}=\overline
{0}$ in this ring.

\item[\textbf{(c)}] Confirm that the quotient ring $R/\left(  17R\right)  $ is
not a field either.
\end{enumerate}

[\textbf{Hint:} For part \textbf{(c)}, find a similar equality as in part
\textbf{(b)}.]
\end{exercise}

\begin{exercise}
\label{exe.quotring.1+i}Let $R=\mathbb{Z}\left[  i\right]  $, as in the above
example. Prove that the quotient ring $R/\left(  \left(  1+i\right)  R\right)
$ has just $2$ elements, and in fact is isomorphic to $\mathbb{Z}/2$. \medskip

[\textbf{Hint:} First, show that both $2$ and $2i$ belong to the principal
ideal $\left(  1+i\right)  R$. Hence, each element of $R$ can be reduced to
the form $a+bi$ with $a,b\in\left\{  0,1\right\}  $ by subtracting an
appropriate element of this ideal. Thus, the only possible residue classes in
$R/\left(  \left(  1+i\right)  R\right)  $ are $\overline{0}$, $\overline{1}$,
$\overline{i}$ and $\overline{1+i}$. But the difference $i-1$ also lies in the
principal ideal $\left(  1+i\right)  R$ (why?), and thus the classes
$\overline{i}$ and $\overline{1}$ are actually identical. So are the classes
$\overline{1+i}$ and $\overline{0}$. We are left with the two classes
$\overline{0}$ and $\overline{1}$, which may or may not be actually distinct.
Prove that they are distinct by arguing that $1$ is not a multiple of $1+i$ in
$R$. This shows that $R/\left(  \left(  1+i\right)  R\right)  $ has exactly
$2$ elements.]
\end{exercise}

\begin{exercise}
Let $R$ be the ring of all real numbers of the form $a+b\sqrt{5}$ with
$a,b\in\mathbb{Z}$. (This is a ring for the same reasons that the $\mathbb{S}$
in Subsection \ref{subsec.rings.def.exas} was a ring; the only difference is
that we now require $a,b\in\mathbb{Z}$ instead of $a,b\in\mathbb{Q}$.)

\begin{enumerate}
\item[\textbf{(a)}] Is the quotient ring $R/\left(  2R\right)  $ a field?

\item[\textbf{(b)}] Is the quotient ring $R/\left(  3R\right)  $ a field?

\item[\textbf{(c)}] Prove that the quotient ring $R/\left(  5R\right)  $ is
not a field, and in fact the residue class $\overline{\sqrt{5}}$ in this
quotient is nilpotent. (See Exercise \ref{exe.21hw1.7} for the meaning of
\textquotedblleft nilpotent\textquotedblright.)
\end{enumerate}
\end{exercise}

We will see more examples of quotient rings soon (and even more in later
chapters, after we introduce polynomial rings). For now, however, let us make
good on our debts and prove Theorem \ref{thm.quotring.welldef}:

\begin{proof}
[Proof of Theorem \ref{thm.quotring.welldef}.]To see that the multiplication
on $R/I$ is well-defined (by the equation (\ref{eq.def.quotring.*})), we must
prove that a product $xy$ with $x,y\in R/I$ does not depend on how exactly we
write $x$ and $y$ as $x=a+I$ and $y=b+I$. In other words, we must show that if
four elements $a,a^{\prime},b,b^{\prime}$ of $R$ satisfy $a+I=a^{\prime}+I$
and $b+I=b^{\prime}+I$, then $ab+I=a^{\prime}b^{\prime}+I$.

So let $a,a^{\prime},b,b^{\prime}\in R$ be such that $a+I=a^{\prime}+I$ and
$b+I=b^{\prime}+I$. From $a+I=a^{\prime}+I$, we obtain $a-a^{\prime}\in I$, so
that $\left(  a-a^{\prime}\right)  b\in I$ (by the second ideal axiom, since
$I$ is an ideal). In other words, $ab-a^{\prime}b\in I$. Hence,
$ab+I=a^{\prime}b+I$. Similarly, we can obtain $a^{\prime}b+I=a^{\prime
}b^{\prime}+I$ (from $b+I=b^{\prime}+I$). Thus, $ab+I=a^{\prime}b+I=a^{\prime
}b^{\prime}+I$, which is just what we need.

So we have shown that the multiplication on $R/I$ is well-defined. Now why is
$R/I$ a ring? This we leave to the reader -- it's a straightforward
consequence of the fact that $R$ is a ring.\footnote{For example, in order to
prove that the multiplication on $R/I$ is associative, we must show that
$x\cdot\left(  y\cdot z\right)  =\left(  x\cdot y\right)  \cdot z$ for any
three elements $x,y,z\in R/I$. But this is straightforward: If we write these
three elements $x,y,z$ as $x=a+I$ and $y=b+I$ and $z=c+I$, then this boils
down to proving that $a\cdot\left(  b\cdot c\right)  =\left(  a\cdot b\right)
\cdot c$, which follows from associativity in $R$.}
\end{proof}

Let me mention some more terminology (some of it informal but fairly popular):

When $R$ is a ring, and $I$ is an ideal of $R$, the quotient ring $R/I$ is
often just called the \textbf{quotient} of $R$ by $I$. It is said to be
obtained by \textbf{quotienting} $R$ by $I$, or by \textbf{quotienting }$I$
\textbf{out of }$R$.

\subsubsection{\label{subsec.rings.quotring.moreexas}More examples of quotient
rings}

Let us give two further hands-on examples of quotient rings. These are not
strictly necessary for the understanding of what follows, but they give some
extra intuition and practice.

\begin{itemize}
\item We recall (from Subsection \ref{subsec.rings.subrings.exas}) that
$\mathbb{Q}^{3\leq3}$ denotes the ring of all upper-triangular $3\times
3$-matrices with rational entries. Thus,%
\[
\mathbb{Q}^{3\leq3}=\left\{  \left(
\begin{array}
[c]{ccc}%
a & b & c\\
0 & d & e\\
0 & 0 & f
\end{array}
\right)  \ \mid\ a,b,c,d,e,f\in\mathbb{Q}\right\}  .
\]
The addition and the multiplication of this ring are matrix addition and
matrix multiplication; its zero is the zero matrix $0_{3\times3}$; its unity
is the identity matrix $I_{3}$.

There is also a ring $\mathbb{Q}^{3=3}$ of all diagonal $3\times3$-matrices
with rational entries. That is,%
\[
\mathbb{Q}^{3=3}=\left\{  \left(
\begin{array}
[c]{ccc}%
a & 0 & 0\\
0 & d & 0\\
0 & 0 & f
\end{array}
\right)  \ \mid\ a,d,f\in\mathbb{Q}\right\}  .
\]
This is a subring of $\mathbb{Q}^{3\leq3}$.

As in Exercise \ref{exe.ideals.triangular-matrices.1} \textbf{(d)}, we
furthermore define $\mathbb{Q}^{3<3}$ to be the set of all matrices in
$\mathbb{Q}^{3\leq3}$ whose diagonal entries are $0$. Thus,%
\[
\mathbb{Q}^{3<3}=\left\{  \left(
\begin{array}
[c]{ccc}%
0 & b & c\\
0 & 0 & e\\
0 & 0 & 0
\end{array}
\right)  \ \mid\ b,c,e\in\mathbb{Q}\right\}  .
\]
The matrices in $\mathbb{Q}^{3<3}$ are known as \textquotedblleft strictly
upper-triangular $3\times3$-matrices\textquotedblright.

We know from Exercise \ref{exe.ideals.triangular-matrices.1} \textbf{(d)} that
$\mathbb{Q}^{3<3}$ is an ideal of $\mathbb{Q}^{3\leq3}$ (and of course, we can
also check this directly\footnote{The only interesting part is to check the
second ideal axiom, i.e., to show that if $A\in\mathbb{Q}^{3<3}$ and
$B\in\mathbb{Q}^{3\leq3}$, then $AB$ and $BA$ belong to $\mathbb{Q}^{3<3}$.
Still, we can do this by direct computation:%
\begin{align*}
\text{If }A  &  =\left(
\begin{array}
[c]{ccc}%
0 & x & y\\
0 & 0 & z\\
0 & 0 & 0
\end{array}
\right)  \text{ and }B=\left(
\begin{array}
[c]{ccc}%
a & b & c\\
0 & d & e\\
0 & 0 & f
\end{array}
\right)  \text{,}\\
\text{then }AB  &  =\left(
\begin{array}
[c]{ccc}%
0 & dx & ex+fy\\
0 & 0 & fz\\
0 & 0 & 0
\end{array}
\right)  \text{ and }BA=\left(
\begin{array}
[c]{ccc}%
0 & ax & ay+bz\\
0 & 0 & dz\\
0 & 0 & 0
\end{array}
\right)  .
\end{align*}
}). What is the quotient ring $\mathbb{Q}^{3\leq3}/\mathbb{Q}^{3<3}$ ?

Any element of this quotient ring is a residue class of the form%
\[
\overline{A}=A+\mathbb{Q}^{3<3}\ \ \ \ \ \ \ \ \ \ \text{for some }%
A\in\mathbb{Q}^{3\leq3}.
\]
In other words, it is a set that consists of some given matrix $A\in
\mathbb{Q}^{3\leq3}$ and all matrices that can be obtained from $A$ by adding
a strictly upper-triangular $3\times3$-matrix. For example, if $A=\left(
\begin{array}
[c]{ccc}%
1 & 2 & 3\\
0 & 4 & 5\\
0 & 0 & 6
\end{array}
\right)  $, then%
\begin{align*}
\overline{A}  &  =A+\mathbb{Q}^{3<3}=\left(
\begin{array}
[c]{ccc}%
1 & 2 & 3\\
0 & 4 & 5\\
0 & 0 & 6
\end{array}
\right)  +\left\{  \left(
\begin{array}
[c]{ccc}%
0 & b & c\\
0 & 0 & e\\
0 & 0 & 0
\end{array}
\right)  \ \mid\ b,c,e\in\mathbb{Q}\right\} \\
&  =\left\{  \left(
\begin{array}
[c]{ccc}%
1 & 2 & 3\\
0 & 4 & 5\\
0 & 0 & 6
\end{array}
\right)  +\left(
\begin{array}
[c]{ccc}%
0 & b & c\\
0 & 0 & e\\
0 & 0 & 0
\end{array}
\right)  \ \mid\ b,c,e\in\mathbb{Q}\right\} \\
&  =\left\{  \left(
\begin{array}
[c]{ccc}%
1 & 2+b & 3+c\\
0 & 4 & 5+e\\
0 & 0 & 6
\end{array}
\right)  \ \mid\ b,c,e\in\mathbb{Q}\right\} \\
&  =\left\{  \left(
\begin{array}
[c]{ccc}%
1 & x & y\\
0 & 4 & z\\
0 & 0 & 6
\end{array}
\right)  \ \mid\ x,y,z\in\mathbb{Q}\right\}
\end{align*}
(here, we have substituted $x,y,z$ for $2+b,3+c,5+e$, because when $b$ ranges
over $\mathbb{Q}$, so does $2+b$, etc.). So this set $\overline{A}$ consists
of all upper-triangular $3\times3$-matrices whose diagonal entries are $1,4,6$
and whose above-diagonal entries are arbitrary rational numbers. We can thus
view this set $\overline{A}$ as a \textquotedblleft partly undetermined
matrix\textquotedblright, in the sense that it is \textquotedblleft a matrix
in which some of the entries can be filled in arbitrarily\textquotedblright%
\ (although, of course, formally speaking, it is not a matrix but a set of
matrices). From this point of view, it makes sense to write $\overline{A}$ as
follows:%
\[
\overline{A}=\left(
\begin{array}
[c]{ccc}%
1 & ? & ?\\
0 & 4 & ?\\
0 & 0 & 6
\end{array}
\right)  ,
\]
where each question mark stands for an undetermined entry (noting that
different question marks are independent, i.e., there are three degrees of
freedom). Formally, such a \textquotedblleft partly undetermined
matrix\textquotedblright\ is meant to be a set of matrices, where each
question mark is a variable that can take any element of $\mathbb{Q}$ as value.

More generally, for any matrix $A=\left(
\begin{array}
[c]{ccc}%
a & b & c\\
0 & d & e\\
0 & 0 & f
\end{array}
\right)  \in\mathbb{Q}^{3\leq3}$, the residue class $\overline{A}\in
\mathbb{Q}^{3\leq3}/\mathbb{Q}^{3<3}$ is%
\begin{equation}
\overline{A}=\overline{\left(
\begin{array}
[c]{ccc}%
a & b & c\\
0 & d & e\\
0 & 0 & f
\end{array}
\right)  }=\left(
\begin{array}
[c]{ccc}%
a & ? & ?\\
0 & d & ?\\
0 & 0 & f
\end{array}
\right)  \label{eq.quotring.R33.1.2}%
\end{equation}
(written as a \textquotedblleft partly undetermined matrix\textquotedblright).
Thus, this class $\overline{A}$ does not depend on the above-diagonal entries
of $A$. (This is not surprising: After all, when we quotient out the ideal
$\mathbb{Q}^{3<3}$, we are equating the strictly upper-triangular matrices
with $0$, which amounts to ignoring the above-diagonal entries.)

Thus, the quotient ring $\mathbb{Q}^{3\leq3}/\mathbb{Q}^{3<3}$ is the set of
\textquotedblleft partly undetermined matrices\textquotedblright\ of the form
$\left(
\begin{array}
[c]{ccc}%
a & ? & ?\\
0 & d & ?\\
0 & 0 & f
\end{array}
\right)  $ (that is, upper-triangular $3\times3$-matrices with fixed entries
on the main diagonal and question marks above it).

According to the formula (\ref{eq.def.quotring.*bar}), the multiplication on
this quotient ring is given by%
\begin{align*}
\overline{\left(
\begin{array}
[c]{ccc}%
a & b & c\\
0 & d & e\\
0 & 0 & f
\end{array}
\right)  }\cdot\overline{\left(
\begin{array}
[c]{ccc}%
x & y & z\\
0 & u & v\\
0 & 0 & w
\end{array}
\right)  }  &  =\overline{\left(
\begin{array}
[c]{ccc}%
a & b & c\\
0 & d & e\\
0 & 0 & f
\end{array}
\right)  \cdot\left(
\begin{array}
[c]{ccc}%
x & y & z\\
0 & u & v\\
0 & 0 & w
\end{array}
\right)  }\\
&  =\overline{\left(
\begin{array}
[c]{ccc}%
ax & bu+ay & az+bv+cw\\
0 & du & dv+ew\\
0 & 0 & fw
\end{array}
\right)  }.
\end{align*}
In terms of \textquotedblleft partly undetermined matrices\textquotedblright,
this can be rewritten in the following simpler form:%
\begin{equation}
\left(
\begin{array}
[c]{ccc}%
a & ? & ?\\
0 & d & ?\\
0 & 0 & f
\end{array}
\right)  \cdot\left(
\begin{array}
[c]{ccc}%
x & ? & ?\\
0 & u & ?\\
0 & 0 & w
\end{array}
\right)  =\left(
\begin{array}
[c]{ccc}%
ax & ? & ?\\
0 & du & ?\\
0 & 0 & fw
\end{array}
\right)  . \label{eq.quotring.R33.1.3}%
\end{equation}
(As we said, the above-diagonal entries don't matter, so we don't need to even
bother computing them.)

Similarly to (\ref{eq.quotring.R33.1.3}), addition in $\mathbb{Q}^{3\leq
3}/\mathbb{Q}^{3<3}$ is given by the formula
\begin{equation}
\left(
\begin{array}
[c]{ccc}%
a & ? & ?\\
0 & d & ?\\
0 & 0 & f
\end{array}
\right)  +\left(
\begin{array}
[c]{ccc}%
x & ? & ?\\
0 & u & ?\\
0 & 0 & w
\end{array}
\right)  =\left(
\begin{array}
[c]{ccc}%
a+x & ? & ?\\
0 & d+u & ?\\
0 & 0 & f+w
\end{array}
\right)  . \label{eq.quotring.R33.1.4}%
\end{equation}

I hope you are disappointed by the formulas (\ref{eq.quotring.R33.1.4}) and
(\ref{eq.quotring.R33.1.3}). After all, what is happening in these formulas is
just entrywise addition and entrywise multiplication of the diagonal entries.
In other words, the \textquotedblleft partly undetermined
matrices\textquotedblright\ $\left(
\begin{array}
[c]{ccc}%
a & ? & ?\\
0 & d & ?\\
0 & 0 & f
\end{array}
\right)  $ in our quotient ring $\mathbb{Q}^{3\leq3}/\mathbb{Q}^{3<3}$ are
behaving (under addition and multiplication) just like the diagonal matrices
$\left(
\begin{array}
[c]{ccc}%
a & 0 & 0\\
0 & d & 0\\
0 & 0 & f
\end{array}
\right)  $ in the ring $\mathbb{Q}^{3=3}$ (since the latter diagonal matrices,
too, are added and multiplied entrywise\footnote{Indeed, addition and
multiplication of diagonal matrices are given by the formulas%
\[
\left(
\begin{array}
[c]{ccc}%
a & 0 & 0\\
0 & d & 0\\
0 & 0 & f
\end{array}
\right)  +\left(
\begin{array}
[c]{ccc}%
x & 0 & 0\\
0 & u & 0\\
0 & 0 & w
\end{array}
\right)  =\left(
\begin{array}
[c]{ccc}%
a+x & 0 & 0\\
0 & d+u & 0\\
0 & 0 & f+w
\end{array}
\right)
\]
and%
\[
\left(
\begin{array}
[c]{ccc}%
a & 0 & 0\\
0 & d & 0\\
0 & 0 & f
\end{array}
\right)  \cdot\left(
\begin{array}
[c]{ccc}%
x & 0 & 0\\
0 & u & 0\\
0 & 0 & w
\end{array}
\right)  =\left(
\begin{array}
[c]{ccc}%
ax & 0 & 0\\
0 & du & 0\\
0 & 0 & fw
\end{array}
\right)  .
\]
These look exactly like the formulas (\ref{eq.quotring.R33.1.4}) and
(\ref{eq.quotring.R33.1.3}) for our \textquotedblleft partly undetermined
matrices\textquotedblright, except that all question marks are replaced by
zeroes.}). To state this more precisely, our quotient ring $\mathbb{Q}%
^{3\leq3}/\mathbb{Q}^{3<3}$ turns out to be isomorphic to the ring
$\mathbb{Q}^{3=3}$ of diagonal $3\times3$-matrices, and the isomorphism is
just the map%
\begin{align*}
\mathbb{Q}^{3\leq3}/\mathbb{Q}^{3<3}  &  \rightarrow\mathbb{Q}^{3=3},\\
\left(
\begin{array}
[c]{ccc}%
a & ? & ?\\
0 & d & ?\\
0 & 0 & f
\end{array}
\right)   &  \mapsto\left(
\begin{array}
[c]{ccc}%
a & 0 & 0\\
0 & d & 0\\
0 & 0 & f
\end{array}
\right)
\end{align*}
that replaces all question marks by zeroes. Thus, our \textquotedblleft partly
undetermined matrices\textquotedblright\ are just diagonal matrices in a
complicated guise, and our quotient ring $\mathbb{Q}^{3\leq3}/\mathbb{Q}%
^{3<3}$ is just an isomorphic copy of the subring $\mathbb{Q}^{3=3}$. This
doesn't feel worth the trouble of defining quotient rings!

If all quotient rings were as boring as this one, then the whole concept
wouldn't be of much use. Fortunately, this is not the case: Not all quotient
rings are subrings in disguise, and not all question marks can just be
replaced by zeroes. We will see this in the next example.

\item We again consider a quotient ring of $\mathbb{Q}^{3\leq3}$, but this
time we quotient out a smaller ideal. Namely, we define the subset%
\[
\mathbb{Q}^{3<<3}:=\left\{  \left(
\begin{array}
[c]{ccc}%
0 & 0 & c\\
0 & 0 & 0\\
0 & 0 & 0
\end{array}
\right)  \ \mid\ c\in\mathbb{Q}\right\}
\]
of $\mathbb{Q}^{3\leq3}$. This set $\mathbb{Q}^{3<<3}$ is an ideal of
$\mathbb{Q}^{3\leq3}$ (again, this can be checked directly\footnote{The only
interesting part is to check the second ideal axiom, i.e., to show that if
$A\in\mathbb{Q}^{3<<3}$ and $B\in\mathbb{Q}^{3\leq3}$, then $AB$ and $BA$
belong to $\mathbb{Q}^{3<<3}$. Still, we can do this by direct computation:%
\begin{align*}
\text{If }A  &  =\left(
\begin{array}
[c]{ccc}%
0 & 0 & y\\
0 & 0 & 0\\
0 & 0 & 0
\end{array}
\right)  \text{ and }B=\left(
\begin{array}
[c]{ccc}%
a & b & c\\
0 & d & e\\
0 & 0 & f
\end{array}
\right)  \text{,}\\
\text{then }AB  &  =\left(
\begin{array}
[c]{ccc}%
0 & 0 & fy\\
0 & 0 & 0\\
0 & 0 & 0
\end{array}
\right)  \text{ and }BA=\left(
\begin{array}
[c]{ccc}%
0 & 0 & ay\\
0 & 0 & 0\\
0 & 0 & 0
\end{array}
\right)  .
\end{align*}
}), and thus there is a quotient ring $\mathbb{Q}^{3\leq3}/\mathbb{Q}^{3<<3}$.
In this quotient ring, the residue class $\overline{A}$ of a matrix $A=\left(
\begin{array}
[c]{ccc}%
a & b & c\\
0 & d & e\\
0 & 0 & f
\end{array}
\right)  \in\mathbb{Q}^{3\leq3}$ is%
\begin{align*}
\overline{A}  &  =\overline{\left(
\begin{array}
[c]{ccc}%
a & b & c\\
0 & d & e\\
0 & 0 & f
\end{array}
\right)  }=\left(
\begin{array}
[c]{ccc}%
a & b & ?\\
0 & d & e\\
0 & 0 & f
\end{array}
\right) \\
&  =\left\{  \left(
\begin{array}
[c]{ccc}%
a & b & x\\
0 & d & e\\
0 & 0 & f
\end{array}
\right)  \ \mid\ x\in\mathbb{Q}\right\}  .
\end{align*}
This is a \textquotedblleft partly undetermined matrix\textquotedblright\ like
those in the previous example, but this time only the northeasternmost entry
is a question mark, since that entry is the only one that can be changed by
adding a matrix in $\mathbb{Q}^{3<<3}$.

According to the formula (\ref{eq.def.quotring.*bar}), the multiplication on
the quotient ring $\mathbb{Q}^{3\leq3}/\mathbb{Q}^{3<<3}$ is given by%
\[
\left(
\begin{array}
[c]{ccc}%
a & b & ?\\
0 & d & e\\
0 & 0 & f
\end{array}
\right)  \cdot\left(
\begin{array}
[c]{ccc}%
x & y & ?\\
0 & u & v\\
0 & 0 & w
\end{array}
\right)  =\left(
\begin{array}
[c]{ccc}%
ax & bu+ay & ?\\
0 & du & dv+ew\\
0 & 0 & fw
\end{array}
\right)  .
\]
Again, the question-mark entry need not be computed. What is new this time is
that \textbf{we cannot replace the question marks by zeroes}. Indeed, the
product%
\[
\left(
\begin{array}
[c]{ccc}%
a & b & 0\\
0 & d & e\\
0 & 0 & f
\end{array}
\right)  \cdot\left(
\begin{array}
[c]{ccc}%
x & y & 0\\
0 & u & v\\
0 & 0 & w
\end{array}
\right)  =\left(
\begin{array}
[c]{ccc}%
ax & bu+ay & bv\\
0 & du & dv+ew\\
0 & 0 & fw
\end{array}
\right)
\]
does not usually have a $0$ in the northeasternmost position, so that the
matrices of the form $\left(
\begin{array}
[c]{ccc}%
a & b & 0\\
0 & d & e\\
0 & 0 & f
\end{array}
\right)  $ do not form a subring of $\mathbb{Q}^{3\leq3}$. Thus, the
multiplication of our \textquotedblleft partly undetermined
matrices\textquotedblright\ does not simply reduce to the multiplication of
regular matrices\footnote{\textquotedblleft Regular\textquotedblright\ as in
\textquotedblleft the usual matrices you know from high
school\textquotedblright, not as in any of the technical meanings of the word
\textquotedblleft regular\textquotedblright.} in a subring (like it did in the
previous example). In other words, $\mathbb{Q}^{3\leq3}/\mathbb{Q}^{3<<3}$ is
not isomorphic to a subring of $\mathbb{Q}^{3\leq3}$ (at least not in an
obvious way as in the previous example), but is a genuinely new ring.
\end{itemize}

The first of the above two examples can be generalized from $3\times
3$-matrices to $n\times n$-matrices (and from rational entries to entries in
an arbitrary ring $R$):

\begin{exercise}
\label{exe.quotring.upper-triangl1}Let $R$ be any ring. Let $n\in\mathbb{N}$.
Recall that $R^{n\leq n}$ denotes the ring of all upper-triangular $n\times
n$-matrices with entries in $R$. In Exercise
\ref{exe.ideals.triangular-matrices.1} \textbf{(d)}, we have seen that
\[
R^{n<n}=\left\{  A\in R^{n\leq n}\ \mid\ \text{all diagonal entries of
}A\text{ equal }0\right\}
\]
is an ideal of this ring $R^{n\leq n}$. The elements of $R^{n<n}$ are called
the strictly upper-triangular matrices.

Let furthermore $R^{n=n}$ denote the set of all diagonal $n\times n$-matrices
in $R^{n\times n}$. That is,%
\begin{align*}
R^{n=n}:=  &  \ \left\{  A\in R^{n\times n}\ \mid\ \text{all off-diagonal
entries of }A\text{ equal }0\right\} \\
=  &  \ \left\{  \left(
\begin{array}
[c]{ccccc}%
a_{1} & 0 & 0 & \cdots & 0\\
0 & a_{2} & 0 & \cdots & 0\\
0 & 0 & a_{3} & \cdots & 0\\
\vdots & \vdots & \vdots & \ddots & \vdots\\
0 & 0 & 0 & \cdots & a_{n}%
\end{array}
\right)  \ \mid\ a_{1},a_{2},\ldots,a_{n}\in R\right\}  .
\end{align*}
It is easy to see that $R^{n=n}$ is a subring of $R^{n\times n}$.

We claim that the quotient ring $R^{n\leq n}/R^{n<n}$ is isomorphic to
$R^{n=n}$. Intuitively, this is reasonable: When you start with all
upper-triangular matrices but \textquotedblleft equate\textquotedblright\ all
the strictly upper-triangular matrices to zero, then you should be left with
the diagonal matrices, since all the off-diagonal entries are
\textquotedblleft being ignored\textquotedblright. Let us make this rigorous.

Define the map $\delta:R^{n\times n}\rightarrow R^{n\times n}$ as in Exercise
\ref{exe.ringmor.matrix-diagonal}. Note that the actual image of this map
$\delta$ is $R^{n=n}$.

\begin{enumerate}
\item[\textbf{(a)}] Prove that the map%
\begin{align*}
R^{n\leq n}/R^{n<n}  &  \rightarrow R^{n=n},\\
\overline{A}  &  \mapsto\delta\left(  A\right)
\end{align*}
is well-defined -- i.e., the value $\delta\left(  A\right)  $ depends not on
the matrix $A\in R^{n\leq n}$ but only on its residue class $\overline{A}\in
R^{n\leq n}/R^{n<n}$. (In other words, prove that if two matrices $A,B\in
R^{n\leq n}$ have the same residue class $\overline{A}=\overline{B}$ in
$R^{n\leq n}/R^{n<n}$, then $\delta\left(  A\right)  =\delta\left(  B\right)
$.)

\item[\textbf{(b)}] Prove that this map is furthermore a ring morphism.

\item[\textbf{(c)}] Prove that this map is invertible.

This shows that the map is a ring isomorphism, and therefore we have $R^{n\leq
n}/R^{n<n}\cong R^{n=n}$ as rings.
\end{enumerate}
\end{exercise}

The ideal $\mathbb{Q}^{3<<3}$ of the second example can also be generalized:

\begin{exercise}
Let $R$ be any ring. Let $n\in\mathbb{N}$. Let $k\in\mathbb{N}$. Recall that
$R^{n\leq n}$ denotes the ring of all upper-triangular $n\times n$-matrices
with entries in $R$.

A matrix $A=\left(
\begin{array}
[c]{cccc}%
a_{1,1} & a_{1,2} & \cdots & a_{1,n}\\
a_{2,1} & a_{2,2} & \cdots & a_{2,n}\\
\vdots & \vdots & \ddots & \vdots\\
a_{n,1} & a_{n,2} & \cdots & a_{n,n}%
\end{array}
\right)  \in R^{n\times n}$ will be called $k$\textbf{-upper-triangular} if
all its entries $a_{i,j}$ with $j<i+k$ are zero (i.e., if it satisfies
$a_{i,j}=0$ for all $i,j\in\left\{  1,2,\ldots,n\right\}  $ satisfying $j<i+k$).

Let $R_{k}^{n\leq n}$ denote the set of all $k$-upper-triangular $n\times
n$-matrices in $R^{n\times n}$. Prove that $R_{k}^{n\leq n}$ is an ideal of
$R^{n\leq n}$.

[The ideals $\mathbb{Q}^{3<3}$ and $\mathbb{Q}^{3<<3}$ from the above two
examples are $\mathbb{Q}_{1}^{3\leq3}$ and $\mathbb{Q}_{2}^{3\leq3}$, respectively.]
\end{exercise}

\subsubsection{The canonical projection}

Ideals of rings are somewhat like normal subgroups of groups: You can
\textquotedblleft quotient them out\textquotedblright\ (this is slang for
\textquotedblleft take a quotient by them\textquotedblright) and get a ring
again (by Theorem \ref{thm.quotring.welldef}).

Now, we are ready to show that any ideal of a ring is the kernel of a ring morphism:

\begin{theorem}
\label{thm.quotring.canproj}Let $R$ be a ring. Let $I$ be an ideal of $R$.
Consider the map
\begin{align*}
\pi:R  &  \rightarrow R/I,\\
r  &  \mapsto r+I.
\end{align*}
Then, $\pi$ is a surjective ring morphism with kernel $I$.
\end{theorem}

\begin{definition}
This morphism $\pi$ is called the \textbf{canonical projection} from $R$ onto
$R/I$.
\end{definition}

\begin{proof}
[Proof of Theorem \ref{thm.quotring.canproj}.]To prove that $\pi$ is a ring
morphism, we need to check that $\pi$ respects addition, multiplication, zero
and unity. All of this is straightforward. For example, in order to see that
$\pi$ respects multiplication, we must show that $\pi\left(  rs\right)
=\pi\left(  r\right)  \cdot\pi\left(  s\right)  $ for all $r,s\in R$; but this
follows from%
\begin{align*}
\pi\left(  r\right)  \cdot\pi\left(  s\right)   &  =\left(  r+I\right)
\cdot\left(  s+I\right)  \ \ \ \ \ \ \ \ \ \ \left(  \text{by the definition
of }\pi\right) \\
&  =rs+I\ \ \ \ \ \ \ \ \ \ \left(
\begin{array}
[c]{c}%
\text{by the definition of the}\\
\text{multiplication in }R/I
\end{array}
\right) \\
&  =\pi\left(  rs\right)  \ \ \ \ \ \ \ \ \ \ \left(  \text{by the definition
of }\pi\right)  .
\end{align*}
Thus, we can show that $\pi$ is a ring morphism.

The surjectivity of $\pi$ is clear, since any element of $R/I$ has the form
$r+I=\pi\left(  r\right)  $ for some $r\in R$.

Finally, we need to show that $\pi$ has kernel $I$. But this, too, is easy:
The kernel of $\pi$ consists of those elements $r\in R$ that satisfy
$\pi\left(  r\right)  =0_{R/I}$. But $\pi\left(  r\right)  =0_{R/I}$ is
equivalent to $r+I=0+I$ (since $\pi\left(  r\right)  $ is the coset $r+I$,
whereas $0_{R/I}$ is the coset $0+I$), which is tantamount to $r\in I$. Thus,
the kernel of $\pi$ is $I$.
\end{proof}

The canonical projection $\pi$ in Theorem \ref{thm.quotring.canproj} can be
viewed as \textquotedblleft putting a bar on each element\textquotedblright,
since it sends each $r\in R$ to the residue class $\overline{r}\in R/I$.

For example, if we take $R=\mathbb{Z}$ and $I=2\mathbb{Z}$ in Theorem
\ref{thm.quotring.canproj}, then the canonical projection $\pi$ is the map%
\begin{align*}
\pi:\mathbb{Z}  &  \rightarrow\mathbb{Z}/2,\\
r  &  \mapsto r+2\mathbb{Z}=\overline{r}.
\end{align*}
This map $\pi$ sends each even integer to $\overline{0}$ and each odd integer
to $\overline{1}$. In other words, $\pi$ assigns to each integer its parity
(as an element of $\mathbb{Z}/2$).

The following two exercises (\cite[homework set \#2, Exercise 10 \textbf{(a)}
and \textbf{(b)}]{21w}) illustrate one of many uses of quotient rings:

\begin{exercise}
\label{exe.21hw2.10a}Let $R$ be a commutative ring, and let $u$ and $n$ be two
nonnegative integers. Let $x,y\in R$ be two elements such that $x-y\in uR$.
(Here, $uR:=\left\{  ur\mid r\in R\right\}  $; this is a principal ideal of
$R$, since $uR=\left(  u1_{R}\right)  R$.)

Prove that
\[
x^{n}-y^{n}\in guR,\ \ \ \ \ \ \ \ \ \ \text{where }g=\gcd\left(  n,u\right)
.
\]

[\textbf{Hint:} Write $x^{n}-y^{n}$ as $\left(  x-y\right)  \left(
x^{n-1}+x^{n-2}y+\cdots+y^{n-1}\right)  $, and show that the second factor
belongs to $gR$. The latter is easiest to do by working in the quotient ring
$R/gR$.]
\end{exercise}

\begin{exercise}
\label{exe.21hw2.10b}Let $\left(  f_{0},f_{1},f_{2},\ldots\right)  $ be the
Fibonacci sequence, defined as in Definition \ref{def.fibonacci.fib}. Prove
that
\[
\gcd\left(  n,f_{d}\right)  \cdot f_{d}\mid f_{dn}%
\ \ \ \ \ \ \ \ \ \ \text{for any }d,n\in\mathbb{N}.
\]

[\textbf{Hint:} Define the matrices $A$ and $B$ and the commutative ring
$\mathcal{F}$ as in Exercise \ref{exe.21hw1.6}, and apply Exercise
\ref{exe.21hw2.10a} to $R=\mathcal{F}$ and $x=A^{d}$ and $y=B^{d}$ and
$u=f_{d}$.]
\end{exercise}

\subsubsection{The universal property of quotient rings: elementwise form}

When trying to understand a quotient ring, it is important to construct ring
morphisms into and out of it. (In the best case scenario, you can find
mutually inverse maps in both directions, and thus obtain an isomorphism to
another ring.)

Constructing morphisms $\alpha:S\rightarrow R/I$ into a quotient ring $R/I$ is
pretty easy (see Theorem \ref{thm.quotring.canproj} for an example).

Constructing morphisms $\beta:R/I\rightarrow S$ out of a quotient ring $R/I$
is trickier. The main problem is to establish that $\beta$ is well-defined: An
element of a quotient ring $R/I$ can be written as $\overline{r}$ for many
different elements $r\in R$, and therefore, when assigning an output value
$\beta\left(  \overline{r}\right)  $ for this element, we need to ensure that
this value depends not on the $r$ but only on the $\overline{r}$. This can be
done by hand (see Exercise \ref{exe.quotring.upper-triangl1} \textbf{(a)} for
an easy example of this), but gets tedious fairly soon. Is there a more
comfortable method?

Yes, and such a method is provided by a theorem known as the \textquotedblleft%
\textbf{universal property of quotient rings}\textquotedblright. This theorem
may appear technical, abstract and pointless at first sight, but it reveals
its usefulness soon after you tire of manually constructing morphisms out of
quotient rings. It provides a mechanical way of constructing a ring morphism
$f^{\prime}:R/I\rightarrow S$ out of a ring morphism $f:R\rightarrow S$, as
soon as you can show that $f$ sends all elements of the ideal $I$ to $0$. The
well-definedness of $f^{\prime}$ and the fact that $f^{\prime}$ is a ring
morphism are automatic consequences of the theorem, once its assumptions have
been satisfied. Here is the precise statement of the theorem (in one of its forms):

\begin{theorem}
[Universal property of quotient rings, elementwise form]%
\label{thm.quotring.uniprop1}Let $R$ be a ring. Let $I$ be an ideal of $R$.

Let $S$ be a ring. Let $f:R\rightarrow S$ be a ring morphism. Assume that
$f\left(  I\right)  =0$ (this is shorthand for saying that $f\left(  a\right)
=0$ for all $a\in I$). Then, the map%
\begin{align*}
f^{\prime}:R/I  &  \rightarrow S,\\
\overline{r}  &  \mapsto f\left(  r\right)  \ \ \ \ \ \ \ \ \ \ \left(
\text{for all }r\in R\right)
\end{align*}
\footnotemark\ is well-defined (i.e., the value $f\left(  r\right)  $ depends
only on the residue class $\overline{r}$, not on $r$ itself) and is a ring morphism.
\end{theorem}

\footnotetext{Recall that $\overline{r}$ means the residue class of $r$ in
$R/I$, that is, the coset $r+I$. Thus, the definition of $f^{\prime}$ can be
rewritten as follows:%
\begin{align*}
f^{\prime}:R/I  &  \rightarrow S,\\
r+I  &  \mapsto f\left(  r\right)  \ \ \ \ \ \ \ \ \ \ \left(  \text{for all
}r\in R\right)  .
\end{align*}
\par
Roughly speaking, the definition of $f^{\prime}$ says that $f^{\prime}$ sends
a residue class to where $f$ would send any element of this residue class. The
(slightly) nontrivial part here is to prove that this is well-defined, i.e.,
that $f$ takes all elements of the given residue class to the same output
value.}Before we prove Theorem \ref{thm.quotring.uniprop1}, let us give an example:

\begin{itemize}
\item Consider the canonical projections%
\begin{align*}
\pi_{6}:\mathbb{Z}  &  \rightarrow\mathbb{Z}/6,\\
r  &  \mapsto r+6\mathbb{Z}%
\end{align*}
and%
\begin{align*}
\pi_{3}:\mathbb{Z}  &  \rightarrow\mathbb{Z}/3,\\
r  &  \mapsto r+3\mathbb{Z}.
\end{align*}
(Each of these two projections sends each integer $r$ to its residue class
$\overline{r}$, but the residue class is a modulo-$6$ class for $\pi_{6}$ and
a modulo-$3$ class for $\pi_{3}$. The notation $\overline{r}$ can mean either
$r+6\mathbb{Z}$ or $r+3\mathbb{Z}$ depending on the context. Thus, pay
attention to the sets to which the elements belong!)

The ideal $6\mathbb{Z}$ of $\mathbb{Z}$ satisfies $\pi_{3}\left(
6\mathbb{Z}\right)  =0$ (because any $j\in6\mathbb{Z}$ is a multiple of $6$,
thus a multiple of $3$, and therefore its residue class $j+3\mathbb{Z}$ is
$\overline{0}$, and thus $\pi_{3}\left(  j\right)  =j+3\mathbb{Z}=\overline
{0}=0_{\mathbb{Z}/3}$). Thus, by Theorem \ref{thm.quotring.uniprop1} (applied
to $R=\mathbb{Z}$, $I=6\mathbb{Z}$, $S=\mathbb{Z}/3$ and $f=\pi_{3}$), we see
that the map%
\begin{align}
\pi_{3}^{\prime}:\mathbb{Z}/6  &  \rightarrow\mathbb{Z}/3,\nonumber\\
\overline{r}  &  \mapsto\pi_{3}\left(  r\right)  \ \ \ \ \ \ \ \ \ \ \left(
\text{that is, }r+6\mathbb{Z}\mapsto r+3\mathbb{Z}\right)
\label{eq.exa.quotring.Z/6-to-Z/3}%
\end{align}
is well-defined and is a ring morphism. Explicitly, this morphism $\pi
_{3}^{\prime}$ sends\footnote{A \textquotedblleft modulo-$6$ residue
class\textquotedblright\ $\overline{r}$ means the residue class $r+6\mathbb{Z}%
$, whereas a \textquotedblleft modulo-$3$ residue class\textquotedblright%
\ $\overline{r}$ means the residue class $r+3\mathbb{Z}$.}%
\begin{align*}
\text{the modulo-}6\text{ residue classes\ \ }  &  \overline{0},\ \ \overline
{1},\ \ \overline{2},\ \ \overline{3},\ \ \overline{4},\ \ \overline{5}\\
\text{to the modulo-}3\text{ residue classes\ \ }  &  \overline{0}%
,\ \ \overline{1},\ \ \overline{2},\ \ \overline{3},\ \ \overline
{4},\ \ \overline{5}.
\end{align*}
In other words, it sends%
\begin{align*}
\text{the modulo-}6\text{ residue classes\ \ }  &  \overline{0},\ \ \overline
{1},\ \ \overline{2},\ \ \overline{3},\ \ \overline{4},\ \ \overline{5}\\
\text{to the modulo-}3\text{ residue classes\ \ }  &  \overline{0}%
,\ \ \overline{1},\ \ \overline{2},\ \ \overline{0},\ \ \overline
{1},\ \ \overline{2}%
\end{align*}
(because in $\mathbb{Z}/3$, we have $\overline{3}=\overline{0}$ and
$\overline{4}=\overline{1}$ and $\overline{5}=\overline{2}$). If you don't
believe in Theorem \ref{thm.quotring.uniprop1} yet, you can easily check by
hand that this is a ring morphism.

More generally, if $n$ and $m$ are two integers such that $m\mid n$, then the
map%
\begin{align}
\mathbb{Z}/n  &  \rightarrow\mathbb{Z}/m,\nonumber\\
\overline{r}  &  \mapsto\overline{r}\ \ \ \ \ \ \ \ \ \ \left(  \text{that is,
}r+n\mathbb{Z}\mapsto r+m\mathbb{Z}\right)  \label{eq.exa.quotring.Z/n-to-Z/m}%
\end{align}
is well-defined and is a ring morphism. This follows from Theorem
\ref{thm.quotring.uniprop1}, applied to $R=\mathbb{Z}$, $I=n\mathbb{Z}$,
$S=\mathbb{Z}/m$ and $f=\pi_{m}$ (the canonical projection from $\mathbb{Z}$
to $\mathbb{Z}/m$), because the condition $m\mid n$ yields $\pi_{m}\left(
n\mathbb{Z}\right)  =0$. The morphism (\ref{eq.exa.quotring.Z/n-to-Z/m}) can
be regarded as reducing a modulo-$n$ residue class \textquotedblleft
further\textquotedblright\ to a modulo-$m$ residue class.

Incidentally, this accounts for all ring morphisms that go between two
quotient rings of $\mathbb{Z}$. That is:

\begin{itemize}
\item If $m$ and $n$ are two integers such that $m\mid n$, then there is only
one ring morphism $\mathbb{Z}/n\rightarrow\mathbb{Z}/m$, and it is the
morphism (\ref{eq.exa.quotring.Z/n-to-Z/m}).

\item If $m$ and $n$ are two integers such that $m\nmid n$, then there is no
ring morphism $\mathbb{Z}/n\rightarrow\mathbb{Z}/m$.
\end{itemize}

Proving this is a nice (and easy) exercise (Exercise
\ref{exe.ring.mors.Z/n-to-Z/m}).
\end{itemize}

Let us now prove the universal property of quotient rings:

\begin{proof}
[Proof of Theorem \ref{thm.quotring.uniprop1}.]We must prove the following two facts:

\begin{enumerate}
\item The map
\begin{align*}
f^{\prime}:R/I  &  \rightarrow S,\\
\overline{r}  &  \mapsto f\left(  r\right)  \ \ \ \ \ \ \ \ \ \ \left(
\text{for all }r\in R\right)
\end{align*}
is well-defined -- i.e., the value $f\left(  r\right)  $ depends only on the
residue class $\overline{r}$ but not on the specific choice of $r$. (This will
ensure that its definition does not give two conflicting output values
$f^{\prime}\left(  x\right)  $ for one and the same residue class $x\in R/I$,
which would spell doom for the map $f^{\prime}$.)

\item The map $f^{\prime}$ is a ring morphism.
\end{enumerate}

Let us prove Fact 1 first. So let $r,r^{\prime}\in R$ be such that
$\overline{r}=\overline{r^{\prime}}$. We must show that $f\left(  r\right)
=f\left(  r^{\prime}\right)  $.

We do what we can: From $\overline{r}=\overline{r^{\prime}}$, we obtain
$r-r^{\prime}\in I$, so that $f\left(  r-r^{\prime}\right)  =0$ because
$f\left(  I\right)  =0$. However, $f$ is a ring morphism and thus respects
differences; hence, $f\left(  r-r^{\prime}\right)  =f\left(  r\right)
-f\left(  r^{\prime}\right)  $. Thus, $f\left(  r\right)  -f\left(  r^{\prime
}\right)  =f\left(  r-r^{\prime}\right)  =0$, so that $f\left(  r\right)
=f\left(  r^{\prime}\right)  $. This proves Fact 1.

Let us now prove Fact 2. We need to show that $f^{\prime}$ is a ring morphism.
There are four axioms to check, but we shall only show one of them (since the
proofs of the other three axioms follow the same mold). Namely, we shall show
that $f^{\prime}$ respects multiplication.

So let $a,b\in R/I$. We must show that $f^{\prime}\left(  ab\right)
=f^{\prime}\left(  a\right)  \cdot f^{\prime}\left(  b\right)  $.

Write the residue classes $a,b\in R/I$ as $a=\overline{r}$ and $b=\overline
{s}$ for some $r,s\in R$. Then, $ab=\overline{r}\cdot\overline{s}%
=\overline{rs}$ by the formula (\ref{eq.def.quotring.*bar}). Hence,
$f^{\prime}\left(  ab\right)  =f^{\prime}\left(  \overline{rs}\right)
=f\left(  rs\right)  $ (by the definition of $f^{\prime}$). On the other hand,
$a=\overline{r}$ and thus $f^{\prime}\left(  a\right)  =f^{\prime}\left(
\overline{r}\right)  =f\left(  r\right)  $ (by the definition of $f^{\prime}%
$). Similarly, $f^{\prime}\left(  b\right)  =f\left(  s\right)  $. Thus,
\begin{align*}
f^{\prime}\left(  ab\right)   &  =f\left(  rs\right)  =\underbrace{f\left(
r\right)  }_{=f^{\prime}\left(  a\right)  }\cdot\underbrace{f\left(  s\right)
}_{=f^{\prime}\left(  b\right)  }\ \ \ \ \ \ \ \ \ \ \left(  \text{since
}f\text{ is a ring morphism}\right) \\
&  =f^{\prime}\left(  a\right)  \cdot f^{\prime}\left(  b\right)  ,
\end{align*}
which is precisely what we wanted to prove. Thus, Fact 2 is proved as well.

So we have shown that our map $f^{\prime}:R/I\rightarrow S$ is well-defined
and is a ring morphism. Thus, Theorem \ref{thm.quotring.uniprop1} is proven.
\end{proof}

As we said, the universal property of quotient rings provides a comfortable
way to construct ring morphisms out of a quotient ring $R/I$. The following
exercises provide some examples of this:

\begin{exercise}
Consider the quotient ring $\mathbb{Q}^{3\leq3}/\mathbb{Q}^{3<<3}$ studied in
Subsection \ref{subsec.rings.quotring.moreexas}.

\begin{enumerate}
\item[\textbf{(a)}] Prove that the map%
\begin{align*}
f:\mathbb{Q}^{3\leq3}  &  \rightarrow\mathbb{Q}^{2\times2},\\
\left(
\begin{array}
[c]{ccc}%
a & b & c\\
0 & d & e\\
0 & 0 & g
\end{array}
\right)   &  \mapsto\left(
\begin{array}
[c]{cc}%
a & b\\
0 & d
\end{array}
\right)
\end{align*}
is a ring morphism.

\item[\textbf{(b)}] Prove that this morphism $f$ satisfies $f\left(
\mathbb{Q}^{3<<3}\right)  =0$.

\item[\textbf{(c)}] Use Theorem \ref{thm.quotring.uniprop1} to conclude that
there is a ring morphism%
\begin{align*}
f^{\prime}:\mathbb{Q}^{3\leq3}/\mathbb{Q}^{3<<3}  &  \rightarrow
\mathbb{Q}^{2\times2},\\
\left(
\begin{array}
[c]{ccc}%
a & b & ?\\
0 & d & e\\
0 & 0 & g
\end{array}
\right)   &  \mapsto\left(
\begin{array}
[c]{cc}%
a & b\\
0 & d
\end{array}
\right)  ,
\end{align*}
where the question mark stands for an undetermined entry (as explained in
Subsection \ref{subsec.rings.quotring.moreexas}).

\item[\textbf{(d)}] Use a similar reasoning to prove the existence of a ring
morphism%
\begin{align*}
F^{\prime}:\mathbb{Q}^{3\leq3}/\mathbb{Q}^{3<<3}  &  \rightarrow
\mathbb{Q}^{4\times4},\\
\left(
\begin{array}
[c]{ccc}%
a & b & ?\\
0 & d & e\\
0 & 0 & g
\end{array}
\right)   &  \mapsto\left(
\begin{array}
[c]{cccc}%
a & b & 0 & 0\\
0 & d & 0 & 0\\
0 & 0 & d & e\\
0 & 0 & 0 & g
\end{array}
\right)  ,
\end{align*}
which is furthermore injective.

\item[\textbf{(e)}] Conclude that the ring $\mathbb{Q}^{3\leq3}/\mathbb{Q}%
^{3<<3}$ is isomorphic to a subring of $\mathbb{Q}^{4\times4}$.
\end{enumerate}
\end{exercise}

\begin{exercise}
Solve parts \textbf{(a)} and \textbf{(b)} of Exercise
\ref{exe.quotring.upper-triangl1} again using the universal property of
quotient rings. (This should be much quicker than the original solutions.)
\end{exercise}

\begin{exercise}
For every integer $m$, define a subring $R_{m}$ of $\mathbb{Q}$ as in Exercise
\ref{exe.21hw1.1}.

Consider the quotient ring $R_{2}/3R_{2}$.

\begin{enumerate}
\item[\textbf{(a)}] Prove that the map%
\begin{align*}
f_{2}:\mathbb{Z}  &  \rightarrow R_{2}/3R_{2},\\
r  &  \mapsto\overline{r}%
\end{align*}
is a ring morphism that satisfies $f_{2}\left(  3\mathbb{Z}\right)  =0$.

\item[\textbf{(b)}] Use the universal property of quotient rings to obtain a
ring morphism%
\begin{align*}
f_{2}^{\prime}:\mathbb{Z}/3  &  \rightarrow R_{2}/3R_{2},\\
\overline{r}  &  \mapsto\overline{r}.
\end{align*}

\item[\textbf{(c)}] Show that this morphism $f_{2}^{\prime}$ is injective.

\item[\textbf{(d)}] Show that this morphism $f_{2}^{\prime}$ is surjective.

\item[\textbf{(e)}] Conclude that $\mathbb{Z}/3\mathbb{Z}\cong R_{2}/3R_{2}$
as rings.

\item[\textbf{(f)}] More generally, let $m$ and $n$ be two coprime integers.
Show that there is a ring isomorphism%
\begin{align*}
\mathbb{Z}/n  &  \rightarrow R_{m}/nR_{m},\\
\overline{r}  &  \mapsto\overline{r}.
\end{align*}

\end{enumerate}

[\textbf{Hint:} Parts \textbf{(a)} and \textbf{(b)} are almost automatic. Part
\textbf{(c)} requires showing that an integer $r$ that is not divisible by $3$
cannot lie in $3R_{2}$ either. Argue this using the coprimality of $2$ and
$3$. Part \textbf{(d)} boils down to showing that the map $f_{2}^{\prime}$
takes $\overline{1/2}$ as a value (why?), but this is easy (why?). Part
\textbf{(f)} requires generalizing the previous parts, using some properties
of coprime integers that we have seen before.]
\end{exercise}

We have also promised another exercise:

\begin{exercise}
\label{exe.ring.mors.Z/n-to-Z/m}Let $m$ and $n$ be two integers. Prove the following:

\begin{enumerate}
\item[\textbf{(a)}] If $m\mid n$, then there is only one ring morphism
$\mathbb{Z}/n\rightarrow\mathbb{Z}/m$, and it is the morphism
(\ref{eq.exa.quotring.Z/n-to-Z/m}).

\item[\textbf{(b)}] If $m\nmid n$, then there is no ring morphism
$\mathbb{Z}/n\rightarrow\mathbb{Z}/m$.
\end{enumerate}

[\textbf{Hint:} For part \textbf{(a)}, show that every ring morphism
$g:\mathbb{Z}/n\rightarrow\mathbb{Z}/m$ must satisfy $g\left(  \overline
{r}\right)  =\overline{r}$ for each $r\in\mathbb{N}$, since $\overline
{r}=r\cdot1_{\mathbb{Z}/n}=\underbrace{1_{\mathbb{Z}/n}+1_{\mathbb{Z}%
/n}+\cdots+1_{\mathbb{Z}/n}}_{r\text{ times}}$. For part \textbf{(b)}, argue
that the existence of a ring morphism $g:\mathbb{Z}/n\rightarrow\mathbb{Z}/m$
forces $n\cdot1_{\mathbb{Z}/m}=0$ because $n\cdot1_{\mathbb{Z}/n}=0$.]
\end{exercise}

\subsubsection{\label{subsec.quotring.uniprop2}The universal property of
quotient rings: abstract form}

For various reasons, it is helpful to have an alternative formulation of
Theorem \ref{thm.quotring.uniprop1}, which does not refer to specific elements
$\overline{r}$ but instead \textquotedblleft implicitly\textquotedblright%
\ describes the morphism $f^{\prime}$ by an equality:

\begin{theorem}
[Universal property of quotient rings, abstract form]%
\label{thm.quotring.uniprop2}Let $R$ be a ring. Let $I$ be an ideal of $R$.
Consider the canonical projection $\pi:R\rightarrow R/I$ (as defined in
Theorem \ref{thm.quotring.canproj}).

Let $S$ be a ring. Let $f:R\rightarrow S$ be a ring morphism. Assume that
$f\left(  I\right)  =0$ (this is shorthand for saying that $f\left(  a\right)
=0$ for all $a\in I$). Then, there is a unique ring morphism $f^{\prime
}:R/I\rightarrow S$ satisfying $f=f^{\prime}\circ\pi$.
\end{theorem}

\begin{proof}
Theorem \ref{thm.quotring.uniprop1} shows that there is a unique ring morphism
$f^{\prime}:R/I\rightarrow S$ that satisfies
\begin{equation}
f^{\prime}\left(  \overline{r}\right)  =f\left(  r\right)
\ \ \ \ \ \ \ \ \ \ \text{for all }r\in R.
\label{pf.thm.quotring.uniprop2.old-cond}%
\end{equation}
(Indeed, the equality (\ref{pf.thm.quotring.uniprop2.old-cond}) clearly
characterizes $f^{\prime}$ uniquely, since every element of $R/I$ can be
written as $\overline{r}$ for some $r\in R$. What Theorem
\ref{thm.quotring.uniprop1} gives us is the \textbf{existence} of such a
morphism $f^{\prime}$.)

We shall now prove that the equality $f=f^{\prime}\circ\pi$ is just an
equivalent restatement of the condition
(\ref{pf.thm.quotring.uniprop2.old-cond}).

Indeed, we have the following chain of equivalences:%
\begin{align*}
&  \ \left(  f=f^{\prime}\circ\pi\right) \\
&  \Longleftrightarrow\ \left(  f\left(  r\right)  =\left(  f^{\prime}\circ
\pi\right)  \left(  r\right)  \text{ for all }r\in R\right)
\ \ \ \ \ \ \ \ \ \ \left(
\begin{array}
[c]{c}%
\text{since two maps are equal}\\
\text{if and only if they}\\
\text{agree on each input}%
\end{array}
\right) \\
&  \Longleftrightarrow\ \left(  f\left(  r\right)  =f^{\prime}\left(
\pi\left(  r\right)  \right)  \text{ for all }r\in R\right)
\ \ \ \ \ \ \ \ \ \ \left(
\begin{array}
[c]{c}%
\text{since }\left(  f^{\prime}\circ\pi\right)  \left(  r\right)  =f^{\prime
}\left(  \pi\left(  r\right)  \right) \\
\text{for each }r
\end{array}
\right) \\
&  \Longleftrightarrow\ \left(  f\left(  r\right)  =f^{\prime}\left(
\overline{r}\right)  \text{ for all }r\in R\right)
\ \ \ \ \ \ \ \ \ \ \left(  \text{since }\pi\left(  r\right)  =\overline
{r}\text{ for each }r\right) \\
&  \Longleftrightarrow\ \left(  f^{\prime}\left(  \overline{r}\right)
=f\left(  r\right)  \text{ for all }r\in R\right)  .
\end{align*}
In other words, the equality $f=f^{\prime}\circ\pi$ is equivalent to the
condition (\ref{pf.thm.quotring.uniprop2.old-cond}).

Now, recall that there is a unique ring morphism $f^{\prime}:R/I\rightarrow S$
that satisfies the condition (\ref{pf.thm.quotring.uniprop2.old-cond}). In
view of the previous sentence, we can reformulate this as follows: There is a
unique ring morphism $f^{\prime}:R/I\rightarrow S$ that satisfies
$f=f^{\prime}\circ\pi$. This proves Theorem \ref{thm.quotring.uniprop2}.
\end{proof}

For example:

\begin{itemize}
\item We can recover the ring morphism $\pi_{3}^{\prime}:\mathbb{Z}%
/6\rightarrow\mathbb{Z}/3$ constructed in (\ref{eq.exa.quotring.Z/6-to-Z/3})
using Theorem \ref{thm.quotring.uniprop2} instead of Theorem
\ref{thm.quotring.uniprop1}. Indeed, applying Theorem
\ref{thm.quotring.uniprop2} to $R=\mathbb{Z}$, $I=6\mathbb{Z}$, $\pi=\pi_{6}$,
$S=\mathbb{Z}/3$ and $f=\pi_{3}$, we see that there is a unique ring morphism
$\pi_{3}^{\prime}:\mathbb{Z}/6\rightarrow\mathbb{Z}/3$ such that $\pi_{3}%
=\pi_{3}^{\prime}\circ\pi_{6}$ (since $\pi_{3}\left(  6\mathbb{Z}\right)
=0$). This morphism $\pi_{3}^{\prime}$ is, of course, the same as the one in
(\ref{eq.exa.quotring.Z/6-to-Z/3}), since the equality $\pi_{3}=\pi
_{3}^{\prime}\circ\pi_{6}$ says precisely that each integer $r$ satisfies
$\pi_{3}\left(  r\right)  =\pi_{3}^{\prime}\left(  \pi_{6}\left(  r\right)
\right)  $, that is, $\overline{r}=\pi_{3}^{\prime}\left(  \overline
{r}\right)  $.
\end{itemize}

A few more remarks are in order.

The equality $f=f^{\prime}\circ\pi$ in Theorem \ref{thm.quotring.uniprop2} is
oftentimes restated as follows: The diagram%
\begin{equation}%
\xymatrix@C=4pc{
R \ar[d]_{\pi} \ar[dr]^f \\
R/I \ar@{-->}[r]_{f'} & S
}
\label{eq.thm.quotring.uniprop2.diagram}%
\end{equation}
commutes. Let me explain what this means: In general, a \textbf{diagram} is a
bunch of sets and a bunch of maps between them, drawn as nodes and arrows.
(Each set is drawn as a node, and each map $g:A\rightarrow B$ is drawn as an
arrow from the $A$-node to the $B$-node.) For instance, the diagram
(\ref{eq.thm.quotring.uniprop2.diagram}) shows the three sets $R$, $R/I$ and
$S$ and the three maps $\pi$, $f$ and $f^{\prime}$. A diagram is said to
\textbf{commute} (or \textbf{be commutative}) if any two ways of going between
two nodes yield the same composed map. In the diagram
(\ref{eq.thm.quotring.uniprop2.diagram}), there are two ways of going from the
$R$-node to the $S$-node: one direct way (which just uses the $f$-arrow), and
one indirect way via the $R/I$-node (using the $\pi$-arrow and the $f^{\prime
}$-arrow). The corresponding composed maps are $f$ (for the direct way) and
$f^{\prime}\circ\pi$ (for the indirect way). This is the only pair of two
different ways that go between the same two nodes in the diagram
(\ref{eq.thm.quotring.uniprop2.diagram}); thus, the diagram commutes if and
only if $f=f^{\prime}\circ\pi$.

Note that we have drawn the map $f^{\prime}$ as a dashed arrow in
(\ref{eq.thm.quotring.uniprop2.diagram}), since this is the map whose
existence is claimed, whereas the other two maps are given and thus drawn as
regular arrows. This is a common convention and helps you distinguish the
things you have from the things you are trying to construct.

In general, diagrams are a good way to visualize situations in which there are
several maps going between the same sets. For example, here is a diagram that
shows the rings $\mathbb{Z}$, $\mathbb{Z}/12$, $\mathbb{Z}/6$, $\mathbb{Z}/4$,
$\mathbb{Z}/3$, $\mathbb{Z}/2$ and $\mathbb{Z}/1$ (the latter ring is trivial)
as well as various morphisms between them:
\[%
\xymatrix{
& \ZZ\ar@/_1pc/[ddl] \ar@/_5pc/[dddl] \ar@/_13pc/[dddd] \ar[d] \ar@
/^1pc/[dddr] \\
& \ZZ/12 \ar[dl] \ar[d] \ar[ddr] \\
\ZZ/4 \ar[d] & \ZZ/6 \ar[dl] \ar[dr] \ar[dd] \\
\ZZ/2 \ar[dr] & & \ZZ/3 \ar[dl] \\
& \ZZ/1
}%
\ \ .
\]
In this diagram, all arrows coming out of the $\mathbb{Z}$-node are canonical
projections $\mathbb{Z}\rightarrow\mathbb{Z}/n$ (sending each $r\in\mathbb{Z}$
to $\overline{r}\in\mathbb{Z}/n$), whereas all the other arrows are instances
of the morphisms (\ref{eq.exa.quotring.Z/n-to-Z/m}) constructed above. Note
that we have not drawn all possible morphisms (e.g., the morphism
$\mathbb{Z}/4\rightarrow\mathbb{Z}/1$ is missing) to avoid crowding the
diagram. This diagram commutes, since each of the arrows sends each residue
class $\overline{r}$ (or, in the case of $\mathbb{Z}$, each integer $r$) to
the corresponding residue class $\overline{r}$ modulo the respective number.

Commutative diagrams become increasingly useful as you go deeper into algebra
(and become ubiquitous when you get to category theory or homological
algebra). For us here, they are just convenient visual and mnemonic devices.

\subsubsection{Injectivity means zero kernel}

Next comes another useful result: a characterization of injectivity in terms
of kernels.

\begin{lemma}
\label{lem.rings.mors.inj}Let $R$ and $S$ be two rings. Let $f:R\rightarrow S$
be a ring morphism. Then, $f$ is injective if and only if $\operatorname*{Ker}%
f=\left\{  0_{R}\right\}  $.
\end{lemma}

\begin{proof}
You probably have seen the analogous results for groups or vector spaces. If
so, then you can just recall the analogous result for groups, and apply it to
the additive groups $\left(  R,+,0\right)  $ and $\left(  S,+,0\right)  $
(since the ring morphism $f$ is clearly a group morphism from $\left(
R,+,0\right)  $ to $\left(  S,+,0\right)  $).

If not, here is the proof:

$\Longrightarrow:$ Assume that $f$ is injective. Then, each $x\in
\operatorname*{Ker}f$ satisfies $f\left(  x\right)  =0_{S}=f\left(
0_{R}\right)  $ (since $f$ is a ring morphism) and thus $x=0_{R}$ because $f$
is injective. In other words, $\operatorname*{Ker}f\subseteq\left\{
0_{R}\right\}  $. But this entails $\operatorname*{Ker}f=\left\{
0_{R}\right\}  $ (since $0_{R}$ always lies in $\operatorname*{Ker}f$). This
proves the \textquotedblleft$\Longrightarrow$\textquotedblright\ direction of
Lemma \ref{lem.rings.mors.inj}.

$\Longleftarrow:$ Assume that $\operatorname*{Ker}f=\left\{  0_{R}\right\}  $.
Now, $f$ is a ring morphism and thus respects differences. Hence, if $a,b\in
R$ satisfy $f\left(  a\right)  =f\left(  b\right)  $, then $f\left(
a-b\right)  =f\left(  a\right)  -f\left(  b\right)  =0$ (since $f\left(
a\right)  =f\left(  b\right)  $) and therefore $a-b\in\operatorname*{Ker}%
f=\left\{  0_{R}\right\}  $, so that $a-b=0$ and thus $a=b$. But this means
that $f$ is injective. This proves the \textquotedblleft$\Longleftarrow
$\textquotedblright\ direction of Lemma \ref{lem.rings.mors.inj}, and thus
completes the proof of the lemma.
\end{proof}

\subsubsection{The First Isomorphism Theorem for sets}

We now approach another important property of quotient rings: the so-called
\textquotedblleft first isomorphism theorem\textquotedblright.

We begin with some basic set theory.

Consider a map $f:R\rightarrow S$ from some set $R$ to some set $S$. Then, I
claim that there is a bijection\footnote{\textquotedblleft
Bijection\textquotedblright\ means the same as \textquotedblleft bijective
map\textquotedblright\ (i.e., a map that is both injective and surjective) and
as \textquotedblleft1-to-1 correspondence\textquotedblright. Also, it is worth
recalling that a map is bijective if and only if it is invertible (i.e., has
an inverse).} hiding inside $f$.

What do I mean by this?

For an example, let $R=\left\{  1,2,3,4,5\right\}  $ and $S=\left\{
1,2,3,4\right\}  $, and let $f:R\rightarrow S$ be the map that sends
$1,2,3,4,5$ to $2,4,2,1,1$, respectively. Here is an illustration of this map
using a standard \textquotedblleft blobs and arrows\textquotedblright%
\ diagram:%
\[%
\begin{tikzpicture}
\filldraw[fill=cyan!20!white] (0, 9) ellipse (1 and 3);
\node(C1) at (0, 11) {$1$};
\node(C2) at (0, 10) {$2$};
\node(C3) at (0, 9) {$3$};
\node(C4) at (0, 8) {$4$};
\node(C5) at (0, 7) {$5$};
\node(Cname) at (0, 12.4) {$R$};
\filldraw[fill=cyan!20!white] (3, 8.5) ellipse (1 and 3);
\node(D1) at (3, 10) {$1$};
\node(D2) at (3, 9) {$2$};
\node(D3) at (3, 8) {$3$};
\node(D4) at (3, 7) {$4$};
\node(Dname) at (3, 11.9) {$S$};
\draw[->, very thick, darkgreen] (C1) -- (D2);
\draw[->, very thick, darkgreen] (C2) -- (D4);
\draw[->, very thick, darkgreen] (C3) -- (D2);
\draw[->, very thick, darkgreen] (C4) -- (D1);
\draw[->, very thick, darkgreen] (C5) -- (D1);
\node(oname) at (1.5, 10.5) {$\color{darkgreen}{f}$};
\end{tikzpicture}%
\ \ .
\]
As you see, this map $f$ is neither injective nor surjective, thus certainly
not bijective. However, I claim that I can \textbf{make} it bijective, by
appropriately tweaking its domain $R$ and its target $S$ as well as the map
$f$ itself. Namely:

\begin{itemize}
\item First, I make $f$ surjective. To do so, I replace the target\footnote{If
$g:U\rightarrow V$ is a map from a set $U$ to a set $V$, then $V$ is called
the \textbf{target} (or \textbf{codomain}) of $g$.} $S$ by the image $f\left(
R\right)  =\left\{  f\left(  r\right)  \ \mid\ r\in R\right\}  $ of the map
$f$. This way, I throw away all elements of $S$ that are not taken as values
by the map $f$. The resulting map%
\begin{align*}
\widetilde{f}:R  &  \rightarrow f\left(  R\right)  ,\\
r  &  \mapsto f\left(  r\right)
\end{align*}
(which differs from $f$ only in its choice of target) is thus surjective.

\item Next, I make $f$ (or, more precisely, $\widetilde{f}$) injective. To do
so, I equate every pair of elements $a,b\in R$ that satisfy $f\left(
a\right)  =f\left(  b\right)  $. The rigorous way to do so is to replace the
elements of $R$ by their equivalence classes with respect to an appropriately
chosen equivalence relation. To wit: We define a binary relation $\sim$ on the
set $R$ by stipulating that two elements $a,b\in R$ should satisfy $a\sim b$
if and only if $f\left(  a\right)  =f\left(  b\right)  $. This relation $\sim$
is an equivalence relation\footnote{For example, it is transitive because if
three elements $a,b,c\in R$ satisfy $f\left(  a\right)  =f\left(  b\right)  $
and $f\left(  b\right)  =f\left(  c\right)  $, then $f\left(  a\right)
=f\left(  c\right)  $.}, and its equivalence classes will be called
$f$\textbf{-classes}. We will use the notation $\overline{r}$ for the
$f$-class that contains a given element $r\in R$. We let $R/f$ denote the set
of all $f$-classes.

Now, consider the map%
\begin{align*}
f^{\prime}:R/f  &  \rightarrow f\left(  R\right)  ,\\
\overline{r}  &  \mapsto f\left(  r\right)  ,
\end{align*}
which sends each $f$-class $\overline{r}$ to the value $f\left(  r\right)  $.
This map $f^{\prime}$ is well-defined, since $f\left(  r\right)  $ depends not
on the element $r$ but only on its $f$-class $\overline{r}$ (because if two
elements $a,b\in R$ have the same $f$-class, then $a\sim b$ and thus $f\left(
a\right)  =f\left(  b\right)  $ by the very definition of $\sim$).

Just like $\widetilde{f}$, the map $f^{\prime}$ is surjective (since every
element of its target $f\left(  R\right)  $ is taken as a value by $f$, and
thus also by $f^{\prime}$). But $f^{\prime}$ is also injective, since any two
elements $a,b$ of $R$ that satisfy $f\left(  a\right)  =f\left(  b\right)  $
have already been merged into the same $f$-class in $R/f$. Thus, $f^{\prime}$
is both injective and surjective, hence bijective.
\end{itemize}

We might call $f^{\prime}$ the \textbf{bijectivization} of $f$ (although there
does not seem to be a standard name for $f^{\prime}$). In our above example,
this map $f^{\prime}$ looks as follows:%
\[%
\begin{tikzpicture}
\filldraw[fill=yellow!20!white] (0, 3) ellipse (1 and 2);
\node(A1) at (0, 4) {$\set{1,3}$};
\node(A2) at (0, 3) {$\set{2}$};
\node(A3) at (0, 2) {$\set{4,5}$};
\filldraw[fill=yellow!20!white] (3, 3) ellipse (1 and 2);
\node(Aname) at (0, 0.6) {$R/f$};
\node(B1) at (3, 4) {$1$};
\node(B2) at (3, 3) {$2$};
\node(B3) at (3, 2) {$4$};
\node(Bname) at (3, 0.6) {$f\tup{R}$};
\draw[->, very thick, darkgreen] (A1) -- (B2);
\draw[->, very thick, darkgreen] (A2) -- (B3);
\draw[->, very thick, darkgreen] (A3) -- (B1);
\node(o'name) at (1.5, 2.1) {$\color{darkgreen}{f'}$};
\end{tikzpicture}%
\ \ .
\]

Moreover, the maps $f$ and $f^{\prime}$ fit together into a nice picture with
two other rather natural maps:%
\[%
\begin{tikzpicture}
\filldraw[fill=yellow!20!white] (0, 3) ellipse (1 and 2);
\node(A1) at (0, 4) {$\set{1,3}$};
\node(A2) at (0, 3) {$\set{2}$};
\node(A3) at (0, 2) {$\set{4,5}$};
\filldraw[fill=yellow!20!white] (3, 3) ellipse (1 and 2);
\node(Aname) at (0, 0.6) {$R/f$};
\node(B1) at (3, 4) {$1$};
\node(B2) at (3, 3) {$2$};
\node(B3) at (3, 2) {$4$};
\node(Bname) at (3, 0.6) {$f\tup{R}$};
\draw[->, very thick, darkgreen] (A1) -- (B2);
\draw[->, very thick, darkgreen] (A2) -- (B3);
\draw[->, very thick, darkgreen] (A3) -- (B1);
\filldraw[fill=cyan!20!white] (0, 9) ellipse (1 and 3);
\node(C1) at (0, 11) {$1$};
\node(C2) at (0, 10) {$2$};
\node(C3) at (0, 9) {$3$};
\node(C4) at (0, 8) {$4$};
\node(C5) at (0, 7) {$5$};
\node(Cname) at (0, 12.4) {$R$};
\filldraw[fill=cyan!20!white] (3, 8.5) ellipse (1 and 3);
\node(D1) at (3, 10) {$1$};
\node(D2) at (3, 9) {$2$};
\node(D3) at (3, 8) {$3$};
\node(D4) at (3, 7) {$4$};
\node(Dname) at (3, 11.9) {$S$};
\draw[->, very thick, darkgreen] (C1) -- (D2);
\draw[->, very thick, darkgreen] (C2) -- (D4);
\draw[->, very thick, darkgreen] (C3) -- (D2);
\draw[->, very thick, darkgreen] (C4) -- (D1);
\draw[->, very thick, darkgreen] (C5) -- (D1);
\draw[->, very thick, red] (C1) edge[bend right=50] (A1);
\draw[->, very thick, red] (C2) edge[bend right=50] (A2);
\draw[->, very thick, red] (C3) edge[bend right=50] (A1);
\draw[->, very thick, red] (C4) edge[bend right=50] (A3);
\draw[->, very thick, red] (C5) edge[bend right=50] (A3);
\draw[->, very thick, red] (B1) edge[bend right=50] (D1);
\draw[->, very thick, red] (B2) edge[bend right=50] (D2);
\draw[->, very thick, red] (B3) edge[bend right=50] (D4);
\node(phinameL) at (-1.8, 5) {$\color{red}{\pi}$};
\node(phinameR) at (4.6, 5) {$\color{red}{\iota}$};
\node(oname) at (1.5, 10.5) {$\color{darkgreen}{f}$};
\node(o'name) at (1.5, 2.1) {$\color{darkgreen}{f'}$};
\end{tikzpicture}%
\ \ .
\]
Here, $\pi:R\rightarrow R/f$ is the \textbf{canonical projection} (i.e., the
map that sends each $r\in R$ to its $f$-class $\overline{r}$), and
$\iota:f\left(  R\right)  \rightarrow S$ is the \textbf{canonical inclusion}
(i.e., the map that sends each $s\in f\left(  R\right)  $ to $s$). These four
maps $f,f^{\prime},\pi,\iota$ satisfy%
\[
f=\iota\circ f^{\prime}\circ\pi,
\]
which means (in the language of Subsection \ref{subsec.quotring.uniprop2})
that the diagram%
\[%
\xymatrix@C=4pc{
R \ar[d]_{\pi} \ar[r]^f & S \\
R/f \ar[r]_{f'} & f\tup{R} \ar[u]^{\iota}
}%
\]
is commutative.

For the sake of completeness, let us state this all as a theorem:

\begin{theorem}
[First Isomorphism Theorem for sets]\label{thm.1it.set}Let $R$ and $S$ be any
two sets, and let $f:R\rightarrow S$ be any map.

Let $\sim$ be the binary relation on the set $R$ defined by requiring that two
elements $a,b\in R$ satisfy $a\sim b$ if and only if $f\left(  a\right)
=f\left(  b\right)  $.

\begin{enumerate}
\item[\textbf{(a)}] This relation $\sim$ is an equivalence relation.
\end{enumerate}

Let us refer to the equivalence classes of this equivalence relation $\sim$ as
the $f$\textbf{-classes}. Let $R/f$ denote the set of all $f$-classes. For any
$r\in R$, we let $\overline{r}$ denote the $f$-class that contains $r$.

\begin{enumerate}
\item[\textbf{(b)}] The image $f\left(  R\right)  :=\left\{  f\left(
r\right)  \ \mid\ r\in R\right\}  $ of $f$ is a subset of $S$.

\item[\textbf{(c)}] The map%
\begin{align*}
f^{\prime}:R/f  &  \rightarrow f\left(  R\right)  ,\\
\overline{r}  &  \mapsto f\left(  r\right)
\end{align*}
is well-defined and bijective.

\item[\textbf{(d)}] Let $\pi:R\rightarrow R/f$ denote the \textbf{canonical
projection} (i.e., the map that sends each $r\in R$ to its $f$-class
$\overline{r}$). Let $\iota:f\left(  R\right)  \rightarrow S$ denote the
\textbf{canonical inclusion} (i.e., the map that sends each $s\in f\left(
R\right)  $ to $s$). Then, the map $f^{\prime}$ defined in part \textbf{(c)}
satisfies
\[
f=\iota\circ f^{\prime}\circ\pi.
\]
In other words, the diagram%
\begin{equation}%
\xymatrix@C=4pc{
R \ar[d]_{\pi} \ar[r]^f & S \\
R/f \ar[r]_{f'} & f\tup{R} \ar[u]^{\iota}
}
\label{eq.thm.1it.set.d.diag}%
\end{equation}
is commutative.
\end{enumerate}
\end{theorem}

\begin{proof}
Part \textbf{(b)} is obvious. We explained the proofs of parts \textbf{(a)}
and \textbf{(c)} before even stating this theorem. Just for the sake of
completeness, we shall now repeat the proof of part \textbf{(c)}, and then
prove part \textbf{(d)}: \medskip

\textbf{(c)} If $a,b\in R$ are two elements satisfying $\overline{a}%
=\overline{b}$, then $f\left(  a\right)  =f\left(  b\right)  $%
\ \ \ \ \footnote{\textit{Proof.} Let $a,b\in R$ be two elements satisfying
$\overline{a}=\overline{b}$. The equality $\overline{a}=\overline{b}$ shows
that $a$ and $b$ belong to the same $f$-class (since $\overline{a}$ denotes
the $f$-class that contains $a$, whereas $\overline{b}$ denotes the $f$-class
that contains $b$). In other words, $a$ and $b$ belong to the same equivalence
class of the equivalence relation $\sim$ (since the $f$-classes are the
equivalence classes of this equivalence relation $\sim$). In other words,
$a\sim b$. In other words, $f\left(  a\right)  =f\left(  b\right)  $ (by the
definition of the relation $\sim$).}. In other words, for any element $r\in
R$, the value $f\left(  r\right)  $ depends only on the $f$-class
$\overline{r}$ and not on $r$ itself. Hence, the map%
\begin{align*}
f^{\prime}:R/f  &  \rightarrow f\left(  R\right)  ,\\
\overline{r}  &  \mapsto f\left(  r\right)
\end{align*}
is well-defined (since each element of $R/f$ can be written as $\overline{r}$
for some $r\in R$). It remains to prove that this map is bijective. To that
purpose, we shall now prove that it is injective and surjective.

\textit{Injectivity:} Let $x,y\in R/f$ be two elements of $R/f$ satisfying
$f^{\prime}\left(  x\right)  =f^{\prime}\left(  y\right)  $. We shall prove
that $x=y$.

We know that $x$ is an element of $R/f$. In other words, $x$ is an $f$-class
(since $R/f$ is the set of all $f$-classes). Thus, we can write $x$ in the
form $x=\overline{a}$ for some $a\in R$. Likewise, we can write $y$ in the
form $y=\overline{b}$ for some $b\in R$. Consider these $a$ and $b$.

From $x=\overline{a}$, we obtain $f^{\prime}\left(  x\right)  =f^{\prime
}\left(  \overline{a}\right)  =f\left(  a\right)  $ (by the definition of
$f^{\prime}$). Similarly, $f^{\prime}\left(  y\right)  =f\left(  b\right)  $
(since $y=\overline{b}$). Hence, $f\left(  a\right)  =f^{\prime}\left(
x\right)  =f^{\prime}\left(  y\right)  =f\left(  b\right)  $. In other words,
$a\sim b$ (by the definition of the relation $\sim$). In other words, $a$ and
$b$ belong to the same equivalence class of the equivalence relation $\sim$.
In other words, $a$ and $b$ belong to the same $f$-class (since the
$f$-classes are the equivalence classes of this equivalence relation $\sim$).
In other words, $\overline{a}=\overline{b}$ (since $\overline{a}$ denotes the
$f$-class that contains $a$, whereas $\overline{b}$ denotes the $f$-class that
contains $b$). In other words, $x=y$ (since $x=\overline{a}$ and
$y=\overline{b}$).

Forget that we fixed $x$ and $y$. We thus have shown that if $x,y\in R/f$ are
two elements of $R/f$ satisfying $f^{\prime}\left(  x\right)  =f^{\prime
}\left(  y\right)  $, then $x=y$. In other words, the map $f^{\prime}$ is injective.

\textit{Surjectivity:} Let $z\in f\left(  R\right)  $. Thus, $z=f\left(
r\right)  $ for some $r\in R$. Now, the $f$-class $\overline{r}\in R/f$
satisfies $f^{\prime}\left(  \overline{r}\right)  =f\left(  r\right)  $ (by
the definition of $f^{\prime}$). Comparing this with $z=f\left(  r\right)  $,
we obtain $z=f^{\prime}\left(  \overline{r}\right)  $. This shows that the map
$f^{\prime}$ takes $z$ as a value.

Forget that we fixed $z$. We thus have shown that the map $f^{\prime}$ takes
each $z\in f\left(  R\right)  $ as a value. In other words, the map
$f^{\prime}$ is surjective.

We now have proved that $f^{\prime}$ is injective and surjective. Thus,
$f^{\prime}$ is bijective, and the proof of Theorem \ref{thm.1it.set}
\textbf{(c)} is complete. \medskip

\textbf{(d)} For each $r\in R$, we have%
\begin{align*}
&  \left(  \iota\circ f^{\prime}\circ\pi\right)  \left(  r\right) \\
&  =\iota\left(  f^{\prime}\left(  \pi\left(  r\right)  \right)  \right) \\
&  =f^{\prime}\left(  \pi\left(  r\right)  \right)
\ \ \ \ \ \ \ \ \ \ \left(  \text{since the definition of }\iota\text{ yields
}\iota\left(  s\right)  =s\text{ for each }s\in f\left(  R\right)  \right) \\
&  =f^{\prime}\left(  \overline{r}\right)  \ \ \ \ \ \ \ \ \ \ \left(
\text{since the definition of }\pi\text{ yields }\pi\left(  r\right)
=\overline{r}\right) \\
&  =f\left(  r\right)  \ \ \ \ \ \ \ \ \ \ \left(  \text{by the definition of
}f^{\prime}\right)  .
\end{align*}
In other words, $\iota\circ f^{\prime}\circ\pi=f$. Thus, $f=\iota\circ
f^{\prime}\circ\pi$. In other words, the diagram (\ref{eq.thm.1it.set.d.diag})
is commutative. This proves Theorem \ref{thm.1it.set} \textbf{(d)}.
\end{proof}

\subsubsection{The First Isomorphism Theorem for rings}

Now, let us extend the First Isomorphism Theorem to rings and ring morphisms
instead of arbitrary sets and maps.

\begin{theorem}
[First Isomorphism Theorem for rings, elementwise form]\label{thm.1it.ring1}%
Let $R$ and $S$ be two rings, and let $f:R\rightarrow S$ be a ring morphism. Then:

\begin{enumerate}
\item[\textbf{(a)}] The kernel $\operatorname*{Ker}f$ is an ideal of $R$.
Thus, $R/\operatorname*{Ker}f$ is a quotient ring of $R$. As a set,
$R/\operatorname*{Ker}f$ is precisely the set $R/f$ defined in Theorem
\ref{thm.1it.set}. The $f$-classes (as defined in Theorem \ref{thm.1it.set})
are precisely the cosets of $\operatorname*{Ker}f$.

\item[\textbf{(b)}] The image $f\left(  R\right)  :=\left\{  f\left(
r\right)  \ \mid\ r\in R\right\}  $ of $f$ is a subring of $S$.

\item[\textbf{(c)}] The map%
\begin{align*}
f^{\prime}:R/\operatorname*{Ker}f  &  \rightarrow f\left(  R\right)  ,\\
\overline{r}  &  \mapsto f\left(  r\right)
\end{align*}
is well-defined and is a ring isomorphism.

\item[\textbf{(d)}] This map $f^{\prime}$ is precisely the map $f^{\prime}$
defined in Theorem \ref{thm.1it.set} \textbf{(c)}.

\item[\textbf{(e)}] Let $\pi:R\rightarrow R/\operatorname*{Ker}f$ denote the
\textbf{canonical projection} (i.e., the map that sends each $r\in R$ to its
coset $\overline{r}$). Let $\iota:f\left(  R\right)  \rightarrow S$ denote the
\textbf{canonical inclusion} (i.e., the map that sends each $s\in f\left(
R\right)  $ to $s$). Then, the map $f^{\prime}$ defined in part \textbf{(c)}
satisfies
\[
f=\iota\circ f^{\prime}\circ\pi.
\]
In other words, the diagram%
\begin{equation}%
\xymatrix@C=4pc{
R \ar[d]_{\pi} \ar[r]^f & S \\
R/\Ker f \ar[r]_{f'} & f\tup{R} \ar[u]^{\iota}
}
\label{eq.thm.1it.ring1.diag}%
\end{equation}
is commutative.

\item[\textbf{(f)}] We have $R/\operatorname*{Ker}f\cong f\left(  R\right)  $
as rings.
\end{enumerate}
\end{theorem}

\begin{proof}
\textbf{(a)} We know that $\operatorname*{Ker}f$ is an ideal of $R$ (by
Theorem \ref{thm.ringmor.ker-ideal}), and therefore $R/\operatorname*{Ker}f$
is a quotient ring of $R$.

Let us next prove that the $f$-classes (as defined in Theorem
\ref{thm.1it.set}) are precisely the cosets of $\operatorname*{Ker}f$.

Indeed, let $\sim$ be the equivalence relation on $R$ defined in Theorem
\ref{thm.1it.set}. Then, the $f$-classes are defined as the equivalence
classes of this relation $\sim$. For any $a,b\in R$, we have $f\left(
b\right)  -f\left(  a\right)  =f\left(  b-a\right)  $ (since $f$ is a ring
morphism and thus respects differences). For any two elements $a,b\in R$, we
have the chain of equivalences%
\begin{align*}
\  &  \ \ \ \ \ \ \ \ \ \left(  a\sim b\right) \\
&  \Longleftrightarrow\ \left(  f\left(  a\right)  =f\left(  b\right)
\right)  \ \ \ \ \ \ \ \ \ \ \left(  \text{by the definition of }\sim\right)
\\
&  \Longleftrightarrow\ \left(  f\left(  b\right)  -f\left(  a\right)
=0\right) \\
&  \Longleftrightarrow\ \left(  f\left(  b-a\right)  =0\right)
\ \ \ \ \ \ \ \ \ \ \left(  \text{since }f\left(  b\right)  -f\left(
a\right)  =f\left(  b-a\right)  \right) \\
&  \Longleftrightarrow\ \left(  b-a\in\operatorname*{Ker}f\right)
\ \ \ \ \ \ \ \ \ \ \left(  \text{by the definition of }\operatorname*{Ker}%
f\right)  .
\end{align*}
However, if $a\in R$ is arbitrary, then%
\begin{align*}
&  \left(  \text{the }f\text{-class that contains }a\right) \\
&  =\left(  \text{the equivalence class of the relation }\sim\text{ that
contains }a\right) \\
&  \ \ \ \ \ \ \ \ \ \ \ \ \ \ \ \ \ \ \ \ \left(  \text{since the
}f\text{-classes are the equivalence classes of }\sim\right) \\
&  =\left\{  b\in R\ \mid\ a\sim b\right\}  \ \ \ \ \ \ \ \ \ \ \left(
\text{by the definition of equivalence classes}\right) \\
&  =\left\{  b\in R\ \mid\ b-a\in\operatorname*{Ker}f\right\} \\
&  \ \ \ \ \ \ \ \ \ \ \ \ \ \ \ \ \ \ \ \ \left(
\begin{array}
[c]{c}%
\text{by the equivalence }\left(  a\sim b\right)  \Longleftrightarrow\left(
b-a\in\operatorname*{Ker}f\right) \\
\text{that we proved above}%
\end{array}
\right) \\
&  =\left\{  b\in R\ \mid\ b\in a+\operatorname*{Ker}f\right\} \\
&  =a+\operatorname*{Ker}f\\
&  =\left(  \text{the coset of }\operatorname*{Ker}f\text{ that contains
}a\right)
\end{align*}
(since the coset of $\operatorname*{Ker}f$ that contains $a$ is
$a+\operatorname*{Ker}f$ by definition). Thus, the $f$-classes are precisely
the cosets of $\operatorname*{Ker}f$.

In other words, the cosets of $\operatorname*{Ker}f$ are precisely the
$f$-classes. Hence, the set $R/\operatorname*{Ker}f$ is precisely the set
$R/f$ (since the former set consists of the cosets of $\operatorname*{Ker}f$,
while the latter set consists of the $f$-classes). This concludes the proof of
Theorem \ref{thm.1it.ring1} \textbf{(a)}. \medskip

\textbf{(b)} This is just Proposition \ref{prop.ringmor.Im-subring}. \medskip

\textbf{(c)} Theorem \ref{thm.1it.ring1} \textbf{(a)} yields that
$R/f=R/\operatorname*{Ker}f$, and that the $f$-classes are precisely the
cosets of $\operatorname*{Ker}f$. Hence, the meaning of the notation
$\overline{r}$ in Theorem \ref{thm.1it.set} is identical with the meaning of
this notation in Theorem \ref{thm.1it.ring1} (indeed, the former denotes the
$f$-class that contains $r$, whereas the latter denotes the coset of
$\operatorname*{Ker}f$ that contains $r$; but as we just said, the $f$-classes
are precisely the cosets of $\operatorname*{Ker}f$). Hence, Theorem
\ref{thm.1it.set} \textbf{(c)} shows that the map%
\begin{align*}
f^{\prime}:R/f  &  \rightarrow f\left(  R\right)  ,\\
\overline{r}  &  \mapsto f\left(  r\right)
\end{align*}
is well-defined and bijective. Since $R/f=R/\operatorname*{Ker}f$ (as sets),
we can restate this as follows: The map
\begin{align*}
f^{\prime}:R/\operatorname*{Ker}f  &  \rightarrow f\left(  R\right)  ,\\
\overline{r}  &  \mapsto f\left(  r\right)
\end{align*}
is well-defined and bijective. It remains to prove that this map $f^{\prime}$
is a ring isomorphism.

We shall first show that $f^{\prime}$ is a ring morphism. Indeed, this is an
easy consequence of Theorem \ref{thm.quotring.uniprop1} (applied to
$I=\operatorname*{Ker}f$), since we have $f\left(  \operatorname*{Ker}%
f\right)  =0$ (by the definition of $\operatorname*{Ker}f$). Alternatively, we
can prove this by hand as follows:

We must show that $f^{\prime}$ is a ring morphism, i.e., that $f^{\prime}$
respects addition, multiplication, zero and unity.

To see that $f^{\prime}$ respects multiplication, we must show that
$f^{\prime}\left(  xy\right)  =f^{\prime}\left(  x\right)  \cdot f^{\prime
}\left(  y\right)  $ for any $x,y\in R/\operatorname*{Ker}f$. So let $x,y\in
R/\operatorname*{Ker}f$ be arbitrary. Then, we can write $x$ and $y$ as
$x=\overline{a}$ and $y=\overline{b}$ for two elements $a,b\in R$. Consider
these $a,b$. From $x=\overline{a}$ and $y=\overline{b}$, we obtain
$xy=\overline{a}\cdot\overline{b}=\overline{ab}$ (by the definition of the
product on $R/\operatorname*{Ker}f$), so that
\begin{align*}
f^{\prime}\left(  xy\right)   &  =f^{\prime}\left(  \overline{ab}\right)
=f\left(  ab\right)  \ \ \ \ \ \ \ \ \ \ \left(  \text{by the definition of
}f^{\prime}\right) \\
&  =f\left(  a\right)  \cdot f\left(  b\right)  \ \ \ \ \ \ \ \ \ \ \left(
\text{since }f\text{ is a ring morphism}\right)  .
\end{align*}
Comparing this with%
\[
f^{\prime}\left(  \underbrace{x}_{=\overline{a}}\right)  \cdot f^{\prime
}\left(  \underbrace{y}_{=\overline{b}}\right)  =\underbrace{f^{\prime}\left(
\overline{a}\right)  }_{\substack{=f\left(  a\right)  \\\text{(by
the}\\\text{definition of }f^{\prime}\text{)}}}\cdot\underbrace{f^{\prime
}\left(  \overline{b}\right)  }_{\substack{=f\left(  b\right)  \\\text{(by
the}\\\text{definition of }f^{\prime}\text{)}}}=f\left(  a\right)  \cdot
f\left(  b\right)  ,
\]
we obtain $f^{\prime}\left(  xy\right)  =f^{\prime}\left(  x\right)  \cdot
f^{\prime}\left(  y\right)  $, just as desired. Thus, we have shown that
$f^{\prime}$ respects multiplication. Similarly, $f^{\prime}$ satisfies all
the other axioms in the definition of a ring morphism.

Thus, we know that $f^{\prime}$ is a ring morphism. Since $f^{\prime}$ is also
invertible (because $f^{\prime}$ is bijective), we conclude that $f^{\prime}$
is an invertible ring morphism. Thus, $f^{\prime}$ is a ring isomorphism
(since Proposition \ref{prop.ringmor.invertible-iso} shows that any invertible
ring morphism is a ring isomorphism). Thus, the proof of Theorem
\ref{thm.1it.ring1} \textbf{(c)} is finished. \medskip

\textbf{(d)} As we have seen in our above proof of Theorem \ref{thm.1it.ring1}
\textbf{(c)}, we have $R/f=R/\operatorname*{Ker}f$, and the meaning of the
notation $\overline{r}$ in Theorem \ref{thm.1it.set} is identical with the
meaning of this notation in Theorem \ref{thm.1it.ring1}. Thus, the map
$f^{\prime}$ in Theorem \ref{thm.1it.ring1} \textbf{(c)} is precisely the map
$f^{\prime}$ defined in Theorem \ref{thm.1it.set} \textbf{(c)}. This proves
Theorem \ref{thm.1it.ring1} \textbf{(d)}. \medskip

\textbf{(e)} As we have seen in our above proof of Theorem \ref{thm.1it.ring1}
\textbf{(c)}, we have $R/f=R/\operatorname*{Ker}f$, and the meaning of the
notation $\overline{r}$ in Theorem \ref{thm.1it.set} is identical with the
meaning of this notation in Theorem \ref{thm.1it.ring1}. Therefore, the maps
$\pi$ and $\iota$ defined in Theorem \ref{thm.1it.ring1} \textbf{(e)} are
precisely the maps $\pi$ and $\iota$ in Theorem \ref{thm.1it.set}
\textbf{(d)}. Hence, the claim of Theorem \ref{thm.1it.ring1} \textbf{(e)}
follows immediately from Theorem \ref{thm.1it.set} \textbf{(d)} (since
$R/f=R/\operatorname*{Ker}f$). \medskip

\textbf{(f)} This follows directly from Theorem \ref{thm.1it.ring1}
\textbf{(c)}.
\end{proof}

As our proof has shown, Theorem \ref{thm.1it.ring1} \textbf{(c)} is merely a
partial improvement on the universal property of quotient rings (Theorem
\ref{thm.quotring.uniprop1}): The latter yields a ring morphism, while the
former produces a ring \textbf{iso}morphism (but in a less general setup:
$R/\operatorname*{Ker}f$ instead of $R/I$). Nevertheless, it is a useful
result, as it can be used to identify certain quotient rings as (isomorphic
copies of) known rings.

\bigskip

Here are some examples for what can be done with the first isomorphism theorem:

\begin{itemize}
\item Consider the map\footnote{See Subsection
\ref{subsec.rings.subrings.exas} for the meaning of the notation
$\mathbb{Q}^{n\leq n}$ (and, more generally, $R^{n\leq n}$ when $R$ is any
ring).}%
\begin{align*}
f:\mathbb{Q}^{4\leq4}  &  \rightarrow\mathbb{Q}^{2\leq2},\\
\left(
\begin{array}
[c]{cccc}%
a & b & c & d\\
0 & u & v & w\\
0 & 0 & x & y\\
0 & 0 & 0 & z
\end{array}
\right)   &  \mapsto\left(
\begin{array}
[c]{cc}%
u & v\\
0 & x
\end{array}
\right)  ,
\end{align*}
which removes the \textquotedblleft outer shell\textquotedblright\ (i.e., the
first and the fourth rows and columns) from an upper-triangular $4\times
4$-matrix. This map $f$ is a ring morphism\footnote{Proving this is a nice
exercise in matrix multiplication! It is obvious that $f$ respects addition,
zero and unity, but you might be skeptical that it respects multiplication.
(And indeed, the analogous map
\begin{align*}
F:\mathbb{Q}^{4\times4}  &  \rightarrow\mathbb{Q}^{2\times2},\\
\left(
\begin{array}
[c]{cccc}%
a & b & c & d\\
a^{\prime} & b^{\prime} & c^{\prime} & d^{\prime}\\
a^{\prime\prime} & b^{\prime\prime} & c^{\prime\prime} & d^{\prime\prime}\\
a^{\prime\prime\prime} & b^{\prime\prime\prime} & c^{\prime\prime\prime} &
d^{\prime\prime\prime}%
\end{array}
\right)   &  \mapsto\left(
\begin{array}
[c]{cc}%
b^{\prime} & c^{\prime}\\
b^{\prime\prime} & c^{\prime\prime}%
\end{array}
\right)  ,
\end{align*}
which removes the \textquotedblleft outer shell\textquotedblright\ from an
arbitrary (not upper-triangular) $4\times4$-matrix, does not respect
multiplication.) You can convince yourself of this property of $f$ by a
straightforward computation:%
\[
\left(
\begin{array}
[c]{cccc}%
a & b & c & d\\
0 & u & v & w\\
0 & 0 & x & y\\
0 & 0 & 0 & z
\end{array}
\right)  \left(
\begin{array}
[c]{cccc}%
a^{\prime} & b^{\prime} & c^{\prime} & d^{\prime}\\
0 & u^{\prime} & v^{\prime} & w^{\prime}\\
0 & 0 & x^{\prime} & y^{\prime}\\
0 & 0 & 0 & z^{\prime}%
\end{array}
\right)  =\left(
\begin{array}
[c]{cccc}%
aa^{\prime} & ab^{\prime}+bu^{\prime} & ac^{\prime}+bv^{\prime}+cx^{\prime} &
ad^{\prime}+bw^{\prime}+cy^{\prime}+dz^{\prime}\\
0 & uu^{\prime} & uv^{\prime}+vx^{\prime} & uw^{\prime}+vy^{\prime}%
+wz^{\prime}\\
0 & 0 & xx^{\prime} & xy^{\prime}+yz^{\prime}\\
0 & 0 & 0 & zz^{\prime}%
\end{array}
\right)
\]
(note the $uu^{\prime}$, $uv^{\prime}+vx^{\prime}$, $0$ and $xx^{\prime}$
entries, which are precisely the entries of $\left(
\begin{array}
[c]{cc}%
u & v\\
0 & x
\end{array}
\right)  \left(
\begin{array}
[c]{cc}%
u^{\prime} & v^{\prime}\\
0 & x^{\prime}%
\end{array}
\right)  $).}.

The kernel of this morphism $f$ is
\begin{align*}
\operatorname*{Ker}f  &  =\left\{  \left(
\begin{array}
[c]{cccc}%
a & b & c & d\\
0 & u & v & w\\
0 & 0 & x & y\\
0 & 0 & 0 & z
\end{array}
\right)  \in\mathbb{Q}^{4\leq4}\ \mid\ \left(
\begin{array}
[c]{cc}%
u & v\\
0 & x
\end{array}
\right)  =0\right\} \\
&  =\left\{  \left(
\begin{array}
[c]{cccc}%
a & b & c & d\\
0 & u & v & w\\
0 & 0 & x & y\\
0 & 0 & 0 & z
\end{array}
\right)  \in\mathbb{Q}^{4\leq4}\ \mid\ u=v=x=0\right\} \\
&  =\left\{  \left(
\begin{array}
[c]{cccc}%
a & b & c & d\\
0 & 0 & 0 & w\\
0 & 0 & 0 & y\\
0 & 0 & 0 & z
\end{array}
\right)  \in\mathbb{Q}^{4\leq4}\right\}  .
\end{align*}
$\allowbreak$

So you can conclude right away that $\operatorname*{Ker}f$ is an ideal of
$\mathbb{Q}^{4\leq4}$. Moreover, the image $f\left(  \mathbb{Q}^{4\leq
4}\right)  $ is the whole $\mathbb{Q}^{2\leq2}$ (that is, the map $f$ is surjective).

The First Isomorphism theorem (Theorem \ref{thm.1it.ring1} \textbf{(c)})
yields a ring isomorphism%
\begin{align*}
f^{\prime}:\mathbb{Q}^{4\leq4}/\operatorname*{Ker}f  &  \rightarrow f\left(
\mathbb{Q}^{4\leq4}\right)  ,\\
\overline{r}  &  \mapsto f\left(  r\right)  .
\end{align*}
In other words, it yields a ring isomorphism%
\begin{align*}
f^{\prime}:\mathbb{Q}^{4\leq4}/\operatorname*{Ker}f  &  \rightarrow
\mathbb{Q}^{2\leq2},\\
\overline{\left(
\begin{array}
[c]{cccc}%
a & b & c & d\\
0 & u & v & w\\
0 & 0 & x & y\\
0 & 0 & 0 & z
\end{array}
\right)  }  &  \mapsto\left(
\begin{array}
[c]{cc}%
u & v\\
0 & x
\end{array}
\right)
\end{align*}
(since $f\left(  \mathbb{Q}^{4\leq4}\right)  =\mathbb{Q}^{2\leq2}$). In
particular, $\mathbb{Q}^{4\leq4}/\operatorname*{Ker}f\cong\mathbb{Q}^{2\leq2}$.

\item We have not properly defined polynomials yet, but once we will, you will
be inundated with good examples for the First Isomorphism Theorem. Many of
these examples will have the form%
\[
\left(  \text{a polynomial ring}\right)  /\left(  \text{an ideal}\right)
\cong\left(  \text{a ring of numbers}\right)  .
\]
For instance, recalling that $\mathbb{R}\left[  x\right]  $ is the ring of all
polynomials in one indeterminate $x$ with real coefficients, we have a ring
morphism%
\begin{align*}
f:\mathbb{R}\left[  x\right]   &  \rightarrow\mathbb{C},\\
p  &  \mapsto p\left(  i\right)
\end{align*}
(which sends each polynomial $p\in\mathbb{R}\left[  x\right]  $ to its value
at the imaginary unit $i=\sqrt{-1}$). This morphism is surjective (that is,
$f\left(  \mathbb{R}\left[  x\right]  \right)  =\mathbb{C}$) and has kernel
$\operatorname*{Ker}f=\left(  x^{2}+1\right)  \mathbb{R}\left[  x\right]  $
(the principal ideal generated by $x^{2}+1$), so that the First Isomorphism
theorem (Theorem \ref{thm.1it.ring1} \textbf{(f)}) yields
\[
\mathbb{R}\left[  x\right]  /\left(  \left(  x^{2}+1\right)  \mathbb{R}\left[
x\right]  \right)  \cong\mathbb{C}.
\]
Informally, this is saying that if you are working with polynomials in an
indeterminate $x$ over $\mathbb{R}$, but you equate the polynomial $x^{2}+1$
to zero (that is, you pretend that $x^{2}=-1$), then you obtain the complex
numbers. This is the rigorous concept behind the classical idea that
\textquotedblleft the complex numbers are what you get if you start with the
real numbers and adjoin a root of the polynomial $x^{2}+1$\textquotedblright.
We will make this precise in a later chapter.
\end{itemize}

\subsubsection{A few remarks on the first isomorphism theorem}

Theorem \ref{thm.1it.ring1} (specifically, its parts \textbf{(c)} and
\textbf{(e)}) is commonly called the \textbf{first isomorphism theorem for
rings}, and is one of the major sources of ring isomorphisms in some parts of
abstract algebra. Due to its importance, a few more comments on it are worth making.

The commutative diagram (\ref{eq.thm.1it.ring1.diag}) in Theorem
\ref{thm.1it.ring1} \textbf{(e)} can be rewritten in a somewhat more
expressive form:%
\[%
\xymatrix@C=5pc{
R \arsurj[ddr]_{\pi} \ar[rrr]^f & & & S \\
\\
& R / \Ker f \ar[r]^\cong_{f'} & f\tup{R} \arinjrev[ruu]_{\iota}
}\ \ .%
\]
Let me explain what you are seeing here: On top is the original ring morphism
$f:R\rightarrow S$. The other four arrows are

\begin{itemize}
\item the canonical projection $\pi:R\rightarrow R/\operatorname*{Ker}f$,
sending each $r\in R$ to its residue class $\overline{r}=r+\operatorname*{Ker}%
f\in R/\operatorname*{Ker}f$;

\item the canonical inclusion $\iota:f\left(  R\right)  \rightarrow S$ (which
just sends each element to itself);

\item the morphism $f^{\prime}:R/\operatorname*{Ker}f\rightarrow f\left(
R\right)  $ claimed by Theorem \ref{thm.1it.ring1} \textbf{(c)}.
\end{itemize}

The special shapes of the arrows signify certain properties:

\begin{itemize}
\item An arrow of shape $\hookrightarrow$ stands for an injective map. (And
indeed, the canonical inclusion $\iota:f\left(  R\right)  \rightarrow S$ is injective.)

\item An arrow of shape $\twoheadrightarrow$ stands for a surjective map. (And
indeed, the canonical projection $\pi$ is surjective.)

\item An arrow with a $\cong$ sign above (or below) it stands for an
isomorphism. (And indeed, our $f^{\prime}$ is an isomorphism.)
\end{itemize}

Note that all four arrows in our diagram are ring morphisms; we thus say that
our diagram is a \textbf{diagram of rings}.

Thus, the first isomorphism theorem for rings shows that each ring morphism
can be decomposed (in a natural way) into a composition of a surjective ring
morphism, a ring isomorphism and an injective ring morphism.

\bigskip

In Theorem \ref{thm.1it.ring1} \textbf{(c)}, we have defined our isomorphism
$f^{\prime}$ explicitly. Alternatively, it can be characterized (uniquely) by
the equation $f=\iota\circ f^{\prime}\circ\pi$ stated in Theorem
\ref{thm.1it.ring1} \textbf{(e)}:

\begin{theorem}
[First isomorphism theorem for rings, abstract form]\label{thm.1it.ring2}Let
$R$ and $S$ be two rings. Let $f:R\rightarrow S$ be a ring morphism. Recall
that $\operatorname*{Ker}f$ is an ideal of $R$, and that $\operatorname{Im}%
f=f\left(  R\right)  $ is a subring of $S$. Let $\iota:f\left(  R\right)
\rightarrow S$ denote the canonical inclusion (i.e., the map that sends each
$s\in f\left(  R\right)  $ to $s$). Then:

\begin{enumerate}
\item[\textbf{(a)}] There is a unique ring morphism $f^{\prime}%
:R/\operatorname*{Ker}f\rightarrow f\left(  R\right)  $ that satisfies the
equation $f=\iota\circ f^{\prime}\circ\pi$ (that is, for which the diagram
(\ref{eq.thm.1it.ring1.diag}) is commutative).

\item[\textbf{(b)}] This morphism $f^{\prime}$ is a ring isomorphism.
\end{enumerate}
\end{theorem}

\begin{proof}
\textbf{(a)} Theorem \ref{thm.1it.ring1} \textbf{(e)} shows that the ring
isomorphism $f^{\prime}:R/\operatorname*{Ker}f\rightarrow f\left(  R\right)  $
constructed in Theorem \ref{thm.1it.ring1} \textbf{(c)} satisfies the equation
$f=\iota\circ f^{\prime}\circ\pi$. Hence, there exists \textbf{at least one}
ring morphism $f^{\prime}:R/\operatorname*{Ker}f\rightarrow f\left(  R\right)
$ that satisfies the equation $f=\iota\circ f^{\prime}\circ\pi$ (namely, the
isomorphism $f^{\prime}$ we were just talking about). It thus remains to prove
that this morphism $f^{\prime}$ is unique.

To show this, we consider an arbitrary ring morphism $f^{\prime}%
:R/\operatorname*{Ker}f\rightarrow f\left(  R\right)  $ that satisfies the
equation $f=\iota\circ f^{\prime}\circ\pi$. Then, for any $r\in R$, we have%
\begin{align*}
f\left(  r\right)   &  =\left(  \iota\circ f^{\prime}\circ\pi\right)  \left(
r\right)  \ \ \ \ \ \ \ \ \ \ \left(  \text{since }f=\iota\circ f^{\prime
}\circ\pi\right) \\
&  =\iota\left(  f^{\prime}\left(  \pi\left(  r\right)  \right)  \right) \\
&  =f^{\prime}\left(  \pi\left(  r\right)  \right)
\ \ \ \ \ \ \ \ \ \ \left(  \text{since the definition of }\iota\text{ yields
}\iota\left(  s\right)  =s\text{ for each }s\in f\left(  R\right)  \right) \\
&  =f^{\prime}\left(  \overline{r}\right)  \ \ \ \ \ \ \ \ \ \ \left(
\text{since the definition of }\pi\text{ yields }\pi\left(  r\right)
=\overline{r}\right)
\end{align*}
and thus $f^{\prime}\left(  \overline{r}\right)  =f\left(  r\right)  $. Thus,
$f^{\prime}$ must be the map%
\begin{align*}
R/\operatorname*{Ker}f  &  \rightarrow f\left(  R\right)  ,\\
\overline{r}  &  \mapsto f\left(  r\right)  ,
\end{align*}
which has been called $f^{\prime}$ in Theorem \ref{thm.1it.ring1}
\textbf{(c)}. In particular, there is only one option for $f^{\prime}$. Thus,
we have shown that $f^{\prime}$ is unique, and the proof of Theorem
\ref{thm.1it.ring2} \textbf{(a)} is complete. \medskip

\textbf{(b)} As we explained above, the unique ring morphism $f^{\prime
}:R/\operatorname*{Ker}f\rightarrow f\left(  R\right)  $ that satisfies the
equation $f=\iota\circ f^{\prime}\circ\pi$ is precisely the map that was
called $f^{\prime}$ in Theorem \ref{thm.1it.ring1} \textbf{(c)}. But the
latter map is a ring isomorphism (by Theorem \ref{thm.1it.ring1}
\textbf{(c)}). Hence, Theorem \ref{thm.1it.ring2} \textbf{(b)} follows.
\end{proof}

There are also second, third and fourth isomorphism theorems. You will meet
them in Section \ref{sec.rings.234it}.

\subsection{\label{sec.rings.dirprods}Direct products of rings (\cite[\S 7.6]%
{DumFoo04})}

\subsubsection{\label{subsec.rings.dirprods.2}Direct products of two rings}

Here is a way of building new rings from old\footnote{There are several other
such ways. We will see a few in this course.}:

\begin{proposition}
\label{prop.rings.dirprod.2.ring} Let $R$ and $S$ be two rings. Then, the
Cartesian product%
\[
R\times S=\left\{  \text{all pairs }\left(  r,s\right)  \text{ with }r\in
R\text{ and }s\in S\right\}
\]
becomes a ring if we endow it with the entrywise addition and multiplication
operations (i.e., addition defined by $\left(  r,s\right)  +\left(  r^{\prime
},s^{\prime}\right)  =\left(  r+r^{\prime},s+s^{\prime}\right)  $, and
multiplication defined by $\left(  r,s\right)  \cdot\left(  r^{\prime
},s^{\prime}\right)  =\left(  rr^{\prime},ss^{\prime}\right)  $) and the zero
$\left(  0_{R},0_{S}\right)  $ and the unity $\left(  1_{R},1_{S}\right)  $.
\end{proposition}

\begin{definition}
\label{def.ring.dirprod.RxS}This ring is denoted by $R\times S$ and is called
the \textbf{direct product} of $R$ and $S$.
\end{definition}

\begin{proof}
[Proof of Proposition \ref{prop.rings.dirprod.2.ring}.]We must check that the
ring axioms are satisfied for $R\times S$. Each one is straightforward to
verify. For example, in order to check the associativity of multiplication, we
need to check that%
\begin{align*}
\left(  r,s\right)  \left(  \left(  r^{\prime},s^{\prime}\right)  \left(
r^{\prime\prime},s^{\prime\prime}\right)  \right)   &  =\left(  \left(
r,s\right)  \left(  r^{\prime},s^{\prime}\right)  \right)  \left(
r^{\prime\prime},s^{\prime\prime}\right) \\
&  \ \ \ \ \ \ \ \ \ \ \text{for all }\left(  r,s\right)  ,\left(  r^{\prime
},s^{\prime}\right)  ,\left(  r^{\prime\prime},s^{\prime\prime}\right)  \in
R\times S
\end{align*}
(because any element of $R\times S$ is a pair). We can do this by computing
both sides and comparing: We have%
\begin{align*}
\left(  r,s\right)  \left(  \left(  r^{\prime},s^{\prime}\right)  \left(
r^{\prime\prime},s^{\prime\prime}\right)  \right)   &  =\left(  r,s\right)
\left(  r^{\prime}r^{\prime\prime},s^{\prime}s^{\prime\prime}\right)  =\left(
r\left(  r^{\prime}r^{\prime\prime}\right)  ,s\left(  s^{\prime}%
s^{\prime\prime}\right)  \right)  \ \ \ \ \ \ \ \ \ \ \text{and}\\
\left(  \left(  r,s\right)  \left(  r^{\prime},s^{\prime}\right)  \right)
\left(  r^{\prime\prime},s^{\prime\prime}\right)   &  =\left(  rr^{\prime
},ss^{\prime}\right)  \left(  r^{\prime\prime},s^{\prime\prime}\right)
=\left(  \left(  rr^{\prime}\right)  r^{\prime\prime},\left(  ss^{\prime
}\right)  s^{\prime\prime}\right)  .
\end{align*}
The right hand sides of these two equalities are equal, since $r\left(
r^{\prime}r^{\prime\prime}\right)  =\left(  rr^{\prime}\right)  r^{\prime
\prime}$ and $s\left(  s^{\prime}s^{\prime\prime}\right)  =\left(  ss^{\prime
}\right)  s^{\prime\prime}$. Thus, the left hand sides are equal as well; this
proves the associativity of multiplication. All other ring axioms follow similarly.
\end{proof}

\subsubsection{\label{subsec.rings.dirprods.n}Direct products of any number of
rings}

More generally, we can define the direct product $R_{1}\times R_{2}%
\times\cdots\times R_{n}$ of any number of rings in the same way (but using
$n$-tuples instead of pairs). Even more generally, we can define the direct
product $\prod_{i\in I}R_{i}$ of any family of rings (including infinite families).

To do so, we recall the notion of a \textquotedblleft family\textquotedblright:

A \textbf{family} is a collection of objects (the \textquotedblleft
entries\textquotedblright\ of the family) indexed by the elements of a given
set $I$ (the \textquotedblleft indexing set\textquotedblright). To be more
specific: A \textbf{family} indexed by a set $I$ means a way of assigning some
object $x_{i}$ to each element $i\in I$. This family is denoted by $\left(
x_{i}\right)  _{i\in I}$ (pronounced \textquotedblleft the family consisting
of $x_{i}$ for each $i\in I$\textquotedblright), and the set $I$ is called its
\textbf{indexing set}, whereas the objects $x_{i}$ are called the
\textbf{entries} of this family.\footnote{Programmers know families under
\href{https://en.wikipedia.org/wiki/Associative_array}{the name
\textquotedblleft dictionaries\textquotedblright\ or \textquotedblleft
associative arrays\textquotedblright} (although the indexing set $I$ is
usually finite in any real-life programming situation).} Two families $\left(
x_{i}\right)  _{i\in I}$ and $\left(  y_{i}\right)  _{i\in I}$ are considered
to be equal if each $i\in I$ satisfies $x_{i}=y_{i}$.

For example, the family $\left(  \mathbb{Z}/n\right)  _{n\in\mathbb{N}}$
consists of all the quotient rings $\mathbb{Z}/n$ of the ring $\mathbb{Z}$
with $n\in\mathbb{N}$. Its indexing set is $\mathbb{N}$, and its entries are
the quotient rings $\mathbb{Z}/n$. For another example, the family $\left(
n^{2}\right)  _{n\in\mathbb{N}}$ consists of the squares of all nonnegative
integers $n$. Its indexing set is $\mathbb{N}$, and its entries are the
squares $n^{2}$.

The notion of a family encompasses several well-known mathematical concepts:

\begin{itemize}
\item $n$-tuples: If $I=\left\{  1,2,\ldots,n\right\}  $, then a family
$\left(  x_{i}\right)  _{i\in I}$ is the $n$-tuple $\left(  x_{1},x_{2}%
,\ldots,x_{n}\right)  $.

\item infinite sequences: If $I=\mathbb{N}=\left\{  0,1,2,\ldots\right\}  $,
then a family $\left(  x_{i}\right)  _{i\in I}$ is the sequence $\left(
x_{0},x_{1},x_{2},\ldots\right)  $.

\item sequences infinite on both sides: If $I=\mathbb{Z}$, then a family
$\left(  x_{i}\right)  _{i\in I}$ is the \textquotedblleft
infinite-on-both-sides sequence\textquotedblright\ $\left(  \ldots
,x_{-2},x_{-1},x_{0},x_{1},x_{2},\ldots\right)  $.

\item maps (i.e., functions): If $I$ and $R$ are any two sets, then a map
$f:I\rightarrow R$ can be viewed as a family $\left(  f\left(  i\right)
\right)  _{i\in I}$ whose entries are the values of $f$. (Fine print: The
family $\left(  f\left(  i\right)  \right)  _{i\in I}$ does not
\textquotedblleft know\textquotedblright\ the set $R$, so it does not fully
represent the map $f$.)
\end{itemize}

Now, if $\left(  X_{i}\right)  _{i\in I}$ is a family of sets (i.e., if
$X_{i}$ is a set for each $i\in I$), then the \textbf{Cartesian product}
$\prod_{i\in I}X_{i}$ of these sets is defined to be the set of all families
$\left(  x_{i}\right)  _{i\in I}$ that satisfy $x_{i}\in X_{i}$ for each $i\in
I$. For instance, an element of the Cartesian product $\prod_{n\in\mathbb{N}%
}\left(  \mathbb{Z}/n\right)  $ is a family $\left(  x_{n}\right)
_{n\in\mathbb{N}}$, where each $x_{n}$ is a residue class in the corresponding
ring $\mathbb{Z}/n$.

We are now ready to define the direct product $\prod_{i\in I}R_{i}$ of an
arbitrary family of rings:

\begin{proposition}
\label{prop.ring.dirprod.Ri}Let $I$ be any set. Let $\left(  R_{i}\right)
_{i\in I}$ be a family of rings (i.e., let $R_{i}$ be a ring for each $i\in
I$). Then, the Cartesian product
\[
\prod_{i\in I}R_{i}=\left\{  \text{all families }\left(  r_{i}\right)  _{i\in
I}\text{ with }r_{i}\in R_{i}\text{ for each }i\in I\right\}
\]
becomes a ring if we endow it with the entrywise addition and multiplication
operations (i.e., addition defined by $\left(  r_{i}\right)  _{i\in I}+\left(
s_{i}\right)  _{i\in I}=\left(  r_{i}+s_{i}\right)  _{i\in I}$, and
multiplication defined by $\left(  r_{i}\right)  _{i\in I}\cdot\left(
s_{i}\right)  _{i\in I}=\left(  r_{i}s_{i}\right)  _{i\in I}$) and the zero
$\left(  0_{R_{i}}\right)  _{i\in I}$ and the unity $\left(  1_{R_{i}}\right)
_{i\in I}$.
\end{proposition}

\begin{definition}
\label{def.ring.dirprod.Ri}The ring defined in Proposition
\ref{prop.ring.dirprod.Ri} is denoted by $\prod_{i\in I}R_{i}$ and is called
the \textbf{direct product} of the rings $R_{i}$. In some special cases, there
are alternative notations for it:

\begin{itemize}
\item If $I=\left\{  1,2,\ldots,n\right\}  $ for some $n\in\mathbb{N}$, then
the ring $\prod_{i\in I}R_{i}$ is also denoted by $R_{1}\times R_{2}%
\times\cdots\times R_{n}$, and we identify each family $\left(  r_{i}\right)
_{i\in I}=\left(  r_{i}\right)  _{i\in\left\{  1,2,\ldots,n\right\}  }$ with
the $n$-tuple $\left(  r_{1},r_{2},\ldots,r_{n}\right)  $. (Thus, the elements
of $R_{1}\times R_{2}\times\cdots\times R_{n}$ are $n$-tuples whose entries
belong to $R_{1},R_{2},\ldots,R_{n}$, respectively.) In particular, for $n=2$,
this recovers the definition of $R\times S$ in Definition
\ref{def.ring.dirprod.RxS}.

\item If all the rings $R_{i}$ are equal to some ring $R$, then their direct
product $\prod_{i\in I}R_{i}=\prod_{i\in I}R$ is also denoted $R^{I}$. Note
that this is the same notation that we previously used for the ring of all
functions from $I$ to $R$ (with pointwise addition and multiplication);
however, these two notations don't really clash, since these two rings are the
same (at least if we identify a function $f:I\rightarrow R$ with the family
$\left(  f\left(  i\right)  \right)  _{i\in I}$). (Pointwise
addition/multiplication of functions corresponds precisely to entrywise
addition/multiplication in the direct product $\prod_{i\in I}R$.)

\item If $n\in\mathbb{N}$, and if $R$ is a ring, then the ring $R^{\left\{
1,2,\ldots,n\right\}  }=\underbrace{R\times R\times\cdots\times R}_{n\text{
times}}$ is also called $R^{n}$.
\end{itemize}
\end{definition}

\subsubsection{\label{subsec.rings.dirprods.exas}Examples}

Here are some examples of direct products:

\begin{itemize}
\item The ring $\mathbb{Z}^{3}=\mathbb{Z}\times\mathbb{Z}\times\mathbb{Z}$
consists of all triples $\left(  r,s,t\right)  $ of integers. They are added
and multiplied entrywise: i.e., we have%
\begin{align*}
\left(  r,s,t\right)  +\left(  r^{\prime},s^{\prime},t^{\prime}\right)   &
=\left(  r+r^{\prime},s+s^{\prime},t+t^{\prime}\right)
\ \ \ \ \ \ \ \ \ \ \text{and}\\
\left(  r,s,t\right)  \cdot\left(  r^{\prime},s^{\prime},t^{\prime}\right)
&  =\left(  rr^{\prime},ss^{\prime},tt^{\prime}\right)  .
\end{align*}

Note that this ring is \textbf{not} an integral domain, since $\left(
0,1,0\right)  \cdot\left(  1,0,0\right)  =\left(  0,0,0\right)  $. Generally,
a direct product of rings is never an integral domain (unless all but one of
its factors are trivial rings and the remaining factor is an integral domain).

\item If $R$, $S$ and $T$ are three rings, then the direct products $R\times
S\times T$ and $\left(  R\times S\right)  \times T$ are not quite the same
(e.g., the former consists of triples $\left(  r,s,t\right)  $, while the
latter consists of nested pairs $\left(  \left(  r,s\right)  ,t\right)  $);
but they are isomorphic through a rather obvious isomorphism: Namely, the map%
\begin{align*}
R\times S\times T  &  \rightarrow\left(  R\times S\right)  \times T,\\
\left(  r,s,t\right)   &  \mapsto\left(  \left(  r,s\right)  ,t\right)
\end{align*}
is a ring isomorphism. This is a quick test of understanding -- if you
understand the definitions, then this should be completely obvious to you.
Similarly, the rings $R\times S\times T$ and $R\times\left(  S\times T\right)
$ are isomorphic. You can easily generalize this to direct products of more
than three rings. We say that the direct product of rings is \textquotedblleft
associative up to isomorphism\textquotedblright.

\item Likewise, the direct product of rings is \textquotedblleft commutative
up to isomorphism\textquotedblright: That is, if $R$ and $S$ are two rings,
then the map%
\begin{align*}
R\times S  &  \rightarrow S\times R,\\
\left(  r,s\right)   &  \mapsto\left(  s,r\right)
\end{align*}
is a ring isomorphism.

\item The ring $\mathbb{C}$ consists of complex numbers, which are defined as
pairs of real numbers (the real part and the imaginary part). Thus,
$\mathbb{C}=\mathbb{R}\times\mathbb{R}$ as sets. Since complex numbers are
added entrywise, we even have $\mathbb{C}=\mathbb{R}\times\mathbb{R}$ as
additive groups (i.e., the additive groups $\left(  \mathbb{C},+,0\right)  $
and $\left(  \mathbb{R}\times\mathbb{R},+,0\right)  $ are identical). However,
$\mathbb{C}$ is \textbf{not} $\mathbb{R}\times\mathbb{R}$ as rings (because
complex numbers are not multiplied entrywise). Even worse, $\mathbb{C}$ is not
even isomorphic to $\mathbb{R}\times\mathbb{R}$ as rings. One way to see this
is by noticing that $\mathbb{C}$ is an integral domain (even a field) whereas
$\mathbb{R}\times\mathbb{R}$ is not (for example, $\left(  1,0\right)
\cdot\left(  0,1\right)  =\left(  0,0\right)  $). Another way to see this is
by noticing that $-1_{\mathbb{C}}$ is a square in $\mathbb{C}$, but
$-1_{\mathbb{R}\times\mathbb{R}}=\left(  -1,-1\right)  $ is not a square in
$\mathbb{R}\times\mathbb{R}$.

Note that these arguments make sense because of the \textquotedblleft
isomorphism principle\textquotedblright\ (which we stated in Subsection
\ref{subsec.ringmor.iso.basic-props}). We recall that this principle says that
isomorphic rings \textquotedblleft behave the same\textquotedblright\ as far
as their properties are concerned -- at least those properties that can be
stated in terms of the ring itself. For example, if $R$ and $S$ are two
isomorphic rings, and if one of $R$ and $S$ is a field, then so is the other.
For yet another example, if $R$ and $S$ are two isomorphic rings, and $R$ has
(say) $15$ units, then so does $S$. For yet another example, if $R$ and $S$
are two isomorphic rings, and $R$ satisfies some property like
\textquotedblleft$x\left(  x+1_{R}\right)  \left(  x-1_{R}\right)  =0$ for all
$x\in R$\textquotedblright, then so does $S$ (with $1_{R}$ replaced by $1_{S}%
$). The only properties of a ring that are not preserved under isomorphism are
properties that refer to specific \textquotedblleft outside\textquotedblright%
\ objects (for example, the rings $R\times S$ and $S\times R$ from the
previous example are isomorphic, but they generally contain different
elements). This all is a general feature of isomorphisms of any sorts of
objects -- not just of rings but also of groups, vector spaces and topological spaces.

\item Let $R$ be any ring. Let $n\in\mathbb{N}$. Let $R^{n=n}$ be the set of
all \textbf{diagonal} matrices in the matrix ring $R^{n\times n}$. In other
words,%
\begin{align*}
R^{n=n}  &  =\left\{  \left(
\begin{array}
[c]{cccc}%
a_{1} & 0 & \cdots & 0\\
0 & a_{2} & \cdots & 0\\
\vdots & \vdots & \ddots & \vdots\\
0 & 0 & \cdots & a_{n}%
\end{array}
\right)  \ \mid\ a_{1},a_{2},\ldots,a_{n}\in R\right\} \\
&  =\left\{  \operatorname*{diag}\left(  a_{1},a_{2},\ldots,a_{n}\right)
\ \mid\ a_{1},a_{2},\ldots,a_{n}\in R\right\}  ,
\end{align*}
where we are using the notation%
\begin{align*}
\operatorname*{diag}\left(  a_{1},a_{2},\ldots,a_{n}\right)   &  =\left(
\text{the diagonal matrix with diagonal }\left(  a_{1},a_{2},\ldots
,a_{n}\right)  \right) \\
&  =\left(
\begin{array}
[c]{cccc}%
a_{1} & 0 & \cdots & 0\\
0 & a_{2} & \cdots & 0\\
\vdots & \vdots & \ddots & \vdots\\
0 & 0 & \cdots & a_{n}%
\end{array}
\right)  .
\end{align*}
(For example,%
\[
R^{2=2}=\left\{  \left(
\begin{array}
[c]{cc}%
a & 0\\
0 & b
\end{array}
\right)  \ \mid\ a,b\in R\right\}  =\left\{  \operatorname*{diag}\left(
a,b\right)  \ \mid\ a,b\in R\right\}  .
\]
)

It is easy to see that $R^{n=n}$ is a subring of $R^{n\times n}$. Moreover,
$R^{n=n}\cong R^{n}$ as rings (where, as we recall, $R^{n}=\underbrace{R\times
R\times\cdots\times R}_{n\text{ times}}=R^{\left\{  1,2,\ldots,n\right\}  }$).
Indeed, the map%
\begin{align*}
R^{n}  &  \rightarrow R^{n=n},\\
\left(  a_{1},a_{2},\ldots,a_{n}\right)   &  \mapsto\operatorname*{diag}%
\left(  a_{1},a_{2},\ldots,a_{n}\right)
\end{align*}
is a ring isomorphism. For example, it respects multiplication, since%
\[
\operatorname*{diag}\left(  a_{1},a_{2},\ldots,a_{n}\right)  \cdot
\operatorname*{diag}\left(  b_{1},b_{2},\ldots,b_{n}\right)
=\operatorname*{diag}\left(  a_{1}b_{1},a_{2}b_{2},\ldots,a_{n}b_{n}\right)
\]
for any $\left(  a_{1},a_{2},\ldots,a_{n}\right)  ,\left(  b_{1},b_{2}%
,\ldots,b_{n}\right)  \in R^{n}$.
\end{itemize}

It is easy to see that a direct product of commutative rings is commutative.

\begin{exercise}
Prove that the direct product $\left(  \mathbb{Z}/2\right)  ^{2}=\left(
\mathbb{Z}/2\right)  \times\left(  \mathbb{Z}/2\right)  $ is isomorphic to the
ring $B_{4}$ from Subsection \ref{subsec.rings.def.exas}.
\end{exercise}

\begin{exercise}
Let $R$ be any ring. Let $I$ be an infinite set. As we recall, the direct
product $R^{I}=\prod_{i\in I}R$ consists of all families $\left(
r_{i}\right)  _{i\in I}$ of elements of $R$. (For instance, if $I=\mathbb{N}$,
then these families are just the infinite sequences $\left(  r_{0},r_{1}%
,r_{2},\ldots\right)  =\left(  r_{i}\right)  _{i\in\mathbb{N}}$ of elements of
$R$.)

\begin{enumerate}
\item[\textbf{(a)}] We say that a family $\left(  r_{i}\right)  _{i\in I}$ is
\textbf{finitary} if it has only finitely many nonzero entries (i.e., if there
are only finitely many $i\in I$ that satisfy $r_{i}\neq0$). Consider the set
of all finitary families in $R^{I}$. Is this set a subring of $R^{I}$ ?

\item[\textbf{(b)}] We say that a family $\left(  r_{i}\right)  _{i\in I}$ is
\textbf{quasifinitary} if all but finitely many of its entries are equal
(i.e., if there exists some $c\in R$ such that only finitely many $i\in I$
satisfy $r_{i}\neq c$). Consider the set of all quasifinitary families in
$R^{I}$. Is this set a subring of $R^{I}$ ?

\item[\textbf{(c)}] We say that a family $\left(  r_{i}\right)  _{i\in I}$ is
\textbf{entry-finite} if it has only finitely many distinct entries (i.e., if
the set $\left\{  r_{i}\ \mid\ i\in I\right\}  $ is finite). Consider the set
of all entry-finite families in $R^{I}$. Is this set a subring of $R^{I}$ ?
\end{enumerate}

(For example, if $R=\mathbb{Z}$ and $I=\mathbb{N}$, then the family $\left(
3,1,0,0,0,0,\ldots\right)  $ (with all entries after the $1$ being zeroes) is
finitary; the family $\left(  3,1,1,1,\ldots\right)  $ (with all entries after
the $3$ being equal to $1$) is quasifinitary; the family $\left(
1,2,1,2,1,2,\ldots\right)  $ (alternating between $1$'s and $2$'s) is entry-finite.)
\end{exercise}

\subsubsection{\label{subsec.rings.dirprods.idp}Direct products and
idempotents}

Direct products of rings are closely related to idempotents. One part of the
connection is easy: If $R$ and $S$ are two rings, then their direct product
$R\times S$ has the two idempotents $\left(  1_{R},0_{S}\right)  $ and
$\left(  0_{R},1_{S}\right)  $ which, in a sense, \textquotedblleft
reveal\textquotedblright\ its two factors; in particular, the multiples of
$\left(  1_{R},0_{S}\right)  $ form a \textquotedblleft copy\textquotedblright%
\ of $R$ (since these multiples are precisely the elements of the form
$\left(  r,0_{S}\right)  $ for $r\in R$), whereas the multiples of $\left(
0_{R},1_{S}\right)  $ form a \textquotedblleft copy\textquotedblright\ of $S$.
The following exercise states this claim precisely:

\begin{exercise}
Let $R$ and $S$ be two rings. Let $a:=\left(  1_{R},0_{S}\right)  \in R\times
S$ and $b:=\left(  0_{R},1_{S}\right)  \in R\times S$. Prove the following:

\begin{enumerate}
\item[\textbf{(a)}] The principal ideal $a\left(  R\times S\right)  $ consists
of all elements of the form $\left(  r,0_{S}\right)  $ with $r\in R$.

\item[\textbf{(b)}] The principal ideal $b\left(  R\times S\right)  $ consists
of all elements of the form $\left(  0_{R},s\right)  $ with $s\in S$.
\end{enumerate}
\end{exercise}

If the ring $R$ is commutative, then this connection has a converse as well:
Any idempotent in a commutative ring $R$ can be used to split $R$ into a
direct product of two rings!\footnote{If the idempotent is $0$ or $1$, then
one of the two factors will be a trivial ring.} Here are the details:

\begin{exercise}
\label{exe.21hw1.3}Let $R$ be a commutative ring, and let $e$ be an idempotent
element of $R$. As we know from Exercise \ref{exe.idemp.basics} \textbf{(a)}
(applied to $a=e$), the element $1-e$ is then idempotent as well.

\begin{enumerate}
\item[\textbf{(a)}] Show that the principal ideal $eR$ is itself a ring, with
addition and multiplication inherited from $R$ and with zero $0_{R}$ and with
unity $e$. (This makes $eR$ a subring of $R$ in the sense of \cite{DumFoo04},
but not in our sense, since its unity is not generally the unity of $R$.)

\item[\textbf{(b)}] Show that the same holds for the principal ideal $\left(
1-e\right)  R$ (except that its unity will be $1-e$ instead of $e$).

\item[\textbf{(c)}] Consider the map
\begin{align*}
f:\left(  eR\right)  \times\left(  \left(  1-e\right)  R\right)   &
\rightarrow R,\\
\left(  a,b\right)   &  \mapsto a+b.
\end{align*}
Prove that this map $f$ is a ring isomorphism.
\end{enumerate}
\end{exercise}

Exercise \ref{exe.21hw1.3} \textbf{(c)} shows that if a commutative ring $R$
has an idempotent element $e$, then $R$ can be decomposed (up to isomorphism)
as a direct product $A\times B$ of two rings $A$ and $B$ (namely, $A=eR$ and
$B=\left(  1-e\right)  R$). If $e$ is not one of the two trivial idempotents
$0$ and $1$, then these two rings $A$ and $B$ will be nontrivial, so the
decomposition really deserves its name.\footnote{As an example, take
$R=\mathbb{Z}/6\mathbb{Z}$, and let $e$ be the idempotent element
$\overline{3}=3+6\mathbb{Z}$ of $R$ (this is idempotent since $3^{2}%
=9\equiv3\operatorname{mod}6$ and thus $\overline{3}^{2}=\overline{3^{2}%
}=\overline{3}$). Then, $eR=\left\{  \overline{0},\overline{3}\right\}
\cong\mathbb{Z}/2\mathbb{Z}$ and $\left(  1-e\right)  R=\left\{  \overline
{0},\overline{2},\overline{4}\right\}  \cong\mathbb{Z}/3\mathbb{Z}$. Hence,
the ring isomorphism $R\cong\left(  eR\right)  \times\left(  \left(
1-e\right)  R\right)  $ becomes a ring isomorphism $\mathbb{Z}/6\mathbb{Z}%
\cong\left(  \mathbb{Z}/2\mathbb{Z}\right)  \times\left(  \mathbb{Z}%
/3\mathbb{Z}\right)  $. We will soon revisit this isomorphism (it is an
instance of the Chinese Remainder Theorem).}

Conversely, as we said above, any direct product of two nontrivial rings has
nontrivial idempotents: If $R$ and $S$ are two rings, then $\left(
1_{R},0_{S}\right)  $ and $\left(  0_{R},1_{S}\right)  $ are two idempotent
elements of the direct product $R\times S$.

Parts \textbf{(a)} and \textbf{(b)} of Exercise \ref{exe.21hw1.3} can be
generalized somewhat: Instead of requiring $R$ to be commutative, it suffices
to require that $er=re$ for all $r\in R$. We cannot, however, drop this
requirement altogether (for instance, the matrix ring $\mathbb{R}^{2\times2}$
has many idempotents, but cannot be written as a direct product of two
nontrivial rings).

\begin{exercise}
\label{exe.circulant.1}Let $R$ be a commutative ring, and $n$ be a positive
integer. An $n\times n$-matrix
\[
A=\left(
\begin{array}
[c]{cccc}%
a_{1,1} & a_{1,2} & \cdots & a_{1,n}\\
a_{2,1} & a_{2,2} & \cdots & a_{2,n}\\
\vdots & \vdots & \ddots & \vdots\\
a_{n,1} & a_{n,2} & \cdots & a_{n,n}%
\end{array}
\right)  \in R^{n\times n}%
\]
is said to be \textbf{circulant} if it has the property that%
\[
a_{i,j}=a_{i^{\prime},j^{\prime}}\text{ whenever }i-j\equiv i^{\prime
}-j^{\prime}\operatorname{mod}n.
\]
In other words, an $n\times n$-matrix $A\in R^{n\times n}$ is circulant if and
only if each row of $A$ equals the preceding row of $A$, cyclically rotated by
$1$ step to the right. For instance, a $4\times4$-matrix is circulant if and
only if it has the form $\left(
\begin{array}
[c]{cccc}%
a & b & c & d\\
d & a & b & c\\
c & d & a & b\\
b & c & d & a
\end{array}
\right)  $ for some $a,b,c,d\in R$.

Let $\operatorname*{Circ}\nolimits_{n}\left(  R\right)  $ denote the set of
all circulant $n\times n$-matrices $A\in R^{n\times n}$.

Let $S\in\operatorname*{Circ}\nolimits_{n}\left(  R\right)  $ be the specific
circulant $n\times n$-matrix whose first row is $\left(  0,1,0,0,0,\ldots
,0\right)  $ (that is, the second entry is $1$ while all the other entries are
$0$). Thus,%
\[
S=\left(
\begin{array}
[c]{cccccc}%
0 & 1 & 0 & 0 & \cdots & 0\\
0 & 0 & 1 & 0 & \cdots & 0\\
0 & 0 & 0 & 1 & \cdots & 0\\
\vdots & \vdots & \vdots & \vdots & \ddots & \vdots\\
0 & 0 & 0 & 0 & \cdots & 1\\
1 & 0 & 0 & 0 & \cdots & 0
\end{array}
\right)  .
\]
(If $n=1$, then $S$ is defined to be the $1\times1$-matrix $\left(
\begin{array}
[c]{c}%
1
\end{array}
\right)  $.)

\begin{enumerate}
\item[\textbf{(a)}] Compute $S^{n}$.

\item[\textbf{(b)}] Prove that every circulant matrix $A\in R^{n\times n}$ can
be written as $a_{0}S^{0}+a_{1}S^{1}+\cdots+a_{n-1}S^{n-1}$, where
$a_{0},a_{1},\ldots,a_{n-1}$ are the entries of the first row of $A$ (from
left to right).

\item[\textbf{(c)}] Prove that
\[
\operatorname*{Circ}\nolimits_{n}\left(  R\right)  =\left\{  a_{0}S^{0}%
+a_{1}S^{1}+\cdots+a_{n-1}S^{n-1}\ \mid\ a_{0},a_{1},\ldots,a_{n-1}\in
R\right\}  .
\]

\item[\textbf{(d)}] Show that $\operatorname*{Circ}\nolimits_{n}\left(
R\right)  $ is a \textbf{commutative} subring of the matrix ring $R^{n\times
n}$.

\item[\textbf{(e)}] Assume that $R$ is nontrivial, that $n\geq2$, and that the
element $n\cdot1_{R}$ of $R$ is invertible. Find an idempotent $e$ in
$\operatorname*{Circ}\nolimits_{n}\left(  R\right)  $ that is distinct from
both the zero matrix $0$ and the identity matrix $I_{n}$.

\item[\textbf{(f)}] Under the same assumptions as in part \textbf{(e)}, prove
that $\operatorname*{Circ}\nolimits_{n}\left(  R\right)  $ is isomorphic to a
direct product of two rings, one of which is $R$.

\item[\textbf{(g)}] Which of the above claims remain true if we no longer
require that $R$ be commutative?
\end{enumerate}
\end{exercise}

\subsubsection{\label{subsec.rings.dirprods.bool}Boolean rings}

Another example of rings comes from (fairly basic) set theory. It rests upon
the notion of \textquotedblleft symmetric difference\textquotedblright:

\begin{definition}
The \textbf{symmetric difference} of two sets $A$ and $B$ is defined to be the
set
\begin{align*}
&  \left(  A\cup B\right)  \setminus\left(  A\cap B\right) \\
&  =\left(  A\setminus B\right)  \cup\left(  B\setminus A\right) \\
&  =\left\{  x\ \mid\ x\text{ belongs to exactly one of the two sets }A\text{
and }B\right\}  .
\end{align*}
This symmetric difference is denoted by $A\mathbin{\bigtriangleup}B$.
\end{definition}

In terms of Venn diagrams, this symmetric difference
$A\mathbin{\bigtriangleup}B$ is the grey zone in the following Venn diagram
(where the two circles are $A$ and $B$):%
\[%
\begin{tikzpicture}
\fill[gray, even odd rule] (0,0) circle (2) (1.5,0) circle (1.7);
\end{tikzpicture}%
\]

Now, let $S$ be any set. Let $\mathcal{P}\left(  S\right)  $ denote the power
set of $S$ (that is, the set of all subsets of $S$). It is easy to check that
the following ten properties hold:
\begin{align*}
A\mathbin{\bigtriangleup} B  &  =B\mathbin{\bigtriangleup} A\qquad\text{for
any sets $A$ and $B$;}\\
A\cap B  &  =B\cap A\qquad\text{for any sets $A$ and $B$;}\\
\left(  A\mathbin{\bigtriangleup} B\right)  \mathbin{\bigtriangleup} C  &
=A\mathbin{\bigtriangleup}\left(  B\mathbin{\bigtriangleup} C\right)
\qquad\text{for any sets $A$, $B$ and $C$;}\\
\left(  A\cap B\right)  \cap C  &  =A\cap\left(  B\cap C\right)
\qquad\text{for any sets $A$, $B$ and $C$;}\\
A\mathbin{\bigtriangleup}\varnothing &  =\varnothing\mathbin{\bigtriangleup}
A=A\qquad\text{for any set $A$;}\\
A\mathbin{\bigtriangleup} A  &  =\varnothing\qquad\text{for any set $A$;}\\
A\cap S  &  =S\cap A=A\qquad\text{for any subset $A$ of $S$;}\\
\varnothing\cap A  &  =A\cap\varnothing=\varnothing\qquad\text{for any set
$A$;}\\
A\cap\left(  B\mathbin{\bigtriangleup} C\right)   &  =\left(  A\cap B\right)
\mathbin{\bigtriangleup}\left(  A\cap C\right)  \qquad\text{for any sets $A$,
$B$ and $C$;}\\
\left(  A\mathbin{\bigtriangleup} B\right)  \cap C  &  =\left(  A\cap
C\right)  \mathbin{\bigtriangleup}\left(  B\cap C\right)  \qquad\text{for any
sets $A$, $B$ and $C$.}%
\end{align*}
Therefore, $\mathcal{P}\left(  S\right)  $ becomes a commutative ring, where
the addition is defined to be the operation $\mathbin{\bigtriangleup}$, the
multiplication is defined to be the operation $\cap$, the zero is defined to
be the set $\varnothing$, and the unity is defined to be the set $S$. (The ten
properties listed above show that the axioms of a commutative ring are
satisfied for $\left(  \mathcal{P}\left(  S\right)
,\mathbin{\bigtriangleup},\cap,\varnothing,S\right)  $. In particular, the
sixth property shows that every subset $A$ of $S$ has an additive inverse --
namely, itself. Of course, it is unusual for an element of a commutative ring
to be its own additive inverse, but in this example it happens all the time!)

The commutative ring $\mathcal{P}\left(  S\right)  $ has the property that
each element $a\in\mathcal{P}\left(  S\right)  $ is idempotent (i.e.,
satisfies $a\cdot a=a$). (Indeed, this simply means that each $A\subseteq S$
satisfies $A\cap A=A$.)

\begin{exercise}
\label{exe.21hw1.4}\ \ 

\begin{enumerate}
\item[\textbf{(a)}] Prove that the ring $\mathcal{P}\left(  S\right)  $ is
isomorphic to the direct product $\left(  \mathbb{Z}/2\mathbb{Z}\right)
^{S}=\prod_{s\in S}\left(  \mathbb{Z}/2\mathbb{Z}\right)  $.

\item[\textbf{(b)}] Let $F$ be the set of all \textbf{finite} subsets of $S$.
Prove that $F$ is an ideal of $\mathcal{P}\left(  S\right)  $.

\item[\textbf{(c)}] Assume that $S$ is infinite. Prove that the ideal $F$ is
not principal.

\item[\textbf{(d)}] Instead, assume that $S$ is finite. Prove that every ideal
of $\mathcal{P}\left(  S\right)  $ is principal.
\end{enumerate}

[\textbf{Hint:} For part \textbf{(d)}, let $I$ be an ideal of $\mathcal{P}%
\left(  S\right)  $, and pick a subset $T\in I$ of largest size. Argue that
each subset of $T$ must also belong to $I$. Conclude that every set in $I$
must be a subset of $T$.]
\end{exercise}

Forget that we fixed $S$. As we noticed, the ring $\mathcal{P}\left(
S\right)  $ that we have just defined has the strange-looking property that
each of its elements is idempotent. Rings with this property are called
\href{https://en.wikipedia.org/wiki/Boolean_ring}{\textbf{Boolean rings}}. (Of
course, $\mathcal{P}\left(  S\right)  $ is the eponymic example for a Boolean
ring; but there are also others.) Let us now study Boolean rings in general:

\begin{definition}
A \textbf{Boolean ring} means a ring $R$ such that every $a\in R$ satisfies
$a^{2}=a$ (that is, every $a\in R$ is idempotent). (Keep in mind that rings
must have a $1$ according to our definition.)
\end{definition}

\begin{exercise}
\label{exe.21hw1.5}Let $R$ be a Boolean ring. Prove the following:

\begin{enumerate}
\item[\textbf{(a)}] We have $2a=0$ for each $a\in R$.

\item[\textbf{(b)}] We have $-a=a$ for each $a\in R$.

\item[\textbf{(c)}] The ring $R$ is commutative.

\item[\textbf{(d)}] If $R$ is finite, then $R\cong\left(  \mathbb{Z}%
/2\mathbb{Z}\right)  ^{n}$ for some $n\in\mathbb{N}$.
\end{enumerate}

[\textbf{Hint:} In part \textbf{(a)}, use $a^{2}=a$ and $\left(  a+1\right)
^{2}=a+1$. In part \textbf{(c)}, expand $\left(  a+b\right)  ^{2}$ (but don't
use the binomial formula, since you don't know yet that $ab=ba$). Finally, for
part \textbf{(d)}, use strong induction on $\left\vert R\right\vert $ as
follows: Pick some $e\in R$ that is distinct from $0$ and $1$ (if no such $e$
exists, the claim is obvious). Then, $e$ is idempotent, so Exercise
\ref{exe.21hw1.3} \textbf{(c)} decomposes the ring $R$ as a direct product of
two smaller rings. You can use without proof that direct products are
associative up to isomorphism (so that $A_{1}\times A_{2}\times\cdots\times
A_{m}\cong\left(  A_{1}\times A_{2}\times\cdots\times A_{k}\right)
\times\left(  A_{k+1}\times A_{k+2}\times\cdots\times A_{m}\right)  $ for any
rings $A_{1},A_{2},\ldots,A_{m}$).]
\end{exercise}

\subsection{\label{sec.rings.ideal-arith}A few operations on ideals
(\cite[\S 7.3]{DumFoo04})}

Next, we shall see three ways to build new ideals of a ring from old. One of
these three ways is intersection: If $I$ and $J$ are two ideals of a ring $R$,
then their intersection $I\cap J$ is easily seen to be an ideal as well (see
Proposition \ref{prop.rings.ideal-arith.laws} \textbf{(a)} below). Let us now
define two other ways:

\begin{definition}
\label{def.rings.ideal-arith.arith}Let $I$ and $J$ be two ideals of a ring $R$.

\begin{enumerate}
\item[\textbf{(a)}] Then, $I+J$ denotes the subset%
\[
\left\{  i+j\ \mid\ i\in I\text{ and }j\in J\right\}  \text{ of }R.
\]

\item[\textbf{(b)}] Next, we define a further subset $IJ$, also denoted
$I\cdot J$. Unlike $I+J$, this will \textbf{not} be defined as $\left\{
i\cdot j\ \mid\ i\in I\text{ and }j\in J\right\}  $. Instead, $IJ=I\cdot J$
will be defined as the set%
\[
\left\{  \text{all finite sums of }\left(  I,J\right)  \text{-products}%
\right\}  ,
\]
where an $\left(  I,J\right)  $\textbf{-product} means a product of the form
$ij$ with $i\in I$ and $j\in J$. In other words,%
\[
IJ=\left\{  i_{1}j_{1}+i_{2}j_{2}+\cdots+i_{k}j_{k}\ \mid\ k\in\mathbb{N}%
\text{ and }i_{1},i_{2},\ldots,i_{k}\in I\text{ and }j_{1},j_{2},\ldots
,j_{k}\in J\right\}  .
\]

\end{enumerate}
\end{definition}

Note that our definition of $IJ$ was more complicated than the one of $I+J$,
as it involved an additional step (viz., taking finite sums). The purpose of
this step is to ensure that $IJ$ is closed under addition (which will later be
used to argue that $IJ$ is an ideal of $R$). It is forced on us if we try to
construct an ideal of $R$ that contains all $\left(  I,J\right)  $-products.
We could have added the same step to our definition of $I+J$, but it would not
have changed anything, since a finite sum of $\left(  I,J\right)  $-sums
(i.e., of sums of the form $i+j$ with $i\in I$ and $j\in J$) can be rewritten
as a single $\left(  I,J\right)  $-sum:%
\begin{align*}
&  \left(  i_{1}+j_{1}\right)  +\left(  i_{2}+j_{2}\right)  +\cdots+\left(
i_{k}+j_{k}\right) \\
&  =\underbrace{\left(  i_{1}+i_{2}+\cdots+i_{k}\right)  }_{\substack{\in
I\\\text{(since }I\text{ is closed under addition)}}}+\underbrace{\left(
j_{1}+j_{2}+\cdots+j_{k}\right)  }_{\substack{\in J\\\text{(since }J\text{ is
closed under addition)}}}.
\end{align*}
For $\left(  I,J\right)  $-products, however, this is not generally the case
(although you won't find a counterexample for $R=\mathbb{Z}$).

Here is an assortment of facts about the above-defined operations on ideals
(see Exercise \ref{exe.21hw1.8} for a proof of the non-obvious
ones):\footnote{Recall that if $R$ is any ring, then the one-element set
$\left\{  0_{R}\right\}  $ and the entire ring $R$ are ideals of $R$. Both of
these ideals are principal ($\left\{  0_{R}\right\}  =0_{R}R$ and $R=1_{R}R$);
they \textquotedblleft bookend\textquotedblright\ all ideals of $R$ (in the
sense that $\left\{  0_{R}\right\}  \subseteq I\subseteq R$ for each ideal $I$
of $R$).
\par
(Here, the ideals $0_{R}R$ and $1_{R}R$ are defined as in Proposition
\ref{prop.ideal.princid}, even though $R$ is not required to be commutative.)}

\begin{proposition}
\label{prop.rings.ideal-arith.laws}Let $R$ be a ring.

\begin{enumerate}
\item[\textbf{(a)}] Let $I$ and $J$ be two ideals of $R$. Then, $I+J$ and
$I\cap J$ and $IJ$ are ideals of $R$ as well.

\item[\textbf{(b)}] Let $I$ and $J$ be two ideals of $R$. Then, $IJ\subseteq
I\cap J\subseteq I\subseteq I+J$ and $IJ\subseteq I\cap J\subseteq J\subseteq
I+J$.

\item[\textbf{(c)}] The set of all ideals of $R$ is a monoid with respect to
the binary operation $+$, with neutral element $\left\{  0_{R}\right\}  =0R$.
That is,%
\begin{align*}
\left(  I+J\right)  +K  &  =I+\left(  J+K\right)
\ \ \ \ \ \ \ \ \ \ \text{for any three ideals }I,J,K\text{ of }R;\\
I+\left\{  0_{R}\right\}   &  =\left\{  0_{R}\right\}
+I=I\ \ \ \ \ \ \ \ \ \ \text{for any ideal }I\text{ of }R.
\end{align*}
Moreover, $I+I=I$ and $R+I=R$ for any ideal $I$ of $R$.

\item[\textbf{(d)}] The set of all ideals of $R$ is a monoid with respect to
the binary operation $\cap$, with neutral element $R=1R$. That is,%
\begin{align*}
\left(  I\cap J\right)  \cap K  &  =I\cap\left(  J\cap K\right)
\ \ \ \ \ \ \ \ \ \ \text{for any three ideals }I,J,K\text{ of }R;\\
I\cap R  &  =R\cap I=I\ \ \ \ \ \ \ \ \ \ \text{for any ideal }I\text{ of }R.
\end{align*}

\item[\textbf{(e)}] The set of all ideals of $R$ is a monoid with respect to
the binary operation $\cdot$, with neutral element $R=1R$. That is,%
\begin{align*}
\left(  IJ\right)  K  &  =I\left(  JK\right)  \ \ \ \ \ \ \ \ \ \ \text{for
any three ideals }I,J,K\text{ of }R;\\
IR  &  =RI=I\ \ \ \ \ \ \ \ \ \ \text{for any ideal }I\text{ of }R.
\end{align*}

\item[\textbf{(f)}] Addition and intersection of ideals are commutative:
\[
I+J=J+I\ \ \ \ \ \ \ \ \ \ \text{and}\ \ \ \ \ \ \ \ \ \ I\cap J=J\cap
I\ \ \ \ \ \ \ \ \ \ \text{for any ideals }I,J\text{ of }R\text{.}%
\]

\item[\textbf{(g)}] If the ring $R$ is commutative, then $IJ=JI$ for any two
ideals $I$ and $J$ of $R$.

\item[\textbf{(h)}] If $I,J,K$ are three ideals of $R$ such that $I\subseteq
J$, then $I+K\subseteq J+K$ and $IK\subseteq JK$ and $KI\subseteq KJ$.
\end{enumerate}
\end{proposition}

Proposition \ref{prop.rings.ideal-arith.laws} shows that the operations $+$,
$\cap$ and $\cdot$ on the set of all ideals of $R$ satisfy a number of laws
similar to the basic laws of arithmetic. This is known as \textbf{ideal
arithmetic}. However, ideals cannot be subtracted (i.e., there is no operation
that undoes addition of ideals\footnote{That is, if $I$ and $J$ are two
ideals, then you cannot recover $I$ from $J$ and $I+J$.
\par
You may, of course, define $I-J$ (in analogy to $I+J$) by setting
$I-J=\left\{  i-j\ \mid\ i\in I\text{ and }j\in J\right\}  $. But this is just
the same ideal as $I+J$ (since $J$ is closed under negation), not a way to
undo addition of ideals!}), and thus the ideals of $R$ do not form an actual
ring. (Likewise, there is no \textquotedblleft division\textquotedblright%
\ operation on ideals that undoes multiplication, although something vaguely
similar is defined in Exercise \ref{exe.ideals.colon} below.)

Here is a diagram showing the inclusions between the ideals $IJ,\ \ I\cap
J,\ \ I+J,\ \ I,\ \ J$:%
\[%
\xymatrix{
& I+J \\
I \arinj[ur] & & J \arinj[ul] \\
& I \cap J \arinj[ul] \arinjrev[ur] \\
& \vphantom{\int_a^b} IJ \arinjrev[u]
}%
\]
(Recall that an arrow of type $X\hookrightarrow Y$ means a canonical inclusion
from $X$ to $Y$, which entails that $X\subseteq Y$.)

\begin{exercise}
\label{exe.21hw1.8}Prove Proposition \ref{prop.rings.ideal-arith.laws}.
\medskip

[\textbf{Hint:} You can be terse here, as there is a lot to show, much of it
straightforward. Part \textbf{(d)} is obvious. For part \textbf{(e)}, I
recommend using the notion of \textquotedblleft$\left(  I,J\right)
$-products\textquotedblright\ from Definition
\ref{def.rings.ideal-arith.arith}; it is often easier to talk abstractly about
sums of $\left(  I,J\right)  $-products than to write them out as $i_{1}%
j_{1}+i_{2}j_{2}+\cdots+i_{k}j_{k}$. For the proof of $\left(  IJ\right)
K=I\left(  JK\right)  $, you can start out by showing that any $\left(
IJ,K\right)  $-product belongs to $I\left(  JK\right)  $.]
\end{exercise}

The following proposition tells us what ideal arithmetic looks like when we
apply it to principal ideals of $\mathbb{Z}$:

\begin{proposition}
\label{prop.rings.ideal-arith.Z}Let $n,m\in\mathbb{Z}$. Let $I=n\mathbb{Z}$
and $J=m\mathbb{Z}$. Then:

\begin{enumerate}
\item[\textbf{(a)}] We have $IJ=nm\mathbb{Z}$.

\item[\textbf{(b)}] We have $I\cap J=\operatorname{lcm}\left(  n,m\right)
\mathbb{Z}$.

\item[\textbf{(c)}] We have $I+J=\gcd\left(  n,m\right)  \mathbb{Z}$.
\end{enumerate}
\end{proposition}

\begin{proof}
\textbf{(a)} From $n\in I$ and $m\in J$, we see that $nm$ is an $\left(
I,J\right)  $-product. Thus, $nm$ is a finite sum of $\left(  I,J\right)
$-products (of just one, to be specific). In other words, $nm\in IJ$. Since
$IJ$ is an ideal of $\mathbb{Z}$, this entails that every multiple of $nm$
also belongs to $IJ$ (by the second ideal axiom); in other words,
$nm\mathbb{Z}\subseteq IJ$.

Conversely: If $i\in I$ and $j\in J$, then $i=nx$ for some $x\in\mathbb{Z}$
(since $i\in I=n\mathbb{Z}$) and $j=my$ for some $y\in\mathbb{Z}$ (since $j\in
J=m\mathbb{Z}$) and therefore $ij=\left(  nx\right)  \left(  my\right)
=nm\left(  xy\right)  \in nm\mathbb{Z}$. Thus, every $\left(  I,J\right)
$-product belongs to $nm\mathbb{Z}$ (because an $\left(  I,J\right)  $-product
always has the form $ij$ for some $i\in I$ and $j\in J$). Hence, any sum of
$\left(  I,J\right)  $-products also belongs to $nm\mathbb{Z}$ (since
$nm\mathbb{Z}$ is closed under addition). In other words, $IJ\subseteq
nm\mathbb{Z}$ (since any element of $IJ$ is a sum of $\left(  I,J\right)
$-products). Therefore, $IJ=nm\mathbb{Z}$ (since we already have seen that
$nm\mathbb{Z}\subseteq IJ$). This proves Proposition
\ref{prop.rings.ideal-arith.Z} \textbf{(a)}. \medskip

\textbf{(b)} We have%
\begin{align*}
I\cap J  &  =\left\{  \text{all elements of }I\text{ that also belong to
}J\right\} \\
&  =\left\{  \text{all multiples of }n\text{ that are also multiples of
}m\right\} \\
&  \ \ \ \ \ \ \ \ \ \ \ \ \ \ \ \ \ \ \ \ \left(
\begin{array}
[c]{c}%
\text{since }I=n\mathbb{Z}=\left\{  \text{all multiples of }n\right\} \\
\text{and }J=m\mathbb{Z}=\left\{  \text{all multiples of }m\right\}
\end{array}
\right) \\
&  =\left\{  \text{all common multiples of }n\text{ and }m\right\} \\
&  =\left\{  \text{all multiples of }\operatorname{lcm}\left(  n,m\right)
\right\} \\
&  \ \ \ \ \ \ \ \ \ \ \ \ \ \ \ \ \ \ \ \ \left(
\begin{array}
[c]{c}%
\text{since a result in elementary number theory}\\
\text{says that the common multiples of }n\text{ and }m\\
\text{are precisely the multiples of }\operatorname{lcm}\left(  n,m\right)
\end{array}
\right) \\
&  =\operatorname{lcm}\left(  n,m\right)  \mathbb{Z}.
\end{align*}
This proves Proposition \ref{prop.rings.ideal-arith.Z} \textbf{(b)}. \medskip

\textbf{(c)} First, we shall show that $I+J\subseteq\gcd\left(  n,m\right)
\mathbb{Z}$. Indeed, any element of $I$ is a multiple of $n$ (since
$I=n\mathbb{Z}$), thus a multiple of $\gcd\left(  n,m\right)  $ (since $n$ is
a multiple of $\gcd\left(  n,m\right)  $). Similarly, any element of $J$ is a
multiple of $\gcd\left(  n,m\right)  $. Thus, an element of $I+J$ is a sum of
two multiples of $\gcd\left(  n,m\right)  $, and therefore itself a multiple
of $\gcd\left(  n,m\right)  $. In other words, any element of $I+J$ belongs to
$\gcd\left(  n,m\right)  \mathbb{Z}$. In other words, $I+J\subseteq\gcd\left(
n,m\right)  \mathbb{Z}$.

Now, we need to prove that $\gcd\left(  n,m\right)  \mathbb{Z}\subseteq I+J$.
For this, it suffices to show that $\gcd\left(  n,m\right)  \in I+J$, because
$I+J$ is an ideal (and thus will contain any multiple of $\gcd\left(
n,m\right)  $ once we know it contains $\gcd\left(  n,m\right)  $). But
Bezout's theorem shows that $\gcd\left(  n,m\right)  =xn+ym$ for some integers
$x$ and $y$. Thus, $\gcd\left(  n,m\right)  \in I+J$ (since $xn=nx\in
n\mathbb{Z}=I$ and $ym=my\in m\mathbb{Z}=J$). This finishes our proof of
$\gcd\left(  n,m\right)  \mathbb{Z}\subseteq I+J$. Combining this with
$I+J\subseteq\gcd\left(  n,m\right)  \mathbb{Z}$, we obtain $I+J=\gcd\left(
n,m\right)  \mathbb{Z}$. This proves Proposition
\ref{prop.rings.ideal-arith.Z} \textbf{(c)}.
\end{proof}

\begin{exercise}
\label{exe.ideals.arithmetic.IJ-not-JI}Let $R$ be any nontrivial ring, and
consider the ideals $I,J,K$ of the upper-triangular matrix ring $R^{2\leq2}$
defined in Exercise \ref{exe.ideals.triangular-matrices.1} \textbf{(a)}. Show
that $IJ=\left\{  0\right\}  $ but $JI=K$. (Thus, $IJ\neq JI$ in this case.)
\end{exercise}

\begin{exercise}
Let $R$ be a ring. Let $I,J,K$ be three ideals of $R$. Prove that
\[
I\left(  J+K\right)  =IJ+IK\ \ \ \ \ \ \ \ \ \ \text{and}%
\ \ \ \ \ \ \ \ \ \ \left(  I+J\right)  K=IK+JK.
\]

\end{exercise}

\begin{exercise}
\label{exe.ideals.abR}Let $R$ be a commutative ring. Let $a$ and $b$ be two
elements of $R$. Prove that $\left(  ab\right)  R=\left(  aR\right)  \left(
bR\right)  $ and $\left(  a+b\right)  R\subseteq aR+bR$.

(Recall that $\left(  ab\right)  R$, $\left(  a+b\right)  R$, $aR$ and $bR$
are principal ideals, and the \textquotedblleft$+$\textquotedblright\ sign in
\textquotedblleft$aR+bR$\textquotedblright\ is a sum of two ideals, not a sum
of two elements of $R$.)
\end{exercise}

\begin{exercise}
\label{exe.ideals.aR+I=bR+I}Let $R$ be a commutative ring. Let $I$ be an ideal
of $R$. Let $a$ and $b$ be two elements of $R$ such that $a-b\in I$. Prove
that $aR+I=bR+I$.
\end{exercise}

The operation that takes an ideal $I$ and a ring element $u$ to the ideal $uI$
(as defined in Proposition \ref{prop.ideal.uI}) is also commonly considered to
be part of ideal arithmetic; the following exercise shows how it interacts
with some of the other operations:

\begin{exercise}
Let $R$ be a commutative ring. Prove the following:

\begin{enumerate}
\item[\textbf{(a)}] For any $u\in R$ and any ideal $I$ of $R$, we have
$uI=\left(  uR\right)  I=I\left(  uR\right)  $.

\item[\textbf{(b)}] For any $u\in R$ and any two ideals $I$ and $J$ of $R$, we
have $u\left(  I+J\right)  =uI+uJ$.

\item[\textbf{(c)}] For any $u,v\in R$ and any ideal $I$ of $R$, we have
$\left(  u+v\right)  I\subseteq uI+vI$, but this is not always an equality.

\item[\textbf{(d)}] For any $u\in R$ and any two ideals $I$ and $J$ of $R$, we
have $u\left(  I\cap J\right)  \subseteq uI\cap uJ$, but this is not always an equality.
\end{enumerate}

[\textbf{Hint:} In part \textbf{(d)}, equality can fail even if $I$ and $J$
are principal ideals. For a counterexample, let $R$ be the ring of all pairs
$\left(  \overline{a},\overline{b}\right)  \in\left(  \mathbb{Z}/8\right)
^{2}$ satisfying $a\equiv b\operatorname{mod}4$.]
\end{exercise}

The next exercise defines yet another (less frequently used) operation on ideals:

\begin{exercise}
\label{exe.ideals.colon}Let $R$ be a ring. Let $I$ and $J$ be two ideals of
$R$.

We say that a given element $a\in R$ \textbf{leads} $J$ into $I$ if and only
if each $j\in J$ satisfies $aj\in I$ and $ja\in I$. In other words, an element
$a\in R$ leads $J$ into $I$ if and only if multiplying this element with any
element of $J$ (from the left or from the right) produces an element of $I$.

We let $\left(  I:J\right)  $ be the set of all elements $a\in R$ that lead
$J$ into $I$.

\begin{enumerate}
\item[\textbf{(a)}] Prove that $\left(  I:J\right)  $ is an ideal of $R$.
(This is called the \textbf{colon ideal} of $I$ and $J$.)

\item[\textbf{(b)}] For $R=\mathbb{Z}$, compute the colon ideals $\left(
6\mathbb{Z}:2\mathbb{Z}\right)  $ and $\left(  6\mathbb{Z}:4\mathbb{Z}\right)
$.

\item[\textbf{(c)}] Let $I$, $J$ and $K$ be three ideals of $R$. Prove that
\[
\left(  I:J\right)  \left(  J:K\right)  \cap\left(  J:K\right)  \left(
I:J\right)  \subseteq\left(  I:K\right)  .
\]

\item[\textbf{(d)}] Let $I$, $J$ and $K$ be three ideals of $R$. Prove that
\[
\left(  I:\left(  J+K\right)  \right)  =\left(  I:J\right)  \cap\left(
I:K\right)  .
\]

\end{enumerate}
\end{exercise}

\subsection{\label{sec.rings.crt}The Chinese Remainder Theorem (\cite[\S 7.6]%
{DumFoo04})}

\subsubsection{Introduction}

In Subsection \ref{subsec.rings.dirprods.exas}, we have seen some examples of
direct products. These examples were not very surprising; they were rings
defined in a way that makes the product structure already quite evident. For
example, the ring of diagonal $n\times n$-matrices was a direct product
because you can easily see that the diagonal entries of diagonal matrices
don't \textquotedblleft interfere\textquotedblright\ with each other when the
matrices are multiplied. Keywords like \textquotedblleft
entrywise\textquotedblright, \textquotedblleft pointwise\textquotedblright%
\ and \textquotedblleft coordinatewise\textquotedblright\ tend to signal that
some structure is a direct product. The $6$-element ring $\mathbb{Z}/6$, on
the other hand, does not look at all like a direct product. Yet, it is
isomorphic to a direct product:%
\[
\mathbb{Z}/6\cong\left(  \mathbb{Z}/2\right)  \times\left(  \mathbb{Z}%
/3\right)  .
\]
Specifically, there is a ring isomorphism%
\[
\mathbb{Z}/6\rightarrow\left(  \mathbb{Z}/2\right)  \times\left(
\mathbb{Z}/3\right)  ,
\]
which sends%
\begin{align*}
\overline{0}  &  \mapsto\left(  \overline{0},\overline{0}\right)
\ \ \ \ \ \ \ \ \ \ \left(  \text{that is, }0+6\mathbb{Z}\mapsto\left(
0+2\mathbb{Z},\ 0+3\mathbb{Z}\right)  \right)  ,\\
\overline{1}  &  \mapsto\left(  \overline{1},\overline{1}\right)  ,\\
\overline{2}  &  \mapsto\left(  \overline{2},\overline{2}\right)  =\left(
\overline{0},\overline{2}\right)  ,\\
\overline{3}  &  \mapsto\left(  \overline{3},\overline{3}\right)  =\left(
\overline{1},\overline{0}\right)  ,\\
\overline{4}  &  \mapsto\left(  \overline{4},\overline{4}\right)  =\left(
\overline{0},\overline{1}\right)  ,\\
\overline{5}  &  \mapsto\left(  \overline{5},\overline{5}\right)  =\left(
\overline{1},\overline{2}\right)  .
\end{align*}
The reason why this works is that $2$ and $3$ are coprime. More generally:

\begin{theorem}
[The Chinese Remainder Theorem for two integers]\label{thm.CRT-2-ints}Let $n$
and $m$ be two coprime integers. Then,%
\[
\mathbb{Z}/\left(  nm\right)  \cong\left(  \mathbb{Z}/n\right)  \times\left(
\mathbb{Z}/m\right)  \ \ \ \ \ \ \ \ \ \ \text{as rings.}%
\]
More concretely, there is a ring isomorphism%
\[
\mathbb{Z}/\left(  nm\right)  \rightarrow\left(  \mathbb{Z}/n\right)
\times\left(  \mathbb{Z}/m\right)
\]
that sends each residue class $\overline{r}$ to the pair $\left(  \overline
{r},\ \overline{r}\right)  $ (or, to use somewhat less ambiguous notation,
sends each residue class $r+nm\mathbb{Z}$ to the pair $\left(  r+n\mathbb{Z}%
,\ r+m\mathbb{Z}\right)  $).
\end{theorem}

Rather than prove this theorem in this form, I will generalize it and then
prove the generalization. After all, this is a course on rings, not just on
$\mathbb{Z}/n$. So I will state and prove a \textquotedblleft Chinese
Remainder Theorem\textquotedblright\ for arbitrary rings. In this theorem,
$\mathbb{Z}$ will be replaced by an arbitrary ring $R$, and the integers $n$
and $m$ will be replaced by two ideals $I$ and $J$ of $R$ (since ideals are
what we can quotient rings by)\footnote{I could also replace the integers $n$
and $m$ by two elements of $R$, but that would be less general: Quotienting by
an element is tantamount to quotienting by a principal ideal, and principal
ideals are just one kind of ideals.}. The condition \textquotedblleft$n$ and
$m$ are coprime\textquotedblright\ will be replaced by the condition
\textquotedblleft$I+J=R$\textquotedblright. Indeed, two integers $n$ and $m$
are coprime if and only if the corresponding principal ideals $I=n\mathbb{Z}$
and $J=m\mathbb{Z}$ of $\mathbb{Z}$ satisfy $I+J=\mathbb{Z}$ (this follows
easily from Proposition \ref{prop.rings.ideal-arith.Z} \textbf{(c)}%
\footnote{\textit{Proof.} Let $n$ and $m$ be two integers. Let $I=n\mathbb{Z}$
and $J=m\mathbb{Z}$ be the corresponding principal ideals of $\mathbb{Z}$.
Then, Proposition \ref{prop.rings.ideal-arith.Z} \textbf{(c)} yields
$I+J=\gcd\left(  n,m\right)  \mathbb{Z}$. If $n$ and $m$ are coprime, then
$\gcd\left(  n,m\right)  =1$, so this rewrites as $I+J=1\mathbb{Z}=\mathbb{Z}%
$. Conversely, if $I+J=\mathbb{Z}$, then $1\in\mathbb{Z}=I+J=\gcd\left(
n,m\right)  \mathbb{Z}$, which shows that $1$ is a multiple of $\gcd\left(
n,m\right)  $; but this entails that $\gcd\left(  n,m\right)  =1$, and
therefore $n$ and $m$ are coprime. Thus, we have shown that $n$ and $m$ are
coprime if and only if $I+J=\mathbb{Z}$.}). Two ideals $I$ and $J$ of a ring
$R$ satisfying $I+J=R$ are said to be \textbf{comaximal}:

\subsubsection{The Chinese Remainder Theorem for two ideals}

\begin{definition}
\label{def.ideals.comax}Let $I$ and $J$ be two ideals of a ring $R$. We say
that $I$ and $J$ are \textbf{comaximal} if $I+J=R$.
\end{definition}

Now we can state the general version of the Chinese Remainder Theorem. We will
state this version in two parts, since they have slightly different
assumptions (the first part requires $R$ to be commutative, while the second
one does not). Both parts will have to be combined to recover Theorem
\ref{thm.CRT-2-ints} later.

\begin{theorem}
[The Chinese Remainder Theorem for two ideals: ideal part]%
\label{thm.CRT-2-ids1}Let $I$ and $J$ be two comaximal ideals of a commutative
ring $R$. (Recall that \textquotedblleft comaximal\textquotedblright\ means
that $I+J=R$.) Then,
\[
I\cap J=IJ.
\]

\end{theorem}

\begin{theorem}
[The Chinese Remainder Theorem for two ideals: quotient part]%
\label{thm.CRT-2-ids2}Let $I$ and $J$ be two comaximal ideals of a ring $R$.
(Recall that \textquotedblleft comaximal\textquotedblright\ means that
$I+J=R$.) Then:

\begin{enumerate}
\item[\textbf{(a)}] We have%
\[
R/\left(  I\cap J\right)  \cong\left(  R/I\right)  \times\left(  R/J\right)
.
\]

\item[\textbf{(b)}] More concretely, there is a ring isomorphism%
\[
R/\left(  I\cap J\right)  \rightarrow\left(  R/I\right)  \times\left(
R/J\right)
\]
that sends each residue class $r+\left(  I\cap J\right)  $ to the pair
$\left(  r+I,\ r+J\right)  $.
\end{enumerate}
\end{theorem}

Let us now prove these theorems. Before we do so, let us agree on a convention
that will save us some parentheses:

\begin{convention}
\label{conv.CRT./-precedence}The \textquotedblleft$/$\textquotedblright\ sign
will have higher precedence than the \textquotedblleft$\times$%
\textquotedblright\ sign, but lower precedence than the \textquotedblleft
implied $\cdot$ sign\textquotedblright. Thus, the expression \textquotedblleft%
$\left(  R/I\right)  \times\left(  R/J\right)  $\textquotedblright\ can be
abbreviated as \textquotedblleft$R/I\times R/J$\textquotedblright\ (without
worrying that it might be misunderstood as \textquotedblleft$R/\left(  I\times
R\right)  /J$\textquotedblright, whatever this would mean), and similarly the
expression \textquotedblleft$R/\left(  IJ\right)  $\textquotedblright\ can be
abbreviated as \textquotedblleft$R/IJ$\textquotedblright\ (without worrying
that it might be misunderstood as \textquotedblleft$\left(  R/I\right)
J$\textquotedblright).
\end{convention}

\begin{proof}
[Proof of Theorem \ref{thm.CRT-2-ids1}.]We have $1\in R=I+J$ (since $I$ and
$J$ are comaximal). In other words,%
\[
\text{there exist }i\in I\text{ and }j\in J\text{ with }1=i+j.
\]
Consider these $i$ and $j$.

Proposition \ref{prop.rings.ideal-arith.laws} \textbf{(b)} yields $IJ\subseteq
I\cap J$. Thus, we only need to show that $I\cap J\subseteq IJ$.

\begin{noncompile}
It is easy to see that the commutativity of $R$ spreads to the ideals of $R$,
too: We have $IJ=JI$.
\end{noncompile}

So let $a\in I\cap J$. Thus, $a\in I$ and $a\in J$. Now,%
\begin{align*}
a  &  =a\cdot\underbrace{1}_{=i+j}=a\cdot\left(  i+j\right)  =\underbrace{ai}%
_{\substack{=ia\\\text{(since }R\text{ is}\\\text{commutative)}}%
}+\,aj=\underbrace{ia}_{\substack{\in IJ\\\text{(since }i\in I\text{ and }a\in
J\text{)}}}+\underbrace{aj}_{\substack{\in IJ\\\text{(since }a\in I\text{ and
}j\in J\text{)}}}\\
&  \in IJ+IJ=IJ.
\end{align*}
(The last equality relied on the fact that $K+K=K$ for any ideal $K$ of $R$.
This is an easy consequence of the fact that $K$ is a subgroup of the additive
group $\left(  R,+,0\right)  $.)

Forget that we fixed $a$. We thus have shown that $a\in IJ$ for each $a\in
I\cap J$. In other words, $I\cap J\subseteq IJ$. As we said above, this
completes the proof of Theorem \ref{thm.CRT-2-ids1}.
\end{proof}

\begin{proof}
[Proof of Theorem \ref{thm.CRT-2-ids2}.]We have $1\in R=I+J$ (since $I$ and
$J$ are comaximal). In other words,%
\[
\text{there exist }i\in I\text{ and }j\in J\text{ with }1=i+j.
\]
Consider these $i$ and $j$.

Consider the map\footnote{Recall that \textquotedblleft$R/I\times
R/J$\textquotedblright\ means \textquotedblleft$\left(  R/I\right)
\times\left(  R/J\right)  $\textquotedblright.}%
\begin{align*}
f:R  &  \rightarrow R/I\times R/J,\\
r  &  \mapsto\left(  r+I,\ r+J\right)  .
\end{align*}
It is straightforward to see that this map $f$ is a ring morphism (from $R$ to
the direct product $R/I\times R/J$).

Moreover, we claim that $\operatorname*{Ker}f=I\cap J$. Indeed, let
$x\in\operatorname*{Ker}f$. Thus, $f\left(  x\right)  =0_{R/I\times
R/J}=\left(  0+I,\ 0+J\right)  $. Since $f\left(  x\right)  $ was defined to
be $\left(  x+I,\ x+J\right)  $, this means that $\left(  x+I,\ x+J\right)
=\left(  0+I,\ 0+J\right)  $. In other words, $x+I=0+I$ and $x+J=0+J$. In
other words, $x\in I$ and $x\in J$. In other words, $x\in I\cap J$.

Forget that we fixed $x$. We thus have shown that $x\in I\cap J$ for each
$x\in\operatorname*{Ker}f$. In other words, $\operatorname*{Ker}f\subseteq
I\cap J$. Reading this argument in reverse shows that $I\cap J\subseteq
\operatorname*{Ker}f$. Thus,
\[
\operatorname*{Ker}f=I\cap J.
\]

Now, we claim that $f$ is surjective. Indeed, $1=i+j$, so that $1-i=j\in J$
and thus $1+J=i+J$. Now, $i+I=0+I$ (since $i\in I$) and $i+J=1+J$ (since
$1+J=i+J$). But the definition of $f$ yields $f\left(  i\right)  =\left(
i+I,\ i+J\right)  =\left(  0+I,\ 1+J\right)  $ (since $i+I=0+I$ and
$i+J=1+J$). Similarly, $f\left(  j\right)  =\left(  1+I,\ 0+J\right)  $. Now,
for every $x\in R$ and $y\in R$, we have%
\begin{align*}
f\left(  yi+xj\right)   &  =\underbrace{f\left(  y\right)  }%
_{\substack{=\left(  y+I,\ y+J\right)  \\\text{(by the definition of
}f\text{)}}}\underbrace{f\left(  i\right)  }_{=\left(  0+I,\ 1+J\right)
}+\underbrace{f\left(  x\right)  }_{\substack{=\left(  x+I,\ x+J\right)
\\\text{(by the definition of }f\text{)}}}\underbrace{f\left(  j\right)
}_{=\left(  1+I,\ 0+J\right)  }\\
&  \ \ \ \ \ \ \ \ \ \ \ \ \ \ \ \ \ \ \ \ \left(  \text{since }f\text{ is a
ring morphism}\right) \\
&  =\underbrace{\left(  y+I,\ y+J\right)  \left(  0+I,\ 1+J\right)
}_{\substack{=\left(  \left(  y+I\right)  \left(  0+I\right)  ,\ \left(
y+J\right)  \left(  1+J\right)  \right)  \\\text{(since the multiplication of
}R/I\times R/J\\\text{is defined to be entrywise)}}}+\underbrace{\left(
x+I,\ x+J\right)  \left(  1+I,\ 0+J\right)  }_{\substack{=\left(  \left(
x+I\right)  \left(  1+I\right)  ,\ \left(  x+J\right)  \left(  0+J\right)
\right)  \\\text{(since the multiplication of }R/I\times R/J\\\text{is defined
to be entrywise)}}}\\
&  =\left(  \underbrace{\left(  y+I\right)  \left(  0+I\right)  }%
_{\substack{=y\cdot0+I\\=0+I}},\ \underbrace{\left(  y+J\right)  \left(
1+J\right)  }_{\substack{=y\cdot1+J\\=y+J}}\right)  +\left(
\underbrace{\left(  x+I\right)  \left(  1+I\right)  }_{\substack{=x\cdot
1+I\\=x+I}},\ \underbrace{\left(  x+J\right)  \left(  0+J\right)
}_{\substack{=x\cdot0+J\\=0+J}}\right) \\
&  =\left(  0+I,\ y+J\right)  +\left(  x+I,\ 0+J\right) \\
&  =\left(  \underbrace{\left(  0+I\right)  +\left(  x+I\right)
}_{\substack{=0+x+I\\=x+I}},\ \underbrace{\left(  y+J\right)  +\left(
0+J\right)  }_{\substack{=y+0+J\\=y+J}}\right) \\
&  \ \ \ \ \ \ \ \ \ \ \ \ \ \ \ \ \ \ \ \ \left(  \text{since the addition of
}R/I\times R/J\text{ is defined to be entrywise}\right) \\
&  =\left(  x+I,\ y+J\right)  ,
\end{align*}
which shows that the pair $\left(  x+I,\ y+J\right)  $ lies in the image of
$f$.

Thus, every element of the form $\left(  x+I,\ y+J\right)  $ for some $x\in R$
and $y\in R$ lies in the image of $f$. Since every element of $R/I\times R/J$
has this form\footnote{because every element of $R/I$ has the form $x+I$ for
some $x\in R$, while every element of $R/J$ has the form $y+J$ for some $y\in
R$}, we thus conclude that every element of $R/I\times R/J$ lies in the image
of $f$. In other words, $f$ is surjective. In other words, $f\left(  R\right)
=R/I\times R/J$.

Now, recall the First isomorphism theorem for rings (Theorem
\ref{thm.1it.ring1} \textbf{(c)}). Applying it to $S=R/I\times R/J$ (and to
our ring morphism $f:R\rightarrow R/I\times R/J$), we see that the map
\begin{align*}
f^{\prime}:R/\operatorname*{Ker}f  &  \rightarrow f\left(  R\right)  ,\\
\overline{r}  &  \mapsto f\left(  r\right)
\end{align*}
is well-defined and is a ring isomorphism.

In our case right now, we have $f\left(  R\right)  =R/I\times R/J$ and
$\operatorname{Ker}f=I\cap J$, so that we can restate this as follows: The
map
\begin{align*}
f^{\prime}:R/\left(  I\cap J\right)   &  \rightarrow R/I\times R/J,\\
\overline{r}  &  \mapsto f\left(  r\right)
\end{align*}
is well-defined and is a ring isomorphism. This map $f^{\prime}$ sends each
residue class $\overline{r}=r+\left(  I\cap J\right)  $ to $f\left(  r\right)
=\left(  r+I,\ r+J\right)  $ (by the definition of $f$). Thus, we have found a
ring isomorphism%
\[
R/\left(  I\cap J\right)  \rightarrow R/I\times R/J
\]
that sends each residue class $r+\left(  I\cap J\right)  $ to the pair
$\left(  r+I,\ r+J\right)  $ (namely, $f^{\prime}$). This proves part
\textbf{(b)} of Theorem \ref{thm.CRT-2-ids2}. Of course, part \textbf{(a)}
thus follows.
\end{proof}

You can get rid of the commutativity requirement on $R$ in Theorem
\ref{thm.CRT-2-ids1} if you replace $IJ$ by $IJ+JI$. (Checking this is a nice
exercise on making sure you understand the above proof.)

\subsubsection{Application to integers}

As a corollary of Theorems \ref{thm.CRT-2-ids1} and \ref{thm.CRT-2-ids2}, we
can now prove the good old number-theoretical Chinese Remainder Theorem
(Theorem \ref{thm.CRT-2-ints}), which we will repeat for convenience:

\begin{theorem}
[The Chinese Remainder Theorem for two integers]\label{thm.CRT-2-ints-rep}Let
$n$ and $m$ be two coprime integers. Then,%
\[
\mathbb{Z}/\left(  nm\right)  \cong\left(  \mathbb{Z}/n\right)  \times\left(
\mathbb{Z}/m\right)  \ \ \ \ \ \ \ \ \ \ \text{as rings.}%
\]
More concretely, there is a ring isomorphism%
\[
\mathbb{Z}/\left(  nm\right)  \rightarrow\left(  \mathbb{Z}/n\right)
\times\left(  \mathbb{Z}/m\right)
\]
that sends each residue class $\overline{r}$ to the pair $\left(  \overline
{r},\ \overline{r}\right)  $ (or, to use somewhat less ambiguous notation,
sends each residue class $r+nm\mathbb{Z}$ to the pair $\left(  r+n\mathbb{Z}%
,\ r+m\mathbb{Z}\right)  $).
\end{theorem}

\begin{proof}
Let $R=\mathbb{Z}$ and $I=n\mathbb{Z}$ and $J=m\mathbb{Z}$. Proposition
\ref{prop.rings.ideal-arith.Z} then yields $IJ=nm\mathbb{Z}$ and $I\cap
J=\operatorname{lcm}\left(  n,m\right)  \mathbb{Z}$ and $I+J=\gcd\left(
n,m\right)  \mathbb{Z}$. Since $n$ and $m$ are coprime, we have $\gcd\left(
n,m\right)  =1$; thus, $I+J=\underbrace{\gcd\left(  n,m\right)  }%
_{=1}\mathbb{Z}=1\mathbb{Z}=\mathbb{Z}$. In other words, the ideals $I$ and
$J$ of $\mathbb{Z}$ are comaximal. Hence, Theorem \ref{thm.CRT-2-ids1} yields
$I\cap J=IJ=nm\mathbb{Z}$. Furthermore, part \textbf{(a)} of Theorem
\ref{thm.CRT-2-ids2} yields $R/\left(  I\cap J\right)  \cong\left(
R/I\right)  \times\left(  R/J\right)  $. In view of $\underbrace{R}%
_{=\mathbb{Z}}/\underbrace{\left(  I\cap J\right)  }_{=nm\mathbb{Z}%
}=\mathbb{Z}/\left(  nm\mathbb{Z}\right)  =\mathbb{Z}/\left(  nm\right)  $ and
$\underbrace{R}_{=\mathbb{Z}}/\underbrace{I}_{=n\mathbb{Z}}=\mathbb{Z}/\left(
n\mathbb{Z}\right)  =\mathbb{Z}/n$ and $\underbrace{R}_{=\mathbb{Z}%
}/\underbrace{J}_{=m\mathbb{Z}}=\mathbb{Z}/\left(  m\mathbb{Z}\right)
=\mathbb{Z}/m$, this rewrites as $\mathbb{Z}/\left(  nm\right)  \cong\left(
\mathbb{Z}/n\right)  \times\left(  \mathbb{Z}/m\right)  $. This proves the
first claim of Theorem \ref{thm.CRT-2-ints-rep}. The \textquotedblleft More
concretely\textquotedblright\ claim likewise follows from part \textbf{(b)} of
Theorem \ref{thm.CRT-2-ids2}.
\end{proof}

\subsubsection{Comaximality for products of ideals}

We shall next prove some auxiliary results about comaximal ideals, which will
later help us generalize Theorem \ref{thm.CRT-2-ids1} and Theorem
\ref{thm.CRT-2-ids2} to $k$ ideals instead of $2$ ideals. These can also serve
as exercises on ideal arithmetic.

\begin{proposition}
\label{prop.comax.prod2}Let $I,J,K$ be three ideals of a ring $R$. Then:

\begin{enumerate}
\item[\textbf{(a)}] We have $\left(  I+K\right)  \left(  J+K\right)  \subseteq
IJ+K$.

\item[\textbf{(b)}] If $I+K=R$ and $J+K=R$, then $IJ+K=R$.
\end{enumerate}
\end{proposition}

\begin{proof}
Using Proposition \ref{prop.rings.ideal-arith.laws} \textbf{(a)}, we easily
see that the sets $I+K$ and $J+K$ and $IJ+K$ are ideals of $R$. Hence, these
sets are closed under addition. \medskip

\textbf{(a)} Let $x\in I+K$ and $y\in J+K$. We shall show that $xy\in IJ+K$.

Indeed, $x\in I+K$. In other words, we can write $x$ in the form $x=i+a$ for
some $i\in I$ and $a\in K$. Consider these $i$ and $a$.

Also, $y\in J+K$. In other words, we can write $y$ in the form $y=j+b$ for
some $j\in J$ and $b\in K$. Consider these $j$ and $b$.

Now, multiplying the two equalities $x=i+a$ and $y=j+b$, we obtain%
\[
xy=\left(  i+a\right)  \left(  j+b\right)  =ij+ib+aj+ab.
\]

From $a\in K$ and $b\in K$, we conclude that each of the three products
$ib,\ aj,\ ab$ belongs to $K$ (by the second ideal axiom, since $K$ is an
ideal). Hence, the sum of these three products also belongs to $K$ (since $K$
is an ideal and thus is closed under addition). In other words, $ib+aj+ab\in
K$. Furthermore, from $i\in I$ and $j\in J$, we see that $ij$ is an $\left(
I,J\right)  $-product and therefore a sum of $\left(  I,J\right)  $-products
(namely, of just one such product). Hence, $ij\in IJ$ (by the definition of
$IJ$). Now,%
\[
xy=\underbrace{ij}_{\in IJ}+\underbrace{ib+aj+ab}_{\in K}\in IJ+K.
\]

Forget that we fixed $x$ and $y$. We thus have shown that $xy\in IJ+K$ for
each $x\in I+K$ and $y\in J+K$. In other words, every $\left(
I+K,\ J+K\right)  $-product belongs to $IJ+K$ (since every $\left(
I+K,\ J+K\right)  $-product has the form $xy$ for some $x\in I+K$ and $y\in
J+K$).

Since $IJ+K$ is closed under addition, we thus conclude that any finite sum of
$\left(  I+K,\ J+K\right)  $-products belongs to $IJ+K$ as well. However, the
definition of $\left(  I+K\right)  \left(  J+K\right)  $ yields%
\[
\left(  I+K\right)  \left(  J+K\right)  =\left\{  \text{all finite sums of
}\left(  I+K,\ J+K\right)  \text{-products}\right\}  \subseteq IJ+K
\]
(since any finite sum of $\left(  I+K,\ J+K\right)  $-products belongs to
$IJ+K$). This proves Proposition \ref{prop.comax.prod2} \textbf{(a)}. \medskip

\textbf{(b)} Assume that $I+K=R$ and $J+K=R$. Proposition
\ref{prop.rings.ideal-arith.laws} \textbf{(e)} yields that $RR=R$. Hence,%
\[
R=\underbrace{R}_{=I+K}\ \ \underbrace{R}_{=J+K}=\left(  I+K\right)  \left(
J+K\right)  \subseteq IJ+K\ \ \ \ \ \ \ \ \ \ \left(  \text{by Proposition
\ref{prop.comax.prod2} \textbf{(a)}}\right)  .
\]
Combined with $IJ+K\subseteq R$ (which is obvious), this yields $IJ+K=R$.
Thus, Proposition \ref{prop.comax.prod2} \textbf{(b)} is proved.
\end{proof}

\begin{exercise}
\ \ 

\begin{enumerate}
\item[\textbf{(a)}] Recall the following fact from elementary number theory:
If $a,b,c$ are three integers such that each of $a$ and $b$ is coprime to $c$,
then $ab$ is also coprime to $c$. How does Proposition \ref{prop.comax.prod2}
\textbf{(b)} generalize this fact?

\item[\textbf{(b)}] What property of greatest common divisors of integers does
Proposition \ref{prop.comax.prod2} \textbf{(a)} generalize?
\end{enumerate}
\end{exercise}

Proposition \ref{prop.comax.prod2} \textbf{(b)} can be extended by replacing
the two ideals $I$ and $J$ by $k$ ideals $I_{1},I_{2},\ldots,I_{k}$:

\begin{proposition}
\label{prop.comax.prodn}Let $I_{1},I_{2},\ldots,I_{k}$ be $k$ ideals of a ring
$R$. Let $K$ be a further ideal of $R$. Assume that
\begin{equation}
I_{i}+K=R\ \ \ \ \ \ \ \ \ \ \text{for each }i\in\left\{  1,2,\ldots
,k\right\}  . \label{eq.prop.comax.prodn.ass}%
\end{equation}
Then, $I_{1}I_{2}\cdots I_{k}+K=R$.
\end{proposition}

\begin{proof}
We proceed by induction on $k$:

\textit{Base case:} The ideal $R$ is the neutral element of the monoid of
ideals of $R$ under multiplication (see Proposition
\ref{prop.rings.ideal-arith.laws} \textbf{(e)}). Thus, the empty product of
ideals of $R$ is defined to be $R$.

Now, in the case $k=0$, the product $I_{1}I_{2}\cdots I_{k}$ is an empty
product of ideals and therefore equals $R$ (by the previous sentence). Hence,
in this case, we have $\underbrace{I_{1}I_{2}\cdots I_{k}}_{=R}+\,K=R+K=R$
(this is very easy to check). Thus, Proposition \ref{prop.comax.prodn} is
proved for $k=0$.

\textit{Induction step:} Let $m$ be a positive integer. Assume (as the
induction hypothesis) that Proposition \ref{prop.comax.prodn} holds for
$k=m-1$. We must prove that Proposition \ref{prop.comax.prodn} holds for $k=m$.

So let $I_{1},I_{2},\ldots,I_{m}$ be $m$ ideals of a ring $R$. Let $K$ be a
further ideal of $R$. Assume that
\begin{equation}
I_{i}+K=R\ \ \ \ \ \ \ \ \ \ \text{for each }i\in\left\{  1,2,\ldots
,m\right\}  . \label{pf.prop.comax.prodn.ass}%
\end{equation}
We must then show that $I_{1}I_{2}\cdots I_{m}+K=R$.

The equality (\ref{pf.prop.comax.prodn.ass}) holds for each $i\in\left\{
1,2,\ldots,m\right\}  $, and thus in particular for each $i\in\left\{
1,2,\ldots,m-1\right\}  $. Hence, we can apply Proposition
\ref{prop.comax.prodn} to $k=m-1$ (since our induction hypothesis says that
Proposition \ref{prop.comax.prodn} holds for $k=m-1$). Thus we obtain
$I_{1}I_{2}\cdots I_{m-1}+K=R$. Since we also have $I_{m}+K=R$ (by
(\ref{pf.prop.comax.prodn.ass}), applied to $i=m$), we can thus apply
Proposition \ref{prop.comax.prod2} \textbf{(b)} to $I=I_{1}I_{2}\cdots
I_{m-1}$ and $J=I_{m}$. We thus obtain%
\[
\left(  I_{1}I_{2}\cdots I_{m-1}\right)  I_{m}+K=R.
\]
In other words, $I_{1}I_{2}\cdots I_{m}+K=R$ (since $\left(  I_{1}I_{2}\cdots
I_{m-1}\right)  I_{m}=I_{1}I_{2}\cdots I_{m}$). Thus, we have proved that
Proposition \ref{prop.comax.prodn} holds for $k=m$. This completes the
induction step, and thus Proposition \ref{prop.comax.prodn} is proved by induction.
\end{proof}

\begin{lemma}
\label{lem.comax.prod-in-cut}Let $I_{1},I_{2},\ldots,I_{k}$ be $k$ ideals of a
ring $R$. Then, $I_{1}I_{2}\cdots I_{k}\subseteq I_{1}\cap I_{2}\cap\cdots\cap
I_{k}$.
\end{lemma}

\begin{proof}
We proceed by induction on $k$:

\textit{Base case:} In the case $k=0$, both $I_{1}I_{2}\cdots I_{k}$ and
$I_{1}\cap I_{2}\cap\cdots\cap I_{k}$ are the ideal $R$ (indeed, this can be
seen as in the proof of Proposition \ref{prop.comax.prodn}, since $R$ is the
neutral element for both the multiplication and the intersection of ideals of
$R$). Thus, in the case $k=0$, we have $I_{1}I_{2}\cdots I_{k}\subseteq
I_{1}\cap I_{2}\cap\cdots\cap I_{k}$ (since $R\subseteq R$). In other words,
Lemma \ref{lem.comax.prod-in-cut} is proved for $k=0$.

\textit{Induction step:} Let $m$ be a positive integer. Assume (as the
induction hypothesis) that Lemma \ref{lem.comax.prod-in-cut} holds for
$k=m-1$. We must prove that Lemma \ref{lem.comax.prod-in-cut} holds for $k=m$.

So let $I_{1},I_{2},\ldots,I_{m}$ be $m$ ideals of a ring $R$. We must show
that $I_{1}I_{2}\cdots I_{m}\subseteq I_{1}\cap I_{2}\cap\cdots\cap I_{m}$.

We can apply Lemma \ref{lem.comax.prod-in-cut} to $k=m-1$ (since our induction
hypothesis says that Lemma \ref{lem.comax.prod-in-cut} holds for $k=m-1$).
Thus we obtain $I_{1}I_{2}\cdots I_{m-1}\subseteq I_{1}\cap I_{2}\cap
\cdots\cap I_{m-1}$.

However, Proposition \ref{prop.rings.ideal-arith.laws} \textbf{(b)} yields
that $IJ\subseteq I\cap J$ for any two ideals $I$ and $J$ of $R$. Applying
this to $I=I_{1}I_{2}\cdots I_{m-1}$ and $J=I_{m}$, we obtain
\begin{align*}
\left(  I_{1}I_{2}\cdots I_{m-1}\right)  I_{m}  &  \subseteq
\underbrace{\left(  I_{1}I_{2}\cdots I_{m-1}\right)  }_{\subseteq I_{1}\cap
I_{2}\cap\cdots\cap I_{m-1}}\cap\,I_{m}\subseteq\left(  I_{1}\cap I_{2}%
\cap\cdots\cap I_{m-1}\right)  \cap I_{m}\\
&  =I_{1}\cap I_{2}\cap\cdots\cap I_{m}.
\end{align*}
Since $\left(  I_{1}I_{2}\cdots I_{m-1}\right)  I_{m}=I_{1}I_{2}\cdots I_{m}$,
we can rewrite this as%
\[
I_{1}I_{2}\cdots I_{m}\subseteq I_{1}\cap I_{2}\cap\cdots\cap I_{m}.
\]
Thus, we have proved that Lemma \ref{lem.comax.prod-in-cut} holds for $k=m$.
This completes the induction step, and thus Lemma \ref{lem.comax.prod-in-cut}
is proved by induction.
\end{proof}

\begin{corollary}
\label{cor.comax.sectn}Let $I_{1},I_{2},\ldots,I_{k}$ be $k$ ideals of a ring
$R$. Let $K$ be a further ideal of $R$. Assume that
\[
I_{i}+K=R\ \ \ \ \ \ \ \ \ \ \text{for each }i\in\left\{  1,2,\ldots
,k\right\}  .
\]
Then, $\left(  I_{1}\cap I_{2}\cap\cdots\cap I_{k}\right)  +K=R$.
\end{corollary}

\begin{proof}
Lemma \ref{lem.comax.prod-in-cut} yields $I_{1}I_{2}\cdots I_{k}\subseteq
I_{1}\cap I_{2}\cap\cdots\cap I_{k}$.

It is easy to see that if $I,J,L$ are three ideals of $R$ such that
$I\subseteq J$, then $I+L\subseteq J+L$. Applying this to $I=I_{1}I_{2}\cdots
I_{k}$, $J=I_{1}\cap I_{2}\cap\cdots\cap I_{k}$ and $L=K$, we obtain%
\[
I_{1}I_{2}\cdots I_{k}+K\subseteq\left(  I_{1}\cap I_{2}\cap\cdots\cap
I_{k}\right)  +K
\]
(since $I_{1}I_{2}\cdots I_{k}\subseteq I_{1}\cap I_{2}\cap\cdots\cap I_{k}$).
Hence,%
\[
\left(  I_{1}\cap I_{2}\cap\cdots\cap I_{k}\right)  +K\supseteq I_{1}%
I_{2}\cdots I_{k}+K=R\ \ \ \ \ \ \ \ \ \ \left(  \text{by Proposition
\ref{prop.comax.prodn}}\right)  .
\]
Combining this with $\left(  I_{1}\cap I_{2}\cap\cdots\cap I_{k}\right)
+K\subseteq R$ (which is obvious), we obtain $\left(  I_{1}\cap I_{2}%
\cap\cdots\cap I_{k}\right)  +K=R$. Thus, Corollary \ref{cor.comax.sectn} is proven.
\end{proof}

\subsubsection{The Chinese Remainder Theorem for $k$ ideals}

As we already suggested, Theorem \ref{thm.CRT-2-ids1} and Theorem
\ref{thm.CRT-2-ids2} can be generalized to $k$ ideals. First, a convention:

\begin{definition}
Let $I_{1},I_{2},\ldots,I_{k}$ be $k$ ideals of a ring $R$. We say that these
$k$ ideals $I_{1},I_{2},\ldots,I_{k}$ are \textbf{mutually comaximal} if
$I_{i}+I_{j}=R$ holds for all $1\leq i<j\leq k$.
\end{definition}

In other words, $k$ ideals $I_{1},I_{2},\ldots,I_{k}$ are mutually comaximal
if $I_{i}$ and $I_{j}$ are comaximal for every $i<j$. When $k>2$, this is a
\textbf{much stronger} statement than $I_{1}+I_{2}+\cdots+I_{k}=R$.

For example, if $n_{1},n_{2},\ldots,n_{k}$ are $k$ arbitrary integers, then
the $k$ principal ideals $n_{1}\mathbb{Z},n_{2}\mathbb{Z},\ldots
,n_{k}\mathbb{Z}$ are mutually comaximal if $n_{1},n_{2},\ldots,n_{k}$ are
mutually coprime (that is, if $n_{i}$ is coprime to $n_{j}$ for all $i<j$).
When $k>2$, this is a \textbf{much stronger} statement than $\gcd\left(
n_{1},n_{2},\ldots,n_{k}\right)  =1$. Be warned! Lots of mistakes have been
made by mistaking \textquotedblleft mutually coprime\textquotedblright\ for
\textquotedblleft gcd of all $k$ numbers is $1$\textquotedblright.

Enough of the warning labels; here are the theorems:

\begin{theorem}
[The Chinese Remainder Theorem for $k$ ideals: ideal part]%
\label{thm.CRT-k-ids1}Let $I_{1},I_{2},\ldots,I_{k}$ be $k$ mutually comaximal
ideals of a commutative ring $R$. Then,%
\[
I_{1}\cap I_{2}\cap\cdots\cap I_{k}=I_{1}I_{2}\cdots I_{k}.
\]

\end{theorem}

\begin{theorem}
[The Chinese Remainder Theorem for $k$ ideals: quotient part]%
\label{thm.CRT-k-ids2}Let $I_{1},I_{2},\ldots,I_{k}$ be $k$ mutually comaximal
ideals of a ring $R$. Then:

\begin{enumerate}
\item[\textbf{(a)}] We have
\[
R/\left(  I_{1}\cap I_{2}\cap\cdots\cap I_{k}\right)  \cong R/I_{1}\times
R/I_{2}\times\cdots\times R/I_{k}.
\]

\item[\textbf{(b)}] More concretely, there is a ring isomorphism%
\[
R/\left(  I_{1}\cap I_{2}\cap\cdots\cap I_{k}\right)  \rightarrow
R/I_{1}\times R/I_{2}\times\cdots\times R/I_{k}%
\]
that sends each residue class $r+\left(  I_{1}\cap I_{2}\cap\cdots\cap
I_{k}\right)  $ to the $k$-tuple $\left(  r+I_{1},\ r+I_{2},\ \ldots
,\ r+I_{k}\right)  $.
\end{enumerate}
\end{theorem}

\begin{proof}
[Proof of Theorem \ref{thm.CRT-k-ids1}.]We proceed by induction on $k$:

\textit{Induction base:} You can take $k=1$ as a base case (it is utterly
trivial), or even $k=0$ if you are brave enough\footnote{Make sure to
understand the empty product of ideals of $R$ to be $R$ itself, since $R$ is
the neutral element of the monoid of ideals of $R$ under multiplication (see
Proposition \ref{prop.rings.ideal-arith.laws} \textbf{(e)}).
\par
Likewise, the empty intersection of ideals of $R$ is $R$ itself, since $R$ is
the neutral element of the monoid of ideals of $R$ under intersection (see
Proposition \ref{prop.rings.ideal-arith.laws} \textbf{(d)}).}.

\textit{Induction step:} Let $n$ be a positive integer. (You can assume $n>1$
if it makes you sleep better.) Assume (as the induction hypothesis) that
Theorem \ref{thm.CRT-k-ids1} holds for $k=n-1$. We must now prove that Theorem
\ref{thm.CRT-k-ids1} holds for $k=n$.

So let $I_{1},I_{2},\ldots,I_{n}$ be $n$ mutually comaximal ideals of a
commutative ring $R$. Then, the induction hypothesis yields that Theorem
\ref{thm.CRT-k-ids1} holds for $I_{1},I_{2},\ldots,I_{n-1}$. In particular, we
have%
\begin{equation}
I_{1}\cap I_{2}\cap\cdots\cap I_{n-1}=I_{1}I_{2}\cdots I_{n-1}.
\label{pf.CRT-k-ids.IH.a}%
\end{equation}

Recall that the ideals $I_{1},I_{2},\ldots,I_{n}$ are mutually comaximal.
Hence, for each $i\in\left\{  1,2,\ldots,n-1\right\}  $, the ideals $I_{i}$
and $I_{n}$ are comaximal, i.e., satisfy $I_{i}+I_{n}=R$. Hence, Corollary
\ref{cor.comax.sectn} (applied to $k=n-1$ and $K=I_{n}$) yields that $\left(
I_{1}\cap I_{2}\cap\cdots\cap I_{n-1}\right)  +I_{n}=R$. In other words, the
two ideals $I_{1}\cap I_{2}\cap\cdots\cap I_{n-1}$ and $I_{n}$ are comaximal.

Hence, we can apply Theorem \ref{thm.CRT-2-ids1} to $I=I_{1}\cap I_{2}%
\cap\cdots\cap I_{n-1}$ and $J=I_{n}$. We thus obtain
\begin{equation}
\left(  I_{1}\cap I_{2}\cap\cdots\cap I_{n-1}\right)  \cap I_{n}=\left(
I_{1}\cap I_{2}\cap\cdots\cap I_{n-1}\right)  I_{n}. \label{pf.CRT-k-ids.In.a}%
\end{equation}

Now,%
\begin{align*}
I_{1}\cap I_{2}\cap\cdots\cap I_{n}  &  =\left(  I_{1}\cap I_{2}\cap\cdots\cap
I_{n-1}\right)  \cap I_{n}\\
&  =\left(  I_{1}\cap I_{2}\cap\cdots\cap I_{n-1}\right)  I_{n}%
\ \ \ \ \ \ \ \ \ \ \left(  \text{by (\ref{pf.CRT-k-ids.In.a})}\right) \\
&  =\left(  I_{1}I_{2}\cdots I_{n-1}\right)  I_{n}\ \ \ \ \ \ \ \ \ \ \left(
\text{by (\ref{pf.CRT-k-ids.IH.a})}\right) \\
&  =I_{1}I_{2}\cdots I_{n}.
\end{align*}
Thus, we have proved that Theorem \ref{thm.CRT-k-ids1} holds for $k=n$. This
completes the induction step. Thus, Theorem \ref{thm.CRT-k-ids1} is proved.
\end{proof}

\begin{proof}
[Proof of Theorem \ref{thm.CRT-k-ids2}.]We proceed by induction on $k$:

\textit{Induction base:} Again, you can use $k=0$ as the base
case\footnote{Just as in the above proof of Theorem \ref{thm.CRT-k-ids1}, the
empty intersection of ideals of $R$ is $R$ by definition. Thus, $R/\left(
I_{1}\cap I_{2}\cap\cdots\cap I_{k}\right)  =R/R$ is a trivial ring for $k=0$.
\par
Also, keep in mind that an empty direct product (i.e., a direct product of $0$
rings) is a trivial ring whose only element is the $0$-tuple $\left(
{}\right)  $.} (or $k=1$ if you want to avoid trivialities).

\textit{Induction step:} Let $n$ be a positive integer. (Again, you can assume
$n>1$ if you prefer.) Assume (as the induction hypothesis) that Theorem
\ref{thm.CRT-k-ids2} holds for $k=n-1$. We must now prove that Theorem
\ref{thm.CRT-k-ids2} holds for $k=n$.

So let $I_{1},I_{2},\ldots,I_{n}$ be $n$ mutually comaximal ideals of a ring
$R$. Then, the induction hypothesis yields that Theorem \ref{thm.CRT-k-ids2}
holds for $I_{1},I_{2},\ldots,I_{n-1}$. In particular, part \textbf{(a)} of
Theorem \ref{thm.CRT-k-ids2} shows that%
\begin{equation}
R/\left(  I_{1}\cap I_{2}\cap\cdots\cap I_{n-1}\right)  \cong R/I_{1}\times
R/I_{2}\times\cdots\times R/I_{n-1}. \label{pf.CRT-k-ids.IH.b}%
\end{equation}
Furthermore, part \textbf{(b)} of Theorem \ref{thm.CRT-k-ids2} shows that
there is a ring isomorphism%
\begin{equation}
R/\left(  I_{1}\cap I_{2}\cap\cdots\cap I_{n-1}\right)  \rightarrow
R/I_{1}\times R/I_{2}\times\cdots\times R/I_{n-1} \label{pf.CRT-k-ids.IH.c}%
\end{equation}
that does what you would expect it to do (viz., sends each residue class
$r+\left(  I_{1}\cap I_{2}\cap\cdots\cap I_{n-1}\right)  $ to the $\left(
n-1\right)  $-tuple $\left(  r+I_{1},\ r+I_{2},\ \ldots,\ r+I_{n-1}\right)  $).

Recall that the ideals $I_{1},I_{2},\ldots,I_{n}$ are mutually comaximal.
Hence, for each $i\in\left\{  1,2,\ldots,n-1\right\}  $, the ideals $I_{i}$
and $I_{n}$ are comaximal, i.e., satisfy $I_{i}+I_{n}=R$. Hence, Corollary
\ref{cor.comax.sectn} (applied to $k=n-1$ and $K=I_{n}$) yields that $\left(
I_{1}\cap I_{2}\cap\cdots\cap I_{n-1}\right)  +I_{n}=R$. In other words, the
two ideals $I_{1}\cap I_{2}\cap\cdots\cap I_{n-1}$ and $I_{n}$ are comaximal.

Hence, we can apply Theorem \ref{thm.CRT-2-ids2} to $I=I_{1}\cap I_{2}%
\cap\cdots\cap I_{n-1}$ and $J=I_{n}$. We thus obtain (from part \textbf{(a)}
of Theorem \ref{thm.CRT-2-ids2}) that%
\begin{align}
&  R/\left(  \left(  I_{1}\cap I_{2}\cap\cdots\cap I_{n-1}\right)  \cap
I_{n}\right) \nonumber\\
&  \cong R/\left(  I_{1}\cap I_{2}\cap\cdots\cap I_{n-1}\right)  \times
R/I_{n}; \label{pf.CRT-k-ids.In.b}%
\end{align}
furthermore, we obtain (from part \textbf{(b)}) that there is a ring
isomorphism
\begin{align}
&  R/\left(  \left(  I_{1}\cap I_{2}\cap\cdots\cap I_{n-1}\right)  \cap
I_{n}\right) \nonumber\\
&  \rightarrow R/\left(  I_{1}\cap I_{2}\cap\cdots\cap I_{n-1}\right)  \times
R/I_{n} \label{pf.CRT-k-ids.In.c}%
\end{align}
that does what you expect (viz., sends each residue class $r+\left(  \left(
I_{1}\cap I_{2}\cap\cdots\cap I_{n-1}\right)  \cap I_{n}\right)  $ to the pair
$\left(  r+\left(  I_{1}\cap I_{2}\cap\cdots\cap I_{n-1}\right)
,\ r+I_{n}\right)  $).

Now, let us combine the isomorphisms that we have found. This requires a
little bit of yak-shaving. We will need the following lemma:

\begin{lemma}
\label{lem.dirprod.AxCisoBxC}Let $A,B,C$ be three rings.

\begin{enumerate}
\item[\textbf{(a)}] If $A\cong B$, then $A\times C\cong B\times C$.

\item[\textbf{(b)}] More specifically: If $f:A\rightarrow B$ is a ring
isomorphism, then the map%
\[
f\times\operatorname*{id}\nolimits_{C}:A\times C\rightarrow B\times C
\]
(this is the map that sends each pair $\left(  a,\ c\right)  \in A\times C$ to
$\left(  f\left(  a\right)  ,\ \operatorname*{id}\nolimits_{C}\left(
c\right)  \right)  =\left(  f\left(  a\right)  ,\ c\right)  \in B\times C$) is
a ring isomorphism, too.
\end{enumerate}
\end{lemma}

This lemma simply says that if you replace a ring in a direct product by an
isomorphic one, then the whole direct product stays isomorphic too. I won't
offend your intellect with the proof of this lemma; it is a purely
paint-by-numbers affair. Such lemmas are a dime a dozen, and you are supposed
to invent one whenever you need it. The idea behind this lemma is simply that
isomorphisms behave like equalities.

So let us go back to our proof of Theorem \ref{thm.CRT-k-ids2}. We have%
\begin{align*}
R/\left(  I_{1}\cap I_{2}\cap\cdots\cap I_{n}\right)   &  =R/\left(  \left(
I_{1}\cap I_{2}\cap\cdots\cap I_{n-1}\right)  \cap I_{n}\right) \\
&  \cong\underbrace{R/\left(  I_{1}\cap I_{2}\cap\cdots\cap I_{n-1}\right)
}_{\substack{\cong R/I_{1}\times R/I_{2}\times\cdots\times R/I_{n-1}%
\\\text{(by (\ref{pf.CRT-k-ids.IH.b}))}}}\times\,R/I_{n}%
\ \ \ \ \ \ \ \ \ \ \left(  \text{by (\ref{pf.CRT-k-ids.In.b})}\right) \\
&  \cong\left(  R/I_{1}\times R/I_{2}\times\cdots\times R/I_{n-1}\right)
\times R/I_{n}\\
&  \ \ \ \ \ \ \ \ \ \ \ \ \ \ \ \ \ \ \ \ \left(  \text{by Lemma
\ref{lem.dirprod.AxCisoBxC} \textbf{(a)}}\right) \\
&  \cong R/I_{1}\times R/I_{2}\times\cdots\times R/I_{n};
\end{align*}
this proves part \textbf{(a)} of Theorem \ref{thm.CRT-k-ids2} for $k=n$.

It remains to prove part \textbf{(b)}. Here we will need Lemma
\ref{lem.dirprod.AxCisoBxC} \textbf{(b)}. Indeed, (\ref{pf.CRT-k-ids.IH.c})
gives us a ring isomorphism $R/\left(  I_{1}\cap I_{2}\cap\cdots\cap
I_{n-1}\right)  \rightarrow R/I_{1}\times R/I_{2}\times\cdots\times R/I_{n-1}%
$, which we can call $f$; thus, Lemma \ref{lem.dirprod.AxCisoBxC} \textbf{(b)}
yields a ring isomorphism%
\begin{align*}
R/\left(  I_{1}\cap I_{2}\cap\cdots\cap I_{n-1}\right)  \times R/I_{n}  &
\rightarrow\left(  R/I_{1}\times R/I_{2}\times\cdots\times R/I_{n-1}\right)
\times R/I_{n},\\
\left(  a,\ c\right)   &  \mapsto\left(  f\left(  a\right)  ,\ c\right)  .
\end{align*}
Now, we compose the arrows in our quiver:%
\begin{align*}
&  R/\left(  I_{1}\cap I_{2}\cap\cdots\cap I_{n}\right) \\
&  =R/\left(  \left(  I_{1}\cap I_{2}\cap\cdots\cap I_{n-1}\right)  \cap
I_{n}\right) \\
&  \rightarrow R/\left(  I_{1}\cap I_{2}\cap\cdots\cap I_{n-1}\right)  \times
R/I_{n}\ \ \ \ \ \ \ \ \ \ \left(  \text{this is the isomorphism from
(\ref{pf.CRT-k-ids.In.c})}\right) \\
&  \rightarrow\left(  R/I_{1}\times R/I_{2}\times\cdots\times R/I_{n-1}%
\right)  \times R/I_{n}\\
&  \ \ \ \ \ \ \ \ \ \ \left(  \text{this is the isomorphism we just
constructed using Lemma \ref{lem.dirprod.AxCisoBxC} \textbf{(b)}}\right) \\
&  \rightarrow R/I_{1}\times R/I_{2}\times\cdots\times R/I_{n}.
\end{align*}
All these arrows are ring isomorphisms; hence, so is their composition. It
remains to show that this isomorphism does what you expect (i.e., sends each
residue class $r+\left(  I_{1}\cap I_{2}\cap\cdots\cap I_{n}\right)  $ to
$\left(  r+I_{1},\ r+I_{2},\ \ldots,\ r+I_{n}\right)  $). This is completely
straightforward, and becomes even more so if you drop the details and just
write $\overline{r}$ for all possible cosets $r+J$ no matter what $J$ is:
Following a coset $\overline{r}=r+\left(  I_{1}\cap I_{2}\cap\cdots\cap
I_{n}\right)  $ through the above arrows, we obtain%
\[
\overline{r}=\overline{r}\mapsto\left(  \overline{r},\ \overline{r}\right)
\mapsto\left(  \left(  \overline{r},\ \overline{r},\ \ldots,\ \overline
{r}\right)  ,\ \overline{r}\right)  \mapsto\left(  \overline{r},\ \overline
{r},\ \ldots,\ \overline{r}\right)  .
\]
While the different $\overline{r}$'s mean different things (namely, they are
cosets for different ideals), we are never in any danger of confusing them for
one another, since we know what sets these maps go between. So the $\left(
\overline{r},\ \overline{r},\ \ldots,\ \overline{r}\right)  $ at the end of
this computation must be $\left(  r+I_{1},\ r+I_{2},\ \ldots,\ r+I_{n}\right)
$, since it is an element of $R/I_{1}\times R/I_{2}\times\cdots\times R/I_{n}%
$. So our isomorphism sends $r+\left(  I_{1}\cap I_{2}\cap\cdots\cap
I_{n}\right)  $ to $\left(  r+I_{1},\ r+I_{2},\ \ldots,\ r+I_{n}\right)  $.
Thus, part \textbf{(b)} of Theorem \ref{thm.CRT-k-ids2} is proved for $k=n$.

Both parts of Theorem \ref{thm.CRT-k-ids2} are thus proved for $k=n$. This
completes the induction step, and thus the proof.
\end{proof}

\subsubsection{Applying to integers again}

We can again apply this to $R=\mathbb{Z}$:

\begin{theorem}
[The Chinese Remainder Theorem for $k$ integers]\label{thm.CRT-k-ints}Let
$n_{1},n_{2},\ldots,n_{k}$ be $k$ mutually coprime integers.
(\textquotedblleft Mutually coprime\textquotedblright\ means that $n_{i}$ is
coprime to $n_{j}$ whenever $i<j$.) Then,%
\[
\mathbb{Z}/\left(  n_{1}n_{2}\cdots n_{k}\right)  \cong\mathbb{Z}/n_{1}%
\times\mathbb{Z}/n_{2}\times\cdots\times\mathbb{Z}/n_{k}.
\]
More concretely, there is a ring isomorphism%
\[
\mathbb{Z}/\left(  n_{1}n_{2}\cdots n_{k}\right)  \rightarrow\mathbb{Z}%
/n_{1}\times\mathbb{Z}/n_{2}\times\cdots\times\mathbb{Z}/n_{k}%
\]
that does what you expect (i.e., sends each residue class $r+n_{1}n_{2}\cdots
n_{k}\mathbb{Z}$ to the $k$-tuple $\left(  r+n_{1}\mathbb{Z},\ r+n_{2}%
\mathbb{Z},\ \ldots,\ r+n_{k}\mathbb{Z}\right)  $).
\end{theorem}

\begin{proof}
This can be derived from Theorem \ref{thm.CRT-k-ids1} and Theorem
\ref{thm.CRT-k-ids2}, in the same way as we derived Theorem
\ref{thm.CRT-2-ints-rep} from Theorem \ref{thm.CRT-2-ids1} and Theorem
\ref{thm.CRT-2-ids2}. Details are LTTR.
\end{proof}

\begin{corollary}
\label{cor.CRT-primepows}Let $p_{1},p_{2},\ldots,p_{k}$ be $k$ distinct
primes. Let $i_{1},i_{2},\ldots,i_{k}$ be $k$ nonnegative integers. Then,%
\[
\mathbb{Z}/\left(  p_{1}^{i_{1}}p_{2}^{i_{2}}\cdots p_{k}^{i_{k}}\right)
\cong\mathbb{Z}/p_{1}^{i_{1}}\times\mathbb{Z}/p_{2}^{i_{2}}\times\cdots
\times\mathbb{Z}/p_{k}^{i_{k}}.
\]
More concretely, there is a ring isomorphism%
\[
\mathbb{Z}/\left(  p_{1}^{i_{1}}p_{2}^{i_{2}}\cdots p_{k}^{i_{k}}\right)
\rightarrow\mathbb{Z}/p_{1}^{i_{1}}\times\mathbb{Z}/p_{2}^{i_{2}}\times
\cdots\times\mathbb{Z}/p_{k}^{i_{k}}%
\]
that does what you expect.
\end{corollary}

\begin{proof}
The prime powers $p_{1}^{i_{1}},p_{2}^{i_{2}},\ldots,p_{k}^{i_{k}}$ are
mutually coprime; thus, we can apply Theorem \ref{thm.CRT-k-ints} to
$n_{j}=p_{j}^{i_{j}}$.
\end{proof}

Note that it is important that the primes be distinct in Corollary
\ref{cor.CRT-primepows}. For example, $\mathbb{Z}/p^{2}$ is not isomorphic to
$\mathbb{Z}/p\times\mathbb{Z}/p$ (not even as additive groups, let alone as rings).

The Chinese Remainder Theorem (and Corollary \ref{cor.CRT-primepows} in
particular) has many down-to-earth consequences. For example, here is one:

\begin{exercise}
Let $n$ be a positive integer. Let $k$ be the number of distinct prime factors
of $n$. (For instance, if $n=360=2^{3}\cdot3^{2}\cdot5$, then $k=3$.)

Show that the ring $\mathbb{Z}/n$ has exactly $2^{k}$ idempotent elements.
\end{exercise}

Let us next use the Chinese Remainder Theorem to revisit Exercise
\ref{exe.21hw0.7}. That exercise asked you to count how many of the numbers
$0,1,\ldots,n-1$ appear as remainders of a perfect square divided by $n$, when
$n$ is $7$ or $14$. Let us now ask ourselves the same question for an
arbitrary positive integer $n$. It is not hard to see that this question is
equivalent to asking how many elements of the ring $\mathbb{Z}/n$ are squares
in this ring. Here I am using the following terminology:

\begin{definition}
\label{def.ring.square}Let $R$ be a ring. An element $r\in R$ is said to be a
\textbf{square} (in $R$) if there exists some $u\in R$ such that $r=u^{2}$.
\end{definition}

For example, the squares in $\mathbb{R}$ are the nonnegative reals, whereas
the squares in $\mathbb{Z}$ are the perfect squares. For another example, the
squares in $\mathbb{Z}/7\mathbb{Z}$ are the four elements $\overline
{0},\overline{1},\overline{2},\overline{4}$. (Indeed, this is equivalent to
the answer to Exercise \ref{exe.21hw0.7} \textbf{(a)}.)

If $n$ is a positive integer, then an element $i\in\left\{  0,1,\ldots
,n-1\right\}  $ is the remainder of some perfect square divided by $n$ if and
only if the element $\overline{i}=i+n\mathbb{Z}$ is a square in $\mathbb{Z}%
/n$. Thus, counting distinct remainders of perfect squares divided by $n$ is
equivalent to counting squares in $\mathbb{Z}/n$.

Now, I claim that the latter can be done easily when the prime factorization
of $n$ is known. The way to do it is in three steps:

\begin{enumerate}
\item Answer the question (i.e., \textquotedblleft how many squares does
$\mathbb{Z}/n$ have?\textquotedblright) when $n$ is prime.

\item Extend the answer to the case when $n$ is a prime power (i.e., a number
of the form $p^{i}$ with $p$ prime and $i\in\mathbb{N}$).

\item Finally, extend the answer to all positive integers $n$.
\end{enumerate}

This three-step program is a standard strategy for answering
number-theoretical questions. Typically, the three steps each have methods
tailored to them:

\begin{enumerate}
\item When $n$ is prime, the ring $\mathbb{Z}/n$ is a field. This makes many
tactics available that would otherwise not work; e.g., Gaussian elimination
works over fields but not generally over arbitrary rings (we will learn more
about this later).

\item There are many tools for \textquotedblleft lifting\textquotedblright%
\ results about primes to analogous results about prime powers.

\item Here, the Chinese Remainder Theorem becomes useful. Any positive integer
$n$ is a product of finitely many mutually coprime prime powers $p_{1}^{a_{1}%
},p_{2}^{a_{2}},\ldots,p_{k}^{a_{k}}$. Thus, the Chinese Remainder Theorem
(more precisely, Corollary \ref{cor.CRT-primepows}) yields%
\begin{equation}
\mathbb{Z}/n\cong\mathbb{Z}/p_{1}^{a_{1}}\times\mathbb{Z}/p_{2}^{a_{2}}%
\times\cdots\times\mathbb{Z}/p_{k}^{a_{k}}. \label{eq.exa.CRT-k-ints.sqs.1}%
\end{equation}

\end{enumerate}

For our specific question (counting squares in $\mathbb{Z}/n$), Step 1 is
quite easy (see Exercise \ref{exe.21hw1.10} \textbf{(c)} below). (More
precisely, that exercise covers the case when $n$ is odd. But the only even
prime is $2$, and you can count the squares in $\mathbb{Z}/2$ on your
hands\footnote{Not fingers, hands.}.) Step 2 is both trickier and more
laborious (see Exercises \ref{exe.count-squares.Z/pk} and
\ref{exe.count-squares.Z/2k} below). Step 3 is now easy (assuming Steps 1 and
2 are done): If $A_{1},A_{2},\ldots,A_{k}$ are rings, then the squares in the
direct product $A_{1}\times A_{2}\times\cdots\times A_{k}$ are just the
$k$-tuples $\left(  a_{1},a_{2},\ldots,a_{k}\right)  $ where each $a_{i}$ is a
square in $A_{i}$; thus,
\begin{align}
&  \left(  \text{the number of squares in }A_{1}\times A_{2}\times\cdots\times
A_{k}\right) \nonumber\\
&  =\prod_{i=1}^{k}\left(  \text{the number of squares in }A_{i}\right)  .
\label{eq.exa.CRT-k-ints.sqs.2}%
\end{align}
Furthermore, isomorphic rings have the same number of squares (since any ring
morphism sends squares to squares, and thus any ring isomorphism restricts to
a bijection between the squares of its domain and the squares of its target).
Thus, (\ref{eq.exa.CRT-k-ints.sqs.1}) yields%
\begin{align*}
&  \left(  \text{the number of squares in }\mathbb{Z}/n\right) \\
&  =\left(  \text{the number of squares in }\mathbb{Z}/p_{1}^{a_{1}}%
\times\mathbb{Z}/p_{2}^{a_{2}}\times\cdots\times\mathbb{Z}/p_{k}^{a_{k}%
}\right) \\
&  =\prod_{i=1}^{k}\left(  \text{the number of squares in }\mathbb{Z}%
/p_{i}^{a_{i}}\right)  \ \ \ \ \ \ \ \ \ \ \left(  \text{by
(\ref{eq.exa.CRT-k-ints.sqs.2})}\right)  .
\end{align*}

As promised, we shall now compute the number of squares in $\mathbb{Z}/p$ when
$p$ is an odd prime. Better even, let us compute the number of squares in any
finite field $F$ that satisfies $2\cdot1_{F}\neq0_{F}$. (This is a more
general question, since $\mathbb{Z}/p$ is a finite field whenever $p$ is a
prime, and it satisfies $2\cdot1_{F}\neq0_{F}$ whenever $p$ is odd.)

\begin{exercise}
\label{exe.21hw1.10}Let $F$ be a field.

\begin{enumerate}
\item[\textbf{(a)}] Prove that if $a,b\in F$ satisfy $a^{2}=b^{2}$, then $a=b$
or $a=-b$.
\end{enumerate}

From now on, assume that $2\cdot1_{F}\neq0_{F}$ (that is, $1_{F}+1_{F}%
\neq0_{F}$). Note that this is satisfied whenever $F=\mathbb{Z}/p\mathbb{Z}$
for a prime $p>2$ (but also for various other finite fields), but fails when
$F=\mathbb{Z}/2\mathbb{Z}$.

\begin{enumerate}
\item[\textbf{(b)}] Prove that $a\neq-a$ for every nonzero $a\in F$.
\end{enumerate}

From now on, assume that $F$ is finite.

\begin{enumerate}
\item[\textbf{(c)}] Prove that the number of squares in $F$ is $\dfrac{1}%
{2}\left(  \left\vert F\right\vert +1\right)  $.

\item[\textbf{(d)}] Conclude that $\left\vert F\right\vert $ is odd.
\end{enumerate}

[\textbf{Hint:} For part \textbf{(c)}, argue that each nonzero square in $F$
can be written as $\alpha^{2}$ for exactly two distinct elements $\alpha\in F$.]
\end{exercise}

Our next step towards counting squares in $\mathbb{Z}/n$ is the following
exercise (\cite[homework set \#2, Exercise 4]{21w}):

\begin{exercise}
\label{exe.21hw2.4}Let $p$ be a prime number.

\begin{enumerate}
\item[\textbf{(a)}] Prove that if $a$ and $b$ are two integers such that
$a^{2}\equiv b^{2}\operatorname{mod}p^{2}$, then $a\equiv b\operatorname{mod}%
p^{2}$ or $a\equiv-b\operatorname{mod}p^{2}$ or $a\equiv b\equiv
0\operatorname{mod}p$.

\item[\textbf{(b)}] Prove that the number of squares in the ring
$\mathbb{Z}/p^{2}$ is $\dfrac{p^{2}-p}{2}+1$.
\end{enumerate}
\end{exercise}

When the prime number $p$ is distinct from $2$, this can be extended to higher
powers of $p$:

\begin{exercise}
\label{exe.count-squares.Z/pk}Let $p>2$ be a prime number. Let $k$ be a
positive integer.

\begin{enumerate}
\item[\textbf{(a)}] Prove that if $a$ and $b$ are two integers such that
$a^{2}\equiv b^{2}\operatorname{mod}p^{k}$, then $a\equiv b\operatorname{mod}%
p^{k}$ or $a\equiv-b\operatorname{mod}p^{k}$ or $a\equiv b\equiv
0\operatorname{mod}p$.

\item[\textbf{(b)}] Prove that the number of squares in the ring
$\mathbb{Z}/p^{k}$ that are units is $\dfrac{p^{k}-p^{k-1}}{2}$.

\item[\textbf{(c)}] Prove that there is a bijection%
\begin{align*}
&  \text{from the set }\left\{  \text{squares in the ring }\mathbb{Z}%
/p^{k}\text{ that are not units}\right\} \\
&  \text{to the set }\left\{  \text{squares in the ring }\mathbb{Z}%
/p^{k-2}\right\}
\end{align*}
whenever $k\geq2$.

\item[\textbf{(d)}] Prove that the number of all squares in the ring
$\mathbb{Z}/p^{k}$ is
\[
\dfrac{1}{2}%
\begin{cases}
p^{k}-p^{k-1}+p^{k-2}-p^{k-3}+p^{k-4}-p^{k-5}\pm\cdots-p^{1}+2, & \text{if
}k\text{ is even};\\
p^{k}-p^{k-1}+p^{k-2}-p^{k-3}+p^{k-4}-p^{k-5}\pm\cdots-p^{0}+2, & \text{if
}k\text{ is odd.}%
\end{cases}
\]

\end{enumerate}
\end{exercise}

Something similar works when $p$ is $2$, but there are some nuances:

\begin{exercise}
\label{exe.count-squares.Z/2k}Let $k\geq2$ be an integer.

\begin{enumerate}
\item[\textbf{(a)}] Prove that if $a$ and $b$ are two integers such that
$a^{2}\equiv b^{2}\operatorname{mod}2^{k}$, then $a\equiv b\operatorname{mod}%
2^{k-1}$ or $a\equiv-b\operatorname{mod}2^{k-1}$ or $a\equiv b\equiv
0\operatorname{mod}2$.

\item[\textbf{(b)}] Conversely, prove that if $a$ and $b$ are two integers
such that $a\equiv b\operatorname{mod}2^{k-1}$ or $a\equiv-b\operatorname{mod}%
2^{k-1}$, then $a^{2}\equiv b^{2}\operatorname{mod}2^{k}$.

\item[\textbf{(c)}] Prove that the number of squares in the ring
$\mathbb{Z}/2^{k}$ that are units is $2^{k-3}$ if $k\geq3$, and is $1$ otherwise.

\item[\textbf{(d)}] Compute the number of all squares in the ring
$\mathbb{Z}/2^{k}$.
\end{enumerate}
\end{exercise}

The Chinese Remainder Theorem can be used to break rings into smaller ones:

\begin{exercise}
\label{exe.CRT-finrings}\ \ 

\begin{enumerate}
\item[\textbf{(a)}] Let $n_{1},n_{2},\ldots,n_{k}$ be $k$ mutually coprime
integers. (\textquotedblleft Mutually coprime\textquotedblright\ means that
$n_{i}$ is coprime to $n_{j}$ whenever $i<j$.) Let $R$ be any ring. For any
integer $m$, we let $mR$ denote the ideal $\left\{  mr\ \mid\ r\in R\right\}
$ of $R$ (this is easily seen to be an ideal, whether or not $R$ is
commutative), and we let $R/m$ denote the corresponding quotient ring $R/mR$.
Prove that%
\[
R/\left(  n_{1}n_{2}\cdots n_{k}\right)  \cong R/n_{1}\times R/n_{2}%
\times\cdots\times R/n_{k}.
\]

\item[\textbf{(b)}] Let $R$ be a finite ring whose size $\left\vert
R\right\vert $ is a product $p_{1}p_{2}\cdots p_{k}$ of distinct primes
$p_{1},p_{2},\ldots,p_{k}$ (such as $3\cdot5\cdot11$). Prove that $R$ is
commutative, and is isomorphic to $\mathbb{Z}/p_{1}\times\mathbb{Z}%
/p_{2}\times\cdots\times\mathbb{Z}/p_{k}$.

\item[\textbf{(c)}] Let $R$ instead be a finite ring such that there exists no
prime $p$ satisfying $p^{3}\mid\left\vert R\right\vert $. (This is a weaker
assumption than in part \textbf{(b)}.) Prove that $R$ is commutative.
\end{enumerate}
\end{exercise}

\subsubsection{Remark on noncommutative rings}

\begin{fineprint}
Theorem \ref{thm.CRT-k-ids1} becomes false if we drop the assumption that $R$
be commutative. Indeed, even Theorem \ref{thm.CRT-2-ids1} becomes false for
noncommutative $R$, as the following exercise (taken from \cite[Example
4]{vanDal06}) shows:

\begin{exercise}
Let $R$ be any nontrivial ring, and consider the ideals $I,J,K$ of the
upper-triangular matrix ring $R^{2\leq2}$ defined in Exercise
\ref{exe.ideals.triangular-matrices.1} \textbf{(a)}.

\begin{enumerate}
\item[\textbf{(a)}] Prove that $I$ and $J$ are comaximal (i.e., we have
$I+J=R^{2\leq2}$).

\item[\textbf{(b)}] Prove that $I\cap J\neq IJ$.
\end{enumerate}
\end{exercise}

However, we can tweak Theorem \ref{thm.CRT-k-ids1} to make it work for
noncommutative rings $R$ as well:

\begin{theorem}
\label{thm.CRT-k-ids-nc1}Let $I_{1},I_{2},\ldots,I_{k}$ be $k$ mutually
comaximal ideals of a (not necessarily commutative) ring $R$. Let $I_{1}\ast
I_{2}\ast\cdots\ast I_{k}$ denote the sum of all the $k!$ products $J_{1}%
J_{2}\cdots J_{k}$, where $J_{1},J_{2},\ldots,J_{k}$ are the $k$ ideals
$I_{1},I_{2},\ldots,I_{k}$ in some order. (For example, if $k=3$, then
$I_{1}\ast I_{2}\ast I_{3}=I_{1}I_{2}I_{3}+I_{1}I_{3}I_{2}+I_{2}I_{1}%
I_{3}+I_{2}I_{3}I_{1}+I_{3}I_{1}I_{2}+I_{3}I_{2}I_{1}$.)

Then,
\[
I_{1}\cap I_{2}\cap\cdots\cap I_{k}=I_{1}\ast I_{2}\ast\cdots\ast I_{k}.
\]

\end{theorem}
\end{fineprint}

\begin{exercise}
\label{exe.thm.CRT-k-ids-nc}Prove Theorem \ref{thm.CRT-k-ids-nc1}. \medskip

[\textbf{Hint:} The proof is a not-too-difficult adaptation of our above proof
of Theorem \ref{thm.CRT-k-ids1}.]
\end{exercise}

\begin{fineprint}
Theorem \ref{thm.CRT-k-ids-nc1} can be improved even further. Namely, instead
of summing all the $k!$ products $J_{1}J_{2}\cdots J_{k}$, we can sum the two
products $I_{1}I_{2}\cdots I_{k}$ and $I_{k}I_{k-1}\cdots I_{1}$ (that is, we
can replace $I_{1}\ast I_{2}\ast\cdots\ast I_{k}$ by the sum $I_{1}I_{2}\cdots
I_{k}+I_{k}I_{k-1}\cdots I_{1}$), and Theorem \ref{thm.CRT-k-ids-nc1} will
remain true! This fascinating result (and a further generalization) is proved
in Birgit van Dalen's nicely written bachelor thesis \cite{vanDal05}, which is
a good reason to learn Dutch\footnote{See \cite{vanDal06} for a summary in
English.}. As an exercise, we suggest proving its $k=3$ case:
\end{fineprint}

\begin{exercise}
Let $R$ be a ring. Let $I,J,K$ be three mutually comaximal ideals of $R$.
Prove that $I\cap J\cap K=IJK+KJI$.
\end{exercise}

\begin{fineprint}
Here are two similar exercises (the first again courtesy of \cite{vanDal05}):
\end{fineprint}

\begin{exercise}
Let $R$ be a ring. Let $I,J,K$ be three mutually comaximal ideals of $R$.
Prove that $I\cap J\cap K=IJK+JKI+KIJ$.
\end{exercise}

\begin{exercise}
Let $R$ be a ring. Let $I,J,K$ be three mutually comaximal ideals of $R$.
Prove that $IJ+JK+KI=R$.
\end{exercise}

\begin{fineprint}
The next exercise shows that comaximality of ideals is passed on from two
ideals to their powers:
\end{fineprint}

\begin{exercise}
Let $R$ be a ring. Let $I$ and $J$ be two comaximal ideals of $R$. Let $n$ be
a positive integer.

\begin{enumerate}
\item[\textbf{(a)}] Prove that the ideals $I$ and $J^{n}$ are comaximal as
well. (Here, $J^{n}$ means $\underbrace{JJ\cdots J}_{n\text{ times}}$, where
we refer to Definition \ref{def.rings.ideal-arith.arith} \textbf{(b)} for the
definition of the product of two ideals.)

\item[\textbf{(b)}] Let $m$ be a further positive integer. Prove that the
ideals $I^{m}$ and $J^{n}$ are comaximal as well.

\item[\textbf{(c)}] By applying this to $R=\mathbb{Z}$, prove that if $a$ and
$b$ are two coprime integers, then their powers $a^{m}$ and $b^{n}$ are also
coprime whenever $m$ and $n$ are two positive integers.
\end{enumerate}
\end{exercise}

\subsection{\label{sec.rings.euclid}Euclidean rings and Euclidean domains
(\cite[\S 8.1]{DumFoo04})}

\subsubsection{All ideals of $\mathbb{Z}$ are principal}

We have talked about ideals of $\mathbb{Z}$ a lot (they give rise to modular
arithmetic), but you might have noticed that all of them were principal. This
is no accident:

\begin{proposition}
\label{prop.eucldom.Z-PID}Any ideal of $\mathbb{Z}$ is principal.
\end{proposition}

\begin{proof}
Let $I$ be an ideal of $\mathbb{Z}$. We must show that $I$ is principal.

If $I=\left\{  0\right\}  $, then this is clear (since $I=0\mathbb{Z}$ in this
case). So we WLOG assume that $I\neq\left\{  0\right\}  $. Since $I$ always
contains $0$, this means that $I$ must contain a nonzero integer as well.
Hence, $I$ contains a positive integer (because if $I$ contains a negative
integer $a$, then $I$ must also contain $\left(  -1\right)  a$, which is
positive). Let $b\in I$ be the \textbf{smallest} positive integer that $I$
contains. Hence, $I$ cannot contain any positive integer smaller than $b$.
However, $I$ contains $b$, and thus contains every multiple of $b$ (since $I$
is an ideal). In other words, $b\mathbb{Z}\subseteq I$.

We will now show that $I\subseteq b\mathbb{Z}$. Indeed, let $a\in I$. Let $r$
be the remainder of $a$ divided by $b$. Then, $r\in\left\{  0,1,\ldots
,b-1\right\}  $ and $r\equiv a\operatorname{mod}b$. Now, from $r\equiv
a\operatorname{mod}b$, we obtain $b\mid r-a$ and thus $r-a\in b\mathbb{Z}%
\subseteq I$. Hence, $r=\underbrace{r-a}_{\in I}+\underbrace{a}_{\in I}\in
I+I=I$ (since $I$ is an ideal of $\mathbb{Z}$). Hence, $r$ cannot be a
positive integer smaller than $b$ (since $I$ cannot contain any positive
integer smaller than $b$). In other words, $r\notin\left\{  1,2,\ldots
,b-1\right\}  $. Contrasting this with $r\in\left\{  0,1,\ldots,b-1\right\}
$, we obtain $r=0$. Thus, $b\mid\underbrace{r}_{=0}-\,a=0-a=-a\mid a$, so that
$a\in b\mathbb{Z}$.

Forget that we fixed $a$. We thus have shown that $a\in b\mathbb{Z}$ for each
$a\in I$. In other words, $I\subseteq b\mathbb{Z}$. Combined with
$b\mathbb{Z}\subseteq I$, this yields $I=b\mathbb{Z}$. Thus, $I$ is principal, qed.
\end{proof}

The key to making this proof work was clearly the concept of division with
remainder. Not every ring has this feature. However, many rings different from
$\mathbb{Z}$ have it; thus, it is worth defining a word for them:

\subsubsection{\label{subsec.rings.euclid.euclid}Euclidean rings and Euclidean
domains}

\begin{definition}
\label{def.rings.euclid}Let $R$ be a commutative ring.

\begin{enumerate}
\item[\textbf{(a)}] A \textbf{norm} on $R$ means a function $N:R\rightarrow
\mathbb{N}$ with $N\left(  0\right)  =0$.

\item[\textbf{(b)}] A norm $N$ on $R$ is said to be \textbf{Euclidean} if for
any $a\in R$ and any nonzero $b\in R$, there exist elements $q,r\in R$ with%
\[
a=qb+r\ \ \ \ \ \ \ \ \ \ \text{and}\ \ \ \ \ \ \ \ \ \ \left(  r=0\text{ or
}N\left(  r\right)  <N\left(  b\right)  \right)  .
\]

\item[\textbf{(c)}] We say that $R$ is a \textbf{Euclidean ring} if $R$ has a
Euclidean norm.

\item[\textbf{(d)}] We say that $R$ is a \textbf{Euclidean domain} if $R$ is a
Euclidean ring and is an integral domain.
\end{enumerate}
\end{definition}

You can think of the norm as a measure of the \textquotedblleft
size\textquotedblright\ of an element of $R$, similar to the absolute value of
an integer or to the degree of a polynomial. (These will indeed be particular
cases.) Note that we are \textbf{not} requiring that the norm have any nice
algebraic properties (such as $N\left(  ab\right)  =N\left(  a\right)
N\left(  b\right)  $, which will be true for some Euclidean norms but not for
others). We are also \textbf{not} requiring the $q$ and the $r$ in Definition
\ref{def.rings.euclid} \textbf{(b)} to be unique. If your familiarity with
norms comes from real analysis, be warned that the concept we have defined
here has nothing in common with the one you know except for the name.

Some examples will help illustrate the definition:

\begin{itemize}
\item Any field $F$ is a Euclidean domain. Indeed, any map $N:F\rightarrow
\mathbb{N}$ with $N\left(  0\right)  =0$ is a Euclidean norm on $F$. (To see
that it satisfies the condition of Definition \ref{def.rings.euclid}
\textbf{(b)}, just set $q=\dfrac{a}{b}$ and $r=0$.)

\item The ring $\mathbb{Z}$ is a Euclidean domain. Indeed, the map%
\begin{align*}
N:\mathbb{Z}  &  \rightarrow\mathbb{N},\\
a  &  \mapsto\left\vert a\right\vert
\end{align*}
is a Euclidean norm on $\mathbb{Z}$. The fact that it is Euclidean follows
from division with remainder\footnote{In more detail: We need to show that for
any $a\in\mathbb{Z}$ and any nonzero $b\in\mathbb{Z}$, there exist elements
$q,r\in\mathbb{Z}$ with%
\[
a=qb+r\ \ \ \ \ \ \ \ \ \ \text{and}\ \ \ \ \ \ \ \ \ \ \left(  r=0\text{ or
}\left\vert r\right\vert <\left\vert b\right\vert \right)  .
\]
To find these $q,r$, we divide $a$ by $\left\vert b\right\vert $ with
remainder. Let $q_{0}$ and $r_{0}$ be the quotient and the remainder that we
obtain. If $b$ is positive, we can then take $q=q_{0}$ and $r=r_{0}$. If $b$
is negative, then we instead take $q=-q_{0}$ and $r=r_{0}$ (because
$b=-\left\vert b\right\vert $).}. However, the $q$ and the $r$ in Definition
\ref{def.rings.euclid} \textbf{(b)} are not unique! For $a=7$ and $b=5$, there
are \textbf{two} pairs $\left(  q,r\right)  \in\mathbb{Z}\times\mathbb{Z}$
with%
\[
a=qb+r\ \ \ \ \ \ \ \ \ \ \text{and}\ \ \ \ \ \ \ \ \ \ \left(  r=0\text{ or
}N\left(  r\right)  <N\left(  b\right)  \right)  .
\]
These two pairs are $\left(  1,2\right)  $ and $\left(  2,-3\right)  $. The
second pair has negative $r$, which is why it does not qualify as a
quotient-remainder pair in the sense of high school arithmetic; but this $r$
nevertheless qualifies for the definition of a Euclidean norm.

\item If $F$ is a field, then the ring $F\left[  x\right]  $ of univariate
polynomials over $F$ is a Euclidean domain. We will discuss this later in more
detail, when we study polynomials. However, polynomial rings in more than $1$
variable are not Euclidean; neither are polynomial rings over nontrivial non-fields.

\item The ring $\mathbb{Z}\left[  i\right]  $ of Gaussian integers is a
Euclidean domain. Indeed, we claim that the map%
\begin{align*}
N:\mathbb{Z}\left[  i\right]   &  \rightarrow\mathbb{N},\\
a+bi  &  \mapsto a^{2}+b^{2}\ \ \ \ \ \ \ \ \ \ \left(  \text{for all }%
a,b\in\mathbb{Z}\right)
\end{align*}
is a Euclidean norm.

To prove this, we must show that for any $\alpha\in\mathbb{Z}\left[  i\right]
$ and any nonzero $\beta\in\mathbb{Z}\left[  i\right]  $, there exist elements
$q,r\in\mathbb{Z}\left[  i\right]  $ with%
\begin{equation}
\alpha=q\beta+r\ \ \ \ \ \ \ \ \ \ \text{and}\ \ \ \ \ \ \ \ \ \ \left(
r=0\text{ or }N\left(  r\right)  <N\left(  \beta\right)  \right)  .
\label{eq.exa.eucdom.Zi.qr1}%
\end{equation}
So let us fix an $\alpha\in\mathbb{Z}\left[  i\right]  $ and a nonzero
$\beta\in\mathbb{Z}\left[  i\right]  $. We are looking for elements
$q,r\in\mathbb{Z}\left[  i\right]  $ that satisfy (\ref{eq.exa.eucdom.Zi.qr1}%
). We can even replace the \textquotedblleft$r=0$ or $N\left(  r\right)
<N\left(  \beta\right)  $\textquotedblright\ condition in
(\ref{eq.exa.eucdom.Zi.qr1}) by the stronger condition \textquotedblleft%
$N\left(  r\right)  <N\left(  \beta\right)  $\textquotedblright.

To find the elements $q,r$ we are seeking, we make the following observation:
The absolute value $\left\vert z\right\vert $ of a complex number $z=a+bi$
(with $a,b\in\mathbb{R}$) is defined as $\left\vert z\right\vert =\sqrt
{a^{2}+b^{2}}$, whereas the norm $N\left(  z\right)  $ of a Gaussian integer
$z=a+bi$ (with $a,b\in\mathbb{Z}$) is defined as $N\left(  z\right)
=a^{2}+b^{2}$. Thus, any $z\in\mathbb{Z}\left[  i\right]  $ satisfies
$N\left(  z\right)  =\left\vert z\right\vert ^{2}$. Hence, we have the
following chain of equivalences:%
\begin{align}
\left(  N\left(  r\right)  <N\left(  \beta\right)  \right)  \  &
\Longleftrightarrow\ \left(  \left\vert r\right\vert ^{2}<\left\vert
\beta\right\vert ^{2}\right)  \ \Longleftrightarrow\ \left(  \left\vert
r\right\vert <\left\vert \beta\right\vert \right)  \ \Longleftrightarrow
\ \left(  \dfrac{\left\vert r\right\vert }{\left\vert \beta\right\vert
}<1\right) \nonumber\\
&  \Longleftrightarrow\ \left(  \left\vert \dfrac{r}{\beta}\right\vert
<1\right)  \label{eq.exa.eucdom.Zi.qr2}%
\end{align}
(since $\dfrac{\left\vert z\right\vert }{\left\vert w\right\vert }=\left\vert
\dfrac{z}{w}\right\vert $ for any two complex numbers $z$ and $w\neq0$).
Moreover, we have the equivalence%
\begin{equation}
\left(  \alpha=q\beta+r\right)  \ \Longleftrightarrow\ \left(  \dfrac{\alpha
}{\beta}=q+\dfrac{r}{\beta}\right)  \ \Longleftrightarrow\ \left(
\dfrac{\alpha}{\beta}-q=\dfrac{r}{\beta}\right)  .
\label{eq.exa.eucdom.Zi.qr3}%
\end{equation}
Now, recall that we are looking for elements $q,r\in\mathbb{Z}\left[
i\right]  $ that satisfy $\alpha=q\beta+r$ and $N\left(  r\right)  <N\left(
\beta\right)  $. In view of (\ref{eq.exa.eucdom.Zi.qr2}) and
(\ref{eq.exa.eucdom.Zi.qr3}), this means that we are looking for elements
$q,r\in\mathbb{Z}\left[  i\right]  $ that satisfy $\dfrac{\alpha}{\beta
}-q=\dfrac{r}{\beta}$ and $\left\vert \dfrac{r}{\beta}\right\vert <1$.
Equivalently, we can look for a Gaussian integer $q\in\mathbb{Z}\left[
i\right]  $ satisfying $\left\vert \dfrac{\alpha}{\beta}-q\right\vert <1$
(because once such a $q$ has been found, we can set $r=\alpha-q\beta$ and
obtain $\dfrac{r}{\beta}=\dfrac{\alpha-q\beta}{\beta}=\dfrac{\alpha}{\beta}%
-q$, so that $\dfrac{\alpha}{\beta}-q=\dfrac{r}{\beta}$ and $\left\vert
\dfrac{r}{\beta}\right\vert =\left\vert \dfrac{\alpha}{\beta}-q\right\vert
<1$). But finding such a $q$ is easy if you remember the geometric meaning of
the Gaussian integers: The Gaussian integers are the lattice points of a
square lattice in the plane:%
\[%
\begin{tikzpicture}[baseline, scale=2.5]
\draw[step=1cm,gray,very thin] (-2.4,-2.4) grid (2.4,2.4);
\draw[very thick,red,->] (-2.5,0) -- (2.5,0);
\draw[very thick,red,->] (0,-2.5) -- (0,2.5);
\foreach\x/\xtext in {-2, -1, 0, 1, 2}
\fill(\x cm, 0) circle (1pt) node[above left=2pt]{$\xtext$};
\foreach\x/\xtext in {-2, -1, 1, 2}
\fill(\x cm, 1 cm) circle (1pt) node[above left=2pt]{$\xtext+ i$};
\fill(0, 1 cm) circle (1pt) node[above left=2pt] {$i$};
\foreach\x/\xtext in {-2, -1, 1, 2}
\fill(\x cm, -1 cm) circle (1pt) node[above left=2pt]{$\xtext- i$};
\fill(0, -1 cm) circle (1pt) node[above left=2pt] {$-i$};
\foreach\x/\xtext in {-2, -1, 1, 2}
\fill(\x cm, 2 cm) circle (1pt) node[above left=2pt]{$\xtext+ 2i$};
\fill(0, 2 cm) circle (1pt) node[above left=2pt] {$2i$};
\foreach\x/\xtext in {-2, -1, 1, 2}
\fill(\x cm, -2 cm) circle (1pt) node[above left=2pt]{$\xtext- 2i$};
\fill(0, -2 cm) circle (1pt) node[above left=2pt] {$-2i$};
\end{tikzpicture}%
\]
(imagine the lattice being extended to infinity in all four directions). So a
Gaussian integer $q\in\mathbb{Z}\left[  i\right]  $ satisfying $\left\vert
\dfrac{\alpha}{\beta}-q\right\vert <1$ simply means a lattice point at a
distance\footnote{The \textbf{distance} between two complex numbers $x$ and
$y$ is defined to be the real number $\left\vert x-y\right\vert $.} less than
$1$ from the point $\dfrac{\alpha}{\beta}$. Geometrically, it is easy to see
that such a lattice point exists (since the point $\dfrac{\alpha}{\beta}$ must
lie in one of the squares of the lattice, and then have distance $\leq
\dfrac{\sqrt{2}}{2}$ from one of the four vertices of the square\footnote{Here
is a close-up picture of the square (with one possible location of
$\dfrac{\alpha}{\beta}$):%
\[%
\begin{tikzpicture}[scale=2.5]
\useasboundingbox(0.7, 0.7) grid (2.3, 2.3);
\draw[step=1cm,gray,very thin] (0.8,0.8) grid (2.2,2.2);
\fill(1, 1) circle (1pt) node[above left=2pt]{$1+i$};
\fill(2, 1) circle (1pt) node[above left=2pt]{$2+i$};
\fill(1, 2) circle (1pt) node[above left=2pt]{$1+2i$};
\fill(2, 2) circle (1pt) node[above left=2pt]{$2+2i$};
\fill[blue] (1.3, 1.6) circle (1pt) node[above right=2pt]{$\alpha/ \beta$};
\end{tikzpicture}%
\]
\par
I am claiming that the point $\dfrac{\alpha}{\beta}$ has distance $\leq
\dfrac{\sqrt{2}}{2}$ from one of the four vertices of the square in which it
lies. The easiest way to see this geometrically is to draw circles of radius
$\dfrac{\sqrt{2}}{2}$ around the vertices of the square, and convince yourself
that these circles cover the entire square:%
\[%
\begin{tikzpicture}[scale=2.5]
\useasboundingbox(0.2, 0.2) grid (2.8, 2.8);
\draw[step=1cm,gray,very thin] (0.8,0.8) grid (2.2,2.2);
\fill(1, 1) circle (1pt) node[above left=2pt]{$1+i$};
\fill(2, 1) circle (1pt) node[above left=2pt]{$2+i$};
\fill(1, 2) circle (1pt) node[above left=2pt]{$1+2i$};
\fill(2, 2) circle (1pt) node[above left=2pt]{$2+2i$};
\fill[blue] (1.3, 1.6) circle (1pt) node[above right=2pt]{$\alpha/ \beta$};
\draw[green] (1, 1) circle [radius=0.7071];
\draw[green] (2, 1) circle [radius=0.7071];
\draw[green] (1, 2) circle [radius=0.7071];
\draw[green] (2, 2) circle [radius=0.7071];
\end{tikzpicture}%
\ \ .
\]
}; but this entails that $\dfrac{\alpha}{\beta}$ has distance $<1$ from this
latter vertex\footnote{since $\dfrac{\sqrt{2}}{2}<1$}). This can also be
proved algebraically\footnote{\textit{Proof.} Write the point $\dfrac{\alpha
}{\beta}$ as $x+yi$, where $x$ and $y$ are real numbers. Each real number $z$
has distance $\leq\dfrac{1}{2}$ from the nearest integer (which is either
$\left\lfloor z\right\rfloor $ or $\left\lceil z\right\rceil $). Thus, $x$ has
distance $\leq\dfrac{1}{2}$ from some integer $n$, and likewise $y$ has
distance $\leq\dfrac{1}{2}$ from some integer $m$. Consider these $n$ and $m$.
Then, $\left\vert x-n\right\vert \leq\dfrac{1}{2}$ and $\left\vert
y-m\right\vert \leq\dfrac{1}{2}$. Since $n$ and $m$ are integers, we have
$n+mi\in\mathbb{Z}\left[  i\right]  $, so that $n+mi$ is a lattice point.
However, $\dfrac{\alpha}{\beta}=x+yi$, so that%
\[
\dfrac{\alpha}{\beta}-\left(  n+mi\right)  =\left(  x+yi\right)  -\left(
n+mi\right)  =\left(  x-n\right)  +\left(  y-m\right)  i.
\]
Hence,%
\begin{align*}
\left\vert \dfrac{\alpha}{\beta}-\left(  n+mi\right)  \right\vert  &
=\left\vert \left(  x-n\right)  +\left(  y-m\right)  i\right\vert \\
&  =\sqrt{\left(  x-n\right)  ^{2}+\left(  y-m\right)  ^{2}}%
\ \ \ \ \ \ \ \ \ \ \left(
\begin{array}
[c]{c}%
\text{by the definition of the absolute}\\
\text{value of a complex number,}\\
\text{since }x-n\text{ and }y-m\text{ are reals}%
\end{array}
\right) \\
&  =\sqrt{\left\vert x-n\right\vert ^{2}+\left\vert y-m\right\vert ^{2}}\\
&  \leq\sqrt{\left(  \dfrac{1}{2}\right)  ^{2}+\left(  \dfrac{1}{2}\right)
^{2}}\ \ \ \ \ \ \ \ \ \ \left(  \text{since }\left\vert x-n\right\vert
\leq\dfrac{1}{2}\text{ and }\left\vert y-m\right\vert \leq\dfrac{1}{2}\right)
\\
&  =\sqrt{\dfrac{1}{2}}=\dfrac{\sqrt{2}}{2}<1.
\end{align*}
Thus, the lattice point $n+mi$ has a distance of $<1$ from the point
$\dfrac{\alpha}{\beta}$.}. Thus, we have found $q$.

(A slightly restated version of this proof can be found in the proof of
Theorem 3.1 in Keith Conrad's \textit{The Gaussian integers} (see
\url{https://kconrad.math.uconn.edu/math5230f12/handouts/Zinotes.pdf} ).)

\item The ring%
\[
\mathbb{Z}\left[  \sqrt{-2}\right]  :=\left\{  a+b\sqrt{-2}\ \mid
\ a,b\in\mathbb{Z}\right\}
\]
(this is another subring of $\mathbb{C}$, since $\sqrt{-2}=\sqrt{2}i$) is
Euclidean, too. (See Exercise \ref{exe.eucldom.Zsqrt-2} for a proof.)

\item The ring%
\[
\mathbb{Z}\left[  \sqrt{-3}\right]  :=\left\{  a+b\sqrt{-3}\ \mid
\ a,b\in\mathbb{Z}\right\}
\]
(this is another subring of $\mathbb{C}$, since $\sqrt{-3}=\sqrt{3}i$) is
\textbf{not} Euclidean. (For a proof, see
\url{https://math.stackexchange.com/questions/115934} or Exercise
\ref{exe.21hw2.7} below.)

However, there is a slightly larger ring that is Euclidean: namely, the
so-called \textbf{ring of Eisenstein integers}, defined as%
\[
\mathbb{Z}\left[  \omega\right]  :=\left\{  a+b\omega\ \mid\ a,b\in
\mathbb{Z}\right\}
\]
for $\omega=\dfrac{-1+\sqrt{-3}}{2}$. (See Exercise
\ref{exe.eucldom.eisenstein} below for the proof that this ring is Euclidean.)

\item The ring
\[
\mathbb{Z}\left[  \sqrt{2}\right]  :=\left\{  a+b\sqrt{2}\ \mid\ a,b\in
\mathbb{Z}\right\}
\]
(this is a subring of $\mathbb{R}$) is Euclidean. A Euclidean norm for it is
the map%
\begin{align*}
\mathbb{Z}\left[  \sqrt{2}\right]   &  \rightarrow\mathbb{N},\\
a+b\sqrt{2}  &  \mapsto\left\vert a^{2}-2b^{2}\right\vert
\ \ \ \ \ \ \ \ \ \ \left(  \text{for }a,b\in\mathbb{Z}\right)  .
\end{align*}
(See Exercise \ref{exe.eucldom.Zsqrt2} below for a proof.)

\item The ring
\[
\mathbb{Z}\left[  \sqrt{14}\right]  :=\left\{  a+b\sqrt{14}\ \mid
\ a,b\in\mathbb{Z}\right\}
\]
is Euclidean. A Euclidean norm for it is notoriously hard to construct (in
particular, it is \textbf{not} the map sending each $a+b\sqrt{14}$ to
$\left\vert a^{2}-14b^{2}\right\vert $). See
\url{https://math.stackexchange.com/questions/1148364} .

\item The ring $\mathbb{Z}\left[  \sqrt{5}\right]  :=\left\{  a+b\sqrt
{5}\ \mid\ a,b\in\mathbb{Z}\right\}  $ is \textbf{not} Euclidean.

\item For each $n\in\mathbb{Z}$, the ring $\mathbb{Z}/n$ is Euclidean (but is
not a Euclidean domain in most cases). A Euclidean norm $N$ on this ring is
easy to construct (e.g., for $n>0$, we can define $N\left(  \overline
{a}\right)  $ to be the smallest nonnegative integer in the residue class
$\overline{a}$).
\end{itemize}

Thus, we have now seen multiple examples and non-examples of Euclidean rings
and Euclidean domains. Now, we claim that all Euclidean domains have a
property that we have previously proved for $\mathbb{Z}$:

\begin{proposition}
\label{prop.eucldom.PID}Let $R$ be a Euclidean ring. Then, any ideal of $R$ is principal.
\end{proposition}

\begin{proof}
The same argument we used for proving Proposition \ref{prop.eucldom.Z-PID} can
easily be adapted to prove Proposition \ref{prop.eucldom.PID}. The main change
is that you now need to take a nonzero $b\in I$ with smallest possible
$N\left(  b\right)  $. (Here, $N$ is a fixed Euclidean norm on $R$.) For
details, see \cite[\S 8.1, proof of Proposition 1]{DumFoo04}.
\end{proof}

\begin{remark}
Euclidean domains are much more well-studied than Euclidean rings. Some
authors go as far as using the word \textquotedblleft Euclidean
ring\textquotedblright\ as a synonym for \textquotedblleft Euclidean
domain\textquotedblright\ (which, of course, conflicts with our definition of
the former).
\end{remark}

See \url{https://kconrad.math.uconn.edu/blurbs/ringtheory/euclideanrk.pdf} for
more about Euclidean domains.

\begin{exercise}
\label{exe.eucldom.Zsqrt-2}Prove that the ring
\[
\mathbb{Z}\left[  \sqrt{-2}\right]  :=\left\{  a+b\sqrt{-2}\ \mid
\ a,b\in\mathbb{Z}\right\}
\]
is Euclidean, and that the map%
\begin{align*}
N:\mathbb{Z}\left[  \sqrt{-2}\right]   &  \rightarrow\mathbb{N},\\
a+b\sqrt{-2}  &  \mapsto a^{2}+2b^{2}\ \ \ \ \ \ \ \ \ \ \left(  \text{for
}a,b\in\mathbb{Z}\right)
\end{align*}
is a Euclidean norm for it. \medskip

[\textbf{Hint:} Imitate the above proof for $\mathbb{Z}\left[  i\right]  $.]
\end{exercise}

\begin{exercise}
\label{exe.eucldom.Zsqrt2}Prove that the ring
\[
\mathbb{Z}\left[  \sqrt{2}\right]  :=\left\{  a+b\sqrt{2}\ \mid\ a,b\in
\mathbb{Z}\right\}
\]
is Euclidean, and that the map%
\begin{align*}
N:\mathbb{Z}\left[  \sqrt{2}\right]   &  \rightarrow\mathbb{N},\\
a+b\sqrt{2}  &  \mapsto\left\vert a^{2}-2b^{2}\right\vert
\ \ \ \ \ \ \ \ \ \ \left(  \text{for }a,b\in\mathbb{Z}\right)
\end{align*}
is a Euclidean norm for it. \medskip

[\textbf{Hint:} First, prove that the latter map $N$ is multiplicative --
i.e., that it satisfies $N\left(  xy\right)  =N\left(  x\right)  \cdot
N\left(  y\right)  $ for all $x,y\in\mathbb{Z}\left[  \sqrt{2}\right]  $.
Moreover, this still holds if we extend $N$ to a map from $\mathbb{Q}\left[
\sqrt{2}\right]  $ to $\mathbb{Q}$ (defined by the same formula).]
\end{exercise}

\begin{exercise}
\label{exe.eucldom.eisenstein}Let $\omega$ denote the complex number
$\dfrac{-1+\sqrt{-3}}{2}\in\mathbb{C}$.

\begin{enumerate}
\item[\textbf{(a)}] Prove that $\omega^{3}=1$ and $\omega^{2}+\omega+1=0$.

\item[\textbf{(b)}] Prove that $\left\vert a+b\omega\right\vert =\sqrt
{a^{2}-ab+b^{2}}$ for any $a,b\in\mathbb{R}$.

\item[\textbf{(c)}] Define a subset $\mathbb{Z}\left[  \omega\right]  $ of
$\mathbb{C}$ by%
\[
\mathbb{Z}\left[  \omega\right]  :=\left\{  a+b\omega\ \mid\ a,b\in
\mathbb{Z}\right\}  .
\]
Prove that $\mathbb{Z}\left[  \omega\right]  $ is a subring of $\mathbb{C}$.
(It is called the ring of \textbf{Eisenstein integers}.)

\item[\textbf{(d)}] Prove that $\mathbb{Z}\left[  \sqrt{-3}\right]  $ is a
subring of $\mathbb{Z}\left[  \omega\right]  $.

\item[\textbf{(e)}] Prove that the ring $\mathbb{Z}\left[  \omega\right]  $ is
Euclidean, and that the map%
\begin{align*}
N:\mathbb{Z}\left[  \omega\right]   &  \rightarrow\mathbb{N},\\
a+b\omega &  \mapsto a^{2}-ab+b^{2}\ \ \ \ \ \ \ \ \ \ \left(  \text{for
}a,b\in\mathbb{Z}\right)
\end{align*}
is a Euclidean norm for it.

\item[\textbf{(f)}] Find all units of the ring $\mathbb{Z}\left[
\omega\right]  $.

\item[\textbf{(g)}] Show that an element $a+b\omega$ of $\mathbb{Z}\left[
\omega\right]  $ (with $a,b\in\mathbb{Z}$) belongs to $\mathbb{Z}\left[
\sqrt{-3}\right]  $ if and only if $b$ is even.

\item[\textbf{(h)}] Consider $\mathbb{Z}\left[  \omega\right]  $ as an
additive group, and $\mathbb{Z}\left[  \sqrt{-3}\right]  $ as a subgroup of
$\mathbb{Z}\left[  \omega\right]  $. Prove that this subgroup $\mathbb{Z}%
\left[  \sqrt{-3}\right]  $ has index $2$ in $\mathbb{Z}\left[  \omega\right]
$ (that is, the quotient group $\mathbb{Z}\left[  \omega\right]
/\mathbb{Z}\left[  \sqrt{-3}\right]  $ has size $2$).

\item[\textbf{(i)}] If $z\in\mathbb{Z}\left[  \omega\right]  $ is any
Eisenstein integer, then at least one of the three numbers $z,\ z\omega
,\ z\omega^{2}$ belongs to $\mathbb{Z}\left[  \sqrt{-3}\right]  $.
\end{enumerate}

[\textbf{Hint:} Geometrically speaking, the three complex numbers
$1,\ \omega,\ \omega^{2}$ are the vertices of an equilateral triangle
inscribed in the unit circle. The elements of $\mathbb{Z}\left[
\omega\right]  $ are the grid points of a triangular lattice that looks as
follows (imagine the picture extended to infinity on all sides):%
\[%
\pgfmathsetmacro{\rows}{5}
\begin{tikzpicture}[baseline, scale=2.5]
\clip(-2.5, -2.5) rectangle (2.5, 2.5);
\foreach\row in {-\rows, ..., \rows} {
\draw[gray, very thin] ($(-\rows,0)+\row*(0.5, {0.5*sqrt(3)})$) -- ($(\rows
,0)+\row*(-0.5, {0.5*sqrt(3)})$);
\draw[gray, very thin] ($({-\rows/2},{-\rows/2*sqrt(3)})+\row
*(0.5,{-0.5*sqrt(3)})$) -- ($(\rows/2,{\rows/2*sqrt(3)})+\row
*(0.5,{-0.5*sqrt(3)})$);
\draw[gray, very thin] ($(\rows/2,{-\rows/2*sqrt(3)})+\row
*(-0.5,{-0.5*sqrt(3)})$) -- ($({-\rows/2},{\rows/2*sqrt(3)})+\row
*(-0.5,{-0.5*sqrt(3)})$);
}
\draw[green!50!black] (0, 0) circle [radius=1];
\fill[red] (0, 0) circle (1pt) node[above left=2pt] {$0$};
\fill[red] (1, 0) circle (1pt) node[above left=2pt] {$1$};
\fill[red] (2, 0) circle (1pt);
\fill[red] (-1, 0) circle (1pt) node[below left=2pt] {$-1$};
\fill[red] (-2, 0) circle (1pt);
\fill[red] (0, {sqrt(3)}) circle (1pt) node[above left=2pt] {$\sqrt{-3}$};
\fill(-1/2, {sqrt(3)/2}) circle (1pt) node[above left=2pt] {$\omega$};
\fill(-1/2, {-sqrt(3)/2}) circle (1pt) node[below left=2pt] {$\omega^2$};
\fill(3/2, {sqrt(3)/2}) circle (1pt) node[above right=2pt] {$\omega+2$};
\fill(5/2, {sqrt(3)/2}) circle (1pt);
\fill(-3/2, {sqrt(3)/2}) circle (1pt);
\fill(-5/2, {sqrt(3)/2}) circle (1pt);
\fill(-5/2, {-sqrt(3)/2}) circle (1pt);
\fill(-3/2, {-sqrt(3)/2}) circle (1pt);
\fill(3/2, {-sqrt(3)/2}) circle (1pt);
\fill(5/2, {-sqrt(3)/2}) circle (1pt);
\fill[red] (-2, {-sqrt(3)}) circle (1pt);
\fill[red] (-1, {-sqrt(3)}) circle (1pt);
\fill[red] (0, {-sqrt(3)}) circle (1pt);
\fill[red] (1, {-sqrt(3)}) circle (1pt);
\fill[red] (2, {-sqrt(3)}) circle (1pt);
\fill[red] (-2, {sqrt(3)}) circle (1pt);
\fill[red] (-1, {sqrt(3)}) circle (1pt);
\fill[red] (1, {sqrt(3)}) circle (1pt);
\fill[red] (2, {sqrt(3)}) circle (1pt);
\fill(1/2, {-sqrt(3)/2}) circle (1pt) node[below=2pt] {$-\omega$};
\fill(1/2, {sqrt(3)/2}) circle (1pt) node[below=2pt] {$-\omega^2$};
\end{tikzpicture}%
\]
(where the red points are the ones that belong to $\mathbb{Z}\left[  \sqrt
{-3}\right]  $).]
\end{exercise}

\begin{exercise}
\label{exe.Rm.norm}Fix an integer $m$. Consider the ring $R_{m}$ defined in
Exercise \ref{exe.21hw1.1}.

Prove that $R_{m}$ is a Euclidean domain. More concretely:

\begin{enumerate}
\item[\textbf{(a)}] For every nonzero $r\in R_{m}$, we let $\widetilde{r}$ be
the smallest positive numerator of $r$. (A \textquotedblleft
numerator\textquotedblright\ of a rational number $r$ means an integer of the
form $dr$ with $d\in\mathbb{Z}$. In other words, if we write $r$ as a ratio of
two integers, then the numerator of this fraction is called a
\textquotedblleft numerator\textquotedblright\ of $r$. For example, $7$ is a
numerator of $\dfrac{7}{9}$, but so are $14$ and $21$ and $-14$ and so on.)

Prove that $\widetilde{r}$ exists.

\item[\textbf{(b)}] Prove that the map%
\begin{align*}
N:R_{m}  &  \rightarrow\mathbb{N},\\
r  &  \mapsto\widetilde{r}\ \ \ \ \ \ \ \ \ \ \text{for }r\neq0,\\
0  &  \mapsto0
\end{align*}
is a Euclidean norm on $R_{m}$.
\end{enumerate}
\end{exercise}

\begin{exercise}
Let $N_{2}:\mathbb{Z}\rightarrow\mathbb{N}$ be the map that sends $0$ to $0$,
while sending each nonzero integer $n$ to $\left\lfloor \log_{2}\left\vert
n\right\vert \right\rfloor $. (Recall that $\left\lfloor x\right\rfloor $
denotes the \textbf{floor} of a real number $x$ -- that is, the largest
integer that is $\leq x$. Thus, $N_{2}\left(  7\right)  =\left\lfloor \log
_{2}\left\vert 7\right\vert \right\rfloor =\left\lfloor 2.807\ldots
\right\rfloor =2$.)

Prove that $N_{2}$ is a Euclidean norm on $\mathbb{Z}$.
\end{exercise}

\begin{exercise}
Prove that the ring $\mathbb{Z}\times\mathbb{Z}$ is not Euclidean. \medskip

[\textbf{Hint:} Assume the contrary, so that $\mathbb{Z}\times\mathbb{Z}$ has
a Euclidean norm $N$. Show that for every $k\in\mathbb{N}$, there are only
finitely many $r\in\mathbb{Z}\times\mathbb{Z}$ satisfying $N\left(  r\right)
=k$. For a contradiction, divide elements by $\left(  1,0\right)  $ and
$\left(  0,1\right)  $.]
\end{exercise}

\begin{exercise}
Let $R$ be a Euclidean ring. Let $I$ be an ideal of $R$. Show that the
quotient ring $R/I$ is again Euclidean. \medskip

[\textbf{Hint:} Let $N$ be a Euclidean norm on $R$. Define a norm
$\overline{N}$ on $R/I$ by setting $\overline{N}\left(  x\right)
=\min\left\{  N\left(  a\right)  \ \mid\ a\in x\right\}  $ for all residue
classes $x\in R/I$.]
\end{exercise}

\begin{exercise}
Prove that if we replace the condition \textquotedblleft$r=0$ or $N\left(
r\right)  <N\left(  b\right)  $\textquotedblright\ by \textquotedblleft%
$N\left(  r\right)  <N\left(  b\right)  $\textquotedblright\ in Definition
\ref{def.rings.euclid}, then the resulting notion of a Euclidean ring will be
equivalent to ours (even though a given Euclidean norm $N$ might no longer
qualify as a Euclidean norm).
\end{exercise}

\subsubsection{\label{subsec.rings.euclid.ea}The (extended) Euclidean
algorithm}

Imagine that you are given some ideal $I$ of $\mathbb{Z}$. Proposition
\ref{prop.eucldom.Z-PID} then guarantees that this ideal $I$ is principal,
i.e., has the form $I=c\mathbb{Z}$ for a single integer $c$. Now, suppose you
want to actually find this $c$.

Our above proof of Proposition \ref{prop.eucldom.Z-PID} is of some help here:
It guarantees that $I=c\mathbb{Z}$, where $c$ is the smallest positive integer
contained in $I$ (or $0$ if no such integer exists)\footnote{This $c$ was
denoted by $b$ in our proof of Proposition \ref{prop.eucldom.Z-PID}.}.
Depending on how much you know about $I$, this can make $c$ easy to find. But
this is not automatic. Indeed, in some cases (e.g., when $I$ is given as the
set of all integers satisfying some complicated uncomputable condition),
finding an integer $c\in\mathbb{Z}$ satisfying $I=c\mathbb{Z}$ is even
algorithmically impossible\footnote{Appreciators of theoretical computer
science will easily concoct an ideal $I$ of $\mathbb{Z}$ such that finding an
integer $c$ satisfying $I=c\mathbb{Z}$ is tantamount to solving
\href{https://en.wikipedia.org/wiki/Halting_problem}{the halting problem}
(which is known to be algorithmically unsolvable).}, even though such a $c$
exists for theoretical reasons.

However, if the ideal $I$ is defined in a sufficiently simple way, then this
problem might be algorithmically solvable. A particularly well-behaved example
is when $I$ is given in the form $I=a\mathbb{Z}+b\mathbb{Z}$ for two
explicitly provided integers $a$ and $b$. In this case, finding an integer $c$
satisfying $I=c\mathbb{Z}$ amounts to finding the greatest common divisor
$\gcd\left(  a,b\right)  $ of $a$ and $b$ (because Proposition
\ref{prop.rings.ideal-arith.Z} \textbf{(c)} yields that $a\mathbb{Z}%
+b\mathbb{Z}=\gcd\left(  a,b\right)  \mathbb{Z}$). The famous Euclidean
algorithm computes this $\gcd\left(  a,b\right)  $, thus solving the problem
of finding $c$. Furthermore, a variant of this algorithm -- known as the
\textbf{extended Euclidean algorithm} -- computes two integers $x$ and $y$
satisfying $\gcd\left(  a,b\right)  =xa+yb$. These integers $x$ and $y$
\textquotedblleft make the equality $a\mathbb{Z}+b\mathbb{Z}=\gcd\left(
a,b\right)  \mathbb{Z}$ explicit\textquotedblright\ (in the sense that they
allow us to actually express an element of $\gcd\left(  a,b\right)
\mathbb{Z}$ in the form \textquotedblleft$a$ times an integer plus $b$ times
an integer\textquotedblright, rather than merely guaranteeing that such an
expression exists).

These two Euclidean algorithms (the usual one and the extended one) can be
found in any textbook on elementary number theory (see, e.g., \cite[Algorithm
2.3.7]{Stein09} for the extended Euclidean algorithm; the usual one can easily
be obtained from it). Here, however, we are interested not in the integers but
in their various generalizations. For what other commutative rings $R$ do such
algorithms (expressing an ideal of the form $aR+bR$ as a principal ideal $cR$)
exist?\footnote{Let me stress that the $aR$, $bR$ and $cR$ here are ideals.
The \textquotedblleft$+$\textquotedblright\ sign in \textquotedblleft%
$aR+bR$\textquotedblright\ stands for a sum of ideals. Thus, the equality
$aR+bR=cR$ has nothing to do with the equality $a+b=c$. (Indeed, the former
equality neither implies nor follows from the latter; the ideals $\left(
a+b\right)  R$ and $aR+bR$ are not the same.)
\par
Translated into the language of elements, the equality $aR+bR=cR$ says that
the elements that can be written as a multiple of $a$ plus a multiple of $b$
are precisely the multiples of $c$.}

Let us make this question more precise: We want an algorithm which, if you
input two elements $a,b\in R$, outputs an element $c\in R$ that satisfies
$aR+bR=cR$. Ideally, this algorithm should also provide \textquotedblleft
evidence\textquotedblright\ for this equality $aR+bR=cR$, that is, a way to
express every element of $cR$ as an element of $aR+bR$ and vice versa. In
order to express every element of $cR$ as an element of $aR+bR$, it suffices
to write $c$ as a sum $xa+yb$ with $x,y\in R$. In order to express every
element of $aR+bR$ as an element of $cR$, it suffices to write $a$ in the form
$a=cu$ for some $u\in R$, and to write $b$ in the form $b=cv$ for some $v\in
R$. So we want our algorithm to output not only $c$ but also $x$, $y$, $u$ and
$v$. In other words, we want it to output the $5$-tuple $\left(
x,y,c,u,v\right)  $ (there is nothing special about the order in which I
listed its entries; I just picked it to put the most important output, $c$, in
the middle).

Let me give this $5$-tuple a name:

\begin{definition}
\label{def.rings.euclid.bezout}Let $a$ and $b$ be two elements of a
commutative ring $R$. Then, a \textbf{Bezout }$5$\textbf{-tuple} for $\left(
a,b\right)  $ shall mean%
\begin{align*}
&  \text{a }5\text{-tuple }\left(  x,y,c,u,v\right)  \in R^{5}\text{ that
satisfies}\\
&  xa+yb=c\text{ and }a=cu\text{ and }b=cv.
\end{align*}

\end{definition}

\begin{example}
Let $R=\mathbb{Z}$ and $a=10$ and $b=6$. Then, $\left(  -1,2,2,5,3\right)  $
is a Bezout $5$-tuple for $\left(  a,b\right)  $, since it satisfies $\left(
-1\right)  a+2b=\left(  -1\right)  10+2\cdot6=2=c$ and $a=10=2\cdot5=cu$ and
$b=6=2\cdot3=cv$. Another Bezout $5$-tuple for $\left(  a,b\right)  $ is
$\left(  1,-2,-2,-5,-3\right)  $. Yet another is $\left(  2,-3,2,5,3\right)
$. There are infinitely many Bezout $5$-tuples $\left(  x,y,c,u,v\right)  $
for $\left(  a,b\right)  $, since we can always replace $x$ and $y$ by $x+3$
and $y-5$ without changing $xa+yb$.
\end{example}

For a general commutative ring $R$, a Bezout $5$-tuple for $\left(
a,b\right)  $ will not always exist. But when it does, it answers all our
questions about the ideal $aR+bR$:

\begin{proposition}
\label{prop.rings.euclid.ea.bezout}Let $R$ be a commutative ring. Let $a,b\in
R$ be arbitrary, and let $\left(  x,y,c,u,v\right)  $ be a Bezout $5$-tuple
for $\left(  a,b\right)  $. Then:

\begin{enumerate}
\item[\textbf{(a)}] We have $aR+bR=cR$.

\item[\textbf{(b)}] Any element $ap+bq$ of $aR+bR$ can be explicitly expressed
as an element of $cR$ by rewriting it in the form $c\left(  up+vq\right)  $.

\item[\textbf{(c)}] Any element $cr$ of $cR$ can be explicitly expressed as an
element of $aR+bR$ by rewriting it as $axr+byr$.
\end{enumerate}
\end{proposition}

\begin{proof}
We know that $\left(  x,y,c,u,v\right)  $ is a Bezout $5$-tuple for $\left(
a,b\right)  $. Thus, the definition of a Bezout $5$-tuple yields that
$xa+yb=c$ and $a=cu$ and $b=cv$.

Now, any element of $aR+bR$ has the form $ap+bq$ for some $p,q\in R$, and
therefore belongs to $cR$, since%
\[
\underbrace{a}_{=cu}p+\underbrace{b}_{=cv}q=cup+cvq=c\underbrace{\left(
up+vq\right)  }_{\in R}\in cR.
\]
This shows that $aR+bR\subseteq cR$.

On the other hand, any element of $cR$ has the form $cr$ for some $r\in R$,
and therefore belongs to $aR+bR$, since%
\[
\underbrace{c}_{=xa+yb}r=\left(  xa+yb\right)  r=a\underbrace{xr}_{\in
R}+b\underbrace{yr}_{\in R}\in aR+bR.
\]
This shows that $cR\subseteq aR+bR$. Combining this with $aR+bR\subseteq cR$,
we obtain $aR+bR=cR$. Therefore, part \textbf{(a)} is proved.

Part \textbf{(b)} was proved above (in the process of proving $aR+bR\subseteq
cR$), and part \textbf{(c)} was proved as well (in the process of showing that
$cR\subseteq aR+bR$). Thus, Proposition \ref{prop.rings.euclid.ea.bezout} is
fully proved.
\end{proof}

Thus, our problem about $aR+bR$ is now reduced to the following: Given two
elements $a$ and $b$ of a commutative ring $R$, when can we algorithmically
find a Bezout $5$-tuple for $\left(  a,b\right)  $ ?

For most rings $R$, the answer is \textquotedblleft no\textquotedblright%
\ already because such a Bezout $5$-tuple doesn't always exist. However, the
answer is \textquotedblleft yes\textquotedblright\ when $R$ is a Euclidean
ring. To be more precise, it is \textquotedblleft yes\textquotedblright\ when
$R$ is a Euclidean ring satisfying certain computability requirements:

\begin{theorem}
\label{thm.rings.euclid.ea.main}Let $R$ be a Euclidean ring.

\begin{enumerate}
\item[\textbf{(a)}] Then, for any two elements $a,b\in R$, there exists a
Bezout $5$-tuple for $\left(  a,b\right)  $.

\item[\textbf{(b)}] Moreover, there exists an algorithm that computes a Bezout
$5$-tuple for any pair $\left(  a,b\right)  \in R^{2}$, provided that the
Euclideanness of $R$ itself is algorithmic\footnotemark.
\end{enumerate}
\end{theorem}

\footnotetext{By \textquotedblleft the Euclideanness of $R$ itself is
algorithmic\textquotedblright, we mean the following:
\par
\begin{itemize}
\item There are algorithms for adding, subtracting and multiplying arbitrary
elements of $R$.{}
\par
\item There is an algorithm for checking whether two given elements of $R$ are
equal.
\par
\item There is a Euclidean norm $N:R\rightarrow\mathbb{N}$ that can be
computed by an algorithm (i.e., there is an algorithm that computes $N\left(
a\right)  $ for each $a\in R$).
\par
\item The $q$ and the $r$ in Definition \ref{def.rings.euclid} \textbf{(b)}
can be computed by an algorithm (i.e., there exists an algorithm that, if you
input an element $a\in R$ and a nonzero element $b\in R$, will output a pair
$\left(  q,r\right)  \in R^{2}$ such that $a=qb+r$ and $\left(  r=0\text{ or
}N\left(  r\right)  <N\left(  b\right)  \right)  $).
\end{itemize}
\par
Actually, the computability of the norm $N$ is not even necessary for our
algorithm.}

\begin{proof}
\textbf{(a)} Since $R$ is Euclidean, there exists a Euclidean norm
$N:R\rightarrow\mathbb{N}$. Consider this norm $N$.

We shall prove Theorem \ref{thm.rings.euclid.ea.main} \textbf{(a)} by strong
induction on the nonnegative integer $N\left(  b\right)  $.

So let $n\in\mathbb{N}$. As the induction hypothesis, we assume that Theorem
\ref{thm.rings.euclid.ea.main} \textbf{(a)} is true for any pair $\left(
a,b\right)  \in R^{2}$ satisfying $N\left(  b\right)  <n$. We must now prove
that Theorem \ref{thm.rings.euclid.ea.main} \textbf{(a)} holds for any pair
$\left(  a,b\right)  \in R^{2}$ satisfying $N\left(  b\right)  =n$.

So let $\left(  a,b\right)  \in R^{2}$ be a pair satisfying $N\left(
b\right)  =n$. Our goal is to prove that Theorem
\ref{thm.rings.euclid.ea.main} \textbf{(a)} holds for this pair, i.e., to
prove that there exists a Bezout $5$-tuple for $\left(  a,b\right)  $.

To construct this $5$-tuple, we distinguish between two cases:

\begin{itemize}
\item \textit{Case 1:} Assume that $b=0$. Then, we set $\left(
x,y,c,u,v\right)  :=\left(  1,0,a,1,0\right)  $. It is easy to see that this
is a Bezout $5$-tuple for $\left(  a,b\right)  $ (since $1a+0b=a$ and
$a=a\cdot1$ and $b=0=a\cdot0$). Hence, we have found a Bezout $5$-tuple for
$\left(  a,b\right)  $ in Case 1.

\item \textit{Case 2:} Assume that $b\neq0$. Since $N$ is a Euclidean norm,
there exist elements $q,r\in R$ with $a=qb+r$ and $\left(  r=0\text{ or
}N\left(  r\right)  <N\left(  b\right)  \right)  $ (by Definition
\ref{def.rings.euclid} \textbf{(b)}). Consider these elements $q,r$. From
$a=qb+r$, we obtain $a-qb=r$. We have $r=0$ or $N\left(  r\right)  <N\left(
b\right)  $; thus, we can break this case into two subcases:

\begin{itemize}
\item \textit{Subcase 2.1:} Assume that $r=0$. Then, $a=qb+\underbrace{r}%
_{=0}=qb=bq$. Hence, $\left(  0,1,b,q,1\right)  $ is a Bezout $5$-tuple for
$\left(  a,b\right)  $ (since $0a+1b=b$ and $a=bq$ and $b=b\cdot1$).

\item \textit{Subcase 2.2:} Assume that $N\left(  r\right)  <N\left(
b\right)  $. Thus, by the induction hypothesis, Theorem
\ref{thm.rings.euclid.ea.main} \textbf{(a)} is true for the pair $\left(
b,r\right)  $ instead of $\left(  a,b\right)  $. In other words,
\[
\text{there exists a Bezout }5\text{-tuple }\left(  x,y,c,u,v\right)  \text{
for this pair }\left(  b,r\right)  .
\]
Consider this $5$-tuple. By the definition of a Bezout $5$-tuple, it satisfies
$xb+yr=c$ and $b=cu$ and $r=cv$.

Now,
\[
\left(  y,\ x-qy,\ c,\ qu+v,\ u\right)  \text{ is a Bezout }5\text{-tuple for
the pair }\left(  a,b\right)  ,
\]
because%
\[
ya+\left(  x-qy\right)  b=ya+xb-qyb=xb+y\underbrace{\left(  a-qb\right)
}_{=r}=xb+yr=c
\]
and%
\[
a=q\underbrace{b}_{=cu}+\underbrace{r}_{=cv}=qcu+cv=c\left(  qu+v\right)
\]
and%
\[
b=cu.
\]

\end{itemize}

Thus, we have found a Bezout $5$-tuple for $\left(  a,b\right)  $ in Case 2
(since we have found such a tuple in each of the two subcases).
\end{itemize}

We have now found a Bezout $5$-tuple for the pair $\left(  a,b\right)  $ in
each of the above two cases. Hence, such a $5$-tuple always exists. This
completes the induction step, and thus Theorem \ref{thm.rings.euclid.ea.main}
\textbf{(a)} is proved. \medskip

\textbf{(b)} Our inductive proof of Theorem \ref{thm.rings.euclid.ea.main}
\textbf{(a)} above gives a recursive algorithm for finding a Bezout $5$-tuple
for the pair $\left(  a,b\right)  $. Indeed, depending on the case and the
subcase, it either solves this problem directly (in Case 1 and in Subcase
2.1), or reduces it to the analogous problem for a different pair $\left(
a,b\right)  $ with a smaller value of $N\left(  b\right)  $ (in Subcase 2.2).
Thus, by a (finite) sequence of such reductions, we eventually arrive at a
pair for which we can directly find a Bezout $5$-tuple, and thus we can do so
for the original pair as well.

Here is this algorithm, spelled out as code (in Python):\footnote{The
\texttt{quo\_rem} function, applied to a pair $\left(  a,b\right)  $, is
assumed to return a pair $\left(  q,r\right)  $ of elements $q,r\in R$ that
satisfy $a=qb+r$ and $\left(  r=0\text{ or }N\left(  r\right)  <N\left(
b\right)  \right)  $. Such a pair exists, since $N$ is a Euclidean
norm.}\medskip

\texttt{\qquad\qquad def bezout\_5tuple(a, b):}

\texttt{\qquad\qquad\qquad\qquad if b == 0: \# Case 1}

\texttt{\qquad\qquad\qquad\qquad\qquad\qquad return (1, 0, a, 1, 0)}

\texttt{\qquad\qquad\qquad\qquad q, r = quo\_rem(a, b)}

\texttt{\qquad\qquad\qquad\qquad if r == 0: \# Subcase 2.1}

\texttt{\qquad\qquad\qquad\qquad\qquad\qquad return (0, 1, b, q, 1)}

\texttt{\qquad\qquad\qquad\qquad\# Subcase 2.2}

\texttt{\qquad\qquad\qquad\qquad(x, y, c, u, v) = bezout\_5tuple(b, r)}

\texttt{\qquad\qquad\qquad\qquad return (y, x - q*y, c, q*u + v, u)}\medskip

This function returns a Bezout $5$-tuple for $\left(  a,b\right)  $.
\end{proof}

Let us give an example for the algorithm that comes out of our above proof of
Theorem \ref{thm.rings.euclid.ea.main} \textbf{(b)}. For simplicity, we use a
very familiar ring -- namely, $\mathbb{Z}$ -- and a very familiar Euclidean norm:

\begin{itemize}
\item Let $R$ be the Euclidean ring $\mathbb{Z}$, and let $N$ be the Euclidean
norm on $\mathbb{Z}$ that sends each integer $a$ to $\left\vert a\right\vert
$. We aim to compute a Bezout $5$-tuple for the pair $\left(  6,14\right)  $.

The existence of such a $5$-tuple is guaranteed by Theorem
\ref{thm.rings.euclid.ea.main} \textbf{(a)} (applied to $a=6$ and $b=14$).
Thus, in order to find such a $5$-tuple, we inspect the proof of Theorem
\ref{thm.rings.euclid.ea.main} \textbf{(a)} we gave above, in the specific
situation when $a=6$ and $b=14$. This situation is an instance of Case 2 in
the above proof (since $b\neq0$), so the first step is to find elements
$q,r\in R$ with $a=qb+r$ and $\left(  r=0\text{ or }N\left(  r\right)
<N\left(  b\right)  \right)  $. Such elements $q,r$ are easily found by
standard division with remainder; we obtain $q=0$ and $r=6$ (since
$a=6=\underbrace{0}_{=q}\cdot\underbrace{14}_{=b}+\underbrace{6}_{=r}$ and
$N\left(  r\right)  =6<14=N\left(  b\right)  $). Note that there is a
different choice as well\footnote{namely, $q=1$ and $r=-8$}, but we pick this one.

Thus, we are in Subcase 2.2 (since $r\neq0$). To continue, we need to find a
Bezout $5$-tuple for the pair $\left(  b,r\right)  =\left(  14,6\right)  $.
How do we find such a $5$-tuple?

Again, we inspect the above proof of Theorem \ref{thm.rings.euclid.ea.main}
\textbf{(a)}, but now for $a=14$ and $b=6$. This situation is again an
instance of Case 2 (since $b\neq0$), so the first step is to find elements
$q,r\in R$ with $a=qb+r$ and $\left(  r=0\text{ or }N\left(  r\right)
<N\left(  b\right)  \right)  $. Such elements $q,r$ are easily found by
standard division with remainder; we obtain $q=2$ and $r=2$ (since $a=2b+2$
and $N\left(  2\right)  <N\left(  b\right)  $). Thus, we are again in Subcase
2.2 (since $r\neq0$). To continue, we need to find a Bezout $5$-tuple for the
pair $\left(  b,r\right)  =\left(  6,2\right)  $. How do we find such a $5$-tuple?

Again, we inspect the above proof of Theorem \ref{thm.rings.euclid.ea.main}
\textbf{(a)}, but now for $a=6$ and $b=2$. This situation is again an instance
of Case 2 (since $b\neq0$), so the first step is to find elements $q,r\in R$
with $a=qb+r$ and $\left(  r=0\text{ or }N\left(  r\right)  <N\left(
b\right)  \right)  $. Such elements $q,r$ are easily found by standard
division with remainder; we obtain $q=3$ and $r=0$. Thus, we are in Subcase
2.1 (since $r=0$), and we conclude that $\left(  0,1,b,q,1\right)  =\left(
0,1,2,3,1\right)  $ is a Bezout $5$-tuple for $\left(  a,b\right)  =\left(
6,2\right)  $.

Having found this Bezout $5$-tuple for $\left(  6,2\right)  $, we can now go
back one step and obtain a Bezout $5$-tuple for $\left(  14,6\right)  $:
Namely, as we learned in Subcase 2.2,%
\begin{align*}
&  \text{if }\left(  x,y,c,u,v\right)  \text{ is a Bezout }5\text{-tuple for
the pair }\left(  b,r\right)  \text{ (where }a=qb+r\text{),}\\
&  \text{then }\left(  y,\ x-qy,\ c,\ qu+v,\ u\right)  \text{ is a Bezout
}5\text{-tuple for the pair }\left(  a,b\right)  \text{.}%
\end{align*}
Thus, knowing that $\left(  0,1,2,3,1\right)  $ is a Bezout $5$-tuple for the
pair $\left(  6,2\right)  $, we conclude (with $a=14$ and $b=6$ and $q=2$ and
$r=2$) that%
\[
\left(  1,\ 0-2\cdot1,\ 2,\ 2\cdot3+1,\ 3\right)  =\left(
1,\ -2,\ 2,\ 7,\ 3\right)
\]
is a Bezout $5$-tuple for the pair $\left(  14,6\right)  $.

Having found this Bezout $5$-tuple for $\left(  14,6\right)  $, we can now go
back one more step and obtain a Bezout $5$-tuple for $\left(  6,14\right)  $:
Namely, as we learned in Subcase 2.2,%
\begin{align*}
&  \text{if }\left(  x,y,c,u,v\right)  \text{ is a Bezout }5\text{-tuple for
the pair }\left(  b,r\right)  \text{ (where }a=qb+r\text{),}\\
&  \text{then }\left(  y,\ x-qy,\ c,\ qu+v,\ u\right)  \text{ is a Bezout
}5\text{-tuple for the pair }\left(  a,b\right)  \text{.}%
\end{align*}
Thus, knowing that $\left(  1,\ -2,\ 2,\ 7,\ 3\right)  $ is a Bezout $5$-tuple
for the pair $\left(  14,6\right)  $, we conclude (with $a=6$ and $b=14$ and
$q=0$ and $r=6$) that%
\[
\left(  -2,\ 1-0\cdot\left(  -2\right)  ,\ 2,\ 0\cdot7+3,\ 7\right)  =\left(
-2,\ 1,\ 2,\ 3,\ 7\right)
\]
is a Bezout $5$-tuple for the pair $\left(  6,14\right)  $. This is precisely
what we were looking for. (The reader can easily verify that this is really a
Bezout $5$-tuple for $\left(  6,14\right)  $, but our proof makes this
verification redundant.)
\end{itemize}

Note that the computation we just showed has a distinctive \textquotedblleft
there-and-back-again\textquotedblright\ pattern: In its first half, we have
been reducing our problem (to find a Bezout $5$-tuple for $\left(
6,14\right)  $) step-by-step to a simpler version of it (to find a Bezout
$5$-tuple for $\left(  6,2\right)  $), which we then solved directly; then, in
its second half, we have been retracing our steps backwards to transform the
solution of the simpler version into a solution of the original problem. This
is typical for recursive algorithms. In our specific case (when $R=\mathbb{Z}%
$), the first half of our computation is just the familiar Euclidean algorithm
for computing the greatest common divisor of two integers $a$ and $b$. Had we
been only looking for this greatest common divisor (which in our case was
$2$), we could have stopped in the middle. However, to find the entire Bezout
$5$-tuple for the original pair $\left(  a,b\right)  $, we had to retrace all
our steps backwards. Thus, the algorithm we used to find a Bezout $5$-tuple is
known as the \textbf{(generalized) extended Euclidean algorithm}. Euclidean
rings owe their name to this very algorithm.

\subsection{\label{sec.rings.pid}Principal ideal domains (\cite[\S 8.1 and
\S 8.2]{DumFoo04})}

\subsubsection{Principal ideal domains}

Proposition \ref{prop.eucldom.PID} is so useful that its conclusion (viz.,
that any ideal of $R$ is principal) has been given its own name:

\begin{definition}
\label{def.PID.def}An integral domain $R$ is said to be a \textbf{principal
ideal domain} (for short, \textbf{PID}) if each ideal of $R$ is principal.
\end{definition}

Thus, Proposition \ref{prop.eucldom.PID} yields the following:

\begin{proposition}
\label{prop.eucldom.PID2}Any Euclidean domain is a PID.
\end{proposition}

The converse is not true, although counterexamples are hard to find. One of
the simplest is the ring $\mathbb{Z}\left[  \alpha\right]  =\left\{
a+b\alpha\ \mid\ a,b\in\mathbb{Z}\right\}  $, where $\alpha=\dfrac
{1+\sqrt{-19}}{2}$ (a complex number). (See \cite[page 282]{DumFoo04} for a
proof that this ring is a PID but not a Euclidean domain.)

\subsubsection{Divisibility in commutative rings}

Much of the basic theory of commutative rings can be viewed as a project to
generalize the classical arithmetic of the integers to wider classes of
\textquotedblleft numbers\textquotedblright. As part of this project, we shall
now define gcds and lcms in commutative rings. Our definition will be stated
for arbitrary commutative rings, but we will soon see that they behave
particularly well when the ring is a PID (which is why we are only doing this
definition now). \footnote{The notions of \textquotedblleft greatest common
divisor\textquotedblright\ and \textquotedblleft lowest common
multiple\textquotedblright\ that we will now introduce are not literal
generalizations of the corresponding notions from classical arithmetic. See
below for the exact relation.}

\begin{definition}
\label{def.gcd-lcm}Let $R$ be a commutative ring.

Let $a\in R$.

\begin{enumerate}
\item[\textbf{(a)}] A \textbf{multiple} of $a$ means an element of the form
$ac$ with $c\in R$. In other words, it means an element of the principal ideal
$aR$.

\item[\textbf{(b)}] A \textbf{divisor} of $a$ means an element $d\in R$ such
that $a$ is a multiple of $d$ (that is, $a\in dR$). We write \textquotedblleft%
$d\mid a$\textquotedblright\ for \textquotedblleft$d$ is a divisor of
$a$\textquotedblright.
\end{enumerate}

Now, let $a\in R$ and $b\in R$.

\begin{enumerate}
\item[\textbf{(c)}] A \textbf{common divisor} of $a$ and $b$ means an element
of $R$ that is a divisor of $a$ and a divisor of $b$ at the same time.

\item[\textbf{(d)}] A \textbf{common multiple} of $a$ and $b$ means an element
of $R$ that is a multiple of $a$ and a multiple of $b$ at the same time.

\item[\textbf{(e)}] A \textbf{greatest common divisor} (short: \textbf{gcd})
of $a$ and $b$ means a common divisor $d$ of $a$ and $b$ such that
\textbf{every} common divisor of $a$ and $b$ is a divisor of $d$.

\item[\textbf{(f)}] A \textbf{lowest common multiple} (also known as
\textbf{least common multiple}; short: \textbf{lcm}) of $a$ and $b$ means a
common multiple $m$ of $a$ and $b$ such that \textbf{every} common multiple of
$a$ and $b$ is a multiple of $m$.
\end{enumerate}
\end{definition}

The concepts of \textquotedblleft multiple\textquotedblright\ and
\textquotedblleft divisor\textquotedblright\ we just introduced are
straightforward generalizations of the corresponding concepts from
arithmetic\footnote{Here I am assuming that you are using the
\textquotedblleft right\textquotedblright\ definitions of the latter concepts.
For example, every integer (including $0$ itself) is a divisor of $0$. Some
authors dislike this and prefer to explicitly require $0$ to not divide $0$;
in that case, of course, my definition does not agree with theirs.}. (You
recover the latter concepts if you set $R=\mathbb{Z}$.) The notions of
\textquotedblleft gcd\textquotedblright\ and \textquotedblleft
lcm\textquotedblright\ are a bit subtler: If $a$ and $b$ are two integers,
then their greatest common divisor $\gcd\left(  a,b\right)  $ in the sense of
classical arithmetic is a gcd of $a$ and $b$ in the sense of Definition
\ref{def.gcd-lcm} \textbf{(e)}; however, so is $-\gcd\left(  a,b\right)  $. So
our new notion of a gcd is slightly more liberal than the classical notion, in
the sense that it allows for negative gcds. The same holds for lcms. Thus,
gcds and lcms in our sense are not literally unique. This is one reason why we
said \textquotedblleft a gcd\textquotedblright\ and \textquotedblleft an
lcm\textquotedblright\ (rather than \textquotedblleft the
gcd\textquotedblright\ and \textquotedblleft the lcm\textquotedblright) in
Definition \ref{def.gcd-lcm}. Another reason is that $a$ and $b$ might not
have any gcd to begin with. (We will later see some examples where this happens.)

Let us first state some basic properties of divisibility:

\begin{proposition}
Let $R$ be a commutative ring. Then:

\begin{enumerate}
\item[\textbf{(a)}] We have $a\mid a$ for each $a\in R$.

\item[\textbf{(b)}] If $a,b,c\in R$ satisfy $a\mid b$ and $b\mid c$, then
$a\mid c$.

\item[\textbf{(c)}] If $a,b,c\in R$ satisfy $a\mid b$ and $a\mid c$, then
$a\mid b+c$.

\item[\textbf{(d)}] If $a,b,c,d\in R$ satisfy $a\mid b$ and $c\mid d$, then
$ac\mid bd$.
\end{enumerate}
\end{proposition}

\begin{proof}
Easy (and analogous to the classical proofs for $R=\mathbb{Z}$).
\end{proof}

\subsubsection{Gcds and lcms}

Before we explore gcds and lcms in arbitrary commutative rings, let us record
the precise relation between them and the classical arithmetic notions:

\begin{proposition}
\label{prop.gcd-lcm.ZZ}Let $a$ and $b$ be two integers. Let $g=\gcd\left(
a,b\right)  $ and $\ell=\operatorname{lcm}\left(  a,b\right)  $, where we are
using the classical arithmetic definitions of $\gcd$ and $\operatorname{lcm}$. Then:

\begin{enumerate}
\item[\textbf{(a)}] The gcds of $a$ and $b$ (in the sense of Definition
\ref{def.gcd-lcm} \textbf{(e)}) are $g$ and $-g$.

\item[\textbf{(b)}] The lcms of $a$ and $b$ (in the sense of Definition
\ref{def.gcd-lcm} \textbf{(f)}) are $\ell$ and $-\ell$.
\end{enumerate}
\end{proposition}

\begin{proof}
\textbf{(a)} It is known from classical arithmetic that $g$ is a common
divisor of $a$ and $b$, and that every common divisor of $a$ and $b$ is a
divisor of $g$. In other words, $g$ is a gcd of $a$ and $b$ in the sense of
Definition \ref{def.gcd-lcm} \textbf{(e)}. It is easy to see that this
property is inherited by $-g$ as well (since divisibilities don't change when
we replace $g$ by $-g$). Thus, both numbers $g$ and $-g$ are gcds of $a$ and
$b$ in the sense of Definition \ref{def.gcd-lcm} \textbf{(e)}. It remains to
show that they are the only gcds of $a$ and $b$ in this sense.

So let $u$ be a gcd of $a$ and $b$ in the sense of Definition
\ref{def.gcd-lcm} \textbf{(e)}. We must show that $u\in\left\{  g,-g\right\}
$.

From the way we introduced $u$, we know that $u$ is a common divisor of $a$
and $b$, and that every common divisor of $a$ and $b$ is a divisor of $u$. The
first of these two facts yields that $u\mid g$ (since any common divisor of
$a$ and $b$ is a divisor of $g$); the second yields that $g\mid u$ (since $g$
is a common divisor of $a$ and $b$, and thus is a divisor of $u$). Combining
$u\mid g$ and $g\mid u$, we find $u=\pm g$. In other words, $u\in\left\{
g,-g\right\}  $. This finishes our proof of part \textbf{(a)}. \medskip

\textbf{(b)} The proof is similar to that for part \textbf{(a)}.
\end{proof}

Now, what about gcds and lcms in other rings? The existence of a gcd is far
from God-given, as the following example shows:

\begin{example}
Let $R$ be the ring%
\[
\mathbb{Z}\left[  \sqrt{-3}\right]  =\left\{  a+b\sqrt{-3}\ \mid
\ a,b\in\mathbb{Z}\right\}  .
\]
Let $a=4$ and $b=2\left(  1+\sqrt{-3}\right)  $. Then, $a$ and $b$ have no gcd
in $R$; nor do they have an lcm in $R$. You will prove this in Exercise
\ref{exe.21hw2.7}.
\end{example}

\subsubsection{Associate elements}

Uniqueness of gcds and lcms is a simpler question: They are rarely unique on
the nose, but they are always unique up to multiplication by a unit when the
ring is an integral domain. Before we show this, let me introduce a word for this:

\begin{definition}
\label{def.ring.associate}Let $R$ be a commutative ring. Let $a,b\in R$. We
say that $a$ is \textbf{associate} to $b$ in $R$ (and we write $a\sim b$) if
there exists a unit $u$ of $R$ such that $a=bu$.

Instead of saying \textquotedblleft$a$ is associate to $b$\textquotedblright,
we shall also say that \textquotedblleft$a$ and $b$ are
associate\textquotedblright. (This is justified by the fact -- which we will
prove in Proposition \ref{prop.ring.associate.eqrel} -- that $\sim$ is an
equivalence relation.)
\end{definition}

For example:

\begin{itemize}
\item Two integers $a$ and $b$ are associate in $\mathbb{Z}$ if and only if
$a=\pm b$ (that is, if and only if $a=b$ or $a=-b$).

\item Any two nonzero elements $a$ and $b$ of a field are associate in that
field (since $\dfrac{a}{b}$ is a unit and satisfies $a=b\cdot\dfrac{a}{b}$).
The element $0$ is associate only to itself.

\item Let $F$ be a field. In the polynomial ring $F\left[  x\right]  $, any
nonzero polynomial $f\in F\left[  x\right]  $ is associate to a monic
polynomial (since its leading coefficient is a unit, and dividing $f$ by this
coefficient results in a monic polynomial).

\item It is not hard to prove that the only units of the ring $\mathbb{Z}%
\left[  i\right]  $ are the four Gaussian integers $1,i,-1,-i$. Thus, two
Gaussian integers $\alpha$ and $\beta$ in $\mathbb{Z}\left[  i\right]  $ are
associate if and only if $\alpha\in\left\{  \beta,i\beta,-\beta,-i\beta
\right\}  $.
\end{itemize}

A general property of associateness is the following:

\begin{proposition}
\label{prop.ring.associate.eqrel}Let $R$ be a commutative ring. The relation
$\sim$ is an equivalence relation.
\end{proposition}

\begin{proof}
This is fairly straightforward. We need to show that the relation $\sim$ is
reflexive, symmetric and transitive.

\textit{Reflexivity:} Any $a\in R$ satisfies $a\sim a$, since the unity
$1_{R}$ is a unit and satisfies $a=a1_{R}$.

\textit{Symmetry:} If $a,b\in R$ satisfy $a\sim b$, then they also satisfy
$b\sim a$. Indeed, $a\sim b$ shows that there is a unit $u$ of $R$ such that
$a=bu$; but this unit $u$ clearly has an inverse $u^{-1}$, which is itself a
unit and satisfies $b=au^{-1}$. But this shows that $b\sim a$.

\textit{Transitivity:} If $a,b,c\in R$ satisfy $a\sim b$ and $b\sim c$, then
they also satisfy $a\sim c$. Indeed, there exist two units $u$ and $v$ of $R$
such that $b=cu$ and $a=bv$ (since $b\sim c$ and $a\sim b$); but the product
$uv$ of these two units is again a unit, and satisfies $a=\underbrace{b}%
_{=cu}v=cuv$, so that $a\sim c$.
\end{proof}

Note that an element $a$ of a ring $R$ is associate to $1$ if and only if $a$
is a unit.

If two elements $a$ and $b$ of a ring $R$ are associate, then each is a
multiple of the other (i.e., we have $a\mid b$ and $b\mid a$). When $R$ is an
integral domain, the converse holds as well:

\begin{proposition}
\label{prop.ring.associate.mutdev}Let $R$ be an integral domain. Let $a,b\in
R$ be such that $a\mid b$ and $b\mid a$. Then, $a\sim b$.
\end{proposition}

\begin{proof}
From $a\mid b$, we see that there exists an $x\in R$ such that $b=ax$.
Consider this $x$.

From $b\mid a$, we see that there exists a $y\in R$ such that $a=by$. Consider
this $y$.

If $a=0$, then the claim is easy (indeed, if $a=0$, then $b=\underbrace{a}%
_{=0}x=0$, so that $a=0=b$ and thus $a\sim b$). Hence, we WLOG assume that
$a\neq0$.

Now, $a=\underbrace{b}_{=ax}y=axy$. In other words, $a\left(  1-xy\right)
=0$. Since $a\neq0$, we thus conclude $1-xy=0$ (since $R$ is an integral
domain). In other words, $xy=1$. Thus, $y$ is a unit (since $R$ is
commutative). Hence, from $a=by$, we obtain $a\sim b$.
\end{proof}

Note that Proposition \ref{prop.ring.associate.mutdev} becomes false if we
drop the \textquotedblleft integral domain\textquotedblright\ condition. Some
sophisticated counterexamples can be found at
\url{https://math.stackexchange.com/questions/14270/} and in Exercise
\ref{exe.polring.mulvar.x+xy} below.

Associate elements \textquotedblleft look the same\textquotedblright\ to
divisibility, by which I mean that a divisibility relation of the form $a\mid
b$ remains equivalent if we replace $a$ by an element associate to $a$ or
replace $b$ by an element associate to $b$. In other words:

\begin{proposition}
\label{prop.ring.associate.div}Let $R$ be a commutative ring. Let
$a,b,a^{\prime},b^{\prime}\in R$ be such that $a\sim a^{\prime}$ and $b\sim
b^{\prime}$. Then, $a\mid b$ if and only if $a^{\prime}\mid b^{\prime}$.
\end{proposition}

\begin{proof}
$\Longrightarrow:$ Assume that $a\mid b$. From $a\sim a^{\prime}$, we see that
$a=a^{\prime}u$ for some unit $u$ of $R$. Hence, $a^{\prime}\mid a$. Also,
from $b\sim b^{\prime}$, we obtain $b^{\prime}\sim b$ (since Proposition
\ref{prop.ring.associate.eqrel} shows that the relation $\sim$ is symmetric).
In other words, $b^{\prime}=bv$ for some unit $v$ of $R$. Thus, $b\mid
b^{\prime}$. Hence, $a^{\prime}\mid a\mid b\mid b^{\prime}$. Thus, we have
proved the \textquotedblleft$\Longrightarrow$\textquotedblright\ direction of
Proposition \ref{prop.ring.associate.div}.

$\Longleftarrow:$ This is analogous to the \textquotedblleft$\Longrightarrow
$\textquotedblright\ direction, since Proposition
\ref{prop.ring.associate.eqrel} shows that the relation $\sim$ is symmetric.
\end{proof}

\begin{exercise}
Let $R$ be a commutative ring. Let $a,b,c,d\in R$ satisfy $a\sim b$ and $c\sim
d$. Prove that $ac\sim bd$.
\end{exercise}

\subsubsection{Uniqueness of gcds and lcms in an integral domain}

We can now state the uniqueness of gcds and lcms in the form in which it does hold:

\begin{proposition}
\label{prop.gcd-lcm.uni}Let $R$ be an integral domain. Let $a,b\in R$. Then:

\begin{enumerate}
\item[\textbf{(a)}] Any two gcds of $a$ and $b$ are associate (i.e., associate
to each other).

\item[\textbf{(b)}] Any two lcms of $a$ and $b$ are associate (i.e., associate
to each other).
\end{enumerate}
\end{proposition}

\begin{proof}
\textbf{(a)} Let $c$ and $d$ be two gcds of $a$ and $b$. We must show that
$c\sim d$.

Any common divisor of $a$ and $b$ is a divisor of $c$ (since $c$ is a gcd of
$a$ and $b$); however, $d$ is a common divisor of $a$ and $b$ (since $d$ is a
gcd of $a$ and $b$). Thus, $d$ is a divisor of $c$. In other words, $d\mid c$.
The same argument, with the roles of $c$ and $d$ swapped, yields $c\mid d$.
Hence, Proposition \ref{prop.ring.associate.mutdev} (applied to $c$ and $d$
instead of $a$ and $b$) yields $c\sim d$.

\textbf{(b)} Analogous to part \textbf{(a)}.
\end{proof}

From Proposition \ref{prop.gcd-lcm.uni}, we recover the fact that gcds and
lcms of integers are unique up to sign (since two integers $a$ and $b$ are
associate in $\mathbb{Z}$ if and only if $a=\pm b$).

\subsubsection{Existence of gcds and lcms in a PID}

We have now talked enough about uniqueness; when do gcds and lcms exist? The
following fact covers one important case:

\begin{theorem}
\label{thm.gcd-lcm.PID}Let $R$ be a PID. Let $a,b\in R$. Then, there exist a
gcd and an lcm of $a$ and $b$.
\end{theorem}

This will follow from the following proposition, which characterizes lcms and
partly characterizes gcds in terms of principal ideals:

\begin{proposition}
\label{prop.gcd-lcm.ideals}Let $R$ be a commutative ring. Let $a,b,c\in R$.

\begin{enumerate}
\item[\textbf{(a)}] If $aR+bR=cR$, then $c$ is a gcd of $a$ and $b$.

\item[\textbf{(b)}] We have $aR\cap bR=cR$ if and only if $c$ is an lcm of $a$
and $b$.
\end{enumerate}
\end{proposition}

Note that $aR+bR=cR$ is an equality between ideals (the $+$ sign on the left
hand side is a sum of ideals); it is \textbf{not} to be confused with $a+b=c$.
Confusingly, $a+b=c$ does \textbf{not} imply $aR+bR=cR$ (since there is no
\textquotedblleft distributivity law\textquotedblright\ that would equate
$\left(  a+b\right)  R$ with $aR+bR$). Instead, it is easy to see that
\textquotedblleft$aR+bR=cR$\textquotedblright\ is equivalent to
\textquotedblleft$a$ and $b$ are multiples of $c$, and there exist two
elements $u,v\in R$ satisfying $c=au+bv$\textquotedblright.

Note the difference between the two parts of Proposition
\ref{prop.gcd-lcm.ideals}: Part \textbf{(b)} is an \textquotedblleft if and
only if\textquotedblright, while part \textbf{(a)} is only an
\textquotedblleft if\textquotedblright. This is no accident: Proposition
\ref{prop.gcd-lcm.ideals} \textbf{(a)} cannot be extended to an
\textquotedblleft if and only if\textquotedblright\ statement. For example, in
the polynomial ring $\mathbb{Q}\left[  x,y\right]  $, the two polynomials $x$
and $y$ have gcd $1$; however, $1$ is not a $\mathbb{Q}\left[  x,y\right]
$-linear combination of $x$ and $y$.

\begin{proof}
[Proof of Proposition \ref{prop.gcd-lcm.ideals}.]\textbf{(a)} Assume that
$aR+bR=cR$. Thus, $c\in cR=aR+bR$. In other words, there exist $x,y\in R$ such
that $c=ax+by$. Hence, if $r$ is a common divisor of $a$ and $b$, then $r\mid
c$\ \ \ \ \footnote{\textit{Proof.} Let $r$ be a common divisor of $a$ and
$b$. Thus, $r\mid a$ and $r\mid b$. In other words, we can write $a$ and $b$
in the forms $a=ra^{\prime}$ and $b=rb^{\prime}$ for some $a^{\prime
},b^{\prime}\in R$. Using these $a^{\prime},b^{\prime}$, we obtain
$c=\underbrace{a}_{=ra^{\prime}}x+\underbrace{b}_{=rb^{\prime}}y=ra^{\prime
}x+rb^{\prime}y=r\left(  a^{\prime}x+b^{\prime}y\right)  $, so that $r\mid c$.
Qed.}. Thus, we have shown that any common divisor of $a$ and $b$ is a divisor
of $c$.

We have $a\in aR\subseteq aR+bR=cR$. In other words, $c\mid a$. Similarly,
$c\mid b$. Hence, $c$ is a common divisor of $a$ and $b$. Combining this
result with the result of the previous paragraph, we conclude that $c$ is a
gcd of $a$ and $b$. This proves Proposition \ref{prop.gcd-lcm.ideals}
\textbf{(a)}.

\textbf{(b)} Recall that an lcm of $a$ and $b$ was defined (in Definition
\ref{def.gcd-lcm} \textbf{(f)}) to be a common multiple $m$ of $a$ and $b$
with the property that every common multiple of $a$ and $b$ is a multiple of
$m$. Hence, we have the following chain of equivalences:%
\begin{align*}
&  \ \left(  c\text{ is an lcm of }a\text{ and }b\right) \\
&  \Longleftrightarrow\ \left(
\begin{array}
[c]{c}%
c\text{ is a common multiple of }a\text{ and }b\text{, and}\\
\text{every common multiple of }a\text{ and }b\text{ is a multiple of }c
\end{array}
\right) \\
&  \Longleftrightarrow\ \left(  c\in aR\cap bR\text{ and every element of
}aR\cap bR\text{ is a multiple of }c\right)
\end{align*}
(since the common multiples of $a$ and $b$ are precisely the elements of
$aR\cap bR$).

Now, let us look a bit closer at the statements on the right hand side. The
statement \textquotedblleft$c\in aR\cap bR$\textquotedblright\ is equivalent
to \textquotedblleft$cR\subseteq aR\cap bR$\textquotedblright\ (indeed, the
set $aR\cap bR$ is an ideal of $R$, and thus it contains the element $c$ if
and only if it contains all multiples of $c$; in other words, it contains the
element $c$ if and only if it contains the subset $cR$). The statement
\textquotedblleft every element of $aR\cap bR$ is a multiple of $c$%
\textquotedblright\ is equivalent to \textquotedblleft$aR\cap bR\subseteq
cR$\textquotedblright\ (since $cR$ is the set of all multiples of $c$). Thus,
our chain of equivalences can be continued as follows:%
\begin{align*}
&  \ \left(  c\text{ is an lcm of }a\text{ and }b\right) \\
&  \Longleftrightarrow\ \left(  \underbrace{c\in aR\cap bR}%
_{\Longleftrightarrow\ cR\subseteq aR\cap bR}\text{ and }%
\underbrace{\text{every element of }aR\cap bR\text{ is a multiple of }%
c}_{\Longleftrightarrow\ aR\cap bR\subseteq cR}\right) \\
&  \Longleftrightarrow\ \left(  cR\subseteq aR\cap bR\text{ and }aR\cap
bR\subseteq cR\right) \\
&  \Longleftrightarrow\ \left(  aR\cap bR=cR\right)  .
\end{align*}
This proves Proposition \ref{prop.gcd-lcm.ideals} \textbf{(b)}.
\end{proof}

\begin{proof}
[Proof of Theorem \ref{thm.gcd-lcm.PID}.]The sum $aR+bR$ is an ideal of $R$,
and thus is a principal ideal (since $R$ is a PID). In other words, $aR+bR=cR$
for some $c\in R$. Consider this $c$. Hence, Proposition
\ref{prop.gcd-lcm.ideals} \textbf{(a)} yields that $c$ is a gcd of $a$ and
$b$. Hence, a gcd of $a$ and $b$ exists.

The intersection $aR\cap bR$ is an ideal of $R$, and thus is a principal ideal
(since $R$ is a PID). In other words, $aR\cap bR=cR$ for some $c\in R$.
Consider this $c$. Hence, Proposition \ref{prop.gcd-lcm.ideals} \textbf{(b)}
yields that $c$ is an lcm of $a$ and $b$. Hence, an lcm of $a$ and $b$ exists.
Theorem \ref{thm.gcd-lcm.PID} is now proven.
\end{proof}

So any two elements of a PID have a gcd and an lcm. If the PID is Euclidean,
then the gcd can be computed by the Euclidean algorithm. Indeed, even more
generally, if a pair $\left(  a,b\right)  $ of two elements of a commutative
ring $R$ has a Bezout $5$-tuple (see Definition \ref{def.rings.euclid.bezout}
for the meaning of this notion), then it has a gcd:

\begin{corollary}
\label{cor.gcd-lcm.eucl-bez}Let $R$ be a commutative ring. Let $a,b\in R$. Let
$\left(  x,y,c,u,v\right)  $ be a Bezout $5$-tuple for $\left(  a,b\right)  $.
Then, $c$ is a gcd of $a$ and $b$.
\end{corollary}

\begin{proof}
Proposition \ref{prop.rings.euclid.ea.bezout} \textbf{(a)} yields that
$aR+bR=cR$. Hence, Proposition \ref{prop.gcd-lcm.ideals} \textbf{(a)} shows
that $c$ is a gcd of $a$ and $b$. This proves Corollary
\ref{cor.gcd-lcm.eucl-bez}.
\end{proof}

In a Euclidean ring $R$, we can use the (generalized) extended Euclidean
algorithm (as explained in the proof of Theorem \ref{thm.rings.euclid.ea.main}
\textbf{(b)}) to compute a Bezout $5$-tuple for any pair $\left(  a,b\right)
\in R\times R$. Thus, by Corollary \ref{cor.gcd-lcm.eucl-bez}, we can compute
a gcd of $a$ and $b$. (See \cite[pages 275--276]{DumFoo04} for an example of
this computation.)

We also note that the converse of Proposition \ref{prop.gcd-lcm.ideals}
\textbf{(a)} is false in general, but true when $R$ is a PID:

\begin{proposition}
\label{prop.gcd-lcm.ideals.pid-a}Let $R$ be a PID. Let $a,b,c\in R$. Then, we
have $aR+bR=cR$ if and only if $c$ is a gcd of $a$ and $b$.
\end{proposition}

\begin{proof}
The \textquotedblleft$\Longrightarrow$\textquotedblright\ direction follows
from Proposition \ref{prop.gcd-lcm.ideals} \textbf{(a)}, so we must only prove
the \textquotedblleft$\Longleftarrow$\textquotedblright\ direction. So let us
assume that $c$ is a gcd of $a$ and $b$. We must show that $aR+bR=cR$.

We have $aR+bR=dR$ for some $d\in R$ (since $R$ is a PID, so the ideal $aR+bR$
must be principal). This $d$ must therefore be a gcd of $a$ and $b$ (by
Proposition \ref{prop.gcd-lcm.ideals} \textbf{(a)}, applied to $d$ instead of
$c$). Hence, Proposition \ref{prop.gcd-lcm.uni} \textbf{(a)} shows that $c$
and $d$ are associate (since both $c$ and $d$ are gcds of $a$ and $b$). That
is, $c\mid d$ and $d\mid c$. Hence, $cR=dR$ (since $c\mid d$ entails
$dR\subseteq cR$, whereas $d\mid c$ entails $cR\subseteq dR$). Therefore,
$aR+bR=dR=cR$, which completes our proof.
\end{proof}

\subsubsection{More about gcds and lcms}

In an integral domain $R$, the gcd and the lcm of two elements $a,b\in R$
determine one another (up to associates) via the formula%
\[
\gcd\left(  a,b\right)  \cdot\operatorname{lcm}\left(  a,b\right)  \sim ab.
\]
This follows from the next exercise (\cite[homework set \#2, Exercise 3]%
{21w}), which also shows that the existence of an lcm implies the existence of
a gcd:

\begin{exercise}
\label{exe.21hw2.3}Let $R$ be an integral domain. Let $a\in R$ and $b\in R$.
Assume that $a$ and $b$ have an lcm $\ell\in R$. Prove that $a$ and $b$ have a
gcd $g\in R$, which furthermore satisfies $g\ell=ab$. \medskip

[\textbf{Hint:} If $u$ and $v$ are two elements of an integral domain $R$,
with $v\neq0$, then you can use the notation $\dfrac{u}{v}$ (or $u/v$) for the
element $w\in R$ satisfying $u=vw$. This element $w$ does not always exist,
but when it does, it is unique, so the notation is unambiguous. It is also
easy to see that standard rules for fractions, such as $\dfrac{u}{v}+\dfrac
{x}{y}=\dfrac{uy+vx}{vy}$ and $\dfrac{u}{v}\cdot\dfrac{x}{y}=\dfrac{ux}{vy}$,
hold as long as the fractions $\dfrac{u}{v}$ and $\dfrac{x}{y}$ exist.]
\end{exercise}

\begin{fineprint}
The converse of Exercise \ref{exe.21hw2.3} is false: The existence of a gcd of
two given elements $a$ and $b$ of an integral domain $R$ does not imply the
existence of an lcm of these two elements. However, if an integral domain $R$
has a gcd for \textbf{each} pair of two elements $a$ and $b$, then it also has
an lcm for each pair. This will be stated as Exercise
\ref{exe.gcd-lcm.gcd-to-lcm} below.

First, we state a more basic exercise, which generalizes the well-known
property $\gcd\left(  am,bm\right)  =\gcd\left(  a,b\right)  \cdot\left\vert
m\right\vert $ that holds for any three integers $a,b,m$:
\end{fineprint}

\begin{exercise}
\label{exe.gcd-lcm.gcdm}Let $R$ be an integral domain. Let $a,b,m\in R$ be
arbitrary. Assume that the elements $a$ and $b$ have a gcd $g$. Assume that
the elements $am$ and $bm$ have a gcd $h$. Prove that $gm\sim h$.
\end{exercise}

\begin{exercise}
\label{exe.gcd-lcm.gcd-to-lcm}Let $R$ be an integral domain. Assume that for
every $a,b\in R$, the elements $a$ and $b$ have a gcd. Prove that for every
$a,b\in R$, the elements $a$ and $b$ have an lcm. \medskip

[\textbf{Hint:} Let $a,b\in R$, and assume WLOG that $a,b\neq0$. Let $g$ be a
gcd of $a$ and $b$. It suffices to show that $\dfrac{ab}{g}$ is an lcm of $a$
and $b$. To this purpose, show first that $\dfrac{ab}{g}$ is a common multiple
of $a$ and $b$. Now, let $m$ be any common multiple of $a$ and $b$. Let $h$ be
a gcd of $am$ and $bm$. Argue that $ab\mid h\mid gm$ (by Exercise
\ref{exe.gcd-lcm.gcdm}). Conclude that $\dfrac{ab}{g}\mid m$.]
\end{exercise}

\begin{exercise}
Let $R$ be a PID. Let $a,b,c\in R$ be arbitrary. Prove the following:

\begin{enumerate}
\item[\textbf{(a)}] If $a\mid bc$, then $a\mid\gcd\left(  a,b\right)  c$.
(Here and in the rest of this exercise, $\gcd\left(  u,v\right)  $ means some
gcd of $u$ and $v$. We don't care which one we choose, since they are all associate.)

\item[\textbf{(b)}] We have $\gcd\left(  a,b\right)  \mid\gcd\left(
a,bc\right)  $.

\item[\textbf{(c)}] We have $\gcd\left(  a,bc\right)  \mid\gcd\left(
a,b\right)  \gcd\left(  a,c\right)  $.

\item[\textbf{(d)}] We have $\gcd\left(  a,b\right)  \gcd\left(  a,c\right)
\mid a\gcd\left(  b,c\right)  $.

\item[\textbf{(e)}] If $\gcd\left(  b,c\right)  \sim1$, then $\gcd\left(
a,bc\right)  \sim\gcd\left(  a,b\right)  \gcd\left(  a,c\right)  $.
\end{enumerate}
\end{exercise}

See \cite{Friedm19} for more about PIDs.

\subsection{\label{sec.rings.ufd}Unique factorization domains (\cite[\S 8.3]%
{DumFoo04})}

\subsubsection{\label{subsec.rings.ufd.irr-prime}Irreducible and prime
elements}

The notions of integral domains, of Euclidean domains and of PIDs are
abstractions for certain properties that hold for the ring $\mathbb{Z}$: The
first one abstracts the fact that products of nonzero integers are nonzero;
the second abstracts division with remainder; the third abstracts the fact
that each ideal of $\mathbb{Z}$ is principal. As we have seen, PIDs are a
weaker form of Euclidean domains. Even weaker is the notion of a \textbf{UFD}
(short for \textbf{Unique Factorization Domain}). This abstracts the existence
and the uniqueness of a prime factorization for integers. How do we define it
in an arbitrary integral domain? What is a good analogue of a prime number in
a general integral domain?

There are at least four such analogues. Let us introduce the first
two:\footnote{The other two are not properties of an element of $R$, but
rather properties of an ideal of $R$. See Definition
\ref{def.ideals.prime-max} for them.}

\begin{definition}
\label{def.rings.prime-irred}Let $R$ be a commutative ring. Let $r\in R$ be
nonzero and not a unit.

\begin{enumerate}
\item[\textbf{(a)}] We say that $r$ is \textbf{irreducible} (in $R$) if it has
the following property: Whenever $a,b\in R$ satisfy $ab=r$, at least one of
$a$ and $b$ is a unit.

\item[\textbf{(b)}] We say that $r$ is \textbf{prime} (in $R$) if it has the
following property: Whenever $a,b\in R$ satisfy $r\mid ab$, we have $r\mid a$
or $r\mid b$.
\end{enumerate}
\end{definition}

Let us see what these concepts mean when $R=\mathbb{Z}$. Both notions
\textquotedblleft irreducible\textquotedblright\ and \textquotedblleft
prime\textquotedblright\ smell like prime numbers, but it is worth being
precise: Not only the prime numbers $2,3,5,7,11,\ldots$ themselves, but also
their negatives $-2,-3,-5,-7,-11,\ldots$ fit both bills (i.e., they are
irreducible and prime in $\mathbb{Z}$). Let us be more explicit:

\begin{proposition}
\label{prop.rings.prime-irred-Z}Let $r\in\mathbb{Z}$. Then, we have the
following equivalences:%
\[
\left(  r\text{ is prime in }\mathbb{Z}\right)  \Longleftrightarrow\left(
r\text{ is irreducible in }\mathbb{Z}\right)  \Longleftrightarrow\left(
\left\vert r\right\vert \text{ is a prime number}\right)  .
\]

\end{proposition}

\begin{proof}
It suffices to prove the three implications%
\begin{align*}
\left(  r\text{ is prime in }\mathbb{Z}\right)   &  \Longrightarrow\left(
r\text{ is irreducible in }\mathbb{Z}\right)  ;\\
\left(  r\text{ is irreducible in }\mathbb{Z}\right)   &  \Longrightarrow
\left(  \left\vert r\right\vert \text{ is a prime number}\right)  ;\\
\left(  \left\vert r\right\vert \text{ is a prime number}\right)   &
\Longrightarrow\left(  r\text{ is prime in }\mathbb{Z}\right)  .
\end{align*}
All of them are LTTR. (The first one is actually a particular case of
Proposition \ref{prop.irred.prime-is-irred} further below. For the other two,
it is recommended to WLOG assume that $r\geq0$, since it is easy to see that
none of the three statements involved changes when $r$ is replaced by $-r$.)
\end{proof}

Thus, in the ring $\mathbb{Z}$, being prime and being irreducible is the same
thing. What about arbitrary integral domains? Here it is not quite the case,
as the following two examples show:

\begin{itemize}
\item In the ring $\mathbb{Z}\left[  \sqrt{-5}\right]  $, the element $3$ is
irreducible but not prime (in $\mathbb{Z}\left[  \sqrt{-5}\right]  $). (See
\cite[\S 8.3]{DumFoo04} for the proof.)

\item Here is an example using polynomials: Define a subset $R$ of the
univariate polynomial ring $\mathbb{Q}\left[  x\right]  $ by\footnote{The
\textquotedblleft$x^{1}$-coefficient\textquotedblright\ of a polynomial $f$
means the coefficient of $f$ before $x^{1}$. For example, the $x^{1}%
$-coefficient of $\left(  x+1\right)  ^{6}$ is $6$, whereas the $x^{1}%
$-coefficient of $x^{2}+1$ is $0$.}%
\begin{align*}
R  &  =\left\{  f\in\mathbb{Q}\left[  x\right]  \ \mid\ \text{the }%
x^{1}\text{-coefficient of }f\text{ is }0\right\} \\
&  =\left\{  f\in\mathbb{Q}\left[  x\right]  \ \mid\ \text{the derivative of
}f\text{ at }0\text{ is }0\right\} \\
&  =\left\{  a_{0}+a_{2}x^{2}+a_{3}x^{3}+\cdots+a_{n}x^{n}\ \mid\ n\geq0\text{
and }a_{0},a_{2},a_{3},\ldots,a_{n}\in\mathbb{Q}\right\}  .
\end{align*}
It is not hard to see that $R$ is a subring of $\mathbb{Q}\left[  x\right]  $.
(Indeed, if $f$ and $g$ are two polynomials whose $x^{1}$-coefficients are
$0$, then the same holds for $f+g$ and $f-g$ and $fg$. This is easiest to see
by computing $f+g$ and $f-g$ and $fg$ and checking that there is no way an
$x^{1}$-monomial can appear in the results.)

When we study polynomials later on, we will prove that $\mathbb{Q}\left[
x\right]  $ is an integral domain.\footnote{This is in fact pretty easy: When
you multiply two nonzero polynomials in $\mathbb{Q}\left[  x\right]  $, their
leading terms get multiplied, so their degrees get added; thus, the product
cannot be $0$.} Thus, the ring $R$ (being a subring of the integral domain
$\mathbb{Q}\left[  x\right]  $) must itself be an integral domain (since a
subring of an integral domain is always itself an integral
domain\footnote{This is obvious if you recall the definition of an integral
domain.}).

Now, the ring $R$ contains no polynomials of degree $1$. However, if
$a,b\in\mathbb{Q}\left[  x\right]  $ are two polynomials satisfying $x^{3}%
=ab$, then $3=\deg\left(  x^{3}\right)  =\deg\left(  ab\right)  =\deg a+\deg
b$, which means that one of the polynomials $a$ and $b$ is either a constant
(and thus a unit in $R$) or has degree $1$ (and thus cannot lie in $R$). This
quickly shows that the element $x^{3}$ of $R$ is irreducible in $R$. However,
this element is not prime in $R$ (since $x^{3}\mid x^{2}x^{4}$ but $x^{3}\nmid
x^{2}$ and $x^{3}\nmid x^{4}$).
\end{itemize}

In each of these two examples, we found an irreducible element that is not
prime. Can we do the opposite? No, as the following fact shows:

\begin{proposition}
\label{prop.irred.prime-is-irred}Let $R$ be an integral domain. Then, any
prime element of $R$ is irreducible.
\end{proposition}

\begin{proof}
Let $r\in R$ be prime. We must show that $r$ is irreducible.

So let $a,b\in R$ satisfy $ab=r$. We must show that at least one of $a$ and
$b$ is a unit.

We have $ab=r$, so that $r\mid ab$. Since $r$ is prime, we thus obtain $r\mid
a$ or $r\mid b$ (by the definition of \textquotedblleft
prime\textquotedblright). Assume WLOG that $r\mid a$ (since otherwise, we have
$r\mid b$, so we can swap $a$ with $b$ to achieve $r\mid a$). Hence, $a=rx$
for some $x\in R$. Consider this $x$. Now, $r=\underbrace{a}_{=rx}b=rxb$, and
therefore $r\left(  1-xb\right)  =r-rxb=0$, and thus $1-xb=0$ (since $r\neq0$
and since $R$ is an integral domain). In other words, $xb=1$. This shows that
$b$ is a unit (since $R$ is commutative). Thus we have shown that at least one
of $a$ and $b$ is a unit. This completes the proof that $r$ is irreducible.
\end{proof}

In a PID, the converse of Proposition \ref{prop.irred.prime-is-irred} also holds:

\begin{proposition}
\label{prop.PID.prime-iff-irred}Let $R$ be a PID. Let $r\in R$. Then, $r$ is
prime if and only if $r$ is irreducible.
\end{proposition}

\begin{proof}
We already showed the \textquotedblleft only if\textquotedblright\ part in
Proposition \ref{prop.irred.prime-is-irred}. We thus only need to prove the
\textquotedblleft if\textquotedblright\ part.

Assume that $r$ is irreducible. We must show that $r$ is prime.

Let $a,b\in R$ satisfy $r\mid ab$. We must prove that $r\mid a$ or $r\mid b$.

Assume the contrary. Thus, we have neither $r\mid a$ nor $r\mid b$.

There is an $h\in R$ such that $ab=rh$ (since $r\mid ab$). Consider this $h$.

Since $R$ is a PID, the ideal $aR+rR$ is principal; in other words, there
exists some $g\in R$ such that $gR=aR+rR$. Consider this $g$. Hence, $a\in
aR\subseteq aR+rR=gR$; in other words, $g\mid a$.

Also, $r\in rR\subseteq aR+rR=gR$; in other words, $g$ is a divisor of $r$.
However, $r$ is irreducible, and thus every divisor of $r$ is either a unit or
associate to $r$\ \ \ \ \footnote{\textit{Proof.} Let $d$ be a divisor of $r$.
We must show that $d$ is either a unit or associate to $r$.
\par
Indeed, there exists some $q\in R$ such that $r=dq$ (since $d$ is a divisor of
$r$). Consider this $q$. Since $r$ is irreducible, at least one of $d$ and $q$
is a unit. Hence, $d$ is either a unit or associate to $r$ (because if $q$ is
a unit, then $d$ is associate to $r$ (since $r=dq$ yields $r\sim d$ and thus
$d\sim r$)).}. Thus, $g$ is either a unit or associate to $r$ (since $g$ is a
divisor of $r$). However, if $g$ was associate to $r$, then we would have
$r\mid g\mid a$, which would contradict the fact that we don't have $r\mid a$.
Thus, $g$ cannot be associate to $r$, and so $g$ must be a unit. Therefore,
$1=gg^{-1}\in gR=aR+rR$. Hence, there exist $u,v\in R$ such that $1=au+rv$.

The same argument (using $b$ instead of $a$) shows that there exist
$u^{\prime},v^{\prime}\in R$ such that $1=bu^{\prime}+rv^{\prime}$.

Now, consider these four elements $u,v,u^{\prime},v^{\prime}$. Multiplying
$1=au+rv$ with $1=bu^{\prime}+rv^{\prime}$ yields%
\begin{align*}
1  &  =\left(  au+rv\right)  \left(  bu^{\prime}+rv^{\prime}\right)
=\underbrace{ab}_{=rh}uu^{\prime}+rvbu^{\prime}+aurv^{\prime}+rvrv^{\prime}\\
&  =rhuu^{\prime}+rvbu^{\prime}+aurv^{\prime}+rvrv^{\prime}%
=r\underbrace{\left(  huu^{\prime}+vbu^{\prime}+auv^{\prime}+vrv^{\prime
}\right)  }_{\in R}\in rR.
\end{align*}
In other words, there exists some $s\in R$ such that $1=rs$. This shows that
$r$ is a unit. This contradicts the fact that $r$ is irreducible. Thus, the
proof of Proposition \ref{prop.PID.prime-iff-irred} is complete.
\end{proof}

\begin{exercise}
Fix a nonzero integer $m$. Consider the ring $R_{m}$ from Exercise
\ref{exe.21hw1.1}.

Let $r\in R_{m}$ be nonzero.

\begin{enumerate}
\item[\textbf{(a)}] Define the $m$\textbf{-core} of $r$ to be the smallest
positive integer that is associate to $r$ in $R_{m}$. Prove that the $m$-core
of $r$ can be obtained as follows: Pick a $k\in\mathbb{N}$ such that
$m^{k}r\in\mathbb{Z}$ (such a $k$ exists, since $r\in R_{m}$). Write the prime
factorization of $\left\vert m^{k}r\right\vert $ as $\left\vert m^{k}%
r\right\vert =p_{1}p_{2}\cdots p_{i}q_{1}q_{2}\cdots q_{j}$, where
$p_{1},p_{2},\ldots,p_{i}$ are primes that divide $m$ and where $q_{1}%
,q_{2},\ldots,q_{j}$ are primes that don't divide $m$. (The primes don't have
to be distinct, and we allow $i=0$ or $j=0$.) Then, the $m$-core of $r$ is
$q_{1}q_{2}\cdots q_{j}$.

\item[\textbf{(b)}] Prove that $r$ is prime in $R_{m}$ if and only if the
$m$-core of $r$ is a prime number (in the usual number-theoretical sense).
\end{enumerate}
\end{exercise}

So we have generalized (in two ways, to boot) the notion of a prime number.
Let us now generalize prime factorization:

\subsubsection{Irreducible factorizations and UFDs}

\begin{definition}
Let $R$ be an integral domain.

\begin{enumerate}
\item[\textbf{(a)}] An \textbf{irreducible factorization} of an element $r\in
R$ means a tuple $\left(  p_{1},p_{2},\ldots,p_{n}\right)  $ of irreducible
elements $p_{1},p_{2},\ldots,p_{n}$ of $R$ such that $r\sim p_{1}p_{2}\cdots
p_{n}$. (Note that this tuple $\left(  p_{1},p_{2},\ldots,p_{n}\right)  $ can
be empty; in this case, the product $p_{1}p_{2}\cdots p_{n}$ is empty and thus
equals $1$. Thus, the empty tuple is an irreducible factorization of any unit
of $R$.)

\item[\textbf{(b)}] We say that $R$ is a \textbf{unique factorization domain}
(or, for short, \textbf{UFD}) if each nonzero $r\in R$ satisfies the following
two statements:

\begin{enumerate}
\item[1.] There exists an irreducible factorization of $r$.

\item[2.] The irreducible factorization of $r$ is unique up to associates.
This means the following: If $\left(  p_{1},p_{2},\ldots,p_{n}\right)  $ and
$\left(  q_{1},q_{2},\ldots,q_{m}\right)  $ are two irreducible factorizations
of $r$ (so that $p_{1},p_{2},\ldots,p_{n}$ and $q_{1},q_{2},\ldots,q_{m}$ are
irreducible elements of $R$ satisfying $r\sim p_{1}p_{2}\cdots p_{n}$ and
$r\sim q_{1}q_{2}\cdots q_{m}$), then we have $n=m$ and there is a bijection
$\alpha:\left\{  1,2,\ldots,n\right\}  \rightarrow\left\{  1,2,\ldots
,m\right\}  $ such that $p_{i}\sim q_{\alpha\left(  i\right)  }$ for each
$i\in\left\{  1,2,\ldots,n\right\}  $.
\end{enumerate}
\end{enumerate}
\end{definition}

My notion of an irreducible factorization differs slightly from that in
\cite{DumFoo04} (in that \cite{DumFoo04} requires $r=p_{1}p_{2}\cdots p_{n}$,
whereas we only require $r\sim p_{1}p_{2}\cdots p_{n}$); I hold mine to be
slightly better-behaved (for example, $-1\in\mathbb{Z}$ would not have an
irreducible factorization in the \cite{DumFoo04} sense). But my definition of
a UFD is equivalent to the one in \cite{DumFoo04}, as can be easily seen.

Soon, we will see that every PID is a UFD, and there are more UFDs than PIDs.
But first, let us see some examples of UFDs:

\begin{itemize}
\item The ring $\mathbb{Z}$ is a UFD. This is, of course, a consequence of
Euclid's famous theorem that says that any positive integer can be uniquely
decomposed into a product of primes. Our definition of an irreducible
factorization differs slightly from the classical notion of a prime
factorization in arithmetic, since our irreducible elements are allowed to be
negative and since we only require $r\sim p_{1}p_{2}\cdots p_{n}$ (rather than
$r=p_{1}p_{2}\cdots p_{n}$); but it is pretty easy to conciliate the two
concepts by replacing all negative factors by their absolute values. For
example, $\left(  -3,-2,2\right)  $ is an irreducible factorization of $-12$,
since $-12\sim\left(  -3\right)  \cdot\left(  -2\right)  \cdot2$; but of
course it corresponds to the classical prime factorization $12=3\cdot2\cdot2$
of the positive integer $12$.

\item Any field is a UFD, since every nonzero element is a unit and thus has
the empty tuple as its only irreducible factorization.

\item We shall soon see that every PID is a UFD.

\item The polynomial rings $\mathbb{Z}\left[  x\right]  $ (consisting of
polynomials in one variable with integer coefficients) and $\mathbb{Q}\left[
x,y\right]  $ (consisting of polynomials in two variables with rational
coefficients) are UFDs, even though they are not PIDs. (Of course,
$\mathbb{Q}\left[  x\right]  $ is a PID and thus a UFD as well.)

\item The rings%
\begin{align*}
\mathbb{Z}\left[  2i\right]   &  =\left\{  a+b\cdot2i\ \mid\ a,b\in
\mathbb{Z}\right\} \\
&  =\left\{  \text{Gaussian integers with an even imaginary part}\right\}
\end{align*}
and
\[
\mathbb{Z}\left[  \sqrt{-5}\right]  =\left\{  a+b\sqrt{-5}\ \mid
\ a,b\in\mathbb{Z}\right\}
\]
are not UFDs. (See Exercise \ref{exe.ufd.Z2i} for more on the former
ring.)\footnote{The ring $\mathbb{Z}\left[  \sqrt{-5}\right]  $ is an example
of a so-called \textbf{half-factorial domain}, i.e., an integral domain
satisfying the following weaker version of UFDness: Even though the
irreducible factorization of a nonzero element $r$ is not always unique up to
associates, the number of factors in this factorization is uniquely determined
by $r$. See \cite[Theorem 1]{Pollac24} for a proof.}
\end{itemize}

Previously (in Proposition \ref{prop.PID.prime-iff-irred}), we proved that an
element of a PID is prime if and only if it is irreducible. We shall now prove
the same result for UFDs (which is stronger, as we will soon see that every
PID is a UFD):

\begin{proposition}
\label{prop.UFD.prime-iff-irred}Let $R$ be a UFD. Let $r\in R$. Then, $r$ is
prime if and only if $r$ is irreducible.
\end{proposition}

\begin{proof}
$\Longrightarrow:$ If $r$ is prime, then $r$ is irreducible (by Proposition
\ref{prop.irred.prime-is-irred}).\footnote{Note that this holds for any
integral domain, not just for any UFD.}

$\Longleftarrow:$ Assume that $r$ is irreducible. We must show that $r$ is prime.

Let $a,b\in R$ satisfy $r\mid ab$. We must prove that $r\mid a$ or $r\mid b$.

Assume the contrary. Thus, neither $a$ nor $b$ is a multiple of $r$. Hence, in
particular, $a$ and $b$ are nonzero (since $0$ is a multiple of $r$). Thus,
$a$ and $b$ have irreducible factorizations (since $R$ is a UFD). Let $\left(
p_{1},p_{2},\ldots,p_{n}\right)  $ and $\left(  q_{1},q_{2},\ldots
,q_{m}\right)  $ be irreducible factorizations of $a$ and $b$. Thus,
$p_{1},p_{2},\ldots,p_{n}$ and $q_{1},q_{2},\ldots,q_{m}$ are irreducible
elements of $R$ satisfying%
\[
a\sim p_{1}p_{2}\cdots p_{n}\ \ \ \ \ \ \ \ \ \ \text{and}%
\ \ \ \ \ \ \ \ \ \ b\sim q_{1}q_{2}\cdots q_{m}.
\]
Multiplying $a\sim p_{1}p_{2}\cdots p_{n}$ with $b\sim q_{1}q_{2}\cdots q_{m}%
$, we see that
\begin{equation}
ab\sim p_{1}p_{2}\cdots p_{n}q_{1}q_{2}\cdots q_{m}
\label{pf.prop.UFD.prime-iff-irred.ab1}%
\end{equation}
(since a product of two units is again a unit).

However, $r\mid ab$. Thus, there exists a $q\in R$ such that $ab=rq$. Consider
this $q$. Note that $ab$ is nonzero (since $a$ and $b$ are nonzero, but $R$ is
an integral domain). Thus, $q$ is nonzero (since $q=0$ would imply
$ab=r\underbrace{q}_{=0}=0$, which would contradict the previous sentence).
Hence, $q$ has an irreducible factorization (since $R$ is a UFD). Let $\left(
s_{1},s_{2},\ldots,s_{k}\right)  $ be an irreducible factorization of $q$.
Thus, $s_{1},s_{2},\ldots,s_{k}$ are irreducible elements of $R$ satisfying
$q\sim s_{1}s_{2}\cdots s_{k}$. From $q\sim s_{1}s_{2}\cdots s_{k}$, we obtain
$rq\sim rs_{1}s_{2}\cdots s_{k}$. Since $ab=rq$, this rewrites as%
\begin{equation}
ab\sim rs_{1}s_{2}\cdots s_{k}. \label{pf.prop.UFD.prime-iff-irred.ab2}%
\end{equation}

Now, we conclude that the two tuples
\[
\left(  p_{1},p_{2},\ldots,p_{n},q_{1},q_{2},\ldots,q_{m}\right)  \text{ and
}\left(  r,s_{1},s_{2},\ldots,s_{k}\right)
\]
are two irreducible factorizations of $ab$ (since all their entries\newline%
$p_{1},p_{2},\ldots,p_{n},q_{1},q_{2},\ldots,q_{m}$ and $r,s_{1},s_{2}%
,\ldots,s_{k}$ are irreducible, and since
(\ref{pf.prop.UFD.prime-iff-irred.ab1}) and
(\ref{pf.prop.UFD.prime-iff-irred.ab2}) hold). Thus, by the uniqueness
condition in the definition of a UFD (which says that the irreducible
factorization of an element is unique up to associates), these two tuples must
be identical up to associates. In particular, every entry of the second tuple
must be associate to some entry of the first. Hence, in particular, the entry
$r$ of the second factorization must be associate to one of the entries
$p_{1},p_{2},\ldots,p_{n},q_{1},q_{2},\ldots,q_{m}$ of the first. In other
words, we must have%
\begin{equation}
r\sim p_{i}\text{ for some }i\in\left\{  1,2,\ldots,n\right\}
\label{pf.prop.UFD.prime-iff-irred.rpi}%
\end{equation}
or%
\begin{equation}
r\sim q_{j}\text{ for some }j\in\left\{  1,2,\ldots,m\right\}  .
\label{pf.prop.UFD.prime-iff-irred.rqj}%
\end{equation}
However, both of these possibilities lead to contradictions: Indeed, if
(\ref{pf.prop.UFD.prime-iff-irred.rpi}) holds, then we have $r\mid a$
(since\footnote{We will use the fact that associates divide each other: i.e.,
if $u$ and $v$ are two elements of $R$ satisfying $u\sim v$, then $u\mid v$.}
$r\sim p_{i}\mid p_{1}p_{2}\cdots p_{n}\sim a$), which contradicts the fact
that $a$ is not a multiple of $r$. Likewise, if
(\ref{pf.prop.UFD.prime-iff-irred.rqj}) holds, then we have $r\mid b$, which
contradicts the fact that $b$ is not a multiple of $r$. Thus, we get a
contradiction in either case, and our proof is complete.
\end{proof}

If $R$ is a UFD, and if $r\in R$ is nonzero, then $r$ is associate to a finite
product $p_{1}p_{2}\cdots p_{n}$ of irreducible elements (by the definition of
a UFD). This product can be simplified by collecting associate factors
together. For example, in $\mathbb{Z}$, we have%
\[
-24=2\cdot\left(  -2\right)  \cdot2\cdot3=-2^{3}\cdot3.
\]
Here is what we get in general:

\begin{proposition}
\label{prop.UFD.prime-power-fac}Let $R$ be a UFD. Let $r\in R$ be nonzero. Then:

\begin{enumerate}
\item[\textbf{(a)}] There exists a list $\left(  q_{1},q_{2},\ldots
,q_{k}\right)  $ of \textbf{mutually non-associate} irreducible elements
$q_{1},q_{2},\ldots,q_{k}\in R$ as well as a list $\left(  e_{1},e_{2}%
,\ldots,e_{k}\right)  $ of \textbf{positive} integers such that%
\[
r\sim q_{1}^{e_{1}}q_{2}^{e_{2}}\cdots q_{k}^{e_{k}}.
\]
We shall refer to these two lists as the \textbf{prime power factorization} of
$r$.

\item[\textbf{(b)}] These two lists are unique up to associates and up to
simultaneous permutation. (That is, any two prime power factorizations of $r$
can be transformed into one another by replacing the irreducible elements
$q_{1},q_{2},\ldots,q_{k}$ by associates, and reordering them while carrying
the exponents $e_{1},e_{2},\ldots,e_{k}$ along with them.)
\end{enumerate}
\end{proposition}

\begin{proof}
[Proof of Proposition \ref{prop.UFD.prime-power-fac}.]\textbf{(a)} Start with
an irreducible factorization of $r$, and collect associate factors together.
For example, if an irreducible factorization of $r$ has the form $\left(
p_{1},p_{2},p_{3},p_{4},p_{5},p_{6}\right)  $ with $p_{1}\sim p_{4}$ and
$p_{2}\sim p_{5}\sim p_{6}$ (and no other associate relations between its
entries), then%
\[
r\sim p_{1}p_{2}p_{3}p_{4}p_{5}p_{6}\sim p_{1}p_{2}p_{3}p_{1}p_{2}p_{2}%
=p_{1}^{2}p_{2}^{3}p_{3},
\]
and this is a prime power factorization of $r$. \medskip

\textbf{(b)} This follows from the uniqueness of an irreducible factorization
(up to associates).
\end{proof}

\begin{proposition}
\label{prop.UFD.gcd-fac-ex}Let $R$ be a UFD. Let $a,b\in R$ be nonzero. Then,
there exists a list $\left(  p_{1},p_{2},\ldots,p_{n}\right)  $ of
\textbf{mutually non-associate} irreducible elements $p_{1},p_{2},\ldots
,p_{n}\in R$ as well as two lists $\left(  e_{1},e_{2},\ldots,e_{n}\right)  $
and $\left(  f_{1},f_{2},\ldots,f_{n}\right)  $ of \textbf{nonnegative}
integers such that%
\[
a\sim p_{1}^{e_{1}}p_{2}^{e_{2}}\cdots p_{n}^{e_{n}}%
\ \ \ \ \ \ \ \ \ \ \text{and}\ \ \ \ \ \ \ \ \ \ b\sim p_{1}^{f_{1}}%
p_{2}^{f_{2}}\cdots p_{n}^{f_{n}}.
\]

\end{proposition}

\begin{proof}
Proposition \ref{prop.UFD.prime-power-fac} shows that $a$ and $b$ have prime
power factorizations%
\[
a\sim q_{1}^{e_{1}}q_{2}^{e_{2}}\cdots q_{k}^{e_{k}}%
\ \ \ \ \ \ \ \ \ \ \text{and}\ \ \ \ \ \ \ \ \ \ b\sim r_{1}^{f_{1}}%
r_{2}^{f_{2}}\cdots r_{m}^{f_{m}}.
\]
All we need now is to reconcile these prime power factorizations so that they
contain the same irreducible elements (albeit possibly with $0$ exponents).
For this purpose, we do the following steps:

\begin{enumerate}
\item If some of the $q_{i}$ are associate to some of the $r_{j}$, then we
replace these $q_{i}$ by the respective $r_{j}$.

\item If some of the $q_{i}$ don't appear among the $r_{j}$, then we insert
$q_{i}^{0}$ factors into the prime power factorization of $b$.

\item If some of the $r_{j}$ don't appear among the $q_{i}$, then we insert
$r_{j}^{0}$ factors into the prime power factorization of $a$.
\end{enumerate}

For example, if $R=\mathbb{Z}$ and $a=12$ and $b=45$, and if we start with the
prime power factorizations $a\sim2^{2}\cdot\left(  -3\right)  ^{1}$ and
$b\sim3^{2}\cdot5^{1}$, then Step 1 transforms the prime power factorization
of $a$ into $a\sim2^{2}\cdot3^{1}$ (since the $-3$ is replaced by the $3$ from
the prime power factorization of $b$); Step 2 then inserts a $2^{0}$ factor
into the prime power factorization of $b$ (so it becomes $b\sim2^{0}\cdot
3^{2}\cdot5^{1}$); Step 3 then inserts a $5^{0}$ factor into the prime power
factorization of $a$ (so it becomes $a\sim2^{2}\cdot3^{1}\cdot5^{0}$). The
resulting factorizations are $a\sim2^{2}\cdot3^{1}\cdot5^{0}$ and $b\sim
2^{0}\cdot3^{2}\cdot5^{1}$, just as promised by Proposition
\ref{prop.UFD.gcd-fac-ex}.
\end{proof}

\subsubsection{Gcds and lcms in a UFD}

\begin{proposition}
\label{prop.UFD.gcd-fac-2}Let $R$ be a UFD. Let $a,b\in R$ be nonzero. Let
$\left(  p_{1},p_{2},\ldots,p_{n}\right)  $, $\left(  e_{1},e_{2},\ldots
,e_{n}\right)  $ and $\left(  f_{1},f_{2},\ldots,f_{n}\right)  $ be as in
Proposition \ref{prop.UFD.gcd-fac-ex}. Then:

\begin{enumerate}
\item[\textbf{(a)}] The element%
\[
p_{1}^{\min\left\{  e_{1},f_{1}\right\}  }p_{2}^{\min\left\{  e_{2}%
,f_{2}\right\}  }\cdots p_{n}^{\min\left\{  e_{n},f_{n}\right\}  }%
\]
is a gcd of $a$ and $b$.

\item[\textbf{(b)}] The element%
\[
p_{1}^{\max\left\{  e_{1},f_{1}\right\}  }p_{2}^{\max\left\{  e_{2}%
,f_{2}\right\}  }\cdots p_{n}^{\max\left\{  e_{n},f_{n}\right\}  }%
\]
is an lcm of $a$ and $b$.
\end{enumerate}
\end{proposition}

\begin{proof}
This is done just as it is commonly done for integers in elementary number
theory. The details are LTTR. (See, e.g., \cite[proof of Proposition
1.11]{Friedm19} for some details on the proof of part \textbf{(a)}; the proof
of part \textbf{(b)} is similar.)
\end{proof}

\begin{corollary}
\label{cor.UFD.gcd-lcm}Any two elements in a UFD have a gcd and an lcm.
\end{corollary}

\begin{proof}
Let $a$ and $b$ be two elements of a UFD $R$. We must show that $a$ and $b$
have a gcd and an lcm.

If $b=0$, then this is easy (just show that $a$ is a gcd of $a$ and $0$, and
that $0$ is an lcm of $a$ and $0$). Thus, we WLOG assume that $b\neq0$. For a
similar reason, we WLOG assume that $a\neq0$. Hence, Proposition
\ref{prop.UFD.gcd-fac-ex} shows that there exists a list $\left(  p_{1}%
,p_{2},\ldots,p_{n}\right)  $ of \textbf{mutually non-associate} irreducible
elements $p_{1},p_{2},\ldots,p_{n}\in R$ as well as two lists $\left(
e_{1},e_{2},\ldots,e_{n}\right)  $ and $\left(  f_{1},f_{2},\ldots
,f_{n}\right)  $ of \textbf{nonnegative} integers such that%
\[
a\sim p_{1}^{e_{1}}p_{2}^{e_{2}}\cdots p_{n}^{e_{n}}%
\ \ \ \ \ \ \ \ \ \ \text{and}\ \ \ \ \ \ \ \ \ \ b\sim p_{1}^{f_{1}}%
p_{2}^{f_{2}}\cdots p_{n}^{f_{n}}.
\]
Thus, Proposition \ref{prop.UFD.gcd-fac-2} shows that $a$ and $b$ have a gcd
and an lcm.
\end{proof}

\subsubsection{Any PID is a UFD}

Finally, as promised, let us state the following theorem, which provides us
many UFDs to apply the above results to:

\begin{theorem}
\label{thm.UFD.PID-is-UFD}Any PID is a UFD.
\end{theorem}

I won't prove Theorem \ref{thm.UFD.PID-is-UFD} here; a proof can be found in
any of \cite[\S 8.3, Theorem 14]{DumFoo04}, \cite[Theorem 9.37]{Grifo26},
\cite[Corollary 12.2.13]{Mileti20} and \cite[Theorem 36.3]{Swanso17}. The
proof of the existence of an irreducible factorization is rather philosophical
and non-constructive; it yields no algorithm for actually finding such a
factorization. (And indeed, there are UFDs in which finding such a
factorization is algorithmically impossible.)\footnote{However, in many
specific PIDs, there are easy (although slow) algorithms for finding
irreducible factorizations. For instance, for $\mathbb{Z}\left[  i\right]  $,
Exercise \ref{exe.21hw2.6bc} asks you to find such an algorithm.} The proof of
the uniqueness of an irreducible factorization is an analogue of the proof you
know from elementary number theory (since we know that irreducible elements
are prime).

\subsubsection{A synopsis}

The following corollary combines several results we have seen above in a
convenient hierarchy:

\begin{corollary}
\label{cor.UFD.hierarchy}We have%
\begin{align*}
\left\{  \text{fields}\right\}   &  \subseteq\left\{  \text{Euclidean
domains}\right\}  \subseteq\left\{  \text{PIDs}\right\}  \subseteq\left\{
\text{UFDs}\right\} \\
&  \subseteq\left\{  \text{integral domains}\right\}  \subseteq\left\{
\text{commutative rings}\right\}  \subseteq\left\{  \text{rings}\right\}  .
\end{align*}

\end{corollary}

Let us illustrate this hierarchy in a symbolic picture:%
\[%
\begin{tikzpicture}[scale=0.6]
\draw(-4, -2) rectangle (16.8, 8);
\draw(-3.5, -1.75) rectangle (14.7, 7);
\draw(-3, -1.5) rectangle (12.6, 6);
\draw(-2.5, -1.25) rectangle (10.5, 5);
\draw(-2, -1) rectangle (8.4, 4);
\draw(-1.5, -0.75) rectangle (6.3, 3);
\draw(-1, -0.5) rectangle (4.2, 2);
\node(A) at (14, 7.5) {rings};
\node(A) at (10, 6.5) {commutative rings};
\node(A) at (8, 5.5) {integral domains};
\node(A) at (8, 4.5) {UFDs};
\node(A) at (6, 3.5) {PIDs};
\node(A) at (2, 2.5) {Euclidean domains};
\node(A) at (2, 1.5) {fields};
\end{tikzpicture}%
\]

All the \textquotedblleft$\subseteq$\textquotedblright\ signs in Corollary
\ref{cor.UFD.hierarchy} are strict inclusions; let us briefly recall some
examples showing this:

\begin{itemize}
\item The rings $\mathbb{Z}$ and $\mathbb{Z}\left[  i\right]  $ and
$\mathbb{Z}\left[  \sqrt{-2}\right]  $ and $\mathbb{Z}\left[  \sqrt{2}\right]
$ and the polynomial ring $\mathbb{Q}\left[  x\right]  $ are Euclidean
domains, but not fields.

\item The ring $\mathbb{Z}\left[  \alpha\right]  $ for $\alpha=\dfrac
{1+\sqrt{-19}}{2}$ is a PID, but not a Euclidean domain.

\item The polynomial rings $\mathbb{Q}\left[  x,y\right]  $ and $\mathbb{Z}%
\left[  x\right]  $ are UFDs, but not PIDs.

\item The rings $\mathbb{Z}\left[  2i\right]  $ and $\mathbb{Z}\left[
\sqrt{-3}\right]  $ are integral domains, but not UFDs.

\item The ring $\mathbb{Z}/6\cong\mathbb{Z}/2\times\mathbb{Z}/3$ is a
commutative ring, but not an integral domain.

\item The matrix ring $\mathbb{Q}^{2\times2}$ and the ring of quaternions
$\mathbb{H}$ are not commutative.
\end{itemize}

\subsection{\label{sec.rings.p=xx+yy}Application: Fermat's $p=x^{2}+y^{2}$
theorem (\cite[\S 8.3]{DumFoo04})}

\subsubsection{Fermat's two-squares theorem}

As an application of some of the above, we will show a result of
Fermat:\footnote{The word \textquotedblleft prime number\textquotedblright\ is
understood as in classical number theory -- i.e., a positive integer $p>1$
whose only positive divisors are $1$ and $p$. In particular, negative numbers
are not allowed as prime numbers, even though they are prime elements of
$\mathbb{Z}$.}

\begin{theorem}
[Fermat's two-squares theorem]\label{thm.fermat.p=xx+yy}Let $p$ be a prime
number such that $p\equiv1\operatorname{mod}4$. Then, $p$ can be written as a
sum of two perfect squares.
\end{theorem}

For example,%
\begin{align*}
5  &  =1^{2}+2^{2};\\
13  &  =2^{2}+3^{2};\\
17  &  =1^{2}+4^{2};\\
29  &  =2^{2}+5^{2}.
\end{align*}
(Note that the prime $2$ can also be written as a sum of two perfect squares:
$2=1^{2}+1^{2}$. But this would distract us from our proof.)

I will prove Theorem \ref{thm.fermat.p=xx+yy} using rings (specifically, using
the ring $\mathbb{Z}/p$ of residue classes and the ring $\mathbb{Z}\left[
i\right]  $ of Gaussian integers). Some of the steps will be left as exercises.

First, we shall show a general curious fact about primes, known as
\textbf{Wilson's theorem}:

\begin{theorem}
[Wilson's theorem]\label{thm.wilson}Let $p$ be a prime. Then, $\left(
p-1\right)  !\equiv-1\operatorname{mod}p$.
\end{theorem}

For example, for $p=5$, this is saying that $4!\equiv-1\operatorname{mod}5$.
And indeed, $4!=24\equiv-1\operatorname{mod}5$.

\begin{proof}
[Proof of Theorem \ref{thm.wilson}.]We must show that $\left(  p-1\right)
!\equiv-1\operatorname{mod}p$. Equivalently, we must show that%
\begin{equation}
\overline{\left(  p-1\right)  !}=\overline{-1}\text{ in }\mathbb{Z}/p.
\label{pf.thm.wilson.goal}%
\end{equation}

However, $\left(  p-1\right)  !=1\cdot2\cdot\cdots\cdot\left(  p-1\right)  $,
so that%
\begin{equation}
\overline{\left(  p-1\right)  !}=\overline{1\cdot2\cdot\cdots\cdot\left(
p-1\right)  }=\overline{1}\cdot\overline{2}\cdot\cdots\cdot\overline{p-1}.
\label{pf.thm.wilson.2}%
\end{equation}

Recall that every ring $R$ has a group of units, which is denoted by
$R^{\times}$. (See Theorem \ref{thm.ringunits.group} for details.)

But $\mathbb{Z}/p$ is a field (as we know, since $p$ is prime) with $p$
elements $\overline{0},\overline{1},\ldots,\overline{p-1}$. Its nonzero
elements $\overline{1},\overline{2},\ldots,\overline{p-1}$ are thus its units.
In other words, its group of units $\left(  \mathbb{Z}/p\right)  ^{\times}$ is
precisely the set $\left\{  \overline{1},\overline{2},\ldots,\overline
{p-1}\right\}  $ (and all the $p-1$ elements $\overline{1},\overline{2}%
,\ldots,\overline{p-1}$ are distinct). Hence,%
\begin{equation}
\prod_{a\in\left(  \mathbb{Z}/p\right)  ^{\times}}a=\overline{1}\cdot
\overline{2}\cdot\cdots\cdot\overline{p-1}. \label{pf.thm.wilson.3}%
\end{equation}

Recall that $\left(  \mathbb{Z}/p\right)  ^{\times}$ is a group. In
particular, any unit has an inverse, which is again a unit. The units
$\overline{1}$ and $\overline{-1}$ are their own inverses (since $\overline
{1}\cdot\overline{1}=\overline{1\cdot1}=\overline{1}$ and $\overline{-1}%
\cdot\overline{-1}=\overline{\left(  -1\right)  \cdot\left(  -1\right)
}=\overline{1}$), and they are the only units that are their own inverses
(this is Exercise \ref{exe.21hw2.5a}). The inverse of the inverse of a unit
$a$ is $a$. Hence, in the product $\prod\limits_{a\in\left(  \mathbb{Z}%
/p\right)  ^{\times}}a$, we can pair up each factor other than $\overline{1}$
and $\overline{-1}$ with its inverse:%
\begin{align}
\prod_{a\in\left(  \mathbb{Z}/p\right)  ^{\times}}a  &  =\underbrace{\left(
a_{1}\cdot a_{1}^{-1}\right)  }_{=\overline{1}}\cdot\underbrace{\left(
a_{2}\cdot a_{2}^{-1}\right)  }_{=\overline{1}}\cdot\cdots\cdot
\underbrace{\left(  a_{k}\cdot a_{k}^{-1}\right)  }_{=\overline{1}}%
\cdot\,\overline{1}\cdot\overline{-1}\nonumber\\
&  =\overline{1}\cdot\overline{1}\cdot\cdots\cdot\overline{1}\cdot\overline
{1}\cdot\overline{-1}=\overline{-1}. \label{pf.thm.wilson.4}%
\end{align}
Now, (\ref{pf.thm.wilson.2}) becomes%
\begin{align*}
\overline{\left(  p-1\right)  !}  &  =\overline{1}\cdot\overline{2}\cdot
\cdots\cdot\overline{p-1}=\prod_{a\in\left(  \mathbb{Z}/p\right)  ^{\times}%
}a\ \ \ \ \ \ \ \ \ \ \left(  \text{by (\ref{pf.thm.wilson.3})}\right) \\
&  =\overline{-1}\ \ \ \ \ \ \ \ \ \ \left(  \text{by (\ref{pf.thm.wilson.4}%
)}\right)  .
\end{align*}
This proves (\ref{pf.thm.wilson.goal}) and thus Theorem \ref{thm.wilson}.

(Caveat: The above was a little bit wrong for $p=2$; in that case, the factors
$\overline{1}$ and $\overline{-1}$ are actually one and the same factor. But
our proof can easily be adapted to this case.)
\end{proof}

\begin{corollary}
\label{cor.wilson.u!2gen}Let $p$ be an odd prime (i.e., a prime distinct from
$2$). Let $u=\dfrac{p-1}{2}\in\mathbb{N}$. Then, $u!^{2}\equiv-\left(
-1\right)  ^{u}\operatorname{mod}p$.
\end{corollary}

\begin{proof}
This is \cite[homework set \#2, Exercise 5 \textbf{(b)}]{21w}. Here is a
sketch of the proof: From $u=\dfrac{p-1}{2}$, we obtain $2u=p-1$ and thus
$2u+1=p$. Also, $p-1=2u$ (since $2u=p-1$) and thus%
\begin{align*}
\left(  p-1\right)  !  &  =\left(  2u\right)  !=\underbrace{1\cdot2\cdot
\cdots\cdot u}_{=u!}\cdot\left(  u+1\right)  \cdot\left(  u+2\right)
\cdot\cdots\cdot\left(  2u\right) \\
&  =u!\cdot\left(  \left(  u+1\right)  \cdot\left(  u+2\right)  \cdot
\cdots\cdot\left(  2u\right)  \right)  .
\end{align*}
In view of%
\begin{align*}
\left(  u+1\right)  \cdot\left(  u+2\right)  \cdot\cdots\cdot\left(
2u\right)   &  =\prod_{k=u+1}^{2u}k=\prod_{j=1}^{u}\underbrace{\left(
2u+1-j\right)  }_{\substack{\equiv-j\operatorname{mod}p\\\text{(since }\left(
2u+1-j\right)  -\left(  -j\right)  =2u+1=p\\\text{is divisible by }p\text{)}%
}}\\
&  \ \ \ \ \ \ \ \ \ \ \ \ \ \ \ \ \ \ \ \ \left(
\begin{array}
[c]{c}%
\text{here, we have substituted }2u+1-j\\
\text{for }k\text{ in the product}%
\end{array}
\right) \\
&  \equiv\prod_{j=1}^{u}\left(  -j\right)  =\left(  -1\right)  ^{u}%
\underbrace{\prod_{j=1}^{u}j}_{=u!}=\left(  -1\right)  ^{u}%
u!\operatorname{mod}p,
\end{align*}
this becomes%
\begin{align*}
\left(  p-1\right)  !  &  =u!\cdot\underbrace{\left(  \left(  u+1\right)
\cdot\left(  u+2\right)  \cdot\cdots\cdot\left(  2u\right)  \right)  }%
_{\equiv\left(  -1\right)  ^{u}u!\operatorname{mod}p}\\
&  \equiv u!\cdot\left(  -1\right)  ^{u}u!\equiv\left(  -1\right)  ^{u}\cdot
u!^{2}\operatorname{mod}p.
\end{align*}
Hence,%
\[
\left(  -1\right)  ^{u}\cdot u!^{2}\equiv\left(  p-1\right)  !\equiv
-1\operatorname{mod}p\ \ \ \ \ \ \ \ \ \ \left(  \text{by Theorem
\ref{thm.wilson}}\right)  .
\]
Multiplying both sides of this congruence by $\left(  -1\right)  ^{u}$, we
obtain $u!^{2}\equiv-\left(  -1\right)  ^{u}\operatorname{mod}p$. This proves
Corollary \ref{cor.wilson.u!2gen}.
\end{proof}

\begin{corollary}
\label{cor.wilson.u!2}Let $p$ be a prime such that $p\equiv1\operatorname{mod}%
4$. Let $u=\dfrac{p-1}{2}\in\mathbb{N}$. Then, $u!^{2}\equiv
-1\operatorname{mod}p$.
\end{corollary}

\begin{proof}
From $p\equiv1\operatorname{mod}4$, we obtain $4\mid p-1$, so that
$2\mid\dfrac{p-1}{2}=u$. Thus, $u$ is even, so that $\left(  -1\right)
^{u}=1$.

The prime $p$ is odd (since $p\equiv1\operatorname{mod}4$). Hence, Corollary
\ref{cor.wilson.u!2gen} yields $u!^{2}\equiv-\underbrace{\left(  -1\right)
^{u}}_{=1}=-1\operatorname{mod}p$. This proves Corollary \ref{cor.wilson.u!2}.
\end{proof}

Now, recall the ring $\mathbb{Z}\left[  i\right]  $ of Gaussian integers. Let
$N:\mathbb{Z}\left[  i\right]  \rightarrow\mathbb{N}$ be the map that sends
each Gaussian integer $a+bi$ (with $a,b\in\mathbb{Z}$) to $a^{2}+b^{2}%
\in\mathbb{N}$. It is straightforward to see:

\begin{proposition}
\label{prop.Zi.norm-mult}We have $N\left(  \alpha\beta\right)  =N\left(
\alpha\right)  N\left(  \beta\right)  $ for any $\alpha,\beta\in
\mathbb{Z}\left[  i\right]  $.
\end{proposition}

\begin{proof}
One way to prove this is by first showing that $N\left(  \gamma\right)
=\gamma\overline{\gamma}$ for each $\gamma\in\mathbb{Z}\left[  i\right]  $
(where $\overline{\gamma}$ denotes the complex conjugate of $\gamma$). Another
is by direct computation: Writing $\alpha$ and $\beta$ as $\alpha=a+bi$ and
$\beta=c+di$, we have $\alpha\beta=\left(  a+bi\right)  \left(  c+di\right)
=\left(  ac-bd\right)  +\left(  ad+bc\right)  i$ and therefore%
\begin{align*}
N\left(  \alpha\beta\right)   &  =N\left(  \left(  ac-bd\right)  +\left(
ad+bc\right)  i\right)  =\left(  ac-bd\right)  ^{2}+\left(  ad+bc\right)
^{2}\\
&  =a^{2}c^{2}-2acbd+b^{2}d^{2}+a^{2}d^{2}+2adbc+b^{2}c^{2}\\
&  =a^{2}c^{2}+b^{2}d^{2}+a^{2}d^{2}+b^{2}c^{2}=\underbrace{\left(
a^{2}+b^{2}\right)  }_{=N\left(  \alpha\right)  }\underbrace{\left(
c^{2}+d^{2}\right)  }_{=N\left(  \beta\right)  }=N\left(  \alpha\right)
N\left(  \beta\right)  .
\end{align*}

\end{proof}

\begin{corollary}
\label{cor.Zi.norm-div}If $z$ and $w$ are two Gaussian integers satisfying
$z\mid w$ in $\mathbb{Z}\left[  i\right]  $, then $N\left(  z\right)  \mid
N\left(  w\right)  $ in $\mathbb{Z}$.
\end{corollary}

\begin{proof}
Exercise (specifically, \cite[homework set \#2, Exercise 6 \textbf{(a)}]{21w}).
\end{proof}

Using this fact, we can characterize the units of $\mathbb{Z}\left[  i\right]
$:

\begin{corollary}
\label{cor.Zi.units}Let $\alpha\in\mathbb{Z}\left[  i\right]  $. Then, we have
the following equivalence:%
\[
\left(  \alpha\text{ is a unit of }\mathbb{Z}\left[  i\right]  \right)
\ \Longleftrightarrow\ \left(  N\left(  \alpha\right)  =1\right)
\ \Longleftrightarrow\ \left(  \alpha\in\left\{  1,i,-1,-i\right\}  \right)
.
\]

\end{corollary}

\begin{proof}
Exercise (specifically, \cite[homework set \#2, Exercise 6 \textbf{(d)}]{21w}).
\end{proof}

The next lemma is also easy to see:

\begin{lemma}
\label{lem.Zi.div-frac}Let $\alpha$ and $\beta$ be Gaussian integers such that
$\alpha\neq0$. Then, $\alpha\mid\beta$ holds in $\mathbb{Z}\left[  i\right]  $
if and only if $\dfrac{\beta}{\alpha}$ is a Gaussian integer.
\end{lemma}

\begin{proof}
This is proved just as the analogous statement for integers is proved.
\end{proof}

\begin{noncompile}
Let $m,a,b$ be integers. Then, $m\mid a+bi$ in $\mathbb{Z}\left[  i\right]  $
if and only if $m\mid a$ and $m\mid b$ in $\mathbb{Z}$.

\begin{proof}
If $m=0$, then this is easy (LTTR). So WLOG assume $m\neq0$. Then, we have the
following equivalences:%
\begin{align*}
\left(  m\mid a+bi\text{ in }\mathbb{Z}\left[  i\right]  \right)  \  &
\Longleftrightarrow\ \left(  \dfrac{a+bi}{m}\in\mathbb{Z}\left[  i\right]
\right)  \ \Longleftrightarrow\ \left(  \dfrac{a}{m}+\dfrac{b}{m}%
i\in\mathbb{Z}\left[  i\right]  \right) \\
&  \Longleftrightarrow\ \left(  \dfrac{a}{m}\text{ and }\dfrac{b}{m}\text{ are
integers}\right)  \ \Longleftrightarrow\ \left(  m\mid a\text{ and }m\mid
b\text{ in }\mathbb{Z}\right)  .
\end{align*}

\end{proof}
\end{noncompile}

Now we can prove Theorem \ref{thm.fermat.p=xx+yy}:

\begin{proof}
[Proof of Theorem \ref{thm.fermat.p=xx+yy}.]Let $u=\dfrac{p-1}{2}$. Then,
$u\in\mathbb{N}$ (actually, $p\equiv1\operatorname{mod}4$ implies that $u$ is
even). Corollary \ref{cor.wilson.u!2} shows that $u!^{2}\equiv
-1\operatorname{mod}p$. That is,%
\[
p\mid u!^{2}-\underbrace{\left(  -1\right)  }_{=i^{2}}=u!^{2}-i^{2}=\left(
u!-i\right)  \left(  u!+i\right)  .
\]
This is a divisibility in $\mathbb{Z}$, thus also in $\mathbb{Z}\left[
i\right]  $.

The number $p$ is a prime number, and thus prime in $\mathbb{Z}$; but this
does \textbf{not} mean that it is prime in $\mathbb{Z}\left[  i\right]  $. And
in fact, we claim that it isn't. Indeed, if $p$ were prime in $\mathbb{Z}%
\left[  i\right]  $, then the divisibility $p\mid\left(  u!-i\right)  \left(
u!+i\right)  $ would entail that $p\mid u!-i$ or $p\mid u!+i$; however,
neither $p\mid u!-i$ nor $p\mid u!+i$ is true\footnote{This is easiest to see
using Lemma \ref{lem.Zi.div-frac}: Indeed, if we had $p\mid u!-i$, then Lemma
\ref{lem.Zi.div-frac} would entail that $\dfrac{u!-i}{p}$ is a Gaussian
integer; however, $\dfrac{u!-i}{p}=\dfrac{u!}{p}+\dfrac{-1}{p}i$ is not a
Gaussian integer (since its imaginary part $\dfrac{-1}{p}$ is not an integer).
Thus, we don't have $p\mid u!-i$. For a similar reason, we don't have $p\mid
u!+i$.}.

Thus, we know that $p$ is not prime in $\mathbb{Z}\left[  i\right]  $. But
$\mathbb{Z}\left[  i\right]  $ is a Euclidean domain (as we proved in
Subsection \ref{subsec.rings.euclid.euclid}), and thus a PID (since
Proposition \ref{prop.eucldom.PID2} says that any Euclidean domain is a PID).
Hence, every irreducible element of $\mathbb{Z}\left[  i\right]  $ is a prime
element of $\mathbb{Z}\left[  i\right]  $ (by Proposition
\ref{prop.PID.prime-iff-irred}). Thus, $p$ cannot be irreducible in
$\mathbb{Z}\left[  i\right]  $ (since $p$ is not prime in $\mathbb{Z}\left[
i\right]  $).

However, $p$ is nonzero and not a unit of $\mathbb{Z}\left[  i\right]  $
(since $\dfrac{1}{p}$ is not a Gaussian integer). Therefore, since $p$ is not
irreducible, there exist two elements $\alpha,\beta\in\mathbb{Z}\left[
i\right]  $ that satisfy $\alpha\beta=p$ but are not units (by the definition
of \textquotedblleft irreducible\textquotedblright). Consider these $\alpha$
and $\beta$.

From $\alpha\beta=p$, we obtain $N\left(  \alpha\beta\right)  =N\left(
p\right)  =N\left(  p+0i\right)  =p^{2}+0^{2}=p^{2}$. Thus, $p^{2}=N\left(
\alpha\beta\right)  =N\left(  \alpha\right)  N\left(  \beta\right)  $ (by
Proposition \ref{prop.Zi.norm-mult}). However, $N\left(  \alpha\right)  $ and
$N\left(  \beta\right)  $ are nonnegative integers (since $N$ is a map
$\mathbb{Z}\left[  i\right]  \rightarrow\mathbb{N}$). Since $p$ is prime, the
only ways to write $p^{2}$ as a product of two nonnegative integers are
$p^{2}=1\cdot p^{2}$ and $p^{2}=p^{2}\cdot1$ and $p^{2}=p\cdot p$ (by the
classical prime factorization theorem from number theory). Hence, the equality
$p^{2}=N\left(  \alpha\right)  N\left(  \beta\right)  $ (with $N\left(
\alpha\right)  $ and $N\left(  \beta\right)  $ being nonnegative integers)
entails that we must be in one of the following two cases:

\textit{Case 1:} One of the two numbers $N\left(  \alpha\right)  $ and
$N\left(  \beta\right)  $ is $1$, and the other is $p^{2}$.

\textit{Case 2:} Both numbers $N\left(  \alpha\right)  $ and $N\left(
\beta\right)  $ are $p$.

Let us consider Case 1. In this case, one of the two numbers $N\left(
\alpha\right)  $ and $N\left(  \beta\right)  $ is $1$. We WLOG assume that
$N\left(  \alpha\right)  =1$ and $N\left(  \beta\right)  =p^{2}$ (since the
other possibility can be transformed into this one by swapping $\alpha$ with
$\beta$). Now, recall that $N\left(  \alpha\right)  =1$ is equivalent to
$\alpha$ being a unit (because of Corollary \ref{cor.Zi.units}). However,
$\alpha$ is not a unit. This is a contradiction. Hence, Case 1 is impossible.

Thus, we must be in Case 2. In other words, $N\left(  \alpha\right)  =p$ and
$N\left(  \beta\right)  =p$.

Now, $\alpha$ is a Gaussian integer, so we can write it as $\alpha=x+yi$ for
some $x,y\in\mathbb{Z}$. Therefore, using these $x,y$, we have $N\left(
\alpha\right)  =x^{2}+y^{2}$. Hence, $x^{2}+y^{2}=N\left(  \alpha\right)  =p$.
Thus, $p$ is a sum of two perfect squares; Theorem \ref{thm.fermat.p=xx+yy} is proven.
\end{proof}

Theorem \ref{thm.fermat.p=xx+yy} is one of the classic results of number
theory; while the proof we just gave is the most famous (variants of it appear
in \cite[\S 6.3]{Aluffi16}, \cite[\S 10.4]{Aluffi21}, \cite[Theorem
9.3.3]{BeaBla19}, \cite[Corollary 8.20]{CuoRot13}, \cite[\S 8.3, Proposition
18 (1)]{DumFoo04}, \cite[Theorem 32.5]{Elman22}, \cite[Theorem 3.5.15]%
{Laurit09}, \cite[Theorem 4.3.7]{Shifri96}, \cite[Theorem 3.8.1]{Weintr22}),
many other proofs are known as well (see, e.g., \cite[Corollary 2]{Bacher22},
\cite{Bachma25}, \cite[Theorem 13.2]{Burton11}, \cite[Theorem 1.2]{Cox22},
\cite[Theorem 2]{Elshol10}, \cite[Theorem 1.1.4]{Lehman19}, \cite[Theorem
100]{Nagell64}, \cite{Shiu21}, \cite[Theorem 5.7.1]{Stein09}). We will
ourselves give a second proof below (see Theorem \ref{thm.qr.jacobs}), using
finite fields instead of Gaussian integers.

We note that Theorem \ref{thm.fermat.p=xx+yy} has a converse as well (which,
however, is rather easy):

\begin{exercise}
\label{exe.fermat.p=xx+yy.converse}Let $p$ be a prime. Prove that $p$ can be
written as a sum of two perfect squares if and only if $p=2$ or $p\equiv
1\operatorname{mod}4$.
\end{exercise}

\subsubsection{Non-primes and numbers of representations}

Theorem \ref{thm.fermat.p=xx+yy} is the beginning of a long sequence of more
and more complex results, whose discovery and proof spanned several centuries.
First of all, one can ask which integers (rather than which primes) can be
written as sums of two perfect squares. The answer is not too hard at this point:

\begin{theorem}
\label{thm.fermat.n=xx+yy}Let $n$ be a positive integer with prime
factorization $n=2^{a}p_{1}^{b_{1}}p_{2}^{b_{2}}\cdots p_{k}^{b_{k}}$, where
$p_{1},p_{2},\ldots,p_{k}$ are distinct primes larger than $2$, and where
$a,b_{1},b_{2},\ldots,b_{k}$ are nonnegative integers. (In particular, if $n$
is odd, then $a=0$.)

Then:

\begin{enumerate}
\item[\textbf{(a)}] The number $n$ can be written as a sum of two perfect
squares if and only if the following condition holds: For each $i\in\left\{
1,2,\ldots,k\right\}  $ satisfying $p_{i}\equiv3\operatorname{mod}4$, the
exponent $b_{i}$ is even.

\item[\textbf{(b)}] If this condition holds, then the number of ways to write
$n$ as a sum of two perfect squares (to be more precise: the number of pairs
$\left(  x,y\right)  \in\mathbb{Z}\times\mathbb{Z}$ satisfying $n=x^{2}+y^{2}%
$) is $4\cdot\prod_{\substack{i\in\left\{  1,2,\ldots,k\right\}
;\\p_{i}\equiv1\operatorname{mod}4}}\left(  b_{i}+1\right)  $.
\end{enumerate}
\end{theorem}

For a proof of this theorem, see \cite[Theorem 4.2.62]{19s} or \cite[\S 8.1,
Corollary 19]{DumFoo04}. (The proof again uses Gaussian integers in a rather
neat way.)

\begin{exercise}
Prove the \textquotedblleft if\textquotedblright\ part of Theorem
\ref{thm.fermat.n=xx+yy} \textbf{(a)}. (Keep in mind that $0$ counts as a
perfect square.)
\end{exercise}

More about decompositions of integers into sums of perfect squares can be found

\begin{itemize}
\item in \cite[\S 8.3]{DumFoo04};

\item in Keith Conrad's
\url{https://kconrad.math.uconn.edu/math5230f12/handouts/Zinotes.pdf} ;

\item in \cite[\S 4.2]{19s}.
\end{itemize}

\subsubsection{Other forms $x^{2}+ny^{2}$}

Instead of writing integers $n$ in the form $n=x^{2}+y^{2}$, we can try to
write them in the form $x^{2}+2y^{2}$ or $x^{2}+3y^{2}$ or $x^{2}+4y^{2}$ or
$x^{2}+5y^{2}$ or $x^{2}+xy+y^{2}$ or $\left\vert x^{2}-2y^{2}\right\vert $ or
many other such forms. Each time, we can ask when this is possible, and how
many ways there are. These questions vary widely in difficulty, and even their
most basic variants (which prime numbers can be written in a given form?) can
be extremely hard. After having answered the $x^{2}+y^{2}$ question for
primes, the $x^{2}+4y^{2}$ question becomes quite easy (see Exercise
\ref{exe.rings.p=xx+4yy} below), but the $x^{2}+2y^{2}$ question will have to
wait until a later chapter (see Exercise \ref{exe.rings.p=xx+2yy} below). A
whole book \cite{Cox22} has been written entirely about the question of
writing prime(!) numbers in the form $x^{2}+ny^{2}$ for positive integers $n$;
just answering these questions for different $n$ requires rather advanced
mathematics. Here is a summary of answers for certain values of $n$ (see
\cite{Cox22} for proofs of these and many more results):

\begin{theorem}
\label{thm.rings.p=xx+nyy}Let $p$ be a prime number.

\begin{enumerate}
\item[\textbf{(a)}] We can write $p$ in the form $p=x^{2}+y^{2}$ with
$x,y\in\mathbb{Z}$ if and only if $p=2$ or $p\equiv1\operatorname{mod}4$.

\item[\textbf{(b)}] We can write $p$ in the form $p=x^{2}+2y^{2}$ with
$x,y\in\mathbb{Z}$ if and only if $p=2$ or $p\equiv1,3\operatorname{mod}8$.
(The notation \textquotedblleft$p\equiv1,3\operatorname{mod}8$%
\textquotedblright\ is shorthand for \textquotedblleft$p$ is congruent to $1$
or to $3$ modulo $8$\textquotedblright. Similar shorthands will be used in the
following parts.)

\item[\textbf{(c)}] We can write $p$ in the form $p=x^{2}+3y^{2}$ with
$x,y\in\mathbb{Z}$ if and only if $p=3$ or $p\equiv1\operatorname{mod}3$.

\item[\textbf{(d)}] We can write $p$ in the form $p=x^{2}+4y^{2}$ with
$x,y\in\mathbb{Z}$ if and only if $p\equiv1\operatorname{mod}4$.

\item[\textbf{(e)}] We can write $p$ in the form $p=x^{2}+5y^{2}$ with
$x,y\in\mathbb{Z}$ if and only if $p=5$ or $p\equiv1,9\operatorname{mod}20$.

\item[\textbf{(f)}] We can write $p$ in the form $p=x^{2}+6y^{2}$ with
$x,y\in\mathbb{Z}$ if and only if $p\equiv1,7\operatorname{mod}24$.

\item[\textbf{(g)}] We can write $p$ in the form $p=x^{2}+14y^{2}$ with
$x,y\in\mathbb{Z}$ if and only if we have $p\equiv
1,9,15,23,25,39\operatorname{mod}56$ and there exists some integer $z$
satisfying $\left(  z^{2}+1\right)  ^{2}\equiv8\operatorname{mod}p$.

\item[\textbf{(h)}] We can write $p$ in the form $p=x^{2}+27y^{2}$ with
$x,y\in\mathbb{Z}$ if and only if we have $p\equiv1\operatorname{mod}3$ and
there exists some integer $z$ satisfying $z^{3}\equiv2\operatorname{mod}p$.
\end{enumerate}
\end{theorem}

Part \textbf{(a)} of Theorem \ref{thm.rings.p=xx+nyy} follows from Theorem
\ref{thm.fermat.p=xx+yy} and Exercise \ref{exe.fermat.p=xx+yy.converse}. As we
said, parts \textbf{(b)} and \textbf{(d)} follow easily from Exercise
\ref{exe.rings.p=xx+2yy} and Exercise \ref{exe.rings.p=xx+4yy}. Part
\textbf{(c)} is similar, but the proof is trickier since $\mathbb{Z}\left[
\sqrt{-3}\right]  $ is not a PID (or even a UFD); nevertheless, fairly
elementary proofs exist, and one such proof is outlined in Exercise
\ref{exe.rings.p=xx+3yy}\footnote{See, e.g.,
\url{https://math.stackexchange.com/a/76917/} for another proof.}. Part
\textbf{(e)} is proved using genus theory of quadratic forms in \cite[(2.22)]%
{Cox22}, and using elementary techniques (quadratic reciprocity) in
\cite{Zhang07}. For part \textbf{(f)}, I am only aware of a genus-theoretical
proof (\cite[Exercise 2.19]{Cox22}). Part \textbf{(g)} requires class field
theory (\cite[Theorem 5.33]{Cox22}). Part \textbf{(h)} can be proved using
cubic reciprocity (\cite[Theorem 4.15]{Cox22}) or using class field theory
(\cite[Theorem 9.8]{Cox22}). Note the additional \textquotedblleft there
exists some integer $z$\textquotedblright\ conditions in parts \textbf{(g)}
and \textbf{(h)}; such conditions can be avoided for $x^{2}+ny^{2}$ when $n$
is small, but eventually become necessary. See
\url{https://mathoverflow.net/questions/79342/} for more about the need for
such conditions, and see \cite[Chapters 2 and 3]{Cox22} for their exact
nature. (Elementary proofs for a few of the more basic parts can be found in
\cite[Theorem 100]{Nagell64}.)

We can also ask about sums of more than two squares. Lagrange proved that
every nonnegative integer can be written as a sum of \textbf{four} squares
(that is, each $n\in\mathbb{N}$ can be written as $n=x^{2}+y^{2}+z^{2}+w^{2}$
for some $x,y,z,w\in\mathbb{Z}$). These days, one of the shortest proofs of
this fact uses the so-called \textit{Hurwitz quaternions} -- a quaternion
analogue of Gaussian integers. See
\url{https://en.wikipedia.org/wiki/Lagrange's_four-square_theorem} or
\cite{Haensc16} or \cite{Schwar14} for the proof.

\subsubsection{Coda}

\begin{exercise}
\label{exe.rings.p=xx+4yy}Let $p$ be a prime. Prove that $p$ can be written in
the form $p=x^{2}+4y^{2}$ for some $x,y\in\mathbb{Z}$ if and only if
$p\equiv1\operatorname{mod}4$.
\end{exercise}

\begin{exercise}
\label{exe.21hw2.6bc}\ \ 

\begin{enumerate}
\item[\textbf{(a)}] Let $z=a+bi\in\mathbb{Z}\left[  i\right]  $ with
$a,b\in\mathbb{Z}$. Assume that $z\neq0$. Let $n=\left\lfloor \left\vert
z\right\vert \right\rfloor =\left\lfloor \sqrt{a^{2}+b^{2}}\right\rfloor $.
Prove that every divisor of $z$ in $\mathbb{Z}\left[  i\right]  $ has the form
$c+di$ with $c,d\in\left\{  -n,-n+1,\ldots,n\right\}  $.

\item[\textbf{(b)}] Without recourse to the general theory of PIDs and UFDs,
prove that every nonzero element of $\mathbb{Z}\left[  i\right]  $ has an
irreducible factorization.
\end{enumerate}
\end{exercise}

\begin{exercise}
\label{exe.ufd.Z2i}Let $R$ be the ring $\mathbb{Z}\left[  i\right]  $ of
Gaussian integers. Let $S$ be the ring%
\begin{align*}
\mathbb{Z}\left[  2i\right]   &  =\left\{  a+b\cdot2i\mid a,b\in
\mathbb{Z}\right\} \\
&  =\left\{  \text{Gaussian integers with an even imaginary part}\right\}  .
\end{align*}
This ring $S$ is a subring of $R$.

Define two elements $x,y\in S$ by $x=2+2i$ and $y=2-2i$.

\begin{enumerate}
\item[\textbf{(a)}] Find the units of $S$.

\item[\textbf{(b)}] Prove that we have $x\sim y$ in $R$, but we don't have
$x\sim y$ in $S$.

\item[\textbf{(c)}] Prove that the ideal $xS+yS$ of $S$ is not principal.

\item[\textbf{(d)}] Conclude that $S$ is not a PID.

\item[\textbf{(e)}] Show that $S$ is not a UFD either.
\end{enumerate}

[\textbf{Hint:} It may be helpful to write $i^{\prime}$ for $2i$ in order to
avoid confusing $i$ for an element of $S$.

For part \textbf{(c)}, argue that if $xS+yS$ was $zS$ for some $z\in S$, then
$xR+yR$ would be $zR$ as well (why?), but this would force $z$ to be associate
to $x$ and $y$ in $R$ (why?), and this would leave only four possibilities for
$z$ (why?).

For part \textbf{(e)}, ponder the equality $xy=2\cdot2\cdot2$. Note that any
divisor of an element $s\in S$ is also a divisor of the same element $s$ in
$R=\mathbb{Z}\left[  i\right]  $ (but not always the other way round).]
\end{exercise}

\begin{exercise}
\label{exe.21hw2.7}Consider the ring
\[
\mathbb{Z}\left[  \sqrt{-3}\right]  =\left\{  a+b\sqrt{-3}\mid a,b\in
\mathbb{Z}\right\}  .
\]
This ring is a subring of $\mathbb{C}$, and thus is an integral domain.

Let $u=2\in\mathbb{Z}\left[  \sqrt{-3}\right]  $ and $v=1+\sqrt{-3}%
\in\mathbb{Z}\left[  \sqrt{-3}\right]  $. Further let $a=2u=4$ and $b=2v$.

\begin{enumerate}
\item[\textbf{(a)}] Prove that both $u$ and $v$ are common divisors of $a$ and
$b$ in $\mathbb{Z}\left[  \sqrt{-3}\right]  $.

\item[\textbf{(b)}] Prove that the only divisors of $4$ in $\mathbb{Z}\left[
\sqrt{-3}\right]  $ are $\pm1$, $\pm2$, $\pm4$, $\pm\left(  1+\sqrt
{-3}\right)  $, and $\pm\left(  1-\sqrt{-3}\right)  $.

\item[\textbf{(c)}] Prove that $a$ and $b$ have no gcd in $\mathbb{Z}\left[
\sqrt{-3}\right]  $.
\end{enumerate}

[\textbf{Remark:} This shows that $\mathbb{Z}\left[  \sqrt{-3}\right]  $ is
not a UFD, thus not a PID and not Euclidean.]
\end{exercise}

\begin{exercise}
\ \ 

\begin{enumerate}
\item[\textbf{(a)}] Prove that there are no ring morphisms from $\mathbb{Z}%
\left[  i\right]  $ to $\mathbb{Z}$.
\end{enumerate}

Now, let $p$ be a prime number. Prove the following:

\begin{enumerate}
\item[\textbf{(b)}] There are no ring morphisms from $\mathbb{Z}\left[
i\right]  $ to $\mathbb{Z}/p$ if $p\equiv3\operatorname{mod}4$.

\item[\textbf{(c)}] There are exactly two ring morphisms from $\mathbb{Z}%
\left[  i\right]  $ to $\mathbb{Z}/p$ if $p\equiv1\operatorname{mod}4$.

\item[\textbf{(d)}] There is a unique ring morphism from $\mathbb{Z}\left[
i\right]  $ to $\mathbb{Z}/p$ if $p=2$.
\end{enumerate}

More generally, prove the following:

\begin{enumerate}
\item[\textbf{(e)}] If $R$ is any ring, then the number of ring morphisms from
$\mathbb{Z}\left[  i\right]  $ to $R$ is the number of all elements $r\in R$
satisfying $r^{2}=-1$.
\end{enumerate}

[\textbf{Hint:} If $f$ is a ring morphism from $\mathbb{Z}\left[  i\right]  $
to $R$, then what equation must $f\left(  i\right)  $ satisfy?]
\end{exercise}

\subsection{\label{sec.rings.234it}More about ideals and quotient rings}

For the sake of completeness, let me mention three more general properties of
quotient rings, known respectively as the \textbf{second}, \textbf{third} and
\textbf{fourth isomorphism theorems for rings}. The second and third
isomorphism theorems claim isomorphisms between quotient rings; the fourth
relates the ideals of a quotient ring to some ideals of the original ring. We
shall not use these theorems, but they are not hard to prove and are part of
an algebraist's culture.

\subsubsection{The second isomorphism theorem for rings}

The \textbf{second isomorphism theorem for rings} is all about the interaction
between an ideal and a subring:

\begin{theorem}
[Second isomorphism theorem for rings]\label{thm.rings.2it}Let $R$ be a ring.
Let $S$ be a subring of $R$. Let $I$ be an ideal of $R$. Define $S+I$ to be
the subset $\left\{  s+i\mid s\in S\text{ and }i\in I\right\}  $ of $R$. Then:

\begin{enumerate}
\item[\textbf{(a)}] This subset $S+I$ is a subring of $R$.

\item[\textbf{(b)}] The set $I$ is an ideal of the ring $S+I$.

\item[\textbf{(c)}] The set $S\cap I$ is an ideal of the ring $S$.

\item[\textbf{(d)}] We have $\left(  S+I\right)  /I\cong S/\left(  S\cap
I\right)  $ (as rings). More concretely, there is a ring isomorphism
$S/\left(  S\cap I\right)  \rightarrow\left(  S+I\right)  /I$ that sends each
residue class $\overline{s}=s+\left(  S\cap I\right)  $ to $\overline{s}=s+I$.
\end{enumerate}
\end{theorem}

\begin{proof}
See \cite[homework set \#2, Exercise 9]{21w}.
\end{proof}

\begin{example}
For an example, we

\begin{itemize}
\item let $R$ be the polynomial ring $\mathbb{Q}\left[  x\right]  $ of all
univariate polynomials with rational coefficients;

\item let $I=\left\{  a_{2}x^{2}+a_{3}x^{3}+\cdots+a_{n}x^{n}\mid n\geq0\text{
and }a_{i}\in\mathbb{Q}\right\}  $ be the ideal consisting of all polynomials
divisible by $x^{2}$ (that is, all polynomials whose $x^{0}$-coefficient and
$x^{1}$-coefficient are $0$);

\item let $S$ be the subring $\mathbb{Q}$ of $R$ (which consists of all
constant polynomials).
\end{itemize}

Then, $S+I=\left\{  a_{0}+a_{2}x^{2}+a_{3}x^{3}+\cdots+a_{n}x^{n}\mid
n\geq0\text{ and }a_{i}\in\mathbb{Q}\right\}  $ is the set of all polynomials
whose $x^{1}$-coefficient is $0$. This is indeed a subring of $R$, as we have
seen in Subsection \ref{subsec.rings.ufd.irr-prime} (where we have used this
subring to find an irreducible element that is not prime).

Other interesting examples of $S+I$ can be obtained using upper-triangular
matrix rings such as $\mathbb{Q}^{3\leq3}$.
\end{example}

\subsubsection{The third isomorphism theorem for rings}

You might have noticed that the quotient rings $R/I$ of a given ring $R$ stand
in an \textquotedblleft antithetical\textquotedblright\ relationship to the
ideals $I$ that produce them: The larger the ideal $I$, the smaller the
quotient ring $R/I$. (In particular, the largest ideal of $R$ is $R$ itself,
and the corresponding quotient ring $R/R$ is trivial, which is as small as a
ring can get.)

Can we make this relationship precise? To some extent, we can. Namely, when
two ideals $I$ and $J$ of a ring $R$ satisfy $I\subseteq J$, the corresponding
quotient rings $R/I$ and $R/J$ are related as well, and specifically, $R/J$ is
isomorphic to a quotient ring of $R/I$. In other words, the quotient $R/J$ by
the \textquotedblleft large\textquotedblright\ ideal $J$ can be obtained by
first quotienting out a \textquotedblleft smaller\textquotedblright\ ideal $I$
and then \textquotedblleft quotienting further\textquotedblright. The
\textbf{third isomorphism theorem for rings} states this relationship in a
more concrete way:

\begin{theorem}
[Third isomorphism theorem for rings]\label{thm.rings.3it}Let $R$ be a ring.
Let $I$ and $J$ be two ideals of $R$ such that $I\subseteq J$. Let $J/I$
denote the set of all cosets $j+I\in R/I$ where $j\in J$. Then:

\begin{enumerate}
\item[\textbf{(a)}] This set $J/I$ is an ideal of $R/I$.

\item[\textbf{(b)}] We have $\left(  R/I\right)  /\left(  J/I\right)  \cong
R/J$ (as rings). More concretely, there is a ring isomorphism $R/J\rightarrow
\left(  R/I\right)  /\left(  J/I\right)  $ that sends each residue class
$\overline{r}=r+J$ to $\overline{r+I}=\left(  r+I\right)  +\left(  J/I\right)
$.
\end{enumerate}
\end{theorem}

\begin{proof}
See \cite[homework set \#2, Exercise 8]{21w}.
\end{proof}

\begin{example}
For an example, take $R=\mathbb{Z}$ and $I=6\mathbb{Z}$ and $J=2\mathbb{Z}$.
In this case, $J/I$ consists of the \textquotedblleft even\textquotedblright%
\ residue classes $\overline{0},\overline{2},\overline{4}$ in $R/I=\mathbb{Z}%
/6$. Theorem \ref{thm.rings.3it} \textbf{(b)} says that if we
\textquotedblleft quotient them out\textquotedblright\ of $\mathbb{Z}/6$, then
we are left with (an isomorphic copy of) $R/J=\mathbb{Z}/2$.
\end{example}

\subsubsection{The inverse image of an ideal}

The following easy fact gives yet another way to construct ideals of rings:

\begin{proposition}
\label{prop.ideals.invimg}Let $R$ and $S$ be two rings. Let $f:R\rightarrow S$
be a ring morphism. Let $K$ be an ideal of $S$.

Then,
\[
f^{-1}\left(  K\right)  :=\left\{  r\in R\mid f\left(  r\right)  \in
K\right\}
\]
is an ideal of $R$ that satisfies $\operatorname{Ker}f\subseteq f^{-1}\left(
K\right)  $.
\end{proposition}

\begin{exercise}
\label{exe.21hw3.1a}Prove Proposition \ref{prop.ideals.invimg}.
\end{exercise}

Here is an example:

\begin{itemize}
\item Let $S$ be the ring $\mathbb{Z}$. Let $R$ be the ring $\mathbb{Z}%
^{2\leq2}$ of all upper-triangular $2\times2$-matrices $\left(
\begin{array}
[c]{cc}%
a & b\\
0 & d
\end{array}
\right)  $ with integer entries $a,b,d\in\mathbb{Z}$. Let $f:R\rightarrow S$
be the map that sends each such matrix $\left(
\begin{array}
[c]{cc}%
a & b\\
0 & d
\end{array}
\right)  $ to its entry $d$. It is easy to see that this map $f$ is a ring
morphism. Let $K$ be the ideal of $S$ consisting of all even integers (that
is, $K=2\mathbb{Z}$). Thus, $f^{-1}\left(  K\right)  $ (as defined in
Proposition \ref{prop.ideals.invimg}) is the set of all upper-triangular
$2\times2$-matrices $\left(
\begin{array}
[c]{cc}%
a & b\\
0 & d
\end{array}
\right)  $ with integer entries $a,b,d\in\mathbb{Z}$ such that $d$ is even.
Proposition \ref{prop.ideals.invimg} shows that this set $f^{-1}\left(
K\right)  $ is an ideal of $R$ (which the reader can easily check).
\end{itemize}

\subsubsection{The fourth isomorphism theorem for rings}

What is the relation between the ideals of a quotient ring $R/I$ and the
ideals of the original ring $R$ ? More generally, what is the relation between
the ideals of two rings $S$ and $R$ when there is a surjective ring morphism
$f:R\rightarrow S$ ? (This is \textquotedblleft more general\textquotedblright%
\ since there is always a surjective morphism $\pi:R\rightarrow R/I$ from a
ring $R$ to a quotient ring $R/I$.) The \textbf{fourth isomorphism theorem for
rings} answers this question:

\begin{theorem}
[Fourth isomorphism theorem for rings]\label{thm.rings.4it}Let $R$ and $S$ be
two rings. Let $f:R\rightarrow S$ be a \textbf{surjective} ring morphism. Then:

\begin{enumerate}
\item[\textbf{(a)}] If $J$ is an ideal of $R$, then $f\left(  J\right)
:=\left\{  f\left(  j\right)  \mid j\in J\right\}  $ is an ideal of $S$.

\item[\textbf{(b)}] The maps
\begin{align*}
\left\{  \text{ideals $J$ of $R$ satisfying $\operatorname{Ker}f\subseteq J$%
}\right\}   &  \rightarrow\left\{  \text{ideals of $S$}\right\}  ,\\
J  &  \mapsto f\left(  J\right)
\end{align*}
and%
\begin{align*}
\left\{  \text{ideals of $S$}\right\}   &  \rightarrow\left\{  \text{ideals
$J$ of $R$ satisfying $\operatorname{Ker}f\subseteq J$}\right\}  ,\\
K  &  \mapsto f^{-1}\left(  K\right)
\end{align*}
(where $f^{-1}\left(  K\right)  $ is defined as in Proposition
\ref{prop.ideals.invimg}) are mutually inverse.
\end{enumerate}
\end{theorem}

\begin{exercise}
\label{exe.21hw3.1b}Prove Theorem \ref{thm.rings.4it}.
\end{exercise}

We note that Theorem \ref{thm.rings.4it} \textbf{(a)} becomes false if we drop
the \textquotedblleft$f$ is surjective\textquotedblright\ assumption.

\subsubsection{Prime and maximal ideals}

The following definition is rather important for the deeper study of ideals in
commutative rings (and, by extension, for algebraic geometry). We will only
touch on it briefly in this little subsection.

\begin{definition}
\label{def.ideals.prime-max}Let $R$ be a commutative ring. Let $I$ be an ideal
of $R$.

\begin{enumerate}
\item[\textbf{(a)}] The ideal $I$ is said to be \textbf{proper} if it is a
proper subset of $R$ (that is, $I\neq R$).

\item[\textbf{(b)}] The ideal $I$ is said to be \textbf{prime} if it is proper
and has the following property: If $a,b\in R$ satisfy $ab\in I$, then $a\in I$
or $b\in I$.

\item[\textbf{(c)}] The ideal $I$ is said to be \textbf{maximal} if it is
proper and the only ideals $J$ of $R$ satisfying $I\subseteq J\subseteq R$ are
$I$ and $R$.
\end{enumerate}
\end{definition}

We note that a principal ideal $aR$ of a commutative ring $R$ (with $a\in R$
nonzero) is prime if and only if $a$ is a prime element of $R$. Thus, the
notion of prime ideals generalizes the notion of prime elements (and,
ultimately, that of prime numbers).

\begin{exercise}
\label{exe.21hw3.2}Let $R$ be a commutative ring. Let $I$ be an ideal of $R$.
Prove the following:

\begin{enumerate}
\item[\textbf{(a)}] The ideal $I$ is prime if and only if the quotient ring
$R/I$ is an integral domain.

\item[\textbf{(b)}] The ideal $I$ is maximal if and only if the quotient ring
$R/I$ is a field.

\item[\textbf{(c)}] Any maximal ideal $I$ of $R$ is prime.
\end{enumerate}

[\textbf{Hint:} Theorem \ref{thm.rings.4it} and Exercise
\ref{exe.ideal.field-has-2} are helpful for part \textbf{(b)}.]
\end{exercise}

\subsection{\label{sec.rings.reg}Regular elements}

As we already mentioned in Section \ref{sec.rings.idomains}, the notion of a
zero divisor is somewhat inconvenient and brittle. It is much better to define
the (essentially) opposite concept:

\begin{definition}
\label{def.rings.reg.reg}Let $R$ be a commutative ring. An element $a\in R$ is
said to be \textbf{regular} (in $R$) if every $b\in R$ satisfying $ab=0$
satisfies $b=0$.
\end{definition}

As we said, regularity (at least for nonzero elements) is just the negation of
being a zero divisor:

\begin{proposition}
\label{prop.rings.reg.equiv}Let $R$ be a commutative ring. Let $a\in R$ be a
nonzero element. Then, $a$ is regular if and only if $a$ is not a zero divisor.
\end{proposition}

\begin{proof}
The negation of the statement \textquotedblleft every $b\in R$ satisfying
$ab=0$ satisfies $b=0$\textquotedblright\ is \textquotedblleft there is a
nonzero $b\in R$ such that $ab=0$\textquotedblright. Thus, the proposition follows.
\end{proof}

The zero $0$ of a commutative ring $R$ is regular if and only if $R$ is
trivial (because if $R$ is nontrivial, then $0\cdot1=0$ but $1\neq0$). Thus,
the regular elements of an integral domain are simply its nonzero elements.

\begin{exercise}
\label{exe.ring.reg.finite}Let $R$ be a \textbf{finite} commutative ring.
Prove that the regular elements of $R$ are precisely the units of $R$.

(Thus, for example, if $n$ is a positive integer, then the regular elements of
$\mathbb{Z}/n$ are the units of $\mathbb{Z}/n$; we have described the latter
in Proposition \ref{prop.ringunits.Z/n}.)
\end{exercise}

We now step to some general properties of regular elements:

\begin{proposition}
\label{prop.rings.reg.123}Let $R$ be a commutative ring. Then:

\begin{enumerate}
\item[\textbf{(a)}] Any unit of $R$ is regular.

\item[\textbf{(b)}] If $u$ and $v$ are two regular elements of $R$, then $uv$
is regular.

\item[\textbf{(c)}] Any product of regular elements of $R$ is again regular.

\item[\textbf{(d)}] Any divisor of a regular element of $R$ is again regular.

\item[\textbf{(e)}] Let $a\in R$. Then, $a$ is regular if and only if the map%
\begin{align*}
R  &  \rightarrow R,\\
x  &  \mapsto ax
\end{align*}
is injective.
\end{enumerate}
\end{proposition}

\begin{proof}
\textbf{(a)} Let $u$ be a unit of $R$. We must show that $u$ is regular.

The element $u$ is a unit, i.e., has an inverse $u^{-1}$. Thus, every $b\in R$
satisfying $ub=0$ satisfies $b=0$ (because $b=u^{-1}\underbrace{ub}_{=0}=0$).
In other words, $u$ is regular. Thus, part \textbf{(a)} is proved. \medskip

\textbf{(b)} Let $u$ and $v$ be two regular elements of $R$. We must show that
$uv$ is regular.

Let $b\in R$ be such that $\left(  uv\right)  b=0$. Then, $u\left(  vb\right)
=\left(  uv\right)  b=0$, whence $vb=0$ (since $u$ is regular). Therefore,
$b=0$ (since $v$ is regular).

Forget that we fixed $b$. We thus have shown that every $b\in R$ satisfying
$\left(  uv\right)  b=0$ satisfies $b=0$. In other words, $uv$ is regular.
This proves part \textbf{(b)}. \medskip

\textbf{(c)} This follows by induction on the number of factors. (The base
case uses part \textbf{(a)}; the induction step uses part \textbf{(b)}.)
\medskip

\textbf{(d)} Let $d$ be a divisor of a regular element $r\in R$. We must show
that $d$ is regular.

We have $r=dx$ for some $x\in R$ (since $d$ is a divisor of $r$). Consider
this $x$. Now, every $b\in R$ satisfying $db=0$ must satisfy $\underbrace{r}%
_{=dx=xd}b=x\underbrace{db}_{=0}=0$ and therefore $b=0$ (since $r$ is
regular). This shows that $d$ is regular. Thus, part \textbf{(d)} is proved.
\medskip

\textbf{(e)} Let $L_{a}$ denote the map
\begin{align*}
R  &  \rightarrow R,\\
x  &  \mapsto ax.
\end{align*}
We must show that $a$ is regular if and only if this map $L_{a}$ is injective.
The \textquotedblleft if\textquotedblright\ direction is easy: If $L_{a}$ is
injective, then every $b\in R$ satisfying $ab=0$ must satisfy $L_{a}\left(
b\right)  =ab=0=a0=L_{a}\left(  0\right)  $ and therefore $b=0$ (since $L_{a}$
is injective), which means that $a$ is regular. To prove the \textquotedblleft
only if\textquotedblright\ direction, we assume that $a$ is regular. Then, for
any two elements $x,y\in R$ satisfying $L_{a}\left(  x\right)  =L_{a}\left(
y\right)  $, we have $a\left(  x-y\right)  =\underbrace{ax}_{=L_{a}\left(
x\right)  }-\underbrace{ay}_{=L_{a}\left(  y\right)  }=L_{a}\left(  x\right)
-L_{a}\left(  y\right)  =0$ (since $L_{a}\left(  x\right)  =L_{a}\left(
y\right)  $) and therefore $x-y=0$ (since $a$ is regular), that is, $x=y$. In
other words, $L_{a}$ is injective. This proves the \textquotedblleft only
if\textquotedblright\ direction. Thus, the proof of part \textbf{(e)} is complete.
\end{proof}

Proposition \ref{prop.rings.reg.123} \textbf{(e)} reveals a different way of
thinking about regularity: It shows that an element $a$ of a commutative ring
$R$ is regular if and only if it can be cancelled from equalities of products
(i.e., if any equality of the form $ax=ay$ with $x,y\in R$ yields $x=y$).

\begin{exercise}
Let $a$ and $b$ be two regular elements of a commutative ring $R$. Prove that
any lcm of $a$ and $b$ is again regular. (Of course, this only says something
if such an lcm exists.)
\end{exercise}

\begin{exercise}
Let $R_{1},R_{2},\ldots,R_{k}$ be $k$ commutative rings. Prove that an element
$\left(  r_{1},r_{2},\ldots,r_{k}\right)  \in R_{1}\times R_{2}\times
\cdots\times R_{k}$ is regular if and only if each $r_{i}$ is regular (in its
respective $R_{i}$).
\end{exercise}

More about regular elements can be found in \cite{regpol}.

Unlike units, regular elements cannot always be divided by. However, when they
can, the quotient is always unique:

\begin{definition}
\label{def.rings.reg.div}Let $R$ be a commutative ring. Let $a,b\in R$ be such
that $b$ is regular. Then, there is \textbf{at most one} element $c\in R$
satisfying $bc=a$ (because Proposition \ref{prop.rings.reg.123} \textbf{(e)}
(applied to $b$ instead of $a$) shows that the map $R\rightarrow R,\ x\mapsto
bx$ is injective, so that $bc_{1}=bc_{2}$ would imply $c_{1}=c_{2}$). If such
an element $c$ exists, then it will be denoted by $\dfrac{a}{b}$. (Its
uniqueness makes this notation unambiguous.)
\end{definition}

This generalizes Definition \ref{def.ringunits.div} (because if $b\in
R^{\times}$, then $b$ is regular by Proposition \ref{prop.rings.reg.123}
\textbf{(a)}, and the element $c\in R$ satisfying $bc=a$ automatically exists
and equals $ab^{-1}$).

\begin{exercise}
Prove the following generalizations of Proposition
\ref{prop.ringunits.div.props} \textbf{(a)} and \textbf{(b)}:

\begin{enumerate}
\item[\textbf{(a)}] Let $R$ be a commutative ring. Let $a,b,c,d\in R$ be such
that $b$ and $d$ are regular and such that $\dfrac{a}{b}$ and $\dfrac{c}{d}$
exist. Then, $\dfrac{ad+bc}{bd}$ and $\dfrac{ac}{bd}$ also exist, and we have%
\[
\dfrac{a}{b}+\dfrac{c}{d}=\dfrac{ad+bc}{bd}%
\]
and%
\[
\dfrac{a}{b}\cdot\dfrac{c}{d}=\dfrac{ac}{bd}.
\]

\item[\textbf{(b)}] If $c$ is regular as well, then%
\[
\dfrac{a}{b}/\dfrac{c}{d}=\dfrac{ad}{bc},
\]
in the sense that if one side of this equality exists, then so does the other,
and the two sides are equal.
\end{enumerate}
\end{exercise}

\subsection{\label{sec.rings.gcdseq}Strong gcd-sequences}

We turn to a more exotic topic (which, in the general form we present it, is
likely novel): the notion of \textbf{strong gcd-sequences}. It will later be
useful to us in the construction of cyclotomic polynomials; but it also
illustrates some techniques for working with ideals in commutative rings
(particularly the use of ideal arithmetic; we recommend the reader to review
Section \ref{sec.rings.ideal-arith} before reading on).

\subsubsection{Examples}

We begin with two classical exercises on elementary number theory
(\cite[Exercise 3.4.1 \textbf{(b)}]{mps} and \cite[Exercise 3.7.2]{mps}, respectively):

\begin{exercise}
\label{exe.gcdseq.an-1}Let $u\in\mathbb{Z}$. Prove that%
\[
\gcd\left(  u^{n}-1,\ u^{m}-1\right)  =\left\vert u^{\gcd\left(  n,m\right)
}-1\right\vert \ \ \ \ \ \ \ \ \ \ \text{for all }n,m\geq0.
\]

\end{exercise}

\begin{exercise}
\label{exe.gcdseq.fib}Consider the Fibonacci sequence $\left(  f_{0}%
,f_{1},f_{2},\ldots\right)  $ from Definition \ref{def.fibonacci.fib}. Prove
that%
\[
\gcd\left(  f_{n},f_{m}\right)  =f_{\gcd\left(  n,m\right)  }%
\ \ \ \ \ \ \ \ \ \ \text{for all }n,m\geq0.
\]

\end{exercise}

These two exercises give two instances of sequences $\left(  a_{1},a_{2}%
,a_{3},\ldots\right)  $ of integers satisfying $\gcd\left(  a_{n}%
,a_{m}\right)  =\left\vert a_{\gcd\left(  n,m\right)  }\right\vert $ for all
$n,m\geq1$ (or even for all $n,m\geq0$ when $a_{0}=0$)\ \ \ \ \footnote{We can
restrict ourselves to the case $n,m\geq1$ because under the condition
$a_{0}=0$ (which is satisfied in both exercises), the statement $\gcd\left(
a_{n},a_{m}\right)  =\left\vert a_{\gcd\left(  n,m\right)  }\right\vert $ is
automatically satisfied when $n=0$ or $m=0$.}. In 1936, Morgan Ward began the
general study of such sequences (\cite{Ward36}, \cite{Ward39}). With a minor
tweak, they can be defined not just for integers, but for elements of any
commutative ring: A sequence $\left(  a_{1},a_{2},a_{3},\ldots\right)  $ of
elements of a commutative ring $R$ is said to be a \textbf{gcd-sequence} if
$a_{\gcd\left(  n,m\right)  }$ is a gcd of $a_{n}$ and $a_{m}$ for all
$n,m\geq1$. Thus, Exercise \ref{exe.gcdseq.an-1} (restricted to $n,m\geq1$)
says that for every $u\in\mathbb{Z}$, the sequence $\left(  u^{1}%
-1,\ u^{2}-1,\ u^{3}-1,\ \ldots\right)  $ of integers is a gcd-sequence.
Likewise, Exercise \ref{exe.gcdseq.fib} (restricted to $n,m\geq1$) says that
$\left(  f_{1},f_{2},f_{3},\ldots\right)  $ is a gcd-sequence. Of course, any
sequence of the form $\left(  1^{k},2^{k},3^{k},\ldots\right)  $ for
$k\in\mathbb{N}$ is a gcd-sequence as well (since $\gcd\left(  n^{k}%
,m^{k}\right)  =\left(  \gcd\left(  n,m\right)  \right)  ^{k}$ for all
$n,m\geq1$). More interesting examples of gcd-sequences can be found in rings
such as $\mathbb{Z}\left[  x\right]  $; one such example will be the key to
our construction of cyclotomic polynomials later on.

As shown in \cite{BFLS13} and \cite{Nowick15}, gcd-sequences have good
properties over UFDs, though not that much can be said over arbitrary
commutative rings. Rather than present these results here, we shall study a
variant of gcd-sequences that is well-behaved over arbitrary commutative
rings, although it is more restrictive in other ways. We call this variant
\textbf{strong gcd-sequences}.

\subsubsection{Definition and first properties}

These sequences are defined as follows:

\begin{definition}
\label{def.gcdseq.gcdseq}Let $R$ be a commutative ring. A sequence $\left(
a_{1},a_{2},a_{3},\ldots\right)  $ of elements of $R$ will be called a
\textbf{strong gcd-sequence} (in $R$) if all integers $n,m\geq1$ satisfy%
\begin{equation}
a_{n}R+a_{m}R=a_{\gcd\left(  n,m\right)  }R. \label{eq.def.gcdseq.gcdseq.a+a}%
\end{equation}

\end{definition}

If $R$ is a PID or the ring $\mathbb{Z}$ in particular, this condition boils
down to the above notions of gcd-sequences:

\begin{proposition}
\label{prop.gcdseq.equiv}\ \ 

\begin{enumerate}
\item[\textbf{(a)}] Let $R$ be a PID. A sequence $\left(  a_{1},a_{2}%
,a_{3},\ldots\right)  $ of elements of $R$ is a strong gcd-sequence if and
only if for all integers $n,m\geq1$, the element $a_{\gcd\left(  n,m\right)
}$ is a gcd of $a_{n}$ and $a_{m}$.

\item[\textbf{(b)}] A sequence $\left(  a_{1},a_{2},a_{3},\ldots\right)  $ of
integers is a strong gcd-sequence (in $\mathbb{Z}$) if and only if for all
integers $n,m\geq1$, we have $\gcd\left(  a_{n},a_{m}\right)  =\left\vert
a_{\gcd\left(  n,m\right)  }\right\vert $.
\end{enumerate}
\end{proposition}

\begin{proof}
\textbf{(a)} Let $\left(  a_{1},a_{2},a_{3},\ldots\right)  $ be a sequence of
elements of $R$. Then, Definition \ref{def.gcdseq.gcdseq} says that $\left(
a_{1},a_{2},a_{3},\ldots\right)  $ is a strong gcd-sequence if and only if for
all integers $n,m\geq1$, the equality $a_{n}R+a_{m}R=a_{\gcd\left(
n,m\right)  }R$ holds. However, the latter equality is equivalent to the
statement that the element $a_{\gcd\left(  n,m\right)  }$ is a gcd of $a_{n}$
and $a_{m}$ (by Proposition \ref{prop.gcd-lcm.ideals.pid-a}, applied to
$a=a_{n}$, $b=a_{m}$ and $c=a_{\gcd\left(  n,m\right)  }$). Thus, we conclude
that $\left(  a_{1},a_{2},a_{3},\ldots\right)  $ is a strong gcd-sequence if
and only if for all integers $n,m\geq1$, the element $a_{\gcd\left(
n,m\right)  }$ is a gcd of $a_{n}$ and $a_{m}$. This proves part \textbf{(a)}.
\medskip

\textbf{(b)} Let $\left(  a_{1},a_{2},a_{3},\ldots\right)  $ be a sequence of
integers. Then, part \textbf{(a)} (applied to $R=\mathbb{Z}$) yields that
$\left(  a_{1},a_{2},a_{3},\ldots\right)  $ is a strong gcd-sequence if and
only if for all integers $n,m\geq1$, the element $a_{\gcd\left(  n,m\right)
}$ is a gcd of $a_{n}$ and $a_{m}$ in $\mathbb{Z}$. Hence, in order to prove
part \textbf{(b)}, it remains to show that (for any given $n,m\geq1$) the
condition \textquotedblleft the element $a_{\gcd\left(  n,m\right)  }$ is a
gcd of $a_{n}$ and $a_{m}$ in $\mathbb{Z}$\textquotedblright\ is equivalent to
\textquotedblleft$\gcd\left(  a_{n},a_{m}\right)  =\left\vert a_{\gcd\left(
n,m\right)  }\right\vert $\textquotedblright. But this is easy: If $n,m\geq1$
are arbitrary, then the gcds of $a_{n}$ and $a_{m}$ in $\mathbb{Z}$ are
$\gcd\left(  a_{n},a_{m}\right)  $ and $-\gcd\left(  a_{n},a_{m}\right)  $ (by
Proposition \ref{prop.gcd-lcm.ZZ} \textbf{(a)}, applied to $a=a_{n}$ and
$b=a_{m}$), and thus we obtain the chain of equivalences%
\begin{align*}
&  \ \left(  \text{the element }a_{\gcd\left(  n,m\right)  }\text{ is a gcd of
}a_{n}\text{ and }a_{m}\text{ in }\mathbb{Z}\right) \\
&  \Longleftrightarrow\ \left(  a_{\gcd\left(  n,m\right)  }\text{ equals
}\gcd\left(  a_{n},a_{m}\right)  \text{ or }-\gcd\left(  a_{n},a_{m}\right)
\right) \\
&  \Longleftrightarrow\ \left(  \left\vert a_{\gcd\left(  n,m\right)
}\right\vert =\left\vert \gcd\left(  a_{n},a_{m}\right)  \right\vert \right)
\\
&  \Longleftrightarrow\ \left(  \left\vert a_{\gcd\left(  n,m\right)
}\right\vert =\gcd\left(  a_{n},a_{m}\right)  \right)
\ \ \ \ \ \ \ \ \ \ \left(
\begin{array}
[c]{c}%
\text{since }\left\vert \gcd\left(  a_{n},a_{m}\right)  \right\vert
=\gcd\left(  a_{n},a_{m}\right) \\
\text{(because }\gcd\left(  a_{n},a_{m}\right)  \geq0\text{)}%
\end{array}
\right) \\
&  \Longleftrightarrow\ \left(  \gcd\left(  a_{n},a_{m}\right)  =\left\vert
a_{\gcd\left(  n,m\right)  }\right\vert \right)  .
\end{align*}
Thus, part \textbf{(b)} is proved.
\end{proof}

We shall next prove a sufficient criterion for a sequence to be a strong
gcd-sequence. This criterion is far from necessary (i.e., it does not
construct all strong gcd-sequences), but it still gives so much freedom that
both Exercise \ref{exe.gcdseq.an-1} and Exercise \ref{exe.gcdseq.fib} (and
several other important examples) are covered by it:

\begin{proposition}
\label{prop.gcdseq.modcrit}Let $R$ be a commutative ring. Let $\left(
a_{1},a_{2},a_{3},\ldots\right)  $ be a sequence of elements of $R$. For any
$m,n\geq1$, let $b_{m,n}$ and $c_{m,n}$ be elements of $R$ such that%
\begin{equation}
a_{m}R+b_{m,n}R+c_{m,n}R=R. \label{eq.prop.gcdseq.modcrit.ass1}%
\end{equation}
Assume furthermore that each $p,q\geq1$ satisfy%
\begin{equation}
a_{p+q}=a_{p}b_{q,p}+a_{q}c_{p,q}\ \ \ \ \ \ \ \ \ \ \text{in }R.
\label{eq.prop.gcdseq.modcrit.ass2}%
\end{equation}
Then, $\left(  a_{1},a_{2},a_{3},\ldots\right)  $ is a strong gcd-sequence.
\end{proposition}

Before reading the following proof, the reader is advised to take a peek below
it, where it is applied to recover Exercise \ref{exe.gcdseq.an-1} and Exercise
\ref{exe.gcdseq.fib}. (Note that the $b_{m,n}$ and $c_{m,n}$ in these
applications only depend on $m$, not on $n$; we are leaving countably many
degrees of freedom unused!)

\begin{proof}
[Proof of Proposition \ref{prop.gcdseq.modcrit}.]We must prove that all
integers $n,m\geq1$ satisfy%
\begin{equation}
a_{n}R+a_{m}R=a_{\gcd\left(  n,m\right)  }R.
\label{pf.prop.gcdseq.modcrit.goal}%
\end{equation}
We shall prove this by strong induction on $n+m$:

Fix a positive integer $N$. Assume that (\ref{pf.prop.gcdseq.modcrit.goal})
holds whenever $n+m<N$. We must prove that (\ref{pf.prop.gcdseq.modcrit.goal})
holds when $n+m=N$.

So let $n,m\geq1$ be two integers such that $n+m=N$. We must prove
(\ref{pf.prop.gcdseq.modcrit.goal}).

Our situation is symmetric in $n$ and $m$ (since $\gcd\left(  n,m\right)
=\gcd\left(  m,n\right)  $ and $a_{n}R+a_{m}R=a_{m}R+a_{n}R$ and $n+m=m+n$).
Thus, we can swap $n$ with $m$ at will. Hence, we can WLOG assume that $n\geq
m$ (otherwise, we swap $n$ with $m$). Assume this.

We must prove (\ref{pf.prop.gcdseq.modcrit.goal}). If $n=m$, then this is
pretty obvious (indeed, in this case, we have $\gcd\left(  n,m\right)
=\gcd\left(  m,m\right)  =m$ and $a_{n}R=a_{m}R$, so that the claim
(\ref{pf.prop.gcdseq.modcrit.goal}) rewrites as $a_{m}R+a_{m}R=a_{m}R$, which
is clearly true since $a_{m}R$ is an ideal). So we WLOG assume that $n\neq m$.
Therefore, $n>m$ (since $n\geq m$). Hence, $n-m\geq1$ and $\left(  n-m\right)
+m=n<n+m$ (since $m\geq1$). In other words, $\left(  n-m\right)  +m<N$ (since
$n+m=N$). Therefore, by our induction hypothesis, we can apply
(\ref{pf.prop.gcdseq.modcrit.goal}) to $n-m$ instead of $n$. As a result, we
obtain%
\begin{equation}
a_{n-m}R+a_{m}R=a_{\gcd\left(  n-m,m\right)  }R=a_{\gcd\left(  n,m\right)  }R
\label{pf.prop.gcdseq.modcrit.IH}%
\end{equation}
(since one of the basic properties of gcds of integers\footnote{specifically,
the fact that the gcd of two integers does not change if we add a multiple of
one to the other} yields $\gcd\left(  n-m,m\right)  =\gcd\left(  n,m\right)  $).

We have $1_{R}\in R=a_{m}R+b_{m,n-m}R+c_{m,n-m}R$ (since
(\ref{eq.prop.gcdseq.modcrit.ass1}) (applied to $n-m$ instead of $n$) yields
$a_{m}R+b_{m,n-m}R+c_{m,n-m}R=R$). That is,
\begin{equation}
1_{R}=a_{m}x+b_{m,n-m}y+c_{m,n-m}z \label{pf.prop.gcdseq.modcrit.1=}%
\end{equation}
for some $x,y,z\in R$. Consider these $x,y,z$.

Let $S$ be the quotient ring $R/a_{m}R$. We shall write $\overline{r}$ for the
residue class of any $r\in R$ in $S$. (That is, $\overline{r}=r+a_{m}R$.)
Projecting the equality (\ref{pf.prop.gcdseq.modcrit.1=}) onto $R/a_{m}R$, we
obtain%
\begin{align}
1_{S}  &  =\overline{a_{m}x+b_{m,n-m}y+c_{m,n-m}z}=\underbrace{\overline
{a_{m}}}_{=0}\overline{x}+\overline{b_{m,n-m}}\,\overline{y}+\overline
{c_{m,n-m}}\,\overline{z}\nonumber\\
&  =\overline{b_{m,n-m}}\,\overline{y}+\overline{c_{m,n-m}}\,\overline{z}.
\label{pf.prop.gcdseq.modcrit.1=S}%
\end{align}

The assumption (\ref{eq.prop.gcdseq.modcrit.ass2}) (applied to $p=n-m$ and
$q=m$) yields%
\[
a_{n}=a_{n-m}b_{m,n-m}+a_{m}c_{n-m,m}.
\]
Projecting this onto $R/a_{m}R$, we obtain%
\begin{equation}
\overline{a_{n}}=\overline{a_{n-m}}\,\overline{b_{m,n-m}}%
+\underbrace{\overline{a_{m}}}_{=0}\overline{c_{n-m,m}}=\overline{a_{n-m}%
}\,\overline{b_{m,n-m}}. \label{pf.prop.gcdseq.modcrit.an1}%
\end{equation}

Moreover, the assumption (\ref{eq.prop.gcdseq.modcrit.ass2}) (applied to $p=m$
and $q=n-m$) yields%
\[
a_{n}=a_{m}b_{n-m,m}+a_{n-m}c_{m,n-m}.
\]
Projecting this onto $R/a_{m}R$, we obtain%
\begin{align}
\overline{a_{n}}  &  =\underbrace{\overline{a_{m}}}_{=0}\overline{b_{n-m,m}%
}+\overline{a_{n-m}}\,\overline{c_{m,n-m}}=\overline{a_{n-m}}\,\overline
{c_{m,n-m}}\label{pf.prop.gcdseq.modcrit.an2}\\
&  \in\overline{a_{n-m}}S.\nonumber
\end{align}
That is, $\overline{a_{n-m}}\mid\overline{a_{n}}$ in $S$. Therefore,%
\begin{equation}
\overline{a_{n}}S\subseteq\overline{a_{n-m}}S.
\label{pf.prop.gcdseq.modcrit.sub1}%
\end{equation}

Multiplying the equality (\ref{pf.prop.gcdseq.modcrit.1=S}) by $\overline
{a_{n-m}}$, we find%
\begin{align*}
\overline{a_{n-m}}  &  =\overline{a_{n-m}}\left(  \overline{b_{m,n-m}%
}\,\overline{y}+\overline{c_{m,n-m}}\,\overline{z}\right)
=\underbrace{\overline{a_{n-m}}\,\overline{b_{m,n-m}}}_{\substack{=\overline
{a_{n}}\\\text{(by (\ref{pf.prop.gcdseq.modcrit.an1}))}}}\overline
{y}+\underbrace{\overline{a_{n-m}}\,\overline{c_{m,n-m}}}%
_{\substack{=\overline{a_{n}}\\\text{(by (\ref{pf.prop.gcdseq.modcrit.an2}))}%
}}\overline{z}\\
&  =\overline{a_{n}}\,\overline{y}+\overline{a_{n}}\,\overline{z}%
=\overline{a_{n}}\left(  \overline{y}+\overline{z}\right)  \in\overline{a_{n}%
}S.
\end{align*}
That is, $\overline{a_{n}}\mid\overline{a_{n-m}}$ in $S$. Therefore,%
\[
\overline{a_{n-m}}S\subseteq\overline{a_{n}}S.
\]
Combining this with (\ref{pf.prop.gcdseq.modcrit.sub1}), we obtain%
\begin{equation}
\overline{a_{n-m}}S=\overline{a_{n}}S. \label{pf.prop.gcdseq.modcrit.eq}%
\end{equation}

However, let $\pi:R\rightarrow R/a_{m}R$ be the canonical projection. Then,
for any $u\in R$, we have%
\[
\pi^{-1}\left(  \overline{u}S\right)  =uR+a_{m}R
\]
(because for any $r\in R$, the following chain of equivalences holds:%
\begin{align*}
\left(  r\in\pi^{-1}\left(  \overline{u}S\right)  \right)  \  &
\Longleftrightarrow\ \left(  \pi\left(  r\right)  \in\overline{u}S\right) \\
&  \Longleftrightarrow\ \left(  \pi\left(  r\right)  =\overline{u}%
\,\overline{s}\text{ for some }s\in R\right) \\
&  \ \ \ \ \ \ \ \ \ \ \ \ \ \ \ \ \ \ \ \ \left(  \text{since }S=\left\{
\overline{s}\ \mid\ s\in R\right\}  \right) \\
&  \Longleftrightarrow\ \left(  r+a_{m}R=us+a_{m}R\text{ for some }s\in
R\right) \\
&  \ \ \ \ \ \ \ \ \ \ \ \ \ \ \ \ \ \ \ \ \left(
\begin{array}
[c]{c}%
\text{since }\pi\left(  r\right)  =\overline{r}=r+a_{m}R\\
\text{and }\overline{u}\,\overline{s}=\overline{us}=us+a_{m}R
\end{array}
\right) \\
&  \Longleftrightarrow\ \left(  r\in us+a_{m}R\text{ for some }s\in R\right)
\\
&  \Longleftrightarrow\ \left(  r\in uR+a_{m}R\right)
\end{align*}
). Thus, applying $\pi^{-1}$ to both sides of the equality
(\ref{pf.prop.gcdseq.modcrit.eq}), we obtain%
\[
a_{n-m}R+a_{m}R=a_{n}R+a_{m}R.
\]
Thus,%
\[
a_{n}R+a_{m}R=a_{n-m}R+a_{m}R=a_{\gcd\left(  n,m\right)  }R
\]
(by (\ref{pf.prop.gcdseq.modcrit.IH})). This completes the induction step.

Thus, (\ref{pf.prop.gcdseq.modcrit.goal}) is proved for all $n,m\geq1$. In
other words, $\left(  a_{1},a_{2},a_{3},\ldots\right)  $ is a strong
gcd-sequence. This proves Proposition \ref{prop.gcdseq.modcrit}.
\end{proof}

Both Exercise \ref{exe.gcdseq.an-1} and Exercise \ref{exe.gcdseq.fib} follow
easily from Proposition \ref{prop.gcdseq.modcrit}:

\begin{fineprint}

\begin{itemize}
\item To solve Exercise \ref{exe.gcdseq.an-1}, we need to show that the
sequence $\left(  u^{1}-1,\ u^{2}-1,\ u^{3}-1,\ \ldots\right)  $ is a strong
gcd-sequence in $\mathbb{Z}$ (because then Proposition \ref{prop.gcdseq.equiv}
\textbf{(b)} will show that $\gcd\left(  u^{n}-1,\ u^{m}-1\right)  =\left\vert
u^{\gcd\left(  n,m\right)  }-1\right\vert $ for all $n,m\geq1$; and it is easy
to verify the same equality for $n=0$ or $m=0$). This follows from Proposition
\ref{prop.gcdseq.modcrit}, applied to $a_{m}=u^{m}-1$ and $b_{m,n}=u^{m}$ and
$c_{m,n}=1$; indeed, these choices of $a_{m}$, $b_{m,n}$ and $c_{m,n}$ satisfy
(\ref{eq.prop.gcdseq.modcrit.ass1}) (since $c_{m,n}=1$ entails $c_{m,n}R=1R=R$
already) and (\ref{eq.prop.gcdseq.modcrit.ass2}) (by an easy computation).

\item Similarly, we can solve Exercise \ref{exe.gcdseq.fib} by applying
Proposition \ref{prop.gcdseq.modcrit} to $a_{m}=f_{m}$ and $b_{m,n}=f_{m+1}$
and $c_{m,n}=f_{m-1}$. Here, the proof of (\ref{eq.prop.gcdseq.modcrit.ass1})
is a bit harder (it suffices to show that $\gcd\left(  f_{m},f_{m+1}\right)
=1$ for each $m\geq0$; this is an easy induction on $m$), whereas the equality
(\ref{eq.prop.gcdseq.modcrit.ass2}) follows from Proposition
\ref{prop.fibonacci.add}.
\end{itemize}

Exercise \ref{exe.gcdseq.an-1} can be generalized further:

\begin{exercise}
\ \ 

\begin{enumerate}
\item[\textbf{(a)}] If $u$ and $v$ are two coprime integers, then prove that%
\[
\gcd\left(  u^{n}-v^{n},\ u^{m}-v^{m}\right)  =\left\vert u^{\gcd\left(
n,m\right)  }-v^{\gcd\left(  n,m\right)  }\right\vert
\ \ \ \ \ \ \ \ \ \ \text{for all }n,m\geq0.
\]

\item[\textbf{(b)}] More generally: Let $R$ be a commutative ring. Let $u$ and
$v$ be two elements of $R$ such that $uR+vR=R$. Prove that $\left(
u^{1}-v^{1},\ u^{2}-v^{2},\ u^{3}-v^{3},\ \ldots\right)  $ is a strong gcd-sequence.
\end{enumerate}
\end{exercise}
\end{fineprint}

We now ask ourselves what can be said about arbitrary strong gcd-sequences.
One property is obvious:

\begin{proposition}
\label{prop.gcdseq.div}Let $R$ be a commutative ring. Let $\left(  a_{1}%
,a_{2},a_{3},\ldots\right)  $ be a strong gcd-sequence in $R$.

\begin{enumerate}
\item[\textbf{(a)}] For any $p,q\geq1$, we have $a_{p+q}\in a_{p}R+a_{q}R$.

\item[\textbf{(b)}] Let $p,q\geq1$ be two integers such that $p\mid q$. Then,
$a_{p}\mid a_{q}$ in $R$.
\end{enumerate}
\end{proposition}

\begin{proof}
\textbf{(a)} Let $p,q\geq1$. Then, $\gcd\left(  p+q,q\right)  =\gcd\left(
p,q\right)  $ (by one of the basic properties of the gcd of two integers). But
$\left(  a_{1},a_{2},a_{3},\ldots\right)  $ is a strong gcd-sequence; hence,
(\ref{eq.def.gcdseq.gcdseq.a+a}) (applied to $n=p+q$ and $m=q$) yields
$a_{p+q}R+a_{q}R=a_{\gcd\left(  p+q,q\right)  }R=a_{\gcd\left(  p,q\right)
}R$ (since $\gcd\left(  p+q,q\right)  =\gcd\left(  p,q\right)  $). On the
other hand, (\ref{eq.def.gcdseq.gcdseq.a+a}) (applied to $n=p$ and $m=q$)
yields $a_{p}R+a_{q}R=a_{\gcd\left(  p,q\right)  }R$. Hence,%
\[
a_{p+q}=a_{p+q}1+a_{q}0\in a_{p+q}R+a_{q}R=a_{\gcd\left(  p,q\right)  }%
R=a_{p}R+a_{q}R.
\]
This proves Proposition \ref{prop.gcdseq.div} \textbf{(a)}. \medskip

\textbf{(b)} From $p\mid q$, we obtain $\gcd\left(  p,q\right)  =p$. But
$\left(  a_{1},a_{2},a_{3},\ldots\right)  $ is a strong gcd-sequence; hence,
(\ref{eq.def.gcdseq.gcdseq.a+a}) (applied to $n=p$ and $m=q$) yields
$a_{p}R+a_{q}R=a_{\gcd\left(  p,q\right)  }R=a_{p}R$ (since $\gcd\left(
p,q\right)  =p$). Hence, $a_{q}=a_{p}0+a_{q}1\in a_{p}R+a_{q}R=a_{p}R$. In
other words, $a_{p}\mid a_{q}$. This proves Proposition \ref{prop.gcdseq.div}
\textbf{(b)}.
\end{proof}

We come to a less obvious property of strong gcd-sequences, originally
observed for integers by Ward in \cite{Ward36}. It says that for any strong
gcd-sequence $a=\left(  a_{1},a_{2},a_{3},\ldots\right)  $, we can construct a
family of \textquotedblleft$a$-binomial coefficients\textquotedblright%
\ $\dbinom{n}{k}_{a}$:

\begin{exercise}
\label{exe.gcdseq.binom}Let $R$ be a commutative ring. Let $\left(
a_{1},a_{2},a_{3},\ldots\right)  $ be a strong gcd-sequence in $R$. Show that
for every nonnegative integers $n\geq k\geq0$, there exists an element
$\dbinom{n}{k}_{a}\in R$ satisfying%
\[
\dbinom{n}{k}_{a}\cdot a_{k}a_{k-1}\cdots a_{1}=a_{n}a_{n-1}\cdots a_{n-k+1}.
\]
(If $R=\mathbb{Z}$ and $a_{m}=m$ for all $m$, then $\dbinom{n}{k}_{a}$ is the
usual binomial coefficient $\dbinom{n}{k}$.) \medskip

[\textbf{Hint:} Define $\dbinom{n}{k}_{a}$ recursively by setting $\dbinom
{n}{0}_{a}=\dbinom{n}{n}_{a}=1$ and $\dbinom{n}{k}_{a}=x_{n,k}\dbinom
{n-1}{k-1}_{a}+y_{n,k}\dbinom{n-1}{k}_{a}$ for all $0<k<n$, where
$x_{n,k},y_{n,k}\in R$ are chosen such that $a_{n}=x_{n,k}a_{k}+y_{n,k}%
a_{n-k}$ (why do these exist?).]
\end{exercise}

\subsubsection{The factorization theorem}

We now shall state the main result of this section: a factorization theorem
for strong gcd-sequences. It is similar to \cite{BFLS13}, but requires more of
the sequence and less of the ring (although it boils down to the same claim
when $R$ is a PID).

\begin{theorem}
[factorization theorem for strong gcd-sequences]\label{thm.gcdseq.fac}Let $R$
be a commutative ring. Let $\left(  a_{1},a_{2},a_{3},\ldots\right)  $ be a
strong gcd-sequence in $R$. Then, there exists a sequence $\left(  b_{1}%
,b_{2},b_{3},\ldots\right)  $ of elements of $R$ such that each $n\geq1$
satisfies%
\begin{equation}
a_{n}=\prod_{d\mid n}b_{d}. \label{eq.thm.gcdseq.fac.a=prod}%
\end{equation}
Here and in the following, the symbol \textquotedblleft$\prod_{d\mid n}%
$\textquotedblright\ means a product over all positive divisors of $n$ (so
that, for example, $\prod_{d\mid6}b_{d}=b_{1}b_{2}b_{3}b_{6}$).
\end{theorem}

\begin{example}
\label{exa.gcdseq.fac.fib}Consider the Fibonacci sequence $\left(  f_{0}%
,f_{1},f_{2},\ldots\right)  $ from Definition \ref{def.fibonacci.fib}. From
Exercise \ref{exe.gcdseq.fib} and Proposition \ref{prop.gcdseq.equiv}
\textbf{(b)}, we conclude that $\left(  f_{1},f_{2},f_{3},\ldots\right)  $ is
a strong gcd-sequence in $\mathbb{Z}$. Hence, Theorem \ref{thm.gcdseq.fac}
(applied to $R=\mathbb{Z}$ and $a_{i}=f_{i}$) says that there exists a
sequence $\left(  b_{1},b_{2},b_{3},\ldots\right)  $ of integers such that
each $n\geq1$ satisfies%
\begin{equation}
f_{n}=\prod_{d\mid n}b_{d}. \label{eq.exa.gcdseq.fac.fib.f=prod}%
\end{equation}
The first entries of this sequence are
\[%
\begin{tabular}
[c]{|c||c|c|c|c|c|c|c|c|c|c|c|c|}\hline
$n$ & $1$ & $2$ & $3$ & $4$ & $5$ & $6$ & $7$ & $8$ & $9$ & $10$ & $11$ &
$12$\\\hline
$b_{n}$ & $1$ & $1$ & $2$ & $3$ & $5$ & $4$ & $13$ & $7$ & $17$ & $11$ & $89$
& $6$\\\hline
\end{tabular}
\ \ \ \ \ \ .
\]
For instance, for $n=10$, the equality (\ref{eq.exa.gcdseq.fac.fib.f=prod})
says $f_{10}=b_{1}b_{2}b_{5}b_{10}$, that is, $55=1\cdot1\cdot5\cdot11$.
\end{example}

The proof of Theorem \ref{thm.gcdseq.fac} will occupy most of the rest of this
section. While there are easier arguments when $R$ is a PID, we shall assume
nothing about $R$, and thus we will have to work essentially using our bare hands.

\subsubsection{A lemma on ideal arithmetic}

Before we get to the proof, we show a simple lemma about ideals, generalizing
Theorem \ref{thm.CRT-2-ids1}:

\begin{lemma}
\label{lem.gcdseq.abcR}Let $R$ be a commutative ring. Let $a\in R$, and let
$I$ and $J$ be two ideals of $R$. If $a\left(  I+J\right)  =aR$, then $aI\cap
aJ=aIJ$.
\end{lemma}

\begin{proof}
Assume that $a\left(  I+J\right)  =aR$. We must prove that $aI\cap aJ=aIJ$.

We have $IJ\subseteq I$ (since $I$ is an ideal) and $IJ\subseteq J$ (since $J$
is an ideal), so that $a\underbrace{IJ}_{\subseteq I}\subseteq aI$ and
$a\underbrace{IJ}_{\subseteq J}\subseteq aJ$. Therefore, $aIJ\subseteq aI\cap
aJ$.

Conversely, let $x\in aI\cap aJ$. Thus, $x\in aI$ and $x\in aJ$. In other
words, $x=ai$ for some $i\in I$ (since $x\in aI$), and $x=aj$ for some $j\in
J$ (since $x\in aJ$). Consider these $i$ and $j$. On the other hand, $a=a1\in
aR=a\left(  I+J\right)  $ (since $a\left(  I+J\right)  =aR$), so that we can
write $a$ as $a=a\left(  u+v\right)  $ for some $u\in I$ and $v\in J$.
Consider these $u$ and $v$. Now, $a=a\left(  u+v\right)  =\left(  u+v\right)
a$ (since $R$ is commutative), and we have%
\begin{align*}
x  &  =\underbrace{a}_{=\left(  u+v\right)  a}i=\left(  u+v\right)
\underbrace{ai}_{=x}=\left(  u+v\right)  x\\
&  =\underbrace{u}_{\in I}\underbrace{x}_{\in aJ}+\underbrace{v}_{\in
J}\underbrace{x}_{\in aI}\in\underbrace{IaJ}_{\substack{=aIJ\\\text{(since
}R\text{ is commutative)}}}+\underbrace{JaI}_{\substack{=aIJ\\\text{(since
}R\text{ is commutative)}}}\\
&  =aIJ+aIJ=aIJ\ \ \ \ \ \ \ \ \ \ \left(  \text{since }aIJ\text{ is an
ideal}\right)  .
\end{align*}

Forget that we fixed $x$. We thus have shown that $x\in aIJ$ for each $x\in
aI\cap aJ$. That is, $aI\cap aJ\subseteq aIJ$. Combining this with
$aIJ\subseteq aI\cap aJ$, we obtain $aI\cap aJ=aIJ$, qed.
\end{proof}

\begin{noncompile}

\begin{lemma}
\label{lem.gcdseq.abcdR}Let $R$ be a commutative ring. Let $a,b\in R$, and let
$I$ and $J$ be two ideals of $R$. If $a\left(  I+J\right)  =aR$ and
$aI\subseteq abR$, then $a\left(  I+bJ\right)  =abR$.
\end{lemma}

\begin{proof}
Assume that $a\left(  I+J\right)  =aR$ and $aI\subseteq abR$.

Clearly, $bI\subseteq I$ (since $I$ is an ideal). But $R$ is commutative, so
that%
\begin{align*}
abR  &  =b\underbrace{aR}_{=a\left(  I+J\right)  }=\underbrace{ba}%
_{=ab}\left(  I+J\right)  =ab\left(  I+J\right) \\
&  =a\underbrace{bI}_{\subseteq I}+\,abJ\subseteq aI+abJ=a\left(  I+bJ\right)
.
\end{align*}
Combining this with%
\[
a\left(  I+bJ\right)  =\underbrace{aI}_{\subseteq abR}+\,ab\underbrace{J}%
_{\subseteq R}\subseteq abR+abR=abR,
\]
we obtain $a\left(  I+bJ\right)  =abR$. This proves Lemma
\ref{lem.gcdseq.abcdR}.
\end{proof}
\end{noncompile}

\subsubsection{Proof of the factorization theorem}

To prove Theorem \ref{thm.gcdseq.fac}, we will construct the sequence $\left(
b_{1},b_{2},b_{3},\ldots\right)  $ by \textquotedblleft iterative
approximation\textquotedblright. First, we introduce some notations. For the
rest of this section, we fix a commutative ring $R$.

\begin{itemize}
\item If $a\in R^{\left\{  1,2,3,\ldots\right\}  }$ is a sequence of elements
of $R$, and if $n$ is a positive integer, then $a_{n}$ shall denote the $n$-th
entry of $a$. (Thus, $a=\left(  a_{1},a_{2},a_{3},\ldots\right)  $.)

\item If $a\in R^{\left\{  1,2,3,\ldots\right\}  }$ is a sequence of elements
of $R$, and if $S$ is a finite set of positive integers, then we set%
\[
a_{S}:=\prod_{d\in S}a_{d}.
\]
(For instance, $a_{\left\{  2,6,7\right\}  }=a_{2}a_{6}a_{7}$.)

\item If $n$ is a positive integer, then $\operatorname*{Div}n$ shall denote
the set of all positive divisors of $n$. (For instance, $\operatorname*{Div}%
6=\left\{  1,2,3,6\right\}  $.) Thus,%
\begin{equation}
a_{\operatorname*{Div}n}=\prod_{d\mid n}a_{d} \label{eq.gcdseq.aDiv}%
\end{equation}
for any sequence $a\in R^{\left\{  1,2,3,\ldots\right\}  }$ (since
\textquotedblleft$\prod_{d\mid n}$\textquotedblright\ means a product over all
positive divisors of $n$, that is, over all elements of $\operatorname*{Div}n$).
\end{itemize}

So Theorem \ref{thm.gcdseq.fac} says that if $a\in R^{\left\{  1,2,3,\ldots
\right\}  }$ is a strong gcd-sequence, then there exists a sequence $b\in
R^{\left\{  1,2,3,\ldots\right\}  }$ such that each $n\geq1$ satisfies
$a_{n}=b_{\operatorname*{Div}n}$.

A crucial ingredient of our proof of this theorem will be the notion of a
\textquotedblleft$k$-factorization\textquotedblright\ of the sequence $a$.
This is a hybrid between the sequence $a$ (which is given) and the sequence
$b$ (which we are trying to construct). This hybrid is not a single sequence,
but a pair $\left(  c,r\right)  $ of two sequences carrying complementary
information (in particular, the first $k$ entries of $r$ are $1$, whereas all
but the first $k$ entries of $c$ are $1$). You can think of $\left(
c,r\right)  $ as a snapshot of a recursive construction of $b$
\textquotedblleft caught in the process\textquotedblright\ (after the first
$k$ entries have been fully processed but the remaining entries are still in
limbo), and this does in fact give a good understanding of its use in our
proof of Theorem \ref{thm.gcdseq.fac}.

\begin{definition}
\label{def.gcdseq.kfac} Given a sequence $a\in R^{\left\{  1,2,3,\ldots
\right\}  }$ and an integer $k\geq0$, we define a $k$\textbf{-factorization}
of $a$ to mean a pair $\left(  c,r\right)  $ of two sequences $c\in
R^{\left\{  1,2,3,\ldots\right\}  }$ and $r\in R^{\left\{  1,2,3,\ldots
\right\}  }$ that satisfy the following four conditions:

\begin{enumerate}
\item We have $c_{n}=1$ for all $n>k$.

\item We have $r_{n}=1$ for all $n\leq k$.

\item We have $a_{n}=c_{\operatorname*{Div}n}r_{n}$ for all $n\geq1$.

\item For any $u,v\geq1$, we have the equality
\begin{equation}
c_{\operatorname*{Div}\ell}\left(  r_{u}R+r_{v}R\right)
=c_{\operatorname*{Div}\ell}r_{\gcd\left(  u,v\right)  }R,
\label{eq.gcdseq.k-fac.4}%
\end{equation}
where $\ell=\operatorname{lcm}\left(  u,v\right)  $.
\end{enumerate}
\end{definition}

\begin{remark}
In some of the most common situations, the fourth condition in Definition
\ref{def.gcdseq.kfac} can be significantly simplified. Namely, assume that all
entries of our sequence $a$ are regular (see Definition
\ref{def.rings.reg.reg}). (For instance, this holds whenever $R$ is an
integral domain and all entries of $a$ are nonzero.) Then, if a pair $\left(
c,r\right)  $ of two sequences satisfies the third of the four conditions
above, then, for any positive integer $n$, the element $c_{\operatorname*{Div}%
n}$ is regular (because the third condition yields $a_{n}%
=c_{\operatorname*{Div}n}r_{n}$, so that $c_{\operatorname*{Div}n}$ is a
divisor of $a_{n}$; but $a_{n}$ is regular by our assumption on $a$, and
therefore Proposition \ref{prop.rings.reg.123} \textbf{(d)} shows that its
divisor $c_{\operatorname*{Div}n}$ is also regular). That is, the map
$R\rightarrow R,\ x\mapsto c_{\operatorname*{Div}n}x$ is injective for each
$n\geq1$ (by Proposition \ref{prop.rings.reg.123} \textbf{(e)}). Hence, we can
cancel $c_{\operatorname*{Div}\ell}$ from (\ref{eq.gcdseq.k-fac.4}), thus
rewriting the condition (\ref{eq.gcdseq.k-fac.4}) as%
\[
r_{u}R+r_{v}R=r_{\gcd\left(  u,v\right)  }R.
\]
Thus, the fourth condition in Definition \ref{def.gcdseq.kfac} is simply
saying that $r$ is a strong gcd-sequence. However, this simplification is not
allowed unless the entries of $a$ are regular; the general case requires the
more complicated equality (\ref{eq.gcdseq.k-fac.4}).
\end{remark}

\begin{example}
Let $f=\left(  f_{1},f_{2},f_{3},\ldots\right)  $ be the Fibonacci sequence
sans its $0$-th entry $f_{0}=0$. Then, $f$ has a $4$-factorization $\left(
c,r\right)  $. The first $10$ entries of $c$ and $r$ (and of $f$, for
comparison) are shown in the following table:%
\[%
\begin{tabular}
[c]{|c||c|c|c|c|c|c|c|c|c|c|}\hline
$n$ & $1$ & $2$ & $3$ & $4$ & $5$ & $6$ & $7$ & $8$ & $9$ & $10$\\\hline\hline
$f_{n}$ & $1$ & $1$ & $2$ & $3$ & $5$ & $8$ & $13$ & $21$ & $34$ &
$55$\\\hline
$c_{n}$ & $1$ & $1$ & $2$ & $3$ & $1$ & $1$ & $1$ & $1$ & $1$ & $1$\\\hline
$r_{n}$ & $1$ & $1$ & $1$ & $1$ & $5$ & $4$ & $13$ & $7$ & $17$ & $55$\\\hline
\end{tabular}
\ \ \ \ \ \ \ .
\]
For instance, the condition $a_{n}=c_{\operatorname*{Div}n}r_{n}$ in the above
definition of a $k$-factorization, applied to $a=f$ and $n=9$, becomes
$f_{9}=c_{\operatorname*{Div}9}r_{9}=c_{1}c_{3}c_{9}r_{9}=1\cdot2\cdot
1\cdot17$.
\end{example}

\begin{exercise}
Let $a\in R^{\left\{  1,2,3,\ldots\right\}  }$ be a sequence and $k\geq0$ an
integer. Assume that all entries of $a$ are regular. Prove that a
$k$-factorization of $a$ is unique if it exists.
\end{exercise}

It is easy to see that any strong gcd-sequence $a\in R^{\left\{
1,2,3,\ldots\right\}  }$ has the $0$-factorization $\left(  \underline{1}%
,a\right)  $, where $\underline{1}\in R^{\left\{  1,2,3,\ldots\right\}  }$ is
the constant sequence $\left(  1,1,1,\ldots\right)  $.

We shall prove Theorem \ref{thm.gcdseq.fac} by starting with this trivial
$0$-factorization and gradually improving it to obtain a $1$-factorization,
then a $2$-factorization, then a $3$-factorization, and so on; the
\textquotedblleft limit\textquotedblright\ of this process will give us the
desired $b_{1},b_{2},b_{3},\ldots$. The method of improvement is encoded in
the following crucial lemma:

\begin{lemma}
\label{lem.gcdseq.1}Let $a\in R^{\left\{  1,2,3,\ldots\right\}  }$ be a strong
gcd-sequence. Let $k$ be a positive integer. Let $\left(  b,r\right)  $ be a
$\left(  k-1\right)  $-factorization of $a$. Then, there exists a
$k$-factorization $\left(  c,s\right)  $ of $a$ that satisfies%
\[
c_{i}=b_{i}\ \ \ \ \ \ \ \ \ \ \text{for all }i\neq k.
\]

\end{lemma}

\begin{proof}
We have assumed that $\left(  b,r\right)  $ is a $\left(  k-1\right)
$-factorization of $a$. By the definition of a $\left(  k-1\right)
$-factorization, this means that we have%
\begin{equation}
b_{n}=1\ \ \ \ \ \ \ \ \ \ \text{for all }n>k-1 \label{pf.lem.gcdseq.1.b1}%
\end{equation}
and%
\begin{equation}
r_{n}=1\ \ \ \ \ \ \ \ \ \ \text{for all }n\leq k-1 \label{pf.lem.gcdseq.1.b2}%
\end{equation}
and%
\begin{equation}
a_{n}=b_{\operatorname*{Div}n}r_{n}\ \ \ \ \ \ \ \ \ \ \text{for all }n\geq1,
\label{pf.lem.gcdseq.1.b3}%
\end{equation}
and that for any $u,v\geq1$, we have the equality
\begin{equation}
b_{\operatorname*{Div}\ell}\left(  r_{u}R+r_{v}R\right)
=b_{\operatorname*{Div}\ell}r_{\gcd\left(  u,v\right)  }R,
\label{pf.lem.gcdseq.1.b4}%
\end{equation}
where $\ell=\operatorname{lcm}\left(  u,v\right)  $.

For any positive integer $n$ satisfying $k\mid n$, we have%
\begin{equation}
b_{\operatorname*{Div}n}r_{n}\in b_{\operatorname*{Div}n}r_{k}R.
\label{pf.lem.gcdseq.1.b5}%
\end{equation}
(\textit{Proof of (\ref{pf.lem.gcdseq.1.b5}):} Let $n$ be a positive integer
such that $k\mid n$. Then, $\gcd\left(  k,n\right)  =k$ and
$\operatorname{lcm}\left(  k,n\right)  =n$. Hence, (\ref{pf.lem.gcdseq.1.b4})
(applied to $u=k$ and $v=n$ and $\ell=n$) yields
\[
b_{\operatorname*{Div}n}\left(  r_{k}R+r_{n}R\right)  =b_{\operatorname*{Div}%
n}r_{\gcd\left(  k,n\right)  }R=b_{\operatorname*{Div}n}r_{k}R
\]
(since $\gcd\left(  k,n\right)  =k$). But $r_{n}=r_{k}0+r_{n}1\in r_{k}%
R+r_{n}R$, so that $b_{\operatorname*{Div}n}r_{n}\in b_{\operatorname*{Div}%
n}\left(  r_{k}R+r_{n}R\right)  =b_{\operatorname*{Div}n}r_{k}R$. This proves
(\ref{pf.lem.gcdseq.1.b5}).)

Now, we define two sequences $c,s\in R^{\left\{  1,2,3,\ldots\right\}  }$ as follows:

\begin{itemize}
\item We define the sequence $c$ by setting%
\begin{equation}
c_{k}=r_{k} \label{pf.lem.gcdseq.1.ck}%
\end{equation}
and%
\begin{equation}
c_{n}=b_{n}\ \ \ \ \ \ \ \ \ \ \text{for each }n\neq k.
\label{pf.lem.gcdseq.1.cn}%
\end{equation}

\item We define the sequence $s$ as follows:

\begin{itemize}
\item We set
\begin{equation}
s_{k}=1. \label{pf.lem.gcdseq.1.sk}%
\end{equation}

\item For each positive integer $n$ satisfying $k\nmid n$, we set%
\begin{equation}
s_{n}=r_{n}. \label{pf.lem.gcdseq.1.sn-no}%
\end{equation}

\item For each positive integer $n$ satisfying $k\mid n$ and $n\neq k$, we
define $s_{n}$ implicitly: Namely, we let $s_{n}$ be an element of $R$ that
satisfies $b_{\operatorname*{Div}n}r_{n}=b_{\operatorname*{Div}n}r_{k}s_{n}$.
(Such an element exists, since (\ref{pf.lem.gcdseq.1.b5}) yields
$b_{\operatorname*{Div}n}r_{n}\in b_{\operatorname*{Div}n}r_{k}R$. If it fails
to be unique, we just pick an arbitrary choice.)
\end{itemize}

Thus, if $n$ is a positive integer satisfying $k\mid n$, then $s_{n}$ is an
element of $R$ that satisfies
\begin{equation}
b_{\operatorname*{Div}n}r_{n}=b_{\operatorname*{Div}n}r_{k}s_{n}.
\label{pf.lem.gcdseq.1.sn-yes}%
\end{equation}
(Indeed, if $n=k$, then this follows from (\ref{pf.lem.gcdseq.1.sk}); on the
other hand, if $n\neq k$, then this follows from the definition of $s_{n}$ as
some element of $R$ that satisfies $b_{\operatorname*{Div}n}r_{n}%
=b_{\operatorname*{Div}n}r_{k}s_{n}$.)
\end{itemize}

We shall now show that $\left(  c,s\right)  $ is a $k$-factorization of $a$.
In order to show this, we must prove the following four claims:

\begin{statement}
\textit{Claim 1:} We have $c_{n}=1$ for all $n>k$.
\end{statement}

\begin{statement}
\textit{Claim 2:} We have $s_{n}=1$ for all $n\leq k$.
\end{statement}

\begin{statement}
\textit{Claim 3:} We have $a_{n}=c_{\operatorname*{Div}n}s_{n}$ for all
$n\geq1$.
\end{statement}

\begin{statement}
\textit{Claim 4:} For any $u,v\geq1$, we have the equality
\[
c_{\operatorname*{Div}\ell}\left(  s_{u}R+s_{v}R\right)
=c_{\operatorname*{Div}\ell}s_{\gcd\left(  u,v\right)  }R,
\]
where $\ell=\operatorname{lcm}\left(  u,v\right)  $.
\end{statement}

The first two of these four claims are particularly easy:

\begin{proof}
[Proof of Claim 1.]Let $n>k$ be an integer. Then, $n\neq k$. Hence,
(\ref{pf.lem.gcdseq.1.cn}) yields $c_{n}=b_{n}=1$ (by
(\ref{pf.lem.gcdseq.1.b1}), since $n>k>k-1$). This proves Claim 1.
\end{proof}

\begin{proof}
[Proof of Claim 2.]Let $n\leq k$ be an integer. We must prove that $s_{n}=1$.

If $n=k$, then this follows immediately from (\ref{pf.lem.gcdseq.1.sk}). Thus,
we WLOG assume that $n\neq k$. Hence, $n<k$ (since $n\leq k$) and therefore
$k\nmid n$ (since $n$ and $k$ are positive integers, so that $k\mid n$ would
entail $k\leq n$). Thus, by (\ref{pf.lem.gcdseq.1.sn-no}), we have
$s_{n}=r_{n}$. However, from $n<k$, we obtain $n\leq k-1$ and therefore
$r_{n}=1$ by (\ref{pf.lem.gcdseq.1.b2}). Hence, $s_{n}=r_{n}=1$. This proves
Claim 2.
\end{proof}

Before we step to the proof of Claim 3, let us prove something basic:

\begin{statement}
\textit{Claim 5:} Let $n\geq1$. Then,%
\[
c_{\operatorname*{Div}n}=%
\begin{cases}
b_{\operatorname*{Div}n}, & \text{if }k\nmid n;\\
r_{k}b_{\operatorname*{Div}n}, & \text{if }k\mid n.
\end{cases}
\]

\end{statement}

\begin{proof}
[Proof of Claim 5.]The equality (\ref{eq.gcdseq.aDiv}) yields
\[
c_{\operatorname*{Div}n}=\prod_{d\mid n}c_{d}\ \ \ \ \ \ \ \ \ \ \text{and}%
\ \ \ \ \ \ \ \ \ \ b_{\operatorname*{Div}n}=\prod_{d\mid n}b_{d}.
\]

We are in one of the following two cases:

\textit{Case 1:} We have $k\nmid n$.

\textit{Case 2:} We have $k\mid n$.

Let us first consider Case 1. In this case, we have $k\nmid n$.

Now, let $d\in\operatorname*{Div}n$. Thus, $d$ is a divisor of $n$, so that
$d\mid n$. Hence, $k\nmid d$ (because otherwise, we would have $k\mid d\mid
n$, contradicting $k\nmid n$) and therefore $d\neq k$ (since otherwise, we
would have $k=d$, contradicting $k\nmid d$). Therefore, $c_{d}=b_{d}$ (by
(\ref{pf.lem.gcdseq.1.cn}), applied to $d$ instead of $n$).

Forget that we fixed $d$. We thus have proved that $c_{d}=b_{d}$ for each
$d\in\operatorname*{Div}n$. Multiplying these equalities over all
$d\in\operatorname*{Div}n$, we obtain $\prod_{d\mid n}c_{d}=\prod_{d\mid
n}b_{d}$ (since \textquotedblleft$\prod_{d\mid n}$\textquotedblright\ is
shorthand for \textquotedblleft$\prod_{d\in\operatorname*{Div}n}%
$\textquotedblright). Hence, $c_{\operatorname*{Div}n}=\prod_{d\mid n}%
c_{d}=\prod_{d\mid n}b_{d}=b_{\operatorname*{Div}n}$. This is exactly what
Claim 5 is saying in the $k\nmid n$ case. Thus, Claim 5 is proved in Case 1.

Let us now consider Case 2. In this case, we have $k\mid n$. Hence, the
product $\prod_{d\mid n}c_{d}$ has a $d=k$ factor. Splitting off this factor,
we obtain%
\begin{equation}
\prod_{d\mid n}c_{d}=\underbrace{c_{k}}_{\substack{=r_{k}\\\text{(by
(\ref{pf.lem.gcdseq.1.ck}))}}}\prod_{\substack{d\mid n;\\d\neq k}%
}\underbrace{c_{d}}_{\substack{=b_{d}\\\text{(by (\ref{pf.lem.gcdseq.1.cn}%
),}\\\text{applied to }d\text{ instead of }n\text{)}}}=r_{k}\prod
_{\substack{d\mid n;\\d\neq k}}b_{d}. \label{pf.lem.gcdseq.1.c5.4}%
\end{equation}
However, (\ref{pf.lem.gcdseq.1.b1}) (applied to $k$ instead of $n$) yields
$b_{k}=1$ (since $k>k-1$). Again, we can split off the $d=k$ factor from the
product $\prod_{d\mid n}b_{d}$ (since $k\mid n$); this yields
\[
\prod_{d\mid n}b_{d}=\underbrace{b_{k}}_{=1}\prod_{\substack{d\mid n;\\d\neq
k}}b_{d}=\prod_{\substack{d\mid n;\\d\neq k}}b_{d}.
\]
Thus,%
\begin{equation}
b_{\operatorname*{Div}n}=\prod_{d\mid n}b_{d}=\prod_{\substack{d\mid n;\\d\neq
k}}b_{d}. \label{pf.lem.gcdseq.1.c5.5}%
\end{equation}
Now,
\begin{align*}
c_{\operatorname*{Div}n}  &  =\prod_{d\mid n}c_{d}=r_{k}\prod_{\substack{d\mid
n;\\d\neq k}}b_{d}\ \ \ \ \ \ \ \ \ \ \left(  \text{by
(\ref{pf.lem.gcdseq.1.c5.4})}\right) \\
&  =r_{k}b_{\operatorname*{Div}n}\ \ \ \ \ \ \ \ \ \ \left(  \text{since
(\ref{pf.lem.gcdseq.1.c5.5}) yields }\prod_{\substack{d\mid n;\\d\neq k}%
}b_{d}=b_{\operatorname*{Div}n}\right)  .
\end{align*}
This is exactly what Claim 5 is saying in the $k\mid n$ case. Thus, Claim 5 is
proved in Case 2.

So the proof of Claim 5 is complete (since we have finished both cases).
\end{proof}

\begin{proof}
[Proof of Claim 3.]Let $n\geq1$. We must prove that $a_{n}%
=c_{\operatorname*{Div}n}s_{n}$. We are in one of the following two cases:

\textit{Case 1:} We have $k\nmid n$.

\textit{Case 2:} We have $k\mid n$.

Let us first consider Case 1. In this case, we have $k\nmid n$. Hence,
(\ref{pf.lem.gcdseq.1.sn-no}) yields $s_{n}=r_{n}$. Furthermore, Claim 5
yields $c_{\operatorname*{Div}n}=b_{\operatorname*{Div}n}$ (since $k\nmid n$).
But (\ref{pf.lem.gcdseq.1.b3}) yields $a_{n}=b_{\operatorname*{Div}n}r_{n}$.
In view of $s_{n}=r_{n}$ and $c_{\operatorname*{Div}n}=b_{\operatorname*{Div}%
n}$, we can rewrite this as $a_{n}=c_{\operatorname*{Div}n}s_{n}$. This proves
Claim 3 in Case 1.

Let us now consider Case 2. In this case, we have $k\mid n$. Hence,
(\ref{pf.lem.gcdseq.1.sn-yes}) yields $b_{\operatorname*{Div}n}r_{n}%
=b_{\operatorname*{Div}n}r_{k}s_{n}$. Furthermore, Claim 5 yields
$c_{\operatorname*{Div}n}=r_{k}b_{\operatorname*{Div}n}$ (since $k\mid n$);
therefore $c_{\operatorname*{Div}n}s_{n}=r_{k}b_{\operatorname*{Div}n}%
s_{n}=b_{\operatorname*{Div}n}r_{k}s_{n}$ (since $R$ is commutative). But
(\ref{pf.lem.gcdseq.1.b3}) yields $a_{n}=b_{\operatorname*{Div}n}%
r_{n}=b_{\operatorname*{Div}n}r_{k}s_{n}$. Comparing these two equalities, we
obtain $a_{n}=c_{\operatorname*{Div}n}s_{n}$. This proves Claim 3 in Case 2.

Thus, we have proved Claim 3 in both cases; so the proof is complete.
\end{proof}

Claim 4 is even more complicated. We begin with three auxiliary claims:

\begin{statement}
\textit{Claim 6:} Let $u,v,w$ be three positive integers such that $u\mid w$
and $v\mid w$. Then,%
\[
b_{\operatorname*{Div}w}\left(  r_{u}R+r_{v}R\right)  =b_{\operatorname*{Div}%
w}r_{\gcd\left(  u,v\right)  }R.
\]

\end{statement}

\begin{proof}
[Proof of Claim 6.]Let $\ell=\operatorname{lcm}\left(  u,v\right)  $. The
integer $w$ is a common multiple of $u$ and $v$ (since $u\mid w$ and $v\mid
w$), and thus is divisible by their lcm. That is, $\operatorname{lcm}\left(
u,v\right)  \mid w$. In other words, $\ell\mid w$ (since $\ell
=\operatorname{lcm}\left(  u,v\right)  $). Hence, every positive divisor of
$\ell$ is a positive divisor of $w$. In other words, $\operatorname*{Div}%
\ell\subseteq\operatorname*{Div}w$.

However, it is easy to see that if $X$ and $Y$ are two finite sets of positive
integers satisfying $X\subseteq Y$, then%
\begin{equation}
b_{Y}=b_{Y\setminus X}b_{X}. \label{pf.lem.gcdseq.1.c6.genAB}%
\end{equation}
(Indeed, the definition of $b_{S}$ for any finite set $S$ of positive integers
shows that
\[
b_{X}=\prod_{d\in X}b_{d}\ \ \ \ \ \ \ \ \ \ \text{and}%
\ \ \ \ \ \ \ \ \ \ b_{Y}=\prod_{d\in Y}b_{d}\ \ \ \ \ \ \ \ \ \ \text{and}%
\ \ \ \ \ \ \ \ \ \ b_{Y\setminus X}=\prod_{d\in Y\setminus X}b_{d}.
\]
Thus, the equality (\ref{pf.lem.gcdseq.1.c6.genAB}) boils down to $\prod_{d\in
Y}b_{d}=\left(  \prod_{d\in Y\setminus X}b_{d}\right)  \left(  \prod_{d\in
X}b_{d}\right)  $, which is a basic property of finite products.)

Applying (\ref{pf.lem.gcdseq.1.c6.genAB}) to $X=\operatorname*{Div}\ell$ and
$Y=\operatorname*{Div}w$, we obtain%
\[
b_{\operatorname*{Div}w}=b_{\left(  \operatorname*{Div}w\right)
\setminus\left(  \operatorname*{Div}\ell\right)  }b_{\operatorname*{Div}\ell
}\ \ \ \ \ \ \ \ \ \ \left(  \text{since }\operatorname*{Div}\ell
\subseteq\operatorname*{Div}w\right)  .
\]
Hence,%
\begin{align*}
\underbrace{b_{\operatorname*{Div}w}}_{=b_{\left(  \operatorname*{Div}%
w\right)  \setminus\left(  \operatorname*{Div}\ell\right)  }%
b_{\operatorname*{Div}\ell}}\left(  r_{u}R+r_{v}R\right)   &  =b_{\left(
\operatorname*{Div}w\right)  \setminus\left(  \operatorname*{Div}\ell\right)
}\underbrace{b_{\operatorname*{Div}\ell}\left(  r_{u}R+r_{v}R\right)
}_{\substack{=b_{\operatorname*{Div}\ell}r_{\gcd\left(  u,v\right)
}R\\\text{(by (\ref{pf.lem.gcdseq.1.b4}))}}}\\
&  =\underbrace{b_{\left(  \operatorname*{Div}w\right)  \setminus\left(
\operatorname*{Div}\ell\right)  }b_{\operatorname*{Div}\ell}}%
_{=b_{\operatorname*{Div}w}}r_{\gcd\left(  u,v\right)  }R\\
&  =b_{\operatorname*{Div}w}r_{\gcd\left(  u,v\right)  }R.
\end{align*}
This proves Claim 6.
\end{proof}

\begin{statement}
\textit{Claim 7:} Let $u,v,w$ be three positive integers such that $u\mid
v\mid w$. Then,%
\[
b_{\operatorname*{Div}w}r_{v}R\subseteq b_{\operatorname*{Div}w}r_{u}R.
\]

\end{statement}

\begin{proof}
[Proof of Claim 7.]We have $\gcd\left(  u,v\right)  =u$ (since $u\mid v$).
But
\begin{align*}
b_{\operatorname*{Div}w}\underbrace{r_{v}R}_{\subseteq r_{u}R+r_{v}R}  &
\subseteq b_{\operatorname*{Div}w}\left(  r_{u}R+r_{v}R\right) \\
&  =b_{\operatorname*{Div}w}r_{\gcd\left(  u,v\right)  }%
R\ \ \ \ \ \ \ \ \ \ \left(  \text{by Claim 6}\right) \\
&  =b_{\operatorname*{Div}w}r_{u}R\ \ \ \ \ \ \ \ \ \ \left(  \text{since
}\gcd\left(  u,v\right)  =u\right)  .
\end{align*}
This proves Claim 7.
\end{proof}

\begin{statement}
\textit{Claim 8:} Let $u$ and $v$ be two positive integers such that $k\mid
u\mid v$. Then,%
\[
b_{\operatorname*{Div}v}r_{u}=c_{\operatorname*{Div}v}s_{u}.
\]

\end{statement}

\begin{proof}
[Proof of Claim 8.]We have $k\mid v$. Hence, Claim 5 (applied to $n=v$) yields
$c_{\operatorname*{Div}v}=r_{k}b_{\operatorname*{Div}v}$. That is,
$r_{k}b_{\operatorname*{Div}v}=c_{\operatorname*{Div}v}$.

But we have $u\mid v$. Hence, every positive divisor of $u$ is a positive
divisor of $v$. In other words, $\operatorname*{Div}u\subseteq
\operatorname*{Div}v$. Hence, applying (\ref{pf.lem.gcdseq.1.c6.genAB}) to
$X=\operatorname*{Div}u$ and $Y=\operatorname*{Div}v$, we obtain%
\[
b_{\operatorname*{Div}v}=b_{\left(  \operatorname*{Div}v\right)
\setminus\left(  \operatorname*{Div}u\right)  }b_{\operatorname*{Div}u}.
\]
But $k\mid u$; hence, (\ref{pf.lem.gcdseq.1.sn-yes}) (applied to $n=u$) yields%
\[
b_{\operatorname*{Div}u}r_{u}=b_{\operatorname*{Div}u}r_{k}s_{u}.
\]
Hence,%
\begin{align*}
\underbrace{b_{\operatorname*{Div}v}}_{=b_{\left(  \operatorname*{Div}%
v\right)  \setminus\left(  \operatorname*{Div}u\right)  }%
b_{\operatorname*{Div}u}}r_{u}  &  =b_{\left(  \operatorname*{Div}v\right)
\setminus\left(  \operatorname*{Div}u\right)  }%
\underbrace{b_{\operatorname*{Div}u}r_{u}}_{=b_{\operatorname*{Div}u}%
r_{k}s_{u}}=\underbrace{b_{\left(  \operatorname*{Div}v\right)  \setminus
\left(  \operatorname*{Div}u\right)  }b_{\operatorname*{Div}u}}%
_{=b_{\operatorname*{Div}v}}r_{k}s_{u}\\
&  =\underbrace{b_{\operatorname*{Div}v}r_{k}}_{=r_{k}b_{\operatorname*{Div}%
v}=c_{\operatorname*{Div}v}}s_{u}=c_{\operatorname*{Div}v}s_{u}.
\end{align*}
This proves Claim 8.
\end{proof}

\begin{proof}
[Proof of Claim 4.]Let $u,v\geq1$ be integers. Let $\ell=\operatorname{lcm}%
\left(  u,v\right)  $. We must prove that
\begin{equation}
c_{\operatorname*{Div}\ell}\left(  s_{u}R+s_{v}R\right)
=c_{\operatorname*{Div}\ell}s_{\gcd\left(  u,v\right)  }R.
\label{pf.lem.gcdseq.1.c4.goal}%
\end{equation}

Note that $\ell=\operatorname{lcm}\left(  u,v\right)  $ entails $u\mid\ell$
and $v\mid\ell$.

Also, set $g=\gcd\left(  u,v\right)  $. Then, $g\mid u$ and $g\mid v$.

We are in one of the following four cases:

\textit{Case 1:} We have $k\nmid u$ and $k\nmid v$.

\textit{Case 2:} We have $k\nmid u$ and $k\mid v$.

\textit{Case 3:} We have $k\mid u$ and $k\nmid v$.

\textit{Case 4:} We have $k\mid u$ and $k\mid v$.

Let us first consider Case 1. In this case, we have $k\nmid u$ and $k\nmid v$.
Thus, $k\nmid\gcd\left(  u,v\right)  $ as well (since otherwise, we would have
$k\mid\gcd\left(  u,v\right)  \mid v$, contradicting $k\nmid v$).

Claim 5 (applied to $n=\ell$) yields
\[
c_{\operatorname*{Div}\ell}=%
\begin{cases}
b_{\operatorname*{Div}\ell}, & \text{if }k\nmid\ell;\\
r_{k}b_{\operatorname*{Div}\ell}, & \text{if }k\mid\ell
\end{cases}
\ \ =%
\begin{cases}
1, & \text{if }k\nmid\ell;\\
r_{k}, & \text{if }k\mid\ell
\end{cases}
\ \ \cdot b_{\operatorname*{Div}\ell}.
\]

Thus, $b_{\operatorname*{Div}\ell}$ divides $c_{\operatorname*{Div}\ell}$ in
$R$. In other words, there exists some $z\in R$ such that
$c_{\operatorname*{Div}\ell}=b_{\operatorname*{Div}\ell}z$. Consider this $z$.

Now, (\ref{pf.lem.gcdseq.1.sn-no}) (applied to $n=\gcd\left(  u,v\right)  $)
yields $s_{\gcd\left(  u,v\right)  }=r_{\gcd\left(  u,v\right)  }$ (since
$k\nmid\gcd\left(  u,v\right)  $). Similarly, $s_{u}=r_{u}$ and $s_{v}=r_{v}$.
In view of these three equalities, we can rewrite (\ref{pf.lem.gcdseq.1.b4})
as%
\[
b_{\operatorname*{Div}\ell}\left(  s_{u}R+s_{v}R\right)
=b_{\operatorname*{Div}\ell}s_{\gcd\left(  u,v\right)  }R.
\]
Multiplying both sides of this equality by $z$, we obtain%
\[
zb_{\operatorname*{Div}\ell}\left(  s_{u}R+s_{v}R\right)
=zb_{\operatorname*{Div}\ell}s_{\gcd\left(  u,v\right)  }R.
\]
Since $c_{\operatorname*{Div}\ell}=b_{\operatorname*{Div}\ell}%
z=zb_{\operatorname*{Div}\ell}$, we can rewrite this as
\[
c_{\operatorname*{Div}\ell}\left(  s_{u}R+s_{v}R\right)
=c_{\operatorname*{Div}\ell}s_{\gcd\left(  u,v\right)  }R.
\]
Thus, Claim 4 is proved in Case 1.

Next, we observe that our situation and goal do not change if we swap $u$ with
$v$ (since $\gcd\left(  u,v\right)  =\gcd\left(  v,u\right)  $ and
$\operatorname{lcm}\left(  u,v\right)  =\operatorname{lcm}\left(  v,u\right)
$ and $s_{u}R+s_{v}R=s_{v}R+s_{u}R$). But such a swap will take us from Case 2
into Case 3 and vice versa. Thus, if we can prove Claim 4 in Case 3, then we
will automatically obtain Claim 4 in Case 2 as well.

Let us therefore consider Case 3 next. In this case, we have $k\mid u$ and
$k\nmid v$.

Thus, $k\nmid g$ as well (since otherwise, we would have $k\mid g\mid v$,
contradicting $k\nmid v$). Hence, (\ref{pf.lem.gcdseq.1.sn-no}) (applied to
$n=g$) yields $s_{g}=r_{g}$. Similarly, $s_{v}=r_{v}$ (since $k\nmid v$).

Let $\beta=b_{\operatorname*{Div}\ell}$. Then, (\ref{pf.lem.gcdseq.1.b4})
rewrites as%
\begin{equation}
\beta\left(  r_{u}R+r_{v}R\right)  =\beta r_{g}R \label{pf.lem.gcdseq.1.c4.2}%
\end{equation}
(since $\beta=b_{\operatorname*{Div}\ell}$ and $g=\gcd\left(  u,v\right)  $).

We have $k\mid u\mid\ell$. Therefore, Claim 5 (applied to $n=\ell$) yields
$c_{\operatorname*{Div}\ell}=r_{k}b_{\operatorname*{Div}\ell}$. In view of
$\beta=b_{\operatorname*{Div}\ell}$, we can rewrite this as
\begin{equation}
c_{\operatorname*{Div}\ell}=r_{k}\beta. \label{pf.lem.gcdseq.1.c4.3}%
\end{equation}

Also, $u\mid\ell$. Hence, each positive divisor of $u$ is a positive divisor
of $\ell$. In other words, $\operatorname*{Div}u\subseteq\operatorname*{Div}%
\ell$. Thus, (\ref{pf.lem.gcdseq.1.c6.genAB}) (applied to
$X=\operatorname*{Div}u$ and $Y=\operatorname*{Div}\ell$) yields
\[
b_{\operatorname*{Div}\ell}=b_{\left(  \operatorname*{Div}\ell\right)
\setminus\left(  \operatorname*{Div}u\right)  }b_{\operatorname*{Div}u}.
\]

Recall that $k\mid u$; thus, (\ref{pf.lem.gcdseq.1.sn-yes}) (applied to $n=u$)
yields $b_{\operatorname*{Div}u}r_{u}=b_{\operatorname*{Div}u}r_{k}s_{u}$.
Thus,%
\begin{align}
\underbrace{\beta}_{\substack{=b_{\operatorname*{Div}\ell}\\=b_{\left(
\operatorname*{Div}\ell\right)  \setminus\left(  \operatorname*{Div}u\right)
}b_{\operatorname*{Div}u}}}r_{u}  &  =b_{\left(  \operatorname*{Div}%
\ell\right)  \setminus\left(  \operatorname*{Div}u\right)  }%
\underbrace{b_{\operatorname*{Div}u}r_{u}}_{=b_{\operatorname*{Div}u}%
r_{k}s_{u}}=\underbrace{b_{\left(  \operatorname*{Div}\ell\right)
\setminus\left(  \operatorname*{Div}u\right)  }b_{\operatorname*{Div}u}%
}_{=b_{\operatorname*{Div}\ell}=\beta}r_{k}s_{u}\nonumber\\
&  =\beta r_{k}s_{u}. \label{pf.lem.gcdseq.1.c4.4}%
\end{align}

Our goal is to show that (\ref{pf.lem.gcdseq.1.c4.goal}) holds, i.e., that%
\begin{equation}
c_{\operatorname*{Div}\ell}\left(  s_{u}R+s_{v}R\right)
=c_{\operatorname*{Div}\ell}s_{\gcd\left(  u,v\right)  }R.
\label{pf.lem.gcdseq.1.c4.5}%
\end{equation}
Since%
\begin{align*}
\underbrace{c_{\operatorname*{Div}\ell}}_{=r_{k}\beta}\left(  s_{u}%
R+s_{v}R\right)   &  =r_{k}\beta\left(  s_{u}R+s_{v}R\right) \\
&  =\underbrace{r_{k}\beta s_{u}}_{\substack{=\beta r_{k}s_{u}=\beta
r_{u}\\\text{(by (\ref{pf.lem.gcdseq.1.c4.4}))}}}R+\underbrace{r_{k}\beta
}_{=\beta r_{k}}\underbrace{s_{v}}_{=r_{v}}R\\
&  =\beta r_{u}R+\beta r_{k}r_{v}R=\beta\left(  r_{u}R+r_{k}r_{v}R\right)
\end{align*}
and%
\[
\underbrace{c_{\operatorname*{Div}\ell}}_{=r_{k}\beta}\underbrace{s_{\gcd
\left(  u,v\right)  }}_{\substack{=s_{g}\\\text{(since }\gcd\left(
u,v\right)  =g\text{)}}}R=r_{k}\beta s_{g}R=\beta r_{k}\underbrace{s_{g}%
}_{=r_{g}}R=\beta r_{k}r_{g}R,
\]
we can rewrite this as
\begin{equation}
\beta\left(  r_{u}R+r_{k}r_{v}R\right)  =\beta r_{k}r_{g}R.
\label{pf.lem.gcdseq.1.c4.6}%
\end{equation}
So this is the equality we must prove.

We have $k\nmid g$ and $\gcd\left(  k,g\right)  \mid g$. Hence, $\gcd\left(
k,g\right)  \neq k$ (since otherwise, we would have $k=\gcd\left(  k,g\right)
\mid g$, contradicting $k\nmid g$). Combining this with $\gcd\left(
k,g\right)  \leq k$ (which follows from $\gcd\left(  k,g\right)  \mid k$), we
obtain $\gcd\left(  k,g\right)  <k$ and therefore $\gcd\left(  k,g\right)
\leq k-1$ (since $\gcd\left(  k,g\right)  $ and $k$ are integers). Thus,
(\ref{pf.lem.gcdseq.1.b2}) (applied to $n=\gcd\left(  k,g\right)  $) yields
$r_{\gcd\left(  k,g\right)  }=1$. But we have $k\mid\ell$ and $g\mid v\mid
\ell$; thus, Claim 6 (applied to $k$, $g$ and $\ell$ instead of $u$, $v$ and
$w$) yields%
\[
b_{\operatorname*{Div}\ell}\left(  r_{k}R+r_{g}R\right)
=b_{\operatorname*{Div}\ell}\underbrace{r_{\gcd\left(  k,g\right)  }}%
_{=1}R=b_{\operatorname*{Div}\ell}R.
\]
In view of $\beta=b_{\operatorname*{Div}\ell}$, we can rewrite this as%
\begin{equation}
\beta\left(  r_{k}R+r_{g}R\right)  =\beta R. \label{pf.lem.gcdseq.1.c4.13}%
\end{equation}

Thus, Lemma \ref{lem.gcdseq.abcR} (applied to $a=\beta$ and $I=r_{k}R$ and
$J=r_{g}R$) yields
\begin{equation}
\beta r_{k}R\cap\beta r_{g}R=\beta\underbrace{\left(  r_{k}R\right)  \left(
r_{g}R\right)  }_{\substack{=r_{k}r_{g}R\\\text{(by the }\left(  ab\right)
R=\left(  aR\right)  \left(  bR\right)  \\\text{part of Exercise
\ref{exe.ideals.abR})}}}=\beta r_{k}r_{g}R. \label{pf.lem.gcdseq.1.c4.10}%
\end{equation}

Furthermore, $k\mid u\mid\ell$. Thus, Claim 7 (applied to $k$, $u$ and $\ell$
instead of $u$, $v$ and $w$) yields $b_{\operatorname*{Div}\ell}%
r_{u}R\subseteq b_{\operatorname*{Div}\ell}r_{k}R$. In other words,%
\begin{equation}
\beta r_{u}R\subseteq\beta r_{k}R \label{pf.lem.gcdseq.1.c4.14}%
\end{equation}
(since $\beta=b_{\operatorname*{Div}\ell}$).

Now, combining the inclusions%
\[
\beta\left(  r_{u}R+r_{k}r_{v}R\right)  =\underbrace{\beta r_{u}%
R}_{\substack{\subseteq\beta r_{k}R\\\text{(by (\ref{pf.lem.gcdseq.1.c4.14}%
))}}}+\,\beta r_{k}\underbrace{r_{v}R}_{\subseteq R}\subseteq\beta
r_{k}R+\beta r_{k}R=\beta r_{k}R
\]
and
\[
\beta\left(  r_{u}R+\underbrace{r_{k}r_{v}R}_{\subseteq r_{v}R}\right)
\subseteq\beta\left(  r_{u}R+r_{v}R\right)  =\beta r_{g}%
R\ \ \ \ \ \ \ \ \ \ \left(  \text{by (\ref{pf.lem.gcdseq.1.c4.2})}\right)  ,
\]
we obtain%
\[
\beta\left(  r_{u}R+r_{k}r_{v}R\right)  \subseteq\beta r_{k}R\cap\beta
r_{g}R=\beta r_{k}r_{g}R\ \ \ \ \ \ \ \ \ \ \left(  \text{by
(\ref{pf.lem.gcdseq.1.c4.10})}\right)  .
\]
Combining this with%
\begin{align*}
\beta r_{k}r_{g}R  &  =r_{k}\underbrace{\beta r_{g}R}_{\substack{=\beta\left(
r_{u}R+r_{v}R\right)  \\\text{(by (\ref{pf.lem.gcdseq.1.c4.2}))}}}=r_{k}%
\beta\left(  r_{u}R+r_{v}R\right)  =\beta r_{k}\left(  r_{u}R+r_{v}R\right) \\
&  =\beta\left(  \underbrace{r_{k}r_{u}R}_{\subseteq r_{u}R}+\,r_{k}%
r_{v}R\right)  \subseteq\beta\left(  r_{u}R+r_{k}r_{v}R\right)  ,
\end{align*}
we obtain%
\[
\beta\left(  r_{u}R+r_{k}r_{v}R\right)  =\beta r_{k}r_{g}R.
\]
This proves (\ref{pf.lem.gcdseq.1.c4.6}). As we explained above, this
completes the proof of Claim 4 in Case 3.

As we already explained, Case 2 can be reduced to Case 3. Thus, the only
remaining case is Case 4.

So let us assume that we are in Case 4. In this case, we have $k\mid u$ and
$k\mid v$. Hence, $k$ is a common divisor of $u$ and $v$, and therefore
divides their gcd $\gcd\left(  u,v\right)  $. In other words, $k\mid
\gcd\left(  u,v\right)  $.

We have $k\mid u\mid\ell$. Hence, Claim 8 (applied to $\ell$ instead of $v$)
yields%
\begin{equation}
b_{\operatorname*{Div}\ell}r_{u}=c_{\operatorname*{Div}\ell}s_{u}.
\label{pf.lem.gcdseq.1.c4.35}%
\end{equation}
Likewise, using $k\mid v\mid\ell$, we obtain%
\begin{equation}
b_{\operatorname*{Div}\ell}r_{v}=c_{\operatorname*{Div}\ell}s_{v}.
\label{pf.lem.gcdseq.1.c4.36}%
\end{equation}
Furthermore, $k\mid\gcd\left(  u,v\right)  \mid u\mid\ell$; thus, Claim 8
(applied to $\gcd\left(  u,v\right)  $ and $\ell$ instead of $u$ and $v$)
yields
\begin{equation}
b_{\operatorname*{Div}\ell}r_{\gcd\left(  u,v\right)  }=c_{\operatorname*{Div}%
\ell}s_{\gcd\left(  u,v\right)  }. \label{pf.lem.gcdseq.1.c4.34}%
\end{equation}
Now,%
\begin{align*}
c_{\operatorname*{Div}\ell}\left(  s_{u}R+s_{v}R\right)   &
=\underbrace{c_{\operatorname*{Div}\ell}s_{u}}%
_{\substack{=b_{\operatorname*{Div}\ell}r_{u}\\\text{(by
(\ref{pf.lem.gcdseq.1.c4.35}))}}}R+\underbrace{c_{\operatorname*{Div}\ell
}s_{v}}_{\substack{=b_{\operatorname*{Div}\ell}r_{v}\\\text{(by
(\ref{pf.lem.gcdseq.1.c4.36}))}}}R\\
&  =b_{\operatorname*{Div}\ell}r_{u}R+b_{\operatorname*{Div}\ell}%
r_{v}R=b_{\operatorname*{Div}\ell}\left(  r_{u}R+r_{v}R\right) \\
&  =\underbrace{b_{\operatorname*{Div}\ell}r_{\gcd\left(  u,v\right)  }%
}_{\substack{=c_{\operatorname*{Div}\ell}s_{\gcd\left(  u,v\right)
}\\\text{(by (\ref{pf.lem.gcdseq.1.c4.34}))}}}R\ \ \ \ \ \ \ \ \ \ \left(
\text{by (\ref{pf.lem.gcdseq.1.b4})}\right) \\
&  =c_{\operatorname*{Div}\ell}s_{\gcd\left(  u,v\right)  }R.
\end{align*}
Thus, Claim 4 is proved in Case 4.

We have now proved Claim 4 in all four cases. Thus, the proof of Claim 4 is complete.
\end{proof}

Now, we have proved all four Claims 1, 2, 3 and 4. Thus, $\left(  c,s\right)
$ is a $k$-factorization of $a$. It remains to show that
\[
c_{i}=b_{i}\ \ \ \ \ \ \ \ \ \ \text{for all }i\neq k.
\]
But this follows immediately from (\ref{pf.lem.gcdseq.1.cn}) (applied to
$n=i$). Thus, the proof of Lemma \ref{lem.gcdseq.1} is complete.
\end{proof}

We can now derive Theorem \ref{thm.gcdseq.fac} easily:

\begin{proof}
[Proof of Theorem \ref{thm.gcdseq.fac}.]By recursion, we shall construct a
sequence
\[
\left(  \left(  b^{\left(  0\right)  },r^{\left(  0\right)  }\right)
,\ \left(  b^{\left(  1\right)  },r^{\left(  1\right)  }\right)  ,\ \left(
b^{\left(  2\right)  },r^{\left(  2\right)  }\right)  ,\ \ldots\right)
\]
of pairs of sequences of elements of $R$ (that is, each $b^{\left(  k\right)
}$ and each $r^{\left(  k\right)  }$ will be a sequence of elements of $R$)
with the following properties:

\begin{enumerate}
\item For each $k\in\mathbb{N}$, the pair $\left(  b^{\left(  k\right)
},r^{\left(  k\right)  }\right)  $ is a $k$-factorization of $a$.

\item For each positive integer $k$, the sequence $b^{\left(  k\right)  }$
differs from $b^{\left(  k-1\right)  }$ only in its $k$-th entry (if it
differs from it at all).
\end{enumerate}

Here is how we perform this construction:

\begin{itemize}
\item We define the first pair $\left(  b^{\left(  0\right)  },r^{\left(
0\right)  }\right)  $ as follows: Let $\underline{1}\in R^{\left\{
1,2,3,\ldots\right\}  }$ be the constant sequence $\left(  1,1,1,\ldots
\right)  $. Note that
\begin{equation}
\underline{1}_{\operatorname*{Div}n}=1\ \ \ \ \ \ \ \ \ \ \text{for each
}n\geq1 \label{pf.thm.gcdseq.fac.1divn}%
\end{equation}
(since $\underline{1}_{\operatorname*{Div}n}=\prod_{d\mid n}%
\underbrace{\underline{1}_{d}}_{=1}=\prod_{d\mid n}1=1$).

We claim that the pair $\left(  \underline{1},a\right)  $ is a $0$%
-factorization of $a$. Indeed, we can check that all four conditions in
Definition \ref{def.gcdseq.kfac} are satisfied for this pair:

\begin{itemize}
\item We have $\underline{1}_{n}=1$ for all $n>0$. (This is clear from the
definition of $\underline{1}$.)

\item We have $a_{n}=1$ for all $n\leq0$. (This is vacuously true, since there
are no positive integers $n\leq0$.)

\item We have $a_{n}=\underline{1}_{\operatorname*{Div}n}a_{n}$ for all
$n\geq1$. (This follows from (\ref{pf.thm.gcdseq.fac.1divn}).)

\item For any $u,v\geq1$, we have the equality
\[
\underline{1}_{\operatorname*{Div}\ell}\left(  a_{u}R+a_{v}R\right)
=\underline{1}_{\operatorname*{Div}\ell}a_{\gcd\left(  u,v\right)  }R,
\]
where $\ell=\operatorname{lcm}\left(  u,v\right)  $. (Indeed,
(\ref{pf.thm.gcdseq.fac.1divn}) shows that $\underline{1}_{\operatorname*{Div}%
\ell}=1$, so this equality boils down to $a_{u}R+a_{v}R=a_{\gcd\left(
u,v\right)  }R$; but the latter holds because $a$ is a strong gcd-sequence.)
\end{itemize}

Thus, $\left(  \underline{1},a\right)  $ is a $0$-factorization of $a$. We
thus define $\left(  b^{\left(  0\right)  },r^{\left(  0\right)  }\right)  $
to be this $0$-factorization $\left(  \underline{1},a\right)  $.

\item Let $k$ be a positive integer. Assume that a $\left(  k-1\right)
$-factorization $\left(  b^{\left(  k-1\right)  },r^{\left(  k-1\right)
}\right)  $ of $a$ has already been defined. Then, Lemma \ref{lem.gcdseq.1}
(applied to $\left(  b,r\right)  =\left(  b^{\left(  k-1\right)  },r^{\left(
k-1\right)  }\right)  $) shows that there exists a $k$-factorization $\left(
c,s\right)  $ of $a$ that satisfies%
\[
c_{i}=\left(  b^{\left(  k-1\right)  }\right)  _{i}%
\ \ \ \ \ \ \ \ \ \ \text{for all }i\neq k.
\]
Fix such a $k$-factorization $\left(  c,s\right)  $, and define $\left(
b^{\left(  k\right)  },r^{\left(  k\right)  }\right)  $ to be this
$k$-factorization. Thus,%
\begin{equation}
\left(  b^{\left(  k\right)  }\right)  _{i}=\left(  b^{\left(  k-1\right)
}\right)  _{i}\ \ \ \ \ \ \ \ \ \ \text{for all }i\neq k.
\label{pf.thm.gcdseq.fac.next}%
\end{equation}
In other words, the sequence $b^{\left(  k\right)  }$ differs from $b^{\left(
k-1\right)  }$ only in its $k$-th entry (if it differs from it at all).
\end{itemize}

Now, define a sequence $b=\left(  b_{1},b_{2},b_{3},\ldots\right)  \in
R^{\left\{  1,2,3,\ldots\right\}  }$ by setting%
\begin{equation}
b_{d}:=\left(  b^{\left(  d\right)  }\right)  _{d}%
\ \ \ \ \ \ \ \ \ \ \text{for each }d\geq1. \label{pf.thm.gcdseq.fac.bdd}%
\end{equation}
Then, we can easily see the following:

\begin{statement}
\textit{Claim 1:} Let $n$ and $d$ be two positive integers such that $n\geq
d$. Then,%
\[
b_{d}=\left(  b^{\left(  n\right)  }\right)  _{d}.
\]

\end{statement}

\begin{proof}
[Proof of Claim 1.]Let $j\geq d$ be any positive integer. Then, $d\leq j<j+1$,
so that $d\neq j+1$. Hence, (\ref{pf.thm.gcdseq.fac.next}) (applied to $k=j+1$
and $i=d$) yields $\left(  b^{\left(  j+1\right)  }\right)  _{d}=\left(
b^{\left(  j+1-1\right)  }\right)  _{d}=\left(  b^{\left(  j\right)  }\right)
_{d}$. That is, $\left(  b^{\left(  j\right)  }\right)  _{d}=\left(
b^{\left(  j+1\right)  }\right)  _{d}$.

Forget that we fixed $j$. We thus have proved the equality $\left(  b^{\left(
j\right)  }\right)  _{d}=\left(  b^{\left(  j+1\right)  }\right)  _{d}$ for
each $j\geq d$. These equalities can be combined into an infinite chain:%
\[
\left(  b^{\left(  d\right)  }\right)  _{d}=\left(  b^{\left(  d+1\right)
}\right)  _{d}=\left(  b^{\left(  d+2\right)  }\right)  _{d}=\left(
b^{\left(  d+3\right)  }\right)  _{d}=\cdots.
\]
In particular, this yields $\left(  b^{\left(  d\right)  }\right)
_{d}=\left(  b^{\left(  n\right)  }\right)  _{d}$ (since $n\geq d$). In view
of (\ref{pf.thm.gcdseq.fac.bdd}), this rewrites as $b_{d}=\left(  b^{\left(
n\right)  }\right)  _{d}$. This proves Claim 1.
\end{proof}

\begin{statement}
\textit{Claim 2:} Each $k\geq1$ satisfies
\[
a_{k}=\prod_{d\mid k}b_{d}.
\]

\end{statement}

\begin{proof}
[Proof of Claim 2.]Let $k\geq1$. Recall that $\left(  b^{\left(  k\right)
},r^{\left(  k\right)  }\right)  $ is a $k$-factorization of $a$. Thus, it
satisfies the four properties of a $k$-factorization demanded in Definition
\ref{def.gcdseq.kfac}. That is:

\begin{enumerate}
\item We have $\left(  b^{\left(  k\right)  }\right)  _{n}=1$ for all $n>k$.

\item We have $\left(  r^{\left(  k\right)  }\right)  _{n}=1$ for all $n\leq
k$.

\item We have $a_{n}=\left(  b^{\left(  k\right)  }\right)
_{\operatorname*{Div}n}\left(  r^{\left(  k\right)  }\right)  _{n}$ for all
$n\geq1$.

\item For any $u,v\geq1$, we have the equality
\[
\left(  b^{\left(  k\right)  }\right)  _{\operatorname*{Div}\ell}\left(
\left(  r^{\left(  k\right)  }\right)  _{u}R+\left(  r^{\left(  k\right)
}\right)  _{v}R\right)  =\left(  b^{\left(  k\right)  }\right)
_{\operatorname*{Div}\ell}\left(  r^{\left(  k\right)  }\right)  _{\gcd\left(
u,v\right)  }R,
\]
where $\ell=\operatorname{lcm}\left(  u,v\right)  $.
\end{enumerate}

We shall only use the second and the third of these properties. The second
property says that $\left(  r^{\left(  k\right)  }\right)  _{n}=1$ for all
$n\leq k$. Applying this to $n=k$, we obtain $\left(  r^{\left(  k\right)
}\right)  _{k}=1$. The third property says that $a_{n}=\left(  b^{\left(
k\right)  }\right)  _{\operatorname*{Div}n}\left(  r^{\left(  k\right)
}\right)  _{n}$ for all $n\geq1$. Applying this to $n=k$, we obtain%
\[
a_{k}=\left(  b^{\left(  k\right)  }\right)  _{\operatorname*{Div}%
k}\underbrace{\left(  r^{\left(  k\right)  }\right)  _{k}}_{=1}=\left(
b^{\left(  k\right)  }\right)  _{\operatorname*{Div}k}=\prod_{d\mid k}\left(
b^{\left(  k\right)  }\right)  _{d}\ \ \ \ \ \ \ \ \ \ \left(  \text{by
(\ref{eq.gcdseq.aDiv})}\right)  .
\]

On the other hand, each positive divisor $d$ of $k$ satisfies $k\geq d$ (since
$d\mid k$) and thus $b_{d}=\left(  b^{\left(  k\right)  }\right)  _{d}$ (by
Claim 1, applied to $n=k$). Thus,%
\[
\prod_{d\mid k}\underbrace{b_{d}}_{\substack{=\left(  b^{\left(  k\right)
}\right)  _{d}}}=\prod_{d\mid k}\left(  b^{\left(  k\right)  }\right)  _{d}.
\]
Comparing this with $a_{k}=\prod_{d\mid k}\left(  b^{\left(  k\right)
}\right)  _{d}$, we obtain $a_{k}=\prod_{d\mid k}b_{d}$. This proves Claim 2.
\end{proof}

Now, the sequence $b=\left(  b_{1},b_{2},b_{3},\ldots\right)  \in R^{\left\{
1,2,3,\ldots\right\}  }$ has the property that each $n\geq1$ satisfies
$a_{n}=\prod_{d\mid n}b_{d}$ (by Claim 2, applied to $k=n$). Hence, a sequence
with this property exists. That is, Theorem \ref{thm.gcdseq.fac} is proved (at last).
\end{proof}

\begin{remark}
Our above proof of Lemma \ref{lem.gcdseq.1} used the
\href{https://en.wikipedia.org/wiki/Axiom_of_dependent_choice}{axiom of
dependent choice}, since we have made infinitely many choices in the
construction of $s$ (recall that we chose an element $s_{n}$ of $R$ that
satisfies $b_{\operatorname*{Div}n}r_{n}=b_{\operatorname*{Div}n}r_{k}s_{n}$;
this element might not be unique). Thus, our proof of Theorem
\ref{thm.gcdseq.fac} depends on this axiom as well. However, this is just an
artifact of the way we have stated these results. To avoid this dependence, we
have to replace the equality (\ref{eq.def.gcdseq.gcdseq.a+a}) in Definition
\ref{def.gcdseq.gcdseq} by an explicit assumption that names \textquotedblleft
witnesses\textquotedblright\ for this equality. That is, instead of requiring
$a_{n}R+a_{m}R=a_{\gcd\left(  n,m\right)  }R$ to hold for all $n,m\geq1$, we
should assume that $\left(  \lambda_{n,m}\right)  _{n,m\geq1}$, $\left(
\mu_{n,m}\right)  _{n,m\geq1}$, $\left(  \nu_{n,m}\right)  _{n,m\geq1}$ and
$\left(  \rho_{n,m}\right)  _{n,m\geq1}$ are four families of elements of $R$
such that all $n,m\geq1$ satisfy
\[
a_{n}=a_{\gcd\left(  n,m\right)  }\lambda_{n,m}\ \ \ \ \ \ \ \ \ \ \text{and}%
\ \ \ \ \ \ \ \ \ \ a_{m}=a_{\gcd\left(  n,m\right)  }\mu_{n,m}%
\]
(these two equalities provide witnesses for $a_{n}R\subseteq a_{\gcd\left(
n,m\right)  }R$ and $a_{m}R\subseteq a_{\gcd\left(  n,m\right)  }R$,
respectively) and%
\[
a_{\gcd\left(  n,m\right)  }=a_{n}\nu_{n,m}+a_{m}\rho_{n,m}%
\]
(this equality provides a witness for $a_{\gcd\left(  n,m\right)  }R\subseteq
a_{n}R+a_{m}R$). Similar changes should be made to the condition
(\ref{eq.gcdseq.k-fac.4}) in Definition \ref{def.gcdseq.kfac}. Once all these
changes are made, the proofs we gave above require no choices.
\end{remark}

\subsubsection{Uniqueness in the regular case}

If the elements $a_{1},a_{2},a_{3},\ldots$ in Theorem \ref{thm.gcdseq.fac} are
regular, then we can obtain two further results:

\begin{proposition}
\label{prop.gcdseq.reg}Let $R$ be a commutative ring. Let $\left(  a_{1}%
,a_{2},a_{3},\ldots\right)  $ be a sequence of regular elements of $R$. Then:

\begin{enumerate}
\item[\textbf{(a)}] There is \textbf{at most one} sequence $\left(
b_{1},b_{2},b_{3},\ldots\right)  $ of elements of $R$ such that each $n\geq1$
satisfies%
\[
a_{n}=\prod_{d\mid n}b_{d}.
\]

\item[\textbf{(b)}] If $\left(  b_{1},b_{2},b_{3},\ldots\right)  $ is such a
sequence, then all its entries $b_{1},b_{2},b_{3},\ldots$ are regular.
\end{enumerate}
\end{proposition}

In particular, Proposition \ref{prop.gcdseq.reg} \textbf{(a)} shows that the
sequence $\left(  b_{1},b_{2},b_{3},\ldots\right)  $ in Theorem
\ref{thm.gcdseq.fac} is uniquely determined if $a_{1},a_{2},a_{3},\ldots$ are regular.

The proof of Proposition \ref{prop.gcdseq.reg} is easy:

\begin{proof}
[Proof of Proposition \ref{prop.gcdseq.reg}.]\textbf{(a)} We must prove that
there is \textbf{at most one} sequence $\left(  b_{1},b_{2},b_{3}%
,\ldots\right)  $ of elements of $R$ such that each $n\geq1$ satisfies%
\begin{equation}
a_{n}=\prod_{d\mid n}b_{d}. \label{pf.prop.gcdseq.reg.a.1}%
\end{equation}

Let $\left(  x_{1},x_{2},x_{3},\ldots\right)  $ and $\left(  y_{1},y_{2}%
,y_{3},\ldots\right)  $ be two such sequences. Thus we must show that $\left(
x_{1},x_{2},x_{3},\ldots\right)  =\left(  y_{1},y_{2},y_{3},\ldots\right)  $.
In other words, we must show that
\begin{equation}
x_{n}=y_{n}\ \ \ \ \ \ \ \ \ \ \text{for each }n\geq1.
\label{pf.prop.gcdseq.reg.a.goal}%
\end{equation}

Recall that $\left(  x_{1},x_{2},x_{3},\ldots\right)  $ is a sequence $\left(
b_{1},b_{2},b_{3},\ldots\right)  $ of elements of $R$ such that each $n\geq1$
satisfies (\ref{pf.prop.gcdseq.reg.a.1}). Thus, each $n\geq1$ satisfies%
\begin{equation}
a_{n}=\prod_{d\mid n}x_{d}=x_{n}\prod_{\substack{d\mid n;\\d\neq n}}x_{d}
\label{pf.prop.gcdseq.reg.a.2}%
\end{equation}
(here, we have split off the factor for $d=n$ from the product). Similarly,
each $n\geq1$ satisfies
\begin{equation}
a_{n}=y_{n}\prod_{\substack{d\mid n;\\d\neq n}}y_{d}.
\label{pf.prop.gcdseq.reg.a.3}%
\end{equation}

Now we shall prove (\ref{pf.prop.gcdseq.reg.a.goal}) by strong induction on
$n$: We fix $n\geq1$, and we assume that%
\begin{equation}
x_{d}=y_{d}\ \ \ \ \ \ \ \ \ \ \text{for all }d<n.
\label{pf.prop.gcdseq.reg.a.4}%
\end{equation}
We must then prove that $x_{n}=y_{n}$.

If $d$ is a positive divisor of $n$ satisfying $d\neq n$, then $d<n$ and
therefore $x_{d}=y_{d}$ (by (\ref{pf.prop.gcdseq.reg.a.4})). Multiplying all
these equalities, we obtain
\[
\prod_{\substack{d\mid n;\\d\neq n}}\underbrace{x_{d}}_{=y_{d}}=\prod
_{\substack{d\mid n;\\d\neq n}}y_{d}.
\]
Set $s:=\prod_{\substack{d\mid n;\\d\neq n}}x_{d}$. Thus,
(\ref{pf.prop.gcdseq.reg.a.2}) rewrites as $a_{n}=x_{n}s$. But
(\ref{pf.prop.gcdseq.reg.a.3}) rewrites as $a_{n}=y_{n}s$ (since
$s=\prod_{\substack{d\mid n;\\d\neq n}}x_{d}=\prod_{\substack{d\mid n;\\d\neq
n}}y_{d}$). Hence, the element $s$ is a divisor of $a_{n}$, and thus is
regular (by Proposition \ref{prop.rings.reg.123} \textbf{(d)}, since $a_{n}$
is regular). Furthermore,
\[
s\left(  x_{n}-y_{n}\right)  =\left(  x_{n}-y_{n}\right)  s=\underbrace{x_{n}%
s}_{=a_{n}}-\underbrace{y_{n}s}_{=a_{n}}=a_{n}-a_{n}=0.
\]
Since $s$ is regular, this entails $x_{n}-y_{n}=0$. In other words,
$x_{n}=y_{n}$. This completes the induction step. Thus,
(\ref{pf.prop.gcdseq.reg.a.goal}) is proved. This completes the proof of
Proposition \ref{prop.gcdseq.reg} \textbf{(a)}. \medskip

\textbf{(b)} Let $\left(  b_{1},b_{2},b_{3},\ldots\right)  $ be such a
sequence. Let $n\geq1$. We must show that $b_{n}$ is regular.

By the assumption on the sequence $\left(  b_{1},b_{2},b_{3},\ldots\right)  $,
we have%
\[
a_{n}=\prod_{d\mid n}b_{d}=b_{n}\prod_{\substack{d\mid n;\\d\neq n}}b_{d}%
\]
(here, we have split off the factor for $d=n$ from the product). Hence,
$b_{n}$ is a divisor of $a_{n}$. Since $a_{n}$ is regular (by assumption), we
thus conclude (by Proposition \ref{prop.rings.reg.123} \textbf{(d)}) that
$b_{n}$ is regular. This proves Proposition \ref{prop.gcdseq.reg} \textbf{(b)}.
\end{proof}

\subsubsection{\label{subsec.gcdseq.frac}The explicit fraction formula for the
factorization}

The sequence $\left(  b_{1},b_{2},b_{3},\ldots\right)  $ in Theorem
\ref{thm.gcdseq.fac} can be expressed explicitly using fractions if the
$a_{1},a_{2},a_{3},\ldots$ are regular. Even without the regularity
assumption, the \textquotedblleft cross-multiplied\textquotedblright\ version
of this expression still holds. We shall state this in Theorem
\ref{thm.gcdseq.fac-pf}; first, we introduce two notations:

\begin{definition}
\label{def.ent.PFn}For any integer $n$, we let $\operatorname*{PF}\left(
n\right)  $ denote the set of all prime divisors of $n$. For example,
$\operatorname*{PF}\left(  20\right)  =\left\{  2,5\right\}  $ and
$\operatorname*{PF}\left(  30\right)  =\left\{  2,3,5\right\}  $. Note that
$\operatorname*{PF}\left(  n\right)  $ is finite unless $n=0$, and is nonempty
unless $\left\vert n\right\vert =1$.
\end{definition}

\begin{definition}
\label{def.ent.prodS}If $S$ is a finite set of positive integers, then
$\operatorname*{prod}S$ denotes the product of all elements of $S$. That is,
$\operatorname*{prod}S=\prod_{s\in S}s$. For instance, $\operatorname*{prod}%
\left(  \left\{  4,7,9\right\}  \right)  =4\cdot7\cdot9$.
\end{definition}

We observe an elementary arithmetical fact:

\begin{lemma}
\label{lem.prod-PF}Let $n$ be an integer, and let $S$ be a finite set of prime
numbers. Then, $S\subseteq\operatorname*{PF}\left(  n\right)  $ if and only if
$\operatorname*{prod}S\mid n$.
\end{lemma}

\begin{fineprint}

\begin{proof}
$\Longrightarrow:$ Assume that $S\subseteq\operatorname*{PF}\left(  n\right)
$. In other words, each $s\in S$ belongs to $\operatorname*{PF}\left(
n\right)  $. Hence, each $s\in S$ divides $n$. Since the elements $s\in S$ are
prime numbers, this entails that their product $\prod_{s\in S}s$ divides $n$
as well (since a product of distinct primes divides an integer if and only if
each of these primes divides the integer). In other words,
$\operatorname*{prod}S$ divides $n$ (since $\operatorname*{prod}S$ is defined
as $\prod_{s\in S}s$). That is, $\operatorname*{prod}S\mid n$. This proves the
\textquotedblleft$\Longrightarrow$\textquotedblright\ direction of the lemma.
\medskip

$\Longleftarrow:$ Assume that $\operatorname*{prod}S\mid n$. Then each $s\in
S$ satisfies $s\mid\operatorname*{prod}S\mid n$ and therefore $s\in
\operatorname*{PF}\left(  n\right)  $ (since $s$ is prime, being an element of
$S$). That is, $S\subseteq\operatorname*{PF}\left(  n\right)  $. Hence, the
\textquotedblleft$\Longleftarrow$\textquotedblright\ direction of the lemma is proved.
\end{proof}
\end{fineprint}

\begin{theorem}
\label{thm.gcdseq.fac-pf}Let $R$ be a commutative ring. Let $\left(
a_{1},a_{2},a_{3},\ldots\right)  $ and $\left(  b_{1},b_{2},b_{3}%
,\ldots\right)  $ be two sequences of elements of $R$ such that each $n\geq1$
satisfies%
\begin{equation}
a_{n}=\prod_{d\mid n}b_{d}. \label{eq.thm.gcdseq.fac-pf.ass}%
\end{equation}
Let $n\geq1$. Then:

\begin{enumerate}
\item[\textbf{(a)}] We have%
\begin{equation}
b_{n}\prod_{\substack{S\subseteq\operatorname*{PF}\left(  n\right)
;\\\left\vert S\right\vert \text{ is odd}}}a_{n/\operatorname*{prod}S}%
=\prod_{\substack{S\subseteq\operatorname*{PF}\left(  n\right)  ;\\\left\vert
S\right\vert \text{ is even}}}a_{n/\operatorname*{prod}S}.
\label{eq.thm.gcdseq.fac-pf.a}%
\end{equation}

\item[\textbf{(b)}] If the element $a_{n}$ is regular in $R$, then (using the
notation of Definition \ref{def.rings.reg.div}) we have%
\begin{equation}
b_{n}=\dfrac{\prod\limits_{\substack{S\subseteq\operatorname*{PF}\left(
n\right)  ;\\\left\vert S\right\vert \text{ is even}}%
}a_{n/\operatorname*{prod}S}}{\prod\limits_{\substack{S\subseteq
\operatorname*{PF}\left(  n\right)  ;\\\left\vert S\right\vert \text{ is odd}%
}}a_{n/\operatorname*{prod}S}}. \label{eq.thm.gcdseq.fac-pf.b}%
\end{equation}

\end{enumerate}
\end{theorem}

\begin{example}
Let $R$, $\left(  a_{1},a_{2},a_{3},\ldots\right)  $ and $\left(  b_{1}%
,b_{2},b_{3},\ldots\right)  $ be as in Theorem \ref{thm.gcdseq.fac-pf}. Then,
$\operatorname*{PF}\left(  20\right)  =\left\{  2,5\right\}  $, so that
Theorem \ref{thm.gcdseq.fac-pf} \textbf{(a)} (applied to $n=20$) yields%
\[
b_{20}a_{20/\operatorname*{prod}\left\{  2\right\}  }%
a_{20/\operatorname*{prod}\left\{  5\right\}  }=a_{20/\operatorname*{prod}%
\varnothing}a_{20/\operatorname*{prod}\left\{  2,5\right\}  }.
\]
In other words,%
\[
b_{20}a_{10}a_{4}=a_{20}a_{2}.
\]
Likewise, Theorem \ref{thm.gcdseq.fac-pf} \textbf{(a)} (applied to $n=30$)
yields%
\[
b_{30}a_{15}a_{10}a_{6}a_{1}=a_{30}a_{2}a_{3}a_{5}.
\]

\end{example}

To prove Theorem \ref{thm.gcdseq.fac-pf}, we need a simple combinatorial lemma
(cf. \cite[\S 6.1]{21s}):

\begin{lemma}
\label{lem.gcdseq.fac-pf.lem-comb}Let $V$ be a finite nonempty set. Then,%
\begin{align*}
&  \left(  \text{number of all }S\subseteq V\text{ such that }\left\vert
S\right\vert \text{ is odd}\right) \\
&  =\left(  \text{number of all }S\subseteq V\text{ such that }\left\vert
S\right\vert \text{ is even}\right)  .
\end{align*}

\end{lemma}

\begin{proof}
We can break up the sum $\sum_{S\subseteq V}\left(  -1\right)  ^{\left\vert
S\right\vert }$ as follows:%
\begin{align}
&  \sum_{S\subseteq V}\left(  -1\right)  ^{\left\vert S\right\vert
}\nonumber\\
&  =\sum_{\substack{S\subseteq V;\\\left\vert S\right\vert \text{ is even}%
}}\underbrace{\left(  -1\right)  ^{\left\vert S\right\vert }}%
_{\substack{=1\\\text{(since }\left\vert S\right\vert \text{ is even)}}%
}+\sum_{\substack{S\subseteq V;\\\left\vert S\right\vert \text{ is odd}%
}}\underbrace{\left(  -1\right)  ^{\left\vert S\right\vert }}%
_{\substack{=-1\\\text{(since }\left\vert S\right\vert \text{ is odd)}%
}}\nonumber\\
&  =\underbrace{\sum_{\substack{S\subseteq V;\\\left\vert S\right\vert \text{
is even}}}1}_{\substack{=\left(  \text{number of all }S\subseteq V\text{ such
that }\left\vert S\right\vert \text{ is even}\right)  \cdot1\\=\left(
\text{number of all }S\subseteq V\text{ such that }\left\vert S\right\vert
\text{ is even}\right)  }}-\underbrace{\sum_{\substack{S\subseteq
V;\\\left\vert S\right\vert \text{ is odd}}}1}_{\substack{=\left(
\text{number of all }S\subseteq V\text{ such that }\left\vert S\right\vert
\text{ is odd}\right)  \cdot1\\=\left(  \text{number of all }S\subseteq
V\text{ such that }\left\vert S\right\vert \text{ is odd}\right)
}}\nonumber\\
&  =\left(  \text{number of all }S\subseteq V\text{ such that }\left\vert
S\right\vert \text{ is even}\right) \nonumber\\
&  \ \ \ \ \ \ \ \ \ \ -\left(  \text{number of all }S\subseteq V\text{ such
that }\left\vert S\right\vert \text{ is odd}\right)  .
\label{pf.lem.gcdseq.fac-pf.lem-comb.1}%
\end{align}
Now, we shall show that the left-hand side of this equality is $0$ (because
this will entail that its right-hand side is also $0$, and thus the lemma will
easily follow). For this purpose, we shall break up the sum $\sum_{S\subseteq
V}\left(  -1\right)  ^{\left\vert S\right\vert }$ in a second way:

Since $V$ is nonempty, we can pick some $v\in V$. Consider this $v$. Then,%
\begin{equation}
\sum_{S\subseteq V}\left(  -1\right)  ^{\left\vert S\right\vert }%
=\sum_{\substack{S\subseteq V;\\v\in S}}\left(  -1\right)  ^{\left\vert
S\right\vert }+\sum_{\substack{S\subseteq V;\\v\notin S}}\left(  -1\right)
^{\left\vert S\right\vert }. \label{pf.lem.gcdseq.fac-pf.lem-comb.5}%
\end{equation}

But we can define two maps%
\begin{align*}
\alpha:\left\{  S\subseteq V\ \mid\ v\notin S\right\}   &  \rightarrow\left\{
S\subseteq V\ \mid\ v\in S\right\}  ,\\
S  &  \mapsto S\cup\left\{  v\right\}
\end{align*}
and%
\begin{align*}
\rho:\left\{  S\subseteq V\ \mid\ v\in S\right\}   &  \rightarrow\left\{
S\subseteq V\ \mid\ v\notin S\right\}  ,\\
S  &  \mapsto S\setminus\left\{  v\right\}  .
\end{align*}
The first of these maps inserts the element $v$ into any subset of $V$ that
does not contain $v$; the second map removes the element $v$ from any subset
of $V$ that does contain $v$. Thus, the two maps $\alpha$ and $\rho$ are
mutually inverse, and therefore are bijections. Hence, we obtain
\begin{equation}
\sum_{\substack{S\subseteq V;\\v\notin S}}\left(  -1\right)  ^{\left\vert
S\right\vert }=\sum_{\substack{S\subseteq V;\\v\in S}}\left(  -1\right)
^{\left\vert \rho\left(  S\right)  \right\vert }
\label{pf.lem.gcdseq.fac-pf.lem-comb.2}%
\end{equation}
(here, we have substituted $\rho\left(  S\right)  $ for $S$ in the sum, since
the map $\rho$ is a bijection).

But each subset $S$ of $V$ satisfying $v\in S$ satisfies $\left\vert
\rho\left(  S\right)  \right\vert =\left\vert S\setminus\left\{  v\right\}
\right\vert =\left\vert S\right\vert -1$ (since $v\in S$) and therefore%
\begin{equation}
\left(  -1\right)  ^{\left\vert \rho\left(  S\right)  \right\vert }=\left(
-1\right)  ^{\left\vert S\right\vert -1}=-\left(  -1\right)  ^{\left\vert
S\right\vert }. \label{pf.lem.gcdseq.fac-pf.lem-comb.3}%
\end{equation}
Hence, (\ref{pf.lem.gcdseq.fac-pf.lem-comb.2}) rewrites as%
\begin{equation}
\sum_{\substack{S\subseteq V;\\v\notin S}}\left(  -1\right)  ^{\left\vert
S\right\vert }=\sum_{\substack{S\subseteq V;\\v\in S}}\left(  -\left(
-1\right)  ^{\left\vert S\right\vert }\right)  =-\sum_{\substack{S\subseteq
V;\\v\in S}}\left(  -1\right)  ^{\left\vert S\right\vert }.
\label{pf.lem.gcdseq.fac-pf.lem-comb.4}%
\end{equation}
Substituting this into (\ref{pf.lem.gcdseq.fac-pf.lem-comb.5}), we find%
\[
\sum_{S\subseteq V}\left(  -1\right)  ^{\left\vert S\right\vert }%
=\sum_{\substack{S\subseteq V;\\v\in S}}\left(  -1\right)  ^{\left\vert
S\right\vert }+\left(  -\sum_{\substack{S\subseteq V;\\v\in S}}\left(
-1\right)  ^{\left\vert S\right\vert }\right)  =0.
\]
Comparing this with (\ref{pf.lem.gcdseq.fac-pf.lem-comb.1}), we obtain%
\begin{align*}
&  \left(  \text{number of all }S\subseteq V\text{ such that }\left\vert
S\right\vert \text{ is even}\right) \\
&  \ \ \ \ \ \ \ \ \ \ -\left(  \text{number of all }S\subseteq V\text{ such
that }\left\vert S\right\vert \text{ is odd}\right) \\
&  =0.
\end{align*}
In other words,
\begin{align*}
&  \left(  \text{number of all }S\subseteq V\text{ such that }\left\vert
S\right\vert \text{ is odd}\right) \\
&  =\left(  \text{number of all }S\subseteq V\text{ such that }\left\vert
S\right\vert \text{ is even}\right)  .
\end{align*}
This proves the lemma.
\end{proof}

\begin{proof}
[Proof of Theorem \ref{thm.gcdseq.fac-pf}.]\textbf{(a)} First, we observe that
each subset $S$ of $\operatorname*{PF}\left(  n\right)  $ is a finite set of
prime numbers (since $\operatorname*{PF}\left(  n\right)  $ consists of prime
numbers) and therefore satisfies $\operatorname*{prod}S\mid n$ (by Lemma
\ref{lem.prod-PF}, since $S\subseteq\operatorname*{PF}\left(  n\right)  $) and
thus $n/\operatorname*{prod}S\in\left\{  1,2,3,\ldots\right\}  $. Hence, all
the expressions $a_{n/\operatorname*{prod}S}$ in (\ref{eq.thm.gcdseq.fac-pf.a}%
) are well-defined. It remains to prove the equality
(\ref{eq.thm.gcdseq.fac-pf.a}).

For each subset $S$ of $\operatorname*{PF}\left(  n\right)  $, we have%
\begin{equation}
a_{n/\operatorname*{prod}S}=\prod_{d\mid n/\operatorname*{prod}S}b_{d}
\label{pf.thm.gcdseq.fac-pf.a.1}%
\end{equation}
(by (\ref{eq.thm.gcdseq.fac-pf.ass}), applied to $n/\operatorname*{prod}S$
instead of $n$). However, we observe another simple arithmetical fact:

\begin{statement}
\textit{Claim 1:} Let $S$ be a subset of $\operatorname*{PF}\left(  n\right)
$. Let $d$ be a positive integer. Then, we have the equivalence
\[
\left(  d\mid n/\operatorname*{prod}S\right)  \ \Longleftrightarrow\ \left(
d\mid n\text{ and }S\subseteq\operatorname*{PF}\left(  n/d\right)  \right)  .
\]

\end{statement}

\begin{fineprint}

\begin{proof}
[Proof of Claim 1.]$\Longrightarrow:$ Assume that $d\mid
n/\operatorname*{prod}S$. Then, $d\mid n/\operatorname*{prod}S\mid n$, so that
$n/d$ is a positive integer. Furthermore, $\dfrac{n/d}{\operatorname*{prod}%
S}=\dfrac{n/\operatorname*{prod}S}{d}\in\mathbb{Z}$ (since $d\mid
n/\operatorname*{prod}S$). In other words, $\operatorname*{prod}S\mid n/d$.
But $S$ is a finite set of prime numbers (since it is a subset of
$\operatorname*{PF}\left(  n\right)  $, which is a finite set of prime
numbers). Thus, Lemma \ref{lem.prod-PF} (applied to $n/d$ instead of $n$)
shows that $S\subseteq\operatorname*{PF}\left(  n/d\right)  $ if and only if
$\operatorname*{prod}S\mid n/d$. Hence, $S\subseteq\operatorname*{PF}\left(
n/d\right)  $ (since $\operatorname*{prod}S\mid n/d$). Thus we have shown that
$d\mid n$ and $S\subseteq\operatorname*{PF}\left(  n/d\right)  $. This proves
the \textquotedblleft$\Longrightarrow$\textquotedblright\ direction of Claim
1. \medskip

$\Longleftarrow:$ Assume that $d\mid n$ and $S\subseteq\operatorname*{PF}%
\left(  n/d\right)  $. Then, $n/d$ is a positive integer (since $d\mid n$).
Also, $S$ is a finite set of prime numbers (since it is a subset of
$\operatorname*{PF}\left(  n/d\right)  $, which is a finite set of prime
numbers). Thus, Lemma \ref{lem.prod-PF} (applied to $n/d$ instead of $n$)
shows that $S\subseteq\operatorname*{PF}\left(  n/d\right)  $ if and only if
$\operatorname*{prod}S\mid n/d$. Hence, $\operatorname*{prod}S\mid n/d$ (since
$S\subseteq\operatorname*{PF}\left(  n/d\right)  $). Hence,
$\operatorname*{prod}S\mid n/d\mid n$, so that $n/\operatorname*{prod}S$ is a
positive integer.

Furthermore, $\dfrac{n/\operatorname*{prod}S}{d}=\dfrac{n/d}%
{\operatorname*{prod}S}\in\mathbb{Z}$ (since $\operatorname*{prod}S\mid n/d$),
so that $d\mid n/\operatorname*{prod}S$. This proves the \textquotedblleft%
$\Longleftarrow$\textquotedblright\ direction of Claim 1.
\end{proof}
\end{fineprint}

Now, if $S$ is a subset of $\operatorname*{PF}\left(  n\right)  $, then Claim
1 shows that%
\[
\prod_{d\mid n/\operatorname*{prod}S}b_{d}=\prod_{\substack{d\mid
n;\\S\subseteq\operatorname*{PF}\left(  n/d\right)  }}b_{d}.
\]
This lets us rewrite (\ref{pf.thm.gcdseq.fac-pf.a.1}) as%
\begin{equation}
a_{n/\operatorname*{prod}S}=\prod_{\substack{d\mid n;\\S\subseteq
\operatorname*{PF}\left(  n/d\right)  }}b_{d}.
\label{pf.thm.gcdseq.fac-pf.a.3}%
\end{equation}

Hence,%
\begin{align}
\prod_{\substack{S\subseteq\operatorname*{PF}\left(  n\right)  ;\\\left\vert
S\right\vert \text{ is odd}}}a_{n/\operatorname*{prod}S}  &  =\prod
_{\substack{S\subseteq\operatorname*{PF}\left(  n\right)  ;\\\left\vert
S\right\vert \text{ is odd}}}\ \ \prod_{\substack{d\mid n;\\S\subseteq
\operatorname*{PF}\left(  n/d\right)  }}b_{d}\nonumber\\
&  =\prod_{d\mid n}\ \ \prod_{\substack{S\subseteq\operatorname*{PF}\left(
n\right)  ;\\\left\vert S\right\vert \text{ is odd;}\\S\subseteq
\operatorname*{PF}\left(  n/d\right)  }}b_{d} \label{pf.thm.gcdseq.fac-pf.a.4}%
\end{align}
(here, we have interchanged the two product signs). However, if $d$ is a
positive divisor of $n$, then $\operatorname*{PF}\left(  n/d\right)
\subseteq\operatorname*{PF}\left(  n\right)  $ (since any prime divisor $p$ of
$n/d$ will satisfy $p\mid n/d\mid n$ and thus be a prime divisor of $n$ as
well), and thus every subset $S\subseteq\operatorname*{PF}\left(  n/d\right)
$ will automatically satisfy $S\subseteq\operatorname*{PF}\left(  n\right)  $.
Hence, the condition \textquotedblleft$S\subseteq\operatorname*{PF}\left(
n\right)  $\textquotedblright\ under the product sign $\prod
_{\substack{S\subseteq\operatorname*{PF}\left(  n\right)  ;\\\left\vert
S\right\vert \text{ is odd;}\\S\subseteq\operatorname*{PF}\left(  n/d\right)
}}$ in (\ref{pf.thm.gcdseq.fac-pf.a.4}) is made redundant by the condition
\textquotedblleft$S\subseteq\operatorname*{PF}\left(  n/d\right)
$\textquotedblright. We can therefore replace this product sign by
$\prod_{\substack{S\subseteq\operatorname*{PF}\left(  n/d\right)
;\\\left\vert S\right\vert \text{ is odd}}}$. Hence,
(\ref{pf.thm.gcdseq.fac-pf.a.4}) rewrites as%
\begin{align}
\prod_{\substack{S\subseteq\operatorname*{PF}\left(  n\right)  ;\\\left\vert
S\right\vert \text{ is odd}}}a_{n/\operatorname*{prod}S}  &  =\prod_{d\mid
n}\ \ \underbrace{\prod_{\substack{S\subseteq\operatorname*{PF}\left(
n/d\right)  ;\\\left\vert S\right\vert \text{ is odd}}}b_{d}}_{=b_{d}^{\left(
\text{number of all }S\subseteq\operatorname*{PF}\left(  n/d\right)  \text{
such that }\left\vert S\right\vert \text{ is odd}\right)  }}\nonumber\\
&  =\prod_{d\mid n}b_{d}^{\left(  \text{number of all }S\subseteq
\operatorname*{PF}\left(  n/d\right)  \text{ such that }\left\vert
S\right\vert \text{ is odd}\right)  }\nonumber\\
&  =b_{n}^{\left(  \text{number of all }S\subseteq\operatorname*{PF}\left(
n/n\right)  \text{ such that }\left\vert S\right\vert \text{ is odd}\right)
}\nonumber\\
&  \ \ \ \ \ \ \ \ \ \ \cdot\prod_{\substack{d\mid n;\\d\neq n}}b_{d}^{\left(
\text{number of all }S\subseteq\operatorname*{PF}\left(  n/d\right)  \text{
such that }\left\vert S\right\vert \text{ is odd}\right)  }%
\ \ \ \ \ \ \ \ \ \ \label{pf.thm.gcdseq.fac-pf.a.5}%
\end{align}
(here, we have split off the factor for $d=n$ from the product). The same
argument (but with the word \textquotedblleft odd\textquotedblright\ replaced
by \textquotedblleft even\textquotedblright) shows that%
\begin{align}
\prod_{\substack{S\subseteq\operatorname*{PF}\left(  n\right)  ;\\\left\vert
S\right\vert \text{ is even}}}a_{n/\operatorname*{prod}S}  &  =b_{n}^{\left(
\text{number of all }S\subseteq\operatorname*{PF}\left(  n/n\right)  \text{
such that }\left\vert S\right\vert \text{ is even}\right)  }\nonumber\\
&  \ \ \ \ \ \ \ \ \ \ \cdot\prod_{\substack{d\mid n;\\d\neq n}}b_{d}^{\left(
\text{number of all }S\subseteq\operatorname*{PF}\left(  n/d\right)  \text{
such that }\left\vert S\right\vert \text{ is even}\right)  }%
.\ \ \ \ \ \ \ \ \ \ \label{pf.thm.gcdseq.fac-pf.a.6}%
\end{align}

However, if $d$ is a positive divisor of $n$ satisfying $d\neq n$, then the
finite set $\operatorname*{PF}\left(  n/d\right)  $ is
nonempty\footnote{\textit{Proof.} Let $d$ be a positive divisor of $n$
satisfying $d\neq n$. Then, $n/d$ is a positive divisor of $n$ satisfying
$n/d\neq1$ (since $d\neq n$), and thus has at least one prime divisor (since
any positive integer $\neq1$ has at least one prime divisor). In other words,
the set $\operatorname*{PF}\left(  n/d\right)  $ is nonempty.}, and therefore
we have%
\begin{align}
&  \left(  \text{number of all }S\subseteq\operatorname*{PF}\left(
n/d\right)  \text{ such that }\left\vert S\right\vert \text{ is odd}\right)
\nonumber\\
&  =\left(  \text{number of all }S\subseteq\operatorname*{PF}\left(
n/d\right)  \text{ such that }\left\vert S\right\vert \text{ is even}\right)
\label{pf.thm.gcdseq.fac-pf.a.7}%
\end{align}
(by Lemma \ref{lem.gcdseq.fac-pf.lem-comb}, applied to $V=\operatorname*{PF}%
\left(  n/d\right)  $). Moreover, the set $\operatorname*{PF}\left(
n/n\right)  =\operatorname*{PF}\left(  1\right)  =\varnothing$ is empty, and
thus has only one subset $S$, namely itself; the size $\left\vert S\right\vert
$ of this one subset is even (in fact, it is $0$). Hence,%
\begin{equation}
\left(  \text{number of all }S\subseteq\operatorname*{PF}\left(  n/n\right)
\text{ such that }\left\vert S\right\vert \text{ is even}\right)  =1
\label{pf.thm.gcdseq.fac-pf.a.8}%
\end{equation}
and%
\begin{equation}
\left(  \text{number of all }S\subseteq\operatorname*{PF}\left(  n/n\right)
\text{ such that }\left\vert S\right\vert \text{ is odd}\right)  =0.
\label{pf.thm.gcdseq.fac-pf.a.9}%
\end{equation}
Now, (\ref{pf.thm.gcdseq.fac-pf.a.5}) becomes%
\begin{align}
\prod_{\substack{S\subseteq\operatorname*{PF}\left(  n\right)  ;\\\left\vert
S\right\vert \text{ is odd}}}a_{n/\operatorname*{prod}S}  &
=\underbrace{b_{n}^{\left(  \text{number of all }S\subseteq\operatorname*{PF}%
\left(  n/n\right)  \text{ such that }\left\vert S\right\vert \text{ is
odd}\right)  }}_{\substack{=b_{n}^{0}\\\text{(by
(\ref{pf.thm.gcdseq.fac-pf.a.9}))}}}\nonumber\\
&  \ \ \ \ \ \ \ \ \ \ \cdot\prod_{\substack{d\mid n;\\d\neq n}%
}\underbrace{b_{d}^{\left(  \text{number of all }S\subseteq\operatorname*{PF}%
\left(  n/d\right)  \text{ such that }\left\vert S\right\vert \text{ is
odd}\right)  }}_{\substack{=b_{d}^{\left(  \text{number of all }%
S\subseteq\operatorname*{PF}\left(  n/d\right)  \text{ such that }\left\vert
S\right\vert \text{ is even}\right)  }\\\text{(by
(\ref{pf.thm.gcdseq.fac-pf.a.7}))}}}\nonumber\\
&  =\underbrace{b_{n}^{0}}_{=1}\cdot\prod_{\substack{d\mid n;\\d\neq n}%
}b_{d}^{\left(  \text{number of all }S\subseteq\operatorname*{PF}\left(
n/d\right)  \text{ such that }\left\vert S\right\vert \text{ is even}\right)
}\nonumber\\
&  =\prod_{\substack{d\mid n;\\d\neq n}}b_{d}^{\left(  \text{number of all
}S\subseteq\operatorname*{PF}\left(  n/d\right)  \text{ such that }\left\vert
S\right\vert \text{ is even}\right)  }. \label{pf.thm.gcdseq.fac-pf.a.10}%
\end{align}
Meanwhile, (\ref{pf.thm.gcdseq.fac-pf.a.6}) becomes%
\begin{align*}
\prod_{\substack{S\subseteq\operatorname*{PF}\left(  n\right)  ;\\\left\vert
S\right\vert \text{ is even}}}a_{n/\operatorname*{prod}S}  &
=\underbrace{b_{n}^{\left(  \text{number of all }S\subseteq\operatorname*{PF}%
\left(  n/n\right)  \text{ such that }\left\vert S\right\vert \text{ is
even}\right)  }}_{\substack{=b_{n}^{1}\\\text{(by
(\ref{pf.thm.gcdseq.fac-pf.a.8}))}}}\\
&  \ \ \ \ \ \ \ \ \ \ \cdot\prod_{\substack{d\mid n;\\d\neq n}}b_{d}^{\left(
\text{number of all }S\subseteq\operatorname*{PF}\left(  n/d\right)  \text{
such that }\left\vert S\right\vert \text{ is even}\right)  }\\
&  =\underbrace{b_{n}^{1}}_{=b_{n}}\cdot\prod_{\substack{d\mid n;\\d\neq
n}}b_{d}^{\left(  \text{number of all }S\subseteq\operatorname*{PF}\left(
n/d\right)  \text{ such that }\left\vert S\right\vert \text{ is even}\right)
}\\
&  =b_{n}\cdot\underbrace{\prod_{\substack{d\mid n;\\d\neq n}}b_{d}^{\left(
\text{number of all }S\subseteq\operatorname*{PF}\left(  n/d\right)  \text{
such that }\left\vert S\right\vert \text{ is even}\right)  }}%
_{\substack{=\prod_{\substack{S\subseteq\operatorname*{PF}\left(  n\right)
;\\\left\vert S\right\vert \text{ is odd}}}a_{n/\operatorname*{prod}%
S}\\\text{(by (\ref{pf.thm.gcdseq.fac-pf.a.10}))}}}\\
&  =b_{n}\prod_{\substack{S\subseteq\operatorname*{PF}\left(  n\right)
;\\\left\vert S\right\vert \text{ is odd}}}a_{n/\operatorname*{prod}S}.
\end{align*}
This proves (\ref{eq.thm.gcdseq.fac-pf.a}). Thus, Theorem
\ref{thm.gcdseq.fac-pf} \textbf{(a)} is proved. \medskip

\textbf{(b)} For each positive divisor $m$ of $n$, we have%
\begin{equation}
a_{m}\mid a_{n} \label{pf.thm.gcdseq.fac-pf.b.div}%
\end{equation}
\footnote{\textit{Proof:} Let $m$ be a positive divisor of $n$. Then, each
positive divisor of $m$ is a positive divisor of $n$ as well. Hence, the
product $\prod_{d\mid m}b_{d}$ is a subproduct of the product $\prod_{d\mid
n}b_{d}$. Therefore, it divides the latter product. In other words,
$\prod_{d\mid m}b_{d}\mid\prod_{d\mid n}b_{d}$.
\par
But (\ref{eq.thm.gcdseq.fac-pf.ass}) yields $a_{n}=\prod_{d\mid n}b_{d}$ and
$a_{m}=\prod_{d\mid m}b_{d}$. Hence, $a_{m}=\prod_{d\mid m}b_{d}\mid
\prod_{d\mid n}b_{d}=a_{n}$. Thus, (\ref{pf.thm.gcdseq.fac-pf.b.div}) is
proved.}.

Now assume that $a_{n}$ is regular in $R$. Then, the element
$a_{n/\operatorname*{prod}S}$ is regular for each $S\subseteq
\operatorname*{PF}\left(  n\right)  $\ \ \ \ \footnote{\textit{Proof.} Let
$S\subseteq\operatorname*{PF}\left(  n\right)  $. Then, $\operatorname*{prod}%
S\mid n$ (by Lemma \ref{lem.prod-PF}), and therefore $n/\operatorname*{prod}S$
is a positive divisor of $n$. Hence, (\ref{pf.thm.gcdseq.fac-pf.b.div})
(applied to $m=n/\operatorname*{prod}S$) yields $a_{n/\operatorname*{prod}%
S}\mid a_{n}$. Hence, $a_{n/\operatorname*{prod}S}$ is a divisor of the
regular element $a_{n}$ of $R$. By Proposition \ref{prop.rings.reg.123}
\textbf{(d)}, this entails that $a_{n/\operatorname*{prod}S}$ is itself
regular. Qed.}. Hence, the product $\prod\limits_{\substack{S\subseteq
\operatorname*{PF}\left(  n\right)  ;\\\left\vert S\right\vert \text{ is odd}%
}}a_{n/\operatorname*{prod}S}$ is regular as well (by Proposition
\ref{prop.rings.reg.123} \textbf{(c)}). This allows us to rewrite
(\ref{eq.thm.gcdseq.fac-pf.a}) as (\ref{eq.thm.gcdseq.fac-pf.b}). Thus,
Theorem \ref{thm.gcdseq.fac-pf} \textbf{(b)} is proved.
\end{proof}

Note that Theorem \ref{thm.gcdseq.fac-pf} \textbf{(b)} proves Proposition
\ref{prop.gcdseq.reg} \textbf{(a)} anew.

\newpage

\section{\label{chp.modules}Modules (\cite[Chapter 10]{DumFoo04})}

We now move on from studying rings themselves to studying \textbf{modules}
over rings. In many ways, modules are even more important than rings, as their
definition offers more freedom (and this freedom is widely used throughout
mathematics). Some would argue that the notion of a ring is merely an
ancillary character to that of a module.

\subsection{\label{sec.modules.def}Definition and examples (\cite[\S 10.1]%
{DumFoo04})}

Before we define modules rigorously, let me give a rough idea of what they
stand for.

We can think of a ring as a system of \textquotedblleft number-like
objects\textquotedblright\ that can be \textquotedblleft
added\textquotedblright\ (with one another) and \textquotedblleft
multiplied\textquotedblright\ (with one another).\footnote{Here I am leaving
the zero and the unity unmentioned, for the sake of brevity.}

In contrast, a \textbf{module} (over a given ring $R$) can be thought of as a
system of \textquotedblleft vector-like objects\textquotedblright\ that can be
\textquotedblleft added\textquotedblright\ (with one another) and
\textquotedblleft scaled\textquotedblright\ (by elements of $R$). In
particular, if $R$ is a field, then the modules over $R$ are just the vector
spaces over $R$ (as defined in any textbook on abstract linear algebra).

To turn this into a proper definition of a module, we just need to decide what
properties of \textquotedblleft adding\textquotedblright\ and
\textquotedblleft scaling\textquotedblright\ we want to require as axioms.

\subsubsection{\label{subsec.modules.def.def}Definition of modules}

For every ring $R$, there are two notions of an \textquotedblleft%
$R$-module\textquotedblright: the \textquotedblleft left $R$%
-modules\textquotedblright\ and the \textquotedblleft right $R$%
-modules\textquotedblright. Let us first define the left ones:

\begin{definition}
\label{def.mod.leftmod}Let $R$ be a ring. A \textbf{left }$R$\textbf{-module}
(or a \textbf{left module over }$R$) means a set $M$ equipped with

\begin{itemize}
\item a binary operation $+$ (that is, a map from $M\times M$ to $M$) that is
called \textbf{addition};

\item an element $0_{M}\in M$ that is called the \textbf{zero element} or the
\textbf{zero vector} or just the \textbf{zero}, and is just denoted by $0$
when there is no ambiguity;

\item a map from $R\times M$ to $M$ that is called the \textbf{action of }$R$
\textbf{on }$M$, and is written as multiplication (i.e., we denote the image
of a pair $\left(  r,m\right)  \in R\times M$ under this map by $rm$ or
$r\cdot m$)
\end{itemize}

\noindent such that the following properties (the \textquotedblleft%
\textbf{module axioms}\textquotedblright) hold:

\begin{itemize}
\item $\left(  M,+,0\right)  $ is an abelian group.

\item The \textbf{right distributivity law} holds: We have $\left(
r+s\right)  m=rm+sm$ for all $r,s\in R$ and $m\in M$.

\item The \textbf{left distributivity law} holds: We have $r\left(
m+n\right)  =rm+rn$ for all $r\in R$ and $m,n\in M$.

\item The \textbf{associativity law} holds: We have $\left(  rs\right)
m=r\left(  sm\right)  $ for all $r,s\in R$ and $m\in M$.

\item We have $0_{R}m=0_{M}$ for every $m\in M$.

\item We have $r\cdot0_{M}=0_{M}$ for every $r\in R$.

\item We have $1m=m$ for every $m\in M$. (Here, \textquotedblleft%
$1$\textquotedblright\ means the unity of $R$.)
\end{itemize}

When $M$ is a left $R$-module, the elements of $M$ are called \textbf{vectors}%
, and the elements of $R$ are called \textbf{scalars}.
\end{definition}

As the name \textquotedblleft left $R$-module\textquotedblright\ suggests,
there is an analogous notion of a \textbf{right }$R$\textbf{-module}:

\begin{definition}
\label{def.mod.rimod}Let $R$ be a ring. A \textbf{right }$R$\textbf{-module}
is defined just as a left $R$-module was defined in Definition
\ref{def.mod.leftmod}, but with the following changes:

\begin{itemize}
\item For a right $R$-module $M$, the action is not a map from $R\times M$ to
$M$, but rather a map from $M\times R$ to $M$.

\item Accordingly, we use the notation $mr$ (rather than $rm$) for the image
of a pair $\left(  m,r\right)  $ under this map.

\item The axioms for a right $R$-module are similar to the module axioms for a
left $R$-module, accounting for the different form of the action. For example,
the associativity law for a right $R$-module is saying that $m\left(
rs\right)  =\left(  mr\right)  s$ for all $r,s\in R$ and $m\in M$.
\end{itemize}
\end{definition}

When $R$ is commutative, any left $R$-module becomes a right $R$-module in a
natural way:

\begin{proposition}
\label{prop.mod.leftmod-is-rightmod}Let $R$ be a commutative ring. Then, we
can make any left $R$-module $M$ into a right $R$-module by setting%
\begin{equation}
mr=rm\ \ \ \ \ \ \ \ \ \ \text{for all }r\in R\text{ and }m\in M.
\label{eq.prop.mod.leftmod-is-rightmod.mr=rm}%
\end{equation}
Similarly, we can make any right $R$-module $M$ into a left $R$-module by
setting%
\begin{equation}
rm=mr\ \ \ \ \ \ \ \ \ \ \text{for all }r\in R\text{ and }m\in M.
\label{eq.prop.mod.leftmod-is-rightmod.rm=mr}%
\end{equation}
These two transformations are mutually inverse, so we shall use them to
identify left $R$-modules with right $R$-modules. Thus, we can use the words
\textquotedblleft left $R$-module\textquotedblright\ and \textquotedblleft
right $R$-module\textquotedblright\ interchangeably, and just speak of
\textquotedblleft$R$\textbf{-modules}\textquotedblright\ instead (without
specifying whether they are left or right). We shall liberally do so in what
follows. (Note that this is not allowed when $R$ is not commutative!)
\end{proposition}

When $R$ is a field, the $R$-modules are also known as the $R$\textbf{-vector
spaces}. These are precisely the vector spaces you have seen in a linear
algebra class. A left $R$-module over an arbitrary ring $R$ is just the
natural generalization of a vector space. But while vector spaces have a very
predictable structure (in particular, a vector space is uniquely determined up
to isomorphism by its dimension), modules can be wild (although the
\textquotedblleft nice\textquotedblright\ families of modules, like $R^{n}$
for $n\in\mathbb{N}$, still exist for every ring). The wilder a ring is, the
more diverse are its modules.

One more remark about Definition \ref{def.mod.leftmod}: The \textquotedblleft%
$0_{R}m=0_{M}$\textquotedblright\ and \textquotedblleft$r\cdot0_{M}=0_{M}%
$\textquotedblright\ axioms are actually redundant (i.e., they follow from the
other axioms). I leave it to you to check this.

\begin{exercise}
Check this!
\end{exercise}

\subsubsection{Submodules and scaling}

We will soon see some examples of $R$-modules; but let us first define
$R$-submodules. If you have seen subspaces of a vector space, this definition
won't surprise you:

\begin{definition}
\label{def.mod.submod}Let $R$ be a ring. Let $M$ be a left $R$-module. An
$R$\textbf{-submodule} (or, to be more precise, a \textbf{left }%
$R$\textbf{-submodule}) of $M$ means a subset $N$ of $M$ such that

\begin{itemize}
\item we have $a+b\in N$ for any $a,b\in N$;

\item we have $ra\in N$ for any $r\in R$ and $a\in N$;

\item we have $0\in N$ (where $0$ means $0_{M}$).
\end{itemize}
\end{definition}

All three axioms in Definition \ref{def.mod.submod} have names: The
\textquotedblleft$a+b\in N$\textquotedblright\ axiom is called
\textquotedblleft$N$ is closed under addition\textquotedblright; the
\textquotedblleft$ra\in N$\textquotedblright\ axiom is called
\textquotedblleft$N$ is closed under scaling\textquotedblright; the
\textquotedblleft$0\in N$\textquotedblright\ axiom is called \textquotedblleft%
$N$ contains the zero vector\textquotedblright. The word \textquotedblleft
scaling\textquotedblright\ that we have just used refers to the following operation:

\begin{definition}
\label{def.mod.scaling}Let $R$ be a ring. Let $M$ be a left $R$-module. Let
$r\in R$ be a scalar. Then, \textbf{scaling} by $r$ (on the module $M$) means
the map%
\begin{align*}
M  &  \rightarrow M,\\
m  &  \mapsto rm.
\end{align*}

\end{definition}

This map is a group homomorphism from the additive group $\left(
M,+,0\right)  $ to itself (since any $m,n\in M$ satisfy $r\left(  m+n\right)
=rm+rn$ and $r\cdot0_{M}=0_{M}$). In particular, scaling by $1$ is the
identity map $\operatorname*{id}\nolimits_{M}:M\rightarrow M$ (since $1m=m$
for each $m\in M$), whereas scaling by $0$ sends each vector $m\in M$ to the
zero vector $0_{M}$.

The \textquotedblleft$ra\in N$\textquotedblright\ axiom in Definition
\ref{def.mod.submod} is saying that $N$ is closed under scaling by every
scalar $r\in R$. We will soon see that an $R$-submodule of $M$ is the same as
a subgroup of the additive group $\left(  M,+,0\right)  $ that is closed under
scaling by every scalar $r\in R$.

Everything we have said about left $R$-modules can be equally said (mutatis
mutandis) for right $R$-modules.

\subsubsection{Examples}

Here are some examples of modules:

\begin{itemize}
\item Let $R$ be any ring. Then, $R$ itself becomes a left $R$-module: Just
define the action to be the multiplication of $R$. In this $R$-module, the
elements of $R$ play both the role of vectors and the role of scalars. Scaling
a vector $m$ by a scalar $r$ just means multiplying $m$ by $r$ (that is,
taking the product $rm$ inside $R$).

The $R$-submodules of this left $R$-module $R$ are the subsets $L$ of $R$ that
are closed under addition and contain $0$ and satisfy $ra\in L$ for all $r\in
R$ and $a\in L$. These subsets are called the \textbf{left ideals} of $R$.
When $R$ is commutative, these are precisely the ideals of $R$. For general
$R$, however, the notion of an ideal is more restrictive than the notion of a
left ideal (since an ideal $L$ has to satisfy not only $ra\in L$ but also
$ar\in L$ for all $r\in R$ and $a\in L$).

For example, if $R$ is the matrix ring $\mathbb{Q}^{2\times2}$, then the only
ideals of $R$ are $\left\{  0_{2\times2}\right\}  $ and $R$ itself, but $R$
has infinitely many left ideals (for example, the set of all matrices of the
form $\left(
\begin{array}
[c]{cc}%
0 & a\\
0 & b
\end{array}
\right)  $ is a left ideal).

\item Let $R$ be any ring. An $R$-module (left or right) is said to be
\textbf{trivial} if it has only one element. The one-element set $\left\{
0\right\}  $ is a trivial left $R$-module (with addition given by $0+0=0$,
action given by $r\cdot0=0$, and zero vector $0$) and a trivial right
$R$-module (likewise).

\item Let $R$ be any ring, and let $n\in\mathbb{N}$. Then, the Cartesian
product%
\[
R^{n}=\left\{  \left(  a_{1},a_{2},\ldots,a_{n}\right)  \ \mid\ \text{all
}a_{i}\text{ belong to }R\right\}
\]
is a left $R$-module, where addition and action are defined entrywise: i.e.,
the addition is defined by
\begin{align*}
\left(  a_{1},a_{2},\ldots,a_{n}\right)  +\left(  b_{1},b_{2},\ldots
,b_{n}\right)   &  =\left(  a_{1}+b_{1},\ a_{2}+b_{2},\ \ldots,\ a_{n}%
+b_{n}\right) \\
\ \ \ \ \ \ \ \ \ \ \text{for all }a_{1},a_{2},\ldots,a_{n}  &  \in R\text{
and }b_{1},b_{2},\ldots,b_{n}\in R,
\end{align*}
and the action is defined by%
\[
r\cdot\left(  a_{1},a_{2},\ldots,a_{n}\right)  =\left(  ra_{1},ra_{2}%
,\ldots,ra_{n}\right)  \text{ for all }r\in R\text{ and }a_{1},a_{2}%
,\ldots,a_{n}\in R.
\]
The zero vector of this $R$-module $R^{n}$ is $\left(  0,0,\ldots,0\right)  $.

If $n=0$, then the left $R$-module $R^{n}$ is trivial, and its only element is
the $0$-tuple $\left(  {}\right)  $.

\item Let $R$ be any ring, and let $n,m\in\mathbb{N}$. Consider the set
$R^{n\times m}$ of all $n\times m$-matrices with entries in $R$. This set
$R^{n\times m}$ is not a ring unless $n=m$ (since two $n\times m$-matrices
cannot be multiplied unless $n=m$). However, it is always a left $R$-module,
where addition and action are defined entrywise: e.g., the action is defined
by%
\begin{align*}
r\cdot\left(
\begin{array}
[c]{cccc}%
a_{1,1} & a_{1,2} & \cdots & a_{1,m}\\
a_{2,1} & a_{2,2} & \cdots & a_{2,m}\\
\vdots & \vdots & \ddots & \vdots\\
a_{n,1} & a_{n,2} & \cdots & a_{n,m}%
\end{array}
\right)   &  =\left(
\begin{array}
[c]{cccc}%
ra_{1,1} & ra_{1,2} & \cdots & ra_{1,m}\\
ra_{2,1} & ra_{2,2} & \cdots & ra_{2,m}\\
\vdots & \vdots & \ddots & \vdots\\
ra_{n,1} & ra_{n,2} & \cdots & ra_{n,m}%
\end{array}
\right) \\
&  \ \ \ \ \ \ \ \ \ \ \ \ \ \ \ \ \ \ \ \ \text{for any }r,a_{i,j}\in R.
\end{align*}
The zero vector of this $R$-module $R^{n\times m}$ is the zero matrix
$0_{n\times m}$.

The set $R^{n\times m}$ is also a right $R$-module (where addition and action
are again defined entrywise, but the action now results in a matrix whose
entries are $a_{i,j}r$ rather than $ra_{i,j}$).

According to Definition \ref{def.mod.leftmod}, the elements of this $R$-module
$R^{n\times m}$ can thus be called \textquotedblleft vectors\textquotedblright%
, even though they are matrices. This shows that our concept of a
\textquotedblleft vector\textquotedblright\ is much more general than the
classical notion of \textquotedblleft vectors\textquotedblright\ from
introductory linear algebra classes (i.e., row vectors and column vectors).
This generality might be an acquired taste, but it is quite useful. For
example, we will soon define linear combinations and linear independence of
vectors; thus we will automatically obtain these notions for matrices.
\end{itemize}

Note that the zero vector of an $R$-module is uniquely determined by its
addition (in fact, this is true for any group); thus, we don't even need to
specify it explicitly when we define an $R$-module.

\begin{exercise}
Let $R$ be a ring, and $n\in\mathbb{N}$. Consider the left $R$-module $R^{n}$.

\begin{enumerate}
\item[\textbf{(a)}] Prove that the set%
\begin{align*}
A:=  &  \left\{  \left(  a_{1},a_{2},\ldots,a_{n}\right)  \in R^{n}%
\ \mid\ a_{1}=a_{2}=\cdots=a_{n}\right\} \\
=  &  \left\{  \left(  \underbrace{a,a,\ldots,a}_{n\text{ times}}\right)
\ \mid\ a\in R\right\}
\end{align*}
is an $R$-submodule of $R^{n}$.

\item[\textbf{(b)}] Prove that the set%
\[
B:=\left\{  \left(  a_{1},a_{2},\ldots,a_{n}\right)  \in R^{n}\ \mid
\ 1a_{1}=2a_{2}=3a_{3}=\cdots=na_{n}\right\}
\]
is an $R$-submodule of $R^{n}$.

\item[\textbf{(c)}] Prove that the set%
\[
C:=\left\{  \left(  a_{1},a_{2},\ldots,a_{n}\right)  \in R^{n}\ \mid
\ a_{1}a_{2}\cdots a_{n}=0\right\}
\]
is an $R$-submodule of $R^{n}$ only if $R$ is trivial or $n=1$.

\item[\textbf{(d)}] Prove that the set
\[
D:=\left\{  \left(  a_{1},a_{2},\ldots,a_{n}\right)  \in R^{n}\ \mid
\ a_{i}=a_{i-1}+a_{i-2}\text{ for all }i\geq3\right\}
\]
is an $R$-submodule of $R^{n}$.
\end{enumerate}
\end{exercise}

\subsubsection{Left vs. right $R$-modules in general}

As we said above, the notion of a left $R$-module is not equivalent to the
notion of a right $R$-module when $R$ is a noncommutative ring. However, the
general notion of a \textquotedblleft left module over a
ring\textquotedblright\ is equivalent to the general notion of a
\textquotedblleft right module over a ring\textquotedblright. Indeed, the next
exercise (\cite[homework set \#2, Exercise 2 \textbf{(d)}]{21w}) provides a
way to convert right modules into left modules (over a different ring):

\begin{exercise}
\label{exe.21hw2.2d}Let $R$ be a ring. Define the opposite ring
$R^{\operatorname*{op}}$ as in Exercise \ref{exe.21hw2.2abc}.

Let $M$ be a right $R$-module. Prove that $M$ becomes a left
$R^{\operatorname{op}}$-module if we define an action of $R^{\operatorname{op}%
}$ on $M$ by
\[
rm=mr\ \ \ \ \ \ \ \ \ \ \text{for all $r\in R^{\operatorname{op}}$ and $m\in
M$}.
\]
(Here, the left hand side is to be understood as the image of $\left(
r,m\right)  $ under the new action of $R^{\operatorname{op}}$ on $M$, whereas
the right hand side is the image of $\left(  m,r\right)  $ under the original
action of $R$ on $M$.)
\end{exercise}

Similarly, we can translate left $R$-modules into right $R^{\operatorname*{op}%
}$-modules. (This is just a generalization of Proposition
\ref{prop.mod.leftmod-is-rightmod} to arbitrary -- i.e., not necessarily
commutative -- rings $R$.)

Thus, for any ring $R$, we can translate left $R$-modules into right
$R^{\operatorname*{op}}$-modules and vice versa. As a consequence, any
property of left $R$-modules can be translated into a property of right
$R^{\operatorname*{op}}$-modules, and vice versa. The same holds with the
words \textquotedblleft left\textquotedblright\ and \textquotedblleft
right\textquotedblright\ interchanged. Thus, we can focus our study on left
$R$-modules, knowing that everything we prove about them will also hold (with
analogous proofs) for right $R^{\operatorname*{op}}$-modules, and thus (if we
replace $R$ by $R^{\operatorname*{op}}$) for right $R$-modules.

\subsection{\label{sec.modules.gens}A couple generalities}

Let us next show a few general properties of modules.

\subsubsection{Negation, subtraction and scaling}

Recall that when a group is written additively (i.e., its operation is denoted
by $+$), the inverse of an element $a$ of this group is denoted by $-a$ (and
is called its additive inverse). The following proposition says that the
additive inverse of a vector in an $R$-module can be obtained by scaling the
vector by $-1$:

\begin{proposition}
\label{prop.submodules.neg}Let $R$ be a ring. Let $M$ be a left $R$-module.
Then, $\left(  -1\right)  m=-m$ for each $m\in M$. (Here, \textquotedblleft%
$-1$\textquotedblright\ stands for $-1_{R}$.)
\end{proposition}

\begin{proof}
Let $m\in M$. Then, $1m=m$ (by one of the module axioms). Thus,%
\begin{align*}
\left(  -1\right)  m+\underbrace{m}_{=1m}  &  =\left(  -1\right)  m+1m\\
&  =\underbrace{\left(  \left(  -1\right)  +1\right)  }_{=0_{R}}%
m\ \ \ \ \ \ \ \ \ \ \left(  \text{by the right distributivity axiom}\right)
\\
&  =0_{R}m=0_{M}\ \ \ \ \ \ \ \ \ \ \left(  \text{by one of the module
axioms}\right)  .
\end{align*}
In other words, $\left(  -1\right)  m$ is an additive inverse of $m$. But the
additive inverse of $m$ is $-m$. Thus, we conclude that $\left(  -1\right)
m=-m$. This proves Proposition \ref{prop.submodules.neg}.
\end{proof}

Further properties of negation and scaling can easily be derived from this.
For example:

\begin{proposition}
\label{prop.submodules.-rm}Let $R$ be a ring. Let $M$ be a left $R$-module.
Let $r\in R$ and $m\in M$. Then,%
\begin{equation}
\left(  -r\right)  m=-\left(  rm\right)  =r\left(  -m\right)
\label{eq.prop.submodules.-rm.1}%
\end{equation}
and
\begin{equation}
\left(  -r\right)  \left(  -m\right)  =rm. \label{eq.prop.submodules.-rm.2}%
\end{equation}

\end{proposition}

\begin{proof}
Left to the reader. (Just as in the proof of Proposition
\ref{prop.submodules.neg}, argue that both $\left(  -r\right)  m$ and
$r\left(  -m\right)  $ are additive inverses of $rm$. This proves
(\ref{eq.prop.submodules.-rm.1}). To get (\ref{eq.prop.submodules.-rm.2}),
apply (\ref{eq.prop.submodules.-rm.1}) to $-m$ instead of $m$.)
\end{proof}

\begin{proposition}
\label{prop.submods.subgp}Let $R$ be a ring. Let $M$ be a left $R$-module.
Then, any $R$-submodule of $M$ is a subgroup of the additive group $\left(
M,+,0\right)  $.
\end{proposition}

\begin{proof}
[Proof of Proposition \ref{prop.submods.subgp}.]Let $N$ be an $R$-submodule of
$M$. Then, $N$ is closed under addition and under scaling and contains the
zero vector. Since $N$ is closed under scaling, we have $\left(  -1\right)
m\in N$ for each $m\in N$. However, each $m\in N$ satisfies $\left(
-1\right)  m=-m$ (by Proposition \ref{prop.submodules.neg}) and thus
$-m=\left(  -1\right)  m\in N$. In other words, $N$ is closed under negation
(= taking additive inverses). Thus, $N$ is a subgroup of $\left(
M,+,0\right)  $.
\end{proof}

\begin{proposition}
\label{prop.submods.via-scaling}Let $R$ be a ring. Let $M$ be a left
$R$-module. Then, an $R$-submodule of $M$ is the same as a subgroup of the
additive group $\left(  M,+,0\right)  $ that is closed under scaling by every
scalar $r\in R$.
\end{proposition}

\begin{proof}
Any $R$-submodule of $M$ is a subgroup of the additive group $\left(
M,+,0\right)  $ (by Proposition \ref{prop.submods.subgp}) that is closed under
scaling by every scalar $r\in R$ (by the definition of a submodule).
Conversely, any subgroup of the additive group $\left(  M,+,0\right)  $ that
is closed under scaling by every scalar $r\in R$ is an $R$-submodule of $M$
(since it satisfies all the axioms for a submodule). Thus, Proposition
\ref{prop.submods.via-scaling} is proved.
\end{proof}

\begin{proposition}
\label{prop.submods.mod}Let $R$ be a ring. Let $M$ be a left $R$-module. Then,
any $R$-submodule of $M$ becomes a left $R$-module in its own right (just like
a subring of a ring becomes a ring).
\end{proposition}

\begin{proof}
Let $N$ be an $R$-submodule of $M$. Then, Proposition \ref{prop.submods.subgp}
shows that $N$ is a subgroup of the additive group $\left(  M,+,0\right)  $.
Hence, $\left(  N,+,0\right)  $ is a group. Since $N$ is closed under scaling,
we furthermore can define an action of $R$ on $N$ in the obvious way (viz.,
inheriting it from $M$). This makes $N$ into a left $R$-module. This proves
Proposition \ref{prop.submods.mod}.
\end{proof}

We also have \textquotedblleft distributivity laws for
subtraction\textquotedblright:

\begin{proposition}
\label{prop.submod.negative-dist}Let $R$ be a ring. Let $M$ be a left
$R$-module. Then:

\begin{enumerate}
\item[\textbf{(a)}] We have $\left(  r-s\right)  m=rm-sm$ for all $r,s\in R$
and $m\in M$.

\item[\textbf{(b)}] We have $r\left(  m-n\right)  =rm-rn$ for all $r\in R$ and
$m,n\in M$.
\end{enumerate}
\end{proposition}

\begin{proof}
LTTR. (The fastest way is to derive these properties from the distributivity
laws by strategic application of (\ref{eq.prop.submodules.-rm.1}).)
\end{proof}

\subsubsection{Finite sums}

Next, let us recall how we defined finite sums $\sum_{s\in S}a_{s}$ of
elements of a ring. In the same way, we can define a finite sum $\sum_{s\in
S}a_{s}$ of elements of any additive group, and thus a finite sum $\sum_{s\in
S}a_{s}$ of elements of any $R$-module (since any $R$-module is an additive
group). Thus, in particular, if $a_{1},a_{2},\ldots,a_{n}$ are $n$ elements of
an $R$-module $M$, then the finite sum $a_{1}+a_{2}+\cdots+a_{n}\in M$ is well-defined.

The following \textquotedblleft generalized distributivity
laws\textquotedblright\ hold in any left $R$-module:

\begin{proposition}
\label{prop.submods.distribs}Let $R$ be a ring. Let $M$ be a left $R$-module. Then:

\begin{enumerate}
\item[\textbf{(a)}] We have%
\[
\left(  r_{1}+r_{2}+\cdots+r_{k}\right)  m=r_{1}m+r_{2}m+\cdots+r_{k}m
\]
for any $r_{1},r_{2},\ldots,r_{k}\in R$ and $m\in M$.

\item[\textbf{(b)}] We have%
\[
r\left(  m_{1}+m_{2}+\cdots+m_{i}\right)  =rm_{1}+rm_{2}+\cdots+rm_{i}%
\]
for any $r\in R$ and $m_{1},m_{2},\ldots,m_{i}\in M$.
\end{enumerate}
\end{proposition}

\begin{proof}
\textbf{(a)} This follows by applying the right distributivity law (one of the
module axioms) many times. (More precisely, this follows by induction on $k$;
the right distributivity law is used in the induction step. The induction base
follows from the $0_{R}m=0_{M}$ axiom.) \medskip

\textbf{(b)} This follows by applying the left distributivity law (one of the
module axioms) many times. (More precisely, this follows by induction on $i$;
the left distributivity law is used in the induction step. The induction base
follows from the $r\cdot0_{M}=0_{M}$ axiom.)
\end{proof}

The following convention is useful when dealing with $R$-modules. Essentially,
it says that (just as with products of multiple elements in a ring or in a
group) we can drop parentheses when we scale an element of an $R$-module by
several elements of $R$:

\begin{convention}
\label{conv.mods.rsm}Let $R$ be a ring. Let $M$ be a left $R$-module. Let
$r,s\in R$ and $m\in M$. Then, $\left(  rs\right)  m$ and $r\left(  sm\right)
$ are the same vector (by the associativity axiom in the definition of a left
$R$-module). We shall denote this vector by $rsm$. Likewise, expressions like
$r_{1}r_{2}\cdots r_{k}m$ (for $r_{1},r_{2},\ldots,r_{k}\in R$ and $m\in M$)
will be understood.
\end{convention}

Everything we said above about left $R$-modules can be adapted to right
$R$-modules in a straightforward way; we leave the details to the reader.

\subsubsection{Some exercises}

\begin{exercise}
\label{exe.mods.modulo}Let $R$ be a ring. Let $M$ be a left $R$-module. Let
$I$ be an $R$-submodule of $M$.

For any two elements $a,b\in M$, we write \textquotedblleft$a\equiv
b\operatorname{mod}I$\textquotedblright\ (and say that \textquotedblleft$a$ is
congruent to $b$ modulo $I$\textquotedblright) if and only if $a-b\in I$.
(This is a generalization of congruence of integers, as it is usually defined
in elementary number theory. Indeed, congruence of integers modulo an integer
$n$ is recovered when $R=\mathbb{Z}$ and $I=n\mathbb{Z}$.)

Prove the following:

\begin{enumerate}
\item[\textbf{(a)}] Each $a\in M$ satisfies $a\equiv a\operatorname{mod}I$.

\item[\textbf{(b)}] If $a,b\in M$ satisfy $a\equiv b\operatorname{mod}I$, then
$b\equiv a\operatorname{mod}I$.

\item[\textbf{(c)}] If $a,b,c\in M$ satisfy $a\equiv b\operatorname{mod}I$ and
$b\equiv c\operatorname{mod}I$, then $a\equiv c\operatorname{mod}I$.

\item[\textbf{(d)}] If $a,b,c,d\in M$ satisfy $a\equiv b\operatorname{mod}I$
and $c\equiv d\operatorname{mod}I$, then $a+c\equiv b+d\operatorname{mod}I$.

\item[\textbf{(e)}] If $a,b\in M$ and $r\in R$ satisfy $a\equiv
b\operatorname{mod}I$, then $ra\equiv rb\operatorname{mod}I$.
\end{enumerate}

Now, we claim a sort of converse:

\begin{enumerate}
\item[\textbf{(f)}] Let us drop the requirement that $I$ be an $R$-submodule
of $M$. Instead, we require that $I$ be any subset of $M$ for which the claims
of parts \textbf{(a)}, \textbf{(c)} and \textbf{(e)} of this exercise hold.
Prove that $I$ is an $R$-submodule of $M$.
\end{enumerate}
\end{exercise}

Exercise \ref{exe.mods.modulo} can be summarized as \textquotedblleft modular
arithmetic modulo a subset $I$ of $M$ works if and only if $I$ is a submodule
of $M$\textquotedblright. In other words, roughly speaking, the submodules of
a module $M$ are precisely the subsets $I$ that allow \textquotedblleft
working modulo $I$\textquotedblright. This is most likely the reason why
modules are called \textquotedblleft modules\textquotedblright\footnote{The
name was coined by Dedekind, although in a less general context.}.

\begin{exercise}
Let $R$ be a ring. Let $M$ be a left $R$-module.

For any subset $K$ of $M$, let $\operatorname*{Ann}K$ denote the subset
$\left\{  r\in R\ \mid\ rk=0\text{ for all }k\in K\right\}  $ of $R$. (This is
called the \textbf{annihilator} of $K$.)

\begin{enumerate}
\item[\textbf{(a)}] Prove that $\operatorname*{Ann}M$ is an ideal of $R$.

\item[\textbf{(b)}] Let $K$ be any subset of $M$. Prove that
$\operatorname*{Ann}K$ is a left ideal of $R$. (Recall that a \textbf{left
ideal} of $R$ means a subset $L$ of $R$ that is closed under addition and
contains $0$ and satisfies $ra\in L$ for all $r\in R$ and $a\in L$.)

\item[\textbf{(c)}] Find an example showing that the $\operatorname*{Ann}K$ in
part \textbf{(b)} is not always an ideal of $R$.
\end{enumerate}
\end{exercise}

\subsection{More operations on modules and submodules}

\subsubsection{\label{subsec.modules.def.dirprod}Direct products and direct
sums}

Fix a ring $R$. In Subsection \ref{subsec.modules.def.def}, we have defined
left $R$-modules (recall: these are essentially additive groups whose elements
can be scaled by elements of $R$), and afterwards we have seen a few examples
of them. Let me briefly repeat the two simplest examples:

\begin{itemize}
\item The ring $R$ itself becomes a left $R$-module: Just define the action to
be the multiplication of $R$. This is called the \textbf{natural left }%
$R$\textbf{-module }$R$. The $R$-submodules of this $R$-module are the left
ideals of $R$. (Every ideal of $R$ is a left ideal of $R$, but usually not
vice versa.)

\item For any $n\in\mathbb{N}$, the set%
\[
R^{n}=\left\{  \left(  a_{1},a_{2},\ldots,a_{n}\right)  \ \mid\ \text{all
}a_{i}\text{ belong to }R\right\}
\]
is a left $R$-module, with addition and action being entrywise\footnote{e.g.,
the action is defined by%
\[
r\cdot\left(  a_{1},a_{2},\ldots,a_{n}\right)  =\left(  ra_{1},ra_{2}%
,\ldots,ra_{n}\right)  \text{ for all }r\in R\text{ and }a_{1},a_{2}%
,\ldots,a_{n}\in R.
\]
} and with the zero vector $\left(  0,0,\ldots,0\right)  $. This generalizes
the Euclidean space $\mathbb{R}^{n}$ from linear algebra, and many of its analogues.
\end{itemize}

Here are some more examples:

\begin{itemize}
\item The left $R$-modules $R^{n}$ (with $n\in\mathbb{N}$) tend to have many
$R$-submodules. When $R$ is a field, this is well-known from linear algebra
(where $R$-submodules are called $R$-vector subspaces); in particular, the
solution set of any given system of homogeneous linear equations in $n$
variables is an $R$-submodule of $R^{n}$. The same applies to any commutative
ring $R$, but here we have even more freedom: Besides equations, our system
can contain congruences too (as long as they are linear and have no constant
term). For instance, for $R=\mathbb{Z}$, the set%
\begin{align*}
&  \left\{  \left(  x,y,z,w\right)  \in\mathbb{Z}^{4}\ \mid\ x\equiv
y\operatorname{mod}2\text{ and }x+y+z+w\equiv0\operatorname{mod}3\right. \\
&
\ \ \ \ \ \ \ \ \ \ \ \ \ \ \ \ \ \ \ \ \ \ \ \ \ \ \ \ \ \ \ \ \ \ \ \ \ \ \ \ \left.
\vphantom{\mathbb{Z}^4}\text{and }x-y+z-w=0\right\}
\end{align*}
is a $\mathbb{Z}$-submodule of $\mathbb{Z}^{4}$. To prove this, you need to
check the axioms (\textquotedblleft closed under addition\textquotedblright,
\textquotedblleft closed under scaling\textquotedblright\ and
\textquotedblleft contains the zero vector\textquotedblright). With a bit of
practice, you can do this all in your head.

If $R$ is noncommutative, you have to be somewhat careful with the side on
which the coefficients stand in your system. If the coefficients are on the
\textbf{right} of the variables, then the solution set is a \textbf{left}
$R$-module (so, e.g., if $a$ and $b$ are two elements of $R$, then $\left\{
\left(  x,y\right)  \in R^{2}\ \mid\ xa+yb=0\right\}  $ is a left $R$-module);
on the other hand, if the coefficients are on the \textbf{left} of the
variables, then the solution set is a \textbf{right} $R$-module. (Again, this
is not hard to check: e.g., the set $\left\{  \left(  x,y\right)  \in
R^{2}\ \mid\ xa+yb=0\right\}  $ is closed under the scaling maps of a left
$R$-module because $xa+yb=0$ implies $rxa+ryb=r\underbrace{\left(
xa+yb\right)  }_{=0}=0$. Meanwhile, in general, this set is not closed under
the scaling maps of a right $R$-module, since $xa+yb=0$ does not imply
$xra+yrb=0$.)

\item Just as we defined the left $R$-module $R^{n}$ consisting of all
$n$-tuples, we can define a left $R$-module \textquotedblleft$R^{\infty}%
$\textquotedblright\ consisting of all infinite sequences. It is commonly
denoted by $R^{\mathbb{N}}$ (since there are different kinds of infinity).
Explicitly, we define the left $R$-module $R^{\mathbb{N}}$ by%
\[
R^{\mathbb{N}}:=\left\{  \left(  a_{0},a_{1},a_{2},\ldots\right)
\ \mid\ \text{all }a_{i}\text{ belong to }R\right\}  ,
\]
where addition and action are defined entrywise.

This left $R$-module $R^{\mathbb{N}}$ has an $R$-submodule%
\begin{align*}
R^{\left(  \mathbb{N}\right)  }:=  &  \ \left\{  \left(  a_{0},a_{1}%
,a_{2},\ldots\right)  \in R^{\mathbb{N}}\ \mid\ \text{only finitely many }%
i\in\mathbb{N}\text{ satisfy }a_{i}\neq0\right\} \\
=  &  \ \left\{  \left(  a_{0},a_{1},a_{2},\ldots\right)  \in R^{\mathbb{N}%
}\ \mid\ \text{all but finitely many }i\in\mathbb{N}\text{ satisfy }%
a_{i}=0\right\}  .
\end{align*}
You can check that this is indeed an $R$-submodule of $R^{\mathbb{N}}$. (For
instance, it is closed under addition, because if only finitely many
$i\in\mathbb{N}$ satisfy $a_{i}\neq0$ and only finitely many $i\in\mathbb{N}$
satisfy $b_{i}\neq0$, then only finitely many $i\in\mathbb{N}$ satisfy
$a_{i}+b_{i}\neq0$.)

For example, if $R=\mathbb{Q}$, then%
\begin{align*}
\left(  1,1,1,\ldots\right)   &  \in R^{\mathbb{N}}\setminus R^{\left(
\mathbb{N}\right)  }\\
\text{and}\ \ \ \ \ \ \ \ \ \ \left(  0,0,0,\ldots\right)   &  \in R^{\left(
\mathbb{N}\right)  }\\
\text{and}\ \ \ \ \ \ \ \ \ \ \left(  1,0,0,0,\ldots\right)   &  \in
R^{\left(  \mathbb{N}\right)  }\\
\text{and}\ \ \ \ \ \ \ \ \ \ \left(  1,0,4,\underbrace{0,0,0,\ldots
}_{\text{zeroes}}\right)   &  \in R^{\left(  \mathbb{N}\right)  }\\
\text{and}\ \ \ \ \ \ \ \ \ \ \left(  \underbrace{1,0,1,0,1,0,\ldots
}_{\text{ones and zeroes in turn}}\right)   &  \in R^{\mathbb{N}}\setminus
R^{\left(  \mathbb{N}\right)  }.
\end{align*}

\end{itemize}

Generalizing $R^{n}$, here is a way to build modules out of other modules:

\begin{definition}
\label{def.mods.dirprod-n}Let $n\in\mathbb{N}$, and let $M_{1},M_{2}%
,\ldots,M_{n}$ be any $n$ left $R$-modules. Then, the Cartesian product
$M_{1}\times M_{2}\times\cdots\times M_{n}$ becomes a left $R$-module itself,
where addition and action are defined entrywise: e.g., the action is defined
by%
\[
r\cdot\left(  m_{1},m_{2},\ldots,m_{n}\right)  =\left(  rm_{1},rm_{2}%
,\ldots,rm_{n}\right)  \text{ for all }r\in R\text{ and }m_{i}\in M_{i}.
\]

This left $R$-module $M_{1}\times M_{2}\times\cdots\times M_{n}$ is called the
\textbf{direct product} of $M_{1},M_{2},\ldots,M_{n}$.
\end{definition}

If all of $M_{1},M_{2},\ldots,M_{n}$ are the natural left $R$-module $R$, then
this direct product is precisely the left $R$-module $R^{n}$ defined above.

This direct product $M_{1}\times M_{2}\times\cdots\times M_{n}$ can be
generalized further, allowing products of infinitely many modules, too. Just
as for rings, the best setting for this is using families, not
lists:\footnote{The proof of Proposition \ref{prop.mods.dirprod-exists} is
easy and LTTR.}

\begin{proposition}
\label{prop.mods.dirprod-exists}Let $I$ be any set. Let $\left(  M_{i}\right)
_{i\in I}$ be any family of left $R$-modules. Then, the Cartesian product%
\[
\prod\limits_{i\in I}M_{i}=\left\{  \text{all families }\left(  m_{i}\right)
_{i\in I}\text{ with }m_{i}\in M_{i}\text{ for all }i\in I\right\}
\]
becomes a left $R$-module if we endow it with the entrywise addition (i.e., we
set $\left(  m_{i}\right)  _{i\in I}+\left(  n_{i}\right)  _{i\in I}=\left(
m_{i}+n_{i}\right)  _{i\in I}$ for any two families $\left(  m_{i}\right)
_{i\in I},\left(  n_{i}\right)  _{i\in I}\in\prod\limits_{i\in I}M_{i}$) and
the entrywise scaling (i.e., we set $r\left(  m_{i}\right)  _{i\in I}=\left(
rm_{i}\right)  _{i\in I}$ for any $r\in R$ and any family $\left(
m_{i}\right)  _{i\in I}\in\prod\limits_{i\in I}M_{i}$) and with the zero
vector $\left(  0\right)  _{i\in I}$.
\end{proposition}

\begin{definition}
\label{def.mods.dirprod}This left $R$-module is denoted by $\prod\limits_{i\in
I}M_{i}$ and called the \textbf{direct product} of the left $R$-modules
$M_{i}$. In some special cases, there are alternative notations for it:

\begin{itemize}
\item If $I=\left\{  1,2,\ldots,n\right\}  $ for some $n\in\mathbb{N}$, then
this left $R$-module is also denoted by $M_{1}\times M_{2}\times\cdots\times
M_{n}$, and we identify a family $\left(  m_{i}\right)  _{i\in I}=\left(
m_{i}\right)  _{i\in\left\{  1,2,\ldots,n\right\}  }$ with the $n$-tuple
$\left(  m_{1},m_{2},\ldots,m_{n}\right)  $. (Thus, $M_{1}\times M_{2}%
\times\cdots\times M_{n}$ is precisely the direct product $M_{1}\times
M_{2}\times\cdots\times M_{n}$ we defined in Definition
\ref{def.mods.dirprod-n}.)

\item If all the left $R$-modules $M_{i}$ are equal to some left $R$-module
$M$, then their direct product $\prod\limits_{i\in I}M_{i}=\prod\limits_{i\in
I}M$ is also denoted $M^{I}$. Note that this generalizes the $R^{\mathbb{N}}$
defined above.

\item We set $M^{n}=M^{\left\{  1,2,\ldots,n\right\}  }$ for each
$n\in\mathbb{N}$ and any left $R$-module $M$. This generalizes the left
$R$-module $R^{n}$ for $n\in\mathbb{N}$ discussed above.
\end{itemize}
\end{definition}

This was quite predictable; but there is more. Indeed, we can generalize not
just $R^{\mathbb{N}}$ but also its submodule $R^{\left(  \mathbb{N}\right)  }%
$, and the result is at least as important:\footnote{The proof of Proposition
\ref{prop.mods.dirsum-exists} is easy and LTTR.}

\begin{proposition}
\label{prop.mods.dirsum-exists}Let $I$ be any set. Let $\left(  M_{i}\right)
_{i\in I}$ be any family of left $R$-modules. Define $\bigoplus\limits_{i\in
I}M_{i}$ to be the subset%
\[
\left\{  \left(  m_{i}\right)  _{i\in I}\in\prod\limits_{i\in I}M_{i}%
\ \mid\text{ only finitely many }i\in I\text{ satisfy }m_{i}\neq0\right\}
\]
of $\prod\limits_{i\in I}M_{i}$. Then, $\bigoplus\limits_{i\in I}M_{i}$ is a
left $R$-submodule of $\prod\limits_{i\in I}M_{i}$, and thus becomes a left
$R$-module itself (by Proposition \ref{prop.submods.mod}).
\end{proposition}

\begin{definition}
\label{def.mods.dirsum}This left $R$-module $\bigoplus\limits_{i\in I}M_{i}$
is called the \textbf{direct sum} of the $R$-modules $M_{i}$.

If $I=\left\{  1,2,\ldots,n\right\}  $ for some $n\in\mathbb{N}$, then this
left $R$-module is also denoted by $M_{1}\oplus M_{2}\oplus\cdots\oplus M_{n}$.
\end{definition}

The last part of this definition might raise some eyebrows. In fact, if the
set $I$ is finite, then $\bigoplus\limits_{i\in I}M_{i}=\prod\limits_{i\in
I}M_{i}$ (since the condition \textquotedblleft only finitely many $i\in I$
satisfy $m_{i}\neq0$\textquotedblright\ is automatically satisfied for any
family $\left(  m_{i}\right)  _{i\in I}$ when $I$ is finite). Thus, in
particular,%
\[
M_{1}\oplus M_{2}\oplus\cdots\oplus M_{n}=M_{1}\times M_{2}\times\cdots\times
M_{n}%
\]
for any left $R$-modules $M_{1},M_{2},\ldots,M_{n}$. So we have introduced two
notations for the same thing (and even worse, one of the notations looks like
a sum, while the other looks like a product!). Nevertheless, both are in use.
Direct sums start differing from direct products when the indexing set $I$ is infinite.

For $I=\mathbb{N}$ and $M_{i}=R$, the direct sum $\bigoplus\limits_{i\in
I}M_{i}=\bigoplus\limits_{i\in\mathbb{N}}R$ is precisely the $R$-module
$R^{\left(  \mathbb{N}\right)  }$ defined above. More generally:

\begin{definition}
\label{def.mods.M(I)}Let $I$ be a set. Let $M$ be any left $R$-module. Then,
$M^{\left(  I\right)  }$ is defined to be the left $R$-module $\bigoplus
\limits_{i\in I}M$.
\end{definition}

\begin{exercise}
\label{exe.mods.dirprod-submods}Let $I$ be any set. Let $\left(  M_{i}\right)
_{i\in I}$ be any family of left $R$-modules. Let $N_{i}$ be an $R$-submodule
of $M_{i}$ for each $i\in I$.

\begin{enumerate}
\item[\textbf{(a)}] Prove that $\prod_{i\in I}N_{i}$ is an $R$-submodule of
the left $R$-module $\prod_{i\in I}M_{i}$.

\item[\textbf{(b)}] Prove that $\bigoplus\limits_{i\in I}N_{i}$ is an
$R$-submodule of the left $R$-module $\bigoplus\limits_{i\in I}M_{i}$.
\end{enumerate}
\end{exercise}

\subsubsection{\label{subsec.modules.def.restr}Restriction of modules}

Here are some more ways to construct modules over rings:

\begin{itemize}
\item If $R$ is a subring of a ring $S$, then $S$ is a left $R$-module (where
the action of $R$ on $S$ is defined by restricting the multiplication map
$S\times S\rightarrow S$ to $R\times S$) and a right $R$-module (in a similar way).

Let me restate this in a more down-to-earth way: If $R$ is a subring of a ring
$S$, then we can multiply any element of $R$ with any element of $S$ (since
both elements lie in the ring $S$); this makes $S$ into a left $R$-module (and
likewise, $S$ becomes a right $R$-module). Explicitly, the action of $R$ on
the left $R$-module $S$ is given by
\[
rs=rs\ \ \ \ \ \ \ \ \ \ \text{for all }r\in R\text{ and }s\in S
\]
(where the \textquotedblleft$rs$\textquotedblright\ on the left hand side
means the image of $\left(  r,s\right)  $ under the action, whereas the
\textquotedblleft$rs$\textquotedblright\ on the right hand side means the
product of $r$ and $s$ in the ring $S$).

Thus, for example, $\mathbb{C}$ is an $\mathbb{R}$-module (since $\mathbb{R}$
is a subring of $\mathbb{C}$) and also a $\mathbb{Q}$-module (for similar
reasons). (In this sentence, you can just as well say \textquotedblleft vector
space\textquotedblright\ instead of \textquotedblleft module\textquotedblright%
, since $\mathbb{R}$ and $\mathbb{Q}$ are fields.)

\item More generally: If $R$ and $S$ are any two rings, and if $f:R\rightarrow
S$ is a ring morphism, then $S$ becomes a left $R$-module (with the action of
$R$ on $S$ being defined by%
\[
rs=f\left(  r\right)  s\ \ \ \ \ \ \ \ \ \ \text{for all }r\in R\text{ and
}s\in S
\]
) and a right $R$-module (in a similar way). The proof of this is easy. These
$R$-module structures are sometimes said to be \textbf{induced} by the
morphism $f$.

Our previous example (in which we made $S$ into an $R$-module whenever $R$ is
a subring of $S$) is the particular case of this construction obtained when
$f$ is the canonical inclusion\footnote{Recall what \textquotedblleft
canonical inclusion\textquotedblright\ means:
\par
If $U$ is a subset of a set $V$, then the map%
\begin{align*}
U  &  \rightarrow V,\\
u  &  \mapsto u
\end{align*}
is called the \textbf{canonical inclusion} of $U$ into $V$.
\par
If $U$ is a subring of a ring $V$, then the canonical inclusion of $U$ into
$V$ is furthermore a ring morphism. (This follows trivially from the
definition of a subring.)} of $R$ into $S$.

Here are some other particular cases:

\begin{itemize}
\item Any quotient ring $R/I$ of a ring $R$ (by some ideal $I$) becomes a left
$R$-module, because the canonical projection $\pi:R\rightarrow R/I$ (which
sends every $r\in R$ to its residue class $\overline{r}\in R/I$) is a ring
morphism. Explicitly, the action of $R$ on $R/I$ is given by%
\[
r\cdot\overline{u}=\underbrace{\pi\left(  r\right)  }_{=\overline{r}}%
\cdot\,\overline{u}=\overline{r}\cdot\overline{u}=\overline{ru}%
\ \ \ \ \ \ \ \ \ \ \text{for all }r,u\in R.
\]
Similarly, $R/I$ becomes a right $R$-module.

\item Here is another particular case (a less transparent one): I claim that
the abelian group $\mathbb{Z}/5$ becomes a $\mathbb{Z}\left[  i\right]
$-module\footnote{As usual, $\mathbb{Z}\left[  i\right]  $ denotes the ring of
the Gaussian integers, with $i=\sqrt{-1}$.}, if we define the action by%
\[
\left(  a+bi\right)  \cdot m=\overline{a+2b}\cdot
m\ \ \ \ \ \ \ \ \ \ \text{for all }a+bi\in\mathbb{Z}\left[  i\right]  \text{
and }m\in\mathbb{Z}/5.
\]
To wit, the map%
\begin{align*}
f:\mathbb{Z}\left[  i\right]   &  \rightarrow\mathbb{Z}/5,\\
a+bi  &  \mapsto\overline{a+2b}%
\end{align*}
is a ring morphism (check this!\footnote{Indeed, it is pretty easy to see that
this map $f$ respects addition, the zero and the unity. It remains to show
that this map respects multiplication. To show this, we fix any $x,y\in
\mathbb{Z}\left[  i\right]  $. We then need to show that $f\left(  xy\right)
=f\left(  x\right)  f\left(  y\right)  $.
\par
Write $x$ and $y$ as $x=a+bi$ and $y=c+di$ for some $a,b,c,d\in\mathbb{Z}$.
Then, $xy=\left(  a+bi\right)  \left(  c+di\right)  =\left(  ac-bd\right)
+\left(  ad+bc\right)  i$ (by the rule for multiplying complex numbers).
Hence,%
\begin{equation}
f\left(  xy\right)  =f\left(  \left(  ac-bd\right)  +\left(  ad+bc\right)
i\right)  =\overline{ac-bd+2\left(  ad+bc\right)  }
\label{pf.exa.Z/5-Zi-mod.pf.1}%
\end{equation}
(by the definition of $f$). On the other hand, $x=a+bi$ entails $f\left(
x\right)  =f\left(  a+bi\right)  =\overline{a+2b}$, and similarly we find
$f\left(  y\right)  =\overline{c+2d}$. Multiplying these two equalities, we
find%
\begin{equation}
f\left(  x\right)  f\left(  y\right)  =\overline{a+2b}\cdot\overline
{c+2d}=\overline{\left(  a+2b\right)  \left(  c+2d\right)  }=\overline
{ac+2^{2}bd+2\left(  ad+bc\right)  } \label{pf.exa.Z/5-Zi-mod.pf.2}%
\end{equation}
(since $\left(  a+2b\right)  \left(  c+2d\right)  =ac+2^{2}bd+2\left(
ad+bc\right)  $). Now, the right hand sides of the equalities
(\ref{pf.exa.Z/5-Zi-mod.pf.1}) and (\ref{pf.exa.Z/5-Zi-mod.pf.2}) are
identical (since $2^{2}\equiv-1\operatorname{mod}5$ and thus $\overline{2}%
^{2}=\overline{-1}$, so that $\overline{2^{2}bd}=\overline{-bd}$); hence, so
are the left hand sides. In other words, $f\left(  xy\right)  =f\left(
x\right)  f\left(  y\right)  $. This completes the proof that the map $f$
respects multiplication; therefore, $f$ is a ring morphism.}); and this can be
used to turn $\mathbb{Z}/5$ into a $\mathbb{Z}\left[  i\right]  $-module by
our above construction; this yields precisely the action I claimed above
(because all $a+bi\in\mathbb{Z}\left[  i\right]  $ and $m\in\mathbb{Z}/5$
satisfy $\left(  a+bi\right)  \cdot m=\underbrace{f\left(  a+bi\right)
}_{=\overline{a+2b}}\cdot\,m=\overline{a+2b}\cdot m$).

This is not the only way to turn $\mathbb{Z}/5$ into a $\mathbb{Z}\left[
i\right]  $-module. We could just as well use the ring morphism%
\begin{align*}
g:\mathbb{Z}\left[  i\right]   &  \rightarrow\mathbb{Z}/5,\\
a+bi  &  \mapsto\overline{a+3b}%
\end{align*}
instead of $f$. This would give us a $\mathbb{Z}\left[  i\right]  $-module
$\mathbb{Z}/5$ with action given by%
\[
\left(  a+bi\right)  \cdot m=\overline{a+3b}\cdot
m\ \ \ \ \ \ \ \ \ \ \text{for all }a+bi\in\mathbb{Z}\left[  i\right]  \text{
and }m\in\mathbb{Z}/5.
\]

Thus, we have obtained two \textbf{different} $\mathbb{Z}\left[  i\right]
$-module structures on $\mathbb{Z}/5$ -- that is, two different $\mathbb{Z}%
\left[  i\right]  $-modules that are equal as sets (and even as additive
groups) but different as $\mathbb{Z}\left[  i\right]  $-modules (and not even
isomorphic as such). None of these two module structures is more natural or
otherwise better than the other. Thus, when you speak of a \textquotedblleft%
$\mathbb{Z}\left[  i\right]  $-module $\mathbb{Z}/5$\textquotedblright, you
need to clarify which one you mean. (Such situations are rather frequent in
algebra. \textquotedblleft Natural\textquotedblright\ $R$-module structures --
i.e., structures that are clearly \textquotedblleft the right
one\textquotedblright\ -- are rare in comparison.)
\end{itemize}

\item Even more generally: If $R$ and $S$ are two rings, and if
$f:R\rightarrow S$ is a ring morphism, then any left $S$-module $M$ (not just
$S$ itself) naturally becomes a left $R$-module, with the action defined by%
\[
rm=f\left(  r\right)  m\ \ \ \ \ \ \ \ \ \ \text{for all }r\in R\text{ and
}m\in M.
\]
(You can think of this as letting $R$ act on $M$ \textquotedblleft by
proxy\textquotedblright: In order to scale a vector $m\in M$ by a scalar $r\in
R$, you just scale it by the scalar $f\left(  r\right)  \in S$.)

This method of turning $S$-modules into $R$-modules is called
\textbf{restriction of scalars} (or, more specifically, \textbf{restricting}
an $S$-module to $R$ via $f$).

If we apply this method to a canonical inclusion (i.e., if $R$ is a subring of
$S$ and if $f:R\rightarrow S$ is the canonical inclusion), then we conclude
that any module over a ring naturally becomes a module over any
subring.\footnote{You can think of it as forgetting how to scale vectors by
scalars that don't belong to the subring.} For example, any $\mathbb{C}%
$-module naturally becomes an $\mathbb{R}$-module (this is known as
\textquotedblleft decomplexification\textquotedblright\ in linear
algebra\footnote{Of course, again, linear algebraists speak of vector spaces
instead of modules.
\par
From linear algebra, you might also know a procedure going in the other
direction: \textquotedblleft complexification\textquotedblright, which turns
an $\mathbb{R}$-vector space into a $\mathbb{C}$-vector space. This, too, can
be generalized to arbitrary ring morphisms, but this would take us rather far
afield.}) and a $\mathbb{Q}$-module and a $\mathbb{Z}$-module.
\end{itemize}

\subsubsection{\label{subsec.modules.def.moreexas}More examples}

Here are some more general constructions of submodules in a given $R$-module
(similar to some of the above constructions for ideals in a given ring):

\begin{proposition}
\label{prop.mods.IcapJ}Let $R$ be a ring. Let $M$ be a left $R$-module. Let
$I$ and $J$ be two $R$-submodules of $M$. Then, $I\cap J$ is an $R$-submodule
of $M$ as well.
\end{proposition}

\begin{proposition}
\label{prop.mods.I+J}Let $R$ be a ring. Let $M$ be a left $R$-module.

\begin{enumerate}
\item[\textbf{(a)}] Let $I$ and $J$ be two $R$-submodules of $M$. Then,%
\[
I+J:=\left\{  i+j\ \mid\ i\in I\text{ and }j\in J\right\}
\]
is an $R$-submodule of $M$ as well.

\item[\textbf{(b)}] If $I$, $J$ and $K$ are three $R$-submodules of $M$, then
$\left(  I+J\right)  +K=I+\left(  J+K\right)  $.
\end{enumerate}
\end{proposition}

\begin{proposition}
\label{prop.mods.IM}Let $R$ be a ring. Let $I$ be an ideal of $R$. Let $M$ be
a left $R$-module. An $\left(  I,M\right)  $\textbf{-product} shall mean a
product of the form $im$ with $i\in I$ and $m\in M$. Then,%
\[
IM:=\left\{  \text{finite sums of }\left(  I,M\right)  \text{-products}%
\right\}
\]
is an $R$-submodule of $M$.
\end{proposition}

\begin{proof}
This is fairly similar to the proof of the fact that the product $IJ$ of two
ideals $I$ and $J$ is again an ideal (see Exercise \ref{exe.21hw1.8}
\textbf{(a)}).
\end{proof}

We can generalize Proposition \ref{prop.ideal.princid} and Proposition
\ref{prop.ideal.uI} from ideals to modules:

\begin{proposition}
\label{prop.mods.aM}Let $R$ be a commutative ring.

\begin{enumerate}
\item[\textbf{(a)}] Let $a\in R$. Let $M$ be an $R$-module. Then,
\[
aM:=\left\{  am\ \mid\ m\in M\right\}
\]
is an $R$-submodule of $M$.

\item[\textbf{(b)}] In particular, $0M=\left\{  0_{M}\right\}  $ and $1M=M$
are $R$-submodules of $M$.
\end{enumerate}
\end{proposition}

\begin{proof}
This is a straightforward generalization of Proposition
\ref{prop.ideal.princid}. The proof is LTTR.
\end{proof}

\begin{exercise}
Prove Propositions \ref{prop.mods.IcapJ}, \ref{prop.mods.I+J},
\ref{prop.mods.IM} and \ref{prop.mods.aM}.
\end{exercise}

Proposition \ref{prop.mods.aM} \textbf{(b)} holds for noncommutative rings
$R$, too: If $M$ is a left $R$-module, then $\left\{  0_{M}\right\}  $ and $M$
are $R$-submodules of $M$. These are the \textquotedblleft
bookends\textquotedblright\ for the $R$-submodules of $M$ (in the sense that
every $R$-submodule $N$ of $M$ satisfies $\left\{  0_{M}\right\}  \subseteq
N\subseteq M$).

Proposition \ref{prop.mods.aM} \textbf{(a)} holds for noncommutative rings $R$
as well, if we assume that $a$ is a central element of $R$. (Of course,
\textquotedblleft$R$-module\textquotedblright\ should then be replaced by
\textquotedblleft left $R$-module\textquotedblright.)

\begin{exercise}
Prove this claim.
\end{exercise}

Here are a few more examples of modules:

\begin{itemize}
\item Let $n\in\mathbb{N}$, and let $R$ be a ring. The set $R^{n}$ is not only
a left $R$-module (as we have seen), but also a right $R^{n\times n}%
$-module\footnote{Recall that $R^{n\times n}$ is the ring of $n\times
n$-matrices over $R$.}, where the action of $R^{n\times n}$ on $R^{n}$ is the
vector-by-matrix multiplication map%
\begin{align*}
R^{n}\times R^{n\times n}  &  \rightarrow R^{n},\\
\left(  v,M\right)   &  \mapsto vM
\end{align*}
(where we identify $n$-tuples $v\in R^{n}$ with row vectors).

\item More generally, for any $n,m\in\mathbb{N}$, the set $R^{n\times m}$ of
all $n\times m$-matrices is a left $R^{n\times n}$-module and a right
$R^{m\times m}$-module\footnote{This is in addition to it being a left
$R$-module and a right $R$-module!} (since an $n\times m$-matrix can be
multiplied by an $n\times n$-matrix from the left and by an $m\times m$-matrix
from the right, and since the module axioms follow from the standard laws of
matrix multiplication such as associativity and distributivity). Even better,
this set is a so-called $\left(  R^{n\times n},R^{m\times m}\right)
$-bimodule\footnote{A \textquotedblleft bimodule\textquotedblright%
\ essentially means a set that has both a left and a right module structure
that fit together well. Here is a precise definition: Let $R$ and $S$ be two
rings. Then, an $\left(  R,S\right)  $\textbf{-bimodule} means a set $M$ that
is simultaneously a left $R$-module and a right $S$-module, where the two
module structures have the same addition (i.e., the meaning of $m+n$ for
$m,n\in M$ does not depend on whether we use the addition of the left
$R$-module $M$ or the addition of the right $S$-module $M$) and satisfy the
compatibility axiom%
\[
\left(  rm\right)  s=r\left(  ms\right)  \ \ \ \ \ \ \ \ \ \ \text{for all
}r\in R\text{, }s\in S\text{ and }m\in M
\]
(which allows us to write \textquotedblleft products\textquotedblright\ like
$rms$ without parentheses).}.

\item Let us study a particular case of this.

Namely, let $R$ be a field $F$, and let $n=2$. So $F^{2}$ is a left
$F$-module, with the action given by%
\[
\lambda\left(  a,b\right)  =\left(  \lambda a,\lambda b\right)
\ \ \ \ \ \ \ \ \ \ \text{for all }\lambda,a,b\in F,
\]
and is a right $F^{2\times2}$-module, with the action given by%
\[
\left(  a,b\right)  \left(
\begin{array}
[c]{cc}%
x & y\\
z & w
\end{array}
\right)  =\left(  ax+bz,ay+bw\right)  \ \ \ \ \ \ \ \ \ \ \text{for all
}a,b,x,y,z,w\in F.
\]

What are the $F$-submodules of $F^{2}$ ? These are precisely the $F$-vector
subspaces of $F^{2}$; as you know from linear algebra, these subspaces are the
whole $F^{2}$ as well as the zero subspace $\left\{  0_{F^{2}}\right\}  $ and
all lines through the origin.

What are the $F^{2\times2}$-submodules of $F^{2}$ ? Only $F^{2}$ and $\left\{
0_{F^{2}}\right\}  $, because any two nonzero vectors in $F^{2}$ can be mapped
to one another by a $2\times2$-matrix.

Now, consider the subring%
\[
F^{2\geq2}:=\left\{  \left(
\begin{array}
[c]{cc}%
x & 0\\
z & w
\end{array}
\right)  \ \mid\ x,z,w\in F\right\}
\]
of $F^{2\times2}$. This is the ring of all lower-triangular $2\times
2$-matrices over $F$. (Yes, it is a subring of $F^{2\times2}$, since the sum
and the product of two lower-triangular matrices are lower-triangular and
since the zero and identity matrices are lower-triangular.) Since $F^{2}$ is a
right $F^{2\times2}$-module, $F^{2}$ must also be a right $F^{2\geq2}$-module
(by restriction). What are the $F^{2\geq2}$-submodules of $F^{2}$ ? Only
$F^{2}$ and $\left\{  0_{F^{2}}\right\}  $ and $\left\{  \left(  a,0\right)
\ \mid\ a\in F\right\}  $. (You might have to prove this on a future homework set.)
\end{itemize}

\begin{exercise}
\label{exe.mods.IM.assoc}Let $I$ and $J$ be two ideals of a ring $R$. Let $M$
be a left $R$-module. Prove that $\left(  IJ\right)  M=I\left(  JM\right)  $.

[Note that this generalizes the identity $\left(  IJ\right)  K=I\left(
JK\right)  $ in Proposition \ref{prop.rings.ideal-arith.laws} \textbf{(e)}.]
\end{exercise}

\subsection{\label{sec.modules.Zmods}Abelian groups as $\mathbb{Z}$-modules
(\cite[\S 10.1]{DumFoo04})}

Now, let us try to understand $\mathbb{Z}$-modules in particular.

\subsubsection{The action of $\mathbb{Z}$ by repeated addition}

Let us recall one of the most basic definitions in elementary mathematics: the
definition of multiplication of integers.

Multiplication of nonnegative integers was defined by repeated addition: If
$n,m\in\mathbb{N}$, then $nm$ means $\underbrace{m+m+\cdots+m}_{n\text{
times}}$. This same formula $nm=\underbrace{m+m+\cdots+m}_{n\text{ times}}$
can be applied to negative integers $m$ as well, but not to negative integers
$n$, since there is no such thing as $\underbrace{m+m+\cdots+m}_{-5\text{
times}}$. Thus, the product $nm$ for negative $n$ had to be defined
differently; one way to define it is by setting $nm=-\left(
\underbrace{m+m+\cdots+m}_{-n\text{ times}}\right)  $ (thus using negation to
reduce the case of negative $n$ to the case of positive $n$). Thus, for
arbitrary integers $n$ and $m$, the product $nm$ is defined by%
\[
nm=%
\begin{cases}
\underbrace{m+m+\cdots+m}_{n\text{ times}}, & \text{if }n\geq0;\\
-\left(  \underbrace{m+m+\cdots+m}_{-n\text{ times}}\right)  , & \text{if
}n<0.
\end{cases}
\]

The same definition can be adapted to any abelian group:

\begin{proposition}
\label{prop.Zmods.1}Let $A$ be an abelian group. Assume that $A$ is written
additively (i.e., the operation of $A$ is denoted by $+$, and the neutral
element by $0$). For any $n\in\mathbb{Z}$ and $a\in A$, define%
\begin{equation}
na=%
\begin{cases}
\underbrace{a+a+\cdots+a}_{n\text{ times}}, & \text{if }n\geq0;\\
-\left(  \underbrace{a+a+\cdots+a}_{-n\text{ times}}\right)  , & \text{if
}n<0.
\end{cases}
\label{eq.prop.Zmods.1.na=}%
\end{equation}
Thus, we have defined a map
\begin{align*}
\mathbb{Z}\times A  &  \rightarrow A,\\
\left(  n,a\right)   &  \mapsto na.
\end{align*}
We shall refer to this map as the \textbf{action of }$\mathbb{Z}$ \textbf{by
repeated addition} (due to the way $na$ was defined in
(\ref{eq.prop.Zmods.1.na=})).

\begin{enumerate}
\item[\textbf{(a)}] The group $A$ becomes a $\mathbb{Z}$-module (where we take
this map as the action of $\mathbb{Z}$ on $A$, and the pre-existing addition
of $A$ as the addition).

\item[\textbf{(b)}] This is the \textbf{only} $\mathbb{Z}$-module structure on
$A$. That is, if $A$ is \textbf{any }$\mathbb{Z}$-module, then the action of
$\mathbb{Z}$ on $A$ is given by the formula (\ref{eq.prop.Zmods.1.na=}) (and
is therefore uniquely determined by the abelian group structure on $A$).

\item[\textbf{(c)}] The $\mathbb{Z}$-submodules of $A$ are precisely the
subgroups of $A$.
\end{enumerate}
\end{proposition}

\begin{proof}
[Proof of Proposition \ref{prop.Zmods.1}.]LTTR. Here are the main ideas:

\textbf{(a)} You have to prove axioms like $\left(  n+m\right)  a=na+ma$ and
$n\left(  a+b\right)  =na+nb$ and $\left(  nm\right)  a=n\left(  ma\right)  $
for all $n,m\in\mathbb{Z}$ and $a,b\in A$. These facts are commonly proved for
$A=\mathbb{Z}$ in standard texts on the construction of the number system; if
you pick the \textquotedblleft right\textquotedblright\ proofs, then you can
adapt them to the general case just by replacing $\mathbb{Z}$ by $A$. The main
idea is \textquotedblleft reduce to the case when $n$ and $m$ are nonnegative,
and then prove them by induction on $n$ and $m$\textquotedblright. The details
are rather laborious, as there are several cases to discuss based on the signs
of $n$, $m$ and $n+m$. \medskip

\textbf{(b)} Given \textbf{any} $\mathbb{Z}$-module structure on $A$, we must
have%
\begin{align*}
na  &  =\underbrace{\left(  1+1+\cdots+1\right)  }_{n\text{ times}%
}a=\underbrace{1a+1a+\cdots+1a}_{n\text{ times}}\ \ \ \ \ \ \ \ \ \ \left(
\text{by Proposition \ref{prop.submods.distribs} \textbf{(a)}}\right) \\
&  =\underbrace{a+a+\cdots+a}_{n\text{ times}}\ \ \ \ \ \ \ \ \ \ \left(
\text{by the }1a=a\text{ axiom}\right)
\end{align*}
for any $n\in\mathbb{N}$ and any $a\in A$. This proves the \textquotedblleft
top half\textquotedblright\ of (\ref{eq.prop.Zmods.1.na=}). It is not hard to
prove the \textquotedblleft bottom half\textquotedblright\ either (use the
right distributivity axiom to see that $na+\left(  -n\right)
a=\underbrace{\left(  n+\left(  -n\right)  \right)  }_{=0}a=0a=0$). \medskip

\textbf{(c)} Proposition \ref{prop.submods.subgp} (applied to $R=\mathbb{Z}$
and $M=A$) shows that any $\mathbb{Z}$-submodule of $A$ is a subgroup of $A$.
Conversely, we must prove that if $B$ is a subgroup of $A$, then $B$ is a
$\mathbb{Z}$-submodule of $A$. So let $B$ be a subgroup of $A$. Then, any
$n\in\mathbb{Z}$ and $b\in B$ satisfy%
\[
nb=%
\begin{cases}
\underbrace{b+b+\cdots+b}_{n\text{ times}}, & \text{if }n\geq0;\\
-\left(  \underbrace{b+b+\cdots+b}_{-n\text{ times}}\right)  , & \text{if }n<0
\end{cases}
\ \ \in B
\]
(since $B$ is closed under addition and negation and contains $0$). In other
words, $B$ is closed under scaling. Hence, $B$ is a $\mathbb{Z}$-submodule of
$A$ (since $B$ is a subgroup of $A$ and therefore closed under addition and
contains $0$), qed.
\end{proof}

Proposition \ref{prop.Zmods.1} reveals what $\mathbb{Z}$-modules really are:
In general, when $R$ is a ring, an $R$-module is an abelian group $A$ with an
extra structure (namely, an action of $R$ on $A$); however, for $R=\mathbb{Z}%
$, this extra structure is redundant (in the sense that it can always be
constructed in a unique way from the abelian group structure), and so a
$\mathbb{Z}$-module is just an abelian group in fancy clothes.\footnote{Don't
get me wrong: The words \textquotedblleft redundant\textquotedblright\ and
\textquotedblleft in fancy clothes\textquotedblright\ don't mean
\textquotedblleft useless\textquotedblright; they just mean that the scaling
is determined by the abelian group structure.} Thus, we shall identify abelian
groups with $\mathbb{Z}$-modules (at least when the abelian groups are written additively).

This has a rather convenient consequence: The theory of $R$-modules is a
generalization of the theory of abelian groups. In particular, anything we
have proved or will prove for $R$-modules can therefore be applied to abelian
groups (by setting $R=\mathbb{Z}$).

\subsubsection{A few words on $\mathbb{Q}$-modules and $\mathbb{R}$-modules}

Thus, we have understood what $\mathbb{Z}$-modules are. What about
$\mathbb{Q}$-modules? Not every abelian group can be made into a $\mathbb{Q}$-module:

\begin{example}
\label{exa.Qmods.notZ/2}There is no $\mathbb{Q}$-module structure on
$\mathbb{Z}/2$ (that is, there is no $\mathbb{Q}$-module whose additive group
is $\mathbb{Z}/2$).
\end{example}

\begin{proof}
This follows from linear algebra (since $\mathbb{Q}$-modules are $\mathbb{Q}%
$-vector spaces and thus have dimensions; but $\mathbb{Z}/2$ is too large to
have dimension $0$ and yet too small to have dimension $>0$). Alternatively,
you can do it by hand: Assume that $\mathbb{Z}/2$ is a $\mathbb{Q}$-module in
some way. Then,%
\[
\dfrac{1}{2}\cdot\left(  2\cdot\overline{1}\right)  =\underbrace{\left(
\dfrac{1}{2}\cdot2\right)  }_{=1}\cdot\overline{1}=1\cdot\overline
{1}=\overline{1},
\]
so that%
\[
\overline{1}=\dfrac{1}{2}\cdot\underbrace{\left(  2\cdot\overline{1}\right)
}_{=\overline{0}}=\dfrac{1}{2}\cdot\overline{0}=\overline{0},
\]
which contradicts $\overline{1}\neq\overline{0}$.
\end{proof}

Thus we see that not every abelian group can be made into a $\mathbb{Q}%
$-module (unlike for $\mathbb{Z}$-modules). However, any abelian group that
can be made into a $\mathbb{Q}$-module can only be made so in one way:

\begin{exercise}
\label{exe.21hw3.3}Let $A$ be an abelian group (written additively). Prove
that there is at most one map $\mathbb{Q}\times A\rightarrow A$ that makes $A$
into a $\mathbb{Q}$-module (where we take this map as the action of
$\mathbb{Q}$ on $A$, and the pre-existing addition of $A$ as the addition).
\medskip

[\textbf{Hint:} Let $\ast_{1}$ and $\ast_{2}$ denote two such maps. Let
$r,s\in\mathbb{Z}$ and $a\in A$ with $s\neq0$. Your goal is to show that
$\dfrac{r}{s}\ast_{1}a=\dfrac{r}{s}\ast_{2}a$. First prove that if two
elements $u,v\in A$ satisfy $s\ast_{1}u=s\ast_{1}v$, then $u=v$.]
\end{exercise}

What about $\mathbb{R}$-modules? Here, we get neither existence nor
uniqueness: There are abelian groups that cannot be made into $\mathbb{R}%
$-modules; there are also abelian groups that can be made into $\mathbb{R}%
$-modules in multiple different ways. So the action of $\mathbb{R}$ on an
$\mathbb{R}$-module cannot be reconstructed from the underlying group of the
latter (unlike for $\mathbb{Z}$ and $\mathbb{Q}$). \textquotedblleft
Most\textquotedblright\ rings behave more like $\mathbb{R}$ than like
$\mathbb{Z}$ and $\mathbb{Q}$ in this regard.

\begin{exercise}
\ \ 

\begin{enumerate}
\item[\textbf{(a)}] Prove that any nontrivial $\mathbb{R}$-module is
uncountable as a set.

\item[\textbf{(b)}] Conclude that $\mathbb{Q}$ cannot be an $\mathbb{R}%
$-module (no matter how we define an action of $\mathbb{R}$ on $\mathbb{Q}$).
\end{enumerate}

[\textbf{Hint:} If $M$ is a nontrivial $\mathbb{R}$-module, and if $m\in M$ is
nonzero, then argue that the elements $rm$ for all $r\in\mathbb{R}$ have to be distinct.]
\end{exercise}

\subsubsection{Repeated addition vs. scaling}

If $R$ is any ring and $M$ is any $R$-module, then $M$ is (in particular) an
additive abelian group, and thus (by Proposition \ref{prop.Zmods.1}
\textbf{(a)}) becomes a $\mathbb{Z}$-module in a natural way (using the action
of $\mathbb{Z}$ by repeated addition). How does this $\mathbb{Z}$-module
structure relate to the original $R$-module structure on $M$ ? The following
proposition shows a certain consistency between the two:

\begin{proposition}
Let $R$ be a ring. Let $M$ be a left $R$-module. Then,%
\[
\left(  nr\right)  a=r\left(  na\right)  =n\left(  ra\right)
\ \ \ \ \ \ \ \ \ \ \text{for all }n\in\mathbb{Z}\text{, }r\in R\text{ and
}a\in M.
\]
Here, we are using both the action of $R$ on $M$ (to make sense of expressions
like $ra$) and the action of $\mathbb{Z}$ by repeated addition (to make sense
of expressions like $na$).
\end{proposition}

\begin{proof}
LTTR. (For $n\geq0$, this is just saying that
\[
\underbrace{\left(  r+r+\cdots+r\right)  }_{n\text{ times}}%
a=r\underbrace{\left(  a+a+\cdots+a\right)  }_{n\text{ times}}%
=\underbrace{ra+ra+\cdots+ra}_{n\text{ times}},
\]
which easily follows from Proposition \ref{prop.submods.distribs}. The case of
$n<0$ can be reduced to the case of $n>0$ using
(\ref{eq.prop.submodules.-rm.1}).)
\end{proof}

\subsection{\label{sec.modules.mors}Module morphisms (\cite[\S 10.2]%
{DumFoo04})}

\subsubsection{\label{subsec.modules.mors.def}Definition}

Module morphisms are defined similarly to ring morphisms, but you probably
already know their definition from linear algebra: they are also known as
linear maps. Let me recall the definition:

\begin{definition}
\label{def.modmor.def}Let $R$ be a ring. Let $M$ and $N$ be two left $R$-modules.

\begin{enumerate}
\item[\textbf{(a)}] A \textbf{left }$R$\textbf{-module homomorphism} (or, for
short, \textbf{left }$R$\textbf{-module morphism}, or \textbf{left }%
$R$\textbf{-linear map}) from $M$ to $N$ means a map $f:M\rightarrow N$ that

\begin{itemize}
\item \textbf{respects addition} (i.e., satisfies $f\left(  a+b\right)
=f\left(  a\right)  +f\left(  b\right)  $ for all $a,b\in M$);

\item \textbf{respects scaling} (i.e., satisfies $f\left(  ra\right)
=rf\left(  a\right)  $ for all $r\in R$ and $a\in M$);

\item \textbf{respects the zero} (i.e., satisfies $f\left(  0_{M}\right)
=0_{N}$).
\end{itemize}

You can drop the word \textquotedblleft left\textquotedblright\ (and, e.g.,
just say \textquotedblleft$R$-module morphism\textquotedblright)\ when it is
clear from the context.

\item[\textbf{(b)}] A \textbf{left }$R$\textbf{-module isomorphism} (or,
informally, \textbf{left }$R$\textbf{-module iso}) from $M$ to $N$ means an
invertible left $R$-module morphism $f:M\rightarrow N$ whose inverse
$f^{-1}:N\rightarrow M$ is also a left $R$-module morphism.

\item[\textbf{(c)}] The left $R$-modules $M$ and $N$ are said to be
\textbf{isomorphic} (this is written $M\cong N$) if there exists a left
$R$-module isomorphism $f:M\rightarrow N$.

\item[\textbf{(d)}] We let $\operatorname*{Hom}\nolimits_{R}\left(
M,N\right)  $ be the set of all left $R$-module morphisms from $M$ to $N$.

\item[\textbf{(e)}] Right $R$-module morphisms (and isomorphisms) are defined similarly.
\end{enumerate}
\end{definition}

It is not hard to show that the \textquotedblleft respects the
zero\textquotedblright\ axiom in Definition \ref{def.modmor.def} \textbf{(a)}
is redundant. (In fact, it is \textquotedblleft doubly
redundant\textquotedblright: It follows from each of the other two axioms!)

\subsubsection{\label{subsec.modules.mors.exas1}Simple examples}

Here are some examples of $R$-module morphisms:

\begin{itemize}
\item You have seen linear maps between vector spaces in linear algebra. These
are precisely the left $R$-module morphisms when $R$ is a field.

\item Let $k\in\mathbb{Z}$. The map
\begin{align*}
\mathbb{Z}  &  \rightarrow\mathbb{Z},\\
a  &  \mapsto ka
\end{align*}
is always a $\mathbb{Z}$-module morphism. (For comparison: It is a ring
morphism only when $k=1$.)

\item More generally: Let $R$ be a \textbf{commutative} ring. Let $k\in R$.
Let $M$ be any $R$-module. Then, the map
\begin{align*}
M  &  \rightarrow M,\\
a  &  \mapsto ka
\end{align*}
is an $R$-module morphism. (This is the map that we have called
\textquotedblleft scaling by $k$\textquotedblright.)

\item Even more generally: Let $R$ be any ring (commutative or not), and let
$k$ be a central element of $R$. Let $M$ be any left $R$-module. Then, the
map
\begin{align*}
M  &  \rightarrow M,\\
a  &  \mapsto ka
\end{align*}
is a left $R$-module morphism. Indeed, this map respects scaling because we
have%
\[
k\left(  ra\right)  =\underbrace{\left(  kr\right)  }%
_{\substack{=rk\\\text{(since }k\text{ is central)}}}a=\left(  rk\right)
a=r\left(  ka\right)  \ \ \ \ \ \ \ \ \ \ \text{for all }r\in R\text{ and
}a\in M.
\]
It respects addition and the zero for similar reasons.

However, if $k$ is not central, then this map is \textbf{not} a left
$R$-module morphism in general.

\item Let $R$ be a ring. Let $n\in\mathbb{N}$. For any $i\in\left\{
1,2,\ldots,n\right\}  $, the map%
\begin{align*}
\pi_{i}:R^{n}  &  \rightarrow R,\\
\left(  a_{1},a_{2},\ldots,a_{n}\right)   &  \mapsto a_{i}%
\end{align*}
is a left $R$-module morphism.

More generally: If $\left(  M_{i}\right)  _{i\in I}$ is a family of left
$R$-modules, and if $j\in I$, then the map%
\begin{align*}
\pi_{j}:\prod_{i\in I}M_{i}  &  \rightarrow M_{j},\\
\left(  m_{i}\right)  _{i\in I}  &  \mapsto m_{j}%
\end{align*}
is a left $R$-module morphism. This follows immediately from the fact that the
structure of $\prod_{i\in I}M_{i}$ (addition, action and zero) is defined entrywise.

\item Let $R$ be a ring. If $M$ and $N$ are two left $R$-modules, then the map%
\begin{align*}
M\times N  &  \rightarrow N\times M,\\
\left(  m,n\right)   &  \mapsto\left(  n,m\right)
\end{align*}
is an $R$-module isomorphism.

\item If $R$ is any ring, and $n,m\in\mathbb{N}$ are arbitrary, then the map%
\begin{align*}
R^{n\times m}  &  \rightarrow R^{m\times n},\\
A  &  \mapsto A^{T}%
\end{align*}
(which sends each matrix $A$ to its transpose) is a left $R$-module
isomorphism. (Indeed, it is an $R$-module morphism because of the formulas
$\left(  A+B\right)  ^{T}=A^{T}+B^{T}$ and $\left(  rA\right)  ^{T}=rA^{T}$
and $0_{n\times m}^{T}=0_{m\times n}$. It is an isomorphism because its
inverse is the analogous map from $R^{m\times n}$ to $R^{n\times m}$.)
\end{itemize}

The $\mathbb{Z}$-module morphisms (i.e., the $\mathbb{Z}$-linear maps) are
simply the group morphisms of additive groups:

\begin{proposition}
\label{prop.modules.mors.Z}Let $M$ and $N$ be two $\mathbb{Z}$-modules. Then,
the $\mathbb{Z}$-module morphisms from $M$ to $N$ are precisely the group
morphisms from $\left(  M,+,0\right)  $ to $\left(  N,+,0\right)  $. In other
words,%
\[
\operatorname*{Hom}\nolimits_{\mathbb{Z}}\left(  M,N\right)  =\left\{
\text{group morphisms }\left(  M,+,0\right)  \rightarrow\left(  N,+,0\right)
\right\}  .
\]

\end{proposition}

\begin{proof}
We have to show that any group morphism $f:\left(  M,+,0\right)
\rightarrow\left(  N,+,0\right)  $ automatically respects the scaling -- i.e.,
that it satisfies $f\left(  na\right)  =nf\left(  a\right)  $ for all
$n\in\mathbb{Z}$ and $a\in M$. This is LTTR.
\end{proof}

\begin{exercise}
Let $R$ be any ring. Consider the map%
\begin{align*}
S:R^{\mathbb{N}}  &  \rightarrow R^{\mathbb{N}},\\
\left(  a_{0},a_{1},a_{2},\ldots\right)   &  \mapsto\left(  a_{0}%
,\ a_{0}+a_{1},\ a_{0}+a_{1}+a_{2},\ a_{0}+a_{1}+a_{2}+a_{3},\ \ldots\right)
\\
&  \ \ \ \ \ \ \ \ \ \ =\left(  b_{0},b_{1},b_{2},\ldots\right)  \text{ where
}b_{i}=a_{0}+a_{1}+\cdots+a_{i}.
\end{align*}
Consider furthermore the map%
\begin{align*}
\Delta:R^{\mathbb{N}}  &  \rightarrow R^{\mathbb{N}},\\
\left(  a_{0},a_{1},a_{2},\ldots\right)   &  \mapsto\left(  a_{0}%
,\ a_{1}-a_{0},\ a_{2}-a_{1},\ a_{3}-a_{2},\ \ldots\right) \\
&  \ \ \ \ \ \ \ \ \ \ =\left(  c_{0},c_{1},c_{2},\ldots\right)  \text{ where
}c_{0}=a_{0}\text{ and }c_{i}=a_{i}-a_{i-1}\text{ for all }i\geq1.
\end{align*}

\begin{enumerate}
\item[\textbf{(a)}] Prove that $S$ and $\Delta$ are $R$-linear maps and are
mutually inverse.

\item[\textbf{(b)}] Recall the $R$-submodule
\[
R^{\left(  \mathbb{N}\right)  }=\left\{  \left(  a_{0},a_{1},a_{2}%
,\ldots\right)  \in R^{\mathbb{N}}\ \mid\ \text{only finitely many }%
i\in\mathbb{N}\text{ satisfy }a_{i}\neq0\right\}
\]
of $R^{\mathbb{N}}$. Define a further $R$-submodule $R^{\left(  \mathbb{N}%
\right)  +}$ of $R^{\mathbb{N}}$ by%
\begin{align*}
R^{\left(  \mathbb{N}\right)  +}  &  :=\left\{  \left(  a_{0},a_{1}%
,a_{2},\ldots\right)  \in R^{\mathbb{N}}\ \mid\ \text{there exists a }c\in
R\text{ such that}\right. \\
&  \ \ \ \ \ \ \ \ \ \ \ \ \ \ \ \ \ \ \ \ \left.  \text{ only finitely many
}i\in\mathbb{N}\text{ satisfy }a_{i}\neq c\right\}  .
\end{align*}
(Thus, a sequence $\left(  a_{0},a_{1},a_{2},\ldots\right)  \in R^{\mathbb{N}%
}$ belongs to $R^{\left(  \mathbb{N}\right)  +}$ if and only if starting from
some point on, all its entries are equal.)

Clearly, $R^{\left(  \mathbb{N}\right)  }$ is a proper subset of $R^{\left(
\mathbb{N}\right)  +}$ (unless $R$ is trivial).

Prove that $R^{\left(  \mathbb{N}\right)  }\cong R^{\left(  \mathbb{N}\right)
+}$ as left $R$-modules, and in fact the restriction of the map $S$ to
$R^{\left(  \mathbb{N}\right)  }$ is a left $R$-module isomorphism from
$R^{\left(  \mathbb{N}\right)  }$ to $R^{\left(  \mathbb{N}\right)  +}$.
\end{enumerate}
\end{exercise}

\subsubsection{\label{subsec.modules.mors.exas2}Ring morphisms as module
morphisms}

Let me give one more, slightly confusing example of module morphisms. Namely,
I claim that any ring morphism is a module morphism, as long as the module
structures are defined correctly (warning: these are often not the module
structures you expect!). To wit:

\begin{itemize}
\item Let $R$ and $S$ be two rings. Let $f:R\rightarrow S$ be a ring morphism.
As we have seen in Subsection \ref{subsec.modules.def.restr}, the ring $S$
then becomes a left $R$-module, with the action of $R$ on $S$ being defined by%
\[
rs=f\left(  r\right)  s\ \ \ \ \ \ \ \ \ \ \text{for all }r\in R\text{ and
}s\in S.
\]
This action is called the action on $S$ induced by $f$. It is now easy to see
that $f$ is a left $R$-module morphism from $R$ to $S$.

Here is an example. There is a ring morphism $f:\mathbb{C}\rightarrow
\mathbb{C}$ that sends each complex number $z=a+bi$ (with $a,b\in\mathbb{R}$)
to its complex conjugate $\overline{z}=a-bi$. Thus, from the previous
paragraph, we can conclude that this morphism $f$ is a $\mathbb{C}$-module
morphism from $\mathbb{C}$ to $\mathbb{C}$. But this is only true if the
$\mathbb{C}$-module structure on the target (but not on the domain) is the one
induced by $f$ (so it is given by $rs=f\left(  r\right)  s=\overline{r}s$ for
all $r\in\mathbb{C}$ and $s\in\mathbb{C}$), which is of course a rather
nonstandard choice of a $\mathbb{C}$-module structure on $\mathbb{C}$. So $f$
is indeed a $\mathbb{C}$-module morphism from $\mathbb{C}$ to $\mathbb{C}$,
but these are two different $\mathbb{C}$-modules $\mathbb{C}$ !

Of course, writing things like this is just inviting confusion. To avoid this
confusion, you need to introduce a new notation for the nonstandard
$\mathbb{C}$-module $\mathbb{C}$ (the one induced by $f$). Namely, let us
denote this new $\mathbb{C}$-module by $\overline{\mathbb{C}}$, while the
unadorned symbol $\mathbb{C}$ will always mean the old, obvious $\mathbb{C}%
$-module structure on $\mathbb{C}$ (in which the action is just the
multiplication). Thus, what we said in the previous paragraph can be restated
as follows: The map $f$ is a $\mathbb{C}$-module morphism from $\mathbb{C}$ to
$\overline{\mathbb{C}}$. Actually, it is easy to see that $f$ is a
$\mathbb{C}$-module \textbf{isomorphism} from $\mathbb{C}$ to $\overline
{\mathbb{C}}$. Thus, the $\mathbb{C}$-modules $\mathbb{C}$ and $\overline
{\mathbb{C}}$ are isomorphic (but still should not be identified to prevent confusion).

More generally, since $f:\mathbb{C}\rightarrow\mathbb{C}$ is a ring morphism,
we can restrict any $\mathbb{C}$-module $M$ to $\mathbb{C}$ via $f$. This
means the following: If $M$ is a $\mathbb{C}$-module, then we define a new
$\mathbb{C}$-module structure on $M$ by%
\[
rm=f\left(  r\right)  m=\overline{r}m\ \ \ \ \ \ \ \ \ \ \text{for all }%
r\in\mathbb{C}\text{ and }m\in M
\]
(where the \textquotedblleft$rm$\textquotedblright\ on the left hand side
refers to the new $\mathbb{C}$-module structure, whereas the \textquotedblleft%
$f\left(  r\right)  m$\textquotedblright\ and \textquotedblleft$\overline{r}%
m$\textquotedblright\ refer to the old one). This new $\mathbb{C}$-module is
called $\overline{M}$ (since calling it $M$ would be asking for trouble). It
is a \textquotedblleft twisted version\textquotedblright\ of $M$: It is
identical to $M$ as an abelian group, but the action of $\mathbb{C}$ on it has
been \textquotedblleft twisted\textquotedblright\ (in the sense that scaling
by $z$ on $\overline{M}$ is the same as scaling by $\overline{z}$ on $M$).

Here is a nice thing about these twisted $\mathbb{C}$-modules: If $V$ and $W$
are two $\mathbb{C}$-modules (i.e., $\mathbb{C}$-vector spaces), then a
$\mathbb{C}$-module morphism $g:V\rightarrow\overline{W}$ is what is known as
an \textbf{antilinear map} from $V$ to $W$ in linear algebra. Thus, antilinear
maps are \textquotedblleft secretly\textquotedblright\ just linear maps, once
you have twisted the vector space structure on the target.
\end{itemize}

\subsubsection{General properties of linearity}

We shall now state a bunch of general facts about module morphisms that are
analogous to some facts we have previously stated for ring morphisms. I won't
distract you with the proofs, as they are all straightforward.

We fix a ring $R$.

\begin{proposition}
\label{prop.modmor.invertible-iso}Let $M$ and $N$ be two left $R$-modules. Let
$f:M\rightarrow N$ be an invertible left $R$-module morphism. Then, $f$ is a
left $R$-module isomorphism.
\end{proposition}

\begin{proposition}
\label{prop.modmor.compose}Let $M$, $N$ and $P$ be three left $R$-modules. Let
$f:N\rightarrow P$ and $g:M\rightarrow N$ be two left $R$-module morphisms.
Then, $f\circ g:M\rightarrow P$ is a left $R$-module morphism.
\end{proposition}

\begin{proposition}
\label{prop.modiso.compose}Let $M$, $N$ and $P$ be three left $R$-modules. Let
$f:N\rightarrow P$ and $g:M\rightarrow N$ be two left $R$-module isomorphisms.
Then, $f\circ g:M\rightarrow P$ is a left $R$-module isomorphism.
\end{proposition}

\begin{proposition}
\label{prop.modiso.inverse}Let $M$ and $N$ be two left $R$-modules. Let
$f:M\rightarrow N$ be a left $R$-module isomorphism. Then, $f^{-1}%
:N\rightarrow M$ is a left $R$-module isomorphism.
\end{proposition}

\begin{corollary}
\label{cor.modiso.equiv}The relation $\cong$ for left $R$-modules is an
equivalence relation.
\end{corollary}

Left $R$-module isomorphisms preserve all \textquotedblleft
intrinsic\textquotedblright\ properties of left $R$-modules (just like ring
isomorphisms do for properties of rings). For example, if $M$ and $N$ are two
isomorphic left $R$-modules, then $M$ has as many $R$-submodules as $N$ does
(and there is a one-to-one correspondence between the $R$-submodules of $M$
and those of $N$).

All of this holds just as well for right $R$-modules; by now this is so
obvious that we don't even need to say it. (Besides, as you have seen in
Exercise \ref{exe.21hw2.2d}, right $R$-modules can be transformed into left
$R^{\operatorname*{op}}$-modules for a certain ring $R^{\operatorname*{op}}$.
This can also be done in reverse, and thus provides a dictionary between left
modules and right modules, which can always be used to translate a statement
about one kind of modules into a statement about the other. Module morphisms
behave as one would expect under this dictionary: When we use this dictionary
to turn two right $R$-modules $M$ and $N$ into left $R^{\operatorname*{op}}%
$-modules, the right $R$-module morphisms from $M$ to $N$ become the left
$R^{\operatorname*{op}}$-module morphisms from $M$ to $N$. This gives you all
excuses you might ever need to ignore right $R$-modules and only work with
left $R$-modules, until you actually need certain \textquotedblleft
hybrid\textquotedblright\ modules with both left and right structures.)

\subsubsection{Adding, subtracting and scaling $R$-linear maps}

In a way, $R$-linear maps (i.e., $R$-module morphisms) behave even better than
ring morphisms: If you add two ring morphisms $f$ and $g$ pointwise (i.e.,
form the map that sends every $r$ to $f\left(  r\right)  +g\left(  r\right)
$), then the resulting map will not usually be a ring morphism. Meanwhile,
$R$-linear maps can be added and sometimes scaled:

\begin{exercise}
\label{exe.21hw3.10ab}Let $R$ be a ring. Let $M$ and $N$ be two left
$R$-modules. Recall that $\operatorname{Hom}_{R}\left(  M,N\right)  $ is the
set of all left $R$-module morphisms from $M$ to $N$.

Prove the following:

\begin{enumerate}
\item[\textbf{(a)}] The map
\begin{align*}
M  &  \rightarrow N,\\
m  &  \mapsto0_{N}%
\end{align*}
is an $R$-module morphism (i.e., is $R$-linear). We shall call it the
\textbf{zero morphism} and denote it by $\mathbf{0}$ (a boldfaced zero).

\item[\textbf{(b)}] If $f\in\operatorname{Hom}_{R}\left(  M,N\right)  $ and
$g\in\operatorname{Hom}_{R}\left(  M,N\right)  $ are two $R$-linear maps from
$M$ to $N$, then the map
\begin{align*}
M  &  \rightarrow N,\\
m  &  \mapsto f\left(  m\right)  +g\left(  m\right)
\end{align*}
is also an $R$-module morphism. We shall denote the latter map by $f+g$, and
we shall call it the \textbf{(pointwise) sum} of $f$ and $g$.

\item[\textbf{(c)}] The set $\operatorname{Hom}_{R}\left(  M,N\right)  $
becomes an additive abelian group if we define addition pointwise (i.e., for
any $f\in\operatorname{Hom}_{R}\left(  M,N\right)  $ and $g\in
\operatorname{Hom}_{R}\left(  M,N\right)  $, we define $f+g$ as in part
\textbf{(b)} of this exercise). Its neutral element is the zero morphism
$\mathbf{0}:M\rightarrow N$ defined in part \textbf{(a)} of this exercise.
This group $\operatorname{Hom}_{R}\left(  M,N\right)  $ is called the
\textbf{Hom group} of $M$ and $N$.

\item[\textbf{(d)}] If $r$ is a central element of $R$ (see Definition
\ref{def.center} \textbf{(a)} for the meaning of \textquotedblleft
central\textquotedblright), and if $f\in\operatorname*{Hom}\nolimits_{R}%
\left(  M,N\right)  $ is an $R$-linear map from $M$ to $N$, then the map
\begin{align*}
M  &  \rightarrow N,\\
m  &  \mapsto rf\left(  m\right)
\end{align*}
is again $R$-linear (i.e., belongs to $\operatorname*{Hom}\nolimits_{R}\left(
M,N\right)  $). We shall denote this latter map by $rf$.

\item[\textbf{(e)}] Find an example where the claim of part \textbf{(d)} can
go wrong if $r$ is not assumed to be central.

\item[\textbf{(f)}] If $R$ is commutative, then the Hom group
$\operatorname{Hom}_{R}\left(  M,N\right)  $ (defined in part \textbf{(c)})
becomes an $R$-module, where the action is defined as follows: For any $r\in
R$ and any $f\in\operatorname{Hom}_{R}\left(  M,N\right)  $, we define
$rf\in\operatorname{Hom}_{R}\left(  M,N\right)  $ as in part \textbf{(d)}.
(This is allowed because in a commutative ring $R$, every element $r$ is central.)
\end{enumerate}
\end{exercise}

\subsubsection{Kernels and images}

Next, we shall study kernels and images of module morphisms.

Again, we fix a ring $R$.

\begin{definition}
Let $M$ and $N$ be two left $R$-modules. Let $f:M\rightarrow N$ be a left
$R$-module morphism. Then, the \textbf{kernel} of $f$ (denoted $\ker f$ or
$\operatorname*{Ker}f$) is defined to be the subset%
\[
\operatorname*{Ker}f:=\left\{  a\in M\ \mid\ f\left(  a\right)  =0_{N}%
\right\}
\]
of $M$.
\end{definition}

Some examples:

\begin{itemize}
\item Let $R$ be a commutative ring. Let $b\in R$. Then, the map
\begin{align*}
R  &  \rightarrow R,\\
r  &  \mapsto br
\end{align*}
is an $R$-module morphism (check this!). The kernel of this map is%
\[
\left\{  r\in R\ \mid\ br=0\right\}  .
\]
We thus conclude that this kernel is $\left\{  0\right\}  $ if and only if $b$
is regular (see Definition \ref{def.rings.reg.reg}).

\item Both $\mathbb{Z}^{3}$ and $\mathbb{Z}\times\left(  \mathbb{Z}/2\right)
$ are $\mathbb{Z}$-modules (since we have seen in Proposition
\ref{prop.Zmods.1} that every additive group is a $\mathbb{Z}$-module). The
map%
\begin{align*}
\mathbb{Z}^{3}  &  \rightarrow\mathbb{Z}\times\left(  \mathbb{Z}/2\right)  ,\\
\left(  a,b,c\right)   &  \mapsto\left(  a-b,\ \overline{b-c}\right)
\end{align*}
is a $\mathbb{Z}$-module morphism. Its kernel is%
\begin{align*}
&  \left\{  \left(  a,b,c\right)  \in\mathbb{Z}^{3}\ \mid\ \left(
a-b,\ \overline{b-c}\right)  =0_{\mathbb{Z}\times\left(  \mathbb{Z}/2\right)
}\right\} \\
&  =\left\{  \left(  a,b,c\right)  \in\mathbb{Z}^{3}\ \mid\ a-b=0\text{ and
}\overline{b-c}=0\right\} \\
&  =\left\{  \left(  a,b,c\right)  \in\mathbb{Z}^{3}\ \mid\ a-b=0\text{ and
}b-c\equiv0\operatorname{mod}2\right\} \\
&  =\left\{  \left(  a,b,c\right)  \in\mathbb{Z}^{3}\ \mid\ a=b\text{ and
}b\equiv c\operatorname{mod}2\right\}  .
\end{align*}

\end{itemize}

Kernels are a standard concept in linear algebra, where they are also called
\textbf{nullspaces}. The following facts should be familiar from abstract
linear algebra (but are also pretty easy to prove):

\begin{theorem}
\label{thm.modmor.ker-ideal}Let $M$ and $N$ be two left $R$-modules. Let
$f:M\rightarrow N$ be a left $R$-module morphism. Then, the kernel
$\operatorname*{Ker}f$ of $f$ is a left $R$-submodule of $M$, whereas the
image $\operatorname{Im}f=f\left(  M\right)  $ of $f$ is a left $R$-submodule
of $N$.
\end{theorem}

\begin{lemma}
\label{lem.modmor.ker-inj}Let $M$ and $N$ be two left $R$-modules. Let
$f:M\rightarrow N$ be a left $R$-module morphism. Then, $f$ is injective if
and only if $\operatorname*{Ker}f=\left\{  0_{M}\right\}  $.
\end{lemma}

Note that Lemma \ref{lem.modmor.ker-inj} is an analogue of Lemma
\ref{lem.rings.mors.inj}.

\begin{exercise}
Prove Theorem \ref{thm.modmor.ker-ideal} and Lemma \ref{lem.modmor.ker-inj}.
\end{exercise}

\subsection{\label{sec.modules.quot}Quotient modules}

We fix a ring $R$ for the entirety of Section \ref{sec.modules.quot}.

\subsubsection{Definition}

We shall next define quotient modules of left $R$-modules, in more or less the
same way as we defined quotient rings of rings (but this time we need to
establish an action instead of a multiplication on the quotient):

\begin{definition}
\label{def.quotmod}Let $I$ be a left $R$-submodule of a left $R$-module $M$.
Thus, $I$ is a subgroup of the additive group $\left(  M,+,0\right)  $, hence
a normal subgroup (since $\left(  M,+,0\right)  $ is abelian). Therefore, the
quotient group $M/I$ itself becomes an abelian group. Its elements are the
cosets $a+I$ of $I$ in $M$.

Note that the addition on $M/I$ is given by
\begin{equation}
\left(  a+I\right)  +\left(  b+I\right)  =\left(  a+b\right)
+I\ \ \ \ \ \ \ \ \ \ \text{for all }a,b\in M. \label{eq.def.quotmod.+}%
\end{equation}

We now define an action of $R$ on $M/I$ by setting%
\begin{equation}
r\left(  a+I\right)  =ra+I\ \ \ \ \ \ \ \ \ \ \text{for all }r\in R\text{ and
}a\in M. \label{eq.def.quotmod.*}%
\end{equation}
(See below for a proof that this is well-defined.)

The set $M/I$, equipped with the addition and the action we just defined and
with the element $0+I$ as zero vector, is a left $R$-module. This left
$R$-module is called the \textbf{quotient left }$R$\textbf{-module} of $M$ by
the submodule $I$; it is also pronounced \textquotedblleft$M$ \textbf{modulo}
$I$\textquotedblright. It is denoted $M/I$ (so when you hear \textquotedblleft
the left $R$-module $M/I$\textquotedblright, it always means the set $M/I$
equipped with the structure just mentioned).

The cosets $a+I$ are called \textbf{residue classes} modulo $I$, and are often
denoted $a\operatorname{mod}I$ or $\left[  a\right]  _{I}$ or $\left[
a\right]  $ or $\overline{a}$. (The last two notations are used when $I$ is
clear from the context.)
\end{definition}

\begin{theorem}
\label{thm.quotmod.welldef}Let $M$ and $I$ be as in Definition
\ref{def.quotmod}. Then, the action of $R$ on $M/I$ is well-defined, and $M/I$
does indeed become a left $R$-module when endowed with the operations and
elements just described.
\end{theorem}

Theorem \ref{thm.quotmod.welldef} is an analogue of Theorem
\ref{thm.quotring.welldef} for modules instead of rings, and its proof is
analogous as well.

Using the notation $\overline{a}$ for the coset $a+I$, we can rewrite the
formulas (\ref{eq.def.quotmod.+}) and (\ref{eq.def.quotmod.*}) as%
\begin{equation}
\overline{a}+\overline{b}=\overline{a+b}\ \ \ \ \ \ \ \ \ \ \text{for all
}a,b\in M \label{eq.def.quotmod.+bar}%
\end{equation}
and%
\begin{equation}
r\cdot\overline{a}=\overline{ra}\ \ \ \ \ \ \ \ \ \ \text{for all }r\in
R\text{ and }a\in M. \label{eq.def.quotmod.*bar}%
\end{equation}
The zero vector $0+I$ of the quotient $R$-module $M/I$ can, of course, be
written as $\overline{0}$.

\begin{fineprint}
\begin{remark}
Note that the residue classes $\overline{a}=a+I$ in Definition
\ref{def.quotmod} are precisely the equivalence classes of the
\textquotedblleft congruent modulo $I$\textquotedblright\ relation defined in
Exercise \ref{exe.mods.modulo}. Thus, the quotient $R$-module $M/I$
generalizes the classical notion of modular arithmetic in $\mathbb{Z}/n$.
\end{remark}
\end{fineprint}

Theorem \ref{thm.quotring.canproj}, too, has an analogue for modules:

\begin{theorem}
\label{thm.quotmod.canproj}Let $I$ be a left $R$-submodule of a left
$R$-module $M$. Consider the map
\begin{align*}
\pi:M  &  \rightarrow M/I,\\
a  &  \mapsto a+I.
\end{align*}
Then, $\pi$ is a surjective $R$-module morphism with kernel $I$.
\end{theorem}

\begin{definition}
This morphism $\pi$ is called the \textbf{canonical projection} from $M$ onto
$M/I$.
\end{definition}

The proof of Theorem \ref{thm.quotmod.canproj} is analogous to the proof of
Theorem \ref{thm.quotring.canproj}.

\subsubsection{Examples}

Examples of quotient modules can be easily created from various sources:

\begin{itemize}
\item Quotients of abelian groups are instances of quotient modules, since
abelian groups are $\mathbb{Z}$-modules.

\item Quotients of vector spaces are instances of quotient modules, since
vector spaces are modules over a field.

For instance, consider the $3$-dimensional vector space (i.e., $\mathbb{R}%
$-module) $\mathbb{R}^{3}$ over the ring $\mathbb{R}$ of real numbers. This
vector space $\mathbb{R}^{3}$ is typically viewed as a model for
three-dimensional space. Define a vector subspace (i.e., $\mathbb{R}%
$-submodule) $I$ of $\mathbb{R}^{3}$ by%
\[
I=\left\{  \left(  x,y,z\right)  \in\mathbb{R}^{3}\ \mid\ x+y+z=0\right\}  .
\]
Geometrically, this is a hyperplane through the origin of $\mathbb{R}^{3}$.
Now, consider the quotient $\mathbb{R}$-module (i.e., quotient vector space)
$\mathbb{R}^{3}/I$. Its elements are residue classes of the form
$\overline{\left(  x,y,z\right)  }$, where two vectors $\left(  x,y,z\right)
$ and $\left(  x^{\prime},y^{\prime},z^{\prime}\right)  $ belong to the same
residue class if and only if their entrywise difference $\left(  x-x^{\prime
},\ y-y^{\prime},\ z-z^{\prime}\right)  $ belongs to $I$ (that is, if we have
$\left(  x-x^{\prime}\right)  +\left(  y-y^{\prime}\right)  +\left(
z-z^{\prime}\right)  =0$). For instance, the two residue classes
$\overline{\left(  3,0,0\right)  }$ and $\overline{\left(  1,1,1\right)  }$
are identical (since $\left(  3-1\right)  +\left(  0-1\right)  +\left(
0-1\right)  =0$), but the two residue classes $\overline{\left(  1,0,0\right)
}$ and $\overline{\left(  2,0,0\right)  }$ are not. It is not hard to see that
each element of $\mathbb{R}^{3}/I$ can be uniquely written in the form
$\overline{\left(  r,0,0\right)  }$ for some $r\in\mathbb{R}$. This shows that
the vector space $\mathbb{R}^{3}/I$ is $1$-dimensional.

\item If $R$ is any ring, and $M$ is any left $R$-module, then the two obvious
$R$-submodules $\left\{  0_{M}\right\}  $ and $M$ of $M$ lead to uninteresting
quotient modules: The quotient module $M/\left\{  0_{M}\right\}  $ is
isomorphic to $M$, whereas the quotient module $M/M$ is trivial (i.e., has
only one element).

\item Let $R$ be a ring. As we recall from Subsection
\ref{subsec.modules.def.dirprod}, the left $R$-module $R^{\mathbb{N}}$ has an
$R$-submodule $R^{\left(  \mathbb{N}\right)  }$. What does the quotient module
$R^{\mathbb{N}}/R^{\left(  \mathbb{N}\right)  }$ look like? Its elements are
residue classes of the form $\overline{\left(  a_{0},a_{1},a_{2}%
,\ldots\right)  }$, where two infinite sequences $\left(  a_{0},a_{1}%
,a_{2},\ldots\right)  $ and $\left(  b_{0},b_{1},b_{2},\ldots\right)  $ belong
to the same residue class if and only if their entrywise difference
\newline$\left(  a_{0}-b_{0},\ a_{1}-b_{1},\ a_{2}-b_{2},\ \ldots\right)  $
belongs to $R^{\left(  \mathbb{N}\right)  }$ (that is, if the two sequences
$\left(  a_{0},a_{1},a_{2},\ldots\right)  $ and $\left(  b_{0},b_{1}%
,b_{2},\ldots\right)  $ agree at all but finitely many positions). Thus, we
can view an element $\overline{\left(  a_{0},a_{1},a_{2},\ldots\right)  }$ of
$R^{\mathbb{N}}/R^{\left(  \mathbb{N}\right)  }$ as an \textquotedblleft
infinite sequence determined up to finite change\textquotedblright\ (where
\textquotedblleft finite change\textquotedblright\ means changing finitely
many entries). This kind of construction is frequent in analysis: For
instance, the limit $\lim\limits_{n\rightarrow\infty}a_{n}$ of a sequence
$\left(  a_{0},a_{1},a_{2},\ldots\right)  $ of real numbers does not depend on
finite changes (i.e., it does not change if we change finitely many entries of
our sequence), and thus (if it exists) can be viewed as a property of the
residue class $\overline{\left(  a_{0},a_{1},a_{2},\ldots\right)  }\in
R^{\mathbb{N}}/R^{\left(  \mathbb{N}\right)  }$ (for $R=\mathbb{R}$).
\end{itemize}

\subsubsection{The universal property of quotient modules}

The universal property of quotient rings (Theorem \ref{thm.quotring.uniprop1}%
), too, has an analogue for modules:

\begin{theorem}
[Universal property of quotient modules, elementwise form]%
\label{thm.quotmod.uniprop1}Let $M$ be a left $R$-module. Let $I$ be a left
$R$-submodule of $M$.

Let $N$ be a left $R$-module. Let $f:M\rightarrow N$ be a left $R$-module
morphism. Assume that $f\left(  I\right)  =0$ (this is shorthand for saying
that $f\left(  a\right)  =0$ for all $a\in I$). Then, the map%
\begin{align*}
f^{\prime}:M/I  &  \rightarrow N,\\
\overline{m}  &  \mapsto f\left(  m\right)  \ \ \ \ \ \ \ \ \ \ \left(
\text{for all }m\in M\right)
\end{align*}
is well-defined (i.e., the value $f\left(  m\right)  $ depends only on the
residue class $\overline{m}$, not on $m$ itself) and is a left $R$-module morphism.
\end{theorem}

The proof of Theorem \ref{thm.quotmod.uniprop1} is analogous to the proof of
Theorem \ref{thm.quotring.uniprop1}.

The abstract form of the universal property of quotient rings (Theorem
\ref{thm.quotring.uniprop2}) has an analogue for modules as well:

\begin{theorem}
[Universal property of quotient modules, abstract form]%
\label{thm.quotmod.uniprop2}Let $M$ be a left $R$-module. Let $I$ be a left
$R$-submodule of $M$. Consider the canonical projection $\pi:M\rightarrow M/I$.

Let $N$ be a left $R$-module. Let $f:M\rightarrow N$ be a left $R$-module
morphism. Assume that $f\left(  I\right)  =0$ (this is shorthand for saying
that $f\left(  a\right)  =0$ for all $a\in I$). Then, there is a unique left
$R$-module morphism $f^{\prime}:M/I\rightarrow N$ satisfying $f=f^{\prime
}\circ\pi$.
\end{theorem}

Just to unravel the abstract definition: This morphism $f^{\prime}$ is exactly
the morphism $f^{\prime}$ from Theorem \ref{thm.quotmod.uniprop1}, i.e., it
sends each coset (= residue class) $\overline{m}=m+I\in M/I$ to $f\left(
m\right)  $.

The proof of Theorem \ref{thm.quotmod.uniprop2} is analogous to the proof of
Theorem \ref{thm.quotring.uniprop2}.

The equality $f=f^{\prime}\circ\pi$ in Theorem \ref{thm.quotmod.uniprop2} is
oftentimes restated as follows: The diagram%
\[%
\xymatrix@C=4pc{
M \ar[d]_{\pi} \ar[dr]^f \\
M/I \ar@{-->}[r]_{f'} & N
}%
\]
commutes.

\subsubsection{The First Isomorphism Theorem for modules}

The First Isomorphism Theorem for rings (Theorem \ref{thm.1it.ring1}) also has
a counterpart for $R$-modules (with an equally straightforward proof):

\begin{theorem}
[First Isomorphism Theorem for modules, elementwise form]\label{thm.1it.mod1}%
Let $M$ and $N$ be two left $R$-modules, and let $f:M\rightarrow N$ be a left
$R$-module morphism. Then:

\begin{enumerate}
\item[\textbf{(a)}] The kernel $\operatorname*{Ker}f$ is an $R$-submodule of
$M$. Thus, $M/\operatorname*{Ker}f$ is a quotient module of $M$. As a set,
$M/\operatorname*{Ker}f$ is precisely the set $M/f$ defined in Theorem
\ref{thm.1it.set} (applied to $M$ and $N$ instead of $R$ and $S$). The
$f$-classes (as defined in Theorem \ref{thm.1it.set}) are precisely the cosets
of $\operatorname*{Ker}f$.

\item[\textbf{(b)}] The image $f\left(  M\right)  :=\left\{  f\left(
m\right)  \ \mid\ m\in M\right\}  $ of $f$ is an $R$-submodule of $N$.

\item[\textbf{(c)}] The map%
\begin{align*}
f^{\prime}:M/\operatorname*{Ker}f  &  \rightarrow f\left(  M\right)  ,\\
\overline{a}  &  \mapsto f\left(  a\right)
\end{align*}
is well-defined and is a left $R$-module isomorphism.

\item[\textbf{(d)}] This map $f^{\prime}$ is precisely the map $f^{\prime}$
defined in Theorem \ref{thm.1it.set} \textbf{(c)} (applied to $M$ and $N$
instead of $R$ and $S$).

\item[\textbf{(e)}] Let $\pi:M\rightarrow M/\operatorname*{Ker}f$ denote the
\textbf{canonical projection} (i.e., the map that sends each $m\in M$ to its
coset $\overline{m}$). Let $\iota:f\left(  M\right)  \rightarrow N$ denote the
\textbf{canonical inclusion} (i.e., the map that sends each $n\in f\left(
M\right)  $ to $n$). Then, the map $f^{\prime}$ defined in part \textbf{(c)}
satisfies
\[
f=\iota\circ f^{\prime}\circ\pi.
\]
In other words, the diagram%
\begin{equation}%
\xymatrix@C=4pc{
M \ar[d]_{\pi} \ar[r]^f & N \\
M/\Ker f \ar[r]_{f'} & f\tup{M} \ar[u]^{\iota}
}
\label{eq.thm.1it.mod1.diag}%
\end{equation}
is commutative.

\item[\textbf{(f)}] We have $M/\operatorname*{Ker}f\cong f\left(  M\right)  $
as left $R$-modules.
\end{enumerate}
\end{theorem}

All results we have stated so far about modules are analogues of known results
about rings. So are their proofs (which is why we have omitted them). The
Second and the Third isomorphism theorem for rings (which you have seen in
Section \ref{sec.rings.234it}) also have analogues for modules.

\begin{remark}
\label{rmk.mods.IT1.rnull1}If you have done some abstract linear algebra, the
formula $M/\operatorname*{Ker}f\cong f\left(  M\right)  $ in Theorem
\ref{thm.1it.mod1} \textbf{(f)} might remind you of something.

Indeed, let $R$ be a field. Thus, $R$-modules are $R$-vector spaces. Let $M$
and $N$ be two finite-dimensional $R$-vector spaces. Let $f:M\rightarrow N$ be
a linear map. Thus, Theorem \ref{thm.1it.mod1} \textbf{(f)} yields that
$M/\operatorname*{Ker}f\cong f\left(  M\right)  $ as $R$-modules (i.e., as
$R$-vector spaces). However, isomorphic vector spaces have equal dimension.
Hence, from $M/\operatorname*{Ker}f\cong f\left(  M\right)  $, we obtain
\begin{equation}
\dim\left(  M/\operatorname*{Ker}f\right)  =\dim\left(  f\left(  M\right)
\right)  . \label{eq.rmk.mods.IT1.rnull1.1}%
\end{equation}

However, it is not hard to see (we will see it soon) that $\dim\left(
M/I\right)  =\dim M-\dim I$ whenever $I$ is a vector subspace of $M$. (The
idea behind this formula is that when you pass from $M$ to $M/I$, you are
\textquotedblleft collapsing\textquotedblright\ the \textquotedblleft
dimensions\textquotedblright\ contained in $I$ (since you are equating any
vector in $I$ with $0$), and thus the dimension of the vector space should go
down by $\dim I$. Formally speaking, this can be shown using bases. We will do
so below.)

As a consequence of the $\dim\left(  M/I\right)  =\dim M-\dim I$ formula, we
have%
\[
\dim\left(  M/\operatorname*{Ker}f\right)  =\dim M-\dim\left(
\operatorname*{Ker}f\right)  .
\]
Hence,%
\[
\dim M-\dim\left(  \operatorname*{Ker}f\right)  =\dim\left(
M/\operatorname*{Ker}f\right)  =\dim\left(  f\left(  M\right)  \right)
\ \ \ \ \ \ \ \ \ \ \left(  \text{by (\ref{eq.rmk.mods.IT1.rnull1.1})}\right)
.
\]
This is the \textbf{rank-nullity formula} from linear algebra (indeed,
$\dim\left(  \operatorname*{Ker}f\right)  $ is called the \textbf{nullity} of
$f$, whereas $\dim\left(  f\left(  M\right)  \right)  $ is called the
\textbf{rank} of $f$).
\end{remark}

Here are some more exercises related to quotient modules:\footnote{Recall once
again that a commutative ring $R$ is an $R$-module itself, and that its
$R$-submodules are precisely its ideals.}

\begin{exercise}
Let $R$ be a commutative ring. Let $N$ be any $R$-module. For any $R$-module
$M$, we define the $R$-module $\operatorname*{Hom}\nolimits_{R}\left(
M,N\right)  $ as in Exercise \ref{exe.21hw3.10ab} \textbf{(f)}.

\begin{enumerate}
\item[\textbf{(a)}] Prove that $\operatorname*{Hom}\nolimits_{R}\left(
R,N\right)  \cong N$ as $R$-modules. More precisely, prove that the map%
\begin{align*}
\operatorname*{Hom}\nolimits_{R}\left(  R,N\right)   &  \rightarrow N,\\
f  &  \mapsto f\left(  1\right)
\end{align*}
(which sends every $R$-linear map $f:R\rightarrow N$ to its value $f\left(
1\right)  $) is an $R$-module isomorphism.

\item[\textbf{(b)}] Let $I$ be an ideal of $R$. Let $N_{I}$ be the subset
$\left\{  a\in N\ \mid\ ia=0\text{ for all }i\in I\right\}  $ of $N$. Prove
that $N_{I}$ is an $R$-submodule of $N$, and that the map%
\begin{align*}
\operatorname*{Hom}\nolimits_{R}\left(  R/I,\ N\right)   &  \rightarrow
N_{I},\\
f  &  \mapsto f\left(  \overline{1}\right)
\end{align*}
is an $R$-module isomorphism.
\end{enumerate}
\end{exercise}

\begin{exercise}
For any two integers $n$ and $m$ with $m\neq0$, prove that
$\operatorname*{Hom}\nolimits_{\mathbb{Z}}\left(  \mathbb{Z}/n,\ \mathbb{Z}%
/m\right)  \cong\mathbb{Z}/\gcd\left(  n,m\right)  $ as $\mathbb{Z}$-modules.
(Here, the $\mathbb{Z}$-module $\operatorname*{Hom}\nolimits_{\mathbb{Z}%
}\left(  \mathbb{Z}/n,\ \mathbb{Z}/m\right)  $ is defined as in Exercise
\ref{exe.21hw3.10ab} \textbf{(f)}.)
\end{exercise}

\begin{exercise}
Let $I$ be any set. Let $\left(  M_{i}\right)  _{i\in I}$ be any family of
left $R$-modules. Let $N_{i}$ be an $R$-submodule of $M_{i}$ for each $i\in I$.

\begin{enumerate}
\item[\textbf{(a)}] Prove that $\left(  \prod_{i\in I}M_{i}\right)  /\left(
\prod_{i\in I}N_{i}\right)  \cong\prod_{i\in I}\left(  M_{i}/N_{i}\right)  $
as left $R$-modules. (The left hand side is well-defined by Exercise
\ref{exe.mods.dirprod-submods} \textbf{(a)}.)

\item[\textbf{(b)}] Prove that $\left(  \bigoplus\limits_{i\in I}M_{i}\right)
/\left(  \bigoplus\limits_{i\in I}N_{i}\right)  \cong\bigoplus\limits_{i\in
I}\left(  M_{i}/N_{i}\right)  $ as left $R$-modules. (The left hand side is
well-defined by Exercise \ref{exe.mods.dirprod-submods} \textbf{(b)}.)
\end{enumerate}
\end{exercise}

An analogue of the Chinese Remainder Theorem (Theorem \ref{thm.CRT-k-ids1} and
Theorem \ref{thm.CRT-k-ids2}) also exists for modules, although it still
involves ideals:

\begin{exercise}
\label{exe.21hw3.4}Prove the \textbf{Chinese Remainder Theorem for Modules}:

Let $R$ be a commutative ring. Let $I_{1},I_{2},\ldots,I_{k}$ be $k$ mutually
comaximal ideals of $R$. Let $M$ be a $R$-module. Then:

\begin{enumerate}
\item[\textbf{(a)}] We have $I_{1}M\cap I_{2}M\cap\cdots\cap I_{k}M=I_{1}%
I_{2}\cdots I_{k}M$. (See Proposition \ref{prop.mods.IM} for the definition of
$IM$ for any ideal $I$ of $R$. The notation \textquotedblleft$I_{1}M\cap
I_{2}M\cap\cdots\cap I_{k}M$\textquotedblright\ means \textquotedblleft%
$\left(  I_{1}M\right)  \cap\left(  I_{2}M\right)  \cap\cdots\cap\left(
I_{k}M\right)  $\textquotedblright.)

\item[\textbf{(b)}] There is an $R$-module isomorphism\footnotemark%
\[
M/\left(  I_{1}I_{2}\cdots I_{k}M\right)  \rightarrow\left(  M/I_{1}M\right)
\times\left(  M/I_{2}M\right)  \times\cdots\times\left(  M/I_{k}M\right)
\]
that sends each coset $m+I_{1}I_{2}\cdots I_{k}M$ to the $k$-tuple $\left(
m+I_{1}M,\ m+I_{2}M,\ \ldots,\ m+I_{k}M\right)  $.
\end{enumerate}
\end{exercise}

\footnotetext{The notation \textquotedblleft$M/IM$\textquotedblright\ (where
$I$ is an ideal of $R$) means \textquotedblleft$M/\left(  IM\right)
$\textquotedblright.}

\subsection{\label{sec.modules.bases}Spanning, linear independence, bases and
free modules (\cite[\S 10.3]{DumFoo04})}

Again, let us fix a ring $R$ for the entirety of Section
\ref{sec.modules.bases}.

\subsubsection{Definitions}

We shall now generalize some classical notions from linear algebra (spanning,
linear independence and bases) to arbitrary $R$-modules.

\begin{definition}
\label{def.mods.lincomb-et-al}Let $M$ be a left $R$-module. Let $m_{1}%
,m_{2},\ldots,m_{n}$ be finitely many vectors in $M$.

\begin{enumerate}
\item[\textbf{(a)}] A \textbf{linear combination} of $m_{1},m_{2},\ldots
,m_{n}$ means a vector of the form%
\[
r_{1}m_{1}+r_{2}m_{2}+\cdots+r_{n}m_{n}\ \ \ \ \ \ \ \ \ \ \text{with }%
r_{1},r_{2},\ldots,r_{n}\in R.
\]

\item[\textbf{(b)}] The set of all linear combinations of $m_{1},m_{2}%
,\ldots,m_{n}$ is called the \textbf{span} of $\left(  m_{1},m_{2}%
,\ldots,m_{n}\right)  $, and is denoted by $\operatorname*{span}\left(
m_{1},m_{2},\ldots,m_{n}\right)  $. (Note that \cite{DumFoo04} calls it
$R\left\{  m_{1},m_{2},\ldots,m_{n}\right\}  $.)

\item[\textbf{(c)}] If the span of $\left(  m_{1},m_{2},\ldots,m_{n}\right)  $
is $M$, then we say that the vectors $m_{1},m_{2},\ldots,m_{n}$ \textbf{span}
$M$ (or \textbf{generate} $M$).

\item[\textbf{(d)}] We say that the vectors $m_{1},m_{2},\ldots,m_{n}$ are
\textbf{linearly independent} if the following holds: If $r_{1},r_{2}%
,\ldots,r_{n}\in R$ satisfy%
\[
r_{1}m_{1}+r_{2}m_{2}+\cdots+r_{n}m_{n}=0,
\]
then $r_{1}=r_{2}=\cdots=r_{n}=0$. (In other words, the vectors $m_{1}%
,m_{2},\ldots,m_{n}$ are said to be linearly independent if the only way to
write $0$ as a linear combination of them is $0=0m_{1}+0m_{2}+\cdots+0m_{n}$.)

\item[\textbf{(e)}] We say that the $n$-tuple $\left(  m_{1},m_{2}%
,\ldots,m_{n}\right)  $ is a \textbf{basis} of the $R$-module $M$ if
$m_{1},m_{2},\ldots,m_{n}$ are linearly independent and span $M$.

\item[\textbf{(f)}] All of this terminology depends on $R$. Thus, if $R$ is
not clear from the context, we will clarify it by saying \textquotedblleft%
$R$-linear combination\textquotedblright\ (or \textquotedblleft linear
combination over $R$\textquotedblright) instead of just \textquotedblleft
linear combination\textquotedblright, and likewise saying \textquotedblleft%
$R$-span\textquotedblright\ or \textquotedblleft$R$-linearly
independent\textquotedblright\ or \textquotedblleft$R$-basis\textquotedblright.
\end{enumerate}
\end{definition}

Fine print: The property of $n$ vectors $m_{1},m_{2},\ldots,m_{n}$ to span $M$
is a joint property (i.e., it is a property of the \textbf{list }$\left(
m_{1},m_{2},\ldots,m_{n}\right)  $, not of each single vector). The same
applies to linear independence. Sometimes, we do say that a single vector $m$
spans $M$ (for example, the vector $1\in\mathbb{Z}$ spans the $\mathbb{Z}%
$-module $\mathbb{Z}$); this means that the one-element list $\left(
m\right)  $ spans $M$.

Definition \ref{def.mods.lincomb-et-al} was tailored to finite lists of
vectors, but we can extend it to arbitrary (possibly infinite) families of vectors:

\begin{definition}
\label{def.mods.lincomb-et-al-inf}Let $M$ be a left $R$-module. Let $\left(
m_{i}\right)  _{i\in I}$ be a family of vectors in $M$ (with $I$ being any set).

\begin{enumerate}
\item[\textbf{(a)}] A \textbf{linear combination} of $\left(  m_{i}\right)
_{i\in I}$ means a vector of the form%
\[
\sum_{i\in I}r_{i}m_{i}%
\]
for some family $\left(  r_{i}\right)  _{i\in I}$ of scalars (i.e., for some
choice of $r_{i}\in R$ for each $i\in I$) with the property that%
\begin{equation}
\text{all but finitely many }i\in I\text{ satisfy }r_{i}=0.
\label{eq.def.mods.lincomb-et-al-inf.a.finit}%
\end{equation}

Here, the sum $\sum_{i\in I}r_{i}m_{i}$ is an infinite sum, but all but
finitely many of its addends are zero (thanks to the condition
(\ref{eq.def.mods.lincomb-et-al-inf.a.finit})). Such a sum is simply defined
to be the sum of the nonzero addends. For example, $3+2+0+0+0+\cdots=3+2=5$.

\item[\textbf{(b)}] The set of all linear combinations of $\left(
m_{i}\right)  _{i\in I}$ is called the \textbf{span} of $\left(  m_{i}\right)
_{i\in I}$, and is denoted by $\operatorname*{span}\left(  m_{i}\right)
_{i\in I}$. (Note that \cite{DumFoo04} calls it $R\left\{  m_{i}\ \mid\ i\in
I\right\}  $.)

\item[\textbf{(c)}] If the span of $\left(  m_{i}\right)  _{i\in I}$ is $M$,
then we say that the family $\left(  m_{i}\right)  _{i\in I}$ \textbf{spans}
$M$ (or \textbf{generates} $M$).

\item[\textbf{(d)}] We say that the family $\left(  m_{i}\right)  _{i\in I}$
is \textbf{linearly independent} if the following holds: If some family
$\left(  r_{i}\right)  _{i\in I}$ of scalars $r_{i}\in R$ has the properties
that%
\begin{equation}
\text{all but finitely many }i\in I\text{ satisfy }r_{i}=0
\label{eq.def.mods.lincomb-et-al-inf.d.finit}%
\end{equation}
and that%
\[
\sum_{i\in I}r_{i}m_{i}=0,
\]
then $r_{i}=0$ for all $i\in I$.

\item[\textbf{(e)}] We say that the family $\left(  m_{i}\right)  _{i\in I}$
is a \textbf{basis} of the $R$-module $M$ if $\left(  m_{i}\right)  _{i\in I}$
is linearly independent and spans $M$.

\item[\textbf{(f)}] All of this terminology depends on $R$. Thus, if $R$ is
not clear from the context, we will clarify it by saying \textquotedblleft%
$R$-linear combination\textquotedblright\ (or \textquotedblleft linear
combination over $R$\textquotedblright) instead of just \textquotedblleft
linear combination\textquotedblright, etc..
\end{enumerate}
\end{definition}

The infinite sums in this definition are a bit of a distraction, but a
necessary one. Fortunately, when studying these notions, it is often
sufficient to work with finite families (i.e., finite sets $I$), since they
are in some sense representative of the general case. To wit:

\begin{proposition}
\label{prop.mods.lincomb-et-all.finitary}Let $M$ be a left $R$-module. Let
$\left(  m_{i}\right)  _{i\in I}$ be a family of vectors in $M$ (with $I$
being any set).

\begin{enumerate}
\item[\textbf{(a)}] Any linear combination of $\left(  m_{i}\right)  _{i\in
I}$ is already a linear combination of some finite subfamily of $\left(
m_{i}\right)  _{i\in I}$. (That is: If $m$ is a linear combination of $\left(
m_{i}\right)  _{i\in I}$, then there exists some finite subset $J$ of $I$ such
that $m$ is a linear combination of $\left(  m_{i}\right)  _{i\in J}$.)

\item[\textbf{(b)}] The family $\left(  m_{i}\right)  _{i\in I}$ is linearly
independent if and only if all its finite subfamilies (i.e., all families of
the form $\left(  m_{i}\right)  _{i\in J}$ with $J$ being a finite subset of
$I$) are linearly independent.
\end{enumerate}
\end{proposition}

\begin{proof}
\textbf{(a)} Let $m$ be a linear combination of $\left(  m_{i}\right)  _{i\in
I}$. Thus, $m$ has the form%
\[
m=\sum_{i\in I}r_{i}m_{i}%
\]
for some family $\left(  r_{i}\right)  _{i\in I}$ of scalars (i.e., for some
choice of $r_{i}\in R$ for each $i\in I$) with the property that%
\[
\text{all but finitely many }i\in I\text{ satisfy }r_{i}=0.
\]
The latter property can be rewritten as follows: There exists a
\textbf{finite} subset $J$ of $I$ such that all $i\in I\setminus J$ satisfy
$r_{i}=0$. Consider this $J$. Then, in the sum $\sum_{i\in I}r_{i}m_{i}$, all
the addends with $i\notin J$ are $0$ (since these addends satisfy $i\notin J$,
thus $i\in I\setminus J$, hence $r_{i}=0$ and therefore $r_{i}m_{i}=0m_{i}%
=0$). Hence, we can throw these addends away and are left with the finite sum
$\sum_{i\in J}r_{i}m_{i}$. Therefore, $\sum_{i\in I}r_{i}m_{i}=\sum_{i\in
J}r_{i}m_{i}$, so that $m=\sum_{i\in I}r_{i}m_{i}=\sum_{i\in J}r_{i}m_{i}$.
This shows that $m$ is a linear combination of the finite subfamily $\left(
m_{i}\right)  _{i\in J}$ of our original family $\left(  m_{i}\right)  _{i\in
I}$. This proves Proposition \ref{prop.mods.lincomb-et-all.finitary}
\textbf{(a)}. \medskip

\textbf{(b)} This is similar to part \textbf{(a)}. The details are left to the
reader. (Again, the key is that the condition
(\ref{eq.def.mods.lincomb-et-al-inf.d.finit}) allows us to restrict ourselves
to a finite subset of $I$.)
\end{proof}

\subsubsection{Spans are submodules}

Next, we show that the span of a family of vectors is always a submodule:

\begin{proposition}
\label{prop.mods.lincomb-et-all.span-submod}Let $M$ be a left $R$-module. Let
$\left(  m_{i}\right)  _{i\in I}$ be a family of vectors in $M$. Then, the
span of this family is an $R$-submodule of $M$.
\end{proposition}

\begin{proof}
You have to show the following three statements:

\begin{enumerate}
\item The sum of two linear combinations of $\left(  m_{i}\right)  _{i\in I}$
is a linear combination of $\left(  m_{i}\right)  _{i\in I}$.

\item Scaling a linear combination of $\left(  m_{i}\right)  _{i\in I}$ by an
$r\in R$ gives a linear combination of $\left(  m_{i}\right)  _{i\in I}$.

\item The zero vector is a linear combination of $\left(  m_{i}\right)  _{i\in
I}$.
\end{enumerate}

All three of these statements are easy. For example, let me show the first
statement: Let $v$ and $w$ be two linear combinations of $\left(
m_{i}\right)  _{i\in I}$. Thus, we can write $v$ and $w$ as%
\begin{equation}
v=\sum_{i\in I}a_{i}m_{i}\ \ \ \ \ \ \ \ \ \ \text{and}%
\ \ \ \ \ \ \ \ \ \ w=\sum_{i\in I}b_{i}m_{i}
\label{pf.prop.mods.lincomb-et-all.span-submod.0}%
\end{equation}
for some two families $\left(  a_{i}\right)  _{i\in I}$ and $\left(
b_{i}\right)  _{i\in I}$ of scalars (i.e., for some choices of $a_{i}\in R$
and $b_{i}\in R$ for each $i\in I$) with the property that%
\begin{equation}
\text{all but finitely many }i\in I\text{ satisfy }a_{i}=0
\label{pf.prop.mods.lincomb-et-all.span-submod.1}%
\end{equation}
and that%
\begin{equation}
\text{all but finitely many }i\in I\text{ satisfy }b_{i}=0.
\label{pf.prop.mods.lincomb-et-all.span-submod.2}%
\end{equation}

Now, adding the two equalities in
(\ref{pf.prop.mods.lincomb-et-all.span-submod.0}) together, we obtain%
\begin{align}
v+w  &  =\sum_{i\in I}a_{i}m_{i}+\sum_{i\in I}b_{i}m_{i}=\sum_{i\in I}\left(
a_{i}m_{i}+b_{i}m_{i}\right) \nonumber\\
&  =\sum_{i\in I}\left(  a_{i}+b_{i}\right)  m_{i}.
\label{pf.prop.mods.lincomb-et-all.span-submod.5}%
\end{align}
Moreover, combining (\ref{pf.prop.mods.lincomb-et-all.span-submod.1}) with
(\ref{pf.prop.mods.lincomb-et-all.span-submod.2}), we see that all but
finitely many $i\in I$ satisfy $a_{i}=0$ and $b_{i}=0$ at the same time (since
the union of two finite sets is still a finite set). Therefore, all but
finitely many $i\in I$ satisfy $a_{i}+b_{i}=0$ (because if $a_{i}=0$ and
$b_{i}=0$, then $a_{i}+b_{i}=0+0=0$). Hence,
(\ref{pf.prop.mods.lincomb-et-all.span-submod.5}) shows that $v+w$ is a linear
combination of $\left(  m_{i}\right)  _{i\in I}$. This proves Statement 1
above. The proofs of Statements 2 and 3 are even easier.
\end{proof}

\subsubsection{Free modules}

\begin{definition}
\label{def.mods.free}\ \ 

\begin{enumerate}
\item[\textbf{(a)}] A left $R$-module is said to be \textbf{free} if it has a basis.

\item[\textbf{(b)}] Let $n\in\mathbb{N}$. A left $R$-module is said to be
\textbf{free of rank }$n$ if it has a basis of size $n$ (i.e., a basis
consisting of $n$ vectors).
\end{enumerate}
\end{definition}

Note that a free $R$-module does not necessarily have a rank, since its basis
could be infinite.\footnote{For some rings $R$, there also exist $R$-modules
that are free of several ranks at the same time -- e.g., an $R$-module can be
free of rank $1$ and free of rank $2$ simultaneously. The simplest such
example is when $R$ is the trivial ring (in which case any $R$-module is
trivial and free of any rank). More interesting examples exist for certain
noncommutative rings -- see, e.g.,
\url{https://math.stackexchange.com/questions/72723/} .}

Let us see some examples of modules that are free and modules that aren't.

You might want to look at $\mathbb{Q}$-modules at first; but they make for
boring examples, because of the following fact:

\begin{theorem}
\label{thm.vs.free}If $F$ is a field, then every $F$-module (= $F$-vector
space) is free.
\end{theorem}

\begin{proof}
This is just the famous fact from linear algebra that every vector space has a
basis. In the most important case (which is when the vector space admits a
finite spanning set -- i.e., there is a finite list $\left(  m_{1}%
,m_{2},\ldots,m_{n}\right)  $ of vectors that spans it\footnote{Such vector
spaces are called \textbf{finite-dimensional}.}), this has fairly neat
elementary proofs (see, e.g., Theorem 2.1 in Keith Conrad's
\url{https://kconrad.math.uconn.edu/blurbs/linmultialg/dimension.pdf} , or
\cite[Theorem 5.3.4]{LaNaSc16} or \cite[Chapter 1, Proposition 2.8]{Treil21}).
In the general case, the proof is tricky and requires the Axiom of Choice (see
Theorem 4.1 in Keith Conrad's
\url{https://kconrad.math.uconn.edu/blurbs/zorn1.pdf}, or \cite[Corollary
212]{Siksek21} or \cite[Theorem 5.23 \textbf{(b)}]{Philip23}).
\end{proof}

For example, Theorem \ref{thm.vs.free} shows that the $\mathbb{Q}$-vector
space $\mathbb{R}$ is free, i.e., has a basis. Such bases are called
\textbf{Hamel bases} and theoretically exist (if you believe in the Axiom of
Choice). Practically, there is no way to construct one.

To find more interesting examples, we need to consider rings that are not
fields. First of all, let us discuss a family of examples that exists for an
arbitrary ring $R$:

\begin{itemize}
\item Consider the left $R$-module%
\[
R^{2}=\left\{  \left(  a,b\right)  \ \mid\ a\in R\text{ and }b\in R\right\}
.
\]

This $R$-module $R^{2}$ is free of rank $2$, since the list $\left(  \left(
1,0\right)  ,\ \ \left(  0,1\right)  \right)  $ is a basis of it. Indeed:

\begin{itemize}
\item The vectors $\left(  1,0\right)  ,\ \left(  0,1\right)  $ span $R^{2}$
(because any vector $\left(  a,b\right)  $ can be written as $a\left(
1,0\right)  +b\left(  0,1\right)  $, and thus is a linear combination of
$\left(  1,0\right)  ,\ \left(  0,1\right)  $).

\item The vectors $\left(  1,0\right)  ,\ \left(  0,1\right)  $ are linearly
independent, since $a\left(  1,0\right)  +b\left(  0,1\right)  =\left(
a,b\right)  $ can only be $0$ if $a=b=0$.
\end{itemize}

\item Likewise, the left $R$-module $R^{3}$ has basis $\left(  \left(
1,0,0\right)  ,\ \left(  0,1,0\right)  ,\ \left(  0,0,1\right)  \right)  $.

\item More generally: If $n\in\mathbb{N}$, then the left $R$-module $R^{n}$
has basis%
\begin{align*}
&  (\left(  1,0,0,\ldots,0\right)  ,\\
&  \ \ \left(  0,1,0,\ldots,0\right)  ,\\
&  \ \ \left(  0,0,1,\ldots,0\right)  ,\\
&  \ \ \ \ \ \ldots,\\
&  \ \ \left(  0,0,0,\ldots,1\right)  ).
\end{align*}
This basis is called the \textbf{standard basis} of $R^{n}$, and its $n$
vectors are called $e_{1},e_{2},\ldots,e_{n}$ (in this order). To make this
more rigorous: For each $i\in\left\{  1,2,\ldots,n\right\}  $, we define
$e_{i}$ to be the vector in $R^{n}$ whose $i$-th entry is $1$ and whose all
remaining entries are $0$ (it is an $n$-tuple, like any vector in $R^{n}$).
Then, the list $\left(  e_{1},e_{2},\ldots,e_{n}\right)  $ is a basis of the
left $R$-module $R^{n}$. Thus, the $R$-module $R^{n}$ is free of rank $n$.

\item As a particular case, the left $R$-module $R^{1}$ is free of rank $1$.
Note that $R^{1}\cong R$, because the map $R\rightarrow R^{1},\ r\mapsto
\left(  r\right)  $ (which merely wraps each scalar into a list to turn it
into a vector) is an $R$-module isomorphism. Hence, the left $R$-module $R$ is
free of rank $1$. Of course, you can see this directly as well: The
one-element list $\left(  1\right)  $ is a basis of it.

Likewise, the left $R$-module $R^{0}$ is free of rank $0$. Note that $R^{0}$
is a trivial $R$-module (it consists of just the zero vector); the empty list
is a basis for it (since the only vector in $R^{0}$ is the zero vector and
thus is a linear combination of nothing). Some authors (e.g., Keith Conrad in
the above-mentioned references) avoid trivial $R$-modules\footnote{A
\textbf{trivial }$R$\textbf{-module} means an $R$-module that consists only of
the zero vector.}, but there is no natural reason to do so except for the
slight weirdness of dealing with empty lists and empty sums.

\item More generally: If $I$ is a set, then\footnote{See Definition
\ref{def.mods.M(I)} for the meaning of the notation $R^{\left(  I\right)  }$
we are using here.}
\[
R^{\left(  I\right)  }=\bigoplus\limits_{i\in I}R=\left\{  \left(
r_{i}\right)  _{i\in I}\in R^{I}\ \mid\ \text{all but finitely many }i\in
I\text{ satisfy }r_{i}=0\right\}
\]
is a free $R$-module. It has a standard basis $\left(  e_{i}\right)  _{i\in
I}$, where each $e_{j}$ is the family that has a $1$ in its $j$-th position
and $0$s in all other positions. (That is, $e_{j}=\left(  \delta_{i,j}\right)
_{i\in I}$, where $\delta_{i,j}=%
\begin{cases}
1, & \text{if }i=j;\\
0, & \text{if }i\neq j
\end{cases}
$\ \ .)

This $R$-module $R^{\left(  I\right)  }$ is an $R$-submodule of
\[
R^{I}=\prod\limits_{i\in I}R=\left\{  \left(  r_{i}\right)  _{i\in I}%
\ \mid\ \text{all }r_{i}\text{ belong to }R\right\}  .
\]
When $I$ is finite, we actually have $R^{\left(  I\right)  }=R^{I}$ (since the
condition \textquotedblleft all but finitely many $i\in I$ satisfy $r_{i}%
=0$\textquotedblright\ is automatically true when $I$ is finite). In general,
however, $R^{\left(  I\right)  }$ is smaller than $R^{I}$, and the $R$-module
$R^{I}=\prod\limits_{i\in I}R$ is usually not free. (For example, the
$\mathbb{Z}$-module $\mathbb{Z}^{\mathbb{N}}$ is not free. This is actually
not easy to prove! A proof is sketched in \cite[\S 10.3, Exercise
24]{DumFoo04}. It is easy to see that the standard basis $\left(
e_{i}\right)  _{i\in\mathbb{N}}$ of $\mathbb{Z}^{\left(  \mathbb{N}\right)  }$
is not a basis of $\mathbb{Z}^{\mathbb{N}}$, since (e.g.) the vector $\left(
1,1,1,1,\ldots\right)  $ is not a linear combination of this
family\footnote{Of course, you could write%
\[
\left(  1,1,1,1,\ldots\right)  =1e_{0}+1e_{1}+1e_{2}+1e_{3}+1e_{4}+\cdots;
\]
however, the sum on the right is properly infinite (with infinitely many
nonzero coefficients) and thus does not count as a linear combination (as it
fails the condition (\ref{eq.def.mods.lincomb-et-al-inf.a.finit}) from
Definition \ref{def.mods.lincomb-et-al-inf}).}. But it is much harder to show
that there is no basis at all.)

\item Let $R$ be a ring. Let $n,m\in\mathbb{N}$. Then, the left $R$-module
$R^{n\times m}$ of all $n\times m$-matrices is free. It has a basis $\left(
E_{i,j}\right)  _{\left(  i,j\right)  \in\left\{  1,2,\ldots,n\right\}
\times\left\{  1,2,\ldots,m\right\}  }$, which consists of the so-called
\textbf{elementary matrices} $E_{i,j}$. For each $i\in\left\{  1,2,\ldots
,n\right\}  $ and $j\in\left\{  1,2,\ldots,m\right\}  $, the respective
elementary matrix $E_{i,j}$ is defined to be the $n\times m$-matrix whose
$\left(  i,j\right)  $-th entry is $1$ while all its other entries are $0$.

For example, for $n=2$ and $m=3$, this basis consists of the six elementary
matrices%
\begin{align*}
E_{1,1}  &  =\left(
\begin{array}
[c]{ccc}%
1 & 0 & 0\\
0 & 0 & 0
\end{array}
\right)  ,\ \ \ \ \ \ \ \ \ \ E_{1,2}=\left(
\begin{array}
[c]{ccc}%
0 & 1 & 0\\
0 & 0 & 0
\end{array}
\right)  ,\ \ \ \ \ \ \ \ \ \ E_{1,3}=\left(
\begin{array}
[c]{ccc}%
0 & 0 & 1\\
0 & 0 & 0
\end{array}
\right)  ,\\
E_{2,1}  &  =\left(
\begin{array}
[c]{ccc}%
0 & 0 & 0\\
1 & 0 & 0
\end{array}
\right)  ,\ \ \ \ \ \ \ \ \ \ E_{2,2}=\left(
\begin{array}
[c]{ccc}%
0 & 0 & 0\\
0 & 1 & 0
\end{array}
\right)  ,\ \ \ \ \ \ \ \ \ \ E_{2,3}=\left(
\begin{array}
[c]{ccc}%
0 & 0 & 0\\
0 & 0 & 1
\end{array}
\right)  .
\end{align*}

There are, of course, many other bases of $R^{n\times m}$ too.

\item Let $R$ be a ring. Let $n\in\mathbb{N}$. The set of all symmetric
$n\times n$-matrices forms a left $R$-submodule $R_{\operatorname*{symm}%
}^{n\times n}$ of the left $R$-module $R^{n\times n}$. It, too, is free. For
example, for $n=2$, it has a basis consisting of the three matrices%
\[
E_{1,1}=\left(
\begin{array}
[c]{cc}%
1 & 0\\
0 & 0
\end{array}
\right)  ,\ \ \ \ \ \ \ \ \ \ E_{1,2}+E_{2,1}=\left(
\begin{array}
[c]{cc}%
0 & 1\\
1 & 0
\end{array}
\right)  ,\ \ \ \ \ \ \ \ \ \ E_{2,2}=\left(
\begin{array}
[c]{cc}%
0 & 0\\
0 & 1
\end{array}
\right)  .
\]

\end{itemize}

Let us now look at $\mathbb{Z}$-modules. Recall that $\mathbb{Z}$-modules are
the same as abelian groups (see Proposition \ref{prop.Zmods.1}), so free
$\mathbb{Z}$-modules are also known as \textbf{free abelian groups} (this is
not the same as free groups).

\begin{itemize}
\item Consider the $\mathbb{Z}$-submodule
\[
U:=\left\{  \left(  a,b,c\right)  \in\mathbb{Z}^{3}\ \mid\ a+b+c=0\right\}
\text{ of }\mathbb{Z}^{3}.
\]
Is $U$ free? Can we find a basis for $U$ ?

So we are trying to find a basis for a submodule of $\mathbb{Z}^{3}$ that is
determined by a set of linear equations (in our case, only one linear equation
-- namely, $a+b+c=0$). If we were using a field (e.g., $\mathbb{Q}$ or
$\mathbb{R}$) instead of $\mathbb{Z}$, then this would be an instance of a
classical problem from linear algebra (solving a system of homogeneous linear
equations\footnote{more precisely: finding a basis for the solution space of
such a system}), which can be solved by Gaussian elimination (see, e.g.,
\cite[\S A.3.2]{LaNaSc16}). If we try to perform Gaussian elimination over
$\mathbb{Z}$, we might run into trouble: Denominators may appear; as a result,
we might not actually get vectors with integer entries. However, for the
submodule $U$ above, this does not happen, and we obtain the basis%
\[
\left(  \left(  -1,1,0\right)  ,\ \left(  -1,0,1\right)  \right)  .
\]
So $U$ is indeed free.\footnote{See Exercise \ref{exe.mods.sum=0} for a
generalization of this.}

What if we have a more complicated submodule and we do run into denominators?
Thus, we do not get a basis using Gaussian elimination. Does this mean that no
basis exists, or does it mean that we have to try something else? We will soon see.

\item The $\mathbb{Z}$-module $\mathbb{Z}/2$ is not free (i.e., does not have
a basis). Indeed, if it had a basis, then this basis would contain at least
one vector (since $\mathbb{Z}/2$ is not trivial), but this vector would not be
linearly independent, since scaling it by $2$ would give $0$.

\item More generally, if $M$ is any \textbf{finite} abelian group of size
larger than $1$, then $M$ is not free (as a $\mathbb{Z}$-module), since a free
$\mathbb{Z}$-module must be either trivial or infinite.

\item The $\mathbb{Z}$-module $\mathbb{Q}$ is not free (i.e., does not have a basis).

\textit{Proof.} Assume the contrary. Thus, there exists a $\mathbb{Z}$-basis
$\left(  m_{i}\right)  _{i\in I}$ of $\mathbb{Q}$. The set $I$ must be
nonempty (since $\mathbb{Q}$ is not trivial); thus, we are in one of the
following two cases:

\begin{itemize}
\item \textit{Case 1:} We have $\left\vert I\right\vert =1$. In this case, $I$
is a $1$-element set, so we can rewrite our basis $\left(  m_{i}\right)
_{i\in I}$ as a list $\left(  m\right)  $ that consists of a single rational
number $m$. This single rational number $m$ must span the entire $\mathbb{Z}%
$-module $\mathbb{Q}$. In other words, every element of $\mathbb{Q}$ must be a
$\mathbb{Z}$-multiple of $m$. But this is absurd (indeed, if $m=0$, then $1$
is not a $\mathbb{Z}$-multiple of $m$; but otherwise, $\dfrac{1}{2}m$ is not a
$\mathbb{Z}$-multiple of $m$).

\item \textit{Case 2:} We have $\left\vert I\right\vert >1$. In this case,
there are at least two vectors $m_{u}$ and $m_{v}$ in this basis $\left(
m_{i}\right)  _{i\in I}$. However, two rational numbers are never $\mathbb{Z}%
$-linearly independent\footnote{Indeed, let $p$ and $q$ be two rational
numbers. We claim that there exist integers $a,b\in\mathbb{Z}$ that are not
both $0$ but still satisfy $ap+bq=0$. (This will clearly prove that $p$ and
$q$ are not $\mathbb{Z}$-linearly independent.)
\par
Indeed, if $p=0$, then we set $a=1$ and $b=0$ and are done. Something similar
works if $q=0$. So we WLOG assume that $p\neq0$ and $q\neq0$. Write $p$ and
$q$ as $p=\dfrac{n}{d}$ and $q=\dfrac{m}{e}$ for some nonzero integers
$n,d,m,e$ (we can do this, since $p$ and $q$ are nonzero rational numbers).
Then, $dmp+\left(  -en\right)  q=0$ (check this!), so we have found our $a$
and $b$ (namely, $a=dm$ and $b=-en$).}. Thus, a fortiori, the whole family
$\left(  m_{i}\right)  _{i\in I}$ cannot be $\mathbb{Z}$-linearly independent
(since a subfamily of a linearly independent family of vectors must always be
linearly independent). This contradicts the assumption that this family is a basis.
\end{itemize}

Thus, in each case, we have found a contradiction, and our proof is complete.

\item Now, consider the $\mathbb{Z}$-submodule%
\[
V:=\left\{  \left(  a,b\right)  \in\mathbb{Z}^{2}\ \mid\ a\equiv
b\operatorname{mod}2\right\}  \text{ of }\mathbb{Z}^{2}.
\]
This $\mathbb{Z}$-submodule $V$ contains the vectors $\left(  0,2\right)  $
and $\left(  1,1\right)  $ and $\left(  1,-1\right)  $ and $\left(
4,-2\right)  $ and many others. Is $V$ free? Can we find a basis for $V$ ?

Let's try the pair $\left(  \left(  2,0\right)  ,\ \ \left(  0,2\right)
\right)  $. Is this pair a basis for $V$ ? Its span is%
\begin{align*}
\operatorname*{span}\left(  \left(  2,0\right)  ,\ \ \left(  0,2\right)
\right)   &  =\left\{  c\left(  2,0\right)  +d\left(  0,2\right)
\ \mid\ c,d\in\mathbb{Z}\right\} \\
&  =\left\{  \left(  2c,2d\right)  \ \mid\ c,d\in\mathbb{Z}\right\} \\
&  =\left\{  \left(  a,b\right)  \in\mathbb{Z}^{2}\ \mid\ a\equiv
b\equiv0\operatorname{mod}2\right\}  .
\end{align*}
This is a $\mathbb{Z}$-submodule of $V$, but not the entire $V$, since (for
example) $\left(  1,1\right)  $ belongs to $V$ but not to
$\operatorname*{span}\left(  \left(  2,0\right)  ,\ \ \left(  0,2\right)
\right)  $. So we have \textquotedblleft undershot\textquotedblright\ our $V$
(by finding a linearly independent family that does not span $V$).

Let's try the triple $\left(  \left(  2,0\right)  ,\ \ \left(  0,2\right)
,\ \ \left(  1,1\right)  \right)  $. This triple does span $V$ (check this!),
but is not linearly independent, since%
\[
1\cdot\left(  2,0\right)  +1\cdot\left(  0,2\right)  +\left(  -2\right)
\cdot\left(  1,1\right)  =0.
\]
So we have \textquotedblleft overshot\textquotedblright\ $V$ now (by finding a
family that spans $V$ but is not linearly independent).

Let us try to correct this by throwing away $\left(  0,2\right)  $. So we are
left with the pair $\left(  \left(  2,0\right)  ,\ \ \left(  1,1\right)
\right)  $. And this pair is indeed a basis of $V$, as can easily be checked.
Indeed, it is linearly independent (you can check this using linear algebra,
since it clearly suffices to prove its $\mathbb{Q}$-linear
independence\footnote{Or you can check this directly: If $a,b\in\mathbb{Z}$
satisfy $a\left(  2,0\right)  +b\left(  1,1\right)  =0$, then $0=a\left(
2,0\right)  +b\left(  1,1\right)  =\left(  2a+b,\ b\right)  $, so that
$\left(  2a+b,\ b\right)  =0=\left(  0,0\right)  $, and therefore $2a+b=0$ and
$b=0$; but this quickly yields $a=b=0$.}), and furthermore spans $V$ because
each $\left(  a,b\right)  \in V$ can be written as a linear combination of
$\left(  2,0\right)  ,\ \ \left(  1,1\right)  $ as follows:%
\[
\left(  a,b\right)  =\underbrace{\dfrac{a-b}{2}}_{\substack{\in\mathbb{Z}%
\\\text{(since }a\equiv b\operatorname{mod}2\text{)}}}\cdot\left(  2,0\right)
+b\cdot\left(  1,1\right)  .
\]

Another basis for $V$ is the pair $\left(  \left(  1,1\right)  ,\ \ \left(
1,-1\right)  \right)  $. Indeed, this pair is linearly independent (check
this!), and it spans $V$, because each $\left(  a,b\right)  \in V$ can be
written as
\[
\left(  a,b\right)  =\underbrace{\dfrac{a+b}{2}}_{\in\mathbb{Z}}\cdot\left(
1,1\right)  +\underbrace{\dfrac{a-b}{2}}_{\in\mathbb{Z}}\cdot\left(
1,-1\right)  \in\operatorname*{span}\left(  \left(  1,1\right)  ,\ \ \left(
1,-1\right)  \right)  .
\]

(See Exercise \ref{exe.mods.all-cong-mod-k} for a generalization of $V$.)
\end{itemize}

\begin{exercise}
\label{exe.mods.sum=0}Let $R$ be any ring. Let $n$ be a positive integer. Let
$U$ be the subset%
\[
\left\{  \left(  a_{1},a_{2},\ldots,a_{n}\right)  \in R^{n}\ \mid\ a_{1}%
+a_{2}+\cdots+a_{n}=0\right\}
\]
of the left $R$-module $R^{n}$.

\begin{enumerate}
\item[\textbf{(a)}] Prove that $U$ is an $R$-submodule of $R^{n}$.

\item[\textbf{(b)}] Show that $U$ is free, and prove that the $n-1$ vectors%
\begin{align*}
&  \left(  1,-1,0,0,\ldots,0\right)  ,\\
&  \left(  0,1,-1,0,\ldots,0\right)  ,\\
&  \ldots,\\
&  \left(  0,0,\ldots,0,1,-1\right)
\end{align*}
(that is, the $n-1$ vectors that consist of a number of $0$'s, followed by a
$1$, followed by a $-1$, followed again by a number of $0$'s) form a basis of
$U$.
\end{enumerate}
\end{exercise}

\begin{exercise}
\label{exe.mods.all-cong-mod-k}Let $n$ and $k$ be two positive integers. Let
$V$ be the subset
\[
\left\{  \left(  a_{1},a_{2},\ldots,a_{n}\right)  \in\mathbb{Z}^{n}%
\ \mid\ a_{1}\equiv a_{2}\equiv\cdots\equiv a_{n}\operatorname{mod}k\right\}
\]
of the $\mathbb{Z}$-module $\mathbb{Z}^{n}$.

\begin{enumerate}
\item[\textbf{(a)}] Prove that $V$ is a $\mathbb{Z}$-submodule of
$\mathbb{Z}^{n}$.

\item[\textbf{(b)}] Show that $V$ is free, and find a basis of $V$.
\end{enumerate}
\end{exercise}

\begin{exercise}
Let $R$ be any ring. Let $n\geq2$ be an integer. Let $W$ be the subset%
\[
\left\{  \left(  a_{1},a_{2},\ldots,a_{n}\right)  \in R^{n}\ \mid
\ a_{i}-a_{i-1}=a_{i+1}-a_{i}\text{ for each }i\in\left\{  2,3,\ldots
,n-1\right\}  \right\}
\]
of the left $R$-module $R^{n}$. (This set $W$ consists of all vectors $\left(
a_{1},a_{2},\ldots,a_{n}\right)  \in R^{n}$ whose entries \textquotedblleft
form an arithmetic progression\textquotedblright, i.e., satisfy $a_{2}%
-a_{1}=a_{3}-a_{2}=a_{4}-a_{3}=\cdots=a_{n}-a_{n-1}$.)

\begin{enumerate}
\item[\textbf{(a)}] Prove that $W$ is an $R$-submodule of $R^{n}$.

\item[\textbf{(b)}] Let $a=\left(  1,1,\ldots,1\right)  \in R^{n}$ and
$b=\left(  1,2,\ldots,n\right)  \in R^{n}$. Prove that $\left(  a,b\right)  $
is a basis of $W$.
\end{enumerate}
\end{exercise}

\begin{exercise}
\ \ 

\begin{enumerate}
\item[\textbf{(a)}] Let $X$ be the $\mathbb{Z}$-submodule
\[
\left\{  \left(  a,b,c\right)  \in\mathbb{Z}^{3}\ \mid\ a\equiv
b\operatorname{mod}2\text{ and }b\equiv c\operatorname{mod}3\right\}
\]
of $\mathbb{Z}^{3}$. Prove that $X$ is free, and find a basis of $X$.

\item[\textbf{(b)}] Let $Y$ be the $\mathbb{Z}$-submodule
\[
\left\{  \left(  a,b,c\right)  \in\mathbb{Z}^{3}\ \mid\ a\equiv
b\operatorname{mod}2\text{ and }b\equiv c\operatorname{mod}3\text{ and
}c\equiv a\operatorname{mod}5\right\}
\]
of $\mathbb{Z}^{3}$. Prove that $Y$ is free, and find a basis of $Y$.
\end{enumerate}
\end{exercise}

\begin{exercise}
\label{exe.mods.asym-matrices}Let $R$ be any ring. A square matrix $A\in
R^{n\times n}$ is said to be \textbf{antisymmetric} (or
\textbf{skew-symmetric}) if $A^{T}=-A$ (where $A^{T}$ is the transpose of
$A$). For each $n\in\mathbb{N}$, we let $R_{\operatorname*{asym}}^{n\times n}$
denote the set of all antisymmetric $n\times n$-matrices in $R^{n\times n}$.

\begin{enumerate}
\item[\textbf{(a)}] Prove that $R_{\operatorname*{asym}}^{n\times n}$ is a
left $R$-submodule of $R^{n\times n}$ for each $n\in\mathbb{N}$.

\item[\textbf{(b)}] Find a basis of the $\mathbb{Q}$-module $\mathbb{Q}%
_{\operatorname*{asym}}^{2\times2}$.

\item[\textbf{(c)}] Find a basis of the $\mathbb{Z}/2$-module $\left(
\mathbb{Z}/2\right)  _{\operatorname*{asym}}^{2\times2}$.

\item[\textbf{(d)}] Prove that the $\mathbb{Z}/4$-module $\left(
\mathbb{Z}/4\right)  _{\operatorname*{asym}}^{2\times2}$ is not free (i.e.,
has no basis).

(This is somewhat surprising when compared to the $R$-module
$R_{\operatorname*{symm}}^{n\times n}$ of symmetric $n\times n$-matrices,
which module always has a basis. In a sense, it shows that antisymmetric
matrices behave more wildly than symmetric matrices.)
\end{enumerate}

[\textbf{Hint:} What exactly does the condition $A^{T}=-A$ say about the
diagonal entries of a matrix $A$ ?]
\end{exercise}

\begin{exercise}
Let $R$ be any ring. A square matrix $A\in R^{n\times n}$ is said to be
\textbf{alternating} if it is antisymmetric (i.e., satisfies $A^{T}=-A$) and
its diagonal entries all equal $0$. For each $n\in\mathbb{N}$, we let
$R_{\operatorname*{alt}}^{n\times n}$ denote the set of all alternating
$n\times n$-matrices in $R^{n\times n}$.

\begin{enumerate}
\item[\textbf{(a)}] Prove that $R_{\operatorname*{alt}}^{n\times n}$ is a left
$R$-submodule of $R^{n\times n}$ for each $n\in\mathbb{N}$.

\item[\textbf{(b)}] Prove that this $R$-module $R_{\operatorname*{alt}%
}^{n\times n}$ is always free, and find a basis of this $R$-module. (This
shows that $R_{\operatorname*{alt}}^{n\times n}$ is a \textquotedblleft
better-behaved\textquotedblright\ variant of the $R$-module
$R_{\operatorname*{asym}}^{n\times n}$ from Exercise
\ref{exe.mods.asym-matrices}.)

\item[\textbf{(c)}] What (fairly simple) condition must $R$ satisfy in order
for this $R$-module $R_{\operatorname*{alt}}^{n\times n}$ to be identical to
the $R$-module $R_{\operatorname*{asym}}^{n\times n}$ from Exercise
\ref{exe.mods.asym-matrices}?
\end{enumerate}
\end{exercise}

\begin{remark}
Let $n\in\mathbb{N}$. A nontrivial theorem says that every $\mathbb{Z}%
$-submodule of $\mathbb{Z}^{n}$ is free, i.e., has a basis, and is isomorphic
to $\mathbb{Z}^{m}$ for some $m\in\left\{  0,1,\ldots,n\right\}  $. For a
proof, see (e.g.) \cite[Theorem 8.25]{Knapp1} (which proves more and in
greater generality). Note that the proof is non-constructive, and there is no
general method for finding a basis (or even the rank) of a given $\mathbb{Z}%
$-submodule of $\mathbb{Z}^{n}$.

More generally, if $R$ is a PID, then any $R$-submodule of $R^{n}$ is free and
isomorphic to $R^{m}$ for some $m\in\left\{  0,1,\ldots,n\right\}  $. This
does not generalize to arbitrary rings, though. For example, if $R$ is the
polynomial ring $\mathbb{Z}\left[  x\right]  $, then $R^{1}$ has an
$R$-submodule consisting of the polynomials with an even constant term. This
$R$-submodule is not free (check this!).
\end{remark}

Let us now return to the general case to state a few theorems:

\begin{theorem}
\label{thm.mods.free-is-iso-1}Let $M$ be a left $R$-module. Let $n\in
\mathbb{N}$. The left $R$-module $M$ is free of rank $n$ if and only if
$M\cong R^{n}$ (as left $R$-modules).
\end{theorem}

More concretely:

\begin{theorem}
\label{thm.mods.free-is-iso-2}Let $M$ be a left $R$-module. Let $m_{1}%
,m_{2},\ldots,m_{n}$ be $n$ vectors in $M$. Consider the map%
\begin{align*}
f:R^{n}  &  \rightarrow M,\\
\left(  r_{1},r_{2},\ldots,r_{n}\right)   &  \mapsto r_{1}m_{1}+r_{2}%
m_{2}+\cdots+r_{n}m_{n}.
\end{align*}
Then:

\begin{enumerate}
\item[\textbf{(a)}] This map $f$ is always a left $R$-module morphism.

\item[\textbf{(b)}] The map $f$ is injective if and only if $m_{1}%
,m_{2},\ldots,m_{n}$ are linearly independent.

\item[\textbf{(c)}] The map $f$ is surjective if and only if $m_{1}%
,m_{2},\ldots,m_{n}$ span $M$.

\item[\textbf{(d)}] The map $f$ is an isomorphism\footnotemark\ if and only if
$\left(  m_{1},m_{2},\ldots,m_{n}\right)  $ is a basis of $M$.
\end{enumerate}
\end{theorem}

\footnotetext{Of course, \textquotedblleft isomorphism\textquotedblright%
\ means \textquotedblleft left $R$-module isomorphism\textquotedblright%
\ here.} Note that the map $f$ in Theorem \ref{thm.mods.free-is-iso-2} takes
an $n$-tuple $\left(  r_{1},r_{2},\ldots,r_{n}\right)  $ of scalars, and uses
these scalars as coefficients to form a linear combination of $m_{1}%
,m_{2},\ldots,m_{n}$. Thus, the values of $f$ are precisely the linear
combinations of $m_{1},m_{2},\ldots,m_{n}$.

\begin{proof}
[Proof of Theorem \ref{thm.mods.free-is-iso-2}.]This is commonly done in
linear algebra texts (albeit usually under the assumption that $R$ is a field,
but the proof is the same); thus I will be brief. \medskip

\textbf{(a)} We must prove that $f$ respects addition, respects scaling and
respects the zero. I will only show that it respects addition, since the other
two statements are analogous.

So we must prove that $f\left(  a+b\right)  =f\left(  a\right)  +f\left(
b\right)  $ for all $a,b\in R^{n}$. Indeed, let $a,b\in R^{n}$. Write $a$ and
$b$ as%
\[
a=\left(  a_{1},a_{2},\ldots,a_{n}\right)  \ \ \ \ \ \ \ \ \ \ \text{and}%
\ \ \ \ \ \ \ \ \ \ b=\left(  b_{1},b_{2},\ldots,b_{n}\right)  .
\]
Then, the definition of $R^{n}$ (as the direct product $\underbrace{R\times
R\times\cdots\times R}_{n\text{ times}}$) yields $a+b=\left(  a_{1}%
+b_{1},a_{2}+b_{2},\ldots,a_{n}+b_{n}\right)  $. Hence, the definition of $f$
yields%
\begin{align*}
f\left(  a+b\right)   &  =\left(  a_{1}+b_{1}\right)  m_{1}+\left(
a_{2}+b_{2}\right)  m_{2}+\cdots+\left(  a_{n}+b_{n}\right)  m_{n}\\
&  =\left(  a_{1}m_{1}+b_{1}m_{1}\right)  +\left(  a_{2}m_{2}+b_{2}%
m_{2}\right)  +\cdots+\left(  a_{n}m_{n}+b_{n}m_{n}\right) \\
&  \ \ \ \ \ \ \ \ \ \ \ \ \ \ \ \ \ \ \ \ \left(  \text{by right
distributivity}\right) \\
&  =\underbrace{\left(  a_{1}m_{1}+a_{2}m_{2}+\cdots+a_{n}m_{n}\right)
}_{\substack{=f\left(  a\right)  \\\text{(by the definition of }f\text{, since
}a=\left(  a_{1},a_{2},\ldots,a_{n}\right)  \text{)}}}+\underbrace{\left(
b_{1}m_{1}+b_{2}m_{2}+\cdots+b_{n}m_{n}\right)  }_{\substack{=f\left(
b\right)  \\\text{(by the definition of }f\text{, since }b=\left(  b_{1}%
,b_{2},\ldots,b_{n}\right)  \text{)}}}\\
&  =f\left(  a\right)  +f\left(  b\right)  ,
\end{align*}
which is what we wanted to show. \medskip

\textbf{(b)} The map $f$ is an $R$-module morphism (by part \textbf{(a)}).
Thus, it is injective if and only if $\operatorname*{Ker}f=\left\{  0_{R^{n}%
}\right\}  $ (by Lemma \ref{lem.modmor.ker-inj}). Hence, we have the following
chain of logical equivalences:%
\begin{align*}
&  \ \left(  f\text{ is injective}\right) \\
&  \Longleftrightarrow\ \left(  \operatorname*{Ker}f=\left\{  0_{R^{n}%
}\right\}  \right) \\
&  \Longleftrightarrow\ \left(  \operatorname*{Ker}f\subseteq\left\{
0_{R^{n}}\right\}  \right)  \ \ \ \ \ \ \ \ \ \ \left(  \text{since }\left\{
0_{R^{n}}\right\}  \text{ is clearly a subset of }\operatorname*{Ker}f\right)
\\
&  \Longleftrightarrow\ \left(  \left\{  a\in R^{n}\ \mid\ f\left(  a\right)
=0\right\}  \subseteq\left\{  0_{R^{n}}\right\}  \right) \\
&  \ \ \ \ \ \ \ \ \ \ \ \ \ \ \ \ \ \ \ \ \left(  \text{since }%
\operatorname*{Ker}f=\left\{  a\in R^{n}\ \mid\ f\left(  a\right)  =0\right\}
\text{ by the definition of }\operatorname*{Ker}f\right) \\
&  \Longleftrightarrow\ \left(  \text{the only }a\in R^{n}\text{ satisfying
}f\left(  a\right)  =0\text{ is }0_{R^{n}}\right) \\
&  \Longleftrightarrow\ \left(
\begin{array}
[c]{c}%
\text{the only }\left(  a_{1},a_{2},\ldots,a_{n}\right)  \in R^{n}\text{
satisfying }f\left(  a_{1},a_{2},\ldots,a_{n}\right)  =0\\
\text{is }\left(  0,0,\ldots,0\right)
\end{array}
\right) \\
&  \ \ \ \ \ \ \ \ \ \ \ \ \ \ \ \ \ \ \ \ \left(
\begin{array}
[c]{c}%
\text{since any }a\in R^{n}\text{ can be written in the form }\left(
a_{1},a_{2},\ldots,a_{n}\right)  \text{,}\\
\text{and since }0_{R^{n}}=\left(  0,0,\ldots,0\right)
\end{array}
\right) \\
&  \Longleftrightarrow\ \left(
\begin{array}
[c]{c}%
\text{the only }\left(  a_{1},a_{2},\ldots,a_{n}\right)  \in R^{n}\text{
satisfying }a_{1}m_{1}+a_{2}m_{2}+\cdots+a_{n}m_{n}=0\\
\text{is }\left(  0,0,\ldots,0\right)
\end{array}
\right) \\
&  \ \ \ \ \ \ \ \ \ \ \ \ \ \ \ \ \ \ \ \ \left(
\begin{array}
[c]{c}%
\text{since }f\left(  a_{1},a_{2},\ldots,a_{n}\right)  =a_{1}m_{1}+a_{2}%
m_{2}+\cdots+a_{n}m_{n}\\
\text{for any }\left(  a_{1},a_{2},\ldots,a_{n}\right)  \in R^{n}%
\end{array}
\right) \\
&  \Longleftrightarrow\ \left(
\begin{array}
[c]{c}%
\text{if }a_{1},a_{2},\ldots,a_{n}\in R\text{ satisfy }a_{1}m_{1}+a_{2}%
m_{2}+\cdots+a_{n}m_{n}=0\text{,}\\
\text{then }a_{1}=a_{2}=\cdots=a_{n}=0
\end{array}
\right) \\
&  \Longleftrightarrow\ \left(  m_{1},m_{2},\ldots,m_{n}\text{ are linearly
independent}\right)
\end{align*}
(by the definition of linear independence). This proves part \textbf{(b)} of
the theorem. \medskip

\textbf{(c)} We have the following chain of logical equivalences:%
\begin{align*}
&  \ \left(  f\text{ is surjective}\right) \\
&  \Longleftrightarrow\ \left(  \text{each }m\in M\text{ can be written as
}f\left(  a\right)  \text{ for some }a\in R^{n}\right) \\
&  \Longleftrightarrow\ \left(
\begin{array}
[c]{c}%
\text{each }m\in M\text{ can be written as }f\left(  a_{1},a_{2},\ldots
,a_{n}\right) \\
\text{for some }\left(  a_{1},a_{2},\ldots,a_{n}\right)  \in R^{n}%
\end{array}
\right) \\
&  \ \ \ \ \ \ \ \ \ \ \ \ \ \ \ \ \ \ \ \ \left(  \text{since any }a\in
R^{n}\text{ can be written in the form }\left(  a_{1},a_{2},\ldots
,a_{n}\right)  \right) \\
&  \Longleftrightarrow\ \left(
\begin{array}
[c]{c}%
\text{each }m\in M\text{ can be written as }a_{1}m_{1}+a_{2}m_{2}+\cdots
+a_{n}m_{n}\\
\text{for some }\left(  a_{1},a_{2},\ldots,a_{n}\right)  \in R^{n}%
\end{array}
\right) \\
&  \ \ \ \ \ \ \ \ \ \ \ \ \ \ \ \ \ \ \ \ \left(
\begin{array}
[c]{c}%
\text{since }f\left(  a_{1},a_{2},\ldots,a_{n}\right)  =a_{1}m_{1}+a_{2}%
m_{2}+\cdots+a_{n}m_{n}\\
\text{for any }\left(  a_{1},a_{2},\ldots,a_{n}\right)  \in R^{n}%
\end{array}
\right) \\
&  \Longleftrightarrow\ \left(  \text{each }m\in M\text{ is a linear
combination of }m_{1},m_{2},\ldots,m_{n}\right) \\
&  \ \ \ \ \ \ \ \ \ \ \ \ \ \ \ \ \ \ \ \ \left(  \text{by the definition of
a linear combination}\right) \\
&  \Longleftrightarrow\ \left(  m_{1},m_{2},\ldots,m_{n}\text{ span }M\right)
.
\end{align*}
This proves part \textbf{(c)} of the theorem. \medskip

\textbf{(d)} We have the following chain of logical equivalences:%
\begin{align*}
&  \ \left(  f\text{ is an }R\text{-module isomorphism}\right) \\
&  \Longleftrightarrow\ \left(  f\text{ is invertible}\right) \\
&  \ \ \ \ \ \ \ \ \ \ \ \ \ \ \ \ \ \ \ \ \left(
\begin{array}
[c]{c}%
\text{since we know from Proposition \ref{prop.modmor.invertible-iso} that
any}\\
\text{invertible }R\text{-module morphism is an isomorphism}%
\end{array}
\right) \\
&  \Longleftrightarrow\ \left(  f\text{ is bijective}\right) \\
&  \Longleftrightarrow\ \underbrace{\left(  f\text{ is injective}\right)
}_{\substack{\Longleftrightarrow\ \left(  m_{1},m_{2},\ldots,m_{n}\text{ are
linearly independent}\right)  \\\text{(by part \textbf{(b)})}}}\wedge
\underbrace{\left(  f\text{ is surjective}\right)  }%
_{\substack{\Longleftrightarrow\ \left(  m_{1},m_{2},\ldots,m_{n}\text{ span
}M\right)  \\\text{(by part \textbf{(c)})}}}\\
&  \Longleftrightarrow\ \left(  m_{1},m_{2},\ldots,m_{n}\text{ are linearly
independent}\right)  \wedge\left(  m_{1},m_{2},\ldots,m_{n}\text{ span
}M\right) \\
&  \Longleftrightarrow\ \left(  \left(  m_{1},m_{2},\ldots,m_{n}\right)
\text{ is a basis of }M\right)
\end{align*}
(by the definition of a basis). This proves part \textbf{(d)} of the theorem.
\end{proof}

\begin{proof}
[Proof of Theorem \ref{thm.mods.free-is-iso-1}.]$\Longrightarrow:$ Assume that
$M$ is free of rank $n$. That is, $M$ has a basis $\left(  m_{1},m_{2}%
,\ldots,m_{n}\right)  $ of size $n$. Consider this basis. Consider the map
$f:R^{n}\rightarrow M$ defined in Theorem \ref{thm.mods.free-is-iso-2}. Thus,
Theorem \ref{thm.mods.free-is-iso-2} \textbf{(d)} yields that $f$ is an
isomorphism. Hence, $R^{n}\cong M$ as left $R$-modules. In other words,
$M\cong R^{n}$ as left $R$-modules. This proves the \textquotedblleft%
$\Longrightarrow$\textquotedblright\ direction of Theorem
\ref{thm.mods.free-is-iso-1}. \medskip

$\Longleftarrow:$ Assume that $M\cong R^{n}$ as left $R$-modules. But the left
$R$-module $R^{n}$ is free of rank $n$ (as we have seen above). Hence, I claim
that the left $R$-module $M$ is also free of rank $n$, since $M\cong R^{n}$.
Indeed, this follows from the \textquotedblleft isomorphism
principle\textquotedblright\ for modules -- i.e., from the \textquotedblleft
meta-theorem\textquotedblright\ that says that module isomorphisms preserve
all \textquotedblleft intrinsic\textquotedblright\ properties of modules (in
this case, this property is \textquotedblleft being free of rank
$n$\textquotedblright).

Here is a more pedestrian way to get to the same conclusion: We have $M\cong
R^{n}$, thus $R^{n}\cong M$. In other words, there exists a left $R$-module
isomorphism $g:R^{n}\rightarrow M$. Consider this $g$. Now, consider the
standard basis $\left(  e_{1},e_{2},\ldots,e_{n}\right)  $ of the left
$R$-module $R^{n}$. Applying $g$ to each vector in this basis, we obtain a
list $\left(  g\left(  e_{1}\right)  ,g\left(  e_{2}\right)  ,\ldots,g\left(
e_{n}\right)  \right)  $ of vectors in $M$. It is straightforward to see that
this new list is a basis of $M$ (indeed, when we apply $g$ to a linear
combination $a_{1}e_{1}+a_{2}e_{2}+\cdots+a_{n}e_{n}$ of the standard basis
$\left(  e_{1},e_{2},\ldots,e_{n}\right)  $ in $R^{n}$, then we obtain
\begin{align*}
g\left(  a_{1}e_{1}+a_{2}e_{2}+\cdots+a_{n}e_{n}\right)   &  =a_{1}g\left(
e_{1}\right)  +a_{2}g\left(  e_{2}\right)  +\cdots+a_{n}g\left(  e_{n}\right)
\\
&  \ \ \ \ \ \ \ \ \ \ \left(  \text{since }g\text{ is }R\text{-linear}%
\right)  ,
\end{align*}
which is the corresponding linear combination of $\left(  g\left(
e_{1}\right)  ,g\left(  e_{2}\right)  ,\ldots,g\left(  e_{n}\right)  \right)
$; thus, linear independence of $\left(  e_{1},e_{2},\ldots,e_{n}\right)  $
translates into linear independence of $\left(  g\left(  e_{1}\right)
,g\left(  e_{2}\right)  ,\ldots,g\left(  e_{n}\right)  \right)  $ (since $g$
sends only $0$ to $0$), and the same holds for spanning (since $g$ is
bijective)). Hence, $M$ has a basis of size $n$. In other words, $M$ is free
of rank $n$.

Either way, the \textquotedblleft$\Longleftarrow$\textquotedblright\ direction
of Theorem \ref{thm.mods.free-is-iso-1} is now proved.
\end{proof}

Theorem \ref{thm.mods.free-is-iso-2} can be generalized to bases of arbitrary size:

\begin{theorem}
\label{thm.mods.free-is-iso-3}Let $M$ be a left $R$-module. Let $\left(
m_{i}\right)  _{i\in I}$ be any family of vectors in $M$. Consider the
map\footnotemark%
\begin{align*}
f:R^{\left(  I\right)  }  &  \rightarrow M,\\
\left(  r_{i}\right)  _{i\in I}  &  \mapsto\sum_{i\in I}r_{i}m_{i}.
\end{align*}
(This is well-defined, since any $\left(  r_{i}\right)  _{i\in I}\in
R^{\left(  I\right)  }$ automatically satisfies the condition
(\ref{eq.def.mods.lincomb-et-al-inf.a.finit}) because of the definition of
$R^{\left(  I\right)  }$.)

Then:

\begin{enumerate}
\item[\textbf{(a)}] This map $f$ is always a left $R$-module morphism.

\item[\textbf{(b)}] The map $f$ is injective if and only if the family
$\left(  m_{i}\right)  _{i\in I}$ is linearly independent.

\item[\textbf{(c)}] The map $f$ is surjective if and only if the family
$\left(  m_{i}\right)  _{i\in I}$ spans $M$.

\item[\textbf{(d)}] The map $f$ is an isomorphism if and only if the family
$\left(  m_{i}\right)  _{i\in I}$ is a basis of $M$.
\end{enumerate}
\end{theorem}

\footnotetext{Recall that $R^{\left(  I\right)  }$ denotes the direct sum
$\bigoplus_{i\in I}R$ here. This is the left $R$-module that consists of all
families $\left(  r_{i}\right)  _{i\in I}\in R^{I}$ such that only finitely
many $i\in I$ satisfy $r_{i}\neq0$.}Note that the map $f$ here has domain
$R^{\left(  I\right)  }$, not $R^{I}$, since the infinite sum $\sum_{i\in
I}r_{i}m_{i}$ is well-defined for all $\left(  r_{i}\right)  _{i\in I}\in
R^{\left(  I\right)  }$ but not (in general) for all $\left(  r_{i}\right)
_{i\in I}\in R^{I}$.

The map $f$ in Theorem \ref{thm.mods.free-is-iso-3} takes a family $\left(
r_{i}\right)  _{i\in I}$ of scalars, and uses it to build a linear combination
of $\left(  m_{i}\right)  _{i\in I}$.

\begin{proof}
[Proof of Theorem \ref{thm.mods.free-is-iso-3}.]Analogous to Theorem
\ref{thm.mods.free-is-iso-2}, with the usual caveats about infinite sums.
\end{proof}

\begin{remark}
As you will have noticed by now, \textquotedblleft free module of rank
$n$\textquotedblright\ is a generalization of \textquotedblleft vector space
of dimension $n$\textquotedblright\ to arbitrary rings.

We have been careful to speak of \textquotedblleft free modules of rank
$n$\textquotedblright, but never of \textquotedblleft the rank of a free
module\textquotedblright. This is due to the somewhat perverse-sounding fact
that there can be modules that are free of several ranks simultaneously (i.e.,
modules that have bases of different sizes). One way to get such modules is by
taking $R$ to be a trivial ring (in which case, any $R$-module is trivial and
is free of every rank simultaneously -- seriously). If this was the only
example, one could discount the issue as a formality, but there are less
trivial (pardon) examples as well: \cite[\S 10.3, exercise 27]{DumFoo04}
constructs a ring $R$ over which $R^{n}\cong R$ as left $R$-modules for each
$n\in\left\{  1,2,3,\ldots\right\}  $ (so $R$ itself is a free $R$-module of
rank $n$ for each $n\in\left\{  1,2,3,\ldots\right\}  $).

If $R$ is a nontrivial commutative ring, then things are nice: The $R$-modules
$R^{0},R^{1},R^{2},\ldots$ are mutually non-isomorphic, so a free $R$-module
can never have two different ranks at the same time. This is not obvious (see
\cite[\S 10.3, exercise 2]{DumFoo04}). We can actually say more: If $R$ is a
nontrivial commutative ring, then an $R$-module morphism $R^{m}\rightarrow
R^{n}$ cannot be injective unless $m\leq n$ (see, e.g.,
\url{https://math.stackexchange.com/questions/106786} or \cite[Theorem
2]{Richma88}), and cannot be surjective unless $m\geq n$ (see, e.g.,
\url{https://math.stackexchange.com/questions/20178} or \cite[Theorem
1]{Richma88}). These facts are in line with the intuition you should have from
linear algebra (injective maps cannot squash dimensions; surjective maps
cannot create dimensions) and also with the Pigeonhole Principles from
combinatorics (a map between two finite sets $M$ and $N$ cannot be injective
unless $\left\vert M\right\vert \leq\left\vert N\right\vert $, and cannot be
surjective unless $\left\vert M\right\vert \geq\left\vert N\right\vert $). But
actually proving them takes real work!
\end{remark}

\subsection{\label{sec.modules.uniprop-free}The universal property of a free
module (\cite[\S 10.3]{DumFoo04})}

As before, we fix a ring $R$.

The next proposition shows that linear maps respect linear combinations (in
the sense that if you apply a linear map to a linear combination of some
vectors, then you get the same linear combination of their images):

\begin{proposition}
\label{prop.mods.lin-map-pres-lc}Let $M$ and $P$ be two left $R$-modules. Let
$f:M\rightarrow P$ be an $R$-linear map. Let $\left(  m_{i}\right)  _{i\in I}$
be any family of vectors in $M$, and let $\left(  r_{i}\right)  _{i\in I}$ be
a family of scalars in $R$ with the property that
\begin{equation}
\text{all but finitely many }i\in I\text{ satisfy }r_{i}=0.
\label{eq.prop.mods.lin-map-pres-lc.fin}%
\end{equation}
Then,
\[
f\left(  \sum_{i\in I}r_{i}m_{i}\right)  =\sum_{i\in I}r_{i}f\left(
m_{i}\right)  .
\]

\end{proposition}

\begin{proof}
[Proof of Proposition \ref{prop.mods.lin-map-pres-lc}.]We give a proof by
example: We assume that $I=\left\{  1,2,3\right\}  $. Thus, the claim we need
to prove is saying that%
\[
f\left(  r_{1}m_{1}+r_{2}m_{2}+r_{3}m_{3}\right)  =r_{1}f\left(  m_{1}\right)
+r_{2}f\left(  m_{2}\right)  +r_{3}f\left(  m_{3}\right)  .
\]
But this is a consequence of the linearity of $f$ (applied several times):%
\begin{align*}
&  f\left(  r_{1}m_{1}+r_{2}m_{2}+r_{3}m_{3}\right) \\
&  =f\left(  r_{1}m_{1}+r_{2}m_{2}\right)  +f\left(  r_{3}m_{3}\right)
\ \ \ \ \ \ \ \ \ \ \left(  \text{since }f\text{ respects addition}\right) \\
&  =f\left(  r_{1}m_{1}\right)  +f\left(  r_{2}m_{2}\right)  +f\left(
r_{3}m_{3}\right)  \ \ \ \ \ \ \ \ \ \ \left(  \text{since }f\text{ respects
addition}\right) \\
&  =r_{1}f\left(  m_{1}\right)  +r_{2}f\left(  m_{2}\right)  +r_{3}f\left(
m_{3}\right)  \ \ \ \ \ \ \ \ \ \ \left(  \text{since }f\text{ respects
scaling}\right)  .
\end{align*}

The same reasoning applies to an arbitrary finite set $I$. (To be fully
rigorous, this is a proof by induction on $\left\vert I\right\vert $.)

The case when $I$ is infinite can be reduced to the case when $I$ is finite
using the assumption (\ref{eq.prop.mods.lin-map-pres-lc.fin}). Indeed, because
of (\ref{eq.prop.mods.lin-map-pres-lc.fin}), there is a finite subset $J$ of
$I$ such that all $i\in I\setminus J$ satisfy $r_{i}=0$. Choosing such a $J$,
we then have
\begin{equation}
\sum_{i\in I}r_{i}m_{i}=\sum_{i\in J}r_{i}m_{i}\ \ \ \ \ \ \ \ \ \ \text{and}%
\ \ \ \ \ \ \ \ \ \ \sum_{i\in I}r_{i}f\left(  m_{i}\right)  =\sum_{i\in
J}r_{i}f\left(  m_{i}\right)  , \label{pf.prop.mods.lin-map-pres-lc.3}%
\end{equation}
since vanishing addends in a sum can be discarded. But since we have already
proved Proposition \ref{prop.mods.lin-map-pres-lc} in the case of a finite set
$I$, we can apply Proposition \ref{prop.mods.lin-map-pres-lc} to $J$, and thus
obtain $f\left(  \sum_{i\in J}r_{i}m_{i}\right)  =\sum_{i\in J}r_{i}f\left(
m_{i}\right)  $. In view of (\ref{pf.prop.mods.lin-map-pres-lc.3}), this
rewrites as $f\left(  \sum_{i\in I}r_{i}m_{i}\right)  =\sum_{i\in I}%
r_{i}f\left(  m_{i}\right)  $, so we are done.
\end{proof}

One useful feature of bases is that they make it easy to define linear maps
out of a free module: Namely, if $M$ is a module with a basis $\left(
m_{i}\right)  _{i\in I}$, and you want to define a linear map $f$ out of $M$,
then it suffices to specify the values $f\left(  m_{i}\right)  $ of the map on
each vector of the basis. These values can be specified arbitrarily; each
possible specification yields a unique linear map $f$. Here is the theorem
that underlies this strategy:

\begin{theorem}
[Universal property of free modules]\label{thm.mods.uniprop-free-1}Let $M$ be
a free left $R$-module with basis $\left(  m_{i}\right)  _{i\in I}$. Let $P$
be a further left $R$-module (not necessarily free). Let $p_{i}\in P$ be a
vector for each $i\in I$. Then, there exists a \textbf{unique} $R$-linear map
$f:M\rightarrow P$ such that%
\begin{equation}
\text{each }i\in I\text{ satisfies }f\left(  m_{i}\right)  =p_{i}.
\label{eq.thm.mods.uniprop-free-1.fmi}%
\end{equation}

\end{theorem}

\begin{proof}
\textbf{Uniqueness:} If $f:M\rightarrow P$ is an $R$-linear map satisfying
(\ref{eq.thm.mods.uniprop-free-1.fmi}), then any $R$-linear combination
$\sum_{i\in I}a_{i}m_{i}$ of $\left(  m_{i}\right)  _{i\in I}$ (where
$a_{i}\in R$ and where all but finitely many $i\in I$ satisfy $a_{i}=0$)
satisfies%
\begin{align}
f\left(  \sum_{i\in I}a_{i}m_{i}\right)   &  =\sum_{i\in I}a_{i}%
\underbrace{f\left(  m_{i}\right)  }_{\substack{=p_{i}\\\text{(by
(\ref{eq.thm.mods.uniprop-free-1.fmi}))}}}\ \ \ \ \ \ \ \ \ \ \left(  \text{by
Proposition \ref{prop.mods.lin-map-pres-lc}}\right) \nonumber\\
&  =\sum_{i\in I}a_{i}p_{i}. \label{pf.thm.mods.uniprop-free-1.uni.1}%
\end{align}
This equality uniquely determines the value of $f$ on each $R$-linear
combination of $\left(  m_{i}\right)  _{i\in I}$. But each element of $M$ can
be written as an $R$-linear combination of $\left(  m_{i}\right)  _{i\in I}$
(since $\left(  m_{i}\right)  _{i\in I}$ is a basis of $M$ and thus spans
$M$). Thus, the equality (\ref{pf.thm.mods.uniprop-free-1.uni.1}) uniquely
determines the value of $f$ on each element of $M$. In other words, it
uniquely determines $f$. Hence, the $R$-linear map $f$ satisfying
(\ref{eq.thm.mods.uniprop-free-1.fmi}) is unique.

\textbf{Existence:} Consider the map%
\begin{align*}
g:R^{\left(  I\right)  }  &  \rightarrow M,\\
\left(  r_{i}\right)  _{i\in I}  &  \mapsto\sum_{i\in I}r_{i}m_{i}.
\end{align*}
This is the map we called $f$ in Theorem \ref{thm.mods.free-is-iso-3} (of
course, we cannot call it $f$ right now, since we need the letter for
something else). Theorem \ref{thm.mods.free-is-iso-3} \textbf{(d)} yields that
the map $g$ is an isomorphism (since the family $\left(  m_{i}\right)  _{i\in
I}$ is a basis of $M$). In particular, this means that $g$ is bijective.
Hence, any element of $M$ can be written as an $R$-linear combination
$\sum_{i\in I}r_{i}m_{i}$ of $\left(  m_{i}\right)  _{i\in I}$ for a
\textbf{unique} family $\left(  r_{i}\right)  _{i\in I}\in R^{\left(
I\right)  }$.

Thanks to this, we can define a map%
\begin{align*}
f:M  &  \rightarrow P,\\
\sum_{i\in I}r_{i}m_{i}  &  \mapsto\sum_{i\in I}r_{i}p_{i}%
\ \ \ \ \ \ \ \ \ \ \left(  \text{for }\left(  r_{i}\right)  _{i\in I}\in
R^{\left(  I\right)  }\right)  .
\end{align*}
Now, it is easy to see that this map $f$ is $R$-linear and satisfies
(\ref{eq.thm.mods.uniprop-free-1.fmi}). Hence, the $R$-linear map $f$
satisfying (\ref{eq.thm.mods.uniprop-free-1.fmi}) exists.

Having proved both existence and uniqueness, we are now done proving Theorem
\ref{thm.mods.uniprop-free-1}.
\end{proof}

In the proof of the \textquotedblleft Uniqueness\textquotedblright\ part
above, we have not used the assumption that the family $\left(  m_{i}\right)
_{i\in I}$ is a basis of $M$; we have only used that it spans $M$. Thus, the
uniqueness of $f$ holds even under this weaker condition. Let us isolate this
into a separate theorem:

\begin{theorem}
[Linear maps are determined on a spanning set]\label{thm.mods.linmap-unidef-1}%
Let $M$ be a left $R$-module. Let $\left(  m_{i}\right)  _{i\in I}$ be a
family of vectors in $M$ that spans $M$. Let $P$ be a further left $R$-module.
Let $f,g:M\rightarrow P$ be two $R$-linear maps such that
\[
\text{each }i\in I\text{ satisfies }f\left(  m_{i}\right)  =g\left(
m_{i}\right)  .
\]
Then, $f=g$.
\end{theorem}

This theorem is often used to prove that two linear maps are equal.

\subsection{\label{sec.modules.bilin}Bilinear maps}

When $R$ is a commutative ring, the addition map%
\[
\operatorname*{add}:R\times R\rightarrow R,\ \ \ \ \ \ \ \ \ \ \left(
a,b\right)  \mapsto a+b
\]
is $R$-linear (where the domain is the direct product of two copies of the
left $R$-module $R$). In fact, if $\left(  a,b\right)  \in R\times R$ and
$\left(  c,d\right)  \in R\times R$ are any two pairs, then%
\begin{align*}
\operatorname*{add}\left(  \underbrace{\left(  a,b\right)  +\left(
c,d\right)  }_{=\left(  a+c,b+d\right)  }\right)   &  =\operatorname*{add}%
\left(  \left(  a+c,b+d\right)  \right)  =\left(  a+c\right)  +\left(
b+d\right)  \ \ \ \ \ \ \ \ \ \ \text{and}\\
\operatorname*{add}\left(  \left(  a,b\right)  \right)  +\operatorname*{add}%
\left(  \left(  c,d\right)  \right)   &  =\left(  a+b\right)  +\left(
c+d\right)  =\left(  a+c\right)  +\left(  b+d\right)
\end{align*}
are clearly the same thing. (This just shows that $\operatorname*{add}$
respects addition; but the other axioms are just as easy.)

In contrast, the multiplication map%
\[
\operatorname*{mul}:R\times R\rightarrow R,\ \ \ \ \ \ \ \ \ \ \left(
a,b\right)  \mapsto ab
\]
is \textbf{not} $R$-linear. However, it is linear in the first argument if we
fix the second. In other words, for any given $b\in R$, the map%
\[
R\rightarrow R,\ \ \ \ \ \ \ \ \ \ a\mapsto ab
\]
is $R$-linear. Likewise, the multiplication map $\operatorname*{mul}:R\times
R\rightarrow R$ is linear in the second argument if we fix the first. Such
maps have a name:

\begin{definition}
\label{def.mods.bilin}Let $R$ be a commutative ring. Let $M$, $N$ and $P$ be
three $R$-modules. A map $f:M\times N\rightarrow P$ is said to be
$R$\textbf{-bilinear} (or just \textbf{bilinear}) if it satisfies the
following two conditions:

\begin{itemize}
\item For any $n\in N$, the map%
\begin{align*}
M  &  \rightarrow P,\\
m  &  \mapsto f\left(  m,n\right)
\end{align*}
is $R$-linear. That is, for any $n\in N$, we have%
\begin{align*}
f\left(  m_{1}+m_{2},n\right)   &  =f\left(  m_{1},n\right)  +f\left(
m_{2},n\right)  \ \ \ \ \ \ \ \ \ \ \text{for all }m_{1},m_{2}\in M;\\
f\left(  rm,n\right)   &  =rf\left(  m,n\right)  \ \ \ \ \ \ \ \ \ \ \text{for
all }r\in R\text{, }m\in M;\\
f\left(  0,n\right)   &  =0.
\end{align*}
This is called \textbf{linearity in the first argument}.

\item For any $m\in M$, the map%
\begin{align*}
N  &  \rightarrow P,\\
n  &  \mapsto f\left(  m,n\right)
\end{align*}
is $R$-linear. That is, for any $m\in M$, we have%
\begin{align*}
f\left(  m,n_{1}+n_{2}\right)   &  =f\left(  m,n_{1}\right)  +f\left(
m,n_{2}\right)  \ \ \ \ \ \ \ \ \ \ \text{for all }n_{1},n_{2}\in N;\\
f\left(  m,rn\right)   &  =rf\left(  m,n\right)  \ \ \ \ \ \ \ \ \ \ \text{for
all }r\in R\text{ and }n\in N;\\
f\left(  m,0\right)   &  =0.
\end{align*}
This is called \textbf{linearity in the second argument}.
\end{itemize}
\end{definition}

Here are some examples of bilinear maps:\footnote{In all these examples, $R$
is assumed to be a commutative ring.}

\begin{itemize}
\item As I said, the multiplication map $R\times R\rightarrow R,\ \left(
a,b\right)  \mapsto ab$ is $R$-bilinear.

\item For any $n\in\mathbb{N}$, the map%
\begin{align*}
R^{n}\times R^{n}  &  \rightarrow R,\\
\left(  \left(  a_{1},a_{2},\ldots,a_{n}\right)  ,\ \left(  b_{1},b_{2}%
,\ldots,b_{n}\right)  \right)   &  \mapsto a_{1}b_{1}+a_{2}b_{2}+\cdots
+a_{n}b_{n}%
\end{align*}
is $R$-bilinear. This map is known as the \textbf{standard scalar product}
(also known as the \textbf{dot product}) on $R^{n}$.

\item Consider the field $\mathbb{C}$ of complex numbers. For any
$n\in\mathbb{N}$, the \textbf{standard inner product}
\begin{align*}
\mathbb{C}^{n}\times\mathbb{C}^{n}  &  \rightarrow\mathbb{C},\\
\left(  \left(  a_{1},a_{2},\ldots,a_{n}\right)  ,\ \left(  b_{1},b_{2}%
,\ldots,b_{n}\right)  \right)   &  \mapsto a_{1}\overline{b_{1}}%
+a_{2}\overline{b_{2}}+\cdots+a_{n}\overline{b_{n}}%
\end{align*}
(where $\overline{z}$ denotes the complex conjugate of a $z\in\mathbb{C}$) is
$\mathbb{R}$-bilinear but not $\mathbb{C}$-bilinear (since it is antilinear
rather than linear in the second argument). However, it becomes $\mathbb{C}%
$-bilinear if you view it as a map $\mathbb{C}^{n}\times\overline{\mathbb{C}%
}^{n}\rightarrow\mathbb{C}$ (with $\overline{\mathbb{C}}$ being the
\textquotedblleft twisted\textquotedblright\ $\mathbb{C}$-module $\mathbb{C}$
from Subsection \ref{subsec.modules.mors.exas2}).

\item The determinant map%
\begin{align*}
\det:R^{2}\times R^{2}  &  \rightarrow R,\\
\left(  \left(  a,b\right)  ,\ \left(  c,d\right)  \right)   &  \mapsto ad-bc
\end{align*}
is $R$-bilinear. (This is called the determinant map because it sends a
$2\times2$-matrix -- encoded as a pair of pairs -- to its determinant.)

\item Matrix multiplication is bilinear. That is: For any $m,n,p\in\mathbb{N}%
$, the map%
\begin{align*}
R^{m\times n}\times R^{n\times p}  &  \rightarrow R^{m\times p},\\
\left(  A,B\right)   &  \mapsto AB
\end{align*}
is $R$-bilinear.

\item The \textbf{cross product} map%
\begin{align*}
R^{3}\times R^{3}  &  \rightarrow R^{3},\\
\left(  \left(  a,b,c\right)  ,\ \left(  x,y,z\right)  \right)   &
\mapsto\left(  bz-cy,\ cx-az,\ ay-bx\right)
\end{align*}
is $R$-bilinear.

\item For any $R$-module $M$, the action map%
\begin{align*}
R\times M  &  \rightarrow M,\\
\left(  r,m\right)   &  \mapsto rm
\end{align*}
is $R$-bilinear. In fact, it is linear in its first argument since every $m\in
M$ satisfies%
\begin{align*}
\left(  r_{1}+r_{2}\right)  m  &  =r_{1}m+r_{2}m\ \ \ \ \ \ \ \ \ \ \text{for
all }r_{1},r_{2}\in R;\\
\left(  rs\right)  m  &  =r\left(  sm\right)  \ \ \ \ \ \ \ \ \ \ \text{for
all }r,s\in R;\\
0_{R}\cdot m  &  =0_{M};
\end{align*}
and it is linear in its second argument since every $r\in R$ satisfies%
\begin{align*}
r\left(  m_{1}+m_{2}\right)   &  =rm_{1}+rm_{2}\ \ \ \ \ \ \ \ \ \ \text{for
all }m_{1},m_{2}\in M;\\
r\left(  sm\right)   &  =s\left(  rm\right)  \ \ \ \ \ \ \ \ \ \ \text{for all
}s\in R\text{ and }m\in M;\\
r\cdot0_{M}  &  =0_{M}.
\end{align*}
(Here, the equality $r\left(  sm\right)  =s\left(  rm\right)  $ follows from
$r\left(  sm\right)  =\underbrace{\left(  rs\right)  }_{=sr}m=\left(
sr\right)  m=s\left(  rm\right)  $. Note how we relied on the commutativity of
$R$ here!)
\end{itemize}

We have always been assuming that $R$ is commutative in this section.
Noncommutative rings $R$ would be a distraction at this point, but might
appear later on.

Theorem \ref{thm.mods.uniprop-free-1} gave us a way to construct linear maps
out of a free module by specifying their values on a basis. We can do the same
for bilinear maps:

\begin{theorem}
[Universal property of free modules wrt bilinear maps]%
\label{thm.mods.uniprop-free-2}Let $R$ be a commutative ring. Let $M$ be a
free $R$-module with basis $\left(  m_{i}\right)  _{i\in I}$. Let $N$ be a
free $R$-module with basis $\left(  n_{j}\right)  _{j\in J}$. Let $P$ be a
further $R$-module (not necessarily free). Let $p_{i,j}\in P$ be a vector for
each pair $\left(  i,j\right)  \in I\times J$. Then, there exists a
\textbf{unique} $R$-bilinear map $f:M\times N\rightarrow P$ such that%
\begin{equation}
\text{each }\left(  i,j\right)  \in I\times J\text{ satisfies }f\left(
m_{i},n_{j}\right)  =p_{i,j}. \label{eq.thm.mods.uniprop-free-2.fmi}%
\end{equation}

\end{theorem}

\begin{proof}
This is similar to the proof of Theorem \ref{thm.mods.uniprop-free-1}. Here
are the details:

\textbf{Uniqueness:} If $f:M\times N\rightarrow P$ is an $R$-bilinear map
satisfying (\ref{eq.thm.mods.uniprop-free-2.fmi}), then any $R$-linear
combination $\sum_{i\in I}a_{i}m_{i}$ of $\left(  m_{i}\right)  _{i\in I}$
(where $a_{i}\in R$ and where all but finitely many $i\in I$ satisfy $a_{i}%
=0$) and any $R$-linear combination $\sum_{j\in J}b_{j}n_{j}$ of $\left(
n_{j}\right)  _{j\in J}$ (where $b_{j}\in R$ and where all but finitely many
$j\in J$ satisfy $b_{j}=0$) satisfy%
\begin{align}
&  f\left(  \sum_{i\in I}a_{i}m_{i},\ \sum_{j\in J}b_{j}n_{j}\right)
\nonumber\\
&  =\sum_{i\in I}a_{i}\underbrace{f\left(  m_{i},\ \sum_{j\in J}b_{j}%
n_{j}\right)  }_{\substack{=\sum_{j\in J}b_{j}f\left(  m_{i},n_{j}\right)
\\\text{(by Proposition \ref{prop.mods.lin-map-pres-lc},}\\\text{since
}f\text{ is }R\text{-linear}\\\text{in its second argument)}}%
}\ \ \ \ \ \ \ \ \ \ \left(
\begin{array}
[c]{c}%
\text{by Proposition \ref{prop.mods.lin-map-pres-lc},}\\
\text{since }f\text{ is }R\text{-linear}\\
\text{in its first argument}%
\end{array}
\right) \nonumber\\
&  =\sum_{i\in I}a_{i}\sum_{j\in J}b_{j}\underbrace{f\left(  m_{i}%
,n_{j}\right)  }_{\substack{=p_{i,j}\\\text{(by
(\ref{eq.thm.mods.uniprop-free-2.fmi}))}}}\ \ \ \ \ \ \ \ \ \ \left(  \text{by
Proposition \ref{prop.mods.lin-map-pres-lc}}\right) \nonumber\\
&  =\sum_{i\in I}a_{i}\sum_{j\in J}b_{j}p_{i,j}=\sum_{\left(  i,j\right)  \in
I\times J}a_{i}b_{j}p_{i,j}. \label{pf.thm.mods.uniprop-free-2.uni.1}%
\end{align}
This equality uniquely determines the value of $f$ on each pair $\left(
x,y\right)  $, where $x$ is an $R$-linear combination of $\left(
m_{i}\right)  _{i\in I}$ and $y$ is an $R$-linear combination of $\left(
n_{j}\right)  _{j\in J}$.

But each element of $M$ can be written as an $R$-linear combination of
$\left(  m_{i}\right)  _{i\in I}$ (since $\left(  m_{i}\right)  _{i\in I}$ is
a basis of $M$ and thus spans $M$), and every element of $N$ can be written as
an $R$-linear combination of $\left(  n_{j}\right)  _{j\in J}$ (for similar
reasons). Thus, the equality (\ref{pf.thm.mods.uniprop-free-2.uni.1}) uniquely
determines the value of $f$ on each pair $\left(  x,y\right)  \in M\times N$.
In other words, it uniquely determines $f$. Hence, the $R$-bilinear map $f$
satisfying (\ref{eq.thm.mods.uniprop-free-2.fmi}) is unique.

\textbf{Existence:} As in the above proof of Theorem
\ref{thm.mods.uniprop-free-1}, we can see that any element of $M$ can be
written as an $R$-linear combination $\sum_{i\in I}r_{i}m_{i}$ of $\left(
m_{i}\right)  _{i\in I}$ for a \textbf{unique} family $\left(  r_{i}\right)
_{i\in I}\in R^{\left(  I\right)  }$. Likewise, any element of $N$ can be
written as an $R$-linear combination $\sum_{j\in J}s_{j}n_{j}$ of $\left(
n_{j}\right)  _{j\in J}$ for a \textbf{unique} family $\left(  s_{j}\right)
_{j\in J}\in R^{\left(  J\right)  }$. Hence, each pair $\left(  x,y\right)
\in M\times N$ can be written in the form $\left(  \sum_{i\in I}r_{i}%
m_{i},\ \sum_{j\in J}s_{j}n_{j}\right)  $ for a \textbf{unique} pair of
families $\left(  r_{i}\right)  _{i\in I}\in R^{\left(  I\right)  }$ and
$\left(  s_{j}\right)  _{j\in J}\in R^{\left(  J\right)  }$.

Thanks to this, we can define a map%
\begin{align*}
f:M\times N  &  \rightarrow P,\\
\left(  \sum_{i\in I}r_{i}m_{i},\ \sum_{j\in J}s_{j}n_{j}\right)   &
\mapsto\sum_{\left(  i,j\right)  \in I\times J}r_{i}s_{j}p_{i,j}%
\ \ \ \ \ \ \ \ \ \ \left(  \text{for }\left(  r_{i}\right)  _{i\in I}\in
R^{\left(  I\right)  }\text{ and }\left(  s_{j}\right)  _{j\in J}\in
R^{\left(  J\right)  }\right)  .
\end{align*}
Now, it is easy to see that this map $f$ is $R$-bilinear and satisfies
(\ref{eq.thm.mods.uniprop-free-2.fmi}). Hence, the $R$-bilinear map $f$
satisfying (\ref{eq.thm.mods.uniprop-free-2.fmi}) exists.

Having proved both existence and uniqueness, we are now done proving Theorem
\ref{thm.mods.uniprop-free-2}.
\end{proof}

\subsection{\label{sec.modules.mulin}Multilinear maps}

Linear and bilinear maps are the first two links in a chain of notions. Here
is the general case:

\begin{definition}
\label{def.mods.mulin}Let $R$ be a commutative ring. Let $M_{1},M_{2}%
,\ldots,M_{n}$ be finitely many $R$-modules. Let $P$ be any $R$-module. A map
$f:M_{1}\times M_{2}\times\cdots\times M_{n}\rightarrow P$ is said to be
$R$\textbf{-multilinear} (or just \textbf{multilinear}) if it satisfies the
following condition:

\begin{itemize}
\item For any $i\in\left\{  1,2,\ldots,n\right\}  $ and any $m_{1}%
,m_{2},\ldots,m_{i-1},m_{i+1},m_{i+2},\ldots,m_{n}$ in the respective modules
(meaning that $m_{k}\in M_{k}$ for each $k\neq i$), the map%
\begin{align*}
M_{i}  &  \rightarrow P,\\
m_{i}  &  \mapsto f\left(  m_{1},m_{2},\ldots,m_{n}\right)
\end{align*}
is $R$-linear. That is, if we fix all arguments of $f$ other than the $i$-th
argument, then $f$ is $R$-linear as a function of the $i$-th argument. This is
called \textbf{linearity in the }$i$\textbf{-th argument}.
\end{itemize}
\end{definition}

Thus, \textquotedblleft bilinear\textquotedblright\ is just \textquotedblleft
multilinear for $n=2$\textquotedblright, whereas \textquotedblleft
linear\textquotedblright\ is \textquotedblleft multilinear for $n=1$%
\textquotedblright.

Here are some examples of multilinear maps:

\begin{itemize}
\item One of the simplest examples of an $R$-multilinear map is the map%
\begin{align*}
\operatorname*{prod}\nolimits_{n}:\underbrace{R\times R\times\cdots\times
R}_{n\text{ times}}  &  \rightarrow R,\\
\left(  a_{1},a_{2},\ldots,a_{n}\right)   &  \mapsto a_{1}a_{2}\cdots a_{n}.
\end{align*}
More generally, for any $r\in R$, the map%
\begin{align*}
\underbrace{R\times R\times\cdots\times R}_{n\text{ times}}  &  \rightarrow
R,\\
\left(  a_{1},a_{2},\ldots,a_{n}\right)   &  \mapsto ra_{1}a_{2}\cdots a_{n}%
\end{align*}
is $R$-multilinear.

\item The most famous example of an $R$-multilinear map is the determinant
function%
\begin{align*}
\det:\underbrace{R^{n}\times R^{n}\times\cdots\times R^{n}}_{n\text{ times}}
&  \rightarrow R,\\
\left(  v_{1},v_{2},\ldots,v_{n}\right)   &  \mapsto\det\left(  v_{1}%
,v_{2},\ldots,v_{n}\right)  ,
\end{align*}
where $\det\left(  v_{1},v_{2},\ldots,v_{n}\right)  $ means the determinant of
the $n\times n$-matrix whose columns are $v_{1},v_{2},\ldots,v_{n}$. (See,
e.g., \cite[Chapter II, Section 7]{Knapp1}, \cite[\S 4.7]{Ford22} or
\cite{Leeb20} for a treatment of determinants based on the concept of
multilinearity\footnote{Note that \cite{Knapp1} is rather terse and abstract,
but it covers the subject almost painlessly if one is familiar with the
advanced viewpoint he is using. More down-to-earth methods are used in
\cite[\S 4.7]{Ford22}, and \cite{Leeb20} is even more elementary (covering
only the case when $R$ is a field; but the same argument can be used in the
general case).}.)
\end{itemize}

There is a universal property of free modules with respect to multilinear maps
(extending Theorem \ref{thm.mods.uniprop-free-1} and Theorem
\ref{thm.mods.uniprop-free-2}), which says that a multilinear map from a
product of free $R$-modules can be defined by specifying its values on all
combinations of basis elements (i.e., on all $n$-tuples whose all entries
belong to the respective bases). I leave it to you to state and prove it.

\subsection{\label{sec.modules.algs}Algebras over commutative rings
(\cite[\S 10.1]{DumFoo04})}

\begin{convention}
In this section, we fix a \textbf{commutative} ring $R$.
\end{convention}

\subsubsection{Definition}

The notion of an $R$\textbf{-algebra} combines the notions of a ring and of an
$R$-module, as well as connecting them by an extra axiom:

\begin{definition}
An $R$\textbf{-algebra} is a set $A$ that is endowed with

\begin{itemize}
\item two binary operations (i.e., maps from $A\times A$ to $A$) that are
called \textbf{addition} and \textbf{multiplication} and denoted by $+$ and
$\cdot$,

\item a map $\cdot$ from $R\times A$ to $A$ that is called \textbf{action} of
$R$ on $A$ (and should not be confused with the multiplication map, which is
also denoted by $\cdot$), and

\item two elements of $A$ that are called \textbf{zero} and \textbf{unity} and
denoted by $0$ and $1$,
\end{itemize}

\noindent such that the following properties (the \textquotedblleft%
\textbf{algebra axioms}\textquotedblright) hold:

\begin{itemize}
\item The addition, the multiplication, the zero and the unity satisfy all the
ring axioms (so that $A$ becomes a ring when equipped with them).

\item The addition, the action and the zero satisfy all the module axioms (so
that $A$ becomes an $R$-module when equipped with them).

\item \textbf{Scale-invariance of multiplication:} We have%
\[
r\left(  ab\right)  =\left(  ra\right)  b=a\left(  rb\right)
\ \ \ \ \ \ \ \ \ \ \text{for all }r\in R\text{ and }a,b\in A.
\]
Here (and in the following), we omit the $\cdot$ signs for multiplication and
action (so \textquotedblleft$ab$\textquotedblright\ means \textquotedblleft%
$a\cdot b$\textquotedblright, and \textquotedblleft$r\left(  ab\right)
$\textquotedblright\ means \textquotedblleft$r\cdot\left(  ab\right)
$\textquotedblright).
\end{itemize}
\end{definition}

Thus, an $R$-algebra is an $R$-module that is also a ring at the same time,
with the same addition (i.e., the addition of the $R$-module must be identical
with the addition of the ring) and the same zero, and satisfying the
\textquotedblleft scale-invariance\textquotedblright\ axiom. In other words,
you get the definition of an $R$-algebra by throwing the definitions of an
$R$-module and a ring together (without duplicating the addition and the zero)
and requiring that the multiplication plays nice with the scaling (in the
sense that scaling a product is equivalent to scaling one of its factors).
Hence, in order to specify an $R$-algebra, it is enough to provide a set with
both a ring structure and an $R$-module structure and show that it satisfies
the \textquotedblleft scale-invariance\textquotedblright\ axiom.

The \textquotedblleft scale-invariance\textquotedblright\ axiom can be
restated as \textquotedblleft the multiplication map
\begin{align*}
A\times A  &  \rightarrow A,\\
\left(  a,b\right)   &  \mapsto ab
\end{align*}
is $R$-bilinear\textquotedblright. More precisely, requiring that the
multiplication map $A\times A\rightarrow A$ be $R$-bilinear is tantamount to
imposing both the scale-invariance axiom and a few of the ring axioms
(distributivity and $0a=a0=0$).

\subsubsection{Examples}

Examples of $R$-algebras include the following:

\begin{itemize}
\item The commutative ring $R$ itself is an $R$-algebra (with both
multiplication and action being the usual multiplication of $R$).

\item The zero ring $\left\{  0\right\}  $ is an $R$-algebra.

\item The matrix ring $R^{n\times n}$ is an $R$-algebra for any $n\in
\mathbb{N}$ (since it is an $R$-module and a ring, and the \textquotedblleft
scale-invariance\textquotedblright\ axiom is easily seen to hold).

\item The ring $\mathbb{C}$ is an $\mathbb{R}$-algebra (since it is an
$\mathbb{R}$-module and a ring, and the \textquotedblleft
scale-invariance\textquotedblright\ axiom is easily seen to hold).

\item The ring $\mathbb{R}$ is a $\mathbb{Q}$-algebra (for similar reasons).

\item More generally: If a commutative ring $R$ is a subring of a
\textbf{commutative} ring $S$, then $S$ becomes an $R$-algebra in a natural
way. In fact, we already know from Subsection \ref{subsec.modules.def.restr}
that $S$ becomes an $R$-module, and it is easy to see that this $R$-module can
be combined with the ring structure on $S$ to form an $R$-algebra.

\item Even more generally: If $R$ and $S$ are two \textbf{commutative} rings,
and if $f:R\rightarrow S$ is a ring morphism, then $S$ becomes an $R$-algebra
in a natural way. In fact, we already know from Subsection
\ref{subsec.modules.def.restr} that $S$ becomes an $R$-module (this is the
$R$-module structure on $S$ induced by $f$), and it is easy to see that this
$R$-module can be combined with the ring structure on $S$ to form an
$R$-algebra.\footnote{In particular, this yields that any quotient ring $R/I$
of $R$ becomes an $R$-algebra, via the canonical projection $\pi:R\rightarrow
R/I$.} This $R$-algebra structure on $S$ is said to be \textbf{induced} by the
morphism $f$.

\item Yet more generally: If $R$ and $S$ are two commutative rings, and if
$f:R\rightarrow S$ is a ring morphism, then any $S$-algebra $A$ becomes an
$R$-algebra in a natural way. In fact, we already know from Subsection
\ref{subsec.modules.def.restr} that $A$ becomes an $R$-module (this is the
$R$-module obtained by restricting the $S$-module $A$ to $R$), and it is easy
to see that this $R$-module can be combined with the ring structure on $A$ to
form an $R$-algebra. This is called the $R$-algebra obtained by
\textbf{restricting} the $S$-algebra $A$ to $R$.

For example, the matrix ring $\mathbb{C}^{2\times2}$ is a $\mathbb{C}%
$-algebra, and thus becomes an $\mathbb{R}$-algebra (since the inclusion map
$\mathbb{R}\rightarrow\mathbb{C}$ is a ring morphism).

\item The quaternion ring $\mathbb{H}$ is an $\mathbb{R}$-algebra. But it is
not a $\mathbb{C}$-algebra, even though it contains $\mathbb{C}$ as a subring.
Indeed, the \textquotedblleft scale-invariance\textquotedblright\ axiom for
$\mathbb{H}$ to be a $\mathbb{C}$-algebra would say that
\[
r\left(  ab\right)  =\left(  ra\right)  b=a\left(  rb\right)
\ \ \ \ \ \ \ \ \ \ \text{for all }r\in\mathbb{C}\text{ and }a,b\in
\mathbb{H};
\]
but this is \textbf{not true} for $r=i$, $a=j$ and $b=1$ because $ij\neq ji$.

This does not contradict the previous bullet points! After all, $\mathbb{H}$
is not commutative.

\item The polynomial ring $R\left[  x\right]  $ (to be defined soon) is an $R$-algebra.
\end{itemize}

\begin{exercise}
\label{exe.algs.Cgen}Recall that the complex numbers were defined as pairs of
real numbers, with entrywise addition and a strange-looking multiplication.
Let us now generalize this construction by replacing real numbers by elements
of the given commutative ring $R$.

Thus, we define an $R$\textbf{-complex number} to be a pair $\left(
a,b\right)  \in R\times R$. We define $\mathbb{C}_{R}$ to be the set of all
$R$-complex numbers. We define an addition $+$ and a multiplication $\cdot$ on
this set $\mathbb{C}_{R}$ by the formulas%
\begin{align*}
\left(  a,b\right)  +\left(  c,d\right)   &  =\left(  a+c,\ b+d\right)
\ \ \ \ \ \ \ \ \ \ \text{and}\\
\left(  a,b\right)  \cdot\left(  c,d\right)   &  =\left(
ac-bd,\ ad+bc\right)
\end{align*}
(these are the same formulas as for the original complex numbers).

\begin{enumerate}
\item[\textbf{(a)}] Prove that the set $\mathbb{C}_{R}$ (equipped with these
two operations, with the zero $\left(  0,0\right)  $ and the unity $\left(
1,0\right)  $) becomes a commutative ring.
\end{enumerate}

We denote this ring by $\mathbb{C}_{R}$, and call it the ring of
$R$\textbf{-complex numbers}. For $R=\mathbb{R}$, we recover the usual complex
numbers: $\mathbb{C}_{\mathbb{R}}=\mathbb{C}$.

We further turn this ring $\mathbb{C}_{R}$ into an $R$-algebra by defining an
action of $R$ on $\mathbb{C}_{R}$ by the equation%
\[
r\left(  a,b\right)  =\left(  ra,rb\right)  \ \ \ \ \ \ \ \ \ \ \text{for any
}r\in R\text{ and any }\left(  a,b\right)  \in\mathbb{C}_{R}.
\]

\begin{enumerate}
\item[\textbf{(b)}] Prove that this does indeed make $\mathbb{C}_{R}$ into an
$R$-algebra.
\end{enumerate}

Let $i$ denote the element $\left(  0,1\right)  $ of $\mathbb{C}_{R}$. (This
is the generalization of the imaginary unit $i$ of $\mathbb{C}$.)

\begin{enumerate}
\item[\textbf{(c)}] Prove that $\left(  a,b\right)  =a\cdot1_{\mathbb{C}_{R}%
}+bi$ for each $\left(  a,b\right)  \in\mathbb{C}_{R}$.
\end{enumerate}

Now, recall that $\mathbb{C}$ is a field. What about its generalized version
$\mathbb{C}_{R}$ ?

\begin{enumerate}
\item[\textbf{(d)}] Prove that $\mathbb{C}_{\mathbb{C}}$ is not a field, but
rather $\mathbb{C}_{\mathbb{C}}\cong\mathbb{C}\times\mathbb{C}$ as rings.

\item[\textbf{(e)}] Prove that $\mathbb{C}_{\mathbb{Z}/2}$ is not a field, and
in fact is isomorphic to the ring $D_{4}$ from Subsection
\ref{subsec.rings.def.exas}.

\item[\textbf{(f)}] Prove that $\mathbb{C}_{\mathbb{Z}/3}$ is a field with $9$ elements.

\item[\textbf{(g)}] Let $S$ be the ring $\mathbb{Z}\left[  i\right]  $ of
Gaussian integers. Prove that $\mathbb{C}_{\mathbb{Z}/n}\cong S/\left(
nS\right)  $ as rings for any integer $n$.
\end{enumerate}
\end{exercise}

\begin{fineprint}
We can similarly generalize the dual numbers (as defined in Exercise
\ref{exe.dualnums.ring}):
\end{fineprint}

\begin{exercise}
\label{exe.algs.dualnums}We define an $R$\textbf{-dual number} to be a pair
$\left(  a,b\right)  \in R\times R$. We define $\mathbb{D}_{R}$ to be the set
of all $R$-dual numbers. We define an addition $+$ and a multiplication
$\cdot$ on this set $\mathbb{D}_{R}$ by the formulas%
\begin{align*}
\left(  a,b\right)  +\left(  c,d\right)   &  =\left(  a+c,\ b+d\right)
\ \ \ \ \ \ \ \ \ \ \text{and}\\
\left(  a,b\right)  \cdot\left(  c,d\right)   &  =\left(  ac,\ ad+bc\right)
\end{align*}
(these are the same formulas as for the original dual numbers).

\begin{enumerate}
\item[\textbf{(a)}] Prove that the set $\mathbb{D}_{R}$ (equipped with these
two operations, with the zero $\left(  0,0\right)  $ and the unity $\left(
1,0\right)  $) becomes a commutative ring.
\end{enumerate}

We denote this ring by $\mathbb{D}_{R}$, and call it the ring of
$R$\textbf{-dual numbers}. For $R=\mathbb{R}$, we recover the usual dual
numbers from Exercise \ref{exe.dualnums.ring}: that is, $\mathbb{D}%
_{\mathbb{R}}=\mathbb{D}$.

We further turn this ring $\mathbb{D}_{R}$ into an $R$-algebra by defining an
action of $R$ on $\mathbb{D}_{R}$ by the equation%
\[
r\left(  a,b\right)  =\left(  ra,rb\right)  \ \ \ \ \ \ \ \ \ \ \text{for any
}r\in R\text{ and any }\left(  a,b\right)  \in\mathbb{D}_{R}.
\]

\begin{enumerate}
\item[\textbf{(b)}] Prove that this does indeed make $\mathbb{D}_{R}$ into an
$R$-algebra.
\end{enumerate}

Let $\varepsilon$ denote the element $\left(  0,1\right)  $ of $\mathbb{D}%
_{R}$.

\begin{enumerate}
\item[\textbf{(c)}] Prove that, for any $a,b\in R$, we have $\left(
a,b\right)  =a\cdot1_{\mathbb{D}_{R}}+b\varepsilon$ in $\mathbb{D}_{R}$.

\item[\textbf{(d)}] Prove that $\varepsilon\in\mathbb{D}_{R}$ is nilpotent,
and in fact $\varepsilon^{2}=0$.

\item[\textbf{(e)}] Prove that $\mathbb{D}_{\mathbb{Z}/2}$ is isomorphic to
the ring $D_{4}$ from Subsection \ref{subsec.rings.def.exas}.
\end{enumerate}
\end{exercise}

\begin{exercise}
\label{exe.21hw3.10cde}Let $R$ be a ring. Let $M$ be an $R$-module. Recall
that the Hom group $\operatorname*{Hom}\nolimits_{R}\left(  M,M\right)  $ is
an additive abelian group (by Exercise \ref{exe.21hw3.10ab} \textbf{(c)}).
Moreover, if $R$ is commutative, then $\operatorname*{Hom}\nolimits_{R}\left(
M,M\right)  $ is also an $R$-module (by Exercise \ref{exe.21hw3.10ab}
\textbf{(f)}). Prove the following:

\begin{enumerate}
\item[\textbf{(a)}] The Hom group $\operatorname{Hom}_{R}\left(  M,M\right)  $
becomes a ring if we define multiplication to be composition (i.e., for any
$f\in\operatorname{Hom}_{R}\left(  M,M\right)  $ and $g\in\operatorname{Hom}%
_{R}\left(  M,M\right)  $, we define $fg$ to be the composition $f\circ g$).
Its unity is the identity map $\operatorname{id}:M\rightarrow M$.

This ring $\operatorname{Hom}_{R}\left(  M,M\right)  $ is also denoted
$\operatorname{End}_{R}\left(  M\right)  $ and known as the
\textbf{endomorphism ring}\footnotemark\ of $M$.

\item[\textbf{(b)}] If $R$ is commutative, then the endomorphism ring
$\operatorname{End}_{R}\left(  M\right)  $ becomes an $R$-algebra (with the
$R$-module structure defined as in Exercise \ref{exe.21hw3.10ab} \textbf{(f)}).

\item[\textbf{(c)}] Let $M=R^{\mathbb{N}}$; this is the left $R$-module of all
infinite sequences $\left(  a_{0},a_{1},a_{2},\ldots\right)  $ of elements of
$R$. Define two left $R$-module morphisms $f:M\rightarrow M$ and
$g:M\rightarrow M$ by
\[
f\left(  a_{0},a_{1},a_{2},\ldots\right)  =\left(  a_{1},a_{2},a_{3}%
,\ldots\right)  \qquad\text{for any $\left(  a_{0},a_{1},a_{2},\ldots\right)
\in M$}%
\]
and
\[
g\left(  a_{0},a_{1},a_{2},\ldots\right)  =\left(  0,a_{0},a_{1},a_{2}%
,\ldots\right)  \qquad\text{for any $\left(  a_{0},a_{1},a_{2},\ldots\right)
\in M$}.
\]
Prove that $fg=1$ but $gf\neq1$ (unless $R$ is trivial) in the ring
$\operatorname*{End}\nolimits_{R}\left(  M\right)  $. (This gives an example
of a left inverse that is not a right inverse.)
\end{enumerate}
\end{exercise}

\footnotetext{An \textbf{endomorphism} is defined to be a morphism from a
structure (e.g., a ring or a module or an algebra) to itself.}

\subsubsection{Rings as $\mathbb{Z}$-algebras}

The most common algebras are the $\mathbb{Z}$-algebras. In fact, every ring is
a $\mathbb{Z}$-algebra in a natural way:

\begin{proposition}
\label{prop.Zalgs.1}Let $A$ be any ring. Then, $A$ is an abelian group (with
respect to addition), so $A$ becomes a $\mathbb{Z}$-module (since we have seen
in Proposition \ref{prop.Zmods.1} that every abelian group naturally becomes a
$\mathbb{Z}$-module). This $\mathbb{Z}$-module structure can be combined with
the ring structure on $A$, turning $A$ into a $\mathbb{Z}$-algebra.
\end{proposition}

\begin{proof}
You have to check \textquotedblleft scale-invariance\textquotedblright. This
is easy and LTTR.
\end{proof}

Thus, every ring becomes a $\mathbb{Z}$-algebra (similarly to how any abelian
group becomes a $\mathbb{Z}$-module). This allows us to equate rings with
$\mathbb{Z}$-algebras. We shall do this whenever convenient.

\subsubsection{The underlying structures}

Every $R$-algebra $A$ has an underlying ring (i.e., the ring obtained from $A$
by forgetting the action) and an underlying $R$-module (i.e., the $R$-module
obtained from $A$ by forgetting the multiplication and the unity); we will
refer to these simply as the \textquotedblleft ring $A$\textquotedblright\ and
the \textquotedblleft$R$-module $A$\textquotedblright. So, for example, if $A$
and $B$ are two $R$-algebras, then the \textquotedblleft ring morphisms from
$A$ to $B$\textquotedblright\ will simply mean the ring morphisms from the
underlying ring of $A$ to the underlying ring of $B$. Similarly the
\textquotedblleft$R$-module morphisms from $A$ to $B$\textquotedblright\ are
to be understood.

\begin{fineprint}
\textbf{Warning:} There is also a notion of \textquotedblleft base
ring\textquotedblright\ which sounds similar to, but has nothing to do with,
the notion of an \textquotedblleft underlying ring\textquotedblright! The
\textbf{base ring} of an $R$-algebra $A$ is defined to be $R$ (unlike the
underlying ring of $A$, which is just $A$ with its action removed).
\end{fineprint}

\subsubsection{Commutative $R$-algebras}

\begin{definition}
An $R$-algebra is said to be \textbf{commutative} if its underlying ring is
commutative (i.e., if its multiplication is commutative).
\end{definition}

\subsubsection{Subalgebras}

Algebras have subalgebras; they are defined exactly as you would expect:

\begin{definition}
\label{def.algs.subalg}Let $A$ be an $R$-algebra. An $R$\textbf{-subalgebra}
of $A$ means a subset of $A$ that is simultaneously a subring and an
$R$-submodule of $A$.
\end{definition}

In pedestrian terms, this means that an $R$-subalgebra of $A$ is a subset of
$A$ that is closed under addition, multiplication and scaling and contains the
zero and the unity. Such an $R$-subalgebra of $A$ clearly becomes an
$R$-algebra in its own right (since we can restrict all relevant operations
from $A$ to this subalgebra).

\subsubsection{$R$-algebra morphisms}

Just as rings have ring morphisms, and $R$-modules have $R$-module morphisms,
there is a notion of $R$-algebra morphisms:

\begin{definition}
\label{def.algs.mor}Let $A$ and $B$ be two $R$-algebras.

\begin{enumerate}
\item[\textbf{(a)}] An $R$\textbf{-algebra morphism} (or, short,
\textbf{algebra morphism}) from $A$ to $B$ means a map $f:A\rightarrow B$ that
is both a ring morphism and an $R$-module morphism (i.e., that respects
addition, multiplication, zero, unity and scaling).

\item[\textbf{(b)}] An $R$\textbf{-algebra isomorphism} (or, informally,
\textbf{algebra iso}) from $A$ to $B$ means an invertible $R$-algebra morphism
$f:A\rightarrow B$ whose inverse $f^{-1}:B\rightarrow A$ is also an
$R$-algebra morphism.

\item[\textbf{(c)}] The $R$-algebras $A$ and $B$ are said to be
\textbf{isomorphic} (this is written $A\cong B$) if there exists an
$R$-algebra isomorphism from $A$ to $B$.
\end{enumerate}
\end{definition}

All the fundamental properties of ring morphisms (stated in Subsection
\ref{subsec.ringmor.mors.basic}) and of ring isomorphisms (stated in
Subsection \ref{subsec.ringmor.iso.basic-props}) have analogues for
$R$-algebra morphisms and isomorphisms, respectively. For example, here is the
analogue of Proposition \ref{prop.ringmor.Im-subring}:

\begin{proposition}
\label{prop.algmor.Im-subalg}Let $A$ and $B$ be two $R$-algebras. Let
$f:A\rightarrow B$ be an $R$-algebra morphism. Then, $\operatorname{Im}%
f=f\left(  A\right)  $ is an $R$-subalgebra of $B$.
\end{proposition}

\begin{proof}
Analogous to the proof of Proposition \ref{prop.ringmor.Im-subring}.
\end{proof}

And here is the analogue of Proposition \ref{prop.ringmor.invertible-iso}:

\begin{proposition}
\label{prop.algmor.invertible-iso}Let $A$ and $B$ be two $R$-algebras. Let
$f:A\rightarrow B$ be an invertible $R$-algebra morphism. Then, $f$ is an
$R$-algebra isomorphism.
\end{proposition}

\begin{proof}
Analogous to the proof of Proposition \ref{prop.ringmor.invertible-iso}.
\end{proof}

An analogue to Proposition \ref{prop.modules.mors.Z} is the following:

\begin{proposition}
\label{prop.algmor.mors.Z}Let $A$ and $B$ be two $\mathbb{Z}$-algebras. Then,
the $\mathbb{Z}$-algebra morphisms from $A$ to $B$ are precisely the ring
morphisms from $A$ to $B$.
\end{proposition}

\begin{proof}
Analogous to the proof of Proposition \ref{prop.modules.mors.Z}.
\end{proof}

\subsubsection{Direct products of algebras}

The direct product of several $R$-algebras is defined just as you would
expect: addition, multiplication and scaling are all entrywise. Just for the
sake of completeness, let me give its precise definition:

\begin{proposition}
\label{prop.algs.dirprod-exists}Let $I$ be any set. Let $\left(  A_{i}\right)
_{i\in I}$ be any family of $R$-algebras. Then, their Cartesian product
$\prod\limits_{i\in I}A_{i}$ becomes an $R$-algebra if we endow it with the
entrywise addition (i.e., we set $\left(  m_{i}\right)  _{i\in I}+\left(
n_{i}\right)  _{i\in I}=\left(  m_{i}+n_{i}\right)  _{i\in I}$ for any two
families $\left(  m_{i}\right)  _{i\in I},\left(  n_{i}\right)  _{i\in I}%
\in\prod\limits_{i\in I}A_{i}$) and the entrywise multiplication (i.e., we set
$\left(  m_{i}\right)  _{i\in I}\cdot\left(  n_{i}\right)  _{i\in I}=\left(
m_{i}\cdot n_{i}\right)  _{i\in I}$ for any two families $\left(
m_{i}\right)  _{i\in I},\left(  n_{i}\right)  _{i\in I}\in\prod\limits_{i\in
I}A_{i}$) and the entrywise scaling (i.e., we set $r\left(  m_{i}\right)
_{i\in I}=\left(  rm_{i}\right)  _{i\in I}$ for any $r\in R$ and any family
$\left(  m_{i}\right)  _{i\in I}\in\prod\limits_{i\in I}A_{i}$) and with the
zero $\left(  0\right)  _{i\in I}$ and the unity $\left(  1\right)  _{i\in I}%
$. The underlying ring of this $R$-algebra $\prod\limits_{i\in I}A_{i}$ is the
direct product of the rings $A_{i}$, whereas the underlying $R$-module of this
$R$-algebra $\prod\limits_{i\in I}A_{i}$ is the direct product of the
$R$-modules $A_{i}$.
\end{proposition}

\begin{definition}
\label{def.algs.dirprod}This $R$-algebra is denoted by $\prod\limits_{i\in
I}A_{i}$ and called the \textbf{direct product} of the $R$-algebras $A_{i}$.

The usual notations apply to these direct products: For example, if
$I=\left\{  1,2,\ldots,n\right\}  $ for some $n\in\mathbb{N}$, then the direct
product $\prod\limits_{i\in I}A_{i}$ is also denoted by $A_{1}\times
A_{2}\times\cdots\times A_{n}$; we further set $A^{I}=\prod\limits_{i\in I}A$
and $A^{n}=A^{\left\{  1,2,\ldots,n\right\}  }$ for each $n\in\mathbb{N}$.
\end{definition}

\subsection{\label{sec.modules.H}Defining algebras: the case of $\mathbb{H}$}

You can think of an $R$-algebra as a ring equipped with an additional piece of
structure -- namely, with an action of $R$ on it. Thus, in order to define an
$R$-algebra, it is natural to start by defining a ring and then defining an
action of $R$ on it (and showing that it satisfies the $R$-module axioms and scale-invariance).

Often, however, it is easier to proceed differently: First, define an
$R$-module, and then define the multiplication and the unity to turn it into
an $R$-algebra. If you do things in this order, you can use the $R$-module
structure as scaffolding for defining the multiplication. This often turns out
to be the simpler way.

Here is an example of how this can work:

Recall the ring $\mathbb{H}$ of quaternions, which were \textquotedblleft
defined\textquotedblright\ to be \textquotedblleft numbers\textquotedblright%
\ of the form $a+bi+cj+dk$ with $a,b,c,d\in\mathbb{R}$ and with the
multiplication rules%
\[
i^{2}=j^{2}=k^{2}%
=-1,\ \ \ \ \ \ \ \ \ \ ij=-ji=k,\ \ \ \ \ \ \ \ \ \ jk=-kj=i,\ \ \ \ \ \ \ \ \ \ ki=-ik=j.
\]
It is clear how to calculate in $\mathbb{H}$ using these rules. But why does
this ring $\mathbb{H}$ exist?

Here is a cautionary tale to show why this is a question: Let's replace
$k^{2}=-1$ by $k^{2}=1$ in our above \textquotedblleft
definition\textquotedblright\ of $\mathbb{H}$ (but still require $i^{2}$ and
$j^{2}$ to be $-1$). Then, $j^{2}\underbrace{k^{2}}_{=1}=j^{2}=-1$, so that%
\[
-1=j^{2}k^{2}=j\underbrace{jk}_{=i}k=j\underbrace{ik}%
_{\substack{=-j\\\text{(since }-ik=j\text{)}}}=j\left(  -j\right)
=-\underbrace{j^{2}}_{=-1}=-\left(  -1\right)  =1.
\]
Adding $1$ to this equality, we find $0=2$, so that $0=1$ (upon division by
$2$). Therefore, the ring is trivial -- i.e., all its elements are $0$.

Ouch. We tried to expand our number system by introducing new
\textquotedblleft numbers\textquotedblright\ $i,j,k$, but instead we ended up
collapsing it (making all numbers equal to $0$).

It should not surprise you that this can happen; after all, the same happens
if you introduce the \textquotedblleft number\textquotedblright\ $\infty
:=\dfrac{1}{0}$ and start doing algebra with it. But why doesn't it happen
with the quaternions? Why is $\mathbb{H}$ actually an extension of our number
system rather than a collapsed version of it?

The simplest way to answer this question is to throw away the wishy-washy
definition of $\mathbb{H}$ we gave above (what does \textquotedblleft numbers
of the form $a+bi+cj+dk$\textquotedblright\ really mean?), and redefine
$\mathbb{H}$ rigorously.

We want $\mathbb{H}$ to be an $\mathbb{R}$-algebra. First, we introduce its
underlying $\mathbb{R}$-module (i.e., $\mathbb{R}$-vector space) structure.
This underlying $\mathbb{R}$-module will be a $4$-dimensional $\mathbb{R}%
$-vector space, i.e., a free $\mathbb{R}$-module of rank $4$. So let me
\textbf{define} $\mathbb{H}$ to be $\mathbb{R}^{4}$ as an $\mathbb{R}$-module.
Let me denote its standard basis by $\left(  \mathbf{e},\mathbf{i}%
,\mathbf{j},\mathbf{k}\right)  $ (so that $\mathbf{e}=\left(  1,0,0,0\right)
$ and $\mathbf{i}=\left(  0,1,0,0\right)  $ and $\mathbf{j}=\left(
0,0,1,0\right)  $ and $\mathbf{k}=\left(  0,0,0,1\right)  $). These four basis
vectors $\mathbf{e},\mathbf{i},\mathbf{j},\mathbf{k}$ will later become the
quaternions $1,i,j,k$, but I'm being cautious for now and avoiding any names
that might be too suggestive. The basis vector $\mathbf{e}$ will be the unity
of $\mathbb{H}$. Next, we define the multiplication of $\mathbb{H}$ to be the
$\mathbb{R}$-bilinear map $\mu:\mathbb{H}\times\mathbb{H}\rightarrow
\mathbb{H}$ that satisfies\footnote{These equations are not chosen at random,
of course; they are simply the equations%
\[
i^{2}=j^{2}=k^{2}%
=-1,\ \ \ \ \ \ \ \ \ \ ij=-ji=k,\ \ \ \ \ \ \ \ \ \ jk=-kj=i,\ \ \ \ \ \ \ \ \ \ ki=-ik=j
\]
(as well as the equations $1\cdot1=1$, $1i=i$, $1j=j$, $1k=k$, $i\cdot1=i$,
$j\cdot1=j$ and $k\cdot1=k$), with $1,i,j,k$ renamed as $\mathbf{e}%
,\mathbf{i},\mathbf{j},\mathbf{k}$.}
\begin{align*}
\mu\left(  \mathbf{e},\mathbf{e}\right)   &  =\mathbf{e}%
,\ \ \ \ \ \ \ \ \ \ \mu\left(  \mathbf{e},\mathbf{i}\right)  =\mathbf{i}%
,\ \ \ \ \ \ \ \ \ \ \mu\left(  \mathbf{e},\mathbf{j}\right)  =\mathbf{j}%
,\ \ \ \ \ \ \ \ \ \ \mu\left(  \mathbf{e},\mathbf{k}\right)  =\mathbf{k},\\
\mu\left(  \mathbf{i},\mathbf{e}\right)   &  =\mathbf{i}%
,\ \ \ \ \ \ \ \ \ \ \mu\left(  \mathbf{i},\mathbf{i}\right)  =-\mathbf{e}%
,\ \ \ \ \ \ \ \ \ \ \mu\left(  \mathbf{i},\mathbf{j}\right)  =\mathbf{k}%
,\ \ \ \ \ \ \ \ \ \ \mu\left(  \mathbf{i},\mathbf{k}\right)  =-\mathbf{j},\\
\mu\left(  \mathbf{j},\mathbf{e}\right)   &  =\mathbf{j}%
,\ \ \ \ \ \ \ \ \ \ \mu\left(  \mathbf{j},\mathbf{i}\right)  =-\mathbf{k}%
,\ \ \ \ \ \ \ \ \ \ \mu\left(  \mathbf{j},\mathbf{j}\right)  =-\mathbf{e}%
,\ \ \ \ \ \ \ \ \ \ \mu\left(  \mathbf{j},\mathbf{k}\right)  =\mathbf{i},\\
\mu\left(  \mathbf{k},\mathbf{e}\right)   &  =\mathbf{k}%
,\ \ \ \ \ \ \ \ \ \ \mu\left(  \mathbf{k},\mathbf{i}\right)  =\mathbf{j}%
,\ \ \ \ \ \ \ \ \ \ \mu\left(  \mathbf{k},\mathbf{j}\right)  =-\mathbf{i}%
,\ \ \ \ \ \ \ \ \ \ \mu\left(  \mathbf{k},\mathbf{k}\right)  =-\mathbf{e}.
\end{align*}
Theorem \ref{thm.mods.uniprop-free-2} guarantees that there is a unique such
bilinear map $\mu$. We set $ab=\mu\left(  a,b\right)  $ for all $a,b\in
\mathbb{H}$.

Why is this a ring? All but two of the ring axioms are obvious (they follow
either from the bilinearity of $\mu$ or from the module axioms for the
$\mathbb{R}$-module $\mathbb{H}=\mathbb{R}^{4}$). The two axioms that are not
obvious are the following:

\begin{enumerate}
\item Associativity of multiplication.

\item Neutrality of $1$ (i.e., the claim that $a\cdot\mathbf{e}=\mathbf{e}%
\cdot a=a$ for each $a\in\mathbb{H}$).
\end{enumerate}

\noindent Fortunately, the bilinearity of $\mu$ will make both of these axioms
straightforward to check. Indeed, let me explain how to check the
associativity of multiplication. In other words, let me prove that the map
$\mu$ is associative -- i.e., that%
\begin{equation}
\mu\left(  \mu\left(  a,b\right)  ,c\right)  =\mu\left(  a,\mu\left(
b,c\right)  \right)  \ \ \ \ \ \ \ \ \ \ \text{for all }a,b,c\in\mathbb{H}.
\label{eq.Hring.assoc-mu}%
\end{equation}

The trick to this is that when a map like $\mu$ is bilinear, its associativity
can be checked on a basis -- or, more generally, on a spanning set:

\begin{lemma}
\label{lem.algebras.assoc-on-basis}Let $R$ be a commutative ring. Let $M$ be
an $R$-module. Let $\left(  m_{i}\right)  _{i\in I}$ be a family of vectors in
$M$ that spans $M$. Let $f:M\times M\rightarrow M$ be an $R$-bilinear map.
Assume that%
\begin{equation}
f\left(  f\left(  m_{i},m_{j}\right)  ,m_{k}\right)  =f\left(  m_{i},f\left(
m_{j},m_{k}\right)  \right)  \ \ \ \ \ \ \ \ \ \ \text{for all }i,j,k\in I.
\label{eq.lem.algebras.assoc-on-basis.ass}%
\end{equation}
Then, we have%
\begin{equation}
f\left(  f\left(  a,b\right)  ,c\right)  =f\left(  a,f\left(  b,c\right)
\right)  \ \ \ \ \ \ \ \ \ \ \text{for all }a,b,c\in M.
\label{eq.lem.algebras.assoc-on-basis.clm}%
\end{equation}

\end{lemma}

\begin{proof}
[Proof of Lemma \ref{lem.algebras.assoc-on-basis}.]Let $a,b,c\in M$. Since the
family $\left(  m_{i}\right)  _{i\in I}$ spans $M$, we can write the three
vectors $a,b,c$ as
\[
a=\sum_{i\in I}a_{i}m_{i},\ \ \ \ \ \ \ \ \ \ b=\sum_{j\in I}b_{j}%
m_{j},\ \ \ \ \ \ \ \ \ \ c=\sum_{k\in I}c_{k}m_{k}%
\]
for some coefficients $a_{i},b_{j},c_{k}\in R$. Consider these coefficients.
Then,\footnote{We will use Proposition \ref{prop.mods.lin-map-pres-lc}
multiple times in this computation.}
\begin{align*}
f\left(  f\left(  a,b\right)  ,c\right)   &  =f\left(  f\left(  \sum_{i\in
I}a_{i}m_{i},\sum_{j\in I}b_{j}m_{j}\right)  ,\sum_{k\in I}c_{k}m_{k}\right)
\\
&  =\sum_{k\in I}c_{k}f\left(  f\left(  \sum_{i\in I}a_{i}m_{i},\sum_{j\in
I}b_{j}m_{j}\right)  ,m_{k}\right) \\
&  \ \ \ \ \ \ \ \ \ \ \ \ \ \ \ \ \ \ \ \ \left(  \text{since }f\text{ is
linear in its second argument}\right) \\
&  =\sum_{k\in I}c_{k}f\left(  \sum_{i\in I}a_{i}f\left(  m_{i},\sum_{j\in
I}b_{j}m_{j}\right)  ,m_{k}\right) \\
&  \ \ \ \ \ \ \ \ \ \ \ \ \ \ \ \ \ \ \ \ \left(  \text{since }f\text{ is
linear in its first argument}\right) \\
&  =\sum_{k\in I}c_{k}f\left(  \sum_{i\in I}a_{i}\sum_{j\in I}b_{j}f\left(
m_{i},m_{j}\right)  ,m_{k}\right) \\
&  \ \ \ \ \ \ \ \ \ \ \ \ \ \ \ \ \ \ \ \ \left(  \text{since }f\text{ is
linear in its second argument}\right) \\
&  =\sum_{k\in I}c_{k}\sum_{i\in I}a_{i}\sum_{j\in I}b_{j}f\left(  f\left(
m_{i},m_{j}\right)  ,m_{k}\right) \\
&  \ \ \ \ \ \ \ \ \ \ \ \ \ \ \ \ \ \ \ \ \left(  \text{since }f\text{ is
linear in its first argument}\right) \\
&  =\sum_{i\in I}\ \ \sum_{j\in I}\ \ \sum_{k\in I}a_{i}b_{j}c_{k}f\left(
f\left(  m_{i},m_{j}\right)  ,m_{k}\right)
\end{align*}
and similarly%
\[
f\left(  a,f\left(  b,c\right)  \right)  =\sum_{i\in I}\ \ \sum_{j\in
I}\ \ \sum_{k\in I}a_{i}b_{j}c_{k}f\left(  m_{i},f\left(  m_{j},m_{k}\right)
\right)  .
\]
The right hand sides of these two equalities are equal by our assumption
(\ref{eq.lem.algebras.assoc-on-basis.ass}). Hence, the left hand sides are
equal. In other words, $f\left(  f\left(  a,b\right)  ,c\right)  =f\left(
a,f\left(  b,c\right)  \right)  $. This proves Lemma
\ref{lem.algebras.assoc-on-basis}.
\end{proof}

Let us now return to $\mathbb{H}$. We want to prove that
\[
\mu\left(  \mu\left(  a,b\right)  ,c\right)  =\mu\left(  a,\mu\left(
b,c\right)  \right)  \ \ \ \ \ \ \ \ \ \ \text{for all }a,b,c\in\mathbb{H}.
\]
By Lemma \ref{lem.algebras.assoc-on-basis} (applied to $R=\mathbb{R}$,
$M=\mathbb{H}$, $\left(  m_{i}\right)  _{i\in I}=\left(  \mathbf{e}%
,\mathbf{i},\mathbf{j},\mathbf{k}\right)  $ and $f=\mu$), it suffices to show
that%
\[
\mu\left(  \mu\left(  a,b\right)  ,c\right)  =\mu\left(  a,\mu\left(
b,c\right)  \right)  \ \ \ \ \ \ \ \ \ \ \text{for all }a,b,c\in\left\{
\mathbf{e},\mathbf{i},\mathbf{j},\mathbf{k}\right\}  .
\]
This is a finite computation: There are only $64$ triples $\left(
a,b,c\right)  $ with $a,b,c\in\left\{  \mathbf{e},\mathbf{i},\mathbf{j}%
,\mathbf{k}\right\}  $, and we can check the equality $\mu\left(  \mu\left(
a,b\right)  ,c\right)  =\mu\left(  a,\mu\left(  b,c\right)  \right)  $ for
each of these triples directly by computation (using the definition of $\mu$).

A computer could do this in the blink of an eye, but we can also do this by
hand. There are some tricks that reduce our work. The first is to notice that
$\mu\left(  \mu\left(  a,b\right)  ,c\right)  =\mu\left(  a,\mu\left(
b,c\right)  \right)  $ is obvious when one of $a,b,c$ is $\mathbf{e}$ (because
$\mu\left(  \mathbf{x},\mathbf{e}\right)  =\mu\left(  \mathbf{e}%
,\mathbf{x}\right)  =\mathbf{x}$ for each $\mathbf{x}\in\left\{
\mathbf{e},\mathbf{i},\mathbf{j},\mathbf{k}\right\}  $). Thus, it suffices to
prove the equality $\mu\left(  \mu\left(  a,b\right)  ,c\right)  =\mu\left(
a,\mu\left(  b,c\right)  \right)  $ in the case when $a,b,c\in\left\{
\mathbf{i},\mathbf{j},\mathbf{k}\right\}  $. This leaves $27$ triples $\left(
a,b,c\right)  $ to check.

The next trick is to observe a cyclic symmetry. Indeed, the definition of
$\mu$ is invariant under cyclic rotation of $\mathbf{i},\mathbf{j},\mathbf{k}%
$, in the sense that if we replace $\mathbf{i},\mathbf{j},\mathbf{k}$ by
$\mathbf{j},\mathbf{k},\mathbf{i}$ (respectively), then the definition remains
unchanged (for example, $\mu\left(  \mathbf{j},\mathbf{i}\right)
=-\mathbf{k}$ becomes $\mu\left(  \mathbf{k},\mathbf{j}\right)  =-\mathbf{i}%
$). Thus, when we are proving $\mu\left(  \mu\left(  a,b\right)  ,c\right)
=\mu\left(  a,\mu\left(  b,c\right)  \right)  $ for all $a,b,c\in\left\{
\mathbf{i},\mathbf{j},\mathbf{k}\right\}  $, we can WLOG assume that
$a=\mathbf{i}$ (since otherwise, we can achieve this by rotating all of
$a,b,c$ until $a$ becomes $\mathbf{i}$). This leaves $9$ triples $\left(
a,b,c\right)  $ to check.

Let me just check one of them: namely, $\left(  a,b,c\right)  =\left(
\mathbf{i},\mathbf{k},\mathbf{k}\right)  $. In this case, we have%
\begin{align*}
\mu\left(  \mu\left(  a,b\right)  ,c\right)   &  =\mu\left(  \underbrace{\mu
\left(  \mathbf{i},\mathbf{k}\right)  }_{=-\mathbf{j}},\mathbf{k}\right)
=\mu\left(  -\mathbf{j},\mathbf{k}\right)  =-\underbrace{\mu\left(
\mathbf{j},\mathbf{k}\right)  }_{=\mathbf{i}}\ \ \ \ \ \ \ \ \ \ \left(
\text{since }\mu\text{ is bilinear}\right) \\
&  =-\mathbf{i}%
\end{align*}
and%
\begin{align*}
\mu\left(  a,\mu\left(  b,c\right)  \right)   &  =\mu\left(  \mathbf{i}%
,\underbrace{\mu\left(  \mathbf{k},\mathbf{k}\right)  }_{=-\mathbf{e}}\right)
=\mu\left(  \mathbf{i},-\mathbf{e}\right)  =-\underbrace{\mu\left(
\mathbf{i},\mathbf{e}\right)  }_{=\mathbf{i}}\ \ \ \ \ \ \ \ \ \ \left(
\text{since }\mu\text{ is bilinear}\right) \\
&  =-\mathbf{i}.
\end{align*}
Thus, $\mu\left(  \mu\left(  a,b\right)  ,c\right)  =\mu\left(  a,\mu\left(
b,c\right)  \right)  $ is proved for this triple. Similarly, the remaining
$9-1=8$ triples can be checked. Thus, associativity of multiplication is
proved for $\mathbb{H}$.

It remains to prove the neutrality of $1$. In other words, it remains to prove
that $a\cdot\mathbf{e}=\mathbf{e}\cdot a=a$ for each $a\in\mathbb{H}$. Once
again, the bilinearity of the multiplication of $\mathbb{H}$ can be used to
reduce this to the case when $a\in\left\{  \mathbf{e},\mathbf{i}%
,\mathbf{j},\mathbf{k}\right\}  $ (here we need to use Theorem
\ref{thm.mods.linmap-unidef-1} instead of Lemma
\ref{lem.algebras.assoc-on-basis}); but in this case, the claim follows from
our definition of $\mu$. The details are LTTR.

\newpage

\section{\label{chp.polys1}Monoid algebras and polynomials (\cite[Chapter
9]{DumFoo04})}

\begin{convention}
Let $R$ be a \textbf{commutative} ring. (This will be a standing assumption
throughout this chapter.)
\end{convention}

In Section \ref{sec.modules.H}, we have learned how to define an $R$-algebra
the \textquotedblleft slick\textquotedblright\ way: Define an $R$-module
first, and then define an $R$-bilinear map on it, which will serve as the
multiplication of the algebra. Then show that the multiplication is
associative (this is best done \textquotedblleft by
linearity\textquotedblright, i.e., using Lemma
\ref{lem.algebras.assoc-on-basis}) and has a unity (this can again be
simplified using linearity).

I illustrated this method on the example of the ring of quaternions (an
$\mathbb{R}$-algebra).

Now let me apply it to define a more important class of algebras: the monoid
algebras, and, as a particular case, the polynomial rings.

\subsection{\label{sec.polys1.monalg}Monoid algebras}

\subsubsection{Definition}

Recall the notion of a \textbf{monoid}: Roughly speaking, it is a
\textquotedblleft group without inverses\textquotedblright. That is, a
\textbf{monoid} is a triple $\left(  M,\cdot,1\right)  $, where $M$ is a set,
$\cdot$ is an associative binary operation on $M$, and $1$ is an element of
$M$ that is neutral for $\cdot$. We will write $mn$ for $m\cdot n$ whenever
$m,n\in M$. Moreover, the element $mn$ will be called the \textbf{product} of
$m$ and $n$ in the monoid $M$. We will write $M$ for the monoid $\left(
M,\cdot,1\right)  $ if $\cdot$ and $1$ are clear from the context. The monoid
$M$ is said to be \textbf{abelian} if $mn=nm$ for all $m,n\in M$. (This
generalizes the notion of an abelian group.) Given a monoid $\left(
M,\cdot,1\right)  $, the binary operation $\cdot$ is called the
\textbf{operation} of $M$, and the element $1$ is called the \textbf{neutral
element} of $M$. We say that the monoid $M$ is \textbf{written
multiplicatively} (or, for short, \textbf{multiplicative}) when its operation
is denoted by $\cdot$, and we say that it is \textbf{written additively} (or,
for short, \textbf{additive}) when its operation is denoted by $+$. Usually,
the neutral element of a multiplicative monoid is denoted by $1$, whereas the
neutral element of an additive monoid is denoted by $0$.

Here is the informal \textbf{idea} behind the notion of a monoid
algebra\textbf{:} The monoid algebra $R\left[  M\right]  $ of a multiplicative
monoid $M$ is the $R$-algebra obtained by \textquotedblleft
adjoining\textquotedblright\ the monoid $M$ to the ring $R$, which means
\textquotedblleft inserting\textquotedblright\ the elements of $M$
\textquotedblleft into\textquotedblright\ $R$. That is, the algebra $R\left[
M\right]  $ consists of \textquotedblleft formal products\textquotedblright%
\ $rm$ with $r\in R$ and $m\in M$, as well as their formal sums. These
products are multiplied using the multiplications of $R$ and $M$:%
\[
\left(  r_{1}m_{1}\right)  \cdot\left(  r_{2}m_{2}\right)  =\left(  r_{1}%
r_{2}\right)  \cdot\left(  m_{1}m_{2}\right)  .
\]

Let us formalize this:\footnote{We recall that $R$ is a commutative ring.
Furthermore, we recall that if $M$ is any set, then $R^{M}$ is the $R$-module%
\[
\left\{  \left(  r_{i}\right)  _{i\in M}\ \mid\ r_{i}\in R\text{ for each
}i\in M\right\}
\]
(consisting of all families $\left(  r_{i}\right)  _{i\in M}$ of elements of
$R$), whereas $R^{\left(  M\right)  }$ is the $R$-submodule%
\[
\bigoplus_{i\in M}R=\left\{  \left(  r_{i}\right)  _{i\in M}\in R^{M}%
\ \mid\ \text{all but finitely many }i\in M\text{ satisfy }r_{i}=0\right\}
\]
of $R^{M}$. If the set $M$ is finite, then $R^{\left(  M\right)  }=R^{M}$.
\par
The $R$-module $R^{\left(  M\right)  }$ is free, and the \textbf{standard
basis} $\left(  e_{m}\right)  _{m\in M}$ of $R^{\left(  M\right)  }$ is
defined as follows: For each $m\in M$, the vector $e_{m}\in R^{\left(
M\right)  }$ is the family whose $m$-th entry is $1$ and whose all other
entries are $0$. (If $M=\left\{  1,2,\ldots,n\right\}  $ for some
$n\in\mathbb{N}$, then this recovers the classical linear-algebraic standard
basis: e.g., if $M=\left\{  1,2,3\right\}  $, then $e_{1}=\left(
1,0,0\right)  $ and $e_{2}=\left(  0,1,0\right)  $ and $e_{3}=\left(
0,0,1\right)  $.)
\par
The standard basis $\left(  e_{m}\right)  _{m\in M}$ of $R^{\left(  M\right)
}$ is, of course, a basis of $R^{\left(  M\right)  }$.}

\begin{definition}
\label{def.monalg.monalg}Let $M$ be a monoid, written multiplicatively (so
that $\cdot$ denotes its operation, and $1$ denotes its neutral element).

The \textbf{monoid algebra} \textbf{of }$M$ \textbf{over }$R$ (also known as
the \textbf{monoid ring of }$M$ \textbf{over} $R$) is the $R$-algebra
$R\left[  M\right]  $ defined as follows: As an $R$-module, it is the free
$R$-module $R^{\left(  M\right)  }$. Its multiplication is defined to be the
unique $R$-bilinear map $\mu:R^{\left(  M\right)  }\times R^{\left(  M\right)
}\rightarrow R^{\left(  M\right)  }$ that satisfies%
\begin{equation}
\mu\left(  e_{m},e_{n}\right)  =e_{mn}\ \ \ \ \ \ \ \ \ \ \text{for all
}m,n\in M. \label{eq.def.monalg.monalg.mu}%
\end{equation}
Here, $\left(  e_{m}\right)  _{m\in M}$ is the standard basis of $R^{\left(
M\right)  }$ (that is, $e_{m}\in R^{\left(  M\right)  }$ is the family whose
$m$-th entry is $1$ and whose all other entries are $0$). The unity of this
$R$-algebra $R\left[  M\right]  $ is $e_{1}$.
\end{definition}

\begin{theorem}
\label{thm.monalg.welldef}This is indeed a well-defined $R$-algebra.
\end{theorem}

\begin{proof}
Theorem \ref{thm.mods.uniprop-free-2} guarantees that there is a unique
$R$-bilinear map $\mu:R^{\left(  M\right)  }\times R^{\left(  M\right)
}\rightarrow R^{\left(  M\right)  }$ that satisfies
(\ref{eq.def.monalg.monalg.mu}). It remains to prove that the $R$-module
$R^{\left(  M\right)  }$ (equipped with the multiplication $\mu$ and the unity
$e_{1}$) is an $R$-algebra.

All we need to show is that $\mu$ is associative, and that $e_{1}$ is a unity.
I will only show the first statement, and leave the second to you.

We need to show that $\mu\left(  \mu\left(  a,b\right)  ,c\right)  =\mu\left(
a,\mu\left(  b,c\right)  \right)  $ for all $a,b,c\in R\left[  M\right]  $.
According to Lemma \ref{lem.algebras.assoc-on-basis}, it suffices to prove
that%
\[
\mu\left(  \mu\left(  e_{m},e_{n}\right)  ,e_{p}\right)  =\mu\left(  e_{m}%
,\mu\left(  e_{n},e_{p}\right)  \right)  \ \ \ \ \ \ \ \ \ \ \text{for all
}m,n,p\in M
\]
(since the family $\left(  e_{m}\right)  _{m\in M}$ is a basis of $R^{\left(
M\right)  }=R\left[  M\right]  $).

Let us do this: If $m,n,p\in M$, then%
\[
\mu\left(  \underbrace{\mu\left(  e_{m},e_{n}\right)  }_{=e_{mn}}%
,e_{p}\right)  =\mu\left(  e_{mn},e_{p}\right)  =e_{\left(  mn\right)
p}=e_{mnp}%
\]
and similarly $\mu\left(  e_{m},\mu\left(  e_{n},e_{p}\right)  \right)
=e_{mnp}$, so we indeed have $\mu\left(  \mu\left(  e_{m},e_{n}\right)
,e_{p}\right)  =\mu\left(  e_{m},\mu\left(  e_{n},e_{p}\right)  \right)  $ as
desired. This completes the proof that $\mu$ is associative. Thus, Theorem
\ref{thm.monalg.welldef} is proven.
\end{proof}

Since the bilinear map $\mu$ in Definition \ref{def.monalg.monalg} is used as
the multiplication of $R\left[  M\right]  $, we can rewrite the equality
(\ref{eq.def.monalg.monalg.mu}) as follows:%
\begin{equation}
e_{m}\cdot e_{n}=e_{mn}\ \ \ \ \ \ \ \ \ \ \text{for all }m,n\in M.
\label{eq.def.monalg.monalg.dot}%
\end{equation}

From the way we defined monoids, it is clear that every group is a monoid.
Monoid algebras of groups have a special name:

\begin{definition}
When a monoid $M$ is a group, its monoid algebra $R\left[  M\right]  $ is
called a \textbf{group algebra} (or \textbf{group ring}).
\end{definition}

\subsubsection{Examples}

Let me show a few examples of monoid algebras.

\begin{example}
\label{exa.monalg.QC2}Consider the order-$2$ cyclic group $C_{2}=\left\{
1,u\right\}  $ with $u^{2}=1$ (of course, we write $C_{2}$ multiplicatively).
This group is better known as $\mathbb{Z}/2$, and its operation is commonly
written as addition, not as multiplication; but we want to write it
multiplicatively here, in order to match the way $M$ is written in Definition
\ref{def.monalg.monalg}.

What does the group algebra $\mathbb{Q}\left[  C_{2}\right]  $ look like? As a
$\mathbb{Q}$-module (i.e., $\mathbb{Q}$-vector space), it has a basis $\left(
e_{m}\right)  _{m\in C_{2}}=\left(  e_{1},e_{u}\right)  $. Thus, any element
of $\mathbb{Q}\left[  C_{2}\right]  $ can be written as $a\underbrace{e_{1}%
}_{=1}+be_{u}=a+be_{u}$ for some unique $a,b\in\mathbb{Q}$. (As usual, we are
writing $1$ for the unity of our ring $\mathbb{Q}\left[  C_{2}\right]  $,
which is $e_{1}$.)

The multiplication on $\mathbb{Q}\left[  C_{2}\right]  $ is $\mathbb{Q}%
$-bilinear and given on the basis by%
\begin{align*}
e_{1}e_{1}  &  =e_{1\cdot1}=e_{1},\ \ \ \ \ \ \ \ \ \ e_{1}e_{u}=e_{1\cdot
u}=e_{u},\\
e_{u}e_{1}  &  =e_{u\cdot1}=e_{u},\ \ \ \ \ \ \ \ \ \ e_{u}e_{u}=e_{u\cdot
u}=e_{u^{2}}=e_{1}.
\end{align*}

Let us use this to compute some products in $\mathbb{Q}\left[  C_{2}\right]
$:%
\begin{align*}
\left(  3+2e_{u}\right)  \left(  1+2e_{u}\right)   &  =3\cdot1+3\cdot
2e_{u}+2e_{u}\cdot1+2e_{u}\cdot2e_{u}\\
&  =3+6e_{u}+2e_{u}+4\underbrace{e_{u}e_{u}}_{=e_{1}=1}\\
&  =3+6e_{u}+2e_{u}+4=7+8e_{u};\\
\left(  1+e_{u}\right)  ^{2}  &  =1+2e_{u}+\underbrace{e_{u}^{2}}_{=e_{u}%
e_{u}=e_{1}=1}=1+2e_{u}+1=2+2e_{u};\\
\left(  1-e_{u}\right)  \left(  1+e_{u}\right)   &  =1-\underbrace{e_{u}^{2}%
}_{=e_{u}e_{u}=e_{1}=1}\ \ \ \ \ \ \ \ \ \ \left(
\begin{array}
[c]{c}%
\text{since }\left(  1-x\right)  \left(  1+x\right)  =1-x^{2}\\
\text{for any }x\text{ in any ring}%
\end{array}
\right) \\
&  =1-1=0.
\end{align*}
The last of these computations shows that $\mathbb{Q}\left[  C_{2}\right]  $
is not an integral domain. In general, for any $a,b,c,d\in\mathbb{Q}$, we have%
\begin{align}
\left(  a+be_{u}\right)  \left(  c+de_{u}\right)   &  =ac+ade_{u}%
+b\underbrace{e_{u}c}_{\substack{=ce_{u}\\\text{(since the}%
\\\text{multiplication of }\mathbb{Q}\left[  C_{2}\right]  \\\text{is
}\mathbb{Q}\text{-bilinear)}}}+\,b\underbrace{e_{u}d}_{\substack{=de_{u}%
\\\text{(since the}\\\text{multiplication of }\mathbb{Q}\left[  C_{2}\right]
\\\text{is }\mathbb{Q}\text{-bilinear)}}}e_{u}\nonumber\\
&  =ac+ade_{u}+bce_{u}+bd\underbrace{e_{u}e_{u}}_{=e_{1}=1}=ac+ade_{u}%
+bce_{u}+bd\nonumber\\
&  =\left(  ac+bd\right)  +\left(  ad+bc\right)  e_{u}.
\label{eq.exa.monalg.QC2.mul}%
\end{align}

What does $\mathbb{Q}\left[  C_{2}\right]  $ \textquotedblleft look
like\textquotedblright? Meaning, which known $\mathbb{Q}$-algebras is
$\mathbb{Q}\left[  C_{2}\right]  $ isomorphic to (if any)?

I \textbf{claim} that
\begin{equation}
\mathbb{Q}\left[  C_{2}\right]  \cong\mathbb{Q}^{2}=\mathbb{Q}\times
\mathbb{Q}\ \ \ \ \ \ \ \ \ \ \left(  \text{as }\mathbb{Q}\text{-algebras}%
\right)  . \label{eq.exa.monalg.QC2.iso}%
\end{equation}
(See Definition \ref{def.algs.dirprod} for the meaning of $\mathbb{Q}^{2}$ and
$\mathbb{Q}\times\mathbb{Q}$.)

[\textit{Proof of (\ref{eq.exa.monalg.QC2.iso}):} First, we observe that
$\mathbb{Q}\left[  C_{2}\right]  $ is commutative (this is easy to check), and
that the element $z:=\dfrac{1+e_{u}}{2}$ of $\mathbb{Q}\left[  C_{2}\right]  $
is idempotent (since an easy computation shows $z^{2}=z$). Hence, Exercise
\ref{exe.21hw1.3} \textbf{(c)} shows that the map%
\begin{align*}
f:\left(  z\mathbb{Q}\left[  C_{2}\right]  \right)  \times\left(  \left(
1-z\right)  \mathbb{Q}\left[  C_{2}\right]  \right)   &  \rightarrow
\mathbb{Q}\left[  C_{2}\right]  ,\\
\left(  a,b\right)   &  \mapsto a+b
\end{align*}
is a ring isomorphism; thus, this map $f$ is invertible. This map $f$ is
furthermore $\mathbb{Q}$-linear and thus is a $\mathbb{Q}$-algebra morphism.
Since $f$ is invertible, we thus conclude that $f$ is a $\mathbb{Q}$-algebra
isomorphism (by Proposition \ref{prop.algmor.invertible-iso}). Now, what are
$z\mathbb{Q}\left[  C_{2}\right]  $ and $\left(  1-z\right)  \mathbb{Q}\left[
C_{2}\right]  $ ? A general element of $\mathbb{Q}\left[  C_{2}\right]  $ has
the form $a+be_{u}$ for some $a,b\in\mathbb{Q}$. Thus, a general element of
$z\mathbb{Q}\left[  C_{2}\right]  $ has the form $z\left(  a+be_{u}\right)  $
for some $a,b\in\mathbb{Q}$. Since%
\begin{align*}
z\left(  a+be_{u}\right)   &  =\dfrac{1+e_{u}}{2}\left(  a+be_{u}\right)
=\dfrac{1}{2}\underbrace{\left(  1+e_{u}\right)  \left(  a+be_{u}\right)
}_{\substack{=a+e_{u}a+be_{u}+e_{u}be_{u}\\=a+ae_{u}+be_{u}+be_{u}e_{u}}}\\
&  =\dfrac{1}{2}\left(  a+ae_{u}+be_{u}+b\underbrace{e_{u}e_{u}}_{=e_{1}%
=1}\right)  =\dfrac{1}{2}\left(  \left(  a+b\right)  +\left(  a+b\right)
e_{u}\right) \\
&  =\left(  a+b\right)  \cdot\underbrace{\dfrac{1+e_{u}}{2}}_{=z}%
=\underbrace{\left(  a+b\right)  }_{\in\mathbb{Q}}z,
\end{align*}
we see that any such element is a \textbf{scalar} multiple of $z$ (that is, an
element of the form $\lambda z$ for some $\lambda\in\mathbb{Q}$, not just a
multiple of $z$ in the ring $\mathbb{Q}\left[  C_{2}\right]  $). In other
words, any such element belongs to the $\mathbb{Q}$-submodule (= $\mathbb{Q}%
$-vector subspace)%
\[
\mathbb{Q}z:=\left\{  \lambda z\ \mid\ \lambda\in\mathbb{Q}\right\}
\ \ \ \ \ \ \ \ \ \ \text{of }\mathbb{Q}\left[  C_{2}\right]  .
\]
Thus, $z\mathbb{Q}\left[  C_{2}\right]  \subseteq\mathbb{Q}z$. Since we also
have $\mathbb{Q}z\subseteq z\mathbb{Q}\left[  C_{2}\right]  $ (since every
$\lambda\in\mathbb{Q}$ satisfies $\lambda z=z\cdot\left(  \lambda
e_{1}\right)  \in z\mathbb{Q}\left[  C_{2}\right]  $), this entails
$z\mathbb{Q}\left[  C_{2}\right]  =\mathbb{Q}z$. Hence, in particular,
$\mathbb{Q}z$ is a $\mathbb{Q}$-algebra with unity $z$. However, the map
\[
\mathbb{Q}\rightarrow\mathbb{Q}z,\ \lambda\mapsto\lambda z
\]
is a $\mathbb{Q}$-algebra morphism (indeed, it is clearly $\mathbb{Q}$-linear;
it respects multiplication since $\left(  \lambda z\right)  \left(  \mu
z\right)  =\lambda\mu\underbrace{z^{2}}_{=z}=\lambda\mu z$ for any
$\lambda,\mu\in\mathbb{Q}$; it respects the unity since $1z=z$ is the unity of
$\mathbb{Q}z$), and thus is a $\mathbb{Q}$-algebra isomorphism (since it is
easily seen to be bijective). Thus, $\mathbb{Q}z\cong\mathbb{Q}$ as
$\mathbb{Q}$-algebras. Combining this with $z\mathbb{Q}\left[  C_{2}\right]
=\mathbb{Q}z$, we obtain $z\mathbb{Q}\left[  C_{2}\right]  =\mathbb{Q}%
z\cong\mathbb{Q}$ as $\mathbb{Q}$-algebras. Similarly, we can prove that
$\left(  1-z\right)  \mathbb{Q}\left[  C_{2}\right]  \cong\mathbb{Q}$ (indeed,
a simple computation shows that $1-z=\dfrac{1-e_{u}}{2}$, and thus we can
mostly repeat our above argument with $1-z$ instead of $z$, with the main
difference being that some plus signs become minus signs).

So the isomorphism $f$ results in%
\[
\mathbb{Q}\left[  C_{2}\right]  \cong\underbrace{\left(  z\mathbb{Q}\left[
C_{2}\right]  \right)  }_{\cong\mathbb{Q}}\times\underbrace{\left(  \left(
1-z\right)  \mathbb{Q}\left[  C_{2}\right]  \right)  }_{\cong\mathbb{Q}}%
\cong\mathbb{Q}\times\mathbb{Q}=\mathbb{Q}^{2}.
\]
This proves (\ref{eq.exa.monalg.QC2.iso}).]

Retracing our proof of (\ref{eq.exa.monalg.QC2.iso}), we actually get an
explicit $\mathbb{Q}$-algebra isomorphism%
\begin{align*}
\mathbb{Q}^{2}  &  \rightarrow\mathbb{Q}\left[  C_{2}\right]  ,\\
\left(  \lambda,\mu\right)   &  \mapsto f\left(  \lambda z,\mu\left(
1-z\right)  \right)  =\lambda z+\mu\left(  1-z\right)  =\lambda\cdot
\dfrac{1+e_{u}}{2}+\mu\cdot\dfrac{1-e_{u}}{2}\\
&  \ \ \ \ \ \ \ \ \ \ \ \ \ \ \ \ \ \ \ \ \ \ \ \ \ \ \ \ \ \ =\dfrac
{\lambda+\mu}{2}+\dfrac{\lambda-\mu}{2}e_{u}.
\end{align*}
The inverse of this isomorphism is the $\mathbb{Q}$-algebra isomorphism%
\begin{align*}
\mathbb{Q}\left[  C_{2}\right]   &  \rightarrow\mathbb{Q}^{2},\\
a+be_{u}  &  \mapsto\left(  a+b,\ a-b\right)  \ \ \ \ \ \ \ \ \ \ \left(
\text{for all }a,b\in\mathbb{Q}\right)  .
\end{align*}

Note that there are only two $\mathbb{Q}$-algebra isomorphisms from
$\mathbb{Q}\left[  C_{2}\right]  $ to $\mathbb{Q}^{2}$: One is the one we just
constructed; the other is%
\begin{align*}
\mathbb{Q}\left[  C_{2}\right]   &  \rightarrow\mathbb{Q}^{2},\\
a+be_{u}  &  \mapsto\left(  a-b,\ a+b\right)  \ \ \ \ \ \ \ \ \ \ \left(
\text{for all }a,b\in\mathbb{Q}\right)
\end{align*}
(which differs from the first only in $a+b$ and $a-b$ being swapped). In
contrast, there are infinitely many $\mathbb{Q}$\textbf{-module} isomorphisms
from $\mathbb{Q}\left[  C_{2}\right]  $ to $\mathbb{Q}^{2}$; the simplest one
just sends each $a+be_{u}$ to $\left(  a,b\right)  $ (for all $a,b\in
\mathbb{Q}$).
\end{example}

\begin{example}
\ \ 

\begin{enumerate}
\item[\textbf{(a)}] We can easily repeat Example \ref{exa.monalg.QC2} using
the field $\mathbb{R}$ (or $\mathbb{C}$) instead of $\mathbb{Q}$. Everything
works just as it did for $\mathbb{Q}$. For example, we get an $\mathbb{R}%
$-algebra isomorphism $\mathbb{R}^{2}\rightarrow\mathbb{R}\left[
C_{2}\right]  $.

\item[\textbf{(b)}] Now, let us try to repeat Example \ref{exa.monalg.QC2}
using the ring $\mathbb{Z}$ instead of $\mathbb{Q}$. The multiplication rule
(\ref{eq.exa.monalg.QC2.mul}) still holds (but now for $a,b,c,d\in\mathbb{Z}%
$). What about the isomorphism (\ref{eq.exa.monalg.QC2.iso}) ? The idempotent
$z$ no longer exists (since we had to divide by $2$ to construct it, but we
cannot divide by $2$ in $\mathbb{Z}$), so our proof of
(\ref{eq.exa.monalg.QC2.iso}) does not work. And indeed,
(\ref{eq.exa.monalg.QC2.iso}) does not hold for $\mathbb{Z}$. The $\mathbb{Z}%
$-algebra%
\[
\mathbb{Z}\left[  C_{2}\right]  =\left\{  a+be_{u}\ \mid\ a,b\in
\mathbb{Z}\right\}
\]
is \textbf{not} isomorphic to any direct product of nontrivial $\mathbb{Z}%
$-algebras. This can be proved by showing that $\mathbb{Z}\left[
C_{2}\right]  $ has no idempotents other than $0$ and $1$. (In fact, if
$a+be_{u}\in\mathbb{Z}\left[  C_{2}\right]  $ is an idempotent, then $\left(
a+be_{u}\right)  ^{2}=a+be_{u}$. But (\ref{eq.exa.monalg.QC2.mul}) yields
$\left(  a+be_{u}\right)  ^{2}=\left(  a^{2}+b^{2}\right)  +2abe_{u}$, so this
idempotency results in $\left(  a^{2}+b^{2}\right)  +2abe_{u}=a+be_{u}$, and
thus $a^{2}+b^{2}=a$ and $2ab=b$ (since $e_{1}=1$ and $e_{u}$ are $\mathbb{Z}%
$-linearly independent). But the only integer solutions $\left(  a,b\right)  $
of this system of two equations are $\left(  0,0\right)  $ and $\left(
1,0\right)  $ (check this!); thus, the only idempotents of $\mathbb{Z}\left[
C_{2}\right]  $ are $0+0e_{u}=0$ and $1+0e_{u}=1$.)
\end{enumerate}
\end{example}

\begin{example}
\label{exa.monalg.QC3}Now, let us take the order-$3$ cyclic group
$C_{3}=\left\{  1,u,v\right\}  $ with $u^{3}=1$ and $v=u^{2}$. (Again, this
group is better known as $\mathbb{Z}/3$, but we write it multiplicatively.)
Then, $\mathbb{Q}\left[  C_{3}\right]  $ is again commutative, and has an
idempotent $z:=\dfrac{1+e_{u}+e_{v}}{3}$; this leads to a $\mathbb{Q}$-algebra
isomorphism%
\[
\mathbb{Q}\left[  C_{3}\right]  \cong\mathbb{Q}\times S,
\]
where the $\mathbb{Q}$ factor is
\[
z\mathbb{Q}\left[  C_{3}\right]  =\mathbb{Q}z=\left\{  a+ae_{u}+ae_{v}%
\ \mid\ a\in\mathbb{Q}\right\}
\]
and where the $S$ factor is%
\[
\left(  1-z\right)  \mathbb{Q}\left[  C_{3}\right]  =\left\{  a+be_{u}%
+ce_{v}\ \mid\ a+b+c=0\right\}  .
\]
The $\mathbb{Q}$ factor is $1$-dimensional (as a $\mathbb{Q}$-vector space),
while the $S$ factor is $2$-dimensional. Can $S$ be decomposed further? What
does $S$ \textquotedblleft look like\textquotedblright? We will later see (in
Exercise \ref{exe.monalg.QC3.fieldext-S} below).
\end{example}

\begin{example}
\label{exa.monalg.H'}Here is a \textbf{non-example}: The ring of quaternions
$\mathbb{H}$ is \textbf{not} a monoid algebra. It is pretty close, in that it
has a basis $\left(  1,i,j,k\right)  $ (over $\mathbb{R}$) with the property
that the product of any two basis elements is either a basis element again
(for example, $ij=k$) or the negative of a basis element (for example,
$ji=-k$). However, for it to be a monoid algebra, it would need a basis such
that the product of any two basis elements is always a basis element (never
the negative of a basis element).\footnotemark\ Such a basis does not exist
for $\mathbb{H}$.

If we remove all the minus signs in the definition of $\mathbb{H}$ (that is,
we replace the multiplication rules by $i^{2}=j^{2}=k^{2}=1$ and $ij=ji=k$ and
$jk=kj=i$ and $ki=ik=j$), then we actually do obtain a monoid algebra (namely,
the group algebra of \href{https://en.wikipedia.org/wiki/Klein_four-group}{the
Klein four-group}).

We can find another group algebra closely related to $\mathbb{H}$. Indeed, we
define the \textbf{quaternion group} $Q_{8}$ to be the subgroup $\left\{
1,i,j,k,-1,-i,-j,-k\right\}  $ of the group of units of $\mathbb{H}$. Then,
consider the group algebra $\mathbb{H}^{\prime}:=\mathbb{R}\left[
Q_{8}\right]  $ of this group $Q_{8}$. This group algebra $\mathbb{H}^{\prime
}$ is $8$-dimensional as an $\mathbb{R}$-vector space, whereas $\mathbb{H}$ is
$4$-dimensional; thus, $\mathbb{H}^{\prime}$ is not quite $\mathbb{H}$ (but
rather close). The main difference between $\mathbb{H}$ and $\mathbb{H}%
^{\prime}$ is that the elements $e_{1}$ and $e_{-1}$ of $\mathbb{H}^{\prime}$
are two different basis elements (thus linearly independent), whereas the
elements $1$ and $-1$ of $\mathbb{H}$ are negatives of each other. Even though
$\mathbb{H}^{\prime}$ is not commutative, we can define a principal ideal
$\left(  e_{1}+e_{-1}\right)  \mathbb{H}^{\prime}$ of $\mathbb{H}^{\prime}$
(by Exercise \ref{exe.ideal.princid-central}, since $e_{1}+e_{-1}$ is a
central element of $\mathbb{H}^{\prime}$), and then it is not hard to show
that the quotient ring $\mathbb{H}^{\prime}/\left(  e_{1}+e_{-1}\right)
\mathbb{H}^{\prime}$ is isomorphic to $\mathbb{H}$. Thus, while $\mathbb{H}$
itself is not a group ring, we can obtain $\mathbb{H}$ from the group ring
$\mathbb{H}^{\prime}=\mathbb{R}\left[  Q_{8}\right]  $ by \textquotedblleft
setting $e_{-1}$ equal to the negative of $e_{1}$\textquotedblright\ (that is,
quotienting out the principal ideal generated by $e_{1}+e_{-1}$).
\end{example}

\footnotetext{In general, you can describe a monoid algebra as an algebra that
has a basis that contains the unity (i.e., the unity of the algebra belongs to
the basis) and is closed under multiplication (i.e., the product of any two
basis elements is again a basis element).}

\begin{exercise}
Let $R$ be a commutative ring, and $n$ be a positive integer. In Exercise
\ref{exe.circulant.1}, we have defined the ring $\operatorname*{Circ}%
\nolimits_{n}$ of circulant $n\times n$-matrices $A\in R^{n\times n}$.

\begin{enumerate}
\item[\textbf{(a)}] Prove that this ring $\operatorname*{Circ}\nolimits_{n}$
is actually an $R$-subalgebra of $R^{n\times n}$.

\item[\textbf{(b)}] Prove that this $R$-algebra $\operatorname*{Circ}%
\nolimits_{n}$ is isomorphic to the group algebra $R\left[  C_{n}\right]  $,
where $C_{n}$ is the order-$n$ cyclic group.
\end{enumerate}
\end{exercise}

\begin{exercise}
Let $R$ be a commutative ring. Let $G$ be a finite group. Let $s$ be the
element $\sum_{g\in G}e_{g}$ of the group algebra $R\left[  G\right]  $.

\begin{enumerate}
\item[\textbf{(a)}] Prove that $s^{2}=\left\vert G\right\vert \cdot s$.

\item[\textbf{(b)}] If $\left\vert G\right\vert \cdot1_{R}$ is invertible in
$R$, then prove that $\dfrac{1}{\left\vert G\right\vert }s\in R\left[
G\right]  $ is idempotent.

[This generalizes the idempotents $z$ in Example \ref{exa.monalg.QC2} and
Example \ref{exa.monalg.QC3}.]
\end{enumerate}
\end{exercise}

\begin{exercise}
\ \ 

\begin{enumerate}
\item[\textbf{(a)}] Prove the claim made in Example \ref{exa.monalg.QC2},
saying that there are only two $\mathbb{Q}$-algebra isomorphisms from
$\mathbb{Q}\left[  C_{2}\right]  $ to $\mathbb{Q}^{2}$.

\item[\textbf{(b)}] More generally: Let $R$ be any integral domain. Let
$n\in\mathbb{N}$. Prove that there are exactly $n!$ many $R$-algebra
isomorphisms from $R^{n}$ to $R^{n}$, and each of them has the form%
\begin{align*}
R^{n}  &  \rightarrow R^{n},\\
\left(  r_{1},r_{2},\ldots,r_{n}\right)   &  \mapsto\left(  r_{\sigma\left(
1\right)  },r_{\sigma\left(  2\right)  },\ldots,r_{\sigma\left(  n\right)
}\right)
\end{align*}
for some permutation $\sigma$ of $\left\{  1,2,\ldots,n\right\}  $. (In other
words, each of these isomorphisms just permutes the entries of the $n$-tuple.)
\end{enumerate}

[\textbf{Hint:} For part \textbf{(b)}, let $f$ be an $R$-algebra isomorphism
$R^{n}\rightarrow R^{n}$. Consider the standard basis vectors $e_{1}%
,e_{2},\ldots,e_{n}$ of $R^{n}$. Form the $n\times n$-matrix whose rows are
the $n$ vectors $f\left(  e_{1}\right)  ,f\left(  e_{2}\right)  ,\ldots
,f\left(  e_{n}\right)  $. Show that each row of this matrix has at least one
nonzero entry, but each column has at most one nonzero entry. Furthermore,
what must these entries be?]
\end{exercise}

\subsubsection{Pretending that the elements of $M$ belong to $R\left[
M\right]  $}

\begin{convention}
\label{conv.monalg.m-as-em}Let $R$ be a commutative ring. Let $M$ be a monoid.
The elements $e_{m}$ of the standard basis $\left(  e_{m}\right)  _{m\in M}$
of $R\left[  M\right]  $ will often be just denoted by $m$ (by abuse of
notation). Thus, for example, the element $a+be_{u}+ce_{v}$ of $\mathbb{Q}%
\left[  C_{3}\right]  $ (from Example \ref{exa.monalg.QC3}) will be written as
$a+bu+cv$. With this notation, an element of $R\left[  M\right]  $ is (at
least notationally) really just a sum of products of elements of $R$ with
elements of $M$.

Do \textbf{not} use this convention when it can create confusion! In
particular, do not use it when some elements of $M$ include minus (or plus)
signs, such as the elements $-1,-i,-j,-k$ in Example \ref{exa.monalg.H'}.
(Indeed, in Example \ref{exa.monalg.H'}, it is crucial that $e_{1}$ and
$e_{-1}$ are two different basis elements of $\mathbb{H}^{\prime}$, not
negatives of each other. Denoting them by $1$ and $-1$ would obscure this and
risk confusing the nonzero element $e_{1}+e_{-1}$ for the zero sum $1+\left(
-1\right)  =0$.)
\end{convention}

\subsubsection{General properties of monoid algebras}

We shall now state some general properties of monoid algebras, and agree on
some conventions that will make it easier to work in those algebras.

\begin{proposition}
\label{prop.monalg.comcm}Let $M$ be an \textbf{abelian} monoid. Then, the
monoid ring $R\left[  M\right]  $ is commutative.
\end{proposition}

\begin{proof}
We must prove that $ab=ba$ for all $a,b\in R\left[  M\right]  $. This is a
typical linearity argument (just as the proof of Lemma
\ref{lem.algebras.assoc-on-basis}): Since $\left(  e_{m}\right)  _{m\in M}$ is
a basis of the $R$-module $R\left[  M\right]  $, we can write $a$ and $b$ as
$R$-linear combinations of this family $\left(  e_{m}\right)  _{m\in M}$. That
is, there exist scalars $a_{m}\in R$ and $b_{m}\in R$ for all $m\in M$ such
that%
\[
a=\sum_{m\in M}a_{m}e_{m}\ \ \ \ \ \ \ \ \ \ \text{and}%
\ \ \ \ \ \ \ \ \ \ b=\sum_{m\in M}b_{m}e_{m}%
\]
(and such that $a_{m}=0$ for all but finitely many $m\in M$, and likewise for
the $b_{m}$). Multiplying these two equalities, we find%
\begin{align*}
ab  &  =\left(  \sum_{m\in M}a_{m}e_{m}\right)  \left(  \sum_{m\in M}%
b_{m}e_{m}\right)  =\left(  \sum_{m\in M}a_{m}e_{m}\right)  \left(  \sum_{n\in
M}b_{n}e_{n}\right) \\
&  \ \ \ \ \ \ \ \ \ \ \ \ \ \ \ \ \ \ \ \ \left(  \text{here, we renamed
}m\text{ as }n\text{ in the second sum}\right) \\
&  =\sum_{m\in M}\ \ \sum_{n\in M}a_{m}b_{n}\underbrace{e_{m}e_{n}%
}_{\substack{=e_{mn}\\\text{(by (\ref{eq.def.monalg.monalg.dot}))}}}\\
&  \ \ \ \ \ \ \ \ \ \ \ \ \ \ \ \ \ \ \ \ \left(  \text{since the
multiplication of the }R\text{-algebra }R\left[  M\right]  \text{ is
}R\text{-bilinear}\right) \\
&  =\sum_{m\in M}\ \ \sum_{n\in M}a_{m}b_{n}\underbrace{e_{mn}}%
_{\substack{=e_{nm}\\\text{(since }M\text{ is abelian,}\\\text{so that
}mn=nm\text{)}}}=\sum_{m\in M}\ \ \sum_{n\in M}a_{m}b_{n}e_{nm}%
\end{align*}
and (if we multiply them in the opposite order)%
\begin{align*}
ba  &  =\left(  \sum_{m\in M}b_{m}e_{m}\right)  \left(  \sum_{m\in M}%
a_{m}e_{m}\right)  =\left(  \sum_{n\in M}b_{n}e_{n}\right)  \left(  \sum_{m\in
M}a_{m}e_{m}\right) \\
&  \ \ \ \ \ \ \ \ \ \ \ \ \ \ \ \ \ \ \ \ \left(  \text{here, we renamed
}m\text{ as }n\text{ in the first sum}\right) \\
&  =\underbrace{\sum_{n\in M}\ \ \sum_{m\in M}}_{=\sum_{m\in M}\ \ \sum_{n\in
M}}\underbrace{b_{n}a_{m}}_{\substack{=a_{m}b_{n}\\\text{(since }R\text{ is
commutative)}}}\underbrace{e_{n}e_{m}}_{\substack{=e_{nm}\\\text{(by
(\ref{eq.def.monalg.monalg.dot}))}}}\\
&  \ \ \ \ \ \ \ \ \ \ \ \ \ \ \ \ \ \ \ \ \left(  \text{since the
multiplication of the }R\text{-algebra }R\left[  M\right]  \text{ is
}R\text{-bilinear}\right) \\
&  =\sum_{m\in M}\ \ \sum_{n\in M}a_{m}b_{n}e_{nm}.
\end{align*}
The right hand sides of these two equalities are equal; thus, so are the left
hand sides. In other words, $ab=ba$. This completes the proof of Proposition
\ref{prop.monalg.comcm}.
\end{proof}

\begin{proposition}
\label{prop.monalg.const}Let $M$ be a monoid with neutral element $1$. Then,
the map%
\begin{align*}
R  &  \rightarrow R\left[  M\right]  ,\\
r  &  \mapsto r\cdot e_{1}%
\end{align*}
is an injective $R$-algebra morphism.
\end{proposition}

\begin{proof}
First of all, this map is clearly injective, because the family $\left(
e_{m}\right)  _{m\in M}$ is a basis of $R\left[  M\right]  $ and thus is
$R$-linearly independent (so $r\cdot e_{1}\neq s\cdot e_{1}$ for any two
distinct $r,s\in R$). It remains to prove that this map is an $R$-algebra
morphism. But this is a particular case of the following general fact: If $A$
is an $R$-algebra, then the map%
\begin{align*}
R  &  \rightarrow A,\\
r  &  \mapsto r\cdot1_{A}%
\end{align*}
is an $R$-algebra morphism. This fact is easy to show (for example, the map
respects multiplication, since any $r,s\in R$ satisfy $\left(  r\cdot
1_{A}\right)  \cdot\left(  s\cdot1_{A}\right)  =rs\cdot1_{A}\cdot1_{A}%
=rs\cdot1_{A}$), and we can apply it to $A=R\left[  M\right]  $ (recalling
that $1_{R\left[  M\right]  }=e_{1}$) to obtain precisely the claim we are
trying to prove.
\end{proof}

\begin{convention}
\label{conv.monalg.const}If $M$ is a monoid, then we shall identify each $r\in
R$ with $r\cdot e_{1}\in R\left[  M\right]  $. This identification is
harmless\footnotemark, and turns $R$ into an $R$-subalgebra of $R\left[
M\right]  $.

An element of $R\left[  M\right]  $ will be called \textbf{constant} if it
lies in this subalgebra (i.e., if it is of the form $r\cdot e_{1}$ for some
$r\in R$). Thus, we have identified each constant element of $R\left[
M\right]  $ with the corresponding element of $R$.
\end{convention}

\footnotetext{Indeed, Proposition \ref{prop.monalg.const} shows that the map
$R\rightarrow R\left[  M\right]  $ sending each $r\in R$ to $r\cdot e_{1}$ is
an injective $R$-algebra morphism. Thus, this map keeps distinct elements of
$R$ distinct in $R\left[  M\right]  $ (since it is injective), and respects
addition and multiplication (since it is an $R$-algebra morphism).}

A \textbf{warning} might be in order: In Example \ref{exa.monalg.QC2}, we have
seen that $\mathbb{Q}\left[  C_{2}\right]  \cong\mathbb{Q}\times\mathbb{Q}$ as
$\mathbb{Q}$-algebras. Now, in Convention \ref{conv.monalg.const}, we have
identified $\mathbb{Q}$ with a $\mathbb{Q}$-subalgebra of $\mathbb{Q}\left[
C_{2}\right]  $. But this subalgebra is not one of the two $\mathbb{Q}$
factors in $\mathbb{Q}\left[  C_{2}\right]  \cong\mathbb{Q}\times\mathbb{Q}$.
Indeed, as a $\mathbb{Q}$-subalgebra, it contains the unity of $\mathbb{Q}%
\left[  C_{2}\right]  $, but none of the two $\mathbb{Q}$ factors does.

\begin{proposition}
\label{prop.monalg.M-into-RM}Let $M$ be a monoid. Then, the map
\begin{align*}
M  &  \rightarrow R\left[  M\right]  ,\\
m  &  \mapsto e_{m}%
\end{align*}
is a monoid morphism from $M$ to $\left(  R\left[  M\right]  ,\cdot,1\right)
$.
\end{proposition}

\begin{proof}
This map respects multiplication (because of (\ref{eq.def.monalg.monalg.dot}))
and sends the neutral element of $M$ to the unity of $R\left[  M\right]  $
(since $e_{1}$ is the unity of $R\left[  M\right]  $). Thus, it is a monoid morphism.
\end{proof}

Note that if we use Convention \ref{conv.monalg.m-as-em}, then the
\textquotedblleft$m\mapsto e_{m}$\textquotedblright\ in Proposition
\ref{prop.monalg.M-into-RM} can be rewritten as \textquotedblleft$m\mapsto
m$\textquotedblright, so the map from Proposition \ref{prop.monalg.M-into-RM}
looks like an inclusion map. This is merely an artefact of our notation. In
truth, the element $m$ of the monoid $M$ is not literally the same as the
corresponding basis element $e_{m}$ of the monoid algebra $R\left[  M\right]
$; we have just agreed to call both of them $m$ for brevity. But Proposition
\ref{prop.monalg.M-into-RM} shows that using the same letter for these two
elements is a mostly harmless abuse of notation. The only possible problem it
can cause is when the map in Proposition \ref{prop.monalg.M-into-RM} fails to
be injective, so we might accidentally equate two distinct elements $m,n$ of
$M$ whose corresponding basis elements $e_{m}$ and $e_{n}$ are equal.
Fortunately, this can only happen if the ring $R$ is trivial (indeed, for any
nontrivial ring $R$, the basis elements $e_{m}$ for $m\in M$ are distinct),
and this is not a very interesting case. (This is also an issue that rarely
comes up in practice. The purpose of Convention \ref{conv.monalg.m-as-em} is
to simplify computations in $R\left[  M\right]  $, not to \textquotedblleft
pull\textquotedblright\ them back into $M$.)

\begin{exercise}
Let $n$ be a positive integer. Consider the symmetric group $S_{n}$ -- that
is, the group of all permutations of the set $\left\{  1,2,\ldots,n\right\}  $.

For any two distinct elements $i$ and $j$ of $\left\{  1,2,\ldots,n\right\}
$, let $t_{i,j}$ be the permutation in $S_{n}$ that swaps $i$ with $j$ while
leaving the remaining elements of $\left\{  1,2,\ldots,n\right\}  $ unchanged.
(This is called a \textbf{transposition}.)

For each $i\in\left\{  1,2,\ldots,n\right\}  $, define an element $Y_{i}%
\in\mathbb{Z}\left[  S_{n}\right]  $ of the group algebra $\mathbb{Z}\left[
S_{n}\right]  $ by%
\[
Y_{i}:=\sum_{j=1}^{i-1}t_{i,j}=t_{i,1}+t_{i,2}+\cdots+t_{i,i-1}%
\]
(where we are using Convention \ref{conv.monalg.m-as-em}, so that $t_{i,j}$
really means $e_{t_{i,j}}$). (Thus, $Y_{1}=0$ and $Y_{2}=t_{2,1}$ and
$Y_{3}=t_{3,1}+t_{3,2}$ and so on.) The $n$ elements $Y_{1},Y_{2},\ldots
,Y_{n}$ are called the \textbf{Young-Jucys-Murphy elements} of $\mathbb{Z}%
\left[  S_{n}\right]  $.

\begin{enumerate}
\item[\textbf{(a)}] Prove that the $n$ elements $Y_{1},Y_{2},\ldots,Y_{n}$
commute (i.e., that we have $Y_{i}Y_{j}=Y_{j}Y_{i}$ for all $i,j\in\left\{
1,2,\ldots,n\right\}  $).

\item[\textbf{(b)}] Prove that the element $Y_{1}+Y_{2}+\cdots+Y_{n}%
=\sum_{1\leq j<i\leq n}t_{i,j}$ belongs to the center of $\mathbb{Z}\left[
S_{n}\right]  $.
\end{enumerate}
\end{exercise}

\begin{exercise}
Let $G$ be a finite group. Let $R$ be a nontrivial commutative ring.

Let $T$ be a subset of $G$. Let $s_{T}$ be the element $\sum_{t\in T}t$ of the
group algebra $R\left[  G\right]  $ (where we use Convention
\ref{conv.monalg.m-as-em}, so that $t$ means $e_{t}$).

Prove that the element $s_{T}$ of $R\left[  G\right]  $ is central if and only
if $T$ is a union of conjugacy classes of $G$ (that is, if every element of
$G$ that is conjugate to an element of $T$ must itself be an element of $T$).
\end{exercise}

\subsection{\label{sec.polys1.polyrings}Polynomial rings}

\subsubsection{\label{subsec.polys1.polyrings.univar}Univariate polynomials}

Now, we can effortlessly define polynomial rings. Recall that $R$ is a
commutative ring. Recall also that $\mathbb{N}=\left\{  0,1,2,\ldots\right\}
$ (so $0\in\mathbb{N}$).

\begin{definition}
\label{def.polring.univar}Let $C$ be the free monoid with a single generator
$x$. This is the monoid whose elements are countably many distinct symbols
named%
\[
x^{0},\ \ x^{1},\ \ x^{2},\ \ x^{3},\ \ \ldots,
\]
and whose operation is defined by%
\[
x^{i}\cdot x^{j}=x^{i+j}\ \ \ \ \ \ \ \ \ \ \text{for all }i,j\in\mathbb{N}.
\]
We write this monoid multiplicatively, but of course it is just the well-known
monoid $\left(  \mathbb{N},+,0\right)  $ in new clothes (we have renamed each
$i\in\mathbb{N}$ as $x^{i}$; we have renamed addition as multiplication). Its
neutral element is $x^{0}$. We set $x=x^{1}$ (so that $x^{i}$ really is the
$i$-th power of $x$).

The elements of $C$ are called \textbf{monomials} (in the variable $x$). The
specific element $x$ is called the \textbf{indeterminate} (or, somewhat
imprecisely, the \textbf{variable}).

Now, the \textbf{univariate polynomial ring} over $R$ is defined to be the
monoid algebra $R\left[  C\right]  $. It is commonly denoted by $R\left[
x\right]  $. Following Convention \ref{conv.monalg.m-as-em}, we simply write
$m$ for $e_{m}$ when $m\in C$ (that is, we write $x^{i}$ for the basis element
$e_{x^{i}}$); thus, $R\left[  x\right]  $ is a free $R$-module with basis%
\[
\left(  x^{0},x^{1},x^{2},x^{3},\ldots\right)  =\left(  1,x,x^{2},x^{3}%
,\ldots\right)  .
\]
Hence, any $p\in R\left[  x\right]  $ can be written as a finite $R$-linear
combination of powers of $x$. In other words, any $p\in R\left[  x\right]  $
can be written in the form%
\[
p=a_{0}x^{0}+a_{1}x^{1}+a_{2}x^{2}+\cdots+a_{n}x^{n}=a_{0}+a_{1}x+a_{2}%
x^{2}+\cdots+a_{n}x^{n}%
\]
for some $n\in\mathbb{N}$ and some $a_{0},a_{1},\ldots,a_{n}\in R$. This
representation is unique up to trailing zeroes (meaning that we can always add
$0x^{n+1}$ addends -- e.g., rewriting $4x^{0}+3x^{1}$ as $4x^{0}+3x^{1}%
+0x^{2}$ --, but other than that it is unique).

Elements of $R\left[  x\right]  $ are called \textbf{polynomials} in $x$ over
$R$.
\end{definition}

Thus, up to notation, the univariate polynomial ring $R\left[  x\right]  $ is
just the monoid ring $R\left[  \mathbb{N}\right]  $ of the abelian monoid
$\mathbb{N}=\left(  \mathbb{N},+,0\right)  $. Hence, this ring $R\left[
x\right]  $ is commutative (by Proposition \ref{prop.monalg.comcm}, since the
monoid $\mathbb{N}$ is abelian). The unity of the $R$-algebra $R\left[
x\right]  $ is $x^{0}=1$.

\begin{example}
\textbf{\ \ }

\begin{enumerate}
\item[\textbf{(a)}] Here is an example of a polynomial:
\[
1+3x^{2}+6x^{3}=1e_{x^{0}}+3e_{x^{2}}+6e_{x^{3}}\in R\left[  x\right]  .
\]

\item[\textbf{(b)}] A non-example: The infinite sum $1+x+x^{2}+x^{3}+\cdots$
is \textbf{not} in $R\left[  x\right]  $. Indeed, polynomials are linear
combinations of powers of $x$, and linear combinations are finite (by
definition); even if you write them as infinite sums, they are de-facto finite
because all but finitely many addends are $0$. Infinite sums such as
$1+x+x^{2}+x^{3}+\cdots$ thus are not polynomials\footnotemark; they are
instead known as \textbf{formal power series}. There is a way to define an
$R$-algebra of formal power series, too, but we won't do so now.
\end{enumerate}
\end{example}

\footnotetext{You might wonder whether such sums are well-defined in the first
place. Yes, they are, if one correctly defines them. For a complete
definition, see (e.g.) \cite[\S 3.2.2]{21s}.}So we have defined
\textbf{univariate} polynomial rings (i.e., polynomial rings in a single
variable). Likewise, we can define \textbf{multivariate} polynomial rings --
i.e., polynomial rings in several variables. For simplicity, let me restrict
myself to finitely many variables.

\subsubsection{\label{subsec.polys1.polyrings.bivar}Bivariate polynomials}

To avoid a barrage of new notations, let me first introduce \textbf{bivariate}
polynomial rings -- i.e., polynomial rings in two variables. Their definition
is just a particular case of the definition of multivariate polynomial rings
(Definition \ref{def.polring.mulvar}) that we will give soon after, but it is
somewhat easier to understand as it involves less complicated notations.

\begin{definition}
\label{def.polring.bivar}Let $C^{\left(  2\right)  }$ be the free abelian
monoid with two generators $x$ and $y$. This is the monoid whose elements are
the distinct symbols%
\begin{align*}
&  x^{0}y^{0},\ \ x^{0}y^{1},\ \ x^{0}y^{2},\ \ x^{0}y^{3},\ \ \ldots,\\
&  x^{1}y^{0},\ \ x^{1}y^{1},\ \ x^{1}y^{2},\ \ x^{1}y^{3},\ \ \ldots,\\
&  x^{2}y^{0},\ \ x^{2}y^{1},\ \ x^{2}y^{2},\ \ x^{2}y^{3},\ \ \ldots,\\
&  \ldots,
\end{align*}
that is, the distinct symbols%
\[
x^{a}y^{b}\ \ \ \ \ \ \ \ \ \ \text{with }a\in\mathbb{N}\text{ and }%
b\in\mathbb{N},
\]
and whose operation is defined by%
\[
\left(  x^{a}y^{b}\right)  \cdot\left(  x^{c}y^{d}\right)  =x^{a+c}%
y^{b+d}\ \ \ \ \ \ \ \ \ \ \text{for all }a,b,c,d\in\mathbb{N}.
\]
We write this monoid multiplicatively, but of course it is just the monoid
$\mathbb{N}^{2}=\left(  \mathbb{N}^{2},+,0\right)  $ in disguise (where the
addition on $\mathbb{N}^{2}$ that we are calling \textquotedblleft%
$+$\textquotedblright\ here is entrywise, and $0$ means the pair $\left(
0,0\right)  $), with each element $\left(  a,b\right)  $ renamed as
$x^{a}y^{b}$ and with addition renamed as multiplication. The elements of
$C^{\left(  2\right)  }$ are called \textbf{monomials} in $x$ and $y$. We
define two specific monomials $x$ and $y$ by%
\[
x=x^{1}y^{0}\ \ \ \ \ \ \ \ \ \ \text{and}\ \ \ \ \ \ \ \ \ \ y=x^{0}y^{1}.
\]
These two monomials $x$ and $y$ are called the \textbf{indeterminates} (or,
somewhat imprecisely, the \textbf{variables}). It is easy to see that any
monomial $x^{a}y^{b}\in C^{\left(  2\right)  }$ is indeed the product of the
powers $x^{a}$ and $y^{b}$ of these indeterminates, just as the notation suggests.

Now, the \textbf{polynomial ring in two variables }$x$\textbf{ and }%
$y$\textbf{ over }$R$ is defined to be the monoid algebra $R\left[  C^{\left(
2\right)  }\right]  $. It is commonly denoted by $R\left[  x,y\right]  $.
Following Convention \ref{conv.monalg.m-as-em}, we simply write $m$ for
$e_{m}$ whenever $m\in C^{\left(  2\right)  }$; thus, $R\left[  x,y\right]  $
is a free $R$-module with basis%
\[
\left(  x^{a}y^{b}\right)  _{\left(  a,b\right)  \in\mathbb{N}^{2}}.
\]
This means that any $p\in R\left[  x,y\right]  $ can be uniquely written as an
$R$-linear combination%
\[
p=\sum_{\left(  a,b\right)  \in\mathbb{N}^{2}}r_{a,b}x^{a}y^{b}%
\]
with $r_{a,b}\in R$ (such that all but finitely many of these coefficients
$r_{a,b}$ are $0$).

Elements of $R\left[  x,y\right]  $ are called \textbf{polynomials} in $x$ and
$y$.
\end{definition}

Thus, up to notation, the multivariate polynomial ring $R\left[  x,y\right]  $
is just the monoid algebra $R\left[  \mathbb{N}^{2}\right]  $ of the abelian
monoid $\mathbb{N}^{2}=\left(  \mathbb{N}^{2},+,0\right)  $. The unity of the
$R$-algebra $R\left[  x,y\right]  $ is $x^{0}y^{0}=1$.

Here are some examples of bivariate polynomials:

\begin{itemize}
\item If $R=\mathbb{Z}$, then%
\[
3x^{2}y^{7}-10x^{1}y^{1}+8x^{0}y^{5}+2x^{0}y^{0}=3x^{2}y^{7}-10xy+8y^{5}+2
\]
is a polynomial in $x$ and $y$, thus belongs to $\mathbb{Z}\left[  x,y\right]
$.

\item If $R=\mathbb{Q}$, then $\dfrac{2}{3}x^{7}y^{2}-\dfrac{1}{2}x^{1}y^{1}$
is a polynomial in $x$ and $y$, thus belongs to $\mathbb{Q}\left[  x,y\right]
$. (Of course, we don't strictly need the coefficients to be non-integers; the
integers are also included in $\mathbb{Q}$. Thus, for example, the polynomial
$x^{2}y^{3}-2x$ also belongs to $\mathbb{Q}\left[  x,y\right]  $.)

\item Note that $x^{2}+7x$ and $y^{2}+7y$ are two (distinct) polynomials in
$\mathbb{Z}\left[  x,y\right]  $. They happen to involve only one of the two
indeterminates each, but this does not make them any less valid. (There are
also constant polynomials in $\mathbb{Z}\left[  x,y\right]  $, such as
$17x^{0}y^{0}=17$.)
\end{itemize}

\subsubsection{\label{subsec.polys1.polyrings.mulvar}Multivariate polynomials}

We shall now define multivariate polynomial rings (in finitely many variables,
which we name $x_{1},x_{2},\ldots,x_{n}$). Their definition is somewhat
notationally dense (subscripts inside superscripts!), but it is just a
straightforward generalization of Definition \ref{def.polring.bivar}, except
that the indeterminates will now be called $x_{1},x_{2},\ldots,x_{n}$ instead
of $x$ and $y$:

\begin{definition}
\label{def.polring.mulvar}Let $n\in\mathbb{N}$. Let $C^{\left(  n\right)  }$
be the free abelian monoid with $n$ generators $x_{1},x_{2},\ldots,x_{n}$.
This is the monoid whose elements are the distinct symbols%
\[
x_{1}^{a_{1}}x_{2}^{a_{2}}\cdots x_{n}^{a_{n}}\ \ \ \ \ \ \ \ \ \ \text{with
}\left(  a_{1},a_{2},\ldots,a_{n}\right)  \in\mathbb{N}^{n},
\]
and whose operation is defined by%
\begin{align*}
&  \left(  x_{1}^{a_{1}}x_{2}^{a_{2}}\cdots x_{n}^{a_{n}}\right)  \cdot\left(
x_{1}^{b_{1}}x_{2}^{b_{2}}\cdots x_{n}^{b_{n}}\right)  =x_{1}^{a_{1}+b_{1}%
}x_{2}^{a_{2}+b_{2}}\cdots x_{n}^{a_{n}+b_{n}}\\
&  \ \ \ \ \ \ \ \ \ \ \ \ \ \ \ \text{for all }\left(  a_{1},a_{2}%
,\ldots,a_{n}\right)  \in\mathbb{N}^{n}\text{ and }\left(  b_{1},b_{2}%
,\ldots,b_{n}\right)  \in\mathbb{N}^{n}.
\end{align*}
We write this monoid multiplicatively, but of course it is just the monoid
$\mathbb{N}^{n}=\left(  \mathbb{N}^{n},+,0\right)  $ in disguise (where the
addition on $\mathbb{N}^{n}$ that we are calling \textquotedblleft%
$+$\textquotedblright\ here is entrywise, and $0$ means the $n$-tuple $\left(
0,0,\ldots,0\right)  $), with each element $\left(  a_{1},a_{2},\ldots
,a_{n}\right)  $ renamed as $x_{1}^{a_{1}}x_{2}^{a_{2}}\cdots x_{n}^{a_{n}}$
and with addition renamed as multiplication. The elements of $C^{\left(
n\right)  }$ are called \textbf{monomials} in $x_{1},x_{2},\ldots,x_{n}$. For
each $i\in\left\{  1,2,\ldots,n\right\}  $, we define a monomial $x_{i}$ by%
\[
x_{i}=x_{1}^{0}x_{2}^{0}\cdots x_{i-1}^{0}x_{i}^{1}x_{i+1}^{0}x_{i+2}%
^{0}\cdots x_{n}^{0}.
\]
These specific elements $x_{1},x_{2},\ldots,x_{n}$ are called the
\textbf{indeterminates}. It is easy to see that any monomial $x_{1}^{a_{1}%
}x_{2}^{a_{2}}\cdots x_{n}^{a_{n}}\in C^{\left(  n\right)  }$ is indeed the
product of the powers $x_{1}^{a_{1}},x_{2}^{a_{2}},\ldots,x_{n}^{a_{n}}$, just
as the notation suggests.

Now, the \textbf{polynomial ring in }$n$ \textbf{variables }$x_{1}%
,x_{2},\ldots,x_{n}$\textbf{ over }$R$ is defined to be the monoid algebra
$R\left[  C^{\left(  n\right)  }\right]  $. It is commonly denoted by
$R\left[  x_{1},x_{2},\ldots,x_{n}\right]  $. Following Convention
\ref{conv.monalg.m-as-em}, we simply write $m$ for $e_{m}$ whenever $m\in
C^{\left(  n\right)  }$; thus, $R\left[  x_{1},x_{2},\ldots,x_{n}\right]  $ is
a free $R$-module with basis%
\[
\left(  x_{1}^{a_{1}}x_{2}^{a_{2}}\cdots x_{n}^{a_{n}}\right)  _{\left(
a_{1},a_{2},\ldots,a_{n}\right)  \in\mathbb{N}^{n}}.
\]
This means that any $p\in R\left[  x_{1},x_{2},\ldots,x_{n}\right]  $ can be
uniquely written as an $R$-linear combination%
\[
p=\sum_{\left(  a_{1},a_{2},\ldots,a_{n}\right)  \in\mathbb{N}^{n}}%
r_{a_{1},a_{2},\ldots,a_{n}}x_{1}^{a_{1}}x_{2}^{a_{2}}\cdots x_{n}^{a_{n}}%
\]
with $r_{a_{1},a_{2},\ldots,a_{n}}\in R$ (such that all but finitely many of
these coefficients $r_{a_{1},a_{2},\ldots,a_{n}}$ are $0$).

Elements of $R\left[  x_{1},x_{2},\ldots,x_{n}\right]  $ are called
\textbf{polynomials} in $x_{1},x_{2},\ldots,x_{n}$.
\end{definition}

Thus, up to notation, the multivariate polynomial ring $R\left[  x_{1}%
,x_{2},\ldots,x_{n}\right]  $ is just the monoid algebra $R\left[
\mathbb{N}^{n}\right]  $ of the abelian monoid $\mathbb{N}^{n}=\left(
\mathbb{N}^{n},+,0\right)  $.

The multivariate polynomial ring $R\left[  x_{1},x_{2},\ldots,x_{n}\right]  $
is commutative (by Proposition \ref{prop.monalg.comcm}, since the monoid
$\mathbb{N}^{n}$ is abelian).

The univariate polynomial ring $R\left[  x\right]  $ can be viewed as a
particular case of the multivariate polynomial ring $R\left[  x_{1}%
,x_{2},\ldots,x_{n}\right]  $ (obtained by taking $n=1$ and renaming $x_{1}$
as $x$)\ \ \ \ \footnote{Strictly speaking, this requires a minor abuse of
notation: We need to identify each nonnegative integer $n\in\mathbb{N}$ with
the $1$-tuple $\left(  n\right)  \in\mathbb{N}^{1}$, so that the free monoid
$C$ becomes identified with $C^{\left(  1\right)  }$, and therefore its monoid
ring $R\left[  C\right]  $ becomes $R\left[  C^{\left(  1\right)  }\right]
$.}. Likewise, the bivariate polynomial ring $R\left[  x,y\right]  $ is a
particular case of the multivariate polynomial ring $R\left[  x_{1}%
,x_{2},\ldots,x_{n}\right]  $ (obtained by taking $n=2$ and renaming $x_{1}$
and $x_{2}$ as $x$ and $y$).

\subsubsection{\label{subsec.polys1.polyrings.eval-uni}Evaluation, aka
substitution for univariate polynomials}

Polynomials as formal linear combinations are already useful and nice. But
they become a much stronger tool once you learn how to evaluate them, i.e.,
substitute things into them. Unlike a function, a univariate polynomial over
$R$ does not have a fixed domain; you can substitute an element of $R$ into
it, but also a square matrix over $R$ or even another polynomial, and more
generally, any element of an $R$-algebra:

\begin{definition}
\label{def.polring.univar-sub}Let $p\in R\left[  x\right]  $ be a univariate
polynomial. Let $A$ be any $R$-algebra. Let $a\in A$.

We define the element $p\left(  a\right)  \in A$ as follows: Write $p$ as
\[
p=\sum\limits_{i\in\mathbb{N}}p_{i}x^{i}%
\]
with $p_{i}\in R$ (where $p_{i}=0$ for all but finitely many $i\in\mathbb{N}%
$), and set
\begin{equation}
p\left(  a\right)  :=\sum\limits_{i\in\mathbb{N}}p_{i}a^{i}.
\label{eq.def.polring.univar-sub.pa=}%
\end{equation}

This element $p\left(  a\right)  $ is called the \textbf{evaluation} of $p$ at
$a$; we also say that it is obtained by \textbf{substituting} $a$ for $x$ in
$p$.

Sometimes I will denote it by $p\left[  a\right]  $ instead of $p\left(
a\right)  $ (for reasons explained below).
\end{definition}

Note that the $p_{i}\in R$ in Definition \ref{def.polring.univar-sub} are
unique, since $\left(  x^{0},x^{1},x^{2},\ldots\right)  $ is a basis of the
$R$-module $R\left[  x\right]  $. Note also that the infinite sum on the right
hand side of (\ref{eq.def.polring.univar-sub.pa=}) is well-defined, since we
have $p_{i}=0$ for all but finitely many $i\in\mathbb{N}$.

As I said, $A$ can be any $R$-algebra in Definition
\ref{def.polring.univar-sub}: for example, $R$ itself, or a matrix ring
$R^{n\times n}$, or the polynomial ring $R\left[  x\right]  $. In particular,
we can substitute $x$ for $x$ in $p$, obtaining $p\left(  x\right)  =p$.

\begin{warning}
The notation $p\left(  a\right)  $ in Definition \ref{def.polring.univar-sub}
has potential for confusion: Is $p\left(  p+1\right)  $ the evaluation of $p$
at $p+1$ or the product of $p$ with $p+1$ ? This is why I prefer the notation
$p\left[  a\right]  $ instead of $p\left(  a\right)  $. I also recommend using
$\cdot$ for products whenever such confusion could arise (thus, write
$p\cdot\left(  p+1\right)  $ if you mean the product of $p$ with $p+1$). When
reading algebra literature, be aware that you will sometimes have to look at
the context and make sanity checks.
\end{warning}

\begin{example}
\label{exa.polring.univar-sub.not-inj}Let $R=\mathbb{Z}/2$, and let $p$ be the
polynomial $x^{2}+x=x\cdot\left(  x+\overline{1}\right)  \in R\left[
x\right]  $. Let us evaluate $p$ at elements of $R$:%
\begin{align*}
p\left(  \overline{0}\right)   &  =\overline{0}^{2}+\overline{0}=\overline
{0};\\
p\left(  \overline{1}\right)   &  =\overline{1}^{2}+\overline{1}=\overline
{1}+\overline{1}=\overline{2}=\overline{0}.
\end{align*}
Thus, the polynomial $p$ gives $\overline{0}$ when evaluated at any element of
$\mathbb{Z}/2$, even though $p$ is not zero as a polynomial. If you want a
nonzero evaluation of $p$, one thing you can do is to evaluate it on a square
matrix:%
\[
p\left(  \left(
\begin{array}
[c]{cc}%
\overline{0} & \overline{1}\\
\overline{1} & \overline{0}%
\end{array}
\right)  \right)  =\left(
\begin{array}
[c]{cc}%
\overline{0} & \overline{1}\\
\overline{1} & \overline{0}%
\end{array}
\right)  ^{2}+\left(
\begin{array}
[c]{cc}%
\overline{0} & \overline{1}\\
\overline{1} & \overline{0}%
\end{array}
\right)  =\left(
\begin{array}
[c]{cc}%
\overline{1} & \overline{1}\\
\overline{1} & \overline{1}%
\end{array}
\right)  \neq0_{2\times2}.
\]
(Or you can evaluate it at $x$, getting $p\left(  x\right)  =p\neq0$.)
\end{example}

Given an $R$-algebra $A$ and an element $a\in A$, we can study the operation
of substituting $a$ for $x$ into polynomials $p\in R\left[  x\right]  $. This
operation is rather well-behaved:

\begin{theorem}
\label{thm.polring.univar-sub-hom}Let $A$ be an $R$-algebra. Let $a\in A$.
Then, the map%
\begin{align*}
R\left[  x\right]   &  \rightarrow A,\\
p  &  \mapsto p\left[  a\right]
\end{align*}
is an $R$-algebra morphism. In particular, for any two polynomials $p,q\in
R\left[  x\right]  $, we have%
\begin{align}
\left(  pq\right)  \left[  a\right]   &  =p\left[  a\right]  \cdot q\left[
a\right]  ;\label{eq.thm.polring.univar-sub-hom.prod}\\
\left(  p+q\right)  \left[  a\right]   &  =p\left[  a\right]  +q\left[
a\right]  . \label{eq.thm.polring.univar-sub-hom.sum}%
\end{align}

\end{theorem}

The proof of this theorem will be easiest to do after showing the following
simple lemma (compare with Lemma \ref{lem.algebras.assoc-on-basis}):

\begin{lemma}
\label{lem.algebras.mor-on-basis}Let $R$ be a commutative ring. Let $A$ and
$B$ be two $R$-algebras. Let $f:A\rightarrow B$ be an $R$-linear map. Let
$\left(  m_{i}\right)  _{i\in I}$ be a family of vectors in $A$ that spans
$A$. If we have%
\begin{equation}
f\left(  m_{i}m_{j}\right)  =f\left(  m_{i}\right)  f\left(  m_{j}\right)
\ \ \ \ \ \ \ \ \ \ \text{for all }i,j\in I,
\label{eq.lem.algebras.mor-on-basis.ass}%
\end{equation}
then we have%
\begin{equation}
f\left(  ab\right)  =f\left(  a\right)  f\left(  b\right)
\ \ \ \ \ \ \ \ \ \ \text{for all }a,b\in A.
\label{eq.lem.algebras.mor-on-basis.clm}%
\end{equation}

\end{lemma}

\begin{proof}
[Proof of Lemma \ref{lem.algebras.mor-on-basis}.]Let $a,b\in A$. Since the
family $\left(  m_{i}\right)  _{i\in I}$ spans $A$, we can write the two
vectors $a$ and $b$ as
\begin{equation}
a=\sum_{i\in I}a_{i}m_{i}\ \ \ \ \ \ \ \ \ \ \text{and}%
\ \ \ \ \ \ \ \ \ \ b=\sum_{j\in I}b_{j}m_{j}
\label{pf.lem.algebras.mor-on-basis.clm.1}%
\end{equation}
for some coefficients $a_{i}$ and $b_{j}$ in $R$. Consider these coefficients.
Hence,%
\[
ab=\left(  \sum_{i\in I}a_{i}m_{i}\right)  \left(  \sum_{j\in I}b_{j}%
m_{j}\right)  =\sum_{i\in I}\ \ \sum_{j\in I}a_{i}b_{j}m_{i}m_{j}%
\]
(since the multiplication of $A$ is $R$-bilinear) and thus%
\begin{align*}
f\left(  ab\right)   &  =f\left(  \sum_{i\in I}\ \ \sum_{j\in I}a_{i}%
b_{j}m_{i}m_{j}\right)  =\sum_{i\in I}\ \ \sum_{j\in I}a_{i}b_{j}%
\underbrace{f\left(  m_{i}m_{j}\right)  }_{\substack{=f\left(  m_{i}\right)
f\left(  m_{j}\right)  \\\text{(by (\ref{eq.lem.algebras.mor-on-basis.ass}))}%
}}\ \ \ \ \ \ \ \ \ \ \left(  \text{since }f\text{ is }R\text{-linear}\right)
\\
&  =\sum_{i\in I}\ \ \sum_{j\in I}a_{i}b_{j}f\left(  m_{i}\right)  f\left(
m_{j}\right)  =\left(  \sum_{i\in I}a_{i}f\left(  m_{i}\right)  \right)
\left(  \sum_{j\in I}b_{j}f\left(  m_{j}\right)  \right)
\end{align*}
(since the multiplication of $B$ is $R$-bilinear). Comparing this with%
\begin{align*}
f\left(  a\right)  f\left(  b\right)   &  =f\left(  \sum_{i\in I}a_{i}%
m_{i}\right)  f\left(  \sum_{j\in I}b_{j}m_{j}\right)
\ \ \ \ \ \ \ \ \ \ \left(  \text{by (\ref{pf.lem.algebras.mor-on-basis.clm.1}%
)}\right) \\
&  =\left(  \sum_{i\in I}a_{i}f\left(  m_{i}\right)  \right)  \left(
\sum_{j\in I}b_{j}f\left(  m_{j}\right)  \right)  \ \ \ \ \ \ \ \ \ \ \left(
\text{since }f\text{ is }R\text{-linear}\right)  ,
\end{align*}
we obtain $f\left(  ab\right)  =f\left(  a\right)  f\left(  b\right)  $. This
proves Lemma \ref{lem.algebras.mor-on-basis}.
\end{proof}

\begin{proof}
[Proof of Theorem \ref{thm.polring.univar-sub-hom}.]Let $f$ be the map%
\begin{align*}
R\left[  x\right]   &  \rightarrow A,\\
p  &  \mapsto p\left[  a\right]  .
\end{align*}
We must show that $f$ is an $R$-algebra morphism.

It is easy to see that $f$ is $R$-linear. (For example, in order to show that
it respects addition, you need to check that $\left(  p+q\right)  \left[
a\right]  =p\left[  a\right]  +q\left[  a\right]  $ for any $p,q\in R\left[
x\right]  $. But this is done exactly as you would think: Write $p$ and $q$ as
$p=\sum_{i\in\mathbb{N}}p_{i}x^{i}$ (with $p_{i}\in R$) and $q=\sum
_{i\in\mathbb{N}}q_{i}x^{i}$ (with $q_{i}\in R$), and conclude that
\[
p+q=\sum_{i\in\mathbb{N}}p_{i}x^{i}+\sum_{i\in\mathbb{N}}q_{i}x^{i}=\sum
_{i\in\mathbb{N}}\left(  p_{i}x^{i}+q_{i}x^{i}\right)  =\sum_{i\in\mathbb{N}%
}\left(  p_{i}+q_{i}\right)  x^{i},
\]
so that
\begin{align*}
\left(  p+q\right)  \left[  a\right]   &  =\sum_{i\in\mathbb{N}}\left(
p_{i}+q_{i}\right)  a^{i}\ \ \ \ \ \ \ \ \ \ \left(  \text{by the definition
of }\left(  p+q\right)  \left[  a\right]  \right) \\
&  =\sum_{i\in\mathbb{N}}p_{i}a^{i}+\sum_{i\in\mathbb{N}}q_{i}a^{i};
\end{align*}
but it is just as easy to see that $p\left[  a\right]  +q\left[  a\right]  $
gives the same result.)

It is furthermore clear that the map $f$ respects the unity; indeed, $f\left(
1\right)  =1\left[  a\right]  =1$ (since substituting $a$ for $x$ in the
polynomial $1=1x^{0}+0x^{1}+0x^{2}+\cdots$ results in $1a^{0}+0a^{1}%
+0a^{2}+\cdots=1$).

All that now remains is to show that $f$ respects multiplication. In other
words, it remains to show that $f\left(  pq\right)  =f\left(  p\right)
f\left(  q\right)  $ for all $p,q\in R\left[  x\right]  $. Lemma
\ref{lem.algebras.mor-on-basis} gives us a shortcut to proving this: Since the
family $\left(  x^{i}\right)  _{i\in\mathbb{N}}$ is a basis of the $R$-module
$R\left[  x\right]  $ (and thus spans this $R$-module), and since we already
know that $f$ is $R$-linear, it suffices to show that%
\begin{equation}
f\left(  x^{i}x^{j}\right)  =f\left(  x^{i}\right)  f\left(  x^{j}\right)
\ \ \ \ \ \ \ \ \ \ \text{for all }i,j\in\mathbb{N}
\label{pf.thm.polring.univar-sub-hom.lastgoal}%
\end{equation}
(because if we can show this, then Lemma \ref{lem.algebras.mor-on-basis} will
yield that $f\left(  pq\right)  =f\left(  p\right)  f\left(  q\right)  $ for
all $p,q\in R\left[  x\right]  $).

So let us prove (\ref{pf.thm.polring.univar-sub-hom.lastgoal}). Fix
$i,j\in\mathbb{N}$. Then, $x^{i}\left[  a\right]  =a^{i}$ (because
substituting $a$ for $x$ in the polynomial $x^{i}=0x^{0}+0x^{1}+\cdots
+0x^{i-1}+1x^{i}+0x^{i+1}+0x^{i+2}+\cdots$ results in $0a^{0}+0a^{1}%
+\cdots+0a^{i-1}+1a^{i}+0a^{i+1}+0a^{i+2}+\cdots=a^{i}$) and similarly
$x^{j}\left[  a\right]  =a^{j}$ and $x^{i+j}\left[  a\right]  =a^{i+j}$. But
$x^{i}x^{j}=x^{i+j}$, so that%
\begin{align*}
f\left(  x^{i}x^{j}\right)   &  =f\left(  x^{i+j}\right)  =x^{i+j}\left[
a\right]  \ \ \ \ \ \ \ \ \ \ \left(  \text{by the definition of }f\right) \\
&  =a^{i+j}=\underbrace{a^{i}}_{\substack{=x^{i}\left[  a\right]  \\=f\left(
x^{i}\right)  \\\text{(by the definition of }f\text{)}}}\underbrace{a^{j}%
}_{\substack{=x^{j}\left[  a\right]  \\=f\left(  x^{j}\right)  \\\text{(by the
definition of }f\text{)}}}=f\left(  x^{i}\right)  f\left(  x^{j}\right)  .
\end{align*}
This proves (\ref{pf.thm.polring.univar-sub-hom.lastgoal}), and thus concludes
the proof of Theorem \ref{thm.polring.univar-sub-hom}.
\end{proof}

\subsubsection{\label{subsec.polys1.polyrings.eval-mul}Evaluation, aka
substitution for multivariate polynomials}

Likewise, we can substitute multiple elements into a multivariate polynomial,
as long as these elements commute:

\begin{definition}
\label{def.polring.mulvar-sub}Let $n\in\mathbb{N}$. Let $p\in R\left[
x_{1},x_{2},\ldots,x_{n}\right]  $ be a multivariate polynomial. Let $A$ be
any $R$-algebra. Let $a_{1},a_{2},\ldots,a_{n}\in A$ be $n$ elements of $A$
that mutually commute (i.e., that satisfy $a_{i}a_{j}=a_{j}a_{i}$ for each
$i,j\in\left\{  1,2,\ldots,n\right\}  $).

We define the element $p\left(  a_{1},a_{2},\ldots,a_{n}\right)  \in A$ as
follows: Write the polynomial $p$ as%
\[
p=\sum\limits_{\left(  i_{1},i_{2},\ldots,i_{n}\right)  \in\mathbb{N}^{n}%
}p_{i_{1},i_{2},\ldots,i_{n}}x_{1}^{i_{1}}x_{2}^{i_{2}}\cdots x_{n}^{i_{n}}%
\]
with $p_{i_{1},i_{2},\ldots,i_{n}}\in R$ (where $p_{i_{1},i_{2},\ldots,i_{n}%
}=0$ for all but finitely many $\left(  i_{1},i_{2},\ldots,i_{n}\right)
\in\mathbb{N}^{n}$), and set
\[
p\left(  a_{1},a_{2},\ldots,a_{n}\right)  :=\sum\limits_{\left(  i_{1}%
,i_{2},\ldots,i_{n}\right)  \in\mathbb{N}^{n}}p_{i_{1},i_{2},\ldots,i_{n}%
}a_{1}^{i_{1}}a_{2}^{i_{2}}\cdots a_{n}^{i_{n}}.
\]

This element $p\left(  a_{1},a_{2},\ldots,a_{n}\right)  $ is called the
\textbf{evaluation} of $p$ at $a_{1},a_{2},\ldots,a_{n}$; we also say that it
is obtained by \textbf{substituting} $a_{1},a_{2},\ldots,a_{n}$ for
$x_{1},x_{2},\ldots,x_{n}$ in $p$.

Sometimes, I will denote it by $p\left[  a_{1},a_{2},\ldots,a_{n}\right]  $
instead of $p\left(  a_{1},a_{2},\ldots,a_{n}\right)  $.
\end{definition}

Needless to say, this definition generalizes Definition
\ref{def.polring.univar-sub}.

It is clear that $p\left(  x_{1},x_{2},\ldots,x_{n}\right)  =p$ for any
polynomial $p\in R\left[  x_{1},x_{2},\ldots,x_{n}\right]  $.

The analogue to Theorem \ref{thm.polring.univar-sub-hom} now is the following:

\begin{theorem}
\label{thm.polring.mulvar-sub-hom}Let $n\in\mathbb{N}$. Let $A$ be an
$R$-algebra. Let $a_{1},a_{2},\ldots,a_{n}\in A$ be $n$ elements of $A$ that
mutually commute. Then, the map%
\begin{align*}
R\left[  x_{1},x_{2},\ldots,x_{n}\right]   &  \rightarrow A,\\
p  &  \mapsto p\left(  a_{1},a_{2},\ldots,a_{n}\right)
\end{align*}
is an $R$-algebra morphism.
\end{theorem}

\begin{proof}
This is similar to the proof of Theorem \ref{thm.polring.univar-sub-hom}, but
more sophisticated due to the presence of multiple variables. Let $f$ denote
the map defined in Theorem \ref{thm.polring.mulvar-sub-hom}. Again, the only
nontrivial task is to show that $f$ respects multiplication. We note that the
definition of $f$ easily yields that%
\begin{equation}
f\left(  x_{1}^{i_{1}}x_{2}^{i_{2}}\cdots x_{n}^{i_{n}}\right)  =a_{1}^{i_{1}%
}a_{2}^{i_{2}}\cdots a_{n}^{i_{n}} \label{pf.thm.polring.mulvar-sub-hom.fm=}%
\end{equation}
for every $\left(  i_{1},i_{2},\ldots,i_{n}\right)  \in\mathbb{N}^{n}$. We
also observe that every $\left(  i_{1},i_{2},\ldots,i_{n}\right)
\in\mathbb{N}^{n}$ and $\left(  j_{1},j_{2},\ldots,j_{n}\right)  \in
\mathbb{N}^{n}$ satisfy%
\begin{align}
&  \left(  a_{1}^{i_{1}}a_{2}^{i_{2}}\cdots a_{n}^{i_{n}}\right)  \cdot\left(
a_{1}^{j_{1}}a_{2}^{j_{2}}\cdots a_{n}^{j_{n}}\right) \nonumber\\
&  =a_{1}^{i_{1}+j_{1}}a_{2}^{i_{2}+j_{2}}\cdots a_{n}^{i_{n}+j_{n}}.
\label{pf.thm.polring.mulvar-sub-hom.aa=a}%
\end{align}
(This follows with a bit of work from the assumption that $a_{1},a_{2}%
,\ldots,a_{n}$ mutually commute\footnote{\textit{Proof sketch.} Fix $\left(
i_{1},i_{2},\ldots,i_{n}\right)  \in\mathbb{N}^{n}$ and $\left(  j_{1}%
,j_{2},\ldots,j_{n}\right)  \in\mathbb{N}^{n}$.
\par
First, show that if a given element $b\in A$ commutes with each of $k$ given
elements $c_{1},c_{2},\ldots,c_{k}\in A$, then it also commutes with their
product $c_{1}c_{2}\cdots c_{k}$. Use this to show that if $b$ and $c$ are two
elements of $A$ that commute, then $b^{i}$ and $c^{j}$ commute for all
$i,j\in\mathbb{N}$. Conclude that each power $a_{k}^{i_{k}}$ commutes with all
the powers $a_{1}^{j_{1}},a_{2}^{j_{2}},\ldots,a_{k-1}^{j_{k-1}}$ and thus
also with their product $a_{1}^{j_{1}}a_{2}^{j_{2}}\cdots a_{k-1}^{j_{k-1}}$.
In other words,%
\begin{equation}
a_{k}^{i_{k}}\cdot\left(  a_{1}^{j_{1}}a_{2}^{j_{2}}\cdots a_{k-1}^{j_{k-1}%
}\right)  =\left(  a_{1}^{j_{1}}a_{2}^{j_{2}}\cdots a_{k-1}^{j_{k-1}}\right)
\cdot a_{k}^{i_{k}} \label{pf.thm.polring.mulvar-sub-hom.aa=a.pf.1}%
\end{equation}
for each $k\in\left\{  1,2,\ldots,n\right\}  $.
\par
Now, you can show that the equality%
\begin{equation}
\left(  a_{1}^{i_{1}}a_{2}^{i_{2}}\cdots a_{k}^{i_{k}}\right)  \cdot\left(
a_{1}^{j_{1}}a_{2}^{j_{2}}\cdots a_{k}^{j_{k}}\right)  =a_{1}^{i_{1}+j_{1}%
}a_{2}^{i_{2}+j_{2}}\cdots a_{k}^{i_{k}+j_{k}}
\label{pf.thm.polring.mulvar-sub-hom.aa=a.pf.2}%
\end{equation}
holds for each $k\in\left\{  0,1,\ldots,n\right\}  $. Indeed, this equality
can be proved by induction on $k$, where the induction step (from $k-1$ to
$k$) proceeds by the following computation:%
\begin{align*}
&  \underbrace{\left(  a_{1}^{i_{1}}a_{2}^{i_{2}}\cdots a_{k}^{i_{k}}\right)
}_{=\left(  a_{1}^{i_{1}}a_{2}^{i_{2}}\cdots a_{k-1}^{i_{k-1}}\right)
a_{k}^{i_{k}}}\cdot\underbrace{\left(  a_{1}^{j_{1}}a_{2}^{j_{2}}\cdots
a_{k}^{j_{k}}\right)  }_{=\left(  a_{1}^{j_{1}}a_{2}^{j_{2}}\cdots
a_{k-1}^{j_{k-1}}\right)  a_{k}^{j_{k}}}\\
&  =\left(  a_{1}^{i_{1}}a_{2}^{i_{2}}\cdots a_{k-1}^{i_{k-1}}\right)
\underbrace{a_{k}^{i_{k}}\cdot\left(  a_{1}^{j_{1}}a_{2}^{j_{2}}\cdots
a_{k-1}^{j_{k-1}}\right)  }_{\substack{=\left(  a_{1}^{j_{1}}a_{2}^{j_{2}%
}\cdots a_{k-1}^{j_{k-1}}\right)  \cdot a_{k}^{i_{k}}\\\text{(by
(\ref{pf.thm.polring.mulvar-sub-hom.aa=a.pf.1}))}}}a_{k}^{j_{k}}\\
&  =\underbrace{\left(  a_{1}^{i_{1}}a_{2}^{i_{2}}\cdots a_{k-1}^{i_{k-1}%
}\right)  \left(  a_{1}^{j_{1}}a_{2}^{j_{2}}\cdots a_{k-1}^{j_{k-1}}\right)
}_{\substack{=a_{1}^{i_{1}+j_{1}}a_{2}^{i_{2}+j_{2}}\cdots a_{k-1}%
^{i_{k-1}+j_{k-1}}\\\text{(by the induction hypothesis)}}}\cdot
\underbrace{a_{k}^{i_{k}}a_{k}^{j_{k}}}_{=a_{k}^{i_{k}+j_{k}}}\\
&  =\left(  a_{1}^{i_{1}+j_{1}}a_{2}^{i_{2}+j_{2}}\cdots a_{k-1}%
^{i_{k-1}+j_{k-1}}\right)  \cdot a_{k}^{i_{k}+j_{k}}=a_{1}^{i_{1}+j_{1}}%
a_{2}^{i_{2}+j_{2}}\cdots a_{k}^{i_{k}+j_{k}}.
\end{align*}
Thus, (\ref{pf.thm.polring.mulvar-sub-hom.aa=a.pf.2}) is proved.
\par
Now, applying (\ref{pf.thm.polring.mulvar-sub-hom.aa=a.pf.2}) to $k=n$, we
obtain
\[
\left(  a_{1}^{i_{1}}a_{2}^{i_{2}}\cdots a_{n}^{i_{n}}\right)  \cdot\left(
a_{1}^{j_{1}}a_{2}^{j_{2}}\cdots a_{n}^{j_{n}}\right)  =a_{1}^{i_{1}+j_{1}%
}a_{2}^{i_{2}+j_{2}}\cdots a_{n}^{i_{n}+j_{n}},
\]
qed.}.)

To prove that $f$ respects multiplication, we again use Lemma
\ref{lem.algebras.mor-on-basis}. This time, instead of proving
(\ref{pf.thm.polring.univar-sub-hom.lastgoal}), we need to prove the equality%
\[
f\left(  \left(  x_{1}^{i_{1}}x_{2}^{i_{2}}\cdots x_{n}^{i_{n}}\right)
\cdot\left(  x_{1}^{j_{1}}x_{2}^{j_{2}}\cdots x_{n}^{j_{n}}\right)  \right)
=f\left(  x_{1}^{i_{1}}x_{2}^{i_{2}}\cdots x_{n}^{i_{n}}\right)  \cdot
f\left(  x_{1}^{j_{1}}x_{2}^{j_{2}}\cdots x_{n}^{j_{n}}\right)
\]
for all $\left(  i_{1},i_{2},\ldots,i_{n}\right)  \in\mathbb{N}^{n}$ and
$\left(  j_{1},j_{2},\ldots,j_{n}\right)  \in\mathbb{N}^{n}$. This equality
follows by comparing%
\begin{align*}
f\left(  \underbrace{\left(  x_{1}^{i_{1}}x_{2}^{i_{2}}\cdots x_{n}^{i_{n}%
}\right)  \cdot\left(  x_{1}^{j_{1}}x_{2}^{j_{2}}\cdots x_{n}^{j_{n}}\right)
}_{=x_{1}^{i_{1}+j_{1}}x_{2}^{i_{2}+j_{2}}\cdots x_{n}^{i_{n}+j_{n}}}\right)
&  =f\left(  x_{1}^{i_{1}+j_{1}}x_{2}^{i_{2}+j_{2}}\cdots x_{n}^{i_{n}+j_{n}%
}\right) \\
&  =a_{1}^{i_{1}+j_{1}}a_{2}^{i_{2}+j_{2}}\cdots a_{n}^{i_{n}+j_{n}%
}\ \ \ \ \ \ \ \ \ \ \left(  \text{by (\ref{pf.thm.polring.mulvar-sub-hom.fm=}%
)}\right)
\end{align*}
with%
\begin{align*}
\underbrace{f\left(  x_{1}^{i_{1}}x_{2}^{i_{2}}\cdots x_{n}^{i_{n}}\right)
}_{\substack{=a_{1}^{i_{1}}a_{2}^{i_{2}}\cdots a_{n}^{i_{n}}\\\text{(by
(\ref{pf.thm.polring.mulvar-sub-hom.fm=}))}}}\cdot\underbrace{f\left(
x_{1}^{j_{1}}x_{2}^{j_{2}}\cdots x_{n}^{j_{n}}\right)  }_{\substack{=a_{1}%
^{j_{1}}a_{2}^{j_{2}}\cdots a_{n}^{j_{n}}\\\text{(by
(\ref{pf.thm.polring.mulvar-sub-hom.fm=}))}}}  &  =\left(  a_{1}^{i_{1}}%
a_{2}^{i_{2}}\cdots a_{n}^{i_{n}}\right)  \cdot\left(  a_{1}^{j_{1}}%
a_{2}^{j_{2}}\cdots a_{n}^{j_{n}}\right) \\
&  =a_{1}^{i_{1}+j_{1}}a_{2}^{i_{2}+j_{2}}\cdots a_{n}^{i_{n}+j_{n}%
}\ \ \ \ \ \ \ \ \ \ \left(  \text{by
(\ref{pf.thm.polring.mulvar-sub-hom.aa=a})}\right)  .
\end{align*}
Thus, we conclude (using Lemma \ref{lem.algebras.mor-on-basis}) that $f$
respects multiplication, and so the proof of Theorem
\ref{thm.polring.mulvar-sub-hom} is easily completed.
\end{proof}

\begin{exercise}
Let $f\in R\left[  x,y\right]  $ be a polynomial in two variables $x$ and $y$.
Let $g=f\left(  y,x\right)  $. (This is the evaluation of $f$ at $y,x$. In
other words, $g$ is the result of replacing each monomial $x^{i}y^{j}$ by
$y^{i}x^{j}$ in $f$. For example, if $f=x^{2}+7xy-y$, then $g=y^{2}+7yx-x$.)

Prove that the difference $f-g$ is divisible by $x-y$ in the ring $R\left[
x,y\right]  $. \medskip

[\textbf{Hint:} Use linearity to reduce the general case to the case when $f$
is a single monomial.]
\end{exercise}

\begin{exercise}
Let $n\in\mathbb{N}$. Let $P$ be the multivariate polynomial ring $R\left[
x_{1},x_{2},\ldots,x_{n}\right]  $.

A polynomial $f\in P$ will be called \textbf{symmetric} if every permutation
$\sigma$ of the set $\left\{  1,2,\ldots,n\right\}  $ satisfies%
\[
f\left(  x_{\sigma\left(  1\right)  },x_{\sigma\left(  2\right)  }%
,\ldots,x_{\sigma\left(  n\right)  }\right)  =f\left(  x_{1},x_{2}%
,\ldots,x_{n}\right)  .
\]
(In other words, a polynomial $f\in P$ is symmetric if and only if $f$ remains
unchanged whenever the indeterminates are permuted. For example, the
polynomials $x_{1}+x_{2}+\cdots+x_{n}$ and $\left(  3+x_{1}\right)  \left(
3+x_{2}\right)  \cdots\left(  3+x_{n}\right)  $ are symmetric, but the
polynomial $x_{1}x_{2}+x_{2}x_{3}+\cdots+x_{n-1}x_{n}$ is not\footnotemark.)

Prove that the set of all symmetric polynomials $f\in P$ is an $R$-subalgebra
of $P$.
\end{exercise}

\footnotetext{unless $R$ is trivial or $n\leq2$}

\subsubsection{Constant polynomials}

Finally, a few more pieces of notation. We recall the notion of a constant
element of a monoid ring (Convention \ref{conv.monalg.const}). Since a
polynomial ring is a monoid ring, we can apply it to polynomial rings, and
obtain the following:

\begin{convention}
\label{conv.polring.constant-pols}Let $n\in\mathbb{N}$. Then, we identify each
$r\in R$ with $r\cdot1\in R\left[  x_{1},x_{2},\ldots,x_{n}\right]  $ (where
$1$ means the monomial $x_{1}^{0}x_{2}^{0}\cdots x_{n}^{0}$, which is the
unity of $R\left[  x_{1},x_{2},\ldots,x_{n}\right]  $). This identification is
harmless, and turns $R$ into an $R$-subalgebra of $R\left[  x_{1},x_{2}%
,\ldots,x_{n}\right]  $.

A polynomial $p\in R\left[  x_{1},x_{2},\ldots,x_{n}\right]  $ is said to be
\textbf{constant} if it lies in this subalgebra (i.e., if it satisfies
$p=r\cdot1$ for some $r\in R$).
\end{convention}

\begin{example}
The polynomial $3=3x^{0}\in R\left[  x\right]  $ is constant, but the
polynomial $3x=3x^{1}$ is not (unless $3=0$ in $R$).
\end{example}

\subsubsection{Coefficients}

By their definition, polynomials are $R$-linear combinations of monomials. Let
us introduce a notation for the coefficients in these $R$-linear combinations:

\begin{definition}
\label{def.polring.coeff}Let $p\in R\left[  x_{1},x_{2},\ldots,x_{n}\right]  $
be a polynomial. Let $\mathfrak{m}=x_{1}^{a_{1}}x_{2}^{a_{2}}\cdots
x_{n}^{a_{n}}$ be a monomial. Then, the \textbf{coefficient} of $\mathfrak{m}$
in $p$ is the element $\left[  \mathfrak{m}\right]  p$ of $R$ defined as
follows: If we write $p$ as%
\[
p=\sum\limits_{\left(  i_{1},i_{2},\ldots,i_{n}\right)  \in\mathbb{N}^{n}%
}p_{i_{1},i_{2},\ldots,i_{n}}x_{1}^{i_{1}}x_{2}^{i_{2}}\cdots x_{n}^{i_{n}}%
\]
with $p_{i_{1},i_{2},\ldots,i_{n}}\in R$, then we set
\[
\left[  \mathfrak{m}\right]  p:=p_{a_{1},a_{2},\ldots,a_{n}}.
\]

\end{definition}

\begin{example}
\label{exa.polring.coeffs}\ \ 

\begin{enumerate}
\item[\textbf{(a)}] For univariate polynomials, we have%
\[
\left[  x^{3}\right]  \left(  \left(  1+x\right)  ^{5}\right)
=10\ \ \ \ \ \ \ \ \ \ \text{and}\ \ \ \ \ \ \ \ \ \ \left[  x^{7}\right]
\left(  \left(  1+x\right)  ^{5}\right)  =0
\]
(since $\left(  1+x\right)  ^{5}=1+5x+10x^{2}+10x^{3}+5x^{4}+x^{5}$).

\item[\textbf{(b)}] For multivariate polynomials, we have%
\[
\left[  x_{1}^{2}x_{2}^{3}\right]  \left(  \left(  x_{1}+x_{2}\right)
^{5}\right)  =10\ \ \ \ \ \ \ \ \ \ \text{and}\ \ \ \ \ \ \ \ \ \ \left[
x_{1}\right]  \left(  \left(  x_{1}+x_{2}\right)  ^{5}\right)  =0
\]
(since $\left(  x_{1}+x_{2}\right)  ^{5}=x_{1}^{5}+5x_{1}^{4}x_{2}+10x_{1}%
^{3}x_{2}^{2}+10x_{1}^{2}x_{2}^{3}+5x_{1}x_{2}^{4}+x_{2}^{5}$).
\end{enumerate}
\end{example}

\subsubsection{Renaming indeterminates}

Often we will want to use symbols other than $x_{1},x_{2},\ldots,x_{n}$ for
indeterminates. Thus, we allow ourselves to rename these indeterminates when
it pleases us. For example, we can rename the indeterminates $x_{1}$ and
$x_{2}$ of the polynomial ring $R\left[  x_{1},x_{2}\right]  $ as $x$ and $y$,
so that the equations in Example \ref{exa.polring.coeffs} \textbf{(b)} become%
\[
\left[  x^{2}y^{3}\right]  \left(  \left(  x+y\right)  ^{5}\right)
=10\ \ \ \ \ \ \ \ \ \ \text{and}\ \ \ \ \ \ \ \ \ \ \left[  x\right]  \left(
\left(  x+y\right)  ^{5}\right)  =0.
\]
When we do this, we shall also rename the ring $R\left[  x_{1},x_{2}\right]  $
as $R\left[  x,y\right]  $. More generally, we can have polynomial rings in
any (finite) set of indeterminates; these rings are written by putting the
names of these indeterminates into the square brackets. For example, $R\left[
a,b,x,y\right]  $ means a polynomial ring in four indeterminates named
$a,b,x,y$.

\begin{remark}
This convention will serve us well in this course, but it eventually reveals
itself to be inconvenient as you move into more advanced territory. In fact,
it is better to think of polynomial rings with differently named
indeterminates as being genuinely distinct rather than merely renamed versions
of one another. For example, the two polynomial rings $R\left[  x\right]  $
and $R\left[  y\right]  $ are best regarded as distinct (even though they are
isomorphic). This allows us to view them both as two \textbf{different}
subrings of $R\left[  x,y\right]  $ (where the first one consists of all
polynomials that don't contain $y$, such as $x^{3}+2x+1$, whereas the second
consists of all polynomials that don't contain $x$, such as $y^{2}-2y$). This
viewpoint is rather natural, but cannot be rigorously justified as long as we
view $y$ as being the same indeterminate as $x$ in all but name. Our
definition of a multivariate polynomial ring $R\left[  x_{1},x_{2}%
,\ldots,x_{n}\right]  $ (Definition \ref{def.polring.mulvar}) depends only on
a ring $R$ and a number $n$, so that it does not support distinguishing
between different polynomial rings with the same $R$ and the same $n$.

Thus, it is advisable to have a more flexible definition, which allows us to
arbitrarily specify the names of the indeterminates. For example, we should be
able to define the polynomial rings $R\left[  x,y\right]  $ and $R\left[
y,z\right]  $, which are each isomorphic to $R\left[  x_{1},x_{2}\right]  $,
but should not be treated as being the same ring (since the former has
indeterminates $x$ and $y$ whereas the latter has indeterminates $y$ and $z$).

Such a definition can be obtained by making some minor changes to our
Definition \ref{def.polring.mulvar}. Namely, let $S$ be any finite set of
symbols, which we want to use as indeterminates (for example, we can have
$S=\left\{  x,y\right\}  $ or $S=\left\{  y,z\right\}  $ or $S=\left\{
\alpha,\mathbf{w},\clubsuit\right\}  $ if we are being silly). Now, instead of
using the monoid
\[
C^{\left(  n\right)  }=\left\{  x_{1}^{a_{1}}x_{2}^{a_{2}}\cdots x_{n}^{a_{n}%
}\ \mid\ \left(  a_{1},a_{2},\ldots,a_{n}\right)  \in\mathbb{N}^{n}\right\}
\]
as the set of monomials, we use the monoid%
\[
C^{\left(  S\right)  }:=\left\{  \prod_{s\in S}s^{a_{s}}\ \mid\ a_{s}%
\in\mathbb{N}\text{ for each }s\in S\right\}  ,
\]
where the \textquotedblleft product\textquotedblright\ $\prod_{s\in S}%
s^{a_{s}}$ is just a formal symbol that encodes a family $\left(
a_{s}\right)  _{s\in S}\in\mathbb{N}^{S}$ of nonnegative integers (one for
each element $s$ of $S$). A monomial in $C^{\left(  S\right)  }$ is thus a
\textquotedblleft formal\textquotedblright\ product of the form $\prod_{s\in
S}s^{a_{s}}$ (with each factor being a formal power of one of our
indeterminates), and two such monomials are multiplied by the rule%
\[
\left(  \prod_{s\in S}s^{a_{s}}\right)  \cdot\left(  \prod_{s\in S}s^{b_{s}%
}\right)  =\prod_{s\in S}s^{a_{s}+b_{s}}.
\]
The polynomial ring in the set $S$ of indeterminates is then defined as the
monoid ring $R\left[  C^{\left(  S\right)  }\right]  $ of this monoid
$C^{\left(  S\right)  }$. We shall refer to such a ring as a
\textbf{multivariate polynomial ring with named variables}, and just call it
$R\left[  S\right]  $.

Thus, for a three-element set $S=\left\{  x,y,z\right\}  $, we obtain the ring
$R\left[  S\right]  =R\left[  C^{\left(  S\right)  }\right]  =R\left[
x,y,z\right]  $, which is the polynomial ring over $R$ in three variables that
are named $x,y,z$. For instance, $x^{2}+7y^{3}z-xyz$ is a polynomial in this
ring $R\left[  x,y,z\right]  $.

Now, of course, this polynomial ring $R\left[  x,y,z\right]  $ is isomorphic
to $R\left[  x_{1},x_{2},x_{3}\right]  $ as an $R$-algebra (via the
isomorphism that sends each monomial $x^{a}y^{b}z^{c}$ to $x_{1}^{a}x_{2}%
^{b}x_{3}^{c}$). More generally, we have $R\left[  S\right]  \cong R\left[
x_{1},x_{2},\ldots,x_{n}\right]  $ whenever $S$ is an $n$-element set. Thus,
the rings $R\left[  x,y\right]  $, $R\left[  y,z\right]  $ and $R\left[
x_{1},x_{2}\right]  $ are all isomorphic, even though they are distinct.
Hence, named variables do not introduce anything genuinely new to our theory
as long as we are studying a single polynomial ring at a time. But their
flexibility is helpful when working with several polynomial rings, e.g., by
allowing us to treat $R\left[  x\right]  $ and $R\left[  y\right]  $ as two
different subrings of $R\left[  x,y\right]  $.

This definition of the polynomial ring $R\left[  S\right]  $ with named
variables can be easily adapted to infinite sets $S$ as well, with a slight
change: The monoid $C^{\left(  S\right)  }$ then needs to be defined as
\[
\left\{  \prod_{s\in S}s^{a_{s}}\ \mid\ a_{s}\in\mathbb{N}\text{ for each
}s\in S\text{, and all but finitely many }s\in S\text{ satisfy }%
a_{s}=0\right\}  .
\]
(Thus, each monomial contains only finitely many indeterminates in nonzero powers.)
\end{remark}

\subsubsection{A remark on noncommutative $R$}

\begin{remark}
We have been requiring that the ring $R$ be commutative. However, it
absolutely is possible to define polynomials over a noncommutative ring, if
you are sufficiently careful. (In particular, this includes polynomials over
matrix rings; these are rather useful in linear algebra. The indeterminates in
such polynomials commute with all elements of $R$.) We have defined the notion
of an $R$-algebra only for commutative rings $R$, but there are ways to adapt
it to the general setup; alternatively, it is possible to redo the
construction of the polynomial ring by hand without using $R$-algebras. See
\cite[\S 1.3]{ChaLoi21} for the latter approach.
\end{remark}

\subsubsection{Polynomials commute with algebra morphisms}

The following proposition is basic and highly useful:

\begin{proposition}
[\textquotedblleft Polynomials commute with algebra
morphisms\textquotedblright]\label{prop.algmor.poly-comm}Let $A$ and $B$ be
two $R$-algebras. Let $f:A\rightarrow B$ be an $R$-algebra morphism. Let $a\in
A$. Let $p\in R\left[  x\right]  $ be a polynomial. Then,%
\[
f\left(  p\left[  a\right]  \right)  =p\left[  f\left(  a\right)  \right]  .
\]

\end{proposition}

\begin{proof}
[Proof of Proposition \ref{prop.algmor.poly-comm}.]Let us give a proof by
example: Set $p=5x^{4}+x^{3}+7x^{1}$. Then, $p\left[  a\right]  =5a^{4}%
+a^{3}+7a^{1}$ and $p\left[  f\left(  a\right)  \right]  =5f\left(  a\right)
^{4}+f\left(  a\right)  ^{3}+7f\left(  a\right)  ^{1}$. Thus, the claim we
have to prove rewrites as%
\[
f\left(  5a^{4}+a^{3}+7a^{1}\right)  =5f\left(  a\right)  ^{4}+f\left(
a\right)  ^{3}+7f\left(  a\right)  ^{1}.
\]
But this follows easily from the fact that $f$ is an $R$-algebra morphism:
Indeed,%
\begin{align*}
f\left(  5a^{4}+a^{3}+7a^{1}\right)   &  =f\left(  5a^{4}\right)  +f\left(
a^{3}\right)  +f\left(  7a^{1}\right)  \ \ \ \ \ \ \ \ \ \ \left(  \text{since
}f\text{ respects addition}\right) \\
&  =5f\left(  a^{4}\right)  +f\left(  a^{3}\right)  +7f\left(  a^{1}\right)
\ \ \ \ \ \ \ \ \ \ \left(  \text{since }f\text{ respects scaling}\right) \\
&  =5f\left(  a\right)  ^{4}+f\left(  a\right)  ^{3}+7f\left(  a\right)
^{1}\ \ \ \ \ \ \ \ \ \ \left(  \text{since }f\text{ respects powers}\right)
.
\end{align*}

The rigorous proof in the general case falls along the same lines: Writing $p$
as $p=\sum_{i=0}^{n}p_{i}x^{i}$ (for $n\in\mathbb{N}$ and $p_{0},p_{1}%
,\ldots,p_{n}\in R$), we can easily see that $p\left[  a\right]  =\sum
_{i=0}^{n}p_{i}a^{i}$ and $p\left[  f\left(  a\right)  \right]  =\sum
_{i=0}^{n}p_{i}\left(  f\left(  a\right)  \right)  ^{i}$, so that the claim we
have to prove rewrites as%
\[
f\left(  \sum_{i=0}^{n}p_{i}a^{i}\right)  =\sum_{i=0}^{n}p_{i}\left(  f\left(
a\right)  \right)  ^{i};
\]
but this follows from
\begin{align*}
f\left(  \sum_{i=0}^{n}p_{i}a^{i}\right)   &  =\sum_{i=0}^{n}p_{i}%
\underbrace{f\left(  a^{i}\right)  }_{\substack{=\left(  f\left(  a\right)
\right)  ^{i}\\\text{(by Proposition \ref{prop.ringmor.respects}
\textbf{(f)},}\\\text{since }f\text{ is a ring morphism)}}%
}\ \ \ \ \ \ \ \ \ \ \left(
\begin{array}
[c]{c}%
\text{by Proposition \ref{prop.mods.lin-map-pres-lc},}\\
\text{since }f\text{ is }R\text{-linear}%
\end{array}
\right) \\
&  =\sum_{i=0}^{n}p_{i}\left(  f\left(  a\right)  \right)  ^{i}.
\end{align*}
Thus, Proposition \ref{prop.algmor.poly-comm} is proved.
\end{proof}

Proposition \ref{prop.algmor.poly-comm} can be generalized to multivariate polynomials:

\begin{proposition}
[\textquotedblleft Polynomials commute with algebra
morphisms\textquotedblright]\label{prop.algmor.poly-comm-multi}Let $A$ and $B$
be two $R$-algebras. Let $f:A\rightarrow B$ be an $R$-algebra morphism. Let
$a_{1},a_{2},\ldots,a_{n}$ be $n$ elements of $A$ that mutually commute. Let
$p\in R\left[  x_{1},x_{2},\ldots,x_{n}\right]  $ be a polynomial. Then, the
elements $f\left(  a_{1}\right)  ,f\left(  a_{2}\right)  ,\ldots,f\left(
a_{n}\right)  $ of $B$ mutually commute as well, and we have%
\begin{equation}
f\left(  p\left[  a_{1},a_{2},\ldots,a_{n}\right]  \right)  =p\left[  f\left(
a_{1}\right)  ,f\left(  a_{2}\right)  ,\ldots,f\left(  a_{n}\right)  \right]
. \label{eq.prop.algmor.poly-comm-multi.eq}%
\end{equation}

\end{proposition}

\begin{proof}
The claim that $f\left(  a_{1}\right)  ,f\left(  a_{2}\right)  ,\ldots
,f\left(  a_{n}\right)  $ mutually commute is easy: For any $i,j\in\left\{
1,2,\ldots,n\right\}  $, we have
\begin{align*}
f\left(  a_{i}\right)  f\left(  a_{j}\right)   &  =f\left(  a_{i}a_{j}\right)
\ \ \ \ \ \ \ \ \ \ \left(  \text{since }f\text{ is a ring morphism}\right) \\
&  =f\left(  a_{j}a_{i}\right)  \ \ \ \ \ \ \ \ \ \ \left(  \text{since }%
a_{1},a_{2},\ldots,a_{n}\text{ commute, so that }a_{i}a_{j}=a_{j}a_{i}\right)
\\
&  =f\left(  a_{j}\right)  f\left(  a_{i}\right)  \ \ \ \ \ \ \ \ \ \ \left(
\text{since }f\text{ is a ring morphism}\right)  .
\end{align*}
It remains to prove (\ref{eq.prop.algmor.poly-comm-multi.eq}). This is
analogous to the proof of Proposition \ref{prop.algmor.poly-comm}, just with
more subscripts: Write $p$ as%
\[
p=\sum\limits_{\left(  i_{1},i_{2},\ldots,i_{n}\right)  \in\mathbb{N}^{n}%
}p_{i_{1},i_{2},\ldots,i_{n}}x_{1}^{i_{1}}x_{2}^{i_{2}}\cdots x_{n}^{i_{n}%
}\ \ \ \ \ \ \ \ \ \ \text{with }p_{i_{1},i_{2},\ldots,i_{n}}\in R.
\]
Then,
\begin{align*}
p\left[  a_{1},a_{2},\ldots,a_{n}\right]   &  =\sum\limits_{\left(
i_{1},i_{2},\ldots,i_{n}\right)  \in\mathbb{N}^{n}}p_{i_{1},i_{2},\ldots
,i_{n}}a_{1}^{i_{1}}a_{2}^{i_{2}}\cdots a_{n}^{i_{n}}%
\ \ \ \ \ \ \ \ \ \ \text{and}\\
p\left[  f\left(  a_{1}\right)  ,f\left(  a_{2}\right)  ,\ldots,f\left(
a_{n}\right)  \right]   &  =\sum\limits_{\left(  i_{1},i_{2},\ldots
,i_{n}\right)  \in\mathbb{N}^{n}}p_{i_{1},i_{2},\ldots,i_{n}}\left(  f\left(
a_{1}\right)  \right)  ^{i_{1}}\left(  f\left(  a_{2}\right)  \right)
^{i_{2}}\cdots\left(  f\left(  a_{n}\right)  \right)  ^{i_{n}},
\end{align*}
so the equality (\ref{eq.prop.algmor.poly-comm-multi.eq}) (which we must
prove) rewrites as%
\begin{align*}
&  f\left(  \sum\limits_{\left(  i_{1},i_{2},\ldots,i_{n}\right)
\in\mathbb{N}^{n}}p_{i_{1},i_{2},\ldots,i_{n}}a_{1}^{i_{1}}a_{2}^{i_{2}}\cdots
a_{n}^{i_{n}}\right) \\
&  =\sum\limits_{\left(  i_{1},i_{2},\ldots,i_{n}\right)  \in\mathbb{N}^{n}%
}p_{i_{1},i_{2},\ldots,i_{n}}\left(  f\left(  a_{1}\right)  \right)  ^{i_{1}%
}\left(  f\left(  a_{2}\right)  \right)  ^{i_{2}}\cdots\left(  f\left(
a_{n}\right)  \right)  ^{i_{n}}.
\end{align*}
But this follows from%
\begin{align*}
&  f\left(  \sum\limits_{\left(  i_{1},i_{2},\ldots,i_{n}\right)
\in\mathbb{N}^{n}}p_{i_{1},i_{2},\ldots,i_{n}}a_{1}^{i_{1}}a_{2}^{i_{2}}\cdots
a_{n}^{i_{n}}\right) \\
&  =\sum\limits_{\left(  i_{1},i_{2},\ldots,i_{n}\right)  \in\mathbb{N}^{n}%
}p_{i_{1},i_{2},\ldots,i_{n}}\underbrace{f\left(  a_{1}^{i_{1}}a_{2}^{i_{2}%
}\cdots a_{n}^{i_{n}}\right)  }_{\substack{=f\left(  a_{1}^{i_{1}}\right)
f\left(  a_{2}^{i_{2}}\right)  \cdots f\left(  a_{n}^{i_{n}}\right)
\\\text{(by Proposition \ref{prop.ringmor.respects} \textbf{(b)}%
,}\\\text{since }f\text{ is a ring morphism)}}}\ \ \ \ \ \ \ \ \ \ \left(
\begin{array}
[c]{c}%
\text{by Proposition \ref{prop.mods.lin-map-pres-lc},}\\
\text{since }f\text{ is }R\text{-linear}%
\end{array}
\right) \\
&  =\sum\limits_{\left(  i_{1},i_{2},\ldots,i_{n}\right)  \in\mathbb{N}^{n}%
}p_{i_{1},i_{2},\ldots,i_{n}}\underbrace{f\left(  a_{1}^{i_{1}}\right)
}_{\substack{=\left(  f\left(  a_{1}\right)  \right)  ^{i_{1}}\\\text{(by
Proposition \ref{prop.ringmor.respects} \textbf{(f)},}\\\text{since }f\text{
is a ring morphism)}}}\underbrace{f\left(  a_{2}^{i_{2}}\right)
}_{\substack{=\left(  f\left(  a_{2}\right)  \right)  ^{i_{2}}%
\\\text{(likewise)}}}\cdots\underbrace{f\left(  a_{n}^{i_{n}}\right)
}_{\substack{=\left(  f\left(  a_{n}\right)  \right)  ^{i_{n}}%
\\\text{(likewise)}}}\\
&  =\sum\limits_{\left(  i_{1},i_{2},\ldots,i_{n}\right)  \in\mathbb{N}^{n}%
}p_{i_{1},i_{2},\ldots,i_{n}}\left(  f\left(  a_{1}\right)  \right)  ^{i_{1}%
}\left(  f\left(  a_{2}\right)  \right)  ^{i_{2}}\cdots\left(  f\left(
a_{n}\right)  \right)  ^{i_{n}}.
\end{align*}

\end{proof}

\subsubsection{\label{subsec.polys1.polyrings.ass-subst}Associativity of
substitution}

Polynomials may not be functions, but they can still be composed: Indeed, the
\textbf{composition} $p\circ q$ of two univariate polynomials $p,q\in R\left[
x\right]  $ is simply defined as the evaluation $p\left[  q\right]  $. It is
tempting to assume that this composition defines the same associative law%
\[
\left(  p\circ q\right)  \circ r=p\circ\left(  q\circ r\right)
\ \ \ \ \ \ \ \ \ \ \text{for }p,q,r\in R\left[  x\right]
\]
as composition of maps. And this is indeed true; even better, this holds if we
let $r$ be any element of any $R$-algebra:

\begin{theorem}
\label{thm.polys1.comp-assoc}Let $p,q\in R\left[  x\right]  $ be two
univariate polynomials. Then:

\begin{enumerate}
\item[\textbf{(a)}] If $A$ is any $R$-algebra, and if $a\in A$ is arbitrary,
then%
\[
\left(  p\left[  q\right]  \right)  \left[  a\right]  =p\left[  q\left[
a\right]  \right]  .
\]

\item[\textbf{(b)}] In particular, if $r\in R\left[  x\right]  $, then
$\left(  p\circ q\right)  \circ r=p\circ\left(  q\circ r\right)  $.
\end{enumerate}
\end{theorem}

\begin{proof}
\textbf{(a)} Let $A$ be any $R$-algebra, and let $a\in A$. We know from
Theorem \ref{thm.polring.univar-sub-hom} that the map%
\begin{align*}
R\left[  x\right]   &  \rightarrow A,\\
p  &  \mapsto p\left[  a\right]
\end{align*}
is an $R$-algebra morphism. Let us denote this morphism by $\operatorname*{ev}%
\nolimits_{a}$. Then, $\operatorname*{ev}\nolimits_{a}\left(  p\left[
q\right]  \right)  =\left(  p\left[  q\right]  \right)  \left[  a\right]  $
and $\operatorname*{ev}\nolimits_{a}\left(  q\right)  =q\left[  a\right]  $.
But Proposition \ref{prop.algmor.poly-comm} (applied to $\operatorname*{ev}%
\nolimits_{a}$, $R\left[  x\right]  $, $A$ and $q$ instead of $f$, $A$, $B$
and $a$) yields $\operatorname*{ev}\nolimits_{a}\left(  p\left[  q\right]
\right)  =p\left[  \operatorname*{ev}\nolimits_{a}\left(  q\right)  \right]
$. In view of $\operatorname*{ev}\nolimits_{a}\left(  p\left[  q\right]
\right)  =\left(  p\left[  q\right]  \right)  \left[  a\right]  $ and
$\operatorname*{ev}\nolimits_{a}\left(  q\right)  =q\left[  a\right]  $, we
can rewrite this as $\left(  p\left[  q\right]  \right)  \left[  a\right]
=p\left[  q\left[  a\right]  \right]  $. This proves part \textbf{(a)}.
\medskip

\textbf{(b)} Let $r\in R\left[  x\right]  $. Then, part \textbf{(a)} (applied
to $A=R\left[  x\right]  $ and $a=r$) yields $\left(  p\left[  q\right]
\right)  \left[  r\right]  =p\left[  q\left[  r\right]  \right]  $. Since the
notation $p\circ q$ denotes $p\left[  q\right]  $, we can rewrite this as
$\left(  p\circ q\right)  \circ r=p\circ\left(  q\circ r\right)  $. Thus, part
\textbf{(b)} is proved.
\end{proof}

Theorem \ref{thm.polys1.comp-assoc} \textbf{(b)} shows that composition of
polynomials (i.e., the operation $\circ$ on the polynomial ring $R\left[
x\right]  $) is associative. It also has a neutral element -- namely, the
polynomial $x$, since each $p\in R\left[  x\right]  $ satisfies $p\circ
x=x\circ p=p$. Thus, $R\left[  x\right]  $ is a monoid with respect to this
operation $\circ$.

There is also a multivariate version of Theorem \ref{thm.polys1.comp-assoc}
\textbf{(a)}:

\begin{theorem}
\label{thm.polys1.comp-assoc-multi}Let $p\in R\left[  x_{1},x_{2},\ldots
,x_{n}\right]  $ be a polynomial in $n$ indeterminates, and let $q_{1}%
,q_{2},\ldots,q_{n}\in R\left[  x_{1},x_{2},\ldots,x_{m}\right]  $ be $n$
polynomials in $m$ indeterminates. Let $A$ be any $R$-algebra. Let
$a_{1},a_{2},\ldots,a_{m}\in A$ be $m$ elements of $A$ that mutually commute.
Then, the elements $q_{1}\left[  a_{1},a_{2},\ldots,a_{m}\right]
,\ q_{2}\left[  a_{1},a_{2},\ldots,a_{m}\right]  ,\ \ldots,\ q_{n}\left[
a_{1},a_{2},\ldots,a_{m}\right]  $ of $A$ commute, and we have%
\begin{align*}
&  \left(  p\left[  q_{1},q_{2},\ldots,q_{n}\right]  \right)  \left[
a_{1},a_{2},\ldots,a_{m}\right] \\
&  =p\left[  q_{1}\left[  a_{1},a_{2},\ldots,a_{m}\right]  ,\ q_{2}\left[
a_{1},a_{2},\ldots,a_{m}\right]  ,\ \ldots,\ q_{n}\left[  a_{1},a_{2}%
,\ldots,a_{m}\right]  \right]  .
\end{align*}

\end{theorem}

\begin{exercise}
Prove Theorem \ref{thm.polys1.comp-assoc-multi}.
\end{exercise}

\subsection{\label{sec.polys1.Rx}Univariate polynomials}

\subsubsection{\label{subsec.polys1.Rx.deg}Degrees and coefficients}

Let us now take a closer look at univariate polynomials (which are the
best-behaved among the polynomials).

As we recall, if $p\in R\left[  x\right]  $ is a polynomial, and if
$i\in\mathbb{N}$, then $\left[  x^{i}\right]  p$ is the coefficient of $x^{i}$
in $p$. That is, if $p$ is written as $p=\sum_{j\in\mathbb{N}}p_{j}x^{j}$ with
$p_{j}\in R$, then $\left[  x^{i}\right]  p=p_{i}$.

\begin{definition}
\label{def.polring.univar-deg}Let $p\in R\left[  x\right]  $ be a univariate polynomial.

\begin{enumerate}
\item[\textbf{(a)}] For any $i\in\mathbb{N}$, the coefficient $\left[
x^{i}\right]  p$ is also called the $x^{i}$\textbf{-coefficient} of $p$.

\item[\textbf{(b)}] If $p\neq0$, then the \textbf{degree} of $p$ is defined to
be the largest $i\in\mathbb{N}$ such that $\left[  x^{i}\right]  p\neq0$. The
degree of the zero polynomial $0\in R\left[  x\right]  $ is defined to be
$-\infty$ (a symbol subject to the rules $-\infty<m$ and $-\infty+m=-\infty$
for any $m\in\mathbb{Z}$).

The degree of $p$ is denoted by $\deg p$.

\item[\textbf{(c)}] If $p\neq0$, then the \textbf{leading coefficient} of $p$
is defined to be the coefficient $\left[  x^{\deg p}\right]  p\in R$.

\item[\textbf{(d)}] The polynomial $p$ is said to be \textbf{monic} (or, as
some say, \textbf{normalized}) if it is nonzero and its leading coefficient is
$1$.
\end{enumerate}
\end{definition}

For example, the polynomial
\[
5x^{3}+2x+1\in\mathbb{Q}\left[  x\right]
\]
has degree $3$ and leading coefficient $5$, and hence is not monic (since
$5\neq1$). The polynomial%
\[
\overline{5}x^{3}+\overline{2}x+\overline{1}\in\left(  \mathbb{Z}/n\right)
\left[  x\right]  \ \ \ \ \ \ \ \ \ \ \left(  \text{for a given integer
}n>0\right)
\]
has

\begin{itemize}
\item degree $3$ if $n>5$;

\item degree $1$ if $n=5$ (because if $n=5$, then the $\overline{5}%
x^{3}=\overline{0}x^{3}$ term disappears);

\item degree $3$ if $n=2,3,4$; and

\item degree $-\infty$ if $n=1$ (since all terms disappear if $n=1$).
\end{itemize}

(Degrees are somewhat degenerate for trivial rings: If $R$ is a trivial ring,
then any polynomial in $R\left[  x\right]  $ is $0$ and has degree $-\infty$.)

The polynomial $\left(  1+x\right)  ^{3}=1+3x+3x^{2}+x^{3}$ is monic (i.e.,
has leading coefficient $1$) and has degree $3$.

Let me stress again that the zero polynomial $0=0x^{0}+0x^{1}+0x^{2}+\cdots$
has degree $-\infty$ by definition. This $-\infty$ is not a number, but we
agree that $-\infty$ is smaller than any integer and does not change if you
add an integer to it (i.e., we have $\left(  -\infty\right)  +m=-\infty$ for
any $m\in\mathbb{Z}$) nor if you add itself to it (i.e., we have $\left(
-\infty\right)  +\left(  -\infty\right)  =-\infty$).

Here are some properties of degrees:

\begin{remark}
\label{rmk.polring.univar-deg-leq}Let $n\in\mathbb{N}$. Then,%
\begin{align*}
&  \left\{  f\in R\left[  x\right]  \ \mid\ \deg f\leq n\right\} \\
&  =\left\{  f\in R\left[  x\right]  \ \mid\ f=a_{0}x^{0}+a_{1}x^{1}%
+\cdots+a_{n}x^{n}\text{ for some }a_{0},a_{1},\ldots,a_{n}\in R\right\} \\
&  =\operatorname*{span}\left(  x^{0},x^{1},\ldots,x^{n}\right)  .
\end{align*}
This is clearly an $R$-submodule of $R\left[  x\right]  $.
\end{remark}

\begin{corollary}
\label{cor.polring.univar-degp+q}Let $p,q\in R\left[  x\right]  $. Then,
\begin{align}
\deg\left(  p+q\right)   &  \leq\max\left\{  \deg p,\deg q\right\}
\ \ \ \ \ \ \ \ \ \ \text{and}\label{eq.cor.polring.univar-degp+q.+}\\
\deg\left(  p-q\right)   &  \leq\max\left\{  \deg p,\deg q\right\}  .
\label{eq.cor.polring.univar-degp+q.-}%
\end{align}

\end{corollary}

\begin{proof}
If $p=q=0$, then the claims are obvious. Thus, WLOG assume that at least one
of $p$ and $q$ is nonzero.

Let $n=\max\left\{  \deg p,\deg q\right\}  $. Let $N$ denote the subset
$\left\{  f\in R\left[  x\right]  \ \mid\ \deg f\leq n\right\}  $ of $R\left[
x\right]  $. Then, we know from Remark \ref{rmk.polring.univar-deg-leq} that
$N$ is an $R$-submodule of $R\left[  x\right]  $. Moreover, the definition of
$n$ shows that $\deg p\leq n$, so that $p\in N$. Similarly, $q\in N$. Hence,
$p+q\in N$ (since $N$ is an $R$-submodule of $R\left[  x\right]  $); in other
words, $\deg\left(  p+q\right)  \leq n$. In other words, $\deg\left(
p+q\right)  \leq\max\left\{  \deg p,\deg q\right\}  $ (since $n=\max\left\{
\deg p,\deg q\right\}  $). Similarly, we can find $\deg\left(  p-q\right)
\leq\max\left\{  \deg p,\deg q\right\}  $. This proves Corollary
\ref{cor.polring.univar-degp+q}.
\end{proof}

\begin{remark}
\label{rmk.polring.univar-deg0}The polynomials of degree $\leq0$ are precisely
the constant polynomials -- i.e., the elements of $R$ (embedded into $R\left[
x\right]  $ as explained in Convention \ref{conv.polring.constant-pols}).
\end{remark}

The following proposition collects some properties of products of univariate polynomials:

\begin{proposition}
\label{prop.polring.univar-degpq}Let $p,q\in R\left[  x\right]  $. Then:

\begin{enumerate}
\item[\textbf{(a)}] We have $\deg\left(  pq\right)  \leq\deg p+\deg q$.

\item[\textbf{(b)}] We have $\deg\left(  pq\right)  =\deg p+\deg q$ if
$p\neq0$ and the leading coefficient of $p$ is a unit.

\item[\textbf{(c)}] We have $\deg\left(  pq\right)  =\deg p+\deg q$ if $R$ is
an integral domain.

\item[\textbf{(d)}] If $n,m\in\mathbb{N}$ satisfy $n\geq\deg p$ and $m\geq\deg
q$, then%
\[
\left[  x^{n+m}\right]  \left(  pq\right)  =\left[  x^{n}\right]  \left(
p\right)  \cdot\left[  x^{m}\right]  \left(  q\right)  .
\]

\item[\textbf{(e)}] If $pq=0$ and $p\neq0$ and if the leading coefficient of
$p$ is a unit, then $q=0$.

\item[\textbf{(f)}] If $p$ and $q$ are monic, then $pq$ is again monic and has
degree $\deg\left(  pq\right)  =\deg p+\deg q$.
\end{enumerate}
\end{proposition}

\begin{corollary}
\label{cor.polring.univar-intdom}If $R$ is an integral domain, then the
polynomial ring $R\left[  x\right]  $ is an integral domain.
\end{corollary}

\begin{proof}
[Proof of Proposition \ref{prop.polring.univar-degpq}.]For parts \textbf{(a)},
\textbf{(b)}, \textbf{(c)} and \textbf{(d)}, we will give an informal
\textquotedblleft proof by example\textquotedblright. Rigorous arguments can
be found in various places\footnote{Actually, in fewer places than one might
expect. Apparently most authors just handwave them away or leave them to the
reader (e.g., Bourbaki writes about part \textbf{(a)} that \textquotedblleft
the proof is immediate\textquotedblright). An explicit proof of part
\textbf{(c)} appears in \cite[Proposition 9.3.6]{Mileti20}, and reading it
between the lines also reveals proofs of \textbf{(a)}, \textbf{(b)} and
\textbf{(d)}. A while ago I have written up proofs for parts \textbf{(a)} and
\textbf{(d)} in \cite{regpol} (where they appear as parts \textbf{(a)} and
\textbf{(b)} of Lemma 3.12), albeit only in the particular case when $p$ is
monic (but the proofs can easily be generalized). A generalization of parts
\textbf{(b)} and \textbf{(c)} also appears in \cite[Proposition (1.3.12)]%
{ChaLoi21}.}.

Let $p$ and $q$ be two polynomials of degrees $\deg p=2$ and $\deg q=3$. Write
$p$ and $q$ as $p=ax^{2}+bx+c$ and $q=dx^{3}+ex^{2}+fx+g$ (with $a,b,c,\ldots
,g\in R$). Then,%
\begin{align}
pq  &  =\left(  ax^{2}+bx+c\right)  \left(  dx^{3}+ex^{2}+fx+g\right)
\nonumber\\
&  =adx^{5}+\left(  \text{lower powers of }x\right)  .
\label{pf.prop.polring.univar-degpq.eq.1}%
\end{align}
Thus, $\deg\left(  pq\right)  \leq5=2+3=\deg p+\deg q$. This proves
Proposition \ref{prop.polring.univar-degpq} \textbf{(a)}.

Moreover, $a\neq0$ (since $\deg p=2$) and $d\neq0$ (since $\deg q=3$). If $R$
is an integral domain, then this entails $ad\neq0$ and therefore $\deg\left(
pq\right)  =5$ (by (\ref{pf.prop.polring.univar-degpq.eq.1})). This proves
Proposition \ref{prop.polring.univar-degpq} \textbf{(c)}. On the other hand,
if $a$ is a unit, then we also have $ad\neq0$ (because otherwise, we would
have $ad=0$ and thus $a^{-1}\underbrace{ad}_{=0}=0$, which would contradict
$a^{-1}ad=d\neq0$) and therefore $\deg\left(  pq\right)  =5$ (by
(\ref{pf.prop.polring.univar-degpq.eq.1})). This proves Proposition
\ref{prop.polring.univar-degpq} \textbf{(b)} (since $a$ is the leading
coefficient of $p$).

The equality (\ref{pf.prop.polring.univar-degpq.eq.1}) shows that the
coefficient of $x^{5}$ in $pq$ is $ad$, and no higher powers of $x$ than
$x^{5}$ appear in $pq$. That is, we have $\left[  x^{5}\right]  \left(
pq\right)  =\underbrace{a}_{=\left[  x^{2}\right]  \left(  p\right)
}\underbrace{d}_{=\left[  x^{3}\right]  \left(  q\right)  }=\left[
x^{2}\right]  \left(  p\right)  \cdot\left[  x^{3}\right]  \left(  q\right)
$, and we have $\left[  x^{i}\right]  \left(  pq\right)  =0$ for all $i>5$.
This quickly yields Proposition \ref{prop.polring.univar-degpq} \textbf{(d)}.

To prove Proposition \ref{prop.polring.univar-degpq} \textbf{(e)}, we assume
the contrary. Thus, $pq=0$ and $p\neq0$ and the leading coefficient of $p$ is
a unit, but $q\neq0$. Then, Proposition \ref{prop.polring.univar-degpq}
\textbf{(b)} yields $\deg\left(  pq\right)  =\underbrace{\deg p}_{\geq
0}+\underbrace{\deg q}_{\geq0}\geq0$. However, $pq=0$, so $\deg\left(
pq\right)  =\deg0=-\infty<0$. These two inequalities clearly contradict each
other, and our proof of Proposition \ref{prop.polring.univar-degpq}
\textbf{(e)} is complete.

Proposition \ref{prop.polring.univar-degpq} \textbf{(f)} follows from parts
\textbf{(a)} and \textbf{(d)}: Indeed, if $p$ and $q$ are monic, then part
\textbf{(d)} (applied to $n=\deg p$ and $m=\deg q$) yields%
\[
\left[  x^{\deg p+\deg q}\right]  \left(  pq\right)  =\underbrace{\left[
x^{\deg p}\right]  \left(  p\right)  }_{\substack{=1\\\text{(since }p\text{ is
monic)}}}\cdot\underbrace{\left[  x^{\deg q}\right]  \left(  q\right)
}_{\substack{=1\\\text{(since }q\text{ is monic)}}}=1,
\]
whereas part \textbf{(a)} yields $\deg\left(  pq\right)  \leq\deg p+\deg q$.
Combining these two facts, we conclude that $pq$ is monic of degree $\deg
p+\deg q$. Thus, Proposition \ref{prop.polring.univar-degpq} \textbf{(f)} is proved.
\end{proof}

\begin{proof}
[Proof of Corollary \ref{cor.polring.univar-intdom}.]Assume that $R$ is an
integral domain. Let $p,q\in R\left[  x\right]  $ be nonzero. Then,
Proposition \ref{prop.polring.univar-degpq} \textbf{(c)} yields $\deg\left(
pq\right)  =\underbrace{\deg p}_{\geq0}+\underbrace{\deg q}_{\geq0}\geq0$, and
thus $pq\neq0$ (since $pq=0$ would yield $\deg\left(  pq\right)
=\deg0=-\infty<0$). Thus, we have shown that $pq\neq0$ for any nonzero $p,q\in
R\left[  x\right]  $. In other words, $R\left[  x\right]  $ is an integral domain.
\end{proof}

\begin{corollary}
\label{cor.polring.univar-degpk}Let $p_{1},p_{2},\ldots,p_{k}\in R\left[
x\right]  $ be $k$ nonzero polynomials. Then:

\begin{enumerate}
\item[\textbf{(a)}] If $p_{1},p_{2},\ldots,p_{k}$ are monic, then the product
$p_{1}p_{2}\cdots p_{k}$ is again monic and has degree $\deg\left(  p_{1}%
p_{2}\cdots p_{k}\right)  =\deg p_{1}+\deg p_{2}+\cdots+\deg p_{k}$.

\item[\textbf{(b)}] If $R$ is an integral domain, then $p_{1}p_{2}\cdots
p_{k}$ is nonzero and has degree $\deg\left(  p_{1}p_{2}\cdots p_{k}\right)
=\deg p_{1}+\deg p_{2}+\cdots+\deg p_{k}$.
\end{enumerate}
\end{corollary}

\begin{proof}
Part \textbf{(a)} follows by induction on $k$ using Proposition
\ref{prop.polring.univar-degpq} \textbf{(f)}.

Part \textbf{(b)} follows by induction on $k$ using Proposition
\ref{prop.polring.univar-degpq} \textbf{(c)}.
\end{proof}

If $R$ is not an integral domain, then degrees of polynomials over $R$ can
behave rather strangely. For example, over $\mathbb{Z}/4$, we have%
\[
\left(  \overline{1}+\overline{2}x\right)  ^{2}=\overline{1}+\overline
{4}x+\overline{4}x^{2}=\overline{1}\ \ \ \ \ \ \ \ \ \ \left(  \text{since
}\overline{4}=\overline{0}\right)  .
\]
So the degree of a polynomial can decrease when it is squared!

\begin{exercise}
\label{exe.21hw3.5}Let $n\in\mathbb{N}$. Prove that we have $x^{2}+x+1\mid
x^{2n}+x^{n}+1$ in the polynomial ring $\mathbb{Z}\left[  x\right]  $ if and
only if $3\nmid n$ in $\mathbb{Z}$. \medskip

[\textbf{Hint:} First show that $x^{3}\equiv1\operatorname{mod}x^{2}+x+1$ in
the ring $\mathbb{Z}\left[  x\right]  $. Here, we are using the notation
$a\equiv b\operatorname{mod}c$ (spoken \textquotedblleft$a$ is congruent to
$b$ modulo $c$\textquotedblright) for $c\mid a-b$ whenever $a,b,c$ are three
elements of a commutative ring $R$. Congruences in $R$ are a straightforward
generalization of congruences of integers (which are known from elementary
number theory), and behave just as nicely; in particular, they can be added,
subtracted and multiplied.]
\end{exercise}

\begin{fineprint}
The equality $\left(  \overline{1}+\overline{2}x\right)  ^{2}=\overline{1}$ in
$\mathbb{Z}/4$ that we observed above is surprising not only because of the
strange \textquotedblleft loss of degree\textquotedblright\ that happens when
$\overline{1}+\overline{2}x$ is squared, but also for another reason: It shows
that the non-constant polynomial $\overline{1}+\overline{2}x$ in $\left(
\mathbb{Z}/4\right)  \left[  x\right]  $ is actually a unit of $\left(
\mathbb{Z}/4\right)  \left[  x\right]  $ ! As Proposition
\ref{prop.polring.univar-degpq} \textbf{(c)} explains, this cannot happen for
polynomials over an integral domain. The following exercise characterizes
precisely when this happens:
\end{fineprint}

\begin{exercise}
\label{exe.21hw4.10}Let $R$ be a commutative ring. Let $f\in R\left[
x\right]  $ be a polynomial. Recall that the notation $\left[  x^{i}\right]
f$ stands for the coefficient of the monomial $x^{i}$ in $f$.

Prove that $f$ is a unit of the ring $R\left[  x\right]  $ if and only if

\begin{itemize}
\item the coefficient $\left[  x^{0}\right]  f$ is a unit of $R$, and

\item all the remaining coefficients $\left[  x^{1}\right]  f,\left[
x^{2}\right]  f,\left[  x^{3}\right]  f,\ldots$ of $f$ are nilpotent.
\end{itemize}

[\textbf{Hint:} Exercise \ref{exe.21hw1.7} \textbf{(c)} shows that the
nilpotent elements of $R\left[  x\right]  $ form an ideal, whereas Exercise
\ref{exe.21hw2.1} \textbf{(b)} shows that the difference of a unit and a
nilpotent element is always a unit (in $R\left[  x\right]  $). This should
help with the \textquotedblleft if\textquotedblright\ direction. For the
\textquotedblleft only if\textquotedblright\ direction, let $f=f_{0}%
x^{0}+f_{1}x^{1}+\cdots+f_{n}x^{n}\in R\left[  x\right]  $ be a unit and
$g=g_{0}x^{0}+g_{1}x^{1}+\cdots+g_{m}x^{m}\in R\left[  x\right]  $ be its
inverse. Assume that $n>0$. Use induction on $r$ to show that $f_{n}%
^{r+1}g_{m-r}=0$ for each $r\in\left\{  0,1,\ldots,m\right\}  $. Use this to
conclude that $f_{n}$ is nilpotent.]
\end{exercise}

\begin{exercise}
Note that Proposition \ref{prop.polring.univar-degpq} \textbf{(e)} can be
restated as follows: Any nonzero polynomial $p\in R\left[  x\right]  $ whose
leading coefficient is a unit is regular (see Definition
\ref{def.rings.reg.reg}). More generally, prove the following:

\begin{enumerate}
\item[\textbf{(a)}] If $p\in R\left[  x\right]  $ is a nonzero polynomial
whose leading coefficient is regular, then $p$ is regular.

\item[\textbf{(b)}] If $p,q\in R\left[  x\right]  $ are two polynomials such
that at least one of $p$ and $q$ has a regular leading coefficient, then
$\deg\left(  pq\right)  =\deg p+\deg q$.
\end{enumerate}
\end{exercise}

\subsubsection{Division with remainder}

The most important feature of univariate polynomials is division with remainder:

\begin{theorem}
[Division-with-remainder theorem for polynomials]%
\label{thm.polring.univar-quorem}Let $b\in R\left[  x\right]  $ be a nonzero
polynomial whose leading coefficient is a unit. Let $a\in R\left[  x\right]  $
be any polynomial.

\begin{enumerate}
\item[\textbf{(a)}] Then, there is a \textbf{unique} pair $\left(  q,r\right)
$ of polynomials in $R\left[  x\right]  $ such that
\[
a=qb+r\ \ \ \ \ \ \ \ \ \ \text{and}\ \ \ \ \ \ \ \ \ \ \deg r<\deg b.
\]

\item[\textbf{(b)}] Moreover, this pair satisfies $\deg q\leq\deg a-\deg b$.
\end{enumerate}
\end{theorem}

The polynomials $q$ and $r$ in Theorem \ref{thm.polring.univar-quorem} are
called the \textbf{quotient} and the \textbf{remainder} obtained when dividing
$a$ by $b$. Note that if $\deg a<\deg b$, then the quotient $q$ is $0$ whereas
the remainder $r$ is $a$. The quotient and the remainder become interesting
when $\deg a\geq\deg b$.

\begin{example}
Let $R=\mathbb{Z}$ and $a=3x^{4}+x^{2}+6x-2$ and $b=x^{2}-3x+1$. Then, Theorem
\ref{thm.polring.univar-quorem} \textbf{(a)} shows that there is a
\textbf{unique} pair $\left(  q,r\right)  $ of polynomials in $R\left[
x\right]  $ such that
\[
a=qb+r\ \ \ \ \ \ \ \ \ \ \text{and}\ \ \ \ \ \ \ \ \ \ \deg r<\deg b.
\]
This pair $\left(  q,r\right)  $ is $\left(  3x^{2}+9x+25,\ 72x-27\right)  $.
Indeed, if you set $q=3x^{2}+9x+25$ and $r=72x-27$, then the equality $a=qb+r$
can be verified by a straightforward computation, whereas the inequality $\deg
r<\deg b$ is obvious. Thus, the quotient obtained when dividing $a$ by $b$ is
$q=3x^{2}+9x+25$, and the remainder is $r=72x-27$.
\end{example}

\begin{example}
Let us give a \textbf{non-example:} Let $R=\mathbb{Z}$ and $b=2$ (a constant
polynomial) and $a=x$. The leading coefficient of $b$ is not a unit (since $2$
is not a unit in $\mathbb{Z}$), so we don't expect Theorem
\ref{thm.polring.univar-quorem} to hold. And indeed: we cannot write $a=qb+r$
with $\deg r<\deg b$. Indeed, this would mean $x=q\cdot2+r$ with $\deg r<0$
(since the constant polynomial $2$ has degree $\deg2=0$); but this is
impossible, since this would entail $x=q\cdot2$, which would contradict the
fact that $x$ is not divisible by $2$.
\end{example}

Instead of proving Theorem \ref{thm.polring.univar-quorem} directly, we will
first show the particular case in which $b$ is required to be monic, and then
use it to derive the general case. The particular case is the following lemma:

\begin{lemma}
\label{lem.polring.univar-quorem.monic}Let $b\in R\left[  x\right]  $ be a
monic polynomial. Let $a\in R\left[  x\right]  $ be any polynomial.

\begin{enumerate}
\item[\textbf{(a)}] Then, there is a \textbf{unique} pair $\left(  q,r\right)
$ of polynomials in $R\left[  x\right]  $ such that
\[
a=qb+r\ \ \ \ \ \ \ \ \ \ \text{and}\ \ \ \ \ \ \ \ \ \ \deg r<\deg b.
\]

\item[\textbf{(b)}] Moreover, this pair satisfies $\deg q\leq\deg a-\deg b$.
\end{enumerate}
\end{lemma}

\begin{proof}
\textbf{(a)} Again, we shall give a proof by example. (For a rigorous proof,
see \cite[Theorem 3.16 and Lemma 3.19]{regpol} or \cite[Theorem 3.6.4]{Ford22}
or \cite[Theorem (1.3.15)]{ChaLoi21} or \cite[Proposition 1.12]{Knapp1} or
\cite[\S 9.2, Theorem 3]{DumFoo04}. Note that some of these sources assume
that $R$ is a field; however, the proofs easily adapt to our general
case.\footnote{Also, \cite[Theorem 9.3.10 part 1]{Mileti20} gives a good proof
of part \textbf{(a)} (even for Theorem \ref{thm.polring.univar-quorem}
\textbf{(a)}).})

We are doing a proof by example, so let us assume that $\deg a=3$ and $\deg
b=2$. Thus, we can write $a$ and $b$ as $a=cx^{3}+dx^{2}+ex+f$ and
$b=x^{2}+gx+h$ for some $c,d,e,f,g,h\in R$ (since $b$ is monic).

\begin{noncompile}
Now, we observe that if a polynomial $p$ satisfies $\deg p\geq\deg b$, then we
can decrease the degree of $p$ by subtracting an appropriate multiple of $b$
from it (specifically: if the leading term of $p$ is $\alpha x^{\deg p}$, then
subtracting $\alpha x^{\deg p-\deg b}b$ from $p$ results in a polynomial
$p-\alpha x^{\deg p-\deg b}b$ that has lower degree than $p$, because the
leading terms of $p$ and $\alpha x^{\deg p-\deg b}b$ cancel each other).
\end{noncompile}

Now, we repeatedly subtract appropriate multiples of $b$ from $a$ in order to
decrease its degree:%
\begin{align*}
&  a=cx^{3}+dx^{2}+ex+f\\
\Longrightarrow\  &  a-\left(  cx\right)  b=\left(  d-cg\right)  x^{2}+\left(
e-ch\right)  x+f\\
&  \ \ \ \ \ \ \ \ \ \ \ \ \ \ \ \ \ \ \ \ \left(  \text{here, we have
subtracted }\left(  cx\right)  b\text{ to kill off the }cx^{3}\text{
term}\right) \\
\Longrightarrow\  &  a-\left(  cx\right)  b-\left(  d-cg\right)  b=\left(
e-ch-\left(  d-cg\right)  g\right)  x+\left(  f-\left(  d-cg\right)  h\right)
\\
&  \ \ \ \ \ \ \ \ \ \ \ \ \ \ \ \ \ \ \ \ \left(
\begin{array}
[c]{c}%
\text{here, we have subtracted }\left(  d-cg\right)  b\text{ to kill off}\\
\text{the }\left(  d-cg\right)  x^{2}\text{ term}%
\end{array}
\right)  .
\end{align*}
Thus,%
\begin{align*}
a  &  =\left(  cx\right)  b+\left(  d-cg\right)  b+\left(  e-ch-\left(
d-cg\right)  g\right)  x+\left(  f-\left(  d-cg\right)  h\right) \\
&  =\left(  cx+\left(  d-cg\right)  \right)  b+\left(  e-ch-\left(
d-cg\right)  g\right)  x+\left(  f-\left(  d-cg\right)  h\right)  .
\end{align*}
Setting $q:=cx+\left(  d-cg\right)  $ and $r:=\left(  e-ch-\left(
d-cg\right)  g\right)  x+\left(  f-\left(  d-cg\right)  h\right)  $, we can
rewrite this as
\[
a=qb+r.
\]
Note that $\deg r<\deg b$ (since any polynomial of degree $\geq\deg b$ could
still be reduced further by subtracting a multiple of $b$ from it).

Thus we have found a pair $\left(  q,r\right)  $ of polynomials satisfying%
\[
a=qb+r\ \ \ \ \ \ \ \ \ \ \text{and}\ \ \ \ \ \ \ \ \ \ \deg r<\deg b.
\]
It remains to prove its uniqueness. In other words, we have to prove that if
$\left(  q_{1},r_{1}\right)  $ and $\left(  q_{2},r_{2}\right)  $ are two
pairs of polynomials satisfying%
\begin{align*}
a  &  =q_{1}b+r_{1}\ \ \ \ \ \ \ \ \ \ \text{and}\ \ \ \ \ \ \ \ \ \ \deg
r_{1}<\deg b\ \ \ \ \ \ \ \ \ \ \text{and}\\
a  &  =q_{2}b+r_{2}\ \ \ \ \ \ \ \ \ \ \text{and}\ \ \ \ \ \ \ \ \ \ \deg
r_{2}<\deg b,
\end{align*}
then $\left(  q_{1},r_{1}\right)  =\left(  q_{2},r_{2}\right)  $. To prove
this, we fix two such pairs $\left(  q_{1},r_{1}\right)  $ and $\left(
q_{2},r_{2}\right)  $. Then, we have%
\[
q_{1}b+r_{1}=a=q_{2}b+r_{2},
\]
so that $r_{1}-r_{2}=q_{2}b-q_{1}b=b\left(  q_{2}-q_{1}\right)  $. Hence,
$b\left(  q_{2}-q_{1}\right)  =r_{1}-r_{2}$, so that%
\begin{align*}
\deg\left(  b\left(  q_{2}-q_{1}\right)  \right)   &  =\deg\left(  r_{1}%
-r_{2}\right)  \leq\max\left\{  \deg r_{1},\deg r_{2}\right\}
\ \ \ \ \ \ \ \ \ \ \left(  \text{by (\ref{eq.cor.polring.univar-degp+q.-}%
)}\right) \\
&  <\deg b\ \ \ \ \ \ \ \ \ \ \left(  \text{since }\deg r_{1}<\deg b\text{ and
}\deg r_{2}<\deg b\right)  .
\end{align*}
However, the leading coefficient of $b$ is a unit\footnote{Indeed, this
coefficient is $1$, since $b$ is monic.}. Hence, if the polynomial
$q_{2}-q_{1}$ was nonzero, then Proposition \ref{prop.polring.univar-degpq}
\textbf{(b)} would entail%
\[
\deg\left(  b\left(  q_{2}-q_{1}\right)  \right)  =\deg b+\underbrace{\deg
\left(  q_{2}-q_{1}\right)  }_{\geq0}\geq\deg b,
\]
which would contradict $\deg\left(  b\left(  q_{2}-q_{1}\right)  \right)
<\deg b$. So $q_{2}-q_{1}$ must be zero. In other words, $q_{2}-q_{1}=0$, so
that $q_{1}=q_{2}$. Moreover, $r_{1}-r_{2}=b\underbrace{\left(  q_{2}%
-q_{1}\right)  }_{=0}=0$, so that $r_{1}=r_{2}$. Hence, $\left(  q_{1}%
,r_{1}\right)  =\left(  q_{2},r_{2}\right)  $. This completes the proof of the
uniqueness of $\left(  q,r\right)  $. Thus, Lemma
\ref{lem.polring.univar-quorem.monic} \textbf{(a)} is proved. \medskip

\textbf{(b)} You can obtain Lemma \ref{lem.polring.univar-quorem.monic}
\textbf{(b)} by a careful analysis of the construction of the pair $\left(
q,r\right)  $ in our proof of part \textbf{(a)}. Indeed, each of the terms of
$q$ was originally a factor that we multiplied with $b$ in order to reduce
$a$; however, the highest power of $x$ in $a$ was $x^{\deg a}$, so the factors
we used did not contain any powers of $x$ higher than $x^{\deg a-\deg b}$.

Alternatively, you can prove Lemma \ref{lem.polring.univar-quorem.monic}
\textbf{(b)} independently of part \textbf{(a)}: Let $\left(  q,r\right)  $ be
a pair of polynomials in $R\left[  x\right]  $ such that
\[
a=qb+r\ \ \ \ \ \ \ \ \ \ \text{and}\ \ \ \ \ \ \ \ \ \ \deg r<\deg b.
\]
We must prove that $\deg q\leq\deg a-\deg b$. Assume the contrary. Thus, $\deg
q>\deg a-\deg b$. Therefore, in particular, $q\neq0$ (since $q=0$ would entail
$\deg q=\deg0=-\infty\leq\deg a-\deg b$), so that $\deg q\geq0$. However, the
leading coefficient of $b$ is a unit\footnote{Indeed, this coefficient is $1$,
since $b$ is monic.}; thus, Proposition \ref{prop.polring.univar-degpq}
\textbf{(b)} yields that
\[
\deg\left(  bq\right)  =\deg b+\deg q>\deg a\ \ \ \ \ \ \ \ \ \ \left(
\text{since }\deg q>\deg a-\deg b\right)  .
\]
Also,%
\[
\deg\left(  bq\right)  =\deg b+\underbrace{\deg q}_{\geq0}\geq\deg b>\deg
r\ \ \ \ \ \ \ \ \ \ \left(  \text{since }\deg r<\deg b\right)  .
\]
Combining these two inequalities, we obtain%
\[
\deg\left(  bq\right)  >\max\left\{  \deg a,\ \deg r\right\}  .
\]

But from $a=qb+r$, we obtain $a-r=qb=bq$, so that $bq=a-r$. Hence,%
\[
\deg\left(  bq\right)  =\deg\left(  a-r\right)  \leq\max\left\{  \deg a,\ \deg
r\right\}  \ \ \ \ \ \ \ \ \ \ \left(  \text{by
(\ref{eq.cor.polring.univar-degp+q.-})}\right)  ,
\]
which contradicts $\deg\left(  bq\right)  >\max\left\{  \deg a,\ \deg
r\right\}  $. This contradiction shows that our assumption was wrong; thus,
Lemma \ref{lem.polring.univar-quorem.monic} \textbf{(b)} is proven.
\end{proof}

We can now prove the general case:

\begin{proof}
[Proof of Theorem \ref{thm.polring.univar-quorem}.]\textbf{(a)} As in our
above proof of Lemma \ref{lem.polring.univar-quorem.monic} \textbf{(a)}, we
can show that the pair $\left(  q,r\right)  $ is unique. It remains to show
that this pair exists.

Let $u$ be the leading coefficient of $b$. Then, $u$ is a unit (by
assumption), and thus has an inverse $u^{-1}$. Scaling the polynomial $b$ by
$u^{-1}$ results in a new polynomial $u^{-1}b$, which has leading coefficient
$u^{-1}u$ (since the leading coefficient $u$ of $b$ gets multiplied by
$u^{-1}$) and thus is monic (since $u^{-1}u=1$). Hence, we can apply Lemma
\ref{lem.polring.univar-quorem.monic} \textbf{(a)} to $u^{-1}b$ instead of
$b$. As a result, we conclude that there is a \textbf{unique} pair $\left(
q,r\right)  $ of polynomials in $R\left[  x\right]  $ such that
\[
a=q\left(  u^{-1}b\right)  +r\ \ \ \ \ \ \ \ \ \ \text{and}%
\ \ \ \ \ \ \ \ \ \ \deg r<\deg\left(  u^{-1}b\right)  .
\]
Let us denote this pair $\left(  q,r\right)  $ by $\left(  \widetilde{q}%
,\widetilde{r}\right)  $. Thus, $\left(  \widetilde{q},\widetilde{r}\right)  $
is a pair of polynomials in $R\left[  x\right]  $ and satisfies%
\[
a=\widetilde{q}\left(  u^{-1}b\right)  +\widetilde{r}%
\ \ \ \ \ \ \ \ \ \ \text{and}\ \ \ \ \ \ \ \ \ \ \deg\widetilde{r}%
<\deg\left(  u^{-1}b\right)  .
\]

Now,
\[
a=\widetilde{q}\left(  u^{-1}b\right)  +\widetilde{r}=\left(  u^{-1}%
\widetilde{q}\right)  b+\widetilde{r}%
\]
and $\deg\widetilde{r}<\deg\left(  u^{-1}b\right)  \leq\deg b$ (because
scaling a polynomial cannot increase its degree). Hence, there is a pair
$\left(  q,r\right)  $ of polynomials in $R\left[  x\right]  $ such that
\[
a=qb+r\ \ \ \ \ \ \ \ \ \ \text{and}\ \ \ \ \ \ \ \ \ \ \deg r<\deg b
\]
(namely, the pair $\left(  q,r\right)  =\left(  u^{-1}\widetilde{q}%
,\widetilde{r}\right)  $). Since we have already shown that such a pair is
unique, we thus have finished proving Theorem \ref{thm.polring.univar-quorem}
\textbf{(a)}. \medskip

\textbf{(b)} Our above proof of Lemma \ref{lem.polring.univar-quorem.monic}
\textbf{(b)} (specifically, the \textquotedblleft
alternative\textquotedblright\ proof that is independent of part \textbf{(a)})
applies to Theorem \ref{thm.polring.univar-quorem} \textbf{(b)} as well.
\end{proof}

We record an automatic corollary of Theorem \ref{thm.polring.univar-quorem}:

\begin{corollary}
\label{cor.polring.univar-quorem.field}Let $F$ be a field. Let $b\in F\left[
x\right]  $ be any nonzero polynomial. Let $a\in F\left[  x\right]  $ be any polynomial.

\begin{enumerate}
\item[\textbf{(a)}] Then, there is a \textbf{unique} pair $\left(  q,r\right)
$ of polynomials in $F\left[  x\right]  $ such that
\[
a=qb+r\ \ \ \ \ \ \ \ \ \ \text{and}\ \ \ \ \ \ \ \ \ \ \deg r<\deg b.
\]

\item[\textbf{(b)}] Moreover, this pair satisfies $\deg q\leq\deg a-\deg b$.
\end{enumerate}
\end{corollary}

\begin{proof}
The polynomial $b$ is nonzero; thus, its leading coefficient is a unit (since
any nonzero element of the field $F$ is a unit). Hence, Theorem
\ref{thm.polring.univar-quorem} applies (to $R=F$).
\end{proof}

The following simple proposition is the polynomial analogue of the classical
fact that a positive integer $b$ divides an integer $a$ if and only if the
remainder that $a$ leaves when divided by $b$ is $0$:

\begin{proposition}
\label{prop.polring.univar-rem0}Let $b\in R\left[  x\right]  $ be a nonzero
polynomial whose leading coefficient is a unit. Let $a\in R\left[  x\right]  $
be any polynomial. Let $q$ and $r$ be the quotient and the remainder obtained
when dividing $a$ by $b$. Then, we have the logical equivalence $\left(
r=0\right)  \ \Longleftrightarrow\ \left(  b\mid a\text{ in }R\left[
x\right]  \right)  $.
\end{proposition}

\begin{proof}
The definition of quotient and remainder yields $a=qb+r$. Hence, if $r=0$,
then $a=qb+\underbrace{r}_{=0}=qb$ and thus $b\mid a$ in $R\left[  x\right]
$. This proves the \textquotedblleft$\Longrightarrow$\textquotedblright%
\ direction of the required equivalence. It thus remains to prove the
\textquotedblleft$\Longleftarrow$\textquotedblright\ direction.

So we assume that $b\mid a$ in $R\left[  x\right]  $. We need to show that
$r=0$.

We have assumed $b\mid a$ in $R\left[  x\right]  $. In other words, there
exists a $c\in R\left[  x\right]  $ such that $a=cb$. Consider this $c$. We
have $a=cb=bc=bc+0$ and $\deg0=-\infty<\deg b$. Thus, $\left(  c,0\right)  $
is a pair $\left(  \widetilde{q},\widetilde{r}\right)  $ of polynomials in
$R\left[  x\right]  $ such that $a=\widetilde{q}b+\widetilde{r}$ and
$\deg\widetilde{r}<\deg b$. But $\left(  q,r\right)  $ is also such a pair (by
the definition of quotient and remainder). However, Theorem
\ref{thm.polring.univar-quorem} \textbf{(a)} shows that there is a
\textbf{unique} such pair. In particular, any two such pairs must be
identical. Thus, the two pairs $\left(  q,r\right)  $ and $\left(  c,0\right)
$ must be identical. That is, we have $q=c$ and $r=0$. In particular, $r=0$;
this completes the proof of the \textquotedblleft$\Longleftarrow
$\textquotedblright\ direction. Proposition \ref{prop.polring.univar-rem0} is
thus proven.
\end{proof}

Note that our above proofs of Lemma \ref{lem.polring.univar-quorem.monic}
\textbf{(a)} and of Theorem \ref{thm.polring.univar-quorem} \textbf{(a)} give
an easy and fairly quick recursive algorithm for computing the $q$ and $r$ in
Theorem \ref{thm.polring.univar-quorem}. This algorithm is known as
\textbf{polynomial division}.

\begin{exercise}
Let $R$ be any commutative ring. Let $n\in\mathbb{N}$.

\begin{enumerate}
\item[\textbf{(a)}] Find the quotient and the remainder obtained when dividing
$\left(  x+1\right)  ^{n}$ by $x$ (in the polynomial ring $R\left[  x\right]
$).

\item[\textbf{(b)}] Find the quotient and the remainder obtained when dividing
$x^{n}$ by $x-1$.

\item[\textbf{(c)}] Find the remainder obtained when dividing $\left(
x+1\right)  ^{n}$ by $x-1$.
\end{enumerate}

[\textbf{Hint:} \textquotedblleft Finding\textquotedblright\ a polynomial here
means computing its coefficients. For instance, in part \textbf{(a)}, the
coefficients of the quotient will be certain binomial coefficients. I am
deliberately not asking for the quotient in part \textbf{(c)}, since I don't
know a closed form for its coefficients that doesn't use summation signs.]
\end{exercise}

\begin{exercise}
Let $R$ be any commutative ring. Let $n\in\mathbb{N}$.

Prove that the remainder obtained when dividing $x^{n}$ by $\left(
x-1\right)  ^{2}$ in the polynomial ring $R\left[  x\right]  $ is $nx-n+1$.
\end{exercise}

\begin{exercise}
Let $R$ be any commutative ring. Let $n$ be a positive integer.

\begin{enumerate}
\item[\textbf{(a)}] Prove that $\left(  x-1\right)  ^{3}\mid x^{2n}%
-nx^{n+1}+nx^{n-1}-1$ in $R\left[  x\right]  $.

\item[\textbf{(b)}] Prove that $x+1\mid x^{2n}-nx^{n+1}+nx^{n-1}-1$ in
$R\left[  x\right]  $.

\item[\textbf{(c)}] Prove that $\left(  x+1\right)  ^{3}\mid x^{2n}%
-nx^{n+1}+nx^{n-1}-1$ in $R\left[  x\right]  $ if $n$ is odd.
\end{enumerate}
\end{exercise}

\begin{exercise}
Let $R$ be any commutative ring. Let $n\in\mathbb{N}$.

In terms of the Fibonacci numbers (Definition \ref{def.fibonacci.fib}), find
the quotient and the remainder obtained when dividing $x^{n}$ by $x^{2}-x-1$.
\medskip

[\textbf{Hint:} Compute them (e.g.) for $n=10$, and prove the pattern you discover.]
\end{exercise}

\begin{exercise}
Let $a\in\mathbb{Z}\left[  x\right]  $ and $b\in\mathbb{Z}\left[  x\right]  $
be two polynomials, with $b$ being nonzero. Without requiring anything about
the leading coefficient of $b$, we cannot apply Theorem
\ref{thm.polring.univar-quorem}, so we don't know whether there exists a pair
$\left(  q,r\right)  $ of polynomials in $\mathbb{Z}\left[  x\right]  $ such
that
\[
a=qb+r\ \ \ \ \ \ \ \ \ \ \text{and}\ \ \ \ \ \ \ \ \ \ \deg r<\deg b.
\]
Nevertheless, such a pair might exist.

\begin{enumerate}
\item[\textbf{(a)}] Does such a pair exist when $a=3x^{2}+1$ and $b=3x-1$ ?

\item[\textbf{(b)}] Does such a pair exist when $a=3x^{3}+1$ and $b=3x-1$ ?

\item[\textbf{(c)}] Does such a pair exist when $a=3x^{2}+2x$ and $b=3x-1$ ?
\end{enumerate}

(Make sure to prove your claims!)
\end{exercise}

The following exercise partially generalizes Theorem
\ref{thm.polring.univar-quorem} to the case when the leading coefficient of
$b$ is not (necessarily) a unit:

\begin{exercise}
Let $R$ be any commutative ring. Let $b\in R\left[  x\right]  $ be a nonzero
polynomial, and let $\lambda\in R$ be its leading coefficient. Let $a\in
R\left[  x\right]  $ be any polynomial such that $\deg a\geq\deg b$. Prove
that there is a pair $\left(  q,r\right)  $ of polynomials in $R\left[
x\right]  $ such that
\begin{align*}
\lambda^{\deg a-\deg b+1}a  &  =qb+r\ \ \ \ \ \ \ \ \ \ \text{and}\\
\deg r  &  <\deg b\ \ \ \ \ \ \ \ \ \ \text{and}\ \ \ \ \ \ \ \ \ \ \deg
q\leq\deg a-\deg b.
\end{align*}
(Note that we can no longer claim that this pair $\left(  q,r\right)  $ is unique.)
\end{exercise}

\subsubsection{Roots}

We shall now discuss roots of polynomials.

\begin{definition}
\label{def.polring.univar-root}Let $A$ be an $R$-algebra. Let $f\in R\left[
x\right]  $. An element $a\in A$ is said to be a \textbf{root} of $f$ if
$f\left(  a\right)  =0$ (that is, $f\left[  a\right]  =0$).
\end{definition}

This is a rather wide notion of roots. For example, the matrix $\left(
\begin{array}
[c]{cc}%
0 & 1\\
0 & 0
\end{array}
\right)  \in\mathbb{Q}^{2\times2}$ is a root of the polynomial $x^{2}%
\in\mathbb{Q}\left[  x\right]  $, since the square of this matrix is $0$. For
another example, the matrix $\left(
\begin{array}
[c]{cc}%
0 & 1\\
1 & 0
\end{array}
\right)  \in\mathbb{Q}^{2\times2}$ is a root of the polynomial $x^{2}%
-1\in\mathbb{Q}\left[  x\right]  $, since
\begin{align*}
\left(  x^{2}-1\right)  \left[  \left(
\begin{array}
[c]{cc}%
0 & 1\\
1 & 0
\end{array}
\right)  \right]   &  =\left(
\begin{array}
[c]{cc}%
0 & 1\\
1 & 0
\end{array}
\right)  ^{2}-1_{\mathbb{Q}^{2\times2}}=\left(
\begin{array}
[c]{cc}%
1 & 0\\
0 & 1
\end{array}
\right)  -\left(
\begin{array}
[c]{cc}%
1 & 0\\
0 & 1
\end{array}
\right) \\
&  =\left(
\begin{array}
[c]{cc}%
0 & 0\\
0 & 0
\end{array}
\right)  =0_{\mathbb{Q}^{2\times2}}.
\end{align*}

The simplest kind of roots, however, are those that lie in the original ring
$R$. Here are some of their properties:

\begin{proposition}
\label{prop.polring.univar-rootdiv}Let $f$ be a polynomial in $R\left[
x\right]  $. Let $a\in R$. Then, we have the following logical equivalence:%
\[
\left(  a\text{ is a root of }f\right)  \Longleftrightarrow\left(  x-a\mid
f\text{ in }R\left[  x\right]  \right)  .
\]

\end{proposition}

\begin{proof}
[Proof of Proposition \ref{prop.polring.univar-rootdiv}.]The polynomial $x-a$
is monic. Hence, Lemma \ref{lem.polring.univar-quorem.monic} \textbf{(a)}
(applied to $f$ and $x-a$ instead of $a$ and $b$) shows that there is a
\textbf{unique} pair $\left(  q,r\right)  $ of polynomials in $R\left[
x\right]  $ such that
\[
f=q\cdot\left(  x-a\right)  +r\ \ \ \ \ \ \ \ \ \ \text{and}%
\ \ \ \ \ \ \ \ \ \ \deg r<\deg\left(  x-a\right)  .
\]
Consider this pair $\left(  q,r\right)  $. From $\deg r<\deg\left(
x-a\right)  =1$, we see that $\deg r\leq0$, which means that $r$ is a
constant. In other words, $r\in R$.

Now, let us substitute $a$ for $x$ on both sides of the equality
$f=q\cdot\left(  x-a\right)  +r$. Thus we get%
\begin{equation}
f\left[  a\right]  =q\left[  a\right]  \cdot\left(  a-a\right)  +r\left[
a\right]  . \label{pf.prop.polring.univar-rootdiv.3}%
\end{equation}
It is worth going through the proof of this equality in some more detail.
Namely, we have $f=q\cdot\left(  x-a\right)  +r$, so that%
\begin{align*}
f\left[  a\right]   &  =\left(  q\cdot\left(  x-a\right)  +r\right)  \left[
a\right] \\
&  =\left(  q\cdot\left(  x-a\right)  \right)  \left[  a\right]  +r\left[
a\right]  \ \ \ \ \ \ \ \ \ \ \left(
\begin{array}
[c]{c}%
\text{by (\ref{eq.thm.polring.univar-sub-hom.sum}), applied to }R\text{,
}q\cdot\left(  x-a\right)  \text{ and }r\\
\text{instead of }A\text{, }p\text{ and }q
\end{array}
\right) \\
&  =q\left[  a\right]  \cdot\underbrace{\left(  x-a\right)  \left[  a\right]
}_{\substack{=a-a\\\text{(by the definition}\\\text{of an evaluation)}%
}}+\,r\left[  a\right]  \ \ \ \ \ \ \ \ \ \ \left(
\begin{array}
[c]{c}%
\text{by (\ref{eq.thm.polring.univar-sub-hom.prod}), applied to }R\text{,
}q\text{ and }x-a\\
\text{instead of }A\text{, }p\text{ and }q
\end{array}
\right) \\
&  =q\left[  a\right]  \cdot\left(  a-a\right)  +r\left[  a\right]  .
\end{align*}
Thus, (\ref{pf.prop.polring.univar-rootdiv.3}) has been proven in detail.

Now, (\ref{pf.prop.polring.univar-rootdiv.3}) becomes%
\[
f\left[  a\right]  =q\left[  a\right]  \cdot\underbrace{\left(  a-a\right)
}_{=0}+\,r\left[  a\right]  =r\left[  a\right]  =r\ \ \ \ \ \ \ \ \ \ \left(
\text{since }r\text{ is a constant}\right)  .
\]
Now, we have the following chain of equivalences:%
\begin{align*}
\left(  a\text{ is a root of }f\right)  \  &  \Longleftrightarrow\ \left(
f\left[  a\right]  =0\right)  \ \ \ \ \ \ \ \ \ \ \ \left(  \text{by the
definition of a root}\right) \\
&  \Longleftrightarrow\ \left(  r=0\right)  \ \ \ \ \ \ \ \ \ \ \ \left(
\text{since }f\left[  a\right]  =r\right) \\
&  \Longleftrightarrow\ \left(  x-a\mid f\text{ in }R\left[  x\right]
\right)
\end{align*}
(by Proposition \ref{prop.polring.univar-rem0}, applied to $x-a$ and $f$
instead of $b$ and $a$). This proves Proposition
\ref{prop.polring.univar-rootdiv}.
\end{proof}

The following theorem is often known as the \textbf{easy half of the FTA
(Fundamental Theorem of Algebra)}:

\begin{theorem}
\label{thm.polring.univar-easyFTA}Let $R$ be an integral domain. Let
$n\in\mathbb{N}$. Then, any nonzero polynomial $f\in R\left[  x\right]  $ of
degree $\leq n$ has at most $n$ roots in $R$. (We are not counting the roots
with multiplicity here.)
\end{theorem}

\begin{proof}
We induct on $n$. The \textit{base case} ($n=0$) is obvious (indeed, a nonzero
polynomial of degree $\leq0$ must be constant, and thus cannot have any roots
to begin with).

\textit{Induction step:} Let $m$ be a positive integer. Assume (as the
induction hypothesis) that Theorem \ref{thm.polring.univar-easyFTA} holds for
$n=m-1$. We must prove that Theorem \ref{thm.polring.univar-easyFTA} holds for
$n=m$.

So let $f\in R\left[  x\right]  $ be a nonzero polynomial of degree $\leq m$.
We must prove that $f$ has at most $m$ roots in $R$.

Indeed, assume the contrary. Thus, $f$ has $m+1$ distinct roots $a_{1}%
,a_{2},\ldots,a_{m+1}$ in $R$ (and possibly more, but we will only need these
$m+1$).

In particular, $a_{m+1}$ is a root of $f$, so that we have $x-a_{m+1}\mid f$
in $R\left[  x\right]  $ (by Proposition \ref{prop.polring.univar-rootdiv},
applied to $a=a_{m+1}$). That is, there exists a polynomial $q\in R\left[
x\right]  $ such that $f=\left(  x-a_{m+1}\right)  \cdot q$. Consider this
$q$. Now, it is easy to see that $a_{1},a_{2},\ldots,a_{m}$ are roots of $q$
(indeed, this uses the fact that $a_{1},a_{2},\ldots,a_{m+1}$ are
\textbf{distinct} roots of $f$ and that $R$ is an integral
domain\footnote{Here is the proof in detail: Let $i\in\left\{  1,2,\ldots
,m\right\}  $. We must show that $a_{i}$ is a root of $q$. Note that $i\neq
m+1$ (since $i\in\left\{  1,2,\ldots,m\right\}  $) and thus $a_{i}\neq
a_{m+1}$ (since $a_{1},a_{2},\ldots,a_{m+1}$ are distinct). Substituting
$a_{i}$ for $x$ in the equality $f=\left(  x-a_{m+1}\right)  \cdot q$, we find%
\[
f\left[  a_{i}\right]  =\left(  a_{i}-a_{m+1}\right)  \cdot q\left[
a_{i}\right]
\]
(formally speaking, this relies on a similar argument as we used to prove
(\ref{pf.prop.polring.univar-rootdiv.3})). Hence,
\[
\left(  a_{i}-a_{m+1}\right)  \cdot q\left[  a_{i}\right]  =f\left[
a_{i}\right]  =0\ \ \ \ \ \ \ \ \ \ \left(  \text{since }a_{i}\text{ is a root
of }f\right)  .
\]
Since $R$ is an integral domain, this entails that we have $a_{i}-a_{m+1}=0$
or $q\left[  a_{i}\right]  =0$. Since $a_{i}-a_{m+1}=0$ is impossible (because
$a_{i}\neq a_{m+1}$), we thus conclude that $q\left[  a_{i}\right]  =0$. In
other words, $a_{i}$ is a root of $q$. Qed.}). Hence, the polynomial $q$ has
at least $m$ roots in $R$ (since these $m$ roots $a_{1},a_{2},\ldots,a_{m}$
are distinct). Also, the polynomial $q$ is nonzero (since otherwise, we would
have $q=0$ and thus $f=\left(  x-a_{m+1}\right)  \cdot\underbrace{q}_{=0}=0$,
contradicting the fact that $f$ is nonzero).

However, Proposition \ref{prop.polring.univar-degpq} \textbf{(c)} (or
Proposition \ref{prop.polring.univar-degpq} \textbf{(b)}, if you wish) yields
\[
\deg\left(  \left(  x-a_{m+1}\right)  \cdot q\right)  =\underbrace{\deg\left(
x-a_{m+1}\right)  }_{=1}+\deg q=1+\deg q,
\]
so that%
\[
\deg q=\deg\left(  \underbrace{\left(  x-a_{m+1}\right)  \cdot q}_{=f}\right)
-1=\underbrace{\deg f}_{\substack{\leq m\\\text{(since }f\text{ has degree
}\leq m\text{)}}}-\,1\leq m-1.
\]
In other words, the polynomial $q$ has degree $\leq m-1$. Hence, by the
induction hypothesis, we can apply Theorem \ref{thm.polring.univar-easyFTA} to
$q$ and $m-1$ instead of $f$ and $n$. We thus conclude that $q$ has at most
$m-1$ roots in $R$. This contradicts the fact that $q$ has at least $m$ roots
in $R$ (which we have shown above). This contradiction completes the induction
step, and so we are done proving Theorem \ref{thm.polring.univar-easyFTA}.
\end{proof}

\begin{remark}
Theorem \ref{thm.polring.univar-easyFTA} can fail if $R$ is not an integral
domain. For instance, the polynomial $x^{2}-\overline{1}\in\left(
\mathbb{Z}/8\right)  \left[  x\right]  $ has degree $2$ but has $4$ roots in
$\mathbb{Z}/8$ (namely, $\overline{1},\ \overline{3},\ \overline
{5},\ \overline{7}$).
\end{remark}

\begin{remark}
Proposition \ref{prop.polring.univar-rootdiv} and Theorem
\ref{thm.polring.univar-easyFTA} are only concerned with roots in the original
ring $R$. As stated, they don't apply to roots that belong to other
$R$-algebras $A$. And indeed, Theorem \ref{thm.polring.univar-easyFTA} fails
quite dramatically if we try to apply it to other $R$-algebras $A$. For
instance, the polynomial $x^{2}+1\in\mathbb{R}\left[  x\right]  $ has
infinitely many roots in the ring of quaternions $\mathbb{H}$ (see Exercise
\ref{exe.quaternions.roots-of--1} \textbf{(b)}). Proposition
\ref{prop.polring.univar-rootdiv} can be extended to commutative $R$-algebras
$A$ (thus allowing $a\in A$ instead of $a\in R$, and \textquotedblleft
converting\textquotedblright\ the polynomial $f$ into a polynomial in
$A\left[  x\right]  $), although this would not make it significantly more
general, since we can \textbf{already} apply it to $A$ instead of $R$ in such
case (see Exercise \ref{exe.polring.univar.Ralg}).
\end{remark}

\begin{remark}
Let us say a few words about the weird name of Theorem
\ref{thm.polring.univar-easyFTA}. The famous \textbf{fundamental theorem of
algebra} (short: \textbf{FTA}) says that any polynomial of degree $n$ in
$\mathbb{C}\left[  x\right]  $ has exactly $n$ roots in $\mathbb{C}$, if we
count the roots with multiplicity. Despite its name, this theorem is not
actually algebraic in nature, since it relies on the analytic structure of the
complex numbers (and the underlying real numbers), and does not hold (e.g.)
for the Gaussian rationals $\mathbb{Q}\left[  i\right]  $. Accordingly, each
proof of the FTA requires at least a little bit of real analysis (and
sometimes far more than a little bit). Various proofs can be found in
\cite[Chapter 3]{LaNaSc16}, \cite[Theorem 7.1]{Aluffi16}, \cite[Chapter IX,
\S 10]{Knapp1}, \cite[Theorem 44.8]{Warner90}, \cite[Theorem 11.6.7]{Steinb06}
and many other places (some of which prove weaker-sounding but equivalent
versions of the result); more exotic proofs are listed in
\url{https://mathoverflow.net/questions/10535} .

However, one \textquotedblleft half\textquotedblright\ of the FTA -- namely,
the claim that a polynomial of degree $n$ in $\mathbb{C}\left[  x\right]  $
always has \textbf{at most} $n$ roots in $\mathbb{C}$ -- actually can be
proved algebraically, and holds not just for $\mathbb{C}$ but also for any
integral domain $R$ in its stead. If we drop the notion of multiplicities,
then this \textquotedblleft half\textquotedblright\ is precisely Theorem
\ref{thm.polring.univar-easyFTA}. Thus, Theorem
\ref{thm.polring.univar-easyFTA} is called the \textquotedblleft easy half of
the FTA\textquotedblright. Despite being the easy half, it is surprisingly
useful, and we will see some of its applications in the following subsections.
In comparison, the \textquotedblleft hard half of the FTA\textquotedblright%
\ (the part that really requires $\mathbb{C}$) is rarely used in abstract
algebra (since algebraists prefer to work in settings more general than
$\mathbb{C}$), but it is important (e.g.) in complex linear algebra, where it
is responsible (e.g.) for the fact that each $n\times n$-matrix over
$\mathbb{C}$ has $n$ eigenvalues (counted with multiplicities). Thus, the name
\textquotedblleft FTA\textquotedblright\ should be regarded as somewhat of a
historical artefact.
\end{remark}

\begin{exercise}
Let $R$ and $S$ be two commutative rings. Let $f:R\rightarrow S$ be a ring
morphism. Let $f\left[  x\right]  $ denote the map%
\begin{align*}
R\left[  x\right]   &  \rightarrow S\left[  x\right]  ,\\
\sum_{i\in\mathbb{N}}r_{i}x^{i}  &  \mapsto\sum_{i\in\mathbb{N}}f\left(
r_{i}\right)  x^{i}\ \ \ \ \ \ \ \ \ \ \left(  \text{for all }r_{i}\in
R\right)  .
\end{align*}
(This map transforms a polynomial by applying $f$ to each of its coefficients.)

Prove that $f\left[  x\right]  $ is a ring morphism.
\end{exercise}

\begin{exercise}
\label{exe.polring.univar.Ralg}Let $R$ be a commutative ring, and let $A$ be a
commutative $R$-algebra. For each polynomial $f=\sum_{i\in\mathbb{N}}%
f_{i}x^{i}\in R\left[  x\right]  $ (with $f_{i}\in R$), we let $f_{A}$ denote
the polynomial $\sum_{i\in\mathbb{N}}\left(  f_{i}\cdot1_{A}\right)  x^{i}\in
A\left[  x\right]  $ (which is simply the polynomial $f$, with each
coefficient \textquotedblleft converted\textquotedblright\ into an element of
$A$ using the standard map $R\rightarrow A,\ r\mapsto r\cdot1_{A}$).

\begin{enumerate}
\item[\textbf{(a)}] Prove that the map%
\begin{align*}
R\left[  x\right]   &  \rightarrow A\left[  x\right]  ,\\
f  &  \mapsto f_{A}%
\end{align*}
is an $R$-algebra morphism.

\item[\textbf{(b)}] Prove that $f_{A}\left[  a\right]  =f\left[  a\right]  $
for any $f\in R\left[  x\right]  $ and any $a\in A$.
\end{enumerate}
\end{exercise}

\begin{fineprint}
The following exercise provides an analogue of Theorem
\ref{thm.polring.univar-easyFTA} for polynomials in two variables. (Analogues
for $n$ variables can be obtained similarly, but require some more cumbersome notation.)
\end{fineprint}

\begin{exercise}
\label{exe.polring.FTA-easy.2var}Let $R$ be an integral domain.

The $x$\textbf{-degree} of a nonzero polynomial $p\in R\left[  x,y\right]  $
is defined to be the largest $i\in\mathbb{N}$ such that there exists some
$j\in\mathbb{N}$ satisfying $\left[  x^{i}y^{j}\right]  p\neq0$. This
$x$-degree is denoted by $\deg_{x}p$. Similarly, the $y$\textbf{-degree} of a
nonzero polynomial $p\in R\left[  x,y\right]  $ is defined to be the largest
$j\in\mathbb{N}$ such that there exists some $i\in\mathbb{N}$ satisfying
$\left[  x^{i}y^{j}\right]  p\neq0$. This $y$-degree is denoted by $\deg_{y}%
p$. (For example, the polynomial $p=2x^{6}y+3xy^{2}-x^{3}+xy\in\mathbb{Z}%
\left[  x,y\right]  $ has $x$-degree $\deg_{x}p=6$ and $y$-degree $\deg
_{y}p=2$.) If $p$ is the zero polynomial, then we set $\deg_{x}p=-\infty$ and
$\deg_{y}p=-\infty$.

Let $n,m\in\mathbb{N}$. Let $p\in R\left[  x,y\right]  $ be any polynomial
satisfying $\deg_{x}p\leq n$ and $\deg_{y}p\leq m$.

Let $a_{0},a_{1},\ldots,a_{n}$ be $n+1$ distinct elements of $R$. Let
$b_{0},b_{1},\ldots,b_{m}$ be $m+1$ distinct elements of $R$. Prove the
following: If%
\[
p\left[  a_{i},b_{j}\right]  =0\ \ \ \ \ \ \ \ \ \ \text{for all }\left(
i,j\right)  \in\left\{  0,1,\ldots,n\right\}  \times\left\{  0,1,\ldots
,m\right\}  ,
\]
then $p=0$ in $R\left[  x,y\right]  $. \medskip

[\textbf{Hint:} Use Theorem \ref{thm.polring.univar-easyFTA} many times.
Specifically, for each $j\in\left\{  0,1,\ldots,m\right\}  $, argue that the
univariate polynomial $p\left[  x,b_{j}\right]  \in R\left[  x\right]  $ has
too many roots to be nonzero. Then, decompose $p$ into $p=\sum_{k=0}^{n}%
p_{k}\left[  y\right]  x^{k}$, where $p_{0},p_{1},\ldots,p_{n}$ are univariate
polynomials of degree $\leq m$. Apply Theorem \ref{thm.polring.univar-easyFTA}
again to each of these polynomials $p_{0},p_{1},\ldots,p_{n}$.]
\end{exercise}

\subsubsection{Application to $\mathbb{Z}/p$: Wilson revisited}

The easy half of the FTA has a surprising plenitude of applications. Let me
show an application to finite fields.

We fix a prime number $p$ for the rest of this subsection.

First, let us reword Fermat's Little Theorem (Proposition
\ref{prop.ent.flt.Z/p}) in the language of polynomials. First, we consider the
polynomial%
\[
x^{p}-x\in\left(  \mathbb{Z}/p\right)  \left[  x\right]  .
\]
Proposition \ref{prop.ent.flt.Z/p} yields that all evaluations of this
polynomial at elements of $\mathbb{Z}/p$ are $0$ (in fact, for each
$u\in\mathbb{Z}/p$, we have $\left(  x^{p}-x\right)  \left[  u\right]
=u^{p}-u=0$, since Proposition \ref{prop.ent.flt.Z/p} yields $u^{p}=u$). The
polynomial itself is not zero, and this is no surprise: It is a degree-$p$
polynomial, so it can afford to have $p$ roots in $\mathbb{Z}/p$ without being
forced by Theorem \ref{thm.polring.univar-easyFTA} to be the zero polynomial.
However, it is \textquotedblleft dangerously close\textquotedblright; if its
degree was even a little bit smaller than $p$, then we would obtain a
contradiction. We can exploit this to extract a nice corollary.

To this end, we define the more sophisticated polynomial%
\[
f:=\left(  x^{p}-x\right)  -\underbrace{\prod_{u\in\mathbb{Z}/p}\left(
x-u\right)  }_{=\left(  x-\overline{0}\right)  \left(  x-\overline{1}\right)
\cdots\left(  x-\overline{\left(  p-1\right)  }\right)  }\in\left(
\mathbb{Z}/p\right)  \left[  x\right]  .
\]
This polynomial $f$ has degree $\leq p-1$ (check
this!\footnote{\textit{Proof.} Both polynomials $x^{p}-x$ and $\prod
_{u\in\mathbb{Z}/p}\left(  x-u\right)  $ have degree $p$ and leading
coefficient $1$. Thus, when you subtract the polynomial $\prod_{u\in
\mathbb{Z}/p}\left(  x-u\right)  $ from $x^{p}-x$, the $x^{p}$ terms of both
polynomials cancel, and what remains is a linear combination of $x^{0}%
,x^{1},\ldots,x^{p-1}$ -- that is, a polynomial of degree $\leq p-1$.}). But
it still has (at least) $p$ roots in $\mathbb{Z}/p$; indeed, all the $p$
elements of $\mathbb{Z}/p$ are roots of $f$, since each $w\in\mathbb{Z}/p$
satisfies%
\[
f\left[  w\right]  =\underbrace{\left(  w^{p}-w\right)  }%
_{\substack{=0\\\text{(since Proposition \ref{prop.ent.flt.Z/p}}\\\text{yields
}w^{p}=w\text{)}}}-\underbrace{\prod_{u\in\mathbb{Z}/p}\left(  w-u\right)
}_{\substack{=0\\\text{(since one of the factors}\\\text{in this product is
}w-w=0\text{)}}}=0-0=0.
\]
If the polynomial $f$ was nonzero, then this would contradict Theorem
\ref{thm.polring.univar-easyFTA} (since $\mathbb{Z}/p$ is a field and thus an
integral domain). Hence, $f$ must be zero. Since we defined $f$ to be the
difference $\left(  x^{p}-x\right)  -\prod_{u\in\mathbb{Z}/p}\left(
x-u\right)  $, we thus conclude that $x^{p}-x=\prod_{u\in\mathbb{Z}/p}\left(
x-u\right)  $. Let us state this as a proposition:

\begin{proposition}
\label{prop.ent.xp-x-over-Fp}Let $p$ be a prime number. Then,%
\[
x^{p}-x=\prod_{u\in\mathbb{Z}/p}\left(  x-u\right)
\ \ \ \ \ \ \ \ \ \ \text{in the polynomial ring }\left(  \mathbb{Z}/p\right)
\left[  x\right]  .
\]

\end{proposition}

Now, let us milk this for consequences. We have%
\begin{align*}
\prod_{u\in\mathbb{Z}/p}\left(  x-u\right)   &  =\left(  x-\overline
{0}\right)  \left(  x-\overline{1}\right)  \cdots\left(  x-\overline{\left(
p-1\right)  }\right) \\
&  \ \ \ \ \ \ \ \ \ \ \ \ \ \ \ \ \ \ \ \ \left(  \text{since }%
\mathbb{Z}/p=\left\{  \overline{0},\overline{1},\ldots,\overline{p-1}\right\}
\right) \\
&  =x\underbrace{\left(  x-\overline{1}\right)  \left(  x-\overline{2}\right)
\cdots\left(  x-\overline{\left(  p-1\right)  }\right)  }_{\substack{=\left(
-\overline{1}\right)  \left(  -\overline{2}\right)  \cdots\left(
-\overline{\left(  p-1\right)  }\right)  \cdot x^{0}+\left(  \text{higher
powers of }x\right)  \\\text{(here, \textquotedblleft higher powers of
}x\text{\textquotedblright\ means}\\\text{\textquotedblleft any powers of
}x\text{ higher than }x^{0}\text{\textquotedblright)}}}\\
&  =x\left(  \left(  -\overline{1}\right)  \left(  -\overline{2}\right)
\cdots\left(  -\overline{\left(  p-1\right)  }\right)  \cdot x^{0}+\left(
\text{higher powers of }x\right)  \right) \\
&  =\left(  -\overline{1}\right)  \left(  -\overline{2}\right)  \cdots\left(
-\overline{\left(  p-1\right)  }\right)  \cdot x^{1}+\left(  \text{higher
powers of }x\right)  .
\end{align*}
Thus, the coefficient of $x^{1}$ in the polynomial $\prod_{u\in\mathbb{Z}%
/p}\left(  x-u\right)  $ is
\begin{align*}
\left(  -\overline{1}\right)  \left(  -\overline{2}\right)  \cdots\left(
-\overline{\left(  p-1\right)  }\right)   &  =\overline{\left(  -1\right)
^{p-1}\cdot\left(  1\cdot2\cdot\cdots\cdot\left(  p-1\right)  \right)  }\\
&  =\overline{\left(  -1\right)  ^{p-1}\cdot\left(  p-1\right)  !}.
\end{align*}
On the other hand, the coefficient of $x^{1}$ in the polynomial $x^{p}-x$ is
$\overline{-1}$ (since $p>1$). But these two coefficients must be equal (since
Proposition \ref{prop.ent.xp-x-over-Fp} says that the polynomials $\prod
_{u\in\mathbb{Z}/p}\left(  x-u\right)  $ and $x^{p}-x$ are equal). Hence,
$\overline{\left(  -1\right)  ^{p-1}\cdot\left(  p-1\right)  !}=\overline{-1}%
$. In other words,%
\[
\left(  -1\right)  ^{p-1}\cdot\left(  p-1\right)  !\equiv-1\operatorname{mod}%
p.
\]
If we multiply this congruence by $\left(  -1\right)  ^{p-1}$, then the left
hand side becomes $\left(  p-1\right)  !$ (since $\left(  -1\right)
^{p-1}\cdot\left(  -1\right)  ^{p-1}=1$), and thus we get
\[
\left(  p-1\right)  !\equiv\left(  -1\right)  ^{p-1}\cdot\left(  -1\right)
=\left(  -1\right)  ^{p}\equiv-1\operatorname{mod}p
\]
(by Theorem \ref{thm.ent.flt}, applied to $a=-1$). Thus, we have proved
Wilson's theorem (Theorem \ref{thm.wilson}) again!

We note that Proposition \ref{prop.ent.xp-x-over-Fp} can be generalized to
arbitrary finite fields:

\begin{exercise}
Let $F$ be a finite field.

\begin{enumerate}
\item[\textbf{(a)}] Prove that%
\[
x^{\left\vert F\right\vert }-x=\prod_{u\in F}\left(  x-u\right)
\ \ \ \ \ \ \ \ \ \ \text{in the polynomial ring }F\left[  x\right]  .
\]

\item[\textbf{(b)}] Prove that the product of all nonzero elements of $F$
equals $-1_{F}$.
\end{enumerate}
\end{exercise}

The method by which we obtained Proposition \ref{prop.ent.xp-x-over-Fp}, too,
can be generalized:

\begin{exercise}
Let $R$ be an integral domain. Let $n\in\mathbb{N}$. Let $f\in R\left[
x\right]  $ be a nonzero polynomial that has degree $\leq n$. Prove the following:

\begin{enumerate}
\item[\textbf{(a)}] If $a_{1},a_{2},\ldots,a_{k}$ are $k$ distinct roots of
$f$ in $R$ for some $k\in\mathbb{N}$, then%
\[
f=q\cdot\left(  x-a_{1}\right)  \left(  x-a_{2}\right)  \cdots\left(
x-a_{k}\right)
\]
for some polynomial $q\in R\left[  x\right]  $ with $\deg q\leq n-k$.

\item[\textbf{(b)}] If $a_{1},a_{2},\ldots,a_{n}$ are $n$ distinct roots of
$f$ in $R$, then%
\[
f=c\left(  x-a_{1}\right)  \left(  x-a_{2}\right)  \cdots\left(
x-a_{n}\right)  ,
\]
where $c=\left[  x^{n}\right]  f$ is the coefficient of $x^{n}$ in $f$.

\item[\textbf{(c)}] More generally, let us replace the assumption
\textquotedblleft Let $R$ be an integral domain\textquotedblright\ by
\textquotedblleft Let $R$ be a commutative ring\textquotedblright. Instead of
requiring that $a_{1},a_{2},\ldots,a_{k}$ (resp., $a_{1},a_{2},\ldots,a_{n}$)
be distinct, we now require that the pairwise differences $a_{i}-a_{j}$ for
$1\leq i<j\leq k$ (resp., $1\leq i<j\leq n$) are units of $R$. Prove that
parts \textbf{(a)} and \textbf{(b)} of this exercise remain valid.
\end{enumerate}
\end{exercise}

Another application of Theorem \ref{thm.polring.univar-easyFTA} is the
following converse to Proposition \ref{prop.finfield.flt}:

\begin{exercise}
\label{exe.21hw4.1b}Let $F$ be a finite field (i.e., a field with finitely
many elements). Let $i>1$ be an integer such that each $u\in F$ satisfies
$u^{i}=u$. Prove that $i\geq\left\vert F\right\vert $.
\end{exercise}

\subsubsection{Application to $\mathbb{Z}/p$: Sum of $k$-th powers}

Here is another surprisingly simple application of Theorem
\ref{thm.polring.univar-easyFTA}:

\begin{proposition}
\label{prop.ent.sum-k-th-powers-mod-p.1}Let $p$ be a prime. Let $k\in\left\{
0,1,\ldots,p-2\right\}  $. Then, the integer $0^{k}+1^{k}+\cdots+\left(
p-1\right)  ^{k}=\sum_{j=0}^{p-1}j^{k}$ is divisible by $p$.
\end{proposition}

For example, for $k=3$ and $p=5$, this says that $0^{3}+1^{3}+2^{3}%
+3^{3}+4^{3}=100$ is divisible by $5$.

\begin{proof}
[Proof of Proposition \ref{prop.ent.sum-k-th-powers-mod-p.1}.]This is obvious
when $k=0$ (because when $k=0$, we have $\sum_{j=0}^{p-1}\underbrace{j^{k}%
}_{=j^{0}=1}=\sum_{j=0}^{p-1}1=p$, which is clearly divisible by $p$). Thus,
we WLOG assume that $k\neq0$.

In $\mathbb{Z}/p$, we have
\begin{equation}
\overline{\sum_{j=0}^{p-1}j^{k}}=\sum_{j=0}^{p-1}\overline{j}^{k}=\sum
_{u\in\mathbb{Z}/p}u^{k} \label{pf.prop.ent.sum-k-th-powers-mod-p.1.2}%
\end{equation}
(since $\mathbb{Z}/p=\left\{  \overline{0},\overline{1},\ldots,\overline
{p-1}\right\}  $). Let us denote the sum $\sum_{u\in\mathbb{Z}/p}u^{k}$ by
$S$. Thus, (\ref{pf.prop.ent.sum-k-th-powers-mod-p.1.2}) becomes
\begin{equation}
\overline{\sum_{j=0}^{p-1}j^{k}}=S.
\label{pf.prop.ent.sum-k-th-powers-mod-p.1.3}%
\end{equation}
If we can show that $S=0$, then (\ref{pf.prop.ent.sum-k-th-powers-mod-p.1.3})
will simplify to $\overline{\sum_{j=0}^{p-1}j^{k}}=0$, which will mean that
$\sum_{j=0}^{p-1}j^{k}$ is divisible by $p$; thus, Proposition
\ref{prop.ent.sum-k-th-powers-mod-p.1} will be proved. Hence, it remains to
prove that $S=0$.

Define a polynomial $f\in\left(  \mathbb{Z}/p\right)  \left[  x\right]  $ by
$f=x^{k+1}-x$. This polynomial is nonzero (since $k\neq0$) and has degree
$k+1\leq p-1$ (since $k\leq p-2$). Thus, Theorem
\ref{thm.polring.univar-easyFTA} (applied to $R=\mathbb{Z}/p$ and $n=p-1$)
yields that it has at most $p-1$ roots in $\mathbb{Z}/p$. Hence, there exists
at least one $a\in\mathbb{Z}/p$ such that $f\left[  a\right]  \neq0$ (because
otherwise, all the $p$ elements $a\in\mathbb{Z}/p$ would be roots of $f$, but
this would give $f$ more than $p-1$ roots). Consider this $a$.

From $f=x^{k+1}-x$, we obtain $f\left[  a\right]  =a^{k+1}-a=a\left(
a^{k}-1\right)  $, so that $a\left(  a^{k}-1\right)  =f\left[  a\right]
\neq0$. Hence, $a\neq0$ and $a^{k}-1\neq0$. The element $a$ of $\mathbb{Z}/p$
is nonzero (since $a\neq0$) and thus a unit (since $\mathbb{Z}/p$ is a field).
It therefore has an inverse $a^{-1}$. Hence, the map%
\begin{align*}
\mathbb{Z}/p  &  \rightarrow\mathbb{Z}/p,\\
u  &  \mapsto au
\end{align*}
is invertible\footnote{Its inverse is the map
\begin{align*}
\mathbb{Z}/p  &  \rightarrow\mathbb{Z}/p,\\
u  &  \mapsto a^{-1}u.
\end{align*}
}, i.e., a bijection. We can thus substitute $au$ for $u$ in the sum
$\sum_{u\in\mathbb{Z}/p}u^{k}$. We obtain%
\[
\sum_{u\in\mathbb{Z}/p}u^{k}=\sum_{u\in\mathbb{Z}/p}\underbrace{\left(
au\right)  ^{k}}_{=a^{k}u^{k}}=\sum_{u\in\mathbb{Z}/p}a^{k}u^{k}=a^{k}%
\sum_{u\in\mathbb{Z}/p}u^{k}.
\]
In view of $\sum_{u\in\mathbb{Z}/p}u^{k}=S$, we can rewrite this equality as
$S=a^{k}S$. Hence, $a^{k}S-S=0$. In other words, $\left(  a^{k}-1\right)
S=0$. Since $\mathbb{Z}/p$ is an integral domain, and since $a^{k}-1\neq0$,
this entails $S=0$ (because otherwise, from $a^{k}-1\neq0$ and $S\neq0$, we
would obtain $\left(  a^{k}-1\right)  S\neq0$). As explained above, this
completes the proof of Proposition \ref{prop.ent.sum-k-th-powers-mod-p.1}.
\end{proof}

\begin{corollary}
\label{cor.ent.sum-f-over-p.1}Let $p$ be a prime. Let $f\in\mathbb{Z}\left[
x\right]  $ be a polynomial (with integer coefficients!) of degree $\leq p-2$.
Then, $\sum_{j=0}^{p-1}f\left(  j\right)  \equiv0\operatorname{mod}p$.
\end{corollary}

\begin{proof}
Write the polynomial $f$ in the form%
\begin{equation}
f=\sum_{k=0}^{p-2}a_{k}x^{k} \label{pf.cor.ent.sum-f-over-p.1.f=}%
\end{equation}
with $a_{0},a_{1},\ldots,a_{p-2}\in\mathbb{Z}$. (This can be done, since
$f\in\mathbb{Z}\left[  x\right]  $ has degree $\leq p-2$.) Thus,%
\[
\sum_{j=0}^{p-1}\underbrace{f\left(  j\right)  }_{\substack{=\sum_{k=0}%
^{p-2}a_{k}j^{k}\\\text{(by (\ref{pf.cor.ent.sum-f-over-p.1.f=}))}}%
}=\sum_{j=0}^{p-1}\ \ \sum_{k=0}^{p-2}a_{k}j^{k}=\sum_{k=0}^{p-2}%
a_{k}\underbrace{\sum_{j=0}^{p-1}j^{k}}_{\substack{\equiv0\operatorname{mod}%
p\\\text{(since Proposition \ref{prop.ent.sum-k-th-powers-mod-p.1}%
}\\\text{yields that }\sum_{j=0}^{p-1}j^{k}\text{ is}\\\text{divisible by
}p\text{)}}}\equiv\sum_{k=0}^{p-2}a_{k}0=0\operatorname{mod}p.
\]
This proves Corollary \ref{cor.ent.sum-f-over-p.1}.
\end{proof}

We can generalize Proposition \ref{prop.ent.sum-k-th-powers-mod-p.1} a little bit:

\begin{exercise}
\label{exe.ent.sum-k-th-powers-mod-p.2}Let $p$ be a prime. Let $k\in
\mathbb{N}$. Prove that the integer $0^{k}+1^{k}+\cdots+\left(  p-1\right)
^{k}=\sum_{j=0}^{p-1}j^{k}$ is divisible by $p$ if and only if $k$ is not a
positive multiple of $p-1$.
\end{exercise}

\begin{exercise}
Let $p$ be a prime. Let $k\in\mathbb{N}$ be not a multiple of $p-1$. Assume
that the rational number $\dfrac{1}{1^{k}}+\dfrac{1}{2^{k}}+\cdots+\dfrac
{1}{\left(  p-1\right)  ^{k}}=\sum_{j=1}^{p-1}\dfrac{1}{j^{k}}$ has been
written as a ratio $\dfrac{u}{v}$ of two integers $u$ and $v$.

\begin{enumerate}
\item[\textbf{(a)}] Prove that $p\mid u$.

\item[\textbf{(b)}] Assume that $p>3$ and $k=1$. Prove that $p^{2}\mid u$.
\end{enumerate}

[\textbf{Hint:} For part \textbf{(a)}, multiply $\sum_{j=1}^{p-1}\dfrac
{1}{j^{k}}$ by $\left(  p-1\right)  !^{k}$ to obtain an integer; then, work in
$\mathbb{Z}/p$. For part \textbf{(b)}, observe that $\dfrac{1}{j}+\dfrac
{1}{p-j}=\dfrac{p}{j\left(  p-j\right)  }$.]
\end{exercise}

\subsubsection{$F\left[  x\right]  $ is a Euclidean domain}

Let us go back to the case of polynomials over a general field. We next record
an abstract consequence of Corollary \ref{cor.polring.univar-quorem.field}
\textbf{(a)}:

\begin{theorem}
\label{thm.polring.univar-Euc}Let $F$ be a field. Then:

\begin{enumerate}
\item[\textbf{(a)}] The polynomial ring $F\left[  x\right]  $ is a Euclidean
domain. The map%
\begin{align*}
N:F\left[  x\right]   &  \rightarrow\mathbb{N},\\
p  &  \mapsto\max\left\{  \deg p,0\right\}  =%
\begin{cases}
\deg p, & \text{if }p\neq0;\\
0, & \text{if }p=0
\end{cases}
\end{align*}
is a Euclidean norm on $F\left[  x\right]  $.

\item[\textbf{(b)}] Thus, the polynomial ring $F\left[  x\right]  $ is a PID,
hence also a UFD.
\end{enumerate}
\end{theorem}

\begin{proof}
\textbf{(a)} Define a map $N:F\left[  x\right]  \rightarrow\mathbb{N}$ by
\[
N\left(  p\right)  =\max\left\{  \deg p,0\right\}
\ \ \ \ \ \ \ \ \ \ \text{for any }p\in F\left[  x\right]  .
\]
Then, Corollary \ref{cor.polring.univar-quorem.field} \textbf{(a)} shows that
$N$ is a Euclidean norm on the ring $F\left[  x\right]  $. Hence, $F\left[
x\right]  $ is a Euclidean domain (since Corollary
\ref{cor.polring.univar-intdom} shows that $F\left[  x\right]  $ is an
integral domain). This proves Theorem \ref{thm.polring.univar-Euc}
\textbf{(a)}. \medskip

\textbf{(b)} We know from Proposition \ref{prop.eucldom.PID2} that every
Euclidean domain is a PID. Hence, $F\left[  x\right]  $ is a PID (since
$F\left[  x\right]  $ is a Euclidean domain), and therefore a UFD (since we
know from Theorem \ref{thm.UFD.PID-is-UFD} that every PID is a UFD).
\end{proof}

Note that the \textquotedblleft UFD\textquotedblright\ part of Theorem
\ref{thm.polring.univar-Euc} \textbf{(b)} is not a very constructive result;
there is no general algorithm for actually finding a prime factorization of a
polynomial (i.e., for factoring a polynomial into irreducible polynomials)
that works over any field. There are reasonably good algorithms for prime
factorization in $\mathbb{Q}\left[  x\right]  $, however (see Section
\ref{sec.polys2.factor} below for one such algorithm, although a very
inefficient one).

Theorem \ref{thm.polring.univar-Euc} entails, in particular, that univariate
polynomials over a field have gcds and lcms (by Theorem \ref{thm.gcd-lcm.PID}%
). Moreover, the analogue of Bezout's theorem holds:

\begin{theorem}
[Bezout's theorem for polynomials]\label{thm.polyring.univar-gcd}Let $F$ be a
field. Let $a,b\in F\left[  x\right]  $ be two polynomials. Then, for any
choice of $\gcd\left(  a,b\right)  $, there exist two polynomials $u,v\in
F\left[  x\right]  $ such that $ua+vb=\gcd\left(  a,b\right)  $.
\end{theorem}

\begin{proof}
This is a general fact that holds in every PID (but not in every UFD). To wit,
let us set $R=F\left[  x\right]  $, and recall that $R$ is a PID (by Theorem
\ref{thm.polring.univar-Euc} \textbf{(b)}). Recall how we proved the existence
of a gcd (the proof of Theorem \ref{thm.gcd-lcm.PID}): Namely, we argued that
there exists a $c\in R$ satisfying $aR+bR=cR$ (since $R$ is a PID, so that the
ideal $aR+bR$ of $R$ must be principal), and then we proved that this $c$ is a
gcd of $a$ and $b$. Now, assume that we have chosen some gcd of $a$ and $b$,
and denoted it by $\gcd\left(  a,b\right)  $. This $\gcd\left(  a,b\right)  $
is not necessarily identical to $c$, but it is clearly associate to $c$, since
Proposition \ref{prop.gcd-lcm.uni} \textbf{(a)} says that any two gcds of $a$
and $b$ are associate. Thus, $\gcd\left(  a,b\right)  =cw$ for some unit $w$
of $R$. Consider this $w$. Now,
\[
\gcd\left(  a,b\right)  =c\underbrace{w}_{\in R}\in cR=aR+bR.
\]
In other words, there exist some $u,v\in R$ such that $\gcd\left(  a,b\right)
=au+bv$. In other words, there exist some $u,v\in R$ such that $\gcd\left(
a,b\right)  =ua+vb$. This proves Theorem \ref{thm.polyring.univar-gcd}.
\end{proof}

A few words about computability are in order. If $F$ is a field, and if
$a,b\in F\left[  x\right]  $ are two polynomials, then we can compute a
$\gcd\left(  a,b\right)  $ by the extended Euclidean algorithm (more
precisely, the algorithm explained in the proof of Theorem
\ref{thm.rings.euclid.ea.main} \textbf{(b)} computes a Bezout $5$-tuple for
$\left(  a,b\right)  $, and then Corollary \ref{cor.gcd-lcm.eucl-bez} reveals
that the third entry of this $5$-tuple is a gcd of $a$ and $b$). This is one
of the most useful features of univariate polynomials. The following exercises
should give you some practice with this algorithm:

\begin{exercise}
Work in the polynomial ring $\mathbb{Q}\left[  x\right]  $.

\begin{enumerate}
\item[\textbf{(a)}] Compute a $\gcd\left(  x^{3}-x,\ x^{5}-3x+2\right)  $.
(The indefinite article \textquotedblleft a\textquotedblright\ refers to the
fact that a gcd is unique only up to multiplying by a unit.)

\item[\textbf{(b)}] Compute a $\gcd\left(  x^{4}+x^{2}+1,\ x^{4}+x^{3}%
+x^{2}+2\right)  $.
\end{enumerate}
\end{exercise}

\begin{exercise}
Let $F$ be a field. We will work in the polynomial ring $F\left[  x\right]  $.

\begin{enumerate}
\item[\textbf{(a)}] Compute a $\gcd\left(  x^{2}-1,\ x^{3}-1\right)  $.

\item[\textbf{(b)}] Compute a $\gcd\left(  x^{2}-1,\ x^{5}-1\right)  $.

\item[\textbf{(c)}] Compute a $\gcd\left(  x^{4}-1,\ x^{6}-1\right)  $.

\item[\textbf{(d)}] Let $m,n\in\mathbb{N}$. Prove that $\gcd\left(
x^{m}-1,\ x^{n}-1\right)  =x^{\gcd\left(  m,n\right)  }-1$. (To be more
precise, show that $x^{\gcd\left(  m,n\right)  }-1$ is a gcd of $x^{m}-1$ and
$x^{n}-1$. Of course, multiplying this gcd by a nonzero scalar in $F$ will
yield another gcd.)
\end{enumerate}
\end{exercise}

\begin{exercise}
Let $F$ be a field. We will work in the polynomial ring $F\left[  x\right]  $.

\begin{enumerate}
\item[\textbf{(a)}] Compute a $\gcd\left(  x^{2}+1,\ x^{3}+1\right)  $ under
the assumption that $2\cdot1_{F}\neq0_{F}$.

\item[\textbf{(b)}] Compute a $\gcd\left(  x^{2}+1,\ x^{3}+1\right)  $ under
the assumption that $2\cdot1_{F}=0_{F}$.

\item[\textbf{(c)}] Compute a $\gcd\left(  x^{3}+1,\ x^{5}+1\right)  $.

\item[\textbf{(d)}] Let $m,n\in\mathbb{N}$ be odd. Prove that $\gcd\left(
x^{m}+1,\ x^{n}+1\right)  =x^{\gcd\left(  m,n\right)  }+1$.

\item[\textbf{(e)}] Let $m$ be an even positive integer, and $n$ an odd
positive integer. Prove that $\gcd\left(  x^{m}+1,\ x^{n}+1\right)  =1$ if
$2\cdot1_{F}\neq0_{F}$.

\item[\textbf{(f)}] More generally, prove the following: If $m$ and $n$ are
positive integers, and if $a,b\in F$ are two elements satisfying $a^{n}\neq
b^{m}$, then $\gcd\left(  x^{m}-a,\ x^{n}-b\right)  =1$.
\end{enumerate}

[\textbf{Hint:} For part \textbf{(d)}, a substitution can be helpful. Part
\textbf{(e)} is easiest to derive from part \textbf{(f)}.]
\end{exercise}

\begin{warning}
Multivariate polynomial rings (like $\mathbb{Q}\left[  x,y\right]  $) are not
PIDs (and thus not Euclidean domains either). For example, if $R=\mathbb{Q}%
\left[  x,y\right]  $, then the ideal $xR+yR$ is not principal. (Check this!
This ideal is easily seen to consist of all polynomials whose constant term (=
coefficient of $x^{0}y^{0}$) is $0$, but these polynomials are not the
multiples of a single polynomial.) However, multivariate polynomial rings over
fields (and, more generally, over UFDs) are still UFDs. This is a deeper
result than the ones we have proved above (see, e.g., \cite[\S 9.3, Corollary
8]{DumFoo04} or \cite[Theorem 0.61]{ClaRos25} or \cite[Theorem 3.7.4]{Ford22}
or \cite[Corollary (2.6.7)]{ChaLoi21} or \cite[Corollary 8.21 and Remark after
it]{Knapp1} or \cite[Theorem 36.11]{Swanso17} for proofs). As a consequence,
polynomials over a field (or a UFD) have gcds; however, they don't generally
satisfy Bezout's theorem unless the polynomials are univariate polynomials
over a field.

Univariate polynomial rings over non-fields (like $\mathbb{Z}\left[  x\right]
$) behave similarly: They are not PIDs, but they are UFDs when the base ring
is a UFD. (That is, if $R$ is a UFD, then so is $R\left[  x\right]  $.)
\end{warning}

\subsubsection{Lagrange interpolation}

Theorem \ref{thm.polring.univar-easyFTA} has the following simple corollary:

\begin{corollary}
[uniqueness of interpolating polynomial]\label{cor.polint.equal}Let $R$ be an
integral domain. Let $n\in\mathbb{N}$. Let $a_{0},a_{1},\ldots,a_{n}$ be $n+1$
distinct elements of $R$. Let $f,g\in R\left[  x\right]  $ be two polynomials
of degree $\leq n$. Assume that
\begin{equation}
f\left[  a_{i}\right]  =g\left[  a_{i}\right]  \ \ \ \ \ \ \ \ \ \ \text{for
all }i\in\left\{  0,1,\ldots,n\right\}  . \label{eq.cor.polint.equal.ass}%
\end{equation}
Then, $f=g$.
\end{corollary}

Corollary \ref{cor.polint.equal} is saying that if two univariate polynomials
of degree $\leq n$ over an integral domain $R$ agree in at least $n+1$
distinct positions $a_{0},a_{1},\ldots,a_{n}$, then they must be equal. In
particular, if two univariate polynomials (of any degree) over an integral
domain $R$ agree at infinitely many distinct positions, then they must be
equal. This is an extremely useful result with applications all over
mathematics.\footnote{For example, its particular case for $R=\mathbb{C}$ is
\cite[Corollary 7.5.7]{mps}, and is subsequently used in \cite[\S 7.5.3]{mps}
to prove identities for binomial coefficients. Many more applications exist in
a similar vein (see, e.g., \cite[Combinatorial proof of Proposition
6.1.1]{21s}); this technique is known as the \textquotedblleft polynomial
identity trick\textquotedblright.}

\begin{proof}
[Proof of Corollary \ref{cor.polint.equal}.]Assume the contrary. Thus, $f\neq
g$, so that $f-g\neq0$. The nonzero polynomial $f-g$ has degree $\leq n$
(since each of $f$ and $g$ has degree $\leq n$). Thus, Theorem
\ref{thm.polring.univar-easyFTA} (applied to $f-g$ instead of $f$) shows that
$f-g$ has at most $n$ roots in $R$.

However, for each $i\in\left\{  0,1,\ldots,n\right\}  $, the element $a_{i}$
is a root of $f-g$, since%
\[
\left(  f-g\right)  \left[  a_{i}\right]  =f\left[  a_{i}\right]  -g\left[
a_{i}\right]  =0\ \ \ \ \ \ \ \ \ \ \left(  \text{by
(\ref{eq.cor.polint.equal.ass})}\right)  .
\]
In other words, the $n+1$ elements $a_{0},a_{1},\ldots,a_{n}$ of $R$ are roots
of $f-g$. Since these elements $a_{0},a_{1},\ldots,a_{n}$ are distinct, this
shows that the polynomial $f-g$ has at least $n+1$ roots in $R$. But this
contradicts the fact that $f-g$ has at most $n$ roots in $R$. This
contradiction shows that our assumption was false. Thus, Corollary
\ref{cor.polint.equal} is proven.
\end{proof}

Corollary \ref{cor.polint.equal} can be restated as follows: A univariate
polynomial $f\in R\left[  x\right]  $ of degree $\leq n$ over an integral
domain $R$ is uniquely determined by any $n+1$ values $f\left[  a_{0}\right]
,\ f\left[  a_{1}\right]  ,\ \ldots,\ f\left[  a_{n}\right]  $ (provided, of
course, that the inputs $a_{0},a_{1},\ldots,a_{n}\in R$ are distinct and
known). This is a nice uniqueness statement. Can we find a matching existence
statement for it? In other words, if we are given $n+1$ distinct elements
$a_{0},a_{1},\ldots,a_{n}$ of an integral domain $R$ and $n+1$ arbitrary
elements $b_{0},b_{1},\ldots,b_{n}$ of $R$, then is there necessarily some
polynomial $f\in R\left[  x\right]  $ of degree $\leq n$ that satisfies%
\[
f\left[  a_{i}\right]  =b_{i}\ \ \ \ \ \ \ \ \ \ \text{for all }i\in\left\{
0,1,\ldots,n\right\}  \text{ ?}%
\]

A bit of thought reveals that the answer is negative. Indeed, there is no
polynomial $f\in\mathbb{Z}\left[  x\right]  $ that satisfies $f\left[
0\right]  =0$ and $f\left[  2\right]  =1$. The reason is a lack of
divisibility: Any polynomial $f\in\mathbb{Z}\left[  x\right]  $ satisfies
$a-b\mid f\left[  a\right]  -f\left[  b\right]  $ for all $a,b\in\mathbb{Z}$
(check this!), but our conditions $f\left[  0\right]  =0$ and $f\left[
2\right]  =1$ would make this divisibility false for $a=2$ and $b=0$. Thus,
our alleged existence statement cannot hold for arbitrary integral domains $R$.

However, it turns out to be true when $R$ is a field. Moreover, the required
polynomial $f$ can be expressed directly:

\begin{theorem}
[Lagrange interpolation]\label{thm.polint.lagrange}Let $F$ be a field.
Consider the univariate polynomial ring $F\left[  x\right]  $. Let
$n\in\mathbb{N}$.

Let $a_{0},a_{1},\ldots,a_{n}$ be $n+1$ distinct elements of $F$. Let
$b_{0},b_{1},\ldots,b_{n}$ be $n+1$ arbitrary elements of $F$. Then:

\begin{enumerate}
\item[\textbf{(a)}] There is a \textbf{unique} polynomial $p\in F\left[
x\right]  $ satisfying $\deg p\leq n$ and
\[
p\left[  a_{i}\right]  =b_{i}\ \ \ \ \ \ \ \ \ \ \text{for all $i\in\left\{
0,1,\ldots,n\right\}  $}.
\]

\item[\textbf{(b)}] This polynomial $p$ is given by
\[
p=\sum_{j=0}^{n}b_{j}\dfrac{\prod_{k\neq j}\left(  x-a_{k}\right)  }%
{\prod_{k\neq j}\left(  a_{j}-a_{k}\right)  }%
\]
(where the \textquotedblleft$\prod_{k\neq j}$\textquotedblright\ sign means a
product over all $k\in\left\{  0,1,\ldots,n\right\}  $ satisfying $k\neq j$).
\end{enumerate}
\end{theorem}

For example, if $n=2$, then

\begin{itemize}
\item Theorem \ref{thm.polint.lagrange} \textbf{(a)} is saying that there is a
unique polynomial $p\in F\left[  x\right]  $ satisfying $\deg p\leq2$ and%
\[
p\left[  a_{0}\right]  =b_{0}\ \ \ \ \ \ \ \ \ \ \text{and}%
\ \ \ \ \ \ \ \ \ \ p\left[  a_{1}\right]  =b_{1}%
\ \ \ \ \ \ \ \ \ \ \text{and}\ \ \ \ \ \ \ \ \ \ p\left[  a_{2}\right]
=b_{2},
\]

\item and Theorem \ref{thm.polint.lagrange} \textbf{(b)} is saying that this
polynomial $p$ is given by%
\[
p=b_{0}\dfrac{\left(  x-a_{1}\right)  \left(  x-a_{2}\right)  }{\left(
a_{0}-a_{1}\right)  \left(  a_{0}-a_{2}\right)  }+b_{1}\dfrac{\left(
x-a_{0}\right)  \left(  x-a_{2}\right)  }{\left(  a_{1}-a_{0}\right)  \left(
a_{1}-a_{2}\right)  }+b_{2}\dfrac{\left(  x-a_{0}\right)  \left(
x-a_{1}\right)  }{\left(  a_{2}-a_{0}\right)  \left(  a_{2}-a_{1}\right)  }.
\]

\end{itemize}

\begin{proof}
[Proof of Theorem \ref{thm.polint.lagrange}.]Define a polynomial $g\in
F\left[  x\right]  $ by%
\begin{equation}
g=\sum_{j=0}^{n}b_{j}\dfrac{\prod_{k\neq j}\left(  x-a_{k}\right)  }%
{\prod_{k\neq j}\left(  a_{j}-a_{k}\right)  }
\label{sol.polring.lagr-interpol.1.g=}%
\end{equation}
(where the \textquotedblleft$\prod_{k\neq j}$\textquotedblright\ signs mean a
product over all $k\in\left\{  0,1,\ldots,n\right\}  $ satisfying $k\neq j$).
Note that $g$ is well-defined; indeed, all the differences $a_{j}-a_{k}$
appearing in the denominators are nonzero (because $a_{0},a_{1},\ldots,a_{n}$
are distinct) and thus are units (since $F$ is a field).

Each of the $n+1$ addends $b_{j}\dfrac{\prod_{k\neq j}\left(  x-a_{k}\right)
}{\prod_{k\neq j}\left(  a_{j}-a_{k}\right)  }$ on the right hand side of
(\ref{sol.polring.lagr-interpol.1.g=}) is a polynomial of degree $\leq n$
(since the numerator $\prod_{k\neq j}\left(  x-a_{k}\right)  $ is a product of
$n$ degree-$1$ polynomials $x-a_{k}$, whereas the remaining pieces $b_{j}$ and
$\prod_{k\neq j}\left(  a_{j}-a_{k}\right)  $ of the expression are elements
of $F$). Hence, their sum must be a polynomial of degree $\leq n$ as well
(since any sum of polynomials of degree $\leq n$ is again a polynomial of
degree $\leq n$). In other words, $g$ is a polynomial of degree $\leq n$
(since (\ref{sol.polring.lagr-interpol.1.g=}) shows that $g$ is their sum).
That is, we have $\deg g\leq n$.

If $i\in\left\{  0,1,\ldots,n\right\}  $ and $j\in\left\{  0,1,\ldots
,n\right\}  $ satisfy $j\neq i$, then we have%
\begin{equation}
\prod_{k\neq j}\left(  a_{i}-a_{k}\right)  =0
\label{sol.polring.lagr-interpol.1.prod=0}%
\end{equation}
(because in this case, the product $\prod_{k\neq j}\left(  a_{i}-a_{k}\right)
$ contains the factor $a_{i}-a_{i}$ (since $i\neq j$), but this factor is $0$,
and therefore the whole product is $0$).

For each $i\in\left\{  0,1,\ldots,n\right\}  $, we have%
\begin{align*}
g\left[  a_{i}\right]   &  =\left(  \sum_{j=0}^{n}b_{j}\dfrac{\prod_{k\neq
j}\left(  x-a_{k}\right)  }{\prod_{k\neq j}\left(  a_{j}-a_{k}\right)
}\right)  \left[  a_{i}\right]  \ \ \ \ \ \ \ \ \ \ \left(  \text{by the
definition of }g\right) \\
&  =\sum_{j=0}^{n}b_{j}\dfrac{\prod_{k\neq j}\left(  a_{i}-a_{k}\right)
}{\prod_{k\neq j}\left(  a_{j}-a_{k}\right)  }\\
&  =b_{i}\underbrace{\dfrac{\prod_{k\neq i}\left(  a_{i}-a_{k}\right)  }%
{\prod_{k\neq i}\left(  a_{i}-a_{k}\right)  }}_{=1}+\sum_{\substack{j\in
\left\{  0,1,\ldots,n\right\}  ;\\j\neq i}}b_{j}\underbrace{\dfrac
{\prod_{k\neq j}\left(  a_{i}-a_{k}\right)  }{\prod_{k\neq j}\left(
a_{j}-a_{k}\right)  }}_{\substack{=0\\\text{(by
(\ref{sol.polring.lagr-interpol.1.prod=0}))}}}\\
&  \ \ \ \ \ \ \ \ \ \ \ \ \ \ \ \ \ \ \ \ \left(  \text{here, we have split
off the addend for }j=i\text{ from the sum}\right) \\
&  =b_{i}+\underbrace{\sum_{\substack{j\in\left\{  0,1,\ldots,n\right\}
;\\j\neq i}}b_{j}\cdot0}_{=0}=b_{i}.
\end{align*}
Hence, $g$ is a polynomial $p\in F\left[  x\right]  $ satisfying $\deg p\leq
n$ and%
\begin{equation}
p\left[  a_{i}\right]  =b_{i}\ \ \ \ \ \ \ \ \ \ \text{for all $i\in\left\{
0,1,\ldots,n\right\}  $} \label{sol.polring.lagr-interpol.1.pcond}%
\end{equation}
(since we already know that $\deg g\leq n$).

\bigskip

\textbf{(a)} We need to prove that there is a unique polynomial $p\in F\left[
x\right]  $ satisfying $\deg p\leq n$ and
(\ref{sol.polring.lagr-interpol.1.pcond}). We already know that such a $p$
exists (because we have just shown that $g$ is such a $p$); thus, it remains
to prove its uniqueness. In other words, we need to prove the following claim:

\begin{statement}
\textit{Claim 1:} Let $p_{1}$ and $p_{2}$ be two polynomials $p\in F\left[
x\right]  $ satisfying $\deg p\leq n$ and
(\ref{sol.polring.lagr-interpol.1.pcond}). Then, $p_{1}=p_{2}$.
\end{statement}

\begin{proof}
[Proof of Claim 1.]We have assumed that $p_{1}$ is a polynomial $p\in F\left[
x\right]  $ satisfying $\deg p\leq n$ and
(\ref{sol.polring.lagr-interpol.1.pcond}). In other words, $p_{1}\in F\left[
x\right]  $ is a polynomial and satisfies $\deg p_{1}\leq n$ and%
\begin{equation}
p_{1}\left[  a_{i}\right]  =b_{i}\ \ \ \ \ \ \ \ \ \ \text{for all }%
i\in\left\{  0,1,\ldots,n\right\}  .
\label{sol.polring.lagr-interpol.1.c1.pf.1}%
\end{equation}
Similarly, $p_{2}\in F\left[  x\right]  $ is a polynomial and satisfies $\deg
p_{2}\leq n$ and
\begin{equation}
p_{2}\left[  a_{i}\right]  =b_{i}\ \ \ \ \ \ \ \ \ \ \text{for all }%
i\in\left\{  0,1,\ldots,n\right\}  .
\label{sol.polring.lagr-interpol.1.c1.pf.2}%
\end{equation}
For each $i\in\left\{  0,1,\ldots,n\right\}  $, we have%
\begin{align*}
p_{1}\left[  a_{i}\right]   &  =b_{i}\ \ \ \ \ \ \ \ \ \ \left(  \text{by
(\ref{sol.polring.lagr-interpol.1.c1.pf.1})}\right) \\
&  =p_{2}\left[  a_{i}\right]  \ \ \ \ \ \ \ \ \ \ \left(  \text{by
(\ref{sol.polring.lagr-interpol.1.c1.pf.2})}\right)  .
\end{align*}
Hence, Corollary \ref{cor.polint.equal} (applied to $R=F$, $f=p_{1}$ and
$g=p_{2}$) yields that $p_{1}=p_{2}$ (since $p_{1}$ and $p_{2}$ both have
degree $\leq n$). This proves Claim 1.
\end{proof}

Now, our proof of Theorem \ref{thm.polint.lagrange} \textbf{(a)} is complete.
\medskip

\textbf{(b)} In our above proof of Theorem \ref{thm.polint.lagrange}
\textbf{(a)}, we have shown not just that there is a unique polynomial $p\in
F\left[  x\right]  $ satisfying $\deg p\leq n$ and
(\ref{sol.polring.lagr-interpol.1.pcond}); we have also shown that $g$ is such
a polynomial. But since this $p$ is unique, this means that $g$ is the
\textbf{only} such polynomial. Thus, the only such polynomial is $g=\sum
_{j=0}^{n}b_{j}\dfrac{\prod_{k\neq j}\left(  x-a_{k}\right)  }{\prod_{k\neq
j}\left(  a_{j}-a_{k}\right)  }$. This proves Theorem
\ref{thm.polint.lagrange} \textbf{(b)}.
\end{proof}

One of many applications of Theorem \ref{thm.polint.lagrange} \textbf{(b)} is
an explicit formula for recovering a univariate polynomial $f$ of degree $\leq
n$ (over a field) from any $n+1$ values $f\left[  a_{0}\right]  ,\ f\left[
a_{1}\right]  ,\ \ldots,\ f\left[  a_{n}\right]  $ of $f$ (provided that the
inputs $a_{0},a_{1},\ldots,a_{n}$ are distinct and known). Let us make this explicit:

\begin{corollary}
\label{cor.polint.lagrange-as-eq}Let $F$ be a field. Consider the univariate
polynomial ring $F\left[  x\right]  $. Let $n\in\mathbb{N}$.

Let $f\in F\left[  x\right]  $ be a polynomial of degree $\leq n$.

Let $a_{0},a_{1},\ldots,a_{n}$ be $n+1$ distinct elements of $F$. Then,%
\[
f=\sum_{j=0}^{n}f\left[  a_{j}\right]  \cdot\dfrac{\prod_{k\neq j}\left(
x-a_{k}\right)  }{\prod_{k\neq j}\left(  a_{j}-a_{k}\right)  }%
\]
(where the \textquotedblleft$\prod_{k\neq j}$\textquotedblright\ sign means a
product over all $k\in\left\{  0,1,\ldots,n\right\}  $ satisfying $k\neq j$).
\end{corollary}

\begin{proof}
Theorem \ref{thm.polint.lagrange} \textbf{(a)} (applied to $b_{i}=f\left[
a_{i}\right]  $) yields that there is a \textbf{unique} polynomial $p\in
F\left[  x\right]  $ satisfying $\deg p\leq n$ and
\[
p\left[  a_{i}\right]  =f\left[  a_{i}\right]  \ \ \ \ \ \ \ \ \ \ \text{for
all $i\in\left\{  0,1,\ldots,n\right\}  $}.
\]
Furthermore, Theorem \ref{thm.polint.lagrange} \textbf{(b)} (applied to
$b_{i}=f\left[  a_{i}\right]  $) yields that this polynomial $p$ is given by
\[
p=\sum_{j=0}^{n}f\left[  a_{j}\right]  \cdot\dfrac{\prod_{k\neq j}\left(
x-a_{k}\right)  }{\prod_{k\neq j}\left(  a_{j}-a_{k}\right)  }.
\]
Combining this, we conclude that if $p\in F\left[  x\right]  $ is a polynomial
satisfying $\deg p\leq n$ and
\[
p\left[  a_{i}\right]  =f\left[  a_{i}\right]  \ \ \ \ \ \ \ \ \ \ \text{for
all $i\in\left\{  0,1,\ldots,n\right\}  ,$}%
\]
then%
\[
p=\sum_{j=0}^{n}f\left[  a_{j}\right]  \cdot\dfrac{\prod_{k\neq j}\left(
x-a_{k}\right)  }{\prod_{k\neq j}\left(  a_{j}-a_{k}\right)  }.
\]
Applying this to $p=f$, we obtain%
\[
f=\sum_{j=0}^{n}f\left[  a_{j}\right]  \cdot\dfrac{\prod_{k\neq j}\left(
x-a_{k}\right)  }{\prod_{k\neq j}\left(  a_{j}-a_{k}\right)  }%
\]
(since $f\in F\left[  x\right]  $ is a polynomial satisfying $\deg f\leq n$
and $f\left[  a_{i}\right]  =f\left[  a_{i}\right]  $ for all $i\in\left\{
0,1,\ldots,n\right\}  $). This proves Corollary
\ref{cor.polint.lagrange-as-eq}.
\end{proof}

\begin{exercise}
\label{exe.polint.sumajk1}Let $F$ be a field. Let $n\in\mathbb{N}$. Let
$a_{0},a_{1},\ldots,a_{n}$ be $n+1$ distinct elements of $F$. Prove that for
each $\ell\in\left\{  0,1,\ldots,n\right\}  $, we have%
\[
\sum_{j=0}^{n}\dfrac{a_{j}^{\ell}}{\prod_{k\neq j}\left(  a_{j}-a_{k}\right)
}=%
\begin{cases}
1, & \text{if }\ell=n;\\
0, & \text{if }\ell<n
\end{cases}
\]
(where the \textquotedblleft$\prod_{k\neq j}$\textquotedblright\ sign means a
product over all $k\in\left\{  0,1,\ldots,n\right\}  $ satisfying $k\neq j$).
\medskip

[\textbf{Hint:} Apply Corollary \ref{cor.polint.lagrange-as-eq} to $f=x^{\ell
}$ and compare the $x^{n}$-coefficients on both sides of the equality.]
\end{exercise}

\begin{exercise}
Let $n\in\mathbb{N}$ and $\ell\in\mathbb{N}$. Prove that
\[
\sum_{j=0}^{n}\left(  -1\right)  ^{n-j}\dbinom{n}{j}j^{\ell}=%
\begin{cases}
n!, & \text{if }\ell=n;\\
0, & \text{if }\ell<n.
\end{cases}
\]

[\textbf{Hint:} Apply Exercise \ref{exe.polint.sumajk1} to $F=\mathbb{Q}$ and
$a_{i}=i$.]
\end{exercise}

\begin{exercise}
Let $a,b\in\mathbb{R}$ and $n\in\mathbb{N}$ and $\ell\in\left\{
0,1,\ldots,n\right\}  $. Prove that%
\[
\sum_{j=0}^{n}\left(  -1\right)  ^{n-j}\dbinom{n}{j}\dbinom{aj+b}{\ell}=%
\begin{cases}
a^{n}, & \text{if }\ell=n;\\
0, & \text{if }\ell<n.
\end{cases}
\]

(Recall that the binomial coefficient $\dbinom{u}{\ell}$ is defined to be
$\dfrac{u\left(  u-1\right)  \left(  u-2\right)  \cdots\left(  u-\ell
+1\right)  }{\ell!}$ for each $u\in\mathbb{R}$.)
\end{exercise}

\begin{exercise}
\label{exe.polint.sumajkf}Let $F$ be a field. Let $n\in\mathbb{N}$. Let
$a_{0},a_{1},\ldots,a_{n}$ be $n+1$ distinct elements of $F$. Let $f\in
F\left[  x\right]  $ be a polynomial of degree $\leq n$. Prove that%
\[
\sum_{j=0}^{n}\dfrac{f\left[  a_{j}\right]  }{\prod_{k\neq j}\left(
a_{j}-a_{k}\right)  }=\left[  x^{n}\right]  f
\]
(where the \textquotedblleft$\prod_{k\neq j}$\textquotedblright\ sign means a
product over all $k\in\left\{  0,1,\ldots,n\right\}  $ satisfying $k\neq j$).
\medskip

[\textbf{Hint:} Apply Corollary \ref{cor.polint.lagrange-as-eq}. Then, take
the $x^{n}$-coefficient.]
\end{exercise}

\begin{exercise}
\label{exe.polint.sumajk2}Let $F$ be a field. Let $n\in\mathbb{N}$. Let
$a_{0},a_{1},\ldots,a_{n}$ be $n+1$ distinct elements of $F$. Prove that%
\[
\sum_{j=0}^{n}\dfrac{a_{j}^{n+1}}{\prod_{k\neq j}\left(  a_{j}-a_{k}\right)
}=a_{0}+a_{1}+\cdots+a_{n}%
\]
(where the \textquotedblleft$\prod_{k\neq j}$\textquotedblright\ sign means a
product over all $k\in\left\{  0,1,\ldots,n\right\}  $ satisfying $k\neq j$).
\medskip

[\textbf{Hint:} Apply Exercise \ref{exe.polint.sumajkf} to $f=x^{n+1}%
-\prod_{k=0}^{n}\left(  x-a_{k}\right)  $.]
\end{exercise}

\begin{exercise}
\label{exe.polint.sumajk4}Let $F$ be a field. Let $n\in\mathbb{N}$. Let
$a_{0},a_{1},\ldots,a_{n}$ be $n+1$ distinct elements of $F$. Prove that%
\[
\sum_{j=0}^{n}a_{j}\dfrac{\prod_{k\neq j}\left(  a_{j}+a_{k}\right)  }%
{\prod_{k\neq j}\left(  a_{j}-a_{k}\right)  }=a_{0}+a_{1}+\cdots+a_{n}%
\]
(where the \textquotedblleft$\prod_{k\neq j}$\textquotedblright\ sign means a
product over all $k\in\left\{  0,1,\ldots,n\right\}  $ satisfying $k\neq j$).
\medskip

[\textbf{Hint:} Let $2_{F}:=2\cdot1_{F}\in F$. Distinguish between the cases
when $2_{F}=0_{F}$ and when $2_{F}\neq0_{F}$. In the former case, $a+b=a-b$
for all $a,b\in F$ (why?), and the claim becomes very easy. Now, consider the
latter case. In this case, $2_{F}$ is a unit of $F$, thus can be cancelled.
Now, apply Exercise \ref{exe.polint.sumajkf} to $f=\prod_{k=0}^{n}\left(
x+a_{k}\right)  -\prod_{k=0}^{n}\left(  x-a_{k}\right)  $. Then, simplify and
cancel $2_{F}$.]

(This exercise is \href{https://cms.math.ca/publications/crux/}{Crux
Mathematicorum} problem \#4762, proposed by Didier Pinchon and George Stoica
in issue 49/2.)
\end{exercise}

\begin{exercise}
\label{exe.polint.sumajk3}Let $F$ be a field. Let $n\in\mathbb{N}$. Let
$a_{0},a_{1},\ldots,a_{n}$ be $n+1$ distinct elements of $F$. Prove that%
\[
\sum_{j=0}^{n}\dfrac{\prod_{k\neq j}\left(  a_{j}+a_{k}\right)  }{\prod_{k\neq
j}\left(  a_{j}-a_{k}\right)  }=%
\begin{cases}
1, & \text{if }n\text{ is even};\\
0, & \text{if }n\text{ is odd}%
\end{cases}
\]
(where the \textquotedblleft$\prod_{k\neq j}$\textquotedblright\ sign means a
product over all $k\in\left\{  0,1,\ldots,n\right\}  $ satisfying $k\neq j$).
\medskip

[\textbf{Hint:} Let $2_{F}:=2\cdot1_{F}\in F$. Deal with the case $2_{F}%
=0_{F}$ as in Exercise \ref{exe.polint.sumajk4}. In the remaining case, show
that the polynomial $\prod_{k=0}^{n}\left(  a_{k}+x\right)  -\prod_{k=0}%
^{n}\left(  a_{k}-x\right)  \in F\left[  x\right]  $ is divisible by $2x$, and
let $f$ be the quotient of this division. Now, apply Exercise
\ref{exe.polint.sumajkf} to this $f$.]
\end{exercise}

\begin{fineprint}
For more identities in the vein of Exercises \ref{exe.polint.sumajk1},
\ref{exe.polint.sumajk2}, \ref{exe.polint.sumajk3} and
\ref{exe.polint.sumajk4}, see \cite{Nica22}.
\end{fineprint}

\begin{exercise}
\label{exe.21hw3.6c}Let $n\in\mathbb{N}$. Let $p\in\mathbb{Q}\left[  x\right]
$ be a polynomial of degree $\leq n$ such that
\[
p\left[  i\right]  =2^{i}\ \ \ \ \ \ \ \ \ \ \text{for all $i\in\left\{
0,1,\ldots,n\right\}  $}.
\]
Find $p\left[  n+1\right]  $.
\end{exercise}

\begin{exercise}
Let $n\in\mathbb{N}$. Let $p\in\mathbb{Q}\left[  x\right]  $ be a polynomial
of degree $\leq n$ such that
\[
p\left[  i\right]  =\dfrac{1}{\dbinom{n+1}{i}}\ \ \ \ \ \ \ \ \ \ \text{for
all $i\in\left\{  0,1,\ldots,n\right\}  $}.
\]
Find $p\left[  n+1\right]  $.
\end{exercise}

\begin{fineprint}
The next exercise gives a two-variable version of Lagrange interpolation (to
be more specific, of Theorem \ref{thm.polint.lagrange}):
\end{fineprint}

\begin{exercise}
\label{exe.polint.lagrange.2var}Let $F$ be a field. Consider the polynomial
ring $F\left[  x,y\right]  $ in two variables $x$ and $y$.

For any polynomial $p\in F\left[  x,y\right]  $, define the $x$-degree
$\deg_{x}p$ and the $y$-degree $\deg_{y}p$ as in Exercise
\ref{exe.polring.FTA-easy.2var}.

Let $n,m\in\mathbb{N}$.

Let $a_{0},a_{1},\ldots,a_{n}$ be $n+1$ distinct elements of $F$. Let
$b_{0},b_{1},\ldots,b_{m}$ be $m+1$ distinct elements of $F$. Let $c_{i,j}$ be
an element of $F$ for each pair $\left(  i,j\right)  \in\left\{
0,1,\ldots,n\right\}  \times\left\{  0,1,\ldots,m\right\}  $. Prove the following:

\begin{enumerate}
\item[\textbf{(a)}] There is a \textbf{unique} polynomial $p\in F\left[
x,y\right]  $ satisfying $\deg_{x}p\leq n$ and $\deg_{y}p\leq m$ and
\[
p\left[  a_{i},b_{j}\right]  =c_{i,j}\ \ \ \ \ \ \ \ \ \ \text{for all
}\left(  i,j\right)  \in\left\{  0,1,\ldots,n\right\}  \times\left\{
0,1,\ldots,m\right\}  .
\]

\item[\textbf{(b)}] This polynomial $p$ is given by
\[
p=\sum_{k=0}^{n}\ \ \sum_{\ell=0}^{m}c_{k,\ell}\dfrac{\prod_{u\neq k}\left(
x-a_{u}\right)  }{\prod_{u\neq k}\left(  a_{k}-a_{u}\right)  }\cdot
\dfrac{\prod_{v\neq\ell}\left(  y-b_{v}\right)  }{\prod_{v\neq\ell}\left(
b_{\ell}-b_{v}\right)  }%
\]
(where the \textquotedblleft$\prod_{u\neq k}$\textquotedblright\ sign means a
product over all $u\in\left\{  0,1,\ldots,n\right\}  $ satisfying $u\neq k$,
and where the \textquotedblleft$\prod_{v\neq\ell}$\textquotedblright\ sign
means a product over all $v\in\left\{  0,1,\ldots,m\right\}  $ satisfying
$v\neq\ell$).
\end{enumerate}
\end{exercise}

\begin{noncompile}
TODO: existence of a polynomial not equivalent to $a_{i}-a_{j}\mid b_{i}%
-b_{j}$

TODO: some apps of pol id trick
\end{noncompile}

\subsection{\label{sec.polys1.quot-of-Ralg}Intermezzo: quotients of
$R$-algebras}

In preparation for the next section, let me quickly introduce quotients of
$R$-algebras. I have previously defined quotients of rings modulo ideals, and
quotients of $R$-modules modulo submodules. These two concepts can be combined
to obtain quotients of $R$-algebras modulo ideals:

\begin{theorem}
\label{thm.algs.quotients}Let $A$ be an $R$-algebra. Let $I$ be an ideal of
$A$. Then:

\begin{enumerate}
\item[\textbf{(a)}] The ideal $I$ is also an $R$-submodule of $A$.

\item[\textbf{(b)}] The quotient ring $A/I$ and the quotient $R$-module $A/I$
fit together to form an $R$-algebra.

\item[\textbf{(c)}] The canonical projection $\pi:A\rightarrow A/I$ (which
sends each $a\in A$ to its residue class $\overline{a}=a+I$) is an $R$-algebra
morphism (from the original $R$-algebra $A$ to the $R$-algebra $A/I$ that we
just constructed in part \textbf{(b)}).
\end{enumerate}
\end{theorem}

\begin{proof}
\textbf{(a)} We already know that $I$ is closed under addition and contains
zero (since $I$ is an ideal). So we must only show that $I$ is closed under
scaling. In other words, we must show that $ri\in I$ for each $r\in R$ and
$i\in I$. But this is easy: If $r\in R$ and $i\in I$, then
\[
r\underbrace{i}_{=1_{A}\cdot i}=r\cdot1_{A}\cdot i=\underbrace{\left(
r\cdot1_{A}\right)  }_{\in A}\cdot\underbrace{i}_{\in I}\in
I\ \ \ \ \ \ \ \ \ \ \left(  \text{since }I\text{ is an ideal of }A\right)  .
\]

\textbf{(b)} LTTR. (You just need to verify the \textquotedblleft
scale-invariance of multiplication\textquotedblright\ axiom, but this is
straightforward.) \medskip

\textbf{(c)} We already know that this canonical projection is a ring morphism
and an $R$-module morphism; thus, it is an $R$-algebra morphism.
\end{proof}

\begin{definition}
Let $A$ and $I$ be as in Theorem \ref{thm.algs.quotients}. Then, the
$R$-algebra $A/I$ constructed in Theorem \ref{thm.algs.quotients} \textbf{(b)}
is called the \textbf{quotient algebra} (or \textbf{quotient }$R$%
\textbf{-algebra}) of $A$ by the ideal $I$.
\end{definition}

Let us next recall the universal property of quotient rings (in its two forms:
Theorem \ref{thm.quotring.uniprop1} and Theorem \ref{thm.quotring.uniprop2}).
This property is the tool of choice for constructing ring morphisms out of a
quotient ring. We can adapt this theorem to $R$-algebras with just trivial
modifications (alas, we have to rename $R$ and $S$ as $A$ and $B$, since $R$
already means something different):

\begin{theorem}
[Universal property of quotient algebras, elementwise form]%
\label{thm.quotalg.uniprop1}Let $A$ be an $R$-algebra. Let $I$ be an ideal of
$A$.

Let $B$ be an $R$-algebra. Let $f:A\rightarrow B$ be an $R$-algebra morphism.
Assume that $f\left(  I\right)  =0$ (this is shorthand for saying that
$f\left(  a\right)  =0$ for all $a\in I$). Then, the map%
\begin{align*}
f^{\prime}:A/I  &  \rightarrow B,\\
\overline{a}  &  \mapsto f\left(  a\right)  \ \ \ \ \ \ \ \ \ \ \left(
\text{for all }a\in A\right)
\end{align*}
is well-defined (i.e., the value $f\left(  a\right)  $ depends only on the
residue class $\overline{a}$, not on $a$ itself) and is an $R$-algebra morphism.
\end{theorem}

\begin{proof}
Adapt the argument that we used to prove Theorem \ref{thm.quotring.uniprop1}.
The only new thing we need to check is that the map $f^{\prime}$ constructed
in the proof is $R$-linear; but this is just as straightforward as showing
that this map is a ring morphism.
\end{proof}

\begin{theorem}
[Universal property of quotient algebras, abstract form]%
\label{thm.quotalg.uniprop2}Let $A$ be an $R$-algebra. Let $I$ be an ideal of
$A$. Consider the canonical projection $\pi:A\rightarrow A/I$.

Let $B$ be an $R$-algebra. Let $f:A\rightarrow B$ be an $R$-algebra morphism.
Assume that $f\left(  I\right)  =0$ (this is shorthand for saying that
$f\left(  a\right)  =0$ for all $a\in I$). Then, there is a unique $R$-algebra
morphism $f^{\prime}:A/I\rightarrow B$ satisfying $f=f^{\prime}\circ\pi$.
\end{theorem}

\begin{proof}
Adapt the argument that we used to prove Theorem \ref{thm.quotring.uniprop2}.
\end{proof}

The First Isomorphism Theorem for rings (Theorem \ref{thm.1it.ring1}) also has
an analogue for $R$-algebras. We leave it to the reader to state it.

\subsection{\label{sec.polys1.adjroot}Adjoining roots}

\subsubsection{Examples}

What is a complex number? Nowadays, the complex numbers are commonly defined
as pairs of real numbers; this is a fairly straightforward process (first you
define addition and multiplication and zero and unity; then you show that the
ring axioms hold) and can be found in many textbooks (e.g., \cite[\S 4.1]%
{19s}). But this is the modern definition. When Girolamo Cardano originally
invented complex numbers back in the 16th century, he had a different vision:
Cardano essentially proposed to \textbf{imagine} that there is a new number
called $i$ that satisfies $i^{2}=-1$ but otherwise behaves like the numbers we
know. Thus, you're allowed to form arbitrary polynomials in $i$, but you have
to equate $i^{2}$ to $-1$, so you never end up getting anything more
complicated than numbers of the form $a+bi$ with $a,b\in\mathbb{R}$ (since any
higher power of $i$ can be reduced to $\pm1$ or $\pm i$ using the $i^{2}=-1$
rule). Thus, it makes sense to encode complex numbers as pairs of real
numbers, but this is merely one way of encoding them.\footnote{I am being
sloppy with the history here. The relevant source is Girolamo Cardano's 1545
book \textit{Ars magna}, specifically its Chapter XXXVII, in which he asks the
reader to \textquotedblleft imagine $\sqrt{-15}$\textquotedblright. Of course,
this is not much different from imagining $i=\sqrt{-1}$, since a square root
of $-15$ could be obtained from a square root of $-1$ by multiplying with the
real number $\sqrt{15}$.
\par
Cardano was writing at a time when even the notion of a negative number was
far from widely accepted in the West (though known in India and Persia).
Cardano himself called negative numbers \textquotedblleft
fictitious\textquotedblright\ (arguably an improvement from previous European
authors, who called them \textquotedblleft absurd\textquotedblright), and did
not quite treat them as first-class numbers. His \textit{Ars magna} is often
considered to be the first serious treatment of negative numbers written in
Europe. That the very same book introduces complex numbers is thus an example
of the \textquotedblleft when it rains, it pours\textquotedblright\ phenomenon
in the history of ideas.
\par
Cardano did not introduce the name $i$ for the imaginary unit $\sqrt{-1}$.
This was done by Euler much later.}

Of course, Cardano's original vision is not a rigorous definition; just as
easily you could introduce a number $j$ satisfying $0j=1$, and thus collapse
the entire number system (since this new number would let you argue that
$1=0j=\left(  0+0\right)  j=0j+0j=1+1=2$). So, if we want to make Cardano's
definition rigorous, we have to rewrite it algebraically. One way to do this
is to define $\mathbb{C}$ as the quotient ring%
\[
\mathbb{R}\left[  x\right]  /\left(  x^{2}+1\right)  \mathbb{R}\left[
x\right]  .
\]
In fact, we start with $\mathbb{R}\left[  x\right]  $ because our complex
numbers should be polynomials in a single symbol $i$ (which will be
represented by the indeterminate $x$ in $\mathbb{R}\left[  x\right]  $); but
then we quotient out the ideal $\left(  x^{2}+1\right)  \mathbb{R}\left[
x\right]  $ since we want $i^{2}+1$ (and thus also each multiple of $i^{2}+1$)
to be $0$ in our complex numbers.

To be on the safe side, let us show that this quotient ring $\mathbb{R}\left[
x\right]  /\left(  x^{2}+1\right)  \mathbb{R}\left[  x\right]  $ is isomorphic
to the complex numbers $\mathbb{C}$ as we know them (i.e., defined in the
modern way, as pairs of real numbers).

First of all, we introduce a shorthand:

\begin{convention}
\label{conv.quotring.R/a}If $R$ is a commutative ring, and if $a\in R$, then
the quotient ring $R/aR$ will be abbreviated as $R/a$. We are already using a
particular case of this notation, as we are writing $\mathbb{Z}/n$ for
$\mathbb{Z}/n\mathbb{Z}$ when $n$ is an integer.

We note that the quotient ring $R/a=R/aR$ is not just a ring, but an
$R$-algebra as well (by Theorem \ref{thm.algs.quotients} \textbf{(b)}).
Furthermore, the $R$-algebra $R/a$ is commutative (since it is a quotient of
$R$). This all will be tacitly used in what follows.
\end{convention}

So we want to prove that $\mathbb{R}\left[  x\right]  /\left(  x^{2}+1\right)
\cong\mathbb{C}$ as rings -- and even better, as $\mathbb{R}$-algebras. Let's
be a little bit more precise:

\begin{proposition}
\label{prop.fieldext.C}We have $\mathbb{R}\left[  x\right]  /\left(
x^{2}+1\right)  \cong\mathbb{C}$ as $\mathbb{R}$-algebras. More concretely:
There is an $\mathbb{R}$-algebra isomorphism%
\begin{align*}
\mathbb{R}\left[  x\right]  /\left(  x^{2}+1\right)   &  \rightarrow
\mathbb{C},\\
\overline{p}  &  \mapsto p\left[  i\right]  .
\end{align*}

\end{proposition}

\begin{proof}
We already know that $\mathbb{C}$ is an $\mathbb{R}$-algebra. Thus, Theorem
\ref{thm.polring.univar-sub-hom} (applied to $R=\mathbb{R}$ and $A=\mathbb{C}$
and $a=i$) yields that the map%
\begin{align*}
f:\mathbb{R}\left[  x\right]   &  \rightarrow\mathbb{C},\\
p  &  \mapsto p\left[  i\right]
\end{align*}
is an $\mathbb{R}$-algebra morphism. This map $f$ sends the principal ideal
$\left(  x^{2}+1\right)  \mathbb{R}\left[  x\right]  $ to $0$, because for
each $q\in\mathbb{R}\left[  x\right]  $, we have%
\[
f\left(  \left(  x^{2}+1\right)  \cdot q\right)  =\left(  \left(
x^{2}+1\right)  \cdot q\right)  \left[  i\right]  =\underbrace{\left(
i^{2}+1\right)  }_{=0}\cdot\,q\left[  i\right]  =0.
\]
Hence, Theorem \ref{thm.quotalg.uniprop1} (applied to $R=\mathbb{R}$,
$A=\mathbb{R}\left[  x\right]  $, $I=\left(  x^{2}+1\right)  \mathbb{R}\left[
x\right]  $ and $B=\mathbb{C}$) shows that the map%
\begin{align*}
f^{\prime}:\mathbb{R}\left[  x\right]  /\left(  x^{2}+1\right)   &
\rightarrow\mathbb{C},\\
\overline{a}  &  \mapsto f\left(  a\right)
\end{align*}
is well-defined and is an $\mathbb{R}$-algebra morphism. Consider this map
$f^{\prime}$. Each $p\in\mathbb{R}\left[  x\right]  $ satisfies
\begin{align}
f^{\prime}\left(  \overline{p}\right)   &  =f\left(  p\right)
\ \ \ \ \ \ \ \ \ \ \left(  \text{by the definition of }f^{\prime}\right)
\nonumber\\
&  =p\left[  i\right]  \ \ \ \ \ \ \ \ \ \ \left(  \text{by the definition of
}f\right)  . \label{pf.prop.fieldext.C.f'p=}%
\end{align}

Now, why is $f^{\prime}$ an isomorphism?

It's not hard to see that $f^{\prime}$ is surjective: Indeed, any
$z\in\mathbb{C}$ can be written as $z=a+bi$ for some $a,b\in\mathbb{R}$, and
then we have $z=a+bi=f^{\prime}\left(  \overline{a+bx}\right)  $ (since
(\ref{pf.prop.fieldext.C.f'p=}) yields $f^{\prime}\left(  \overline
{a+bx}\right)  =\left(  a+bx\right)  \left[  i\right]  =a+bi$).

Now, how can we prove that $f^{\prime}$ is injective? Since $f^{\prime}$ is
$\mathbb{R}$-linear, it suffices to show that $\operatorname*{Ker}\left(
f^{\prime}\right)  =\left\{  0\right\}  $ (by Lemma \ref{lem.modmor.ker-inj}).

Let $u\in\operatorname*{Ker}\left(  f^{\prime}\right)  $. Thus, $u\in
\mathbb{R}\left[  x\right]  /\left(  x^{2}+1\right)  $, so that $u=\overline
{p}$ for some $p\in\mathbb{R}\left[  x\right]  $. Consider this $p$.

However, Theorem \ref{thm.polring.univar-quorem} \textbf{(a)} (applied to
$R=\mathbb{R}$, $b=x^{2}+1$ and $a=p$) yields that there is a unique pair
$\left(  q,r\right)  $ of polynomials in $\mathbb{R}\left[  x\right]  $ such
that
\[
p=q\cdot\left(  x^{2}+1\right)  +r\ \ \ \ \ \ \ \ \ \ \text{and}%
\ \ \ \ \ \ \ \ \ \ \deg r<\deg\left(  x^{2}+1\right)  .
\]
Consider this pair $\left(  q,r\right)  $. From $\deg r<\deg\left(
x^{2}+1\right)  =2$, we see that the polynomial $r$ can be written as $a+bx$
for some $a,b\in\mathbb{R}$. Consider these $a,b$. From $p=q\cdot\left(
x^{2}+1\right)  +r$, we obtain $p-r=q\cdot\left(  x^{2}+1\right)  \in\left(
x^{2}+1\right)  \mathbb{R}\left[  x\right]  $; thus, $\overline{p}%
=\overline{r}$ in the quotient ring $\mathbb{R}\left[  x\right]  /\left(
x^{2}+1\right)  $. Now,
\begin{align*}
u  &  =\overline{p}=\overline{r}=\overline{a+bx}\ \ \ \ \ \ \ \ \ \ \left(
\text{since }r=a+bx\right)  ,\ \ \ \ \ \ \ \ \ \ \text{so that}\\
f^{\prime}\left(  u\right)   &  =f^{\prime}\left(  \overline{a+bx}\right)
=\left(  a+bx\right)  \left[  i\right]  \ \ \ \ \ \ \ \ \ \ \left(  \text{by
(\ref{pf.prop.fieldext.C.f'p=})}\right) \\
&  =a+bi.
\end{align*}
Hence, $a+bi=f^{\prime}\left(  u\right)  =0$ (since $u\in\operatorname*{Ker}%
\left(  f^{\prime}\right)  $). Since $a,b\in\mathbb{R}$, this entails $a=b=0$
(since the complex numbers $1$ and $i$ are $\mathbb{R}$-linearly independent).
Thus, $u=\overline{a+bx}$ rewrites as $u=\overline{0+0x}=0\in\left\{
0\right\}  $.

Forget that we fixed $u$. We thus have shown that $u\in\left\{  0\right\}  $
for each $u\in\operatorname*{Ker}\left(  f^{\prime}\right)  $. In other words,
$\operatorname*{Ker}\left(  f^{\prime}\right)  \subseteq\left\{  0\right\}  $.
Since the reverse inclusion $\left\{  0\right\}  \subseteq\operatorname*{Ker}%
\left(  f^{\prime}\right)  $ is obvious, we thus conclude that
$\operatorname*{Ker}\left(  f^{\prime}\right)  =\left\{  0\right\}  $. As we
have said, this entails that $f^{\prime}$ is injective.

Now we know that the map $f^{\prime}$ is injective and surjective. Hence,
$f^{\prime}$ is bijective, i.e., invertible. Since every invertible
$\mathbb{R}$-algebra morphism is an $\mathbb{R}$-algebra isomorphism (by
Proposition \ref{prop.algmor.invertible-iso}), we thus conclude that
$f^{\prime}$ is an $\mathbb{R}$-algebra isomorphism. This proves Proposition
\ref{prop.fieldext.C} (since the map $f^{\prime}$ satisfies $f^{\prime}\left(
\overline{p}\right)  =p\left[  i\right]  $ for each $p\in\mathbb{R}\left[
x\right]  $, and thus is precisely the alleged isomorphism claimed in
Proposition \ref{prop.fieldext.C}).
\end{proof}

Note the use of polynomial division (with remainder) in our above proof of
Proposition \ref{prop.fieldext.C}. It has a natural usefulness in the study of
quotient rings of $\mathbb{R}\left[  x\right]  $, just as integer division
(with remainder) is crucial to the study of quotient rings of $\mathbb{Z}$.

Similarly to Proposition \ref{prop.fieldext.C}, we can reveal further quotient
rings of polynomial rings as certain rings we know:

\begin{proposition}
\label{prop.fieldext.Zi}\ \ 

\begin{enumerate}
\item[\textbf{(a)}] Recall the ring $\mathbb{Z}\left[  i\right]  $ of Gaussian
integers. We have $\mathbb{Z}\left[  x\right]  /\left(  x^{2}+1\right)
\cong\mathbb{Z}\left[  i\right]  $ as $\mathbb{Z}$-algebras. More concretely:
There is a $\mathbb{Z}$-algebra isomorphism%
\begin{align*}
\mathbb{Z}\left[  x\right]  /\left(  x^{2}+1\right)   &  \rightarrow
\mathbb{Z}\left[  i\right]  ,\\
\overline{p}  &  \mapsto p\left[  i\right]  .
\end{align*}

\item[\textbf{(b)}] Recall the ring $\mathbb{S}=\mathbb{Q}\left[  \sqrt
{5}\right]  =\left\{  a+b\sqrt{5}\ \mid\ a,b\in\mathbb{Q}\right\}  $ (a
subring of $\mathbb{R}$). We have $\mathbb{Q}\left[  x\right]  /\left(
x^{2}-5\right)  \cong\mathbb{S}$ as $\mathbb{Q}$-algebras. More concretely:
There is a $\mathbb{Q}$-algebra isomorphism%
\begin{align*}
\mathbb{Q}\left[  x\right]  /\left(  x^{2}-5\right)   &  \rightarrow
\mathbb{S},\\
\overline{p}  &  \mapsto p\left[  \sqrt{5}\right]  .
\end{align*}

\end{enumerate}
\end{proposition}

\begin{proof}
\textbf{(a)} Analogous to the proof of Proposition \ref{prop.fieldext.C}.

\textbf{(b)} Analogous to the proof of Proposition \ref{prop.fieldext.C}.
\end{proof}

Proposition \ref{prop.fieldext.C} and Proposition \ref{prop.fieldext.Zi}
suggest that when we start with a ring $R$ and a polynomial $b\in R\left[
x\right]  $, then the quotient ring $R\left[  x\right]  /b$ is (in some way)
an \textquotedblleft extension\textquotedblright\ of $R$ by a root of $b$, in
the sense that it contains $R$ as a subring (at least up to isomorphism) but
also contains a root of $b$ (namely, $\overline{x}$). Thus, we can hope that
by taking the quotient ring $R\left[  x\right]  /b$, we can \textquotedblleft
adjoin\textquotedblright\ a root of $b$ to the ring $R$ even if $b$ has no
root in $R$ (just as Cardano defined the complex numbers by \textquotedblleft
adjoining\textquotedblright\ a root of $x^{2}+1$ to $\mathbb{R}$). The word
\textquotedblleft adjoin\textquotedblright\ here means something like
\textquotedblleft insert\textquotedblright, \textquotedblleft
attach\textquotedblright\ or \textquotedblleft throw in\textquotedblright.

This is a good intuition, but there are nuances: In the process of
\textquotedblleft adjoining\textquotedblright\ our root to our ring $R$, we
may end up making $R$ \textquotedblleft smaller\textquotedblright, in the
sense that different elements of $R$ become equal when the root is
\textquotedblleft adjoined\textquotedblright\ (and thus the resulting ring is
not really an \textquotedblleft extension\textquotedblright\ of $R$). The
following example (in which we take a quotient of $\mathbb{Z}\left[  x\right]
$ by a constant polynomial) demonstrates this:

\begin{proposition}
\label{prop.fieldext.Zx/2}\ \ 

\begin{enumerate}
\item[\textbf{(a)}] We have $\left(  \mathbb{Z}\left[  x\right]  \right)
/m\cong\left(  \mathbb{Z}/m\right)  \left[  x\right]  $ as $\mathbb{Z}%
$-algebras (i.e., as rings) for any integer $m$.

\item[\textbf{(b)}] The ring $\left(  \mathbb{Z}\left[  x\right]  \right)  /1$
is trivial.
\end{enumerate}
\end{proposition}

\begin{proof}
[Proof sketch.]\textbf{(a)} Let $m$ be an integer. Then, the principal ideal
$m\mathbb{Z}\left[  x\right]  $ of $\mathbb{Z}\left[  x\right]  $ consists of
all polynomials whose all coefficients are multiples of $m$. Thus, it is easy
to see that the map
\begin{align*}
f:\left(  \mathbb{Z}\left[  x\right]  \right)  /m  &  \rightarrow\left(
\mathbb{Z}/m\right)  \left[  x\right]  ,\\
\overline{a_{0}x^{0}+a_{1}x^{1}+a_{2}x^{2}+\cdots}  &  \mapsto\overline{a_{0}%
}x^{0}+\overline{a_{1}}x^{1}+\overline{a_{2}}x^{2}+\cdots
\end{align*}
is well-defined and is a $\mathbb{Z}$-algebra isomorphism. This proves
Proposition \ref{prop.fieldext.Zx/2} \textbf{(a)}. \medskip

\textbf{(b)} More generally: If $R$ is any commutative ring, then the ring
$R/1$ is trivial. This is because the principal ideal $1R$ of $R$ is the whole
ring $R$, so there is only one coset modulo this ideal.
\end{proof}

Proposition \ref{prop.fieldext.Zx/2} \textbf{(a)} (applied to $m=2$) shows
that if we take the quotient ring of $\mathbb{Z}\left[  x\right]  $ modulo
(the principal ideal generated by) the constant polynomial $2$, then we don't
get an \textquotedblleft extension\textquotedblright\ of $\mathbb{Z}$; what we
instead get is the polynomial ring $\left(  \mathbb{Z}/2\right)  \left[
x\right]  $, in which (unlike in $\mathbb{Z}$) we have $1+1=0$ (so it
certainly cannot contain a copy of $\mathbb{Z}$ as a subring). But if you
think about this carefully, you will realize that this perfectly agrees with
the idea of \textquotedblleft adjoining a root\textquotedblright. Indeed, to
\textquotedblleft adjoin\textquotedblright\ a root of the constant polynomial
$2$ to $\mathbb{Z}$ means to introduce a new \textquotedblleft
number\textquotedblright\ $x$ satisfying $2=0$. The equation $2=0$ tells us
nothing about the \textquotedblleft number\textquotedblright\ $x$ (so it
remains completely unconstrained), but collapses all even integers to $0$,
thus leaving us with the ring $\left(  \mathbb{Z}/2\right)  \left[  x\right]
$. This is precisely what Proposition \ref{prop.fieldext.Zx/2} \textbf{(a)}
told us. Likewise, \textquotedblleft adjoining\textquotedblright\ a root of
$1$ to $\mathbb{Z}$ causes $1=0$, which renders the ring trivial (since any
element of a ring is a multiple of $1$); this agrees with Proposition
\ref{prop.fieldext.Zx/2} \textbf{(b)}.

The examples so far have taught us that -- yes -- we can \textquotedblleft
adjoin\textquotedblright\ a root of any polynomial to a commutative ring $R$,
but we don't always get an extension of $R$ (although we do always get an
$R$-algebra). In Theorem \ref{thm.fieldext.basis-monic} \textbf{(c)}, we will
see a (sufficient) criterion for when we do. \medskip

Here is another natural question: What happens if we \textquotedblleft
adjoin\textquotedblright\ a root of a polynomial $b$ that already has a root
in $R$ ? For example, let us take the polynomial $x^{2}-1$ over $\mathbb{Q}$
(which has $1$ and $-1$ as roots). It turns out that the resulting quotient
ring $\mathbb{Q}\left[  x\right]  /\left(  x^{2}-1\right)  $ is a good friend
of ours by now:

\begin{proposition}
\label{prop.fieldext.xx-1}Recall the group algebra $\mathbb{Q}\left[
C_{2}\right]  $ of the cyclic group $C_{2}$ from Example \ref{exa.monalg.QC2}.
Then,%
\[
\mathbb{Q}\left[  x\right]  /\left(  x^{2}-1\right)  \cong\mathbb{Q}\left[
C_{2}\right]  \cong\mathbb{Q}\times\mathbb{Q}\ \ \ \ \ \ \ \ \ \ \text{as
}\mathbb{Q}\text{-algebras.}%
\]

\end{proposition}

\begin{proof}
In Example \ref{exa.monalg.QC2}, we have seen that the group algebra
$\mathbb{Q}\left[  C_{2}\right]  $ has a basis $\left(  e_{1},e_{u}\right)  $
(as a $\mathbb{Q}$-module). By Convention \ref{conv.monalg.m-as-em}, we can
write $1$ and $u$ for $e_{1}$ and $e_{u}$, so that this basis becomes $\left(
1,u\right)  $. We also know (from Example \ref{exa.monalg.QC2}) that
$\mathbb{Q}\left[  C_{2}\right]  \cong\mathbb{Q}\times\mathbb{Q}$ as
$\mathbb{Q}$-algebras. It thus remains to prove that $\mathbb{Q}\left[
x\right]  /\left(  x^{2}-1\right)  \cong\mathbb{Q}\left[  C_{2}\right]  $.

Note the similarity between $\mathbb{Q}\left[  C_{2}\right]  $ and
$\mathbb{C}$:

\begin{itemize}
\item The $\mathbb{Q}$-module $\mathbb{Q}\left[  C_{2}\right]  $ has basis
$\left(  1,u\right)  $, with $u^{2}=1$.

\item The $\mathbb{R}$-module $\mathbb{C}$ has basis $\left(  1,i\right)  $,
with $i^{2}=-1$.
\end{itemize}

This suggests that we just copypaste our above proof of Proposition
\ref{prop.fieldext.C}, replacing $\mathbb{R}$, $\mathbb{C}$ and $i$ by
$\mathbb{Q}$, $\mathbb{Q}\left[  C_{2}\right]  $ and $u$ and occasionally
flipping signs (and renaming the element $u\in\operatorname*{Ker}\left(
f^{\prime}\right)  $ as $v$ to avoid clashing with the unrelated generator $u$
of $C_{2}$). This is precisely what we are now going to do (but in a smaller
font, to avoid wasting paper).

\begin{fineprint}
Theorem \ref{thm.polring.univar-sub-hom} (applied to $R=\mathbb{Q}$ and
$A=\mathbb{Q}\left[  C_{2}\right]  $ and $a=u$) yields that the map%
\begin{align*}
f:\mathbb{Q}\left[  x\right]   &  \rightarrow\mathbb{Q}\left[  C_{2}\right]
,\\
p  &  \mapsto p\left[  u\right]
\end{align*}
is a $\mathbb{Q}$-algebra morphism. This map $f$ sends the principal ideal
$\left(  x^{2}-1\right)  \mathbb{Q}\left[  x\right]  $ to $0$, because for
each $q\in\mathbb{Q}\left[  x\right]  $, we have%
\[
f\left(  \left(  x^{2}-1\right)  \cdot q\right)  =\left(  \left(
x^{2}-1\right)  \cdot q\right)  \left[  u\right]  =\underbrace{\left(
u^{2}-1\right)  }_{\substack{=0\\\text{(since }u^{2}=1\text{)}}}\cdot
\,q\left[  u\right]  =0.
\]
Hence, Theorem \ref{thm.quotalg.uniprop1} (applied to $R=\mathbb{Q}$,
$A=\mathbb{Q}\left[  x\right]  $, $I=\left(  x^{2}-1\right)  \mathbb{Q}\left[
x\right]  $ and $B=\mathbb{Q}\left[  C_{2}\right]  $) shows that the map%
\begin{align*}
f^{\prime}:\mathbb{Q}\left[  x\right]  /\left(  x^{2}-1\right)   &
\rightarrow\mathbb{Q}\left[  C_{2}\right]  ,\\
\overline{a}  &  \mapsto f\left(  a\right)
\end{align*}
is well-defined and is a $\mathbb{Q}$-algebra morphism. Consider this
$f^{\prime}$. Each $p\in\mathbb{Q}\left[  x\right]  $ satisfies%
\begin{align}
f^{\prime}\left(  \overline{p}\right)   &  =f\left(  p\right)
\ \ \ \ \ \ \ \ \ \ \left(  \text{by the definition of }f^{\prime}\right)
\nonumber\\
&  =p\left[  u\right]  \ \ \ \ \ \ \ \ \ \ \left(  \text{by the definition of
}f\right)  . \label{pf.prop.fieldext.xx-1.f'p=}%
\end{align}

Now, why is $f^{\prime}$ an isomorphism?

It's not hard to see that $f^{\prime}$ is surjective: Indeed, any
$z\in\mathbb{Q}\left[  C_{2}\right]  $ can be written as $z=a+bu$ for some
$a,b\in\mathbb{Q}$, and then we have $z=a+bu=f^{\prime}\left(  \overline
{a+bx}\right)  $ (since (\ref{pf.prop.fieldext.xx-1.f'p=}) yields $f^{\prime
}\left(  \overline{a+bx}\right)  =\left(  a+bx\right)  \left[  u\right]
=a+bu$).

Now, how can we prove that $f^{\prime}$ is injective? Since $f^{\prime}$ is
$\mathbb{Q}$-linear, it suffices to show that $\operatorname*{Ker}\left(
f^{\prime}\right)  =\left\{  0\right\}  $ (by Lemma \ref{lem.modmor.ker-inj}).

Let $v\in\operatorname*{Ker}\left(  f^{\prime}\right)  $. Thus, $v\in
\mathbb{Q}\left[  x\right]  /\left(  x^{2}-1\right)  $, so that $v=\overline
{p}$ for some $p\in\mathbb{Q}\left[  x\right]  $. Consider this $p$.

However, Theorem \ref{thm.polring.univar-quorem} \textbf{(a)} (applied to
$R=\mathbb{Q}$, $b=x^{2}-1$ and $a=p$) yields that there is a unique pair
$\left(  q,r\right)  $ of polynomials in $\mathbb{Q}\left[  x\right]  $ such
that
\[
p=q\cdot\left(  x^{2}-1\right)  +r\ \ \ \ \ \ \ \ \ \ \text{and}%
\ \ \ \ \ \ \ \ \ \ \deg r<\deg\left(  x^{2}-1\right)  .
\]
Consider this pair $\left(  q,r\right)  $. From $\deg r<\deg\left(
x^{2}-1\right)  =2$, we see that the polynomial $r$ can be written as $a+bx$
for some $a,b\in\mathbb{Q}$. Consider these $a,b$. From $p=q\cdot\left(
x^{2}-1\right)  +r$, we obtain $p-r=q\cdot\left(  x^{2}-1\right)  \in\left(
x^{2}-1\right)  \mathbb{Q}\left[  x\right]  $; thus, $\overline{p}%
=\overline{r}$ in the quotient ring $\mathbb{Q}\left[  x\right]  /\left(
x^{2}-1\right)  $. Now,
\begin{align*}
v  &  =\overline{p}=\overline{r}=\overline{a+bx}\ \ \ \ \ \ \ \ \ \ \left(
\text{since }r=a+bx\right)  ,\ \ \ \ \ \ \ \ \ \ \text{so that}\\
f^{\prime}\left(  v\right)   &  =f^{\prime}\left(  \overline{a+bx}\right)
=\left(  a+bx\right)  \left[  u\right]  \ \ \ \ \ \ \ \ \ \ \left(  \text{by
(\ref{pf.prop.fieldext.xx-1.f'p=})}\right) \\
&  =a+bu.
\end{align*}
Hence, $a+bu=f^{\prime}\left(  v\right)  =0$ (since $v\in\operatorname*{Ker}%
\left(  f^{\prime}\right)  $). Since $a,b\in\mathbb{Q}$, this entails $a=b=0$
(since the vectors $1$ and $u$ in $\mathbb{Q}\left[  C_{2}\right]  $ are
$\mathbb{Q}$-linearly independent). Thus, $v=\overline{a+bx}$ rewrites as
$v=\overline{0+0x}=0\in\left\{  0\right\}  $.

Forget that we fixed $v$. We thus have shown that $v\in\left\{  0\right\}  $
for each $v\in\operatorname*{Ker}\left(  f^{\prime}\right)  $. In other words,
$\operatorname*{Ker}\left(  f^{\prime}\right)  \subseteq\left\{  0\right\}  $.
Since the reverse inclusion $\left\{  0\right\}  \subseteq\operatorname*{Ker}%
\left(  f^{\prime}\right)  $ is obvious, we thus conclude that
$\operatorname*{Ker}\left(  f^{\prime}\right)  =\left\{  0\right\}  $. As we
have said, this entails that $f^{\prime}$ is injective.

Now we know that the map $f^{\prime}$ is injective and surjective. Hence,
$f^{\prime}$ is bijective, i.e., invertible. Since every invertible
$\mathbb{Q}$-algebra morphism is a $\mathbb{Q}$-algebra isomorphism (by
Proposition \ref{prop.algmor.invertible-iso}), we thus conclude that
$f^{\prime}$ is an $\mathbb{Q}$-algebra isomorphism. Hence, $\mathbb{Q}\left[
x\right]  /\left(  x^{2}-1\right)  \cong\mathbb{Q}\left[  C_{2}\right]  $. As
we said, this proves Proposition \ref{prop.fieldext.xx-1}.
\end{fineprint}
\end{proof}

\begin{exercise}
\label{exe.fieldext.dualnums}Let $R$ be a commutative ring. Recall the
$R$-algebra $\mathbb{D}_{R}$ defined in Exercise \ref{exe.algs.dualnums},
along with the element $\varepsilon=\left(  0,1\right)  \in\mathbb{D}_{R}$
defined ibidem. Prove that there is an $R$-algebra isomorphism%
\begin{align*}
R\left[  x\right]  /\left(  x^{2}\right)   &  \rightarrow\mathbb{D}_{R},\\
\overline{p}  &  \mapsto p\left[  \varepsilon\right]  .
\end{align*}

\end{exercise}

\begin{exercise}
\label{exe.21hw3.9}Let $\varphi$ be the \textbf{golden ratio} -- i.e., the
real number $\dfrac{1+\sqrt{5}}{2}\approx1.618\ldots$.

Let $\mathbb{Z}\left[  \varphi\right]  $ be the set of all reals of the form
$a+b\varphi$ with $a,b\in\mathbb{Z}$.

\begin{enumerate}
\item[\textbf{(a)}] Prove that $\mathbb{Z}\left[  \varphi\right]  $ is a
subring of $\mathbb{R}$.

\item[\textbf{(b)}] Prove that
\[
\mathbb{Z}\left[  \varphi\right]  \cong\mathcal{F}\cong\mathbb{Z}\left[
x\right]  /\left(  x^{2}-x-1\right)  \qquad\text{as rings,}%
\]
where $\mathcal{F}$ is the ring defined in Exercise \ref{exe.21hw1.6}.
\end{enumerate}
\end{exercise}

\begin{exercise}
\label{exe.monalg.QC3.fieldext-S}Let $R$ be the commutative group algebra
$\mathbb{Q}\left[  C_{3}\right]  $ discussed in Example \ref{exa.monalg.QC3}.
Consider its idempotent element $z=\dfrac{1+e_{u}+e_{v}}{3}=\dfrac{1+u+v}{3}$.
Let $S$ be the principal ideal $\left(  1-z\right)  R$ of $R$. As we know from
Exercise \ref{exe.21hw1.3} \textbf{(b)} (applied to $e=z$), this principal
ideal $S$ is itself a ring, with addition and multiplication inherited from
$R$ and with zero $0_{R}$ and with unity $1-z$.

\begin{enumerate}
\item[\textbf{(a)}] Prove that this ring $S$ is isomorphic to $\mathbb{Q}%
\left[  x\right]  /\left(  x^{2}+x+1\right)  $ as rings.

\item[\textbf{(b)}] Prove that $S$ is furthermore isomorphic to the subring
\[
\mathbb{Q}\left[  \sqrt{-3}\right]  =\left\{  a+b\sqrt{-3}\ \mid
\ a,b\in\mathbb{Q}\right\}  \ \ \ \ \ \ \ \ \ \ \text{of }\mathbb{C}\text{.}%
\]

\end{enumerate}

[\textbf{Hint:} For part \textbf{(a)}, first show that
\[
S=\left\{  a+bu+cv\ \mid\ a,b,c\in\mathbb{Q}\text{ with }a+b+c=0\right\}  .
\]
Then, show that the $\mathbb{Q}$-algebra morphism%
\begin{align*}
f:\mathbb{Q}\left[  x\right]  /\left(  x^{2}+x+1\right)   &  \rightarrow S,\\
\overline{p}  &  \mapsto p\left[  u\left(  1-z\right)  \right]
\end{align*}
is well-defined and invertible. The quickest way to verify invertibility is
using linear algebra over $\mathbb{Q}$, as $f$ is a $\mathbb{Q}$-linear map
between two $2$-dimensional $\mathbb{Q}$-vector spaces.

For part \textbf{(b)}, find a root of the polynomial $x^{2}+x+1$ in
$\mathbb{Q}\left[  \sqrt{-3}\right]  $.]
\end{exercise}

In our proofs of Propositions \ref{prop.fieldext.C}, \ref{prop.fieldext.xx-1}
and \ref{prop.fieldext.Zi} (even though I left the latter to the reader), we
used that the leading coefficients of the polynomials we were quotienting out
were units. Indeed, this is what allowed us to apply Theorem
\ref{thm.polring.univar-quorem} \textbf{(a)}, which was a crucial step in
proving that $f^{\prime}$ is injective. Describing quotient rings becomes much
more complicated when the leading coefficient of the polynomial is not a unit.
Sometimes it is nevertheless possible. Here is a particularly well-behaved example:

\begin{proposition}
\label{prop.fieldext.mx-1}Fix a nonzero integer $m$. Define the ring $R_{m}$
as in Exercise \ref{exe.21hw1.1}; that is, $R_{m}$ is the subring%
\[
\left\{  r\in\mathbb{Q}\ \mid\ \text{there exists }k\in\mathbb{N}\text{
satisfying }m^{k}r\in\mathbb{Z}\right\}
\]
of $\mathbb{Q}$. Then,%
\[
\mathbb{Z}\left[  x\right]  /\left(  mx-1\right)  \cong R_{m}%
\ \ \ \ \ \ \ \ \ \ \text{as }\mathbb{Z}\text{-algebras (i.e., as rings).}%
\]
More concretely: There is a $\mathbb{Z}$-algebra isomorphism%
\begin{align*}
\mathbb{Z}\left[  x\right]  /\left(  mx-1\right)   &  \rightarrow R_{m},\\
\overline{p}  &  \mapsto p\left[  \dfrac{1}{m}\right]  .
\end{align*}

\end{proposition}

\begin{proof}
[Proof sketch.]Intuitively, this should be exactly what you expect: According
to our \textquotedblleft adjoining roots\textquotedblright\ philosophy, the
ring $\mathbb{Z}\left[  x\right]  /\left(  mx-1\right)  $ is what you get if
you \textquotedblleft adjoin\textquotedblright\ a root of the polynomial
$mx-1$ to $\mathbb{Z}$. But such a root would behave like the rational number
$\dfrac{1}{m}$; so it is no surprise that the resulting ring would be
isomorphic to $R_{m}$ (since $R_{m}$ is really just \textquotedblleft the
numbers you can get if you start with the integers and also allow multiplying
by $\dfrac{1}{m}$\textquotedblright). This, of course, is not a proof.

An actual proof can be done along the following lines:

\begin{enumerate}
\item Show that a $\mathbb{Z}$-algebra morphism%
\begin{align*}
\alpha:\mathbb{Z}\left[  x\right]  /\left(  mx-1\right)   &  \rightarrow
R_{m},\\
\overline{p}  &  \mapsto p\left[  \dfrac{1}{m}\right]
\end{align*}
exists. This is similar to the corresponding part of the proof of Proposition
\ref{prop.fieldext.C} (where we called the corresponding morphism $f^{\prime}$
rather than $\alpha$); the main roles are played by Theorem
\ref{thm.polring.univar-sub-hom} and Theorem \ref{thm.quotalg.uniprop1}.

\item (Optional:) Show that this morphism $\alpha$ is surjective. (In fact,
each element of $R_{m}$ has the form $\dfrac{a}{m^{k}}$ for some
$a\in\mathbb{Z}$ and some $k\in\mathbb{N}$, and thus equals $\alpha\left(
\overline{ax^{k}}\right)  $.)

\item Don't waste your time trying to show that $\alpha$ is injective; there
is no quick way to prove this directly.

\item Show that there is a map
\begin{align*}
\beta:R_{m}  &  \rightarrow\mathbb{Z}\left[  x\right]  /\left(  mx-1\right)
,\\
\dfrac{a}{m^{k}}  &  \mapsto\overline{ax^{k}}\ \ \ \ \ \ \ \ \ \ \left(
\text{where }a\in\mathbb{Z}\text{ and }k\in\mathbb{N}\right)  .
\end{align*}
(You need to show that this is well-defined -- i.e., that if an element of
$R_{m}$ has been written in the form $\dfrac{a}{m^{k}}$ in two different ways,
then the resulting residue classes $\overline{ax^{k}}$ will be equal.)

\item Show that $\beta$ is a $\mathbb{Z}$-algebra morphism. (This is an
exercise in bringing fractions to a common denominator.)

\item Show that $\beta\circ\alpha=\operatorname*{id}$. (Indeed, $\beta
\circ\alpha$ is a $\mathbb{Z}$-algebra morphism, since $\beta$ and $\alpha$
are $\mathbb{Z}$-algebra morphisms. Moreover, it is easy to show that $\left(
\beta\circ\alpha\right)  \left(  \overline{x}\right)  =\overline{x}$. Hence,
$\left(  \beta\circ\alpha\right)  \left(  \sum_{i=0}^{n}c_{i}\overline{x}%
^{i}\right)  =\sum_{i=0}^{n}c_{i}\overline{x}^{i}$ for each $n\in\mathbb{N}$
and any coefficients $c_{0},c_{1},\ldots,c_{n}\in\mathbb{Z}$ (since
$\beta\circ\alpha$ is a $\mathbb{Z}$-algebra morphism). But this is saying
that $\beta\circ\alpha=\operatorname*{id}$, since every element of
$\mathbb{Z}\left[  x\right]  /\left(  mx-1\right)  $ can be written as
$\sum_{i=0}^{n}c_{i}\overline{x}^{i}$ for some $n\in\mathbb{N}$ and some
coefficients $c_{0},c_{1},\ldots,c_{n}\in\mathbb{Z}$.)

\item Show that $\alpha\circ\beta=\operatorname*{id}$. (Indeed, if you have
done Step 2, then this follows from $\beta\circ\alpha=\operatorname*{id}$.
Otherwise, show it directly.)

\item Conclude from Steps 6 and 7 that the maps $\alpha$ and $\beta$ are
mutually inverse, and thus $\alpha$ is invertible. Since $\alpha$ is a
$\mathbb{Z}$-algebra morphism, this entails that $\alpha$ is a $\mathbb{Z}%
$-algebra isomorphism, and you are done. \qedhere

\end{enumerate}
\end{proof}

\subsubsection{The general construction}

In the previous subsection, we have seen a few examples of the construction in
which we start with a commutative ring $R$ and a polynomial $b\in R\left[
x\right]  $, and construct the quotient ring $R\left[  x\right]  /b$. To
recall, the bottom line of this construction is \textquotedblleft throw a new
root of $b$ into the ring $R$ and see what happens\textquotedblright. Often,
this produces a ring extension of $R$ -- i.e., a larger ring that contains $R$
as a subring. (For example, this happens if $R=\mathbb{R}$ and $b=x^{2}+1$;
this is how Cardano defined the complex numbers.) However, this doesn't always
go well. Sometimes, what happens instead is that the ring $R$ collapses to a
trivial ring (e.g., if $b=1$) or at least becomes smaller (e.g., we have
$\left(  \mathbb{Z}/6\right)  \left[  x\right]  /\left(  2x-1\right)
\cong\mathbb{Z}/3$). Sometimes, the ring loses some of its properties: e.g.,
if we throw a new root of $x^{2}-1$ into the field $\mathbb{Q}$, then the
resulting ring $\mathbb{Q}\left[  x\right]  /\left(  x^{2}-1\right)  $ not
only fails to be a field, but even fails to be an integral domain (indeed, we
have seen that this ring is isomorphic to $\mathbb{Q}\times\mathbb{Q}$).

Let us put these things in order. First, let us show that the residue class
$\overline{x}$ in $R\left[  x\right]  /b$ is a root of $b$, so that our
construction really creates a root of $b$:

\begin{proposition}
\label{prop.fieldext.proj-morphism}Let $b\in R\left[  x\right]  $ be a
polynomial. (Recall that $R$ is still a fixed commutative ring.)

\begin{enumerate}
\item[\textbf{(a)}] The projection map%
\begin{align*}
\pi_{b}:R\left[  x\right]   &  \rightarrow R\left[  x\right]  /b,\\
p  &  \mapsto\overline{p}%
\end{align*}
is an $R\left[  x\right]  $-algebra morphism, and thus an $R$-algebra morphism.

\item[\textbf{(b)}] The map\footnotemark%
\begin{align*}
R  &  \rightarrow R\left[  x\right]  /b,\\
r  &  \mapsto\overline{r}%
\end{align*}
is an $R$-algebra morphism.

\item[\textbf{(c)}] For any $p\in R\left[  x\right]  $, we have $p\left[
\overline{x}\right]  =\overline{p}$ in $R\left[  x\right]  /b$.

\item[\textbf{(d)}] The element $\overline{x}\in R\left[  x\right]  /b$ is a
root of $b$.
\end{enumerate}
\end{proposition}

\footnotetext{Note the difference between the maps in part \textbf{(a)} and in
part \textbf{(b)}: The map in part \textbf{(a)} takes as input a polynomial
$p\in R\left[  x\right]  $, whereas the map in part \textbf{(b)} takes as
input a scalar $r\in R$ (and treats it as a constant polynomial, i.e., as
$rx^{0}\in R\left[  x\right]  $). If you regard $R$ as a subring of $R\left[
x\right]  $, you can thus view the map in part \textbf{(b)} as a restriction
of the map in part \textbf{(a)}.}

\begin{proof}
[Proof of Proposition \ref{prop.fieldext.proj-morphism}.]\textbf{(a)} This
follows from the general fact (Theorem \ref{thm.algs.quotients} \textbf{(c)})
that the canonical projection from an $R$-algebra to its quotient is an
$R$-algebra morphism. Note that we need to apply this fact to $R\left[
x\right]  $ instead of $R$ here, in order to conclude that the map in question
is an $R\left[  x\right]  $-algebra morphism. \medskip

\textbf{(b)} The map%
\begin{align*}
R  &  \rightarrow R\left[  x\right]  /b,\\
r  &  \mapsto\overline{r}%
\end{align*}
is the composition of the projection map $\pi_{b}$ from part \textbf{(a)} with
the inclusion map%
\begin{align*}
R  &  \rightarrow R\left[  x\right]  ,\\
r  &  \mapsto r=rx^{0}.
\end{align*}
Thus, it is a composition of two $R$-algebra morphisms (since both $\pi_{b}$
and the inclusion map are $R$-algebra morphisms). Hence, it is an $R$-algebra
morphism itself\footnote{Indeed, there is an easy fact (which we never stated,
but which is completely straightforward to prove after what we have seen) that
any composition of two $R$-algebra morphisms is itself an $R$-algebra
morphism.}. This proves Proposition \ref{prop.fieldext.proj-morphism}
\textbf{(b)}. \medskip

\textbf{(c)} Here is an abstract argument: Let $p\in R\left[  x\right]  $. The
projection map $\pi_{b}$ from Proposition \ref{prop.fieldext.proj-morphism}
\textbf{(a)} is an $R$-algebra morphism (by Proposition
\ref{prop.fieldext.proj-morphism} \textbf{(a)}). Hence, Proposition
\ref{prop.algmor.poly-comm} (applied to $A=R\left[  x\right]  $ and
$B=R\left[  x\right]  /b$ and $a=x$ and $f=\pi_{b}$) yields
\begin{equation}
\pi_{b}\left(  p\left[  x\right]  \right)  =p\left[  \pi_{b}\left(  x\right)
\right]  . \label{pf.prop.fieldext.proj-morphism.c.1}%
\end{equation}
However, the definition of $\pi_{b}$ yields $\pi_{b}\left(  p\left[  x\right]
\right)  =\overline{p\left[  x\right]  }=\overline{p}$ (since $p\left[
x\right]  =p$) and $\pi_{b}\left(  x\right)  =\overline{x}$. Hence,
(\ref{pf.prop.fieldext.proj-morphism.c.1}) rewrites as $\overline{p}=p\left[
\overline{x}\right]  $. This proves Proposition
\ref{prop.fieldext.proj-morphism} \textbf{(c)}.

Alternatively, you can prove it directly by writing $p$ as $p=\sum_{i=0}%
^{n}p_{i}x^{i}$ with $p_{i}\in R$. (Indeed, if you do this, then the claim
rewrites as $\sum_{i=0}^{n}p_{i}\overline{x}^{i}=\overline{\sum_{i=0}^{n}%
p_{i}x^{i}}$; but this is an easy consequence of how the quotient $R\left[
x\right]  /b$ was defined.) \medskip

\textbf{(d)} Proposition \ref{prop.fieldext.proj-morphism} \textbf{(c)}
(applied to $p=b$) yields $b\left[  \overline{x}\right]  =\overline
{b}=\overline{0}$ (since $b\in bR\left[  x\right]  $). In other words,
$\overline{x}$ is a root of $b$. This proves Proposition
\ref{prop.fieldext.proj-morphism} \textbf{(d)}.
\end{proof}

Next, for a large class of polynomials $b\in R\left[  x\right]  $ (including
the monic ones, and all the nonzero polynomials over a field), we are going to
show what $R\left[  x\right]  /b$ looks like as an $R$-module:

\begin{theorem}
\label{thm.fieldext.basis-monic}Let $m\in\mathbb{N}$. Let $b\in R\left[
x\right]  $ be a polynomial of degree $m$ such that its leading coefficient
$\left[  x^{m}\right]  b$ is a unit of $R$. Then:

\begin{enumerate}
\item[\textbf{(a)}] Each element of $R\left[  x\right]  /b$ can be uniquely
written in the form%
\[
a_{0}\overline{x^{0}}+a_{1}\overline{x^{1}}+\cdots+a_{m-1}\overline{x^{m-1}%
}\ \ \ \ \ \ \ \ \ \ \text{with }a_{0},a_{1},\ldots,a_{m-1}\in R.
\]

\item[\textbf{(b)}] The $m$ vectors $\overline{x^{0}},\overline{x^{1}}%
,\ldots,\overline{x^{m-1}}$ form a basis of the $R$-module $R\left[  x\right]
/b$. Thus, this $R$-module $R\left[  x\right]  /b$ is free of rank $m=\deg b$.

\item[\textbf{(c)}] Assume that $m>0$. Then, the $R$-algebra
morphism\footnotemark%
\begin{align*}
R  &  \rightarrow R\left[  x\right]  /b,\\
r  &  \mapsto\overline{r}%
\end{align*}
is injective. Therefore, $R$ can be viewed as an $R$-subalgebra (thus a
subring) of $R\left[  x\right]  /b$ if we identify each $r\in R$ with its
image $\overline{r}\in R\left[  x\right]  /b$.

\item[\textbf{(d)}] In particular, under the assumption that $m>0$, there
exists a commutative ring that contains $R$ as a subring and that contains a
root of $b$.
\end{enumerate}
\end{theorem}

\footnotetext{This is the map from Proposition
\ref{prop.fieldext.proj-morphism} \textbf{(b)}.}

\begin{proof}
\textbf{(a)} Let $\alpha\in R\left[  x\right]  /b$. Then, $\alpha=\overline
{a}$ for some polynomial $a\in R\left[  x\right]  $. Consider this $a$. The
division-with-remainder theorem for polynomials (Theorem
\ref{thm.polring.univar-quorem} \textbf{(a)}) tells us that there is a unique
pair $\left(  q,r\right)  $ of polynomials in $R\left[  x\right]  $ such that%
\[
a=qb+r\ \ \ \ \ \ \ \ \ \ \text{and}\ \ \ \ \ \ \ \ \ \ \deg r<\deg b.
\]
Consider this pair $\left(  q,r\right)  $. Then, in $R\left[  x\right]  /b$,
we have $\overline{a}=\overline{r}$ (since $a=qb+r$ entails $a-r=qb=bq\in
bR\left[  x\right]  $).

We have $\deg r<\deg b=m$; thus, we can write $r$ in the form $r=r_{0}%
x^{0}+r_{1}x^{1}+\cdots+r_{m-1}x^{m-1}$ for some $r_{0},r_{1},\ldots
,r_{m-1}\in R$. Consider these $r_{0},r_{1},\ldots,r_{m-1}$. We have%
\begin{align*}
\alpha &  =\overline{a}=\overline{r}=\overline{r_{0}x^{0}+r_{1}x^{1}%
+\cdots+r_{m-1}x^{m-1}}\\
&  \ \ \ \ \ \ \ \ \ \ \ \ \ \ \ \ \ \ \ \ \left(  \text{since }r=r_{0}%
x^{0}+r_{1}x^{1}+\cdots+r_{m-1}x^{m-1}\right) \\
&  =r_{0}\overline{x^{0}}+r_{1}\overline{x^{1}}+\cdots+r_{m-1}\overline
{x^{m-1}}%
\end{align*}
(since the scaling and the addition of the quotient algebra $R\left[
x\right]  /b$ were inherited from $R\left[  x\right]  $).

Thus, we have represented our $\alpha\in R\left[  x\right]  /b$ in the form%
\[
a_{0}\overline{x^{0}}+a_{1}\overline{x^{1}}+\cdots+a_{m-1}\overline{x^{m-1}%
}\ \ \ \ \ \ \ \ \ \ \text{with }a_{0},a_{1},\ldots,a_{m-1}\in R
\]
(namely, for $a_{i}=r_{i}$). It remains to show that this representation is unique.

This can be shown by walking the above proof backwards and using the
uniqueness part of the division-with-remainder theorem. Here are the details:
Assume that%
\[
\alpha=b_{0}\overline{x^{0}}+b_{1}\overline{x^{1}}+\cdots+b_{m-1}%
\overline{x^{m-1}}\ \ \ \ \ \ \ \ \ \ \text{with }b_{0},b_{1},\ldots
,b_{m-1}\in R
\]
is some representation of $\alpha$ in the above form. We must then show that
this representation is actually the representation that we constructed above
-- i.e., that we have $b_{i}=r_{i}$ for each $i\in\left\{  0,1,\ldots
,m-1\right\}  $.

Indeed, define a polynomial $s\in R\left[  x\right]  $ by $s=b_{0}x^{0}%
+b_{1}x^{1}+\cdots+b_{m-1}x^{m-1}$. Then, $\deg s\leq m-1<m=\deg b$. Also,%
\[
\overline{a}=\alpha=b_{0}\overline{x^{0}}+b_{1}\overline{x^{1}}+\cdots
+b_{m-1}\overline{x^{m-1}}=\overline{b_{0}x^{0}+b_{1}x^{1}+\cdots
+b_{m-1}x^{m-1}}=\overline{s}%
\]
(since $b_{0}x^{0}+b_{1}x^{1}+\cdots+b_{m-1}x^{m-1}=s$). In other words,
$a-s\in bR\left[  x\right]  $. In other words,%
\[
a-s=bd\ \ \ \ \ \ \ \ \ \ \text{for some }d\in R\left[  x\right]  .
\]
Consider this $d$. Thus, $a=bd+s=db+s$. Now, the pair $\left(  d,s\right)  $
is a pair of polynomials in $R\left[  x\right]  $ satisfying $a=db+s$ and
$\deg s<\deg b$. This means that it satisfies the exact conditions that the
pair $\left(  q,r\right)  $ was asked to satisfy. However, the
division-with-remainder theorem for polynomials said that the pair $\left(
q,r\right)  $ satisfying those conditions was unique. Hence, we must have
$\left(  d,s\right)  =\left(  q,r\right)  $ (since $\left(  d,s\right)  $
satisfies the same conditions as $\left(  q,r\right)  $). Thus, $d=q$ and
$s=r$.

Now,%
\[
b_{0}x^{0}+b_{1}x^{1}+\cdots+b_{m-1}x^{m-1}=s=r=r_{0}x^{0}+r_{1}x^{1}%
+\cdots+r_{m-1}x^{m-1}.
\]
Comparing coefficients in these polynomials, we conclude that $b_{i}=r_{i}$
for each $i\in\left\{  0,1,\ldots,m-1\right\}  $ (since $\left(  x^{0}%
,x^{1},x^{2},\ldots\right)  $ is a basis of the $R$-module $R\left[  x\right]
$). This is what we needed to show. Theorem \ref{thm.fieldext.basis-monic}
\textbf{(a)} is thus proved. \medskip

\textbf{(b)} This is just Theorem \ref{thm.fieldext.basis-monic} \textbf{(a)},
rewritten in terms of modules and bases.

In some more detail:

\begin{itemize}
\item Each element of $R\left[  x\right]  /b$ can be written in the form%
\[
a_{0}\overline{x^{0}}+a_{1}\overline{x^{1}}+\cdots+a_{m-1}\overline{x^{m-1}%
}\ \ \ \ \ \ \ \ \ \ \text{with }a_{0},a_{1},\ldots,a_{m-1}\in R
\]
(according to Theorem \ref{thm.fieldext.basis-monic} \textbf{(a)}). In other
words, each element of $R\left[  x\right]  /b$ is an $R$-linear combination of
$\overline{x^{0}},\overline{x^{1}},\ldots,\overline{x^{m-1}}$. Thus, the list
$\left(  \overline{x^{0}},\overline{x^{1}},\ldots,\overline{x^{m-1}}\right)  $
spans the $R$-module $R\left[  x\right]  /b$.

\item Each element of $R\left[  x\right]  /b$ can be \textbf{uniquely}
represented in the form%
\[
a_{0}\overline{x^{0}}+a_{1}\overline{x^{1}}+\cdots+a_{m-1}\overline{x^{m-1}%
}\ \ \ \ \ \ \ \ \ \ \text{with }a_{0},a_{1},\ldots,a_{m-1}\in R
\]
(according to Theorem \ref{thm.fieldext.basis-monic} \textbf{(a)}). Hence, in
particular, the zero vector $\overline{0}\in R\left[  x\right]  /b$ can be
\textbf{uniquely} represented in this form. But it is clear how to represent
$\overline{0}$ in this form: We just write%
\[
\overline{0}=0\overline{x^{0}}+0\overline{x^{1}}+\cdots+0\overline{x^{m-1}}.
\]
Since we have just said that $\overline{0}$ can be \textbf{uniquely}
represented in this form, we thus conclude that this is the \textbf{only} way
to represent $\overline{0}$ in this form. In other words, if $\overline{0}$
has been represented in the form $a_{0}\overline{x^{0}}+a_{1}\overline{x^{1}%
}+\cdots+a_{m-1}\overline{x^{m-1}}$ with $a_{0},a_{1},\ldots,a_{m-1}\in R$,
then we must have $a_{0}=a_{1}=\cdots=a_{m-1}=0$. In other words, if
$a_{0},a_{1},\ldots,a_{m-1}\in R$ satisfy $a_{0}\overline{x^{0}}%
+a_{1}\overline{x^{1}}+\cdots+a_{m-1}\overline{x^{m-1}}=\overline{0}$, then
$a_{0}=a_{1}=\cdots=a_{m-1}=0$. But this is saying precisely that the list
$\left(  \overline{x^{0}},\overline{x^{1}},\ldots,\overline{x^{m-1}}\right)  $
is $R$-linearly independent.
\end{itemize}

Thus, we have shown that the list $\left(  \overline{x^{0}},\overline{x^{1}%
},\ldots,\overline{x^{m-1}}\right)  $ is $R$-linearly independent and spans
$R\left[  x\right]  /b$. In other words, this list is a basis of $R\left[
x\right]  /b$. This proves Theorem \ref{thm.fieldext.basis-monic}
\textbf{(b)}. \medskip

\textbf{(c)} We know (from Proposition \ref{prop.fieldext.proj-morphism}
\textbf{(b)}) that the map%
\begin{align*}
R  &  \rightarrow R\left[  x\right]  /b,\\
r  &  \mapsto\overline{r}%
\end{align*}
is an $R$-algebra morphism. We only need to show that it is injective. It
clearly suffices to show that its kernel is $\left\{  0\right\}  $ (because we
know that an $R$-module morphism is injective if and only if its kernel is
$\left\{  0\right\}  $).

So let $r$ be in the kernel of this morphism. We must prove that $r=0$.

Since $r$ is in the kernel of the above morphism, we have $\overline{r}=0$ in
$R\left[  x\right]  /b$. In other words, $r$ is a multiple of $b$. In other
words, $r=bc$ for some polynomial $c\in R\left[  x\right]  $. Consider this
$c$. From $r=bc$, we obtain $\deg r=\deg\left(  bc\right)  =\deg b+\deg c$ (by
Proposition \ref{prop.polring.univar-degpq} \textbf{(b)}, since the leading
coefficient of $b$ is a unit). Thus, $\deg b+\deg c=\deg r\leq0$ (since $r$ is
constant). However, $\deg b=m>0$ by assumption. Hence, $\deg b>0\geq\deg
b+\deg c$. This entails $\deg c<0$. This means that $c=0$, whence
$r=b\underbrace{c}_{=0}=0$.

Forget that we fixed $r$. We thus have proved that if $r$ is in the kernel of
our morphism, then $r=0$. Hence, the kernel of our morphism is $\left\{
0\right\}  $ (since $0$ is clearly in its kernel). Thus, the morphism is
injective, and Theorem \ref{thm.fieldext.basis-monic} \textbf{(c)} is proven.
\medskip

\textbf{(d)} Assume that $m>0$. The ring $R\left[  x\right]  /b$ contains a
root of $b$ (namely, $\overline{x}$, according to Proposition
\ref{prop.fieldext.proj-morphism} \textbf{(d)}), and also contains
\textquotedblleft a copy of $R$\textquotedblright, in the sense that there is
an injective ring morphism from $R$ to $R\left[  x\right]  /b$ (namely, the
one we constructed in Theorem \ref{thm.fieldext.basis-monic} \textbf{(c)}). If
we replace this copy of $R$ by the original $R$ (by replacing each
$\overline{r}\in R\left[  x\right]  /b$ with the corresponding $r\in R$), then
we obtain a ring that contains $R$ as a subring but also contains a root of
$b$. This proves Theorem \ref{thm.fieldext.basis-monic} \textbf{(d)}.
\end{proof}

Let us summarize: We have generalized the construction of $\mathbb{C}$.
Namely, we have found a way to \textquotedblleft adjoin\textquotedblright\ a
root of a polynomial $b\in R\left[  x\right]  $ to a commutative ring $R$ by
forming the quotient ring $R\left[  x\right]  /b$. This latter ring is always
a commutative ring and an $R$-algebra. Moreover, if $b$ is \textquotedblleft
nice\textquotedblright\ (that is, we have $\deg b>0$, and the leading
coefficient of $b$ is a unit), then this latter ring $R\left[  x\right]  /b$
will contain $R$ as a subring (by Theorem \ref{thm.fieldext.basis-monic}
\textbf{(c)}) and also will be a free $R$-module of rank $\deg b$ (by Theorem
\ref{thm.fieldext.basis-monic} \textbf{(b)}). If $b$ is not as
\textquotedblleft nice\textquotedblright, then the ring $R\left[  x\right]
/b$ may fail to contain $R$ as a subring (even though it still is an
$R$-algebra), and may be smaller than $R$ or even trivial.

\subsection{\label{sec.polys1.fieldexts}Field extensions from adjoining roots}

Let $F$ be a field. Then, any non-constant univariate polynomial $b\in
F\left[  x\right]  $ is \textquotedblleft nice\textquotedblright\ in the sense
of the preceding paragraph, so that $F\left[  x\right]  /b$ is a commutative
ring that contains $F$ as a subring and that contains a root of $b$. When will
this ring $F\left[  x\right]  /b$ be a field?

We first state a simple fact about the units of $F\left[  x\right]  $:

\begin{proposition}
\label{prop.fieldext.pol-units}Let $F$ be a field. The units of the polynomial
ring $F\left[  x\right]  $ are precisely the nonzero constant polynomials.
\end{proposition}

\begin{proof}
Any nonzero constant polynomial is a unit of $F\left[  x\right]  $ (since it
is a unit of $F$). Conversely, any unit of $F\left[  x\right]  $ must be a
nonzero constant polynomial\footnote{\textit{Proof.} Let $u$ be a unit of
$F\left[  x\right]  $. We must show that $u$ is a nonzero constant polynomial.
\par
We know that $u$ is a unit of $F\left[  x\right]  $; hence, there exists some
$v\in F\left[  x\right]  $ satisfying $uv=1$. Consider this $v$. From
$uv=1\neq0$, we obtain $u\neq0$, so that $u$ is nonzero. Hence, $\deg\left(
uv\right)  =\deg u+\deg v$ (by Proposition \ref{prop.polring.univar-degpq}
\textbf{(c)}, since $F$ is an integral domain). Moreover, from $uv=1$, we
obtain $\deg\left(  uv\right)  =\deg1=0$, so that $0=\deg\left(  uv\right)
=\deg u+\underbrace{\deg v}_{\geq0}\geq\deg u$, which entails that $u$ is
constant. Thus, $u$ is a nonzero constant polynomial, qed.}.
\end{proof}

Recall (from Theorem \ref{thm.polring.univar-Euc}) that $F\left[  x\right]  $
is a Euclidean domain, hence a PID (by Proposition \ref{prop.eucldom.PID2}),
hence a UFD (by Theorem \ref{thm.UFD.PID-is-UFD}). Furthermore, an element
$p\in F\left[  x\right]  $ is prime\footnote{See Definition
\ref{def.rings.prime-irred} for the definitions of prime and irreducible
elements of an integral domain.} if and only if it is irreducible (by
Proposition \ref{prop.PID.prime-iff-irred}, since $F\left[  x\right]  $ is a
PID). The notion of \textquotedblleft irreducible\textquotedblright\ in
$F\left[  x\right]  $\ is precisely the classical concept of an irreducible polynomial:

\begin{proposition}
\label{prop.fieldext.pol-irr}Let $F$ be a field. Let $p\in F\left[  x\right]
$. Then, $p$ is irreducible if and only if $p$ is non-constant and cannot be
written as a product of two non-constant polynomials.
\end{proposition}

\begin{proof}
The definition of \textquotedblleft irreducible\textquotedblright\ says that
$p$ is irreducible if and only if $p$ is nonzero and not a unit and has the
property that whenever $a,b\in F\left[  x\right]  $ satisfy $ab=p$, at least
one of $a$ and $b$ must be a unit.

In view of Proposition \ref{prop.fieldext.pol-units}, this can be rewritten as
follows: $p$ is irreducible if and only if $p$ is nonzero and not a nonzero
constant polynomial and has the property that whenever $a,b\in F\left[
x\right]  $ satisfy $ab=p$, at least one of $a$ and $b$ must be a nonzero
constant polynomial.

We can declutter this statement (e.g., \textquotedblleft nonzero and not a
nonzero constant polynomial\textquotedblright\ can be shortened to
\textquotedblleft non-constant\textquotedblright), and thus obtain the
following: $p$ is irreducible if and only if $p$ is non-constant and has the
property that whenever $a,b\in F\left[  x\right]  $ satisfy $ab=p$, at least
one of $a$ and $b$ must be constant. In other words, $p$ is irreducible if and
only if $p$ is non-constant and cannot be written as a product of two
non-constant polynomials.
\end{proof}

Now, we can characterize when a quotient ring of the form $F\left[  x\right]
/p$ is a field:

\begin{theorem}
\label{thm.fieldext.field-irr}Let $F$ be a field. Let $p\in F\left[  x\right]
$ be nonzero. Then, the ring $F\left[  x\right]  /p$ is a field if and only if
$p$ is irreducible.
\end{theorem}

For example, the irreducible polynomial $x^{2}+1$ over the field $\mathbb{R}$
yields the field $\mathbb{R}\left[  x\right]  /\left(  x^{2}+1\right)  $
(which is $\cong\mathbb{C}$), but the non-irreducible polynomial $x^{2}-1$
over the field $\mathbb{R}$ yields the non-field $\mathbb{R}\left[  x\right]
/\left(  x^{2}-1\right)  \cong\mathbb{R}\times\mathbb{R}$.

Theorem \ref{thm.fieldext.field-irr} is analogous to the fact that
$\mathbb{Z}/n$ is a field (for a positive integer $n$) if and only if $n$ is
prime. Just like the latter fact, it is a particular case of the following
general property of PIDs:

\begin{theorem}
\label{thm.PID.quot-field}Let $R$ be a PID. Let $p\in R$ be nonzero. Then, the
ring $R/p$ is a field if and only if $p$ is irreducible.
\end{theorem}

\begin{proof}
$\Longrightarrow:$ LTTR.

$\Longleftarrow:$ Assume that $p$ is irreducible. We must show that $R/p$ is a field.

First of all, $p$ is not a unit (since $p$ is irreducible), so that $1$ is not
a multiple of $p$. Hence, $\overline{1}\neq\overline{0}$ in $R/p$. In other
words, the ring $R/p$ is not trivial. This ring is furthermore commutative
(since $R$ is commutative).

Now, let $\alpha\in R/p$ be a nonzero element. We shall prove that $\alpha$ is
a unit.

Write $\alpha$ as $\overline{a}$ for some $a\in R$. Then, $\overline{a}%
=\alpha\neq0$ in $R/p$ (since $\alpha$ is nonzero), so that $p\nmid a$.

Now, recall that $R$ is a PID, so that any ideal of $R$ is principal. In
particular, this entails that the ideal $aR+pR$ is principal. In other words,
there exists some $g\in R$ such that $aR+pR=gR$. Consider this $g$. According
to Proposition \ref{prop.gcd-lcm.ideals} \textbf{(a)}, we can conclude from
$aR+pR=gR$ that $g$ is a gcd of $a$ and $p$. Thus, $g\mid a$ and $g\mid p$.

However, $p$ is irreducible; hence, every divisor of $p$ is either a unit or
associate to $p$ (indeed, this is easily seen to be a consequence of the
definition of \textquotedblleft irreducible\textquotedblright\footnote{Indeed:
If $d$ is a divisor of $p$, then there exists an $e\in R$ such that $p=de$.
Consider this $e$. From $p=de$, we conclude that $d$ or $e$ is a unit (since
$p$ is irreducible). In the first case, $d$ is a unit; in the second case, $d$
is associate to $p$.}). Thus, $g$ is either a unit or associate to $p$ (since
$g\mid p$). However, $g$ cannot be associate to $p$ (because if $g$ was
associate to $p$, then we would have $p\mid g\mid a$, which would contradict
$p\nmid a$). Hence, $g$ must be a unit. Thus, it has an inverse $g^{-1}$.

But $g=g\cdot1\in gR=aR+pR$. In other words, there exist two elements $u,v\in
R$ such that $g=au+pv$. Consider these $u,v$. Then, $g=au+pv=ua+pv$, so that
\[
\overline{g}=\overline{ua+pv}=\overline{ua}\ \ \ \ \ \ \ \ \ \ \left(
\text{since }pv\in pR\right)
\]
in $R/p$. Therefore,%
\[
\overline{g^{-1}u}\cdot\overline{a}=\overline{g^{-1}}\cdot\overline{u}%
\cdot\overline{a}=\overline{g^{-1}}\cdot\underbrace{\overline{ua}}%
_{=\overline{g}}=\overline{g^{-1}}\cdot\overline{g}=\overline{g^{-1}%
g}=\overline{1}.
\]
But this equality shows that $\overline{g^{-1}u}$ is an inverse of
$\overline{a}$ in the ring $R/p$ (because we know that $R/p$ is commutative,
so that we don't need to check $\overline{a}\cdot\overline{g^{-1}u}%
=\overline{1}$ as well). Thus, $\overline{a}$ is a unit. In other words,
$\alpha$ is a unit (since $\alpha=\overline{a}$).

Forget that we fixed $\alpha$. We thus have shown that any nonzero $\alpha\in
R/p$ is a unit. In other words, $R/p$ is a field (since $R/p$ is a nontrivial
commutative ring).
\end{proof}

Theorem \ref{thm.fieldext.field-irr} is a particular case of Theorem
\ref{thm.PID.quot-field} (since $F\left[  x\right]  $ is a PID when $F$ is a field).

As a consequence of Theorem \ref{thm.fieldext.field-irr}, we can now
\textquotedblleft adjoin\textquotedblright\ a root of an irreducible
polynomial to a field without destroying its field-ness: Namely, if we have a
field $F$ and some irreducible polynomial $b\in F\left[  x\right]  $, then the
quotient ring $F\left[  x\right]  /b$ will be a field that contains $F$ as a
subring and that contains a root of $b$. This generalizes Cardano's definition
of $\mathbb{C}$, but can also be applied to adjoin roots to fields other than
$\mathbb{R}$.

\begin{noncompile}
(OLD EXAMPLE, DOESN'T FIT, WRONG NOTATIONS)

Consider the ring
\[
\mathbb{Q}\left[  \sqrt[3]{2}\right]  :=\left\{  a+b\sqrt[3]{2}+c\sqrt[3]%
{4}\ \mid\ a,b,c\in\mathbb{Q}\right\}  .
\]
Then,
\[
\mathbb{Q}\left[  \sqrt[3]{2}\right]  \cong\mathbb{Q}\left[  x\right]
/\left(  x^{3}-2\right)  .
\]
More concretely, there is a $\mathbb{Q}$-algebra isomorphism%
\begin{align*}
\mathbb{Q}\left[  x\right]  /\left(  x^{3}-2\right)   &  \rightarrow
\mathbb{Q}\left[  \sqrt[3]{2}\right]  ,\\
\overline{p}  &  \mapsto p\left(  \sqrt[3]{2}\right)  .
\end{align*}
Since the polynomial $x^{3}-2\in\mathbb{Q}\left[  x\right]  $ is irreducible
(this is not hard to check), Theorem \ref{thm.fieldext.field-irr} predicts
that $\mathbb{Q}\left[  x\right]  /\left(  x^{3}-2\right)  $ is a field. Thus,
$\mathbb{Q}\left[  \sqrt[3]{2}\right]  $ is a field. Hence, any nonzero
element $a+b\sqrt[3]{2}+c\sqrt[3]{4}$ of $\mathbb{Q}\left[  \sqrt[3]%
{2}\right]  $ has an inverse $a^{\prime}+b^{\prime}\sqrt[3]{2}+c^{\prime
}\sqrt[3]{4}$. Let us compute it on an example: Let's compute the inverse of
$1+\sqrt[3]{2}+\sqrt[3]{4}$ (so $a=b=c=1$).

Let's follow the proof of the above theorem. We want to find an inverse for
$1+\sqrt[3]{2}+\sqrt[3]{4}$. In other words, we want to find an inverse for%
\[
\overline{1+x+x^{2}}\in\mathbb{Q}\left[  x\right]  /\left(  x^{3}-2\right)  .
\]
To do so, we apply the Euclidean algorithm to compute the gcd of $x^{3}-2$ and
$1+x+x^{2}$:

Dividing $x^{3}-2$ by $1+x+x^{2}$ with remainder, we find%
\[
x^{3}-2=\left(  x-1\right)  \left(  1+x+x^{2}\right)  -1.
\]

This rewrites as
\[
\underbrace{\left(  -1\right)  }_{=u}\cdot\underbrace{\left(  x^{3}-2\right)
}_{=p}+\underbrace{\left(  x-1\right)  }_{=v}\cdot\underbrace{\left(
1+x+x^{2}\right)  }_{=a}=\underbrace{1}_{=\gcd\left(  p,a\right)  }.
\]
Recall that $\overline{\left(  \gcd\left(  p,a\right)  \right)  ^{-1}}%
\cdot\overline{v}$ is an inverse to $\overline{a}$. Thus, in our case, the
inverse is%
\[
\overline{1^{-1}}\cdot\overline{x-1}=\overline{x-1}.
\]
Sending it back through the isomorphism, we get $\sqrt[3]{2}-1$.

So the inverse of $1+\sqrt[3]{2}+\sqrt[3]{4}$ is $\sqrt[3]{2}-1$. And indeed,%
\[
\left(  1+\sqrt[3]{2}+\sqrt[3]{4}\right)  \left(  \sqrt[3]{2}-1\right)  =1.
\]

\end{noncompile}

\begin{example}
\label{exa.fieldext.F9}The polynomial $x^{2}+1\in\left(  \mathbb{Z}/3\right)
\left[  x\right]  $ is irreducible. (Indeed, $\mathbb{Z}/3$ being a finite
field, we could verify this by going through all nonconstant polynomials of
degree $<2$ and checking that none of them divides $x^{2}+1$.)

Thus, Theorem \ref{thm.fieldext.field-irr} yields that $\left(  \mathbb{Z}%
/3\right)  \left[  x\right]  /\left(  x^{2}+1\right)  $ is a field. This field
is a free $\mathbb{Z}/3$-module of rank $2$ (by Theorem
\ref{thm.fieldext.basis-monic} \textbf{(b)}), and thus is isomorphic to
$\left(  \mathbb{Z}/3\right)  ^{2}=\left(  \mathbb{Z}/3\right)  \times\left(
\mathbb{Z}/3\right)  $ as a $\mathbb{Z}/3$-module (but not as a ring, of
course). Hence, the size of this field is $\left\vert \left(  \mathbb{Z}%
/3\right)  ^{2}\right\vert =\left\vert \mathbb{Z}/3\right\vert ^{2}=3^{2}=9$.

Thus, we have found a finite field of size $9$. We have obtained it from
$\mathbb{Z}/3$ in the same way as $\mathbb{C}$ was obtained from $\mathbb{R}$:
by adjoining a square root of $-1$.

Incidentally, this field can also be constructed as $\mathbb{Z}\left[
i\right]  /3$.
\end{example}

\newpage

\section{\label{chp.finfields}Finite fields}

\subsection{\label{sec.finfields.basics}Basics}

Example \ref{exa.fieldext.F9} may make you wonder: what finite fields can we
find? We know that for each prime $p$, the quotient ring $\mathbb{Z}/p$ is a
field of size $p$; thus, we know a finite field of any prime size. Now we have
found a finite field of size $9$, too. What other finite fields exist?

Let's first grab the low-hanging fruit:

\begin{proposition}
\label{prop.finfield.p^2}Let $p$ be a prime number. Then:

\begin{enumerate}
\item[\textbf{(a)}] There exists an irreducible polynomial $b\in\left(
\mathbb{Z}/p\right)  \left[  x\right]  $ of degree $2$ over $\mathbb{Z}/p$.

\item[\textbf{(b)}] There exists a finite field of size $p^{2}$.
\end{enumerate}
\end{proposition}

\begin{example}
\ \ 

\begin{enumerate}
\item[\textbf{(a)}] Let $p=2$. Then, Proposition \ref{prop.finfield.p^2}
\textbf{(a)} yields that there exists an irreducible polynomial $b\in\left(
\mathbb{Z}/2\right)  \left[  x\right]  $ of degree $2$ over $\mathbb{Z}/2$. It
is easy to see that this polynomial $b$ is actually unique; it is
$b=x^{2}+x+1$. Furthermore, Proposition \ref{prop.finfield.p^2} \textbf{(b)}
yields that there exists a finite field of size $p^{2}=2^{2}=4$. As we will
see in the proof of the proposition, this latter field can be obtained as
$\left(  \mathbb{Z}/2\right)  \left[  x\right]  /b$ for the afore-mentioned
polynomial $b$. It has four elements $\overline{0},\overline{1},\overline
{x},\overline{x+1}$. It is easy to see that this field is actually isomorphic
to the four-element ring $F_{4}$ constructed in Subsection
\ref{subsec.rings.def.exas}. (One possible isomorphism sends its four elements
$\overline{0},\overline{1},\overline{x},\overline{x+1}$ to the four elements
$0,1,a,b$ of the latter.) This explains why the $F_{4}$ from Subsection
\ref{subsec.rings.def.exas} is a ring (and even a field)!

\item[\textbf{(b)}] Let us now take $p=3$ instead. Then, Proposition
\ref{prop.finfield.p^2} \textbf{(a)} yields that there exists an irreducible
polynomial $b\in\left(  \mathbb{Z}/3\right)  \left[  x\right]  $ of degree $2$
over $\mathbb{Z}/3$. Actually, there are three such polynomials that are
monic:
\[
b_{1}=x^{2}+1;\ \ \ \ \ \ \ \ \ \ b_{2}=x^{2}+x+2;\ \ \ \ \ \ \ \ \ \ b_{3}%
=x^{2}+2x+2.
\]
We can use them to construct a finite field of size $p^{2}=3^{2}=9$. Actually,
the three fields obtained turn out to be mutually isomorphic.
\end{enumerate}
\end{example}

\begin{proof}
[Proof of Proposition \ref{prop.finfield.p^2}.]We write $F$ for $\mathbb{Z}%
/p$. Thus, $F$ is a field and satisfies $\left\vert F\right\vert =\left\vert
\mathbb{Z}/p\right\vert =p$. \medskip

\textbf{(a)} If $p=2$, then we can take $b=x^{2}+x+1$; it is easy to check
that this $b$ is irreducible.

Thus, WLOG assume that $p\neq2$. Hence, $p>2$. Thus, $\overline{1}%
\neq\overline{-1}$ in $\mathbb{Z}/p$. In other words, $\overline{1}%
\neq\overline{-1}$ in $F$ (since $\mathbb{Z}/p=F$). The map%
\begin{align*}
F  &  \rightarrow F,\\
a  &  \mapsto a^{2}%
\end{align*}
is not injective (since $\overline{1}^{2}=\overline{-1}^{2}$ but $\overline
{1}\neq\overline{-1}$), and thus cannot be surjective (by the pigeonhole
principle). Thus, there exists some $u\in F$ that is not in the image of this
map. In other words, there exists some $u\in F$ that is not a square. Consider
such a $u$. Then, the polynomial $x^{2}-u$ has no roots in $F$.

Now it is not hard to prove that the polynomial $x^{2}-u$ is
irreducible.\footnote{\textit{Proof.} Assume that we have written $x^{2}-u$ as
a product $fg$ of two non-constant polynomials $f,g\in F\left[  x\right]  $.
We shall derive a contradiction.
\par
Indeed, we have assumed that $x^{2}-u=fg$; hence, $\deg\left(  x^{2}-u\right)
=\deg\left(  fg\right)  =\deg f+\deg g$ (since $F$ is an integral domain).
Thus, $\deg f+\deg g=\deg\left(  x^{2}-u\right)  =2$. Since $\deg f$ and $\deg
g$ are positive integers (because $f$ and $g$ are non-constant), this entails
that $\deg f$ and $\deg g$ must equal $1$ (since the only pair of positive
integers that add up to $2$ is $\left(  1,1\right)  $). Thus, in particular,
$\deg f=1$. Hence, $f=ax+b$ for some $a,b\in F$ with $a\neq0$. Consider these
$a,b$. From $f=ax+b$, we obtain $f\left[  \dfrac{-b}{a}\right]  =a\cdot
\dfrac{-b}{a}+b=0$. Thus, the polynomial $f$ has a root in $F$ (namely,
$\dfrac{-b}{a}$). Hence, the polynomial $x^{2}-u$ has a root in $F$ as well
(indeed, $f\mid fg=x^{2}-u$, so that every root of $f$ is also a root of
$x^{2}-u$). This contradicts the fact that the polynomial $x^{2}-u$ has no
roots in $F$.
\par
Thus, we have found a contradiction stemming from our assumption that
$x^{2}-u$ is a product $fg$ of two non-constant polynomials $f,g\in F\left[
x\right]  $. Hence, $x^{2}-u$ cannot be written as such a product. In other
words, $x^{2}-u$ is irreducible (since $x^{2}-u$ is a non-constant
polynomial). Qed.} This proves Proposition \ref{prop.finfield.p^2}
\textbf{(a)} (since $x^{2}-u\in\left(  \mathbb{Z}/p\right)  \left[  x\right]
$ is an irreducible polynomial of degree $2$). \medskip

\textbf{(b)} Proposition \ref{prop.finfield.p^2} \textbf{(a)} yields that
there exists an irreducible polynomial $b\in F\left[  x\right]  $ of degree
$2$ over $F$ (since $\mathbb{Z}/p=F$). Consider this $b$. Theorem
\ref{thm.fieldext.field-irr} (applied to $b$ instead of $p$) then yields that
the ring $F\left[  x\right]  /b$ is a field. Moreover, $F\left[  x\right]  /b$
is a free $F$-module of rank $2$ (by Theorem \ref{thm.fieldext.basis-monic}
\textbf{(b)}), and thus is isomorphic to $F^{2}$ as a $F$-module, and
therefore has size $\left\vert F^{2}\right\vert =\left\vert F\right\vert
^{2}=p^{2}$ (since $\left\vert F\right\vert =p$). Hence, $F\left[  x\right]
/b$ is a finite field of size $p^{2}$. This proves Proposition
\ref{prop.finfield.p^2} \textbf{(b)}.
\end{proof}

By more complicated but somewhat similar arguments\footnote{Not too similar!
It is not true that the map%
\begin{align*}
F  &  \rightarrow F,\\
a  &  \mapsto a^{3}%
\end{align*}
is always non-surjective when $F=\mathbb{Z}/p$ for $p>3$. Instead, you have to
argue the existence of an irreducible polynomial $b\in\left(  \mathbb{Z}%
/p\right)  \left[  x\right]  $ of degree $3$ over $\mathbb{Z}/p$ by a counting
argument: Show that the total number of monic degree-$3$ polynomials in
$\left(  \mathbb{Z}/p\right)  \left[  x\right]  $ is $p^{3}$, whereas the
total number of monic degree-$3$ polynomials in $\left(  \mathbb{Z}/p\right)
\left[  x\right]  $ that can be written as a product of a degree-$1$ and a
degree-$2$ polynomial is smaller than $p^{3}$; thus, at least one monic
degree-$3$ polynomial cannot be written as such a product.}, we can also see
that there exists a finite field of size $p^{3}$ for any prime $p$. This
suggests generalizing to $p^{m}$; but this is much harder. Indeed, a
nonconstant polynomial over $F$ of degree $\leq3$ will always be irreducible
if it has no roots in $F$ (check this!\footnote{See Exercise
\ref{exe.polring.xxx-r} \textbf{(a)}.}); however, for polynomials of degree
$\geq4$, this is no longer the case (fun exercise: prove that the polynomial
$x^{4}+4\in\mathbb{Q}\left[  x\right]  $ is not irreducible, despite, of
course, not having any roots over $\mathbb{Q}$). Thus, our trick for finding
irreducible polynomials will no longer work for degrees $>3$. We can still
find a field of size $p^{4}$ by applying our trick twice (first get a finite
field of size $p^{2}$, then proceed to find an irreducible polynomial of
degree $2$ over that field), and by induction we can find fields of sizes
$p^{8},p^{16},p^{32},\ldots$. But we don't get a field of size $p^{5}$ this way.

So do such fields exist?

\subsection{\label{sec.finfields.char}The characteristic of a field}

\begin{noncompile}
In the previous section, we studied the possible sizes of a finite field. We
showed how to find fields of sizes $p,\ p^{2},\ p^{3},\ p^{4}$ and various
others whenever $p$ is a prime, but we could not answer the question for size
$p^{5}$ (or $p^{7}$, or many others).
\end{noncompile}

Leaving prime powers aside for a moment, what about fields of size $6$ ? It
turns out that such fields don't exist, for a fairly simple reason. Fields
have an important invariant, the so-called \textbf{characteristic}:

\begin{definition}
\label{def.finfield.char}Let $F$ be a field. The \textbf{characteristic} of
$F$ is an integer called $\operatorname*{char}F$, which is defined as follows:

\begin{itemize}
\item If there exists a positive integer $n$ such that $n\cdot1_{F}=0_{F}$,
then $\operatorname*{char}F$ is defined to be the \textbf{smallest} such $n$.

\item If such an $n$ does not exist, then $\operatorname*{char}F$ is defined
to be $0$.
\end{itemize}
\end{definition}

Roughly speaking, $\operatorname*{char}F$ is \textquotedblleft how often you
have to add $1_{F}$ to itself to obtain $0_{F}$\textquotedblright\ (with the
caveat that we define it to be $0$ if you never obtain $0_{F}$ by adding
$1_{F}$ to itself). Here are some examples:

\begin{itemize}
\item We have $\operatorname*{char}\mathbb{Q}=0$, since there exists no
positive integer $n$ such that $n\cdot1_{\mathbb{Q}}=0_{\mathbb{Q}}$. For the
same reason, $\operatorname*{char}\mathbb{R}=0$ and $\operatorname*{char}%
\mathbb{C}=0$.

\item For any prime $p$, we have $\operatorname*{char}\left(  \mathbb{Z}%
/p\right)  =p$. Indeed, $p\cdot1_{\mathbb{Z}/p}=p\cdot\overline{1}%
=\overline{p\cdot1}=\overline{p}=\overline{0}$ in $\mathbb{Z}/p$, but every
positive integer $n<p$ satisfies $n\cdot1_{\mathbb{Z}/p}=n\cdot\overline
{1}=\overline{n\cdot1}=\overline{n}\neq\overline{0}$ in $\mathbb{Z}/p$.

\item For our fields $F$ of size $p^{2}$ or $p^{3}$, we also have
$\operatorname*{char}F=p$, since they contain $\mathbb{Z}/p$ as subrings.
\end{itemize}

Let us now see what the characteristic of a field satisfies in general, and
what it can tell us about the field.

\begin{theorem}
[Properties of characteristics]\label{thm.finfield.char.basics}Let $F$ be a
field. Let $p=\operatorname*{char}F$. Then:

\begin{enumerate}
\item[\textbf{(a)}] The field $F$ is a $\mathbb{Z}/p$-algebra. (Remember:
$\mathbb{Z}/0\cong\mathbb{Z}$.)

\item[\textbf{(b)}] We have $pa=0$ for each $a\in F$.

\item[\textbf{(c)}] The number $p$ is either prime or $0$.

\item[\textbf{(d)}] If $F$ is finite, then $p$ is a prime.

\item[\textbf{(e)}] If $F$ is finite, then $\left\vert F\right\vert =p^{m}$
for some positive integer $m$.

\item[\textbf{(f)}] If $p$ is a prime, then $F$ contains \textquotedblleft a
copy of $\mathbb{Z}/p$\textquotedblright\ (meaning: a subring isomorphic to
$\mathbb{Z}/p$).

\item[\textbf{(g)}] If $p=0$, then $F$ contains \textquotedblleft a copy of
$\mathbb{Q}$\textquotedblright\ (meaning: a subring isomorphic to $\mathbb{Q}%
$): namely, the map%
\begin{align*}
\mathbb{Q}  &  \rightarrow F,\\
\dfrac{a}{b}  &  \mapsto\dfrac{a\cdot1_{F}}{b\cdot1_{F}}%
\ \ \ \ \ \ \ \ \ \ \left(  \text{for }a,b\in\mathbb{Z}\text{ with }%
b\neq0\right)
\end{align*}
is an injective ring morphism.
\end{enumerate}
\end{theorem}

\begin{proof}
We have $p\cdot1_{F}=0_{F}$. Indeed, if $p=0$, then this is obvious; but
otherwise it follows from the definition of $\operatorname*{char}F$. \medskip

\textbf{(b)} Let $a\in F$. Then, $a=1_{F}\cdot a$. Thus,%
\[
pa=p\left(  1_{F}\cdot a\right)  =\underbrace{\left(  p\cdot1_{F}\right)
}_{=0_{F}}\cdot\,a=0_{F}\cdot a=0_{F}=0.
\]
This proves Theorem \ref{thm.finfield.char.basics} \textbf{(b)}. \medskip

\textbf{(a)} We define an action of the ring $\mathbb{Z}/p$ on $F$ by
\[
\overline{k}\cdot a=ka\ \ \ \ \ \ \ \ \ \ \text{for all }k\in\mathbb{Z}\text{
and }a\in F.
\]
Why is this well-defined? In other words, why is it true that if two integers
$k$ and $\ell$ satisfy $\overline{k}=\overline{\ell}$, then $ka=\ell a$ for
all $a\in F$ ?

Let us check this directly: Let $k$ and $\ell$ be two integers satisfying
$\overline{k}=\overline{\ell}$ in $\mathbb{Z}/p$. This means $k\equiv
\ell\operatorname{mod}p$, so that $k-\ell$ is a multiple of $p$. That is,
$k-\ell=pu$ for some $u\in\mathbb{Z}$. Consider this $u$. Now,%
\[
ka-\ell a=\underbrace{\left(  k-\ell\right)  }_{=pu}a=pua=0
\]
(by Theorem \ref{thm.finfield.char.basics} \textbf{(b)}, applied to $ua$
instead of $a$). Thus, $ka=\ell a$, which is precisely what we wanted to prove.

Thus, the action of $\mathbb{Z}/p$ on $F$ is well-defined. Now, it remains to
show that $F$ is a $\mathbb{Z}/p$-module, and that the \textquotedblleft
scale-invariance\textquotedblright\ axiom is satisfied. All of this is easy
and LTTR\footnote{For example, let us prove the associativity law, which says
that $\left(  rs\right)  m=r\left(  sm\right)  $ for all $r,s\in\mathbb{Z}/p$
and $m\in F$. Indeed, let $r,s\in\mathbb{Z}/p$ and $m\in F$. Write $r$ and $s$
as $\overline{k}$ and $\overline{\ell}$ for some integers $k$ and $\ell$.
Then, $rs=\overline{k}\cdot\overline{\ell}=\overline{k\ell}$, so that $\left(
rs\right)  m=\overline{k\ell}\cdot m=k\ell m$ (by our definition of the action
of $\mathbb{Z}/p$ on $F$). Also, from $r=\overline{k}$ and $s=\overline{\ell}%
$, we obtain%
\begin{align*}
r\left(  sm\right)   &  =\overline{k}\cdot\left(  \overline{\ell}\cdot
m\right)  =k\left(  \overline{\ell}\cdot m\right)  \ \ \ \ \ \ \ \ \ \ \left(
\text{by our definition of the action}\right) \\
&  =k\left(  \ell m\right)  \ \ \ \ \ \ \ \ \ \ \left(  \text{since our
definition of the action yields }\overline{\ell}\cdot m=\ell m\right) \\
&  =k\ell m.
\end{align*}
Comparing this with $\left(  rs\right)  m=k\ell m$, we obtain $\left(
rs\right)  m=r\left(  sm\right)  $, qed.}. Thus, $F$ becomes a $\mathbb{Z}%
/p$-algebra. This proves Theorem \ref{thm.finfield.char.basics} \textbf{(a)}.
\medskip

\textbf{(c)} Assume the contrary. Thus, $p$ is neither a prime nor $0$. Hence,
$p$ is either $1$ or a composite\footnote{A positive integer is said to be
\textbf{composite} if it can be written as a product of two integers each
larger than $1$.} positive integer (since $p$ is always a nonnegative integer).

Since $F$ is a field, we have $1\neq0$ in $F$. In other words, $1_{F}\neq
0_{F}$. If we had $p=1$, then we would thus have $\underbrace{p}_{=1}%
\cdot\,1_{F}=1\cdot1_{F}=1_{F}\neq0_{F}$, which would contradict $p\cdot
1_{F}=0_{F}$. Thus, we cannot have $p=1$. Hence, $p$ must be composite (since
$p$ is either $1$ or composite). In other words, $p=uv$ for some integers
$u>1$ and $v>1$. Consider these integers $u$ and $v$.

From $u>1$ and $v>1$ and $p=uv$, we see that both integers $u$ and $v$ are
smaller than $p$. Hence, neither $u\cdot1_{F}$ nor $v\cdot1_{F}$ can be
$0_{F}$ (since $p=\operatorname*{char}F$ was defined to be the
\textbf{smallest} positive integer $n$ such that $n\cdot1_{F}=0_{F}$). Since
$F$ is an integral domain (because $F$ is a field), this yields that the
product $\left(  u\cdot1_{F}\right)  \cdot\left(  v\cdot1_{F}\right)  $ is
also nonzero.

Now, $p\cdot1_{F}=0_{F}$, so%
\[
0_{F}=\underbrace{p}_{=uv}\cdot\,1_{F}=uv\cdot1_{F}=\left(  u\cdot
1_{F}\right)  \cdot\left(  v\cdot1_{F}\right)  .
\]
This contradicts the fact that the product $\left(  u\cdot1_{F}\right)
\cdot\left(  v\cdot1_{F}\right)  $ is nonzero. This proves Theorem
\ref{thm.finfield.char.basics} \textbf{(c)}. \medskip

\textbf{(d)} Assume that $F$ is finite. We must show that $p$ is a prime.

According to Theorem \ref{thm.finfield.char.basics} \textbf{(c)}, it suffices
to show that $p\neq0$. So let us show this. Assume the contrary. Then, $p=0$.
Hence, none of the elements $1\cdot1_{F}$, $2\cdot1_{F}$, $3\cdot1_{F}$,
$\ldots$ of $F$ is $0_{F}$ (by the definition of $\operatorname*{char}F$). But
$F$ is finite, so two of these elements must be equal (by the Pigeonhole
Principle). In other words, there exist positive integers $u<v$ such that
$u\cdot1_{F}=v\cdot1_{F}$. Consider these $u$ and $v$. Then, $v-u$ is a
positive integer, and we have $\left(  v-u\right)  \cdot1_{F}=v\cdot
1_{F}-u\cdot1_{F}=0_{F}$ (since $u\cdot1_{F}=v\cdot1_{F}$). But $\left(
v-u\right)  \cdot1_{F}$ is one of the elements $1\cdot1_{F}$, $2\cdot1_{F}$,
$3\cdot1_{F}$, $\ldots$ (since $u<v$), and we just said that none of these
elements is $0_{F}$. This contradicts $\left(  v-u\right)  \cdot1_{F}=0_{F}$.
Thus, our assumption was false; hence, Theorem \ref{thm.finfield.char.basics}
\textbf{(d)} is proven. \medskip

\textbf{(e)} \textit{First proof of part \textbf{(e)}:} Assume that $F$ is
finite. Thus, by Theorem \ref{thm.finfield.char.basics} \textbf{(d)}, we know
that $p$ is prime.

Since $F$ is a field, we have $1\neq0$ in $F$. Hence, $\left\vert F\right\vert
>1$.

From Theorem \ref{thm.finfield.char.basics} \textbf{(a)}, we know that $F$ is
a $\mathbb{Z}/p$-algebra. Thus, in particular, $F$ is a $\mathbb{Z}/p$-module.
But since $p$ is prime, $\mathbb{Z}/p$ is a field.

Now, recall that a module over a field is nothing but a vector space. In
particular, every module over a field is free (since any vector space has a
basis\footnote{This fact is Theorem \ref{thm.vs.free}.
\par
Once again, I haven't actually proved this fact in this course, but you can
easily bridge this gap yourself or look it up in any text on linear algebra
(or in Keith Conrad's
\url{https://kconrad.math.uconn.edu/blurbs/linmultialg/dimension.pdf} ). Our
situation is simpler than the general case, since we know that $F$ is finite,
so it is clear that there is a finite list of vectors in $F$ that span $F$
(because you can just take a list of \textbf{all} elements of $F$). In order
to obtain a basis from such a list, you only need to successively remove
vectors that are linear combinations of other vectors; once no such vectors
remain, the list will be a basis (make sure you understand why!).}). Thus, in
particular, the $\mathbb{Z}/p$-module $F$ is free. In other words, the
$\mathbb{Z}/p$-module $F$ has a basis. This basis must be finite (since $F$
itself is finite). Thus, $F\cong\left(  \mathbb{Z}/p\right)  ^{m}$ as
$\mathbb{Z}/p$-modules for some $m\in\mathbb{N}$. Consider this $m$. From
$F\cong\left(  \mathbb{Z}/p\right)  ^{m}$, we obtain $\left\vert F\right\vert
=\left\vert \left(  \mathbb{Z}/p\right)  ^{m}\right\vert =\left\vert
\mathbb{Z}/p\right\vert ^{m}=p^{m}$. It remains to prove that $m$ is positive.
But this is easy: If $m$ was $0$, then $\left\vert F\right\vert =p^{m}$ would
imply $\left\vert F\right\vert =p^{0}=1$, which would contradict $\left\vert
F\right\vert >1$. Thus, the proof of Theorem \ref{thm.finfield.char.basics}
\textbf{(e)} is complete. \medskip

\textit{Second proof of part \textbf{(e)}:} There is an alternative proof of
Theorem \ref{thm.finfield.char.basics} \textbf{(e)}, which avoids any use of
linear algebra but instead uses some group theory. Specifically, we will use
\textbf{Cauchy's theorem}, which says the following: If $G$ is a finite group,
and if $q$ is a prime number that divides the size $\left\vert G\right\vert $,
then $G$ has an element of order $q$. Proofs of this theorem can be found,
e.g., in \url{https://kconrad.math.uconn.edu/blurbs/grouptheory/cauchypf.pdf}
or in \cite[Theorem 4.9.4]{Sharif22} or \cite[Remark after Theorem
4.59]{Knapp1} or \cite[Theorem 21.22]{Elman22} or \cite[Corollary
2.4.14]{Ford22} or \cite[Lemma 17.2.2]{Hudso25}, just to mention some freely
available sources. (It is also an easy consequence of the first Sylow theorem.)

Now, assume that $F$ is finite. Thus, by Theorem
\ref{thm.finfield.char.basics} \textbf{(d)}, we know that $p$ is prime.

Since $F$ is a field, we have $1\neq0$ in $F$. Hence, $\left\vert F\right\vert
>1$.

Assume (for the sake of contradiction) that the integer $\left\vert
F\right\vert $ has a prime divisor $q$ distinct from $p$. Thus, Cauchy's
theorem (applied to $G$ being the additive group $\left(  F,+,0\right)  $)
yields that the additive group $\left(  F,+,0\right)  $ has an element of
order $q$. Let $a$ be this element. Then, $a\neq0$ but $qa=0$ (since $a$ has
order $q$). However, Theorem \ref{thm.finfield.char.basics} \textbf{(b)}
yields $pa=0$.

However, $p$ and $q$ are two distinct primes, and thus are coprime. Thus,
$\gcd\left(  p,q\right)  =1$. But Bezout's theorem shows that there exist two
integers $x$ and $y$ such that $xp+yq=\gcd\left(  p,q\right)  $. Consider
these $x$ and $y$. Then,
\[
\underbrace{\left(  xp+yq\right)  }_{=\gcd\left(  p,q\right)  =1}a=1a=a\neq0
\]
contradicts%
\[
\left(  xp+yq\right)  a=x\underbrace{pa}_{=0}+\,y\underbrace{qa}%
_{=0}=x0+y0=0.
\]
This contradiction shows that our assumption (that the integer $\left\vert
F\right\vert $ has a prime divisor $q$ distinct from $p$) is false.

Hence, the integer $\left\vert F\right\vert $ has no prime divisor distinct
from $p$. Thus, the only prime divisor of $\left\vert F\right\vert $ is $p$.
Therefore, $\left\vert F\right\vert =p^{m}$ for some $m\in\mathbb{N}$. This
$m$ must furthermore be positive (since $\left\vert F\right\vert >1$). Thus,
Theorem \ref{thm.finfield.char.basics} \textbf{(e)} is proven again. \medskip

\textbf{(f)} Assume that $p$ is a prime. Then, $F$ is a $\mathbb{Z}/p$-algebra
(by Theorem \ref{thm.finfield.char.basics} \textbf{(a)}), so we can define a
map%
\begin{align*}
\mathbb{Z}/p  &  \rightarrow F,\\
\alpha &  \mapsto\alpha\cdot1_{F}.
\end{align*}
It is straightforward to check that this map is a ring morphism\footnote{For
instance, it respects multiplication because $\left(  \alpha\cdot1_{F}\right)
\left(  \beta\cdot1_{F}\right)  =\alpha\beta\cdot\underbrace{1_{F}1_{F}%
}_{=1_{F}}=\alpha\beta\cdot1_{F}$ for any $\alpha,\beta\in\mathbb{Z}/p$.};
furthermore, it is easily seen to be injective\footnote{This is actually best
understood as a particular case of the following general fact: \textbf{Any}
ring morphism from a field to a nontrivial ring is injective!
\par
The proof of this general fact is pretty easy: Let $f:K\rightarrow R$ be a
ring morphism from a field $K$ to a nontrivial ring $R$. We must show that $f$
is injective.
\par
Assume that $a$ is a nonzero element of $\operatorname*{Ker}f$. Then, $a$ is a
unit of $K$ (since $K$ is a field, and thus any nonzero element of $K$ is a
unit), and thus $a^{-1}$ exists. Since $f$ is a ring morphism, we have%
\[
f\left(  aa^{-1}\right)  =\underbrace{f\left(  a\right)  }%
_{\substack{=0\\\text{(since }a\in\operatorname*{Ker}f\text{)}}}\cdot
\,f\left(  a^{-1}\right)  =0,
\]
so that $0=f\left(  \underbrace{aa^{-1}}_{=1}\right)  =f\left(  1\right)  =1$
(since $f$ is a ring morphism). This entails that the ring $R$ is trivial,
which contradicts our assumption that $R$ is nontrivial.
\par
Forget that we fixed $a$. We thus obtained a contradiction for any nonzero
element $a$ of $\operatorname*{Ker}f$. Hence, no such element exists. In other
words, any $a\in\operatorname*{Ker}f$ must be $0$. Thus, $\operatorname*{Ker}%
f\subseteq\left\{  0\right\}  $, so that $\operatorname*{Ker}f=\left\{
0\right\}  $, and therefore $f$ is injective (by Lemma
\ref{lem.rings.mors.inj}). This completes our proof.}. Hence, its image is a
subring of $F$ that is isomorphic to $\mathbb{Z}/p$ (by Proposition
\ref{prop.ringmor.inj-image} \textbf{(a)}). This proves Theorem
\ref{thm.finfield.char.basics} \textbf{(f)}. \medskip

\textbf{(g)} We will be very brief, since we won't use Theorem
\ref{thm.finfield.char.basics} \textbf{(g)} in what follows.

Assume that $p=0$. Then, for any nonzero integer $b$, the element $b\cdot
1_{F}$ of $F$ is nonzero (why?) and therefore a unit of $F$ (since $F$ is a
field). Hence, for any rational number $\dfrac{a}{b}\in\mathbb{Q}$ (written in
such a way that $a,b\in\mathbb{Z}$ and $b\neq0$), the element $\dfrac
{a\cdot1_{F}}{b\cdot1_{F}}\in F$ is well-defined. Now, of course, the
representation of a rational number as $\dfrac{a}{b}$ with $a,b\in\mathbb{Z}$
is not unique (for instance, $\dfrac{6}{4}$ and $\dfrac{3}{2}$ are the same
rational number); however, it is not hard to show that $\dfrac{a\cdot1_{F}%
}{b\cdot1_{F}}$ is uniquely determined by $\dfrac{a}{b}$ (meaning that if
$a,b,c,d\in\mathbb{Z}$ satisfy $\dfrac{a}{b}=\dfrac{c}{d}$, then we also have
$\dfrac{a\cdot1_{F}}{b\cdot1_{F}}=\dfrac{c\cdot1_{F}}{d\cdot1_{F}}$). Thus,
the map
\begin{align*}
\mathbb{Q}  &  \rightarrow F,\\
\dfrac{a}{b}  &  \mapsto\dfrac{a\cdot1_{F}}{b\cdot1_{F}}%
\ \ \ \ \ \ \ \ \ \ \left(  \text{for }a,b\in\mathbb{Z}\text{ with }%
b\neq0\right)
\end{align*}
is well-defined. Next, it can be shown that this map is a ring morphism and is
injective\footnote{The injectivity follows just as in part \textbf{(f)}.}.
Hence, its image is a subring of $F$ that is isomorphic to $\mathbb{Q}$. This
proves Theorem \ref{thm.finfield.char.basics} \textbf{(g)}.
\end{proof}

Parts \textbf{(f)} and \textbf{(g)} of Theorem \ref{thm.finfield.char.basics}
show that any field $F$ has at its \textquotedblleft core\textquotedblright\ a
\textquotedblleft small\textquotedblright\ field: either (a copy of)
$\mathbb{Z}/p$ (if its characteristic is a prime $p$) or (a copy of)
$\mathbb{Q}$ (if its characteristic is $0$).

Parts \textbf{(d)} and \textbf{(e)} of Theorem \ref{thm.finfield.char.basics}
(in combination) show that the size of any finite field is a power of a prime.
Thus, there are no finite fields of size $6$ or $10$ or $12$.

Hence, we can limit our search for finite fields to those of size $p^{m}$ for
$p$ prime and $m>0$. We have already found such fields for $m=1$ and for $m=2$
(for all $p$), and briefly hinted at the cases $m=3$ and $m=4$, but we are
still missing the case of general $m$.

\begin{exercise}
\label{exe.polring.xxx-r}\ \ \ 

\begin{enumerate}
\item[\textbf{(a)}] Let $F$ be any field, and let $p\in F\left[  x\right]  $
be a polynomial of degree $2$ or $3$. Prove that $p$ is irreducible if and
only if $p$ has no root in $F$.

\item[\textbf{(b)}] Let $r\in\mathbb{Q}$ be a rational number that is not the
cube of any rational number. Prove that the polynomial $x^{3}-r\in
\mathbb{Q}\left[  x\right]  $ is irreducible.

\item[\textbf{(c)}] Prove that there is an injective ring morphism%
\begin{align*}
\mathbb{Q}\left[  x\right]  /\left(  x^{3}-5\right)   &  \rightarrow
\mathbb{R},\\
\overline{p}  &  \mapsto p\left(  \sqrt[3]{5}\right)  .
\end{align*}

\item[\textbf{(d)}] In Subsection \ref{subsec.rings.def.exas}, we claimed that
there are no $a,b\in\mathbb{Q}$ that satisfy $1+2\sqrt[3]{5}+\left(
\sqrt[3]{5}\right)  ^{2}=a+b\sqrt[3]{5}$. Prove this claim.
\end{enumerate}
\end{exercise}

\subsection{\label{sec.finfields.tools}Tools}

\subsubsection{Splitting polynomials}

We are still trying to find finite fields of sizes $p^{m}$ for $p$ prime and
$m>0$. We will approach the general case indirectly (no easy and direct proofs
are known). We will need a bunch of tools. The first is the notion of a
\textbf{splitting field}. We begin with a definition:

\begin{definition}
\label{def.finfield.splits}Let $R$ be a commutative ring. Let $b\in R\left[
x\right]  $ be a polynomial over $R$. We say that $b$ \textbf{splits} over $R$
if there exist elements $r_{1},r_{2},\ldots,r_{m}$ of $R$ such that%
\[
b=\left(  x-r_{1}\right)  \left(  x-r_{2}\right)  \cdots\left(  x-r_{m}%
\right)  .
\]

\end{definition}

Note that in this definition, we must necessarily have $\deg b=m$ (unless $R$
is trivial). Also, a polynomial cannot split unless it is monic. This might
differ from how other authors define the notion of \textquotedblleft
splitting\textquotedblright, but it is sufficient for what we will do with it.

\begin{example}
\label{exa.finfield.splits.1}\ \ 

\begin{enumerate}
\item[\textbf{(a)}] The polynomial $x^{2}-1$ splits over $\mathbb{Q}$, since%
\[
x^{2}-1=\left(  x-1\right)  \left(  x+1\right)  =\left(  x-1\right)  \left(
x-\left(  -1\right)  \right)  .
\]

\item[\textbf{(b)}] The polynomial $x^{2}+1$ does not split over $\mathbb{R}$
(since it has no roots in $\mathbb{R}$), but it splits over $\mathbb{C}$,
since%
\[
x^{2}+1=\left(  x-i\right)  \left(  x+i\right)  =\left(  x-i\right)  \left(
x-\left(  -i\right)  \right)  .
\]

\item[\textbf{(c)}] The polynomial $x^{2}$ splits over $\mathbb{Q}$, since
$x^{2}=xx=\left(  x-0\right)  \left(  x-0\right)  $.

\item[\textbf{(d)}] The polynomial $x^{4}-9$ does not split over $\mathbb{R}$.
Indeed, it has a factorization%
\[
x^{4}-9=\left(  x-\sqrt{3}\right)  \left(  x+\sqrt{3}\right)  \left(
x^{2}+3\right)  ,
\]
but the $x^{2}+3$ factor is still not of the form $x-r$ and cannot be factored
further over $\mathbb{R}$. However, this polynomial does split over
$\mathbb{C}$, since%
\[
x^{4}-9=\left(  x-\sqrt{3}\right)  \left(  x+\sqrt{3}\right)  \left(
x-\sqrt{3}i\right)  \left(  x+\sqrt{3}i\right)  .
\]

\item[\textbf{(e)}] Any monic polynomial of degree $1$ automatically splits
over whatever commutative ring it is defined over. So does the constant
polynomial $1$ (since it is an empty product).
\end{enumerate}
\end{example}

When a polynomial splits over a field, its roots can be read off directly from
the splitting:

\begin{proposition}
\label{prop.finfield.splits.roots}Let $F$ be a field. Let $r_{1},r_{2}%
,\ldots,r_{m}\in F$. Then,
\begin{align*}
&  \left\{  \text{the roots of the polynomial }\left(  x-r_{1}\right)  \left(
x-r_{2}\right)  \cdots\left(  x-r_{m}\right)  \in F\left[  x\right]  \text{ in
}F\right\} \\
&  =\left\{  r_{1},r_{2},\ldots,r_{m}\right\}  .
\end{align*}

\end{proposition}

\begin{proof}
The ring $F$ is a field, thus an integral domain. Thus, a product $uv$ of two
elements $u,v\in F$ is zero if and only if one of its factors is zero. Hence,
a finite product $u_{1}u_{2}\cdots u_{k}$ of elements of $F$ is zero if and
only if one of its factors is zero\footnote{Indeed, this follows easily by
induction on $k$, using the preceding sentence in the induction step.}.

Now, we have%
\begin{align*}
&  \left\{  \text{the roots of the polynomial }\left(  x-r_{1}\right)  \left(
x-r_{2}\right)  \cdots\left(  x-r_{m}\right)  \in F\left[  x\right]  \text{ in
}F\right\} \\
&  =\left\{  a\in F\ \mid\ \left(  \left(  x-r_{1}\right)  \left(
x-r_{2}\right)  \cdots\left(  x-r_{m}\right)  \right)  \left[  a\right]
=0\right\} \\
&  \ \ \ \ \ \ \ \ \ \ \ \ \ \ \ \ \ \ \ \ \left(  \text{by the definition of
a \textquotedblleft root\textquotedblright}\right) \\
&  =\left\{  a\in F\ \mid\ \left(  a-r_{1}\right)  \left(  a-r_{2}\right)
\cdots\left(  a-r_{m}\right)  =0\right\} \\
&  \ \ \ \ \ \ \ \ \ \ \ \ \ \ \ \ \ \ \ \ \left(
\begin{array}
[c]{c}%
\text{since the evaluation }\left(  \left(  x-r_{1}\right)  \left(
x-r_{2}\right)  \cdots\left(  x-r_{m}\right)  \right)  \left[  a\right] \\
\text{equals }\left(  a-r_{1}\right)  \left(  a-r_{2}\right)  \cdots\left(
a-r_{m}\right)
\end{array}
\right) \\
&  =\left\{  a\in F\ \mid\ \text{one of }a-r_{1},a-r_{2},\ldots,a-r_{m}\text{
is zero}\right\} \\
&  \ \ \ \ \ \ \ \ \ \ \ \ \ \ \ \ \ \ \ \ \left(
\begin{array}
[c]{c}%
\text{since a finite product }u_{1}u_{2}\cdots u_{k}\text{ of elements of
}F\text{ is zero}\\
\text{if and only if one of its factors is zero}%
\end{array}
\right) \\
&  =\left\{  a\in F\ \mid\ a=r_{1}\text{ or }a=r_{2}\text{ or }\cdots\text{ or
}a=r_{m}\right\} \\
&  =\left\{  r_{1},r_{2},\ldots,r_{m}\right\}  .
\end{align*}
This proves Proposition \ref{prop.finfield.splits.roots}.
\end{proof}

\begin{remark}
It is worth noting that Proposition \ref{prop.finfield.splits.roots} still
holds if we replace \textquotedblleft field\textquotedblright\ by
\textquotedblleft integral domain\textquotedblright\ (and the same proof
applies); but it does not hold when $F$ is just a general commutative ring.
For example, if $F=\mathbb{Z}/4$, then the polynomial $\left(  x-0\right)
\left(  x-0\right)  \left(  x-1\right)  \left(  x-3\right)  \in F\left[
x\right]  $ has roots $0,1,2,3$, rather than just $0,1,3$ as Proposition
\ref{prop.finfield.splits.roots} would predict. A similar construction works
for $\mathbb{Z}/n$ where $n$ is any composite integer $>1$ (see Exercise
\ref{exe.prop.finfield.splits.roots.Z/n}).
\end{remark}

\begin{exercise}
\label{exe.prop.finfield.splits.roots.Z/n}Let $n>1$ be a composite integer
(i.e., an integer that is not prime). Prove the following:

\begin{enumerate}
\item[\textbf{(a)}] Each element of $\mathbb{Z}/n$ is a root of the polynomial
$\prod_{i=1}^{n-1}\left(  x-\overline{i}\right)  ^{2}\in\left(  \mathbb{Z}%
/n\right)  \left[  x\right]  $.

\item[\textbf{(b)}] If $n>4$, then each element of $\mathbb{Z}/n$ is a root of
the polynomial $\prod_{i=1}^{n-1}\left(  x-\overline{i}\right)  \in\left(
\mathbb{Z}/n\right)  \left[  x\right]  $.
\end{enumerate}
\end{exercise}

\subsubsection{Splitting fields}

The Fundamental Theorem of Algebra says that each monic univariate polynomial
over $\mathbb{C}$ splits over $\mathbb{C}$. This is not actually a theorem of
algebra, since it relies on the definition of $\mathbb{C}$ (which is
analytic); however, it explains some of the significance of $\mathbb{C}$. In
general, a field $F$ is said to be \textbf{algebraically closed} if each monic
univariate polynomial over $F$ splits over $F$. The field $\mathbb{C}$ is not
the only algebraically closed field, but it is perhaps the best-known.

We won't need algebraically closed fields in this course; we will need a more
\textquotedblleft local\textquotedblright\ notion: that of a splitting field.
To introduce it, we make a simple observation, which we have already (tacitly)
used in Example \ref{exa.finfield.splits.1} (as we have been treating the same
polynomial $x^{2}+1$ first as a polynomial in $\mathbb{R}\left[  x\right]  $
and then as a polynomial in $\mathbb{C}\left[  x\right]  $):

\begin{proposition}
\label{prop.polring.subring}Let $S$ be a commutative ring. Let $R$ be a
subring of $S$. Then, any polynomial over $R$ automatically is a polynomial
over $S$ as well (since its coefficients lie in $R$ and therefore also lie in
$S$), and thus the polynomial ring $R\left[  x\right]  $ becomes a subring of
$S\left[  x\right]  $.
\end{proposition}

For example, $\mathbb{R}\left[  x\right]  $ is a subring of $\mathbb{C}\left[
x\right]  $. Polynomials like $x^{2}+1$ might not split over $\mathbb{R}$, but
they split over $\mathbb{C}$. This suggests that if a monic polynomial does
not split over a ring, we might fix this by making the ring larger
(\textquotedblleft extending\textquotedblright\ the ring, possibly by
\textquotedblleft adjoining\textquotedblright\ some roots), just as
$\mathbb{C}$ was constructed from $\mathbb{R}$ in order to make $x^{2}+1$
split. Thus we make the following definition:

\begin{definition}
\label{def.finfield.splitfield}Let $F$ be a field. Let $b\in F\left[
x\right]  $ be a monic polynomial over $F$. Then, a \textbf{splitting field}
of $b$ (over $F$) means a field $S$ such that

\begin{itemize}
\item $F$ is a subring of $S$;

\item the polynomial $b$ (regarded as a polynomial in $S\left[  x\right]  $)
splits over $S$.
\end{itemize}
\end{definition}

Examples:

\begin{itemize}
\item $\mathbb{C}$ is a splitting field of $x^{2}+1$ over $\mathbb{R}$.

\item $\mathbb{C}$ is a splitting field of $x^{2}-2$ over $\mathbb{Q}$, but so
is $\mathbb{R}$ (since $x^{2}-2$ already splits over $\mathbb{R}$) or even the
smaller field $\mathbb{Q}\left[  \sqrt{2}\right]  =\left\{  a+b\sqrt{2}%
\ \mid\ a,b\in\mathbb{Q}\right\}  $.

\item $\mathbb{Q}$ itself is a splitting field of $x^{2}-1$ over $\mathbb{Q}$.
\end{itemize}

(Be careful with the literature: Many authors have a more restrictive concept
of a \textquotedblleft splitting field\textquotedblright, which requires not
only that the polynomial split over it, but also that the field -- in some
reasonable way -- is minimal with this property. For example, these authors do
not accept $\mathbb{R}$ as a splitting field of $x^{2}-2$ over $\mathbb{Q}$,
since the much smaller field $\mathbb{Q}\left[  \sqrt{2}\right]  $ suffices to
split the polynomial. But our definition suffices for our purposes.)

The most important fact about splitting fields is that they always exist:

\begin{theorem}
\label{thm.finfield.splitfield}Let $F$ be a field. Let $b\in F\left[
x\right]  $ be a monic polynomial over $F$. Then:

\begin{enumerate}
\item[\textbf{(a)}] We can write $b$ as a product $b=c_{1}c_{2}\cdots c_{k}$
of monic irreducible polynomials $c_{1},c_{2},\ldots,c_{k}\in F\left[
x\right]  $.

\item[\textbf{(b)}] If $\deg b>0$, then there is a field that contains $F$ as
a subring and that contains a root of $b$.

\item[\textbf{(c)}] There exists a splitting field of $b$ over $F$.
\end{enumerate}
\end{theorem}

\begin{proof}
\textbf{(a)} Any nonzero polynomial in $F\left[  x\right]  $ can be made monic
by scaling it with a nonzero scalar (namely, if $g\in F\left[  x\right]  $ is
a nonzero polynomial, and if $c$ is its leading coefficient, then $c^{-1}g$ is
a monic polynomial). This scaling does not interfere with its divisibility
properties; thus, if $g$ is irreducible, then it remains so after the scaling.

Hence, it suffices to show that we can write $b$ as a product $b=c_{1}%
c_{2}\cdots c_{k}$ of irreducible polynomials $c_{1},c_{2},\ldots,c_{k}\in
F\left[  x\right]  $.

Abstractly, this follows easily from the fact that $F\left[  x\right]  $ is a
UFD. In a more down-to-earth manner, this can be shown just like the classical
fact that each positive integer can be written as a product of primes. The
proof proceeds by strong induction on $\deg b$; the main idea is
\textquotedblleft either $b$ is itself irreducible, in which case we are done;
or $b$ can be written as a product of two polynomials of smaller degree, in
which case the induction hypothesis applies\textquotedblright.

(Note that this proof is constructive when $F$ is finite, since we can
actually try out all polynomials of degree smaller than $\deg b$ and check
which of them divide $b$.)

Theorem \ref{thm.finfield.splitfield} \textbf{(a)} is thus proved.

\textbf{(b)} Assume that $\deg b>0$. We must find a field that contains $F$ as
a subring and that contains a root of $b$.

It is tempting to take $F\left[  x\right]  /b$, but this might fail to be a
field (since $b$ might fail to be irreducible).

Instead, we use Theorem \ref{thm.finfield.splitfield} \textbf{(a)} to write
$b$ as a product $b=c_{1}c_{2}\cdots c_{k}$ of monic irreducible polynomials
$c_{1},c_{2},\ldots,c_{k}\in F\left[  x\right]  $, and then we take the field
$F\left[  x\right]  /c_{1}$ (which is indeed a field, because $c_{1}$ is
irreducible\footnote{We are using Theorem \ref{thm.fieldext.field-irr}
here.}). This field will contain a root of $c_{1}$, and thus also contain a
root of $b$ (since a root of $c_{1}$ is always a root of $b$). So Theorem
\ref{thm.finfield.splitfield} \textbf{(b)} is proved.

(Where did I use the assumption $\deg b>0$ in this proof? Hint: Why is there a
$c_{1}$ ?)

\textbf{(c)} Here is a proof by example: Assume that $\deg b=3$.

Theorem \ref{thm.finfield.splitfield} \textbf{(b)} says that there is a field
$F^{\prime}$ that contains $F$ as a subring and that contains a root of $b$.
Consider this $F^{\prime}$, and let $r_{1}$ be the root of $b$ that it
contains. Thus, $x-r_{1}\mid b$ in $F^{\prime}\left[  x\right]  $ (since
$r_{1}$ is a root of $b$). Hence, the polynomial $\dfrac{b}{x-r_{1}}\in
F^{\prime}\left[  x\right]  $ is well-defined. Moreover, this polynomial
$\dfrac{b}{x-r_{1}}$ has degree $3-1=2$ and is monic\footnote{Here, we are
using the fact that when we divide a monic polynomial $f$ by a monic
polynomial $g$ with $\deg g \leq\deg f$, the quotient will again be monic.
(The proof is LTTR. Note that this holds even if there is a remainder!)}.

Now, we apply Theorem \ref{thm.finfield.splitfield} \textbf{(b)} again, but
this time to the field $F^{\prime}$ and the monic polynomial $\dfrac
{b}{x-r_{1}}$ over it. Thus we conclude that there is a field $F^{\prime
\prime}$ that contains $F^{\prime}$ as a subring and that contains a root of
$\dfrac{b}{x-r_{1}}$. Consider this $F^{\prime\prime}$, and let $r_{2}$ be the
root of $\dfrac{b}{x-r_{1}}$ that it contains. Thus, $x-r_{2}\mid\dfrac
{b}{x-r_{1}}$ in $F^{\prime\prime}\left[  x\right]  $ (since $r_{2}$ is a root
of $\dfrac{b}{x-r_{1}}$). Hence, the polynomial $\dfrac{b}{x-r_{1}}/\left(
x-r_{2}\right)  \in F^{\prime\prime}\left[  x\right]  $ is well-defined. In
other words, the polynomial $\dfrac{b}{\left(  x-r_{1}\right)  \left(
x-r_{2}\right)  }\in F^{\prime\prime}\left[  x\right]  $ is well-defined.
Moreover, this polynomial $\dfrac{b}{\left(  x-r_{1}\right)  \left(
x-r_{2}\right)  }$ has degree $3-2=1$ and is monic.

Now, we apply Theorem \ref{thm.finfield.splitfield} \textbf{(b)} again, but
this time to the field $F^{\prime\prime}$ and the monic polynomial $\dfrac
{b}{\left(  x-r_{1}\right)  \left(  x-r_{2}\right)  }$ over it. Thus we
conclude that there is a field $F^{\prime\prime\prime}$ that contains
$F^{\prime\prime}$ as a subring and that contains a root of $\dfrac{b}{\left(
x-r_{1}\right)  \left(  x-r_{2}\right)  }$. Consider this $F^{\prime
\prime\prime}$, and let $r_{3}$ be the root of $\dfrac{b}{\left(
x-r_{1}\right)  \left(  x-r_{2}\right)  }$ that it contains. Thus,
$x-r_{3}\mid\dfrac{b}{\left(  x-r_{1}\right)  \left(  x-r_{2}\right)  }$ in
$F^{\prime\prime\prime}\left[  x\right]  $. Hence, the polynomial
\newline$\dfrac{b}{\left(  x-r_{1}\right)  \left(  x-r_{2}\right)  \left(
x-r_{3}\right)  }\in F^{\prime\prime\prime}\left[  x\right]  $ is
well-defined. Furthermore, this polynomial has degree $3-3=0$ and is monic. In
other words, this polynomial equals $1$. In other words, $b=\left(
x-r_{1}\right)  \left(  x-r_{2}\right)  \left(  x-r_{3}\right)  $ in
$F^{\prime\prime\prime}\left[  x\right]  $. This shows that $b$ splits over
$F^{\prime\prime\prime}$. Moreover, by construction, $F^{\prime\prime\prime}$
is a field that contains $F$ as a subring (since $F\subseteq F^{\prime
}\subseteq F^{\prime\prime}\subseteq F^{\prime\prime\prime}$, and each of
these \textquotedblleft$\subseteq$\textquotedblright\ signs is not just a
subset but actually a subring).

Thus, we have proved Theorem \ref{thm.finfield.splitfield} \textbf{(c)} in our
example. Proving it in the general case is just a matter of formalizing what
we did as an induction on $\deg b$.
\end{proof}

\subsubsection{\label{subsec.finfields.tools.frob}The Idiot's Binomial Formula
and the Frobenius endomorphism}

Next, to something different. A rather surprising property of fields of
positive characteristic is the following theorem (often called
\textbf{Freshman's Dream} or \textbf{Idiot's Binomial Formula} due to its
similarity to a popular student mistake):

\begin{theorem}
[Idiot's Binomial Formula, aka Freshman's Dream]\label{thm.finfield.idiot}Let
$p$ be a prime number. Let $F$ be a field of characteristic $p$, or, more
generally, any commutative $\mathbb{Z}/p$-algebra. Then:

\begin{enumerate}
\item[\textbf{(a)}] We have $\left(  a+b\right)  ^{p}=a^{p}+b^{p}$ for any
$a,b\in F$.

\item[\textbf{(b)}] We have $\left(  a+b\right)  ^{p^{m}}=a^{p^{m}}+b^{p^{m}}$
for any $a,b\in F$ and $m\in\mathbb{N}$.

\item[\textbf{(c)}] We have $\left(  a-b\right)  ^{p}=a^{p}-b^{p}$ for any
$a,b\in F$.

\item[\textbf{(d)}] We have $\left(  a-b\right)  ^{p^{m}}=a^{p^{m}}-b^{p^{m}}$
for any $a,b\in F$ and $m\in\mathbb{N}$.
\end{enumerate}
\end{theorem}

For example, for $p=3$, Theorem \ref{thm.finfield.idiot} \textbf{(a)} says
that $\left(  a+b\right)  ^{3}=a^{3}+b^{3}$. And indeed, we can show this
directly: Theorem \ref{thm.finfield.char.basics} \textbf{(b)} shows that
$3u=0$ for any $u\in F$ (for $p=3$), and the binomial formula yields%
\[
\left(  a+b\right)  ^{3}=a^{3}+\underbrace{3a^{2}b}%
_{\substack{=0\\\text{(since }3u=0\\\text{for any }u\in F\text{)}%
}}+\underbrace{3ab^{2}}_{\substack{=0\\\text{(since }3u=0\\\text{for any }u\in
F\text{)}}}+\,b^{3}=a^{3}+b^{3}.
\]
We note that Theorem \ref{thm.finfield.idiot} says nothing about powers other
than $p$-th or $p^{m}$-th powers. So it is not a replacement for the binomial formula!

To prove Theorem \ref{thm.finfield.idiot} \textbf{(a)} in general, we will
argue in the same way as in the $p=3$ example we just showed; we will just
need to know that all but the leftmost and rightmost addends in the binomial
formula vanish. This is a consequence of the following property of binomial coefficients:

\begin{lemma}
\label{lem.ent.pdivpchoosek}Let $p$ be a prime number. Let $k\in\left\{
1,2,\ldots,p-1\right\}  $. Then, $p\mid\dbinom{p}{k}$.
\end{lemma}

Note that this does indeed depend on $p$ being a prime. For example, the
number $4$ is not prime, and we \textbf{do not} have $4\mid\dbinom{4}{2}$
(since $\dbinom{4}{2}=6$).

\begin{proof}
[Proof of Lemma \ref{lem.ent.pdivpchoosek}.]There is an easy-to-prove formula
(see, e.g., \cite[Proposition 3.22]{detnotes}) saying that%
\[
\dbinom{p}{k}=\dfrac{p}{k}\cdot\dbinom{p-1}{k-1}.
\]
Hence,%
\[
k\dbinom{p}{k}=p\dbinom{p-1}{k-1}.
\]

Hence, $k\dbinom{p}{k}$ is divisible by $p$. But $k$ is coprime to $p$ (since
$p$ is prime), so we can cancel $k$ from this divisibility, and conclude that
$\dbinom{p}{k}$ is divisible by $p$. Lemma \ref{lem.ent.pdivpchoosek} is
proved. (See \cite[discussion of Exercise 9.1.6]{mps} for another proof.)
\end{proof}

\begin{proof}
[Proof of Theorem \ref{thm.finfield.idiot}.]\textbf{(a)} Let $a,b\in F$. Then,
$ab=ba$ (since $F$ is commutative); thus, the Binomial Formula yields%
\begin{equation}
\left(  a+b\right)  ^{p}=\sum_{k=0}^{p}\dbinom{p}{k}a^{k}b^{p-k}=a^{p}%
+\sum_{k=1}^{p-1}\dbinom{p}{k}a^{k}b^{p-k}+b^{p}.
\label{pf.thm.finfield.idiot.a.1}%
\end{equation}

Now, we claim that all the addends in the sum $\sum_{k=1}^{p-1}\dbinom{p}%
{k}a^{k}b^{p-k}$ vanish. Indeed, let $k\in\left\{  1,2,\ldots,p-1\right\}  $.
Then, Lemma \ref{lem.ent.pdivpchoosek} tells us that $\dbinom{p}{k}=mp$ for
some $m\in\mathbb{Z}$. Consider this $m$. Clearly, each $u\in F$ satisfies
$pu=0$ (indeed, $F$ is a $\mathbb{Z}/p$-algebra, so that
$pu=\underbrace{\overline{p}}_{=0\text{ in }\mathbb{Z}/p}u=0u=0$). Thus, each
$u\in F$ satisfies $\dbinom{p}{k}u=m\underbrace{pu}_{\substack{=0\\\text{(as
we just saw)}}}=m\cdot0=0$. Hence, in particular, we have%
\begin{equation}
\dbinom{p}{k}a^{k}b^{p-k}=0. \label{pf.thm.finfield.idiot.a.2}%
\end{equation}

Now, forget that we fixed $k$. We thus have shown that
(\ref{pf.thm.finfield.idiot.a.2}) holds for each $k\in\left\{  1,2,\ldots
,p-1\right\}  $. Hence, (\ref{pf.thm.finfield.idiot.a.1}) becomes%
\[
\left(  a+b\right)  ^{p}=a^{p}+\sum_{k=1}^{p-1}\underbrace{\dbinom{p}{k}%
a^{k}b^{p-k}}_{\substack{=0\\\text{(by (\ref{pf.thm.finfield.idiot.a.2}))}%
}}+\,b^{p}=a^{p}+b^{p}.
\]
This proves Theorem \ref{thm.finfield.idiot} \textbf{(a)}. \medskip

\textbf{(b)} This follows by induction on $m$ using Theorem
\ref{thm.finfield.idiot} \textbf{(a)}, since any $u\in F$ and any $m>0$
satisfy $u^{p^{m}}=\left(  u^{p^{m-1}}\right)  ^{p}$. \medskip

\textbf{(c)} Let $a,b\in F$. Applying Theorem \ref{thm.finfield.idiot}
\textbf{(a)} to $a-b$ instead of $a$, we get%
\[
\left(  \left(  a-b\right)  +b\right)  ^{p}=\left(  a-b\right)  ^{p}+b^{p}.
\]
Solving this for $\left(  a-b\right)  ^{p}$, we get%
\[
\left(  a-b\right)  ^{p}=\left(  \underbrace{\left(  a-b\right)  +b}%
_{=a}\right)  ^{p}-b^{p}=a^{p}-b^{p}.
\]
This proves Theorem \ref{thm.finfield.idiot} \textbf{(c)}. \medskip

\textbf{(d)} This follows by induction on $m$ using Theorem
\ref{thm.finfield.idiot} \textbf{(c)}.
\end{proof}

\begin{corollary}
\label{cor.finfield.frobenius}Let $p$ be a prime number. Let $F$ be a field of
characteristic $p$, or, more generally, any commutative $\mathbb{Z}%
/p$-algebra. Then, the map%
\begin{align*}
F  &  \rightarrow F,\\
a  &  \mapsto a^{p}%
\end{align*}
is a ring morphism.
\end{corollary}

\begin{proof}
Theorem \ref{thm.finfield.idiot} \textbf{(a)} says that this map respects
addition. But it is also clear that this map respects multiplication (since
$\left(  ab\right)  ^{p}=a^{p}b^{p}$ for any $a,b\in F$) and respects zero and
unity (since $0^{p}=0$ and $1^{p}=1$). Thus, it is a ring morphism.
\end{proof}

The ring morphism in Corollary \ref{cor.finfield.frobenius} is known as the
\textbf{Frobenius endomorphism} of $F$. It exists for arbitrary commutative
$\mathbb{Z}/p$-algebras, but it is particularly well-behaved for finite
fields. In particular, it is bijective when $F$ is a finite field:

\begin{exercise}
\label{exe.21hw4.2}Let $p$ be a prime number. Let $F$ be a finite field of
characteristic $p$. Let $f$ be the Frobenius endomorphism of $F$ (that is, the
map $F\rightarrow F,\ a\mapsto a^{p}$). Recall that $f$ is a ring morphism (by
Corollary \ref{cor.finfield.frobenius}).

\begin{enumerate}
\item[\textbf{(a)}] Prove that $f$ is a ring isomorphism from $F$ to $F$ (so
it is invertible).

\item[\textbf{(b)}] Now, replace the words \textquotedblleft field of
characteristic $p$\textquotedblright\ by the (more general) \textquotedblleft
commutative $\mathbb{Z}/p$-algebra\textquotedblright\ in the above. Find an
example where the claim of part \textbf{(a)} becomes false.
\end{enumerate}
\end{exercise}

We will use the Idiot's Binomial Formula to construct finite fields, but it
has many other applications. Here is just one:

\begin{exercise}
\label{exe.21hw4.7}Let $p$ be a prime number.

\begin{enumerate}
\item[\textbf{(a)}] Prove that $\left(  1+x\right)  ^{ap+c}=\left(
1+x^{p}\right)  ^{a}\left(  1+x\right)  ^{c}$ in the polynomial ring $\left(
\mathbb{Z}/p\right)  \left[  x\right]  $ for any $a,c\in\mathbb{N}$.

\item[\textbf{(b)}] Prove \textbf{Lucas's congruence}: Any $a,b\in\mathbb{N}$
and any $c,d\in\left\{  0,1,\ldots,p-1\right\}  $ satisfy
\[
\dbinom{ap+c}{bp+d}\equiv\dbinom{a}{b}\dbinom{c}{d}\operatorname{mod}p.
\]

\item[\textbf{(c)}] Explain the similarity between the modulo-$2$ Pascal's
triangle and Sierpinski's triangle (see Exercise \ref{exe.21hw0.9}).
\end{enumerate}
\end{exercise}

\begin{fineprint}
The following exercise is a converse to Lemma \ref{lem.ent.pdivpchoosek}:
\end{fineprint}

\begin{exercise}
Let $n>1$ be an integer that is not prime. Prove that there exists some
$k\in\left\{  1,2,\ldots,n-1\right\}  $ such that $n\nmid\dbinom{n}{k}$.
\medskip

[\textbf{Hint:} Let $k$ be a prime divisor of $n$.]
\end{exercise}

\subsubsection{\label{subsec.finfields.tools.deriv}The derivative of a
polynomial}

Our last tool is optional, as we will use it in one proof but avoid it in an
alternative proof of the same claim. Nevertheless, it is worth seeing, since
its usefulness is not limited to the cursory use we will put it to. This tool
is the notion of the \textquotedblleft formal\textquotedblright\ (i.e., purely
algebraic) \textbf{derivative} of a polynomial (which is defined over an
arbitrary commutative ring, not just over the real or complex numbers):

\begin{definition}
\label{def.polring.derivative}Let $R$ be a commutative ring. Let $f\in
R\left[  x\right]  $ be a polynomial. The \textbf{derivative} $f^{\prime}$ of
$f$ is a polynomial in $R\left[  x\right]  $ defined as follows: Writing $f$
in the form
\[
f=\sum\limits_{k\in\mathbb{N}}f_{k}x^{k}=f_{0}x^{0}+f_{1}x^{1}+f_{2}%
x^{2}+f_{3}x^{3}+\cdots
\]
for some $f_{0},f_{1},f_{2},\ldots\in R$, we set%
\[
f^{\prime}=\sum\limits_{k>0}f_{k}kx^{k-1}=f_{1}\cdot1x^{0}+f_{2}\cdot
2x^{1}+f_{3}\cdot3x^{2}+f_{4}\cdot4x^{3}+\cdots
\]

\end{definition}

For example, if $f=7x^{4}+2x+3$, then $f^{\prime}=7\cdot4x^{3}+2\cdot
1x^{0}=28x^{3}+2$ (where we have, of course, ignored zero coefficients).

Definition \ref{def.polring.derivative} is obviously inspired by the formula
for the derivative of a polynomial function in calculus. Unlike in calculus,
we are not wasting our time with little $\varepsilon$s and convergence issues;
instead, we are just defining $f^{\prime}$ using the explicit formula that
probably took you some time to prove back in calculus. There is no free lunch
here -- with this definition, you cannot re-use anything you have learned
about derivatives in your analysis classes (already because you are working in
a much more general setting now, with a commutative ring $R$ instead of the
real numbers); thus, a host of basic properties of derivatives need to be
proven before the notion becomes useful. In particular, the following needs to
be shown:

\begin{proposition}
\label{prop.polring.derivative.leibniz}Let $R$ be a commutative ring. Let
$f,g\in R\left[  x\right]  $. Then: \medskip

\begin{enumerate}
\item[\textbf{(a)}] We have $\left(  f+g\right)  ^{\prime}=f^{\prime
}+g^{\prime}$.

\item[\textbf{(b)}] We have $\left(  fg\right)  ^{\prime}=f^{\prime
}g+fg^{\prime}$. (This is called the \textbf{Leibniz rule}.)
\end{enumerate}
\end{proposition}

\begin{exercise}
\label{exe.21hw3.7bc}Prove Proposition \ref{prop.polring.derivative.leibniz}.
\medskip

[\textbf{Hint:} For part \textbf{(b)}, it is easiest to first prove it in the
particular case when $f=x^{a}$ and $g=x^{b}$ for some $a$ and $b$, and then
obtain the general case by interchanging summations.]
\end{exercise}

The following corollary is an algebraic analogue of the well-known fact
\textquotedblleft a double root of a polynomial is a root of its
derivative\textquotedblright:

\begin{corollary}
\label{cor.polring.derivative.doubleroot}Let $R$ be a commutative ring. Let
$f\in R\left[  x\right]  $ and $r\in R$. If $\left(  x-r\right)  ^{2}\mid f$
in $R\left[  x\right]  $, then $x-r\mid f^{\prime}$ in $R\left[  x\right]  $.
\end{corollary}

\begin{proof}
Assume that $\left(  x-r\right)  ^{2}\mid f$. Thus, we can write $f$ as
$f=\left(  x-r\right)  ^{2}g$ for some $g\in R\left[  x\right]  $. Consider
this $g$. From $f=\left(  x-r\right)  ^{2}g$, we obtain%
\begin{align*}
f^{\prime}  &  =\left(  \left(  x-r\right)  ^{2}g\right)  ^{\prime
}=\underbrace{\left(  \left(  x-r\right)  ^{2}\right)  ^{\prime}%
}_{\substack{=2\left(  x-r\right)  \\\text{(this is easy to}\\\text{check
directly)}}}g+\left(  x-r\right)  ^{2}g^{\prime}\ \ \ \ \ \ \ \ \ \ \left(
\text{by the Leibniz rule}\right) \\
&  =2\left(  x-r\right)  g+\left(  x-r\right)  ^{2}g^{\prime}=\left(
x-r\right)  \left(  2g+\left(  x-r\right)  g^{\prime}\right)  .
\end{align*}
Thus, $x-r\mid f^{\prime}$, so that Corollary
\ref{cor.polring.derivative.doubleroot} is proven.
\end{proof}

Corollary \ref{cor.polring.derivative.doubleroot} can be applied (in its
contrapositive) as a sufficient criterion for a polynomial to have distinct
roots. This is exactly how we will soon apply it. We note that Corollary
\ref{cor.polring.derivative.doubleroot} also has a converse:

\begin{exercise}
Let $R$ be a commutative ring. Let $f\in R\left[  x\right]  $ and $r\in R$. If
$r$ is a root of both $f$ and $f^{\prime}$, then prove that $\left(
x-r\right)  ^{2}\mid f$ in $R\left[  x\right]  $.
\end{exercise}

Here are some further properties of derivatives of polynomials:

\begin{exercise}
Let $R$ be a commutative ring. Let $f\in R\left[  x\right]  $ be any
polynomial. Prove the following:

\begin{enumerate}
\item[\textbf{(a)}] We have $\deg\left(  f^{\prime}\right)  \leq\deg f-1$
(where we understand that $-\infty-1=-\infty$).

\item[\textbf{(b)}] If $R$ is a $\mathbb{Q}$-algebra and $f$ is not constant,
then $\deg\left(  f^{\prime}\right)  =\deg f-1$.
\end{enumerate}
\end{exercise}

\begin{exercise}
\ \ 

\begin{enumerate}
\item[\textbf{(a)}] Let $R$ be a commutative ring. Let $f,g\in R\left[
x\right]  $ be any two polynomials. Prove that $\left(  f\left[  g\right]
\right)  ^{\prime}=f^{\prime}\left[  g\right]  \cdot g^{\prime}$. (This is
called the \textbf{chain rule for polynomials}, and is the algebraic analogue
of the chain rule in calculus.)

\item[\textbf{(b)}] Use this chain rule to give a different proof of the
equality $\left(  \left(  x-r\right)  ^{2}\right)  ^{\prime}=2\left(
x-r\right)  $ in the above proof of Corollary
\ref{cor.polring.derivative.doubleroot}.
\end{enumerate}
\end{exercise}

\begin{exercise}
Let $R$ be a commutative ring. Let $f\in R\left[  x\right]  $ be a polynomial.
Consider also the polynomial ring $R\left[  x,y\right]  $ in two
indeterminates $x$ and $y$.

\begin{enumerate}
\item[\textbf{(a)}] Prove that there is a unique polynomial $g\in R\left[
x,y\right]  $ satisfying $f\left[  y\right]  -f\left[  x\right]  =\left(
y-x\right)  \cdot g$.

\item[\textbf{(b)}] Prove that this polynomial $g$ satisfies $g\left[
x,x\right]  =f^{\prime}$.
\end{enumerate}

(Because of its definition, the polynomial $g$ can be written as
$\dfrac{f\left[  y\right]  -f\left[  x\right]  }{y-x}$, even though $y-x$ is
not a unit of $R\left[  x,y\right]  $. However, when computing $g\left[
x,x\right]  $, we cannot just set $y$ to $x$ in the fraction $\dfrac{f\left[
y\right]  -f\left[  x\right]  }{y-x}$; instead, we must first expand this
fraction into a polynomial. Thus, the claim of part \textbf{(b)} is the
algebraic analogue of the formula $f^{\prime}\left(  x\right)  =\lim
\limits_{y\rightarrow x}\dfrac{f\left(  y\right)  -f\left(  x\right)  }{y-x}$
from calculus.)
\end{exercise}

\begin{exercise}
\label{exe.dualnums.deriv}Let $R$ be a commutative ring. Recall the
$R$-algebra $\mathbb{D}_{R}$ of $R$-dual numbers from Exercise
\ref{exe.algs.dualnums}. Let $\varepsilon$ be the element $\left(  0,1\right)
$ of $\mathbb{D}_{R}$. Let $f\in R\left[  x\right]  $ be an arbitrary
polynomial. Let $a,b\in R$. Prove that%
\[
f\left[  \left(  a,b\right)  \right]  =\left(  f\left[  a\right]
,\ bf^{\prime}\left[  a\right]  \right)  \ \ \ \ \ \ \ \ \ \ \text{in
}\mathbb{D}_{R}.
\]
(In other words, prove that $f\left[  a+b\varepsilon\right]  =f\left[
a\right]  +bf^{\prime}\left[  a\right]  \varepsilon$, where we are identifying
each element $r\in R$ with the $R$-dual number $\left(  r,0\right)
\in\mathbb{D}_{R}$.)
\end{exercise}

\begin{exercise}
\label{exe.21hw3.7bdef}Let $R$ be a commutative ring.

Let $D:R\left[  x\right]  \rightarrow R\left[  x\right]  $ be the map sending
each polynomial $f$ to its derivative $f^{\prime}$. We refer to $D$ as
\textbf{(formal) differentiation}. As usual, for any $n\in\mathbb{N}$, we let
$D^{n}$ denote $\underbrace{D\circ D\circ\cdots\circ D}_{n\text{ times}}$
(which means $\operatorname{id}$ if $n=0$).

Prove the following:

\begin{enumerate}
\item[\textbf{(a)}] The map $D:R\left[  x\right]  \rightarrow R\left[
x\right]  $ is $R$-linear.

\item[\textbf{(b)}] We have $D^{n}\left(  x^{k}\right)  =n!\dbinom{k}%
{n}x^{k-n}$ for all $n\in\mathbb{N}$ and $k\in\mathbb{N}$. Here, the
expression \textquotedblleft$\dbinom{k}{n}x^{k-n}$\textquotedblright\ is to be
understood as $0$ when $k<n$.

\item[\textbf{(c)}] If $\mathbb{Q}$ is a subring of $R$, then every polynomial
$f\in R\left[  x\right]  $ satisfies
\[
f\left[  x+a\right]  =\sum_{n\in\mathbb{N}}\dfrac{1}{n!}\left(  D^{n}\left(
f\right)  \right)  \left[  a\right]  \cdot x^{n}\qquad\text{for all $a\in R$%
}.
\]
(The infinite sum on the right hand side has only finitely many nonzero addends.)

\item[\textbf{(d)}] If $p$ is a prime such that $p\cdot1_{R}=0$ (for example,
this happens if $R=\mathbb{Z}/p$), then $D^{p}\left(  f\right)  =0$ for each
$f\in R\left[  x\right]  $.
\end{enumerate}
\end{exercise}

\begin{exercise}
\label{exe.21hw3.7ghi}Let $R$ be a commutative ring such that $\mathbb{Q}$ is
a subring of $R$. For each polynomial
\[
f=\sum_{k\in\mathbb{N}}f_{k}x^{k}=f_{0}x^{0}+f_{1}x^{1}+f_{2}x^{2}+\cdots\in
R\left[  x\right]  \ \ \ \ \ \ \ \ \ \ \left(  \text{where $f_{i}\in R$%
}\right)  \text{,}%
\]
we define the \textbf{(formal) integral} $\int f$ of $f$ to be the polynomial
\[
\sum_{k\in\mathbb{N}}\dfrac{1}{k+1}f_{k}x^{k+1}=\dfrac{1}{1}f_{0}x^{1}%
+\dfrac{1}{2}f_{1}x^{2}+\dfrac{1}{3}f_{2}x^{3}+\cdots\in R\left[  x\right]  .
\]
(This definition imitates the standard procedure for integrating power series
in analysis, but again works for any commutative ring $R$ that contains
$\mathbb{Q}$ as subring.)

Let $J:R\left[  x\right]  \rightarrow R\left[  x\right]  $ be the map sending
each polynomial $f$ to its integral $\int f$. Prove the following:

\begin{enumerate}
\item[\textbf{(a)}] The map $J:R\left[  x\right]  \rightarrow R\left[
x\right]  $ is $R$-linear.

\item[\textbf{(b)}] We have $D\circ J=\operatorname{id}$ (where $D$ is as in
Exercise \ref{exe.21hw3.7bdef}).

\item[\textbf{(c)}] We have $J\circ D\neq\operatorname{id}$.
\end{enumerate}
\end{exercise}

\subsection{\label{sec.finfields.exist}Existence of finite fields}

Now we are in walking distance of the existence of fields of size $p^{m}$:

\begin{theorem}
\label{thm.finfield.exist}Let $p$ be a prime number. Let $m$ be a positive
integer. Then, there exists a finite field of size $p^{m}$.
\end{theorem}

\begin{proof}
From $p>1$ and $m>0$, we obtain $p^{m}>1$. Hence, the polynomial $x^{p^{m}}-x$
is monic. Thus, by Theorem \ref{thm.finfield.splitfield} \textbf{(c)}, there
exists a splitting field of this polynomial over $\mathbb{Z}/p$. Let $S$ be
such a splitting field. Thus, the polynomial $x^{p^{m}}-x$ splits over $S$. In
other words, there exist elements $r_{1},r_{2},\ldots,r_{p^{m}}$ of $S$ such
that%
\begin{equation}
x^{p^{m}}-x=\left(  x-r_{1}\right)  \left(  x-r_{2}\right)  \cdots\left(
x-r_{p^{m}}\right)  . \label{pf.thm.finfield.exist.1}%
\end{equation}
Consider these $r_{1},r_{2},\ldots,r_{p^{m}}$.

Let%
\[
L=\left\{  r_{1},r_{2},\ldots,r_{p^{m}}\right\}  .
\]

Our goal will be to show that $L$ is a finite field of size $p^{m}$.

\textbf{Everything} in this statement needs proof!\footnote{Except for the
\textquotedblleft finite\textquotedblright\ part, which is obvious but not
overly helpful by itself.} Even the size is not obvious, let alone that $L$ is
a field. Let us start with the size:

\begin{statement}
\textit{Claim 1:} We have $\left\vert L\right\vert =p^{m}$.
\end{statement}

Let us give two proofs of Claim 1:

[\textit{First proof of Claim 1:} This amounts to showing that $r_{1}%
,r_{2},\ldots,r_{p^{m}}$ are distinct (since this will immediately yield that
$L=\left\{  r_{1},r_{2},\ldots,r_{p^{m}}\right\}  $ is a $p^{m}$-element set).
Let us thus do this. Indeed, assume the contrary. Then, $r_{i}=r_{j}$ for some
$i<j$. Hence, the $x-r_{i}$ and $x-r_{j}$ factors on the right hand side of
(\ref{pf.thm.finfield.exist.1}) are identical. Thus, $x-r_{i}$ appears twice
as a factor on this right hand side; consequently,
(\ref{pf.thm.finfield.exist.1}) entails that $\left(  x-r_{i}\right)  ^{2}\mid
x^{p^{m}}-x$. Hence, Corollary \ref{cor.polring.derivative.doubleroot}
(applied to $R=S$ and $f=x^{p^{m}}-x$ and $r=r_{i}$) yields $x-r_{i}%
\mid\left(  x^{p^{m}}-x\right)  ^{\prime}$. But Definition
\ref{def.polring.derivative} yields%
\[
\left(  x^{p^{m}}-x\right)  ^{\prime}=\underbrace{p^{m}x^{p^{m}-1}%
}_{\substack{=0\\\text{(since }pu=0\\\text{for any }u\in S\text{)}}}-\,1=-1.
\]
Thus, $x-r_{i}\mid\left(  x^{p^{m}}-x\right)  ^{\prime}=-1\mid1$. But it is
impossible for the degree-$1$ polynomial $x-r_{i}$ to divide the degree-$0$
polynomial $1$ (for degree reasons). So we have found a contradiction.]

[\textit{Second proof of Claim 1:} The following proof of Claim 1 (which I
have learnt from \cite[Chapter 5, Theorem 3.3]{Shifri96}) avoids the use of derivatives.

Again, it suffices to show that $r_{1},r_{2},\ldots,r_{p^{m}}$ are distinct.
Again, assume the contrary. Then, $r_{i}=r_{j}$ for some $i<j$. As before, we
then find that $\left(  x-r_{i}\right)  ^{2}\mid x^{p^{m}}-x$. In other words,
$x^{p^{m}}-x=\left(  x-r_{i}\right)  ^{2}\cdot g$ for some polynomial $g\in
S\left[  x\right]  $. Consider this $g$. Substituting $r_{i}$ for $x$ in the
equality $x^{p^{m}}-x=\left(  x-r_{i}\right)  ^{2}\cdot g$, we obtain
$r_{i}^{p^{m}}-r_{i}=\underbrace{\left(  r_{i}-r_{i}\right)  ^{2}}_{=0}%
\cdot\,g\left[  r_{i}\right]  =0$. In other words, $r_{i}^{p^{m}}=r_{i}$.

Substituting $x+r_{i}$ for $x$ in the equality $x^{p^{m}}-x=\left(
x-r_{i}\right)  ^{2}\cdot g$, we obtain%
\[
\left(  x+r_{i}\right)  ^{p^{m}}-\left(  x+r_{i}\right)  =\left(
\underbrace{x+r_{i}-r_{i}}_{=x}\right)  ^{2}\cdot g\left[  x+r_{i}\right]
=x^{2}\cdot g\left[  x+r_{i}\right]  .
\]
In view of%
\[
\underbrace{\left(  x+r_{i}\right)  ^{p^{m}}}_{\substack{=x^{p^{m}}%
+r_{i}^{p^{m}}\\\text{(by Theorem \ref{thm.finfield.idiot} \textbf{(b)})}%
}}-\left(  x+r_{i}\right)  =x^{p^{m}}+\underbrace{r_{i}^{p^{m}}}_{=r_{i}%
}-\left(  x+r_{i}\right)  =x^{p^{m}}+r_{i}-\left(  x+r_{i}\right)  =x^{p^{m}%
}-x,
\]
we can rewrite this as $x^{p^{m}}-x=x^{2}\cdot g\left[  x+r_{i}\right]  $.
Thus, $x^{2}\mid x^{p^{m}}-x$. But this is absurd, since the $x^{1}%
$-coefficient of $x^{p^{m}}-x$ is nonzero. Thus, we have found a contradiction
again.] \medskip

Next, let us characterize $L$ somewhat differently:

\begin{statement}
\textit{Claim 2:} We have%
\begin{align*}
L  &  =\left\{  \text{the roots of }x^{p^{m}}-x\text{ in }S\right\}  =\left\{
a\in S\ \mid\ a^{p^{m}}-a=0\right\} \\
&  =\left\{  a\in S\ \mid\ a^{p^{m}}=a\right\}  .
\end{align*}

\end{statement}

[\textit{Proof of Claim 2:} The equation (\ref{pf.thm.finfield.exist.1})
yields that%
\begin{align*}
&  \left\{  \text{the roots of }x^{p^{m}}-x\text{ in }S\right\} \\
&  =\left\{  \text{the roots of }\left(  x-r_{1}\right)  \left(
x-r_{2}\right)  \cdots\left(  x-r_{p^{m}}\right)  \text{ in }S\right\} \\
&  =\left\{  r_{1},r_{2},\ldots,r_{p^{m}}\right\}  \ \ \ \ \ \ \ \ \ \ \left(
\text{by Proposition \ref{prop.finfield.splits.roots}}\right) \\
&  =L.
\end{align*}
Hence,%
\begin{align*}
L  &  =\left\{  \text{the roots of }x^{p^{m}}-x\text{ in }S\right\}  =\left\{
a\in S\ \mid\ a^{p^{m}}-a=0\right\} \\
&  =\left\{  a\in S\ \mid\ a^{p^{m}}=a\right\}  .
\end{align*}
This proves Claim 2.]

Now, why is $L$ a field? First, let us check that $L$ is a ring:

\begin{statement}
\textit{Claim 3:} The set $L$ is a subring of $S$.
\end{statement}

[\textit{Proof of Claim 3:} Claim 2 yields that
\begin{equation}
L=\left\{  a\in S\ \mid\ a^{p^{m}}=a\right\}  .
\label{pf.thm.finfield.exist.c3.pf.1}%
\end{equation}
Hence, $0\in L$ (since $0^{p^{m}}=0$) and $1\in L$ (since $1^{p^{m}}=1$).
Furthermore, I claim that $L$ is closed under addition. Indeed, if $a,b\in L$,
then (\ref{pf.thm.finfield.exist.c3.pf.1}) yields $a^{p^{m}}=a$ and $b^{p^{m}%
}=b$, so that%
\begin{align*}
\left(  a+b\right)  ^{p^{m}}  &  =\underbrace{a^{p^{m}}}_{=a}%
+\underbrace{b^{p^{m}}}_{=b}\ \ \ \ \ \ \ \ \ \ \left(  \text{by Theorem
\ref{thm.finfield.idiot} \textbf{(b)}}\right) \\
&  =a+b,
\end{align*}
and this means $a+b\in L$ (again because of
(\ref{pf.thm.finfield.exist.c3.pf.1})). This shows that $L$ is closed under
addition. For a similar reason, $L$ is closed under subtraction\footnote{Use
Theorem \ref{thm.finfield.idiot} \textbf{(d)} instead of Theorem
\ref{thm.finfield.idiot} \textbf{(b)} here.}, so that $L$ is closed under
negation. Finally, $L$ is closed under multiplication, since $\left(
ab\right)  ^{p^{m}}=a^{p^{m}}b^{p^{m}}$ for any $a,b\in L$. Hence, $L$ is a
subring of $S$.]

Thus, in particular, $L$ is a commutative ring (since $S$ is a field, thus a
commutative ring). Now, let us see that $L$ is a field:

\begin{statement}
\textit{Claim 4:} The ring $L$ is a field.
\end{statement}

[\textit{Proof of Claim 4:} We know that $S$ is a field, so that $0\neq1$ in
$S$, and this of course means that $0\neq1$ in $L$. It thus remains to show
that every nonzero element of $L$ is a unit.

Let $a\in L$ be nonzero. Then, $a$ has an inverse in $S$, since $S$ is a
field. This inverse $a^{-1}$ satisfies $\left(  a^{-1}\right)  ^{p^{m}%
}=\left(  a^{p^{m}}\right)  ^{-1}$ (indeed, this is a particular case of the
identity $\left(  g^{-1}\right)  ^{k}=\left(  g^{k}\right)  ^{-1}$, which
holds whenever $g$ is an element of a group and $k$ is an integer). But $a\in
L$ and thus $a^{p^{m}}=a$ (by (\ref{pf.thm.finfield.exist.c3.pf.1})). Hence,
\[
\left(  a^{-1}\right)  ^{p^{m}}=\left(  \underbrace{a^{p^{m}}}_{=a}\right)
^{-1}=a^{-1},
\]
so that $a^{-1}\in L$ (by (\ref{pf.thm.finfield.exist.c3.pf.1}) again). Thus,
$a$ has an inverse in $L$; in other words, $a$ is a unit of $L$.

Thus, we have shown that every nonzero element of $L$ is a unit. As we said,
this finishes the proof of Claim 4.]

Combining Claims 1 and 4, we conclude that $L$ is a field of size $p^{m}$.
Thus, such a field exists. This proves Theorem \ref{thm.finfield.exist}.
\end{proof}

So we are done with the first deep result of this course! But some further
questions suggest themselves:

\begin{itemize}
\item We have obtained $L$ rather indirectly: First, we took a splitting field
$S$ of the huge polynomial $x^{p^{m}}-x$; then we carved $L$ out of it by
taking the roots of this polynomial. Could we get $L$ more directly? For
example, if there is an irreducible polynomial $f$ of degree $m$ over
$\mathbb{Z}/p$, then we can just take the field $\left(  \mathbb{Z}/p\right)
\left[  x\right]  /f$. Is there such an $f$ ?

\item Can there be several non-isomorphic fields of size $p^{m}$ (for fixed
$p$ and $m$)? For example, can there be two non-isomorphic fields of size
$p^{2}$ ? It is not hard to see that any field of size $p^{2}$ can be obtained
(up to isomorphism) by adjoining a root of an irreducible quadratic polynomial
to $\mathbb{Z}/p$; thus, the question is whether different such polynomials
can lead to different fields.

If we were working with infinite fields, examples of such behavior would be
easy to find. For example, adjoining a root of $x^{2}-2$ to $\mathbb{Q}$
yields the field $\mathbb{Q}\left[  \sqrt{2}\right]  $, whereas adjoining a
root of $x^{2}-3$ to $\mathbb{Q}$ yields the field $\mathbb{Q}\left[  \sqrt
{3}\right]  $. It is not hard to see that $\mathbb{Q}\left[  \sqrt{2}\right]
$ is not isomorphic to $\mathbb{Q}\left[  \sqrt{3}\right]  $ (for example, you
can show that $2$ is a square in $\mathbb{Q}\left[  \sqrt{2}\right]  $ but not
in $\mathbb{Q}\left[  \sqrt{3}\right]  $). Can this happen with $\mathbb{Z}/p$
instead of $\mathbb{Q}$ ?
\end{itemize}

These questions will be answered in the next section.

\subsection{Uniqueness of finite fields}

Theorem \ref{thm.finfield.exist} shows that for any prime power $p^{m}$, there
exists a finite field of size $p^{m}$. The next natural question is: How many
such fields are there? Literally, of course, there are infinitely many, since
each one has infinitely many isomorphic (but not literally identical) copies.
Of course, the right question to ask is how many non-isomorphic finite fields
there are of a given size.

The answer is surprisingly simple: There is only one. Proving this will take
us some work. (This section is more abstract and notationally dense than many
others, and can be skipped.)

\subsubsection{Annihilating polynomials and minimal polynomials}

We begin with two fundamental concepts of Galois theory: the notions of
\textquotedblleft annihilating polynomials\textquotedblright\ and of
\textquotedblleft minimal polynomials\textquotedblright. We will not delve
deeper than this into Galois theory, but we will explore these notions in some detail.

Roughly speaking, the main problem of classical algebra is solving polynomial
equations: Given a polynomial $f$, what are its roots? Our abstract viewpoint
has allowed us to extend the field over which the polynomial is defined, and
in such an extension we can always find a root for any non-constant univariate
polynomial (see Theorem \ref{thm.finfield.splitfield} \textbf{(b)}).

We now turn the question around: Given an element $a$ of a field $F$, what are
the polynomials $f$ that have $a$ as a root? If we want $f$ to belong to
$F\left[  x\right]  $, then this is an easy question, and the answer is
\textquotedblleft exactly those polynomials that are divisible by
$x-a$\textquotedblright\ (see Proposition \ref{prop.polring.univar-rootdiv}).
But this question becomes more interesting if we require $f$ to belong to
$S\left[  x\right]  $, where $S$ is a smaller field than $F$. For instance,
$a$ could be the imaginary unit $i=\sqrt{-1}\in\mathbb{C}$, and $S$ could be
the field $\mathbb{R}$, so we would be asking about the real polynomials that
have $i$ as a root. Clearly, $x^{2}+1$ is one of these, and thus any multiple
of $x^{2}+1$ qualifies as well. Are there any others?

To study this kind of question, we introduce the notions of annihilating and
minimal polynomials:\footnote{Mnemonic: The letters $S$ and $F$ refer to
\textquotedblleft subfield\textquotedblright\ and \textquotedblleft
field\textquotedblright.}

\begin{definition}
\label{def.minpol.def}Let $S$ and $F$ be two fields such that $S$ is a subring
of $F$. (For example, we can take $S=\mathbb{Q}$ and $F=\mathbb{R}$.)

Let $a\in F$ be an arbitrary element.

\begin{enumerate}
\item[\textbf{(a)}] An \textbf{annihilating polynomial} of \emph{$a$} shall
mean a polynomial $f\in S\left[  x\right]  $ such that $a$ is a root of $f$.
For instance:

\begin{itemize}
\item If $a\in S$, then $x-a$ is an annihilating polynomial of $a$.

\item If $a$ is a square root of an element $v\in S$, then $x^{2}-v$ is an
annihilating polynomial of $a$.

\item If $S=\mathbb{Q}$ and $F=\mathbb{R}$, then $x^{4}-10x^{2}+1$ is an
annihilating polynomial of $\sqrt{2}+\sqrt{3}$.

\item The real number $\pi$ is known to be transcendental; this means that
there exists no nonzero annihilating polynomial of $\pi$ (for $S=\mathbb{Q}$
and $F=\mathbb{R}$).
\end{itemize}

\item[\textbf{(b)}] The \textbf{minimal polynomial} of $a$ (over the subfield
$S$) is defined to be the monic annihilating polynomial of $a$ of smallest
possible degree (if such a polynomial exists).\footnotemark
\end{enumerate}
\end{definition}

\footnotetext{Thus, the minimal polynomial of $a$ is defined in the exact same
way as the minimal polynomial of a square matrix was defined in linear
algebra. However, the minimal polynomial of a matrix is not always
irreducible, whereas in our case the minimal polynomial will be irreducible
(see below).}It is not yet clear that the minimal polynomial of $a$ really
deserves the definite article! Couldn't there be several monic annihilating
polynomials of $a$ of smallest possible degree? Which of them deserves to be
called \textquotedblleft the\textquotedblright\ minimal polynomial?

Fortunately, this question never arises, since there is only one:

\begin{theorem}
\label{thm.minpol.basics}Let $S$ and $F$ be two fields such that $S$ is a
subring of $F$. Let $a\in F$. Then:

\begin{enumerate}
\item[\textbf{(a)}] The minimal polynomial of $a$ is unique if it exists. That
is, if there is at least one monic annihilating polynomial of $a$, then only
one of these polynomials has smallest possible degree.

\item[\textbf{(b)}] The minimal polynomial of $a$ is always irreducible if it exists.

\item[\textbf{(c)}] If $f$ is the minimal polynomial of $a$, then the map
\begin{align*}
S\left[  x\right]  /f  &  \rightarrow F,\\
\overline{g}  &  \mapsto g\left[  a\right]
\end{align*}
is a (well-defined) $S$-algebra morphism, and is injective.

\item[\textbf{(d)}] Assume that $a$ has a minimal polynomial. Let $f$ be the
minimal polynomial of $a$. Then, the annihilating polynomials of $a$ are
precisely the polynomials $g\in S\left[  x\right]  $ that are divisible by $f$.
\end{enumerate}
\end{theorem}

Before we prove this, let us show a lemma:

\begin{lemma}
\label{lem.minpol.div}Let $S$ and $F$ be two fields such that $S$ is a subring
of $F$. Let $a\in F$.

Let $f\in S\left[  x\right]  $ be a minimal polynomial of $a$. (We could say
\textquotedblleft\textbf{the} minimal polynomial of $a$\textquotedblright\ if
we knew that it is unique, but we don't know this yet; we will prove this soon.)

Let $g\in S\left[  x\right]  $ be any polynomial. Then:

\begin{enumerate}
\item[\textbf{(a)}] If $g\left[  a\right]  =0$, then $f\mid g$.

\item[\textbf{(b)}] If $f\mid g$, then $g\left[  a\right]  =0$.
\end{enumerate}
\end{lemma}

\begin{proof}
[Proof of Lemma \ref{lem.minpol.div}.]We know that $f$ is a minimal polynomial
of $a$. Thus, $f$ is a monic annihilating polynomial of $a$. Hence, $a$ is a
root of $f$; in other words, $f\left[  a\right]  =0$. Also, the leading
coefficient of $f$ is $1$ (since $f$ is monic), and thus is a unit of $S$.
Hence, Theorem \ref{thm.polring.univar-quorem} \textbf{(a)} (applied to $S$,
$f$ and $g$ instead of $R$, $b$ and $a$) yields that there is a unique pair
$\left(  q,r\right)  $ of polynomials in $S\left[  x\right]  $ such that
\[
g=qf+r\ \ \ \ \ \ \ \ \ \ \text{and}\ \ \ \ \ \ \ \ \ \ \deg r<\deg f.
\]
Consider this pair $\left(  q,r\right)  $. \medskip

\textbf{(a)} Assume that $g\left[  a\right]  =0$. Thus,%
\begin{align*}
0  &  =g\left[  a\right]  =\left(  qf+r\right)  \left[  a\right]
\ \ \ \ \ \ \ \ \ \ \left(  \text{since }g=qf+r\right) \\
&  =q\left[  a\right]  \cdot\underbrace{f\left[  a\right]  }_{=0}+\,r\left[
a\right]  =r\left[  a\right]  .
\end{align*}
Hence, $r\left[  a\right]  =0$.

We assume (for the sake of contradiction) that $r\neq0$. Hence, $r$ has a
leading coefficient $\lambda$. Consider this $\lambda$. This $\lambda$ is
nonzero, and thus is a unit of $S$ (since $S$ is a field). Hence, the
polynomial $\lambda^{-1}r$ is well-defined. Moreover, since $\lambda^{-1}$ is
a nonzero constant, we have $\deg\left(  \lambda^{-1}r\right)  =\deg r<\deg f$.

However, $\left(  \lambda^{-1}r\right)  \left[  a\right]  =\lambda^{-1}%
\cdot\underbrace{r\left[  a\right]  }_{=0}=0$. In other words, $a$ is a root
of $\lambda^{-1}r$. Thus, $\lambda^{-1}r$ is an annihilating polynomial of
$a$. Furthermore, this polynomial $\lambda^{-1}r$ is monic (since $\lambda$ is
the leading coefficient of $r$, so that scaling $r$ by $\lambda^{-1}$ turns
the leading coefficient into $1$).

Now, recall that $f$ is a minimal polynomial of $a$. In other words, $f$ is a
monic annihilating polynomial of $a$ of smallest possible degree. Therefore,
$\deg\left(  \lambda^{-1}r\right)  \geq\deg f$ (since $\lambda^{-1}r$, too, is
a monic annihilating polynomial of $a$). This contradicts $\deg\left(
\lambda^{-1}r\right)  <\deg f$. This contradiction shows that our assumption
(that $r\neq0$) was false. Hence, $r=0$.

Now, $g=qf+\underbrace{r}_{=0}=qf$, so that $f\mid g$. This proves Lemma
\ref{lem.minpol.div} \textbf{(a)}. \medskip

\textbf{(b)} Assume that $f\mid g$. Thus, $g=fq$ for some $q\in S\left[
x\right]  $. Consider this $q$. From $g=fq$, we obtain $g\left[  a\right]
=\left(  fq\right)  \left[  a\right]  =\underbrace{f\left[  a\right]  }%
_{=0}\cdot\,q\left[  a\right]  =0$. Thus, Lemma \ref{lem.minpol.div}
\textbf{(b)} is proved.
\end{proof}

\begin{proof}
[Proof of Theorem \ref{thm.minpol.basics}.]\textbf{(a)} Assume that $a$ has a
minimal polynomial. We must show that the minimal polynomial of $a$ is unique.

Assume the contrary. Thus, there exist two distinct minimal polynomials $f$
and $g$ of $a$. Consider these $f$ and $g$. Both $f$ and $g$ are minimal
polynomials of $a$, and thus are monic annihilating polynomials of $a$. In
particular, we have $g\left[  a\right]  =0$ (since $g$ is an annihilating
polynomial of $a$). Hence, Lemma \ref{lem.minpol.div} \textbf{(a)} yields
$f\mid g$. The same argument (but with the roles of $f$ and $g$ interchanged)
yields $g\mid f$.

We have $f\mid g$. In other words, there exists a polynomial $q\in S\left[
x\right]  $ such that $g=fq$. Consider this $q$. We have $fq=g\neq0$ (since
$g$ is monic) and thus $q\neq0$. Also, $f\neq0$ (since $f$ is monic). But $F$
is a field and thus an integral domain. Hence, Proposition
\ref{prop.polring.univar-degpq} \textbf{(c)} yields $\deg\left(  fq\right)
=\deg f+\deg q$. In view of $g=fq$, this rewrites as $\deg g=\deg f+\deg q$.
Hence, $\deg g\geq\deg f$ (since $q\neq0$ entails $\deg q\geq0$). Similarly,
$\deg f\geq\deg g$. Combining these two inequalities, we find $\deg f=\deg g$.

Thus, $\deg f=\deg g=\deg f+\deg q$, so that $\deg q=0$ (since $f\neq0$
entails $\deg f\geq0$). In other words, $q$ is a nonzero constant.

The leading coefficient of $f$ is $1$ (since $f$ is monic). Thus, the leading
coefficient of $fq$ is $q$ (since $q$ is a nonzero constant). In view of
$fq=g$, this rewrites as follows: The leading coefficient of $g$ is $q$.
However, the leading coefficient of $g$ is $1$ (since $g$ is monic). Comparing
the previous two sentences, we conclude that $q=1$. Hence, $g=f\underbrace{q}%
_{=1}=f$. This contradicts our assumption that $f$ and $g$ are distinct. This
contradiction shows that our assumption was wrong. Theorem
\ref{thm.minpol.basics} \textbf{(a)} is thus proved. \medskip

\textbf{(b)} Assume that $a$ has a minimal polynomial. Let $f$ be the minimal
polynomial of $a$. We must show that $f$ is irreducible.

Recall that $f$ is the minimal polynomial of $a$. In other words, $f$ is a
monic annihilating polynomial of $a$ of smallest possible degree. Thus, $a$ is
a root of $f$ (since $f$ is an annihilating polynomial of $a$); in other
words, $f\left[  a\right]  =0$. Thus, $f$ cannot be a nonzero constant
(because this would entail $f\left[  a\right]  =f\neq0$, contradicting
$f\left[  a\right]  =0$). Hence, $f$ is not a unit of the ring $S\left[
x\right]  $. (This is something that needs to be checked if you want to show
that $f$ is irreducible. Don't forget about this!)

Now, let $u,v\in S\left[  x\right]  $ be such that $uv=f$. We shall show that
at least one of $u$ and $v$ is a unit of $S\left[  x\right]  $.

Indeed, assume the contrary. Thus, neither $u$ nor $v$ is a unit of $S\left[
x\right]  $. Moreover, neither $u$ nor $v$ equals $0$ (since $uv=f\neq0$).
Hence, neither $u$ nor $v$ is constant (since a constant polynomial is either
a unit of $S\left[  x\right]  $ or equals $0$). Thus, $\deg u\geq1$ and $\deg
v\geq1$. However, $f=uv$ and thus $\deg f=\deg\left(  uv\right)  =\deg u+\deg
v$ (since $S$ is a field and thus an integral domain). Hence, $\deg f=\deg
u+\underbrace{\deg v}_{\geq1>0}>\deg u$, so that $\deg u<\deg f$.

Now, $f=uv$, so that $f\left[  a\right]  =\left(  uv\right)  \left[  a\right]
=u\left[  a\right]  \cdot v\left[  a\right]  $ and therefore $u\left[
a\right]  \cdot v\left[  a\right]  =f\left[  a\right]  =0$. Since $F$ is a
field and thus an integral domain, this entails that $u\left[  a\right]  =0$
or $v\left[  a\right]  =0$. We WLOG assume that $u\left[  a\right]  =0$ (since
otherwise, we can simply swap $u$ with $v$). Let $\lambda$ denote the leading
coefficient of $u$ (this is well-defined, since $u$ does not equal $0$). Then,
the polynomial $\lambda^{-1}u$ is monic, and is an annihilating polynomial of
$a$ (since $\left(  \lambda^{-1}u\right)  \left[  a\right]  =\lambda^{-1}%
\cdot\underbrace{u\left[  a\right]  }_{=0}=0$). Thus, the degree of this
polynomial $\lambda^{-1}u$ must be at least as large as that of $f$ (since $f$
is a monic annihilating polynomial of $a$ of smallest possible degree). In
other words, $\deg\left(  \lambda^{-1}u\right)  \geq\deg f$. This contradicts
$\deg\left(  \lambda^{-1}u\right)  =\deg u<\deg f$.

This contradiction shows that our assumption was wrong. Hence, at least one of
$u$ and $v$ is a unit of $S\left[  x\right]  $.

Forget that we fixed $u,v$. We thus have shown that whenever $u,v\in S\left[
x\right]  $ satisfy $uv=f$, at least one of $u$ and $v$ is a unit of $S\left[
x\right]  $. Thus, $f$ is irreducible (since $f$ is not a unit of $S\left[
x\right]  $). This proves Theorem \ref{thm.minpol.basics} \textbf{(b)}.
\medskip

\textbf{(c)} Assume that $a$ has a minimal polynomial. Let $f$ be the minimal
polynomial of $a$. Thus, $f$ is a monic annihilating polynomial of $a$. Hence,
$f\left[  a\right]  =0$.

We know that $F$ is an $S$-algebra (since $F$ is a commutative ring, and $S$
is a subring of $F$). The map%
\begin{align*}
\varphi:S\left[  x\right]   &  \rightarrow F,\\
g  &  \mapsto g\left[  a\right]
\end{align*}
is an $S$-algebra morphism (by Theorem \ref{thm.polring.univar-sub-hom}). This
map $\varphi$ sends the principal ideal $fS\left[  x\right]  $ to $0$, because
for each $q\in S\left[  x\right]  $, we have%
\[
\varphi\left(  fq\right)  =\left(  fq\right)  \left[  a\right]
=\underbrace{f\left[  a\right]  }_{=0}\cdot\,q\left[  a\right]  =0.
\]
Hence, the universal property of quotient algebras (Theorem
\ref{thm.quotalg.uniprop1}, applied to $S$, $S\left[  x\right]  $, $fS\left[
x\right]  $, $F$ and $\varphi$ instead of $R$, $A$, $I$, $B$ and $f$) yields
that the map%
\begin{align*}
\varphi^{\prime}:S\left[  x\right]  /f  &  \rightarrow F,\\
\overline{g}  &  \mapsto\varphi\left(  g\right)  \ \ \ \ \ \ \ \ \ \ \left(
\text{for all }g\in S\left[  x\right]  \right)
\end{align*}
is well-defined and is an $S$-algebra morphism. Consider this $\varphi
^{\prime}$. Thus, each $g\in S\left[  x\right]  $ satisfies
\begin{align}
\varphi^{\prime}\left(  \overline{g}\right)   &  =\varphi\left(  g\right)
\ \ \ \ \ \ \ \ \ \ \left(  \text{by the definition of }\varphi^{\prime
}\right) \nonumber\\
&  =g\left[  a\right]  \ \ \ \ \ \ \ \ \ \ \left(  \text{by the definition of
}\varphi\right)  . \label{sol.minpols.1.c.3}%
\end{align}
Thus, $\varphi^{\prime}$ is precisely the map
\begin{align*}
S\left[  x\right]  /f  &  \rightarrow F,\\
\overline{g}  &  \mapsto g\left[  a\right]
\end{align*}
that Theorem \ref{thm.minpol.basics} \textbf{(c)} is talking about. In
particular, we now know that this map is well-defined and is an $S$-algebra
morphism. It remains to prove that this map $\varphi^{\prime}$ is injective.

Since $\varphi^{\prime}$ is $S$-linear, it suffices to show that
$\operatorname*{Ker}\left(  \varphi^{\prime}\right)  =\left\{  0\right\}  $
(by Lemma \ref{lem.modmor.ker-inj}).

Let $u\in\operatorname*{Ker}\left(  \varphi^{\prime}\right)  $. We shall prove
that $u=0$.

Indeed, we have $u\in\operatorname*{Ker}\left(  \varphi^{\prime}\right)
\subseteq S\left[  x\right]  /f$, so that $u=\overline{g}$ for some $g\in
S\left[  x\right]  $. Consider this $g$. We have $\varphi^{\prime}\left(
\overline{g}\right)  =0$ (since $\overline{g}=u\in\operatorname*{Ker}\left(
\varphi^{\prime}\right)  $), so that $0=\varphi^{\prime}\left(  \overline
{g}\right)  =g\left[  a\right]  $ (by (\ref{sol.minpols.1.c.3})). Hence, Lemma
\ref{lem.minpol.div} \textbf{(a)} yields $f\mid g$. In other words,
$\overline{g}=0$ in $S\left[  x\right]  /f$. Hence, $u=\overline{g}%
=0\in\left\{  0\right\}  $.

Forget that we fixed $u$. We thus have shown that $u\in\left\{  0\right\}  $
for each $u\in\operatorname*{Ker}\left(  \varphi^{\prime}\right)  $. In other
words, $\operatorname*{Ker}\left(  \varphi^{\prime}\right)  \subseteq\left\{
0\right\}  $. Since the reverse inclusion $\left\{  0\right\}  \subseteq
\operatorname*{Ker}\left(  \varphi^{\prime}\right)  $ is obvious, we thus
conclude that $\operatorname*{Ker}\left(  \varphi^{\prime}\right)  =\left\{
0\right\}  $. As we have said, this entails that $\varphi^{\prime}$ is
injective. This completes the proof of Theorem \ref{thm.minpol.basics}
\textbf{(c)}. \medskip

\textbf{(d)} We must show that the annihilating polynomials of $a$ are
precisely the polynomials $g\in S\left[  x\right]  $ that are divisible by
$f$. In other words, we must prove the following two statements:

\begin{statement}
\textit{Statement 1:} Any annihilating polynomial of $a$ is divisible by $f$.
\end{statement}

\begin{statement}
\textit{Statement 2:} Any polynomial $g\in S\left[  x\right]  $ that is
divisible by $f$ is an annihilating polynomial of $a$.
\end{statement}

Fortunately, Statement 1 is just a restatement of Lemma \ref{lem.minpol.div}
\textbf{(a)} (since an annihilating polynomial of $a$ is precisely a
polynomial $g\in S\left[  x\right]  $ such that $g\left[  a\right]  =0$), and
Statement 2 is likewise a restatement of Lemma \ref{lem.minpol.div}
\textbf{(b)}. Thus, both Statements 1 and 2 have already been proven, and
Theorem \ref{thm.minpol.basics} \textbf{(d)} is proved.
\end{proof}

\begin{corollary}
\label{cor.minpol.irred->min}Let $S$ and $F$ be two fields such that $S$ is a
subring of $F$. Let $a\in F$.

Let $g\in S\left[  x\right]  $ be a monic irreducible polynomial such that
$g\left[  a\right]  =0$. Then, $g$ is the minimal polynomial of $a$ (over $S$).
\end{corollary}

\begin{proof}
We have $g\left[  a\right]  =0$. In other words, $a$ is a root of $g$. In
other words, $g$ is an annihilating polynomial of $a$. We thus conclude that
$a$ has a monic annihilating polynomial in $S\left[  x\right]  $ (namely,
$g$). Hence, there exists at least one monic annihilating polynomial of $a$ in
$S\left[  x\right]  $. Thus, there also exists such a polynomial of smallest
degree. In other words, $a$ has a minimal polynomial.

Let $f$ be this minimal polynomial. Thus, $f$ is monic. Moreover, $f$ is
irreducible (by Theorem \ref{thm.minpol.basics} \textbf{(b)}), and hence
cannot be a unit of $S\left[  x\right]  $.

Furthermore, Lemma \ref{lem.minpol.div} \textbf{(a)} yields that $f\mid g$
(since $g\left[  a\right]  =0$). In other words, there exists a polynomial
$q\in S\left[  x\right]  $ such that $g=fq$. Consider this $q$. Since $g$ is
irreducible, we conclude that one of $f$ and $q$ is a unit (since $g=fq$).
Since $f$ cannot be a unit, we thus conclude that $q$ is a unit. In other
words, $q$ is a nonzero scalar. This scalar $q$ must be $1$ (since both $f$
and $fq=g$ are monic). Thus, $g=f\underbrace{q}_{=1}=f$. Hence, $g$ is the
minimal polynomial of $a$ (since $f$ is the minimal polynomial of $a$). This
proves Corollary \ref{cor.minpol.irred->min}.
\end{proof}

\begin{exercise}
\label{exe.minpol.sqrtd}Let $S=\mathbb{Q}$ and $F=\mathbb{C}$. Let
$d\in\mathbb{Q}$. Prove the following:

\begin{enumerate}
\item[\textbf{(a)}] The minimal polynomial of $\sqrt{d}$ (over the field
$S=\mathbb{Q}$) is $%
\begin{cases}
x-\sqrt{d}, & \text{if }\sqrt{d}\in\mathbb{Q};\\
x^{2}-d, & \text{otherwise.}%
\end{cases}
$

(In particular, the minimal polynomial of the imaginary unit $i=\sqrt{-1}$ is
$x^{2}+1$.)

\item[\textbf{(b)}] The minimal polynomial of $\sqrt[3]{d}$ (over the field
$S=\mathbb{Q}$) is $%
\begin{cases}
x-\sqrt[3]{d}, & \text{if }\sqrt[3]{d}\in\mathbb{Q};\\
x^{3}-d, & \text{otherwise.}%
\end{cases}
$

\item[\textbf{(c)}] If $d>0$ and $\sqrt{d}\notin\mathbb{Q}$, then the minimal
polynomial of $\sqrt[4]{d}$ (over the field $S=\mathbb{Q}$) is $x^{4}-d$.

\item[\textbf{(d)}] What is the minimal polynomial of $\sqrt[4]{-4}$ (over the
field $S=\mathbb{Q}$) ? (Interpret $\sqrt[4]{-4}$ as the $4$-th root of $-4$
with argument $\pi/4$.)

\item[\textbf{(e)}] Describe the minimal polynomial of $\sqrt[4]{d}$ in general.
\end{enumerate}

[\textbf{Hint:} In part \textbf{(d)}, the answer is neither a degree-$1$
polynomial nor a degree-$4$ polynomial. In part \textbf{(e)}, there are three cases.]
\end{exercise}

\begin{exercise}
\label{exe.21hw4.4}Let $S=\mathbb{Q}$ and $F=\mathbb{R}$. Let $p$ and $q$ be
two positive integers such that none of $p$, $q$ and $pq$ is a perfect square
(i.e., a square in $\mathbb{Z}$). (For example, we can take $p=5$ and $q=8$.)
Let $a=\sqrt{p}+\sqrt{q}\in F$.

Let $f$ denote the polynomial
\[
\left(  x^{2}-p-q\right)  ^{2}-4pq=x^{4}-2\left(  p+q\right)  x^{2}+\left(
p-q\right)  ^{2}\in S\left[  x\right]  .
\]

\begin{enumerate}
\item[\textbf{(a)}] Show that $f$ is an annihilating polynomial of $a$ (that
is, $f\left[  a\right]  =0$).

\item[\textbf{(b)}] Show that $f$ has no rational root.

\item[\textbf{(c)}] Show that $f$ is irreducible (in $S\left[  x\right]  $).

\item[\textbf{(d)}] Conclude that $f$ is the minimal polynomial of $a$ (over
the field $S=\mathbb{Q}$).
\end{enumerate}

[\textbf{Hint:} Part \textbf{(c)} is the tricky one. One way to prove it is to
decompose $f$ as%
\[
f=\left(  x-\left(  \sqrt{p}+\sqrt{q}\right)  \right)  \left(  x-\left(
\sqrt{p}-\sqrt{q}\right)  \right)  \left(  x-\left(  -\sqrt{p}+\sqrt
{q}\right)  \right)  \left(  x-\left(  -\sqrt{p}-\sqrt{q}\right)  \right)
\]
in $\mathbb{R}\left[  x\right]  $, and show that no two of the four factors
here yield a polynomial in $\mathbb{Q}\left[  x\right]  $ when multiplied.
Another approach uses the fact that the polynomial $f$ is \textbf{even} --
meaning that $f\left[  -x\right]  =f$, or, equivalently (since $S$ has
characteristic $0$) that no odd powers of $x$ appear in $f$. This can be used
to argue that if $f=g_{1}g_{2}\cdots g_{k}$ is the factorization of $f$ into
monic irreducible polynomials, then substituting $-x$ for $x$ into it must
yield another factorization $f=f\left[  -x\right]  =g_{1}\left[  -x\right]
g_{2}\left[  -x\right]  \cdots g_{k}\left[  -x\right]  $ of $f$ into monic
irreducible polynomials (why are they still monic?). Since $S\left[  x\right]
$ is a UFD, the two factorizations must be identical (up to the order of the
factors). This narrows down the possibilities for $g_{1},g_{2},\ldots,g_{k}$ substantially.]
\end{exercise}

\begin{exercise}
Let $S$ and $F$ be the subrings
\[
S:=\left\{  \left(
\begin{array}
[c]{cc}%
u & 0\\
0 & u
\end{array}
\right)  \ \mid\ u\in\mathbb{Z}/4\right\}  \ \ \ \ \ \ \ \ \ \ \text{and}%
\ \ \ \ \ \ \ \ \ \ F:=\left\{  \left(
\begin{array}
[c]{cc}%
u & v\\
0 & u
\end{array}
\right)  \ \mid\ u,v\in\mathbb{Z}/4\right\}
\]
of the matrix ring $\left(  \mathbb{Z}/4\right)  ^{2\times2}$. It is easy to
see that $F$ is a commutative ring, and that $S$ is a subring of $F$.

Let $a$ be the element $\left(
\begin{array}
[c]{cc}%
0 & 2\\
0 & 0
\end{array}
\right)  $ of $F$.

\begin{enumerate}
\item[\textbf{(a)}] Prove that $a^{2}=2a=0_{2\times2}$.

\item[\textbf{(b)}] Let us extend Definition \ref{def.minpol.def} from fields
to commutative rings in the obvious way (i.e., replacing \textquotedblleft
field\textquotedblright\ by \textquotedblleft commutative
ring\textquotedblright\ throughout this definition). Prove that the element
$a$ of $F$ has two minimal polynomials, namely $x^{2}$ and $x^{2}+2x$.

\item[\textbf{(c)}] Conclude that Theorem \ref{thm.minpol.basics} \textbf{(a)}
no longer holds if we generalize it from fields to commutative rings.
\end{enumerate}
\end{exercise}

\subsubsection{Minimal polynomials in finite fields}

Next, we apply the notion of minimal polynomials to finite fields:

\begin{proposition}
\label{prop.minpol.pm}Let $p$ be a prime number. Let $S$ be the field
$\mathbb{Z}/p$.

Let $F$ be a finite field of characteristic $p$. Assume that $S=\mathbb{Z}/p$
is a subring of $F$.

We shall use the terminology from Definition \ref{def.minpol.def}.

Let $a\in F$ be arbitrary. Then:

\begin{enumerate}
\item[\textbf{(a)}] The minimal polynomial of $a$ (over $S$) exists (i.e.,
there is always at least one monic annihilating polynomial of $a$).

\item[\textbf{(b)}] If the minimal polynomial of $a$ has degree $k$, then
$a^{p^{k}}=a$.

\item[\textbf{(c)}] Let $m$ be a positive integer satisfying $\left\vert
F\right\vert \leq p^{m}$. If the minimal polynomial of $a$ has degree $k$,
then $k\leq m$.
\end{enumerate}
\end{proposition}

\begin{proof}
First, we note that $0\neq1$ in $F$ (since $F$ is a field). Thus, $F$ has at
least two distinct elements; that is, we have $\left\vert F\right\vert >1$.

Note also that $\left\vert S\right\vert =p$ (since $S=\mathbb{Z}/p$). \medskip

\textbf{(a)} Proposition \ref{prop.finfield.flt} yields $a^{\left\vert
F\right\vert }=a$. Thus, $a$ is a root of the polynomial $x^{\left\vert
F\right\vert }-x\in S\left[  x\right]  $. In other words, $x^{\left\vert
F\right\vert }-x\in S\left[  x\right]  $ is an annihilating polynomial of $a$.
Since this polynomial is monic, we thus conclude that $a$ has a monic
annihilating polynomial in $S\left[  x\right]  $ (namely, $x^{\left\vert
F\right\vert }-x$). Hence, there exists at least one monic annihilating
polynomial of $a$ in $S\left[  x\right]  $. Thus, there also exists such a
polynomial of smallest degree. In other words, $a$ has a minimal polynomial.
This proves Proposition \ref{prop.minpol.pm} \textbf{(a)}. \medskip

\textbf{(b)} Let $f$ be the minimal polynomial of $a$. Let $k=\deg f$. We must
show that $a^{p^{k}}=a$.

Theorem \ref{thm.minpol.basics} \textbf{(b)} shows that the polynomial $f$ is
irreducible. Hence, the quotient ring $S\left[  x\right]  /f$ is a field (by
Theorem \ref{thm.fieldext.field-irr}). On the other hand, the leading
coefficient of $f$ is a unit (since $f$ is monic). Thus, as an $S$-module,
$S\left[  x\right]  /f$ is free of rank $k=\deg f$ (by Theorem
\ref{thm.fieldext.basis-monic} \textbf{(b)}). Hence, $S\left[  x\right]
/f\cong S^{k}$ as an $S$-module. Hence, $\left\vert S\left[  x\right]
/f\right\vert =\left\vert S^{k}\right\vert =\left\vert S\right\vert ^{k}%
=p^{k}$ (since $\left\vert S\right\vert =p$). Therefore, the field $S\left[
x\right]  /f$ is finite.

Now, $\overline{x}\in S\left[  x\right]  /f$ is an element of this finite
field $S\left[  x\right]  /f$. Hence, Proposition \ref{prop.finfield.flt}
(applied to $S\left[  x\right]  /f$ and $\overline{x}$ instead of $F$ and $u$)
yields $\overline{x}^{\left\vert S\left[  x\right]  /f\right\vert }%
=\overline{x}$. In view of $\left\vert S\left[  x\right]  /f\right\vert
=p^{k}$, this rewrites as $\overline{x}^{p^{k}}=\overline{x}$. In other words,
$\overline{x^{p^{k}}}=\overline{x}$ (since $\overline{x^{p^{k}}}=\overline
{x}^{p^{k}}$).

Theorem \ref{thm.minpol.basics} \textbf{(c)} shows that the map
\begin{align*}
S\left[  x\right]  /f  &  \rightarrow F,\\
\overline{g}  &  \mapsto g\left[  a\right]
\end{align*}
is a (well-defined) $S$-algebra morphism, and is injective. Applying this map
to both sides of the equality $\overline{x^{p^{k}}}=\overline{x}$, we obtain
$x^{p^{k}}\left[  a\right]  =x\left[  a\right]  $. In other words, $a^{p^{k}%
}=a$. This proves Proposition \ref{prop.minpol.pm} \textbf{(b)}. \medskip

\textbf{(c)} Let $f$ be the minimal polynomial of $a$. Let $k=\deg f$. We must
show that $k\leq m$.

We have already shown (in the above proof of part \textbf{(b)}) that
$\left\vert S\left[  x\right]  /f\right\vert =p^{k}$. Moreover, we have
already shown (in the above proof of part \textbf{(b)}) that the map%
\begin{align*}
S\left[  x\right]  /f  &  \rightarrow F,\\
\overline{g}  &  \mapsto g\left[  a\right]
\end{align*}
is injective. Hence, $\left\vert S\left[  x\right]  /f\right\vert
\leq\left\vert F\right\vert $ (because if $U$ and $V$ are two finite sets, and
if there exists an injective map from $U$ to $V$, then $\left\vert
U\right\vert \leq\left\vert V\right\vert $). In view of $\left\vert S\left[
x\right]  /f\right\vert =p^{k}$, this rewrites as $p^{k}\leq\left\vert
F\right\vert $. Hence, $p^{k}\leq\left\vert F\right\vert \leq p^{m}$, and
therefore $k\leq m$ (since $p>1$). This proves Proposition
\ref{prop.minpol.pm} \textbf{(c)}.
\end{proof}

\begin{fineprint}

\begin{remark}
Proposition \ref{prop.minpol.pm} can be proved in many other ways. For
example, part \textbf{(a)} can also be obtained from the pigeonhole principle
(to wit: the principle shows that two of the $\left\vert F\right\vert +1$
elements $a^{0},a^{1},a^{2},\ldots,a^{\left\vert F\right\vert }$ are equal;
thus, $a$ has an annihilating polynomial of the form $x^{i}-x^{j}$ for some
$i>j\geq0$). For parts \textbf{(b)} and \textbf{(c)}, it helps to look at the
subset%
\[
A:=\left\{  g\left[  a\right]  \ \mid\ g\in S\left[  x\right]  \right\}
\]
of $F$. This subset $A$ is easily seen to be a subring of $F$, and thus a
field (indeed, any subring of $F$ is a finite integral domain and therefore a
field). Moreover, it has size $p^{k}$ (since the division-with-remainder
theorem for polynomials shows that it is a free $\mathbb{Z}/p$-module with
basis $\left(  a^{0},a^{1},\ldots,a^{k-1}\right)  $). Thus, $p^{k}=\left\vert
A\right\vert \leq\left\vert F\right\vert \leq p^{m}$, so that $k\leq m$, and
this yields part \textbf{(c)} of Proposition \ref{prop.minpol.pm}. Moreover,
part \textbf{(b)} follows by applying Proposition \ref{prop.finfield.flt} to
$A$ instead of $F$ (since $a\in A$).
\end{remark}
\end{fineprint}

\begin{exercise}
Let $p$, $S$, $F$ and $a$ be as in Proposition \ref{prop.minpol.pm}.

Assume that the minimal polynomial of $a$ has degree $k$. Prove that this
minimal polynomial is%
\[
\left(  x-a^{p^{0}}\right)  \left(  x-a^{p^{1}}\right)  \cdots\left(
x-a^{p^{k-1}}\right)  =\prod_{i=0}^{k-1}\left(  x-a^{p^{i}}\right)  .
\]

[\textbf{Hint:} First prove that not only $a$, but all the $k$ powers
$a^{p^{0}},a^{p^{1}},\ldots,a^{p^{k-1}}$ are roots of the minimal polynomial.
Then prove that these $k$ powers are distinct. To do the latter, consider two
integers $i,j$ with $0\leq i<j<k$. Then, show that the set
\[
F_{i,j}:=\left\{  u\in F\ \mid\ u^{p^{i}}=u^{p^{j}}\right\}
\]
is a subring of $F$ (why?) and thus is a field (why?), but has size
$\left\vert F_{i,j}\right\vert \leq p^{j}$ (why?). Now apply Proposition
\ref{prop.minpol.pm} \textbf{(c)} to $F_{i,j}$ and $j$ instead of $F$ and $m$
to obtain a contradiction if $a\in F_{i,j}$.]
\end{exercise}

\begin{lemma}
\label{lem.minpol.pm2}Let $p$ be a prime number. Let $F$ be a finite field of
characteristic $p$. Let $m$ be the positive integer satisfying $\left\vert
F\right\vert =p^{m}$. (We know from Theorem \ref{thm.finfield.char.basics}
\textbf{(e)} that this $m$ exists.)

Then, there exists at least one $a\in F$ such that none of the $m-1$ powers
$a^{p^{1}},a^{p^{2}},\ldots,a^{p^{m-1}}$ equals $a$.
\end{lemma}

\begin{proof}
Assume the contrary. Thus, each $a\in F$ satisfies at least one of the $m-1$
equations%
\[
a^{p^{1}}=a,\ \ \ \ \ \ \ \ \ \ a^{p^{2}}=a,\ \ \ \ \ \ \ \ \ \ \ldots
,\ \ \ \ \ \ \ \ \ \ a^{p^{m-1}}=a.
\]
In other words, each $a\in F$ is a root of the polynomial%
\[
\left(  x^{p^{1}}-x\right)  \left(  x^{p^{2}}-x\right)  \cdots\left(
x^{p^{m-1}}-x\right)  \in F\left[  x\right]  .
\]
But this polynomial has degree $p^{1}+p^{2}+\cdots+p^{m-1}$, and thus has at
most $p^{1}+p^{2}+\cdots+p^{m-1}$ many roots (by Theorem
\ref{thm.polring.univar-easyFTA}). Hence, it has fewer than $p^{m}$ roots
(since it is easy to see that $p^{1}+p^{2}+\cdots+p^{m-1}=\dfrac{p^{m}-1}%
{p-1}-1<\dfrac{p^{m}-1}{p-1}\leq p^{m}-1<p^{m}$). However, we just found out
that each $a\in F$ is a root of this polynomial; thus, this polynomial has at
least $p^{m}$ roots (since $\left\vert F\right\vert =p^{m}$). The preceding
two sentences contradict each other. This contradiction shows that our
assumption was wrong; hence, Lemma \ref{lem.minpol.pm2} is proved.
\end{proof}

\subsubsection{Each finite field is obtained from $\mathbb{Z}/p$ by a single
root adjunction}

The above results lead to a first interesting property of finite fields:

\begin{theorem}
\label{thm.finfields.one-root}Let $p$ be a prime number. Let $S$ be the field
$\mathbb{Z}/p$.

Let $F$ be a finite field of characteristic $p$.

Let $m$ be the positive integer satisfying $\left\vert F\right\vert =p^{m}$.
(We know from Theorem \ref{thm.finfield.char.basics} \textbf{(e)} that this
$m$ exists.)

Then, there exists at least one monic irreducible polynomial $f\in S\left[
x\right]  =\left(  \mathbb{Z}/p\right)  \left[  x\right]  $ of degree $m$ that
satisfies $F\cong S\left[  x\right]  /f$.
\end{theorem}

This theorem shows that any finite field can be constructed (up to
isomorphism) by adjoining a (single) root of an irreducible polynomial to a
field of the form $\mathbb{Z}/p$. It also shows that for any prime $p$ and any
positive integer $m$, there exists an irreducible polynomial of degree $m$
over $\mathbb{Z}/p$ (because Theorem \ref{thm.finfield.exist} says that there
exists a finite field $F$ of size $\left\vert F\right\vert =p^{m}$).

\begin{proof}
[Proof of Theorem \ref{thm.finfields.one-root}.]Theorem
\ref{thm.finfield.char.basics} \textbf{(f)} shows that the field $F$ contains
\textquotedblleft a copy of $\mathbb{Z}/p$\textquotedblright\ (that is, a
subring isomorphic to $\mathbb{Z}/p$). By renaming the elements of $F$
accordingly, we WLOG assume that this copy is $\mathbb{Z}/p$ itself, i.e.,
that $\mathbb{Z}/p$ is a subring of $F$. In other words, $S$ is a subring of
$F$ (since $S=\mathbb{Z}/p$).

Lemma \ref{lem.minpol.pm2} shows that there exists at least one $a\in F$ such
that none of the $m-1$ powers $a^{p^{1}},a^{p^{2}},\ldots,a^{p^{m-1}}$ equals
$a$. Consider this $a$. Proposition \ref{prop.minpol.pm} \textbf{(a)} shows
that $a$ has a minimal polynomial (over $S$). Let $f\in S\left[  x\right]  $
be this polynomial, and let $k=\deg f$ be its degree. Theorem
\ref{thm.minpol.basics} \textbf{(b)} shows that this minimal polynomial $f$ is
irreducible. Thus, $f$ is not constant; hence, its degree $k$ is positive.
That is, we have $k\neq0$.

Proposition \ref{prop.minpol.pm} \textbf{(b)} shows that $a^{p^{k}}=a$;
therefore, $k\notin\left\{  1,2,\ldots,m-1\right\}  $ (since none of the $m-1$
powers $a^{p^{1}},a^{p^{2}},\ldots,a^{p^{m-1}}$ equals $a$). Combining this
with $k\neq0$, we thus obtain $k\notin\left\{  0,1,\ldots,m-1\right\}  $, so
that $k\geq m$. However, Proposition \ref{prop.minpol.pm} \textbf{(c)} shows
that $k\leq m$. Combined with $k\geq m$, this yields $k=m$. Thus, $f$ has
degree $m$ (since $f$ has degree $\deg f=k=m$).

It remains to show that $F\cong S\left[  x\right]  /f$ (as $S$-algebras).

As in the proof of Proposition \ref{prop.minpol.pm} \textbf{(b)}, we can see
that $\left\vert S\left[  x\right]  /f\right\vert =p^{k}$. In view of $k=m$,
this rewrites as $\left\vert S\left[  x\right]  /f\right\vert =p^{m}$.
Compared with $\left\vert F\right\vert =p^{m}$, this yields $\left\vert
S\left[  x\right]  /f\right\vert =\left\vert F\right\vert $. Hence, $S\left[
x\right]  /f$ and $F$ are two finite sets of the same size.

Theorem \ref{thm.minpol.basics} \textbf{(c)} yields that the map
\begin{align*}
S\left[  x\right]  /f  &  \rightarrow F,\\
\overline{g}  &  \mapsto g\left[  a\right]
\end{align*}
is a (well-defined) $S$-algebra morphism, and is injective. This map is thus
an injective map between two finite sets of the same size (since $S\left[
x\right]  /f$ and $F$ are two finite sets of the same size), and therefore is
bijective (since the pigeonhole principle shows that any injective map between
two finite sets of the same size is bijective). In other words, this map is
invertible. Since it is an $S$-algebra morphism, it is thus an $S$-algebra
isomorphism (by Proposition \ref{prop.algmor.invertible-iso}). Hence, $F\cong
S\left[  x\right]  /f$ (as $S$-algebras). This completes the proof of Theorem
\ref{thm.finfields.one-root}.
\end{proof}

We also observe the following curious fact:

\begin{theorem}
\label{thm.finfields.irr-div-xpm-x}Let $p$ be a prime number. Let $S$ be the
field $\mathbb{Z}/p$.

Let $m$ be a positive integer.

Then, any irreducible polynomial $f\in S\left[  x\right]  $ of degree $m$
divides $x^{p^{m}}-x\in S\left[  x\right]  $.
\end{theorem}

\begin{proof}
Let $f\in S\left[  x\right]  $ be an irreducible polynomial of degree $m$. We
must show that $f\mid x^{p^{m}}-x$ in $S\left[  x\right]  $.

The quotient ring $S\left[  x\right]  /f$ is a field (by Theorem
\ref{thm.fieldext.field-irr}, since $f$ is irreducible). On the other hand,
the leading coefficient of $f$ is a unit (since $S$ is a field, so that every
nonzero element of $S$ is a unit). Thus, as an $S$-module, $S\left[  x\right]
/f$ is free of rank $m=\deg f$ (by Theorem \ref{thm.fieldext.basis-monic}
\textbf{(b)}). Hence, $S\left[  x\right]  /f\cong S^{m}$ as an $S$-module.
Thus, $\left\vert S\left[  x\right]  /f\right\vert =\left\vert S^{m}%
\right\vert =\left\vert S\right\vert ^{m}=p^{m}$ (since $\left\vert
S\right\vert =p$). Therefore, the field $S\left[  x\right]  /f$ is finite.

Now, let $F$ be this finite field $S\left[  x\right]  /f$. Let $a$ be the
element $\overline{x}\in S\left[  x\right]  /f=F$. Then, Proposition
\ref{prop.finfield.flt} yields $a^{\left\vert F\right\vert }=a$. However,
$F=S\left[  x\right]  /f$, so that $\left\vert F\right\vert =\left\vert
S\left[  x\right]  /f\right\vert =p^{m}$. Hence,
\begin{align*}
a^{\left\vert F\right\vert }  &  =a^{p^{m}}=\overline{x}^{p^{m}}%
\ \ \ \ \ \ \ \ \ \ \left(  \text{since }a=\overline{x}\right) \\
&  =\overline{x^{p^{m}}},
\end{align*}
so that $\overline{x^{p^{m}}}=a^{\left\vert F\right\vert }=a=\overline{x}$. In
other words, $x^{p^{m}}-x\in fS\left[  x\right]  $. In other words, $f\mid
x^{p^{m}}-x$ in $S\left[  x\right]  $. This proves Theorem
\ref{thm.finfields.irr-div-xpm-x}.
\end{proof}

As a consequence of this theorem, we can show that Theorem
\ref{thm.finfields.one-root} can be strengthened: Not only is there
\textbf{some} monic irreducible polynomial $f\in S\left[  x\right]  =\left(
\mathbb{Z}/p\right)  \left[  x\right]  $ of degree $m$ that satisfies $F\cong
S\left[  x\right]  /f$, but actually \textbf{any} monic irreducible polynomial
of this degree will do! Let us state this more carefully:

\begin{corollary}
\label{cor.finfields.pm-already-root}Let $p$ be a prime number. Let $S$ be the
field $\mathbb{Z}/p$.

Let $m$ be a positive integer. Let $F$ be a finite field of characteristic $p$
such that $\left\vert F\right\vert =p^{m}$.

Let $f\in S\left[  x\right]  $ be a monic irreducible polynomial of degree
$m$. Then, $F\cong S\left[  x\right]  /f$ (as $S$-algebras).
\end{corollary}

\begin{proof}
Theorem \ref{thm.finfield.char.basics} \textbf{(f)} shows that the field $F$
contains \textquotedblleft a copy of $\mathbb{Z}/p$\textquotedblright\ (that
is, a subring isomorphic to $\mathbb{Z}/p$). By renaming the elements of $F$
accordingly, we WLOG assume that this copy is $\mathbb{Z}/p$ itself, i.e.,
that $\mathbb{Z}/p$ is a subring of $F$. In other words, $S$ is a subring of
$F$ (since $S=\mathbb{Z}/p$).

We have $\deg f=m>0$. Hence, Theorem \ref{thm.finfield.splitfield}
\textbf{(b)} (applied to $b=f$) shows that there is a field that contains $F$
as a subring and that contains a root of $f$. Let $F^{\prime}$ be this field,
and let $a$ be this root. Thus, $a\in F^{\prime}$ and $f\left[  a\right]  =0$.

Our first goal is to prove that $a$ actually belongs to $F$.

Indeed, Theorem \ref{thm.finfields.irr-div-xpm-x} shows that $f$ divides
$x^{p^{m}}-x$ in $S\left[  x\right]  $. Hence, $a$ is a root of the polynomial
$x^{p^{m}}-x$ (since $a$ is a root of $f$). In other words, $a^{p^{m}}-a=0$.

On the other hand, Proposition \ref{prop.finfield.flt} yields that
$u^{\left\vert F\right\vert }=u$ for each $u\in F$. In other words, $u^{p^{m}%
}=u$ for each $u\in F$ (since $\left\vert F\right\vert =p^{m}$). In other
words, $u^{p^{m}}-u=0$ holds for each $u\in F$. This equality $u^{p^{m}}-u=0$
holds for $u=a$ as well (since $a^{p^{m}}-a=0$). Thus, we conclude that
$u^{p^{m}}-u=0$ holds for each $u\in F\cup\left\{  a\right\}  $. In other
words, each $u\in F\cup\left\{  a\right\}  $ is a root of the polynomial
$x^{p^{m}}-x$. Therefore, the polynomial $x^{p^{m}}-x$ has at least
$\left\vert F\cup\left\{  a\right\}  \right\vert $ many roots in $F^{\prime}$.
But this polynomial $x^{p^{m}}-x$ has degree $p^{m}$ (since $p^{m}>1$), and
thus has at most $p^{m}$ roots in $F^{\prime}$ (by Theorem
\ref{thm.polring.univar-easyFTA}, applied to $F^{\prime}$, $p^{m}$ and
$x^{p^{m}}-x$ instead of $R$, $n$ and $f$). Confronting the preceding two
sentences with each other, we obtain $\left\vert F\cup\left\{  a\right\}
\right\vert \leq p^{m}$. In view of $p^{m}=\left\vert F\right\vert $, this
rewrites as $\left\vert F\cup\left\{  a\right\}  \right\vert \leq\left\vert
F\right\vert $. Since $F$ is a finite set, this inequality yields $a\in F$
(since otherwise, we would have $\left\vert F\cup\left\{  a\right\}
\right\vert =\left\vert F\right\vert +1>\left\vert F\right\vert $). Thus, we
have shown that $a$ belongs to $F$. We can now forget about $F^{\prime}$.

We have $\deg f=m$, and the leading coefficient of $f$ is a unit (since $S$ is
a field, so that every nonzero element of $S$ is a unit). Thus, as an
$S$-module, $S\left[  x\right]  /f$ is free of rank $\deg f=m$ (by Theorem
\ref{thm.fieldext.basis-monic} \textbf{(b)}). Hence, $S\left[  x\right]
/f\cong S^{m}$ as an $S$-module. Hence, $\left\vert S\left[  x\right]
/f\right\vert =\left\vert S^{m}\right\vert =\left\vert S\right\vert ^{m}%
=p^{m}$ (since $\left\vert S\right\vert =p$). Comparing this with $\left\vert
F\right\vert =p^{m}$, we obtain $\left\vert S\left[  x\right]  /f\right\vert
=\left\vert F\right\vert $. In other words, $S\left[  x\right]  /f$ and $F$
are two finite sets of the same size.

The polynomial $f$ is monic and irreducible, and satisfies $f\left[  a\right]
=0$. Hence, Corollary \ref{cor.minpol.irred->min} (applied to $g=f$) shows
that $f$ is the minimal polynomial of $a$ (over $S$). Thus, Theorem
\ref{thm.minpol.basics} \textbf{(c)} yields that the map
\begin{align*}
S\left[  x\right]  /f  &  \rightarrow F,\\
\overline{g}  &  \mapsto g\left[  a\right]
\end{align*}
is a (well-defined) $S$-algebra morphism, and is injective. This map is thus
an injective map between two finite sets of the same size (since $S\left[
x\right]  /f$ and $F$ are two finite sets of the same size), and therefore is
bijective (since the pigeonhole principle shows that any injective map between
two finite sets of the same size is bijective). In other words, this map is
invertible. Since it is an $S$-algebra morphism, it is thus an $S$-algebra
isomorphism (by Proposition \ref{prop.algmor.invertible-iso}). Hence, $F\cong
S\left[  x\right]  /f$ (as $S$-algebras). This proves Corollary
\ref{cor.finfields.pm-already-root}.
\end{proof}

\begin{exercise}
\label{exe.finfields.xpm-x.1}Let $p$ be a prime number. Let $S$ be the field
$\mathbb{Z}/p$.

Let $m$ be a positive integer. Let $f\in S\left[  x\right]  $ be an
irreducible polynomial.

Prove the following:

\begin{enumerate}
\item[\textbf{(a)}] If $\deg f\mid m$, then $f\mid x^{p^{m}}-x$ in $S\left[
x\right]  $.

(Note that this generalizes Theorem \ref{thm.finfields.irr-div-xpm-x}.)

\item[\textbf{(b)}] Conversely, if $f\mid x^{p^{m}}-x$ in $S\left[  x\right]
$, then $\deg f\mid m$.
\end{enumerate}
\end{exercise}

\begin{exercise}
Let $p$ be a prime number. Let $S$ be the field $\mathbb{Z}/p$. Let $m$ be a
positive integer. Prove that the polynomial $x^{p^{m}}-x\in S\left[  x\right]
$ equals the product of all monic irreducible polynomials $f\in S\left[
x\right]  $ satisfying $\deg f\mid m$.

(For example, for $p=2$ and $m=3$, this yields%
\[
x^{8}-x=\underbrace{x\left(  x+1\right)  }_{\substack{\text{irreducible}%
\\\text{polynomials}\\\text{of degree }1}}\underbrace{\left(  x^{3}%
+x+1\right)  \left(  x^{3}+x^{2}+1\right)  }_{\substack{\text{irreducible}%
\\\text{polynomials}\\\text{of degree }3}}
\]
and%
\[
x^{16}-x=\underbrace{x\left(  x+1\right)  }_{\substack{\text{irreducible}%
\\\text{polynomials}\\\text{of degree }1}}\underbrace{\left(  x^{2}%
+x+1\right)  }_{\substack{\text{irreducible}\\\text{polynomials}\\\text{of
degree }2}}\underbrace{\left(  x^{4}+x+1\right)  \left(  x^{4}+x^{3}+1\right)
\left(  x^{4}+x^{3}+x^{2}+x+1\right)  }_{\substack{\text{irreducible}%
\\\text{polynomials}\\\text{of degree }4}}
\]
in $\left(  \mathbb{Z}/2\right)  \left[  x\right]  $.) \medskip

[\textbf{Hint:} First argue that the factorization of $x^{p^{m}}-x$ into monic
irreducible polynomials contains no factor more than once. Then use Exercise
\ref{exe.finfields.xpm-x.1}.]
\end{exercise}

\subsubsection{Proof of the uniqueness}

We can now easily prove the uniqueness of a finite field of given size (up to isomorphism):

\begin{theorem}
\label{thm.finfield.unique}Any two finite fields that have the same size are isomorphic.
\end{theorem}

\begin{proof}
Let $F$ and $F^{\prime}$ be two finite fields that have the same size. Thus,
$\left\vert F\right\vert =\left\vert F^{\prime}\right\vert $. Our goal is to
prove that $F\cong F^{\prime}$.

Let $p=\operatorname*{char}F$. Then, Theorem \ref{thm.finfield.char.basics}
\textbf{(d)} shows that $p$ is a prime. Also, Theorem
\ref{thm.finfield.char.basics} \textbf{(e)} shows that $\left\vert
F\right\vert =p^{m}$ for some positive integer $m$. Consider this $m$.
Comparing $\left\vert F\right\vert =p^{m}$ with $\left\vert F\right\vert
=\left\vert F^{\prime}\right\vert $, we find $\left\vert F^{\prime}\right\vert
=p^{m}$.

Let $q=\operatorname*{char}\left(  F^{\prime}\right)  $. Then, Theorem
\ref{thm.finfield.char.basics} \textbf{(d)} shows that $q$ is a prime. Also,
Theorem \ref{thm.finfield.char.basics} \textbf{(e)} shows that $\left\vert
F^{\prime}\right\vert =q^{n}$ for some positive integer $n$. Consider this
$n$. Now, $p^{m}=\left\vert F\right\vert =\left\vert F^{\prime}\right\vert
=q^{n}$. Since $p$ and $q$ are primes (and $m$ and $n$ are positive integers),
this can only happen if $p=q$ and $m=n$. Thus, we obtain $p=q$ and $m=n$.
Hence, $F^{\prime}$ has characteristic $p$ (since $p=q=\operatorname*{char}%
\left(  F^{\prime}\right)  $).

Set $S=\mathbb{Z}/p$. Theorem \ref{thm.finfields.one-root} yields that there
exists at least one monic irreducible polynomial $f\in S\left[  x\right]
=\left(  \mathbb{Z}/p\right)  \left[  x\right]  $ of degree $m$ that satisfies
$F\cong S\left[  x\right]  /f$. Consider this $f$. Recall that $\left\vert
F^{\prime}\right\vert =p^{m}$. Hence, Corollary
\ref{cor.finfields.pm-already-root} (applied to $F^{\prime}$ instead of $F$)
yields that $F^{\prime}\cong S\left[  x\right]  /f$ (as $S$-algebras).
Combining this with $F\cong S\left[  x\right]  /f$, we obtain $F\cong S\left[
x\right]  /f\cong F^{\prime}$. As we said, this proves Theorem
\ref{thm.finfield.unique}.
\end{proof}

\begin{fineprint}
As we observed after Theorem \ref{thm.finfields.one-root}, there exists an
irreducible polynomial of any positive degree $m$ over $\mathbb{Z}/p$ for any
prime $p$. Explicitly finding such polynomials is not at all easy, unless you
are content with the brute-force method \textquotedblleft try all the
polynomials and check for irreducibility\textquotedblright\ (a finite
algorithm, although a rather laborious one).\footnote{A popular method for
choosing these polynomials is known as
\href{https://en.wikipedia.org/wiki/Conway_polynomial_(finite_fields)}{\textquotedblleft
Conway polynomials\textquotedblright}.} However, for specific degrees, there
are better methods. In particular, for $m=p$, there is an explicit choice:

\begin{exercise}
Let $p$ be a prime number. Let $S$ be the field $\mathbb{Z}/p$. Let $a\in S$
be nonzero. Prove that the polynomial $x^{p}-x+a\in S\left[  x\right]  $ is
irreducible. \medskip

[\textbf{Hint:} Assume that $x^{p}-x+a=fg$ for two non-units $f,g\in S\left[
x\right]  $ with $f$ irreducible. Argue that $\overline{x}^{p}=\overline{x}-a$
in the quotient ring $S\left[  x\right]  /f$. Use this to show that
$\overline{x}^{p^{i}}=\overline{x}-ia$ for all $i\in\mathbb{N}$. On the other
hand, let $k=\deg f\in\left\{  1,2,\ldots,p-1\right\}  $, and argue that
$\overline{x}^{p^{k}}=\overline{x}$ in $S\left[  x\right]  /f$. Can these
formulas coexist?]
\end{exercise}
\end{fineprint}

\subsection{\label{sec.finfields.morelem}Lemmas on $p$-th powers}

The following lemma will be used twice in the next section:

\begin{lemma}
\label{lem.finfield.ap=a-only}Let $p$ be a prime. Let $F$ be a field such that
$\mathbb{Z}/p$ is a subring of $F$. Then,%
\[
\left\{  a\in F\ \mid\ a^{p}=a\right\}  =\mathbb{Z}/p.
\]

\end{lemma}

This lemma gives a criterion for showing that an element of $F$ lies in
$\mathbb{Z}/p$: namely, just show that $a^{p}=a$.

\begin{proof}
[Proof of Lemma \ref{lem.finfield.ap=a-only}.]For each $u\in\mathbb{Z}/p$, we
have $u^{p}=u$ (by Proposition \ref{prop.ent.flt.Z/p}) and thus $u\in\left\{
a\in F\ \mid\ a^{p}=a\right\}  $. In other words, $\mathbb{Z}/p\subseteq
\left\{  a\in F\ \mid\ a^{p}=a\right\}  $.

Now, I claim that $\left\vert \left\{  a\in F\ \mid\ a^{p}=a\right\}
\right\vert \leq p$. Indeed, $F$ is an integral domain. Thus, the easy half of
the FTA (Theorem \ref{thm.polring.univar-easyFTA}) yields that if $n$ is a
nonnegative integer, then any nonzero polynomial of degree $\leq n$ over $F$
has at most $n$ roots in $F$. Applying this to the polynomial $x^{p}-x$ (which
is nonzero and has degree $p$), we conclude that the polynomial $x^{p}-x$ has
at most $p$ roots in $F$. But the set of all roots of this polynomial
$x^{p}-x$ in $F$ is $\left\{  a\in F\ \mid\ a^{p}=a\right\}  $; hence, the
preceding sentence says that $\left\vert \left\{  a\in F\ \mid\ a^{p}%
=a\right\}  \right\vert \leq p$. Thus, in particular, the set $\left\{  a\in
F\ \mid\ a^{p}=a\right\}  $ is finite.

However, an easy and fundamental fact in combinatorics says that if $X$ and
$Y$ are two finite sets with $X\subseteq Y$ and $\left\vert Y\right\vert
\leq\left\vert X\right\vert $, then $X=Y$. Applying this to $X=\mathbb{Z}/p$
and $Y=\left\{  a\in F\ \mid\ a^{p}=a\right\}  $, we obtain $\mathbb{Z}%
/p=\left\{  a\in F\ \mid\ a^{p}=a\right\}  $ (since $\mathbb{Z}/p\subseteq
\left\{  a\in F\ \mid\ a^{p}=a\right\}  $ and $\left\vert \left\{  a\in
F\ \mid\ a^{p}=a\right\}  \right\vert \leq p=\left\vert \mathbb{Z}%
/p\right\vert $). This proves Lemma \ref{lem.finfield.ap=a-only}.
\end{proof}

\begin{exercise}
Let $p$ be a prime. Let $F$ be a finite field of size $p^{m}$, where $m$ is a
positive integer. Assume that $\mathbb{Z}/p$ is a subring of $F$. Let $a\in F$.

\begin{enumerate}
\item[\textbf{(a)}] Let $r=p^{0}+p^{1}+\cdots+p^{m-1}$. Prove that $a^{r}%
\in\mathbb{Z}/p$.

\item[\textbf{(b)}] Let $b=a^{p^{0}}+a^{p^{1}}+\cdots+a^{p^{m-1}}$. Prove that
$b\in\mathbb{Z}/p$.
\end{enumerate}
\end{exercise}

Another useful lemma says (in terms of Subsection
\ref{subsec.finfields.tools.frob}) that the Frobenius endomorphism of a field
of characteristic $p$ is always injective:

\begin{lemma}
\label{lem.finfield.frob-inj}Let $p$ be a prime. Let $F$ be a field of
characteristic $p$. Let $a,b\in F$ satisfy $a\neq b$. Then, $a^{p}\neq b^{p}$.
\end{lemma}

Note that this would fail for $F=\mathbb{R}$ and $p=2$ (because, for example,
$1\neq-1$ but $1^{2}=\left(  -1\right)  ^{2}$), and also fail for
$F=\mathbb{C}$ and any $p>1$. Thus, this marks one more of the situations
where fields of prime characteristic $p$ behave better than fields of
characteristic $0$.

\begin{proof}
[Proof of Lemma \ref{lem.finfield.frob-inj}.]From $a\neq b$, we see that the
element $a-b$ of $F$ is nonzero. But $F$ is a field, and thus an integral
domain. Hence, it is easy to see (by induction on $k$) that any finite product
$u_{1}u_{2}\cdots u_{k}$ of nonzero elements of $F$ is nonzero. Thus, in
particular, the product $\underbrace{\left(  a-b\right)  \left(  a-b\right)
\cdots\left(  a-b\right)  }_{p\text{ times}}$ is nonzero (since $a-b$ is
nonzero). In other words, $\left(  a-b\right)  ^{p}$ is nonzero.

However, Theorem \ref{thm.finfield.idiot} \textbf{(c)} yields $\left(
a-b\right)  ^{p}=a^{p}-b^{p}$, so that $a^{p}-b^{p}=\left(  a-b\right)
^{p}\neq0$ (since $\left(  a-b\right)  ^{p}$ is nonzero). In other words,
$a^{p}\neq b^{p}$. This proves Lemma \ref{lem.finfield.frob-inj}.
\end{proof}

At this point, we end (at least for a while) our theoretical study of finite
fields, and instead focus on some of their applications. The reader can learn
more about finite fields from texts such as \cite{LidNie00} and
\cite{MulMum07}.

\subsection{\label{sec.finfields.approot}An application of root adjunction}

What are finite fields (particularly the ones that are not just $\mathbb{Z}%
/p$) good for? Known applications include error-correcting codes
(\href{https://en.wikipedia.org/wiki/BCH_code}{BCH codes}), group theory (many
\href{https://en.wikipedia.org/wiki/Classification_of_finite_simple_groups}{finite
simple groups} can be constructed as
\href{https://en.wikipedia.org/wiki/Group_of_Lie_type}{matrix groups over
finite fields}), \href{https://en.wikipedia.org/wiki/Block_design}{block
designs} (roughly speaking, finite structures with symmetries that resemble
geometries) and, of course, number theory (not unexpectedly; number theory
uses everything). Various applications along these lines can be found in
\cite{MulMum07}. Let me show a more humble -- but also more self-contained --
application. Namely, by adjoining roots of polynomials to $\mathbb{Z}/p$, we
will prove a curious fact about Fibonacci numbers (\cite[\S 25]{Vorobi02}):

\begin{theorem}
\label{thm.fib.p-divides}Let $\left(  f_{0},f_{1},f_{2},\ldots\right)  $ be
the Fibonacci sequence (as defined in Definition \ref{def.fibonacci.fib}).

Let $p$ be a prime. Then:

\begin{enumerate}
\item[\textbf{(a)}] If $p\equiv\pm1\operatorname{mod}5$ (meaning that $p$ is
congruent to one of $1$ and $-1$ modulo $5$), then $p\mid f_{p-1}$.

\item[\textbf{(b)}] If $p\equiv\pm2\operatorname{mod}5$ (meaning that $p$ is
congruent to one of $2$ and $-2$ modulo $5$), then $p\mid f_{p+1}$.
\end{enumerate}
\end{theorem}

For example:

\begin{itemize}
\item For $p=2$, Theorem \ref{thm.fib.p-divides} \textbf{(b)} says that $2\mid
f_{3}$ (since $2\equiv2\operatorname{mod}5$), and indeed we have $f_{3}=2$.

\item For $p=7$, Theorem \ref{thm.fib.p-divides} \textbf{(b)} says that $7\mid
f_{8}$ (since $7\equiv2\operatorname{mod}5$), and indeed we have
$f_{8}=21=3\cdot7$.

\item For $p=11$, Theorem \ref{thm.fib.p-divides} \textbf{(a)} says that
$11\mid f_{10}$ (since $11\equiv1\operatorname{mod}5$), and indeed we have
$f_{10}=55=5\cdot11$.
\end{itemize}

Our proof of Theorem \ref{thm.fib.p-divides} will be inspired by the famous
\textbf{Binet formula} for Fibonacci numbers:

\begin{theorem}
[Binet formula for Fibonacci numbers]\label{thm.fib.binet}Let
\[
\varphi=\dfrac{1+\sqrt{5}}{2}\approx1.618\ \ \ \ \ \ \ \ \ \ \text{and}%
\ \ \ \ \ \ \ \ \ \ \psi=\dfrac{1-\sqrt{5}}{2}\approx-0.618
\]
be the two roots of the quadratic polynomial $x^{2}-x-1$ in $\mathbb{R}$. Let
$\left(  f_{0},f_{1},f_{2},\ldots\right)  $ be the Fibonacci sequence. Then,%
\[
f_{n}=\dfrac{1}{\sqrt{5}}\varphi^{n}-\dfrac{1}{\sqrt{5}}\psi^{n}%
\ \ \ \ \ \ \ \ \ \ \text{for each }n\in\mathbb{N}.
\]

\end{theorem}

This is somewhat mysterious -- why should irrational numbers like $\sqrt{5}$
appear in a formula for an integer sequence like $\left(  f_{0},f_{1}%
,f_{2},\ldots\right)  $ ? Proving Theorem \ref{thm.fib.binet} is an easy
exercise in strong induction\footnote{See \cite[Proof of Theorem 2.3.1]{mps}
for the proof in detail.}. Finding it is trickier -- the matrix approach from
Exercise \ref{exe.21hw1.6} can help here. Indeed, once you know that the
matrix $A=\left(
\begin{array}
[c]{cc}%
0 & 1\\
1 & 1
\end{array}
\right)  \in\mathbb{R}^{2\times2}$ satisfies $A^{n}=f_{n}A+f_{n-1}I_{2}$ for
each $n\geq1$ (this was proven in Exercise \ref{exe.21hw1.6}), you can boil
down the computation of $f_{n}$ to the computation of $A^{n}$. But there is a
famous trick for computing powers of a matrix: namely, you diagonalize the
matrix and take the powers of its diagonal entries\footnote{Namely: If
$A=QDQ^{-1}$, then $A^{n}=QD^{n}Q^{-1}$ for any $n\in\mathbb{N}$. If the
matrix $D$ is diagonal, then $D^{n}$ is easily computed by taking its diagonal
entries to the $n$-th powers; thus, $A^{n}$ can be obtained as well.}. This
trick only works if the matrix is diagonalizable; but fortunately, our matrix
$A$ is diagonalizable, so we can compute $A^{n}$ using this trick, ultimately
obtaining Theorem \ref{thm.fib.binet} stated above. This demystifies the
formula: $x^{2}-x-1$ is just the characteristic polynomial of the matrix $A$,
and $\varphi$ and $\psi$ are its eigenvalues.

Anyway, how does this help us prove Theorem \ref{thm.fib.p-divides}? The Binet
formula involves irrational numbers and division; we thus cannot directly draw
any conclusions about divisibility from it.

We can, however, use it as an inspiration. To wit, we shall introduce
analogues of $\varphi$ and $\psi$ in \textquotedblleft characteristic
$p$\textquotedblright. These should be roots of the same polynomial
$x^{2}-x-1$, but regarded as a polynomial over $\mathbb{Z}/p$ instead of
$\mathbb{R}$. Depending on $p$, this polynomial may or may not have roots in
$\mathbb{Z}/p$, but we can always construct a splitting field in which it will
have roots (see Theorem \ref{thm.finfield.splitfield} \textbf{(c)}). Let us
use this to attempt a proof of Theorem \ref{thm.fib.p-divides}:

\begin{proof}
[Proof of Theorem \ref{thm.fib.p-divides}, part 1.]First, we WLOG assume that
$p\neq5$ (since Theorem \ref{thm.fib.p-divides} makes no statement about
$p=5$). Hence, $p\nmid5$ (since $p$ is prime), so that $\overline{5}%
\neq\overline{0}$ in $\mathbb{Z}/p$. Furthermore, from $p\neq5$, we obtain
$5\nmid p$ (since $p$ is prime); thus, the remainder of $p$ upon division by
$5$ must be $1$, $2$, $3$ or $4$. Therefore, $p$ must satisfy one of the
conditions $p\equiv\pm1\operatorname{mod}5$ and $p\equiv\pm2\operatorname{mod}%
5$.

Let $F$ be a splitting field of the polynomial $x^{2}-x-1$ over $\mathbb{Z}%
/p$. (We know from Theorem \ref{thm.finfield.splitfield} \textbf{(c)} that
such an $F$ exists, since the polynomial is monic.) Thus,%
\[
x^{2}-x-1=\left(  x-\varphi\right)  \left(  x-\psi\right)
\ \ \ \ \ \ \ \ \ \ \text{for some }\varphi,\psi\in F.
\]
Consider these $\varphi,\psi$. Comparing coefficients in front of the
monomials $x^{1}$ and $x^{0}$ in the polynomial identity
\[
x^{2}-x-1=\left(  x-\varphi\right)  \left(  x-\psi\right)  =x^{2}-\left(
\varphi+\psi\right)  x+\varphi\psi
\]
yields\footnote{This is perhaps a good time to recall the warnings about
evaluating polynomials over finite fields. Two polynomials $f$ and $g$ over a
finite field $F$ do not need to be identical just because their evaluations at
all elements of $F$ are identical (for example, the polynomials $x^{2}$ and
$x$ over $\mathbb{Z}/2$ are not identical, but their evaluations on both
elements $\overline{0}$ and $\overline{1}$ of $\mathbb{Z}/2$ are identical).
However, our two polynomials $x^{2}-x-1$ and $x^{2}-\left(  \varphi
+\psi\right)  x+\varphi\psi$ (whose coefficients we are comparing here) are
known to be identical (not just their evaluations but the polynomials
themselves); thus, we can compare their coefficients.}%
\[
-1=-\left(  \varphi+\psi\right)  \ \ \ \ \ \ \ \ \ \ \text{and}%
\ \ \ \ \ \ \ \ \ \ -1=\varphi\psi.
\]
(Of course, the \textquotedblleft$1$\textquotedblright\ here stands for
$1_{F}$.) In other words,%
\[
\varphi+\psi=1\ \ \ \ \ \ \ \ \ \ \text{and}\ \ \ \ \ \ \ \ \ \ \varphi
\psi=-1.
\]

Define an element $\sqrt{5}$ of $F$ by $\sqrt{5}=\varphi-\psi$. This is
certainly a strange notation (this $\sqrt{5}$ is not the actual number
$\sqrt{5}$ but just an analogue of it in our field $F$), but it is harmless
(as we won't deal with the actual number $\sqrt{5}$ in this proof, but only
with the element $\sqrt{5}=\varphi-\psi$ that we just introduced). Moreover,
it is justified because%
\[
\left(  \varphi-\psi\right)  ^{2}=\varphi^{2}-2\varphi\psi+\psi^{2}=\left(
\underbrace{\varphi+\psi}_{=1}\right)  ^{2}-4\underbrace{\varphi\psi}%
_{=-1}=1^{2}-4\left(  -1\right)  =\overline{5}.
\]
As a consequence, $\left(  \sqrt{5}\right)  ^{2}=\overline{5}\neq\overline{0}%
$, so that $\sqrt{5}\neq\overline{0}$. Thus, $\sqrt{5}$ is a unit of $F$
(since $F$ is a field), so we can divide by $\sqrt{5}$.

Now, we claim that an analogue of the Binet formula holds in $F$: Namely, we
have%
\begin{equation}
\overline{f_{n}}=\dfrac{1}{\sqrt{5}}\varphi^{n}-\dfrac{1}{\sqrt{5}}\psi
^{n}\ \ \ \ \ \ \ \ \ \ \text{for each }n\in\mathbb{N}.
\label{pf.thm.fib.p-divides.binet}%
\end{equation}
This can be proved by the same strong induction argument as the original Binet
formula (Theorem \ref{thm.fib.binet}).

Now, we want to show that $p\mid f_{p-1}$ for some primes $p$ and that $p\mid
f_{p+1}$ for other primes $p$ (remember: we have already gotten rid of the
$p=5$ case). In other words, we want to show that $\overline{f_{p-1}}=0$ for
some primes $p$, and that $\overline{f_{p+1}}=0$ for other primes $p$. For
now, let us ignore the question of which primes $p$ satisfy which of these.

Here comes a trick that will look magical, but is actually an instance of a
general method. We have $\varphi^{2}-\varphi-1=0$ (since $\varphi$ is a root
of the polynomial $x^{2}-x-1$), so that $\varphi^{2}=\varphi+1$. Taking this
equality to the $p$-th power, we obtain%
\begin{align*}
\varphi^{2p}  &  =\left(  \varphi+1\right)  ^{p}=\varphi^{p}+1^{p}%
\ \ \ \ \ \ \ \ \ \ \left(  \text{by Theorem \ref{thm.finfield.idiot}
\textbf{(a)}}\right) \\
&  =\varphi^{p}+1.
\end{align*}
In other words, $\varphi^{2p}-\varphi^{p}-1=0$. Thus, $\varphi^{p}$ is a root
of the polynomial $x^{2}-x-1=\left(  x-\varphi\right)  \left(  x-\psi\right)
$. In other words, $\left(  \varphi^{p}-\varphi\right)  \left(  \varphi
^{p}-\psi\right)  =0$. Since $F$ is an integral domain, this entails
$\varphi^{p}-\varphi=0$ or $\varphi^{p}-\psi=0$. In other words, $\varphi
^{p}=\varphi$ or $\varphi^{p}=\psi$. In other words, $\varphi^{p}\in\left\{
\varphi,\psi\right\}  $. Similarly, $\psi^{p}\in\left\{  \varphi,\psi\right\}
$.

Moreover, $\varphi-\psi=\sqrt{5}\neq\overline{0}$, so that $\varphi\neq\psi$
and therefore $\varphi^{p}\neq\psi^{p}$ (by Lemma \ref{lem.finfield.frob-inj},
since $F$ has characteristic $p$). Combining this with $\varphi^{p}\in\left\{
\varphi,\psi\right\}  $ and $\psi^{p}\in\left\{  \varphi,\psi\right\}  $, we
conclude that $\varphi^{p}$ and $\psi^{p}$ are two \textbf{distinct} elements
of the set $\left\{  \varphi,\psi\right\}  $. Thus, $\left\{  \varphi^{p}%
,\psi^{p}\right\}  =\left\{  \varphi,\psi\right\}  $. So we are in one of the
following two cases:

\textit{Case 1:} We have $\varphi^{p}=\varphi$ and $\psi^{p}=\psi$.

\textit{Case 2:} We have $\varphi^{p}=\psi$ and $\psi^{p}=\varphi$.

Let us consider Case 1. In this case, we have $\varphi^{p}=\varphi$ and
$\psi^{p}=\psi$. Now, $\varphi\neq0$ (since $\varphi^{2}=\varphi+1$ would turn
into the absurd equality $0=1$ if $\varphi$ were $0$); thus, we can cancel
$\varphi$ from the equality $\varphi^{p}=\varphi$ (since $F$ is a field). As a
result, we obtain $\varphi^{p-1}=1$. Similarly, $\psi^{p-1}=1$. Now,
(\ref{pf.thm.fib.p-divides.binet}) yields%
\[
\overline{f_{p-1}}=\dfrac{1}{\sqrt{5}}\underbrace{\varphi^{p-1}}_{=1}%
-\,\dfrac{1}{\sqrt{5}}\underbrace{\psi^{p-1}}_{=1}=\dfrac{1}{\sqrt{5}}%
\cdot1-\dfrac{1}{\sqrt{5}}\cdot1=0.
\]
Thus, we have shown that $\overline{f_{p-1}}=0$ (that is, $p\mid f_{p-1}$) in
Case 1.

Let us next consider Case 2. In this case, we have $\varphi^{p}=\psi$ and
$\psi^{p}=\varphi$. Thus, $\varphi^{p+1}=\underbrace{\varphi^{p}}_{=\psi
}\varphi=\psi\varphi=\varphi\psi=-1$ and similarly $\psi^{p+1}=-1$. Now,
(\ref{pf.thm.fib.p-divides.binet}) yields%
\[
\overline{f_{p+1}}=\dfrac{1}{\sqrt{5}}\underbrace{\varphi^{p+1}}%
_{=-1}-\,\dfrac{1}{\sqrt{5}}\underbrace{\psi^{p+1}}_{=-1}=\dfrac{1}{\sqrt{5}%
}\left(  -1\right)  -\dfrac{1}{\sqrt{5}}\left(  -1\right)  =0.
\]
Thus, we have shown that $\overline{f_{p+1}}=0$ (that is, $p\mid f_{p+1}$) in
Case 2.

So we have shown that we always have $p\mid f_{p-1}$ or $p\mid f_{p+1}$. But
why does the former hold for $p\equiv\pm1\operatorname{mod}5$ and the latter
for $p\equiv\pm2\operatorname{mod}5$ ? In other words, why does our Case 1
correspond to $p\equiv\pm1\operatorname{mod}5$ and our Case 2 to $p\equiv
\pm2\operatorname{mod}5$ ?

This will take some more work. We have the following chain of equivalences:%
\begin{align}
&  \ \left(  \text{we are in Case 1}\right) \nonumber\\
&  \Longleftrightarrow\ \left(  \varphi^{p}=\varphi\text{ and }\psi^{p}%
=\psi\right) \nonumber\\
&  \Longleftrightarrow\ \left(  \varphi^{p}=\varphi\right)
\ \ \ \ \ \ \ \ \ \ \left(
\begin{array}
[c]{c}%
\text{because if }\varphi^{p}=\varphi\text{,}\\
\text{then }\psi^{p}\text{ cannot be }\varphi\text{ (since }\varphi^{p}%
\neq\psi^{p}\text{)}\\
\text{and thus must be }\psi\text{ (since }\psi^{p}\in\left\{  \varphi
,\psi\right\}  \text{)}%
\end{array}
\right) \nonumber\\
&  \Longleftrightarrow\ \left(  \varphi\in\left\{  a\in F\ \mid\ a^{p}%
=a\right\}  \right) \nonumber\\
&  \Longleftrightarrow\ \left(  \varphi\in\mathbb{Z}/p\right)
\ \ \ \ \ \ \ \ \ \ \left(  \text{by Lemma \ref{lem.finfield.ap=a-only}%
}\right) \nonumber\\
&  \Longleftrightarrow\ \left(  \text{the polynomial }x^{2}-x-1\text{ has a
root in }\mathbb{Z}/p\right)  . \label{pf.thm.fib.p-divides.equiv1}%
\end{align}
(In the last equivalence sign, the \textquotedblleft$\Longrightarrow
$\textquotedblright\ part is obvious (since $\varphi$ is a root of $x^{2}%
-x-1$). The \textquotedblleft$\Longleftarrow$\textquotedblright\ part can be
proved as follows: If the polynomial $x^{2}-x-1$ has a root in $\mathbb{Z}/p$,
then this root must be either $\varphi$ or $\psi$ (because $x^{2}-x-1=\left(
x-\varphi\right)  \left(  x-\psi\right)  $); however, in either of these
cases, we obtain $\varphi\in\mathbb{Z}/p$ (because if $\psi\in\mathbb{Z}/p$,
then $\varphi=\underbrace{\left(  \varphi+\psi\right)  }_{=1\in\mathbb{Z}%
/p}-\underbrace{\psi}_{\in\mathbb{Z}/p}\in\mathbb{Z}/p$).)

Thus, our question is reduced to asking when the polynomial $x^{2}-x-1$ has a
root in $\mathbb{Z}/p$. In other words, when can we find our $\varphi$ and
$\psi$ in $\mathbb{Z}/p$, and when do we have to go into a larger field to
find them?

We WLOG assume that $p\neq2$ (since the case $p=2$ is trivial to do by hand).
Thus, $\overline{2}\in\mathbb{Z}/p$ is nonzero and thus has an inverse. This
allows us to complete the square (just as in high school, but over the field
$\mathbb{Z}/p$ now):%
\begin{equation}
x^{2}-x-1=\left(  x-\dfrac{\overline{1}}{\overline{2}}\right)  ^{2}%
-\dfrac{\overline{5}}{\overline{4}}. \label{pf.thm.fib.p-divides.cts}%
\end{equation}
Thus, the polynomial $x^{2}-x-1$ has a root in $\mathbb{Z}/p$ if and only if
$\dfrac{\overline{5}}{\overline{4}}$ is a square in $\mathbb{Z}/p$. Obviously,
$\dfrac{\overline{5}}{\overline{4}}$ is a square in $\mathbb{Z}/p$ if and only
if $\overline{5}$ is a square in $\mathbb{Z}/p$ (since $\overline{4}%
=\overline{2}^{2}$ is always a square in $\mathbb{Z}/p$). Thus, in order to
prove Theorem \ref{thm.fib.p-divides}, it remains to prove the following:
\end{proof}

\begin{theorem}
\label{thm.qr.5/p}Let $p$ be a prime such that $p\neq2$. Then:

\begin{enumerate}
\item[\textbf{(a)}] If $p\equiv\pm1\operatorname{mod}5$ (meaning that $p$ is
congruent to one of $1$ and $-1$ modulo $5$), then $\overline{5}$ is a square
in $\mathbb{Z}/p$.

\item[\textbf{(b)}] If $p\equiv\pm2\operatorname{mod}5$ (meaning that $p$ is
congruent to one of $2$ and $-2$ modulo $5$), then $\overline{5}$ is not a
square in $\mathbb{Z}/p$.
\end{enumerate}
\end{theorem}

For example, $\overline{5}\in\mathbb{Z}/p$ is not a square for $p=7$, but is a
square for $p=11$ (namely, $\overline{5}=\overline{4}^{2}$).

I will now prove Theorem \ref{thm.qr.5/p}; then, I will explain how it helps
complete the above proof of Theorem \ref{thm.fib.p-divides}, and afterwards
(perhaps most interestingly) discuss how to generalize it to other numbers
instead of $5$.

\begin{proof}
[Proof of Theorem \ref{thm.qr.5/p}.]The following proof (due to Gauss) will
again use field extensions. We WLOG assume that $p\neq5$ (since Theorem
\ref{thm.qr.5/p} makes no claim about the case $p=5$).

An element $z$ of a field $F$ is said to be a \textbf{primitive 5-th root of
unity} if it satisfies $z^{5}=1$ but $z\neq1$. In other words, the element $z$
is a primitive 5-th root of unity if it is nonzero and its order in the group
$F^{\times}$ (this is the group of units of $F$) is $5$.

For example, $\mathbb{R}$ has no primitive 5-th roots of unity (since a real
number $z$ satisfying $z^{5}=1$ must necessarily satisfy $z=1$), but
$\mathbb{C}$ has four of them: namely, $e^{2\pi ik/5}$ for $k\in\left\{
1,2,3,4\right\}  $. (See
\url{https://upload.wikimedia.org/wikipedia/commons/4/40/One5Root.svg} for an
illustration of the latter on the Argand diagram: The $5$ blue points, which
are the vertices of a regular pentagon, all satisfy $z^{5}=1$, and all but one
of them are primitive 5-th roots of unity.)

Does $\mathbb{Z}/p$ have any primitive 5-th roots of unity? Sometimes yes
(e.g., for $p=11$); sometimes no (e.g., for $p=7$). We don't care -- we shall
just adjoin one.

To see how, we notice the following: If $F$ is a field of characteristic $p$,
then a primitive 5-th root of unity in $F$ is just an element $z\in F$ that
satisfies $z^{4}+z^{3}+z^{2}+z+1=0$.\ \ \ \ \footnote{\textit{Proof.} If $z$
is a primitive 5-th root of unity in $F$, then $z^{5}=1$ but $z\neq1$, so that
$\dfrac{z^{5}-1}{z-1}=0$ (since the numerator $z^{5}-1$ is $0$ but the
denominator $z-1$ is nonzero), and therefore $z^{4}+z^{3}+z^{2}+z+1=0$ (since
$z^{4}+z^{3}+z^{2}+z+1=\dfrac{z^{5}-1}{z-1}$).
\par
Conversely, assume that $z^{4}+z^{3}+z^{2}+z+1=0$. Then, $z^{5}-1=\left(
z-1\right)  \underbrace{\left(  z^{4}+z^{3}+z^{2}+z+1\right)  }_{=0}=0$, so
that $z^{5}=1$. However, if we had $z=1$, then we would have $z^{4}%
+z^{3}+z^{2}+z+1=1^{4}+1^{3}+1^{2}+1+1=\overline{5}\neq\overline{0}$, which
would contradict $z^{4}+z^{3}+z^{2}+z+1=0=\overline{0}$. Hence, we must have
$z\neq1$. Thus we have shown that $z^{5}=1$ and $z\neq1$; in other words, $z$
is a primitive 5-th root of unity.} Knowing this, we can easily adjoin a
primitive 5-th root of unity to $\mathbb{Z}/p$: Namely, $x^{4}+x^{3}%
+x^{2}+x+1\in\left(  \mathbb{Z}/p\right)  \left[  x\right]  $ is a monic
polynomial of degree $4$ over $\mathbb{Z}/p$. Thus, by Theorem
\ref{thm.finfield.splitfield} \textbf{(b)}, there exists a field that contains
$\mathbb{Z}/p$ as a subring and that contains a root of this polynomial. Let
$S$ be such a field, and let $z$ be this root. Thus, $z\in S$ satisfies
$z^{4}+z^{3}+z^{2}+z+1=0$, and therefore is a primitive 5-th root of unity (by
what we have just said). That is, we have $z^{5}=1$ and $z\neq1$.

Now comes the magic: Set $\tau=z-z^{2}-z^{3}+z^{4}\in S$. Then,%
\begin{align*}
\tau^{2}  &  =\left(  z-z^{2}-z^{3}+z^{4}\right)  ^{2}\\
&  =z^{2}+z^{4}+z^{6}+z^{8}-2zz^{2}-2zz^{3}+2zz^{4}+2z^{2}z^{3}-2z^{2}%
z^{4}-2z^{3}z^{4}\\
&  \ \ \ \ \ \ \ \ \ \ \ \ \ \ \ \ \ \ \ \ \left(  \text{by expanding the
square}\right) \\
&  =z^{2}+z^{4}+z^{6}+z^{8}-2z^{3}-2z^{4}+2z^{5}+2z^{5}-2z^{6}-2z^{7}\\
&  =z^{2}+z^{4}+z+z^{3}-2z^{3}-2z^{4}+\overline{2}+\overline{2}-2z-2z^{2}\\
&  \ \ \ \ \ \ \ \ \ \ \ \ \ \ \ \ \ \ \ \ \left(  \text{since }z^{5}=1\text{
and thus }z^{6}=z\text{ and }z^{7}=z^{2}\right) \\
&  =\overline{4}-\left(  z+z^{2}+z^{3}+z^{4}\right)  =\overline{5}%
-\underbrace{\left(  z^{4}+z^{3}+z^{2}+z+1\right)  }_{=0}=\overline{5}.
\end{align*}
Thus, $\tau$ is a \textquotedblleft square root\textquotedblright\ of
$\overline{5}$ in $S$ (meaning: an element of $S$ whose square is
$\overline{5}$). Hence, the only \textquotedblleft square
roots\textquotedblright\ of $\overline{5}$ in $S$ are $\tau$ and $-\tau
$\ \ \ \ \footnote{This is a particular case of the following general fact: If
$R$ is an integral domain, and if $u,v\in R$ satisfy $u^{2}=v$, then the only
\textquotedblleft square roots\textquotedblright\ of $v$ in $R$ are $u$ and
$-u$. (To check this, argue as follows: If $w$ is a square root of $v$ in $R$,
then $\left(  w-u\right)  \left(  w+u\right)  =\underbrace{w^{2}}%
_{=v}-\underbrace{u^{2}}_{=v}=v-v=0$, so that $w-u=0$ or $w+u=0$ (since $R$ is
an integral domain), so that $w=u$ or $w=-u$.)}.

This suggests that studying $\tau$ should help understand whether
$\overline{5}$ is a square in $\mathbb{Z}/p$. Indeed, if $\tau$ belongs to
$\mathbb{Z}/p$, then $\overline{5}$ is a square in $\mathbb{Z}/p$ (since
$\tau^{2}=\overline{5}$). Conversely (but less obviously), if $\tau$ does
\textbf{not} belong to $\mathbb{Z}/p$, then $\overline{5}$ is not a square in
$\mathbb{Z}/p$ (because the only \textquotedblleft square
roots\textquotedblright\ of $\overline{5}$ in $S$ are $\tau$ and $-\tau$, and
neither of them belongs to $\mathbb{Z}/p$\ \ \ \ \footnote{Indeed, from
$\tau\notin\mathbb{Z}/p$, we obtain $-\tau\notin\mathbb{Z}/p$ (since
otherwise, $\tau=-\left(  -\tau\right)  $ would yield $\tau\in\mathbb{Z}/p$%
).}). Now, how can we tell whether $\tau$ belongs to $\mathbb{Z}/p$ ?

Inspired by Lemma \ref{lem.finfield.ap=a-only}, we compute $\tau^{p}$. From
$\tau=z-z^{2}-z^{3}+z^{4}$, we obtain%
\[
\tau^{p}=\left(  z-z^{2}-z^{3}+z^{4}\right)  ^{p}=z^{p}-z^{2p}-z^{3p}+z^{4p}%
\]
(by parts \textbf{(a)} and \textbf{(c)} of Theorem \ref{thm.finfield.idiot},
applied several times). The right hand side of this can be greatly simplified
if you know the remainder of $p$ upon division by $5$. Indeed, we have
$z^{5}=1$, so that $z^{6}=z$ and $z^{7}=z^{2}$ and more generally
$z^{k}=z^{\ell}$ for any two integers $k$ and $\ell$ satisfying $k\equiv
\ell\operatorname{mod}5$. Hence, in order to simplify the right hand side, we
distinguish the following four cases:

\textit{Case 1:} We have $p\equiv1\operatorname{mod}5$.

\textit{Case 2:} We have $p\equiv2\operatorname{mod}5$.

\textit{Case 3:} We have $p\equiv3\operatorname{mod}5$.

\textit{Case 4:} We have $p\equiv4\operatorname{mod}5$.

(There is no Case 0, since $5\nmid p$ entails $p\not \equiv
0\operatorname{mod}5$.)

In Case 2, we have%
\begin{align*}
\tau^{p}  &  =\underbrace{z^{p}}_{\substack{=z^{2}\\\text{(since }%
p\equiv2\operatorname{mod}5\text{)}}}-\underbrace{z^{2p}}_{\substack{=z^{4}%
\\\text{(since }2p\equiv4\operatorname{mod}5\text{)}}}-\underbrace{z^{3p}%
}_{\substack{=z^{1}\\\text{(since }3p\equiv1\operatorname{mod}5\text{)}%
}}+\underbrace{z^{4p}}_{\substack{=z^{3}\\\text{(since }4p\equiv
3\operatorname{mod}5\text{)}}}\\
&  =z^{2}-z^{4}-z^{1}+z^{3}=-\underbrace{\left(  z-z^{2}-z^{3}+z^{4}\right)
}_{=\tau}=-\tau.
\end{align*}
Similarly, we get $\tau^{p}=-\tau$ in Case 3, and we get $\tau^{p}=\tau$ in
Cases 1 and 4.

Thus, in Cases 1 and 4, we have $\tau^{p}=\tau$ and therefore $\tau\in\left\{
a\in S\ \mid\ a^{p}=a\right\}  =\mathbb{Z}/p$ (by Lemma
\ref{lem.finfield.ap=a-only}), and thus $\overline{5}$ is a square in
$\mathbb{Z}/p$ (since $\tau^{2}=\overline{5}$). On the other hand, in Cases 2
and 3, we have $\tau^{p}=-\tau\neq\tau$ (since $2\tau\neq0$%
\ \ \ \ \footnote{This can be shown as follows: From $\tau^{2}=\overline
{5}\neq0$, we obtain $\tau\neq0$. Moreover, $p\neq2$ shows that $\overline
{2}\neq0$ in $\mathbb{Z}/p$. Now, $S$ is an integral domain; hence, from
$\overline{2}\neq0$ and $\tau\neq0$, we obtain $\overline{2}\tau\neq0$. In
other words, $2\tau\neq0$.}) and therefore $\tau\notin\left\{  a\in
S\ \mid\ a^{p}=a\right\}  =\mathbb{Z}/p$ (by Lemma
\ref{lem.finfield.ap=a-only}), and thus $\overline{5}$ is not a square in
$\mathbb{Z}/p$ (as explained above). This proves Theorem \ref{thm.qr.5/p}.
\end{proof}

The \textquotedblleft magical\textquotedblright\ use of $z$ (a primitive 5-th
root of unity) to construct a square root of $5$ is connected to the ubiquity
of $\sqrt{5}$ in
\href{https://en.wikipedia.org/wiki/Pentagon#Regular_pentagons}{the geometry
of regular pentagons}. But it is not specific to the number $5$: Gauss showed
that for any prime $p>2$, we can similarly construct the square root
$\sqrt{\left(  -1\right)  ^{\left(  p-1\right)  /2}p}$ from a primitive $p$-th
root of unity. See Theorem \ref{thm.qr.gausssum2} \textbf{(c)} below for this
construction.\footnote{Another exposition can be found in \cite[\S 4.4]%
{Stein09} (where it is only stated for the field $\mathbb{C}$, but most other
fields can be handled similarly).
\par
Just to whet your appetite a bit more: Recall that the equilateral triangle
inscribed in the unit circle has sidelength $\sqrt{3}$. This is the $p=3$ case
of this connection.}

Next, let us use Theorem \ref{thm.qr.5/p} to complete our above proof of
Theorem \ref{thm.fib.p-divides}:

\begin{proof}
[Proof of Theorem \ref{thm.fib.p-divides}, part 2.]Recall the two Cases 1 and
2 that appeared in part 1 of this proof. We extend the equivalence
(\ref{pf.thm.fib.p-divides.equiv1}) as follows:%
\begin{align}
&  \ \ \left(  \text{we are in Case 1}\right) \nonumber\\
&  \Longleftrightarrow\ \left(  \text{the polynomial }x^{2}-x-1\text{ has a
root in }\mathbb{Z}/p\right) \nonumber\\
&  \Longleftrightarrow\ \left(  \text{the polynomial }\left(  x-\dfrac
{\overline{1}}{\overline{2}}\right)  ^{2}-\dfrac{\overline{5}}{\overline{4}%
}\text{ has a root in }\mathbb{Z}/p\right) \nonumber\\
&  \ \ \ \ \ \ \ \ \ \ \ \ \ \ \ \ \ \ \ \ \left(  \text{by
(\ref{pf.thm.fib.p-divides.cts})}\right) \nonumber\\
&  \Longleftrightarrow\ \left(  \dfrac{\overline{5}}{\overline{4}}\text{ is a
square in }\mathbb{Z}/p\right) \nonumber\\
&  \Longleftrightarrow\ \left(  \overline{5}\text{ is a square in }%
\mathbb{Z}/p\right) \nonumber\\
&  \ \ \ \ \ \ \ \ \ \ \ \ \ \ \ \ \ \ \ \ \left(  \text{since }%
\dfrac{\overline{5}}{\overline{4}}=a^{2}\text{ is equivalent to }\overline
{5}=\left(  2a\right)  ^{2}\right) \nonumber\\
&  \Longleftrightarrow\ \left(  p\equiv\pm1\operatorname{mod}5\right)
\label{pf.thm.fib.p-divides.equiv2}%
\end{align}
(by Theorem \ref{thm.qr.5/p}, since $p$ must satisfy one of the conditions
$p\equiv\pm1\operatorname{mod}5$ and $p\equiv\pm2\operatorname{mod}5$). But we
have shown that if we are in Case 1, then $p\mid f_{p-1}$. Thus, we conclude
using (\ref{pf.thm.fib.p-divides.equiv2}) that if $p\equiv\pm
1\operatorname{mod}5$, then $p\mid f_{p-1}$. This proves Theorem
\ref{thm.fib.p-divides} \textbf{(a)}. Likewise, if $p\equiv\pm
2\operatorname{mod}5$, then we do \textbf{not} have $p\equiv\pm
1\operatorname{mod}5$, so that we are \textbf{not} in Case 1 (by
(\ref{pf.thm.fib.p-divides.equiv2})), and thus we are in Case 2; hence, as we
proved above, we must have $p\mid f_{p+1}$ in this case. Thus, Theorem
\ref{thm.fib.p-divides} \textbf{(b)} is proved.
\end{proof}

\subsection{\label{sec.finfields.qr1}Quadratic residues: an introduction}

\subsubsection{Definitions and examples}

We have touched upon an interesting subject, so let us delve deeper. Theorem
\ref{thm.qr.5/p} answers the question for which primes $p$ the residue class
$\overline{5}\in\mathbb{Z}/p$ is a square in $\mathbb{Z}/p$; but we can ask
the same question about the residue class $\overline{a}$ of any $a\in
\mathbb{Z}$.

\begin{definition}
\label{def.qr.qr}Let $p$ be a prime. Let $a$ be an integer not divisible by
$p$.

Then, $a$ is said to be a \textbf{quadratic residue modulo }$p$ (short: a
\textbf{QR mod }$p$) if the residue class $\overline{a}\in\mathbb{Z}/p$ is a
square (or, equivalently, if there is an integer $b$ such that $a\equiv
b^{2}\operatorname{mod}p$).

Otherwise, $a$ is said to be a \textbf{quadratic nonresidue modulo }$p$
(short: a \textbf{QNR mod }$p$).
\end{definition}

For instance, if $p=7$, then the three integers $1$, $2$ and $4$ are QRs mod
$7$ (since $\overline{1}=\overline{1}^{2}$ and $\overline{2}=\overline{3}^{2}$
and $\overline{4}=\overline{2}^{2}$ in $\mathbb{Z}/7$), and therefore any
integers that are congruent to any of these three integers modulo $7$ are QRs
mod $7$ as well. The integers $3$, $5$ and $6$ (and any integers congruent to
them modulo $7$) are QNRs mod $7$. Integers divisible by $7$ (such as $0$)
count neither as QRs nor as QNRs mod $7$.

\begin{definition}
\label{def.qr.legendre}Let $p$ be a prime. Let $a$ be an integer. The
\textbf{Legendre symbol} $\left(  \dfrac{a}{p}\right)  $ (do not mistake this
for a fraction! this is not a fraction!) is the integer defined as follows:%
\[
\left(  \dfrac{a}{p}\right)  =%
\begin{cases}
0, & \text{if }p\mid a;\\
1, & \text{if }a\text{ is a QR mod }p\text{;}\\
-1, & \text{if }a\text{ is a QNR mod }p.
\end{cases}
\]

\end{definition}

Note that the Legendre symbol $\left(  \dfrac{a}{p}\right)  $ depends only on
the prime $p$ and on the residue class $\overline{a}\in\mathbb{Z}/p$, but not
on the integer $a$ itself. In other words, if a prime $p$ and two integers $a$
and $b$ satisfy $a\equiv b\operatorname{mod}p$, then%
\begin{equation}
\left(  \dfrac{a}{p}\right)  =\left(  \dfrac{b}{p}\right)  .
\label{eq.qr.legendre.congruent}%
\end{equation}
For example, $\left(  \dfrac{-1}{p}\right)  =\left(  \dfrac{p-1}{p}\right)  $
for any prime $p$. Here are some more examples of Legendre symbols:

\begin{itemize}
\item $2$ is a QR mod $7$, since $2\equiv3^{2}\operatorname{mod}7$. Thus,
$\left(  \dfrac{2}{7}\right)  =1$.

\item $2$ is a QNR mod $5$, since the squares in $\mathbb{Z}/5$ are
$\overline{0},\overline{1},\overline{4}$. Thus, $\left(  \dfrac{2}{5}\right)
=-1$.

\item $-1$ is a QR mod $5$, since $-1\equiv2^{2}\operatorname{mod}5$. Thus,
$\left(  \dfrac{-1}{5}\right)  =1$.

\item $-1$ is a QNR mod $3$. Thus, $\left(  \dfrac{-1}{3}\right)  =-1$.

\item Theorem \ref{thm.qr.5/p} says that every prime $p\neq2$ satisfies%
\[
\left(  \dfrac{5}{p}\right)  =%
\begin{cases}
0, & \text{if }p=5;\\
1, & \text{if }p\equiv\pm1\operatorname{mod}5;\\
-1, & \text{if }p\equiv\pm2\operatorname{mod}5.
\end{cases}
\]

\end{itemize}

This might whet an appetite: can we find similarly simple expressions for
$\left(  \dfrac{2}{p}\right)  $ or $\left(  \dfrac{3}{p}\right)  $ or $\left(
\dfrac{-1}{p}\right)  $ ? What can we say about Legendre symbols in general?

\subsubsection{Counting squares}

We begin by counting the squares in $\mathbb{Z}/p$:

\begin{proposition}
\label{prop.qr.count}Let $p\neq2$ be a prime. Then:

\begin{enumerate}
\item[\textbf{(a)}] The number of nonzero squares in $\mathbb{Z}/p$ is
$\left(  p-1\right)  /2$.

\item[\textbf{(b)}] The number of squares in $\mathbb{Z}/p$ is $\left(
p+1\right)  /2$.

\item[\textbf{(c)}] The number of elements of $\mathbb{Z}/p$ that are not
squares is $\left(  p-1\right)  /2$.
\end{enumerate}
\end{proposition}

\begin{proof}
This is a particular case of something that was proved in the solution of
Exercise \ref{exe.21hw1.10} \textbf{(c)}; but let us show the proof here to
keep this section self-contained: \medskip

\textbf{(a)} Given any element $v\in\mathbb{Z}/p$, we define a \textbf{square
root} of $v$ to be an element $u\in\mathbb{Z}/p$ that satisfies $u^{2}=v$. It
is easy to see that each nonzero square $v\in\mathbb{Z}/p$ has exactly two
(distinct) square roots\footnote{\textit{Proof.} Let $v\in\mathbb{Z}/p$ be a
nonzero square. Then, $v=c^{2}$ for some $c\in\mathbb{Z}/p$. Consider this
$c$. Since $c^{2}=v$ is nonzero, we see that $c$ is nonzero. However,
$\overline{2}\in\mathbb{Z}/p$ is also nonzero (since $p\neq2$ is a prime), and
thus $\overline{2}\cdot c\neq0$ (since $\mathbb{Z}/p$ is an integral domain,
and since $\overline{2}$ and $c$ are nonzero). Thus, $c+c=2\underbrace{c}%
_{=1_{\mathbb{Z}/p}c}=\underbrace{2\cdot1_{\mathbb{Z}/p}}_{=\overline{2}%
}c=\overline{2}\cdot c\neq0$. Subtracting $c$ from both sides of this
non-equality, we obtain $c\neq-c$.
\par
Now, both $c$ and $-c$ are square roots of $v$ (since $c^{2}=v$ and $\left(
-c\right)  ^{2}=c^{2}=v$). Conversely, any square root $d$ of $v$ must be one
of $c$ and $-c$ (because it satisfies $d^{2}=v$ and thus $\left(  d-c\right)
\left(  d+c\right)  =\underbrace{d^{2}}_{=v}-\underbrace{c^{2}}_{=v}=v-v=0$,
which entails (since $\mathbb{Z}/p$ is an integral domain) that $d-c$ or $d+c$
must be $0$, which in turn means that $d$ is either $c$ or $-c$). This shows
that $c$ and $-c$ are the only square roots of $v$. Since $c\neq-c$, we thus
conclude that $v$ has exactly two square roots. Qed.}, and these two square
roots are themselves nonzero (since $\overline{0}^{2}=\overline{0}$). Thus,
each nonzero square $v\in\mathbb{Z}/p$ has exactly two nonzero square roots.

If $c\in\left\{  \overline{1},\overline{2},\ldots,\overline{p-1}\right\}  $,
then $c\neq0$ and thus $c^{2}\neq0$ (since $\mathbb{Z}/p$ is an integral
domain), so that $c^{2}$ is a nonzero square. Hence, there is a map%
\begin{align*}
\left\{  \overline{1},\overline{2},\ldots,\overline{p-1}\right\}   &
\rightarrow\left\{  \text{nonzero squares }v\in\mathbb{Z}/p\right\}  ,\\
c  &  \mapsto c^{2}.
\end{align*}
This map is a 2-to-1 correspondence (i.e., each element of the set
\newline$\left\{  \text{nonzero squares }v\in\mathbb{Z}/p\right\}  $ has
exactly two preimages under this map), because each nonzero square
$v\in\mathbb{Z}/p$ has exactly two nonzero square roots. Thus,
\[
\left\vert \left\{  \overline{1},\overline{2},\ldots,\overline{p-1}\right\}
\right\vert =2\cdot\left\vert \left\{  \text{nonzero squares }v\in
\mathbb{Z}/p\right\}  \right\vert .
\]
Therefore,%
\[
\left\vert \left\{  \text{nonzero squares }v\in\mathbb{Z}/p\right\}
\right\vert =\dfrac{1}{2}\cdot\underbrace{\left\vert \left\{  \overline
{1},\overline{2},\ldots,\overline{p-1}\right\}  \right\vert }_{=p-1}=\dfrac
{1}{2}\cdot\left(  p-1\right)  =\left(  p-1\right)  /2.
\]
This proves Proposition \ref{prop.qr.count} \textbf{(a)}. \medskip

\textbf{(b)} The squares in $\mathbb{Z}/p$ are of two kinds: the zero squares
and the nonzero squares. Of the former kind, there is only one (namely,
$\overline{0}$, which is a square because $\overline{0}=\overline{0}^{2}$). Of
the latter kind, there are $\left(  p-1\right)  /2$ (by Proposition
\ref{prop.qr.count} \textbf{(a)}). Thus, in total, there are $1+\left(
p-1\right)  /2=\left(  p+1\right)  /2$ squares in $\mathbb{Z}/p$. This proves
Proposition \ref{prop.qr.count} \textbf{(b)}. \medskip

\textbf{(c)} The ring $\mathbb{Z}/p$ has $p$ elements, and exactly $\left(
p+1\right)  /2$ of them are squares (by Proposition \ref{prop.qr.count}
\textbf{(b)}). Hence, exactly $p-\left(  p+1\right)  /2=\left(  p-1\right)
/2$ of its elements are not squares. This proves Proposition
\ref{prop.qr.count} \textbf{(c)}.
\end{proof}

\begin{noncompile}
Note that Proposition \ref{prop.qr.count} \textbf{(b)} is a particular case of
Exercise \ref{exe.21hw1.10} \textbf{(c)}.
\end{noncompile}

\subsubsection{Euler's QR criterion}

Next, we will show a simple yet surprising rule, which was discovered by Euler:

\begin{theorem}
[Euler's QR criterion]\label{thm.qr.euler} Let $p\neq2$ be a prime. Let $a$ be
an integer. Then,%
\[
\left(  \dfrac{a}{p}\right)  \equiv a^{\left(  p-1\right)  /2}%
\operatorname{mod}p.
\]

\end{theorem}

\begin{proof}
Since $p$ is prime and satisfies $p\neq2$, we see that $p$ is odd and $\geq3$.
Hence, $\left(  p-1\right)  /2$ is a positive integer. Thus, $0^{\left(
p-1\right)  /2}=0$.

We must prove that $\left(  \dfrac{a}{p}\right)  \equiv a^{\left(  p-1\right)
/2}\operatorname{mod}p$. If $p\mid a$, then this boils down to showing that
$0\equiv0\operatorname{mod}p$ (since $0^{\left(  p-1\right)  /2}=0$). Thus, we
WLOG assume that $p\nmid a$.

Let $u=\overline{a}\in\mathbb{Z}/p$; thus, $u$ is nonzero (since $p\nmid a$).
Hence, the definition of $\left(  \dfrac{a}{p}\right)  $ yields%
\[
\left(  \dfrac{a}{p}\right)  =%
\begin{cases}
1, & \text{if }a\text{ is a QR mod }p\text{;}\\
-1, & \text{if }a\text{ is a QNR mod }p
\end{cases}
\ \ =%
\begin{cases}
1, & \text{if }u\text{ is a square};\\
-1, & \text{if }u\text{ is not a square}%
\end{cases}
\]
(by the definition of QRs and QNRs) and thus%
\begin{equation}
\overline{\left(  \dfrac{a}{p}\right)  }=%
\begin{cases}
\overline{1}, & \text{if }u\text{ is a square};\\
\overline{-1}, & \text{if }u\text{ is not a square.}%
\end{cases}
\label{pf.thm.qr.euler.1}%
\end{equation}
Also,
\begin{equation}
\overline{a^{\left(  p-1\right)  /2}}=\overline{a}^{\left(  p-1\right)
/2}=u^{\left(  p-1\right)  /2} \label{pf.thm.qr.euler.2}%
\end{equation}
(since $\overline{a}=u$). Now, we have the following chain of equivalences.%
\begin{align*}
&  \ \left(  \left(  \dfrac{a}{p}\right)  \equiv a^{\left(  p-1\right)
/2}\operatorname{mod}p\right)  \ \ \ \ \ \ \ \ \ \ \left(  \text{this is the
claim we are proving}\right) \\
&  \Longleftrightarrow\ \left(  \overline{\left(  \dfrac{a}{p}\right)
}=\overline{a^{\left(  p-1\right)  /2}}\right)  \ \ \Longleftrightarrow
\ \ \left(  \overline{a^{\left(  p-1\right)  /2}}=\overline{\left(  \dfrac
{a}{p}\right)  }\right) \\
&  \Longleftrightarrow\ \left(  u^{\left(  p-1\right)  /2}=%
\begin{cases}
\overline{1}, & \text{if }u\text{ is a square};\\
\overline{-1}, & \text{if }u\text{ is not a square}%
\end{cases}
\right)
\end{align*}
(by (\ref{pf.thm.qr.euler.1}) and (\ref{pf.thm.qr.euler.2})). Hence, it
remains to prove that

\begin{itemize}
\item we have $u^{\left(  p-1\right)  /2}=\overline{1}$ if $u$ is a square;

\item we have $u^{\left(  p-1\right)  /2}=\overline{-1}$ if $u$ is not a square.
\end{itemize}

Equivalently, we shall prove the following three claims:

\begin{statement}
\textit{Claim 1:} Any nonzero element $v\in\mathbb{Z}/p$ satisfies $v^{\left(
p-1\right)  /2}=\overline{1}$ or $v^{\left(  p-1\right)  /2}=\overline{-1}$.
\end{statement}

\begin{statement}
\textit{Claim 2:} Any nonzero square $v\in\mathbb{Z}/p$ satisfies $v^{\left(
p-1\right)  /2}=\overline{1}$.
\end{statement}

\begin{statement}
\textit{Claim 3:} Any element $v\in\mathbb{Z}/p$ that is not a square
satisfies $v^{\left(  p-1\right)  /2}\neq\overline{1}$.
\end{statement}

This will prove the two bullet points we claimed above: The first bullet point
will follow from Claim 2, while the second will follow from Claims 1 and 3. So
it remains to prove the three Claims 1, 2 and 3.

\begin{proof}
[Proof of Claim 1.]Let $v\in\mathbb{Z}/p$ be a nonzero element. Then,
Proposition \ref{prop.ent.flt.Z/p} (applied to $v$ instead of $u$) yields
$v^{p}=v$. We can cancel $v$ from this equality (since $v$ is nonzero and
$\mathbb{Z}/p$ is a field), and thus obtain $v^{p-1}=1$. Since $p-1$ is even,
we have $\left(  v^{\left(  p-1\right)  /2}\right)  ^{2}=v^{p-1}=1$, so that
$\left(  v^{\left(  p-1\right)  /2}\right)  ^{2}-1=0$. In view of $\left(
v^{\left(  p-1\right)  /2}\right)  ^{2}-1=\left(  v^{\left(  p-1\right)
/2}-1\right)  \left(  v^{\left(  p-1\right)  /2}+1\right)  $, this rewrites as
$\left(  v^{\left(  p-1\right)  /2}-1\right)  \left(  v^{\left(  p-1\right)
/2}+1\right)  =0$. Since $\mathbb{Z}/p$ is an integral domain, this entails
$v^{\left(  p-1\right)  /2}-1=0$ or $v^{\left(  p-1\right)  /2}+1=0$. In other
words, $v^{\left(  p-1\right)  /2}=\overline{1}$ or $v^{\left(  p-1\right)
/2}=\overline{-1}$. This proves Claim 1.
\end{proof}

\begin{proof}
[Proof of Claim 2.]Let $v\in\mathbb{Z}/p$ be a nonzero square. Thus, $v=w^{2}$
for some $w\in\mathbb{Z}/p$. Consider this $w$. Now, $w\neq0$ (since $w^{2}=v$
is nonzero). But Proposition \ref{prop.ent.flt.Z/p} (applied to $w$ instead of
$u$) yields $w^{p}=w$. We can cancel $w$ from this equality (since $w\neq0$
and $\mathbb{Z}/p$ is a field), and thus obtain $w^{p-1}=1$. Now, from
$v=w^{2}$, we obtain $v^{\left(  p-1\right)  /2}=\left(  w^{2}\right)
^{\left(  p-1\right)  /2}=w^{p-1}=1=\overline{1}$. This proves Claim 2.
\end{proof}

\begin{proof}
[Proof of Claim 3.]Here we take a bird's eye view (as in our above proof of
Lemma \ref{lem.finfield.ap=a-only}), rather than treating a single element
$v$. Indeed, $\mathbb{Z}/p$ is an integral domain. Thus, the easy half of the
FTA (Theorem \ref{thm.polring.univar-easyFTA}) yields that if $n$ is a
nonnegative integer, then any nonzero polynomial of degree $\leq n$ over
$\mathbb{Z}/p$ has at most $n$ roots in $\mathbb{Z}/p$. Applying this to the
polynomial $x^{\left(  p-1\right)  /2}-1$ (which is nonzero and has degree
$\left(  p-1\right)  /2$), we conclude that the polynomial $x^{\left(
p-1\right)  /2}-1$ has at most $\left(  p-1\right)  /2$ roots in
$\mathbb{Z}/p$. But the set of all roots of this polynomial $x^{\left(
p-1\right)  /2}-1$ in $\mathbb{Z}/p$ is $\left\{  v\in\mathbb{Z}%
/p\ \mid\ v^{\left(  p-1\right)  /2}=\overline{1}\right\}  $; hence, the
preceding sentence says that $\left\vert \left\{  v\in\mathbb{Z}%
/p\ \mid\ v^{\left(  p-1\right)  /2}=\overline{1}\right\}  \right\vert
\leq\left(  p-1\right)  /2$.

On the other hand, $\left\{  \text{nonzero squares }v\in\mathbb{Z}/p\right\}
\subseteq\left\{  v\in\mathbb{Z}/p\ \mid\ v^{\left(  p-1\right)  /2}%
=\overline{1}\right\}  $ (by Claim 2) and $\left\vert \left\{  \text{nonzero
squares }v\in\mathbb{Z}/p\right\}  \right\vert =\left(  p-1\right)  /2$ (by
Proposition \ref{prop.qr.count} \textbf{(a)}).

However, an easy and fundamental fact in combinatorics says that if $X$ and
$Y$ are two finite sets with $X\subseteq Y$ and $\left\vert Y\right\vert
\leq\left\vert X\right\vert $, then $X=Y$. Applying this to $X=\left\{
\text{nonzero squares }v\in\mathbb{Z}/p\right\}  $ and $Y=\left\{
v\in\mathbb{Z}/p\ \mid\ v^{\left(  p-1\right)  /2}=\overline{1}\right\}  $, we
obtain
\[
\left\{  \text{nonzero squares }v\in\mathbb{Z}/p\right\}  =\left\{
v\in\mathbb{Z}/p\ \mid\ v^{\left(  p-1\right)  /2}=\overline{1}\right\}
\]
(since $\left\{  \text{nonzero squares }v\in\mathbb{Z}/p\right\}
\subseteq\left\{  v\in\mathbb{Z}/p\ \mid\ v^{\left(  p-1\right)  /2}%
=\overline{1}\right\}  $ and \newline$\left\vert \left\{  v\in\mathbb{Z}%
/p\ \mid\ v^{\left(  p-1\right)  /2}=\overline{1}\right\}  \right\vert
\leq\left(  p-1\right)  /2=\left\vert \left\{  \text{nonzero squares }%
v\in\mathbb{Z}/p\right\}  \right\vert $). Thus, every $v\in\mathbb{Z}/p$
satisfying $v^{\left(  p-1\right)  /2}=\overline{1}$ must be a nonzero square.
By taking the contrapositive of this statement, we obtain Claim 3.
\end{proof}

Having proved Claims 1, 2 and 3, we thus have completed the proof of Theorem
\ref{thm.qr.euler}.
\end{proof}

We note that Theorem \ref{thm.qr.euler} has a generalization to squares in
arbitrary finite fields:

\begin{exercise}
\label{exe.qr.euler-ff}Let $F$ be a finite field of size $q$, where $q$ is
odd. Let $u\in F$ be a nonzero element. Prove that:

\begin{enumerate}
\item[\textbf{(a)}] If $u$ is a square in $F$ (that is, if $u=v^{2}$ for some
$v\in F$), then $u^{\left(  q-1\right)  /2}=1$.

\item[\textbf{(b)}] If $u$ is not a square in $F$, then $u^{\left(
q-1\right)  /2}=-1$.
\end{enumerate}
\end{exercise}

\subsubsection{The arithmetic of Legendre symbols}

Euler's criterion has a surprising corollary:

\begin{corollary}
[Multiplicativity of the Legendre symbol]\label{cor.qr.mult} Let $p\neq2$ be a
prime. Let $a,b\in\mathbb{Z}$. Then,%
\[
\left(  \dfrac{ab}{p}\right)  =\left(  \dfrac{a}{p}\right)  \left(  \dfrac
{b}{p}\right)  .
\]

\end{corollary}

To prove this, we will need the following near-trivial lemma (yes, we will;
just wait):

\begin{lemma}
\label{lem.qr.01-1}Let $p\neq2$ be a prime. Let $u,v\in\left\{
0,1,-1\right\}  $ be two integers satisfying $u\equiv v\operatorname{mod}p$.
Then, $u=v$.
\end{lemma}

\begin{proof}
We have $p>2$ (since $p\neq2$ and since $p$ is a prime). Thus, neither $1$ nor
$2$ nor $-1$ nor $-2$ is divisible by $p$. Therefore, the three integers
$0,1,-1$ are pairwise incongruent\footnote{\textquotedblleft
Incongruent\textquotedblright\ means \textquotedblleft not
congruent\textquotedblright.} modulo $p$. In other words, if two of them are
congruent modulo $p$, then these two integers are just equal. Hence, from
$u\equiv v\operatorname{mod}p$, we obtain $u=v$ (since $u,v\in\left\{
0,1,-1\right\}  $). This proves Lemma \ref{lem.qr.01-1}.
\end{proof}

\begin{proof}
[Proof of Corollary \ref{cor.qr.mult}.]Theorem \ref{thm.qr.euler} yields%
\begin{equation}
\left(  \dfrac{ab}{p}\right)  \equiv\left(  ab\right)  ^{\left(  p-1\right)
/2}=a^{\left(  p-1\right)  /2}b^{\left(  p-1\right)  /2}\operatorname{mod}p.
\label{pf.cor.qr.mult.1}%
\end{equation}
But Theorem \ref{thm.qr.euler} also yields $\left(  \dfrac{a}{p}\right)
\equiv a^{\left(  p-1\right)  /2}\operatorname{mod}p$ and $\left(  \dfrac
{b}{p}\right)  \equiv b^{\left(  p-1\right)  /2}\operatorname{mod}p$.
Multiplying these two congruences, we obtain%
\[
\left(  \dfrac{a}{p}\right)  \left(  \dfrac{b}{p}\right)  \equiv a^{\left(
p-1\right)  /2}b^{\left(  p-1\right)  /2}\operatorname{mod}p.
\]
Comparing this congruence with (\ref{pf.cor.qr.mult.1}), we find
\begin{equation}
\left(  \dfrac{ab}{p}\right)  \equiv\left(  \dfrac{a}{p}\right)  \left(
\dfrac{b}{p}\right)  \operatorname{mod}p. \label{pf.cor.qr.mult.3}%
\end{equation}

But we want an equality, not a congruence! Luckily, the congruence
(\ref{pf.cor.qr.mult.3}) turns out to entail the equality. Indeed, both sides
of the congruence (\ref{pf.cor.qr.mult.3}) equal $0$ or $1$ or $-1$ (since any
Legendre symbol is either $0$ or $1$ or $-1$, and the same holds for a product
of Legendre symbols). Hence, their congruence implies their equality (by Lemma
\ref{lem.qr.01-1}). This proves Corollary \ref{cor.qr.mult}.
\end{proof}

Corollary \ref{cor.qr.mult} has two nice corollaries of its own:

\begin{corollary}
\label{cor.qr.mult2}Let $p\neq2$ be a prime. The map%
\begin{align*}
\left(  \mathbb{Z}/p\right)  ^{\times}  &  \rightarrow\left\{  1,-1\right\}
,\\
\overline{a}  &  \mapsto\left(  \dfrac{a}{p}\right)
\end{align*}
is a group morphism (i.e., a homomorphism of groups).
\end{corollary}

\begin{proof}
The map is well-defined, since (as we have explained above) $\left(  \dfrac
{a}{p}\right)  $ depends only on $p$ and on $\overline{a}\in\mathbb{Z}/p$ (but
not on $a$ itself). Let us now show that this map is a group morphism.

In order to show that a map between two groups is a group morphism, it
suffices to show that this map respects multiplication (this is well-known).
Thus, it suffices to show that our map respects multiplication. In other
words, it suffices to show that $\left(  \dfrac{ab}{p}\right)  =\left(
\dfrac{a}{p}\right)  \left(  \dfrac{b}{p}\right)  $ for any $a,b\in\mathbb{Z}$
that are not divisible by $p$ (since any two elements of the group $\left(
\mathbb{Z}/p\right)  ^{\times}$ can be written in the forms $\overline{a}$ and
$\overline{b}$ for two such $a,b\in\mathbb{Z}$, and then their product will be
$\overline{ab}$). But this follows from Corollary \ref{cor.qr.mult}.
\end{proof}

\begin{corollary}
\label{cor.qr.mult3}Let $p\neq2$ be a prime. Let $u,v\in\mathbb{Z}/p$ be two
nonzero residue classes. Then:

\begin{enumerate}
\item[\textbf{(a)}] If $u$ and $v$ are squares, then $uv$ is a square.

\item[\textbf{(b)}] If only one of $u$ and $v$ is a square, then $uv$ is not a square.

\item[\textbf{(c)}] If none of $u$ and $v$ is a square, then $uv$ is a square.
\end{enumerate}
\end{corollary}

Note that Corollary \ref{cor.qr.mult3} \textbf{(c)} would fail if we replaced
$\mathbb{Z}/p$ by $\mathbb{Q}$. For example, none of the rational numbers $2$
and $3$ is a square, but neither is $2\cdot3$. But it does hold in
$\mathbb{Z}/p$ (as we shall now show), and (more generally) in finite fields,
as well as in $\mathbb{R}$ (since the non-squares in $\mathbb{R}$ are
precisely the negative reals, but a product of two negative reals is always positive).

\begin{proof}
[Proof of Corollary \ref{cor.qr.mult3}.]We shall only prove part \textbf{(c)},
for two reasons: First of all, parts \textbf{(a)} and \textbf{(b)} hold for
any field (unlike part \textbf{(c)}, as we just discussed), and can easily be
proved using nothing but the field axioms. Also, the proof we will give for
part \textbf{(c)} can easily be adapted to the other two parts.

\textbf{(c)} Assume that none of $u$ and $v$ is a square. Write $u$ and $v$ in
the form $u=\overline{a}$ and $v=\overline{b}$ for some integers $a$ and $b$.
Then, $a$ is a QNR mod $p$ (since $\overline{a}=u$ is not a square and thus
nonzero), and thus $\left(  \dfrac{a}{p}\right)  =-1$ (by the definition of
the Legendre symbol). Similarly, $\left(  \dfrac{b}{p}\right)  =-1$. Hence,
Corollary \ref{cor.qr.mult} yields $\left(  \dfrac{ab}{p}\right)
=\underbrace{\left(  \dfrac{a}{p}\right)  }_{=-1}\underbrace{\left(  \dfrac
{b}{p}\right)  }_{=-1}=\left(  -1\right)  \left(  -1\right)  =1$. In other
words, $ab$ is a QR mod $p$ (by the definition of the Legendre symbol). In
other words, $\overline{ab}$ is nonzero and a square. In view of
$\overline{ab}=\overline{a}\cdot\overline{b}=uv$ (since $\overline{a}=u$ and
$\overline{b}=v$), this yields that $uv$ is a square. Thus, Corollary
\ref{cor.qr.mult3} \textbf{(c)} is proven.
\end{proof}

\begin{exercise}
Let $F$ be a finite field. Prove that the polynomial $x^{4}+1\in F\left[
x\right]  $ is not irreducible. \medskip

[\textbf{Hint:} It suffices to consider the case $F=\mathbb{Z}/p$ for a prime
$p$. In this case, assume further that $p\neq2$, since the $p=2$ case is
easily done by hand. Use Corollary \ref{cor.qr.mult3} \textbf{(c)} to show
that at least one of the residue classes $\overline{-1}$, $\overline{2}$ and
$\overline{-2}$ is a square in $\mathbb{Z}/p$. In each of these cases, factor
$x^{4}+1$.]
\end{exercise}

\subsubsection{When $-1$ is a QR}

Let us now return to the computation of Legendre symbols. Thanks to Corollary
\ref{cor.qr.mult}, we have (for example) $\left(  \dfrac{6}{p}\right)
=\left(  \dfrac{2}{p}\right)  \left(  \dfrac{3}{p}\right)  $ and $\left(
\dfrac{-6}{p}\right)  =\left(  \dfrac{-1}{p}\right)  \left(  \dfrac{2}%
{p}\right)  \left(  \dfrac{3}{p}\right)  $ for any prime $p$. But how do we
compute $\left(  \dfrac{-1}{p}\right)  $, $\left(  \dfrac{2}{p}\right)  $ and
$\left(  \dfrac{3}{p}\right)  $ ?

We begin with $\left(  \dfrac{-1}{p}\right)  $, which is probably the easiest one:

\begin{theorem}
\label{thm.qr.-1}Let $p\neq2$ be a prime. Then, $-1$ is a QR mod $p$ (that is,
$\overline{-1}\in\mathbb{Z}/p$ is a square) if and only if $p\equiv
1\operatorname{mod}4$. In other words,%
\[
\left(  \dfrac{-1}{p}\right)  =%
\begin{cases}
1, & \text{if }p\equiv1\operatorname{mod}4;\\
-1, & \text{if }p\equiv3\operatorname{mod}4.
\end{cases}
\]

\end{theorem}

\begin{proof}
Since $p$ is a prime satisfying $p\neq2$, the number $p$ is odd. Hence,
$\left(  p-1\right)  /2\in\mathbb{Z}$.

Theorem \ref{thm.qr.euler} yields the congruence%
\[
\left(  \dfrac{-1}{p}\right)  \equiv\left(  -1\right)  ^{\left(  p-1\right)
/2}\operatorname{mod}p.
\]
Lemma \ref{lem.qr.01-1} shows that this congruence must actually be an
equality, since both of its sides are $0$ or $1$ or $-1$. In other words,%
\begin{equation}
\left(  \dfrac{-1}{p}\right)  =\left(  -1\right)  ^{\left(  p-1\right)  /2}.
\label{pf.thm.qr.-1.2}%
\end{equation}
Now, $p$ must satisfy $p\equiv1\operatorname{mod}4$ or $p\equiv
3\operatorname{mod}4$ (since $p$ is odd). In the former case, $\left(
-1\right)  ^{\left(  p-1\right)  /2}$ is $1$; in the latter, $-1$. Hence,
(\ref{pf.thm.qr.-1.2}) can be rewritten as%
\[
\left(  \dfrac{-1}{p}\right)  =%
\begin{cases}
1, & \text{if }p\equiv1\operatorname{mod}4;\\
-1, & \text{if }p\equiv3\operatorname{mod}4.
\end{cases}
\]

\end{proof}

Theorem \ref{thm.qr.-1} has some unexpected applications. For instance, it can
be used to prove a strengthening of Euclid's famous result that there are
infinitely many prime numbers:

\begin{exercise}
\label{exe.dirichlet.1mod4}The \textbf{Fermat numbers} are the positive
integers $F_{0},F_{1},F_{2},\ldots$ defined by%
\[
F_{n}:=2^{2^{n}}+1\ \ \ \ \ \ \ \ \ \ \text{for each }n\in\mathbb{N}.
\]

\begin{enumerate}
\item[\textbf{(a)}] Prove that $\gcd\left(  F_{i},F_{j}\right)  =1$ for all
$i\neq j$.

\item[\textbf{(b)}] Prove that every prime $p$ that divides some Fermat number
$F_{n}$ with $n\geq1$ must satisfy $p\equiv1\operatorname{mod}4$.

\item[\textbf{(c)}] Conclude that there are infinitely many prime numbers $p$
satisfying $p\equiv1\operatorname{mod}4$.
\end{enumerate}
\end{exercise}

More generally,
\href{https://en.wikipedia.org/wiki/Dirichlet's_theorem_on_arithmetic_progressions}{\textbf{Dirichlet's
theorem on prime numbers in arithmetic progressions}} says that whenever $u$
and $v$ are two coprime integers satisfying $u>0$, then there are infinitely
many prime numbers $p$ satisfying $p\equiv v\operatorname{mod}u$. However,
this is far harder to prove; all known proofs require analysis (see, e.g.,
\cite[Chapter 1]{Klazar10a} and \cite[Chapter 2]{Klazar10b}). We will give a
proof of the $v=1$ case further below (see Theorem \ref{thm.dirichlet-1modn}).

\subsubsection{Quadratic reciprocity}

So we have simple formulas for $\left(  \dfrac{-1}{p}\right)  $ and $\left(
\dfrac{5}{p}\right)  $. What about $\left(  \dfrac{a}{p}\right)  $ for a
general $a$ ? We only need to know a formula for $\left(  \dfrac{q}{p}\right)
$ for each prime $q$ (because, as per Corollary \ref{cor.qr.mult} above, we
can then get a general formula for $\left(  \dfrac{a}{p}\right)  $ by
decomposing $a$ into a product of primes and possibly $-1$, and multiplying).
Here is one:

\begin{theorem}
[Quadratic Reciprocity Law]\label{thm.qr.qr}\ \ 

\begin{enumerate}
\item[\textbf{(a)}] Let $p\neq2$ be a prime. Then,%
\[
\left(  \dfrac{2}{p}\right)  =\left(  -1\right)  ^{\left(  p^{2}-1\right)
/8}=%
\begin{cases}
1, & \text{if }p\equiv\pm1\operatorname{mod}8;\\
-1, & \text{if }p\equiv\pm3\operatorname{mod}8.
\end{cases}
\]

\item[\textbf{(b)}] Let $p$ and $q$ be two distinct primes distinct from $2$.
Then,%
\[
\left(  \dfrac{q}{p}\right)  =\left(  -1\right)  ^{\left(  p-1\right)  \left(
q-1\right)  /4}\left(  \dfrac{p}{q}\right)  .
\]

\end{enumerate}
\end{theorem}

For example, if $p\neq2$ is any prime distinct from $5$, then Theorem
\ref{thm.qr.qr} \textbf{(b)} (applied to $q=5$) yields%
\begin{align*}
\left(  \dfrac{5}{p}\right)   &  =\underbrace{\left(  -1\right)  ^{\left(
p-1\right)  \left(  5-1\right)  /4}}_{\substack{=1\\\text{(since }5-1=4\text{
and thus}\\\left(  p-1\right)  \left(  5-1\right)  /4=p-1\text{ is even)}%
}}\left(  \dfrac{p}{5}\right)  =\left(  \dfrac{p}{5}\right) \\
&  =%
\begin{cases}
1, & \text{if }\overline{p}\in\mathbb{Z}/5\text{ is a square;}\\
-1, & \text{if }\overline{p}\in\mathbb{Z}/5\text{ is not a square}%
\end{cases}
\ \ =%
\begin{cases}
1, & \text{if }p\equiv\pm1\operatorname{mod}5;\\
-1, & \text{if }p\equiv\pm2\operatorname{mod}5
\end{cases}
\end{align*}
(the last equality follows from the fact that the nonzero squares in
$\mathbb{Z}/5$ are $\overline{1}$ and $\overline{-1}$); this recovers the
claim of Theorem \ref{thm.qr.5/p}. So Theorem \ref{thm.qr.5/p} was merely the
tip of an iceberg.

Theorem \ref{thm.qr.qr} is known as the \textbf{law of quadratic
reciprocity}\footnote{Some authors refer only to Theorem \ref{thm.qr.qr}
\textbf{(b)} as the law of quadratic reciprocity; they correspondingly call
Theorem \ref{thm.qr.qr} \textbf{(a)} the \textquotedblleft second
supplementary law of quadratic reciprocity\textquotedblright. (The name
\textquotedblleft first supplementary law\textquotedblright\ then refers to
Theorem \ref{thm.qr.-1}.)}, and is one of the most classical theorems in
mathematics -- discovered by Euler, proved by Gauss. By now, it has received
over 300 proofs (see
\url{https://web.archive.org/web/20250106010310/https://www.mathi.uni-heidelberg.de/~flemmermeyer/qrg_proofs.html}
for a list), and new proofs keep getting published. You'll get to prove its
part \textbf{(a)} in the following exercise, inspired by the $q=5$ case we
proved above:

\begin{exercise}
\label{exe.21hw4.8}Let $p$ be an odd prime. Let $\zeta$ be a root of the
polynomial $x^{4}+1\in\left(  \mathbb{Z}/p\right)  \left[  x\right]  $ in a
commutative $\mathbb{Z}/p$-algebra $A$. Thus, $\zeta^{4}=-1$, so that $\zeta$
is a unit (with inverse $-\zeta^{3}$). Let $\tau\in A$ be defined by
$\tau=\zeta+\zeta^{-1}$. Prove the following:

\begin{enumerate}
\item[\textbf{(a)}] We have $\tau^{2}=2$. (Here, $2$ stands for $2\cdot
1_{A}\in A$.)

\item[\textbf{(b)}] We have $\tau^{p}=\left(  \dfrac{2}{p}\right)  \tau$,
where $\left(  \dfrac{2}{p}\right)  $ means a Legendre symbol.

\item[\textbf{(c)}] If $p\equiv\pm1\operatorname{mod}8$ (that is, if $p$ is
congruent to $1$ or to $-1$ modulo $8$), then $\tau^{p}=\tau$.

\item[\textbf{(d)}] If $p\equiv\pm3\operatorname{mod}8$ (that is, if $p$ is
congruent to $3$ or to $-3$ modulo $8$), then $\tau^{p}=-\tau$.

\item[\textbf{(e)}] Prove Theorem \ref{thm.qr.qr} \textbf{(a)}.
\end{enumerate}

[\textbf{Hint:} For part \textbf{(b)}, start out by writing $\tau^{p}=\left(
\tau^{2}\right)  ^{\left(  p-1\right)  /2}\tau$.]
\end{exercise}

Hopefully, Exercise \ref{exe.21hw4.8} sheds some light on the strange
definition of $\tau$.

The rest of this section will be devoted to proving Theorem \ref{thm.qr.qr}
\textbf{(b)}. Before we embark on the actual proof, we shall show a few
auxiliary results, which are themselves of some interest.

\subsubsection{A sum of Legendre symbols}

The first auxiliary result is a classical formula for a certain sum of
Legendre symbols:

\begin{proposition}
\label{prop.qr.sumii-k}Let $p\neq2$ be a prime. Let $k\in\mathbb{Z}$. Then:

\begin{enumerate}
\item[\textbf{(a)}] If $p\mid k$, then $\sum_{i=0}^{p-1}\left(  \dfrac
{i\left(  i-k\right)  }{p}\right)  =p-1$.

\item[\textbf{(b)}] If $p\nmid k$, then $\sum_{i=0}^{p-1}\left(
\dfrac{i\left(  i-k\right)  }{p}\right)  =-1$.
\end{enumerate}
\end{proposition}

\begin{proof}
We first observe that $\left(  \dfrac{0\left(  0-k\right)  }{p}\right)
=\left(  \dfrac{0}{p}\right)  =0$ (since $p\mid0$). In other words, the $i=0$
addend of the sum $\sum_{i=0}^{p-1}\left(  \dfrac{i\left(  i-k\right)  }%
{p}\right)  $ is $0$. Hence, we can remove this addend from the sum, and
obtain%
\begin{equation}
\sum_{i=0}^{p-1}\left(  \dfrac{i\left(  i-k\right)  }{p}\right)  =\sum
_{i=1}^{p-1}\left(  \dfrac{i\left(  i-k\right)  }{p}\right)  .
\label{pf.prop.qr.sumii-k.-0}%
\end{equation}

We also note that $p$ is odd (since $p\neq2$ is a prime), and thus $\left(
p-1\right)  /2\in\mathbb{N}$. \medskip

\textbf{(a)} Assume that $p\mid k$. Then, $k\equiv0\operatorname{mod}p$, so
that $\overline{k}=\overline{0}$ in $\mathbb{Z}/p$.

Let $i\in\left\{  1,2,\ldots,p-1\right\}  $. Then, $p\nmid i$, so that
$\overline{i}\neq0$ in $\mathbb{Z}/p$. Hence, $\overline{i}^{2}\neq0$ in
$\mathbb{Z}/p$ as well (since $\mathbb{Z}/p$ is an integral domain). However,
in $\mathbb{Z}/p$, we have%
\[
\overline{i\left(  i-k\right)  }=\overline{i}\left(  \overline{i}%
-\underbrace{\overline{k}}_{=\overline{0}}\right)  =\overline{i}\left(
\overline{i}-\overline{0}\right)  =\overline{i}^{2}\neq0,
\]
so that $i\left(  i-k\right)  $ is not divisible by $p$. Furthermore, the
element $\overline{i\left(  i-k\right)  }$ of $\mathbb{Z}/p$ is a square
(since $\overline{i\left(  i-k\right)  }=\overline{i}^{2}$), so that $i\left(
i-k\right)  $ is a QR mod $p$ (since $i\left(  i-k\right)  $ is not divisible
by $p$). Therefore, $\left(  \dfrac{i\left(  i-k\right)  }{p}\right)  =1$.

Forget that we fixed $i$. We thus have proved the equality $\left(
\dfrac{i\left(  i-k\right)  }{p}\right)  =1$ for each $i\in\left\{
1,2,\ldots,p-1\right\}  $. Summing this equality over all $i\in\left\{
1,2,\ldots,p-1\right\}  $, we obtain $\sum_{i=1}^{p-1}\left(  \dfrac{i\left(
i-k\right)  }{p}\right)  =\sum_{i=1}^{p-1}1=p-1$. Hence,
(\ref{pf.prop.qr.sumii-k.-0}) can be rewritten as $\sum_{i=0}^{p-1}\left(
\dfrac{i\left(  i-k\right)  }{p}\right)  =p-1$. This proves Proposition
\ref{prop.qr.sumii-k} \textbf{(a)}. \medskip

\textbf{(b)} Assume that $p\nmid k$. We consider the polynomial
\[
f:=\left(  x\left(  x-k\right)  \right)  ^{\left(  p-1\right)  /2}-x^{p-1}%
\in\mathbb{Z}\left[  x\right]
\]
(this is well-defined, since $\left(  p-1\right)  /2\in\mathbb{N}$). We claim
that this polynomial $f$ has degree $\leq p-2$. Indeed, $f$ is the difference
of the two polynomials $\left(  x\left(  x-k\right)  \right)  ^{\left(
p-1\right)  /2}$ and $x^{p-1}$, both of which have degree $p-1$ and leading
coefficient $1$ (this is obvious for $x^{p-1}$, and for $\left(  x\left(
x-k\right)  \right)  ^{\left(  p-1\right)  /2}$ it follows from Proposition
\ref{prop.polring.univar-degpq} \textbf{(d)}\footnote{In more detail: An easy
consequence of Proposition \ref{prop.polring.univar-degpq} \textbf{(d)} is
that if $g$ is a monic polynomial of degree $i$, then $g^{m}$ is a monic
polynomial of degree $mi$ for each $m\in\mathbb{N}$. Applying this to
$g=x\left(  x-k\right)  $ and $i=2$ and $m=\left(  p-1\right)  /2$, we
conclude that $\left(  x\left(  x-k\right)  \right)  ^{\left(  p-1\right)
/2}$ is a monic polynomial of degree $\left(  p-1\right)  /2\cdot2=p-1$. In
other words, $\left(  x\left(  x-k\right)  \right)  ^{\left(  p-1\right)  /2}$
is a polynomial with degree $p-1$ and leading coefficient $1$.}). When we
subtract these two polynomials, the leading terms cancel out (since the
leading coefficients are equal), and thus we are left with a polynomial of
degree $\leq p-2$. In other words, the polynomial $f$ has degree $\leq p-2$.

Hence, Corollary \ref{cor.ent.sum-f-over-p.1} yields $\sum_{j=0}^{p-1}f\left(
j\right)  \equiv0\operatorname{mod}p$. In other words,%
\[
\sum_{j=0}^{p-1}\left(  \left(  j\left(  j-k\right)  \right)  ^{\left(
p-1\right)  /2}-j^{p-1}\right)  \equiv0\operatorname{mod}p
\]
(since the definition of $f$ yields $f\left(  j\right)  =\left(  j\left(
j-k\right)  \right)  ^{\left(  p-1\right)  /2}-j^{p-1}$ for each
$j\in\mathbb{Z}$). In other words,%
\[
\sum_{j=0}^{p-1}\left(  j\left(  j-k\right)  \right)  ^{\left(  p-1\right)
/2}-\sum_{j=0}^{p-1}j^{p-1}\equiv0\operatorname{mod}p.
\]
Hence,%
\begin{align*}
\sum_{j=0}^{p-1}\left(  j\left(  j-k\right)  \right)  ^{\left(  p-1\right)
/2}  &  \equiv\sum_{j=0}^{p-1}j^{p-1}=\underbrace{0^{p-1}}%
_{\substack{=0\\\text{(since }p-1>0\text{)}}}+\sum_{j=1}^{p-1}%
\underbrace{j^{p-1}}_{\substack{\equiv1\operatorname{mod}p\\\text{(by
Corollary \ref{cor.ent.flt.p-1},}\\\text{since }p\nmid j\text{)}}}\\
&  \equiv0+\sum_{j=1}^{p-1}1=p-1\equiv-1\operatorname{mod}p.
\end{align*}

Now,
\begin{align}
\sum_{i=0}^{p-1}\underbrace{\left(  \dfrac{i\left(  i-k\right)  }{p}\right)
}_{\substack{\equiv\left(  i\left(  i-k\right)  \right)  ^{\left(  p-1\right)
/2}\operatorname{mod}p\\\text{(by Theorem \ref{thm.qr.euler},}\\\text{applied
to }a=i\left(  i-k\right)  \text{)}}}  &  \equiv\sum_{i=0}^{p-1}\left(
i\left(  i-k\right)  \right)  ^{\left(  p-1\right)  /2}=\sum_{j=0}%
^{p-1}\left(  j\left(  j-k\right)  \right)  ^{\left(  p-1\right)
/2}\nonumber\\
&  \equiv-1\operatorname{mod}p. \label{pf.prop.qr.sumii-k.b.cong}%
\end{align}

This is very close to the thing we want to prove, but we are not quite there
yet: We want to prove that the two sides of (\ref{pf.prop.qr.sumii-k.b.cong})
are identical, not just congruent modulo $p$. Thus, we need an extra argument.
Let $s:=\sum_{i=0}^{p-1}\left(  \dfrac{i\left(  i-k\right)  }{p}\right)  $.
This $s$ is a sum of $p$ Legendre symbols, each of which is either $0$ or $1$
or $-1$ (since any Legendre symbol is either $0$ or $1$ or $-1$). Hence, $s$
is an integer between $-p$ and $p$ (inclusive). However, we can show something
slightly better: We can show that $s$ is an integer between $-\left(
p-2\right)  $ and $p-2$.

In order to show this, we let $\ell$ be the remainder obtained when dividing
$k$ by $p$. Then, $\ell\in\left\{  0,1,\ldots,p-1\right\}  $ and $\ell\neq0$
(since $p\nmid k$) and $\ell\equiv k\operatorname{mod}p$. Hence, $\ell-k$ is
divisible by $p$, so that $\ell\left(  \ell-k\right)  $ is divisible by $p$.
Thus, $\left(  \dfrac{\ell\left(  \ell-k\right)  }{p}\right)  =0$ by
Definition \ref{def.qr.legendre}.

Now, the sum
\begin{equation}
s=\sum_{i=0}^{p-1}\left(  \dfrac{i\left(  i-k\right)  }{p}\right)
\label{pf.prop.qr.sumii-k.b.s=}%
\end{equation}
has at least two addends that equal $0$: namely, the $i=0$ addend (which is
$\left(  \dfrac{0\left(  0-k\right)  }{p}\right)  =0$) and the $i=\ell$ addend
(which is $\left(  \dfrac{\ell\left(  \ell-k\right)  }{p}\right)  =0$). These
are two distinct addends (since $\ell\neq0$). The remaining $p-2$ addends in
the sum $s$ are Legendre symbols, so each of them is either $0$ or $1$ or
$-1$. Hence, the entire sum $s$ is an integer between $-\left(  p-2\right)  $
and $p-2$ (inclusive).

We have now encircled $s$ from all sides: On the one hand, we know that
$s=\sum_{i=0}^{p-1}\left(  \dfrac{i\left(  i-k\right)  }{p}\right)
\equiv-1\operatorname{mod}p$ (by (\ref{pf.prop.qr.sumii-k.b.cong})), so that
$s$ is congruent to $-1$ modulo $p$. On the other hand, we know that $s$ is an
integer between $-\left(  p-2\right)  $ and $p-2$ (inclusive). However, the
only integer between $-\left(  p-2\right)  $ and $p-2$ that is congruent to
$-1$ modulo $p$ is the integer $-1$ itself (since $p-1$ is too large, while
$-p-1$ is too small). Thus, $s$ must be $-1$. In other words, $\sum
_{i=0}^{p-1}\left(  \dfrac{i\left(  i-k\right)  }{p}\right)  $ must be $-1$
(by (\ref{pf.prop.qr.sumii-k.b.s=})). This proves Proposition
\ref{prop.qr.sumii-k} \textbf{(b)}.
\end{proof}

The following exercise (a result of Jacobsthal from 1907, see \cite{Jacobs07})
generalizes Proposition \ref{prop.qr.sumii-k}:

\begin{exercise}
Let $p\neq2$ be a prime. Let $a,b\in\mathbb{Z}$. Prove that
\[
\sum_{i=0}^{p-1}\left(  \dfrac{i^{2}+ai+b}{p}\right)  =%
\begin{cases}
p-1, & \text{if }a^{2}\equiv4b\operatorname{mod}p;\\
-1, & \text{else.}%
\end{cases}
\]

\end{exercise}

\subsubsection{Gaussian sums}

In our above proof of Theorem \ref{thm.qr.5/p}, we did some seemingly
unmotivated things: We adjoined an element $z$ satisfying $z^{4}+z^{3}%
+z^{2}+z+1=0$ to our field (which was $\mathbb{Z}/p$, but this is not
important); then, we defined $\tau=z-z^{2}-z^{3}+z^{4}$; then, we showed (by
computation) that $\tau^{2}=\overline{5}$.

This element $\tau$ was not chosen at random; its definition can be rewritten
using Legendre symbols as
\[
\tau=\sum_{i=0}^{4}\left(  \dfrac{i}{5}\right)  z^{i}%
\]
(since $\sum_{i=0}^{4}\left(  \dfrac{i}{5}\right)  z^{i}=\underbrace{\left(
\dfrac{0}{5}\right)  }_{=0}z^{0}+\underbrace{\left(  \dfrac{1}{5}\right)
}_{=1}z^{1}+\underbrace{\left(  \dfrac{2}{5}\right)  }_{=-1}z^{2}%
+\underbrace{\left(  \dfrac{3}{5}\right)  }_{=-1}z^{3}+\underbrace{\left(
\dfrac{4}{5}\right)  }_{=1}z^{4}=0z^{0}+1z^{1}+\left(  -1\right)
z^{2}+\left(  -1\right)  z^{3}+1z^{4}=z-z^{2}-z^{3}+z^{4}$). This suggests a
generalization: Replacing $5$ by an arbitrary prime $p\neq2$, we can adjoin an
element $z$ satisfying $z^{p-1}+z^{p-2}+\cdots+z^{0}=0$ to our field (which,
as we said, can be arbitrary), and define an element $\tau=\sum_{i=0}%
^{p-1}\left(  \dfrac{i}{p}\right)  z^{i}$. The square $\tau^{2}$ of this
element will then be $\overline{\left(  -1\right)  ^{\left(  p-1\right)  /2}%
p}$ instead of $\overline{5}$. This was found by Gauss, and we shall prove
this as part of the following theorem:

\begin{theorem}
[square of the Gaussian sum]\label{thm.qr.gausssum2}Let $p$ be a prime. Let
$A$ be a commutative ring. Let $z\in A$ be an element satisfying
\begin{equation}
z^{p-1}+z^{p-2}+\cdots+z^{0}=0. \label{eq.thm.qr.gausssum2.1}%
\end{equation}
Then:

\begin{enumerate}
\item[\textbf{(a)}] We have $z^{p}=1_{A}$.

\item[\textbf{(b)}] If $u,v\in\mathbb{N}$ satisfy $u\equiv v\operatorname{mod}%
p$, then $z^{u}=z^{v}$.

\item[\textbf{(c)}] Assume that $p\neq2$. Define an element $\tau\in A$ by%
\[
\tau=\sum_{i=0}^{p-1}\left(  \dfrac{i}{p}\right)  z^{i}.
\]
Then,
\[
\tau^{2}=\left(  -1\right)  ^{\left(  p-1\right)  /2}p\cdot1_{A}.
\]

\end{enumerate}
\end{theorem}

The $\tau$ defined in Theorem \ref{thm.qr.gausssum2} \textbf{(c)} is known as
a \textbf{Gaussian sum}. The best-known particular case of Theorem
\ref{thm.qr.gausssum2} is when $A=\mathbb{C}$ and $z=e^{2\pi i/p}=\cos
\dfrac{2\pi}{p}+i\sin\dfrac{2\pi}{p}$; in this case, the claim of Theorem
\ref{thm.qr.gausssum2} \textbf{(c)} is saying that $\tau^{2}=\left(
-1\right)  ^{\left(  p-1\right)  /2}p$. This implies that $\tau=\pm
\sqrt{\left(  -1\right)  ^{\left(  p-1\right)  /2}p}$ in this case, and one
can wonder whether the $\pm$ sign is a $+$ or a $-$. This question has been
answered (the sign is always a $+$ sign), but the proof is tricky and would
take us further afield than I'd like. You can find it in \cite[Theorem
G.2]{Elman22}. To us, this question is not important, since it only concerns
the case $A=\mathbb{C}$, whereas we will apply Theorem \ref{thm.qr.gausssum2}
to a different ring $A$.

\begin{proof}
[Proof of Theorem \ref{thm.qr.gausssum2}.]\textbf{(a)} Multiplying both sides
of the equality (\ref{eq.thm.qr.gausssum2.1}) by $z$, we obtain $z\cdot\left(
z^{p-1}+z^{p-2}+\cdots+z^{0}\right)  =0$. In other words,%
\[
z^{p}+z^{p-1}+\cdots+z^{1}=0.
\]
Subtracting this equality from the original equality
(\ref{eq.thm.qr.gausssum2.1}) (and cancelling all the addends that appear in
both), we obtain $z^{0}-z^{p}=0$. In other words, $z^{0}=z^{p}$. Hence,
$z^{p}=z^{0}=1_{A}$. This proves Theorem \ref{thm.qr.gausssum2} \textbf{(a)}.
\medskip

\textbf{(b)} Let $u,v\in\mathbb{N}$ satisfy $u\equiv v\operatorname{mod}p$. We
must prove that $z^{u}=z^{v}$.

We WLOG assume that $u\geq v$ (otherwise, swap $u$ with $v$). Thus,
$u-v\in\mathbb{N}$. Since $u\equiv v\operatorname{mod}p$, we also have $p\mid
u-v$, so that $\dfrac{u-v}{p}\in\mathbb{Z}$ and therefore $\dfrac{u-v}{p}%
\in\mathbb{N}$ (since $u-v\in\mathbb{N}$). Let us denote this number
$\dfrac{u-v}{p}$ by $m$. Thus, $m=\dfrac{u-v}{p}\in\mathbb{N}$. Also, solving
the equation $m=\dfrac{u-v}{p}$ for $u$, we find $u=mp+v$. Hence,
$z^{u}=z^{mp+v}=\left(  z^{p}\right)  ^{m}z^{v}$. However, Theorem
\ref{thm.qr.gausssum2} \textbf{(a)} yields $z^{p}=1_{A}$. Thus, $\left(
z^{p}\right)  ^{m}=1_{A}^{m}=1_{A}$ and therefore $z^{u}=\underbrace{\left(
z^{p}\right)  ^{m}}_{=1_{A}}z^{v}=z^{v}$. Thus, Theorem \ref{thm.qr.gausssum2}
\textbf{(b)} is proved. \medskip

\textbf{(c)} From (\ref{eq.thm.qr.gausssum2.1}), we obtain%
\begin{align*}
0  &  =z^{p-1}+z^{p-2}+\cdots+z^{0}=\left(  z^{p-1}+z^{p-2}+\cdots
+z^{1}\right)  +\underbrace{z^{0}}_{=1_{A}}\\
&  =\left(  z^{p-1}+z^{p-2}+\cdots+z^{1}\right)  +1_{A}.
\end{align*}
In other words,
\begin{equation}
z^{p-1}+z^{p-2}+\cdots+z^{1}=-1_{A}. \label{pf.thm.qr.gausssum2.c.-1}%
\end{equation}

We have $\tau=\sum_{i=0}^{p-1}\left(  \dfrac{i}{p}\right)  z^{i}$ and
$\tau=\sum_{i=0}^{p-1}\left(  \dfrac{i}{p}\right)  z^{i}=\sum_{j=0}%
^{p-1}\left(  \dfrac{j}{p}\right)  z^{j}$. Multiplying these two equalities,
we obtain%
\begin{align}
\tau\tau &  =\left(  \sum_{i=0}^{p-1}\left(  \dfrac{i}{p}\right)
z^{i}\right)  \left(  \sum_{j=0}^{p-1}\left(  \dfrac{j}{p}\right)
z^{j}\right)  =\sum_{i=0}^{p-1}\ \ \sum_{j=0}^{p-1}\left(  \dfrac{i}%
{p}\right)  z^{i}\left(  \dfrac{j}{p}\right)  z^{j}\nonumber\\
&  =\underbrace{\sum_{i=0}^{p-1}\ \ \sum_{j=0}^{p-1}}_{=\sum_{\left(
i,j\right)  \in\left\{  0,1,\ldots,p-1\right\}  ^{2}}}\left(  \dfrac{i}%
{p}\right)  \left(  \dfrac{j}{p}\right)  z^{i}z^{j}=\sum_{\left(  i,j\right)
\in\left\{  0,1,\ldots,p-1\right\}  ^{2}}\left(  \dfrac{i}{p}\right)  \left(
\dfrac{j}{p}\right)  z^{i}z^{j}\nonumber\\
&  =\sum_{k=0}^{p-1}\ \ \sum_{\substack{\left(  i,j\right)  \in\left\{
0,1,\ldots,p-1\right\}  ^{2};\\i+j\equiv k\operatorname{mod}p}}\left(
\dfrac{i}{p}\right)  \left(  \dfrac{j}{p}\right)  \underbrace{z^{i}z^{j}%
}_{\substack{=z^{i+j}=z^{k}\\\text{(by Theorem \ref{thm.qr.gausssum2}
\textbf{(b)},}\\\text{since }i+j\equiv k\operatorname{mod}p\text{)}%
}}\nonumber\\
&  \ \ \ \ \ \ \ \ \ \ \ \ \ \ \ \ \ \ \ \ \left(
\begin{array}
[c]{c}%
\text{because for each }\left(  i,j\right)  \in\left\{  0,1,\ldots
,p-1\right\}  ^{2}\text{, there exists a}\\
\text{unique }k\in\left\{  0,1,\ldots,p-1\right\}  \text{ that satisfies
}i+j\equiv k\operatorname{mod}p\\
\text{(namely, this }k\text{ is the remainder of }i+j\text{ upon division by
}p\text{)}%
\end{array}
\right) \nonumber\\
&  =\sum_{k=0}^{p-1}\ \ \sum_{\substack{\left(  i,j\right)  \in\left\{
0,1,\ldots,p-1\right\}  ^{2};\\i+j\equiv k\operatorname{mod}p}}\left(
\dfrac{i}{p}\right)  \left(  \dfrac{j}{p}\right)  z^{k}.
\label{pf.thm.qr.gausssum2.c.1}%
\end{align}

Now, let $k\in\left\{  0,1,\ldots,p-1\right\}  $ be a number, and let $\left(
i,j\right)  \in\left\{  0,1,\ldots,p-1\right\}  ^{2}$ be a pair satisfying
$i+j\equiv k\operatorname{mod}p$. Then, Corollary \ref{cor.qr.mult} yields
\begin{equation}
\left(  \dfrac{ij}{p}\right)  =\left(  \dfrac{i}{p}\right)  \left(  \dfrac
{j}{p}\right)  . \label{pf.thm.qr.gausssum2.c.2}%
\end{equation}

However, from $i+j\equiv k\operatorname{mod}p$, we obtain $j\equiv
k-i\operatorname{mod}p$ and thus $ij\equiv i\left(  k-i\right)  =\left(
-1\right)  i\left(  i-k\right)  \operatorname{mod}p$. Hence, $\left(
\dfrac{ij}{p}\right)  =\left(  \dfrac{\left(  -1\right)  i\left(  i-k\right)
}{p}\right)  $ (by (\ref{eq.qr.legendre.congruent}), applied to $a=ij$ and
$b=\left(  -1\right)  i\left(  i-k\right)  $). Comparing this with
(\ref{pf.thm.qr.gausssum2.c.2}), we obtain%
\begin{align}
\left(  \dfrac{i}{p}\right)  \left(  \dfrac{j}{p}\right)   &  =\left(
\dfrac{\left(  -1\right)  i\left(  i-k\right)  }{p}\right)
=\underbrace{\left(  \dfrac{-1}{p}\right)  }_{\substack{=\left(  -1\right)
^{\left(  p-1\right)  /2}\\\text{(by (\ref{pf.thm.qr.-1.2}))}}}\left(
\dfrac{i\left(  i-k\right)  }{p}\right) \nonumber\\
&  \ \ \ \ \ \ \ \ \ \ \ \ \ \ \ \ \ \ \ \ \left(  \text{by Corollary
\ref{cor.qr.mult}, applied to }a=-1\text{ and }b=i\left(  i-k\right)  \right)
\nonumber\\
&  =\left(  -1\right)  ^{\left(  p-1\right)  /2}\left(  \dfrac{i\left(
i-k\right)  }{p}\right)  . \label{pf.thm.qr.gausssum2.c.3}%
\end{align}

Forget that we fixed $k$ and $\left(  i,j\right)  $. We thus have proved
(\ref{pf.thm.qr.gausssum2.c.3}) for each $k\in\left\{  0,1,\ldots,p-1\right\}
$ and each pair $\left(  i,j\right)  \in\left\{  0,1,\ldots,p-1\right\}  ^{2}$
satisfying $i+j\equiv k\operatorname{mod}p$. Therefore,
(\ref{pf.thm.qr.gausssum2.c.1}) becomes%
\begin{align}
\tau\tau &  =\sum_{k=0}^{p-1}\ \ \sum_{\substack{\left(  i,j\right)
\in\left\{  0,1,\ldots,p-1\right\}  ^{2};\\i+j\equiv k\operatorname{mod}%
p}}\underbrace{\left(  \dfrac{i}{p}\right)  \left(  \dfrac{j}{p}\right)
}_{\substack{=\left(  -1\right)  ^{\left(  p-1\right)  /2}\left(
\dfrac{i\left(  i-k\right)  }{p}\right)  \\\text{(by
(\ref{pf.thm.qr.gausssum2.c.3}))}}}z^{k}\nonumber\\
&  =\left(  -1\right)  ^{\left(  p-1\right)  /2}\sum_{k=0}^{p-1}%
\ \ \underbrace{\sum_{\substack{\left(  i,j\right)  \in\left\{  0,1,\ldots
,p-1\right\}  ^{2};\\i+j\equiv k\operatorname{mod}p}}}_{\substack{=\sum
_{i=0}^{p-1}\ \ \sum_{\substack{j\in\left\{  0,1,\ldots,p-1\right\}
;\\i+j\equiv k\operatorname{mod}p}}\\=\sum_{i=0}^{p-1}\ \ \sum_{\substack{j\in
\left\{  0,1,\ldots,p-1\right\}  ;\\j\equiv k-i\operatorname{mod}%
p}}\\\text{(since the congruence }i+j\equiv k\operatorname{mod}p\\\text{is
equivalent to }j\equiv k-i\operatorname{mod}p\text{)}}}\left(  \dfrac{i\left(
i-k\right)  }{p}\right)  z^{k}\nonumber\\
&  =\left(  -1\right)  ^{\left(  p-1\right)  /2}\sum_{k=0}^{p-1}\ \ \sum
_{i=0}^{p-1}\ \ \sum_{\substack{j\in\left\{  0,1,\ldots,p-1\right\}
;\\j\equiv k-i\operatorname{mod}p}}\left(  \dfrac{i\left(  i-k\right)  }%
{p}\right)  z^{k}. \label{pf.thm.qr.gausssum2.c.4}%
\end{align}

Now, fix two numbers $k,i\in\left\{  0,1,\ldots,p-1\right\}  $. Then, there is
a unique number $j\in\left\{  0,1,\ldots,p-1\right\}  $ that satisfies
$j\equiv k-i\operatorname{mod}p$ (namely, the remainder of $k-i$ upon division
by $p$). Hence, the sum $\sum_{\substack{j\in\left\{  0,1,\ldots,p-1\right\}
;\\j\equiv k-i\operatorname{mod}p}}\left(  \dfrac{i\left(  i-k\right)  }%
{p}\right)  $ has exactly $1$ addend, and therefore simplifies to $\left(
\dfrac{i\left(  i-k\right)  }{p}\right)  $.

Forget that we fixed $k,i$. We thus have shown that%
\begin{equation}
\sum_{\substack{j\in\left\{  0,1,\ldots,p-1\right\}  ;\\j\equiv
k-i\operatorname{mod}p}}\left(  \dfrac{i\left(  i-k\right)  }{p}\right)
=\left(  \dfrac{i\left(  i-k\right)  }{p}\right)
\label{pf.thm.qr.gausssum2.c.5}%
\end{equation}
for each $k,i\in\left\{  0,1,\ldots,p-1\right\}  $. Thus,
\begin{align*}
&  \sum_{k=0}^{p-1}\ \ \sum_{i=0}^{p-1}\ \ \underbrace{\sum_{\substack{j\in
\left\{  0,1,\ldots,p-1\right\}  ;\\j\equiv k-i\operatorname{mod}p}}\left(
\dfrac{i\left(  i-k\right)  }{p}\right)  }_{\substack{=\left(  \dfrac{i\left(
i-k\right)  }{p}\right)  \\\text{(by (\ref{pf.thm.qr.gausssum2.c.5}))}}%
}z^{k}\\
&  =\sum_{k=0}^{p-1}\ \ \sum_{i=0}^{p-1}\left(  \dfrac{i\left(  i-k\right)
}{p}\right)  z^{k}\\
&  =\underbrace{\sum_{i=0}^{p-1}\left(  \dfrac{i\left(  i-0\right)  }%
{p}\right)  }_{\substack{=p-1\\\text{(by Proposition \ref{prop.qr.sumii-k}
\textbf{(a)},}\\\text{applied to }k=0\text{)}}}\underbrace{z^{0}}_{=1_{A}%
}+\sum_{k=1}^{p-1}\ \ \underbrace{\sum_{i=0}^{p-1}\left(  \dfrac{i\left(
i-k\right)  }{p}\right)  }_{\substack{=-1\\\text{(by Proposition
\ref{prop.qr.sumii-k} \textbf{(b)},}\\\text{since }k\in\left\{  1,2,\ldots
,p-1\right\}  \text{ entails }p\nmid k\text{)}}}z^{k}\\
&  \ \ \ \ \ \ \ \ \ \ \ \ \ \ \ \ \ \ \ \ \left(  \text{here, we have split
off the addend for }k=0\text{ from the sum}\right) \\
&  =\left(  p-1\right)  \cdot1_{A}+\sum_{k=1}^{p-1}\left(  -1\right)  z^{k}\\
&  =\left(  p-1\right)  \cdot1_{A}-\underbrace{\sum_{k=1}^{p-1}z^{k}%
}_{\substack{=z^{p-1}+z^{p-2}+\cdots+z^{1}=-1_{A}\\\text{(by
(\ref{pf.thm.qr.gausssum2.c.-1}))}}}\\
&  =\left(  p-1\right)  \cdot1_{A}-\left(  -1_{A}\right)  =p\cdot1_{A}.
\end{align*}
Therefore, we can rewrite (\ref{pf.thm.qr.gausssum2.c.4}) as%
\[
\tau\tau=\left(  -1\right)  ^{\left(  p-1\right)  /2}p\cdot1_{A}.
\]
In other words, $\tau^{2}=\left(  -1\right)  ^{\left(  p-1\right)  /2}%
p\cdot1_{A}$. This proves Theorem \ref{thm.qr.gausssum2} \textbf{(c)}.
\end{proof}

\subsubsection{Proof of quadratic reciprocity for two odd primes}

We are now ready to prove the Quadratic Reciprocity Law for two odd primes
(i.e., part \textbf{(b)} of Theorem \ref{thm.qr.qr}):

\begin{proof}
[Proof of Theorem \ref{thm.qr.qr} \textbf{(b)}.]Let $F$ be the field
$\mathbb{Z}/q$. (This is a field, since $q$ is prime.) We have $q\nmid p$
(since $p$ and $q$ are two distinct primes); thus, the residue class
$\overline{p}$ in $\mathbb{Z}/q$ is nonzero. Since $\mathbb{Z}/q$ is a field,
this shows that $\overline{p}$ is a unit.

Now, we shall extend the field $F$ to a ring by adjoining an element $z$ that
satisfies $z^{p-1}+z^{p-2}+\cdots+z^{0}=0$.

Indeed, the polynomial $x^{p-1}+x^{p-2}+\cdots+x^{0}\in F\left[  x\right]  $
is a monic polynomial of degree $p-1$ over $F$, so that its leading
coefficient is a unit. Thus, by Theorem \ref{thm.fieldext.basis-monic}
\textbf{(d)} (applied to $R=F$ and $b=x^{p-1}+x^{p-2}+\cdots+x^{0}$ and
$m=p-1$), there exists a commutative ring that contains $F$ as a subring and
that contains a root of $x^{p-1}+x^{p-2}+\cdots+x^{0}$.

Let $A$ be this ring, and let $z$ be this root that it contains. Thus, $z\in
A$ satisfies $z^{p-1}+z^{p-2}+\cdots+z^{0}=0$. (Note that we could use Theorem
\ref{thm.finfield.splitfield} \textbf{(b)} to obtain a field instead of this
ring $A$, but we will have no need for this, so we have applied the less
advanced Theorem \ref{thm.fieldext.basis-monic} \textbf{(d)}.)

Note that $A$ is a commutative $\mathbb{Z}/q$-algebra (since $A$ is a
commutative ring that contains $F=\mathbb{Z}/q$ as a subring).

Define an element $\tau\in A$ by%
\begin{equation}
\tau=\sum_{i=0}^{p-1}\left(  \dfrac{i}{p}\right)  z^{i}.
\label{pf.thm.qr.qr.b.tau=}%
\end{equation}
Then, Theorem \ref{thm.qr.gausssum2} \textbf{(c)} yields\footnote{Here and in
the rest of this proof, the notation $\overline{r}$ always means the residue
class of an integer $r$ in $\mathbb{Z}/q$.}
\begin{align}
\tau^{2}  &  =\left(  -1\right)  ^{\left(  p-1\right)  /2}p\cdot
\underbrace{1_{A}}_{\substack{=\overline{1}\\\text{(in }\mathbb{Z}/q\text{)}%
}}=\left(  -1\right)  ^{\left(  p-1\right)  /2}p\cdot\overline{1}\nonumber\\
&  =\overline{\left(  -1\right)  ^{\left(  p-1\right)  /2}p}%
\label{pf.thm.qr.qr.b.tau2=}\\
&  =\overline{\left(  -1\right)  ^{\left(  p-1\right)  /2}}\cdot\overline
{p}.\nonumber
\end{align}
This shows that $\tau^{2}$ is a unit of $\mathbb{Z}/q$ (since $\overline
{\left(  -1\right)  ^{\left(  p-1\right)  /2}}$ and $\overline{p}$ are units
of $\mathbb{Z}/q$).

We shall now compute $\tau^{q-1}$.

First, $q$ is odd (since $q\neq2$ is a prime), and thus $\left(  q-1\right)
/2\in\mathbb{N}$. Taking the equality (\ref{pf.thm.qr.qr.b.tau2=}) to the
$\left(  q-1\right)  /2$-th power, we obtain%
\begin{align*}
\left(  \tau^{2}\right)  ^{\left(  q-1\right)  /2}  &  =\left(  \overline
{\left(  -1\right)  ^{\left(  p-1\right)  /2}p}\right)  ^{\left(  q-1\right)
/2}=\overline{\left(  \left(  -1\right)  ^{\left(  p-1\right)  /2}p\right)
^{\left(  q-1\right)  /2}}\\
&  =\overline{\left(  \left(  -1\right)  ^{\left(  p-1\right)  /2}\right)
^{\left(  q-1\right)  /2}p^{\left(  q-1\right)  /2}}=\overline{\left(  \left(
-1\right)  ^{\left(  p-1\right)  /2}\right)  ^{\left(  q-1\right)  /2}}%
\cdot\overline{p^{\left(  q-1\right)  /2}}.
\end{align*}
In view of $\left(  \tau^{2}\right)  ^{\left(  q-1\right)  /2}=\tau^{q-1}$ and
$\left(  \left(  -1\right)  ^{\left(  p-1\right)  /2}\right)  ^{\left(
q-1\right)  /2}=\left(  -1\right)  ^{\left(  p-1\right)  \left(  q-1\right)
/4}$, we can rewrite this as%
\begin{equation}
\tau^{q-1}=\overline{\left(  -1\right)  ^{\left(  p-1\right)  \left(
q-1\right)  /4}}\cdot\overline{p^{\left(  q-1\right)  /2}}.
\label{pf.thm.qr.qr.b.4}%
\end{equation}

However, Theorem \ref{thm.qr.euler} (applied to $q$ and $p$ instead of $p$ and
$a$) yields $\left(  \dfrac{p}{q}\right)  \equiv p^{\left(  q-1\right)
/2}\operatorname{mod}q$. In other words, $\overline{\left(  \dfrac{p}%
{q}\right)  }=\overline{p^{\left(  q-1\right)  /2}}$ in $\mathbb{Z}/q$. Hence,
we can rewrite (\ref{pf.thm.qr.qr.b.4}) as%
\begin{align}
\tau^{q-1}  &  =\overline{\left(  -1\right)  ^{\left(  p-1\right)  \left(
q-1\right)  /4}}\cdot\overline{\left(  \dfrac{p}{q}\right)  }\nonumber\\
&  =\overline{\left(  -1\right)  ^{\left(  p-1\right)  \left(  q-1\right)
/4}\left(  \dfrac{p}{q}\right)  }. \label{pf.thm.qr.qr.b.5}%
\end{align}

On the other hand, define a map%
\begin{align*}
\mathbf{f}:A  &  \rightarrow A,\\
a  &  \mapsto a^{q}.
\end{align*}
Then, Corollary \ref{cor.finfield.frobenius} (applied to $A$ and $q$ instead
of $F$ and $p$) shows that $\mathbf{f}$ is a ring morphism (since $A$ is a
commutative $\mathbb{Z}/q$-algebra). Applying this map $\mathbf{f}$ to both
sides of (\ref{pf.thm.qr.qr.b.tau=}), we find%
\begin{align*}
\mathbf{f}\left(  \tau\right)   &  =\mathbf{f}\left(  \sum_{i=0}^{p-1}\left(
\dfrac{i}{p}\right)  z^{i}\right) \\
&  =\sum_{i=0}^{p-1}\left(  \dfrac{i}{p}\right)  \underbrace{\mathbf{f}\left(
z^{i}\right)  }_{\substack{=\left(  z^{i}\right)  ^{q}\\\text{(by the
definition of }\mathbf{f}\text{)}}}\ \ \ \ \ \ \ \ \ \ \left(  \text{since
}\mathbf{f}\text{ is a ring morphism}\right) \\
&  =\sum_{i=0}^{p-1}\left(  \dfrac{i}{p}\right)  \underbrace{\left(
z^{i}\right)  ^{q}}_{=z^{iq}}=\sum_{i=0}^{p-1}\left(  \dfrac{i}{p}\right)
z^{iq}.
\end{align*}
Since $\mathbf{f}\left(  \tau\right)  =\tau^{q}$ (by the definition of
$\mathbf{f}$), we can rewrite this as
\begin{equation}
\tau^{q}=\sum_{i=0}^{p-1}\left(  \dfrac{i}{p}\right)  z^{iq}.
\label{pf.thm.qr.qr.b.6}%
\end{equation}

We shall now rewrite $\tau$ in a different way. To this purpose, we consider
the map%
\begin{align*}
\left\{  0,1,\ldots,p-1\right\}   &  \rightarrow\left\{  0,1,\ldots
,p-1\right\}  ,\\
i  &  \mapsto\left(  iq\right)  \%p,
\end{align*}
where $\left(  iq\right)  \%p$ denotes the remainder that $iq$ leaves upon
division by $p$. This map is easily seen to be
bijective\footnote{\textit{Proof.} It suffices to show that this map is
injective (because an injective map between two finite sets of the same size
is automatically bijective). So let us show this.
\par
Let $i_{1},i_{2}$ be two distinct elements of $\left\{  0,1,\ldots
,p-1\right\}  $. We must prove that $\left(  i_{1}q\right)  \%p\neq\left(
i_{2}q\right)  \%p$.
\par
Assume the contrary. Thus, $\left(  i_{1}q\right)  \%p=\left(  i_{2}q\right)
\%p$. Hence, $i_{1}q\equiv i_{2}q\operatorname{mod}p$ (because two integers
leave the same remainder upon division by $p$ if and only if they are
congruent modulo $p$). In other words, $p\mid i_{1}q-i_{2}q$. In other words,
$p\mid\left(  i_{1}-i_{2}\right)  q$. Since $p$ is a prime, this entails that
either $p\mid i_{1}-i_{2}$ or $p\mid q$ (or both). Since $p\mid q$ is
impossible (because $p$ and $q$ are two distinct primes), we thus conclude
that $p\mid i_{1}-i_{2}$. In other words, $i_{1}\equiv i_{2}\operatorname{mod}%
p$. However, since $i_{1},i_{2}\in\left\{  0,1,\ldots,p-1\right\}  $, this
entails $i_{1}=i_{2}$, which contradicts the fact that $i_{1},i_{2}$ are
distinct. This contradiction shows that our assumption was false, qed.}. Thus,
we can substitute $\left(  iq\right)  \%p$ for $i$ in the sum $\sum
_{i=0}^{p-1}\left(  \dfrac{i}{p}\right)  z^{i}$. We thus obtain%
\begin{align*}
\sum_{i=0}^{p-1}\left(  \dfrac{i}{p}\right)  z^{i}  &  =\sum_{i=0}%
^{p-1}\underbrace{\left(  \dfrac{\left(  iq\right)  \%p}{p}\right)
}_{\substack{=\left(  \dfrac{iq}{p}\right)  \\\text{(by
(\ref{eq.qr.legendre.congruent}),}\\\text{since }\left(  iq\right)  \%p\equiv
iq\operatorname{mod}p\text{)}}}\underbrace{z^{\left(  iq\right)  \%p}%
}_{\substack{=z^{iq}\\\text{(by Theorem \ref{thm.qr.gausssum2} \textbf{(b)}%
,}\\\text{since }\left(  iq\right)  \%p\equiv iq\operatorname{mod}p\text{)}%
}}\\
&  =\sum_{i=0}^{p-1}\underbrace{\left(  \dfrac{iq}{p}\right)  }%
_{\substack{=\left(  \dfrac{i}{p}\right)  \left(  \dfrac{q}{p}\right)
\\\text{(by Corollary \ref{cor.qr.mult})}}}z^{iq}=\left(  \dfrac{q}{p}\right)
\cdot\underbrace{\sum_{i=0}^{p-1}\left(  \dfrac{i}{p}\right)  z^{iq}%
}_{\substack{=\tau^{q}\\\text{(by (\ref{pf.thm.qr.qr.b.6}))}}}=\left(
\dfrac{q}{p}\right)  \cdot\tau^{q}.
\end{align*}
In view of (\ref{pf.thm.qr.qr.b.tau=}), we can rewrite this as%
\[
\tau=\left(  \dfrac{q}{p}\right)  \cdot\tau^{q}.
\]
Multiplying both sides of this equality by $\tau$, we obtain%
\begin{align*}
\tau^{2}  &  =\left(  \dfrac{q}{p}\right)  \cdot\underbrace{\tau^{q+1}}%
_{=\tau^{q-1}\cdot\tau^{2}}=\left(  \dfrac{q}{p}\right)  \cdot\tau^{q-1}%
\cdot\tau^{2}\\
&  =\left(  \dfrac{q}{p}\right)  \cdot\overline{\left(  -1\right)  ^{\left(
p-1\right)  \left(  q-1\right)  /4}\left(  \dfrac{p}{q}\right)  }\cdot\tau
^{2}\ \ \ \ \ \ \ \ \ \ \left(  \text{by (\ref{pf.thm.qr.qr.b.5})}\right)  .
\end{align*}
Since $\tau^{2}$ is a unit of $\mathbb{Z}/q$, we can cancel $\tau^{2}$ from
both sides of this equality (by multiplying by its inverse). Thus, we obtain%
\begin{align*}
\overline{1}  &  =\left(  \dfrac{q}{p}\right)  \cdot\overline{\left(
-1\right)  ^{\left(  p-1\right)  \left(  q-1\right)  /4}\left(  \dfrac{p}%
{q}\right)  }=\overline{\left(  \dfrac{q}{p}\right)  \cdot\left(  -1\right)
^{\left(  p-1\right)  \left(  q-1\right)  /4}\left(  \dfrac{p}{q}\right)  }\\
&  =\overline{\left(  -1\right)  ^{\left(  p-1\right)  \left(  q-1\right)
/4}\left(  \dfrac{p}{q}\right)  \left(  \dfrac{q}{p}\right)  }.
\end{align*}
In other words,%
\[
1\equiv\left(  -1\right)  ^{\left(  p-1\right)  \left(  q-1\right)  /4}\left(
\dfrac{p}{q}\right)  \left(  \dfrac{q}{p}\right)  \operatorname{mod}q.
\]
Thus, Lemma \ref{lem.qr.01-1} (applied to $q$, $1$ and $\left(  -1\right)
^{\left(  p-1\right)  \left(  q-1\right)  /4}\left(  \dfrac{p}{q}\right)
\left(  \dfrac{q}{p}\right)  $ instead of $p$, $u$ and $v$) shows that%
\begin{equation}
1=\left(  -1\right)  ^{\left(  p-1\right)  \left(  q-1\right)  /4}\left(
\dfrac{p}{q}\right)  \left(  \dfrac{q}{p}\right)  \label{pf.thm.qr.qr.b.at}%
\end{equation}
(because both $1$ and $\left(  -1\right)  ^{\left(  p-1\right)  \left(
q-1\right)  /4}\left(  \dfrac{p}{q}\right)  \left(  \dfrac{q}{p}\right)  $ are
elements of $\left\{  0,1,-1\right\}  $).

However, $p\nmid q$ (since $p$ and $q$ are two distinct primes), and thus
$\left(  \dfrac{q}{p}\right)  $ is either $1$ or $-1$. Hence, in either case,
we have $\left(  \dfrac{q}{p}\right)  ^{2}=1$. Now, multiplying both sides of
the equality (\ref{pf.thm.qr.qr.b.at}) by $\left(  \dfrac{q}{p}\right)  $, we
obtain%
\[
\left(  \dfrac{q}{p}\right)  =\left(  -1\right)  ^{\left(  p-1\right)  \left(
q-1\right)  /4}\left(  \dfrac{p}{q}\right)  \underbrace{\left(  \dfrac{q}%
{p}\right)  ^{2}}_{=1}=\left(  -1\right)  ^{\left(  p-1\right)  \left(
q-1\right)  /4}\left(  \dfrac{p}{q}\right)  .
\]
This proves Theorem \ref{thm.qr.qr} \textbf{(b)}.
\end{proof}

See \cite[Chapter 9]{Burton11} or \cite[Chapter 4]{Stein09} (or almost any
text on elementary number theory) for more about quadratic residues. A
collection of proofs of Theorem \ref{thm.qr.qr} has also been published as a
book (\cite{Baumga15}); one of the most elementary proofs is presented in
\cite[\S 3.12]{KeeGui20}. See also \cite[particularly Chapter 16]{Schroe09}
for an application of quadratic residues to the acoustics of concert halls.
Inside mathematics, there are many other applications:

\begin{exercise}
\label{exe.rings.p=xx+2yy}Let $p\neq2$ be a prime. Prove the following:

\begin{enumerate}
\item[\textbf{(a)}] We have%
\[
\left(  \dfrac{-2}{p}\right)  =%
\begin{cases}
1, & \text{if }p\equiv1\operatorname{mod}8\text{ or }p\equiv
3\operatorname{mod}8;\\
-1, & \text{if }p\equiv-1\operatorname{mod}8\text{ or }p\equiv
-3\operatorname{mod}8.
\end{cases}
\]

\item[\textbf{(b)}] If $p\equiv1\operatorname{mod}8$ or $p\equiv
3\operatorname{mod}8$, then $p$ can be written in the form $x^{2}+2y^{2}$ for
some $x,y\in\mathbb{Z}$. \medskip
\end{enumerate}

[\textbf{Hint:} For part \textbf{(b)}, recall Exercise
\ref{exe.eucldom.Zsqrt-2} and Section \ref{sec.rings.p=xx+yy}.]
\end{exercise}

\begin{exercise}
\label{exe.rings.p=xx+3yy}Let $p>3$ be a prime. Prove the following:

\begin{enumerate}
\item[\textbf{(a)}] We have%
\[
\left(  \dfrac{-3}{p}\right)  =%
\begin{cases}
1, & \text{if }p\equiv1\operatorname{mod}3;\\
-1, & \text{if }p\equiv-1\operatorname{mod}3.
\end{cases}
\]

\item[\textbf{(b)}] If $p\equiv1\operatorname{mod}3$, then $p$ can be written
in the form $u^{2}-uv+v^{2}$ for some $u,v\in\mathbb{Z}$.

\item[\textbf{(c)}] If $p\equiv1\operatorname{mod}3$, then $p$ can be written
in the form $x^{2}+3y^{2}$ for some $x,y\in\mathbb{Z}$. \medskip
\end{enumerate}

[\textbf{Hint:} For parts \textbf{(b)} and \textbf{(c)}, recall Exercise
\ref{exe.eucldom.eisenstein} and Section \ref{sec.rings.p=xx+yy}. This time,
$\mathbb{Z}\left[  \sqrt{-3}\right]  $ is the wrong ring to be working in,
since it is not Euclidean; but $\mathbb{Z}\left[  \omega\right]  $ acts its
part well enough for part \textbf{(b)}. In order to obtain part \textbf{(c)},
use part \textbf{(b)} and then rewrite $u^{2}-uv+v^{2}$ as $\left(  \dfrac
{u}{2}-v\right)  ^{2}+3\left(  \dfrac{u}{2}\right)  ^{2}$ if $u$ is even, or
as $\left(  \dfrac{v}{2}-u\right)  ^{2}+3\left(  \dfrac{v}{2}\right)  ^{2}$ if
$v$ is even, or as $\left(  \dfrac{u+v}{2}\right)  ^{2}+3\left(  \dfrac
{u-v}{2}\right)  ^{2}$ in the remaining case.]
\end{exercise}

\begin{exercise}
\label{exe.rings.p=xx-2yy}Let $p$ be a prime such that $p\equiv
1\operatorname{mod}8$ or $p\equiv7\operatorname{mod}8$. Prove that $p$ can be
written in the form $x^{2}-2y^{2}$ for some $x,y\in\mathbb{Z}$. \medskip

[\textbf{Hint:} Recall Exercise \ref{exe.eucldom.Zsqrt2}. Furthermore, observe
that if $a^{2}-2b^{2}=-p$, then $\left(  a-2b\right)  ^{2}-2\left(
a-b\right)  ^{2}=p$.]
\end{exercise}

\subsubsection{Jacobsthal's explicit formulas for $p=x^{2}+y^{2}$}

\begin{fineprint}
Legendre symbols are useful not only for studying squares in $\mathbb{Z}/p$. A
surprising application was found by Jacobsthal in 1907 \cite[Seite
240]{Jacobs07}: Recall that Theorem \ref{thm.fermat.p=xx+yy} says that each
prime $p$ satisfying $p\equiv1\operatorname{mod}4$ can be written as a sum of
two perfect squares. Jacobsthal used Legendre symbols not only to prove this
theorem in a new way, but also to give \textquotedblleft
explicit\textquotedblright\ formulas for two perfect squares that sum to $p$.
We are using scare quotes around the word \textquotedblleft
explicit\textquotedblright, since using these formulas to compute the squares
is much slower than searching for the squares by brute force (let alone than
an actually efficient algorithm, such as the one given in \cite[\S 5.7]%
{Stein09}), but the formulas are fascinating in their own right.

Jacobsthal's theorem can be stated as follows:

\begin{theorem}
[Jacobsthal's formulas]\label{thm.qr.jacobs}Let $p$ be a prime such that
$p\equiv1\operatorname{mod}4$. For any integer $h$, define%
\begin{equation}
W\left(  h\right)  :=\sum_{i=0}^{p-1}\left(  \dfrac{i\left(  i^{2}+h\right)
}{p}\right)  \in\mathbb{Z}. \label{eq.thm.qr.jacobs.defWh}%
\end{equation}
Let $m$ be a QNR mod $p$ (that is, an integer such that $\overline{m}%
\in\mathbb{Z}/p$ is not a square). Let $a=W\left(  1\right)  $ and $b=W\left(
m\right)  $. Then:

\begin{enumerate}
\item[\textbf{(a)}] The integers $a$ and $b$ are even.

\item[\textbf{(b)}] We have $p=\left(  a/2\right)  ^{2}+\left(  b/2\right)
^{2}$.
\end{enumerate}
\end{theorem}

Take a moment to appreciate this theorem as a miracle, before we somewhat
dispel the mystery through the proof. First, an example:

\begin{itemize}
\item Let $p=13$ (a prime that satisfies $p\equiv1\operatorname{mod}4$), and
let $m=5$ (a QNR mod $p$). Then, using the notation of Theorem
\ref{thm.qr.jacobs}, we have%
\begin{align*}
a  &  =W\left(  1\right)  =\sum_{i=0}^{p-1}\left(  \dfrac{i\left(
i^{2}+1\right)  }{p}\right) \\
&  =\left(  \dfrac{0\left(  0^{2}+1\right)  }{13}\right)  +\left(
\dfrac{1\left(  1^{2}+1\right)  }{13}\right)  +\cdots+\left(  \dfrac{12\left(
12^{2}+1\right)  }{13}\right) \\
&  =0+\left(  -1\right)  +1+1+1+0+1+1+0+1+1+1+\left(  -1\right) \\
&  =6
\end{align*}
and%
\begin{align*}
b  &  =W\left(  m\right)  =W\left(  5\right)  =\sum_{i=0}^{p-1}\left(
\dfrac{i\left(  i^{2}+5\right)  }{p}\right) \\
&  =\left(  \dfrac{0\left(  0^{2}+5\right)  }{13}\right)  +\left(
\dfrac{1\left(  1^{2}+5\right)  }{13}\right)  +\cdots+\left(  \dfrac{12\left(
12^{2}+5\right)  }{13}\right) \\
&  =0+\left(  -1\right)  +\left(  -1\right)  +1+\left(  -1\right)  +\left(
-1\right)  +1+1\\
&  \ \ \ \ \ \ \ \ \ \ +\left(  -1\right)  +\left(  -1\right)  +1+\left(
-1\right)  +\left(  -1\right) \\
&  =-4,
\end{align*}
so that $a$ and $b$ are indeed even and we indeed have $p=\left(  a/2\right)
^{2}+\left(  b/2\right)  ^{2}$ (since $\left(  a/2\right)  ^{2}+\left(
b/2\right)  ^{2}=\left(  6/2\right)  ^{2}+\left(  -4/2\right)  ^{2}=9+4=13=p$).
\end{itemize}

The numbers $W\left(  h\right)  $ in Theorem \ref{thm.qr.jacobs} (and several
similarly defined numbers) are known as \textbf{Jacobsthal sums}.

\bigskip

A sequence of lemmas will pave our way to the proof of Theorem
\ref{thm.qr.jacobs}. First, however, we introduce a simple piece of notation:

\begin{definition}
\label{def.qr.Kpu}Let $p\neq2$ be a prime. Let $u\in\mathbb{Z}/p$. Then,
$K_{p}\left(  u\right)  $ shall denote the Legendre symbol $\left(  \dfrac
{a}{p}\right)  $, where $a$ is an integer satisfying $\overline{a}=u$. This
Legendre symbol is well-defined, i.e., it depends only on $u$ (not on $a$),
because if $a$ and $b$ are two integers satisfying $\overline{a}=\overline{b}%
$, then $a\equiv b\operatorname{mod}p$ and therefore $\left(  \dfrac{a}%
{p}\right)  =\left(  \dfrac{b}{p}\right)  $ (by
(\ref{eq.qr.legendre.congruent})).
\end{definition}

For example, $K_{5}\left(  \overline{7}\right)  =\left(  \dfrac{7}{5}\right)
=-1$ and $K_{5}\left(  \overline{10}\right)  =\left(  \dfrac{10}{5}\right)
=0$.

We fix a prime $p\neq2$ for the rest of this subsection (but we don't require
that $p\equiv1\operatorname{mod}4$). From Definition \ref{def.qr.Kpu}, it
follows that%
\begin{equation}
K_{p}\left(  \overline{a}\right)  =\left(  \dfrac{a}{p}\right)
\ \ \ \ \ \ \ \ \ \ \text{for any }a\in\mathbb{Z}. \label{eq.def.qr.Kpu.Kpa}%
\end{equation}

As a warm-up, we shall prove some facts about the $K_{p}\left(  u\right)  $,
many of which are mere restatements of known properties of Legendre symbols
using our new notation for them:

\begin{lemma}
\label{lem.qr.jacobs.Kpdef}Let $u\in\mathbb{Z}/p$. Then:

\begin{enumerate}
\item[\textbf{(a)}] If $u=0$, then $K_{p}\left(  u\right)  =0$.

\item[\textbf{(b)}] If $u\neq0$ but $u$ is a square, then $K_{p}\left(
u\right)  =1$.

\item[\textbf{(c)}] If $u$ is not a square, then $K_{p}\left(  u\right)  =-1$.

\item[\textbf{(d)}] If $u\neq0$, then $K_{p}\left(  u\right)  \equiv
1\operatorname{mod}2$.

\item[\textbf{(e)}] We always have $\left(  K_{p}\left(  u\right)  \right)
^{3}=K_{p}\left(  u\right)  $.

\item[\textbf{(f)}] If $u\neq0$, then $\left(  K_{p}\left(  u\right)  \right)
^{2}=1$.
\end{enumerate}
\end{lemma}

\begin{proof}
Write the residue class $u$ in the form $u=\overline{a}$ for some
$a\in\mathbb{Z}$. Then, $K_{p}\left(  u\right)  =\left(  \dfrac{a}{p}\right)
$ (by Definition \ref{def.qr.Kpu}).

\textbf{(a)} We must prove that $K_{p}\left(  0\right)  =0$. However, the zero
of $\mathbb{Z}/p$ is the residue class $\overline{0}$. That is, $0=\overline
{0}$. Hence, $K_{p}\left(  0\right)  =K_{p}\left(  \overline{0}\right)
=\left(  \dfrac{0}{p}\right)  $ (by (\ref{eq.def.qr.Kpu.Kpa})), so that
$K_{p}\left(  0\right)  =\left(  \dfrac{0}{p}\right)  =0$ (by Definition
\ref{def.qr.legendre}, since $p\mid0$). This proves Lemma
\ref{lem.qr.jacobs.Kpdef} \textbf{(a)}.

\textbf{(b)} Assume that $u\neq0$ and that $u$ is a square. Then,
$\overline{a}=u\neq0$ in $\mathbb{Z}/p$, so that $a$ is not divisible by $p$.
Moreover, $\overline{a}=u$ is a square in $\mathbb{Z}/p$. Hence, $a$ is a QR
mod $p$ (by the definition of a QR). Therefore, $\left(  \dfrac{a}{p}\right)
=1$ (by the definition of the Legendre symbol). Now, $K_{p}\left(  u\right)
=\left(  \dfrac{a}{p}\right)  =1$. This proves Lemma \ref{lem.qr.jacobs.Kpdef}
\textbf{(b)}.

\textbf{(c)} Assume that $u$ is not a square. Then, $\overline{a}=u\neq0$ in
$\mathbb{Z}/p$ (since $0$ is a square but $u$ is not), so that $a$ is not
divisible by $p$. Moreover, $\overline{a}=u$ is not a square in $\mathbb{Z}%
/p$. Hence, $a$ is a QNR mod $p$ (by the definition of a QNR). Therefore,
$\left(  \dfrac{a}{p}\right)  =-1$ (by the definition of the Legendre symbol).
Now, $K_{p}\left(  u\right)  =\left(  \dfrac{a}{p}\right)  =-1$. This proves
Lemma \ref{lem.qr.jacobs.Kpdef} \textbf{(c)}.

\textbf{(d)} Assume that $u\neq0$. Then, parts \textbf{(b)} and \textbf{(c)}
of Lemma \ref{lem.qr.jacobs.Kpdef} show that $K_{p}\left(  u\right)  $ is
either $1$ or $-1$. In either case, $K_{p}\left(  u\right)  $ is odd, i.e., we
have $K_{p}\left(  u\right)  \equiv1\operatorname{mod}2$. This proves Lemma
\ref{lem.qr.jacobs.Kpdef} \textbf{(d)}.

\textbf{(e)} The number $K_{p}\left(  u\right)  $ is a Legendre symbol (by its
definition) and thus belongs to the set $\left\{  0,1,-1\right\}  $ (since any
Legendre symbol belongs to this set). In other words, $K_{p}\left(  u\right)
\in\left\{  0,1,-1\right\}  $.

However, each number $x\in\left\{  0,1,-1\right\}  $ satisfies $x^{3}=x$ (just
check this for each of the values $0$, $1$ and $-1$). Applying this to
$x=K_{p}\left(  u\right)  $, we obtain $\left(  K_{p}\left(  u\right)
\right)  ^{3}=K_{p}\left(  u\right)  $ (since $K_{p}\left(  u\right)
\in\left\{  0,1,-1\right\}  $). This proves Lemma \ref{lem.qr.jacobs.Kpdef}
\textbf{(e)}.

\textbf{(f)} Assume that $u\neq0$. Then, parts \textbf{(b)} and \textbf{(c)}
of Lemma \ref{lem.qr.jacobs.Kpdef} show that $K_{p}\left(  u\right)  $ is
either $1$ or $-1$. In either case, $\left(  K_{p}\left(  u\right)  \right)
^{2}=1$. This proves Lemma \ref{lem.qr.jacobs.Kpdef} \textbf{(f)}.
\end{proof}

\begin{lemma}
\label{lem.qr.jacobs.Kp=num}Let $u\in\mathbb{Z}/p$. Then,\footnotemark%
\[
K_{p}\left(  u\right)  =\left(  \text{\# of elements }x\in\mathbb{Z}/p\text{
such that }x^{2}=u\right)  -1.
\]

\end{lemma}

\footnotetext{The symbol \textquotedblleft\#\textquotedblright\ means
\textquotedblleft number\textquotedblright. For instance, $\left(  \text{\# of
odd integers }i\in\left\{  0,1,\ldots,10\right\}  \right)  =5$.}

\begin{proof}
Write the residue class $u$ as $u=\overline{a}$ for some integer $a$. We are
in one of the following three cases:

\textit{Case 1:} We have $u=0$.

\textit{Case 2:} We have $u\neq0$, but $u$ is a square in $\mathbb{Z}/p$.

\textit{Case 3:} The element $u$ is not a square in $\mathbb{Z}/p$.

Let us first consider Case 1. In this case, $u=0$. Thus, Lemma
\ref{lem.qr.jacobs.Kpdef} \textbf{(a)} yields $K_{p}\left(  u\right)  =0$. On
the other hand, $\mathbb{Z}/p$ is an integral domain. Thus, there is only one
element $x\in\mathbb{Z}/p$ such that $x^{2}=0$ (namely, $0$). Thus,%
\[
\left(  \text{\# of elements }x\in\mathbb{Z}/p\text{ such that }%
x^{2}=0\right)  =1.
\]
In other words,%
\[
\left(  \text{\# of elements }x\in\mathbb{Z}/p\text{ such that }%
x^{2}=u\right)  =1
\]
(since $u=0$). Hence,%
\[
\left(  \text{\# of elements }x\in\mathbb{Z}/p\text{ such that }%
x^{2}=u\right)  -1=0=K_{p}\left(  u\right)
\]
(since $K_{p}\left(  u\right)  =0$). Thus, Lemma \ref{lem.qr.jacobs.Kp=num} is
proved in Case 1.

Let us now consider Case 2. In this case, we have $u\neq0$, but $u$ is a
square in $\mathbb{Z}/p$. Hence, Lemma \ref{lem.qr.jacobs.Kpdef} \textbf{(b)}
yields $K_{p}\left(  u\right)  =1$.

We have $u=y^{2}$ for some $y\in\mathbb{Z}/p$ (since $u$ is a square).
Consider this $y$. Then, $y\neq0$ (since $y^{2}=u\neq0$). Also, $\overline
{2}\neq0$ in $\mathbb{Z}/p$ (since $p\neq2$ is a prime). Now,
$2y=\underbrace{\left(  2\cdot1_{\mathbb{Z}/p}\right)  }_{=\overline{2}%
}y=\overline{2}y\neq0$ (since $\overline{2}\neq0$ and $y\neq0$, and since
$\mathbb{Z}/p$ is an integral domain). Subtracting $y$ from this non-equation,
we obtain $y\neq-y$.

Now,%
\begin{align}
\left\{  x\in\mathbb{Z}/p\ \mid\ x^{2}=y^{2}\right\}   &  =\left\{
x\in\mathbb{Z}/p\ \mid\ x^{2}-y^{2}=0\right\} \nonumber\\
&  =\left\{  x\in\mathbb{Z}/p\ \mid\ \left(  x-y\right)  \left(  x+y\right)
=0\right\} \nonumber\\
&  \ \ \ \ \ \ \ \ \ \ \ \ \ \ \ \ \ \ \ \ \left(  \text{since }x^{2}%
-y^{2}=\left(  x-y\right)  \left(  x+y\right)  \right) \nonumber\\
&  =\left\{  x\in\mathbb{Z}/p\ \mid\ x-y=0\text{ or }x+y=0\right\} \nonumber\\
&  \ \ \ \ \ \ \ \ \ \ \ \ \ \ \ \ \ \ \ \ \left(  \text{since }%
\mathbb{Z}/p\text{ is an integral domain}\right) \nonumber\\
&  =\left\{  x\in\mathbb{Z}/p\ \mid\ x=y\text{ or }x=-y\right\} \nonumber\\
&  =\left\{  y,-y\right\}  . \label{pf.lem.qr.jacobs.Kp=num.3}%
\end{align}
However, $u=y^{2}$, and thus%
\begin{align*}
&  \left(  \text{\# of elements }x\in\mathbb{Z}/p\text{ such that }%
x^{2}=u\right) \\
&  =\left(  \text{\# of elements }x\in\mathbb{Z}/p\text{ such that }%
x^{2}=y^{2}\right) \\
&  =\left\vert \left\{  x\in\mathbb{Z}/p\ \mid\ x^{2}=y^{2}\right\}
\right\vert =\left\vert \left\{  y,-y\right\}  \right\vert
\ \ \ \ \ \ \ \ \ \ \left(  \text{by (\ref{pf.lem.qr.jacobs.Kp=num.3})}\right)
\\
&  =2\ \ \ \ \ \ \ \ \ \ \left(  \text{since }y\neq-y\right)  ,
\end{align*}
so that%
\[
\left(  \text{\# of elements }x\in\mathbb{Z}/p\text{ such that }%
x^{2}=u\right)  -1=2-1=1=K_{p}\left(  u\right)
\]
(since $K_{p}\left(  u\right)  =1$). Thus, Lemma \ref{lem.qr.jacobs.Kp=num} is
proved in Case 2.

Finally, let us consider Case 3. In this case, the element $u$ is not a square
in $\mathbb{Z}/p$. Thus, Lemma \ref{lem.qr.jacobs.Kpdef} \textbf{(c)} yields
$K_{p}\left(  u\right)  =-1$. Also, there exist no elements $x\in\mathbb{Z}/p$
such that $x^{2}=u$ (since $u$ is not a square). Hence,%
\[
\left(  \text{\# of elements }x\in\mathbb{Z}/p\text{ such that }%
x^{2}=u\right)  =0,
\]
so that%
\[
\left(  \text{\# of elements }x\in\mathbb{Z}/p\text{ such that }%
x^{2}=u\right)  -1=0-1=-1=K_{p}\left(  u\right)
\]
(since $K_{p}\left(  u\right)  =-1$). Thus, Lemma \ref{lem.qr.jacobs.Kp=num}
is proved in Case 3.

We have now proved Lemma \ref{lem.qr.jacobs.Kp=num} in all three cases, so
that Lemma \ref{lem.qr.jacobs.Kp=num} is really proved.
\end{proof}

\begin{lemma}
\label{lem.qr.jacobs.sumKp}We have%
\[
\sum_{u\in\mathbb{Z}/p}K_{p}\left(  u\right)  =0.
\]

\end{lemma}

\begin{proof}
By Lemma \ref{lem.qr.jacobs.Kp=num}, we have%
\begin{align*}
\sum_{u\in\mathbb{Z}/p}K_{p}\left(  u\right)   &  =\sum_{u\in\mathbb{Z}%
/p}\left(  \left(  \text{\# of elements }x\in\mathbb{Z}/p\text{ such that
}x^{2}=u\right)  -1\right) \\
&  =\underbrace{\sum_{u\in\mathbb{Z}/p}\left(  \text{\# of elements }%
x\in\mathbb{Z}/p\text{ such that }x^{2}=u\right)  }_{\substack{=\left(
\text{\# of all elements }x\in\mathbb{Z}/p\right)  \\\text{(since each }%
x\in\mathbb{Z}/p\text{ satisfies }x^{2}=u\text{ for exactly one }%
u\in\mathbb{Z}/p\text{)}}}-\underbrace{\sum_{u\in\mathbb{Z}/p}1}_{=\left\vert
\mathbb{Z}/p\right\vert }\\
&  =\underbrace{\left(  \text{\# of all elements }x\in\mathbb{Z}/p\right)
}_{=\left\vert \mathbb{Z}/p\right\vert }-\left\vert \mathbb{Z}/p\right\vert
=\left\vert \mathbb{Z}/p\right\vert -\left\vert \mathbb{Z}/p\right\vert =0.
\end{align*}
This proves Lemma \ref{lem.qr.jacobs.sumKp}.
\end{proof}

\begin{lemma}
\label{lem.qr.jacobs.Kpprod}Let $u,v\in\mathbb{Z}/p$. Then, $K_{p}\left(
uv\right)  =K_{p}\left(  u\right)  \cdot K_{p}\left(  v\right)  $.
\end{lemma}

\begin{proof}
This is just a restatement of Corollary \ref{cor.qr.mult}. In more detail:

Write the residue classes $u$ and $v$ in the forms $u=\overline{a}$ and
$v=\overline{b}$ for some $a,b\in\mathbb{Z}$. Then, $uv=\overline{a}%
\cdot\overline{b}=\overline{ab}$. Hence,
\begin{align*}
K_{p}\left(  uv\right)   &  =K_{p}\left(  \overline{ab}\right)  =\left(
\dfrac{ab}{p}\right)  \ \ \ \ \ \ \ \ \ \ \left(  \text{by
(\ref{eq.def.qr.Kpu.Kpa}), applied to }ab\text{ instead of }a\right) \\
&  =\left(  \dfrac{a}{p}\right)  \left(  \dfrac{b}{p}\right)
\ \ \ \ \ \ \ \ \ \ \left(  \text{by Corollary \ref{cor.qr.mult}}\right)  .
\end{align*}
Comparing this with%
\[
K_{p}\left(  \underbrace{u}_{=\overline{a}}\right)  K_{p}\left(
\underbrace{v}_{=\overline{b}}\right)  =\underbrace{K_{p}\left(  \overline
{a}\right)  }_{\substack{=\left(  \dfrac{a}{p}\right)  \\\text{(by
(\ref{eq.def.qr.Kpu.Kpa}))}}}\ \ \underbrace{K_{p}\left(  \overline{b}\right)
}_{\substack{=\left(  \dfrac{b}{p}\right)  \\\text{(by
(\ref{eq.def.qr.Kpu.Kpa}))}}}=\left(  \dfrac{a}{p}\right)  \left(  \dfrac
{b}{p}\right)  ,
\]
we find $K_{p}\left(  uv\right)  =K_{p}\left(  u\right)  \cdot K_{p}\left(
v\right)  $. This proves Lemma \ref{lem.qr.jacobs.Kpprod}.
\end{proof}

\begin{lemma}
\label{lem.qr.jacobs.Kpxy}Let $y\in\mathbb{Z}/p$ be nonzero. Then,%
\[
\sum_{\substack{x\in\mathbb{Z}/p;\\x^{2}=y^{2}}}K_{p}\left(  xy\right)
=1+\left(  -1\right)  ^{\left(  p-1\right)  /2}.
\]

\end{lemma}

\begin{proof}
Let $u=y^{2}$. Then, $u=y^{2}=yy$ is nonzero (since $y$ is nonzero, and since
$\mathbb{Z}/p$ is an integral domain). Thus, as in the proof of Lemma
\ref{lem.qr.jacobs.Kp=num} (specifically, in Case 2), we can see that
$\left\{  x\in\mathbb{Z}/p\ \mid\ x^{2}=y^{2}\right\}  =\left\{  y,-y\right\}
$ and $y\neq-y$. This shows that the sum $\sum_{\substack{x\in\mathbb{Z}%
/p;\\x^{2}=y^{2}}}K_{p}\left(  xy\right)  $ has only two addends, namely the
addends for $x=y$ and for $x=-y$. Hence,%
\begin{align}
\sum_{\substack{x\in\mathbb{Z}/p;\\x^{2}=y^{2}}}K_{p}\left(  xy\right)   &
=K_{p}\left(  \underbrace{yy}_{=y^{2}}\right)  +K_{p}\left(
\underbrace{\left(  -y\right)  y}_{=\left(  -1\right)  y^{2}}\right)
=K_{p}\left(  y^{2}\right)  +\underbrace{K_{p}\left(  \left(  -1\right)
y^{2}\right)  }_{\substack{=K_{p}\left(  -1\right)  K_{p}\left(  y^{2}\right)
\\\text{(by Lemma \ref{lem.qr.jacobs.Kpprod})}}}\nonumber\\
&  =K_{p}\left(  y^{2}\right)  +K_{p}\left(  -1\right)  K_{p}\left(
y^{2}\right)  =\left(  1+K_{p}\left(  -1\right)  \right)  K_{p}\left(
\underbrace{y^{2}}_{=u}\right) \nonumber\\
&  =\left(  1+K_{p}\left(  -1\right)  \right)  K_{p}\left(  u\right)  .
\label{pf.lem.qr.jacobs.Kpxy.1}%
\end{align}

However, $u$ is a square (since $u=y^{2}$) and satisfies $u\neq0$ (since $u$
is nonzero). Thus, Lemma \ref{lem.qr.jacobs.Kpdef} \textbf{(b)} yields
$K_{p}\left(  u\right)  =1$. Moreover, in $\mathbb{Z}/p$, we have
$-1=\overline{-1}$, so that
\begin{align*}
K_{p}\left(  -1\right)   &  =K_{p}\left(  \overline{-1}\right)  =\left(
\dfrac{-1}{p}\right)  \ \ \ \ \ \ \ \ \ \ \left(  \text{by
(\ref{eq.def.qr.Kpu.Kpa})}\right) \\
&  =%
\begin{cases}
1, & \text{if }p\equiv1\operatorname{mod}4;\\
-1, & \text{if }p\equiv3\operatorname{mod}4
\end{cases}
\ \ \ \ \ \ \ \ \ \ \left(  \text{by Theorem \ref{thm.qr.-1}}\right) \\
&  =\left(  -1\right)  ^{\left(  p-1\right)  /2}.
\end{align*}

Thus, (\ref{pf.lem.qr.jacobs.Kpxy.1}) becomes%
\[
\sum_{\substack{x\in\mathbb{Z}/p;\\x^{2}=y^{2}}}K_{p}\left(  xy\right)
=\left(  1+K_{p}\left(  -1\right)  \right)  \underbrace{K_{p}\left(  u\right)
}_{=1}=1+\underbrace{K_{p}\left(  -1\right)  }_{=\left(  -1\right)  ^{\left(
p-1\right)  /2}}=1+\left(  -1\right)  ^{\left(  p-1\right)  /2}.
\]
This proves Lemma \ref{lem.qr.jacobs.Kpxy}.
\end{proof}

For each $d\in\mathbb{Z}/p$, we define an integer $W\left(  d\right)  $ by%
\[
W\left(  d\right)  =\sum_{u\in\mathbb{Z}/p}K_{p}\left(  u\left(
u^{2}+d\right)  \right)  .
\]
Thus, for each $h\in\mathbb{Z}$, we have%
\begin{align}
W\left(  \overline{h}\right)   &  =\sum_{u\in\mathbb{Z}/p}K_{p}\left(
u\left(  u^{2}+\overline{h}\right)  \right)  =\sum_{i=0}^{p-1}K_{p}\left(
\underbrace{\overline{i}\left(  \overline{i}^{2}+\overline{h}\right)
}_{=\overline{i\left(  i^{2}+h\right)  }}\right) \nonumber\\
&  \ \ \ \ \ \ \ \ \ \ \ \ \ \ \ \ \ \ \ \ \left(
\begin{array}
[c]{c}%
\text{here, we have substituted }\overline{i}\text{ for }u\text{ in the
sum,}\\
\text{since the map }\left\{  0,1,\ldots,p-1\right\}  \rightarrow
\mathbb{Z}/p\\
\text{that sends each }i\text{ to }\overline{i}\text{ is a bijection}%
\end{array}
\right) \nonumber\\
&  =\sum_{i=0}^{p-1}\underbrace{K_{p}\left(  \overline{i\left(  i^{2}%
+h\right)  }\right)  }_{\substack{=\left(  \dfrac{i\left(  i^{2}+h\right)
}{p}\right)  \\\text{(by (\ref{eq.def.qr.Kpu.Kpa}))}}}=\sum_{i=0}^{p-1}\left(
\dfrac{i\left(  i^{2}+h\right)  }{p}\right)  . \label{eq.def.qr.Kpu.Wh}%
\end{align}
The right hand side here is precisely what was called $W\left(  h\right)  $ in
Theorem \ref{thm.qr.jacobs}.

Now, we shall prove the following lemmas about these values $W\left(
d\right)  $:

\begin{lemma}
\label{lem.qr.jacobs.0}The residue class $\overline{0}\in\mathbb{Z}/p$
satisfies $W\left(  \overline{0}\right)  =0$.
\end{lemma}

\begin{proof}
By definition of $W\left(  \overline{0}\right)  $, we have%
\begin{align*}
W\left(  \overline{0}\right)   &  =\sum_{u\in\mathbb{Z}/p}K_{p}\left(
u\underbrace{\left(  u^{2}+\overline{0}\right)  }_{=u^{2}}\right)  =\sum
_{u\in\mathbb{Z}/p}\underbrace{K_{p}\left(  uu^{2}\right)  }_{\substack{=K_{p}%
\left(  u\right)  \cdot K_{p}\left(  u^{2}\right)  \\\text{(by Lemma
\ref{lem.qr.jacobs.Kpprod})}}}=\sum_{u\in\mathbb{Z}/p}K_{p}\left(  u\right)
\cdot\underbrace{K_{p}\left(  u^{2}\right)  }_{\substack{=K_{p}\left(
uu\right)  \\=K_{p}\left(  u\right)  \cdot K_{p}\left(  u\right)  \\\text{(by
Lemma \ref{lem.qr.jacobs.Kpprod})}}}\\
&  =\sum_{u\in\mathbb{Z}/p}\underbrace{K_{p}\left(  u\right)  \cdot
K_{p}\left(  u\right)  \cdot K_{p}\left(  u\right)  }_{\substack{=\left(
K_{p}\left(  u\right)  \right)  ^{3}\\=K_{p}\left(  u\right)  \\\text{(by
Lemma \ref{lem.qr.jacobs.Kpdef} \textbf{(e)})}}}=\sum_{u\in\mathbb{Z}/p}%
K_{p}\left(  u\right)  =0
\end{align*}
(by Lemma \ref{lem.qr.jacobs.sumKp}). This proves Lemma \ref{lem.qr.jacobs.0}.
\end{proof}

\begin{lemma}
\label{lem.qr.jacobs.even}For any $d\in\mathbb{Z}/p$, the integer $W\left(
d\right)  $ is even.
\end{lemma}

\begin{proof}
We know that $p\neq2$ is a prime; hence, $p$ is odd. In other words,
$p\equiv1\operatorname{mod}2$.

Let $d\in\mathbb{Z}/p$. The definition of $W\left(  d\right)  $ yields%
\begin{align}
W\left(  d\right)   &  =\sum_{u\in\mathbb{Z}/p}K_{p}\left(  u\left(
u^{2}+d\right)  \right) \nonumber\\
&  =\sum_{\substack{u\in\mathbb{Z}/p;\\u\left(  u^{2}+d\right)  =0}%
}\underbrace{K_{p}\left(  u\left(  u^{2}+d\right)  \right)  }%
_{\substack{=0\\\text{(by Lemma \ref{lem.qr.jacobs.Kpdef} \textbf{(a)}%
,}\\\text{applied to }u\left(  u^{2}+d\right)  \text{ instead of }u\text{)}%
}}+\sum_{\substack{u\in\mathbb{Z}/p;\\u\left(  u^{2}+d\right)  \neq
0}}\underbrace{K_{p}\left(  u\left(  u^{2}+d\right)  \right)  }%
_{\substack{\equiv1\operatorname{mod}2\\\text{(by Lemma
\ref{lem.qr.jacobs.Kpdef} \textbf{(d)},}\\\text{applied to }u\left(
u^{2}+d\right)  \text{ instead of }u\text{)}}}\nonumber\\
&  \equiv\underbrace{\sum_{\substack{u\in\mathbb{Z}/p;\\u\left(
u^{2}+d\right)  =0}}0}_{=0}+\sum_{\substack{u\in\mathbb{Z}/p;\\u\left(
u^{2}+d\right)  \neq0}}1=\sum_{\substack{u\in\mathbb{Z}/p;\\u\left(
u^{2}+d\right)  \neq0}}1\nonumber\\
&  =\left\vert \left\{  u\in\mathbb{Z}/p\ \mid\ u\left(  u^{2}+d\right)
\neq0\right\}  \right\vert \cdot1\nonumber\\
&  =\left\vert \left\{  u\in\mathbb{Z}/p\ \mid\ u\left(  u^{2}+d\right)
\neq0\right\}  \right\vert \nonumber\\
&  =\underbrace{\left\vert \mathbb{Z}/p\right\vert }_{\substack{=p\\\equiv
1\operatorname{mod}2}}-\left\vert \left\{  u\in\mathbb{Z}/p\ \mid\ u\left(
u^{2}+d\right)  =0\right\}  \right\vert \nonumber\\
&  \ \ \ \ \ \ \ \ \ \ \ \ \ \ \ \ \ \ \ \ \left(
\begin{array}
[c]{c}%
\text{since the set }\left\{  u\in\mathbb{Z}/p\ \mid\ u\left(  u^{2}+d\right)
\neq0\right\}  \text{ is the}\\
\text{complement of }\left\{  u\in\mathbb{Z}/p\ \mid\ u\left(  u^{2}+d\right)
=0\right\}  \text{ within }\mathbb{Z}/p
\end{array}
\right) \nonumber\\
&  \equiv1-\left\vert \left\{  u\in\mathbb{Z}/p\ \mid\ u\left(  u^{2}%
+d\right)  =0\right\}  \right\vert \operatorname{mod}2.
\label{pf.lem.qr.jacobs.even.1}%
\end{align}

We shall now prove that the number $\left\vert \left\{  u\in\mathbb{Z}%
/p\ \mid\ u\left(  u^{2}+d\right)  =0\right\}  \right\vert $ is odd.

Indeed, assume the contrary. Thus, $\left\vert \left\{  u\in\mathbb{Z}%
/p\ \mid\ u\left(  u^{2}+d\right)  =0\right\}  \right\vert $ is even. In other
words, the number of all $u\in\mathbb{Z}/p$ satisfying $u\left(
u^{2}+d\right)  =0$ is even. In other words, the number of roots of the
polynomial $x\left(  x^{2}+d\right)  \in\left(  \mathbb{Z}/p\right)  \left[
x\right]  $ in $\mathbb{Z}/p$ is even.

The polynomial $x\left(  x^{2}+d\right)  \in\left(  \mathbb{Z}/p\right)
\left[  x\right]  $ has degree $3$, and thus has $\leq3$ roots in
$\mathbb{Z}/p$ (by Theorem \ref{thm.polring.univar-easyFTA}, since
$\mathbb{Z}/p$ is an integral domain). In other words, its number of roots is
$\leq3$. Since we also know that this number is even, we thus conclude that
this number is $0$ or $2$ (since the only even nonnegative integers $\leq3$
are $0$ and $2$). In other words, the polynomial $x\left(  x^{2}+d\right)
\in\left(  \mathbb{Z}/p\right)  \left[  x\right]  $ has either $0$ or $2$
roots in $\mathbb{Z}/p$. Since it cannot have $0$ roots in $\mathbb{Z}/p$
(because $0$ is clearly a root of this polynomial), we thus conclude that it
has $2$ roots in $\mathbb{Z}/p$. One of these $2$ roots is $0$ (since $0$ is
clearly a root of this polynomial); let $r$ be the other root. Thus, $r\neq0$
and $r\left(  r^{2}+d\right)  =0$. Hence, $\left(  -r\right)  \left(  \left(
-r\right)  ^{2}+d\right)  =-\underbrace{r\left(  r^{2}+d\right)  }_{=0}=0$.
This shows that $-r$ is a root of the polynomial $x\left(  x^{2}+d\right)
\in\left(  \mathbb{Z}/p\right)  \left[  x\right]  $ in $\mathbb{Z}/p$. Since
the only roots of this polynomial are $0$ and $r$, we thus conclude that $-r$
must be either $0$ or $r$. Since $-r$ cannot be $0$ (because $r\neq0$), we
thus conclude that $-r=r$. Adding $r$ to both sides of this equality, we find
$0=2r=\underbrace{\left(  2\cdot1_{\mathbb{Z}/p}\right)  }_{=\overline{2}%
}r=\overline{2}r$. However, $\overline{2}\neq0$ in $\mathbb{Z}/p$ (since
$p\neq2$ is a prime) and $r\neq0$. Since $\mathbb{Z}/p$ is an integral domain,
these entail that $\overline{2}r\neq0$. This contradicts $0=\overline{2}r$.

This contradiction shows that our assumption was wrong. Hence, we have shown
that $\left\vert \left\{  u\in\mathbb{Z}/p\ \mid\ u\left(  u^{2}+d\right)
=0\right\}  \right\vert $ is odd. In other words,%
\[
\left\vert \left\{  u\in\mathbb{Z}/p\ \mid\ u\left(  u^{2}+d\right)
=0\right\}  \right\vert \equiv1\operatorname{mod}2.
\]
Thus, (\ref{pf.lem.qr.jacobs.even.1}) becomes%
\[
W\left(  d\right)  \equiv1-\underbrace{\left\vert \left\{  u\in\mathbb{Z}%
/p\ \mid\ u\left(  u^{2}+d\right)  =0\right\}  \right\vert }_{\equiv
1\operatorname{mod}2}\equiv1-1=0\operatorname{mod}2.
\]
In other words, $W\left(  d\right)  $ is even. This proves Lemma
\ref{lem.qr.jacobs.even}.
\end{proof}

\begin{lemma}
\label{lem.qr.jacobs.uuh}Let $d,c\in\mathbb{Z}/p$. Then, $W\left(
c^{2}d\right)  =K_{p}\left(  c\right)  W\left(  d\right)  $.
\end{lemma}

\begin{proof}
We have $K_{p}\left(  \overline{0}\right)  =0$ (by Lemma
\ref{lem.qr.jacobs.Kpdef} \textbf{(a)}, applied to $u=\overline{0}$).
Furthermore, $\overline{0}^{2}d=\overline{0}$, and therefore%
\begin{align*}
W\left(  \overline{0}^{2}d\right)   &  =W\left(  \overline{0}\right)
=0\ \ \ \ \ \ \ \ \ \ \left(  \text{by Lemma \ref{lem.qr.jacobs.0}}\right) \\
&  =K_{p}\left(  \overline{0}\right)  W\left(  d\right)
\ \ \ \ \ \ \ \ \ \ \left(  \text{since }\underbrace{K_{p}\left(  \overline
{0}\right)  }_{=0}W\left(  d\right)  =0\right)  .
\end{align*}
Thus, Lemma \ref{lem.qr.jacobs.uuh} is proved in the case when $c=\overline
{0}$. For the rest of this proof, we thus WLOG assume that $c\neq\overline{0}%
$. Hence, $c$ is a nonzero element of $\mathbb{Z}/p$, and thus is a unit
(since $\mathbb{Z}/p$ is a field). Hence, the map%
\begin{align}
\mathbb{Z}/p  &  \rightarrow\mathbb{Z}/p,\nonumber\\
u  &  \mapsto cu \label{pf.lem.qr.jacobs.uuh.bij}%
\end{align}
is a bijection. Furthermore, $c\neq\overline{0}=0$ entails $c^{2}\neq0$ (since
$\mathbb{Z}/p$ is an integral domain), and clearly $c^{2}$ is a square. Hence,
Lemma \ref{lem.qr.jacobs.Kpdef} \textbf{(b)} (applied to $u=c^{2}$) shows that
$K_{p}\left(  c^{2}\right)  =1$.

Now, the definition of $W\left(  c^{2}d\right)  $ yields%
\begin{align*}
W\left(  c^{2}d\right)   &  =\sum_{u\in\mathbb{Z}/p}K_{p}\left(  u\left(
u^{2}+c^{2}d\right)  \right)  =\sum_{u\in\mathbb{Z}/p}K_{p}\left(
\underbrace{\left(  cu\right)  \left(  \left(  cu\right)  ^{2}+c^{2}d\right)
}_{=c^{3}u\left(  u^{2}+d\right)  }\right) \\
&  \ \ \ \ \ \ \ \ \ \ \ \ \ \ \ \ \ \ \ \ \left(
\begin{array}
[c]{c}%
\text{here, we have substituted }cu\text{ for }u\text{ in the sum,}\\
\text{since the map (\ref{pf.lem.qr.jacobs.uuh.bij}) is a bijection}%
\end{array}
\right) \\
&  =\sum_{u\in\mathbb{Z}/p}\underbrace{K_{p}\left(  c^{3}u\left(
u^{2}+d\right)  \right)  }_{\substack{=K_{p}\left(  c^{3}\right)  K_{p}\left(
u\left(  u^{2}+d\right)  \right)  \\\text{(by Lemma \ref{lem.qr.jacobs.Kpprod}%
)}}}=K_{p}\left(  \underbrace{c^{3}}_{=cc^{2}}\right)  \underbrace{\sum
_{u\in\mathbb{Z}/p}K_{p}\left(  u\left(  u^{2}+d\right)  \right)
}_{\substack{=W\left(  d\right)  \\\text{(by the definition of }W\left(
d\right)  \text{)}}}\\
&  =\underbrace{K_{p}\left(  cc^{2}\right)  }_{\substack{=K_{p}\left(
c\right)  K_{p}\left(  c^{2}\right)  \\\text{(by Lemma
\ref{lem.qr.jacobs.Kpprod})}}}W\left(  d\right)  =K_{p}\left(  c\right)
\underbrace{K_{p}\left(  c^{2}\right)  }_{=1}W\left(  d\right)  =K_{p}\left(
c\right)  W\left(  d\right)  .
\end{align*}
This proves Lemma \ref{lem.qr.jacobs.uuh}.
\end{proof}

\begin{lemma}
\label{lem.qr.jacobs.dg}Let $g\in\mathbb{Z}/p$ be an element of $\mathbb{Z}/p$
that is not a square. Let $d\in\mathbb{Z}/p$ be a further element. Then:

\begin{enumerate}
\item[\textbf{(a)}] If $d\in\mathbb{Z}/p$ is a nonzero square, then $\left(
W\left(  d\right)  \right)  ^{2}=\left(  W\left(  \overline{1}\right)
\right)  ^{2}$.

\item[\textbf{(b)}] If $d\in\mathbb{Z}/p$ is not a square, then $\left(
W\left(  d\right)  \right)  ^{2}=\left(  W\left(  g\right)  \right)  ^{2}$.
\end{enumerate}
\end{lemma}

\begin{proof}
\textbf{(a)} Assume that $d\in\mathbb{Z}/p$ is a nonzero square. Thus,
$d=c^{2}$ for some $c\in\mathbb{Z}/p$. Consider this $c$. Then, $c\neq0$
(since $c^{2}=d\neq0$). Thus, Lemma \ref{lem.qr.jacobs.Kpdef} \textbf{(f)}
(applied to $u=c$) shows that $\left(  K_{p}\left(  c\right)  \right)  ^{2}%
=1$. However, Lemma \ref{lem.qr.jacobs.uuh} (applied to $\overline{1}$ instead
of $d$) yields $W\left(  c^{2}\cdot\overline{1}\right)  =K_{p}\left(
c\right)  W\left(  \overline{1}\right)  $. In view of $c^{2}\cdot\overline
{1}=c^{2}=d$, we can rewrite this as $W\left(  d\right)  =K_{p}\left(
c\right)  W\left(  \overline{1}\right)  $. Squaring this equality, we obtain%
\[
\left(  W\left(  d\right)  \right)  ^{2}=\left(  K_{p}\left(  c\right)
W\left(  \overline{1}\right)  \right)  ^{2}=\underbrace{\left(  K_{p}\left(
c\right)  \right)  ^{2}}_{=1}\left(  W\left(  \overline{1}\right)  \right)
^{2}=\left(  W\left(  \overline{1}\right)  \right)  ^{2}.
\]
This proves Lemma \ref{lem.qr.jacobs.dg} \textbf{(a)}.

\textbf{(b)} Assume that $d\in\mathbb{Z}/p$ is not a square. Thus, $d\neq0$.
Lemma \ref{lem.qr.jacobs.Kpdef} \textbf{(c)} yields $K_{p}\left(  d\right)
=-1$ (since $d$ is not a square). Also, Lemma \ref{lem.qr.jacobs.Kpdef}
\textbf{(c)} yields $K_{p}\left(  g\right)  =-1$ (since $g$ is not a square).

We have $g\neq0$ (since $g$ is not a square). Hence, $\dfrac{d}{g}%
\in\mathbb{Z}/p$ is well-defined (since $\mathbb{Z}/p$ is a field). We have
$d=\dfrac{d}{g}\cdot g$, so that%
\[
K_{p}\left(  d\right)  =K_{p}\left(  \dfrac{d}{g}\cdot g\right)  =K_{p}\left(
\dfrac{d}{g}\right)  \cdot K_{p}\left(  g\right)  \ \ \ \ \ \ \ \ \ \ \left(
\text{by Lemma \ref{lem.qr.jacobs.Kpprod}}\right)  .
\]
In view of $K_{p}\left(  d\right)  =-1$, this rewrites as $-1=K_{p}\left(
\dfrac{d}{g}\right)  \cdot\underbrace{K_{p}\left(  g\right)  }_{=-1}%
=-K_{p}\left(  \dfrac{d}{g}\right)  $, so that $K_{p}\left(  \dfrac{d}%
{g}\right)  =1$.

However, if $\dfrac{d}{g}$ was not a square, then Lemma
\ref{lem.qr.jacobs.Kpdef} \textbf{(c)} would yield $K_{p}\left(  \dfrac{d}%
{g}\right)  =-1$, which would contradict $K_{p}\left(  \dfrac{d}{g}\right)
=1$. Thus, $\dfrac{d}{g}$ must be a square. In other words, $\dfrac{d}%
{g}=c^{2}$ for some $c\in\mathbb{Z}/p$. Consider this $c$. Then, $c^{2}%
g=d\neq0$, so that $c\neq0$. Thus, Lemma \ref{lem.qr.jacobs.Kpdef}
\textbf{(f)} (applied to $u=c$) shows that $\left(  K_{p}\left(  c\right)
\right)  ^{2}=1$. However, Lemma \ref{lem.qr.jacobs.uuh} (applied to $g$
instead of $d$) yields $W\left(  c^{2}g\right)  =K_{p}\left(  c\right)
W\left(  g\right)  $. In view of $c^{2}g=d$, we can rewrite this as $W\left(
d\right)  =K_{p}\left(  c\right)  W\left(  g\right)  $. Squaring this
equality, we obtain%
\[
\left(  W\left(  d\right)  \right)  ^{2}=\left(  K_{p}\left(  c\right)
W\left(  g\right)  \right)  ^{2}=\underbrace{\left(  K_{p}\left(  c\right)
\right)  ^{2}}_{=1}\left(  W\left(  g\right)  \right)  ^{2}=\left(  W\left(
g\right)  \right)  ^{2}.
\]
This proves Lemma \ref{lem.qr.jacobs.dg} \textbf{(b)}.
\end{proof}

\begin{lemma}
\label{lem.qr.jacobs.sumsq}Let $g\in\mathbb{Z}/p$ be an element of
$\mathbb{Z}/p$ that is not a square. Then,%
\[
\sum_{d\in\mathbb{Z}/p}\left(  W\left(  d\right)  \right)  ^{2}=\dfrac{p-1}%
{2}\left(  \left(  W\left(  \overline{1}\right)  \right)  ^{2}+\left(
W\left(  g\right)  \right)  ^{2}\right)  .
\]

\end{lemma}

\begin{proof}
We have%
\begin{align*}
\sum_{d\in\mathbb{Z}/p}\left(  W\left(  d\right)  \right)  ^{2}  &
=\underbrace{\left(  W\left(  \overline{0}\right)  \right)  ^{2}%
}_{\substack{=0\\\text{(since Lemma \ref{lem.qr.jacobs.0}}\\\text{yields
}W\left(  \overline{0}\right)  =0\text{)}}}+\sum_{\substack{d\in
\mathbb{Z}/p;\\d\neq\overline{0}}}\left(  W\left(  d\right)  \right)
^{2}=\sum_{\substack{d\in\mathbb{Z}/p;\\d\neq\overline{0}}}\left(  W\left(
d\right)  \right)  ^{2}\\
&  =\sum_{\substack{d\in\mathbb{Z}/p;\\d\neq\overline{0};\\d\text{ is a
square}}}\ \ \underbrace{\left(  W\left(  d\right)  \right)  ^{2}%
}_{\substack{=\left(  W\left(  \overline{1}\right)  \right)  ^{2}\\\text{(by
Lemma \ref{lem.qr.jacobs.dg} \textbf{(a)})}}}+\underbrace{\sum_{\substack{d\in
\mathbb{Z}/p;\\d\neq\overline{0};\\d\text{ is not a square}}}}%
_{\substack{=\sum_{\substack{d\in\mathbb{Z}/p;\\d\text{ is not a square}%
}}\\\text{(since the condition \textquotedblleft}d\neq\overline{0}%
\text{\textquotedblright}\\\text{follows from the}\\\text{condition
\textquotedblleft}d\text{ is not a square\textquotedblright)}}%
}\ \ \underbrace{\left(  W\left(  d\right)  \right)  ^{2}}_{\substack{=\left(
W\left(  g\right)  \right)  ^{2}\\\text{(by Lemma \ref{lem.qr.jacobs.dg}
\textbf{(b)})}}}\\
&  =\sum_{\substack{d\in\mathbb{Z}/p;\\d\neq\overline{0};\\d\text{ is a
square}}}\left(  W\left(  \overline{1}\right)  \right)  ^{2}+\sum
_{\substack{d\in\mathbb{Z}/p;\\d\text{ is not a square}}}\left(  W\left(
g\right)  \right)  ^{2}\\
&  =\underbrace{\left\vert \left\{  d\in\mathbb{Z}/p\ \mid\ d\neq\overline
{0}\text{, and }d\text{ is a square}\right\}  \right\vert }%
_{\substack{=\left(  \text{number of nonzero squares in }\mathbb{Z}/p\right)
\\=\left(  p-1\right)  /2\\\text{(by Proposition \ref{prop.qr.count}
\textbf{(a)})}}}\cdot\left(  W\left(  \overline{1}\right)  \right)  ^{2}\\
&  \ \ \ \ \ \ \ \ \ \ +\underbrace{\left\vert \left\{  d\in\mathbb{Z}%
/p\ \mid\text{ }d\text{ is not a square}\right\}  \right\vert }%
_{\substack{=\left(  \text{number of elements of }\mathbb{Z}/p\text{ that are
not squares}\right)  \\=\left(  p-1\right)  /2\\\text{(by Proposition
\ref{prop.qr.count} \textbf{(c)})}}}\cdot\left(  W\left(  g\right)  \right)
^{2}\\
&  =\left(  \left(  p-1\right)  /2\right)  \cdot\left(  W\left(  \overline
{1}\right)  \right)  ^{2}+\left(  \left(  p-1\right)  /2\right)  \cdot\left(
W\left(  g\right)  \right)  ^{2}\\
&  =\dfrac{p-1}{2}\left(  \left(  W\left(  \overline{1}\right)  \right)
^{2}+\left(  W\left(  g\right)  \right)  ^{2}\right)  .
\end{align*}
This proves Lemma \ref{lem.qr.jacobs.sumsq}.
\end{proof}

\begin{lemma}
\label{lem.qr.jacobs.sumxy}Let $x,y\in\mathbb{Z}/p$. Then,%
\[
\sum_{d\in\mathbb{Z}/p}K_{p}\left(  \left(  x^{2}+d\right)  \left(
y^{2}+d\right)  \right)  =-1+%
\begin{cases}
p, & \text{if }x^{2}=y^{2};\\
0, & \text{otherwise.}%
\end{cases}
\]

\end{lemma}

\begin{proof}
Write the residue classes $x$ and $y$ as $x=\overline{a}$ and $y=\overline{b}$
for some $a,b\in\mathbb{Z}$. Let $k:=b^{2}-a^{2}\in\mathbb{Z}$. Thus, in
$\mathbb{Z}/p$, we have $\overline{k}=\overline{b^{2}-a^{2}}=\overline{b}%
^{2}-\overline{a}^{2}=y^{2}-x^{2}$ (since $\overline{a}=x$ and $\overline
{b}=y$).

The map%
\begin{align*}
\mathbb{Z}/p  &  \rightarrow\mathbb{Z}/p,\\
v  &  \mapsto v-y^{2}%
\end{align*}
is a bijection. Hence, we can substitute $v-y^{2}$ for $d$ in the sum
$\sum_{d\in\mathbb{Z}/p}K_{p}\left(  \left(  x^{2}+d\right)  \left(
y^{2}+d\right)  \right)  $. Thus, we find%
\begin{align}
&  \sum_{d\in\mathbb{Z}/p}K_{p}\left(  \left(  x^{2}+d\right)  \left(
y^{2}+d\right)  \right) \nonumber\\
&  =\sum_{v\in\mathbb{Z}/p}K_{p}\left(  \underbrace{\left(  x^{2}+\left(
v-y^{2}\right)  \right)  }_{=v-\left(  y^{2}-x^{2}\right)  }%
\underbrace{\left(  y^{2}+\left(  v-y^{2}\right)  \right)  }_{=v}\right)
\nonumber\\
&  =\sum_{v\in\mathbb{Z}/p}K_{p}\left(  \underbrace{\left(  v-\left(
y^{2}-x^{2}\right)  \right)  v}_{=v\left(  v-\left(  y^{2}-x^{2}\right)
\right)  }\right)  =\sum_{v\in\mathbb{Z}/p}K_{p}\left(  v\left(
v-\underbrace{\left(  y^{2}-x^{2}\right)  }_{=\overline{k}}\right)  \right)
\nonumber\\
&  =\sum_{v\in\mathbb{Z}/p}K_{p}\left(  v\left(  v-\overline{k}\right)
\right)  =\sum_{i=0}^{p-1}K_{p}\left(  \underbrace{\overline{i}\left(
\overline{i}-\overline{k}\right)  }_{=\overline{i\left(  i-k\right)  }}\right)
\nonumber\\
&  \ \ \ \ \ \ \ \ \ \ \ \ \ \ \ \ \ \ \ \ \left(
\begin{array}
[c]{c}%
\text{here, we have substituted }\overline{i}\text{ for }v\text{ in the
sum,}\\
\text{since the map }\left\{  0,1,\ldots,p-1\right\}  \rightarrow
\mathbb{Z}/p\text{ that}\\
\text{sends each }i\text{ to }\overline{i}\text{ is a bijection}%
\end{array}
\right) \nonumber\\
&  =\sum_{i=0}^{p-1}\underbrace{K_{p}\left(  \overline{i\left(  i-k\right)
}\right)  }_{\substack{=\left(  \dfrac{i\left(  i-k\right)  }{p}\right)
\\\text{(by (\ref{eq.def.qr.Kpu.Kpa}))}}}=\sum_{i=0}^{p-1}\left(
\dfrac{i\left(  i-k\right)  }{p}\right) \nonumber\\
&  =%
\begin{cases}
p-1, & \text{if }p\mid k;\\
-1, & \text{if }p\nmid k
\end{cases}
\ \ \ \ \ \ \ \ \ \ \left(  \text{by Proposition \ref{prop.qr.sumii-k}}\right)
\nonumber\\
&  =%
\begin{cases}
p-1, & \text{if }p\mid k;\\
-1, & \text{otherwise}%
\end{cases}
\nonumber\\
&  =-1+%
\begin{cases}
p, & \text{if }p\mid k;\\
0, & \text{otherwise.}%
\end{cases}
\label{pf.lem.qr.jacobs.sumxy.at}%
\end{align}

However, the statement \textquotedblleft$p\mid k$\textquotedblright\ is
equivalent to \textquotedblleft$x^{2}=y^{2}$\textquotedblright\ (because we
have the following chain of logical equivalences:%
\begin{align*}
\left(  p\mid k\right)  \  &  \Longleftrightarrow\ \left(  \overline
{k}=\overline{0}\text{ in }\mathbb{Z}/p\right) \\
&  \Longleftrightarrow\ \left(  y^{2}-x^{2}=\overline{0}\text{ in }%
\mathbb{Z}/p\right)  \ \ \ \ \ \ \ \ \ \ \left(  \text{since }\overline
{k}=y^{2}-x^{2}\right) \\
&  \Longleftrightarrow\ \left(  y^{2}=x^{2}\right)  \ \Longleftrightarrow
\ \left(  x^{2}=y^{2}\right)
\end{align*}
). Thus, we can rewrite (\ref{pf.lem.qr.jacobs.sumxy.at}) as
\[
\sum_{d\in\mathbb{Z}/p}K_{p}\left(  \left(  x^{2}+d\right)  \left(
y^{2}+d\right)  \right)  =-1+%
\begin{cases}
p, & \text{if }x^{2}=y^{2};\\
0, & \text{otherwise.}%
\end{cases}
\]
This proves Lemma \ref{lem.qr.jacobs.sumxy}.
\end{proof}

\begin{lemma}
\label{lem.qr.jacobs.sumsq2}We have
\[
\sum_{d\in\mathbb{Z}/p}\left(  W\left(  d\right)  \right)  ^{2}=p\left(
p-1\right)  \left(  1+\left(  -1\right)  ^{\left(  p-1\right)  /2}\right)  .
\]

\end{lemma}

\begin{proof}
For each $d\in\mathbb{Z}/p$, we have%
\begin{align*}
\left(  W\left(  d\right)  \right)  ^{2}  &  =\left(  \sum_{u\in\mathbb{Z}%
/p}K_{p}\left(  u\left(  u^{2}+d\right)  \right)  \right)  ^{2}%
\ \ \ \ \ \ \ \ \ \ \left(  \text{by the definition of }W\left(  d\right)
\right) \\
&  =\left(  \sum_{u\in\mathbb{Z}/p}K_{p}\left(  u\left(  u^{2}+d\right)
\right)  \right)  \left(  \sum_{u\in\mathbb{Z}/p}K_{p}\left(  u\left(
u^{2}+d\right)  \right)  \right) \\
&  =\left(  \sum_{x\in\mathbb{Z}/p}K_{p}\left(  x\left(  x^{2}+d\right)
\right)  \right)  \left(  \sum_{y\in\mathbb{Z}/p}K_{p}\left(  y\left(
y^{2}+d\right)  \right)  \right) \\
&  \ \ \ \ \ \ \ \ \ \ \ \ \ \ \ \ \ \ \ \ \left(  \text{here, we have renamed
both summation indices}\right) \\
&  =\sum_{\left(  x,y\right)  \in\left(  \mathbb{Z}/p\right)  ^{2}%
}\underbrace{K_{p}\left(  x\left(  x^{2}+d\right)  \right)  \cdot K_{p}\left(
y\left(  y^{2}+d\right)  \right)  }_{\substack{=K_{p}\left(  x\left(
x^{2}+d\right)  \cdot y\left(  y^{2}+d\right)  \right)  \\\text{(by Lemma
\ref{lem.qr.jacobs.Kpprod},}\\\text{applied to }u=x\left(  x^{2}+d\right)
\text{ and }v=y\left(  y^{2}+d\right)  \text{)}}}\\
&  =\sum_{\left(  x,y\right)  \in\left(  \mathbb{Z}/p\right)  ^{2}}%
K_{p}\left(  \underbrace{x\left(  x^{2}+d\right)  \cdot y\left(
y^{2}+d\right)  }_{=\left(  xy\right)  \left(  \left(  x^{2}+d\right)  \left(
y^{2}+d\right)  \right)  }\right) \\
&  =\sum_{\left(  x,y\right)  \in\left(  \mathbb{Z}/p\right)  ^{2}%
}\underbrace{K_{p}\left(  \left(  xy\right)  \left(  \left(  x^{2}+d\right)
\left(  y^{2}+d\right)  \right)  \right)  }_{\substack{=K_{p}\left(
xy\right)  \cdot K_{p}\left(  \left(  x^{2}+d\right)  \left(  y^{2}+d\right)
\right)  \\\text{(by Lemma \ref{lem.qr.jacobs.Kpprod})}}}\\
&  =\sum_{\left(  x,y\right)  \in\left(  \mathbb{Z}/p\right)  ^{2}}%
K_{p}\left(  xy\right)  \cdot K_{p}\left(  \left(  x^{2}+d\right)  \left(
y^{2}+d\right)  \right)  .
\end{align*}
Summing this equality over all $d\in\mathbb{Z}/p$, we obtain%
\begin{align*}
\sum_{d\in\mathbb{Z}/p}\left(  W\left(  d\right)  \right)  ^{2}  &
=\sum_{d\in\mathbb{Z}/p}\ \ \sum_{\left(  x,y\right)  \in\left(
\mathbb{Z}/p\right)  ^{2}}K_{p}\left(  xy\right)  \cdot K_{p}\left(  \left(
x^{2}+d\right)  \left(  y^{2}+d\right)  \right) \\
&  =\sum_{\left(  x,y\right)  \in\left(  \mathbb{Z}/p\right)  ^{2}}%
K_{p}\left(  xy\right)  \underbrace{\sum_{d\in\mathbb{Z}/p}K_{p}\left(
\left(  x^{2}+d\right)  \left(  y^{2}+d\right)  \right)  }_{\substack{=-1+%
\begin{cases}
p, & \text{if }x^{2}=y^{2};\\
0, & \text{otherwise}%
\end{cases}
\\\text{(by Lemma \ref{lem.qr.jacobs.sumxy})}}}\\
&  =\sum_{\left(  x,y\right)  \in\left(  \mathbb{Z}/p\right)  ^{2}}%
K_{p}\left(  xy\right)  \left(  -1+%
\begin{cases}
p, & \text{if }x^{2}=y^{2};\\
0, & \text{otherwise}%
\end{cases}
\right) \\
&  =\sum_{\left(  x,y\right)  \in\left(  \mathbb{Z}/p\right)  ^{2}}%
K_{p}\left(  xy\right)  \left(  -1\right)  +\sum_{\left(  x,y\right)
\in\left(  \mathbb{Z}/p\right)  ^{2}}K_{p}\left(  xy\right)
\begin{cases}
p, & \text{if }x^{2}=y^{2};\\
0, & \text{otherwise.}%
\end{cases}
\end{align*}

Let us now simplify the two sums on the right hand side separately.

First, we note that
\begin{align*}
\sum_{\left(  x,y\right)  \in\left(  \mathbb{Z}/p\right)  ^{2}}K_{p}\left(
xy\right)  \left(  -1\right)   &  =-\sum_{\left(  x,y\right)  \in\left(
\mathbb{Z}/p\right)  ^{2}}K_{p}\left(  xy\right) \\
&  =-\sum_{x\in\mathbb{Z}/p}\ \ \sum_{y\in\mathbb{Z}/p}\underbrace{K_{p}%
\left(  xy\right)  }_{\substack{=K_{p}\left(  x\right)  \cdot K_{p}\left(
y\right)  \\\text{(by Lemma \ref{lem.qr.jacobs.Kpprod})}}}\\
&  =-\sum_{x\in\mathbb{Z}/p}\ \ \sum_{y\in\mathbb{Z}/p}K_{p}\left(  x\right)
\cdot K_{p}\left(  y\right) \\
&  =-\sum_{x\in\mathbb{Z}/p}K_{p}\left(  x\right)  \underbrace{\sum
_{y\in\mathbb{Z}/p}K_{p}\left(  y\right)  }_{\substack{=\sum_{u\in
\mathbb{Z}/p}K_{p}\left(  u\right)  =0\\\text{(by Lemma
\ref{lem.qr.jacobs.sumKp})}}}\\
&  =-\underbrace{\sum_{x\in\mathbb{Z}/p}K_{p}\left(  x\right)  0}_{=0}=0.
\end{align*}

The second sum is not much trickier by now. We have%
\begin{align*}
&  \sum_{\left(  x,y\right)  \in\left(  \mathbb{Z}/p\right)  ^{2}}K_{p}\left(
xy\right)
\begin{cases}
p, & \text{if }x^{2}=y^{2};\\
0, & \text{otherwise}%
\end{cases}
\\
&  =\sum_{\substack{\left(  x,y\right)  \in\left(  \mathbb{Z}/p\right)
^{2};\\x^{2}=y^{2}}}K_{p}\left(  xy\right)  p\ \ \ \ \ \ \ \ \ \ \left(
\begin{array}
[c]{c}%
\text{since all the addends in this sum are }0\\
\text{except for those with }x^{2}=y^{2}%
\end{array}
\right) \\
&  =\sum_{\substack{\left(  x,y\right)  \in\left(  \mathbb{Z}/p\right)
^{2};\\x^{2}=y^{2};\\y=0}}K_{p}\left(  \underbrace{xy}%
_{\substack{=0\\\text{(since }y=0\text{)}}}\right)  p+\sum_{\substack{\left(
x,y\right)  \in\left(  \mathbb{Z}/p\right)  ^{2};\\x^{2}=y^{2};\\y\neq0}%
}K_{p}\left(  xy\right)  p\\
&  =\sum_{\substack{\left(  x,y\right)  \in\left(  \mathbb{Z}/p\right)
^{2};\\x^{2}=y^{2};\\y=0}}\underbrace{K_{p}\left(  0\right)  }%
_{\substack{=0\\\text{(by Lemma \ref{lem.qr.jacobs.Kpdef} \textbf{(a)})}%
}}p+\sum_{\substack{\left(  x,y\right)  \in\left(  \mathbb{Z}/p\right)
^{2};\\x^{2}=y^{2};\\y\neq0}}K_{p}\left(  xy\right)  p\\
&  =\underbrace{\sum_{\substack{\left(  x,y\right)  \in\left(  \mathbb{Z}%
/p\right)  ^{2};\\x^{2}=y^{2};\\y=0}}0p}_{=0}+\sum_{\substack{\left(
x,y\right)  \in\left(  \mathbb{Z}/p\right)  ^{2};\\x^{2}=y^{2};\\y\neq0}%
}K_{p}\left(  xy\right)  p\\
&  =\sum_{\substack{\left(  x,y\right)  \in\left(  \mathbb{Z}/p\right)
^{2};\\x^{2}=y^{2};\\y\neq0}}K_{p}\left(  xy\right)  p=\sum_{\substack{y\in
\mathbb{Z}/p;\\y\neq0}}\ \ \underbrace{\sum_{\substack{x\in\mathbb{Z}%
/p;\\x^{2}=y^{2}}}K_{p}\left(  xy\right)  }_{\substack{=1+\left(  -1\right)
^{\left(  p-1\right)  /2}\\\text{(by Lemma \ref{lem.qr.jacobs.Kpxy})}}}p\\
&  =\sum_{\substack{y\in\mathbb{Z}/p;\\y\neq0}}\left(  1+\left(  -1\right)
^{\left(  p-1\right)  /2}\right)  p=\left(  p-1\right)  \left(  1+\left(
-1\right)  ^{\left(  p-1\right)  /2}\right)  p
\end{align*}
(since this sum has $p-1$ addends).

Now, we combine all we have found: We have%
\begin{align*}
\sum_{d\in\mathbb{Z}/p}\left(  W\left(  d\right)  \right)  ^{2}  &
=\underbrace{\sum_{\left(  x,y\right)  \in\left(  \mathbb{Z}/p\right)  ^{2}%
}K_{p}\left(  xy\right)  \left(  -1\right)  }_{=0}+\underbrace{\sum_{\left(
x,y\right)  \in\left(  \mathbb{Z}/p\right)  ^{2}}K_{p}\left(  xy\right)
\begin{cases}
p, & \text{if }x^{2}=y^{2};\\
0, & \text{otherwise}%
\end{cases}
}_{=\left(  p-1\right)  \left(  1+\left(  -1\right)  ^{\left(  p-1\right)
/2}\right)  p}\\
&  =\left(  p-1\right)  \left(  1+\left(  -1\right)  ^{\left(  p-1\right)
/2}\right)  p=p\left(  p-1\right)  \left(  1+\left(  -1\right)  ^{\left(
p-1\right)  /2}\right)  .
\end{align*}
This proves Lemma \ref{lem.qr.jacobs.sumsq2}.
\end{proof}

\begin{proof}
[Proof of Theorem \ref{thm.qr.jacobs}.]We have $4\mid p-1$ (since
$p\equiv1\operatorname{mod}4$), so that $\left(  p-1\right)  /2$ is even.
Thus, $\left(  -1\right)  ^{\left(  p-1\right)  /2}=1$.

Comparing the equalities (\ref{eq.thm.qr.jacobs.defWh}) and
(\ref{eq.def.qr.Kpu.Wh}), we see that%
\[
W\left(  \overline{h}\right)  =W\left(  h\right)
\ \ \ \ \ \ \ \ \ \ \text{for each }h\in\mathbb{Z}\text{.}%
\]
Thus, in particular, $W\left(  \overline{1}\right)  =W\left(  1\right)  =a$
and $W\left(  \overline{m}\right)  =W\left(  m\right)  =b$. Moreover,
$\overline{m}\in\mathbb{Z}/p$ is not a square (by the definition of $m$).

Lemma \ref{lem.qr.jacobs.sumsq2} yields%
\[
\sum_{d\in\mathbb{Z}/p}\left(  W\left(  d\right)  \right)  ^{2}=p\left(
p-1\right)  \left(  1+\underbrace{\left(  -1\right)  ^{\left(  p-1\right)
/2}}_{=1}\right)  =p\left(  p-1\right)  \left(  1+1\right)  =2p\left(
p-1\right)  .
\]
Hence,%
\[
2p\left(  p-1\right)  =\sum_{d\in\mathbb{Z}/p}\left(  W\left(  d\right)
\right)  ^{2}=\dfrac{p-1}{2}\left(  \left(  W\left(  \overline{1}\right)
\right)  ^{2}+\left(  W\left(  \overline{m}\right)  \right)  ^{2}\right)
\]
(by Lemma \ref{lem.qr.jacobs.sumsq}, applied to $g=\overline{m}$). Dividing
both sides of this equality by $2\left(  p-1\right)  $, we obtain%
\[
p=\dfrac{1}{4}\left(  \left(  \underbrace{W\left(  \overline{1}\right)  }%
_{=a}\right)  ^{2}+\left(  \underbrace{W\left(  \overline{m}\right)  }%
_{=b}\right)  ^{2}\right)  =\dfrac{1}{4}\left(  a^{2}+b^{2}\right)  =\left(
a/2\right)  ^{2}+\left(  b/2\right)  ^{2}.
\]
This proves Theorem \ref{thm.qr.jacobs} \textbf{(b)}.

\textbf{(a)} We have $a=W\left(  \overline{1}\right)  $, which is even (by
Lemma \ref{lem.qr.jacobs.even}). Also, we have $b=W\left(  \overline
{m}\right)  $, which is even (by Lemma \ref{lem.qr.jacobs.even}). Thus,
Theorem \ref{thm.qr.jacobs} \textbf{(a)} is proven.
\end{proof}

Of course, Theorem \ref{thm.qr.jacobs} gives a new proof of Theorem
\ref{thm.fermat.p=xx+yy} (because Theorem \ref{thm.qr.jacobs} \textbf{(a)}
shows that $a/2$ and $b/2$ are integers, and Theorem \ref{thm.qr.jacobs}
\textbf{(b)} shows that $p$ is the sum of the squares of these two integers).
\medskip

Some curious variants of Theorem \ref{thm.qr.jacobs} can be found in
\cite{Whitem52}, \cite{ChLoYa11}, \cite{Willia23} and \cite[Theorem
3.16]{Nica25}; I suspect that there is more to be discovered. See also
\cite{BeEvWi98} and \cite{Nica25} for a wealth of related material.

\begin{exercise}
Make the same assumptions and definitions as in Theorem \ref{thm.qr.jacobs}.

\begin{enumerate}
\item[\textbf{(a)}] Prove that $a/2$ is odd, but $b/2$ is even.

\item[\textbf{(b)}] Now assume that $p\equiv1\operatorname{mod}8$. Prove that
$b/2$ is divisible by $4$, and thus we have $p=x^{2}+16y^{2}$ for $x=a/2$ and
$y=b/8$.
\end{enumerate}
\end{exercise}

\begin{exercise}
Make the same assumptions and definitions as in Theorem \ref{thm.qr.jacobs}.
Prove that%
\[
W\left(  h\right)  \equiv-h^{\left(  p-1\right)  /4}\dbinom{\left(
p-1\right)  /2}{\left(  p-1\right)  /4}\operatorname{mod}%
p\ \ \ \ \ \ \ \ \ \ \text{for any }h\in\mathbb{Z}\text{.}%
\]

\end{exercise}
\end{fineprint}

\subsection{\label{sec.finfield.primroot}Primitive roots}

We shall now state a curious (and actually rather useful) property of finite
fields that we won't prove until later.

We begin with a fun observation:

\begin{example}
The sequence of residue classes $\overline{1},\overline{2},\ldots,\overline
{6}$ in $\mathbb{Z}/7$ is an arithmetic sequence (in the sense that there
exists some \textquotedblleft difference\textquotedblright\ $d\in\mathbb{Z}/7$
such that each entry of this sequence equals the preceding entry plus $d$).

I claim that you can permute this sequence so that it becomes a geometric
sequence (in the sense that there exists some \textquotedblleft
quotient\textquotedblright\ $q\in\mathbb{Z}/7$ such that each entry of the
permuted sequence equals the preceding entry times $q$) !

Namely, $\overline{1},\overline{3},\overline{2},\overline{6},\overline
{4},\overline{5}$ is a geometric sequence. Its \textquotedblleft
quotient\textquotedblright\ is $\overline{3}$, meaning that each entry equals
the preceding entry times $\overline{3}$:%
\[
\overline{3}=\overline{1}\cdot\overline{3},\ \ \ \ \ \ \ \ \ \ \overline
{2}=\overline{3}\cdot\overline{3},\ \ \ \ \ \ \ \ \ \ \overline{6}%
=\overline{2}\cdot\overline{3},\ \ \ \ \ \ \ \ \ \ \ldots.
\]

\end{example}

This can be generalized: For any prime $p$, we can arrange the residue classes
$\overline{1},\overline{2},\ldots,\overline{p-1}$ in a geometric sequence.
Here is another way to put it:

\begin{theorem}
[Gauss]\label{thm.finfield.gauss1}Let $p$ be a prime. Then, there exists some
$g\in\left(  \mathbb{Z}/p\right)  ^{\times}$ such that its $p-1$ powers
$g^{0},g^{1},\ldots,g^{p-2}$ are distinct and satisfy%
\[
\left(  \mathbb{Z}/p\right)  ^{\times}=\left\{  g^{0},g^{1},\ldots
,g^{p-2}\right\}  .
\]
Such a $g$ is called a \textbf{primitive root} modulo $p$.
\end{theorem}

More generally:

\begin{theorem}
\label{thm.finfield.gauss2}Let $F$ be any finite field. Then, the group
$F^{\times}=F\setminus\left\{  0\right\}  $ is cyclic.
\end{theorem}

(This is easily seen to generalize Theorem \ref{thm.finfield.gauss1}: Indeed,
applying Theorem \ref{thm.finfield.gauss2} to $F=\mathbb{Z}/p$ yields that the
group $\left(  \mathbb{Z}/p\right)  ^{\times}$ is cyclic; since this group has
order $\left\vert \left(  \mathbb{Z}/p\right)  ^{\times}\right\vert =p-1$,
this yields further that this group has the form described in Theorem
\ref{thm.finfield.gauss1}. We will revisit this later on after proving Theorem
\ref{thm.finfield.gauss2}.)

Theorem \ref{thm.finfield.gauss2} can be generalized even further:

\begin{theorem}
\label{thm.finfield.gauss3}Let $F$ be any field. Let $G$ be a \textbf{finite}
subgroup of its group $F^{\times}=F\setminus\left\{  0\right\}  $ of units.
Then, $G$ is cyclic.
\end{theorem}

We will prove Theorem \ref{thm.finfield.gauss3} in two ways further below:
once using cyclotomic polynomials in Subsection \ref{subsec.cyclotomic.Fxcyc},
and once again using the classification of finite abelian groups in Section
\ref{sec.modPID.Fxcyc}. There are more elementary proofs as well. Keith
Conrad's note
\url{https://kconrad.math.uconn.edu/blurbs/grouptheory/cyclicmodp.pdf} and the
MathOverflow thread \url{https://mathoverflow.net/questions/54735/} contain
multiple proofs of Theorem \ref{thm.finfield.gauss1}, many of which generalize
to Theorem \ref{thm.finfield.gauss2} and to Theorem \ref{thm.finfield.gauss3}.

\newpage

\section{\label{chp.polys2}Polynomials II}

We shall now resume the study of polynomials.

\begin{convention}
We fix a commutative ring $R$. This convention will remain in force for the
entire chapter.
\end{convention}

\subsection{\label{sec.polys2.mulvar-again}Multivariate polynomials again}

Let us repeat Theorem \ref{thm.fieldext.basis-monic} (in a slightly shortened version):

\begin{theorem}
\label{thm.fieldext.basis-monic.repeat}Let $m\in\mathbb{N}$. Let $b\in
R\left[  x\right]  $ be a polynomial of degree $m$ such that its leading
coefficient $\left[  x^{m}\right]  b$ is a unit. Then, each element of
$R\left[  x\right]  /b$ can be uniquely written in the form%
\[
a_{0}\overline{x^{0}}+a_{1}\overline{x^{1}}+\cdots+a_{m-1}\overline{x^{m-1}%
}\ \ \ \ \ \ \ \ \ \ \text{with }a_{0},a_{1},\ldots,a_{m-1}\in R.
\]
Equivalently, the $m$ vectors $\overline{x^{0}},\overline{x^{1}}%
,\ldots,\overline{x^{m-1}}$ form a basis of the $R$-module $R\left[  x\right]
/b$. Thus, this $R$-module $R\left[  x\right]  /b$ is free of rank $m=\deg b$.
If $m>0$, then the ring $R\left[  x\right]  /b$ contains \textquotedblleft a
copy of $R$\textquotedblright.
\end{theorem}

Thus we understand quotients of univariate polynomial rings rather well when
the leading coefficient is a unit. They are less predictable when it is not a
unit. If $R$ is a field, however, then the leading coefficient of a nonzero
polynomial $b\in R\left[  x\right]  $ is always a unit, so we don't need to
worry about this issue.

But can we do this with multivariate polynomials?

Consider, for example, the two-variable polynomial ring $R\left[  x,y\right]
$. What does $R\left[  x,y\right]  /b$ look like for a polynomial $b\in
R\left[  x,y\right]  $ ? Keep in mind that the \textquotedblleft
idea\textquotedblright\ behind quotienting out $b$ is that we are setting $b$
to $0$. So $R\left[  x,y\right]  /b$ is \textquotedblleft the ring of
polynomials in $x$ and $y$ subject to the assumption that $b\left(
x,y\right)  =0$\textquotedblright.

Let us first try to answer this question for some special polynomials $b$; we
will then look for a pattern. There is a lot to be learned from the examples.

\subsubsection{Example 1: $R\left[  x,y\right]  /y$}

What is $R\left[  x,y\right]  /y$ ? We expect this to be isomorphic to
$R\left[  x\right]  $, because setting $y$ to $0$ in a polynomial $f\left(
x,y\right)  $ should give $f\left(  x,0\right)  \in R\left[  x\right]  $.

This is indeed true, and the formal proof is essentially just a formalization
of this informal argument:

\begin{proposition}
\label{prop.polring.Rxy/y}We have $R\left[  x,y\right]  /y\cong R\left[
x\right]  $ as $R$-algebras.
\end{proposition}

\begin{proof}
Define a map%
\begin{align*}
\alpha:R\left[  x,y\right]  /y  &  \rightarrow R\left[  x\right]  ,\\
\overline{f}  &  \mapsto f\left(  x,0\right)  .
\end{align*}
First, we need to check that this map $\alpha$ is well-defined. In other
words, we need to check the following:

\begin{statement}
\textit{Claim 1:} If $f,g\in R\left[  x,y\right]  $ are two polynomials
satisfying $\overline{f}=\overline{g}$ in $R\left[  x,y\right]  /y$, then
$f\left(  x,0\right)  =g\left(  x,0\right)  $.
\end{statement}

[\textit{Proof of Claim 1:} Let $f,g\in R\left[  x,y\right]  $ be two
polynomials satisfying $\overline{f}=\overline{g}$ in $R\left[  x,y\right]
/y$. Then, $\overline{f}=\overline{g}$ means that $f-g\in yR\left[
x,y\right]  $; in other words, $f-g=yp$ for some polynomial $p\in R\left[
x,y\right]  $. Consider this $p$. Now, evaluating both sides of the equality
$f-g=yp$ at $\left(  x,0\right)  $ (that is, substituting $0$ for $y$) yields
$f\left(  x,0\right)  -g\left(  x,0\right)  =0p\left(  x,0\right)  =0$ and
thus $f\left(  x,0\right)  =g\left(  x,0\right)  $. This proves Claim 1.]

Having proved Claim 1, we thus know that the map $\alpha$ is well-defined. It
is straightforward to see that $\alpha$ is an $R$-algebra morphism (because
the map $R\left[  x,y\right]  \rightarrow R\left[  x\right]  ,\ f\mapsto
f\left(  x,0\right)  $ is an $R$-algebra morphism\footnote{This is a
particular case of Theorem \ref{thm.polring.mulvar-sub-hom}.}).

In the opposite direction, define a map%
\begin{align*}
\beta:R\left[  x\right]   &  \rightarrow R\left[  x,y\right]  /y,\\
g  &  \mapsto\overline{g\left[  x\right]  }.
\end{align*}
It is again clear that this is an $R$-algebra morphism.

Now, we shall show that the maps $\alpha$ and $\beta$ are mutually inverse. To
prove this, we need to check that $\alpha\circ\beta=\operatorname*{id}$ and
$\beta\circ\alpha=\operatorname*{id}$. Checking $\alpha\circ\beta
=\operatorname*{id}$ is the easy part. The \textquotedblleft hard
part\textquotedblright\ is showing that $\beta\circ\alpha=\operatorname*{id}$.
There are two ways to do this:

[\textit{First proof of }$\beta\circ\alpha=\operatorname*{id}$\textit{:} To
show this, we need to prove that $\left(  \beta\circ\alpha\right)  \left(
\overline{f}\right)  =\overline{f}$ for each $f\in R\left[  x,y\right]  $. So
let us fix an $f\in R\left[  x,y\right]  $. Then,%
\begin{align*}
\left(  \beta\circ\alpha\right)  \left(  \overline{f}\right)   &
=\beta\left(  \alpha\left(  \overline{f}\right)  \right) \\
&  =\beta\left(  f\left(  x,0\right)  \right)  \ \ \ \ \ \ \ \ \ \ \left(
\text{since }\alpha\left(  \overline{f}\right)  \text{ was defined to be
}f\left(  x,0\right)  \right) \\
&  =\overline{\left(  f\left(  x,0\right)  \right)  \left[  x\right]
}\ \ \ \ \ \ \ \ \ \ \left(  \text{by the definition of }\beta\right) \\
&  =\overline{f\left(  x,0\right)  }\ \ \ \ \ \ \ \ \ \ \left(  \text{since
}\left(  f\left(  x,0\right)  \right)  \left[  x\right]  =f\left(  x,0\right)
\right)  .
\end{align*}
Thus, it remains to show that $\overline{f\left(  x,0\right)  }=\overline{f}$
(because we want to show that $\left(  \beta\circ\alpha\right)  \left(
\overline{f}\right)  =\overline{f}$). In other words, it remains to show that
$f-f\left(  x,0\right)  \in yR\left[  x,y\right]  $. We do this directly:
Write $f$ in the form $f=\sum\limits_{i,j\in\mathbb{N}}a_{i,j}x^{i}y^{j}$
(with $a_{i,j}\in R$). Then,%
\begin{align*}
f\left(  x,0\right)   &  =\sum_{i,j\in\mathbb{N}}a_{i,j}x^{i}0^{j}=\sum
_{i\in\mathbb{N}}a_{i,0}x^{i}\underbrace{0^{0}}_{=1}+\sum_{\substack{i,j\in
\mathbb{N};\\j>0}}a_{i,j}x^{i}\underbrace{0^{j}}_{\substack{=0\\\text{(since
}j>0\text{)}}}\\
&  \ \ \ \ \ \ \ \ \ \ \ \ \ \ \ \ \ \ \ \ \left(
\begin{array}
[c]{c}%
\text{here, we have split the sum into two parts:}\\
\text{one that contains all terms with }j=0\\
\text{and one that contains all terms with }j>0
\end{array}
\right) \\
&  =\sum_{i\in\mathbb{N}}a_{i,0}x^{i}=\sum_{\substack{i,j\in\mathbb{N}%
;\\j=0}}a_{i,j}x^{i}y^{j}\ \ \ \ \ \ \ \ \ \ \left(  \text{since }%
y^{j}=1\text{ and }a_{i,j}=a_{i,0}\text{ for }j=0\right)  .
\end{align*}
Subtracting this from $f=\sum\limits_{i,j\in\mathbb{N}}a_{i,j}x^{i}y^{j}$, we
find%
\begin{align*}
f-f\left(  x,0\right)   &  =\sum\limits_{i,j\in\mathbb{N}}a_{i,j}x^{i}%
y^{j}-\sum_{\substack{i,j\in\mathbb{N};\\j=0}}a_{i,j}x^{i}y^{j}=\sum
_{\substack{i,j\in\mathbb{N};\\j>0}}a_{i,j}x^{i}\underbrace{y^{j}%
}_{\substack{=yy^{j-1}\\\text{(we can do this}\\\text{because }j>0\text{)}}}\\
&  =\sum_{\substack{i,j\in\mathbb{N};\\j>0}}a_{i,j}x^{i}yy^{j-1}%
=y\sum_{\substack{i,j\in\mathbb{N};\\j>0}}a_{i,j}x^{i}y^{j-1}\in yR\left[
x,y\right]  ,
\end{align*}
as we wanted to prove. Thus, $\overline{f\left(  x,0\right)  }=\overline{f}$,
so that $\left(  \beta\circ\alpha\right)  \left(  \overline{f}\right)
=\overline{f\left(  x,0\right)  }=\overline{f}$. This proves $\beta\circ
\alpha=\operatorname*{id}$.]

[\textit{Second proof of }$\beta\circ\alpha=\operatorname*{id}$\textit{:} Here
is a more \textquotedblleft cultured\textquotedblright\ proof. We know that
$\beta$ and $\alpha$ are $R$-algebra morphisms, hence are $R$-linear maps.
Thus, $\beta\circ\alpha$ and $\operatorname*{id}$ are two $R$-linear maps from
$R\left[  x,y\right]  /y$ to $R\left[  x,y\right]  /y$. Our goal is to prove
that these two $R$-linear maps $\beta\circ\alpha$ and $\operatorname*{id}$ are
equal. As we have learned in Theorem \ref{thm.mods.linmap-unidef-1}, there is
a shortcut for proving that two $R$-linear maps are equal: It suffices to pick
a family of vectors that spans the domain (in our case, the $R$-module
$R\left[  x,y\right]  /y$), and to show that the two maps agree on the vectors
of this family. In our case, there is a rather natural choice of such a
family: the family of monomials, or rather of their cosets. That is, we choose
the family $\left(  \overline{x^{i}y^{j}}\right)  _{i,j\in\mathbb{N}}$. This
family spans the $R$-module $R\left[  x,y\right]  /y$ (since the family
$\left(  x^{i}y^{j}\right)  _{i,j\in\mathbb{N}}$ spans the $R$-module
$R\left[  x,y\right]  $, and since the canonical projection onto $R\left[
x,y\right]  /y$ clearly preserves their spanning property). Thus, we only need
to show that the two maps $\beta\circ\alpha$ and $\operatorname*{id}$ agree on
the vectors of this family -- i.e., to show that%
\[
\left(  \beta\circ\alpha\right)  \left(  \overline{x^{i}y^{j}}\right)
=\operatorname*{id}\left(  \overline{x^{i}y^{j}}\right)
\ \ \ \ \ \ \ \ \ \ \text{for any }i,j\in\mathbb{N}.
\]
But this is straightforward: We fix $i,j\in\mathbb{N}$, and set out to show
that $\left(  \beta\circ\alpha\right)  \left(  \overline{x^{i}y^{j}}\right)
=\operatorname*{id}\left(  \overline{x^{i}y^{j}}\right)  $. If $j>0$, then
$\overline{x^{i}y^{j}}=0$ (since $x^{i}y^{j}\in yR\left[  x,y\right]  $ in
this case) and therefore both $\left(  \beta\circ\alpha\right)  \left(
\overline{x^{i}y^{j}}\right)  $ and $\operatorname*{id}\left(  \overline
{x^{i}y^{j}}\right)  $ must be $0$ in this case (since $R$-linear maps always
send $0$ to $0$). If, on the other hand, $j=0$, then $\overline{x^{i}y^{j}%
}=\overline{x^{i}y^{0}}=\overline{x^{i}}$ and therefore $\alpha\left(
\overline{x^{i}y^{j}}\right)  =\alpha\left(  \overline{x^{i}}\right)  =x^{i}$
(since substituting $0$ for $y$ does not change the monomial $x^{i}$) and thus
$\left(  \beta\circ\alpha\right)  \left(  \overline{x^{i}y^{j}}\right)
=\beta\left(  x^{i}\right)  =\overline{x^{i}}=\overline{x^{i}y^{j}%
}=\operatorname*{id}\left(  \overline{x^{i}y^{j}}\right)  $. Hence, in both
cases, we have shown that $\left(  \beta\circ\alpha\right)  \left(
\overline{x^{i}y^{j}}\right)  =\operatorname*{id}\left(  \overline{x^{i}y^{j}%
}\right)  $. This completes the proof of $\beta\circ\alpha=\operatorname*{id}$.]

Either way, we have now shown that $\beta\circ\alpha=\operatorname*{id}$.
Combined with $\alpha\circ\beta=\operatorname*{id}$, this yields that the two
maps $\alpha$ and $\beta$ are mutually inverse. Thus, $\alpha$ is an
invertible $R$-algebra morphism, hence an $R$-algebra isomorphism. This proves
Proposition \ref{prop.polring.Rxy/y}.
\end{proof}

We can easily generalize this to multiple variables:

\begin{proposition}
\label{prop.polring.Rxmodxn}For any $n>0$, we have%
\[
R\left[  x_{1},x_{2},\ldots,x_{n}\right]  /x_{n}\cong R\left[  x_{1}%
,x_{2},\ldots,x_{n-1}\right]  \text{ as }R\text{-algebras.}%
\]

\end{proposition}

\begin{proof}
Same idea as for Proposition \ref{prop.polring.Rxy/y}, but requiring more
subscripts to juggle.
\end{proof}

\subsubsection{Example 2: $R\left[  x,y\right]  /\left(  x^{2}+y^{2}-1\right)
$}

What does $R\left[  x,y\right]  /\left(  x^{2}+y^{2}-1\right)  $ look like?

This is a fairly useful $R$-algebra; it can be viewed as the algebra of
polynomial functions on the unit circle. Indeed, any element $\overline{f}\in
R\left[  x,y\right]  /\left(  x^{2}+y^{2}-1\right)  $ can be \textquotedblleft
evaluated\textquotedblright\ at a point $\left(  a,b\right)  $ on the unit
circle (meaning, a pair of elements $a,b\in R$ with $a^{2}+b^{2}=1$). (Take
this geometric interpretation with a grain of salt: It takes some artistic
license to even speak of a \textquotedblleft circle\textquotedblright\ when
$R=\mathbb{Z}/2$.)

There are various interesting ring-theoretical questions to be asked about the
quotient ring $R\left[  x,y\right]  /\left(  x^{2}+y^{2}-1\right)  $; however,
let us restrict ourselves to studying it as an $R$-module. As an $R$-module,
is $R\left[  x,y\right]  /\left(  x^{2}+y^{2}-1\right)  $ free? What is a
basis? This boils down to asking whether (and how) we can divide polynomials
with remainder by $x^{2}+y^{2}-1$.

Here we will be helped by the following fact:

\begin{proposition}
\label{prop.polring.Rxy-stepwise}We have%
\[
R\left[  x,y\right]  \cong\left(  R\left[  x\right]  \right)  \left[
y\right]  \ \ \ \ \ \ \ \ \ \ \text{as }R\text{-algebras.}%
\]
More concretely, the map%
\begin{align*}
\varphi:R\left[  x,y\right]   &  \rightarrow\left(  R\left[  x\right]
\right)  \left[  y\right]  ,\\
\sum_{i,j\in\mathbb{N}}a_{i,j}x^{i}y^{j}  &  \mapsto\sum_{j\in\mathbb{N}%
}\left(  \sum_{i\in\mathbb{N}}a_{i,j}x^{i}\right)  y^{j}%
\ \ \ \ \ \ \ \ \ \ \left(  \text{where }a_{i,j}\in R\right)
\end{align*}
is an $R$-algebra isomorphism.
\end{proposition}

\begin{proof}
First of all, you are excused for wondering what the deal is: Isn't the above
map $\varphi$ just the identity map, since $\sum_{j\in\mathbb{N}}\left(
\sum_{i\in\mathbb{N}}a_{i,j}x^{i}\right)  y^{j}$ is the same polynomial as
$\sum_{i,j\in\mathbb{N}}a_{i,j}x^{i}y^{j}$ (just rewritten)?

Essentially yes, but there is a technical difference between the rings
$R\left[  x,y\right]  $ and $\left(  R\left[  x\right]  \right)  \left[
y\right]  $. The former is a polynomial ring in two indeterminates $x,y$ over
$R$, whereas the latter is a polynomial ring in one indeterminate $y$ over the
ring $R\left[  x\right]  $. Hence,

\begin{itemize}
\item the elements of $R\left[  x,y\right]  $ are polynomials in two variables
$x,y$ with coefficients in $R$, whereas

\item the elements of $\left(  R\left[  x\right]  \right)  \left[  y\right]  $
are polynomials in one variable $y$ with coefficients in $R\left[  x\right]  $
(that is, their coefficients themselves are polynomials in one variable $x$
over $R$).
\end{itemize}

Thus, even if a polynomial in $R\left[  x,y\right]  $ and a polynomial in
$\left(  R\left[  x\right]  \right)  \left[  y\right]  $ look exactly the same
(such as, for example, the polynomials $2x^{2}y^{3}$ in both rings), they are
technically different. (The polynomial $2x^{2}y^{3}$ in $R\left[  x,y\right]
$ has the monomial $x^{2}y^{3}$ appear in it with coefficient $2$, whereas the
polynomial $2x^{2}y^{3}$ in $\left(  R\left[  x\right]  \right)  \left[
y\right]  $ has the monomial $y^{3}$ appear in it with coefficient $2x^{2}$.)
The map $\varphi$ thus sends each polynomial in $R\left[  x,y\right]  $ to the
identically-looking polynomial in $\left(  R\left[  x\right]  \right)  \left[
y\right]  $.

This being said, the claim we are proving is saying precisely that the
difference between $R\left[  x,y\right]  $ and $\left(  R\left[  x\right]
\right)  \left[  y\right]  $ is only a technicality; in essence the two rings
are the same. The proof is rather straightforward. The simplest way is as
follows: The map $\varphi$ defined in the proposition is easily seen to be
well-defined and an $R$\textbf{-module} isomorphism. Thus, it remains to prove
that this map $\varphi$ respects multiplication and respects the unity. It is
clear enough that $\varphi$ respects the unity (since the unities of both
rings equal $x^{0}y^{0}$), so we only need to check that $\varphi$ respects
multiplication. According to Lemma \ref{lem.algebras.mor-on-basis}, it
suffices to prove this on a family of vectors that spans the $R$-module
$R\left[  x,y\right]  $; in other words, we only need to find a family
$\left(  m_{i}\right)  _{i\in I}$ of vectors in $R\left[  x,y\right]  $ that
spans $R\left[  x,y\right]  $, and show that%
\[
\varphi\left(  m_{i}m_{j}\right)  =\varphi\left(  m_{i}\right)  \varphi\left(
m_{j}\right)  \ \ \ \ \ \ \ \ \ \ \text{for all }i,j\in I.
\]
Fortunately, the family of monomials $\left(  x^{i}y^{j}\right)  _{\left(
i,j\right)  \in\mathbb{N}^{2}}$ is such a family of vectors (even better, it
is a basis of the $R$-module $R\left[  x,y\right]  $); thus, we only need to
prove that%
\[
\varphi\left(  x^{i}y^{u}\cdot x^{j}y^{v}\right)  =\varphi\left(  x^{i}%
y^{u}\right)  \cdot\varphi\left(  x^{j}y^{v}\right)
\ \ \ \ \ \ \ \ \ \ \text{for all }\left(  i,u\right)  ,\left(  j,v\right)
\in\mathbb{N}^{2}.
\]
But this is easy (the left and right hand sides both equal $x^{i+j}y^{u+v}%
\in\left(  R\left[  x\right]  \right)  \left[  y\right]  $). Thus, we conclude
that $\varphi$ respects multiplication; as we said above, this completes the
proof of Proposition \ref{prop.polring.Rxy-stepwise}.
\end{proof}

Now, in view of Proposition \ref{prop.polring.Rxy-stepwise}, we have the
$R$-algebra isomorphism%
\begin{equation}
R\left[  x,y\right]  /\left(  x^{2}+y^{2}-1\right)  \cong\left(  R\left[
x\right]  \right)  \left[  y\right]  /\left(  x^{2}+y^{2}-1\right)
\label{eq.polring.Rxy/xx+yy-1.1}%
\end{equation}
(since the isomorphism $\varphi$ from Proposition
\ref{prop.polring.Rxy-stepwise} sends the polynomial $x^{2}+y^{2}-1\in
R\left[  x,y\right]  $ to the identically-looking polynomial $x^{2}+y^{2}%
-1\in\left(  R\left[  x\right]  \right)  \left[  y\right]  $).

The ring on the right hand side of (\ref{eq.polring.Rxy/xx+yy-1.1}) is a
quotient ring of the \textbf{univariate} polynomial ring $\left(  R\left[
x\right]  \right)  \left[  y\right]  $ modulo the \textbf{monic} polynomial
$x^{2}+y^{2}-1=y^{2}+\underbrace{\left(  x^{2}-1\right)  }_{\text{constant
term in }R\left[  x\right]  }$ in the variable $y$. Thus, Theorem
\ref{thm.fieldext.basis-monic.repeat} (applied to $2$, $R\left[  x\right]  $,
$y$ and $x^{2}+y^{2}-1$ instead of $m$, $R$, $x$ and $b$) shows that this
quotient ring $\left(  R\left[  x\right]  \right)  \left[  y\right]  /\left(
x^{2}+y^{2}-1\right)  $ has a basis $\left(  \overline{y^{0}},\overline{y^{1}%
}\right)  $ as an $R\left[  x\right]  $-module. This means that any element of
$\left(  R\left[  x\right]  \right)  \left[  y\right]  /\left(  x^{2}%
+y^{2}-1\right)  $ can be uniquely written as%
\[
\alpha\overline{y^{0}}+\beta\overline{y^{1}}\ \ \ \ \ \ \ \ \ \ \text{for some
}\alpha,\beta\in R\left[  x\right]  .
\]
Since elements of $R\left[  x\right]  $ themselves can be uniquely written as
$R$-linear combinations of powers of $x$, we thus conclude that any element of
$\left(  R\left[  x\right]  \right)  \left[  y\right]  /\left(  x^{2}%
+y^{2}-1\right)  $ can be uniquely written as%
\begin{align*}
&  \left(  \alpha_{0}x^{0}+\alpha_{1}x^{1}+\alpha_{2}x^{2}+\cdots\right)
\overline{y^{0}}+\left(  \beta_{0}x^{0}+\beta_{1}x^{1}+\beta_{2}x^{2}%
+\cdots\right)  \overline{y^{1}}\\
&  =\overline{\left(  \alpha_{0}x^{0}+\alpha_{1}x^{1}+\alpha_{2}x^{2}%
+\cdots\right)  y^{0}+\left(  \beta_{0}x^{0}+\beta_{1}x^{1}+\beta_{2}%
x^{2}+\cdots\right)  y^{1}}\\
&  =\alpha_{0}\overline{x^{0}y^{0}}+\alpha_{1}\overline{x^{1}y^{0}}+\alpha
_{2}\overline{x^{2}y^{0}}+\cdots+\beta_{0}\overline{x^{0}y^{1}}+\beta
_{1}\overline{x^{1}y^{1}}+\beta_{2}\overline{x^{2}y^{1}}+\cdots
\end{align*}
for some $\alpha_{0},\alpha_{1},\alpha_{2},\ldots,\beta_{0},\beta_{1}%
,\beta_{2},\ldots\in R$ (with all but finitely many of these coefficients
$\alpha_{0},\alpha_{1},\alpha_{2},\ldots,\beta_{0},\beta_{1},\beta_{2},\ldots$
being $0$).

Thus, as an $R$-module, $\left(  R\left[  x\right]  \right)  \left[  y\right]
/\left(  x^{2}+y^{2}-1\right)  $ has a basis%
\[
\left(  \overline{x^{0}y^{0}},\overline{x^{1}y^{0}},\overline{x^{2}y^{0}%
},\ldots,\overline{x^{0}y^{1}},\overline{x^{1}y^{1}},\overline{x^{2}y^{1}%
},\ldots\right)  .
\]
In view of the $R$-algebra isomorphism (\ref{eq.polring.Rxy/xx+yy-1.1}) (which
sends each $\overline{x^{i}y^{j}}$ to $\overline{x^{i}y^{j}}$), we can thus
conclude that, as an $R$-module, $R\left[  x,y\right]  /\left(  x^{2}%
+y^{2}-1\right)  $ has a basis%
\begin{equation}
\left(  \overline{x^{0}y^{0}},\overline{x^{1}y^{0}},\overline{x^{2}y^{0}%
},\ldots,\overline{x^{0}y^{1}},\overline{x^{1}y^{1}},\overline{x^{2}y^{1}%
},\ldots\right)  . \label{eq.polring.Rxy/xx+yy-1.6}%
\end{equation}

\subsubsection{Indeterminates one at a time}

We digress from our series of examples in order to make a few comments about
Proposition \ref{prop.polring.Rxy-stepwise}. We first observe the following:

\begin{proposition}
\label{prop.polring.Rxy-stepwise-Rx}The map $\varphi:R\left[  x,y\right]
\rightarrow\left(  R\left[  x\right]  \right)  \left[  y\right]  $ from
Proposition \ref{prop.polring.Rxy-stepwise} is not just an $R$-algebra
isomorphism, but also an $R\left[  x\right]  $-algebra isomorphism. Here, we
view $R\left[  x,y\right]  $ as an $R\left[  x\right]  $-algebra via the ring
morphism%
\begin{align*}
R\left[  x\right]   &  \rightarrow R\left[  x,y\right]  ,\\
f  &  \mapsto f\left[  x\right]  .
\end{align*}

\end{proposition}

\begin{proof}
LTTR. (It only needs to be shown that $\varphi\left(  fg\right)
=f\varphi\left(  g\right)  $ for any $f\in R\left[  x\right]  $ and $g\in
R\left[  x,y\right]  $.)
\end{proof}

The order in which we list the variables doesn't matter much in a polynomial
ring; thus, Proposition \ref{prop.polring.Rxy-stepwise} has the following
analogue (which is proved similarly):

\begin{proposition}
\label{prop.polring.Rxy-stepwise2}We have%
\[
R\left[  x,y\right]  \cong\left(  R\left[  y\right]  \right)  \left[
x\right]  \ \ \ \ \ \ \ \ \ \ \text{as }R\text{-algebras.}%
\]
More concretely, the map%
\begin{align*}
\varphi:R\left[  x,y\right]   &  \rightarrow\left(  R\left[  y\right]
\right)  \left[  x\right]  ,\\
\sum_{i,j\in\mathbb{N}}a_{i,j}x^{i}y^{j}  &  \mapsto\sum_{i\in\mathbb{N}%
}\left(  \sum_{j\in\mathbb{N}}a_{i,j}y^{j}\right)  x^{i}%
\ \ \ \ \ \ \ \ \ \ \left(  \text{where }a_{i,j}\in R\right)
\end{align*}
is an $R$-algebra isomorphism.
\end{proposition}

Again, this isomorphism is an $R\left[  y\right]  $-algebra isomorphism
(similarly to Proposition \ref{prop.polring.Rxy-stepwise-Rx}). \medskip

Proposition \ref{prop.polring.Rxy-stepwise} can also be generalized to more
than $2$ variables:

\begin{proposition}
\label{prop.polring.Rxn-stepwise}For any $n>0$, we have%
\[
R\left[  x_{1},x_{2},\ldots,x_{n}\right]  \cong\left(  R\left[  x_{1}%
,x_{2},\ldots,x_{n-1}\right]  \right)  \left[  x_{n}\right]  \text{ as
}R\text{-algebras.}%
\]
More concretely, the map%
\begin{align*}
\varphi:R\left[  x_{1},x_{2},\ldots,x_{n}\right]   &  \rightarrow\left(
R\left[  x_{1},x_{2},\ldots,x_{n-1}\right]  \right)  \left[  x_{n}\right]  ,\\
\sum_{\left(  i_{1},i_{2},\ldots,i_{n}\right)  \in\mathbb{N}^{n}}%
a_{i_{1},i_{2},\ldots,i_{n}}x_{1}^{i_{1}}x_{2}^{i_{2}}\cdots x_{n}^{i_{n}}  &
\mapsto\sum_{j\in\mathbb{N}}\left(  \sum_{\left(  i_{1},i_{2},\ldots
,i_{n-1}\right)  \in\mathbb{N}^{n-1}}a_{i_{1},i_{2},\ldots,i_{n-1},j}%
x_{1}^{i_{1}}x_{2}^{i_{2}}\cdots x_{n-1}^{i_{n-1}}\right)  x_{n}^{j}\\
&  \ \ \ \ \ \ \ \ \ \ \ \ \ \ \ \ \ \ \ \ \left(  \text{where }a_{i_{1}%
,i_{2},\ldots,i_{n}}\in R\right)
\end{align*}
is an $R$-algebra isomorphism.
\end{proposition}

\begin{proof}
Generalize the proof of Proposition \ref{prop.polring.Rxy-stepwise} (same
idea, more subscripts).
\end{proof}

As in Proposition \ref{prop.polring.Rxy-stepwise-Rx}, the map $\varphi$ in
Proposition \ref{prop.polring.Rxn-stepwise} is not just an $R$-algebra
isomorphism but also an $R\left[  x_{1},x_{2},\ldots,x_{n-1}\right]  $-algebra isomorphism.

\subsubsection{More examples?}

Having understood the $R$-modules $R\left[  x,y\right]  /y$ and $R\left[
x,y\right]  /\left(  x^{2}+y^{2}-1\right)  $, we move on to further examples.

What does $R\left[  x,y\right]  /\left(  xy\right)  $ look like? We cannot
answer this using the methods used above, since the polynomial $xy$ is neither
monic in $y$ when considered as a polynomial in $\left(  R\left[  x\right]
\right)  \left[  y\right]  $ nor monic in $x$ when considered as a polynomial
in $\left(  R\left[  y\right]  \right)  \left[  x\right]  $.

What about $R\left[  x,y\right]  /\left(  xy\left(  x-y\right)  \right)  $ ?
Can we divide $\left(  x+y\right)  ^{3}$ by $xy\left(  x-y\right)  $ with
remainder? What is the remainder? Should we replace $x^{2}y$ by $xy^{2}$ or
vice versa?

To make things more complicated (but also more useful), let's not forget that
we can quotient a ring by an ideal, not just by a single element. Even if $R$
is a field, the polynomial ring $R\left[  x,y\right]  $ is not a PID (unlike
$R\left[  x\right]  $ for a field $R$), so not every ideal is principal.

The following shorthand will be useful:

\begin{definition}
\label{def.ideal.gen-by}Let $S$ be a commutative ring. Let $a_{1},a_{2}%
,\ldots,a_{k}$ be elements of $S$. Then, the ideal $a_{1}S+a_{2}S+\cdots
+a_{k}S$ (this is the set of all $S$-linear combinations of $a_{1}%
,a_{2},\ldots,a_{k}$) is called \textbf{the ideal generated by }$a_{1}%
,a_{2},\ldots,a_{k}$. The quotient ring $S/\left(  a_{1}S+a_{2}S+\cdots
+a_{k}S\right)  $ will be denoted by $S/\left(  a_{1},a_{2},\ldots
,a_{k}\right)  $.
\end{definition}

(Many authors actually write $\left(  a_{1},a_{2},\ldots,a_{k}\right)  $ for
the ideal $a_{1}S+a_{2}S+\cdots+a_{k}S$, but this risks confusion since
$\left(  a_{1},a_{2},\ldots,a_{k}\right)  $ also means the $k$-tuple.)

Informally, $S/\left(  a_{1},a_{2},\ldots,a_{k}\right)  $ is what is obtained
from $S$ if you set all of $a_{1},a_{2},\ldots,a_{k}$ to $0$.

For an example, we can look at $R\left[  x,y\right]  /\left(  x+y,x-y\right)
$. This behaves differently depending on $R$:

\begin{itemize}
\item If $R=\mathbb{Q}$, then
\[
R\left[  x,y\right]  /\left(  x+y,x-y\right)  =\mathbb{Q}\left[  x,y\right]
/\left(  x+y,x-y\right)  =\mathbb{Q}\left[  x,y\right]  /\left(  x,y\right)
\]
(since it is easy to see that the $\mathbb{Q}\left[  x,y\right]  $-linear
combinations of $x+y$ and $x-y$ are precisely the $\mathbb{Q}\left[
x,y\right]  $-linear combinations of $x$ and $y$), and thus%
\[
R\left[  x,y\right]  /\left(  x+y,x-y\right)  =\mathbb{Q}\left[  x,y\right]
/\left(  x,y\right)  \cong\mathbb{Q}.
\]

\item If $R=\mathbb{Z}/2$, then%
\begin{align*}
R\left[  x,y\right]  /\left(  x+y,x-y\right)   &  =\left(  \mathbb{Z}%
/2\right)  \left[  x,y\right]  /\left(  \underbrace{x+y}%
_{\substack{=x-y\\\text{(since we are in}\\\text{characteristic }2\text{)}%
}},x-y\right) \\
&  =\left(  \mathbb{Z}/2\right)  \left[  x,y\right]  /\left(  x-y,x-y\right)
\\
&  =\left(  \mathbb{Z}/2\right)  \left[  x,y\right]  /\left(  x-y\right)
\cong\left(  \mathbb{Z}/2\right)  \left[  x\right]  .
\end{align*}

\end{itemize}

We can easily come up with more complicated examples:

\begin{itemize}
\item What is $R\left[  x,y,z\right]  /\left(  x^{2}-yz,y^{2}-zx,z^{2}%
-xy\right)  $ ? What lies in the ideal \newline$\left(  x^{2}-yz\right)
R\left[  x,y,z\right]  +\left(  y^{2}-zx\right)  R\left[  x,y,z\right]
+\left(  z^{2}-xy\right)  R\left[  x,y,z\right]  $ ?

\item What is $R\left[  x,y,z\right]  /\left(  x^{2}+xy,y^{2}+yz,z^{2}%
+zx\right)  $ ? What lies in the ideal \newline$\left(  x^{2}+xy\right)
R\left[  x,y,z\right]  +\left(  y^{2}+yz\right)  R\left[  x,y,z\right]
+\left(  z^{2}+zx\right)  R\left[  x,y,z\right]  $ ? For example, I claim that
$z^{4}$ lies in this ideal when $R=\mathbb{Q}$, but $z^{3}$ does not. How do I
know? How can you tell?
\end{itemize}

In theory, you could imagine that there are ideals that do not even have a
finite list of elements generating them. There are rings that have such
ideals. For example, the polynomial ring $\mathbb{Z}\left[  x_{1},x_{2}%
,x_{3},\ldots\right]  $ in infinitely many variables has such ideals (see
Exercise \ref{exe.polring.inf-many-vars-inf-ideal} below). But polynomial
rings in finitely many variables over a field are not this bad. Indeed:

\begin{theorem}
[Hilbert's basis theorem]\label{thm.polring.hilbert}Let $F$ be a field. Let
$S$ be the polynomial ring $F\left[  x_{1},x_{2},\ldots,x_{n}\right]  $ for
some $n\in\mathbb{N}$. Then, any ideal $I$ of $S$ is finitely generated (this
means that there is a finite list $\left(  a_{1},a_{2},\ldots,a_{k}\right)  $
of elements of $I$ such that $I=a_{1}S+a_{2}S+\cdots+a_{k}S$).
\end{theorem}

\begin{proof}
See \cite[\S 9.6, Corollary 22]{DumFoo04} or \cite[Corollary 5.4.8]{Laurit09}
or (for a more general result) \cite[Theorem 36.12]{Swanso17}.
\end{proof}

\begin{warning}
If $n=1$, then the ideal $I$ in Theorem \ref{thm.polring.hilbert} is principal
(since $F\left[  x_{1}\right]  $ is a PID), so you can get by with a
length-$1$ list (i.e., with $k=1$). However, if $n=2$, then the list can be
arbitrarily large. You cannot always find a length-$2$ list. For example, let
$d\in\mathbb{N}$ be arbitrary; then, in the polynomial ring $F\left[
x,y\right]  $, the ideal generated by all monomials of degree $d$ (that is, by
$x^{d},x^{d-1}y,x^{d-2}y^{2},\ldots,y^{d}$) cannot be generated by $d$ or
fewer elements.
\end{warning}

\begin{exercise}
\label{exe.polring.inf-many-vars-inf-ideal}Let $S_{\infty}$ be the polynomial
ring $\mathbb{Z}\left[  x_{1},x_{2},x_{3},\ldots\right]  $ in infinitely many
variables. Strictly speaking, we have never defined this ring, but you can
easily produce its definition: It still is a monoid ring, but each monomial
now has the form $x_{1}^{a_{1}}x_{2}^{a_{2}}x_{3}^{a_{3}}\cdots$ for some
infinite sequence $\left(  a_{1},a_{2},a_{3},\ldots\right)  $ of nonnegative
integers with the property that only finitely many of the exponents $a_{k}$
are nonzero. (Thus, there are infinitely many indeterminates, but each single
monomial can only use finitely many of them. For instance, infinite monomials
like $x_{1}x_{2}x_{3}\cdots$ are not allowed. As a consequence, a polynomial
in $S_{\infty}$ must also use only a finite set of indeterminates.)

Let $J$ be the set of all polynomials in $S_{\infty}$ whose constant term
(i.e., coefficient of the monomial $x_{1}^{0}x_{2}^{0}x_{3}^{0}\cdots$) is $0$.

\begin{enumerate}
\item[\textbf{(a)}] Show that $J$ is an ideal of $S_{\infty}$.

\item[\textbf{(b)}] Show that $J$ is not an ideal generated by any finite list
of elements of $S_{\infty}$.
\end{enumerate}
\end{exercise}

\subsection{\label{sec.polys2.deg-lex}Degrees and the deg-lex order}

Let us now attempt a more general approach.

\begin{convention}
From now on, for the rest of this chapter, we fix a commutative ring $R$ and
an $n\in\mathbb{N}$.

We let $P$ denote the polynomial ring $R\left[  x_{1},x_{2},\ldots
,x_{n}\right]  $.
\end{convention}

As we recall, a \textbf{monomial} is an element of the free abelian monoid
$C^{\left(  n\right)  }$ with $n$ generators $x_{1},x_{2},\ldots,x_{n}$; it
has the form $x_{1}^{a_{1}}x_{2}^{a_{2}}\cdots x_{n}^{a_{n}}$ for some
$\left(  a_{1},a_{2},\ldots,a_{n}\right)  \in\mathbb{N}^{n}$.

\subsubsection{Degrees}

Our first goal is to define the degree of a polynomial in $n$ variables. We
begin by defining the degree of a monomial:

\begin{definition}
\label{def.polring.monomial-deg}The \textbf{degree }of a monomial
$\mathfrak{m}=x_{1}^{a_{1}}x_{2}^{a_{2}}\cdots x_{n}^{a_{n}}\in C^{\left(
n\right)  }$ is defined to be the number $a_{1}+a_{2}+\cdots+a_{n}%
\in\mathbb{N}$. It is denoted by $\deg\mathfrak{m}$.
\end{definition}

For example, the monomial $x_{1}^{5}x_{2}x_{4}^{2}=x_{1}^{5}x_{2}^{1}x_{3}%
^{0}x_{4}^{2}$ has degree $5+1+0+2=8$.

\begin{definition}
\label{def.polring.monomial-appear}A monomial $\mathfrak{m}$ is said to
\textbf{appear} in a polynomial $f\in P$ if $\left[  \mathfrak{m}\right]
f\neq0$. (Recall that $\left[  \mathfrak{m}\right]  f$ means the coefficient
of $\mathfrak{m}$ in $f$.)
\end{definition}

For example, the monomial $x^{2}y$ appears in $\left(  x+y\right)  ^{3}\in
R\left[  x,y\right]  $ (if $3\neq0$ in $R$), but the monomial $xy$ does not.

\begin{definition}
\label{def.polring.mulvar-deg}The \textbf{degree} (or \textbf{total degree})
of a nonzero polynomial $f\in P$ is the largest degree of a monomial that
appears in $f$. We also define the degree of the zero polynomial $0\in P$ to
be $-\infty$.
\end{definition}

For example:

\begin{itemize}
\item The polynomial $\left(  x+y+1\right)  ^{3}\in\mathbb{Q}\left[
x,y\right]  $ has degree $3$.

\item The polynomial $\left(  x+y+1\right)  ^{3}-\left(  x+y\right)  ^{3}%
\in\mathbb{Q}\left[  x,y\right]  $ has degree $2$, since it equals
$3x^{2}+3y^{2}+6xy+3x+3y+1$.

\item The polynomial $\left(  x+y+\overline{1}\right)  ^{3}-\left(
x+y\right)  ^{3}\in\left(  \mathbb{Z}/3\right)  \left[  x,y\right]  $ has
degree $0$, since it equals $\overline{1}$.
\end{itemize}

Definition \ref{def.polring.mulvar-deg} generalizes our old definition of
degree for nonzero univariate polynomials.

The following proposition generalizes a fact that we previously proved for
univariate polynomials (parts \textbf{(a)} and \textbf{(c)} of Proposition
\ref{prop.polring.univar-degpq}):

\begin{proposition}
[Degree-of-a-product formula]\label{prop.polring.deg-pq-mulvar}Let $R$ be a
commutative ring. Let $p,q\in P$ be nonzero.

\begin{enumerate}
\item[\textbf{(a)}] We have $\deg\left(  pq\right)  \leq\deg p+\deg q$.

\item[\textbf{(b)}] We have $\deg\left(  pq\right)  =\deg p+\deg q$ if $R$ is
an integral domain.
\end{enumerate}
\end{proposition}

Part \textbf{(a)} of this proposition is pretty clear. (The reason is that
$\deg\left(  \mathfrak{mn}\right)  =\deg\mathfrak{m}+\deg\mathfrak{n}$ for any
monomials $\mathfrak{m},\mathfrak{n}$.)

What about part \textbf{(b)}? We proved this for univariate polynomials using
leading coefficients. What is a leading coefficient when several monomials can
have the same degree? In order to define it, we need to break ties (i.e.,
establish an ordering on monomials of equal degrees) in a way that will be
compatible with products\footnote{I will explain what this means later.}. To
that aim, we shall introduce a total order on the set $C^{\left(  n\right)  }$
of all monomials.

\subsubsection{The deg-lex order}

Recall that a \textbf{total order} (or, to be more precise, a \textbf{strict
total order}) on a set $S$ is a binary relation $\prec$ on $S$ that is

\begin{itemize}
\item \textbf{asymmetric} (meaning that no two elements $a$ and $b$ of $S$
satisfy $a\prec b$ and $b\prec a$ at the same time);

\item \textbf{transitive} (meaning that if $a,b,c\in S$ satisfy $a\prec b$ and
$b\prec c$, then $a\prec c$);

\item \textbf{trichotomous} (meaning that for any two elements $a$ and $b$ of
$S$, we have $a\prec b$ or $a=b$ or $b\prec a$).
\end{itemize}

Here are some examples of total orders:

\begin{itemize}
\item The relation $<$ on the set $\mathbb{N}$ or on the set $\mathbb{Z}$ or
on the set $\mathbb{R}$ is a total order.

\item So is the relation $>$ on each of these three sets.

\item If $S$ is a finite set, and if $\left(  s_{1},s_{2},\ldots,s_{k}\right)
$ is a list of all elements of $S$, with each element of $S$ appearing exactly
once in this list, then we can define a total order $\prec$ on $S$ as follows:
We declare that two elements $u,v\in S$ satisfy $u\prec v$ if and only if $u$
appears prior to $v$ in this list $\left(  s_{1},s_{2},\ldots,s_{k}\right)  $
(that is, if $u=s_{i}$ and $v=s_{j}$ for some $i<j$).

\item On the other hand, the relation $\subsetneq$ (\textquotedblleft proper
subset of\textquotedblright) on the power set $\mathcal{P}\left(  X\right)  $
of a set $X$ is not a total order unless $\left\vert X\right\vert \leq1$.
(Indeed, it is asymmetric and transitive, but not trichotomous, because if
$\alpha$ and $\beta$ are two distinct elements of $X$, then we have neither
$\left\{  \alpha\right\}  \subsetneq\left\{  \beta\right\}  $ nor $\left\{
\alpha\right\}  =\left\{  \beta\right\}  $ nor $\left\{  \beta\right\}
\subsetneq\left\{  \alpha\right\}  $.)
\end{itemize}

If $\prec$ is a total order on a set $S$, then we view relations of the form
$a\prec b$ as saying that $a$ is in some sense smaller than $b$. We will use
the notations $\preccurlyeq$, $\succ$ and $\succcurlyeq$ accordingly; this
means that

\begin{itemize}
\item we write \textquotedblleft$a\preccurlyeq b$\textquotedblright\ for
\textquotedblleft$a\prec b$ or $a=b$\textquotedblright.

\item we write \textquotedblleft$a\succ b$\textquotedblright\ for
\textquotedblleft$b\prec a$\textquotedblright.

\item we write \textquotedblleft$a\succcurlyeq b$\textquotedblright\ for
\textquotedblleft$a\succ b$ or $a=b$\textquotedblright.
\end{itemize}

So we all know a total order on the set $\mathbb{R}$ of all real numbers. But
what about other sets? For example, how can we find a total order on the set
of words in the English language? A long time ago, creators of dictionaries
and encyclopedias were faced with this very problem, because it would be hard
to look a word up in a dictionary if there was no well-known total order in
which the words appeared in the dictionary. The total order commonly used in
dictionaries is known as the \textbf{lexicographic order} (or
\textbf{dictionary order}): Words are ordered by their first letter (e.g.,
\textquotedblleft ant\textquotedblright\ $\prec$ \textquotedblleft
bear\textquotedblright); ties are broken using the second letter
(\textquotedblleft ant\textquotedblright\ $\prec$ \textquotedblleft
armadillo\textquotedblright); remaining ties are broken using the third letter
(\textquotedblleft camel\textquotedblright\ $\prec$ \textquotedblleft
cat\textquotedblright); and so on; absent letters at the end are treated as
being smaller than present letters (e.g., \textquotedblleft
ant\textquotedblright\ $\prec$ \textquotedblleft anteater\textquotedblright).
We use this as an inspiration for defining a total order on $C^{\left(
n\right)  }$, but we shall use the degree as the first level of comparison.

\begin{definition}
\label{def.polring.deglex}We define a total order $\prec$ (called the
\textbf{degree-lexicographic order}, or -- for short -- the \textbf{deg-lex
order}) on the set $C^{\left(  n\right)  }$ of all monomials as follows:

For two monomials $\mathfrak{m}=x_{1}^{a_{1}}x_{2}^{a_{2}}\cdots x_{n}^{a_{n}%
}$ and $\mathfrak{n}=x_{1}^{b_{1}}x_{2}^{b_{2}}\cdots x_{n}^{b_{n}}$, we
declare that $\mathfrak{m}\prec\mathfrak{n}$ if and only if

\begin{itemize}
\item \textbf{either} $\deg\mathfrak{m}<\deg\mathfrak{n}$;

\item \textbf{or} $\deg\mathfrak{m}=\deg\mathfrak{n}$ and the following holds:
There is an $i\in\left\{  1,2,\ldots,n\right\}  $ such that $a_{i}\neq b_{i}$,
and the \textbf{smallest} such $i$ satisfies $a_{i}<b_{i}$.
\end{itemize}
\end{definition}

In words:

\begin{itemize}
\item If two monomials have different degrees, then we declare the monomial
with smaller degree to be the smaller one.

\item If they have equal degrees, then we look at the first variable that has
different exponents in the two monomials, and we declare the monomial with the
smaller exponent on this variable to be smaller.
\end{itemize}

For example:

\begin{itemize}
\item We have $x_{1}^{2}\prec x_{2}x_{3}^{2}$, since $\deg\left(  x_{1}%
^{2}\right)  =2<3=\deg\left(  x_{2}x_{3}^{2}\right)  $.

\item We have $x_{1}^{5}x_{2}x_{3}x_{4}^{2}\prec x_{1}^{5}x_{2}x_{3}^{2}x_{4}%
$, since the two monomials have the same degree, and the first variable that
has different exponents in these two monomials is $x_{3}$, and this variable
appears with a smaller exponent in $x_{1}^{5}x_{2}x_{3}x_{4}^{2}$ (namely,
with exponent $1$) than in $x_{1}^{5}x_{2}x_{3}^{2}x_{4}$ (namely, with
exponent $2$).

\item We have $x_{1}x_{3}^{2}\prec x_{1}x_{2}x_{3}$, since the first variable
that has different exponents in these two monomials is $x_{2}$, and this
variable appears with a smaller exponent in $x_{1}x_{3}^{2}$ (namely, with
exponent $0$) than in $x_{1}x_{2}x_{3}$ (namely, with exponent $1$).

\item The reader may easily check that $x_{3}^{3}\prec x_{1}x_{2}x_{3}%
^{2}\prec x_{1}x_{2}^{2}x_{3}\prec x_{1}^{2}x_{2}x_{3}\prec x_{3}^{5}\prec
x_{1}^{2}x_{2}^{2}x_{3}^{2}\prec x_{1}^{6}$.
\end{itemize}

You can intuitively think of the deg-lex order as follows: A monomial becomes
larger (in this order) if you increase its degree, and also becomes larger if
you replace an $x_{i}$ factor by an $x_{j}$ factor with $j<i$.

The deg-lex order has several good properties:

\begin{proposition}
\label{prop.polring.deglex.basics}\ \ 

\begin{enumerate}
\item[\textbf{(a)}] The deg-lex order really is a total order on $C^{\left(
n\right)  }$.

\item[\textbf{(b)}] If $\mathfrak{m},\mathfrak{n},\mathfrak{p}\in C^{\left(
n\right)  }$ satisfy $\mathfrak{m}\prec\mathfrak{n}$, then $\mathfrak{mp}%
\prec\mathfrak{np}$.

\item[\textbf{(c)}] We have $1\preccurlyeq\mathfrak{m}$ for any $\mathfrak{m}%
\in C^{\left(  n\right)  }$.

\item[\textbf{(d)}] Let $\mathfrak{m}\in C^{\left(  n\right)  }$ be any
monomial. Then, there are only finitely many monomials $\mathfrak{p}$ such
that $\mathfrak{p}\prec\mathfrak{m}$.

\item[\textbf{(e)}] There are no infinite decreasing chains $\mathfrak{m}%
_{0}\succ\mathfrak{m}_{1}\succ\mathfrak{m}_{2}\succ\cdots$ of monomials.

\item[\textbf{(f)}] If $T$ is a nonempty finite set of monomials, then $T$ has
a largest element with respect to $\prec$ (that is, an element $\mathfrak{t}%
\in T$ such that $\mathfrak{m}\preccurlyeq\mathfrak{t}$ for all $\mathfrak{m}%
\in T$).

\item[\textbf{(g)}] If $T$ is a nonempty set of monomials, then $T$ has a
smallest element with respect to $\prec$ (that is, an element $\mathfrak{t}\in
T$ such that $\mathfrak{m}\succcurlyeq\mathfrak{t}$ for all $\mathfrak{m}\in
T$).
\end{enumerate}
\end{proposition}

Note that we require $T$ to be finite in Proposition
\ref{prop.polring.deglex.basics} \textbf{(f)} but not in Proposition
\ref{prop.polring.deglex.basics} \textbf{(g)}. This is similar to the
situation for sets of nonnegative integers (viz., any nonempty set of
nonnegative integers has a smallest element, but only finite nonempty sets of
nonnegative integers have largest elements).

\begin{proof}
[Hints to the proof of Proposition \ref{prop.polring.deglex.basics}%
.]\textbf{(a)}, \textbf{(b)}, \textbf{(c)}, \textbf{(d)} LTTR.

\textbf{(e)} This follows from \textbf{(d)}.

\textbf{(f)} This holds for any total order on any set.

\textbf{(g)} This is easily proved using \textbf{(d)} (or, less easily, using
\textbf{(e)}). LTTR.
\end{proof}

(Proposition \ref{prop.polring.deglex.basics} \textbf{(b)} is what I meant
when I said that the deg-lex order is \textquotedblleft compatible with
products\textquotedblright.)

\subsubsection{Leading coefficients, monomials and terms}

Now, we can define leading coefficients of multivariate polynomials:

\begin{definition}
\label{def.polring.leading}Let $f\in P$ be a nonzero polynomial.

\begin{enumerate}
\item[\textbf{(a)}] The \textbf{leading monomial} of $f$ means the largest
(with respect to $\prec$) monomial that appears in $f$. It is denoted by
$\operatorname*{LM}f$.

\item[\textbf{(b)}] The \textbf{leading coefficient} of $f$ means the
coefficient $\left[  \operatorname*{LM}f\right]  f$. It is denoted by
$\operatorname*{LC}f$.

\item[\textbf{(c)}] The \textbf{leading term} of $f$ means the product
$\operatorname*{LC}f\cdot\operatorname*{LM}f$. It is denoted by
$\operatorname*{LT}f$.
\end{enumerate}
\end{definition}

For example, if $3\neq0$ in $R$, then%
\begin{align*}
\operatorname*{LM}\left(  \left(  x_{1}+x_{2}+1\right)  ^{3}-x_{1}^{3}\right)
&  =x_{1}^{2}x_{2};\\
\operatorname*{LC}\left(  \left(  x_{1}+x_{2}+1\right)  ^{3}-x_{1}^{3}\right)
&  =3;\\
\operatorname*{LT}\left(  \left(  x_{1}+x_{2}+1\right)  ^{3}-x_{1}^{3}\right)
&  =3x_{1}^{2}x_{2}.
\end{align*}

Two simple consequences of this definition are:

\begin{proposition}
\label{prop.polring.leading.subtract}Let $f\in P$ be a nonzero polynomial.
Then, $f-\operatorname*{LT}f=0$ or else $\operatorname*{LM}\left(
f-\operatorname*{LT}f\right)  \prec\operatorname*{LM}f$.
\end{proposition}

\begin{proof}
By Definition \ref{def.polring.leading}, we have%
\[
f=\operatorname*{LT}f+\left(  \text{an }R\text{-linear combination of
monomials }\mathfrak{m}\text{ with }\mathfrak{m}\prec\operatorname*{LM}%
f\right)  .
\]
Hence, $f-\operatorname*{LT}f$ is an $R$-linear combination of monomials
$\mathfrak{m}$ with $\mathfrak{m}\prec\operatorname*{LM}f$. Therefore,
$f-\operatorname*{LT}f=0$ or else $\operatorname*{LM}\left(
f-\operatorname*{LT}f\right)  \prec\operatorname*{LM}f$.
\end{proof}

\begin{proposition}
\label{prop.polring.leading.prod}Let $f,g\in P$ be nonzero polynomials such
that $\operatorname*{LC}f$ is regular in $R$ (see Definition
\ref{def.rings.reg.reg}). Then,%
\[
\operatorname*{LM}\left(  fg\right)  =\operatorname*{LM}f\cdot
\operatorname*{LM}g\ \ \ \ \ \ \ \ \ \ \text{and}%
\ \ \ \ \ \ \ \ \ \ \operatorname*{LC}\left(  fg\right)  =\operatorname*{LC}%
f\cdot\operatorname*{LC}g.
\]

\end{proposition}

\begin{proof}
LTTR. (Use Proposition \ref{prop.polring.deglex.basics} \textbf{(b)}.)
\end{proof}

Now we can easily prove Proposition \ref{prop.polring.deg-pq-mulvar}
\textbf{(b)}. (The details are LTTR.)

From Proposition \ref{prop.polring.deg-pq-mulvar} \textbf{(b)}, we obtain the following:

\begin{corollary}
If $R$ is an integral domain, then the polynomial ring $P=R\left[  x_{1}%
,x_{2},\ldots,x_{n}\right]  $ is an integral domain.
\end{corollary}

(Alternatively, this can also be proved by induction on $n$, using Proposition
\ref{prop.polring.Rxn-stepwise}.)

\subsection{\label{sec.polys2.groebner}Division with remainder and Gr\"{o}bner
bases}

By defining leading monomials and leading coefficients, we have recovered one
piece of the nice theory of univariate polynomials in the multivariate case.
Can we do more? Can we define division with remainder?

\subsubsection{The case of principal ideals}

We \textbf{can} divide with remainder by a single polynomial\footnote{Recall
that $P=R\left[  x_{1},x_{2},\ldots,x_{n}\right]  $.}:

\begin{theorem}
[Division-with-remainder theorem for multivariate polynomials]%
\label{thm.polring.quorem-mulvar}Let $b\in P$ be a nonzero polynomial whose
leading coefficient $\operatorname*{LC}b$ is a unit of $R$. Let $a\in P$ be
any polynomial.

Then, there is a \textbf{unique} pair $\left(  q,r\right)  $ of polynomials in
$P$ such that%
\[
a=qb+r\ \ \ \ \ \ \ \ \ \ \text{and}\ \ \ \ \ \ \ \ \ \ r\text{ is
}\operatorname*{LM}b\text{-reduced.}%
\]
Here, a polynomial $r\in P$ is said to be $\mathfrak{m}$\textbf{-reduced}
(where $\mathfrak{m}$ is a monomial) if no monomial divisible by
$\mathfrak{m}$ appears in $r$.
\end{theorem}

This generalizes the division-with-remainder theorem for univariate
polynomials (Theorem \ref{thm.polring.univar-quorem} \textbf{(a)}); indeed, if
$n=1$, then the condition \textquotedblleft$r$ is $\operatorname*{LM}%
b$-reduced\textquotedblright\ is equivalent to \textquotedblleft$\deg r<\deg
b$\textquotedblright\ (which is familiar from the case of univariate
polynomials). The entries $q$ and $r$ of the pair $\left(  q,r\right)  $ in
Theorem \ref{thm.polring.quorem-mulvar} will be called the \textbf{quotient}
and the \textbf{remainder} of the division of $a$ by $b$.

Let us illustrate Theorem \ref{thm.polring.quorem-mulvar} on an example:

\begin{itemize}
\item Let $n=2$ and $R=\mathbb{Z}$, and let us rename the indeterminates
$x_{1},x_{2}$ as $x,y$. Thus, $P=\mathbb{Z}\left[  x,y\right]  $. Let
$b=xy\left(  x-y\right)  \in P$. Thus, $\operatorname*{LM}b=x^{2}y$ and
$\operatorname*{LC}b=1$.

Let $a=\left(  x+y\right)  ^{4}$. We want to divide $a$ by $b$ with remainder.
That is, we want to find the pair $\left(  q,r\right)  $ in Theorem
\ref{thm.polring.quorem-mulvar}.

Theorem \ref{thm.polring.quorem-mulvar} says that $a$ can be written as a
multiple of $b$ plus some $\operatorname*{LM}b$-reduced polynomial. In other
words, it says that by subtracting an appropriate multiple of $b$ from $a$, we
can obtain an $\operatorname*{LM}b$-reduced polynomial. How do we find the
right multiple to subtract?

In the univariate case, \textquotedblleft$\operatorname*{LM}b$%
-reduced\textquotedblright\ was simply saying that $\deg r<\deg b$, and we
achieved this by repeatedly subtracting multiples of $b$ from $a$ in order to
chip away at the leading term (reducing the degree by at least $1$ in each
step). We can do this similarly in the multivariate case: We simply check
whether $a$ is already $\operatorname*{LM}b$-reduced. As long as it isn't, we
find some monomial divisible by $\operatorname*{LM}b$ that appears in $a$, and
we clear it out by subtracting an appropriate multiple of $b$ (so that this
monomial no longer appears in $a$). More precisely, we clear out the highest
such monomial that appears in $a$. We keep doing this until no such monomials
remain (which means that $a$ has become $\operatorname*{LM}b$-reduced).

Let us actually do this in our above example: We start with%
\[
a=\left(  x+y\right)  ^{4}=x^{4}+4x^{3}y+6x^{2}y^{2}+4xy^{3}+y^{4}.
\]
Two monomials that are multiples of $\operatorname*{LM}b=x^{2}y$ appear on the
right hand side: $x^{3}y$ and $x^{2}y^{2}$. The highest of them is $x^{3}y$,
so we clear it out by subtracting an appropriate multiple of $b$. This
appropriate multiple is $4xb$, since we want to clear out a $4x^{3}y$ term. So
we get%
\begin{align*}
a-4xb  &  =\left(  x^{4}+4x^{3}y+6x^{2}y^{2}+4xy^{3}+y^{4}\right)  -4x\cdot
xy\left(  x-y\right) \\
&  =x^{4}+10x^{2}y^{2}+4xy^{3}+y^{4}.
\end{align*}

Now we still have one monomial left that is a multiple of $\operatorname*{LM}%
b=x^{2}y$, namely $x^{2}y^{2}$. We clear it out by subtracting $10yb$, and we
end up with%
\begin{align*}
a-4xb-10yb  &  =\left(  x^{4}+10x^{2}y^{2}+4xy^{3}+y^{4}\right)  -10y\cdot
xy\left(  x-y\right) \\
&  =x^{4}+14xy^{3}+y^{4}.
\end{align*}
The right hand side of this equality is $\operatorname*{LM}b$-reduced, so it
is the remainder we were looking for. That is, the $r$ in our pair $\left(
q,r\right)  $ is $x^{4}+14xy^{3}+y^{4}$. The $q$ in this pair we find by
collecting the multiples of $b$ that we have subtracted; thus, we get
$q=4x+10y$ (since we have subtracted $4xb$ and $10yb$). Hence, our pair
$\left(  q,r\right)  $ is%
\[
\left(  q,r\right)  =\left(  4x+10y,\ x^{4}+14xy^{3}+y^{4}\right)  .
\]

This example was somewhat simplistic. In more complicated cases, it can happen
that subtracting a multiple of $b$ will create new monomials that are not
$\operatorname*{LM}b$-reduced. However, if we keep following our method, all
those new monomials will eventually get removed as well.
\end{itemize}

\begin{proof}
[Hints to the proof of Theorem \ref{thm.polring.quorem-mulvar}.]The existence
of the pair $\left(  q,r\right)  $ is proved by the same idea as in the
example we just did. All we need to do is to explain why our procedure
terminates (i.e., doesn't keep running forever). This is not hard: We observe
that, as we keep subtracting appropriate multiples of $b$ from $a$, the
\textbf{highest} monomial that is a multiple of $\operatorname*{LM}b$ and
appears in $a$ becomes smaller and smaller (because each subtraction clears
out the highest such monomial, and can only introduce lower such monomials).
Thus, if our procedure were to run forever, then we would obtain an infinite
decreasing chain $\mathfrak{m}_{0}\succ\mathfrak{m}_{1}\succ\mathfrak{m}%
_{2}\succ\cdots$ of monomials; but this would contradict Proposition
\ref{prop.polring.deglex.basics} \textbf{(e)}. Thus, the algorithm eventually
terminates, and this proves the existence of $\left(  q,r\right)  $.

To prove the uniqueness of $\left(  q,r\right)  $, it suffices to show that no
nonzero multiple of $b$ is $\operatorname*{LM}b$-reduced\footnote{Indeed, if
$\left(  q_{1},r_{1}\right)  $ and $\left(  q_{2},r_{2}\right)  $ are two
pairs $\left(  q,r\right)  $ satisfying the claim of Theorem
\ref{thm.polring.quorem-mulvar}, then $r_{1}-r_{2}=\left(  q_{2}-q_{1}\right)
b$ is a multiple of $b$ that is $\operatorname*{LM}b$-reduced (since $r_{1}$
and $r_{2}$ are $\operatorname*{LM}b$-reduced).}. But this follows easily from
Proposition \ref{prop.polring.leading.prod}.
\end{proof}

As a consequence of Theorem \ref{thm.polring.quorem-mulvar} (or, more
precisely, of the algorithm for the construction of $\left(  q,r\right)  $
that we demonstrated in the above example), we obtain an algorithmic way to
tell whether a polynomial $a\in P$ is divisible by $b$ or not (whenever $b\in
P$ is a nonzero polynomial whose leading coefficient $\operatorname*{LC}b$ is
a unit of $R$). Namely, we compute the pair $\left(  q,r\right)  $ from
Theorem \ref{thm.polring.quorem-mulvar}, and check whether $r=0$. The
uniqueness of this pair easily yields that $b\mid a$ if and only if $r=0$.

Another consequence of Theorem \ref{thm.polring.quorem-mulvar} is the
following corollary that explicitly constructs a basis of the $R$-module $P/b$:

\begin{corollary}
\label{cor.polring.mod-mulvar}Let $b\in P$ be a nonzero polynomial whose
leading coefficient $\operatorname*{LC}b$ is a unit of $R$. Then, each element
of $P/b$ can be uniquely written in the form%
\[
\sum_{\substack{\mathfrak{m}\text{ is a monomial}\\\text{not divisible by
}\operatorname*{LM}b}}a_{\mathfrak{m}}\overline{\mathfrak{m}}%
\ \ \ \ \ \ \ \ \ \ \text{with }a_{\mathfrak{m}}\in R
\]
(where all but finitely many $\mathfrak{m}$ satisfy $a_{\mathfrak{m}}=0$).
Equivalently, the family $\left(  \overline{\mathfrak{m}}\right)
_{\mathfrak{m}\text{ is a monomial not divisible by }\operatorname*{LM}b}$ is
a basis of the $R$-module $P/b$. If $b$ is not constant, then the ring $P/b$
contains \textquotedblleft a copy of $R$\textquotedblright.
\end{corollary}

Corollary \ref{cor.polring.mod-mulvar} generalizes Theorem
\ref{thm.fieldext.basis-monic.repeat} (and is proved in the same way, except
that we use Theorem \ref{thm.polring.quorem-mulvar} instead of the univariate
division-with-remainder theorem). Here are some examples:

\begin{itemize}
\item Let us take $P=R\left[  x,y\right]  $ and $b=y$ in Corollary
\ref{cor.polring.mod-mulvar}. Then, $\operatorname*{LM}b=y$, so that Corollary
\ref{cor.polring.mod-mulvar} yields that the family $\left(  \overline
{\mathfrak{m}}\right)  _{\mathfrak{m}\text{ is a monomial not divisible by }%
y}$ is a basis of the $R$-module $P/b=R\left[  x,y\right]  /y$. Since the
monomials not divisible by $y$ are precisely the powers of $x$ (that is,
$x^{0},x^{1},x^{2},\ldots$), we can rewrite this as follows: The family
$\left(  \overline{x^{i}}\right)  _{i\in\mathbb{N}}=\left(  \overline{x^{0}%
},\overline{x^{1}},\overline{x^{2}},\ldots\right)  $ is a basis of the
$R$-module $P/b=R\left[  x,y\right]  /y$. This is in line with Proposition
\ref{prop.polring.Rxy/y} (indeed, the isomorphism $R\left[  x,y\right]
/y\rightarrow R\left[  x\right]  $ sends this family to the standard basis
$\left(  x^{0},x^{1},x^{2},\ldots\right)  $ of $R\left[  x\right]  $).

\item Let us take $P=R\left[  x,y\right]  $ and $b=x^{2}+y^{2}-1$ in Corollary
\ref{cor.polring.mod-mulvar}. Then, $\operatorname*{LM}b=x^{2}$, so that
Corollary \ref{cor.polring.mod-mulvar} yields that the family \newline$\left(
\overline{\mathfrak{m}}\right)  _{\mathfrak{m}\text{ is a monomial not
divisible by }x^{2}}$ is a basis of the $R$-module \newline$P/b=R\left[
x,y\right]  /\left(  x^{2}+y^{2}-1\right)  $. Since the monomials not
divisible by $x^{2}$ are precisely the monomials $x^{i}y^{j}$ with $i<2$, we
can rewrite this as follows: The family
\[
\left(  \overline{x^{i}y^{j}}\right)  _{\left(  i,j\right)  \in\mathbb{N}%
^{2};\ i<2}=\left(  \overline{x^{0}y^{0}},\overline{x^{0}y^{1}},\overline
{x^{0}y^{2}},\ldots,\overline{x^{1}y^{0}},\overline{x^{1}y^{1}},\overline
{x^{1}y^{2}},\ldots\right)
\]
is a basis of the $R$-module $P/b=R\left[  x,y\right]  /\left(  x^{2}%
+y^{2}-1\right)  $. This is not the basis that we obtained back in
(\ref{eq.polring.Rxy/xx+yy-1.6}), but rather is obtained from the latter by
interchanging $x$ and $y$. Of course, it is no surprise that interchanging $x$
and $y$ turns a basis into a basis; indeed, the variables $x$ and $y$ clearly
play symmetric roles in $R\left[  x,y\right]  /\left(  x^{2}+y^{2}-1\right)
$, so every basis that treats them unequally has a \textquotedblleft
mirror\textquotedblright\ version with $x$ and $y$ interchanged.

\item Let us take $P=R\left[  x,y\right]  $ and $b=xy$ in Corollary
\ref{cor.polring.mod-mulvar}. Then, $\operatorname*{LM}b=xy$, so that
Corollary \ref{cor.polring.mod-mulvar} yields that the family $\left(
\overline{\mathfrak{m}}\right)  _{\mathfrak{m}\text{ is a monomial not
divisible by }xy}$ is a basis of the $R$-module $P/b=R\left[  x,y\right]
/\left(  xy\right)  $. Since the monomials not divisible by $xy$ are precisely
the monomials $1,x^{1},x^{2},x^{3},\ldots,y^{1},y^{2},y^{3},\ldots$ (that is,
the monomials that are powers of a single indeterminate), we can rewrite this
as follows: The family
\[
\left(  \overline{1},\overline{x^{1}},\overline{x^{2}},\overline{x^{3}}%
,\ldots,\overline{y^{1}},\overline{y^{2}},\overline{y^{3}},\ldots\right)
\]
is a basis of the $R$-module $P/b=R\left[  x,y\right]  /\left(  xy\right)  $.
This can be obtained in more direct ways, too.

\item Likewise, applying Corollary \ref{cor.polring.mod-mulvar} to $P=R\left[
x,y\right]  $ and $b=xy\left(  x-y\right)  $ yields that the family%
\begin{align*}
&  \left(  \overline{\mathfrak{m}}\right)  _{\mathfrak{m}\text{ is a monomial
not divisible by }x^{2}y}\\
&  =\left(  \overline{1},\overline{x^{1}},\overline{x^{2}},\overline{x^{3}%
},\ldots,\overline{y^{1}},\overline{y^{2}},\overline{y^{3}},\ldots
,\overline{xy^{1}},\overline{xy^{2}},\overline{xy^{3}},\ldots\right)
\end{align*}
is a basis of the $R$-module $R\left[  x,y\right]  /\left(  xy\left(
x-y\right)  \right)  $.
\end{itemize}

\begin{exercise}
Consider the setting of Theorem \ref{thm.polring.quorem-mulvar}. Prove that
the remainder of the division of $a$ by $b$ is the unique $\operatorname*{LM}%
b$-reduced polynomial $p\in P$ that satisfies $a-p\in bP$.
\end{exercise}

\begin{exercise}
Let $R=\mathbb{Z}$, and let us rename the variables $x_{1},x_{2},x_{3},x_{4}$
as $x,y,z,w$. Let $n$ be a positive integer.

\begin{enumerate}
\item[\textbf{(a)}] Find the remainder of the division of $\left(  x+y\right)
^{n}$ by $xy\left(  x-y\right)  $.

\item[\textbf{(b)}] Find the remainder of the division of $x^{n}$ by $\left(
x-y\right)  ^{2}$.

\item[\textbf{(c)}] Find the remainder of the division of $\left(  xz\right)
^{n}$ by $\left(  x-y\right)  \left(  z-w\right)  $.
\end{enumerate}
\end{exercise}

\begin{exercise}
\label{exe.polring.mulvar.x+xy}Let $R$ be any commutative ring. Let $S$ be the
ring $R\left[  x,y\right]  /\left(  xy^{2}\right)  $.

\begin{enumerate}
\item[\textbf{(a)}] Prove that the family%
\begin{align*}
&  \left(  \overline{\mathfrak{m}}\right)  _{\mathfrak{m}\text{ is a monomial
not divisible by }xy^{2}}\\
&  =\left(  \overline{1},\overline{x},\overline{x^{2}},\overline{x^{3}}%
,\ldots,\overline{y},\overline{xy},\overline{x^{2}y},\overline{x^{3}y}%
,\ldots,\overline{y^{2}},\overline{y^{3}},\overline{y^{4}},\ldots\right)
\end{align*}
is a basis of the $R$-module $S$.

\item[\textbf{(b)}] Prove that the maps%
\begin{align*}
S  &  \rightarrow R\left[  x\right]  ,\\
\overline{f}  &  \mapsto f\left[  x,0\right]  \ \ \ \ \ \ \ \ \ \ \left(
\text{\textquotedblleft substituting }0\text{ for }y\text{\textquotedblright%
}\right)
\end{align*}
and%
\begin{align*}
S  &  \rightarrow R\left[  y\right]  ,\\
\overline{f}  &  \mapsto f\left[  0,y\right]  \ \ \ \ \ \ \ \ \ \ \left(
\text{\textquotedblleft substituting }0\text{ for }x\text{\textquotedblright%
}\right)
\end{align*}
are $R$-algebra morphisms.

\item[\textbf{(c)}] Assume that $R$ is a field. Prove that the units of $S$
are precisely the elements of the form $\overline{\lambda+xyf}$ for $f\in
R\left[  x,y\right]  $ and $\lambda\in R^{\times}$.
\end{enumerate}

Now, define two elements $a=\overline{x}$ and $b=\overline{x+xy}$ in this ring
$S$.

\begin{enumerate}
\item[\textbf{(d)}] Prove that $aS=bS$.

\item[\textbf{(e)}] Prove that $a$ is not associate to $b$ in $S$ if $R$ is a
field. (The notion of \textquotedblleft associate\textquotedblright\ is
defined in Definition \ref{def.ring.associate}.)

\item[\textbf{(f)}] Conclude that Proposition \ref{prop.ring.associate.mutdev}
becomes false if we don't require $R$ to be an integral domain.
\end{enumerate}

[\textbf{Hint:} For part \textbf{(c)}, use Exercise \ref{exe.21hw2.1}
(observing that $\overline{xy}\in S$ is nilpotent) and then show that an
element of the form $\overline{f}$ for an $xy$-reduced polynomial $f$ can only
be a unit if $f$ is constant (because otherwise, one of the two morphisms from
part \textbf{(b)} would send this element to a non-constant univariate polynomial).

For part \textbf{(d)}, compute $\overline{x+xy}\cdot\overline{1-y}$ in $S$.]
\end{exercise}

\subsubsection{The case of arbitrary ideals}

Now what if we want to know how $P/I$ looks like for a non-principal ideal
$I$, say $I=b_{1}P+b_{2}P+\cdots+b_{k}P$ for some $b_{1},b_{2},\ldots,b_{k}\in
P$ ? Can we divide a polynomial by $I$ with remainder? Can we check whether a
polynomial belongs to $I$ ? (Remember: If $I=bP$ is a principal ideal, then
this means checking whether the polynomial is divisible by $b$. We have seen
how to do this using Theorem \ref{thm.polring.quorem-mulvar}.)

We can try to replicate the above \textquotedblleft division with
remainder\textquotedblright\ logic.

\begin{example}
\label{exa.polring.quorem-mulvar-I.ex1}Let $n=2$, and let us write $x,y$ for
the indeterminates $x_{1},x_{2}$. Let $R=\mathbb{Q}$ (just to be specific),
and let $I=b_{1}P+b_{2}P$ with%
\begin{align*}
b_{1}  &  =xy+1,\\
b_{2}  &  =y+1.
\end{align*}

Let $a\in P$ be any polynomial. We try to divide $a$ by $I$ with remainder.
This means writing $a$ in the form $a=i+r$ where $i\in I$ and $r$ is a
\textquotedblleft remainder\textquotedblright. Here, a \textquotedblleft%
\textbf{remainder}\textquotedblright\ (modulo $b_{1}$ and $b_{2}$) means a
polynomial that is both $\operatorname*{LM}b_{1}$-reduced and
$\operatorname*{LM}b_{2}$-reduced, i.e., that contains neither multiples of
$\operatorname*{LM}b_{1}$ nor multiples of $\operatorname*{LM}b_{2}$ among its
monomials. We can achieve this by subtracting multiples of $b_{1}$ and
multiples of $b_{2}$ from $a$ until no such remain. In more detail: Whenever
some monomial that is a multiple of $\operatorname*{LM}b_{1}$ appears in our
polynomial, we can subtract an appropriate multiple of $b_{1}$ from our
polynomial to remove this monomial. (Namely, the multiple of $b_{1}$ that we
choose is the one whose leading term would cancel the multiple of
$\operatorname*{LM}b_{1}$ we want to remove from our polynomial.) Similarly we
get rid of multiples of $\operatorname*{LM}b_{2}$. When no more monomials that
are multiples of $\operatorname*{LM}b_{1}$ or multiples of $\operatorname*{LM}%
b_{2}$ remain in our polynomial, then we have found our \textquotedblleft
remainder\textquotedblright.

We refer to this procedure as the \textbf{division-with-remainder algorithm}.
Note that this is a nondeterministic algorithm, in the sense that you often
have a choice of which step you make. For instance, if your polynomial
contains a monomial that is a multiple of both $\operatorname*{LM}b_{1}$ and
$\operatorname*{LM}b_{2}$ at the same time, do you remove it by subtracting a
multiple of $b_{1}$ or by subtracting a multiple of $b_{2}$ ? Thus, the
\textquotedblleft remainder\textquotedblright\ at the end might fail to be unique.

Let us check this on a specific example. Let $a=xy-y\in P$. Here is one way to
perform our division-with-remainder algorithm:%
\begin{align*}
a=xy-y\ \  &  \underset{\text{to get rid of the }xy\text{ monomial}%
}{\overset{\text{subtract }1b_{1}}{\longrightarrow}}\ \ \left(  xy-y\right)
-\left(  xy+1\right)  =-y-1\\
&  \underset{\text{to get rid of the }y\text{ monomial}%
}{\overset{\text{subtract }-1b_{2}}{\longrightarrow}}\ \ \left(  -y-1\right)
-\left(  -1\right)  \left(  y+1\right)  =0.
\end{align*}

Here is another way to do it:%
\begin{align*}
a=xy-y\ \  &  \underset{\text{to get rid of the }xy\text{ monomial}%
}{\overset{\text{subtract }xb_{2}}{\longrightarrow}}\ \ \left(  xy-y\right)
-x\left(  y+1\right)  =-x-y\\
&  \underset{\text{to get rid of the }y\text{ monomial}%
}{\overset{\text{subtract }-1b_{2}}{\longrightarrow}}\ \ \left(  -x-y\right)
-\left(  -1\right)  \left(  y+1\right)  =-x+1.
\end{align*}

Both results we have obtained are both $\operatorname*{LM}b_{1}$-reduced and
$\operatorname*{LM}b_{2}$-reduced, so they qualify as \textquotedblleft
remainders\textquotedblright\ of $a$ modulo $b_{1}$ and $b_{2}$. However, they
are not equal! So the remainder is not unique this time (unlike in Theorem
\ref{thm.polring.quorem-mulvar}). In particular, the first remainder we
obtained was $0$, which showed that $a\in I$ (because this remainder was
obtained from $a$ by subtracting multiples of $b_{1}$ and $b_{2}$, and of
course these multiples all belong to $I$); but the second remainder was not
$0$, thus allowing no such conclusion. So we don't have a sure way of telling
whether a polynomial belongs to $I$ or not; if we are unlucky, we get a
nonzero remainder even for a polynomial that does belong to $I$.
\end{example}

This is bad! But it gets even worse:

\begin{example}
\label{exa.polring.quorem-mulvar-I.ex2}A simpler example: Let $n=2$ and
$I=b_{1}P+b_{2}P$ with
\begin{align*}
b_{1}  &  =xy+x,\\
b_{2}  &  =xy+y.
\end{align*}
The polynomial $x-y$ lies in $I$ (since $x-y=b_{1}-b_{2}$), but it is both
$\operatorname*{LM}b_{1}$-reduced and $\operatorname*{LM}b_{2}$-reduced, so we
cannot see this from our division-with-remainder algorithm no matter what
choices we make (because the algorithm does nothing: $x-y$ already is a
\textquotedblleft remainder\textquotedblright). We could, of course, subtract
$b_{1}$ from $x-y$ (to obtain $\left(  x-y\right)  -\left(  xy+x\right)
=-y-xy$), but this would be a \textquotedblleft step
backwards\textquotedblright, as it would increase the leading monomial (and
even the degree) of our polynomial. The idea of the division-with-remainder
algorithm is to reduce the polynomial step by step, always \textquotedblleft
walking downhill\textquotedblright, rather than having to \textquotedblleft
cross a mountain\textquotedblright\ first (temporarily increasing the leading monomial).
\end{example}

Example \ref{exa.polring.quorem-mulvar-I.ex2} might give you an idea of what
is standing in our way here: It is the fact that when we compute $b_{1}-b_{2}%
$, the leading terms $xy$ cancel. It means, in a sense, that our $b_{1}$ and
$b_{2}$ are \textquotedblleft unnecessarily convoluted\textquotedblright; we
should perhaps fix this by replacing $b_{2}$ by the smaller polynomial
$b_{2}-b_{1}=y-x$. This simplifies $b_{2}$ but does not change $I$ (since
$b_{1}P+b_{2}P=b_{1}P+\left(  b_{2}-b_{1}\right)  P$\ \ \ \ \footnote{This is
a consequence of Exercise \ref{exe.ideals.aR+I=bR+I} (applied to $b_{2}$,
$b_{2}-b_{1}$ and $b_{1}P$ instead of $a$, $b$ and $I$).}). This is similar to
one of the row-reduction steps involved in bringing a matrix to row echelon form.

What does it mean in general that a list $\left(  b_{1},b_{2},\ldots
,b_{k}\right)  $ of polynomials is \textquotedblleft unnecessarily
convoluted\textquotedblright? The $xy$ cancellation in $b_{1}-b_{2}$ above was
easy to see; what other cancellations can lurk in a list of polynomials?

Let me formalize this question. The following definition will be a bit
long-winded but it is just giving names to the kind of observations you would
have made when trying to discuss the above algorithm:

\begin{definition}
\label{def.polring.grobas}Let $\mathbf{b}=\left(  b_{1},b_{2},\ldots
,b_{k}\right)  $ be a list of nonzero polynomials in $P$ whose leading
coefficients are units of $R$.

\begin{enumerate}
\item[\textbf{(a)}] Given two polynomials $c,d\in P$, we write
$c\underset{\mathbf{b}}{\longrightarrow}d$ (and say \textquotedblleft$c$ can
be \textbf{reduced} to $d$ in a \textbf{single step} using $\mathbf{b}%
$\textquotedblright) if

\begin{itemize}
\item some monomial $\mathfrak{m}$ appearing in $c$ is a multiple of
$\operatorname*{LM}b_{i}$ for some $i\in\left\{  1,2,\ldots,k\right\}  $;

\item we have%
\[
d=c-\dfrac{\left[  \mathfrak{m}\right]  c}{\operatorname*{LC}b_{i}}\cdot
\dfrac{\mathfrak{m}}{\operatorname*{LM}b_{i}}\cdot b_{i}.
\]
(This equation essentially says that we obtain $d$ from $c$ by subtracting the
appropriate multiple of $b_{i}$ to get rid of the monomial $\mathfrak{m}$. The
multiple is $\dfrac{\left[  \mathfrak{m}\right]  c}{\operatorname*{LC}b_{i}%
}\cdot\dfrac{\mathfrak{m}}{\operatorname*{LM}b_{i}}\cdot b_{i}$, since the
$\dfrac{\mathfrak{m}}{\operatorname*{LM}b_{i}}$ factor is needed to turn the
leading monomial of $b_{i}$ into $\mathfrak{m}$, whereas the $\dfrac{\left[
\mathfrak{m}\right]  c}{\operatorname*{LC}b_{i}}$ factor serves to make the
coefficient of this monomial the same as that in $c$. Note that the fraction
$\dfrac{\left[  \mathfrak{m}\right]  c}{\operatorname*{LC}b_{i}}\in R$ is
well-defined since $\operatorname*{LC}b_{i}$ is a unit, whereas the fraction
$\dfrac{\mathfrak{m}}{\operatorname*{LM}b_{i}}\in C^{\left(  n\right)  }$ is
well-defined since $\mathfrak{m}$ is a multiple of $\operatorname*{LM}b_{i}$.)
\end{itemize}

For instance, using the notations of Example
\ref{exa.polring.quorem-mulvar-I.ex1} and setting $\mathbf{b}=\left(
b_{1},b_{2}\right)  $, we have%
\[
xy-y\underset{\mathbf{b}}{\longrightarrow}-x-y,
\]
because we obtain $-x-y$ from $xy-y$ by subtracting the multiple $xb_{2}$ of
$b_{2}$ (which kills the $xy$ monomial). Likewise, for the same $\mathbf{b}$,
we have%
\[
5x^{2}y^{3}\underset{\mathbf{b}}{\longrightarrow}-5xy^{2},
\]
because we obtain $-5xy^{2}$ from $5x^{2}y^{3}$ by subtracting the multiple
$5xy^{2}b_{1}$ of $b_{1}$ (which kills the $x^{2}y^{3}$ monomial).

\item[\textbf{(b)}] Given two polynomials $c,d\in P$, we write $c\overset{\ast
}{\underset{\mathbf{b}}{\longrightarrow}}d$ (and say \textquotedblleft$c$ can
be \textbf{reduced} to $d$ in \textbf{many steps} using $\mathbf{b}%
$\textquotedblright) if there is a sequence $\left(  c_{0},c_{1},\ldots
,c_{m}\right)  $ of polynomials in $P$ such that $c_{0}=c$ and $c_{m}=d$ and%
\[
c_{i}\underset{\mathbf{b}}{\longrightarrow}c_{i+1}%
\ \ \ \ \ \ \ \ \ \ \text{for each }i\in\left\{  0,1,\ldots,m-1\right\}  .
\]

Note that this sequence can be trivial (i.e., we can have $m=0$), in which
case of course we have $c=d$. Thus, $c\overset{\ast}{\underset{\mathbf{b}%
}{\longrightarrow}}c$ for any $c\in P$. (Like any true algebraists, we
understand \textquotedblleft many steps\textquotedblright\ to allow
\textquotedblleft zero steps\textquotedblright.) We also can have $m=1$; thus,
$c\overset{\ast}{\underset{\mathbf{b}}{\longrightarrow}}d$ holds if
$c\underset{\mathbf{b}}{\longrightarrow}d$. (That is, \textquotedblleft many
steps\textquotedblright\ allows \textquotedblleft one step\textquotedblright.)

As an example of a nontrivial many-steps reduction, we observe that using the
notations of Example \ref{exa.polring.quorem-mulvar-I.ex1} and setting
$\mathbf{b}=\left(  b_{1},b_{2}\right)  $, we have%
\[
5x^{2}y^{3}\underset{\mathbf{b}}{\longrightarrow}-5xy^{2}\underset{\mathbf{b}%
}{\longrightarrow}5y\underset{\mathbf{b}}{\longrightarrow}-5
\]
and thus $5x^{2}y^{3}\overset{\ast}{\underset{\mathbf{b}}{\longrightarrow}}-5$.

\item[\textbf{(c)}] We say that a polynomial $r\in P$ is $\mathbf{b}%
$\textbf{-reduced} if it is $\operatorname*{LM}b_{i}$-reduced for all
$i\in\left\{  1,2,\ldots,k\right\}  $. This is equivalent to saying that there
exists no polynomial $s\in P$ with $r\underset{\mathbf{b}}{\longrightarrow}s$
(that is, \textquotedblleft$r$ cannot be reduced any further using
$\mathbf{b}$\textquotedblright).

\item[\textbf{(d)}] A \textbf{remainder} of a polynomial $a\in P$ modulo
$\mathbf{b}$ means a $\mathbf{b}$-reduced polynomial $r\in P$ such that
$a\overset{\ast}{\underset{\mathbf{b}}{\longrightarrow}}r$. Such a remainder
always exists (this is not hard to show), but is not always unique (as we have
seen in Example \ref{exa.polring.quorem-mulvar-I.ex1}).

\item[\textbf{(e)}] We say that the list $\mathbf{b}$ is a \textbf{Gr\"{o}bner
basis} if any $a\in P$ has a \textbf{unique} remainder modulo $\mathbf{b}$.
\end{enumerate}
\end{definition}

(Don't take the word \textquotedblleft basis\textquotedblright\ in
\textquotedblleft Gr\"{o}bner basis\textquotedblright\ to heart. It is closer
to \textquotedblleft generating set\textquotedblright\ or \textquotedblleft
spanning set\textquotedblright\ than to any sort of \textquotedblleft
basis\textquotedblright\ in linear algebra. In particular, a Gr\"{o}bner basis
can be $R$-linearly dependent or even contain the same polynomial twice.)

So we have seen that not every list of nonzero polynomials is a Gr\"{o}bner
basis. This leads to the following two questions:

\begin{itemize}
\item Can we \textbf{tell} whether a list of nonzero polynomials is a
Gr\"{o}bner basis? (We cannot afford to check every $a\in P$ and every way of
reducing it modulo $\mathbf{b}$.)

\item If a list is not a Gr\"{o}bner basis, can we at least \textbf{find} a
Gr\"{o}bner basis that generates the same ideal as the list?
\end{itemize}

If $R$ is not a field, then the answer to the second question is
\textquotedblleft no\textquotedblright\ for reasons that should be familiar
from the univariate case (non-unit leading coefficients).

When $R$ is a field, Bruno Buchberger has answered both questions in the
positive in the 1960s (and his answer to the first question applies just as
well to arbitrary $R$). The algorithms he found are one of the pillars of
modern computer algebra. I will state the main results without proof, but you
can find proofs in the literature (e.g., \cite[\S 9.6]{DumFoo04} or
\cite[Chapter 1]{deGraaf-CA}).

We will need the notion of an \textbf{S-polynomial}:

\begin{definition}
Let $f,g\in P$ be nonzero polynomials whose leading coefficients are units of
$R$.

Let $\mathfrak{p}=x_{1}^{p_{1}}x_{2}^{p_{2}}\cdots x_{n}^{p_{n}}%
=\operatorname*{LM}f$ and $\mathfrak{q}=x_{1}^{q_{1}}x_{2}^{q_{2}}\cdots
x_{n}^{q_{n}}=\operatorname*{LM}g$ be their leading monomials, and let
$\lambda=\operatorname*{LC}f$ and $\mu=\operatorname*{LC}g$ be their leading
coefficients. So%
\begin{align*}
f  &  =\lambda\mathfrak{p}+\left(  \text{smaller terms}\right)  ;\\
g  &  =\mu\mathfrak{q}+\left(  \text{smaller terms}\right)  .
\end{align*}

Let
\[
\mathfrak{m}=x_{1}^{\max\left\{  p_{1},q_{1}\right\}  }x_{2}^{\max\left\{
p_{2},q_{2}\right\}  }\cdots x_{n}^{\max\left\{  p_{n},q_{n}\right\}  }.
\]
(This is the lcm of $\mathfrak{p}$ and $\mathfrak{q}$ among the monomials; it
is the smallest-degree monomial that is divisible by both $\mathfrak{p}$ and
$\mathfrak{q}$.) Note that $\dfrac{\mathfrak{m}}{\mathfrak{p}}$ and
$\dfrac{\mathfrak{m}}{\mathfrak{q}}$ are well-defined monomials (since
$\mathfrak{p}\mid\mathfrak{m}$ and $\mathfrak{q}\mid\mathfrak{m}$).

The \textbf{S-polynomial} (short for \textbf{syzygy polynomial}) of $f$ and
$g$ is defined to be the polynomial%
\[
S\left(  f,g\right)  :=\dfrac{1}{\lambda}\cdot\dfrac{\mathfrak{m}%
}{\mathfrak{p}}\cdot f-\dfrac{1}{\mu}\cdot\dfrac{\mathfrak{m}}{\mathfrak{q}%
}\cdot g\in P.
\]

Here is the intuition behind this: $S\left(  f,g\right)  $ is the simplest way
to form a $P$-linear combination of $f$ and $g$ in which the leading terms of
$f$ and $g$ cancel. Namely, in order to obtain such a $P$-linear combination,
we must first rescale $f$ and $g$ so that their leading coefficients become
equal (this can be achieved by scaling $f$ by $\dfrac{1}{\lambda}$ and scaling
$g$ by $\dfrac{1}{\mu}$); then we must multiply them with appropriate
monomials to make their leading monomials equal (this can be achieved by
multiplying $f$ by $\dfrac{\mathfrak{m}}{\mathfrak{p}}$ and multiplying $g$ by
$\dfrac{\mathfrak{m}}{\mathfrak{q}}$, so that both leading monomials become
$\mathfrak{m}$). The resulting two polynomials have equal leading terms
(namely, $\mathfrak{m}$), so their leading terms cancel out when we subtract
them. The result of this subtraction is $S\left(  f,g\right)  $. To be more
specific, when we multiplied $f$ and $g$ with appropriate monomials to make
their leading monomials equal, we made sure to choose the latter monomials as
low-degree as possible; this is why we took $\mathfrak{m}$ to be the lcm of
$\mathfrak{p}$ and $\mathfrak{q}$ and not some other monomial divisible by
$\mathfrak{p}$ and $\mathfrak{q}$ (such as the product $\mathfrak{pq}$).
\end{definition}

\begin{example}
For $n=2$ (and denoting $x_{1},x_{2}$ by $x,y$ as usual), we have%
\[
S\left(  x^{2}y+1,\ xy^{2}+1\right)  =y\left(  x^{2}y+1\right)  -x\left(
xy^{2}+1\right)  =y-x
\]
and%
\[
S\left(  xy+1,\ 2x\right)  =1\left(  xy+1\right)  -\dfrac{1}{2}\cdot
y\cdot2x=1
\]
and%
\begin{align*}
S\left(  x^{2}+2x+y,\ y^{2}+x+2y\right)   &  =y^{2}\left(  x^{2}+2x+y\right)
-x^{2}\left(  y^{2}+x+2y\right) \\
&  =-x^{3}-2x^{2}y+2xy^{2}+y^{3}.
\end{align*}

\end{example}

Note that the cancellation of the leading terms in the construction of
$S\left(  f,g\right)  $ is precisely the sort of cancellation that prevented
us from having a unique remainder in our above examples.

The following crucial theorem says that these cancellations are a canary in
the mine: If they don't happen, then the list is a Gr\"{o}bner basis.

\begin{theorem}
[Buchberger's criterion]\label{thm.polring.grobas.buch-crit}Let $\mathbf{b}%
=\left(  b_{1},b_{2},\ldots,b_{k}\right)  $ be a list of nonzero polynomials
in $P$ whose leading coefficients are units of $R$.

Then, $\mathbf{b}$ is a Gr\"{o}bner basis if and only if every $i<j$ satisfy%
\[
S\left(  b_{i},b_{j}\right)  \overset{\ast}{\underset{\mathbf{b}%
}{\longrightarrow}}0.
\]

\end{theorem}

The idea behind this theorem is that a list of polynomials (whose leading
coefficients are units) is a Gr\"{o}bner basis if and only if any S-polynomial
of two polynomials in the list reduces to $0$ modulo the list. Note that
\textquotedblleft reduces to $0$ modulo the list\textquotedblright\ means that
there is some way to get the remainder $0$ when applying the
division-with-remainder algorithm to this S-polynomial; we are not requiring
that \textbf{every} way of applying the division-with-remainder algorithm to
it will give $0$. (But this will follow automatically if we have shown that
$\mathbf{b}$ is a Gr\"{o}bner basis.)

\begin{example}
Let $n=2$, and write $x,y$ for $x_{1},x_{2}$. Let $I=b_{1}P+b_{2}P$, where%
\begin{align*}
b_{1}  &  =xy+1,\\
b_{2}  &  =y+1.
\end{align*}
We already know from Example \ref{exa.polring.quorem-mulvar-I.ex1} that
$\left(  b_{1},b_{2}\right)  $ is not a Gr\"{o}bner basis, but let us now see
this using Buchberger's criterion:%
\[
S\left(  b_{1},b_{2}\right)  =1\left(  xy+1\right)  -x\left(  y+1\right)
=1-x.
\]
This polynomial $1-x$ is already $\mathbf{b}$-reduced (where $\mathbf{b}%
=\left(  b_{1},b_{2}\right)  $), and it is not $0$, so we \textbf{don't} have
$S\left(  b_{1},b_{2}\right)  \overset{\ast}{\underset{\mathbf{b}%
}{\longrightarrow}}0$. Thus, Theorem \ref{thm.polring.grobas.buch-crit}
confirms again that our $\mathbf{b}$ is not a Gr\"{o}bner basis.
\end{example}

\begin{example}
\label{exa.polring.grobas.xx-yz.1}Let $n=3$, and write $x,y,z$ for
$x_{1},x_{2},x_{3}$. Let $I=b_{1}P+b_{2}P+b_{3}P$, where%
\begin{align*}
b_{1}  &  =x^{2}-yz,\\
b_{2}  &  =y^{2}-zx,\\
b_{3}  &  =z^{2}-xy.
\end{align*}
Is $\mathbf{b}:=\left(  b_{1},b_{2},b_{3}\right)  $ a Gr\"{o}bner basis? We
check this using Buchberger's criterion. First, we rewrite $b_{1},b_{2},b_{3}$
in a way that their leading terms are up front:%
\begin{align*}
b_{1}  &  =x^{2}-yz,\\
b_{2}  &  =-zx+y^{2},\\
b_{3}  &  =-xy+z^{2}.
\end{align*}
(It is generally advised to always write the terms of a polynomial in the
deg-lex order, from highest to lowest, when performing division-with-remainder
or computing S-polynomials. Otherwise, it is too easy to get confused about
which terms are leading!)

Now, we compute remainders of $S\left(  b_{i},b_{j}\right)  $ modulo
$\mathbf{b}$ for all $i<j$:

\begin{itemize}
\item We have%
\begin{align*}
S\left(  b_{1},b_{2}\right)   &  =S\left(  x^{2}-yz,\ -zx+y^{2}\right) \\
&  =z\left(  x^{2}-yz\right)  -\left(  -x\right)  \left(  -zx+y^{2}\right)
=xy^{2}-yz^{2}\\
&  \underset{\mathbf{b}}{\longrightarrow}\left(  xy^{2}-yz^{2}\right)
-\left(  -y\right)  \left(  -xy+z^{2}\right) \\
&  \ \ \ \ \ \ \ \ \ \ \ \ \ \ \ \ \ \ \ \ \left(
\begin{array}
[c]{c}%
\text{here, we subtracted }-yb_{3}\\
\text{in order to remove the }xy^{2}\text{ monomial}%
\end{array}
\right) \\
&  =0,
\end{align*}
so that $S\left(  b_{1},b_{2}\right)  \overset{\ast}{\underset{\mathbf{b}%
}{\longrightarrow}}0$.

\item We have%
\begin{align*}
S\left(  b_{1},b_{3}\right)   &  =S\left(  x^{2}-yz,\ -xy+z^{2}\right) \\
&  =y\left(  x^{2}-yz\right)  -\left(  -x\right)  \left(  -xy+z^{2}\right)
=xz^{2}-y^{2}z\\
&  \underset{\mathbf{b}}{\longrightarrow}\left(  xz^{2}-y^{2}z\right)
-\left(  -z\right)  \left(  -zx+y^{2}\right) \\
&  \ \ \ \ \ \ \ \ \ \ \ \ \ \ \ \ \ \ \ \ \left(
\begin{array}
[c]{c}%
\text{here, we subtracted }-zb_{2}\\
\text{in order to remove the }xz^{2}\text{ monomial}%
\end{array}
\right) \\
&  =0,
\end{align*}
so that $S\left(  b_{1},b_{3}\right)  \overset{\ast}{\underset{\mathbf{b}%
}{\longrightarrow}}0$.

\item We have%
\begin{align*}
S\left(  b_{2},b_{3}\right)   &  =S\left(  -zx+y^{2},\ -xy+z^{2}\right) \\
&  =\dfrac{1}{-1}y\left(  -zx+y^{2}\right)  -\dfrac{1}{-1}z\left(
-xy+z^{2}\right) \\
&  =-y^{3}+z^{3}\text{ is }\mathbf{b}\text{-reduced and not }0.
\end{align*}
Thus, we \textbf{do not} have $S\left(  b_{2},b_{3}\right)  \overset{\ast
}{\underset{\mathbf{b}}{\longrightarrow}}0$. This shows that $\left(
b_{1},b_{2},b_{3}\right)  $ is \textbf{not} a Gr\"{o}bner basis.
\end{itemize}

(This example was a bit unusual in that our many-step reductions were actually
one-step reductions. But it is certainly not unusual in that we have wasted a
lot of work before getting the answer \textquotedblleft no\textquotedblright.)
\end{example}

Buchberger's criterion is proved (e.g.) in any of \cite[p. 324]{DumFoo04},
\cite[Theorem 5.6.8]{Laurit09}, \cite[proof of Theorem 1.1.33]{deGraaf-CA},
\cite[Theorem 5.48, (iii) $\Longleftrightarrow$ (i)]{BecWei98}, \cite[Theorem
2.14]{EneHer12}, \cite[V.3 \textbf{i)} $\Longleftrightarrow$ \textbf{ii)}%
]{MalBlo15} and (in a slight variation) in \cite[Chapter 2, \S 6, Theorem
6]{CoLiOs15}\footnote{\textbf{Caveat:} Different authors define the relation
$\overset{\ast}{\underset{\mathbf{b}}{\longrightarrow}}$ somewhat differently.
The proofs usually apply to all versions.}. The \textquotedblleft only
if\textquotedblright\ part is obvious; the \textquotedblleft
if\textquotedblright\ part is interesting.

Gr\"{o}bner bases help us better understand ideals of $P$:

\begin{definition}
Let $I$ be an ideal of $P$. A \textbf{Gr\"{o}bner basis of }$I$ means a
Gr\"{o}bner basis $\left(  b_{1},b_{2},\ldots,b_{k}\right)  $ that generates
$I$ (that is, that satisfies $I=b_{1}P+b_{2}P+\cdots+b_{k}P$).
\end{definition}

\begin{corollary}
[Macaulay's basis theorem]\label{cor.polring.mod-grobas}Let $\mathbf{b}%
=\left(  b_{1},b_{2},\ldots,b_{k}\right)  $ be a list of nonzero polynomials
in $P$ whose leading coefficients are units of $R$. Assume that $\mathbf{b}$
is a Gr\"{o}bner basis.

Let $I$ be the ideal $b_{1}P+b_{2}P+\cdots+b_{k}P$ of $P$. Then, each element
of $P/I$ can be uniquely written in the form%
\[
\sum_{\substack{\mathfrak{m}\text{ is a }\mathbf{b}\text{-reduced}%
\\\text{monomial}}}a_{\mathfrak{m}}\overline{\mathfrak{m}}%
\ \ \ \ \ \ \ \ \ \ \text{with }a_{\mathfrak{m}}\in R
\]
(where all but finitely many $\mathfrak{m}$ satisfy $a_{\mathfrak{m}}=0$).
Equivalently, the family $\left(  \overline{\mathfrak{m}}\right)
_{\mathfrak{m}\text{ is a }\mathbf{b}\text{-reduced monomial}}$ is a basis of
the $R$-module $P/I$. If none of the polynomials $b_{1},b_{2},\ldots,b_{k}$ is
constant, then the ring $P/I$ contains \textquotedblleft a copy of
$R$\textquotedblright.
\end{corollary}

\begin{proof}
LTTR.
\end{proof}

To summarize: If we know a Gr\"{o}bner basis of an ideal $I$ of $P$, then we
know a lot about $I$ (in particular, we can tell when a polynomial belongs to
$I$, and we can find a basis for $P/I$). But how do we find a Gr\"{o}bner
basis of an ideal? Is there always one?

Not for arbitrary $R$. But if $R$ is a field, then yes:

\begin{theorem}
[Buchberger's theorem]\label{thm.polring.grobas.buch-alg}Let $R$ be a field.
Let $I$ be an ideal of the polynomial ring $P=R\left[  x_{1},x_{2}%
,\ldots,x_{n}\right]  $. Then, $I$ has a Gr\"{o}bner basis.

Moreover, if $b_{1},b_{2},\ldots,b_{k}$ are nonzero polynomials such that
$I=b_{1}P+b_{2}P+\cdots+b_{k}P$, then we can construct a Gr\"{o}bner basis of
$I$ by the following algorithm (\textbf{Buchberger's algorithm}):

\begin{itemize}
\item Initially, let $\mathbf{b}$ be the list $\left(  b_{1},b_{2}%
,\ldots,b_{k}\right)  $.

\item As long as there exist two entries of $\mathbf{b}$ whose S-polynomial
has a nonzero remainder modulo $\mathbf{b}$, we append this remainder to the
list. (It is enough to compute one remainder for each pair of entries of
$\mathbf{b}$.)

\item Once no such two entries exist any more, we are done: $\mathbf{b}$ is a
Gr\"{o}bner basis of $I$.
\end{itemize}

This algorithm always terminates after finitely many steps (i.e., we don't
keep adding new entries to $\mathbf{b}$ forever).
\end{theorem}

We won't prove this (see, e.g., \cite[Proposition 1.1.35]{deGraaf-CA} or
\cite[Theorem 5.7.2]{Laurit09} for proofs), but we will give examples:

\begin{example}
Let $n=3$, and write $x,y,z$ for $x_{1},x_{2},x_{3}$. Let $I=b_{1}%
P+b_{2}P+b_{3}P$, where%
\begin{align*}
b_{1}  &  =x^{2}-yz,\\
b_{2}  &  =y^{2}-zx,\\
b_{3}  &  =z^{2}-xy.
\end{align*}
We want to find a Gr\"{o}bner basis of this ideal $I$.

As we have seen in Example \ref{exa.polring.grobas.xx-yz.1}, the list
$\mathbf{b}:=\left(  b_{1},b_{2},b_{3}\right)  $ is not a Gr\"{o}bner basis,
since $S\left(  b_{2},b_{3}\right)  =-y^{3}+z^{3}$ does not have remainder $0$
modulo $\mathbf{b}$. Its remainder is $-y^{3}+z^{3}$ itself. Thus, following
Buchberger's algorithm, we append this remainder to the list. That is, we set
$b_{4}=-y^{3}+z^{3}$, and continue with the list $\left(  b_{1},b_{2}%
,b_{3},b_{4}\right)  $. We call this list $\mathbf{b}$ again.

Since $\mathbf{b}$ has grown, we must now also check whether the new
S-polynomials%
\[
S\left(  b_{1},b_{4}\right)  ,\ \ S\left(  b_{2},b_{4}\right)  ,\ \ S\left(
b_{3},b_{4}\right)
\]
reduce to $0$ modulo $\mathbf{b}$. Fortunately, they do. Thus, our new list
$\mathbf{b}=\left(  b_{1},b_{2},b_{3},b_{4}\right)  $ is a Gr\"{o}bner basis
of $I$.
\end{example}

\begin{example}
Let $n=3$, and write $x,y,z$ for $x_{1},x_{2},x_{3}$. Assume that
$R=\mathbb{Q}$. Let $I=b_{1}P+b_{2}P+b_{3}P$, where%
\begin{align*}
b_{1}  &  =x^{2}+xy,\\
b_{2}  &  =y^{2}+yz,\\
b_{3}  &  =z^{2}+zx.
\end{align*}
Then, again, it is not hard to see that $\left(  b_{1},b_{2},b_{3}\right)  $
is not a Gr\"{o}bner basis of $I$. Using Buchberger's algorithm, we can easily
compute its Gr\"{o}bner basis. For example, $I$ has Gr\"{o}bner basis%
\[
\left(  x^{2}+xy,\ \ y^{2}+yz,\ \ xz+z^{2},\ \ yz^{2}-z^{3},\ \ z^{4}\right)
.
\]
(Note that the Gr\"{o}bner basis of an ideal is not unique, so you might get a
different one if you perform Buchberger's algorithm differently. When there
are several pairs $\left(  b_{i},b_{j}\right)  $ whose S-polynomial does not
reduce to $0$, you have a choice of which of these pairs you handle first.)

This Gr\"{o}bner basis reveals that $z^{4}\in I$ but $z^{3}\notin I$ (since
$z^{3}$ is reduced modulo the above Gr\"{o}bner basis). Just working from the
original definition of $I$, this would be far from obvious!

You can do Gr\"{o}bner basis computations with most computer algebra systems
(e.g., SageMath, Mathematica, Singular, SymPy). For example,
\href{https://sagecell.sagemath.org/?z=eJwdjLEKgCAUAPegf3g0mYlDc7W3VR9QKD5EMAVbnn591nbD3W1yIgFZQFlghi36HOLtlD9csGzfBcRkMM2dQeuRur5t1s-TzqDyjM4RBiBe-_xj5qWufiycPlvaFFEHTJdWj3tY_wKyCiIZ&lang=sage&interacts=eJyLjgUAARUAuQ==}{here}
is SageMath code for the Gr\"{o}bner basis of the above ideal. Note that we
took $R=\mathbb{Q}$ in this computation (the \textquotedblleft\texttt{QQ}%
\textquotedblright\ means the field of rational numbers, and the
\textquotedblleft\texttt{order="deglex"}\textquotedblright\ argument ensures
that the deg-lex order is used), but the same computation works over any field
$R$ of characteristic $\neq2$ (and, because our ideal is rather nice, even
over any commutative ring $R$ in which $2$ is invertible; this is not automatic).
\end{example}

\begin{exercise}
Let $n=3$, and let us rename the indeterminates $x_{1},x_{2},x_{3}$ as
$x,y,z$. Define two polynomials $b_{1}$ and $b_{2}$ in $P$ by%
\[
b_{1}=x^{2}-y\ \ \ \ \ \ \ \ \ \ \text{and}\ \ \ \ \ \ \ \ \ \ b_{2}=x^{3}-z.
\]
Let $I=b_{1}P+b_{2}P$. Find a Gr\"{o}bner basis of $I$. (Feel free to assume
that $R=\mathbb{Q}$ for simplicity.)
\end{exercise}

\begin{exercise}
Let $n=3$, and write $x,y,z$ for $x_{1},x_{2},x_{3}$. Let $I=b_{1}%
P+b_{2}P+b_{3}P$, where%
\begin{align*}
b_{1}  &  =x+xyz,\\
b_{2}  &  =y+xyz,\\
b_{3}  &  =z+xyz.
\end{align*}
Find a Gr\"{o}bner basis of $I$. (Feel free to assume that $R=\mathbb{Q}$ for simplicity.)
\end{exercise}

\begin{exercise}
\label{exe.polring.grobas.xx+yz.1}Let $n=3$, and write $x,y,z$ for
$x_{1},x_{2},x_{3}$. Let $I=b_{1}P+b_{2}P+b_{3}P$, where%
\begin{align*}
b_{1}  &  =x^{2}+yz,\\
b_{2}  &  =y^{2}+zx,\\
b_{3}  &  =z^{2}+xy.
\end{align*}

\begin{enumerate}
\item[\textbf{(a)}] Find a Gr\"{o}bner basis of $I$, assuming that
$R=\mathbb{Q}$.

\item[\textbf{(b)}] Prove that $I$ has no Gr\"{o}bner basis when
$R=\mathbb{Z}$.

[\textbf{Hint:} Show that $P/I$ is not a free $R$-module when $R=\mathbb{Z}$.]
\end{enumerate}
\end{exercise}

\subsubsection{Monomial orders}

We have so far been using the deg-lex order on the monomials. There are many
other total orders that share most of its nice properties and are often more
suited for specific problems.

Let me only mention the \textbf{lexicographic order}, which is defined just as
the deg-lex order but without taking the degree into account. That is:

\begin{definition}
We define a total order $\prec$ (called the \textbf{lexicographic order}, or
-- for short -- the \textbf{lex order}) on the set $C^{\left(  n\right)  }$ of
all monomials as follows:

For two monomials $\mathfrak{m}=x_{1}^{a_{1}}x_{2}^{a_{2}}\cdots x_{n}^{a_{n}%
}$ and $\mathfrak{n}=x_{1}^{b_{1}}x_{2}^{b_{2}}\cdots x_{n}^{b_{n}}$, we
declare that $\mathfrak{m}\prec\mathfrak{n}$ if and only if

\begin{itemize}
\item there is an $i\in\left\{  1,2,\ldots,n\right\}  $ such that $a_{i}\neq
b_{i}$, and the \textbf{smallest} such $i$ satisfies $a_{i}<b_{i}$.
\end{itemize}
\end{definition}

Several properties of the deg-lex order were collected in Proposition
\ref{prop.polring.deglex.basics}. All of those properties except for
Proposition \ref{prop.polring.deglex.basics} \textbf{(d)} hold for the lex
order as well. Proposition \ref{prop.polring.deglex.basics} \textbf{(d)} fails
for the lex order (for $n>1$), because $x_{1}$ is larger (with respect to the
lex order) than \textbf{any} power of $x_{2}$ (and, of course, there are
infinitely many powers of $x_{2}$). Proposition
\ref{prop.polring.deglex.basics} \textbf{(e)} is still true for the lex order,
but its proof is harder. However, the theory of Gr\"{o}bner bases does not use
Proposition \ref{prop.polring.deglex.basics} \textbf{(d)}, so it still can be
done with the lex order. This yields new (in general, different) Gr\"{o}bner bases.

\begin{example}
\label{exa.polring.grobas.lex.xx-y}Let $n=3$ and let $I=\left(  x^{2}%
-y\right)  P+\left(  y^{2}-z\right)  P+\left(  z^{2}-x\right)  P$ (where we
write $x,y,z$ for $x_{1},x_{2},x_{3}$). Then, a Gr\"{o}bner basis of $I$ with
respect to the deg-lex order is%
\[
\left(  x^{2}-y,\ \ y^{2}-z,\ \ z^{2}-x\right)
\]
(this is precisely the list of generators that we started with). But this is
not a Gr\"{o}bner basis with respect to the lex order. Instead, a Gr\"{o}bner
basis of $I$ with respect to the lex order is%
\[
\left(  x-z^{2},\ \ y-z^{4},\ \ z^{8}-z\right)  .
\]

\end{example}

\begin{example}
Let $n=3$ and let $I=\left(  x^{2}-y^{3}\right)  P+\left(  y^{4}-z^{2}\right)
P+\left(  z^{2}-x^{5}\right)  P$ (where we write $x,y,z$ for $x_{1}%
,x_{2},x_{3}$). Then, a Gr\"{o}bner basis of $I$ with respect to the deg-lex
order is%
\begin{align*}
&  \left(  z^{6}-yz^{2},\ \ x^{3}z^{2}-yz^{2},\ \ xy^{2}z^{2}-z^{2}%
,\ \ xz^{4}-y^{2}z^{2},\right. \\
&  \ \ \ \ \ \ \ \ \ \ \left.  yz^{4}-xz^{2},\ \ x^{4}-y^{2}z^{2}%
,\ \ x^{2}y-z^{2},\ \ y^{3}-x^{2}\right)  .
\end{align*}
But a Gr\"{o}bner basis of $I$ with respect to the lex order is%
\[
\left(  x^{2}-y^{3},\ \ xz^{2}-z^{8},\ \ y^{4}-z^{2},\ \ yz^{2}-z^{6}%
,\ \ z^{16}-z^{2}\right)  .
\]

In the SageMath computer algebra system, you can signal the use of the lex
order (as opposed to the deg-lex order) by replacing \textquotedblleft%
\texttt{PolynomialRing(QQ, order="deglex")}\textquotedblright\ by
\textquotedblleft\texttt{PolynomialRing(QQ, order="lex")}\textquotedblright.
\end{example}

This last example illustrates one reason to vary the total order on monomials:
Gr\"{o}bner bases can be rather long, even if the ideal is easy to write down.
The size of a Gr\"{o}bner basis can be doubly exponential in the size of the
input (see \cite{MoSrSa24}). In real life, this worst-case scenario doesn't
happen very often, but when it does, switching to a different monomial order
will often ameliorate it. (Think of it as a way to re-roll the dice if you got
an unlucky roll.)

\begin{exercise}
Let $n=3$, and let us rename the indeterminates $x_{1},x_{2},x_{3}$ as
$x,y,z$. Define two polynomials $b_{1}$ and $b_{2}$ in $P$ by%
\[
b_{1}=x^{2}-yz\ \ \ \ \ \ \ \ \ \ \text{and}\ \ \ \ \ \ \ \ \ \ b_{2}%
=x-y^{2}.
\]
Let $I=b_{1}P+b_{2}P$.

\begin{enumerate}
\item[\textbf{(a)}] Find a Gr\"{o}bner basis of $I$ with respect to the
deg-lex order.

\item[\textbf{(b)}] Find a Gr\"{o}bner basis of $I$ with respect to the lex order.
\end{enumerate}
\end{exercise}

The following exercise generalizes Example \ref{exa.polring.grobas.lex.xx-y}
from $n=3$ to arbitrary $n$:

\begin{exercise}
Let $n$ be a positive integer. Define $n$ polynomials $b_{1},b_{2}%
,\ldots,b_{n}$ in $P$ by setting%
\begin{align*}
b_{i}  &  =x_{i}^{2}-x_{i+1}\ \ \ \ \ \ \ \ \ \ \text{for each }i\in\left\{
1,2,\ldots,n-1\right\}  \ \ \ \ \ \ \ \ \ \ \text{and}\\
b_{n}  &  =x_{n}^{2}-x_{1}.
\end{align*}
Let $I=b_{1}P+b_{2}P+\cdots+b_{n}P$.

\begin{enumerate}
\item[\textbf{(a)}] Find a Gr\"{o}bner basis of $I$ with respect to the
deg-lex order. (This should go very quick!)

\item[\textbf{(b)}] Find a Gr\"{o}bner basis of $I$ with respect to the lex order.
\end{enumerate}

[\textbf{Hint:} For part \textbf{(b)}, compute the answer for some small
values of $n$ and spot the pattern.]
\end{exercise}

\subsection{\label{sec.polys2.polsys}Solving polynomial systems using
Gr\"{o}bner bases}

Another occasion to use Gr\"{o}bner bases (and the lex order in particular) is
solving systems of polynomial equations. Polynomial equations are closely
connected to ideals:

\begin{definition}
Let $b_{1},b_{2},\ldots,b_{k}$ be $k$ polynomials in $P$, and let $A$ be a
commutative $R$-algebra. Then, a \textbf{root} (or, alternatively, a
\textbf{common root}) of $\left(  b_{1},b_{2},\ldots,b_{k}\right)  $ in $A$
means an $n$-tuple $\left(  a_{1},a_{2},\ldots,a_{n}\right)  \in A^{n}$ such
that%
\[
b_{i}\left(  a_{1},a_{2},\ldots,a_{n}\right)  =0\ \ \ \ \ \ \ \ \ \ \text{for
all }i\in\left\{  1,2,\ldots,k\right\}  .
\]

\end{definition}

This definition generalizes the standard notion of a root of a polynomial to
multiple variables and multiple polynomials.

Thus, solving systems of polynomial equations means finding roots of lists of
polynomials. It turns out that the list of polynomials doesn't really matter;
what does is the ideal it generates:

\begin{proposition}
\label{prop.polysys.ideal}Let $b_{1},b_{2},\ldots,b_{k}$ be $k$ polynomials in
$P$, and let $A$ be a commutative $R$-algebra.

Then, the roots of $\left(  b_{1},b_{2},\ldots,b_{k}\right)  $ in $A$ depend
only on the ideal generated by $b_{1},b_{2},\ldots,b_{k}$, rather than on the
polynomials $b_{1},b_{2},\ldots,b_{k}$ themselves.

More concretely: If $I=b_{1}P+b_{2}P+\cdots+b_{k}P$ is the ideal of $P$
generated by $b_{1},b_{2},\ldots,b_{k}$, then the roots of $\left(
b_{1},b_{2},\ldots,b_{k}\right)  $ are precisely the $n$-tuples $\left(
a_{1},a_{2},\ldots,a_{n}\right)  \in A^{n}$ such that%
\[
f\left(  a_{1},a_{2},\ldots,a_{n}\right)  =0\ \ \ \ \ \ \ \ \ \ \text{for all
}f\in I.
\]

\end{proposition}

\begin{proof}
Easy, LTTR. (You have to prove that if $\left(  a_{1},a_{2},\ldots
,a_{n}\right)  \in A^{n}$ is a root of $\left(  b_{1},b_{2},\ldots
,b_{k}\right)  $, then $f\left(  a_{1},a_{2},\ldots,a_{n}\right)  =0$ for all
$f\in I$. But this is easy: Each $f\in I$ is a $P$-linear combination
$c_{1}b_{1}+c_{2}b_{2}+\cdots+c_{k}b_{k}$ of $\left(  b_{1},b_{2},\ldots
,b_{k}\right)  $, and therefore satisfies%
\begin{align*}
&  f\left(  a_{1},a_{2},\ldots,a_{n}\right) \\
&  =c_{1}\left(  a_{1},a_{2},\ldots,a_{n}\right)  \underbrace{b_{1}\left(
a_{1},a_{2},\ldots,a_{n}\right)  }_{=0}+\,c_{2}\left(  a_{1},a_{2}%
,\ldots,a_{n}\right)  \underbrace{b_{2}\left(  a_{1},a_{2},\ldots
,a_{n}\right)  }_{=0}\\
&  \ \ \ \ \ \ \ \ \ \ +\cdots+c_{k}\left(  a_{1},a_{2},\ldots,a_{n}\right)
\underbrace{b_{k}\left(  a_{1},a_{2},\ldots,a_{n}\right)  }_{=0}\\
&  =0.
\end{align*}
The converse is even more obvious, since the polynomials $b_{1},b_{2}%
,\ldots,b_{k}$ all belong to $I$.)
\end{proof}

Thus, if we want to find the roots of $\left(  b_{1},b_{2},\ldots
,b_{k}\right)  $, we can replace $\left(  b_{1},b_{2},\ldots,b_{k}\right)  $
by any other tuple of polynomials that generates the same ideal of $P$. (This
is just the polynomial analogue of the classical \textquotedblleft
addition\textquotedblright\ strategy for solving systems of linear equations.)

One of the most useful ways to do this is to replace $\left(  b_{1}%
,b_{2},\ldots,b_{k}\right)  $ by a Gr\"{o}bner basis of the ideal it generates
-- particularly, by a Gr\"{o}bner basis with respect to the lex order. Let us
see how this helps on an example:

\begin{example}
Recall Exercise \ref{exe.21hw0.8}:

\begin{statement}
Solve the following system of equations:%
\begin{align*}
a^{2}+b+c  &  =1,\\
b^{2}+c+a  &  =1,\\
c^{2}+a+b  &  =1
\end{align*}
for three complex numbers $a,b,c$.
\end{statement}

Let us formalize this in terms of polynomials and roots. We set $R=\mathbb{Q}$
and $n=3$, and we write $x,y,z$ for $x_{1},x_{2},x_{3}$. Thus, the exercise is
asking for the roots of%
\[
\left(  x^{2}+y+z-1,\ \ y^{2}+z+x-1,\ \ z^{2}+x+y-1\right)
\]
in the $\mathbb{Q}$-algebra $\mathbb{C}$.

Let $I$ be the ideal of $P=\mathbb{Q}\left[  x,y,z\right]  $ generated by the
three polynomials $x^{2}+y+z-1,\ \ y^{2}+z+x-1,\ \ z^{2}+x+y-1$. Using a
computer (or a lot of patience), we can easily find a Gr\"{o}bner basis of $I$
with respect to the lex order. We get%
\[
\left(  x+y+z^{2}-1,\ \ y^{2}-y-z^{2}+z,\ \ yz^{2}+\dfrac{1}{2}z^{4}-\dfrac
{1}{2}z^{2},\ \ z^{6}-4z^{4}+4z^{3}-z^{2}\right)  .
\]

We observe that the last polynomial in this Gr\"{o}bner basis only involves
the variable $z$ ! Thus, the $c$ entry in each of the solutions $\left(
a,b,c\right)  $ of our system must be a root of this polynomial $z^{6}%
-4z^{4}+4z^{3}-z^{2}$. We can therefore find all possibilities for $c$ by
finding the roots of this polynomial (I am here assuming that you can solve
univariate polynomials; we will learn a bit more about this in Section
\ref{sec.polys2.factor} perhaps). In our concrete case, we can easily do
this:
\[
z^{6}-4z^{4}+4z^{3}-z^{2}=z^{2}\left(  z-1\right)  ^{2}\left(  z^{2}%
+2z-1\right)  .
\]
Thus, the options for $c$ are $0,1,\sqrt{2}-1,-\sqrt{2}-1$.

Now let us find $b$. Either we use the symmetry of the original system to
argue that the options for $b$ must be the same as for $c$; or we use the
second-to-last polynomial in our Gr\"{o}bner basis (or the second one) to
compute $b$ now that $c$ is known. At last, we get to $a$ in a similar way.

In the end, we get finitely many options for $\left(  a,b,c\right)  $. We need
to check which of these options actually are solutions of the original system.
This is a lot of work, but a computer can do it.
\end{example}

Of course, there are more elegant ways to solve the above exercise (otherwise,
I would not have posed it). However, the way we just showed is generalizable.
In general, if a system of polynomial equations over $\mathbb{C}$ has only
finitely many solutions, then we can find them all in this way (provided that
we have an algorithm for finding all roots of a univariate
polynomial).\footnote{If a system of polynomial equations has infinitely many
solutions, then this strategy usually will not work. For example, if we try to
use it to solve the system%
\begin{align*}
ab  &  =0,\\
bc  &  =0,\\
ca  &  =0,
\end{align*}
then we find the Gr\"{o}bner basis $\left(  xy,yz,xz\right)  $, which doesn't
get us any closer to the solutions. Blame this on the problem, not on the
Gr\"{o}bner basis: The system has a more complicated combinatorial structure
(its solution set is the union of the three axes in 3D space; there are
infinitely many options for each of $a,b,c$).} Thus, using Gr\"{o}bner bases
with respect to the lex order, we can (often) reduce solving systems of
polynomial equations in multiple variables to solving polynomial equations in
a single variable. (See \cite[\S 5.9]{Laurit09} for more details on this.)

Some things don't look like systems of polynomial equations, but yet boil down
to such systems. Here is an example:

\begin{example}
\label{exa.21hw0.4}Recall Exercise \ref{exe.21hw0.4}:

\begin{statement}
Simplify $\sqrt[3]{2+\sqrt{5}}+\sqrt[3]{2-\sqrt{5}}$.
\end{statement}

There are various ways of solving this using some creativity or lucky ideas.
Let us try to be more methodical here. We set
\[
a=\sqrt[3]{2+\sqrt{5}},\ \ \ \ \ \ \ \ \ \ b=\sqrt[3]{2-\sqrt{5}%
},\ \ \ \ \ \ \ \ \ \ c=\sqrt[3]{2+\sqrt{5}}+\sqrt[3]{2-\sqrt{5}}.
\]
Thus, we want to find a simpler expression for $c$. A good first step would be
to find a polynomial whose root $c$ is (since we would then have a chance of
finding $c$ by root-finding techniques). We see that $a,b,c$ satisfy the
following system of equations:%
\begin{align*}
\left(  a^{3}-2\right)  ^{2}-5  &  =0,\\
\left(  b^{3}-2\right)  ^{2}-5  &  =0,\\
a+b-c  &  =0.
\end{align*}
(Indeed, the first equation comes from \textquotedblleft
unraveling\textquotedblright\ $a=\sqrt[3]{2+\sqrt{5}}$, and likewise for the
second; the third comes from the obvious fact that $c=a+b$.)

We try to solve this system using Gr\"{o}bner bases. Thus, we consider the
ideal%
\[
I:=\left(  \left(  x^{3}-2\right)  ^{2}-5\right)  P+\left(  \left(
y^{3}-2\right)  ^{2}-5\right)  P+\left(  x+y-z\right)  P
\]
of the polynomial ring $P=\mathbb{Q}\left[  x,y,z\right]  $. Using SageMath,
we can easily find a Gr\"{o}bner basis of this ideal $I$ with respect to the
lex order. Its last entry is a polynomial that involves only the variable $z$,
so we can narrow down the options for $c$ to the roots of this polynomial.

This looks nice in theory, but in practice you will realize that this last
entry is%
\[
z^{21}-40z^{18}+218z^{15}-72z^{12}-9931z^{9}-5216z^{6}+19136z^{3}-4096.
\]
Eek. With a good computer algebra system, you can factor this polynomial, but
there will be some degree-$4$ factors irreducible over $\mathbb{Q}$. The
polynomial has $5$ real roots, so $c$ must be one of them, but we need some
harder work to find out which one. This is all not very convenient.

But our approach can be salvaged. We have been \textquotedblleft throwing
away\textquotedblright\ information about our $a,b,c$; no wonder that we got
so many options for $c$. Indeed, the equation $\left(  a^{3}-2\right)
^{2}-5=0$ doesn't really mean $a=\sqrt[3]{2+\sqrt{5}}$; it only means that $a$
is \textbf{some} cube root of ($2$ plus \textbf{some} square root of $5$).
Here, we are using the word \textquotedblleft root\textquotedblright\ in the
wider sense, so a nonzero complex number has two square roots and three cube
roots; thus, there are $6$ possibilities in total for $a$. Likewise for $b$.
Our system of equations above allows $c$ to be a sum of any of the $6$
possible $a$'s with any of the $6$ possible $b$'s. Unsurprisingly, this leaves
lots of different options for $c$.

Thus, we need to integrate a bit more information about the actual values of
$a,b$ into our system. Of course, we know that $a$ is the \textbf{real} cube
root of $2+\sqrt{5}$, where $\sqrt{5}$ is the \textbf{positive} square root of
$5$. But this is not the kind of information we can easily integrate into a
system of equations.

However, we can observe that
\begin{align*}
ab  &  =\sqrt[3]{2+\sqrt{5}}\cdot\sqrt[3]{2-\sqrt{5}}=\sqrt[3]{\left(
2+\sqrt{5}\right)  \cdot\left(  2-\sqrt{5}\right)  }\\
&  \ \ \ \ \ \ \ \ \ \ \left(  \text{since }\sqrt[3]{u}\cdot\sqrt[3]%
{v}=\sqrt[3]{uv}\text{ for any }u,v\in\mathbb{R}\right) \\
&  =\sqrt[3]{-1}=-1.
\end{align*}
Thus, we can extend our system to%
\begin{align*}
\left(  a^{3}-2\right)  ^{2}-5  &  =0,\\
\left(  b^{3}-2\right)  ^{2}-5  &  =0,\\
a+b-c  &  =0,\\
ab+1  &  =0.
\end{align*}
This is a different system and has a smaller set of solutions than the
previous one, but that's good news, since the solution we are looking for is
one of its solutions.

Now, solving this new system using the Gr\"{o}bner basis technique, we find
that $c$ is a root of the polynomial $z^{3}+3z-4$ (since this polynomial is
the last entry of the Gr\"{o}bner basis we find). But the roots of this
polynomial are easy to find: The factorization%
\[
z^{3}+3z-4=\left(  z-1\right)  \underbrace{\left(  z^{2}+z+4\right)
}_{\text{always positive on }\mathbb{R}}%
\]
shows that its only real root is $1$, so that $c$ must be $1$ (since $c$ is
real by definition). Thus our exercise is solved.
\end{example}

\begin{exercise}
Given nonzero constants $a,b,c\in\mathbb{R}$, solve the system $\left\{
\begin{array}
[c]{l}%
x^{2}-yz=a;\\
y^{2}-zx=b;\\
z^{2}-xy=c
\end{array}
\right.  $ in real unknowns $x,y,z$. (The choice of field actually does not
matter here; we could just as well use $\mathbb{C}$ or $\mathbb{Q}$.) \medskip

[\textbf{Hint:} Compute a Gr\"{o}bner basis of the ideal generated by
$x^{2}-yz-a$ etc.]
\end{exercise}

See \cite{CoLiOs15} for more about solving systems of polynomial equations,
and for further applications of Gr\"{o}bner bases.

\subsection{\label{sec.polys2.factor}Factorization of polynomials}

In Example \ref{exa.21hw0.4}, we used a computer to factor a polynomial. Let
me say some words about the algorithms that are used for this (or, at least,
about an algorithm that could theoretically be used for this, but is too slow
in practice; actual computers use faster algorithms). (This was also asked in
Exercise \ref{exe.21hw0.3}.)

\subsubsection{Factoring univariate polynomials}

Let $F$ be a field.

Recall that the ring $F\left[  x\right]  $ is a UFD; thus, each nonzero
polynomial in $F\left[  x\right]  $ has an essentially unique factorization
into irreducible polynomials. (\textquotedblleft Essentially\textquotedblright%
\ means \textquotedblleft up to order and up to associates\textquotedblright.
Keep in mind that the units of $F\left[  x\right]  $ are precisely the nonzero
constant polynomials, so that two polynomials $f,g\in F\left[  x\right]  $ are
associate if and only if there exists some $\lambda\in F\setminus\left\{
0\right\}  $ satisfying $g=\lambda f$.)

How do we find this factorization (into irreducible polynomials)?

When $F$ is finite, we can just check all possibilities by brute force.
Indeed, any factor in the factorization of a nonzero polynomial $f$ must be a
polynomial of degree $\leq\deg f$, and this leaves finitely many options for
it when $F$ is finite.

For general fields $F$, there is no algorithm that finds the factorization of
every polynomial.\footnote{See \url{https://mathoverflow.net/a/350877/} for an
outline of the proof.} But what about well-known fields like $\mathbb{Q}$,
$\mathbb{R}$ and $\mathbb{C}$ ?

Over $\mathbb{R}$ and $\mathbb{C}$ you cannot \textquotedblleft
really\textquotedblright\ factor polynomials, because this is not a
numerically stable problem. For example, the polynomial $x^{2}-2x+1$ factors
over $\mathbb{R}$ (as $\left(  x-1\right)  ^{2}$), but $x^{2}-1.999x+1$ does
not (nontrivially at least). Approximate algorithms that work for sufficiently
non-singular inputs exist, but this is more a question of numerical analysis
than of algebra.

What about polynomials over $\mathbb{Q}$ ? There is an algorithm, whose main
ingredient is the following fact:

\begin{proposition}
[Gauss's lemma in one of its forms]\label{prop.polfactor.gausslemZ1}Let
$f\in\mathbb{Z}\left[  x\right]  $ be a non-constant polynomial. If $f$ is
irreducible in $\mathbb{Z}\left[  x\right]  $, then $f$ is irreducible in
$\mathbb{Q}\left[  x\right]  $.
\end{proposition}

\begin{proof}
Assume the contrary. Thus, $f=gh$ for some nonconstant polynomials
$g,h\in\mathbb{Q}\left[  x\right]  $ (since the units of $\mathbb{Q}\left[
x\right]  $ are precisely the nonzero constant polynomials). By multiplying
the two polynomials $g$ and $h$ with the lowest common denominators of their
coefficients, we obtain two nonconstant polynomials $u$ and $v$ in
$\mathbb{Z}\left[  x\right]  $. These two polynomials $u$ and $v$ satisfy
$uv=Ngh$ for some positive integer $N$ (since $u$ and $v$ are positive integer
multiples of $g$ and $h$). Consider this $N$. We have $uv=N\underbrace{gh}%
_{=f}=Nf$, so that $Nf=uv$.

Thus, we have found two nonconstant polynomials $u,v\in\mathbb{Z}\left[
x\right]  $ and a positive integer $N$ such that
\begin{equation}
Nf=uv. \label{pf.prop.polfactor.gausslemZ1.Nf=uv}%
\end{equation}

We WLOG assume that $N$ is \textbf{minimal} with the property that such $u,v$
exist. (In other words, among all triples $\left(  u,v,N\right)  $ of two
nonconstant polynomials $u,v\in\mathbb{Z}\left[  x\right]  $ and a positive
integer $N$ satisfying (\ref{pf.prop.polfactor.gausslemZ1.Nf=uv}), we pick one
in which $N$ is minimal. This might not be the one that we obtained from $g$
and $h$ above.)

If $N=1$, then (\ref{pf.prop.polfactor.gausslemZ1.Nf=uv}) rewrites as $f=uv$,
which contradicts the assumption that $f$ is irreducible (since $u$ and $v$
are nonconstant and thus non-units). Hence, we cannot have $N=1$. Thus, there
exists a prime $p$ that divides $N$. Consider such a $p$. Recall that
$\mathbb{Z}/p$ is a field (since $p$ is prime). Therefore, $\mathbb{Z}/p$ is
an integral domain, so that the polynomial ring $\left(  \mathbb{Z}/p\right)
\left[  x\right]  $ is an integral domain as well (by Corollary
\ref{cor.polring.univar-intdom}).

We shall now show a way to turn any polynomial $s\in\mathbb{Z}\left[
x\right]  $ into a polynomial $\overline{s}\in\left(  \mathbb{Z}/p\right)
\left[  x\right]  $. It is as simple as you can imagine: We simply replace
each coefficient by its residue class modulo $p$. In other words, if
$s=s_{0}x^{0}+s_{1}x^{1}+\cdots+s_{n}x^{n}$ is a polynomial in $\mathbb{Z}%
\left[  x\right]  $ (with $s_{i}\in\mathbb{Z}$), then we define a polynomial
$\overline{s}:=\overline{s_{0}}x^{0}+\overline{s_{1}}x^{1}+\cdots
+\overline{s_{n}}x^{n}\in\left(  \mathbb{Z}/p\right)  \left[  x\right]  $
(where $\overline{s_{i}}$ means the residue class of $s_{i}$ modulo $p$). For
example, if $p=5$, then $\overline{2x^{3}+7}=\overline{2}x^{3}+\overline
{7}=\overline{2}x^{3}+\overline{2}$. It is easy to see that the map%
\begin{align*}
\mathbb{Z}\left[  x\right]   &  \rightarrow\left(  \mathbb{Z}/p\right)
\left[  x\right]  ,\\
s  &  \mapsto\overline{s}%
\end{align*}
is a ring morphism (since the rules for adding and multiplying polynomials are
the same over $\mathbb{Z}$ and over $\mathbb{Z}/p$). Thus, $\overline
{uv}=\overline{u}\cdot\overline{v}$.

Now, $f\in\mathbb{Z}\left[  x\right]  $; hence, all coefficients of the
polynomial $Nf$ are divisible by $N$, and thus also divisible by $p$ (since
$p$ divides $N$). Thus, $\overline{Nf}=0$ in $\left(  \mathbb{Z}/p\right)
\left[  x\right]  $. However, (\ref{pf.prop.polfactor.gausslemZ1.Nf=uv})
entails $\overline{Nf}=\overline{uv}=\overline{u}\cdot\overline{v}$. Thus,
$\overline{u}\cdot\overline{v}=\overline{Nf}=0$. Since $\left(  \mathbb{Z}%
/p\right)  \left[  x\right]  $ is an integral domain, this shows that
$\overline{u}=0$ or $\overline{v}=0$. We WLOG assume that $\overline{u}=0$
(since otherwise, we can simply swap $u$ with $v$).

Now, $\overline{u}=0$ means that all coefficients of $u$ are multiples of $p$.
In other words, $\dfrac{1}{p}u\in\mathbb{Z}\left[  x\right]  $. Now, the
equality (\ref{pf.prop.polfactor.gausslemZ1.Nf=uv}) yields%
\[
\dfrac{N}{p}f=\left(  \dfrac{1}{p}u\right)  v.
\]
Since $\dfrac{N}{p}$ is a positive integer (because $p$ divides $N$) and since
$\dfrac{1}{p}u\in\mathbb{Z}\left[  x\right]  $, this equality shows that
$\left(  \dfrac{1}{p}u,v,\dfrac{N}{p}\right)  $ is a triple of two nonconstant
polynomials $\dfrac{1}{p}u,v\in\mathbb{Z}\left[  x\right]  $ and a positive
integer $\dfrac{N}{p}$ satisfying (\ref{pf.prop.polfactor.gausslemZ1.Nf=uv})
(with $u$ and $N$ replaced by $\dfrac{1}{p}u$ and $\dfrac{N}{p}$). But recall
that among all such triples, we chose $\left(  u,v,N\right)  $ to be one with
minimal $N$. Thus, $N\leq\dfrac{N}{p}$. Therefore, $p\leq1$ (since $N$ is a
positive integer). This contradicts the assumption that $p$ is prime. This
contradiction completes the proof.
\end{proof}

Let us now address two computational problems for polynomials with integer or
rational coefficients.

\begin{statement}
\textbf{Problem 1:} Let $f,g\in\mathbb{Z}\left[  x\right]  $ be two
polynomials with $g\neq0$. Check whether $g$ divides $f$ in $\mathbb{Z}\left[
x\right]  $.
\end{statement}

\begin{proof}
[Solution (sketched).]The leading coefficient of $g$ may or may not be a unit
of $\mathbb{Z}$; however, it is always a unit of $\mathbb{Q}$. Thus, we can
use division with remainder to check whether $g$ divides $f$ in the (larger)
ring $\mathbb{Q}\left[  x\right]  $. If the answer is \textquotedblleft
no\textquotedblright, then (a fortiori) $g$ cannot divide $f$ in
$\mathbb{Z}\left[  x\right]  $ (since $\mathbb{Z}\left[  x\right]  $ is a
subring of $\mathbb{Q}\left[  x\right]  $). If the answer is \textquotedblleft
yes\textquotedblright, then we compute the quotient $\dfrac{f}{g}\in
\mathbb{Q}\left[  x\right]  $ and check whether it belongs to $\mathbb{Z}%
\left[  x\right]  $ (that is, whether its coefficients are integers). If yes,
then the answer is \textquotedblleft yes\textquotedblright; if no, then the
answer is \textquotedblleft no\textquotedblright. Problem 1 is thus solved.
\end{proof}

\begin{statement}
\textbf{Problem 2:} Let $f\in\mathbb{Z}\left[  x\right]  $ be a nonzero
polynomial. Construct a list of all divisors of $f$ in $\mathbb{Z}\left[
x\right]  $.
\end{statement}

\begin{proof}
[Solution (sketched).]Let $n=\deg f$. Pick $n+1$ distinct integers
$a_{0},a_{1},\ldots,a_{n}$ that are \textbf{not} roots of $f$. (Such $n+1$
integers can always be found, since $f$ is a nonzero polynomial of degree $n$
and thus has at most $n$ roots in the integral domain $\mathbb{Z}$. Thus, for
example, among the $2n+1$ numbers $-n,-n+1,\ldots,n$, at least $n+1$ many are
not roots of $f$.)

For each $i\in\left\{  0,1,\ldots,n\right\}  $, let $D_{i}$ be the set of all
divisors of the integer $f\left[  a_{i}\right]  $ (recall that $f\left[
a\right]  $ is our notation for the evaluation of $f$ at $a$; this is commonly
denoted $f\left(  a\right)  $). This set $D_{i}$ is finite (since $f\left[
a_{i}\right]  \neq0$), and its elements can be explicitly listed. Hence, the
set $D_{0}\times D_{1}\times\cdots\times D_{n}$ is finite as well, and its
elements can be explicitly listed.

Now, let $g$ be a divisor of $f$ in $\mathbb{Z}\left[  x\right]  $. Then,
$g\in\mathbb{Z}\left[  x\right]  $, and there exists a further polynomial
$h\in\mathbb{Z}\left[  x\right]  $ such that $f=gh$. Consider this $h$. From
$f=gh$, we obtain $\deg f=\deg\left(  gh\right)  =\deg g+\underbrace{\deg
h}_{\geq0}\geq\deg g$, so that $\deg g\leq\deg f=n$. In other words, the
polynomial $g$ must have degree $\leq n$.

For each $i\in\left\{  0,1,\ldots,n\right\}  $, we have
\[
f\left[  a_{i}\right]  =g\left[  a_{i}\right]  h\left[  a_{i}\right]
\ \ \ \ \ \ \ \ \ \ \left(  \text{since }f=gh\right)
\]
and thus $g\left[  a_{i}\right]  \mid f\left[  a_{i}\right]  $, so that
$g\left[  a_{i}\right]  \in D_{i}$. Hence,%
\[
\left(  g\left[  a_{0}\right]  ,\ g\left[  a_{1}\right]  ,\ \ldots,\ g\left[
a_{n}\right]  \right)  \in D_{0}\times D_{1}\times\cdots\times D_{n}.
\]
Thus, for each divisor $g$ of $f$ in $\mathbb{Z}\left[  x\right]  $, we know
that the $\left(  n+1\right)  $-tuple \newline$\left(  g\left[  a_{0}\right]
,\ g\left[  a_{1}\right]  ,\ \ldots,\ g\left[  a_{n}\right]  \right)  $
belongs to the finite set $D_{0}\times D_{1}\times\cdots\times D_{n}$ (which
does not depend on $g$ and can be explicitly found). Hence, we have finitely
many options for this $\left(  n+1\right)  $-tuple.

However, given the $\left(  n+1\right)  $-tuple $\left(  g\left[
a_{0}\right]  ,\ g\left[  a_{1}\right]  ,\ \ldots,\ g\left[  a_{n}\right]
\right)  $, we can uniquely reconstruct the polynomial $g$. (Indeed, we know
that $g$ has degree $\leq n$, so that Corollary
\ref{cor.polint.lagrange-as-eq} (applied to $g$ instead of $f$) yields
$g=\sum_{j=0}^{n}g\left[  a_{j}\right]  \cdot\dfrac{\prod_{k\neq j}\left(
x-a_{k}\right)  }{\prod_{k\neq j}\left(  a_{j}-a_{k}\right)  }$. This is an
explicit formula for $g$ in terms of the $\left(  n+1\right)  $-tuple
\newline$\left(  g\left[  a_{0}\right]  ,\ g\left[  a_{1}\right]
,\ \ldots,\ g\left[  a_{n}\right]  \right)  $). Therefore, given the $\left(
n+1\right)  $-tuple \newline$\left(  g\left[  a_{0}\right]  ,\ g\left[
a_{1}\right]  ,\ \ldots,\ g\left[  a_{n}\right]  \right)  $ for a divisor $g$
of $f$, we can uniquely reconstruct $g$.)

Thus, we have finitely many options for $g$ (since we have finitely many
options for this $\left(  n+1\right)  $-tuple). Usually, few of these options
will actually produce a polynomial $g\in\mathbb{Z}\left[  x\right]  $ that
divides $f$ (indeed, many of them will produce polynomials with non-integer
coefficients; and even among the polynomials that do have integer
coefficients, many will fail to divide $f$). However, we can check which of
these options do produce a polynomial $g\in\mathbb{Z}\left[  x\right]  $ that
divides $f$ (our above solution to Problem 1 helps here). Thus, we end up with
a list of all divisors of $f$ in $\mathbb{Z}\left[  x\right]  $. This solves
Problem 2.
\end{proof}

\begin{statement}
\textbf{Problem 3:} Let $f\in\mathbb{Q}\left[  x\right]  $ be a nonzero
polynomial. Find a factorization of $f$ into a product of irreducible polynomials.
\end{statement}

\begin{proof}
[Solution sketch.]WLOG assume that $f$ is not constant (otherwise, $f$ is a
unit and thus associate to the empty product). Also, WLOG assume that
$f\in\mathbb{Z}\left[  x\right]  $ (otherwise, multiply $f$ with the lowest
common denominator of its coefficients). Furthermore, WLOG assume that the gcd
of the coefficients of $f$ is $1$ (otherwise, divide $f$ by this gcd). We find
a list of all divisors of $f$ in $\mathbb{Z}\left[  x\right]  $ (using the
solution to Problem 2). If the only such divisors are $\pm1$ and $\pm f$, then
$f$ is irreducible in $\mathbb{Z}\left[  x\right]  $ and thus also irreducible
in $\mathbb{Q}\left[  x\right]  $ (by Proposition
\ref{prop.polfactor.gausslemZ1}), so we are done. Else, we find a divisor $g$
of $f$ that is neither $\pm1$ nor $\pm f$, and thus we can decompose $f$ as a
product $gh$ of two nonconstant polynomials $g,h\in\mathbb{Z}\left[  x\right]
$. In that case, we have reduced the problem to the same problem with the
(lower-degree) polynomials $g$ and $h$. Thus, recursively iterating the
procedure, we end up with a factorization of $f$ into a product of irreducible polynomials.
\end{proof}

Our solution to Problem 3 is a theoretical algorithm for factoring a
polynomial in $\mathbb{Q}\left[  x\right]  $ into irreducible polynomials. The
algorithm is too computationally intensive to be viable in practice, so
computers use different methods (often using $\mathbb{Z}/p$ as a stand-in for
$\mathbb{Z}$ and using the Chinese Remainder Theorem to \textquotedblleft
glue\textquotedblright\ the factorizations over different $\mathbb{Z}/p$'s
together). See \cite[Chapters 14--16]{GatGer13} for a serious introduction to
these methods.

\subsubsection{Factoring multivariate polynomials}

Factoring multivariate polynomials over $\mathbb{Q}$ can be done similarly
using multivariate Lagrange interpolation\footnote{See Exercise
\ref{exe.polint.lagrange.2var} for Lagrange interpolation in the case of $2$
variables. The case of $n$ variables is conceptually similar, though there are
many more subscripts to deal with.}. (The word \textquotedblleft
similarly\textquotedblright\ is doing some heavy duty here.) Alternatively, it
can be reduced to the univariate case by the following trick: If
$f\in\mathbb{Q}\left[  x,y\right]  $ is a polynomial of degree $<N$ (for some
$N\in\mathbb{N}$), then the univariate polynomial $f\left(  x,x^{N}\right)  $
\textquotedblleft carries all the information of $f$\textquotedblright\ (in
the sense that no two different terms of $f$ get merged when we substitute
$x^{N}$ for $y$). For example, if $f=x^{2}+xy+y^{2}$ and $N=5$, then%
\[
f\left(  x,x^{N}\right)  =f\left(  x,x^{5}\right)  =x^{2}+xx^{5}+\left(
x^{5}\right)  ^{2}=x^{2}+x^{6}+x^{10}.
\]
Thus, in order to factor $f$, it suffices to factor $f\left(  x,x^{N}\right)
$ (a univariate polynomial), and then try to lift the factorization back by
\textquotedblleft substituting $y$ for $x^{N}$\textquotedblright. This trick
is justified by the following exercise:

\begin{exercise}
\label{exe.21hw4.9}Let $R$ be a commutative ring. Let $P$ be the polynomial
ring $R\left[  x,y\right]  $.

Fix $N\in\mathbb{N}$. Let $P_{N}$ be the $R$-submodule
\[
\left\{  f\in P\ \mid\ f=0\text{ or }\deg f<N\right\}
\]
of $P$. (This is an $R$-submodule, since it is the span of the family $\left(
x^{i}y^{j}\right)  _{\left(  i,j\right)  \in\mathbb{N}^{2};\ i+j<N}$.)

\begin{enumerate}
\item[\textbf{(a)}] Consider the $R$-algebra morphism
\begin{align*}
S:P  &  \rightarrow R\left[  x\right]  ,\\
f  &  \mapsto f\left(  x,x^{N}\right)  .
\end{align*}
(This is the map that substitutes $x^{N}$ for $y$ in any polynomial $f\in P$.
It is an $R$-algebra morphism, as we know from Theorem
\ref{thm.polring.mulvar-sub-hom}.)

Prove that the restriction of $S$ to $P_{N}$ is injective.

\item[\textbf{(b)}] Assume that $R$ is a field. Let $f\in P_{N}$ be such that
the polynomial $S\left(  f\right)  \in R\left[  x\right]  $ is irreducible.
Show that $f\in P=R\left[  x,y\right]  $ is irreducible.
\end{enumerate}

[\textbf{Remark:} The converse of part \textbf{(b)} does not hold. For
example, if $R=\mathbb{Q}$ and $N=2$, then the polynomial $f:=1+2x+y\in P$ is
irreducible, but the polynomial $S\left(  f\right)  =1+2x+x^{2}=\left(
1+x\right)  ^{2}\in R\left[  x\right]  $ is not.]
\end{exercise}

This trick (of substituting $x^{N}$ for $y$) is easily generalized to
polynomials in more than two variables. For example, a polynomial in $4$
variables $x,y,z,w$ can be transformed into a polynomial in $3$ variables
$x,y,z$ by substituting $z^{N}$ for $w$.

\subsection{\label{sec.cyclotomic}Cyclotomic polynomials}

\subsubsection{Definition}

We shall now study a sequence of special polynomials in the univariate
polynomial ring $\mathbb{Z}\left[  x\right]  $ that are highly important in
algebra and number theory\footnote{We will see some of their applications
later in this section.}: the \textbf{cyclotomic polynomials}. Their definition
relies on the following theorem, which will be proved later in this same section:

\begin{theorem}
\label{thm.cyclotomic.EC}There is a unique sequence $\left(  \Phi_{1},\Phi
_{2},\Phi_{3},\ldots\right)  $ of polynomials in $\mathbb{Z}\left[  x\right]
$ such that each $n\geq1$ satisfies%
\begin{equation}
x^{n}-1=\prod_{d\mid n}\Phi_{d}\ \ \ \ \ \ \ \ \ \ \text{in }\mathbb{Z}\left[
x\right]  . \label{eq.thm.cyclotomic.EC.eq}%
\end{equation}
Here and in the following, the symbol \textquotedblleft$\prod_{d\mid n}%
$\textquotedblright\ means a product over all positive divisors of $n$ (so
that, for example, $\prod_{d\mid6}\Phi_{d}=\Phi_{1}\Phi_{2}\Phi_{3}\Phi_{6}$).
\end{theorem}

We will prove Theorem \ref{thm.cyclotomic.EC} by showing that the sequence
\newline$\left(  x^{1}-1,\ x^{2}-1,\ x^{3}-1,\ \ldots\right)  $ is a strong
gcd-sequence (see Definition \ref{def.gcdseq.gcdseq}); this will yield the
existence part of Theorem \ref{thm.cyclotomic.EC} using Theorem
\ref{thm.gcdseq.fac}, while the uniqueness part will easily follow from
Proposition \ref{prop.gcdseq.reg}. First, however, let us give these
polynomials a name, compute some examples and state a few more of their properties.

\begin{definition}
The polynomials $\Phi_{n}$ that constitute the sequence $\left(  \Phi_{1}%
,\Phi_{2},\Phi_{3},\ldots\right)  $ in Theorem \ref{thm.cyclotomic.EC} are
known as the \textbf{cyclotomic polynomials}. Specifically, each $\Phi_{n}$ is
called the $n$\textbf{-th cyclotomic polynomial}.
\end{definition}

The equality (\ref{eq.thm.cyclotomic.EC.eq}) can be rewritten as\footnote{The
product $\prod_{\substack{d\mid n;\\d\neq n}}\Phi_{d}$ is empty when $n=1$; as
usual, an empty product should be understood as $1$.}%
\begin{equation}
x^{n}-1=\prod_{d\mid n}\Phi_{d}=\Phi_{n}\prod_{\substack{d\mid n;\\d\neq
n}}\Phi_{d} \label{eq.thm.cyclotomic.EC.eq-n}%
\end{equation}
(here, we have split off the factor for $d=n$ from the product); this allows
us to compute $\Phi_{n}$ recursively by polynomial division using the
\textquotedblleft smaller\textquotedblright\ cyclotomic polynomials $\Phi
_{1},\Phi_{2},\ldots,\Phi_{n-1}$ (since all the factors $\Phi_{d}$ of the
product $\prod_{\substack{d\mid n;\\d\neq n}}\Phi_{d}$ have $d<n$). Thus, we
obtain the following table of the first $16$ cyclotomic polynomials:%
\begin{align*}
\Phi_{1}  &  =x-1;\\
\Phi_{2}  &  =x+1;\\
\Phi_{3}  &  =x^{2}+x+1;\\
\Phi_{4}  &  =x^{2}+1;\\
\Phi_{5}  &  =x^{4}+x^{3}+x^{2}+x+1;\\
\Phi_{6}  &  =x^{2}-x+1;\\
\Phi_{7}  &  =x^{6}+x^{5}+x^{4}+x^{3}+x^{2}+x+1;\\
\Phi_{8}  &  =x^{4}+1;\\
\Phi_{9}  &  =x^{6}+x^{3}+1;\\
\Phi_{10}  &  =x^{4}-x^{3}+x^{2}-x+1;\\
\Phi_{11}  &  =x^{10}+x^{9}+x^{8}+x^{7}+x^{6}+x^{5}+x^{4}+x^{3}+x^{2}+x+1;\\
\Phi_{12}  &  =x^{4}-x^{2}+1;\\
\Phi_{13}  &  =x^{12}+x^{11}+x^{10}+x^{9}+x^{8}+x^{7}+x^{6}+x^{5}+x^{4}%
+x^{3}+x^{2}+x+1;\\
\Phi_{14}  &  =x^{6}-x^{5}+x^{4}-x^{3}+x^{2}-x+1;\\
\Phi_{15}  &  =x^{8}-x^{7}+x^{5}-x^{4}+x^{3}-x+1;\\
\Phi_{16}  &  =x^{8}+1.
\end{align*}
These polynomials exhibit some rather peculiar behavior. For one, their
degrees are not increasing; they fluctuate (while still going to $\infty$). It
might appear that their coefficients are always $1$, $0$ and $-1$; but this is
not the case, as the example of
\begin{align*}
\Phi_{105}  &  =x^{48}+x^{47}+x^{46}-x^{43}-x^{42}-2x^{41}-x^{40}%
-x^{39}+x^{36}+x^{35}+x^{34}+x^{33}+x^{32}\\
&  \ \ \ \ \ \ \ \ \ \ +x^{31}-x^{28}-x^{26}-x^{24}-x^{22}-x^{20}%
+x^{17}+x^{16}+x^{15}+x^{14}+x^{13}+x^{12}\\
&  \ \ \ \ \ \ \ \ \ \ -x^{9}-x^{8}-2x^{7}-x^{6}-x^{5}+x^{2}+x+1
\end{align*}
shows (this is the smallest $n$ for which $\Phi_{n}$ has a coefficient
distinct from $1$, $0$ and $-1$). Despite these vagaries, an attentive reader
will likely have spotted a few actual patterns, such as the following:

\begin{proposition}
\label{prop.cyclotomic.Phip}For any prime number $p$, we have%
\[
\Phi_{p}=x^{0}+x^{1}+\cdots+x^{p-1}.
\]

\end{proposition}

\begin{proof}
Let $p$ be a prime number. Then, (\ref{eq.thm.cyclotomic.EC.eq}) (applied to
$n=p$) yields
\[
x^{p}-1=\prod_{d\mid p}\Phi_{d}=\Phi_{1}\Phi_{p}%
\]
(since the only positive divisors of $p$ are $1$ and $p$). But $\Phi_{1}=x-1$
(this follows easily by applying (\ref{eq.thm.cyclotomic.EC.eq}) to $n=1$).
Hence, $\Phi_{1}\Phi_{p}=\left(  x-1\right)  \Phi_{p}$, so that%
\[
\left(  x-1\right)  \Phi_{p}=\Phi_{1}\Phi_{p}=x^{p}-1=\left(  x-1\right)
\left(  x^{0}+x^{1}+\cdots+x^{p-1}\right)  .
\]
We can cancel $x-1$ from this equality (e.g., because Corollary
\ref{cor.polring.univar-intdom} shows that $\mathbb{Z}\left[  x\right]  $ is
an integral domain, so the nonzero polynomial $x-1$ cannot be a zero
divisor)\footnote{This can also be argued more tastefully using Proposition
\ref{prop.polring.univar-degpq} \textbf{(e)}.}, and thus obtain%
\[
\Phi_{p}=x^{0}+x^{1}+\cdots+x^{p-1}.
\]
This proves Proposition \ref{prop.cyclotomic.Phip}.
\end{proof}

Another simple pattern is observed at $\Phi_{2},\Phi_{4},\Phi_{8},\Phi_{16}$:

\begin{exercise}
Let $m$ be a positive integer. Prove that $\Phi_{2^{m}}=x^{2^{m-1}}+1$.
\end{exercise}

A more-or-less explicit formula for $\Phi_{n}$ can be obtained from Theorem
\ref{thm.gcdseq.fac-pf} \textbf{(b)}:

\begin{theorem}
\label{thm.cyclotomic.frac}For any integer $n\geq1$, we have%
\begin{equation}
\Phi_{n}=\dfrac{\prod\limits_{\substack{S\subseteq\operatorname*{PF}\left(
n\right)  ;\\\left\vert S\right\vert \text{ is even}}}\left(
x^{n/\operatorname*{prod}S}-1\right)  }{\prod\limits_{\substack{S\subseteq
\operatorname*{PF}\left(  n\right)  ;\\\left\vert S\right\vert \text{ is odd}%
}}\left(  x^{n/\operatorname*{prod}S}-1\right)  }.
\label{eq.thm.cyclotomic.frac.eq}%
\end{equation}
Here, the notations $\operatorname*{PF}\left(  n\right)  $ and
$\operatorname*{prod}S$ are defined as in Definition \ref{def.ent.PFn} and
Definition \ref{def.ent.prodS}, and we are using Definition
\ref{def.rings.reg.div} to make sense of the fraction on the right-hand side.
\end{theorem}

\begin{exercise}
\ \ 

\begin{enumerate}
\item[\textbf{(a)}] Prove that every integer $n\geq1$ satisfies%
\[
x^{0}+x^{1}+\cdots+x^{n-1}=\prod_{\substack{d\mid n;\\d\neq1}}\Phi_{d}.
\]

\item[\textbf{(b)}] Prove that every integer $n>1$ satisfies $\Phi_{n}\left(
0\right)  =1$.

\item[\textbf{(c)}] Prove that every integer $n>1$ satisfies%
\[
\Phi_{n}\left(  1\right)  =%
\begin{cases}
p, & \text{if }n\text{ is a power of a prime }p;\\
1, & \text{otherwise.}%
\end{cases}
\]

\end{enumerate}
\end{exercise}

\subsubsection{Cyclotomic polynomials and roots of unity}

Most authors usually don't define cyclotomic polynomials the way we do.
Instead, they define them by the explicit formula%
\begin{equation}
\Phi_{n}=\prod_{\substack{k\in\left\{  1,2,\ldots,n\right\}  ;\\\gcd\left(
k,n\right)  =1}}\left(  x-\zeta^{k}\right)  , \label{eq.thm.cyclotomic.ngon}%
\end{equation}
where $\zeta=e^{2\pi i/n}$ is the \textquotedblleft standard\textquotedblright%
\ primitive $n$-th root of unity in $\mathbb{C}$ (that is, the point on the
unit circle having argument $2\pi/n$; thus its $n$ powers $\zeta^{0},\zeta
^{1},\ldots,\zeta^{n-1}$ form the vertices of a regular $n$-gon). This
definition has the advantage of being explicit, but it relies on the complex
numbers and requires nontrivial work to show that $\Phi_{n}$ is actually a
polynomial in $\mathbb{Z}\left[  x\right]  $ (as opposed to merely in
$\mathbb{C}\left[  x\right]  $). Thus, we are not using it as our definition
of $\Phi_{n}$. However, we shall show that it is equivalent; in fact, we will
prove the following more general result (which we shall prove further below,
in the same section):

\begin{theorem}
\label{thm.cyclotomic.RC}Let $F$ be a field. Consider its group of units
$F^{\times}$ (see Theorem \ref{thm.ringunits.group}). Let $n$ be a positive
integer. Let $u$ be an element of the group $F^{\times}$ that has order $n$
(that is, satisfies $u^{n}=1$ but $u^{k}\neq1$ for all positive integers
$k<n$). Then,%
\[
\Phi_{n}=\prod_{\substack{k\in\left\{  1,2,\ldots,n\right\}  ;\\\gcd\left(
k,n\right)  =1}}\left(  x-u^{k}\right)  \ \ \ \ \ \ \ \ \ \ \text{in }F\left[
x\right]  .
\]

\end{theorem}

An element $u$ satisfying the conditions of Theorem \ref{thm.cyclotomic.RC} is
known as a \textbf{primitive }$n$\textbf{-th root of unity} in $F$. For
example, the primitive $n$-th roots of unity in $\mathbb{C}$ are the complex
numbers $\zeta^{k}$ for $\zeta=e^{2\pi i/n}$ and $k\in\left\{  1,2,\ldots
,n\right\}  $ satisfying $\gcd\left(  k,n\right)  =1$. On the other hand, the
field $\mathbb{R}$ contains the primitive $1$-st root of unity $1$ and the
primitive $2$-nd root of unity $-1$, but no primitive $n$-th roots of unity
for any $n>2$.

\begin{exercise}
Let $p$ be a prime. Let $n$ be a positive integer such that $p\mid n$. Prove
that a field of characteristic $p$ cannot contain any primitive $n$-th root of unity.
\end{exercise}

\subsubsection{The Euler totient function and the degree of $\Phi_{n}$}

One of the most natural questions about the cyclotomic polynomials is to
compute their degrees. There is a simple answer, which however requires a
classical number-theoretical function (or integer sequence):

\begin{definition}
\label{def.ent.phi.phi}The \textbf{Euler totient function} (aka the
\textbf{Euler }$\phi$\textbf{-function}) is the function $\phi:\left\{
1,2,3,\ldots\right\}  \rightarrow\mathbb{N}$ defined by%
\[
\phi\left(  n\right)  =\left(  \text{number of all }k\in\left\{
1,2,\ldots,n\right\}  \text{ satisfying }\gcd\left(  k,n\right)  =1\right)
\]
for all positive integers $n$. It is always denoted by $\phi$ in number theory.
\end{definition}

Here is a table of the early values of $\phi$:%
\[%
\begin{tabular}
[c]{|c||c|c|c|c|c|c|c|c|c|c|c|c|c|}\hline
$n$ & $1$ & $2$ & $3$ & $4$ & $5$ & $6$ & $7$ & $8$ & $9$ & $10$ & $11$ & $12$
& $13$\\\hline
$\phi\left(  n\right)  $ & $1$ & $1$ & $2$ & $2$ & $4$ & $2$ & $6$ & $4$ & $6$
& $4$ & $10$ & $4$ & $12$\\\hline
\end{tabular}
\ \ \ \ \ \ .
\]
Several properties of $\phi$ can be found in places such as \cite[\S 2.14]%
{19s}, \cite[\S 9.1.7]{mps}, \cite[\S 1.7]{Laurit09}, and essentially any
textbook on number theory. We recall only two of them (both proved in
\cite[\S 2.14]{19s}). One is the explicit formula for $\phi\left(  n\right)  $
when the prime factors of $n$ are known:

\begin{theorem}
\label{thm.ent.phi.explicit}Let $n$ be a positive integer.

\begin{enumerate}
\item[\textbf{(a)}] For any prime $p$, let $v_{p}\left(  n\right)  $ denote
the largest $m\in\mathbb{N}$ such that $p^{m}\mid n$. Then,%
\[
\phi\left(  n\right)  =\prod_{\substack{p\text{ prime;}\\p\mid n}}\left(
\left(  p-1\right)  p^{v_{p}\left(  n\right)  -1}\right)  =n\cdot
\prod_{\substack{p\text{ prime;}\\p\mid n}}\left(  1-\dfrac{1}{p}\right)  .
\]

\item[\textbf{(b)}] Use the notations $\operatorname*{PF}\left(  n\right)  $
and $\operatorname*{prod}S$ defined in Definition \ref{def.ent.PFn} and
Definition \ref{def.ent.prodS}. Then,%
\[
\phi\left(  n\right)  =\sum_{S\subseteq\operatorname*{PF}\left(  n\right)
}\left(  -1\right)  ^{\left\vert S\right\vert }\cdot\dfrac{n}%
{\operatorname*{prod}S}.
\]

\end{enumerate}
\end{theorem}

\begin{proof}
\textbf{(a)} This is \cite[Theorem 2.14.5]{19s}. \medskip

\textbf{(b)} This is an intermediate result in \cite[Proof of Theorem
6.2.7]{21s}. But it also is a restatement of part \textbf{(a)}: Indeed, there
is a general identity (a particular case of \cite[Exercise 6.1 \textbf{(a)}%
]{detnotes}) saying that if $\alpha_{1},\alpha_{2},\ldots,\alpha_{k}$ are
finitely many elements of a commutative ring $R$, then%
\[
\left(  1+\alpha_{1}\right)  \left(  1+\alpha_{2}\right)  \cdots\left(
1+\alpha_{k}\right)  =\sum_{S\subseteq\left\{  1,2,\ldots,k\right\}
}\ \ \prod_{p\in S}\alpha_{p}.
\]
Replacing the indexes $1,2,\ldots,k$ here by the elements of an arbitrary
finite set $U$, we can rewrite this as follows: If $U$ is a finite set, and
$\alpha_{p}\in R$ for each $p\in U$ (where $R$ is a commutative ring), then%
\begin{equation}
\prod_{p\in U}\left(  1+\alpha_{p}\right)  =\sum_{S\subseteq U}\ \ \prod_{p\in
S}\alpha_{p}. \label{pf.thm.ent.phi.explicit.b.gen}%
\end{equation}

Now, by the definition of $\operatorname*{PF}\left(  n\right)  $, we have%
\begin{align*}
\prod_{\substack{p\text{ prime;}\\p\mid n}}\left(  1-\dfrac{1}{p}\right)   &
=\prod_{p\in\operatorname*{PF}\left(  n\right)  }\underbrace{\left(
1-\dfrac{1}{p}\right)  }_{=1+\dfrac{-1}{p}}=\prod_{p\in\operatorname*{PF}%
\left(  n\right)  }\left(  1+\dfrac{-1}{p}\right) \\
&  =\sum_{S\subseteq\operatorname*{PF}\left(  n\right)  }\ \ \underbrace{\prod
_{p\in S}\dfrac{-1}{p}}_{\substack{=\dfrac{\left(  -1\right)  ^{\left\vert
S\right\vert }}{\prod_{p\in S}p}=\dfrac{\left(  -1\right)  ^{\left\vert
S\right\vert }}{\operatorname*{prod}S}\\\text{(since }\prod_{p\in
S}p=\operatorname*{prod}S\text{)}}}\ \ \ \ \ \ \ \ \ \ \left(
\begin{array}
[c]{c}%
\text{by (\ref{pf.thm.ent.phi.explicit.b.gen}), applied to }R=\mathbb{Q}\\
\text{and }U=\operatorname*{PF}\left(  n\right)  \text{ and }\alpha_{p}%
=\dfrac{-1}{p}%
\end{array}
\right) \\
&  =\sum_{S\subseteq\operatorname*{PF}\left(  n\right)  }\dfrac{\left(
-1\right)  ^{\left\vert S\right\vert }}{\operatorname*{prod}S}.
\end{align*}
But part \textbf{(a)} yields%
\[
\phi\left(  n\right)  =n\cdot\underbrace{\prod_{\substack{p\text{
prime;}\\p\mid n}}\left(  1-\dfrac{1}{p}\right)  }_{=\sum_{S\subseteq
\operatorname*{PF}\left(  n\right)  }\dfrac{\left(  -1\right)  ^{\left\vert
S\right\vert }}{\operatorname*{prod}S}}=n\cdot\sum_{S\subseteq
\operatorname*{PF}\left(  n\right)  }\dfrac{\left(  -1\right)  ^{\left\vert
S\right\vert }}{\operatorname*{prod}S}=\sum_{S\subseteq\operatorname*{PF}%
\left(  n\right)  }\left(  -1\right)  ^{\left\vert S\right\vert }\cdot
\dfrac{n}{\operatorname*{prod}S}.
\]
This proves part \textbf{(b)}.
\end{proof}

Another is a beautiful identity due to Gauss, which we will use.

\begin{theorem}
\label{thm.ent.phi.sum-div}Let $n$ be a positive integer. Then,%
\[
\sum_{d\mid n}\phi\left(  d\right)  =n.
\]
Here and in the following, the symbol \textquotedblleft$\sum_{d\mid n}%
$\textquotedblright\ stands for \textquotedblleft sum over all positive
divisors $d$ of $n$\textquotedblright\ (so that, for example, $\sum_{d\mid
6}\phi\left(  d\right)  =\phi\left(  1\right)  +\phi\left(  2\right)
+\phi\left(  3\right)  +\phi\left(  6\right)  =1+1+2+2$).
\end{theorem}

We do recall the proof of this theorem, since we will need all of its ideas
more than once. It rests on two basic lemmas:

\begin{lemma}
\label{lem.ent.phi.phi.sum-div.1}Let $n$ be a positive integer. Then, the map
\begin{align}
\left\{  \text{positive divisors of }n\right\}   &  \rightarrow\left\{
\text{positive divisors of }n\right\}  ,\nonumber\\
d  &  \mapsto n/d \label{eq.lem.ent.phi.phi.sum-div.1.bij}%
\end{align}
is a bijection.
\end{lemma}

\begin{fineprint}

\begin{proof}
This map is inverse to itself (since $n/\left(  n/d\right)  =d$ for each
positive divisor $d$ of $n$). Thus, it is a bijection.
\end{proof}
\end{fineprint}

\begin{lemma}
\label{lem.ent.phi.phi.sum-div.2}Let $n$ be a positive integer. Let $d$ be a
positive divisor of $n$. Then, the map
\begin{align}
\left\{  i\in\left\{  1,2,\ldots,\dfrac{n}{d}\right\}  \ \mid\ \gcd\left(
i,\dfrac{n}{d}\right)  =1\right\}   &  \rightarrow\left\{  k\in\left\{
1,2,\ldots,n\right\}  \ \mid\ \gcd\left(  k,n\right)  =d\right\}  ,\nonumber\\
i  &  \mapsto di \label{eq.lem.ent.phi.phi.sum-div.2.bij}%
\end{align}
is well-defined and is a bijection.
\end{lemma}

\begin{fineprint}

\begin{proof}
Clearly, $\dfrac{n}{d}$ is a positive integer (since $d$ is a positive divisor
of $n$).

If $k\in\left\{  1,2,\ldots,n\right\}  $ satisfies $\gcd\left(  k,n\right)
=d$, then $k$ is a multiple of $d$ (since $d=\gcd\left(  k,n\right)  \mid k$)
and thus belongs to $\left\{  1d,2d,\ldots,\dfrac{n}{d}d\right\}  $ (since
$k>0$ and $k\leq n=\dfrac{n}{d}d$); therefore, we have $\dfrac{k}{d}%
\in\left\{  1,2,\ldots,\dfrac{n}{d}\right\}  $. Furthermore, in this
situation, we have $\gcd\left(  \dfrac{k}{d},\dfrac{n}{d}\right)  =1$ (since
$d\cdot\gcd\left(  \dfrac{k}{d},\dfrac{n}{d}\right)  =\gcd\left(
\underbrace{d\cdot\dfrac{k}{d}}_{=k},\underbrace{d\cdot\dfrac{n}{d}}%
_{=n}\right)  =\gcd\left(  k,n\right)  =d$). Thus, the map%
\begin{align*}
\Phi:\left\{  k\in\left\{  1,2,\ldots,n\right\}  \ \mid\ \gcd\left(
k,n\right)  =d\right\}   &  \rightarrow\left\{  i\in\left\{  1,2,\ldots
,\dfrac{n}{d}\right\}  \ \mid\ \gcd\left(  i,\dfrac{n}{d}\right)  =1\right\}
,\\
k  &  \mapsto\dfrac{k}{d}%
\end{align*}
is well-defined.

On the other hand, if $i\in\left\{  1,2,\ldots,\dfrac{n}{d}\right\}  $
satisfies $\gcd\left(  i,\dfrac{n}{d}\right)  =1$, then
\[
di\in\left\{  d\cdot1,d\cdot2,\ldots,d\cdot\dfrac{n}{d}\right\}  =\left\{
1d,2d,\ldots,\dfrac{n}{d}d\right\}  \subseteq\left\{  1,2,\ldots,n\right\}
\]
and $\gcd\left(  di,n\right)  =\gcd\left(  di,d\cdot\dfrac{n}{d}\right)
=d\cdot\underbrace{\gcd\left(  i,\dfrac{n}{d}\right)  }_{=1}=d$. Hence, the
map%
\begin{align*}
\Psi:\left\{  i\in\left\{  1,2,\ldots,\dfrac{n}{d}\right\}  \ \mid
\ \gcd\left(  i,\dfrac{n}{d}\right)  =1\right\}   &  \rightarrow\left\{
k\in\left\{  1,2,\ldots,n\right\}  \ \mid\ \gcd\left(  k,n\right)  =d\right\}
,\\
i  &  \mapsto di
\end{align*}
is well-defined.

The maps $\Phi$ and $\Psi$ are mutually inverse (since $\Phi$ divides its
input by $d$, whereas $\Psi$ multiplies it by $d$); thus, they are bijections.
In particular, $\Psi$ is a bijection. This proves Lemma
\ref{lem.ent.phi.phi.sum-div.2}, since the map described in the lemma is
precisely $\Psi$.
\end{proof}
\end{fineprint}

\begin{fineprint}

\begin{proof}
[Proof of Theorem \ref{thm.ent.phi.sum-div}.]We have%
\begin{align*}
n  &  =\left(  \text{number of all }k\in\left\{  1,2,\ldots,n\right\}  \right)
\\
&  =\sum_{d\mid n}\left(  \text{number of all }k\in\left\{  1,2,\ldots
,n\right\}  \text{ satisfying }\gcd\left(  k,n\right)  =d\right)
\end{align*}
(since $\gcd\left(  k,n\right)  $ is always a positive divisor of $n$).

However, if $d$ is a positive divisor of $n$, then%
\begin{align}
&  \left(  \text{number of all }k\in\left\{  1,2,\ldots,n\right\}  \text{
satisfying }\gcd\left(  k,n\right)  =d\right) \nonumber\\
&  =\left\vert \left\{  k\in\left\{  1,2,\ldots,n\right\}  \ \mid\ \gcd\left(
k,n\right)  =d\right\}  \right\vert \nonumber\\
&  =\left\vert \left\{  i\in\left\{  1,2,\ldots,\dfrac{n}{d}\right\}
\ \mid\ \gcd\left(  i,\dfrac{n}{d}\right)  =1\right\}  \right\vert \nonumber\\
&  \ \ \ \ \ \ \ \ \ \ \ \ \ \ \ \ \ \ \ \ \left(
\begin{array}
[c]{c}%
\text{by the bijection principle, since Lemma \ref{lem.ent.phi.phi.sum-div.2}
shows that}\\
\text{there is a bijection from }\left\{  i\in\left\{  1,2,\ldots,\dfrac{n}%
{d}\right\}  \ \mid\ \gcd\left(  i,\dfrac{n}{d}\right)  =1\right\} \\
\text{to }\left\{  k\in\left\{  1,2,\ldots,n\right\}  \ \mid\ \gcd\left(
k,n\right)  =d\right\}
\end{array}
\right) \nonumber\\
&  =\left(  \text{number of all }i\in\left\{  1,2,\ldots,\dfrac{n}{d}\right\}
\text{ satisfying }\gcd\left(  i,\dfrac{n}{d}\right)  =1\right) \nonumber\\
&  =\phi\left(  \dfrac{n}{d}\right)  \ \ \ \ \ \ \ \ \ \ \left(  \text{by the
definition of }\phi\left(  \dfrac{n}{d}\right)  \right)  .
\label{pf.thm.ent.phi.sum-div.4}%
\end{align}
Thus,
\[
n=\sum_{d\mid n}\underbrace{\left(  \text{number of all }k\in\left\{
1,2,\ldots,n\right\}  \text{ satisfying }\gcd\left(  k,n\right)  =d\right)
}_{\substack{=\phi\left(  \dfrac{n}{d}\right)  \\\text{(by
(\ref{pf.thm.ent.phi.sum-div.4}))}}}=\sum_{d\mid n}\phi\left(  \dfrac{n}%
{d}\right)  =\sum_{d\mid n}\phi\left(  d\right)
\]
(here, we have substituted $d$ for $\dfrac{n}{d}$ in the sum, since the map
(\ref{eq.lem.ent.phi.phi.sum-div.1.bij}) is a bijection). Thus, Theorem
\ref{thm.ent.phi.sum-div} is proved.
\end{proof}
\end{fineprint}

We can now answer the question about the degree of $\Phi_{n}$:

\begin{theorem}
\label{thm.cyclotomic.DC}Let $n$ be a positive integer. Then, $\Phi_{n}$ is a
monic polynomial of degree $\phi\left(  n\right)  $.
\end{theorem}

Again, we shall prove this further below.

\subsubsection{Proofs of Theorems \ref{thm.cyclotomic.EC} and
\ref{thm.cyclotomic.frac}}

We now step to the proofs of the theorems shown above. We begin with the
existence and uniqueness of cyclotomic polynomials, so that we stand on firm
ground when we prove their other properties:

\begin{proof}
[Proof of Theorem \ref{thm.cyclotomic.EC}.]For any $n\geq1$, set%
\begin{equation}
a_{n}:=x^{n}-1\ \ \ \ \ \ \ \ \ \ \text{in }\mathbb{Z}\left[  x\right]  .
\label{pf.thm.cyclotomic.EC.an=}%
\end{equation}
Recall the notion of a strong gcd-sequence (Definition \ref{def.gcdseq.gcdseq}%
). We claim the following:

\begin{statement}
\textit{Claim 1:} The sequence $\left(  a_{1},a_{2},a_{3},\ldots\right)  $ of
elements of $\mathbb{Z}\left[  x\right]  $ is a strong gcd-sequence.
\end{statement}

\begin{proof}
[Proof of Claim 1.]We will obtain this by applying Proposition
\ref{prop.gcdseq.modcrit} to $R=\mathbb{Z}\left[  x\right]  $ and
$b_{m,n}=x^{m}$ and $c_{m,n}=1$. For this to work, we must verify that the
assumptions of Proposition \ref{prop.gcdseq.modcrit} are satisfied -- i.e.,
that all $m,n\geq1$ satisfy%
\begin{equation}
a_{m}\mathbb{Z}\left[  x\right]  +x^{m}\mathbb{Z}\left[  x\right]
+1\mathbb{Z}\left[  x\right]  =\mathbb{Z}\left[  x\right]  ,
\label{pf.thm.cyclotomic.EC.1}%
\end{equation}
and that each $p,q\geq1$ satisfy%
\begin{equation}
a_{p+q}=a_{p}x^{q}+a_{q}1\ \ \ \ \ \ \ \ \ \ \text{in }\mathbb{Z}\left[
x\right]  . \label{pf.thm.cyclotomic.EC.2}%
\end{equation}
But both of these verifications are easy: The equality
(\ref{pf.thm.cyclotomic.EC.1}) follows because the principal ideal
$1\mathbb{Z}\left[  x\right]  $ on the left-hand side is already the whole
ring $\mathbb{Z}\left[  x\right]  $ (and, of course, adding further ideals to
$\mathbb{Z}\left[  x\right]  $ still results in $\mathbb{Z}\left[  x\right]
$), whereas the equality (\ref{pf.thm.cyclotomic.EC.2}) follows by expanding
and canceling terms\footnote{In more detail: From
(\ref{pf.thm.cyclotomic.EC.an=}), we obtain $a_{p+q}=x^{p+q}-1$ and
$a_{p}=x^{p}-1$ and $a_{q}=x^{q}-1$. Hence,%
\[
a_{p}x^{q}+a_{q}1=\left(  x^{p}-1\right)  x^{q}+\left(  x^{q}-1\right)
1=x^{p}x^{q}-x^{q}+x^{q}-1=x^{p}x^{q}-1=x^{p+q}-1=a_{p+q},
\]
so that (\ref{pf.thm.cyclotomic.EC.2}) is proved.}. Thus, Proposition
\ref{prop.gcdseq.modcrit} can be applied, and we conclude that $\left(
a_{1},a_{2},a_{3},\ldots\right)  $ is a strong gcd-sequence. This proves Claim 1.
\end{proof}

\begin{statement}
\textit{Claim 2:} For each $n\geq1$, the polynomial $a_{n}$ is a regular
element of $\mathbb{Z}\left[  x\right]  $.
\end{statement}

\begin{proof}
[Proof of Claim 2.]Let $n\geq1$. The polynomial $a_{n}=x^{n}-1\in
\mathbb{Z}\left[  x\right]  $ is monic. Thus, it is nonzero and its leading
coefficient is $1$, which is a unit. Hence, Proposition
\ref{prop.polring.univar-degpq} \textbf{(e)} (applied to $R=\mathbb{Z}$ and
$p=a_{n}$) shows that if some $q\in\mathbb{Z}\left[  x\right]  $ satisfies
$a_{n}q=0$, then $q=0$. In other words, the element $a_{n}$ of $\mathbb{Z}%
\left[  x\right]  $ is regular. This proves Claim 2.
\end{proof}

Now, we know from Claim 1 that $\left(  a_{1},a_{2},a_{3},\ldots\right)  $ is
a strong gcd-sequence in $\mathbb{Z}\left[  x\right]  $. Hence, Theorem
\ref{thm.gcdseq.fac} (applied to $R=\mathbb{Z}\left[  x\right]  $) shows that
there exists a sequence $\left(  b_{1},b_{2},b_{3},\ldots\right)  $ of
elements of $\mathbb{Z}\left[  x\right]  $ such that each $n\geq1$ satisfies%
\[
a_{n}=\prod_{d\mid n}b_{d}.
\]
Moreover, Proposition \ref{prop.gcdseq.reg} \textbf{(a)} (applied to
$R=\mathbb{Z}\left[  x\right]  $) shows that there is \textbf{at most one}
such sequence $\left(  b_{1},b_{2},b_{3},\ldots\right)  $ (since Claim 2 shows
that $\left(  a_{1},a_{2},a_{3},\ldots\right)  $ is a sequence of regular
elements of $\mathbb{Z}\left[  x\right]  $). Combining these two facts, we
conclude that there is a \textbf{unique} sequence $\left(  b_{1},b_{2}%
,b_{3},\ldots\right)  $ of elements of $\mathbb{Z}\left[  x\right]  $ such
that each $n\geq1$ satisfies%
\[
a_{n}=\prod_{d\mid n}b_{d}.
\]
Renaming the $b_{i}$ as $\Phi_{i}$ in this sentence, and recalling that
$a_{n}=x^{n}-1$, we can rewrite this as follows: There is a unique sequence
$\left(  \Phi_{1},\Phi_{2},\Phi_{3},\ldots\right)  $ of elements of
$\mathbb{Z}\left[  x\right]  $ such that each $n\geq1$ satisfies%
\[
x^{n}-1=\prod_{d\mid n}\Phi_{d}.
\]
This proves Theorem \ref{thm.cyclotomic.EC}.
\end{proof}

Theorem \ref{thm.cyclotomic.frac} is an easy consequence:

\begin{proof}
[Proof of Theorem \ref{thm.cyclotomic.frac}.]Define $a_{1},a_{2},a_{3}%
,\ldots\in\mathbb{Z}\left[  x\right]  $ as in the above proof of Theorem
\ref{thm.cyclotomic.EC}. Then, each $n\geq1$ satisfies%
\[
a_{n}=x^{n}-1=\prod_{d\mid n}\Phi_{d}\ \ \ \ \ \ \ \ \ \ \left(  \text{by
(\ref{eq.thm.cyclotomic.EC.eq})}\right)  .
\]
Thus, the assumption (\ref{eq.thm.gcdseq.fac-pf.ass}) of Theorem
\ref{thm.gcdseq.fac-pf} is satisfied for $R=\mathbb{Z}\left[  x\right]  $ and
$b_{m}=\Phi_{m}$.

Now, let $n\geq1$ be an integer. Then, the element $a_{n}$ is regular in
$\mathbb{Z}\left[  x\right]  $ (as we proved in Claim 2 of our above proof of
Theorem \ref{thm.cyclotomic.EC}). Thus, Theorem \ref{thm.gcdseq.fac-pf}
\textbf{(b)} (applied to $R=\mathbb{Z}\left[  x\right]  $ and $b_{m}=\Phi_{m}%
$) yields
\[
\Phi_{n}=\dfrac{\prod\limits_{\substack{S\subseteq\operatorname*{PF}\left(
n\right)  ;\\\left\vert S\right\vert \text{ is even}}}\left(
x^{n/\operatorname*{prod}S}-1\right)  }{\prod\limits_{\substack{S\subseteq
\operatorname*{PF}\left(  n\right)  ;\\\left\vert S\right\vert \text{ is odd}%
}}\left(  x^{n/\operatorname*{prod}S}-1\right)  }.
\]
This proves Theorem \ref{thm.cyclotomic.frac}.
\end{proof}

\subsubsection{Proof of Theorem \ref{thm.cyclotomic.DC}}

We prove Theorem \ref{thm.cyclotomic.DC} next; this will be an easy strong
induction, once the following easy lemma is available:

\begin{lemma}
\label{lem.polring.univar-backmonic}Let $R$ be a commutative ring. Let
$p,q,r\in R\left[  x\right]  $ be three polynomials such that $pq=r$. Assume
that both $p$ and $r$ are monic. Then, $q$ is monic as well, and has degree
$\deg q=\deg r-\deg p$.
\end{lemma}

\begin{proof}
The polynomial $p$ is monic; thus, we have $p\neq0$ and the leading
coefficient of $p$ is $1$, which is a unit. Hence, Proposition
\ref{prop.polring.univar-degpq} \textbf{(b)} yields that $\deg\left(
pq\right)  =\deg p+\deg q$. In other words, $\deg r=\deg p+\deg q$ (since
$r=pq$). Hence, $\deg q=\deg r-\deg p$.

It remains to prove that $q$ is monic. Let $n=\deg p$ and $m=\deg q$. Then,
Proposition \ref{prop.polring.univar-degpq} \textbf{(d)} yields%
\begin{equation}
\left[  x^{n+m}\right]  \left(  pq\right)  =\left[  x^{n}\right]  \left(
p\right)  \cdot\left[  x^{m}\right]  \left(  q\right)  .
\label{pf.lem.polring.univar-backmonic.1}%
\end{equation}
But $n=\deg p$ and $m=\deg q$ lead to $n+m=\deg p+\deg q=\deg\left(
pq\right)  $, so that $\left[  x^{n+m}\right]  \left(  pq\right)  =\left[
x^{\deg\left(  pq\right)  }\right]  \left(  pq\right)  =1$ (since $pq=r$ is
monic). Hence, (\ref{pf.lem.polring.univar-backmonic.1}) rewrites as
\[
1=\underbrace{\left[  x^{n}\right]  \left(  p\right)  }%
_{\substack{=1\\\text{(since }p\text{ is monic}\\\text{of degree }\deg
p=n\text{)}}}\cdot\left[  x^{m}\right]  \left(  q\right)  =\left[
x^{m}\right]  \left(  q\right)  =\left[  x^{\deg q}\right]  \left(  q\right)
\]
(since $m=\deg q$). In other words, $\left[  x^{\deg q}\right]  \left(
q\right)  =1$. That is, $q$ is monic. This completes the proof of Lemma
\ref{lem.polring.univar-backmonic}.
\end{proof}

\begin{proof}
[Proof of Theorem \ref{thm.cyclotomic.DC}.]We proceed by strong induction on
$n$. Thus, we fix a positive integer $n$, and we assume (as induction
hypothesis) that Theorem \ref{thm.cyclotomic.DC} holds for each of the
positive integers $1,2,\ldots,n-1$ instead of $n$. We must then prove that
Theorem \ref{thm.cyclotomic.DC} holds for $n$ as well. In other words, we must
show that $\Phi_{n}$ is a monic polynomial of degree $\phi\left(  n\right)  $.

From (\ref{eq.thm.cyclotomic.EC.eq-n}), we have%
\begin{equation}
x^{n}-1=\Phi_{n}\prod_{\substack{d\mid n;\\d\neq n}}\Phi_{d}=\left(
\prod_{\substack{d\mid n;\\d\neq n}}\Phi_{d}\right)  \Phi_{n}.
\label{pf.thm.cyclotomic.DC.1}%
\end{equation}
Moreover, Theorem \ref{thm.ent.phi.sum-div} yields%
\begin{equation}
n=\sum_{d\mid n}\phi\left(  d\right)  =\phi\left(  n\right)  +\sum
_{\substack{d\mid n;\\d\neq n}}\phi\left(  d\right)
\label{pf.thm.cyclotomic.DC.2}%
\end{equation}
(here, we have split off the addend for $d=n$ from the sum).

However, let $d$ be a proper divisor of $n$ satisfying $d\neq n$. Then, $d<n$
and thus $d\leq n-1$, so that $d\in\left\{  1,2,\ldots,n-1\right\}  $. Hence,
our induction hypothesis shows that Theorem \ref{thm.cyclotomic.DC} holds for
$d$ instead of $n$. In other words, $\Phi_{d}$ is a monic polynomial of degree
$\phi\left(  d\right)  $.

Forget that we fixed $d$. We thus have shown that if $d$ is a proper divisor
of $n$ satisfying $d\neq n$, then $\Phi_{d}$ is a monic polynomial of degree
$\phi\left(  d\right)  $. Hence, $\prod_{\substack{d\mid n;\\d\neq n}}\Phi
_{d}$ is a product of monic polynomials of respective degrees $\deg\Phi
_{d}=\phi\left(  d\right)  $. Thus, Corollary \ref{cor.polring.univar-degpk}
\textbf{(a)} shows that this product $\prod_{\substack{d\mid n;\\d\neq n}%
}\Phi_{d}$ is itself monic and has degree%
\[
\deg\prod_{\substack{d\mid n;\\d\neq n}}\Phi_{d}=\sum_{\substack{d\mid
n;\\d\neq n}}\underbrace{\deg\Phi_{d}}_{=\phi\left(  d\right)  }%
=\sum_{\substack{d\mid n;\\d\neq n}}\phi\left(  d\right)  =n-\phi\left(
n\right)  \ \ \ \ \ \ \ \ \ \ \left(  \text{by (\ref{pf.thm.cyclotomic.DC.2}%
)}\right)  .
\]

Now, we know that both polynomials $\prod_{\substack{d\mid n;\\d\neq n}%
}\Phi_{d}$ and $x^{n}-1$ are monic (the former we just proved; the latter is
obvious) and satisfy $\left(  \prod_{\substack{d\mid n;\\d\neq n}}\Phi
_{d}\right)  \Phi_{n}=x^{n}-1$ (by (\ref{pf.thm.cyclotomic.DC.1})). Hence,
Lemma \ref{lem.polring.univar-backmonic} (applied to $p=\prod_{\substack{d\mid
n;\\d\neq n}}\Phi_{d}$ and $q=\Phi_{n}$ and $r=x^{n}-1$) shows that $\Phi_{n}$
is monic as well, and has degree%
\[
\deg\left(  \Phi_{n}\right)  =\underbrace{\deg\left(  x^{n}-1\right)  }%
_{=n}-\underbrace{\deg\prod_{\substack{d\mid n;\\d\neq n}}\Phi_{d}}%
_{=n-\phi\left(  n\right)  }=n-\left(  n-\phi\left(  n\right)  \right)
=\phi\left(  n\right)  .
\]
In other words, $\Phi_{n}$ is a monic polynomial of degree $\phi\left(
n\right)  $. This completes the induction step. Thus, Theorem
\ref{thm.cyclotomic.DC} is proved.
\end{proof}

\subsubsection{Proof of Theorem \ref{thm.cyclotomic.RC}}

It remains to prove Theorem \ref{thm.cyclotomic.RC}. Before we do this, we
show something simpler:

\begin{theorem}
\label{thm.cyclotomic.RCxn-1}Let $F$ be a field. Consider its group of units
$F^{\times}$ (see Theorem \ref{thm.ringunits.group}). Let $n$ be a positive
integer. Let $u$ be an element of the group $F^{\times}$ that has order $n$
(that is, satisfies $u^{n}=1$ but $u^{k}\neq1$ for all positive integers
$k<n$). Then,%
\begin{equation}
x^{n}-1=\prod_{k\in\left\{  1,2,\ldots,n\right\}  }\left(  x-u^{k}\right)
\ \ \ \ \ \ \ \ \ \ \text{in }F\left[  x\right]  .
\label{eq.thm.cyclotomic.RCxn-1.eq}%
\end{equation}

\end{theorem}

\begin{proof}
Consider the cyclic subgroup $\left\langle u\right\rangle $ of $F^{\times}$
generated by $u$. The order of this subgroup $\left\langle u\right\rangle $
equals the order of $u$, which is $n$ (since $u$ has order $n$). Thus,
$\left\langle u\right\rangle $ is a cyclic group of order $n$. By the
classical theory of cyclic groups, this yields that it consists of the $n$
distinct elements $u^{1},u^{2},\ldots,u^{n}$. Hence,
\[
\prod_{\alpha\in\left\langle u\right\rangle }\left(  x-\alpha\right)
=\prod_{k\in\left\{  1,2,\ldots,n\right\}  }\left(  x-u^{k}\right)  .
\]

Now, define the polynomial%
\[
p:=\underbrace{x^{n}-1}_{\substack{\text{a monic polynomial}\\\text{of degree
}n}}-\underbrace{\prod_{\alpha\in\left\langle u\right\rangle }\left(
x-\alpha\right)  }_{\substack{\text{a monic polynomial}\\\text{of degree
}\left\vert \left\langle u\right\rangle \right\vert =n}}\in F\left[  x\right]
.
\]
This polynomial $p$ has degree $\deg p<n$ (since it is a difference of two
monic polynomials of degree $n$; their leading terms $x^{n}$ cancel when they
are subtracted). Moreover, by the definition of $p$, we know that each
$i\in\mathbb{Z}$ satisfies
\begin{align*}
p\left(  u^{i}\right)   &  =\underbrace{\left(  u^{i}\right)  ^{n}%
}_{\substack{=u^{in}=\left(  u^{n}\right)  ^{i}=1^{i}\\\text{(since }%
u^{n}=1\text{)}}}-\,1-\underbrace{\prod_{\alpha\in\left\langle u\right\rangle
}\left(  u^{i}-\alpha\right)  }_{\substack{=0\\\text{(since }u^{i}%
\in\left\langle u\right\rangle \text{, so that one of the factors}\\\text{of
this product is }u^{i}-u^{i}=0\text{)}}}\\
&  =1^{i}-1-0=1^{i}-1=0\ \ \ \ \ \ \ \ \ \ \left(  \text{since }%
1^{i}=1\right)  .
\end{align*}
In other words, for each $i\in\mathbb{Z}$, the element $u^{i}$ is a root of
$p$. Thus, all the $n$ distinct elements $u^{1},u^{2},\ldots,u^{n}$ are roots
of $p$. Hence, $p$ has at least $n$ roots in $F$. Since $n>\deg p$ (because
$\deg p<n$), this entails that $p$ has more than $\deg p$ roots in $F$. But
Theorem \ref{thm.polring.univar-easyFTA} (applied to $p$ and $\deg p$ instead
of $f$ and $n$) shows that if $p$ is nonzero, then $p$ has at most $\deg p$
roots in $F$. Hence, $p$ cannot be nonzero (since $p$ has more than $\deg p$
roots in $F$). That is, $p=0$. Therefore,%
\begin{align*}
x^{n}-1  &  =\prod_{\alpha\in\left\langle u\right\rangle }\left(
x-\alpha\right)  \ \ \ \ \ \ \ \ \ \ \left(  \text{since }x^{n}-1-\prod
_{\alpha\in\left\langle u\right\rangle }\left(  x-\alpha\right)  =p=0\right)
\\
&  =\prod_{k\in\left\{  1,2,\ldots,n\right\}  }\left(  x-u^{k}\right)  .
\end{align*}
This proves Theorem \ref{thm.cyclotomic.RCxn-1}.
\end{proof}

\begin{proof}
[Proof of Theorem \ref{thm.cyclotomic.RC}.]We proceed by strong induction on
$n$. Thus, we fix a positive integer $n$, and assume (as the induction
hypothesis) that Theorem \ref{thm.cyclotomic.RC} holds for each of the
positive integers $1,2,\ldots,n-1$ instead of $n$. We must then prove that
Theorem \ref{thm.cyclotomic.RC} holds for $n$ as well.

Our induction hypothesis is saying that Theorem \ref{thm.cyclotomic.RC} holds
for each of the positive integers $1,2,\ldots,n-1$ instead of $n$. In other
words, if $m\in\left\{  1,2,\ldots,n-1\right\}  $, then Theorem
\ref{thm.cyclotomic.RC} holds for $m$ instead of $n$. In other words, if
$m\in\left\{  1,2,\ldots,n-1\right\}  $, and if $u$ is an element of the group
$F^{\times}$ that has order $m$, then%
\begin{equation}
\Phi_{m}=\prod_{\substack{k\in\left\{  1,2,\ldots,m\right\}  ;\\\gcd\left(
k,m\right)  =1}}\left(  x-u^{k}\right)  . \label{pf.thm.cyclotomic.RC.IH}%
\end{equation}

Now, let us prove Theorem \ref{thm.cyclotomic.RC} for $n$. Let $u$ be an
element of the group $F^{\times}$ that has order $n$.

For any positive divisor $d$ of $n$, we define the polynomial%
\[
p_{d}:=\prod_{\substack{k\in\left\{  1,2,\ldots,n\right\}  ;\\\gcd\left(
k,n\right)  =d}}\left(  x-u^{k}\right)  \in F\left[  x\right]  .
\]
For each $k\in\left\{  1,2,\ldots,n\right\}  $, the number $\gcd\left(
k,n\right)  $ is a positive divisor of $n$. Thus,%
\[
\prod_{k\in\left\{  1,2,\ldots,n\right\}  }\left(  x-u^{k}\right)
=\prod_{d\mid n}\ \ \underbrace{\prod_{\substack{k\in\left\{  1,2,\ldots
,n\right\}  ;\\\gcd\left(  k,n\right)  =d}}\left(  x-u^{k}\right)
}_{\substack{=p_{d}\\\text{(by the definition of }p_{d}\text{)}}}=\prod_{d\mid
n}p_{d}.
\]
In view of (\ref{eq.thm.cyclotomic.RCxn-1.eq}), we can rewrite this as%
\begin{equation}
x^{n}-1=\prod_{d\mid n}p_{d}. \label{pf.thm.cyclotomic.RC.1}%
\end{equation}

Recall that $u$ has order $n$ in the group $F^{\times}$. That is, $u^{n}=1$,
but%
\begin{equation}
u^{j}\neq1\ \ \ \ \ \ \ \ \ \ \text{for each positive integer }j<n.
\label{pf.thm.cyclotomic.RC.ujn1}%
\end{equation}

Now, let $d$ be a positive divisor of $n$ satisfying $d\neq1$. Then, $\left(
u^{d}\right)  ^{n/d}=u^{d\left(  n/d\right)  }=u^{n}=1$, whereas each positive
integer $k<n/d$ satisfies
\[
\left(  u^{d}\right)  ^{k}=u^{dk}\neq1
\]
(by (\ref{pf.thm.cyclotomic.RC.ujn1}), since $dk$ is a positive integer
satisfying $dk<n$ (because $k<n/d$)). Thus, the element $u^{d}$ has order
$n/d$ in the group $F^{\times}$. Furthermore, $n/d\neq n$ (since $d\neq1$), so
that $n/d<n$ (since $n/d$ is a divisor of $n$ and thus satisfies $n/d\leq n$).
Since $n/d$ is a positive integer, this entails $n/d\in\left\{  1,2,\ldots
,n-1\right\}  $. Hence, (\ref{pf.thm.cyclotomic.RC.IH}) (applied to $n/d$ and
$u^{d}$ instead of $m$ and $u$) yields%
\begin{align}
\Phi_{n/d}  &  =\prod_{\substack{k\in\left\{  1,2,\ldots,n/d\right\}
;\\\gcd\left(  k,n/d\right)  =1}}\left(  x-\left(  u^{d}\right)  ^{k}\right)
=\prod_{\substack{i\in\left\{  1,2,\ldots,n/d\right\}  ;\\\gcd\left(
i,n/d\right)  =1}}\left(  x-\underbrace{\left(  u^{d}\right)  ^{i}}_{=u^{di}%
}\right) \nonumber\\
&  =\prod_{\substack{i\in\left\{  1,2,\ldots,n/d\right\}  ;\\\gcd\left(
i,n/d\right)  =1}}\left(  x-u^{di}\right)  . \label{pf.thm.cyclotomic.RC.4}%
\end{align}

But there is a bijection%
\begin{align*}
\left\{  i\in\left\{  1,2,\ldots,\dfrac{n}{d}\right\}  \ \mid\ \gcd\left(
i,\dfrac{n}{d}\right)  =1\right\}   &  \rightarrow\left\{  k\in\left\{
1,2,\ldots,n\right\}  \ \mid\ \gcd\left(  k,n\right)  =d\right\}  ,\\
i  &  \mapsto di
\end{align*}
(by Lemma \ref{lem.ent.phi.phi.sum-div.2}). Thus, we can substitute $di$ for
$k$ in the product $\prod_{\substack{k\in\left\{  1,2,\ldots,n\right\}
;\\\gcd\left(  k,n\right)  =d}}\left(  x-u^{k}\right)  $. Thus, we obtain%
\[
\prod_{\substack{k\in\left\{  1,2,\ldots,n\right\}  ;\\\gcd\left(  k,n\right)
=d}}\left(  x-u^{k}\right)  =\prod_{\substack{i\in\left\{  1,2,\ldots
,n/d\right\}  ;\\\gcd\left(  i,n/d\right)  =1}}\left(  x-u^{di}\right)
=\Phi_{n/d}%
\]
(by (\ref{pf.thm.cyclotomic.RC.4})). Hence, by the definition of $p_{d}$, we
have%
\begin{equation}
p_{d}=\prod_{\substack{k\in\left\{  1,2,\ldots,n\right\}  ;\\\gcd\left(
k,n\right)  =d}}\left(  x-u^{k}\right)  =\Phi_{n/d}.
\label{pf.thm.cyclotomic.RC.6}%
\end{equation}

Now, forget that we fixed $d$. We thus have proved
(\ref{pf.thm.cyclotomic.RC.6}) for every positive divisor $d$ of $n$
satisfying $d\neq1$.

Now, (\ref{pf.thm.cyclotomic.RC.1}) becomes%
\begin{align*}
x^{n}-1  &  =\prod_{d\mid n}p_{d}=p_{1}\prod_{\substack{d\mid n;\\d\neq
1}}\underbrace{p_{d}}_{\substack{=\Phi_{n/d}\\\text{(by
(\ref{pf.thm.cyclotomic.RC.6}))}}}\ \ \ \ \ \ \ \ \ \ \left(
\begin{array}
[c]{c}%
\text{here, we have split off the factor}\\
\text{for }d=1\text{ from the product}%
\end{array}
\right) \\
&  =p_{1}\prod_{\substack{d\mid n;\\d\neq1}}\Phi_{n/d}=p_{1}\underbrace{\prod
_{\substack{d\mid n;\\n/d\neq1}}}_{\substack{=\prod_{\substack{d\mid n;\\d\neq
n}}\\\text{(since }n/d\neq1\text{ is}\\\text{equivalent to }d\neq n\text{)}%
}}\underbrace{\Phi_{n/\left(  n/d\right)  }}_{\substack{=\Phi_{d}%
\\\text{(since }n/\left(  n/d\right)  =d\text{)}}}\\
&  \ \ \ \ \ \ \ \ \ \ \ \ \ \ \ \ \ \ \ \ \left(
\begin{array}
[c]{c}%
\text{here, we have substituted }n/d\text{ for }d\text{ in the product,}\\
\text{since the map (\ref{eq.lem.ent.phi.phi.sum-div.1.bij}) is a bijection}%
\end{array}
\right) \\
&  =p_{1}\prod_{\substack{d\mid n;\\d\neq n}}\Phi_{d}.
\end{align*}
But (\ref{eq.thm.cyclotomic.EC.eq-n}) yields%
\begin{equation}
x^{n}-1=\Phi_{n}\prod_{\substack{d\mid n;\\d\neq n}}\Phi_{d}.
\label{pf.thm.cyclotomic.RC.7}%
\end{equation}
Comparing these two equalities, we find%
\begin{equation}
\Phi_{n}\prod_{\substack{d\mid n;\\d\neq n}}\Phi_{d}=p_{1}\prod
_{\substack{d\mid n;\\d\neq n}}\Phi_{d}. \label{pf.thm.cyclotomic.RC.8}%
\end{equation}
Moreover, the polynomial $x^{n}-1\in F\left[  x\right]  $ is nonzero (since it
is a monic polynomial of degree $n$). The equality
(\ref{pf.thm.cyclotomic.RC.7}) shows that $\prod_{\substack{d\mid n;\\d\neq
n}}\Phi_{d}$ is a divisor of this nonzero polynomial $x^{n}-1$, and thus is
nonzero itself. Hence, we can cancel $\prod_{\substack{d\mid n;\\d\neq n}%
}\Phi_{d}$ from the equality (\ref{pf.thm.cyclotomic.RC.8}) (since $F\left[
x\right]  $ is an integral domain), obtaining%
\[
\Phi_{n}=p_{1}=\prod_{\substack{k\in\left\{  1,2,\ldots,n\right\}
;\\\gcd\left(  k,n\right)  =1}}\left(  x-u^{k}\right)
\]
(by the definition of $p_{1}$). Thus, Theorem \ref{thm.cyclotomic.RC} holds
for our $n$. This completes the induction step, and -- with it -- the proof of
Theorem \ref{thm.cyclotomic.RC}.
\end{proof}

\subsubsection{Further formulas}

An explicit formula for $\Phi_{n}$ that does not rely on division or complex
numbers is not known, and unlikely to exist. However, some identities exist
that help compute $\Phi_{n}$ for specific values of $n$:

\begin{proposition}
\label{prop.cyclotomic.2m}Let $m>1$ be an odd positive integer. Then,
$\Phi_{2m}=\Phi_{m}\left[  -x\right]  $.
\end{proposition}

\begin{fineprint}

\begin{proof}
We shall first prove that
\begin{equation}
\prod_{d\mid m}\Phi_{2d}=-\prod_{d\mid m}\Phi_{d}\left[  -x\right]  .
\label{pf.prop.cyclotomic.2m.goal1}%
\end{equation}

\begin{proof}
[Proof of (\ref{pf.prop.cyclotomic.2m.goal1}).]Indeed,
(\ref{eq.thm.cyclotomic.EC.eq}) (applied to $n=m$) yields%
\begin{equation}
x^{m}-1=\prod_{d\mid m}\Phi_{d}. \label{pf.prop.cyclotomic.2m.0a}%
\end{equation}
Substituting $-x$ for $x$ in this equality, we obtain%
\[
\left(  -x\right)  ^{m}-1=\prod_{d\mid m}\Phi_{d}\left[  -x\right]  .
\]
Since $\left(  -x\right)  ^{m}=-x^{m}$ (because $m$ is odd), we can rewrite
this as%
\[
-x^{m}-1=\prod_{d\mid m}\Phi_{d}\left[  -x\right]  .
\]
Multiplying this equality by $-1$, we obtain%
\begin{equation}
x^{m}+1=-\prod_{d\mid m}\Phi_{d}\left[  -x\right]  .
\label{pf.prop.cyclotomic.2m.0b}%
\end{equation}
Furthermore, (\ref{eq.thm.cyclotomic.EC.eq}) (applied to $n=2m$) yields%
\begin{equation}
x^{2m}-1=\prod_{d\mid2m}\Phi_{d}=\left(  \prod_{\substack{d\mid2m;\\d\text{ is
odd}}}\Phi_{d}\right)  \left(  \prod_{\substack{d\mid2m;\\d\text{ is even}%
}}\Phi_{d}\right)  \label{pf.prop.cyclotomic.2m.1}%
\end{equation}
(since each $d$ is either odd or even but never both). But $m$ is odd. Thus,
the \textbf{odd} positive divisors of $2m$ are precisely the positive divisors
of $m$, whereas the \textbf{even} positive divisors of $2m$ are precisely $2r$
for $r$ being the positive divisors of $m$. The first of these two statements
yields%
\[
\prod_{\substack{d\mid2m;\\d\text{ is odd}}}\Phi_{d}=\prod_{d\mid m}\Phi_{d},
\]
whereas the second statement yields%
\[
\prod_{\substack{d\mid2m;\\d\text{ is even}}}\Phi_{d}=\prod_{r\mid m}\Phi
_{2r}.
\]
In view of these two equalities, we can rewrite (\ref{pf.prop.cyclotomic.2m.1}%
) as%
\[
x^{2m}-1=\underbrace{\left(  \prod_{d\mid m}\Phi_{d}\right)  }%
_{\substack{=x^{m}-1\\\text{(by (\ref{pf.prop.cyclotomic.2m.0a})))}}}\left(
\prod_{r\mid m}\Phi_{2r}\right)  =\left(  x^{m}-1\right)  \left(  \prod_{r\mid
m}\Phi_{2r}\right)  .
\]
Thus,%
\[
\left(  x^{m}-1\right)  \left(  \prod_{r\mid m}\Phi_{2r}\right)
=x^{2m}-1=\left(  x^{m}\right)  ^{2}-1=\left(  x^{m}-1\right)  \left(
x^{m}+1\right)  .
\]
We can cancel $x^{m}-1$ from this equality (since $x^{m}-1$ is a nonzero
element of the integral domain $\mathbb{Z}\left[  x\right]  $), and obtain%
\[
\prod_{r\mid m}\Phi_{2r}=x^{m}+1=-\prod_{d\mid m}\Phi_{d}\left[  -x\right]
\ \ \ \ \ \ \ \ \ \ \left(  \text{by (\ref{pf.prop.cyclotomic.2m.0b})}\right)
.
\]
Upon renaming the index $r$ as $d$ on the left-hand side, this becomes
(\ref{pf.prop.cyclotomic.2m.goal1}). Thus, (\ref{pf.prop.cyclotomic.2m.goal1})
is proven.
\end{proof}

Now, let us prove Proposition \ref{prop.cyclotomic.2m}. We proceed by strong
induction on $m$. So we fix an odd $m>1$, and we assume (as induction
hypothesis) that Proposition \ref{prop.cyclotomic.2m} is already proved for
all the odd numbers among $2,3,\ldots,m-1$ instead of $m$. We must prove
Proposition \ref{prop.cyclotomic.2m} for $m$.

We have
\[
\prod_{d\mid m}\Phi_{2d}=\Phi_{2\cdot1}\cdot\Phi_{2m}\cdot\prod
_{\substack{d\mid m;\\d\neq1;\\d\neq m}}\Phi_{2d}%
\]
(here, we have split off the factors for $d=1$ and for $d=m$ from the
product\footnote{These two factors are distinct, since $m>1$.}) and
\[
\prod_{d\mid m}\Phi_{d}\left[  -x\right]  =\Phi_{1}\left[  -x\right]
\cdot\Phi_{m}\left[  -x\right]  \cdot\prod_{\substack{d\mid m;\\d\neq1;\\d\neq
m}}\Phi_{d}\left[  -x\right]
\]
(similarly). Thus, the equality (\ref{pf.prop.cyclotomic.2m.goal1}) (which we
proved above) can be rewritten as%
\[
\Phi_{2\cdot1}\cdot\Phi_{2m}\cdot\prod_{\substack{d\mid m;\\d\neq1;\\d\neq
m}}\Phi_{2d}=-\Phi_{1}\left[  -x\right]  \cdot\Phi_{m}\left[  -x\right]
\cdot\prod_{\substack{d\mid m;\\d\neq1;\\d\neq m}}\Phi_{d}\left[  -x\right]
.
\]
In view of $\Phi_{2\cdot1}=\Phi_{2}=x+1$ and $\Phi_{1}\left[  -x\right]
=-x-1$ (since $\Phi_{1}=x-1$), we can rewrite this further as
\begin{align*}
\left(  x+1\right)  \cdot\Phi_{2m}\cdot\prod_{\substack{d\mid m;\\d\neq
1;\\d\neq m}}\Phi_{2d}  &  =\underbrace{-\left(  -x-1\right)  }_{=x+1}%
\cdot\,\Phi_{m}\left[  -x\right]  \cdot\prod_{\substack{d\mid m;\\d\neq
1;\\d\neq m}}\Phi_{d}\left[  -x\right] \\
&  =\left(  x+1\right)  \cdot\Phi_{m}\left[  -x\right]  \cdot\prod
_{\substack{d\mid m;\\d\neq1;\\d\neq m}}\Phi_{d}\left[  -x\right]  .
\end{align*}
We can cancel the factor $x+1$ from this equality (since $x+1$ is a nonzero
polynomial in the integral domain $\mathbb{Z}\left[  x\right]  $), and thus
obtain%
\begin{equation}
\Phi_{2m}\cdot\prod_{\substack{d\mid m;\\d\neq1;\\d\neq m}}\Phi_{2d}=\Phi
_{m}\left[  -x\right]  \cdot\prod_{\substack{d\mid m;\\d\neq1;\\d\neq m}%
}\Phi_{d}\left[  -x\right]  . \label{pf.prop.cyclotomic.2m.goal4}%
\end{equation}

Now, let $d$ be a positive divisor of $m$ satisfying $d\neq1$ and $d\neq m$.
Then, from $d\mid m$, we obtain $d\leq m$, so that $d<m$ (since $d\neq m$).
Thus, $d\leq m-1$, so that $d\in\left\{  2,3,\ldots,m-1\right\}  $ (since $d$
is positive and $d\neq1$). Furthermore, $d$ is odd (since $d$ is a divisor of
the odd integer $m$). So $d$ is an odd number among $2,3,\ldots,m-1$. Hence,
our induction hypothesis shows that Proposition \ref{prop.cyclotomic.2m} holds
for $d$ instead of $m$; that is, we have $\Phi_{2d}=\Phi_{d}\left[  -x\right]
$. That is, $\Phi_{d}\left[  -x\right]  =\Phi_{2d}$.

Forget that we fixed $d$. We thus have proved the equality $\Phi_{d}\left[
-x\right]  =\Phi_{2d}$ for each positive divisor $d$ of $m$ satisfying
$d\neq1$ and $d\neq m$. Multiplying all these equalities, we find%
\[
\prod_{\substack{d\mid m;\\d\neq1;\\d\neq m}}\Phi_{d}\left[  -x\right]
=\prod_{\substack{d\mid m;\\d\neq1;\\d\neq m}}\Phi_{2d}.
\]
Hence, we can rewrite (\ref{pf.prop.cyclotomic.2m.goal4}) as%
\[
\Phi_{2m}\cdot\prod_{\substack{d\mid m;\\d\neq1;\\d\neq m}}\Phi_{2d}=\Phi
_{m}\left[  -x\right]  \cdot\prod_{\substack{d\mid m;\\d\neq1;\\d\neq m}%
}\Phi_{2d}.
\]
We can cancel $\prod_{\substack{d\mid m;\\d\neq1;\\d\neq m}}\Phi_{2d}$ from
this equality (since $\prod_{\substack{d\mid m;\\d\neq1;\\d\neq m}}\Phi_{2d}$
is a nonzero polynomial in the integral domain $\mathbb{Z}\left[  x\right]
$\ \ \ \ \footnote{Indeed, each cyclotomic polynomial $\Phi_{2d}$ is monic (by
Theorem \ref{thm.cyclotomic.DC}) and thus nonzero. Hence, $\prod
_{\substack{d\mid m;\\d\neq1;\\d\neq m}}\Phi_{2d}$ is a product of nonzero
polynomials in $\mathbb{Z}\left[  x\right]  $ and thus itself nonzero (since
$\mathbb{Z}\left[  x\right]  $ is an integral domain).}), and thus obtain%
\[
\Phi_{2m}=\Phi_{m}\left[  -x\right]  .
\]
That is, Proposition \ref{prop.cyclotomic.2m} holds for our $m$. This
completes the induction step; thus Proposition \ref{prop.cyclotomic.2m} is proved.
\end{proof}
\end{fineprint}

Further identities between cyclotomic polynomials can be proved in a similar vein:

\begin{exercise}
\label{exe.cyclotomic.pn}Let $p$ be a prime. Let $n$ be a positive integer.
Prove the following:

\begin{enumerate}
\item[\textbf{(a)}] If $p\mid n$, then $\Phi_{pn}=\Phi_{n}\left[
x^{p}\right]  $.

\item[\textbf{(b)}] If $p\nmid n$, then $\Phi_{n}\cdot\Phi_{pn}=\Phi
_{n}\left[  x^{p}\right]  $.
\end{enumerate}
\end{exercise}

\begin{exercise}
\label{exe.cyclotomic.nk}Let $n$ and $k$ be two positive integers such that
each prime divisor of $k$ divides $n$. (For instance, we can have
$n=2^{2}\cdot3\cdot7$ and $k=2\cdot3^{4}$.) Prove that $\Phi_{nk}=\Phi
_{n}\left[  x^{k}\right]  $.
\end{exercise}

Exercise \ref{exe.cyclotomic.nk} shows that the computation of arbitrary
cyclotomic polynomials $\Phi_{m}$ can be reduced to the computation of
$\Phi_{n}$ for squarefree positive integers $n$ (where \textquotedblleft
squarefree\textquotedblright\ means \textquotedblleft product of distinct
primes\textquotedblright: e.g., the integer $70=2\cdot5\cdot7$ is squarefree,
while $18=2\cdot3^{2}$ is not). Indeed, any positive integer $m$ can be
written as $m=nk$ for two positive integers $n$ and $k$ such that $n$ is
squarefree and each prime divisor of $k$ divides $n$. The simplest cases are
$\Phi_{p}$ for a single prime $p$ (see Proposition \ref{prop.cyclotomic.Phip}
above) and $\Phi_{pq}$ for two distinct primes $p$ and $q$ (see
\cite{LamLeu96}, \cite{Thanga09} or \cite[(4)]{ZhuQi08}).

\begin{exercise}
Let $n$ and $m$ be two distinct positive integers. Prove that $\Phi_{n}$ and
$\Phi_{m}$ are coprime in $\mathbb{Z}\left[  x\right]  $ (that is, $1$ is a
gcd of $\Phi_{n}$ and $\Phi_{m}$ in $\mathbb{Z}\left[  x\right]  $).
\end{exercise}

We shall not prove the following result of Gauss, which becomes important in
number theory:

\begin{theorem}
\label{thm.cyclotomic.irred}Let $n$ be a positive integer. Then, the
cyclotomic polynomial $\Phi_{n}$ is irreducible in $\mathbb{Z}\left[
x\right]  $.
\end{theorem}

For proofs of this theorem, see any of \cite[Theorem 3.3]{DiLoQu25},
\cite[\S 4.5, Proposition 11]{Bosch18}, \cite[Theorem 31.3]{Gallia21},
\cite[Proposition 9.1]{Klopsc08}, \cite[Lemma 9.42]{Knapp1}, \cite[\S 47]%
{Nagell64}, \cite[\S 8.4]{Waerde91}, \cite[Theorem 50.3]{Warner90},
\cite{Weintr13}, \cite[Theorem 4.14.7]{Weintr22}.

In contrast, over finite fields, $\Phi_{n}$ is not always irreducible. Some
cyclotomic polynomials are actually \textbf{never} irreducible over finite fields:

\begin{exercise}
Prove that the cyclotomic polynomials $\Phi_{8}$ and $\Phi_{15}$ are reducible
(i.e., not irreducible) in $\mathbb{F}_{p}\left[  x\right]  $ for
\textbf{every} prime $p$ (and thus also in $F\left[  x\right]  $ for every
field $F$ of positive characteristic). \medskip

[\textbf{Hint:} For $\Phi_{8}$, argue that if $-1$ is not a square in
$\mathbb{F}_{p}$, then one of the two square roots of $-1$ in $\mathbb{F}%
_{p^{2}}$ is a square in $\mathbb{F}_{p^{2}}$. For $\Phi_{15}$, argue that
$\mathbb{F}_{p^{4}}$ (for $p>5$) contains both a primitive $3$-rd root of
unity and a primitive $5$-th root of unity. Denote these two roots by $\alpha$
and $\beta$, and argue that $\alpha\beta\in\mathbb{F}_{p^{4}}$ is a primitive
$15$-th root of unity.]
\end{exercise}

\subsubsection{Roots of unity vs. roots of cyclotomic polynomials}

As Theorem \ref{thm.cyclotomic.RC} already showed, cyclotomic polynomials are
closely related to the order of elements of $F^{\times}$. Let us take this
connection a bit further:

\begin{proposition}
\label{prop.cyclotomic.roots}Let $F$ be a field, and $n$ be a positive
integer. Then:

\begin{enumerate}
\item[\textbf{(a)}] Each element of order $n$ in the group $F^{\times}$ is a
root of $\Phi_{n}$.

\item[\textbf{(b)}] Assume that $n\cdot1_{F}\neq0$. Then, the elements of
order $n$ in the group $F^{\times}$ are precisely the roots of $\Phi_{n}$ in
$F$.
\end{enumerate}
\end{proposition}

To prove this proposition, we need a simple lemma:

\begin{lemma}
\label{lem.cyclotomic.xd-1phin}Let $n$ be a positive integer. Let $d$ be a
positive divisor of $n$ such that $d\neq n$. Then,%
\[
\left(  x^{d}-1\right)  \cdot\Phi_{n}\mid x^{n}-1\ \ \ \ \ \ \ \ \ \ \text{in
}\mathbb{Z}\left[  x\right]  .
\]

\end{lemma}

\begin{proof}
[Proof of Lemma \ref{lem.cyclotomic.xd-1phin}.]We have%
\begin{equation}
x^{n}-1=\Phi_{n}\prod_{\substack{k\mid n;\\k\neq n}}\Phi_{k}
\label{pf.lem.cyclotomic.xd-1phin.1}%
\end{equation}
(indeed, this is just the equality (\ref{eq.thm.cyclotomic.EC.eq-n}), with the
index $d$ renamed as $k$ in the product) and%
\begin{equation}
x^{d}-1=\prod_{k\mid d}\Phi_{k} \label{pf.lem.cyclotomic.xd-1phin.2}%
\end{equation}
(indeed, this is just the equality (\ref{eq.thm.cyclotomic.EC.eq}), applied to
$d$ instead of $n$, and with the index $d$ renamed as $k$ in the product). But
$d$ is a positive divisor of $n$ such that $d\neq n$. Hence, $d\mid n$ and
$d<n$. Furthermore, any positive divisor $k$ of $d$ is a positive divisor of
$n$ as well (since $k\mid d\mid n$) and satisfies $k\neq n$ (since $k\leq
d<n$). Thus, the product $\prod_{k\mid d}\Phi_{k}$ is a subproduct of
$\prod_{\substack{k\mid n;\\k\neq n}}\Phi_{k}$. In particular, the former
product divides the latter. That is,
\begin{equation}
\prod_{\substack{k\mid n;\\k\neq n}}\Phi_{k}=\left(  \prod_{k\mid d}\Phi
_{k}\right)  w \label{pf.lem.cyclotomic.xd-1phin.3}%
\end{equation}
for some $w\in\mathbb{Z}\left[  x\right]  $. Using this $w$, we have%
\begin{align*}
x^{n}-1  &  =\Phi_{n}\prod_{\substack{k\mid n;\\k\neq n}}\Phi_{k}%
\ \ \ \ \ \ \ \ \ \ \left(  \text{by (\ref{pf.lem.cyclotomic.xd-1phin.1}%
)}\right) \\
&  =\Phi_{n}\underbrace{\left(  \prod_{k\mid d}\Phi_{k}\right)  }%
_{\substack{=x^{d}-1\\\text{(by (\ref{pf.lem.cyclotomic.xd-1phin.2}))}%
}}w\ \ \ \ \ \ \ \ \ \ \left(  \text{by (\ref{pf.lem.cyclotomic.xd-1phin.3}%
)}\right) \\
&  =\Phi_{n}\cdot\left(  x^{d}-1\right)  \cdot w=\left(  x^{d}-1\right)
\cdot\Phi_{n}\cdot w.
\end{align*}
Thus, $\left(  x^{d}-1\right)  \cdot\Phi_{n}\mid x^{n}-1$, which proves the lemma.
\end{proof}

\begin{proof}
[Proof of Proposition \ref{prop.cyclotomic.roots}.]\textbf{(a)} Let $u$ be an
element of order $n$ in the group $F^{\times}$. We must show that $u$ is a
root of $\Phi_{n}$. But Theorem \ref{thm.cyclotomic.RC} yields%
\begin{equation}
\Phi_{n}=\prod_{\substack{k\in\left\{  1,2,\ldots,n\right\}  ;\\\gcd\left(
k,n\right)  =1}}\left(  x-u^{k}\right)  . \label{pf.prop.cyclotomic.roots.3}%
\end{equation}
One of the factors in the product on the right-hand side is $x-u^{1}$ (since
$1\in\left\{  1,2,\ldots,n\right\}  $ and $\gcd\left(  1,n\right)  =1$).
Hence, this product is divisible by $x-u^{1}=x-u$. Since this product is
$\Phi_{n}$ (by (\ref{pf.prop.cyclotomic.roots.3})), this shows that $\Phi_{n}$
is divisible by $x-u$. That is, $x-u\mid\Phi_{n}$ in $F\left[  x\right]  $. By
Proposition \ref{prop.polring.univar-rootdiv}, this means that $u$ is a root
of $\Phi_{n}$. Thus, Proposition \ref{prop.cyclotomic.roots} \textbf{(a)} is
proved. \medskip

\textbf{(b)} Part \textbf{(a)} shows that each element of order $n$ in the
group $F^{\times}$ is a root of $\Phi_{n}$ in $F$. It remains to prove the
converse -- i.e., that each root of $\Phi_{n}$ in $F$ is an element of order
$n$ in the group $F^{\times}$.

So let $a$ be a root of $\Phi_{n}$ in $F$. We must prove that $a$ is an
element of order $n$ in the group $F^{\times}$.

Since $a$ is a root of $\Phi_{n}$, we have $x-a\mid\Phi_{n}$ in $F\left[
x\right]  $ (by Proposition \ref{prop.polring.univar-rootdiv}).

From (\ref{eq.thm.cyclotomic.EC.eq-n}), we obtain $\Phi_{n}\mid x^{n}-1$ in
$\mathbb{Z}\left[  x\right]  $, thus also in $F\left[  x\right]  $. Hence,
$x-a\mid\Phi_{n}\mid x^{n}-1$ in $F\left[  x\right]  $. By Proposition
\ref{prop.polring.univar-rootdiv}, this entails that $a$ is a root of
$x^{n}-1$. That is, $a^{n}-1=0$. In other words, $a^{n}=1$. Hence, $a$ is a
unit of $F$ (with inverse $a^{n-1}$), so that $a\in F^{\times}$.

It remains to prove that $a$ has order $n$ in $F^{\times}$.

Assume the contrary. Thus, the order of $a$ in $F^{\times}$ is not $n$; but
this order is a positive divisor of $n$ (since $a^{n}=1$). Hence, it must be a
positive divisor $d$ of $n$ such that $d\neq n$. Consider this $d$. Thus,
$a^{d}=1$, so that $a^{d}-1=0$. In other words, $a$ is a root of the
polynomial $x^{d}-1$. By Proposition \ref{prop.polring.univar-rootdiv}, this
entails that $x-a\mid x^{d}-1$ in $F\left[  x\right]  $. Multiplying this
divisibility by the divisibility $x-a\mid\Phi_{n}$, we obtain%
\[
\left(  x-a\right)  \cdot\left(  x-a\right)  \mid\left(  x^{d}-1\right)
\cdot\Phi_{n}\mid x^{n}-1\ \ \ \ \ \ \ \ \ \ \left(  \text{by Lemma
\ref{lem.cyclotomic.xd-1phin}}\right)  .
\]
That is, $\left(  x-a\right)  ^{2}\mid x^{n}-1$. Hence, Corollary
\ref{cor.polring.derivative.doubleroot} (applied to $f=x^{n}-1$ and $r=a$)
yields%
\[
x-a\mid\left(  x^{n}-1\right)  ^{\prime}=nx^{n-1}\ \ \ \ \ \ \ \ \ \ \left(
\text{by the definition of a derivative}\right)  .
\]
By Proposition \ref{prop.polring.univar-rootdiv}, this means that $a$ is a
root of the polynomial $nx^{n-1}$. That is, $na^{n-1}=0$. Multiplying this by
$a$, we obtain $na^{n}=0$. Since $a^{n}=1$, we can rewrite this as
$n\cdot1_{F}=0$; but this contradicts $n\cdot1_{F}\neq0$. This contradiction
shows that our assumption was false. Hence, the proof of Proposition
\ref{prop.cyclotomic.roots} \textbf{(b)} is complete.
\end{proof}

\begin{remark}
Proposition \ref{prop.cyclotomic.roots} \textbf{(b)} does not hold without the
$n\cdot1_{F}\neq0$ assumption. For example, if $p$ is a prime, then $1$ is a
root of $\Phi_{p}$ in $\mathbb{F}_{p}$ (in fact, it is easy to see that
$\Phi_{p}=\left(  x-1\right)  ^{p-1}$ in $\mathbb{F}_{p}\left[  x\right]  $),
but does not have order $p$ in $\mathbb{F}_{p}^{\times}$.
\end{remark}

Using Proposition \ref{prop.cyclotomic.roots} \textbf{(b)}, we obtain a useful
result: A primitive $n$-th root of unity can always be adjoined to a field
$F$, provided that $n\cdot1_{F}\neq0$:

\begin{corollary}
\label{cor.cyclotomic.create-root}Let $F$ be a field, and $n$ be a positive
integer such that $n\cdot1_{F}\neq0$. Then, there exists a field $L$ that
contains $F$ as a subring and such that the group $L^{\times}$ contains an
element of order $n$.
\end{corollary}

\begin{proof}
Consider $\Phi_{n}\in\mathbb{Z}\left[  x\right]  $ as a polynomial in
$F\left[  x\right]  $. Theorem \ref{thm.cyclotomic.DC} shows that this
polynomial $\Phi_{n}\left(  x\right)  $ is monic of degree $\deg\Phi_{n}%
=\phi\left(  n\right)  $ (this clearly does not change as we pass from
$\mathbb{Z}\left[  x\right]  $ to $F\left[  x\right]  $, since its leading
coefficient $1$ remains $1$ in $F\left[  x\right]  $). But $\phi\left(
n\right)  $ is defined as the number of all $k\in\left\{  1,2,\ldots
,n\right\}  $ satisfying $\gcd\left(  k,n\right)  =1$. Since there is clearly
at least one such $k$ (namely, $1$), we thus have $\phi\left(  n\right)
\geq1$. Thus, $\deg\Phi_{n}=\phi\left(  n\right)  \geq1>0$. Hence, Theorem
\ref{thm.finfield.splitfield} \textbf{(b)} (applied to $b=\Phi_{n}$) shows
that there is a field that contains $F$ as a subring and that contains a root
of $\Phi_{n}$. Let $L$ be this field. Then, $L$ contains a root of $\Phi_{n}$.
Proposition \ref{prop.cyclotomic.roots} \textbf{(b)} (applied to $L$ instead
of $F$) shows that this root is an element of order $n$ in $L^{\times}$ (since
$n\cdot1_{L}=n\cdot1_{F}\neq0$). Hence, the group $L^{\times}$ contains an
element of order $n$. This proves Corollary \ref{cor.cyclotomic.create-root}.
\end{proof}

\subsubsection{Application: Infinitely many primes $\equiv1\operatorname{mod}%
n$}

The famous \textbf{Dirichlet theorem on arithmetic progressions} says that for
any two coprime integers $a$ and $n\geq1$, there are infinitely many primes
that are $\equiv a\operatorname{mod}n$. All known proofs of this theorem use
some form of analysis (real or complex) and are fairly deep. However, in the
case of $a=1$, there is a short algebraic proof using cyclotomic polynomials.
Let us show this proof, following \cite[Theorem 96]{Nagell64}:

\begin{theorem}
\label{thm.dirichlet-1modn}Let $n$ be a positive integer. Then, there are
infinitely many primes $p$ satisfying $p\equiv1\operatorname{mod}n$.
\end{theorem}

\begin{proof}
We proceed as in Euclid's famous proof of the infinitude of primes (which is
actually the $n=1$ case of Theorem \ref{thm.dirichlet-1modn}). So we assume
that $p_{1},p_{2},\ldots,p_{k}$ are finitely many primes $p$ that satisfy
$p\equiv1\operatorname{mod}n$. We want to find another such prime.

Theorem \ref{thm.cyclotomic.DC} shows that the cyclotomic polynomial $\Phi
_{n}\in\mathbb{Z}\left[  x\right]  $ is monic and has degree $\phi\left(
n\right)  \geq1$ (since there is always at least one $k\in\left\{
1,2,\ldots,n\right\}  $ satisfying $\gcd\left(  k,n\right)  =1$, namely
$k=1$). Hence, this polynomial grows unboundedly when regarded as a function
from $\mathbb{Z}$ to $\mathbb{Z}$. Thus, $\Phi_{n}\left(  N\right)  >1$ for
all sufficiently large integers $N$. In particular, we can thus pick a
positive integer $m$ large enough that $\Phi_{n}\left(  mnp_{1}p_{2}\cdots
p_{k}\right)  >1$ (since $np_{1}p_{2}\cdots p_{k}$ is a positive integer, so
that $mnp_{1}p_{2}\cdots p_{k}\rightarrow\infty$ as $m\rightarrow\infty$).
Pick such an $m$.

Set $N:=mnp_{1}p_{2}\cdots p_{k}$. Thus, $\Phi_{n}\left(  N\right)  =\Phi
_{n}\left(  mnp_{1}p_{2}\cdots p_{k}\right)  >1$. Hence, the integer $\Phi
_{n}\left(  N\right)  $ has a prime divisor $q$.

Consider this $q$. Then, $q\mid\Phi_{n}\left(  N\right)  $, so that $\Phi
_{n}\left(  N\right)  \equiv0\operatorname{mod}q$. Thus, the residue class
$\overline{N}$ of $N$ in $\mathbb{F}_{q}$ satisfies $\Phi_{n}\left(
\overline{N}\right)  =0$ in $\mathbb{F}_{q}$. In other words, this class
$\overline{N}$ is a root of $\Phi_{n}$ in $\mathbb{F}_{q}$.

Moreover, (\ref{eq.thm.cyclotomic.EC.eq-n}) yields $\Phi_{n}\mid x^{n}-1$ in
$\mathbb{Z}\left[  x\right]  $. Substituting $N$ for $x$ yields $\Phi
_{n}\left(  N\right)  \mid N^{n}-1$, so that $q\mid\Phi_{n}\left(  N\right)
\mid N^{n}-1$. Therefore, $N^{n}\equiv1\not \equiv 0\operatorname{mod}q$, so
that $N\not \equiv 0\operatorname{mod}q$. That is, $q\nmid N=mnp_{1}%
p_{2}\cdots p_{k}$. Therefore, $q$ cannot be any of $p_{1},p_{2},\ldots,p_{k}%
$. From $q\nmid mnp_{1}p_{2}\cdots p_{k}$, we also obtain $q\nmid n$, so that
$n\neq0$ in $\mathbb{F}_{q}$. Hence, Proposition \ref{prop.cyclotomic.roots}
\textbf{(b)} (applied to $F=\mathbb{F}_{q}$) yields that the elements of order
$n$ in the group $\mathbb{F}_{q}^{\times}$ are precisely the roots of
$\Phi_{n}$ in $\mathbb{F}_{q}$. Therefore, $\overline{N}$ is an element of
order $n$ in the group $\mathbb{F}_{q}^{\times}$ (since $\overline{N}$ is a
root of $\Phi_{n}$ in $\mathbb{F}_{q}$). Thus, the group $\mathbb{F}%
_{q}^{\times}$ has an element of order $n$. By Lagrange's theorem (in the form
\textquotedblleft the order of an element of a finite group always divides the
order of the group\textquotedblright), this entails that $n$ divides the order
of the group $\mathbb{F}_{q}^{\times}$. Since the latter order is $\left\vert
\mathbb{F}_{q}^{\times}\right\vert =q-1$, this is saying that $n\mid q-1$.
That is, $q\equiv1\operatorname{mod}n$.

So we have found a prime $q$ that satisfies $q\equiv1\operatorname{mod}n$ but
cannot be any of $p_{1},p_{2},\ldots,p_{k}$. This is precisely what we were
looking for. So Theorem \ref{thm.dirichlet-1modn} is proven.
\end{proof}

The above proof of Theorem \ref{thm.dirichlet-1modn} contains a generalizable idea:

\begin{exercise}
Let $f\in\mathbb{Z}\left[  x\right]  $ be a non-constant polynomial. Prove
that there are infinitely many primes $p$ such that there exists a positive
integer $n$ satisfying $p\mid f\left(  n\right)  $.
\end{exercise}

\subsubsection{\label{subsec.cyclotomic.Fxcyc}Application: Finite subgroups of
$F^{\times}$ are cyclic}

Let us next use cyclotomic polynomials to prove Theorem
\ref{thm.finfield.gauss3} (and, with it, Theorem \ref{thm.finfield.gauss2} and
Theorem \ref{thm.finfield.gauss1}).

Recall what Theorem \ref{thm.finfield.gauss3} says: It says that if $F$ is any
field, then any \textbf{finite} subgroup $G$ of the group $F^{\times
}=F\setminus\left\{  0\right\}  $ is cyclic.

\begin{proof}
[First proof of Theorem \ref{thm.finfield.gauss3}.]Let $n=\left\vert
G\right\vert $ (this is well-defined, since $G$ is finite); thus, $G$ is a
group of order $n$. Recall the equality (\ref{eq.thm.cyclotomic.EC.eq-n}). Let
$f\in F\left[  x\right]  $ be the product $\prod_{\substack{d\mid n;\\d\neq
n}}\Phi_{d}$, regarded as a polynomial in $F\left[  x\right]  $. Then, $f$ is
divisible by all the factors of this product. That is, for each positive
divisor $d$ of $n$ satisfying $d\neq n$, we have%
\begin{equation}
\Phi_{d}\mid f. \label{pf.thm.finfield.gauss3.a1}%
\end{equation}

Also, the equality (\ref{eq.thm.cyclotomic.EC.eq-n}) (transported into
$F\left[  x\right]  $) rewrites as%
\begin{equation}
x^{n}-1=\Phi_{n}f. \label{pf.thm.finfield.gauss3.a2}%
\end{equation}
Hence,
\begin{equation}
\deg\left(  x^{n}-1\right)  =\deg\left(  \Phi_{n}f\right)  =\deg\Phi_{n}+\deg
f. \label{pf.thm.finfield.gauss3.a3}%
\end{equation}

But we can easily see that $\deg\Phi_{n}\geq1$ (indeed, we already showed this
in our above proof of Corollary \ref{cor.cyclotomic.create-root}). Thus,
(\ref{pf.thm.finfield.gauss3.a3}) becomes%
\[
\deg\left(  x^{n}-1\right)  =\underbrace{\deg\Phi_{n}}_{\geq1>0}+\deg f>\deg
f,
\]
so that%
\[
\deg f<\deg\left(  x^{n}-1\right)  =n.
\]
Moreover, the polynomial $f$ is nonzero (since $\Phi_{n}f=x^{n}-1$ is nonzero).

But the easy half of the FTA (Theorem \ref{thm.polring.univar-easyFTA}) says
that a nonzero polynomial of degree $\leq m$ over an integral domain $R$ has
at most $m$ roots in $R$. That is, any nonzero polynomial $g$ over an integral
domain $R$ has at most $\deg g$ roots in $R$. Applying this to $R=F$ and
$g=f$, we conclude that $f$ has at most $\deg f$ roots in $F$. Hence, $f$ has
fewer than $n$ roots in $F$ (since $\deg f<n$). Therefore, at least one
element $u\in G$ is not a root of $f$ (since $G$ has $\left\vert G\right\vert
=n$ elements, but $f$ has fewer than $n$ roots). Consider this $u$. Thus,
$f\left[  u\right]  \neq0$ (since $u$ is not a root of $f$).

Lagrange's theorem says that any element $g$ of an order-$n$ group satisfies
$g^{n}=1$. Thus, any $g\in G$ satisfies $g^{n}=1$ (since the group $G$ has
order $\left\vert G\right\vert =n$). Applying this to $g=u$, we find $u^{n}%
=1$. Hence, $u$ is a unit of $F$ (with inverse $u^{n-1}$); that is, $u\in
F^{\times}$. From $u^{n}=1$, we now conclude that $u$ has a finite order in
the group $F^{\times}$, and this order must be a positive divisor $d$ of $n$.
Consider this $d$. Proposition \ref{prop.cyclotomic.roots} \textbf{(a)}
(applied to $d$ instead of $n$) shows that $u$ is a root of $\Phi_{d}$ (since
$u$ is an element of order $d$ in the group $F^{\times}$). That is, $\Phi
_{d}\left[  u\right]  =0$.

Now, if we had $d\neq n$, then (\ref{pf.thm.finfield.gauss3.a1}) would yield
$\Phi_{d}\mid f$, and thus $\Phi_{d}\left[  u\right]  =0$ would entail
$f\left[  u\right]  =0$, which would contradict $f\left[  u\right]  \neq0$.
Hence, $d=n$. Thus, the order of $u$ in the group $F^{\times}$ is $n$ (since
this order is $d$). In other words, the order of $u$ in the group $G$ is $n$
(since $G$ is a subgroup of $F^{\times}$ that contains $u$). Thus, the cyclic
subgroup of $G$ generated by $u$ has order $n$. Since the whole group $G$ has
order $n$ as well, this shows that this subgroup must be the whole $G$. Hence,
$G$ is a cyclic group. This proves Theorem \ref{thm.finfield.gauss3}.
\end{proof}

\begin{proof}
[Proof of Theorem \ref{thm.finfield.gauss2}.]Apply Theorem
\ref{thm.finfield.gauss3} to $G=F^{\times}$.
\end{proof}

\begin{proof}
[Proof of Theorem \ref{thm.finfield.gauss1}.]Apply Theorem
\ref{thm.finfield.gauss2} to $F=\mathbb{Z}/p$. This yields that the group
$\left(  \mathbb{Z}/p\right)  ^{\times}$ is cyclic. In other words, there
exists some $g\in\left(  \mathbb{Z}/p\right)  ^{\times}$ such that its powers
$g^{0},g^{1},\ldots,g^{\left\vert \left(  \mathbb{Z}/p\right)  ^{\times
}\right\vert -1}$ are distinct and satisfy%
\[
\left(  \mathbb{Z}/p\right)  ^{\times}=\left\{  g^{0},g^{1},\ldots
,g^{\left\vert \left(  \mathbb{Z}/p\right)  ^{\times}\right\vert -1}\right\}
.
\]
In view of $\left\vert \left(  \mathbb{Z}/p\right)  ^{\times}\right\vert
=p-1$, this rewrites as follows: There exists some $g\in\left(  \mathbb{Z}%
/p\right)  ^{\times}$ such that its $p-1$ powers $g^{0},g^{1},\ldots,g^{p-2}$
are distinct and satisfy%
\[
\left(  \mathbb{Z}/p\right)  ^{\times}=\left\{  g^{0},g^{1},\ldots
,g^{p-2}\right\}  .
\]
This proves Theorem \ref{thm.finfield.gauss1}.
\end{proof}

\subsubsection{Application: Wedderburn's Little Theorem}

\begin{fineprint}
Another application of cyclotomic polynomials is perhaps even more unexpected.

Let us recall that we have initially constructed the finite field
$\mathbb{F}_{9}$ as a finite analogue of the complex numbers $\mathbb{C}$, by
adjoining a square root of $-1$ to $\mathbb{F}_{3}$ instead of to $\mathbb{R}%
$. Do the Hamiltonian quaternions have a similar finite analogue? Their
construction can certainly be reenacted over any commutative ring instead of
$\mathbb{R}$ (see, e.g., \cite{whyquat}); but their most striking property --
namely, that every nonzero quaternion is invertible -- will never work over a
finite field. This is not a coincidence. Indeed, starting with a finite field,
we will only obtain finite rings, and a finite ring cannot be a division ring
(i.e., a nontrivial ring whose nonzero elements are all invertible) without
being commutative. This fact is known as \textbf{Wedderburn's Little Theorem},
and we will give a proof using cyclotomic polynomials.

First, let us state it properly:

\begin{definition}
Let $R$ be a ring. Assume that $0\neq1$ in $R$. We say that $R$ is a
\textbf{division ring} (aka a \textbf{skew-field}) if every nonzero element of
$R$ is a unit.
\end{definition}

\begin{theorem}
\label{thm.little-wedderburn}Let $R$ be a finite division ring. Then, $R$ is a field.
\end{theorem}

The proof we shall give is due to Witt \cite{Witt31}, whose 1-page paper
leaves a lot more unsaid than it explains. This proof is rather widely known
(appearing, e.g., in \cite[\S 2]{Asgarl14}, \cite[Theorem 8.5.6]{BeaBla19},
\cite[part 2, Theorem C-2.31]{Rotman3e}, \cite[Theorem 67.7]{Elman22} and
various other texts), though we put a little twist on it by keeping it
entirely algebraic.

We will use a few lemmas. The first is a general and simple property of the
center of a ring (recall Definition \ref{def.center}):

\begin{lemma}
\label{lem.little-wedderburn.center}Let $R$ be a ring.

\begin{enumerate}
\item[\textbf{(a)}] If $a\in Z\left(  R\right)  $ is a unit of $R$, then
$a^{-1}$ also belongs to $Z\left(  R\right)  $.

\item[\textbf{(b)}] Assume that $R$ is a division ring. Then, the center
$Z\left(  R\right)  $ of $R$ is a field.
\end{enumerate}
\end{lemma}

\begin{proof}
\textbf{(a)} Assume that $a\in Z\left(  R\right)  $ is a unit of $R$. Thus,
its inverse $a^{-1}$ exists. Moreover, $a$ is central in $R$ (since $a\in
Z\left(  R\right)  $); that is, each $b\in R$ satisfies $ab=ba$. Therefore,
each $b\in R$ satisfies $a^{-1}b=ba^{-1}$ (because multiplying the equality
$ab=ba$ with $a^{-1}$ from the right, we obtain $aba^{-1}=b\underbrace{aa^{-1}%
}_{=1}=b$; and multiplying this latter equality with $a^{-1}$ from the left,
we obtain $a^{-1}aba^{-1}=a^{-1}b$, so that $a^{-1}b=\underbrace{a^{-1}a}%
_{=1}ba^{-1}=ba^{-1}$). In other words, $a^{-1}$ is central in $R$. That is,
$a^{-1}$ belongs to $Z\left(  R\right)  $. This proves Lemma
\ref{lem.little-wedderburn.center} \textbf{(a)}. \medskip

\textbf{(b)} We know from Exercise \ref{exe.center.subring} \textbf{(a)} that
$Z\left(  R\right)  $ is a commutative subring of $R$. Moreover, $0\neq1$ in
$R$ (since $R$ is a division ring), thus in $Z\left(  R\right)  $ (since $0$
and $1$ do belong to $Z\left(  R\right)  $). Hence, it remains to prove that
every nonzero element of $Z\left(  R\right)  $ is a unit -- which means here
that it is a unit of $Z\left(  R\right)  $.

It is certainly clear that every nonzero element $a$ of $Z\left(  R\right)  $
is a nonzero element of $R$, and thus is a unit of $R$ (since $R$ is a
division ring); thus, Lemma \ref{lem.little-wedderburn.center} \textbf{(a)}
shows that its inverse $a^{-1}$ also belongs to $Z\left(  R\right)  $, and
therefore $a$ is a unit of $Z\left(  R\right)  $ (since $a^{-1}$ is its
inverse in $Z\left(  R\right)  $). This is precisely what we wanted to prove.
Thus, Lemma \ref{lem.little-wedderburn.center} \textbf{(b)} is proved.
\end{proof}

Lemma \ref{lem.little-wedderburn.center} can be extended from centers to
centralizers (see Definition \ref{def.centralizer.centralizer}):

\begin{lemma}
\label{lem.little-wedderburn.centralizer}Let $R$ be a ring. Let $S$ be a
subset of $R$.

\begin{enumerate}
\item[\textbf{(a)}] If $a\in Z_{R}\left(  S\right)  $ is a unit of $R$, then
$a^{-1}$ also belongs to $Z_{R}\left(  S\right)  $.

\item[\textbf{(b)}] Assume that $R$ is a division ring. Then, the centralizer
$Z_{R}\left(  S\right)  $ is a division ring.
\end{enumerate}
\end{lemma}

\begin{proof}
Repeat our above proof of Lemma \ref{lem.little-wedderburn.center}, but now
use $b\in S$ instead of $b\in R$ (and replace $Z\left(  R\right)  $ by
$Z_{R}\left(  S\right)  $). Of course, in part \textbf{(b)}, the ring
$Z_{R}\left(  S\right)  $ is no longer commutative (in general), so it will be
only a division ring, not a field.
\end{proof}

We are now almost ready for Witt's proof of Wedderburn's Little Theorem:

\begin{proof}
[Proof of Theorem \ref{thm.little-wedderburn}.]Lemma
\ref{lem.little-wedderburn.center} \textbf{(b)} says that the center $Z\left(
R\right)  $ of $R$ is a field. Let us denote this field by $Z$. Thus,
$Z=Z\left(  R\right)  $ is a finite field. Let $q=\left\vert Z\right\vert $ be
its size. Of course, $q>1$ (since $Z$ is a field).

The ring $R$ contains the field $Z$ as a subring. Thus, $R$ becomes a left
$Z$-module (i.e., a $Z$-vector space) where we define the action of $Z$ on $R$
to be simply multiplication (i.e., the map $Z\times R\rightarrow R$ sending
each pair $\left(  z,r\right)  $ to the product $z\cdot r$%
).\ \ \ \ \footnote{Actually, we can say even more: Each $a\in Z$ and $r\in R$
satisfy $ar=ra$ (since $a\in Z=Z\left(  R\right)  $, so that $a$ is central in
$R$). Thus, the ring $R$ is a $Z$-algebra, where the action of $Z$ on $R$ is
simply multiplication (left or right -- it does not matter because of the
$ar=ra$ equality from the preceding sentence). (The scale-invariance of
multiplication is an easy consequence of centrality as well.)
\par
However, we will only need the left $Z$-module structure in what follows.}
This $Z$-vector space $R$ is finite-dimensional (since it is finite as a set),
and thus has a finite basis (by Theorem \ref{thm.vs.free}), i.e., is free of
rank $n$ for some $n\in\mathbb{N}$. That is, $R$ is isomorphic to $Z^{n}$ for
some $n\in\mathbb{N}$ (by Theorem \ref{thm.mods.free-is-iso-1}). Consider this
$n$. Thus, $\left\vert R\right\vert =\left\vert Z^{n}\right\vert =\left\vert
Z\right\vert ^{n}=q^{n}$ (since $\left\vert Z\right\vert =q$). Since $0\neq1$
in $R$, we have $\left\vert R\right\vert >1$ and thus $q^{n}=\left\vert
R\right\vert >1$, so that $n>0$. Moreover, $R$ is a division ring, so all its
nonzero elements are units. That is, $R^{\times}=R\setminus\left\{  0\right\}
$. Hence, $\left\vert R^{\times}\right\vert =\left\vert R\right\vert
-1=q^{n}-1$ (since $\left\vert R\right\vert =q^{n}$).

Two elements $a$ and $b$ of $R$ are said to be \textbf{conjugate} if there
exists a unit $u\in R^{\times}$ satisfying $a=ubu^{-1}$. The relation of being
\textquotedblleft conjugate\textquotedblright\ is a binary relation on $R$,
and is easily seen to be an equivalence relation (for example, it is symmetric
because $a=ubu^{-1}$ implies $b=u^{-1}a\left(  u^{-1}\right)  ^{-1}$; and it
is transitive because $a=ubu^{-1}$ and $b=vcv^{-1}$ imply $a=\left(
uv\right)  c\left(  uv\right)  ^{-1}$). Its equivalence classes are called the
\textbf{conjugacy classes} of $R$. Clearly, they are disjoint, and their union
is $R$; thus,%
\begin{equation}
\left\vert R\right\vert =\sum_{C\text{ is a conjugacy class of }R}\left\vert
C\right\vert . \label{pf.thm.little-wedderburn.classeq1}%
\end{equation}
In the following, we shall simply call them \textquotedblleft conjugacy
classes\textquotedblright, omitting the words \textquotedblleft of
$R$\textquotedblright. Now, let us take a closer look at these conjugacy classes:

\begin{statement}
\textit{Claim 1:} Let $z\in Z$. Then, the $1$-element set $\left\{  z\right\}
$ is a conjugacy class of $R$.
\end{statement}

\begin{proof}
[Proof of Claim 1.]We have $z\in Z=Z\left(  R\right)  $; that is, $z$ is
central in $R$. Now, let $C$ be the conjugacy class of $R$ that contains $z$.
Then, any element of $C$ is conjugate to $z$; in other words, it has the form
$uzu^{-1}$ for some unit $u\in R^{\times}$. Since $uz=zu$ (because $z$ is
central in $R$), it thus equals $\underbrace{uz}_{=zu}u^{-1}=zuu^{-1}=z$. So
we have shown that any element of $C$ equals $z$. In other words,
$C\subseteq\left\{  z\right\}  $. Since $z\in C$, we thus conclude $C=\left\{
z\right\}  $. Hence, $\left\{  z\right\}  $ is a conjugacy class of $R$ (since
$C$ is a conjugacy class of $R$). This proves Claim 1.
\end{proof}

\begin{statement}
\textit{Claim 2:} Let $C$ be a conjugacy class of $R$ that has size $1$ (that
is, satisfies $\left\vert C\right\vert =1$). Then, $C=\left\{  z\right\}  $
for some $z\in Z$.
\end{statement}

\begin{proof}
[Proof of Claim 2.]Certainly, $C=\left\{  z\right\}  $ for some $z\in R$
(since $\left\vert C\right\vert =1$). Consider this $z$. We must then show
that $z\in Z$. In other words, we must show that $z$ is central in $R$ (since
$Z=Z\left(  R\right)  $). In other words, we must show that $zb=bz$ for each
$b\in R$.

So let us fix $b\in R$. If $b=0$, then $zb=bz$ is clearly true (since
$z0=0=0z$). Thus, we WLOG assume that $b\neq0$. But $R$ is a division ring, so
that every nonzero element of $R$ is a unit. Hence, $b$ is a unit of $R$
(since $b\neq0$). Thus, $bzb^{-1}$ is conjugate to $z$, and thus lies in the
same conjugacy class as $z$. That is, $bzb^{-1}\in C$ (since $C$ is a
conjugacy class containing $z$). Therefore, $bzb^{-1}\in C=\left\{  z\right\}
$, so that $bzb^{-1}=z$ and therefore $bz=zb$. In other words, $zb=bz$.

Forget that we have fixed $b$. We thus have shown that $zb=bz$ for each $b\in
R$. In other words, $z$ is central in $R$. That is, $z\in Z\left(  R\right)
=Z$. This proves Claim 2 (since $C=\left\{  z\right\}  $).
\end{proof}

Claim 1 shows that each element $z\in Z$ gives rise to a conjugacy class
$\left\{  z\right\}  $ of $R$ that has size $1$. Since all these classes
$\left\{  z\right\}  $ are distinct, their total number is $\left\vert
Z\right\vert =q$. Thus, we have found $q$ conjugacy classes of $R$ that have
size $1$. Claim 2 shows that these $q$ conjugacy classes are the only
conjugacy classes of $R$ that have size $1$. Thus, there are exactly $q$
conjugacy classes of $R$ that have size $1$.

\begin{statement}
\textit{Claim 3:} Let $C$ be a conjugacy class of $R$. Then, $\left\vert
C\right\vert =\dfrac{q^{n}-1}{q^{d}-1}$ for some positive divisor $d$ of $n$.
\end{statement}

\begin{proof}
[Proof of Claim 3.]Pick any $z\in C$. Thus, $C$ is the conjugacy class that
contains $z$.

Let $U$ be the centralizer of $\left\{  z\right\}  $ in $R$ (see Definition
\ref{def.centralizer.centralizer}). That is, $U=Z_{R}\left(  \left\{
z\right\}  \right)  $. Then, $U$ is itself a subring of $R$ (by Exercise
\ref{exe.centralizer.subring} \textbf{(a)}). Moreover, it is easy to see that
$Z\subseteq U$\ \ \ \ \footnote{\textit{Proof.} Let $x\in Z$. Then, $x$ is
central in $R$ (since $x\in Z=Z\left(  R\right)  $). In other words, $x$
commutes with all elements of $R$, and thus, in particular, commutes with all
elements of $\left\{  z\right\}  $. In other words, $x$ centralizes $\left\{
z\right\}  $. Hence, $x\in Z_{R}\left(  \left\{  z\right\}  \right)  =U$.
Since we have proved this for all $x\in Z$, we thus have shown that
$Z\subseteq U$.}. Thus, $U$ is a subring of $R$ that contains $Z$ as a
subring. Hence, $U$ becomes a left $Z$-module (i.e., a $Z$-vector space) where
we define the action of $Z$ on $U$ to be simply multiplication (i.e., the map
$Z\times U\rightarrow U$ sending each pair $\left(  z,r\right)  $ to the
product $z\cdot r$). This $Z$-vector space $U$ is finite-dimensional (since it
is finite as a set), and thus has a finite basis (by Theorem \ref{thm.vs.free}%
), i.e., is free of rank $d$ for some $d\in\mathbb{N}$. That is, $U$ is
isomorphic to $Z^{d}$ for some $d\in\mathbb{N}$ (by Theorem
\ref{thm.mods.free-is-iso-1}). Consider this $d$. Thus, $\left\vert
U\right\vert =\left\vert Z^{d}\right\vert =\left\vert Z\right\vert ^{d}=q^{d}$
(since $\left\vert Z\right\vert =q$). Since $0\neq1$ in $U$ (because $U$ is a
subring of $R$), we have $\left\vert U\right\vert >1$ and thus $q^{d}%
=\left\vert U\right\vert >1$, so that $d>0$. Moreover, $U$ is a subset of $R$,
so that $\left\vert U\right\vert \leq\left\vert R\right\vert $ and thus
$q^{d}=\left\vert U\right\vert \leq\left\vert R\right\vert =q^{n}$; therefore,
$d\leq n$ (since $q>1$). Combining this with $d>0$, we obtain $d\in\left\{
1,2,\ldots,n\right\}  $.

The ring $U=Z_{R}\left(  \left\{  z\right\}  \right)  $ is a division ring by
Lemma \ref{lem.little-wedderburn.centralizer} \textbf{(b)}. Hence, all its
nonzero elements are units. That is, $U^{\times}=U\setminus\left\{  0\right\}
$, and thus $\left\vert U^{\times}\right\vert =\left\vert U\right\vert
-1=q^{d}-1$ (since $\left\vert U\right\vert =q^{d}$). But any unit of $U$ is
obviously a unit of $R$. That is, $U^{\times}$ is a subgroup of $R^{\times}$.
Therefore, $\left\vert U^{\times}\right\vert \mid\left\vert R^{\times
}\right\vert $ by Lagrange's theorem. That is, $q^{d}-1\mid q^{n}-1$ (since
$\left\vert U^{\times}\right\vert =q^{d}-1$ and $\left\vert R^{\times
}\right\vert =q^{n}-1$).

Now we claim that $d\mid n$. Indeed, Exercise \ref{exe.gcdseq.an-1} (applied
to $q$, $d$ and $n$ instead of $u$, $n$ and $m$) shows that
\[
\gcd\left(  q^{d}-1,\ q^{n}-1\right)  =\left\vert q^{\gcd\left(  d,n\right)
}-1\right\vert =q^{\gcd\left(  d,n\right)  }-1\ \ \ \ \ \ \ \ \ \ \left(
\text{since }q^{\gcd\left(  d,n\right)  }\geq1\right)  .
\]
Comparing this with $\gcd\left(  q^{d}-1,\ q^{n}-1\right)  =q^{d}-1$ (since
$q^{d}-1\mid q^{n}-1$ and $q^{d}-1\geq0$), we obtain $q^{d}-1=q^{\gcd\left(
d,n\right)  }-1$, whence $q^{d}=q^{\gcd\left(  d,n\right)  }$ and thus
$d=\gcd\left(  d,n\right)  $ (since $q>1$). Thus, $d=\gcd\left(  d,n\right)
\mid n$.

So we know that $d$ is a positive divisor of $n$ (since $d\mid n$ and $d>0$).

For every $b\in R^{\times}$, the element $bzb^{-1}$ of $R$ is conjugate to $z$
(by the definition of \textquotedblleft conjugate\textquotedblright) and thus
belongs to the conjugacy class $C$ (since $C$ is the conjugacy class that
contains $z$). Thus, we can define a map%
\begin{align*}
\alpha:R^{\times}  &  \rightarrow C,\\
b  &  \mapsto bzb^{-1}.
\end{align*}
We shall now show that for each $c\in C$, there are precisely $\left\vert
U^{\times}\right\vert $ many elements $b\in R^{\times}$ that satisfy
$\alpha\left(  b\right)  =c$.

Readers familiar with group actions will identify this as a particular case of
the orbit-stabilizer formula (where the group $R^{\times}$ acts on itself by
conjugation, and $C$ is the orbit of $z$ under this action, while $U^{\times}$
is the stabilizer of $z$). However, let us give a self-contained elementary
reformulation of this proof.

Fix $c\in C$. Thus, $c$ is conjugate to $z$ (since $C$ is the conjugacy class
that contains $z$). In other words, $c=uzu^{-1}$ for some $u\in R^{\times}$.
Consider this $u$.

Now, for any $b\in R^{\times}$, we have the following chain of equivalences:%
\begin{align*}
&  \ \left(  \alpha\left(  b\right)  =c\right) \\
&  \Longleftrightarrow\ \left(  bzb^{-1}=uzu^{-1}\right)
\ \ \ \ \ \ \ \ \ \ \left(  \text{since }\alpha\left(  b\right)
=bzb^{-1}\text{ and }c=uzu^{-1}\right) \\
&  \Longleftrightarrow\ \left(  u^{-1}bzb^{-1}=zu^{-1}\right)
\ \ \ \ \ \ \ \ \ \ \left(
\begin{array}
[c]{c}%
\text{here, we multiplied the equality}\\
\text{by }u^{-1}\text{ from the left}%
\end{array}
\right) \\
&  \Longleftrightarrow\ \left(  u^{-1}bz=zu^{-1}b\right)
\ \ \ \ \ \ \ \ \ \ \left(
\begin{array}
[c]{c}%
\text{here, we multiplied the equality by }b\\
\text{from the right}%
\end{array}
\right) \\
&  \Longleftrightarrow\ \left(  u^{-1}b\text{ commutes with }z\right) \\
&  \Longleftrightarrow\ \left(  u^{-1}b\text{ commutes with every element of
}\left\{  z\right\}  \right) \\
&  \Longleftrightarrow\ \left(  u^{-1}b\in Z_{R}\left(  \left\{  z\right\}
\right)  \right)  \ \ \ \ \ \ \ \ \ \ \left(  \text{by the definition of the
centralizer }Z_{R}\left(  \left\{  z\right\}  \right)  \right) \\
&  \Longleftrightarrow\ \left(  u^{-1}b\in U\right)
\ \ \ \ \ \ \ \ \ \ \left(  \text{since }Z_{R}\left(  \left\{  z\right\}
\right)  =U\right) \\
&  \Longleftrightarrow\ \left(  u^{-1}b\in U\setminus\left\{  0\right\}
\right)  \ \ \ \ \ \ \ \ \ \ \left(  \text{since }u^{-1}b\text{ cannot be
}0\text{ (because }u,b\in R^{\times}\text{)}\right) \\
&  \Longleftrightarrow\ \left(  u^{-1}b\in U^{\times}\right)
\ \ \ \ \ \ \ \ \ \ \left(  \text{since }U\setminus\left\{  0\right\}
=U^{\times}\right) \\
&  \Longleftrightarrow\ \left(  b\in uU^{\times}\right)
\end{align*}
(where $uU^{\times}$ is a coset of the subgroup $U^{\times}$ in the group
$R^{\times}$). Hence, the elements $b\in R^{\times}$ that satisfy
$\alpha\left(  b\right)  =c$ are precisely the elements $b\in R^{\times}$ that
satisfy $b\in uU^{\times}$; in other words, they are precisely the elements of
the coset $uU^{\times}$. Hence, their number is $\left\vert uU^{\times
}\right\vert =\left\vert U^{\times}\right\vert $ (since a coset of a subgroup
always has the same size as the subgroup itself). In other words, there are
precisely $\left\vert U^{\times}\right\vert $ many elements $b\in R^{\times}$
that satisfy $\alpha\left(  b\right)  =c$.

Forget that we fixed $c$. We have thus shown that for each $c\in C$, there are
precisely $\left\vert U^{\times}\right\vert $ many elements $b\in R^{\times}$
that satisfy $\alpha\left(  b\right)  =c$. In other words, for each $c\in C$,
we have%
\begin{equation}
\left(  \text{number of all }b\in R^{\times}\text{ that satisfy }\alpha\left(
b\right)  =c\right)  =\left\vert U^{\times}\right\vert =q^{d}%
-1.\ \ \ \ \ \ \ \ \ \ \label{pf.thm.little-wedderburn.c3.pf.9}%
\end{equation}

Now,%
\begin{align*}
q^{n}-1  &  =\left\vert R^{\times}\right\vert =\sum_{c\in C}%
\underbrace{\left(  \text{number of all }b\in R^{\times}\text{ that satisfy
}\alpha\left(  b\right)  =c\right)  }_{\substack{=q^{d}-1\\\text{(by
(\ref{pf.thm.little-wedderburn.c3.pf.9}))}}}\\
&  =\sum_{c\in C}\left(  q^{d}-1\right)  =\left\vert C\right\vert \cdot\left(
q^{d}-1\right)  ,
\end{align*}
so that
\[
\left\vert C\right\vert =\dfrac{q^{n}-1}{q^{d}-1}.
\]
Since $d$ is a positive divisor of $n$, we thus have proved Claim 3.
\end{proof}

\begin{statement}
\textit{Claim 4:} Let $C$ be a conjugacy class of $R$ such that $\left\vert
C\right\vert \neq1$. Then, $\left\vert C\right\vert \equiv0\operatorname{mod}%
\Phi_{n}\left(  q\right)  $.
\end{statement}

\begin{proof}
[Proof of Claim 4.]Claim 3 shows that $\left\vert C\right\vert =\dfrac
{q^{n}-1}{q^{d}-1}$ for some positive divisor $d$ of $n$. Consider this $d$.
Thus, $\dfrac{q^{n}-1}{q^{d}-1}=\left\vert C\right\vert \neq1$, so that
$q^{n}-1\neq q^{d}-1$ and therefore $n\neq d$. In other words, $d\neq n$.
Hence, Lemma \ref{lem.cyclotomic.xd-1phin} yields $\left(  x^{d}-1\right)
\cdot\Phi_{n}\mid x^{n}-1$ in $\mathbb{Z}\left[  x\right]  $. Evaluating both
sides of this divisibility at the integer $q$, we obtain $\left(
q^{d}-1\right)  \cdot\Phi_{n}\left(  q\right)  \mid q^{n}-1$ in $\mathbb{Z}$.
Hence, $\Phi_{n}\left(  q\right)  \mid\dfrac{q^{n}-1}{q^{d}-1}=\left\vert
C\right\vert $. In other words, $\left\vert C\right\vert \equiv
0\operatorname{mod}\Phi_{n}\left(  q\right)  $. This proves Claim 4.
\end{proof}

Furthermore, (\ref{eq.thm.cyclotomic.EC.eq-n}) yields $\Phi_{n}\mid x^{n}-1$
in $\mathbb{Z}\left[  x\right]  $. Evaluating both sides of this divisibility
at the integer $q$, we obtain
\begin{equation}
\Phi_{n}\left(  q\right)  \mid q^{n}-1\ \ \ \ \ \ \ \ \ \ \text{in }%
\mathbb{Z}. \label{pf.thm.little-wedderburn.R=1}%
\end{equation}

Now, (\ref{pf.thm.little-wedderburn.classeq1}) becomes%
\begin{align*}
\left\vert R\right\vert  &  =\sum_{C\text{ is a conjugacy class of }%
R}\left\vert C\right\vert \\
&  =\sum_{\substack{C\text{ is a conjugacy class of }R\text{;}\\\left\vert
C\right\vert =1}}\underbrace{\left\vert C\right\vert }_{=1}+\sum
_{\substack{C\text{ is a conjugacy class of }R\text{;}\\\left\vert
C\right\vert \neq1}}\underbrace{\left\vert C\right\vert }_{\substack{\equiv
0\operatorname{mod}\Phi_{n}\left(  q\right)  \\\text{(by Claim 4)}}}\\
&  \equiv\sum_{\substack{C\text{ is a conjugacy class of }R\text{;}%
\\\left\vert C\right\vert =1}}1+\underbrace{\sum_{\substack{C\text{ is a
conjugacy class of }R\text{;}\\\left\vert C\right\vert \neq1}}0}_{=0}\\
&  =\sum_{\substack{C\text{ is a conjugacy class of }R\text{;}\\\left\vert
C\right\vert =1}}1\\
&  \equiv\left(  \text{number of conjugacy classes }C\text{ of }R\text{ such
that }\left\vert C\right\vert =1\right) \\
&  =\left(  \text{number of conjugacy classes of }R\text{ that have size
}1\right) \\
&  =q\operatorname{mod}\Phi_{n}\left(  q\right)
\end{align*}
(since there are exactly $q$ conjugacy classes of $R$ that have size $1$).
Hence,%
\[
q\equiv\left\vert R\right\vert =q^{n}\equiv1\operatorname{mod}\Phi_{n}\left(
q\right)  \ \ \ \ \ \ \ \ \ \ \left(  \text{by
(\ref{pf.thm.little-wedderburn.R=1})}\right)  .
\]
In other words, $\Phi_{n}\left(  q\right)  \mid q-1$. Since $q-1$ is a
positive integer (because $q>1$), this entails that%
\begin{equation}
\Phi_{n}\left(  q\right)  \leq q-1 \label{pf.thm.little-wedderburn.ineq1}%
\end{equation}
(because any divisor of a positive integer is at most as large as the integer).

We now claim that $n=1$. Indeed, assume the contrary. Then, $n>1$. As we will
soon show (Lemma \ref{lem.little-wedderburn.ineq2} below), this entails that%
\begin{equation}
\Phi_{n}\left(  m\right)  >\left(  m-1\right)  ^{\phi\left(  n\right)
}\ \ \ \ \ \ \ \ \ \ \text{for each real }m>1.
\label{pf.thm.little-wedderburn.ineq2a}%
\end{equation}
Applying this to $m=q$, we obtain $\Phi_{n}\left(  q\right)  >\left(
q-1\right)  ^{\phi\left(  n\right)  }$ (since $q>1$). Since $\phi\left(
n\right)  \geq1$ (as we saw in the proof of Corollary
\ref{cor.cyclotomic.create-root}) and $q-1\geq1$ (since $q>1$ is an integer),
we have $\left(  q-1\right)  ^{\phi\left(  n\right)  }\geq\left(  q-1\right)
^{1}=q-1$, whence $\Phi_{n}\left(  q\right)  >\left(  q-1\right)
^{\phi\left(  n\right)  }\geq q-1$. This contradicts
(\ref{pf.thm.little-wedderburn.ineq1}).

This contradiction shows that our assumption was false. Thus, $n=1$.
Therefore, our above equality $\left\vert R\right\vert =\left\vert
Z\right\vert ^{n}$ rewrites as $\left\vert R\right\vert =\left\vert
Z\right\vert ^{1}=\left\vert Z\right\vert $. Hence, $R=Z$ (since $Z\subseteq
R$). Since $Z$ is a field, this shows that $R$ is a field. Theorem
\ref{thm.little-wedderburn} is thus proved...
\end{proof}

... up to the proof of an inequality which we promised: the inequality
(\ref{pf.thm.little-wedderburn.ineq2a}). Before we redeem this debt, let us
show a simpler auxiliary inequality for Euler's totient function $\phi$:

\begin{lemma}
\label{lem.little-wedderburn.ineq1}Let $n$ be a positive integer. Then,%
\[
\phi\left(  n\right)  \geq2^{\left\vert \operatorname*{PF}\left(  n\right)
\setminus\left\{  2\right\}  \right\vert }.
\]
(See Definition \ref{def.ent.PFn} for the meaning of $\operatorname*{PF}%
\left(  n\right)  $.)
\end{lemma}

\begin{proof}
Theorem \ref{thm.ent.phi.explicit} \textbf{(a)} yields
\begin{align*}
\phi\left(  n\right)   &  =\underbrace{\prod_{\substack{p\text{ prime;}\\p\mid
n}}}_{=\prod_{p\in\operatorname*{PF}\left(  n\right)  }}\left(  \left(
p-1\right)  \underbrace{p^{v_{p}\left(  n\right)  -1}}_{\geq1}\right) \\
&  \geq\prod_{p\in\operatorname*{PF}\left(  n\right)  }\left(  p-1\right)
=\prod_{\substack{p\in\operatorname*{PF}\left(  n\right)  ;\\p=2}%
}\underbrace{\left(  p-1\right)  }_{\substack{=1\\\text{(since }p=2\text{)}%
}}\cdot\underbrace{\prod_{\substack{p\in\operatorname*{PF}\left(  n\right)
;\\p\neq2}}}_{=\prod_{p\in\operatorname*{PF}\left(  n\right)  \setminus
\left\{  2\right\}  }}\underbrace{\left(  p-1\right)  }_{\substack{\geq
2\\\text{(since }p\text{ is a prime }\neq2\text{,}\\\text{so that }%
p\geq3\text{)}}}\\
&  \geq\underbrace{\prod_{\substack{p\in\operatorname*{PF}\left(  n\right)
;\\p=2}}1}_{=1}\cdot\underbrace{\prod_{p\in\operatorname*{PF}\left(  n\right)
\setminus\left\{  2\right\}  }2}_{=2^{\left\vert \operatorname*{PF}\left(
n\right)  \setminus\left\{  2\right\}  \right\vert }}=2^{\left\vert
\operatorname*{PF}\left(  n\right)  \setminus\left\{  2\right\}  \right\vert
}.
\end{align*}
This proves Lemma \ref{lem.little-wedderburn.ineq1}.
\end{proof}

And now we can finish our above proof of Theorem \ref{thm.little-wedderburn}
by proving (\ref{pf.thm.little-wedderburn.ineq2a}):

\begin{lemma}
\label{lem.little-wedderburn.ineq2}Let $n>1$ be an integer, and let $m>1$ be a
real. Then, $\Phi_{n}\left(  m\right)  >\left(  m-1\right)  ^{\phi\left(
n\right)  }$.
\end{lemma}

\begin{proof}
Use the notations $\operatorname*{PF}\left(  n\right)  $ and
$\operatorname*{prod}S$ defined in Definition \ref{def.ent.PFn} and Definition
\ref{def.ent.prodS}. The integer $n$ has at least one prime divisor (since
$n>1$); in other words, the set $\operatorname*{PF}\left(  n\right)  $ is nonempty.

If $V$ is any finite nonempty set, then exactly half of all subsets of $V$
have even size (since Lemma \ref{lem.gcdseq.fac-pf.lem-comb} says that the
number of subsets of $V$ that have odd size equals the number of subsets of
$V$ that have even size); in other words, we have%
\begin{align}
\left(  \text{number of all }S\subseteq V\text{ such that }\left\vert
S\right\vert \text{ is even}\right)   &  =\dfrac{1}{2}\cdot\underbrace{\left(
\text{number of all }S\subseteq V\right)  }_{=2^{\left\vert V\right\vert }%
}=\dfrac{1}{2}\cdot2^{\left\vert V\right\vert }\nonumber\\
&  =2^{\left\vert V\right\vert -1}.
\label{pf.lem.little-wedderburn.ineq2.comb1}%
\end{align}
Applying this to $V=\operatorname*{PF}\left(  n\right)  $, we obtain%
\begin{equation}
\left(  \text{number of all }S\subseteq\operatorname*{PF}\left(  n\right)
\text{ such that }\left\vert S\right\vert \text{ is even}\right)
=2^{\left\vert \operatorname*{PF}\left(  n\right)  \right\vert -1}%
\ \ \ \ \ \ \ \ \ \ \label{pf.lem.little-wedderburn.ineq2.comb2}%
\end{equation}
(since $\operatorname*{PF}\left(  n\right)  $ is nonempty).

Theorem \ref{thm.ent.phi.explicit} \textbf{(b)} yields%
\[
\phi\left(  n\right)  =\sum_{S\subseteq\operatorname*{PF}\left(  n\right)
}\left(  -1\right)  ^{\left\vert S\right\vert }\cdot\dfrac{n}%
{\operatorname*{prod}S}=\sum\limits_{\substack{S\subseteq\operatorname*{PF}%
\left(  n\right)  ;\\\left\vert S\right\vert \text{ is even}}}\dfrac
{n}{\operatorname*{prod}S}-\sum\limits_{\substack{S\subseteq\operatorname*{PF}%
\left(  n\right)  ;\\\left\vert S\right\vert \text{ is odd}}}\dfrac
{n}{\operatorname*{prod}S},
\]
so that%
\begin{align}
m^{\phi\left(  n\right)  }  &  =m^{\sum\limits_{\substack{S\subseteq
\operatorname*{PF}\left(  n\right)  ;\\\left\vert S\right\vert \text{ is
even}}}\dfrac{n}{\operatorname*{prod}S}-\sum\limits_{\substack{S\subseteq
\operatorname*{PF}\left(  n\right)  ;\\\left\vert S\right\vert \text{ is odd}%
}}\dfrac{n}{\operatorname*{prod}S}}\nonumber\\
&  =\dfrac{\prod\limits_{\substack{S\subseteq\operatorname*{PF}\left(
n\right)  ;\\\left\vert S\right\vert \text{ is even}}%
}m^{n/\operatorname*{prod}S}}{\prod\limits_{\substack{S\subseteq
\operatorname*{PF}\left(  n\right)  ;\\\left\vert S\right\vert \text{ is odd}%
}}m^{n/\operatorname*{prod}S}}. \label{pf.lem.little-wedderburn.bound.mphin}%
\end{align}

Note that each subset $S$ of $\operatorname*{PF}\left(  n\right)  $ satisfies
$\operatorname*{prod}S\mid n$ (by Lemma \ref{lem.prod-PF}) and thus
$\operatorname*{prod}S\leq n$, so that $n/\operatorname*{prod}S\geq1$, thus
$n/\operatorname*{prod}S-1\geq0$ and therefore
\begin{equation}
m^{n/\operatorname*{prod}S-1}\geq1\ \ \ \ \ \ \ \ \ \ \left(  \text{since
}m>1\right)  , \label{pf.lem.little-wedderburn.bound.2}%
\end{equation}
so that%
\begin{align}
m^{n/\operatorname*{prod}S}-1  &  \geq m^{n/\operatorname*{prod}%
S}-m^{n/\operatorname*{prod}S-1}=m^{n/\operatorname*{prod}S}\cdot\left(
1-\dfrac{1}{m}\right) \nonumber\\
&  =m^{n/\operatorname*{prod}S}\cdot\dfrac{m-1}{m}.
\label{pf.lem.little-wedderburn.bound.3}%
\end{align}
Moreover, if $S$ is the empty set $\varnothing$, then this inequality is
strict (since $\operatorname*{prod}S=\operatorname*{prod}\varnothing=1$ leads
to $n/\underbrace{\operatorname*{prod}S}_{=1}=n>1$ and thus
$n/\operatorname*{prod}S-1>0$).

Evaluating both sides of the equality (\ref{eq.thm.cyclotomic.frac.eq}) at
$m$, we find%
\begin{equation}
\Phi_{n}\left(  m\right)  =\dfrac{\prod\limits_{\substack{S\subseteq
\operatorname*{PF}\left(  n\right)  ;\\\left\vert S\right\vert \text{ is
even}}}\left(  m^{n/\operatorname*{prod}S}-1\right)  }{\prod
\limits_{\substack{S\subseteq\operatorname*{PF}\left(  n\right)  ;\\\left\vert
S\right\vert \text{ is odd}}}\left(  m^{n/\operatorname*{prod}S}-1\right)  }.
\label{pf.lem.little-wedderburn.bound.1}%
\end{equation}
\footnote{If evaluating a fraction sounds suspicious to you, you can
cross-multiply (\ref{eq.thm.cyclotomic.frac.eq}), thus obtaining the equality%
\[
\Phi_{n}\cdot\prod\limits_{\substack{S\subseteq\operatorname*{PF}\left(
n\right)  ;\\\left\vert S\right\vert \text{ is odd}}}\left(
x^{n/\operatorname*{prod}S}-1\right)  =\prod\limits_{\substack{S\subseteq
\operatorname*{PF}\left(  n\right)  ;\\\left\vert S\right\vert \text{ is
even}}}\left(  x^{n/\operatorname*{prod}S}-1\right)  ;
\]
then evaluate this equality at $m$ and finally divide by $\prod
\limits_{\substack{S\subseteq\operatorname*{PF}\left(  n\right)  ;\\\left\vert
S\right\vert \text{ is odd}}}\left(  m^{n/\operatorname*{prod}S}-1\right)  $
(which is nonzero because each subset $S$ of $\operatorname*{PF}\left(
n\right)  $ satisfies $n/\operatorname*{prod}S>0$ and therefore
$m^{n/\operatorname*{prod}S}>1$).} However, if we multiply the inequalities
(\ref{pf.lem.little-wedderburn.bound.3}) over all subsets $S$ of
$\operatorname*{PF}\left(  n\right)  $ whose size $\left\vert S\right\vert $
is even, then we obtain
\[
\prod\limits_{\substack{S\subseteq\operatorname*{PF}\left(  n\right)
;\\\left\vert S\right\vert \text{ is even}}}\left(  m^{n/\operatorname*{prod}%
S}-1\right)  >\prod\limits_{\substack{S\subseteq\operatorname*{PF}\left(
n\right)  ;\\\left\vert S\right\vert \text{ is even}}}\left(
m^{n/\operatorname*{prod}S}\cdot\dfrac{m-1}{m}\right)
\]
(this really is a strict inequality, because at least one of the inequalities
being multiplied is strict: namely, for the subset $S=\varnothing$, whose size
$\left\vert S\right\vert =\left\vert \varnothing\right\vert =0$ is obviously
even, the inequality (\ref{pf.lem.little-wedderburn.bound.3}) is strict, as we
observed above). On the other hand,%
\[
\prod\limits_{\substack{S\subseteq\operatorname*{PF}\left(  n\right)
;\\\left\vert S\right\vert \text{ is odd}}}\underbrace{\left(
m^{n/\operatorname*{prod}S}-1\right)  }_{\leq m^{n/\operatorname*{prod}S}}%
\leq\prod\limits_{\substack{S\subseteq\operatorname*{PF}\left(  n\right)
;\\\left\vert S\right\vert \text{ is odd}}}m^{n/\operatorname*{prod}S}.
\]
Dividing the first of these two inequalities by the second, we find%
\begin{align}
\dfrac{\prod\limits_{\substack{S\subseteq\operatorname*{PF}\left(  n\right)
;\\\left\vert S\right\vert \text{ is even}}}\left(  m^{n/\operatorname*{prod}%
S}-1\right)  }{\prod\limits_{\substack{S\subseteq\operatorname*{PF}\left(
n\right)  ;\\\left\vert S\right\vert \text{ is odd}}}\left(
m^{n/\operatorname*{prod}S}-1\right)  }  &  >\dfrac{\prod
\limits_{\substack{S\subseteq\operatorname*{PF}\left(  n\right)  ;\\\left\vert
S\right\vert \text{ is even}}}\left(  m^{n/\operatorname*{prod}S}\cdot
\dfrac{m-1}{m}\right)  }{\prod\limits_{\substack{S\subseteq\operatorname*{PF}%
\left(  n\right)  ;\\\left\vert S\right\vert \text{ is odd}}%
}m^{n/\operatorname*{prod}S}}\nonumber\\
&  =\dfrac{\prod\limits_{\substack{S\subseteq\operatorname*{PF}\left(
n\right)  ;\\\left\vert S\right\vert \text{ is even}}%
}m^{n/\operatorname*{prod}S}\cdot\prod\limits_{\substack{S\subseteq
\operatorname*{PF}\left(  n\right)  ;\\\left\vert S\right\vert \text{ is
even}}}\dfrac{m-1}{m}}{\prod\limits_{\substack{S\subseteq\operatorname*{PF}%
\left(  n\right)  ;\\\left\vert S\right\vert \text{ is odd}}%
}m^{n/\operatorname*{prod}S}}\nonumber\\
&  =\underbrace{\dfrac{\prod\limits_{\substack{S\subseteq\operatorname*{PF}%
\left(  n\right)  ;\\\left\vert S\right\vert \text{ is even}}%
}m^{n/\operatorname*{prod}S}}{\prod\limits_{\substack{S\subseteq
\operatorname*{PF}\left(  n\right)  ;\\\left\vert S\right\vert \text{ is odd}%
}}m^{n/\operatorname*{prod}S}}}_{\substack{=m^{\phi\left(  n\right)
}\\\text{(by (\ref{pf.lem.little-wedderburn.bound.mphin}))}}}\cdot
\underbrace{\prod\limits_{\substack{S\subseteq\operatorname*{PF}\left(
n\right)  ;\\\left\vert S\right\vert \text{ is even}}}\dfrac{m-1}{m}%
}_{\substack{=\left(  \dfrac{m-1}{m}\right)  ^{2^{\left\vert
\operatorname*{PF}\left(  n\right)  \right\vert -1}}\\\text{(since this
product has }2^{\left\vert \operatorname*{PF}\left(  n\right)  \right\vert
-1}\\\text{factors (by (\ref{pf.lem.little-wedderburn.ineq2.comb2})))}%
}}\nonumber\\
&  =m^{\phi\left(  n\right)  }\cdot\left(  \dfrac{m-1}{m}\right)
^{2^{\left\vert \operatorname*{PF}\left(  n\right)  \right\vert -1}%
}\nonumber\\
&  =m^{\phi\left(  n\right)  -2^{\left\vert \operatorname*{PF}\left(
n\right)  \right\vert -1}}\cdot\left(  m-1\right)  ^{2^{\left\vert
\operatorname*{PF}\left(  n\right)  \right\vert -1}}%
.\ \ \ \ \ \ \ \ \ \ \label{pf.lem.little-wedderburn.bound.5}%
\end{align}
But Lemma \ref{lem.little-wedderburn.ineq1} yields
\[
\phi\left(  n\right)  \geq2^{\left\vert \operatorname*{PF}\left(  n\right)
\setminus\left\{  2\right\}  \right\vert }\geq2^{\left\vert \operatorname*{PF}%
\left(  n\right)  \right\vert -1}\ \ \ \ \ \ \ \ \ \ \left(  \text{since
}\left\vert \operatorname*{PF}\left(  n\right)  \setminus\left\{  2\right\}
\right\vert \geq\left\vert \operatorname*{PF}\left(  n\right)  \right\vert
-1\right)  ,
\]
and thus $\phi\left(  n\right)  -2^{\left\vert \operatorname*{PF}\left(
n\right)  \right\vert -1}\geq0$. Hence, we can take the inequality $m\geq m-1$
to the $\left(  \phi\left(  n\right)  -2^{\left\vert \operatorname*{PF}\left(
n\right)  \right\vert -1}\right)  $-th power and obtain $m^{\phi\left(
n\right)  -2^{\left\vert \operatorname*{PF}\left(  n\right)  \right\vert -1}%
}\geq\left(  m-1\right)  ^{\phi\left(  n\right)  -2^{\left\vert
\operatorname*{PF}\left(  n\right)  \right\vert -1}}$. Now,
(\ref{pf.lem.little-wedderburn.bound.1}) becomes%
\begin{align*}
\Phi_{n}\left(  m\right)   &  =\dfrac{\prod\limits_{\substack{S\subseteq
\operatorname*{PF}\left(  n\right)  ;\\\left\vert S\right\vert \text{ is
even}}}\left(  m^{n/\operatorname*{prod}S}-1\right)  }{\prod
\limits_{\substack{S\subseteq\operatorname*{PF}\left(  n\right)  ;\\\left\vert
S\right\vert \text{ is odd}}}\left(  m^{n/\operatorname*{prod}S}-1\right)
}>\underbrace{m^{\phi\left(  n\right)  -2^{\left\vert \operatorname*{PF}%
\left(  n\right)  \right\vert -1}}}_{\geq\left(  m-1\right)  ^{\phi\left(
n\right)  -2^{\left\vert \operatorname*{PF}\left(  n\right)  \right\vert -1}}%
}\cdot\left(  m-1\right)  ^{2^{\left\vert \operatorname*{PF}\left(  n\right)
\right\vert -1}}\ \ \ \ \ \ \ \ \ \ \left(  \text{by
(\ref{pf.lem.little-wedderburn.bound.5})}\right) \\
&  \geq\left(  m-1\right)  ^{\phi\left(  n\right)  -2^{\left\vert
\operatorname*{PF}\left(  n\right)  \right\vert -1}}\cdot\left(  m-1\right)
^{2^{\left\vert \operatorname*{PF}\left(  n\right)  \right\vert -1}}=\left(
m-1\right)  ^{\phi\left(  n\right)  },
\end{align*}
and thus Lemma \ref{lem.little-wedderburn.ineq2} is proved.
\end{proof}

Note that most texts prove Lemma \ref{lem.little-wedderburn.ineq2}
geometrically, using (\ref{eq.thm.cyclotomic.ngon}) and the position of roots
of unity on the unit circle (see, e.g., \cite[\S 2]{Asgarl14}, \cite[Theorem
8.5.6]{BeaBla19}, \cite[Theorem 67.7]{Elman22}, \cite[part 2, Theorem
C-2.31]{Rotman3e}). Our proof is trickier but purely algebraic.
\end{fineprint}

\subsubsection{Aside: Dynatomic polynomials}

We take this occasion to introduce a more obscure family of polynomials, known
as the \textbf{dynatomic polynomials}, whose existence too follows from the
factorization theorem for gcd-sequences (Theorem \ref{thm.gcdseq.fac}).

We first observe an easy consequence of Proposition \ref{prop.gcdseq.modcrit}:

\begin{corollary}
\label{cor.gcdseq.divcrit}Let $R$ be a commutative ring. Let $\left(
a_{1},a_{2},a_{3},\ldots\right)  $ be a sequence of elements of $R$. Assume
that%
\begin{equation}
a_{n}\mid a_{n+m}-a_{m}\ \ \ \ \ \ \ \ \ \ \text{for all }n,m\geq1.
\label{eq.cor.gcdseq.divcrit.ass}%
\end{equation}
Then, $\left(  a_{1},a_{2},a_{3},\ldots\right)  $ is a strong gcd-sequence.
\end{corollary}

\begin{fineprint}

\begin{proof}
Let $m,n\geq1$. Choose an element $b_{m,n}\in R$ satisfying
\begin{equation}
a_{n+m}-a_{m}=a_{n}b_{m,n} \label{pf.cor.gcdseq.divcrit.b}%
\end{equation}
(such an element exists because of (\ref{eq.cor.gcdseq.divcrit.ass})). Also,
define an element $c_{m,n}\in R$ by $c_{m,n}=1$.

Forget that we fixed $m,n$. We thus have defined two elements $b_{m,n}$ and
$c_{m,n}$ of $R$ for all $m,n\geq1$. These elements satisfy
(\ref{pf.cor.gcdseq.divcrit.b}) and $c_{m,n}=1$. Thus, for every $m,n\geq1$,
we have%
\[
a_{m}R+b_{m,n}R+\underbrace{c_{m,n}}_{=1}R=a_{m}R+b_{m,n}R+R=R
\]
(since the principal ideal $R$ does not change if we add further ideals to
it). Moreover, each $p,q\geq1$ satisfy%
\begin{equation}
a_{p+q}=a_{p}b_{q,p}+a_{q}c_{p,q}\ \ \ \ \ \ \ \ \ \ \text{in }R
\label{pf.cor.gcdseq.divcrit.4}%
\end{equation}
(because applying (\ref{pf.cor.gcdseq.divcrit.b}) to $n=p$ and $m=q$ yields
$a_{p+q}-a_{q}=a_{p}b_{q,p}$ and therefore
\[
a_{p+q}=a_{p}b_{q,p}+\underbrace{a_{q}}_{\substack{=a_{q}c_{p,q}\\\text{(since
}c_{p,q}=1\text{)}}}=a_{p}b_{q,p}+a_{q}c_{p,q},
\]
which is precisely (\ref{pf.cor.gcdseq.divcrit.4})). Thus, Proposition
\ref{prop.gcdseq.modcrit} shows that $\left(  a_{1},a_{2},a_{3},\ldots\right)
$ is a strong gcd-sequence. This proves Corollary \ref{cor.gcdseq.divcrit}.
\end{proof}
\end{fineprint}

The next proposition is another classical property of polynomials:

\begin{proposition}
\label{prop.polring.univar-x-y-div}Let $f$ be a polynomial in $R\left[
x\right]  $. Let $A$ be a commutative $R$-algebra. Let $a,b\in A$. Then,%
\[
a-b\mid f\left[  a\right]  -f\left[  b\right]  \ \ \ \ \ \ \ \ \ \ \text{in
}A.
\]

\end{proposition}

\begin{fineprint}

\begin{proof}
Let $g$ be the polynomial $f\left[  x\right]  -f\left[  b\right]  \in A\left[
x\right]  $. Then, $g\left[  b\right]  =f\left[  b\right]  -f\left[  b\right]
=0$, so that $b$ is a root of $g$. But Proposition
\ref{prop.polring.univar-rootdiv} (applied to $A$, $b$ and $g$ instead of $R$,
$a$ and $f$) yields that we have the following logical equivalence:%
\[
\left(  b\text{ is a root of }g\right)  \Longleftrightarrow\left(  x-b\mid
g\text{ in }A\left[  x\right]  \right)  .
\]
Hence, $x-b\mid g$ in $A\left[  x\right]  $ (since $b$ is a root of $g$). In
other words, there exists a polynomial $h\in A\left[  x\right]  $ such that
$g=\left(  x-b\right)  \cdot h$. Consider this $h$. Evaluating both sides of
$g=\left(  x-b\right)  \cdot h$ at $a$ (using Theorem
\ref{thm.polring.univar-sub-hom}), we find $g\left[  a\right]  =\left(
a-b\right)  \cdot h\left[  a\right]  $, so that
\[
a-b\mid g\left[  a\right]  =f\left[  a\right]  -f\left[  b\right]
\ \ \ \ \ \ \ \ \ \ \left(  \text{since }g=f\left[  x\right]  -f\left[
b\right]  \right)  .
\]
This proves Proposition \ref{prop.polring.univar-x-y-div}.
\end{proof}
\end{fineprint}

Using Corollary \ref{cor.gcdseq.divcrit}, we can easily see the following:

\begin{theorem}
\label{thm.gcdseq.dynatom}Let $P\in R\left[  x\right]  $ be a univariate
polynomial over a commutative ring $R$. For each $n\in\mathbb{N}$, let
$P^{\circ n}$ denote the polynomial $\underbrace{P\circ P\circ\cdots\circ
P}_{n\text{ times }P}\in R\left[  x\right]  $, where the symbol $\circ$
denotes composition of polynomials (see Subsection
\ref{subsec.polys1.polyrings.ass-subst}). (In particular, $P^{\circ0}=x$,
since $x$ is the neutral element of the operation $\circ$ on $R\left[
x\right]  $; furthermore, $P^{\circ1}=P$ and $P^{\circ2}=P\circ P$).

Let $A$ be a commutative $R$-algebra, and let $a\in A$. Let $a_{n}:=P^{\circ
n}\left[  a\right]  -a\in A$ for each $n\geq1$. Then, the sequence $\left(
a_{1},a_{2},a_{3},\ldots\right)  $ is a strong gcd-sequence in $A$.
\end{theorem}

\begin{example}
Let $P=x^{2}-2$ in $\mathbb{Z}\left[  x\right]  $. Let $A$ be a commutative
$\mathbb{Z}$-algebra (i.e., a commutative ring), and let $a\in A$. Then, the
sequence $\left(  a_{1},a_{2},a_{3},\ldots\right)  $ defined in Theorem
\ref{thm.gcdseq.dynatom} starts as follows:%
\begin{align*}
a_{1}  &  =P^{\circ1}\left[  a\right]  -a=P\left[  a\right]  -a=\left(
a^{2}-2\right)  -a=a^{2}-a-2;\\
a_{2}  &  =P^{\circ2}\left[  a\right]  -a=P\left[  P\left[  a\right]  \right]
-a=\left(  \left(  a^{2}-2\right)  ^{2}-2\right)  -a=a^{4}-4a^{2}-a+2;\\
a_{3}  &  =P^{\circ3}\left[  a\right]  -a=a^{8}-8a^{6}+20a^{4}-16a^{2}-a+2;\\
a_{4}  &  =P^{\circ4}\left[  a\right]  -a=a^{16}-16a^{14}+104a^{12}%
-352a^{10}+660a^{8}\\
&  \qquad\qquad\qquad\qquad-672a^{6}+336a^{4}-64a^{2}-a+2.
\end{align*}
The polynomials $P^{\circ n}$ clearly grow fast (their degrees get doubled
from $n$ to $n+1$).
\end{example}

\begin{fineprint}

\begin{proof}
[Proof of Theorem \ref{thm.gcdseq.dynatom}.]By Corollary
\ref{cor.gcdseq.divcrit} (applied to $A$ instead of $R$), it will suffice to
show that
\[
a_{n}\mid a_{n+m}-a_{m}\ \ \ \ \ \ \ \ \ \ \text{for all }n,m\geq1.
\]
So let us show this. Let $n,m\geq1$. Then, $a_{n+m}=P^{\circ\left(
n+m\right)  }\left[  a\right]  -a$ and $a_{m}=P^{\circ m}\left[  a\right]
-a$, so that%
\begin{align}
a_{n+m}-a_{m}  &  =\left(  P^{\circ\left(  n+m\right)  }\left[  a\right]
-a\right)  -\left(  P^{\circ m}\left[  a\right]  -a\right)
=\underbrace{P^{\circ\left(  n+m\right)  }}_{\substack{=P^{\circ\left(
m+n\right)  }\\=P^{\circ m}\circ P^{\circ n}\\=P^{\circ m}\left[  P^{\circ
n}\right]  }}\left[  a\right]  -P^{\circ m}\left[  a\right] \nonumber\\
&  =\underbrace{\left(  P^{\circ m}\left[  P^{\circ n}\right]  \right)
\left[  a\right]  }_{\substack{=P^{\circ m}\left[  P^{\circ n}\left[
a\right]  \right]  \\\text{(by Theorem \ref{thm.polys1.comp-assoc}
\textbf{(a)})}}}-\,P^{\circ m}\left[  a\right] \nonumber\\
&  =P^{\circ m}\left[  P^{\circ n}\left[  a\right]  \right]  -P^{\circ
m}\left[  a\right]  . \label{pf.thm.gcdseq.dynatom.1}%
\end{align}
Furthermore, the definition of $a_{n}$ shows that
\begin{equation}
a_{n}=P^{\circ n}\left[  a\right]  -a. \label{pf.thm.gcdseq.dynatom.2}%
\end{equation}
But Proposition \ref{prop.polring.univar-x-y-div} (applied to $P^{\circ
n}\left[  a\right]  $, $a$ and $P^{\circ m}$ instead of $a$, $b$ and $f$)
yields
\[
P^{\circ n}\left[  a\right]  -a\mid P^{\circ m}\left[  P^{\circ n}\left[
a\right]  \right]  -P^{\circ m}\left[  a\right]  \ \ \ \ \ \ \ \ \ \ \text{in
}A.
\]
In view of (\ref{pf.thm.gcdseq.dynatom.1}) and (\ref{pf.thm.gcdseq.dynatom.2}%
), we can rewrite this as%
\[
a_{n}\mid a_{n+m}-a_{m}\ \ \ \ \ \ \ \ \ \ \text{in }A.
\]
This is precisely what we needed to show. Thus, the proof of Theorem
\ref{thm.gcdseq.dynatom} is complete.
\end{proof}
\end{fineprint}

In the situation of Theorem \ref{thm.gcdseq.dynatom}, we can thus apply
Theorem \ref{thm.gcdseq.fac} to the strong gcd-sequence $\left(  a_{1}%
,a_{2},a_{3},\ldots\right)  $ (using $A$ instead of $R$), obtaining a sequence
$\left(  b_{1},b_{2},b_{3},\ldots\right)  $ of elements of $A$ such that each
$n\geq1$ satisfies%
\[
a_{n}=\prod_{d\mid n}b_{d}.
\]
These elements $b_{1},b_{2},b_{3},\ldots$ are not necessarily unique, but (by
Proposition \ref{prop.gcdseq.reg} \textbf{(a)}) they are so when $a_{1}%
,a_{2},a_{3},\ldots$ are regular in $A$. This regularity condition holds in
particular when $A=R\left[  x\right]  $ and $a=x$ (so that $a_{n}=P^{\circ
n}\left[  a\right]  -a$ simplifies to $a_{n}=P^{\circ n}-x$) and $P$ is a
monic polynomial of degree $\geq2$:

\begin{exercise}
\label{exe.gcdseq.dynatom.reg}Let $P\in R\left[  x\right]  $ be a monic
polynomial of degree $\geq2$ over a commutative ring $R$. Let $P^{\circ n}$ be
defined as in Theorem \ref{thm.gcdseq.dynatom}. Prove that the polynomials
$P^{\circ n}-x\in R\left[  x\right]  $ are monic and therefore regular for all
$n\geq1$. Conclude that there is a unique sequence $\left(  b_{1},b_{2}%
,b_{3},\ldots\right)  $ of elements of $R\left[  x\right]  $ such that each
$n\geq1$ satisfies%
\begin{equation}
P^{\circ n}-x=\prod_{d\mid n}b_{d}. \label{eq.exe.gcdseq.dynatom.reg.1}%
\end{equation}
Moreover, conclude that these $b_{n}$ are given by%
\begin{equation}
b_{n}=\dfrac{\prod\limits_{\substack{S\subseteq\operatorname*{PF}\left(
n\right)  ;\\\left\vert S\right\vert \text{ is even}}}\left(  P^{\circ\left(
n/\operatorname*{prod}S\right)  }-x\right)  }{\prod
\limits_{\substack{S\subseteq\operatorname*{PF}\left(  n\right)  ;\\\left\vert
S\right\vert \text{ is odd}}}\left(  P^{\circ\left(  n/\operatorname*{prod}%
S\right)  }-x\right)  }\ \ \ \ \ \ \ \ \ \ \text{for all }n\geq1
\label{eq.exe.gcdseq.dynatom.reg.2}%
\end{equation}
(where we are using Definition \ref{def.ent.PFn} and Definition
\ref{def.ent.prodS}). These polynomials $b_{n}$ are called the
\textbf{dynatomic polynomials} of $P$, studied (with a definition equivalent
to (\ref{eq.exe.gcdseq.dynatom.reg.2})) e.g. in \cite[Definition 3]{BFLS13}
and \cite{DoFiHy21} (see also \cite[\S 4.1]{Silver07}). They have applications
in arithmetic dynamics (e.g., the study of what happens when a polynomial is
applied to a number over and over again).
\end{exercise}

\begin{exercise}
Let $k$ be a positive integer.

\begin{enumerate}
\item[\textbf{(a)}] Let $\left(  a_{1},a_{2},a_{3},\ldots\right)  $ be the
sequence of integers defined recursively by $a_{1}=a_{2}=\cdots=a_{k}=1$ and
$a_{n}=1+a_{n-1}a_{n-2}\cdots a_{n-k}$ for all $n>k$. Prove that $\left(
a_{1},a_{2},a_{3},\ldots\right)  $ is a strong gcd-sequence.

\item[\textbf{(b)}] Generalize this, e.g., by replacing the recursion
$a_{n}=1+a_{n-1}a_{n-2}\cdots a_{n-k}$ by $a_{n}=1+a_{n-1}^{i_{1}}%
a_{n-2}^{i_{2}}\cdots a_{n-k}^{i_{k}}$ for fixed positive integers
$i_{1},i_{2},\ldots,i_{k}$.
\end{enumerate}
\end{exercise}

\newpage

\section{Modules over a PID (specifically, over $\mathbb{Z}$)}

Modules over a field are rather well-behaved: they are all free, i.e., they
have bases and thus are isomorphic to \textquotedblleft direct sum
powers\textquotedblright\ of the field.

Modules over an arbitrary ring can be rather wild.

Studying modules over a PID is a middle ground: they are not that wild, but
still sufficiently frequent in \textquotedblleft real life\textquotedblright.

I will just give a taste of their theory. The only PID I will work with is
$\mathbb{Z}$, and the only modules I will classify are finite, but you will
see some germs of more general arguments in my brief treatment of this rather
special case.

\subsection{Classifying finite abelian groups}

\subsubsection{The classification theorem}

Classifying finite groups is notoriously hard. Even the so-called
\textquotedblleft simple\textquotedblright\ groups have
\href{https://en.wikipedia.org/wiki/Classification_of_finite_simple_groups}{a
classification that spans a page} (and takes a dozen books to prove). The
finite \textbf{abelian} groups, on the other hand, do have a rather manageable classification:

\begin{theorem}
[Classification of finite abelian groups]\label{thm.modPID.class-finabgr}Let
$G$ be a finite abelian group.

\begin{enumerate}
\item[\textbf{(a)}] Then, $G$ is isomorphic to a direct product of finitely
many finite cyclic groups.

In other words, there exist positive integers $n_{1},n_{2},\ldots,n_{k}$ such
that%
\[
G\cong\left(  \mathbb{Z}/n_{1}\right)  \times\left(  \mathbb{Z}/n_{2}\right)
\times\cdots\times\left(  \mathbb{Z}/n_{k}\right)  .
\]

\item[\textbf{(b)}] Moreover, we can choose these $n_{1},n_{2},\ldots,n_{k}$
in such a way that they are $>1$ and satisfy%
\[
n_{1}\mid n_{2}\mid\cdots\mid n_{k}.
\]

\item[\textbf{(c)}] Finally, if we choose them in such a way, then they are unique.
\end{enumerate}
\end{theorem}

I will outline a proof of parts \textbf{(a)} and \textbf{(b)} of this theorem
using modules over $\mathbb{Z}$. (There are other proofs, e.g., using group
theory -- see, e.g., \cite[Theorem 14.11]{Elman22}.) A proof of Theorem
\ref{thm.modPID.class-finabgr} \textbf{(c)} is outlined in Exercise
\ref{exe.modPID.class-finabgr.c.3} below.

\subsubsection{On modules and matrices}

How do modules come into play here in the first place? Recall from Proposition
\ref{prop.Zmods.1} that abelian groups are $\mathbb{Z}$-modules\footnote{And
we know from Proposition \ref{prop.modules.mors.Z} that the group morphisms
between these abelian groups are exactly the $\mathbb{Z}$-module morphisms.};
thus, classifying finite abelian groups is the same as classifying finite
$\mathbb{Z}$-modules.

One other thing that will be crucial is good old matrices. Recall from linear
algebra that matrices over a field $F$ correspond to linear maps between
$F$-vector spaces. Likewise, matrices over an arbitrary commutative ring $R$
correspond to linear maps between free $R$-modules. Specifically:

\begin{convention}
For any commutative ring $R$ and any $n\in\mathbb{N}$, we identify the
$n$-tuples $\left(  a_{1},a_{2},\ldots,a_{n}\right)  \in R^{n}$ with the
column vectors $\left(
\begin{array}
[c]{c}%
a_{1}\\
a_{2}\\
\vdots\\
a_{n}%
\end{array}
\right)  \in R^{n\times1}$. Thus, $R^{n}$ becomes the $R$-module $R^{n\times
1}$ of column vectors of size $n$.
\end{convention}

\begin{proposition}
\label{prop.modPID.la.mat-Rn}Let $R$ be a commutative ring.

\begin{enumerate}
\item[\textbf{(a)}] If $A\in R^{n\times m}$ is an $n\times m$-matrix over $R$,
then the map%
\begin{align}
R^{m}  &  \rightarrow R^{n},\nonumber\\
v  &  \mapsto Av \label{eq.prop.modPID.la.mat-Rn.map}%
\end{align}
is an $R$-linear map.

\item[\textbf{(b)}] Moreover, any $R$-linear map from $R^{m}$ to $R^{n}$ has
the form (\ref{eq.prop.modPID.la.mat-Rn.map}) for a unique $n\times m$-matrix
$A\in R^{n\times m}$.

\item[\textbf{(c)}] Thus, there is a 1-to-1 correspondence between $n\times
m$-matrices over $R$ and linear maps from $R^{m}$ to $R^{n}$.
\end{enumerate}
\end{proposition}

\begin{proof}
As in linear algebra. (See \cite[Theorem 6.8.4]{19s} for part \textbf{(a)},
and see \cite[Proposition 6.8.5]{19s} for part \textbf{(b)}.)
\end{proof}

\begin{definition}
\label{def.modPID.la.ColA}Let $R$ be a commutative ring. Let $A\in R^{n\times
m}$ be an $n\times m$-matrix over $R$.

\begin{enumerate}
\item[\textbf{(a)}] We set%
\begin{align*}
\operatorname*{Col}A:  &  =\left\{  Av\ \mid\ v\in R^{m}\right\} \\
&  =\left(  \text{the image of the linear map
(\ref{eq.prop.modPID.la.mat-Rn.map})}\right) \\
&  =\left(  \text{the span of the columns of }A\right)  .
\end{align*}
This is an $R$-submodule of $R^{n}$, and is called the \textbf{column space}
of $A$. (This is all exactly as in linear algebra.)

\item[\textbf{(b)}] The \textbf{cokernel} of $A$ is defined to be the quotient
$R$-module $R^{n}/\operatorname*{Col}A$.
\end{enumerate}
\end{definition}

\begin{definition}
\label{def.modPID.la.finpres}Let $R$ be a commutative ring. An $R$-module is
said to be \textbf{finitely presented} if it is isomorphic to the cokernel of
some matrix over $R$.
\end{definition}

\begin{remark}
This latter definition might appear somewhat random. Here is some intuition
for those who know a bit about groups, specifically
\href{https://en.wikipedia.org/wiki/Presentation_of_a_group}{about their
presentations}. An $R$-module is finitely presented if it can be
\textquotedblleft defined by finitely many generators and finitely many
relations\textquotedblright. For example, recall that the $R$-module $R^{4}$
can be viewed as the $R$-module consisting of all \textquotedblleft
formal\textquotedblright\ $R$-linear combinations $ax+by+cz+dw$ of four
independent symbols $x,y,z,w$. Likewise, the $R$-module%
\[
R^{4}/\operatorname*{Col}A\ \ \ \ \ \ \ \ \ \ \text{for }A=\left(
\begin{array}
[c]{cc}%
3 & 2\\
4 & 7\\
-5 & 0\\
-6 & -4
\end{array}
\right)
\]
can be expressed as the $R$-module consisting of all \textquotedblleft
formal\textquotedblright\ $R$-linear combinations $ax+by+cz+dw$ but subject to
the relations $3x+4y=5z+6w$ and $2x+7y=4w$. Here, the \textquotedblleft
generators\textquotedblright\ $x,y,z,w$ are the cosets $e_{1}%
+\operatorname*{Col}A,\ \ e_{2}+\operatorname*{Col}A,\ \ e_{3}%
+\operatorname*{Col}A,\ \ e_{4}+\operatorname*{Col}A$ of the four standard
basis elements $e_{1},e_{2},e_{3},e_{4}$ of $R^{4}$; they satisfy the
relations $3x+4y=5z+6w$ and $2x+7y=4w$ because we have factored out the
submodule
\begin{align*}
\operatorname*{Col}A  &  =\operatorname*{span}\left(  \left(
\begin{array}
[c]{c}%
3\\
4\\
-5\\
-6
\end{array}
\right)  ,\left(
\begin{array}
[c]{c}%
2\\
7\\
0\\
-4
\end{array}
\right)  \right) \\
&  =\operatorname*{span}\left(  3e_{1}+4e_{2}-5e_{3}-6e_{4},2e_{1}%
+7e_{2}-4e_{4}\right)  .
\end{align*}

\end{remark}

\subsubsection{Every finite $\mathbb{Z}$-module is finitely presented}

Our first step towards the classification theorem is the following:

\begin{lemma}
\label{lem.modPID.fin-Z-is-fp}Let $G$ be a finite $\mathbb{Z}$-module.
(\textquotedblleft Finite\textquotedblright\ means that the set $G$ is
finite.) Then, $G$ is finitely presented.
\end{lemma}

\begin{proof}
The set $G$ is finite and nonempty (since it contains $0$); thus, its size
$\left\vert G\right\vert $ is a positive integer. Let us denote this positive
integer by $n$.

The abelian group $\left(  G,+,0\right)  $ is finite; thus, Lagrange's theorem
yields that $\left\vert G\right\vert \cdot a=0$ for each $a\in G$. In other
words,
\begin{equation}
na=0\ \ \ \ \ \ \ \ \ \ \text{for each }a\in G
\label{pf.lem.modPID.fin-Z-is-fp.na=0}%
\end{equation}
(since $n=\left\vert G\right\vert $).

Let $\left(  m_{1},m_{2},\ldots,m_{n}\right)  $ be a list of all the $n$
elements of $G$ (each listed exactly once). Thus, $G=\left\{  m_{1}%
,m_{2},\ldots,m_{n}\right\}  $.

Consider the free $\mathbb{Z}$-module $\mathbb{Z}^{n}$ with its standard basis
$\left(  e_{1},e_{2},\ldots,e_{n}\right)  $. The map%
\begin{align*}
f:\mathbb{Z}^{n}  &  \rightarrow G,\\
\left(  r_{1},r_{2},\ldots,r_{n}\right)   &  \mapsto r_{1}m_{1}+r_{2}%
m_{2}+\cdots+r_{n}m_{n}%
\end{align*}
is a $\mathbb{Z}$-module morphism (according to Theorem
\ref{thm.mods.free-is-iso-2} \textbf{(a)}). Moreover, this map $f$ satisfies
$f\left(  e_{i}\right)  =m_{i}$ for each $i\in\left\{  1,2,\ldots,n\right\}
$, and thus its image contains all of $m_{1},m_{2},\ldots,m_{n}$; thus, this
map $f$ is surjective (since $G=\left\{  m_{1},m_{2},\ldots,m_{n}\right\}  $).
The First isomorphism theorem for modules (Theorem \ref{thm.1it.mod1}
\textbf{(f)}, applied to $M=\mathbb{Z}^{n}$ and $N=G$) yields%
\begin{equation}
\mathbb{Z}^{n}/\operatorname*{Ker}f\cong f\left(  \mathbb{Z}^{n}\right)
=G\ \ \ \ \ \ \ \ \ \ \left(  \text{since }f\text{ is surjective}\right)  .
\label{pf.lem.modPID.fin-Z-is-fp.3}%
\end{equation}

Now, we shall construct an $n\times k$-matrix (for some $k\in\mathbb{N}$)
satisfying $\operatorname*{Ker}f=\operatorname*{Col}A$.

Indeed, we consider the following two kinds of vectors in $\mathbb{Z}^{n}$:

\begin{itemize}
\item The $n$\textbf{-stretched basis vectors} shall mean the $n$ vectors
$ne_{1},ne_{2},\ldots,ne_{n}$. These $n$ vectors belong to
$\operatorname*{Ker}f$, since each $i\in\left\{  1,2,\ldots,n\right\}  $
satisfies%
\begin{align*}
f\left(  ne_{i}\right)   &  =nm_{i}\ \ \ \ \ \ \ \ \ \ \left(  \text{by the
definition of }f\right) \\
&  =0\ \ \ \ \ \ \ \ \ \ \left(  \text{by
(\ref{pf.lem.modPID.fin-Z-is-fp.na=0}), applied to }a=m_{i}\right)
\end{align*}
and thus $ne_{i}\in\operatorname*{Ker}f$.

\item The \textbf{reduced kernel vectors} shall mean the vectors
\[
\left(  r_{1},r_{2},\ldots,r_{n}\right)  \in\left\{  0,1,\ldots,n-1\right\}
^{n}%
\]
that belong to $\operatorname*{Ker}f$. There are finitely many such vectors,
since the set $\left\{  0,1,\ldots,n-1\right\}  ^{n}$ is finite.
\end{itemize}

We have just shown that all $n$-stretched basis vectors and all reduced kernel
vectors belong to $\operatorname*{Ker}f$. Hence, any $\mathbb{Z}$-linear
combination of $n$-stretched basis vectors and reduced kernel vectors belongs
to $\operatorname*{Ker}f$ (because $\operatorname*{Ker}f$ is a $\mathbb{Z}%
$-submodule of $\mathbb{Z}^{n}$, and thus is closed under linear combination).
Conversely, using division with remainder, it is not hard to see that any
vector in $\operatorname*{Ker}f$ is a $\mathbb{Z}$-linear combination of
$n$-stretched basis vectors and reduced kernel
vectors\footnote{\textit{Proof.} Let $v=\left(  v_{1},v_{2},\ldots
,v_{n}\right)  $ be a vector in $\operatorname*{Ker}f$. We must show that $v$
is a $\mathbb{Z}$-linear combination of $n$-stretched basis vectors and
reduced kernel vectors.
\par
For each $i\in\left\{  1,2,\ldots,n\right\}  $, we write $v_{i}=q_{i}n+r_{i}$,
where $q_{i}$ and $r_{i}$ are the quotient and the remainder obtained when
dividing $v_{i}$ by $n$. Then,
\begin{align*}
v  &  =\left(  v_{1},v_{2},\ldots,v_{n}\right)  =\left(  q_{1}n+r_{1}%
,\ \ q_{2}n+r_{2},\ \ \ldots,\ \ q_{n}n+r_{n}\right) \\
&  =q_{1}ne_{1}+q_{2}ne_{2}+\cdots+q_{n}ne_{n}+\left(  r_{1},r_{2}%
,\ldots,r_{n}\right)  ,
\end{align*}
so that%
\[
\left(  r_{1},r_{2},\ldots,r_{n}\right)  =v-\left(  q_{1}ne_{1}+q_{2}%
ne_{2}+\cdots+q_{n}ne_{n}\right)  \in\operatorname*{Ker}f
\]
(since the vector $v$ as well as all the $n$ vectors $ne_{1},ne_{2}%
,\ldots,ne_{n}$ belong to $\operatorname*{Ker}f$, and since
$\operatorname*{Ker}f$ is a $\mathbb{Z}$-submodule of $\mathbb{Z}^{n}$). Thus,
$\left(  r_{1},r_{2},\ldots,r_{n}\right)  $ is a reduced kernel vector (since
the definition of the $r_{i}$ as remainders ensures that $r_{i}\in\left\{
0,1,\ldots,n-1\right\}  $ for all $i$, and thus $\left(  r_{1},r_{2}%
,\ldots,r_{n}\right)  \in\left\{  0,1,\ldots,n-1\right\}  ^{n}$). Thus, from%
\[
v=q_{1}ne_{1}+q_{2}ne_{2}+\cdots+q_{n}ne_{n}+\left(  r_{1},r_{2},\ldots
,r_{n}\right)  ,
\]
we conclude that $v$ is a $\mathbb{Z}$-linear combination of $n$-stretched
basis vectors and reduced kernel vectors. Qed.}. Hence, $\operatorname*{Ker}f$
is precisely the set of all $\mathbb{Z}$-linear combinations of $n$-stretched
basis vectors and reduced kernel vectors. In other words, $\operatorname*{Ker}%
f$ is the span of the vectors we just mentioned.

Now, let $A$ be the matrix whose columns are precisely the $n$-stretched basis
vectors and the reduced kernel vectors. (This is well-defined, since there are
only finitely many of these vectors.) Then, $\operatorname*{Col}A$ is the span
of the vectors we just mentioned. But we have seen in the previous paragraph
that $\operatorname*{Ker}f$ is the span of these vectors. Comparing these two
results, we conclude that $\operatorname*{Ker}f=\operatorname*{Col}A$. Hence,
(\ref{pf.lem.modPID.fin-Z-is-fp.3}) rewrites as%
\[
\mathbb{Z}^{n}/\operatorname*{Col}A\cong G.
\]
In other words, $G$ is isomorphic to the cokernel of $A$. Hence, $G$ is
finitely presented. This proves Lemma \ref{lem.modPID.fin-Z-is-fp}.
\end{proof}

\subsubsection{Understanding cokernels of diagonal matrices}

Recall that we still want to prove Theorem \ref{thm.modPID.class-finabgr}
\textbf{(a)}, which claims that every finite $\mathbb{Z}$-module $G$ is
isomorphic to a direct product of finitely many finite cyclic groups. Lemma
\ref{lem.modPID.fin-Z-is-fp} shows that $G$ is finitely presented. How does
this help us?

Well, $G$ is finitely presented, i.e., isomorphic to the cokernel of a matrix.
If this matrix happens to be diagonal, then we are basically done! Indeed, for
example, here is how the cokernel of a diagonal $3\times3$-matrix looks like:%
\begin{align*}
&  \mathbb{Z}^{3}/\operatorname*{Col}\left(
\begin{array}
[c]{ccc}%
a & 0 & 0\\
0 & b & 0\\
0 & 0 & c
\end{array}
\right) \\
&  =\mathbb{Z}^{3}/\operatorname*{span}\left(  \left(
\begin{array}
[c]{c}%
a\\
0\\
0
\end{array}
\right)  ,\left(
\begin{array}
[c]{c}%
0\\
b\\
0
\end{array}
\right)  ,\left(
\begin{array}
[c]{c}%
0\\
0\\
c
\end{array}
\right)  \right) \\
&  =\mathbb{Z}^{3}/\operatorname*{span}\left(  ae_{1},be_{2},ce_{3}\right) \\
&  \ \ \ \ \ \ \ \ \ \ \ \ \ \ \ \ \ \ \ \ \left(  \text{where }e_{1}%
,e_{2},e_{3}\text{ are the standard basis vectors of }\mathbb{Z}^{3}\right) \\
&  \cong\left(  \mathbb{Z}/a\right)  \times\left(  \mathbb{Z}/b\right)
\times\left(  \mathbb{Z}/c\right)  .
\end{align*}
(The \textquotedblleft$\cong$\textquotedblright\ sign here is a nice exercise
in understanding quotients of modules. Explicitly, it stems from the map%
\begin{align*}
\mathbb{Z}^{3}/\operatorname*{span}\left(  ae_{1},be_{2},ce_{3}\right)   &
\rightarrow\left(  \mathbb{Z}/a\right)  \times\left(  \mathbb{Z}/b\right)
\times\left(  \mathbb{Z}/c\right)  ,\\
\overline{\left(  u,v,w\right)  }  &  \mapsto\left(  \overline{u},\overline
{v},\overline{w}\right)  ,
\end{align*}
which is easily seen to be a $\mathbb{Z}$-module isomorphism. The intuition is
simply that when we take the quotient of the free $\mathbb{Z}$-module
$\mathbb{Z}^{3}$ by its submodule $\operatorname*{span}\left(  ae_{1}%
,be_{2},ce_{3}\right)  $, we end up identifying any two vectors $\left(
u,v,w\right)  $ and $\left(  u^{\prime},v^{\prime},w^{\prime}\right)  $ that
satisfy $u\equiv u^{\prime}\operatorname{mod}a$ and $v\equiv v^{\prime
}\operatorname{mod}b$ and $w\equiv w^{\prime}\operatorname{mod}c$; but this is
tantamount to replacing the first entry of our vector by a residue class
modulo $a$, the second by a residue class modulo $b$, and the third by a
residue class modulo $c$.)

Usually, the matrix whose cokernel we need will be rectangular, not square;
however, even for rectangular matrices there is a notion of being diagonal:

\begin{definition}
\label{def.modPID.la.rect-diag}Let $R$ be a ring. A rectangular matrix $A\in
R^{n\times m}$ is said to be \textbf{diagonal} if its $\left(  i,j\right)
$-th entry is $0$ whenever $i\neq j$.
\end{definition}

This is a looser notion of \textquotedblleft diagonal\textquotedblright\ than
the one you learnt in linear algebra, since we are not requiring that $n=m$.
For example, a diagonal $2\times4$-matrix looks like $\left(
\begin{array}
[c]{cccc}%
a & 0 & 0 & 0\\
0 & b & 0 & 0
\end{array}
\right)  $, whereas a diagonal $4\times2$-matrix looks like $\left(
\begin{array}
[c]{cc}%
a & 0\\
0 & b\\
0 & 0\\
0 & 0
\end{array}
\right)  $.

\begin{proposition}
\label{prop.modPID.la.cok-diag}Let $A\in\mathbb{Z}^{n\times m}$ be diagonal.
Then, its cokernel $\mathbb{Z}^{n}/\operatorname*{Col}A$ is isomorphic to a
direct product of finitely many cyclic groups (which, however, are not
necessarily finite).
\end{proposition}

\begin{proof}
[Proof of Proposition \ref{prop.modPID.la.cok-diag} (sketched).]We give a
\textquotedblleft proof by example\textquotedblright, or rather a proof by two
(hopefully representative) examples:%
\begin{align*}
\mathbb{Z}^{2}/\operatorname*{Col}\left(
\begin{array}
[c]{cccc}%
a & 0 & 0 & 0\\
0 & b & 0 & 0
\end{array}
\right)   &  =\mathbb{Z}^{2}/\operatorname*{span}\left(  ae_{1},be_{2}%
,0,0\right)  =\mathbb{Z}^{2}/\operatorname*{span}\left(  ae_{1},be_{2}\right)
\\
&  \cong\left(  \mathbb{Z}/a\right)  \times\left(  \mathbb{Z}/b\right)
\end{align*}
and%
\[
\mathbb{Z}^{4}/\operatorname*{Col}\left(
\begin{array}
[c]{cc}%
a & 0\\
0 & b\\
0 & 0\\
0 & 0
\end{array}
\right)  =\mathbb{Z}^{4}/\operatorname*{span}\left(  ae_{1},be_{2}\right)
\cong\left(  \mathbb{Z}/a\right)  \times\left(  \mathbb{Z}/b\right)
\times\mathbb{Z}\times\mathbb{Z}.
\]

\end{proof}

\subsubsection{The proof strategy}

This suggests a somewhat daring strategy for proving parts \textbf{(a)} and
\textbf{(b)} of Theorem \ref{thm.modPID.class-finabgr}:

\begin{enumerate}
\item Let $G$ be a finite abelian group. Thus, $G$ is a finite $\mathbb{Z}$-module.

\item By Lemma \ref{lem.modPID.fin-Z-is-fp}, the $\mathbb{Z}$-module $G$ is
finitely presented. In other words, there is a matrix $A\in\mathbb{Z}^{n\times
m}$ (for some $m\in\mathbb{N}$) such that $G\cong\mathbb{Z}^{n}%
/\operatorname*{Col}A$.

\item Tweaking this matrix $A$ in a strategic way, we can make it diagonal
without changing $\mathbb{Z}^{n}/\operatorname*{Col}A$ too much (to be
precise: $\mathbb{Z}^{n}/\operatorname*{Col}A$ stays isomorphic to $G$).

\item Then, we use Proposition \ref{prop.modPID.la.cok-diag} to argue that
$\mathbb{Z}^{n}/\operatorname*{Col}A$ is isomorphic to a direct product of
finitely many cyclic groups (which are not necessarily finite).

\item We notice that these cyclic groups must be finite, because their direct
product is finite (after all, this direct product is isomorphic to $G$, which
is finite).

\item Thus, $G\cong\left(  \mathbb{Z}/n_{1}\right)  \times\left(
\mathbb{Z}/n_{2}\right)  \times\cdots\times\left(  \mathbb{Z}/n_{k}\right)  $
for some positive integers $n_{1},n_{2},\ldots,n_{k}$. (This proves Theorem
\ref{thm.modPID.class-finabgr} \textbf{(a)}.)

\item We WLOG assume that $n_{1},n_{2},\ldots,n_{k}$ are $>1$, since any
$n_{i}$ that equals $1$ only contributes a trivial factor $\mathbb{Z}/1$ to
the direct product $\left(  \mathbb{Z}/n_{1}\right)  \times\left(
\mathbb{Z}/n_{2}\right)  \times\cdots\times\left(  \mathbb{Z}/n_{k}\right)  $
(and of course such a factor can simply be removed from the product).

\item Finally, by fudging the $n_{1},n_{2},\ldots,n_{k}$ appropriately, we
ensure that $n_{1}\mid n_{2}\mid\cdots\mid n_{k}$. (This proves Theorem
\ref{thm.modPID.class-finabgr} \textbf{(b)}.)
\end{enumerate}

Steps 1, 2, 4, 5, 6, 7 should be rather clear by now. But Steps 3 and 8 sound
rather ambitious. How can we turn an arbitrary matrix into a diagonal one? How
can we pull $n_{1}\mid n_{2}\mid\cdots\mid n_{k}$ out of thin air?

\subsubsection{Row and column operations and equivalent matrices}

To make Step 3 a reality, the tools of choice are \textbf{row operations} and
\textbf{column operations}. These are a mild generalization of the row and
column operations that you know from linear algebra. Here is one way to define them:

\begin{definition}
\label{def.modPID.la.rowops}\ \ 

\begin{enumerate}
\item[\textbf{(a)}] A square matrix $A\in\mathbb{Z}^{k\times k}$ is said to be
\textbf{invertible} if it has an inverse matrix in $\mathbb{Z}^{k\times k}$
(that is, an inverse matrix with integer entries). In other words, it is said
to be invertible if it is a unit of the matrix ring $\mathbb{Z}^{k\times k}$.

For example, $\left(
\begin{array}
[c]{cc}%
1 & 1\\
1 & -1
\end{array}
\right)  \in\mathbb{Z}^{2\times2}$ is not invertible. It has an inverse in
$\mathbb{Q}^{2\times2}$, but that doesn't count!

\item[\textbf{(b)}] A \textbf{row operation} means an operation transforming a
matrix $A\in\mathbb{Z}^{n\times m}$ into $BA$, where $B\in\mathbb{Z}^{n\times
n}$ is some invertible $n\times n$-matrix.

\item[\textbf{(c)}] A \textbf{column operation} means an operation
transforming a matrix $A\in\mathbb{Z}^{n\times m}$ into $AC$, where
$C\in\mathbb{Z}^{m\times m}$ is some invertible $m\times m$-matrix.

\item[\textbf{(d)}] Two matrices $A,A^{\prime}\in\mathbb{Z}^{n\times m}$ are
said to be \textbf{equivalent} if there exist invertible matrices
$B\in\mathbb{Z}^{n\times n}$ and $C\in\mathbb{Z}^{m\times m}$ such that
$A^{\prime}=BAC$. In other words, $A,A^{\prime}$ are said to be equivalent if
$A$ can be transformed into $A^{\prime}$ using row and column operations.
\end{enumerate}
\end{definition}

You know all these notions in the case of a field; we are just adapting them
to the case of $\mathbb{Z}$.

\begin{remark}
\label{rmk.modPID.la.rowops1}\ \ 

\begin{enumerate}
\item[\textbf{(a)}] Any row operation can be undone by another row operation.

\item[\textbf{(b)}] Adding a multiple of a row to another row is a row operation.

\item[\textbf{(c)}] Swapping two rows is a row operation.

\item[\textbf{(d)}] Scaling a row by $-1$ is a row operation. (But scaling a
row by $2$ is not!)

\item[\textbf{(e)}] The analogues of all these statements for columns instead
of rows hold.
\end{enumerate}
\end{remark}

\begin{proof}
Part \textbf{(a)} is obvious (since multiplying a matrix by an invertible
matrix $B$ can be undone by multiplying it by $B^{-1}$). The other parts are
proved as in linear algebra.
\end{proof}

\begin{proposition}
\label{prop.modPID.la.congruent-cok-iso}If two matrices $A,A^{\prime}%
\in\mathbb{Z}^{n\times m}$ are equivalent, then their cokernels $\mathbb{Z}%
^{n}/\operatorname*{Col}A$ and $\mathbb{Z}^{n}/\operatorname*{Col}A^{\prime}$
are isomorphic.
\end{proposition}

\begin{proof}
Let $A,A^{\prime}\in\mathbb{Z}^{n\times m}$ be two matrices that are
equivalent. Thus, there exist invertible matrices $B\in\mathbb{Z}^{n\times n}$
and $C\in\mathbb{Z}^{m\times m}$ such that $A^{\prime}=BAC$. Consider these
$B$ and $C$.

I claim that the map%
\begin{align*}
f:\mathbb{Z}^{n}/\operatorname*{Col}A  &  \rightarrow\mathbb{Z}^{n}%
/\operatorname*{Col}A^{\prime},\\
\overline{v}  &  \mapsto\overline{Bv}%
\end{align*}
is well-defined and is a $\mathbb{Z}$-module isomorphism.

First of all, let me prove that $f$ is well-defined. Indeed, let
$v,w\in\mathbb{Z}^{n}$ be such that $\overline{v}=\overline{w}$ in
$\mathbb{Z}^{n}/\operatorname*{Col}A$. We must prove that $\overline
{Bv}=\overline{Bw}$ in $\mathbb{Z}^{n}/\operatorname*{Col}A^{\prime}$.

From $\overline{v}=\overline{w}$ in $\mathbb{Z}^{n}/\operatorname*{Col}A$, we
obtain $v-w\in\operatorname*{Col}A$. In other words, $v-w=Au$ for some
$u\in\mathbb{Z}^{m}$ (since $\operatorname*{Col}A=\left\{  Au\ \mid
\ u\in\mathbb{Z}^{m}\right\}  $). Consider this $u$. We have $C^{-1}%
\in\mathbb{Z}^{m\times m}$ (since $C$ is invertible) and thus $C^{-1}%
u\in\mathbb{Z}^{m}$. Now,%
\[
Bv-Bw=B\underbrace{\left(  v-w\right)  }_{=Au}=\underbrace{BA}%
_{\substack{=A^{\prime}C^{-1}\\\text{(since }BAC=A^{\prime}\text{)}%
}}u=A^{\prime}\underbrace{C^{-1}u}_{\in\mathbb{Z}^{m}}\in\operatorname*{Col}%
A^{\prime}%
\]
(since $\operatorname*{Col}A^{\prime}=\left\{  A^{\prime}z\ \mid
\ z\in\mathbb{Z}^{m}\right\}  $). In other words, $\overline{Bv}=\overline
{Bw}$ in $\mathbb{Z}^{n}/\operatorname*{Col}A^{\prime}$, which is precisely
what we wanted to show.

Thus, we have shown that $f$ is well-defined.

It is straightforward to see that $f$ is a $\mathbb{Z}$-module morphism. Next,
in order to show that $f$ is invertible, I will construct an inverse.

Indeed, I claim that the map%
\begin{align*}
g:\mathbb{Z}^{n}/\operatorname*{Col}A^{\prime}  &  \rightarrow\mathbb{Z}%
^{n}/\operatorname*{Col}A,\\
\overline{v}  &  \mapsto\overline{B^{-1}v}%
\end{align*}
is well-defined and is inverse to $f$. The \textquotedblleft
well-defined\textquotedblright\ part of this claim is left to the reader (the
proof is analogous to the proof that $f$ is well-defined, since $A^{\prime
}=BAC$ entails $A=B^{-1}A^{\prime}C^{-1}$). The \textquotedblleft inverse to
$f$\textquotedblright\ part is straightforward (we have $BB^{-1}v=v$ and
$B^{-1}Bv=v$ for any $v$).

Now, $f$ is invertible (since $g$ is inverse to $f$), and thus is a
$\mathbb{Z}$-module isomorphism (since $f$ is a $\mathbb{Z}$-module morphism).
Hence, the $\mathbb{Z}$-modules $\mathbb{Z}^{n}/\operatorname*{Col}A$ and
$\mathbb{Z}^{n}/\operatorname*{Col}A^{\prime}$ are isomorphic. This proves
Proposition \ref{prop.modPID.la.congruent-cok-iso}.
\end{proof}

\subsubsection{The Smith normal form algorithm}

The following theorem will be crucial for Step 3:

\begin{theorem}
[Smith normal form, weak version]\label{thm.modPID.SNF-Z}Each rectangular
matrix $A\in\mathbb{Z}^{n\times m}$ is equivalent to a diagonal matrix (i.e.,
can be transformed into a diagonal matrix via row and column operations).
\end{theorem}

This theorem, combined with Proposition \ref{prop.modPID.la.congruent-cok-iso}%
, suffices to complete Step 3 of our plan. Thus, we need to prove Theorem
\ref{thm.modPID.SNF-Z}. Here is a very rough outline of the proof:

\begin{proof}
[Proof of Theorem \ref{thm.modPID.SNF-Z} (sketched).]Again, we give a
\textquotedblleft proof by example\textquotedblright. We start with the matrix
$\left(
\begin{array}
[c]{cc}%
4 & 6\\
3 & 2\\
2 & 2
\end{array}
\right)  \in\mathbb{Z}^{3\times2}$, and we try to transform it into a diagonal
matrix by a sequence of row operations and column operations. Note that this
is in some sense a subtler version of Gaussian elimination (subtler because we
are not allowed to scale rows or columns by any numbers other than $-1$, and
because we can only add $\mathbb{Z}$-multiples of rows/columns to other
rows/columns, rather than $\mathbb{Q}$-multiples). We shall use the
\textquotedblleft$\overset{\operatorname*{R}}{\longrightarrow}$%
\textquotedblright\ arrow for \textquotedblleft row
operation\textquotedblright\ and the \textquotedblleft%
$\overset{\operatorname*{C}}{\longrightarrow}$\textquotedblright\ arrow for
\textquotedblleft column operation\textquotedblright.%
\begin{align*}
\left(
\begin{array}
[c]{cc}%
4 & 6\\
3 & 2\\
2 & 2
\end{array}
\right)   &  \overset{\operatorname*{C}}{\longrightarrow}\left(
\begin{array}
[c]{cc}%
4 & 2\\
3 & -1\\
2 & 0
\end{array}
\right)  \ \ \ \ \ \ \ \ \ \ \left(  \text{here we subtracted column 1 from
column 2}\right) \\
&  \overset{\operatorname*{C}}{\longrightarrow}\left(
\begin{array}
[c]{cc}%
0 & 2\\
5 & -1\\
2 & 0
\end{array}
\right)  \ \ \ \ \ \ \ \ \ \ \left(  \text{here we subtracted }2\cdot
\text{column 2 from column 1}\right) \\
&  \overset{\operatorname*{C}}{\longrightarrow}\left(
\begin{array}
[c]{cc}%
2 & 0\\
-1 & 5\\
0 & 2
\end{array}
\right)  \ \ \ \ \ \ \ \ \ \ \left(  \text{here we swapped columns 1 and
2}\right) \\
&  \overset{\operatorname*{R}}{\longrightarrow}\left(
\begin{array}
[c]{cc}%
2 & 0\\
1 & -5\\
0 & 2
\end{array}
\right)  \ \ \ \ \ \ \ \ \ \ \left(  \text{here we scaled row 2 by }-1\right)
\\
&  \overset{\operatorname*{R}}{\longrightarrow}\left(
\begin{array}
[c]{cc}%
0 & 10\\
1 & -5\\
0 & 2
\end{array}
\right)  \ \ \ \ \ \ \ \ \ \ \left(  \text{here we subtracted }2\cdot\text{row
2 from row 1}\right) \\
&  \overset{\operatorname*{R}}{\longrightarrow}\left(
\begin{array}
[c]{cc}%
1 & -5\\
0 & 10\\
0 & 2
\end{array}
\right)  \ \ \ \ \ \ \ \ \ \ \left(  \text{here we swapped rows 1 and
2}\right) \\
&  \overset{\operatorname*{C}}{\longrightarrow}\left(
\begin{array}
[c]{cc}%
1 & 0\\
0 & 10\\
0 & 2
\end{array}
\right)  \ \ \ \ \ \ \ \ \ \ \left(  \text{here we added }5\cdot\text{column 1
to column 2}\right) \\
&  \overset{\operatorname*{R}}{\longrightarrow}\left(
\begin{array}
[c]{cc}%
1 & 0\\
0 & 0\\
0 & 2
\end{array}
\right)  \ \ \ \ \ \ \ \ \ \ \left(  \text{here we subtracted }5\cdot\text{row
3 from row 2}\right) \\
&  \overset{\operatorname*{R}}{\longrightarrow}\left(
\begin{array}
[c]{cc}%
1 & 0\\
0 & 2\\
0 & 0
\end{array}
\right)  \ \ \ \ \ \ \ \ \ \ \left(  \text{here we swapped rows 2 and
3}\right)  ,
\end{align*}
and this is a diagonal matrix.

The general procedure is as follows (you can check that this is precisely what
we have done in the example above):

\begin{itemize}
\item We first \textquotedblleft clear out\textquotedblright\ the 1st row;
this means turning it into $\left(  g,0,0,\ldots,0\right)  $, where $g$ is the
gcd of its entries. This is achieved as follows: We first ensure that all
entries in the 1st row are nonnegative by appropriate column operations
(namely, whenever an entry is negative, we scale the respective column by
$-1$). Then, as long as the 1st row contains at least two nonzero entries, we
subtract the column that contains the smaller one (or, better, an appropriate
multiple of this column) from the column that contains the larger
one\footnote{If the two entries are equal, then we subtract the column that
lies further left from the column that lies further right.}. Note that this is
essentially the Euclidean algorithm (or, to be more precise, a variant thereof
for multiple integers). Finally, when there is only one nonzero entry left in
the 1st row, we move this entry into the position $\left(  1,1\right)  $ by
another column operation (swapping its column with the first column).

\item Then, we use the same method (but using row operations instead of column
operations) to clear out the 1st column (i.e., to ensure that its only nonzero
entry is the $\left(  1,1\right)  $-entry).

Note that this might mess up the 1st row again (i.e., some entries of the 1st
row that were previously $0$ might become nonzero again); in this case, we
again clear out the 1st row, then again clear out the 1st column, and so on,
until neither the 1st row nor the 1st column contain any nonzero entries
except for the $\left(  1,1\right)  $-entry.

I claim that this loop cannot go on forever. To see why, you should note that
each of the \textquotedblleft clear out the 1st row\textquotedblright\ and
\textquotedblleft clear out the 1st column\textquotedblright\ subroutines
causes the $\left(  1,1\right)  $-entry to be replaced by a gcd of several
entries, one of which is the $\left(  1,1\right)  $-entry. Clearly, such a
replacement cannot make the $\left(  1,1\right)  $-entry larger (at least in
absolute value), since $\left\vert \gcd\left(  a_{1},a_{2},\ldots
,a_{k}\right)  \right\vert \leq\left\vert a_{1}\right\vert $ for any integers
$a_{1},a_{2},\ldots,a_{k}$ with $a_{1}$ nonzero\footnote{There is one
exception: If the $\left(  1,1\right)  $-entry was $0$, then it can become
larger during a \textquotedblleft clearing-out\textquotedblright\ subroutine
(for example, if the 1st row was $\left(  0,51\right)  $, then its $\left(
1,1\right)  $-entry $0$ will get replaced with $\gcd\left(  0,\ 51\right)
=51$). But this can happen only once during our algorithm, because once the
$\left(  1,1\right)  $-entry is nonzero, it will never become $0$ again (since
the gcd of multiple numbers cannot be $0$ unless \textbf{all} these numbers
were $0$). Thus, this exception still cannot make our loop go on forever.}.
Moreover, it will make this entry strictly smaller, unless the $\left(
1,1\right)  $-entry was the gcd of all the entries in its row/column to begin
with. Let us refer to the latter case as the \textquotedblleft degenerate
case\textquotedblright; in this case, the $\left(  1,1\right)  $-entry does
not change as we clear out the 1st row/column. Sooner or later, we will
necessarily encounter this \textquotedblleft degenerate case\textquotedblright%
, since otherwise the $\left(  1,1\right)  $-entry would keep decreasing
indefinitely (but a nonnegative integer cannot do that).

What happens when we encounter the \textquotedblleft degenerate
case\textquotedblright? Let us say that we encounter it as we are clearing out
the 1st row, after the 1st column has been already cleared out (this is always
the case after the second time we apply our \textquotedblleft clearing-out
subroutine\textquotedblright). Say that the 1st row is $\left(  a_{1}%
,a_{2},\ldots,a_{m}\right)  $ at the moment we start clearing it out, and
becomes $\left(  a_{1},0,0,\ldots,0\right)  $ after we clear it out (the first
entry is $a_{1}$ because we are in the \textquotedblleft degenerate
case\textquotedblright). This entails that $a_{1}=\gcd\left(  a_{1}%
,a_{2},\ldots,a_{m}\right)  $, so that all the entries $a_{1},a_{2}%
,\ldots,a_{m}$ of the 1st row (before the \textquotedblleft clearing
out\textquotedblright) must be multiples of $a_{1}$. This property is clearly
preserved during the \textquotedblleft clearing-out\textquotedblright%
\ subroutine as we clear out the 1st row (since the gcd of its entries does
not change, thus remains $a_{1}$ throughout it). As a consequence, we never
have to change the 1st column during this subroutine (since the 1st column
only has to change if its top entry is larger than another entry of the 1st
row or if its top entry is $\leq0$, but neither of these can happen when all
entries of the 1st row are multiples of the $\left(  1,1\right)  $-entry). So
the 1st column does not change.

Thus, after clearing out the 1st row, we are left with a matrix whose 1st row
is clear (i.e., contains no nonzero entries except for the $\left(
1,1\right)  $-entry) and whose 1st column is also clear (since it has not
changed during the \textquotedblleft clearing out the 1st
row\textquotedblright\ subroutine, but was itself cleared just before it). In
other words, we are left with a matrix whose 1st row and whose 1st column only
contain a single nonzero entry (if any!), which is the $\left(  1,1\right)  $-entry.

At this point, we forget about the 1st row and the 1st column, and play the
same game with the rest of the matrix. (So we are working recursively. Note
that whatever operations we do with the rest of the matrix, the 1st row and
the 1st column will be unaffected, because they are filled with $0$s
everywhere apart from the $\left(  1,1\right)  $-entry. Thus, we won't ever
have to clear them up again.)

\item At the end of the procedure, the matrix will be diagonal.
\end{itemize}

\noindent Thus, after a sequence of row operations and column operations, our
matrix has become diagonal. This proves Theorem \ref{thm.modPID.SNF-Z}.
\end{proof}

This completes Step 3 of our plan.

\subsubsection{A few words on arbitrary rings}

Before I move on to Step 8, let me say a few words about generalizing Theorem
\ref{thm.modPID.SNF-Z} to other rings. In our proof of Theorem
\ref{thm.modPID.SNF-Z}, we seemingly used the fact that the entries of our
matrix are integers (since we argued that a nonnegative integer cannot keep
decreasing indefinitely). However, the proof is easily adapted to any
Euclidean domain instead of $\mathbb{Z}$ (we just need to argue that the
Euclidean norm of the $\left(  1,1\right)  $-th entry decreases, instead of
that entry itself). In truth, Theorem \ref{thm.modPID.SNF-Z} holds even more
generally, with $\mathbb{Z}$ replaced by a PID. This level of generality is a
tad too advanced for us, but proofs of this version of Theorem
\ref{thm.modPID.SNF-Z} can be found in various algebra texts (e.g., in
\cite[Theorem (5.3.10)]{ChaLoi21}). Note that the diagonal matrix in Theorem
\ref{thm.modPID.SNF-Z} is not unique.

\begin{remark}
When the base ring is a field, the diagonal matrix in Theorem
\ref{thm.modPID.SNF-Z} becomes the rank normal form of $A$ (see, e.g.,
\url{https://math.stackexchange.com/questions/371497/} ).
\end{remark}

\subsubsection{Solving systems of linear equations over $\mathbb{Z}$}

\begin{remark}
Incidentally, Theorem \ref{thm.modPID.SNF-Z} also helps solve systems of
linear equations in integer unknowns (as in Exercise \ref{exe.21hw0.5}). To
wit, if two matrices $A,A^{\prime}\in\mathbb{Z}^{n\times m}$ are equivalent,
and if $B\in\mathbb{Z}^{n\times n}$ and $C\in\mathbb{Z}^{m\times m}$ are two
invertible matrices satisfying $A^{\prime}=BAC$, and if $v\in\mathbb{Z}^{n}$
is any vector, then there is a bijection
\begin{align*}
\left\{  w\in\mathbb{Z}^{m}\ \mid\ Aw=v\right\}   &  \rightarrow\left\{
y\in\mathbb{Z}^{m}\ \mid\ A^{\prime}y=Bv\right\}  ,\\
w  &  \mapsto C^{-1}w
\end{align*}
(check this!). Thus, solving the equation $Aw=v$ for an unknown vector
$w\in\mathbb{Z}^{m}$ is tantamount to solving the equation $A^{\prime}y=Bv$
for an unknown vector $y\in\mathbb{Z}^{m}$. But Theorem \ref{thm.modPID.SNF-Z}
tells us that we can choose $A^{\prime}$ to be diagonal, and then the equation
$A^{\prime}y=Bv$ is rather easy to solve. Thus, we obtain an algorithm for
solving a vector equation of the form $Aw=v$ for an unknown vector
$w\in\mathbb{Z}^{m}$; that is, we obtain an algorithm for solving systems of
linear equations in integer unknowns.
\end{remark}

\subsubsection{Step 8: streamlining direct products of $\mathbb{Z}/n$'s}

Let us return to our multi-step plan for proving Theorem
\ref{thm.modPID.class-finabgr}. Step 8 is fun. Let me first discuss it in the
case when $k=2$. In this case, I need to explain how a direct product of the
form $\left(  \mathbb{Z}/n_{1}\right)  \times\left(  \mathbb{Z}/n_{2}\right)
$ with two positive integers $n_{1}$ and $n_{2}$ can be rewritten (up to
isomorphism) as a direct product of the form $\left(  \mathbb{Z}/n_{1}%
^{\prime}\right)  \times\left(  \mathbb{Z}/n_{2}^{\prime}\right)  $ with
$n_{1}^{\prime}\mid n_{2}^{\prime}$. For simplicity, let me rename $n_{1}$ and
$n_{2}$ as $n$ and $m$; then I claim that $n_{1}^{\prime}$ and $n_{2}^{\prime
}$ can be chosen to be $\gcd\left(  n,m\right)  $ and $\operatorname{lcm}%
\left(  n,m\right)  $, respectively (these clearly satisfy $n_{1}^{\prime}\mid
n_{2}^{\prime}$, since $\gcd\left(  n,m\right)  \mid n\mid\operatorname{lcm}%
\left(  n,m\right)  $). In order to prove this claim, I need to show the
following lemma:

\begin{lemma}
\label{lem.modPID.SNF.Zg-Zl}Let $n,m\in\mathbb{Z}$. Let $g=\gcd\left(
n,m\right)  $ and $\ell=\operatorname{lcm}\left(  n,m\right)  $.

\begin{enumerate}
\item[\textbf{(a)}] Then, the matrices%
\[
\left(
\begin{array}
[c]{cc}%
n & 0\\
0 & m
\end{array}
\right)  \ \ \ \ \ \ \ \ \ \ \text{and}\ \ \ \ \ \ \ \ \ \ \left(
\begin{array}
[c]{cc}%
g & 0\\
0 & \ell
\end{array}
\right)
\]
in $\mathbb{Z}^{2\times2}$ are equivalent.

\item[\textbf{(b)}] As a consequence,%
\[
\left(  \mathbb{Z}/n\right)  \times\left(  \mathbb{Z}/m\right)  \cong\left(
\mathbb{Z}/g\right)  \times\left(  \mathbb{Z}/\ell\right)
\]
as groups.
\end{enumerate}
\end{lemma}

\begin{proof}
This is so enjoyable that you should probably try to prove this on your own!
Read on at your own spoiler risk.

\textbf{(a)} We WLOG assume that $g\neq0$ (since otherwise, we have $n=m=0$,
and thus the two matrices in question both equal the zero matrix).

Bezout's theorem shows that there exist integers $x,y$ such that $g=xn+ym$
(since $g=\gcd\left(  n,m\right)  $). Consider these $x,y$. Moreover, there
exists some $u\in\mathbb{Z}$ such that $n=gu$ (since $g\mid n$). Likewise,
there exists some $v\in\mathbb{Z}$ such that $m=gv$ (since $g\mid m$).
Consider these $u$ and $v$.

Furthermore, it is known that $\gcd\left(  n,m\right)  \cdot\operatorname{lcm}%
\left(  n,m\right)  =\left\vert nm\right\vert $. In other words,
$g\ell=\left\vert nm\right\vert $. Thus, $g\ell=\pm\underbrace{n}_{=gu}m=\pm
gum$. Cancelling $g$ from this equality, we find $\ell=\pm um$ (since $g\neq
0$). Thus, $um=\pm\ell$, so that $-um=-\left(  \pm\ell\right)  =\mp\ell$.

Now, we transform the matrix $\left(
\begin{array}
[c]{cc}%
n & 0\\
0 & m
\end{array}
\right)  $ as follows (using the \textquotedblleft$\overset{\operatorname*{R}%
}{\longrightarrow}$\textquotedblright\ arrow for \textquotedblleft row
operation\textquotedblright\ and the \textquotedblleft%
$\overset{\operatorname*{C}}{\longrightarrow}$\textquotedblright\ arrow for
\textquotedblleft column operation\textquotedblright):%
\begin{align*}
\left(
\begin{array}
[c]{cc}%
n & 0\\
0 & m
\end{array}
\right)   &  \overset{\operatorname*{C}}{\longrightarrow}\left(
\begin{array}
[c]{cc}%
n & xn\\
0 & m
\end{array}
\right)  \ \ \ \ \ \ \ \ \ \ \left(  \text{here we added }x\cdot\text{column 1
to column 2}\right) \\
&  \overset{\operatorname*{R}}{\longrightarrow}\left(
\begin{array}
[c]{cc}%
n & xn+ym\\
0 & m
\end{array}
\right)  \ \ \ \ \ \ \ \ \ \ \left(  \text{here we added }y\cdot\text{row 2 to
row 1}\right) \\
&  =\left(
\begin{array}
[c]{cc}%
gu & g\\
0 & m
\end{array}
\right)  \ \ \ \ \ \ \ \ \ \ \left(  \text{since }n=gu\text{ and
}xn+ym=g\right) \\
&  \overset{\operatorname*{C}}{\longrightarrow}\left(
\begin{array}
[c]{cc}%
0 & g\\
-um & m
\end{array}
\right)  \ \ \ \ \ \ \ \ \ \ \left(  \text{here we subtracted }u\cdot
\text{column 2 from column 1}\right) \\
&  =\left(
\begin{array}
[c]{cc}%
0 & g\\
-um & gv
\end{array}
\right)  \ \ \ \ \ \ \ \ \ \ \left(  \text{since }m=gv\right) \\
&  \overset{\operatorname*{R}}{\longrightarrow}\left(
\begin{array}
[c]{cc}%
0 & g\\
-um & 0
\end{array}
\right)  \ \ \ \ \ \ \ \ \ \ \left(  \text{here we subtracted }v\cdot\text{row
1 from row 2}\right) \\
&  \overset{\operatorname*{C}}{\longrightarrow}\left(
\begin{array}
[c]{cc}%
g & 0\\
0 & -um
\end{array}
\right)  \ \ \ \ \ \ \ \ \ \ \left(  \text{here, we swapped column 1 with
column 2}\right) \\
&  =\left(
\begin{array}
[c]{cc}%
g & 0\\
0 & \mp\ell
\end{array}
\right)  \ \ \ \ \ \ \ \ \ \ \left(  \text{since }-um=\mp\ell\right)  .
\end{align*}
If the $\mp\ell$ here is a $+\ell$, then we have thus obtained the matrix
$\left(
\begin{array}
[c]{cc}%
g & 0\\
0 & \ell
\end{array}
\right)  $, so that we conclude that the two matrices $\left(
\begin{array}
[c]{cc}%
n & 0\\
0 & m
\end{array}
\right)  $ and $\left(
\begin{array}
[c]{cc}%
g & 0\\
0 & \ell
\end{array}
\right)  $ are equivalent, as we wanted to show. If it is a $-\ell$ instead,
then we need one more column operation (viz., scaling the second column by
$-1$) in order to get to the same result and therefore to the same conclusion.
Thus, Lemma \ref{lem.modPID.SNF.Zg-Zl} \textbf{(a)} is proved.

\textbf{(b)} Lemma \ref{lem.modPID.SNF.Zg-Zl} \textbf{(a)} yields that the
matrices%
\[
\left(
\begin{array}
[c]{cc}%
n & 0\\
0 & m
\end{array}
\right)  \ \ \ \ \ \ \ \ \ \ \text{and}\ \ \ \ \ \ \ \ \ \ \left(
\begin{array}
[c]{cc}%
g & 0\\
0 & \ell
\end{array}
\right)
\]
in $\mathbb{Z}^{2\times2}$ are equivalent. Hence, Proposition
\ref{prop.modPID.la.congruent-cok-iso} yields that their cokernels%
\[
\mathbb{Z}^{2}/\operatorname*{Col}\left(
\begin{array}
[c]{cc}%
n & 0\\
0 & m
\end{array}
\right)  \ \ \ \ \ \ \ \ \ \ \text{and}\ \ \ \ \ \ \ \ \ \ \mathbb{Z}%
^{2}/\operatorname*{Col}\left(
\begin{array}
[c]{cc}%
g & 0\\
0 & \ell
\end{array}
\right)
\]
are isomorphic. In view of
\[
\mathbb{Z}^{2}/\operatorname*{Col}\left(
\begin{array}
[c]{cc}%
n & 0\\
0 & m
\end{array}
\right)  =\mathbb{Z}^{2}/\operatorname*{span}\left(  ne_{1},me_{2}\right)
\cong\left(  \mathbb{Z}/n\right)  \times\left(  \mathbb{Z}/m\right)
\]
and%
\[
\mathbb{Z}^{2}/\operatorname*{Col}\left(
\begin{array}
[c]{cc}%
g & 0\\
0 & \ell
\end{array}
\right)  =\mathbb{Z}^{2}/\operatorname*{span}\left(  ge_{1},\ell e_{2}\right)
\cong\left(  \mathbb{Z}/g\right)  \times\left(  \mathbb{Z}/\ell\right)  ,
\]
this means that $\left(  \mathbb{Z}/n\right)  \times\left(  \mathbb{Z}%
/m\right)  $ and $\left(  \mathbb{Z}/g\right)  \times\left(  \mathbb{Z}%
/\ell\right)  $ are isomorphic (as $\mathbb{Z}$-modules, and thus as groups).
This proves Lemma \ref{lem.modPID.SNF.Zg-Zl} \textbf{(b)}.
\end{proof}

Lemma \ref{lem.modPID.SNF.Zg-Zl} \textbf{(b)} is sufficient to complete Step 8
in the case when $k=2$ (that is, when $G$ is a direct product of two cyclic
groups). In the general case, we can try to use Lemma
\ref{lem.modPID.SNF.Zg-Zl} \textbf{(b)} multiple times; in fact, applying
Lemma \ref{lem.modPID.SNF.Zg-Zl} \textbf{(b)} to any pair of consecutive
factors $\mathbb{Z}/n_{i}$ and $\mathbb{Z}/n_{i+1}$ in the direct product
$\left(  \mathbb{Z}/n_{1}\right)  \times\left(  \mathbb{Z}/n_{2}\right)
\times\cdots\times\left(  \mathbb{Z}/n_{k}\right)  $ will replace these two
factors by $\mathbb{Z}/n_{i}^{\prime}$ and $\mathbb{Z}/n_{i+1}^{\prime}$ with
$n_{i}^{\prime}\mid n_{i+1}^{\prime}$. For example, if $k=3$, then we can thus
construct the following chain of isomorphisms:%
\begin{align*}
&  \underbrace{\left(  \mathbb{Z}/n_{1}\right)  \times\left(  \mathbb{Z}%
/n_{2}\right)  }_{\substack{\cong\left(  \mathbb{Z}/n_{1}^{\prime}\right)
\times\left(  \mathbb{Z}/n_{2}^{\prime}\right)  \\\text{for }n_{1}^{\prime
}=\gcd\left(  n_{1},n_{2}\right)  \text{ and }n_{2}^{\prime}%
=\operatorname{lcm}\left(  n_{1},n_{2}\right)  \\\text{(by Lemma
\ref{lem.modPID.SNF.Zg-Zl} \textbf{(b)})}}}\times\left(  \mathbb{Z}%
/n_{3}\right) \\
&  \cong\left(  \mathbb{Z}/n_{1}^{\prime}\right)  \times\underbrace{\left(
\mathbb{Z}/n_{2}^{\prime}\right)  \times\left(  \mathbb{Z}/n_{3}\right)
}_{\substack{\cong\left(  \mathbb{Z}/n_{2}^{\prime\prime}\right)
\times\left(  \mathbb{Z}/n_{3}^{\prime\prime}\right)  \\\text{for }%
n_{2}^{\prime\prime}=\gcd\left(  n_{2}^{\prime},n_{3}\right)  \text{ and
}n_{3}^{\prime\prime}=\operatorname{lcm}\left(  n_{2}^{\prime},n_{3}\right)
\\\text{(by Lemma \ref{lem.modPID.SNF.Zg-Zl} \textbf{(b)})}}}\\
&  \cong\underbrace{\left(  \mathbb{Z}/n_{1}^{\prime}\right)  \times\left(
\mathbb{Z}/n_{2}^{\prime\prime}\right)  }_{\substack{\cong\left(
\mathbb{Z}/n_{1}^{\prime\prime\prime}\right)  \times\left(  \mathbb{Z}%
/n_{2}^{\prime\prime\prime}\right)  \\\text{for }n_{1}^{\prime\prime\prime
}=\gcd\left(  n_{1}^{\prime},n_{2}^{\prime\prime}\right)  \text{ and }%
n_{2}^{\prime\prime\prime}=\operatorname{lcm}\left(  n_{1}^{\prime}%
,n_{2}^{\prime\prime}\right)  \\\text{(by Lemma \ref{lem.modPID.SNF.Zg-Zl}
\textbf{(b)})}}}\times\left(  \mathbb{Z}/n_{3}^{\prime\prime}\right) \\
&  \cong\left(  \mathbb{Z}/n_{1}^{\prime\prime\prime}\right)  \times\left(
\mathbb{Z}/n_{2}^{\prime\prime\prime}\right)  \times\left(  \mathbb{Z}%
/n_{3}^{\prime\prime}\right)  .
\end{align*}
It takes some thought to confirm that the resulting numbers $n_{1}%
^{\prime\prime\prime},n_{2}^{\prime\prime\prime},n_{3}^{\prime\prime}$ really
do satisfy $n_{1}^{\prime\prime\prime}\mid n_{2}^{\prime\prime\prime}\mid
n_{3}^{\prime\prime}$. (Indeed, $n_{1}^{\prime\prime\prime}\mid n_{2}%
^{\prime\prime\prime}$ follows from the definitions of $n_{1}^{\prime
\prime\prime}$ and $n_{2}^{\prime\prime\prime}$ as gcd and lcm of one and the
same pair of integers. As for proving $n_{2}^{\prime\prime\prime}\mid
n_{3}^{\prime\prime}$, you have to first argue that combining
\begin{align*}
n_{1}^{\prime}  &  =\gcd\left(  n_{1},n_{2}\right)  \mid\operatorname{lcm}%
\left(  n_{1},n_{2}\right)  =n_{2}^{\prime}\mid\operatorname{lcm}\left(
n_{2}^{\prime},n_{3}\right)  =n_{3}^{\prime\prime}%
\ \ \ \ \ \ \ \ \ \ \text{and}\\
n_{2}^{\prime\prime}  &  =\gcd\left(  n_{2}^{\prime},n_{3}\right)
\mid\operatorname{lcm}\left(  n_{2}^{\prime},n_{3}\right)  =n_{3}%
^{\prime\prime}%
\end{align*}
leads to $\operatorname{lcm}\left(  n_{1}^{\prime},n_{2}^{\prime\prime
}\right)  \mid n_{3}^{\prime\prime}$, so that $n_{2}^{\prime\prime\prime
}=\operatorname{lcm}\left(  n_{1}^{\prime},n_{2}^{\prime\prime}\right)  \mid
n_{3}^{\prime\prime}$.) It might not be obvious, but this generalizes to
arbitrary $k$:

\begin{itemize}
\item First apply Lemma \ref{lem.modPID.SNF.Zg-Zl} \textbf{(b)} to the first
two factors of the direct product, then to the second and third factors, then
to the third and fourth factors, and so on, until you have reached the right
end of the direct product. After this, the numbers $n_{1},n_{2},\ldots
,n_{k-1}$ will all divide $n_{k}$.

\item Then do the same, but stop just before the last factor (i.e., leave the
last factor untouched). After this, the numbers $n_{1},n_{2},\ldots,n_{k-2}$
will all divide $n_{k-1}$, but the numbers $n_{1},n_{2},\ldots,n_{k-1}$ will
still all divide $n_{k}$.

\item Then do the same, but stop just before the second-to-last factor (i.e.,
leave the last two factors untouched). After this, the numbers $n_{1}%
,n_{2},\ldots,n_{k-3}$ will all divide $n_{k-2}$, but the previously mentioned
divisibilities will remain intact.

\item And so on, until at the end there are no more factors left to apply
Lemma \ref{lem.modPID.SNF.Zg-Zl} \textbf{(b)} to. At that point, you will have
$n_{1}\mid n_{2}\mid\cdots\mid n_{k}$.
\end{itemize}

\begin{fineprint}
(Fun fact: There are many other ways to achieve $n_{1}\mid n_{2}\mid\cdots\mid
n_{k}$ by a sequence of moves of the form \textquotedblleft replace $n_{i}$
and $n_{i+1}$ by $\gcd\left(  n_{i},n_{i+1}\right)  $ and $\operatorname{lcm}%
\left(  n_{i},n_{i+1}\right)  $\textquotedblright. Indeed, any sufficiently
long sequence of such moves will eventually come to a halt -- at least if we
make sure to only apply such a move to pairs $\left(  n_{i},n_{i+1}\right)  $
that don't already satisfy $n_{i}\mid n_{i+1}$ -- and the resulting numbers
will satisfy $n_{1}\mid n_{2}\mid\cdots\mid n_{k}$. Moreover, the resulting
numbers will not depend on the sequence of moves. Proving this all is Problem
A3 on the Putnam contest 2008:
\href{https://kskedlaya.org/putnam-archive/2008.pdf}{problem statements} and
\href{https://kskedlaya.org/putnam-archive/2008s.pdf}{solutions}. Our specific
choreographed sequence above was merely the easiest one to analyze.)
\end{fineprint}

Thus we have outlined a proof of parts \textbf{(a)} and \textbf{(b)} of
Theorem \ref{thm.modPID.class-finabgr}. \medskip

\begin{fineprint}
Did Lemma \ref{lem.modPID.SNF.Zg-Zl} \textbf{(b)} remind you of the Chinese
Remainder Theorem? There is indeed a connection, although it is not as obvious
as one might expect:

\begin{exercise}
\label{exe.modPID.SNF.Zg-Zl-rings}Let $n,m\in\mathbb{Z}$. Let $g=\gcd\left(
n,m\right)  $ and $\ell=\operatorname{lcm}\left(  n,m\right)  $.

Lemma \ref{lem.modPID.SNF.Zg-Zl} \textbf{(b)} shows that $\left(
\mathbb{Z}/n\right)  \times\left(  \mathbb{Z}/m\right)  \cong\left(
\mathbb{Z}/g\right)  \times\left(  \mathbb{Z}/\ell\right)  $ as groups.

\begin{enumerate}
\item[\textbf{(a)}] Prove that $\left(  \mathbb{Z}/n\right)  \times\left(
\mathbb{Z}/m\right)  \cong\left(  \mathbb{Z}/g\right)  \times\left(
\mathbb{Z}/\ell\right)  $ as rings as well.

\item[\textbf{(b)}] Use this to prove Theorem \ref{thm.CRT-2-ints-rep} again.
\end{enumerate}

[\textbf{Hint:} Part \textbf{(a)} is easiest to solve using the prime
factorizations of $n$ and $m$. In particular, this yields a new proof of Lemma
\ref{lem.modPID.SNF.Zg-Zl} \textbf{(b)}.]
\end{exercise}

Note that Exercise \ref{exe.modPID.SNF.Zg-Zl-rings} cannot be generalized to
the extent Theorem \ref{thm.CRT-2-ints-rep} was generalized to Theorem
\ref{thm.CRT-2-ids2}. In general, if $I$ and $J$ are two ideals of a
commutative ring $R$, then the rings $\left(  R/I\right)  \times\left(
R/J\right)  $ and $\left(  R/\left(  I+J\right)  \right)  \times\left(
R/\left(  I\cap J\right)  \right)  $ are usually not isomorphic. However,
Exercise \ref{exe.modPID.SNF.Zg-Zl-rings} \textbf{(a)} can be generalized
slightly by replacing $\mathbb{Z}$ by an arbitrary PID.
\end{fineprint}

\subsubsection{Uniqueness of the SNF}

The next two exercises contain a do-it-yourself proof of Theorem
\ref{thm.modPID.class-finabgr} \textbf{(c)} (see also \cite[last claim of
Corollary (5.4.4)]{ChaLoi21} for a more general result).

\begin{exercise}
\label{exe.modPID.class-finabgr.c.1}Let $G$ be a $\mathbb{Z}$-module. Let
$n_{1},n_{2},\ldots,n_{k}$ be $k$ positive integers such that $G\cong\left(
\mathbb{Z}/n_{1}\right)  \times\left(  \mathbb{Z}/n_{2}\right)  \times
\cdots\times\left(  \mathbb{Z}/n_{k}\right)  $ and $n_{1}\mid n_{2}\mid
\cdots\mid n_{k}$. Let $p$ be a prime number, and let $i\in\mathbb{N}$. Prove
that
\[
\left\vert p^{i}G/p^{i+1}G\right\vert =p^{\left(  \text{the number of all
}j\in\left\{  1,2,\ldots,k\right\}  \text{ such that }p^{i+1}\mid
n_{j}\right)  }%
\]
(where we are regarding $G$ as a $\mathbb{Z}$-module, so that $p^{i}G=\left\{
p^{i}g\ \mid\ g\in G\right\}  $ and $p^{i+1}G=\left\{  p^{i+1}g\ \mid\ g\in
G\right\}  $).
\end{exercise}

\begin{exercise}
\label{exe.modPID.class-finabgr.c.3}Prove Theorem
\ref{thm.modPID.class-finabgr} \textbf{(c)}. \medskip

[\textbf{Hint:} Assume that $G\cong\left(  \mathbb{Z}/n_{1}\right)
\times\left(  \mathbb{Z}/n_{2}\right)  \times\cdots\times\left(
\mathbb{Z}/n_{k}\right)  $ and $n_{1}\mid n_{2}\mid\cdots\mid n_{k}$ and
$n_{1},n_{2},\ldots,n_{k}>1$. Prove that knowing the numbers%
\[
\left(  \text{the number of all }j\in\left\{  1,2,\ldots,k\right\}  \text{
such that }p^{i+1}\mid n_{j}\right)
\]
for all primes $p$ and all $i\in\mathbb{N}$ uniquely characterizes
$n_{1},n_{2},\ldots,n_{k}$. Now use Exercise
\ref{exe.modPID.class-finabgr.c.1}.]
\end{exercise}

\subsection{\label{sec.modPID.Fxcyc}Application: Primitive roots}

Theorem \ref{thm.modPID.class-finabgr} allows us to prove Theorem
\ref{thm.finfield.gauss3} (and therefore Theorem \ref{thm.finfield.gauss1} and
Theorem \ref{thm.finfield.gauss2}) again:

\begin{proof}
[Second proof of Theorem \ref{thm.finfield.gauss3}.]The group $G$ is finite
and abelian. Thus, by Theorem \ref{thm.modPID.class-finabgr} (parts
\textbf{(a)} and \textbf{(b)}), we have%
\begin{equation}
G\cong\left(  \mathbb{Z}/n_{1}\right)  \times\left(  \mathbb{Z}/n_{2}\right)
\times\cdots\times\left(  \mathbb{Z}/n_{k}\right)
\label{pf.thm.finfield.gauss3.1}%
\end{equation}
for some positive integers $n_{1},n_{2},\ldots,n_{k}>1$ satisfying $n_{1}\mid
n_{2}\mid\cdots\mid n_{k}$. Consider these $n_{1},n_{2},\ldots,n_{k}$.

Our goal is to show that $k\leq1$ (because then,
(\ref{pf.thm.finfield.gauss3.1}) will show that $G$ is cyclic). In order to
prove this, we assume the contrary. Thus, $k>1$, so $k\geq2$.

Now, $G$ is not just a random abelian group. It has a peculiar property:
Namely, for any positive integer $d$, the group $G$ has no more than $d$
elements $g$ satisfying $g^{d}=1$. (Indeed, all such elements $g$ must be
roots of the degree-$d$ polynomial $x^{d}-1\in F\left[  x\right]  $, but we
know that a degree-$d$ polynomial over a field has no more than $d$ roots.)

Applying this to $d=n_{1}$, we conclude that $G$ has no more than $n_{1}$
elements $g$ satisfying $g^{n_{1}}=1$.

However, the $\mathbb{Z}/n_{1}$ factor on the right hand side of
(\ref{pf.thm.finfield.gauss3.1}) contributes $n_{1}$ such elements (indeed,
each element $g$ of $\mathbb{Z}/n_{1}$ becomes $0$ when multiplied by $n_{1}$,
and thus -- if we rewrite the group multiplicatively -- satisfies $g^{n_{1}%
}=1$), and the $\mathbb{Z}/n_{2}$ factor also contributes $n_{1}$ such
elements (since $n_{1}\mid n_{2}$, so that each of the $n_{1}$ multiples of
$\overline{n_{2}/n_{1}}$ in $\mathbb{Z}/n_{2}$ is such an element). These two
factors overlap only in the identity element. Thus, we have found at least
$2n_{1}-1$ many elements $g\in G$ satisfying $g^{n_{1}}=1$. But there are no
more than $n_{1}$ such elements, as we have seen above. Thus, $2n_{1}-1\leq
n_{1}$, or, equivalently, $n_{1}\leq1$. This contradicts $n_{1}>1$. This
contradiction shows that our assumption was wrong, and this completes the
proof of Theorem \ref{thm.finfield.gauss3}.
\end{proof}

See Subsection \ref{subsec.cyclotomic.Fxcyc} for how Theorem
\ref{thm.finfield.gauss1} and Theorem \ref{thm.finfield.gauss2} follow from
Theorem \ref{thm.finfield.gauss3}.

\begin{exercise}
Let $p$ be a prime. Prove that the group $\left(  \mathbb{Z}/\left(
2p\right)  \right)  ^{\times}$ is cyclic (even though $\mathbb{Z}/\left(
2p\right)  $ is not a field). \medskip

[\textbf{Hint:} There is a ring morphism $\pi:\mathbb{Z}/\left(  2p\right)
\rightarrow\mathbb{Z}/p$ that sends each $\overline{a}$ to $\overline{a}$.
Like any ring morphism, this morphism sends units to units, and thus its
restriction to $\left(  \mathbb{Z}/\left(  2p\right)  \right)  ^{\times}$ is a
group morphism $\left(  \mathbb{Z}/\left(  2p\right)  \right)  ^{\times
}\rightarrow\left(  \mathbb{Z}/p\right)  ^{\times}$. Prove that the latter
group morphism is an isomorphism whenever $p>2$. Deal with the $p=2$ case by
hand. Alternatively, use the Chinese Remainder Theorem.]
\end{exercise}

\begin{exercise}
\ \ 

\begin{enumerate}
\item[\textbf{(a)}] Prove that the group $\left(  \mathbb{Z}/8\right)
^{\times}$ is \textbf{not} cyclic.

\item[\textbf{(b)}] More generally: Let $k$ be a positive integer. Prove that
the group $\left(  \mathbb{Z}/2^{k}\right)  ^{\times}$ is cyclic if and only
if $k\leq2$.
\end{enumerate}

[\textbf{Hint:} For part \textbf{(b)}, if $k>3$, find at least two distinct
elements of $\left(  \mathbb{Z}/2^{k}\right)  ^{\times}$ that have order $2$.]
\end{exercise}

\begin{exercise}
\label{exe.primroot.Z/p2}Let $p\neq2$ be a prime. Prove that the group
$\left(  \mathbb{Z}/p^{2}\right)  ^{\times}$ is cyclic (even though
$\mathbb{Z}/p^{2}$ is not a field). \medskip

[\textbf{Hint:} This group has size $p^{2}-p=p\left(  p-1\right)  $. Thus, it
suffices to find an element of this group whose order is $p\left(  p-1\right)
$.

Pick $a\in\mathbb{Z}$ such that $\overline{a}$ is a generator of the cyclic
group $\left(  \mathbb{Z}/p\right)  ^{\times}$. Show first that $\left(
a+p\right)  ^{p-1}\equiv a^{p-1}-pa^{p-2}\operatorname{mod}p^{2}$. Conclude
that at least one of the integers $a^{p-1}$ and $\left(  a+p\right)  ^{p-1}$
is not congruent to $1$ modulo $p^{2}$. In other words, there exists a
$b\in\left\{  a,\ a+p\right\}  $ such that $b^{p-1}\not \equiv
1\operatorname{mod}p^{2}$. Now, show that the corresponding element
$\overline{b}$ of $\mathbb{Z}/p^{2}$ belongs to the group $\left(
\mathbb{Z}/p^{2}\right)  ^{\times}$ and has order $p\left(  p-1\right)  $ in
this group.]
\end{exercise}

\begin{exercise}
\label{exe.primroot.Z/pk}Let $p\neq2$ be a prime. Let $k$ be a positive
integer. Prove that the group $\left(  \mathbb{Z}/p^{k}\right)  ^{\times}$ is
cyclic (even though $\mathbb{Z}/p^{k}$ is not a field when $k>1$). \medskip

[\textbf{Hint:} Induct on $k$. Exercise \ref{exe.primroot.Z/p2} yields that
there exists an $a\in\mathbb{Z}$ such that $\overline{a}$ is a generator of
the cyclic group $\left(  \mathbb{Z}/p^{2}\right)  ^{\times}$. Consider such
an $a$. Show that $a^{p-1}\equiv1\operatorname{mod}p$ but $a^{p-1}%
\not \equiv 1\operatorname{mod}p^{2}$, and thus $a^{p-1}=1+pu$ for some
integer $u$ not divisible by $p$. Use this to show that $a^{p^{k-2}\left(
p-1\right)  }\not \equiv 1\operatorname{mod}p^{k}$ whenever $k\geq2$. Also
show that $a^{p^{k-1}v}\not \equiv 1\operatorname{mod}p^{k}$ for any proper
divisor $v$ of $p-1$. Use these non-congruences to conclude that $\overline
{a}\in\mathbb{Z}/p^{k}$ has order $p^{k-1}\left(  p-1\right)  $ in the group
$\left(  \mathbb{Z}/p^{k}\right)  ^{\times}$, and therefore the latter group
is cyclic.]
\end{exercise}

With all these exercises, we can finally characterize for which numbers $n$
the group $\left(  \mathbb{Z}/n\right)  ^{\times}$ is cyclic:

\begin{exercise}
\label{exe.primroot.n}Let $n$ be a positive integer. Prove that the group
$\left(  \mathbb{Z}/n\right)  ^{\times}$ is cyclic if and only if
\begin{align*}
n  &  \in\left\{  1,2,4\right\}  \cup\left\{  p^{k}\ \mid\ p\neq2\text{ is a
prime, and }k\text{ is a positive integer}\right\} \\
&  \ \ \ \ \ \ \ \ \ \ \cup\left\{  2p^{k}\ \mid\ p\neq2\text{ is a prime, and
}k\text{ is a positive integer}\right\}  .
\end{align*}

[\textbf{Hint:} For the \textquotedblleft if\textquotedblright\ direction:
What is the connection between the groups $\left(  \mathbb{Z}/m\right)
^{\times}$ and $\left(  \mathbb{Z}/\left(  2m\right)  \right)  ^{\times}$ when
$m$ is odd?

For the \textquotedblleft only if\textquotedblright\ direction: Argue that
$\left(  \mathbb{Z}/n\right)  ^{\times}$ has more than $2$ elements of order
$2$ if $n$ does not belong to the given set.]
\end{exercise}

\end{document}